%% file: main_190414.tex
\documentclass[10pt, a4paper]{amsart}
\pdfoutput=1 
\newcommand{\mydate}{14 Apr 2019}
\makeindex

\usepackage{fancyhdr}
\newcommand{\thisfile}{main\_180925}
\usepackage{enumitem}
\setlist{noitemsep}
\setlist[itemize]{topsep=0ex}
\setlist[enumerate]{topsep=0ex}

\title{Derived Categories} 
\dedicatory{Dedicated to Alexander Grothendieck, in Memoriam}

\usepackage{amscd}
\usepackage{amssymb}
\usepackage{graphicx}
\usepackage[all, pdf]{xy}
\usepackage{lmodern}
\usepackage[T1]{fontenc}

\usepackage{comment}
\usepackage{hyperref}
\hypersetup{colorlinks=false}

\author{Amnon Yekutieli}

\address{Department of  Mathematics
Ben Gurion University, Be'er Sheva 84105, Israel. \newline  \indent 
\textup{\textit{Email}: 
\href{mailto:amyekut@math.bgu.ac.il}{\nolinkurl{amyekut@math.bgu.ac.il}}, 
\textit{Web}: 
\url{http://www.math.bgu.ac.il/~amyekut}.}}

\newtheorem{thm}[equation]{Theorem}
\newtheorem{cor}[equation]{Corollary}
\newtheorem{prop}[equation]{Proposition}
\newtheorem{lem}[equation]{Lemma}
\theoremstyle{definition}
\newtheorem{dfn}[equation]{Definition}
\newtheorem{rem}[equation]{Remark}
\newtheorem{exa}[equation]{Example}

\newtheorem{exer}[equation]{Exercise}
\newtheorem{que}[equation]{Question}

\newtheorem{setup}[equation]{Setup}

\newtheorem{notation}[equation]{Notation}
\newtheorem{conv}[equation]{Convention}
\numberwithin{equation}{subsection}
\newcommand{\iso}{\xrightarrow{\simeq}}
\newcommand{\inj}{\rightarrowtail}  
\newcommand{\surj}{\twoheadrightarrow}
\newcommand{\xar}{\xrightarrow}
\newcommand{\xlar}{\xleftarrow}
\newcommand{\lar}{\leftarrow}
\newcommand{\opn}{\operatorname}

\newcommand{\ul}{\underline}
\newcommand{\ol}{\overline}
\newcommand{\rmitem}[1]{\item[\text{\textup{(#1)}}]}
\newcommand{\mfrak}[1]{\mathfrak{#1}}
\newcommand{\mcal}[1]{\mathcal{#1}}

\newcommand{\mbf}[1]{\mathbf{#1}}
\newcommand{\mrm}[1]{\mathrm{#1}}
\newcommand{\mbb}[1]{\mathbb{#1}}

\newcommand{\tup}[1]{\textup{#1}}
\newcommand{\bsym}[1]{\boldsymbol{#1}}
\newcommand{\boplus}{\bigoplus\nolimits}

\newcommand{\ot}{\otimes}

\newcommand{\til}[1]{\tilde{#1}}

\newcommand{\what}[1]{\widehat{#1}}
\newcommand{\bra}[1]{\langle #1 \rangle}

\newcommand{\K}{\mbb{K}}
\newcommand{\R}{\mbb{R}}
\newcommand{\N}{\mbb{N}}
\newcommand{\Z}{\mbb{Z}}
\newcommand{\Q}{\mbb{Q}}

\newcommand{\OO}{\mcal{O}}

\renewcommand{\AA}{\mcal{A}}
\newcommand{\BB}{\mcal{B}}

\newcommand{\LL}{\mcal{L}}
\newcommand{\MM}{\mcal{M}}
\newcommand{\NN}{\mcal{N}}
\newcommand{\KK}{\mcal{K}}
\newcommand{\PP}{\mcal{P}}

\newcommand{\II}{\mcal{I}}

\newcommand{\g}{\mfrak{g}}
\newcommand{\m}{\mfrak{m}}

\newcommand{\n}{\mfrak{n}}
\renewcommand{\a}{\mfrak{a}}
\renewcommand{\b}{\mfrak{b}}

\newcommand{\p}{\mfrak{p}}
\newcommand{\q}{\mfrak{q}}
\renewcommand{\r}{\mfrak{r}}
\newcommand{\kk}{\bsym{k}}
\newcommand{\om}{\omega}
\newcommand{\al}{\alpha}
\newcommand{\be}{\beta} 
\newcommand{\ga}{\gamma}
\newcommand{\la}{\lambda}
\newcommand{\de}{\delta}
\newcommand{\ze}{\zeta}
\newcommand{\ep}{\epsilon}
\renewcommand{\th}{\theta}
\newcommand{\si}{\sigma}
\newcommand{\Ga}{\Gamma}
\newcommand{\Om}{\Omega}

\newcommand{\bwedge}{{\textstyle \bigwedge}}

\renewcommand{\d}{\mrm{d}}
\newcommand{\pa}{\partial}
\newcommand{\cd}{\mspace{1.5mu}{\cdot}\mspace{1.5mu}} 

\newcommand{\sbmat}[1]{\left[ \begin{smallmatrix} #1
\end{smallmatrix} \right]}
\newcommand{\bmat}[1]{\begin{bmatrix} #1 \end{bmatrix}}
\newcommand{\twoto}{\Rightarrow}
\newcommand{\twoiso}{\stackrel{\simeq\ }{\Longrightarrow}}
\newcommand{\xtwoto}[1]{\stackrel{#1 \ }{\Longrightarrow}}

\newcommand{\cat}[1]{\opn{\mathsf{#1}}}
\newcommand{\catt}[1]{\mrm{\mathsf{#1}}}
\newcommand{\dcat}[1]{\boldsymbol{\mathsf{#1}}}
\newcommand{\bcat}[1]{\operatorname{\boldsymbol{\mathsf{#1}}}}
\newcommand{\sub}{\subseteq}
\newcommand{\subneq}{\varsubsetneq}
\newcommand{\Hom}{\mcal{H}om}

\newcommand{\lb}{\linebreak}
\newcommand{\nl}{\newline}
\newcommand{\centover}{/_{\mspace{-2mu} \mrm{c}} \mspace{2mu}}

\newcommand{\feftover}{/_{\mspace{-2mu} \mrm{feft}} \mspace{2mu}}
\newcommand{\eftover}{/_{\mspace{-2mu} \mrm{eft}} \mspace{2mu}}  
\newcommand{\ringcfeft}[1]{\catt{Rng}_{\mrm{c}} \feftover #1}

\newcommand{\fcentover}{/_{\! \mrm{fc}}\,}
\newcommand{\abs}[1]{\lvert #1 \rvert}
\newcommand{\ttev}{\opn{ev}}

\newcommand{\hnote}[1]{}
\newcommand{\lms}{}

\newcommand{\mysection}{\section}
\newcommand{\mysubsection}{\subsection}
\newcommand{\indash}{--}

\newcommand{\AYcopyright}{
\footnote[0]{This material will be published by Cambridge University Press
as {\em Derived Categories} by Amnon Yekutieli. This prepublication version
is free to view and download for personal use only. Not for redistribution,
resale or use in derivative works. \textcopyright Amnon Yekutieli, 2019.}}

\begin{document}

\begin{abstract}
This is the fourth (and last) prepublication version of a book on {\em derived 
categories}, that will be published by Cambridge University Press. 
The purpose of the book is to provide solid foundations for the theory of 
derived categories, and to present several applications of this theory 
in {\em commutative and noncommutative algebra}. The emphasis is on 
constructions and examples, rather than on axiomatics.  
 
Here is a brief description of the book. After a review of standard facts on 
abelian categories, we study {\em differential graded algebra} in depth 
(including DG rings, DG modules, DG categories and DG functors). 

We then move to {\em triangulated categories}  
and {\em triangulated functors} between them, and explain how they arise from 
the DG background. Specifically, we talk about the homotopy category 
$\dcat{K}(A, \cat{M})$ of DG $A$-modules in $\cat{M}$, where $A$ is a DG 
ring and $\cat{M}$ is an abelian category. 
The {\em derived category} $\dcat{D}(A, \cat{M})$ 
is the localization of $\dcat{K}(A, \cat{M})$ with respect to the 
quasi-isomorphisms. We define the left and right {\em derived functors} of a 
triangulated functor, and prove their uniqueness. Existence of derived 
functors relies on having enough {\em K-injective}, {\em 
K-projective} or {\em K-flat} DG modules, as the case may be. 
We give constructions of resolutions by such DG modules in 
several important algebraic situations. 

There is a detailed study of the derived Hom and tensor bifunctors, and the 
related adjunction formulas. We talk about cohomological dimensions of 
functors and how they are used. 

The last sections of the book are more specialized. 
There is a section on {\em dualizing and residue complexes} over commutative 
noetherian rings, as defined by Grothendieck. We introduce {\em Van den Bergh 
rigidity} in this context. The algebro-geometric applications of commutative 
rigid dualizing complexes are outlined. 

A section is devoted to {\em perfect DG modules} and {\em tilting DG bimodules} 
over NC (noncommutative) DG rings. This includes theorems on {\em derived 
Morita theory} and on the structure of the {\em NC derived Picard group}. 
 
Three sections are on {\em NC connected graded rings}. In the first of them we 
give some basic definitions and results on {\em algebraically graded rings}, 
including {\em Artin-Schelter regular rings}. The second section is on {\em 
derived torsion} for NC connected graded rings, it relation to 
the {\em $\chi$ condition} of Artin-Zhang, and the {\em NC MGM Equivalence}.
In the third section we introduce {\em balanced dualizing complexes}, and we 
prove their uniqueness, existence and trace functoriality.  
The final section of the book is on {\em NC rigid dualizing complexes}, 
following Van den Bergh. We prove the uniqueness and existence of 
these complexes, give a few examples, and discuss their relation to {\em 
Calabi-Yau rings}. 

Feedback from readers (corrections and suggestions) is most welcome. 
\end{abstract}

\AYcopyright

\maketitle

\thispagestyle{fancy}

\cleardoublepage
\tableofcontents

\include{block0_190413}

\include{block1_190413}

\include{block2_190413}

\include{block3_190413}

\include{block4_190413}

\include{block5_190413}

\include{block6_190413}

\include{block7_190413}

\printindex

\end{document}

%% file: block0_190413.tex
\renewcommand{\thisfile}{block0\_190328}
 
\cleardoublepage
\setcounter{section}{-1}
\mysection{Introduction}
\AYcopyright

\mysubsection{On the Subject} \label{subsec:on-the-sub}
{\em Derived categories} were introduced by A. Grothendieck and J.-L.Verdier 
around 1960, and were first published in the book \cite{RD} by R. 
Hartshorne. The basic idea was as follows. They had realized that the derived 
functors of classical homological algebra, namely the functors 
$\mrm{R}^q F, \, \mrm{L}_q F : \cat{M} \to \cat{N}$ 
derived from an additive functor 
$F : \cat{M} \to \cat{N}$
between abelian categories, are too limited to allow several rather natural 
manipulations. Perhaps the most important operation that was lacking was the 
composition of derived functors; the best approximation of it was a spectral 
sequence. 

The solution to the problem was to invent a new category, starting from a given 
abelian category $\cat{M}$. The objects of this new category are the complexes 
of objects of $\cat{M}$. These are the same complexes that play an auxiliary 
role in classical homological algebra, as resolutions of objects of $\cat{M}$. 
The complexes form a category $\dcat{C}(\cat{M})$, but this category is not 
sufficiently intricate to carry in it the information of derived 
functors. So it must be modified. 

A morphism $\phi : M \to N$ in $\dcat{C}(\cat{M})$
is called a {\em quasi-isomorphism}
if in each degree $q$ the cohomology morphism 
$\opn{H}^q(\phi) : \opn{H}^q(M) \to \opn{H}^q(N)$
in $\cat{M}$ is an isomorphism. The modification that is needed is to make the 
quasi-isomorphisms invertible. This is done by a formal 
localization procedure, and the resulting category (with the same objects as 
$\dcat{C}(\cat{M})$) is the derived category $\dcat{D}(\cat{M})$.
There is a functor 
$\opn{Q} : \dcat{C}(\cat{M}) \to \dcat{D}(\cat{M})$,
which is the identity on objects, and it has a universal property (it is 
initial among the functors that send the quasi-isomorphisms to isomorphisms).  
A theorem (analogous to Ore localization in noncommutative ring theory) says 
that every morphism $\th$ in $\dcat{D}(\cat{M})$ can be written as a simple 
left 
or right fraction:
\begin{equation} \label{eqn:4994}
\th = \opn{Q}(\psi_0)^{-1} \circ \opn{Q}(\phi_0) = 
\opn{Q}(\phi_1) \circ \opn{Q}(\psi_1)^{-1} ,
\end{equation}
where $\phi_i$ and $\psi_i$ are morphisms in $\dcat{C}(\cat{M})$, and 
$\psi_i$ are quasi-isomorphisms. 

The cohomology functors 
$\opn{H}^q : \dcat{D}(\cat{M}) \to \cat{M}$,
for all $q \in \Z$, are still defined. 
It turns out that the functor $\cat{M} \to \dcat{D}(\cat{M})$,
that sends an object $M$ to the complex $M$ concentrated in degree $0$, is 
fully faithful. 

The next step is to say what is a left or a right derived functor of an 
additive functor $F : \cat{M} \to \cat{N}$. The functor $F$ can be extended 
in an obvious manner to a functor on complexes
$F : \dcat{C}(\cat{M}) \to \dcat{C}(\cat{N})$.
A {\em right derived functor} of $F$ is a functor 
\begin{equation} \label{eqn:4990}
\mrm{R} F : \dcat{D}(\cat{M}) \to \dcat{D}(\cat{N}) ,
\end{equation}
together with a morphism of functors 
$\eta^{\mrm{R}} : \opn{Q}_{\cat{N}} \circ \, F \to \mrm{R} F \circ 
\opn{Q}_{\cat{M}}$.
The pair $(\mrm{R} F, \eta^{\mrm{R}})$ has to be {\em initial among all such 
pairs}. The uniqueness of such a functor $\mrm{R} F$, up to a unique 
isomorphism, is relatively easy to prove (using the language of 
$2$-categories). As for existence of $\mrm{R} F$, it relies on the 
existence of suitable resolutions (like the injective resolutions in the 
classical situation). If these resolutions exist, 
and if the original functor $F$ is left exact,
then there is a canonical isomorphism of functors 
\begin{equation} \label{eqn:4991}
\mrm{R}^q F \cong \opn{H}^q \circ \, \mrm{R} F : \cat{M} \to \cat{N} 
\end{equation}
for every $q \geq 0$. 

The left derived functor 
\begin{equation} \label{eqn:4992}
\mrm{L} F : \dcat{D}(\cat{M}) \to \dcat{D}(\cat{N}) 
\end{equation}
is defined similarly. When suitable resolutions exist, and when 
$F$ is right exact, there is a canonical isomorphism of functors
\begin{equation} \label{eqn:4993}
\mrm{L}_q F \cong \opn{H}^{-q} \circ \, \mrm{L} F : \cat{M} \to \cat{N} 
\end{equation}
for every $q \geq 0$. 

There are several variations: $F$ could be a contravariant additive 
functor; or it could be an additive bifunctor, contravariant in one or two of 
its arguments. In all these situations the derived (bi)functors 
$\mrm{R} F$ and $\mrm{L} F$ can be defined. 

The derived category $\dcat{D}(\cat{M})$ is additive, but it is not abelian. 
The notion of short exact sequence (in $\cat{M}$ and in $\dcat{C}(\cat{M})$) is 
replaced by that of {\em distinguished 
triangle}, and thus $\dcat{D}(\cat{M})$ is a {\em triangulated category}. 
The derived functors $\mrm{R} F$ and $\mrm{L} F$ are {\em triangulated 
functors}, which means that they send distinguished triangles in 
$\dcat{D}(\cat{M})$ to distinguished triangles in $\dcat{D}(\cat{N})$. 

Already in classical homological algebra we are interested in the {\em 
bifunctors} $\opn{Hom}(-, -)$ and $(- \ot -)$. These bifunctors can also be 
derived. To simplify matters, let's assume that $A$ is a commutative ring, and 
$\cat{M} = \cat{N} = \cat{Mod} A$, the category of $A$-modules. 
We then have bifunctors 
\[ \opn{Hom}_A(-, -) \, : \, (\cat{Mod} A)^{\mrm{op}} \, \times \, \cat{Mod} A
\, \to \, \cat{Mod} A \]
and 
\[ (- \ot_A -) \, : \, \cat{Mod} A \, \times \, \cat{Mod} A \, \to \, 
\cat{Mod} A \, , \]
where the superscript ``op'' denotes the opposite category, that encodes the
contravariance in the first argument of Hom. 
In this situation all resolutions exist, and we have the right derived 
bifunctor 
\begin{equation} \label{eqn:4995}
\opn{RHom}_A(-, -) \, : \, \dcat{D}(\cat{Mod} A)^{\mrm{op}} \, \times \, 
\dcat{D}(\cat{Mod} A) \, \to \, \dcat{D}(\cat{Mod} A) 
\end{equation}
and the left derived bifunctor 
\begin{equation} \label{eqn:4996}
(- \ot^{\mrm{L}}_{A} -) \, : \, \dcat{D}(\cat{Mod} A) \, \times \, 
\dcat{D}(\cat{Mod} A) \, \to \, \dcat{D}(\cat{Mod} A) \, .
\end{equation}
The compatibility with the classical derived bifunctors is this: there are 
canonical isomorphisms 
\begin{equation} \label{eqn:4997}
\opn{Ext}^q_A(M, N) \, \cong \,
\opn{H}^q \bigl( \opn{RHom}_A(M, N) \bigr) 
\end{equation}
and
\begin{equation} \label{eqn:4998}
\opn{Tor}_q^A(M, N) \, \cong  \, \opn{H}^{-q} (M \ot^{\mrm{L}}_{A} N) 
\end{equation}
for all $M, N \in \cat{Mod} A$ and $q \geq 0$.  

This is what derived categories and derived functors are. As to what can be 
done with them, here are some of the things we will explore in our book:
\begin{itemize}
\item {\em Dualizing complexes} and {\em residue complexes} over noetherian 
commutative rings. Besides the original treatment from \cite{RD}, that we 
present in detail here, we also include {\em Van den Bergh rigidity} in the 
commutative setting, that gives rise to {\em rigid residue complexes}. 

\item {\em Perfect DG modules} and {\em tilting DG bimodules}  over 
noncommutative DG rings, and a few variants of {\em derived 
Morita Theory}, including the {\em Rickard-Keller Theorem}. 

\item {\em Derived torsion} and {\em balanced dualizing complexes}
over connected graded NC rings, and {\em rigid 
dualizing complexes} over NC rings, including a full proof of the {\em Van den 
Bergh Existence Theorem} for NC rigid dualizing complexes. 
\end{itemize}

A topic that is beyond the scope of this book, but of which we 
provide an outline here, is:
\begin{itemize}
\item The {\em rigid approach to Grothendieck Duality} on 
noetherian schemes and Deligne-Mumford stacks. 
\end{itemize}

Derived categories have important roles in several areas of mathematics; below 
is a partial list. We will not be able to talk about any of these 
topics in our book, so instead we give some references alongside each 
topic. 
\begin{itemize}
\item[$\triangleright$] $\mcal{D}$-modules, perverse sheaves, and 
representations of algebraic groups and Lie algebras.  
See \cite{BBD} and \cite{Bor}. 
More recently the focus in this area is on the {\em Geometric Langlands 
Correspondence}, that can only be stated in terms of derived categories (see 
the survey \cite{Gai}).

\item[$\triangleright$] 
Algebraic analysis, including differential, 
microdifferential and DQ modules (see \cite{Ka}, \cite{SKK}, \cite{KaSc3}) 
and microlocal sheaf theory (see \cite{KaSc1}), with its application to 
symplectic topology (see \cite{Ta}, \cite{NaZa}). 

\item[$\triangleright$] Representations of finite groups and quivers, 
including {\em cluster algebras} and the {\em Brou\'e Conjecture}. 
See \cite{Hap}, \cite{Kel2}, \cite{CrRo}. 

\item[$\triangleright$] Birational algebraic geometry. This includes {\em 
Fourier-Mukai transforms} and {\em semi-orthogonal decompositions}. 
See the surveys \cite{HiVdB} and \cite{MacSt}, and the book \cite{Huy}. 

\item[$\triangleright$] Homological mirror symmetry. It relates the derived 
category of coherent sheaves on a complex algebraic variety $X$ to the 
{\em derived Fukaya category} of the mirror partner $Y$, which is a symplectic 
manifold. See Remark \ref{rem:4960} and the online reference \cite{HMS}.  

\item[$\triangleright$] Derived algebraic geometry. Here not only the category 
of sheaves is derived, but also the underlying geometric objects (schemes or 
stacks). See Example \ref{exa:4205}, Remark \ref{rem:4970}, and the references 
\cite{Lur} and \cite{To2}.  
\end{itemize}

\mysubsection{A Motivating Discussion: Duality} \label{subsec:2160}
Let us now approach derived categories from another perspective,
very different from the one taken in the previous subsection, by 
considering the idea of {\em duality in algebra}. 

We begin with something elementary: linear algebra. 
Take a field $\K$. Given a $\K$-module $M$ (i.e.\  a vector space), let
$D(M) := \opn{Hom}_{\K}(M, \K)$,
the dual module. There is a canonical homomorphism 
\begin{equation} \label{eqn:3250}
\opn{ev}_M : M \to D(D(M)) ,
\end{equation}
called {\em Hom-evaluation}, whose formula is 
$\opn{ev}_M(m)(\phi) := \phi(m)$ for 
$m \in M$ and $\phi \in D(M)$. 
If $M$ is finitely generated then $\opn{ev}_M$ is an isomorphism (actually this 
is ``if and only if''). 

To formalize this situation, let $\cat{Mod} \K$ denote the category of
$\K$-modules. Then 
$D : \cat{Mod} \K \to \cat{Mod} \K$
is a contravariant functor, and 
$\opn{ev} : \opn{Id} \to D \circ D$
is a morphism of functors (i.e.\ a natural transformation). Here $\opn{Id}$ is 
the identity functor of $\cat{Mod} \K$. 

Now let us replace $\K$ by some nonzero commutative ring $A$. Again we can
define a contravariant functor 
\begin{equation} \label{eqn:3241}
D : \cat{Mod} A \to \cat{Mod} A , \quad
D(M) := \opn{Hom}_A(M, A) ,
\end{equation}
and a morphism of functors $\opn{ev} : \opn{Id} \to D \circ D$.
It is easy to see that $\opn{ev}_M : M \to D(D(M))$ is an isomorphism if $M$
is a finitely generated free $A$-module. 
Of course we can't expect reflexivity (i.e.\ $\opn{ev}_M$ being an isomorphism) 
if $M$ is not finitely generated; but what about a finitely generated module 
that is not free?

In order to understand this better, let us concentrate on the ring $A = \Z$.
Since $\Z$-modules are just abelian groups, the category $\cat{Mod} \Z$ is 
often denoted by $\cat{Ab}$. 
Let $\cat{Ab}_{\mrm{f}}$ be the full subcategory of finitely generated abelian 
groups. Every finitely generated abelian group is of the form 
$M \cong T \oplus F$, with $T$ finite and $F$ free. (The letters  ``T'' and 
``F''  stand for ``torsion'' and ``free'' respectively.) It is important to 
note 
that this is {\em not a canonical 
isomorphism}. There is a canonical short exact sequence
\begin{equation} \label{eqn:1000}
0 \to T \xar{\phi} M \xar{\psi} F \to 0 
\end{equation}
in $\cat{Ab}_{\mrm{f}}$, but the decomposition $M \cong T \oplus F$ comes from 
{\em choosing a splitting} $\sigma : F \to M$ of this sequence. 

\begin{exer}
Prove that the exact sequence (\ref{eqn:1000}) is functorial; 
namely there are functors 
$T, F : \cat{Ab}_{\mrm{f}} \to \cat{Ab}_{\mrm{f}}$,
and natural transformations
$\phi : T \to \opn{Id}$ and  
$\psi :  \opn{Id} \to F$, such that for each $M \in \cat{Ab}_{\mrm{f}}$ 
the group $T(M)$ is finite, the group $F(M)$ is free, and 
the sequence of homomorphisms
\begin{equation} \label{eqn:1002}
0 \to T(M) \xar{\phi_M} M  \xar{\psi_M} F(M) \to  0
\end{equation}
is exact. 

Next, prove that there does not exist a {\em functorial decomposition} of a
finitely generated abelian group into a free part and a finite part. 
Namely, there is no natural transformation
$\sigma : F \to \opn{Id}$, such that for every $M$ the homomorphism 
$\sigma_M : F(M) \to M$ 
splits the sequence (\ref{eqn:1002}). (Hint: find a counterexample.)
\end{exer}

We know that for a free finitely generated abelian group $F$ there is 
reflexivity, i.e.\ $\opn{ev}_F : F \to D(D(F))$ is an isomorphism. But for a
finite abelian group $T$ we have
$D(T) = \opn{Hom}_{\Z}(T, \Z) = 0$.
Thus, for a group $M \in \cat{Ab}_{\mrm{f}}$  with a nonzero torsion subgroup
$T$, reflexivity fails: $\opn{ev}_M : M \to D(D(M))$ is not an isomorphism.

On the other hand, for an abelian group $M$ we can define another sort of dual:
$D'(M) := \opn{Hom}_{\Z}(M, \mbb{Q} / \Z)$.
There is a morphism of functors 
$\opn{ev}' : \opn{Id} \to D' \circ D'$. 
For a finite abelian group $T$ the homomorphism 
$\opn{ev}'_T : T \to D'(D'(T))$ is an isomorphism; this can be seen by 
decomposing $T$ into cyclic groups, and for a finite cyclic group it is clear. 
So $D'$ is a duality for finite abelian groups.
(We may view the abelian group $\mbb{Q} / \Z$ as the group of roots of $1$ in
$\mbb{C}$, via the exponential function; and then $D'$ becomes {\em 
Pontryagin Duality}.)

But for a finitely generated free abelian group $F$ we get
$D'(D'(F)) = \what{F}$, the profinite completion of $F$. So once more this is
not a good duality for all finitely generated abelian groups.

This is where the {\em derived category} enters. 
For every commutative ring $A$ there is the derived category
$\dcat{D}(\cat{Mod} A)$. Here is a very quick explanation of it,
in concrete terms -- as opposed to the abstract point of view taken in the 
previous subsection. 

Recall that a {\em complex} of $A$-modules is a diagram
\begin{equation} \label{eqn:1010}
M = \bigl(  \cdots \to M^{-1} \xar{\d_M^{-1}} M^{0} 
\xar{\d_M^0} M^{1} \to \cdots \bigr) 
\end{equation}
in the category $\cat{Mod} A$.
Namely the $M^i$ are $A$-modules, and the $\d_M^{i}$ are homomorphisms.
The condition is that $\d_M^{i + 1} \circ \d_M^{i} = 0$.
We sometimes write 
$M = \{ M^i \}_{i \in \Z}$. 
The collection $\d_M = \{ \d_M^{i} \}_{i \in \Z}$ is called the 
{\em differential} of $M$.

Given a second complex 
\[ N = \bigl(  \cdots \to N^{-1} \xar{\d_N^{-1}} N^{0} 
\xar{\d_N^0} N^{1} \to \cdots \bigr) , \]
a {\em homomorphism of complexes} 
$\phi^{} : M^{} \to N^{}$
is a collection 
$\phi = \{ \phi^{i} \}_{i \in \Z}$ 
of homomorphisms 
$\phi^i : M^i \to N^i$ 
in $\cat{Mod} A$ satisfying 
$\phi^{i + 1} \circ \d_M^{i} =  \d_N^{i} \circ \phi^i$.
The resulting category is denoted by $\dcat{C}(\cat{Mod} A)$. 

The $i$-th {\em cohomology} of the complex $M$ is 
\begin{equation} \label{eqn:5000}
\opn{H}^i(M^{}) := 
\opn{Ker}(\d_M^{i}) \, / \, \opn{Im}(\d_M^{i - 1}) \in \cat{Mod} A . 
\end{equation}
A homomorphism  $\phi^{} : M^{} \to N^{}$
in $\dcat{C}(\cat{Mod} A)$ induces homomorphisms 
$\opn{H}^i(\phi^{}) : \opn{H}^i(M^{}) \to \opn{H}^i(N^{})$ 
in $\cat{Mod} A$. We call $\phi^{}$ a {\em quasi-isomorphism} if all the 
homomorphisms $\opn{H}^i(\phi^{})$ are isomorphisms. 

The derived category $\dcat{D}(\cat{Mod} A)$ is the localization of 
$\dcat{C}(\cat{Mod} A)$ with respect to the quasi-isomorphisms. 
This means that $\dcat{D}(\cat{Mod} A)$ has the same objects as 
$\dcat{C}(\cat{Mod} A)$. There is a functor 
\begin{equation} \label{eqn:5001}
\opn{Q} : \dcat{C}(\cat{Mod} A) \to \dcat{D}(\cat{Mod} A)
\end{equation}
which is the identity on objects, it sends quasi-isomorphisms to 
isomorphisms, and it is universal for this property. 

A single $A$-module $M^0$ can be viewed as a complex $M$ concentrated in degree 
$0$:
\[ M^{} = \bigl(  \cdots \to 0 \xar{0} M^0 \xar{0} 0 \to \cdots \bigr) . \]
This turns out to be a fully faithful embedding 
\begin{equation} \label{eqn:1004}
\cat{Mod} A \to \dcat{D}(\cat{Mod} A) .
\end{equation}
The essential image of this embedding is the full subcategory 
of $\dcat{D}(\cat{Mod} A)$ on the complexes $M^{}$ whose cohomology is 
concentrated in degree $0$. In this way we have {\em enlarged} the category of 
$A$-modules. All this is explained in Sections 
\ref{sec:loc-cats}-\ref{sec:der-cat} of the book. 

Here is a very important kind of quasi-isomorphism. Suppose $M$ is an 
$A$-module and 
\begin{equation} \label{eqn:1020}
\cdots \to P^{-2} \xar{\d_P^{-2}}  P^{-1} \xar{\d_P^{-1}}  P^{0} 
\xar{\rho} M \to 0 
\end{equation}
is a projective resolution of it.
We can view $M$ as a complex concentrated in degree $0$, by the embedding 
(\ref{eqn:1004}). Define the complex 
\begin{equation} \label{eqn:3246}
P := \bigl(  \cdots \to P^{-2} \xar{\d_P^{-2}} 
P^{-1} \xar{\d_P^{-1}} P^0 \to 0  \to \cdots \bigr) ,
\end{equation}
concentrated in nonpositive degrees. 
Then $\rho$ becomes a morphism of complexes 
$\rho : P \to  M$.
The exactness of the sequence (\ref{eqn:1020}) says that $\rho$ is actually a 
quasi-iso\-morphism. 
Thus $\opn{Q}(\rho) : P^{} \to M$
is an isomorphism in $\dcat{D}(\cat{Mod} A)$.

Let us fix a complex $R \in \dcat{C}(\cat{Mod} A)$.
For every complex $M \in \dcat{C}(\cat{Mod} A)$ we can form the complex
\[ D(M) := \opn{Hom}_{A}(M, R) \in \dcat{C}(\cat{Mod} A) . \]
This is the usual Hom complex (we recall it in 
Subsection \ref{subsec:complexes}). 
As $M$ changes we get a contravariant functor 
\[ D : \dcat{C}(\cat{Mod} A) \to \dcat{C}(\cat{Mod} A) . \]
The functor $D$ has a {\em contravariant right derived functor} 
\begin{equation} \label{eqn:5002}
\mrm{R} D : \dcat{D}(\cat{Mod} A) \to \dcat{D}(\cat{Mod} A) .
\end{equation}
If $P$ is a bounded above complex of projective modules (like in formula 
(\ref{eqn:3246})), or more generally a {\em K-projective complex}
(see Subsection \ref{subsec:K-proj}), then there is a canonical isomorphism
\begin{equation} \label{eqn:3242}
\mrm{R} D(P) \cong D(P) = \opn{Hom}_A(P, R) . 
\end{equation}
Every complex $M$ admits a K-projective resolution
$\rho : P \to M$, and this allows us to calculate $\mrm{R} D(M)$.
Indeed, because the morphism 
$\opn{Q}(\rho) : P \to M$ 
is an isomorphism in 
$\dcat{D}(\cat{Mod} A)$, it follows that 
$\mrm{R} D(\opn{Q}(\rho)) : \mrm{R} D(M) \to \mrm{R} D(P)$
is an isomorphism in $\dcat{D}(\cat{Mod} A)$. 
And the complex $\mrm{R} D(P)$ is known by the canonical isomorphism 
(\ref{eqn:3242}). All this is explained in Sections \ref{sec:der-funcs}, 
\ref{sec:resol} and \ref{sec:exist-resol} of the book. 

It turns out that there is a canonical morphism 
\begin{equation} \label{eqn:3247}
\opn{ev}^{\mrm{R}} : \opn{Id} \to \mrm{R} D \circ \mrm{R} D
\end{equation}
of functors from $\dcat{D}(\cat{Mod} A)$ to itself, called {\em derived 
Hom-evaluation}. See Subsection 
\ref{subsec:du-cplxs}. 

Let us now return to the ring $A = \Z$ and the complex $R = \Z$. 
So the functor $D$ is the same one we had in (\ref{eqn:3241}).
Given a finitely generated abelian group $M$, we want to calculate the 
complexes $\mrm{R} D(M)$ and $\mrm{R} D(\mrm{R} D(M))$, and the morphism 
\begin{equation} \label{eqn:3248}
\opn{ev}^{\mrm{R}}_M : M \to \mrm{R} D(\mrm{R} D(M)) 
\end{equation}
in $\dcat{D}(\cat{Mod} A)$. As explained above, for this we choose a 
projective resolution $\rho : P \to M$,
and then we calculate  the complexes $\mrm{R} D(P)$ and 
$\mrm{R} D(\mrm{R} D(P))$, and the morphism $\opn{ev}^{\mrm{R}}_P$. 
For convenience we choose a projective resolution $P$ of this shape:
\[ \begin{aligned}
& P^{} = 
\bigl(  \cdots \to 0 \to P^{-1} \xar{\d_P^{-1}} P^0 \to 0 \to \cdots \bigr)
\\ & \quad \ 
= \bigl( \cdots \to 0 \xar{}  \Z^{r_1} \xar{\bsym{a} \cd (-)} \Z^{r_0} \xar{} 0 
\cdots 
\bigr) \, , 
\end{aligned} \]
where $r_0, r_1 \in \N$ and $\bsym{a}$ is a matrix of integers. 
The complex $\mrm{R} D(P)$ is this: 
\[ \mrm{R} D(P) \cong D(P) = \opn{Hom}_{\Z}(P^{}, \Z) =
\bigl( \cdots \to 0 \xar{}  \Z^{r_0} \xar{\bsym{a}^{\mrm{t}} \cd (-)}  
\Z^{r_1} \xar{} 0 \cdots \bigr) , \]
a complex of free modules concentrated in degrees $0$ and $1$, with the 
transpose matrix $\bsym{a}^{\mrm{t}}$ as its differential. 

Because $\mrm{R} D(P^{}) \cong D(P^{})$ is itself a bounded complex of 
free modules, its derived dual is 
\begin{equation} \label{eqn:5003}
\mrm{R} D(\mrm{R} D(P^{})) \cong D(D(P^{})) = 
\opn{Hom}_{\Z} \bigl( \opn{Hom}_{\Z}(P^{}, \Z), \Z \bigr) .
\end{equation}
Under the isomorphism (\ref{eqn:5003}), the derived Hom-evaluation morphism 
$\opn{ev}^{\mrm{R}}_P$ in this case is 
just the naive Hom-evaluation homomorphism 
$\opn{ev}_P : P^{} \to D(D(P^{}))$ 
in $\dcat{C}(\cat{Mod} \Z)$ from (\ref{eqn:3250}); see Exercise 
\ref{exer:3245}. Because $P^0$ and $P^{-1}$ are finite rank free modules, it 
follows that $\opn{ev}_P$ is an isomorphism in $\dcat{C}(\cat{Mod} \Z)$. 
Therefore the morphism $\opn{ev}^{\mrm{R}}_M$ in $\dcat{D}(\cat{Mod} \Z)$
is an isomorphism. We see that  {\em $\mrm{R} D$ is a duality that holds for 
all finitely generated $\Z$-modules} $M$~!

Actually, much more is true. Let us denote by 
$\dcat{D}_{\mrm{f}}(\cat{Mod} \Z)$
the full subcategory of $\dcat{D}(\cat{Mod} \Z)$ on the complexes $M$
such that $\opn{H}^i(M)$ is finitely generated for all $i$. Then, according to 
Theorem \ref{thm:2155}, $\opn{ev}^{\mrm{R}}_M$ is an isomorphism for every 
$M \in \dcat{D}_{\mrm{f}}(\cat{Mod} \Z)$. 
It follows that 
\begin{equation} \label{eqn:5004}
\mrm{R} D : \dcat{D}_{\mrm{f}}(\cat{Mod} \Z) \to 
\dcat{D}_{\mrm{f}}(\cat{Mod} \Z)
\end{equation}
is a duality (a contravariant equivalence). This is the celebrated {\em 
Grothendieck Duality}. 

Here is the connection between the derived duality $\mrm{R} D$ and the 
classical dualities $D$ and $D'$. Take a finitely generated abelian group 
$M$, with short exact sequence (\ref{eqn:1000}). There are
canonical isomorphisms
\[ \mrm{H}^0 (\mrm{R} D(M)) \cong \opn{Ext}^0_{\Z}(M, \Z) \cong 
\opn{Hom}_{\Z}(M, \Z) \cong \opn{Hom}_{\Z}(F, \Z) = D(F)  \]
and 
\[ \mrm{H}^1 (\mrm{R} D(M)) \cong \opn{Ext}^1_{\Z}(M, \Z) 
\cong \opn{Ext}^1_{\Z}(T, \Z) \cong D'(T)  .  \]
The cohomologies $\mrm{H}^i (\mrm{R} D(M))$ vanish for $i \neq 0, 1$. 
We see that if $M$ is neither free nor finite, then  
$\mrm{H}^0 (\mrm{R} D(M))$ and $\mrm{H}^1 (\mrm{R} D(M))$ are both nonzero; 
so that the complex $\mrm{R} D(M)$ is not isomorphic
in $\dcat{D}(\cat{Mod} \Z)$ to an object of 
$\cat{Mod} \Z$, under the embedding (\ref{eqn:1004}). 

Grothendieck Duality holds for many noetherian commutative rings $A$. 
A sufficient condition is that $A$ is a finitely generated ring over a 
regular noetherian ring $\K$ (e.g.\ $\K = \Z$ or a field). 
A complex $R \in \dcat{D}(\cat{Mod} A)$ for which the contravariant functor 
\begin{equation} \label{eqn:5010}
\mrm{R} D = \opn{RHom}_{A}(-, R) : \dcat{D}_{\mrm{f}}(\cat{Mod} A) \to 
\dcat{D}_{\mrm{f}}(\cat{Mod} A)
\end{equation}
is an equivalence is called a  {\em dualizing complex}. 
(This is not quite accurate -- see Definition \ref{dfn:2155} for the precise 
technical conditions on $R$.) 
A dualizing complex $R$ over $A$ is unique (up to a degree translation and 
tensoring with an invertible module). 
See Theorems \ref{thm:2155}, \ref{thm:2170} and  \ref{thm:2175}. 

Interestingly, the structure of the dualizing complex $R$ depends on the
geometry of the ring $A$ (i.e.\ of the affine scheme $\opn{Spec}(A)$). 
If $A$ is a regular ring (like $\Z$) then $R = A$ is dualizing. 
If $A$ is a Cohen-Macaulay ring (and $\opn{Spec}(A)$ is connected) then $R$ is 
a single $A$-module (up to a shift in degrees). But if $A$ is a more 
singular ring, then $R$ must live in several degrees, as the next example 
shows. 

\begin{figure}[!htb]
\centering
\includegraphics[scale=0.46]{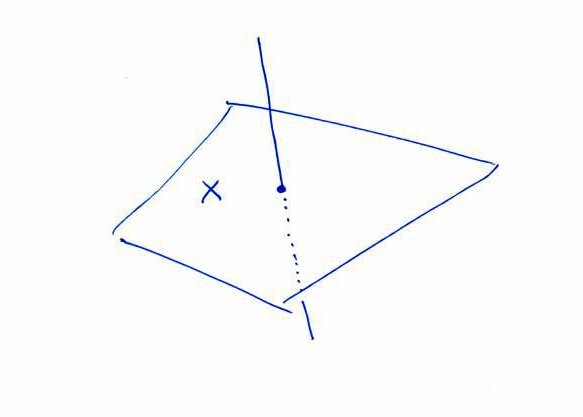}
\caption{An algebraic variety $X$ that is connected but not equidimensional, 
and hence not Cohen-Macaulay.} 
\label{fig:1}
\end{figure}

\begin{exa} \label{exa:2252}
Consider the affine algebraic variety 
$X \subseteq \mbf{A}^3_{\mbb{R}}$
which is the union of a plane and a line that meet at a point, with 
coordinate ring
\[ A = \mbb{R}[t_1, t_2, t_3] / (t_3 \cd t_1, \, t_3 \cd t_2) . \]
See figure \ref{fig:1}. 
A dualizing complex $R$ over $A$ must live in two adjacent degrees;
namely there is some $i$ such that both 
$\mrm{H}^i(R)$ and $\mrm{H}^{i+1}(R)$ are nonzero. 
This calculation is worked out in full in Example \ref{exa:2245}.
\end{exa}

One can also talk about dualizing complexes over {\em noncommutative rings}. 
We will do this in Sections \ref{sec:BDC} and \ref{sec:rigid-DC-NC}.

\mysubsection{On the Book}

This book develops the theory of derived categories, starting from the 
foundations, and going all the way to applications in commutative and 
noncommutative algebra. 
The emphasis is on explicit constructions (with examples), as opposed to 
axiomatics. The most abstract concept we use is probably that of an abelian 
category (which seems indispensable). 

A special feature of this book is that most of the theory deals with 
the category $\dcat{C}(A, \cat{M})$ of {\em DG $A$-modules in $\cat{M}$}, 
where $A$ is a DG  ring and $\cat{M}$ is an 
abelian category. Here ``DG'' is short for ``differential graded'', and our 
DG rings are more commonly known as unital associative DG algebras.
See Subsections \ref{subsec:DG-rng-mod} and 
\ref{subsec:DGModinM} for the definitions. 
The notion $\dcat{C}(A, \cat{M})$ covers most important examples that 
arise in algebra and geometry: 
\begin{itemize}
\item The category $\dcat{C}(A)$ of DG $A$-modules, for any DG ring 
$A$. This includes unbounded complexes of modules over an ordinary ring $A$. 

\item The category $\dcat{C}(\cat{M})$ of unbounded complexes in any 
abelian category $\cat{M}$. This includes 
$\cat{M} = \cat{Mod} \AA$, the category of sheaves of $\AA$-modules on a ringed 
space $(X, \AA)$. 
\end{itemize}

The category $\dcat{C}(A, \cat{M})$ is a {\em DG category}, and 
its DG structure determines the {\em homotopy category} 
$\dcat{K}(A, \cat{M})$ with its {\em triangulated structure}. 
We prove that every {\em DG functor} 
$F : \dcat{C}(A, \cat{M}) \to \dcat{C}(B, \cat{N})$
induces a {\em triangulated functor} 
\begin{equation} \label{eqn:3205}
F : \dcat{K}(A, \cat{M}) \to \dcat{K}(B, \cat{N})
\end{equation}
between the homotopy categories. 

We can now reveal that in the previous subsections we were a bit imprecise (for 
the sake of simplifying the exposition): what we referred to there as  
$\dcat{C}(\cat{M})$ was actually the {\em strict subcategory}
$\dcat{C}_{\mrm{str}}(\cat{M})$, whose morphisms are the degree $0$ cocycles in 
the DG category $\dcat{C}(\cat{M})$. For the same reason the homotopy 
category $\dcat{K}(\cat{M})$ was suppressed there. 

The {\em derived category} $\dcat{D}(A, \cat{M})$ is obtained 
from $\dcat{K}(A, \cat{M})$ by inverting the quasi-isomorphisms.
If suitable resolutions exist, the triangulated functor (\ref{eqn:3205}) 
can be derived, on the left or on the right, to yield triangulated functors
\begin{equation} \label{eqn:5011}
\mrm{L} F, \, \mrm{R} F : 
\dcat{D}(A, \cat{M}) \to \dcat{D}(B, \cat{N}) .
\end{equation}
We prove existence of {\em K-injective}, {\em K-projective} and {\em K-flat} 
resolutions in $\dcat{K}(A, \cat{M})$ in several important contexts,
and explain their roles in deriving functors and in presenting morphisms
in $\dcat{D}(A, \cat{M})$.

In the last six sections of the book we discuss a few key applications of 
derived categories to commutative and noncommutative algebra. 

There is a section-by-section synopsis of the book in subsection 
\ref{subsec:synopsis} of the Introduction. 

The book is based on notes for advanced courses given at Ben Gurion 
University, in the academic years 2011-12, 2015-16 and 2016-17. 
The main sources for Sections 1-12 of the book are \cite{RD} and \cite{KaSc1}; 
but the DG theory component is absent from those earlier texts, and is pretty 
much our own interpretation of folklore results. The material covered in 
Sections 13-18 is adapted from research papers.

\mysubsection{Synopsis of the Book} \label{subsec:synopsis}
Here is a section-by-section description of the material in the book. 

\medskip \noindent 
{\bf Sections \ref{sec:basic}-\ref{sec:abelian}.} These sections are pretty 
much a review of the standard material on categories and functors, especially 
abelian categories and additive functors, that is needed for the book. 

There is one new topic here -- it is our {\em sheaf tricks}, in Subsection 
\ref{subsec:sheaf-tricks}, that are a ``cheap substitute'' for the 
Freyd-Mitchell Theorem, in the sense that they facilitate proofs of various 
results in abstract abelian categories, yet they are themselves easy to prove 
(and we provide these proofs). 

A reader who is familiar with this material can skip these sections; yet we 
do recommend looking at our notational conventions, that are spelled out in 
Conventions \ref{conv:2490} and \ref{conv:3170}. 

\medskip \noindent 
{\bf Section \ref{sec:DG-algebra}.} A good understanding of {\em DG algebra} 
(``DG'' is short for ``differential graded'') is essential in our approach to 
derived categories. By DG algebra we mean DG rings, DG modules, DG categories 
and DG functors. There do not exist (to our knowledge) detailed textbook 
references for DG algebra. Therefore we have included a lot of basic material 
in this section.

We work over a fixed nonzero commutative base ring $\K$ (e.g.\ a field or the 
ring of integers $\Z$). All linear operations (rings, categories, functors, 
etc.) are assumed to be $\K$-linear. A {\em DG module} 
is a module $M$ with a direct sum decomposition 
$M = \boplus_{i \in \Z} M^i$ into submodules, with a differential $\d_M$ of 
degree $1$ satisfying $\d_M \circ \d_M = 0$. The grading on $M$ is called  
{\em cohomological grading}. 
Tensor products of DG modules have {\em signed braiding isomorphisms}: for 
DG modules $M$ and $N$, and for homogeneous elements 
$m \in M^i$ and $n \in N^j$, we define 
\[ \opn{br}_{M, N}(m \ot n) := (-1)^{i \cd j} \cd n \ot m \in N \ot_{\K} M . \]
This braiding is usually referred to as the {\em Koszul sign rule}. 

A {\em DG category} is a $\K$-linear
category $\cat{C}$ in which the morphism sets 
$\opn{Hom}_{\cat{C}}(M, N)$ are endowed with DG module structures, 
and the composition functions are DG bilinear.  
There are two important sources of DG categories: the category 
$\dcat{C}(\cat{M})$ of complexes in a $\K$-linear abelian category $\cat{M}$,
and the category of DG modules over a central DG $\K$-ring $A$
(traditionally called a unital associative DG $\K$-algebra).
Since we want to consider both setups, but we wish to avoid repetition, we 
have devised a new concept that combines both: the DG category 
$\dcat{C}(A, \cat{M})$ of {\em DG $A$-modules in $\cat{M}$}. 
By definition, a DG $A$-module in $\cat{M}$ is a pair $(M, f)$, consisting of a 
complex $M \in \dcat{C}(\cat{M})$, together with a DG $\K$-ring homomorphism 
$f : A \to \opn{End}_{\cat{M}}(M)$.
The morphisms in $\dcat{C}(A, \cat{M})$ are the morphism in 
$\dcat{C}(\cat{M})$ that respect the action of $A$. 
When $A = \K$ we are in the case 
$\dcat{C}(A, \cat{M}) = \dcat{C}(\cat{M})$, and when 
$\cat{M} = \cat{Mod} \K$ we are in the case 
$\dcat{C}(A, \cat{M}) = \dcat{C}(A)$. 

A morphism $\phi : M \to N$ in $\dcat{C}(A, \cat{M})$
is called a {\em strict morphism} if it is a degree $0$ cocycle. 
The {\em strict subcategory} $\dcat{C}_{\mrm{str}}(A, \cat{M})$
is the subcategory of $\dcat{C}(A, \cat{M})$ on all the objects, 
and its morphisms are the strict morphisms. The strict category is abelian. 

It is important to mention that $\dcat{C}(A, \cat{M})$ has 
more structure than just a DG category. Here the objects have DG structures 
too, and there is the cohomology functor 
$\opn{H} : \dcat{C}_{\mrm{str}}(A, \cat{M}) \to 
\dcat{G}_{\mrm{str}}(\cat{M})$,
where the latter is the category of graded objects in $\cat{M}$. 
The cohomology functor determines the {\em quasi-isomorphisms}: these are the 
morphisms $\psi$ in $\dcat{C}_{\mrm{str}}(A, \cat{M})$ s.t.\ 
$\opn{H}(\psi)$ is an isomorphism. 
The set of quasi-isomorphisms is denoted by 
$\dcat{S}(A, \cat{M})$.  

\medskip \noindent 
{\bf Section \ref{sec:stand-tri}.} 
We talk about the {\em translation functor} and {\em standard triangles}
in $\dcat{C}(A, \cat{M})$.
This section consists mostly of new material, some of it implicit in the paper 
\cite{BoKa} on {\em pretriangulated DG categories}. 

The translation $\opn{T}(M)$ of a DG module $M$ is the usual one (a shift by 
$1$ in degree, and differential $-\d_M$). A calculation shows that 
$\opn{T} : \dcat{C}(A, \cat{M}) \to \dcat{C}(A, \cat{M})$
is a DG functor. We introduce the ``little t operator'', which is an invertible 
degree $-1$ morphism 
$\opn{t} : \opn{Id} \to \opn{T}$
of DG functors from $\dcat{C}(A, \cat{M})$ to itself.
The operator $\opn{t}$ facilitates many calculations.  

A morphism $\phi : M \to N$ in $\dcat{C}_{\mrm{str}}(A, \cat{M})$
gives rise to the standard triangle
\[ M \xar{\phi} N \xar{ e_{\phi}} \opn{Cone}(\phi) \xar{p_{\phi}} \opn{T}(M) \]
in $\dcat{C}_{\mrm{str}}(A, \cat{M})$.
As a graded module the standard cone is 
\[ \opn{Cone}(\phi) := N \oplus \opn{T}(M) =
\bmat{N \\ \opn{T}(M)}  , \]
written as a column module. The differential of $\opn{Cone}(\phi)$ is left 
multiplication by the matrix of degree $1$ operators 
\[ \d_{\mrm{Cone}} :=
\bmat{\d_N & \phi \circ \opn{t}_M^{-1} \\[0.1em] 0 & \d_{\opn{T}(M)}} \]

Consider a DG functor
\begin{equation} \label{eqn:1310}
F : \dcat{C}(A, \cat{M}) \to \dcat{C}(B, \cat{N}) , 
\end{equation}
where $B$ is another DG ring and  $\cat{N}$ is another abelian 
category. In Theorem \ref{thm:1150} we show that there is a canonical 
isomorphism
\begin{equation} \label{eqn:1720}
\tau_F : F \circ \opn{T} \iso \opn{T} \circ \, F
\end{equation}
of functors 
$\dcat{C}_{\mrm{str}}(A, \cat{M}) \to \dcat{C}_{\mrm{str}}(B, \cat{N})$,
called the {\em translation isomorphism}. The pair 
$(F, \tau_F)$ is called a {\em T-additive functor}. 
Then, in Theorem \ref{thm:1185}, 
we prove that $F$ sends standard triangles in  
$\dcat{C}_{\mrm{str}}(A, \cat{M})$
to standard triangles in $\dcat{C}_{\mrm{str}}(B, \cat{N})$. 

We end this section with several examples of DG functors.
These examples are prototypes -- they can be easily extended to other setups.

 \medskip \noindent 
{\bf Section \ref{sec:triangulated}.}
We start with the theory of {\em triangulated categories} and 
{\em triangulated functors}, mainly following \cite{RD}. Because the  
octahedral axiom plays no role in our book, we give it minimal attention. 

The {\em homotopy category} $\dcat{K}(A, \cat{M})$ has the same objects as 
$\dcat{C}(A, \cat{M})$, and its morphisms are 
\[ \opn{Hom}_{\dcat{K}(A, \cat{M})} (M, N) := 
\opn{H}^0 \bigl( \opn{Hom}_{\dcat{C}(A, \cat{M})} (M, N) \bigr) . \]
There is a functor 
$\opn{P} : \dcat{C}_{\mrm{str}}(A, \cat{M}) \to \dcat{K}(A, \cat{M})$, 
which is the identity on objects and surjective on morphisms. 

In Subsection \ref{subsec:K-is-triang} we prove that the homotopy 
category $\dcat{K}(A, \cat{M})$ is triangulated. Its {\em distinguished 
triangles} are the triangles that are isomorphic to the images under the 
functor $\opn{P}$ of the standard triangles in 
$\dcat{C}_{\mrm{str}}(A, \cat{M})$. 

Theorem \ref{thm:1265} says that given a DG 
functor $F$ as in (\ref{eqn:1310}), with translation isomorphism $\tau_F$
from (\ref{eqn:1720}), the T-additive functor 
\begin{equation} \label{eqn:5012}
(F, \tau_F) : \dcat{K}(A, \cat{M}) \to \dcat{K}(B, \cat{N})
\end{equation}
is triangulated. This seems to be a new result, unifying well-known yet 
disparate examples. 

In Subsection \ref{subsec:opp-hom-triang} we put a triangulated structure on 
the opposite homotopy category 
$\dcat{K}(A, \cat{M})^{\mrm{op}}$.
We prove that a contravariant DG functor $F$, like (\ref{eqn:1310}) but
with source $\dcat{C}(A, \cat{M})^{\mrm{op}}$, induces a 
contravariant triangulated functor $(F, \tau_F)$, like (\ref{eqn:5012}) but
with source $\dcat{K}(A, \cat{M})^{\mrm{op}}$.

\medskip \noindent 
{\bf Section \ref{sec:loc-cats}.}
In this section we take a close look at {\em localization of categories}. 
We give a detailed proof of the theorem on {\em Ore localization} (also known 
as noncommutative localization). We then prove that the localization 
$\cat{K}_{\cat{S}}$ of a linear category $\cat{K}$ at a denominator set 
$\cat{S}$ is a linear category too, and the localization functor 
$\opn{Q} : \cat{K} \to \cat{K}_{\cat{S}}$ is linear.

\medskip \noindent 
{\bf Section \ref{sec:der-cat}.}
We begin by proving that if $\cat{K}$ is a triangulated category and 
$\cat{S} \sub \cat{K}$ is a {\em denominator set of cohomological origin}, then 
the localized category $\cat{K}_{\cat{S}}$ is triangulated, 
and the localization functor 
$\opn{Q} : \cat{K} \to \cat{K}_{\cat{S}}$ is triangulated.

In the case of the triangulated category $\dcat{K}(A, \cat{M})$, and the 
quasi-iso\-morph\-isms $\dcat{S}(A, \cat{M})$ in it, we get the {\em derived 
category}
\begin{equation} \label{eqn:5013}
\dcat{D}(A, \cat{M}) := \dcat{K}(A, \cat{M})_{\dcat{S}(A, \cat{M})} ,
\end{equation}
and the triangulated localization functor 
\begin{equation} \label{eqn:5014}
\opn{Q} : \dcat{K}(A, \cat{M}) \to \dcat{D}(A, \cat{M}) . 
\end{equation}

We look at the full subcategory $\dcat{K}^{\star}(A, \cat{M})$
of $\dcat{K}(A, \cat{M})$ corresponding to a boundedness condition 
$\star \in \{ +, -, \mrm{b} \}$. 
We prove that when the DG ring $A$ is nonpositive, the localization
$\dcat{D}^{\star}(A, \cat{M})$ of $\dcat{K}^{\star}(A, \cat{M})$
with respect to the quasi-iso\-morphisms in it embeds
fully faithfully in $\dcat{D}(A, \cat{M})$. 
We also prove that the obvious functor 
$\cat{M} \to \dcat{D}(\cat{M})$ is fully faithful.

The section ends with a study of the triangulated structure of the opposite 
derived category $\dcat{D}(A, \cat{M})^{\mrm{op}}$.

\medskip \noindent 
{\bf Section \ref{sec:der-funcs}.}
In this section we talk about {\em derived functors}. To make the definitions 
of derived functors precise, we introduce some {\em $2$-categorical 
notation} here (in Subsections \ref{subsec:2-cat-not}-\ref{subsec:fun-cats}).

In Subsection \ref{subsec:abstr-der-funs} we look at {\em abstract derived 
functors}. Namely, $\cat{K}$ and $\cat{E}$ are categories (without any extra 
structure), $F : \cat{K} \to \cat{E}$ is a functor, and 
$\cat{S} \sub \cat{K}$ is a multiplicatively closed set of morphisms. 
The localization of $\cat{K}$ w.r.t.\ $\cat{S}$ is 
$\opn{Q} : \cat{K} \to \cat{K}_{\cat{S}}$. 
A {\em right derived functor of $F$ w.r.t.\ $\cat{S}$} is a pair 
$(\mrm{R} F ,\eta^{\mrm{R}})$, where 
$\mrm{R} F : \cat{K}_{\cat{S}} \to \cat{E}$
is a functor, and 
$\eta^{\mrm{R}} : F \to \mrm{R} F \circ \opn{Q}$
is a morphism of triangulated functors. The pair $(\mrm{R} F ,\eta^{\mrm{R}})$ 
has a universal property (it is initial among all such pairs), making it unique 
up to a unique isomorphism. The  {\em left derived functor}
$(\mrm{L} F ,\eta^{\mrm{L}})$ is defined similarly.

We provide a general existence theorem for derived functors. For the right 
derived functor $\mrm{R} F$, existence is proved when $\cat{S}$ is a left 
denominator set, and there exists a full 
subcategory $\cat{J} \sub \cat{K}$ which is $F$-acyclic and 
right-resolves all objects of $\cat{K}$.  Likewise, for the left derived functor
$\mrm{L} F$ we prove existence when $\cat{S}$ is a right denominator set, and 
there exists a full subcategory $\cat{P} \sub \cat{K}$ which is $F$-acyclic and 
left-resolves all objects of $\cat{K}$. 

In Subsection \ref{subsec:tri-der-funs} we specialize to {\em triangulated 
derived functors}. Here $\cat{K}$ and $\cat{E}$ are triangulated categories, 
$F : \cat{K} \to \cat{E}$ is a triangulated functor, and 
$\cat{S} \sub \cat{K}$ is a multiplicatively closed set of cohomological 
origin (i.e.\ $\cat{S}$ consists of the quasi-isomorphisms w.r.t.\ 
some cohomological functor $H : \cat{K} \to \cat{M}$). The right and left 
derived functors 
$\mrm{R} F, \mrm{L} F : \cat{K}_{\cat{S}} \to \cat{E}$
are defined like in the abstract setting, and their uniqueness 
is also proved the same way. Existence requires resolving subcategories 
$\cat{J}, \cat{P} \sub \cat{K}$ as above, that are also full triangulated 
subcategories. 

The section is concluded with a discussion of {\em contravariant triangulated 
derived functors}. 

\medskip \noindent 
{\bf Section \ref{sec:DG-tri-bifun}.}
This section is devoted to  {\em DG and triangulated bifunctors}. 
We prove that a DG bifunctor 
\[ F : \dcat{C}(A_1, \cat{M}_1) \times \dcat{C}(A_2, \cat{M}_2)
\to \dcat{C}(B, \cat{N}) \]
induces a triangulated bifunctor 
\[ (F, \tau_1, \tau_2) :  \dcat{K}(A_1, \cat{M}_1) \times 
\dcat{K}(A_2, \cat{M}_2) \to \dcat{K}(B, \cat{N}) . \]

Then we define left and right derived bifunctors, and prove their uniqueness 
and existence, under suitable conditions. Bifunctors that are contravariant in 
one or two of the arguments are also studied.

\medskip \noindent 
{\bf Section \ref{sec:resol}.}
Here we define {\em K-injective} and {\em K-projective} objects in 
$\dcat{K}(A, \cat{M})$, and also  {\em K-flat} objects in 
$\dcat{K}(A)$. These constitute full triangulated subcategories of
$\dcat{K}(A, \cat{M})$, and we refer to them as {\em resolving subcategories}.
The category $\dcat{K}^{\star}(A, \cat{M})_{\mrm{inj}}$ of K-injectives in 
$\dcat{K}^{\star}(A, \cat{M})$, for a boundedness condition $\star$, 
plays the role of the category $\cat{J}$ in Section \ref{sec:der-funcs} 
above, and the category of K-projectives 
$\dcat{K}^{\star}(A, \cat{M})_{\mrm{prj}}$ 
plays the role of the category $\cat{P}$ there.
The K-flat DG modules are acyclic for the tensor functor. 

Furthermore, we prove that the functors 
$\opn{Q} : \dcat{K}^{\star}(A, \cat{M})_{\mrm{inj}} \to 
\dcat{D}^{\star}(A, \cat{M})$
and 
$\opn{Q} : \dcat{K}^{\star}(A, \cat{M})_{\mrm{prj}} \to 
\dcat{D}^{\star}(A, \cat{M})$,
see equation (\ref{eqn:5014}), are fully faithful.

\medskip \noindent 
{\bf Section \ref{sec:exist-resol}.}
In this section we prove {\em existence} of K-injective, K-projective and 
K-flat resolutions in several important cases of 
$\dcat{C}^{\star}(A, \cat{M})$~:
\begin{itemize}
\item K-projective resolutions in $\dcat{C}^{-}(\cat{M})$, where $\cat{M}$ is 
an abelian category with enough projectives. This is classical (i.e.\ it is 
already in \cite{RD}).  

\item K-projective resolutions in $\dcat{C}^{}(A)$, where $A$ is any DG ring. 
This includes the category of unbounded complexes of modules over a ring $A$. 

\item K-injective resolutions in $\dcat{C}^{+}(\cat{M})$, where $\cat{M}$ is 
an abelian category with enough injectives. This is classical too. 

\item K-injective resolutions in $\dcat{C}^{}(A)$, where $A$ is any DG ring. 
\end{itemize}
Our proofs are explicit, and we use limits of complexes cautiously (since 
this is known to be a pitfall). We follow several sources: \cite{RD}, 
\cite{Spa}, \cite{Kel1}, \cite{KaSc1} and private notes provided by B. Keller.

\medskip \noindent 
{\bf Section \ref{sec:adj-equ-cohdim}.}
This section is quite varied. 
In Subsections \ref{subsec:RHom}-\ref{subsec:LTors} there is
a detailed look at the derived bifunctors 
$\opn{RHom}(-,-)$ and $(- \ot^{\mrm{L}} -)$. 

Next, in Subsection \ref{subsec:way-out}, we study {\em cohomological 
dimensions of functors}. This is a refinement of the notion of {\em way-out 
functor} from \cite{RD}. It is used in Subsection \ref{subsec:thm-funcs-fin} to 
prove a few theorems about triangulated functors, such as a sufficient condition 
for a morphism $\ze : F \to G$ of triangulated functors to be an isomorphism. 

In Subsections \ref{subsec:adjs-NC}-\ref{subsec:DGR-quisoms} we study 
several adjunction formulas that involve 
the bifunctors $\opn{RHom}(-,-)$ and $(- \ot^{\mrm{L}} -)$.
We define {\em derived forward adjunction} and  
{\em derived backward adjunction}. These adjunction formulas hold for 
arbitrary DG rings (without any commutativity or boundedness conditions). 
We prove that if $A \to B$ is a {\em quasi-isomorphism of DG rings}, then 
the restriction functor 
$\opn{Rest}_{B / A} : \dcat{D}(B) \to \dcat{D}(A)$
is an equivalence, and it respects the derived bifunctors 
$\opn{RHom}(-,-)$ and $(- \ot^{\mrm{L}} -)$.

Resolutions of DG rings are important in several contexts. In Subsection 
\ref{subsec:DG-ring-res} we provide a full proof that given a 
DG $\K$-ring $A$, there exists a {\em noncommutative semi-free DG ring 
resolution} $\til{A} \to A$ over $\K$.

In Subsection \ref{subsec:der-tens-eval} there is a very general 
theorem providing sufficient conditions for the {\em derived tensor-evaluation 
morphism}
\begin{equation} \label{eqn:4610}
\opn{ev}^{\mrm{R, L}}_{L, M, N} :
\opn{RHom}_{A}(L, M) \ot^{\mrm{L}}_{B} N \to 
\opn{RHom}_{A}(L, M \ot^{\mrm{L}}_{B} N)
\end{equation}
in $\dcat{D}(\K)$ to be an isomorphism. Here $A$ and $B$ are DG $\K$-rings, 
$L \in \dcat{D}(A)$, $M \in \dcat{D}(A \ot_{\K} B^{\mrm{op}})$ and 
$N \in \dcat{D}(B)$.

Lastly, in Subsection \ref{subsec:adjs-comm} we present some 
adjunction formulas that pertain only to weakly commutative DG rings.

\medskip \noindent 
{\bf Section \ref{sec:dual-cplx-comm-rng}.}
This section starts (in the first three subsections) by retelling the material 
in \cite{RD} on {\em dualizing complexes} and {\em residue complexes} over 
noetherian commutative rings.

A complex $R \in \dcat{D}(A)$ is called {\em dualizing} if it has finitely 
generated cohomology modules, finite injective dimension, and the {\em derived 
Morita property}, which says that the {\em derived homothety morphism}
\begin{equation} \label{eqn:4670}
\opn{hm}^{\mrm{R}}_{R} : A \to \opn{RHom}_{A}(R, R)
\end{equation}
in $\dcat{D}(A)$ is an isomorphism. 
(This is a variant of the derived Hom-evaluation morphism 
$\opn{ev}^{\mrm{R}}_M$ from equation (\ref{eqn:3248}), with 
$M = A$.)

To a dualizing complex $R$ we associate the {\em duality functor}
\[ D := \opn{RHom}_{A}(-, R) : \dcat{D}(A)^{\mrm{op}}
\to \dcat{D}(A) . \]
The duality functor induces an equivalence of triangulated categories
$D : \dcat{D}_{\mrm{f}}(A)^{\mrm{op}} \lb \to \dcat{D}_{\mrm{f}}(A)$.
Here $\dcat{D}_{\mrm{f}}(A)$ is the full subcategory of $\dcat{D}(A)$
on the complexes with finitely generated cohomology modules. 

A residue complex is a dualizing complex 
$\KK$ that consists of injective modules of the correct multiplicity in each 
degree. If $R = \KK$ is a residue complex, then the associated duality functor 
is $D = \opn{Hom}_{A}(-, \KK)$.

We prove uniqueness of dualizing complexes over a noetherian commutative ring 
$A$ (up to the obvious twists), and existence when $A$ is {\em essentially 
finite type over a regular noetherian ring} $\K$. In this book the adjective 
``regular noetherian'' includes the condition that $\K$ has finite Krull 
dimension. 

There is a stronger uniqueness property for residue complexes. 
Residue complexes exist whenever dualizing complexes exist: given a dualizing 
complex $R$, its minimal injective resolution $\KK$ is a residue complex. 
To understand residue complexes we review the {\em Matlis classification of 
injective $A$-modules}.  

In remarks we provide sketches of Matlis Duality, Local Duality and the 
interpretation of Cohen-Macaulay complexes as perverse modules.  

In the last two subsections we talk about {\em Van den Bergh rigidity}.
Let $\K$ be a regular noetherian ring, and let $A$ be a flat essentially finite 
type (FEFT) $\K$-ring. (The flatness condition is just to simplify matters; see 
Remark  \ref{rem:2285}.)  Given a complex $M \in \dcat{D}(A)$, its {\em square 
relative to $\K$} is the complex
\[ \opn{Sq}_{A / \K}(M) := 
\opn{RHom}_{A \ot_{\K} A}(A, M \ot^{\mrm{L}}_{\K} M) \in \dcat{D}(A) . \]
We prove that $\opn{Sq}_{A / \K}$ is a {\em quadratic functor}. 
A {\em rigid complex over $A$ relative to $\K$} is a pair 
$(M, \rho)$, where $M \in \dcat{D}(A)$ and 
$\rho : M \iso \opn{Sq}_{A / \K}(M)$
is an isomorphism in $\dcat{D}(A)$, called a {\em rigidifying isomorphism}.
Given another rigid complex $(N, \si)$, a {\em rigid morphism} between 
them is a morphism $\phi : M \to N$ in $\dcat{D}(A)$ such that 
the diagram 
\begin{equation} \label{eqn:5015}
\UseTips \xymatrix @C=8ex @R=6ex {
M
\ar[d]_{\phi}
\ar[r]^(0.4){\rho}
&
\opn{Sq}_{A / \K}(M)
\ar[d]^{\opn{Sq}_{A / \K}(\phi)}
\\
N
\ar[r]^(0.4){\si}
&
\opn{Sq}_{A / \K}(N)
}
\end{equation}
in $\dcat{D}(A)$ is commutative. 
We prove that if $(M, \rho)$ is a rigid complex such that $M$ has the derived 
Morita property, then the only rigid automorphism of $(M, \rho)$ is the 
identity. 

A {\em rigid dualizing complex over $A$ relative to $\K$} is a rigid complex 
$(R, \rho)$ such that $R$ is a dualizing complex. We prove that if $A$ has a 
rigid dualizing complex, then it is unique up to a unique rigid isomorphism.
Existence of a rigid dualizing complex is harder to prove, and we just give a 
reference to it. 

A {\em rigid residue complex over $A$ relative to $\K$} is a rigid complex 
$(\KK, \rho)$ such that $\KK$ is a residue complex. These always exist, 
and they are unique in the following very strong sense: if $(\KK', \rho')$ 
is another rigid residue complex, then there is a 
unique isomorphism $\phi : \KK \iso \KK'$ in $\dcat{C}_{\mrm{str}}(A)$, such 
that $\opn{Q}(\phi)$ is a rigid isomorphism in $\dcat{D}(A)$.

We end this section 
with two remarks (Remark \ref{rem:4192} and \ref{rem:4193}) 
that explain how rigid residue complexes allow a new approach to {\em residues 
and duality on schemes and Deligne-Mumford stacks}, with references.

\medskip \noindent 
{\bf Section \ref{sec:perf-tilt-NC}.}
We begin (in Subsection \ref{subsec:perf-dgmods}) with a systematic 
study of {\em algebraically perfect DG modules} over a NC DG ring $A$. 
The abbreviation ``NC'' stands for ``noncommutative''. By 
definition, a DG module $L \in \dcat{D}(A)$ is called algebraically perfect if 
it belongs to the saturated full triangulated subcategory of $\dcat{D}(A)$ 
generated by the DG module $A$. We give several characterizations of 
algebraically perfect DG modules; one of them is that $L$ is a {\em compact 
object} of $\dcat{D}(A)$. When $A$ is a ring, we prove that $L$ is  
algebraically perfect if and only if it is isomorphic, in $\dcat{D}(A)$, to a 
bounded complex of finitely generated projective $A$-modules. 

In Subsection \ref{subsec:der-morita} we prove the following {\em 
Derived Morita Theorem}, Theorem \ref{thm:3430}:
Let $\cat{E} \sub \dcat{D}(A)$ be a full triangulated subcategory which is 
closed under infinite direct sums, let $P$ be a compact generator of 
$\cat{E}$ which is either K-projective or K-injective, and let
$B := \opn{End}_{A}(P)^{\mrm{op}}$. Then the functor 
$\opn{RHom}_{A}(P, -) : \cat{E} \to \dcat{D}(B)$
is an equivalence of triangulated categories, with quasi-inverse 
$P \ot^{\mrm{L}}_{B} (-)$. 

From Subsection \ref{subsec:flat-dg-rng} to the end of this section 
we assume that the DG rings in question are K-flat over the base ring $\K$. 
Note that every DG $\K$-ring $A$ admits a NC semi-free DG ring resolution 
$\til{A} \to A$. The DG ring $\til{A}$ is K-flat over $\K$, and the 
restriction functor $\dcat{D}(A) \to \dcat{D}(\til{A})$ is an equivalence of 
triangulated categories; so the flatness restriction can be  
circumvented. Subsection \ref{subsec:flat-dg-rng} contains some basic 
constructions of derived functors between categories of DG bimodules (that 
require the flatness over the base ring). 

Next, in Subsection \ref{subsec:tilting}, we define {\em tilting DG 
bimodules}. A DG bimodule $T \in \dcat{D}(B \ot_{\K} A^{\mrm{op}})$ is called
a tilting DG $B$-$A$-bimodule if there exists some 
$S \in \lb \dcat{D}(A \ot_{\K} B^{\mrm{op}})$ 
with isomorphisms $S \ot^{\mrm{L}}_{B} T \cong A$ in 
$\dcat{D}(A^{\mrm{en}})$ and $T \ot^{\mrm{L}}_{A} S \cong B$ in 
$\dcat{D}(B^{\mrm{en}})$. Here 
$A^{\mrm{en}} := A \ot_{\K} A^{\mrm{op}}$ is the {\em enveloping DG ring of 
$A$}, and likewise $B^{\mrm{en}}$. Among other results, we prove that 
$T \in \dcat{D}(B \ot_{\K} A^{\mrm{op}})$
is tilting if and only if $T$ is a compact generator on the $B$ side
(i.e.\ it is a compact generator of $ \dcat{D}(B)$), and it 
has the {\em NC derived Morita property} on the $B$ side, namely the canonical
morphism $A \to \opn{RHom}_{B}(T, T)$
in $\dcat{D}(A^{\mrm{en}})$ is an isomorphism.
We also prove the {\em Rickard-Keller Theorem}, asserting that if $A$ and $B$ 
are rings, and there exists a $\K$-linear equivalence of 
triangulated categories
$\dcat{D}(A) \to \dcat{D}(B)$, then there exists a tilting DG $B$-$A$-bimodule. 

Lastly, in Subsection \ref{subsec:tilting-rings}, we introduce the {\em NC
derived Picard group} $\opn{DPic}_{\K}(A)$ of a flat $\K$-ring $A$. The 
structure of this group is calculated when $A$ is either local or commutative.

\medskip \noindent 
{\bf Section \ref{sec:alg-gra-rings}.}
This section, as well as Sections \ref{sec:der-tors-conn} and \ref{sec:BDC}, 
are on {\em algebraically graded rings}, which is our name for 
$\Z$-graded rings that have lower indices and do not involve the Koszul sign 
rule.  (This is in contrast with the cohomologically graded rings mentioned 
above, that underly DG rings). Simply put, these are the usual graded rings 
that one encounters in textbooks on commutative and noncommutative algebra. 
With few exceptions, the base ring $\K$ in the four final sections of the book 
is a field. 

Let $A = \boplus_{i \in \Z} A_i$ be an algebraically graded central $\K$-ring. 
In Section \ref{sec:alg-gra-rings} we define the category 
of algebraically graded $A$-modules $\dcat{M}(A, \mrm{gr})$. 
Its objects are the algebraically graded (left) $A$-modules
$M = \boplus_{i \in \Z} M_i$, and the homomorphisms are the degree $0$
$A$-linear homomorphisms. It is a $\K$-linear abelian category. 
We talk about finiteness in the algebraically graded context, and about 
various kinds of homological properties, such as {\em 
graded-injectivity}. 

The category of complexes with entries in $\dcat{M}(A, \mrm{gr})$ is the DG 
category \lb 
$\dcat{C}_{}(A, \mrm{gr}) := 
\dcat{C}_{} \bigl( \dcat{M}(A, \mrm{gr}) \bigr)$.
Its objects are bigraded: a complex 
$M \in \dcat{C}_{}(A, \mrm{gr})$ has a direct sum decomposition 
$M = \boplus_{i, j} M^i_j$ into $\K$-modules. Here $i$ is the cohomological 
degree and $j$ is the algebraic degree. The differential goes like this:
$\d_M : M^i_j \to M^{i + 1}_j$. 

Just like for any other abelian category, we have the derived category 
$\dcat{D}_{}(A, \mrm{gr}) := 
\dcat{D}_{} \bigl( \dcat{M}(A, \mrm{gr}) \bigr)$.
This is a triangulated category. We present (quickly) the algebraically graded 
variants of K-projective resolutions etc., and the relevant derived functors.

Special emphasis is given to {\em connected graded rings}. An algebraically 
graded ring $A$ is called connected if $A = \boplus_{i \in \N} A_i$, 
$A_0 = \K$, and each $A_i$ is a finite $\K$-module.
In a connected graded ring $A$ there is the {\em augmentation ideal} 
$\m := \boplus_{i > 0} A_i$. We view $A / \m \cong \K$ as a graded 
$A$-bimodule. 

Among the connected graded rings we single out the {\em Artin-Schelter regular 
graded rings}. A noetherian connected graded ring $A$ is called AS regular if 
it has finite graded global dimension, and if there are isomorphisms
\begin{equation} \label{eqn:4621}
\opn{RHom}_{A}(\K, A) \cong \opn{RHom}_{A^{\mrm{op}}}(\K, A) \cong
\K(l)[-n]
\end{equation}
in $\dcat{D}(\K)$ for some integers $l$ and $n$. One of the results is that if 
$A$ is a noetherian connected graded ring, $a \in A$ is a regular central 
homogeneous element of positive degree, and the ring $B := A / (a)$ is AS 
regular, then $A$ is also AS regular.

\medskip \noindent 
{\bf Section \ref{sec:der-tors-conn}.}
Let $A$ be a connected graded ring over the base field $\K$, with augmentation 
ideal $\m$. In this section we study {\em derived 
$\m$-torsion}, both for complexes of graded $A$-modules and for complexes of 
graded bimodules. By this we mean that, taking a second graded ring $B$, we 
look at the triangulated functor 
\[ \mrm{R} \Ga_{\m} : \dcat{D}(A \ot_{\K} B^{\mrm{op}}, \mrm{gr}) \to
\dcat{D}(A \ot_{\K} B^{\mrm{op}}, \mrm{gr}) .  \]

One of the main results (Theorem \ref{thm:3750}) says that if $A$ is 
a left noetherian connected graded ring, and if the functor 
$\mrm{R} \Ga_{\m}$ has finite cohomological dimension, then there is a 
functorial isomorphism 
\begin{equation} \label{eqn:5016}
\opn{ev}^{\mrm{R, L}}_{\m, M} : 
P \ot^{\mrm{L}}_{A} M \iso \mrm{R} \Ga_{\m}(M)  
\end{equation}
for $M \in \dcat{D}(A \ot_{\K} B^{\mrm{op}}, \mrm{gr})$,
where
\begin{equation} \label{eqn:4620}
P := \mrm{R} \Ga_{\m}(A) \in \dcat{D}(A^{\mrm{en}}, \mrm{gr}) . 
\end{equation}
We also prove the {\em NC MGM Equivalence} in the connected graded context (see 
Theorem \ref{thm:3792}).

The {\em $\chi$ condition} of M. Artin and J.J. Zhang is introduced in 
Subsection \ref{subsec:sym-der-tor}. We study how this condition 
interacts with derived torsion.
To a complex of graded $A$-bimodules, namely an object of 
$\dcat{D}(A^{\mrm{en}}, \mrm{gr})$, we can also apply derived $\m$-torsion from 
the right side, thus obtaining a complex 
$\mrm{R} \Ga_{\m^{\mrm{op}}}(M) \in \dcat{D}(A^{\mrm{en}}, \mrm{gr})$. 
Theorem \ref{thm:4051} says that if $A$ is a 
noetherian connected graded ring of finite local cohomological dimension (i.e.\ 
both functors $\mrm{R} \Ga_{\m}$ and $\mrm{R} \Ga_{\m^{\mrm{op}}}$ have finite 
cohomological dimension), that satisfies the $\chi$ condition, then there is a 
functorial isomorphism 
\begin{equation} \label{eqn:5017}
\ep_M : \mrm{R} \Ga_{\m}(M) \iso \mrm{R} \Ga_{\m^{\mrm{op}}}(M) 
\end{equation}
in $\dcat{D}(A^{\mrm{en}}, \mrm{gr})$ for all complexes 
$M \in \dcat{D}(A^{\mrm{en}}, \mrm{gr})$ whose cohomologies 
are finite $A$-modules on both sides. We call this phenomenon {\em symmetric 
derived $\m$-torsion}.

\medskip \noindent 
{\bf Section \ref{sec:BDC}.}
The focus of this section is on {\em balanced NC dualizing complexes}, 
following \cite{Ye1}.
Let $A$ be a noetherian connected graded ring over the base field $\K$. A 
complex $R \in \dcat{D}^{\mrm{b}}(A^{\mrm{en}}, \mrm{gr})$
is called a {\em graded NC dualizing complex} if it satisfies these three 
conditions:
\begin{enumerate}
\rmitem{i} The bimodules $\opn{H}^q(R)$ are finitely generated 
$A$-modules on both sides.

\rmitem{ii} The complex $R$ has finite graded injective 
dimension on both sides.

\rmitem{iii} The complex $R$ has the NC derived Morita property on both 
sides; namely the canonical morphisms 
$A \to \opn{RHom}_{A}(R, R)$ and 
$A \to \opn{RHom}_{A^{\mrm{op}}}(R, R)$ 
in $\dcat{D}(A^{\mrm{en}}, \mrm{gr})$ are isomorphisms. 
\end{enumerate}

A balanced NC dualizing complex over $A$ is a pair $(R, \be)$,
where $R$ is a graded NC dualizing complex over $A$ with symmetric derived 
$\m$-torsion, and 
$\be : \mrm{R} \Ga_{\m}(R) \iso A^*$
is an isomorphism in $\dcat{D}(A^{\mrm{en}}, \mrm{gr})$, called a {\em 
balancing isomorphism}. Here and later we write 
$M^* := \opn{Hom}_{\K}(M, \K)$, the graded $\K$-linear dual of a graded module 
$M$. The graded bimodule $A^*$ is a graded-injective $A$-module on both sides. 

To a balanced dualizing complex $R$ we associate the duality functors 
\[ D_A := \opn{RHom}_{A}(-, R) : \dcat{D}(A, \mrm{gr})^{\mrm{op}}
\to \dcat{D}(A^{\mrm{op}}, \mrm{gr}) \]
and 
\[ D_{A^{\mrm{op}}} := \opn{RHom}_{A^{\mrm{op}}}(-, R) : 
\dcat{D}(A^{\mrm{op}}, \mrm{gr})^{\mrm{op}}
\to \dcat{D}(A, \mrm{gr}) . \]
When restricted to the subcategories of complexes with finitely generated 
cohomology modules, this gives an equivalence 
\begin{equation} \label{eqn:5018}
D_A :  \dcat{D}_{\mrm{f}}(A, \mrm{gr})^{\mrm{op}}
\to \dcat{D}_{\mrm{f}}(A^{\mrm{op}}, \mrm{gr})
\end{equation}
with quasi-inverse $D_{A^{\mrm{op}}}$. 

A balanced dualizing complex $(R, \be)$ is unique up to a unique isomorphism 
(Theorem \ref{thm:3710}), and it satisfies the NC Local Duality Theorem 
\ref{thm:4616}.  Results of Yekutieli, Zhang and M. Van den Bergh say that 
a noetherian connected graded ring $A$ has a balanced dualizing complex if and 
only if $A$ satisfies the $\chi$ condition and has finite local cohomological 
dimension. The formula for the balanced dualizing complex is 
$R := P^*$, where $P$ is the complex from (\ref{eqn:4620}). 
See Corollary \ref{cor:4585}.

If $A$ is an AS regular (or more generally AS Gorenstein) graded ring, then it 
has a balanced dualizing complex $R = A(\phi, -l)[n]$. Here $l$ and $n$ are the 
integers from formula (\ref{eqn:4621}), and $\phi$ is a graded ring 
automorphism 
of $A$. See Corollary \ref{cor:3780}.

\medskip \noindent 
{\bf Section \ref{sec:rigid-DC-NC}.}
The final section of the book deals with {\em NC rigid dualizing complexes}. 
Let $A$ be a NC central $\K$-ring (where $\K$ is the base field). A {\em NC 
dualizing complex} over $A$ is a complex 
$R \in \dcat{D}^{\mrm{b}}(A^{\mrm{en}})$
whose cohomology bimodules are finitely generated on both sides; it has finite 
injective dimension on both sides; and it has the NC derived Morita property on 
both sides. This is a modification of the graded definition given above.

Let $M \in \dcat{D}(A^{\mrm{en}})$. The {\em NC square} of $M$
is the complex 
\begin{equation} \label{eqn:4660}
\opn{Sq}_{A / \K}(M) := \opn{RHom}_{A^{\mrm{en}}}(A, R \ot_{\K} R) \in 
\dcat{D}(A^{\mrm{en}}) .
\end{equation}
This formula is ambiguous, because  
the complex of $\K$-modules $R \ot_{\K} R$ has on it four distinct 
commuting actions by the ring $A$; but formula (\ref{eqn:4660}) is made precise 
at the beginning of Subsection \ref{subsec:RNCDC-uniq}. 
A {\em rigidifying isomorphism} for $M$ is an isomorphism 
$\rho : M \iso  \opn{Sq}_{A / \K}(M)$
in $\dcat{D}(A^{\mrm{en}})$, and the pair $(M, \rho)$ is a {\em NC rigid 
complex}. A NC rigid dualizing complex is a NC rigid complex $(R, \rho)$ 
such that $R$ is a NC dualizing complex. 

The definition of NC rigid dualizing complex was introduced by Van den Bergh in 
his paper \cite{VdB}. He also proved that a rigid NC dualizing complex is 
unique up to isomorphism. Later, in \cite{Ye4}, it was proved that this 
isomorphism is itself unique. We reproduce these results in Subsection 
\ref{subsec:RNCDC-uniq}.  

We also give Van den Bergh's theorem on the existence of 
NC rigid dualizing complexes from \cite{VdB}, with a full proof. 
This is Theorem \ref{thm:4396}. Here is what it says: 
Suppose the ring $A$ admits a filtration 
$F = \{ F_j(A) \}_{j \geq -1}$
such that the associated graded ring $\opn{Gr}^F(A)$ 
is noetherian connected and has a balanced dualizing complex. Then $A$ has a 
rigid dualizing complex. 

An important special case of the Van den Bergh Existence Theorem is  
Theorem \ref{thm:4490}. It says that if the graded ring  $\opn{Gr}^F(A)$ from 
the previous paragraph is AS regular, then the rigid NC dualizing complex of 
$A$ is $R = A(\mu)[n]$,
where $\mu$ is a $\K$-ring automorphism of $A$ that respects the filtration 
$F$, and $n$ is an integer. The automorphism $\nu := \mu^{-1}$ is called the 
{\em Nakayama automorphism} of $A$. In modern terminology the ring $A$ is 
called an {\em $n$-dimensional twisted Calabi-Yau ring}. 

Finally we state and prove the {\em Van den Bergh Duality Theorem} for 
Hoch\-schild (co)homology, and give an example of a {\em Calabi-Yau category} 
of fractional dimension.

\mysubsection{What Is Not in the Book}
A very important aspect of the theory of derived categories is {\em geometric}. 
Unfortunately our book does not discuss this aspect (except in passing).
We had hoped to include the geometric aspect, but as the 
book grew longer this became impractical. 

Geometric derived categories come in two kinds. The first kind is the derived 
category $\dcat{D}(X) = \dcat{D}(\cat{Mod} \OO_X)$, 
where $(X, \OO_X)$ is a scheme, and $\cat{Mod} \OO_X$ is the abelian category 
of sheaves of $\OO_X$-modules. This kind of derived category is the subject of 
the original book \cite{RD}. For recent treatments see \cite{Huy}, 
\cite{SP} or \cite{Ols}. The last two references also treat derived 
category $\dcat{D}(\mfrak{X})$ of sheaves of modules on an algebraic stack 
$\mfrak{X}$. We address this aspect of derived categories in Remarks 
\ref{rem:4675}, \ref{rem:4676}, \ref{rem:4192} and \ref{rem:4193}. 

The second kind of geometric derived category is 
$\dcat{D}(\K_X) = \dcat{D}(\cat{Mod} \K_X)$, where $X$ is a topological space
(or a site, such as the \'etale site of a scheme), $\K_X$ is the constant 
sheaf of rings $\K$ on $X$, and 
$\cat{Mod} \K_X$ is the abelian category of sheaves of $\K_X$-modules on $X$.
For this kind of derived categories we recommend the book \cite{KaSc1}. 

One should also mention, in this context, the new and important theories of 
{\em derived algebraic geometry} (DAG), in which schemes are replaced by
{\em derived stacks}. These new geometric objects carry derived categories of 
modules. See the references \cite{Lur}, \cite{To2} and \cite{nLab}, 
Remark \ref{rem:4960} and Example \ref{exa:4205}.
The {\em commutative DG rings} that we discuss briefly in Subsection 
\ref{subsec:DG-rng-mod} are {\em affine derived schemes} from the point of view 
of DAG.

\mysubsection{Prerequisites and Recommended Bibliography}

In preparing the book, the assumptions were that the reader is already familiar 
with these topics:
\begin{itemize}
\item Categories and functors, and classical homological 
algebra, namely the derived functors $\mrm{R}^q F$ and $\mrm{L}_q F$
of an additive functor $F : \cat{M} \to \cat{N}$ between abelian  
categories. For these topics we recommend the books 
\cite{Mac1}, \cite{HiSt} and \cite{Rot}. 

\item Commutative and noncommutative ring theory, e.g.\ from the books
\cite{Eis}, \cite{Mats}, \cite{AlKl}, \cite{Row} and \cite{Rot}.  
\end{itemize}

For the topics above we merely review the material, and point to references 
when necessary.

There are a few earlier books that deal with derived categories, in varying 
degrees of detail and depth. Some of them -- \cite{RD} (the original 
reference), \cite{KaSc1} and \cite{KaSc2} -- served as sources for us when 
writing the present book. Other books, such as \cite{GeMa},  
\cite{We} and \cite{Huy}, are somewhat sketchy in their treatment of derived 
categories, but they might be useful for a reader who wants another perspective 
on the subject. 

Finally we want to mention the evolving online reference \cite{SP}, that 
contains a huge amount of information on all the topics listed above. 

\mysubsection{Credo,  Writing Style and Goals}
Since its inception around 1960, there has been very little literature on the 
theory of derived categories. In some respect, the only detailed account for 
many years was the original book \cite{RD}, written by R. Hartshorne, following 
notes by J.-L. Verdier (for Chapter I) and by A. Grothendieck (for Chapters 
II-VII). Several accounts appeared later as parts of the books 
\cite{We}, \cite{GeMa}, \cite{KaSc1}, \cite{KaSc2}, \cite{Huy}, and maybe a few 
others -- but none of these accounts provided enough detailed content to make 
it possible for a mathematician to learn  how to work with derived categories, 
beyond a rather superficial level. {\em The theory thus remained mysterious}. 

A personal belief of mine is that {\em mathematics should not be mysterious}. 
Some mathematics is very easy to explain. However, a few branches of 
mathematics are truly hard; among them are algebraic geometry and derived 
categories. My moral goal in this book is to demonstrate 
that {\em the theory of derived categories is difficult, but not mysterious}. 

The  series of books \cite{EGA} by Grothendieck and J. Dieudonn\'e, and then 
the book \cite{Har} by Hartshorne, have shown us that {\em algebraic geometry 
is difficult but not mysterious}. The definitions and the statements are 
precise, and the proofs are available (to be read or to be taken on trust, as 
the reader prefers). I hope that the present book will do the same for derived 
categories. (Although I doubt I can match the excellent writing talent of the 
aforementioned authors!).

In more practical terms, the goal of this book is to develop the theory of 
derived categories in a systematic fashion, with full details, and with several 
important applications. The expectation is that our book will open the doors 
for researchers in algebra and geometry (and related disciplines such as 
mathematical physics) to productive work using derived categories. Doors 
that have been to a large extent locked until now (or at least hidden by the 
shrouds of mystery, to use the prior metaphor). 

This book is by far too lengthy for a one-semester graduate course.
(In my rather sluggish lecturing style, Sections 1-13 of the book took about 
four semesters.) The book is intended to be used as a reference, or for 
personal learning. Perhaps a lecturer who has the ability to concentrate the 
material sufficiently, or to choose only certain key aspects, might extract a 
one-semester course from this book; if so, please let me know!

\mysubsection{Acknowledgments}

As already mentioned, the book originated in a course on derived categories, 
that was held at Ben Gurion University in Spring 2012. I want to thank the 
participants of this course for correcting many of my mistakes
(both in real time during the lectures, and afterwards when writing the notes 
\cite{Ye6}).  Thanks also to Joseph Lipman, Pierre Schapira, Amnon Neeman 
and Charles Weibel for helpful discussions during the preparation of that 
course. Vincent Beck, Yang Han and Lucas Simon sent me corrections and useful 
comments on the material in \cite{Ye6}.  

I started writing the book itself while teaching a four
semester course on the subject at Ben Gurion University, spanning the 
academic years 2015-16 and 2016-17. I wish to thank the participants of this 
course, and especially Rishi Vyas, Stephan Snigerov, Asaf Yekutieli, S.D. 
Banerjee, Alberto Boix and William Woods, who contributed material and 
corrected numerous  mistakes. 

Ben Gurion University was generous enough to permit this long project. The 
project was supported by the Israel Science Foundation grant no.\ 253/13. 

Stephan Snigerov, William Woods and Rishi Vyas helped me prepare the manuscript 
for publication, and I wish to thank them for that. Bernhard Keller kindly sent 
me private notes containing clean proofs of several theorems on the 
existence of resolutions of DG modules, and he also helped me with numerous 
suggestions. 

In the process of writing the book I have also benefited from the advice of 
Pierre Schapira, Robin Hartshorne, Rodney Sharp, 
Manuel Saorin, Suresh Nayak, Liran Shaul, Johan de Jong, Steven Kleiman, 
Louis Rowen, Amnon Neeman, Amiram Braun, Jesse Wolfson, 
Bjarn de Jong, James Zhang, Jun-Ichi Miyachi, Vladimir Hinich, 
Brooke Shipley, Jeremy Rickard, Peter J{\o}rgensen, 
Michael Sharpe (on LaTeX), Roland Berger, Mattia Ornaghi, Georges Maltsiniotis, 
Harry Gindi, Goncalo Tabuada, Sefi Ladkani, Saurabh Singh,
and the anonymous referees for Cambridge University Press. 

Special thanks to Thomas Harris, my editor at Cambridge University Press,  
whose initiative made the publication of this book possible; and to 
Clare Dennison and Libby Haynes at CUP, for their assistance.

%% file: block1_190413.tex
\renewcommand{\thisfile}{block1\_190328}

\cleardoublepage
\mysection{Basic Facts on Categories} \label{sec:basic}
\AYcopyright

It is assumed that the reader has a working knowledge of categories and 
functors. References for this material are \cite{Mac1}, \cite{Mac2}, 
\cite{Rot}, \cite{HiSt}, \cite{We}, \cite{KaSc1} and \cite{KaSc2}. In this 
section we review the relevant material, and establish notation. 

\mysubsection{Set Theory} \label{subsec:set-theor}
In this book we will not try to be precise about issues of set theory. The
blanket assumption is that we are given a {\em Grothendieck universe}
$\cat{U}$. This is a ``large'' infinite set.
A {\em small set}, or a $\cat{U}$-small set, is a set $S$ that is an element of 
$\cat{U}$. We want all the products $\prod_{i \in I} S_i$ and disjoint unions 
$\coprod_{i \in I} S_i$, with $I$ and $S_i$ small sets, 
to be small sets too. 
(This requirement is not crucial for us, and it is more a matter of convenience.
When dealing with higher categories, one usually needs a hierarchy of 
universes anyhow.) We assume that the axiom of choice holds in $\cat{U}$. 

A {\em $\cat{U}$-category} is a category $\cat{C}$ whose
set of objects $\opn{Ob}(\cat{C})$ is a subset of $\cat{U}$, 
and for every $C, D \in \opn{Ob}(\cat{C})$ the set of morphisms 
$\opn{Hom}_{\cat{C}}(C, D)$ is small. 
If $\opn{Ob}(\cat{C})$ is also small, then $\cat{C}$ is called a {\em small 
category}. See \cite{SGA-4}, \cite{Mac2} or \cite[Section 1.1]{KaSc2}. Another 
approach, involving ``sets'' vs ``classes'',  can be found in \cite{Ne1}. 

We denote by $\cat{Set}$ the category of all small sets. So 
$\opn{Ob}(\cat{Set}) = \cat{U}$, and $\cat{Set}$ is a $\cat{U}$-category.
A group (or a ring, etc.) is called small if its underlying set is
small. We denote by $\cat{Grp}$, $\cat{Ab}$, $\cat{Rng}$ and 
$\cat{Rng_{\mrm{c}}}$ the categories of small groups, small abelian groups, 
small rings  and small commutative rings respectively.  
For a small ring $A$ we denote by $\cat{Mod} A$ the category of 
small left $A$-modules. 

By default we work with $\cat{U}$-categories, and from now on $\cat{U}$ will
remain implicit. There are several places in which we shall encounter 
set theoretical issues (regarding functor categories and localization of 
categories); but these problems can be solved by introducing a bigger universe
$\cat{V}$ such that $\cat{U} \in \cat{V}$.

\mysubsection{Notation and Conventions} \label{subsec:notation}

Let $\cat{C}$ be a category. 
We often write $C \in \cat{C}$ as an abbreviation for 
$C \in \opn{Ob}(\cat{C})$. For an object $C$, its identity automorphism is 
denoted by $\opn{id}_C$. The identity functor of $\cat{C}$ is denoted by 
$\opn{Id}_{\cat{C}}$. 

The opposite category of $\cat{C}$ is $\cat{C}^{\mrm{op}}$.
It has the same objects as $\cat{C}$, but the morphism sets are 
$\opn{Hom}_{\cat{C}^{\mrm{op}}}(C_0, C_1) :=
\opn{Hom}_{\cat{C}}(C_1, C_0)$,
and composition is reversed. Of course 
$(\cat{C}^{\mrm{op}})^{\mrm{op}} = \cat{C}$. 
The identity functor of $\cat{C}$ can be viewed as a contravariant functor 
$\opn{Op} : \cat{C} \to \cat{C}^{\mrm{op}}$. 
To be explicit, on objects we take 
$\opn{Op}(C) := C$. As for morphisms, given a morphism 
$\phi : C_0 \to C_1$ in $\cat{C}$, we let 
$\opn{Op}(\phi) : \opn{Op}(C_1) \to \opn{Op}(C_0)$
be the morphism $\opn{Op}(\phi) := \phi$ in $\cat{C}^{\mrm{op}}$.
The inverse functor $\cat{C}^{\mrm{op}} \to \cat{C}$ is also denoted by 
$\opn{Op}$.
Thus $\opn{Op} \circ \opn{Op} = \opn{Id}$. 

A contravariant functor $F : \cat{C} \to \cat{D}$
is the same as a covariant functor 
$F \circ \opn{Op} : \cat{C}^{\mrm{op}} \to \cat{D}$.
By default all functors will be covariant, unless explicitly mentioned 
otherwise. Contravariant functors will almost always we dealt with by replacing 
the source 
category with its opposite.

\begin{dfn} \label{dfn:4125}
Let $\K$ be a commutative ring. By a 
{\em central $\K$-ring}%
\index{Central $\K$-ring} 
we mean a ring 
$A$, together with a ring homomorphism $\K \to A$, called the structural 
homomorphism, such that the image of 
$\K$ is inside the center of $A$. 

The category of central $\K$-rings, whose morphisms are the ring homomorphisms
$f : A \to B$ that respect the structural homomorphisms from $\K$, is denoted 
by $\catt{Rng} \centover \K$.  
\end{dfn}

Traditionally, a central $\K$-ring was called a ``unital associative 
$\K$-algebra''. Of course all rings and ring homomorphisms are unital.
When $\K = \Z$, a central $\K$-ring is just a ring, and then we 
sometimes use the notation $\catt{Rng}$. 

\begin{exa}  \label{exa:1035}
Let $\K$ be a nonzero commutative ring, and let $n$ be a positive integer. Then 
$\opn{Mat}_n(\K)$, the ring of $n \times n$ matrices with entries in $\K$, is a 
central $\K$-ring. 
\end{exa}

\begin{dfn} \label{dfn:4126}
Let $A$ be a ring. 
We denote by $\cat{Mod} A$, or by the abbreviated notation $\dcat{M}(A)$, the 
category of left $A$-modules.
\end{dfn}

Rings and modules are very important for us, so let us also put forth the next 
convention.

\begin{conv} \label{conv:2490} 
Below are the default implicit assumptions for linear structures and 
operations. 
\begin{enumerate}
\item There is a nonzero commutative base ring $\K$ (e.g.\ the ring of integers 
$\Z$ of a field). 

\item The unadorned tensor symbol $\ot$ means $\ot_{\K}$. 

\item All rings are central $\K$-rings (see Definition 
\ref{dfn:4125}), all ring homomorphisms are over $\K$, and all bimodules are 
$\K$-central. 

\item Generalizing (3), all linear categories are $\K$-linear (see 
Definition \ref{dfn:1083}), and all linear functors are $\K$-linear (see 
Definition \ref{dfn:1031}).

\item For a ring $A$, by all $A$-modules are left $A$-modules, unless 
explicitly stated otherwise. 
\end{enumerate}
\end{conv}

Right $A$-modules are left modules over the opposite ring $A^{\mrm{op}}$, and 
this is the way we shall most often deal with them. 
Morphisms in the categories of groups, rings, $A$-modules, etc.\ will 
usually be called group homomorphism, ring homomorphisms, $A$-module 
homomorphisms, etc., respectively. 

\begin{conv} \label{conv:3170}
We will try to keep the following font and letter conventions: 
\begin{itemize}
\item $f : C \to D$ is a morphism between objects in a category. 
\item $F : \cat{C} \to \cat{D}$ is a functor between categories. 
\item $\eta : F \to G$ is morphism of functors (i.e.\ a natural transformation) 
between functors $F, G : \cat{C} \to \cat{D}$.
\item $\phi, \psi, \phi_i : M \to N$ are morphisms between objects in an 
abelian category $\cat{M}$. 
\item $F : \cat{M} \to \cat{N}$ is a linear functor between abelian 
categories.
\item The category of complexes in an abelian category $\cat{M}$ is 
$\dcat{C}(\cat{M})$.
\item If $\cat{M}$ is a module category, and $M \in \opn{Ob}(\cat{M})$,
then elements of $M$ will be denoted by $m, n, m_i, \ldots$.
\end{itemize}
\end{conv}

\mysubsection{Epimorphisms and Monomorphisms} \label{subsec:zero}
Let $\cat{C}$ be a category. Recall that a morphism $f : C \to D$ in $\cat{C}$ 
is called an {\em isomorphism} if there is a morphism 
$g : D \to C$ such that $f \circ g = \opn{id}_D$ and 
$g \circ f = \opn{id}_C$. The morphism $g$ is called the {\em inverse} of $f$, 
it is unique (if it exists), and it is denoted by $f^{-1}$. 
An isomorphism is often denoted by this shape of arrow: 
$f : C \iso D$. 

A morphism $f : C \to D$ in $\cat{C}$ is called an
{\em epimorphism} if it has the right cancellation property: for every
$g, g' : D \to E$, $g \circ f = g' \circ f$ implies $g = g'$. 
An epimorphism is often denoted by this shape of arrow: 
$f : C \surj D$. 

A morphism $f : C \to D$ is called a {\em monomorphism} if it has the left
cancellation property: for every
$g, g' : E \to C$, $f \circ g = f \circ g'$ implies $g = g'$. 
A monomorphism is often denoted by this shape of arrow: 
$f : C \inj D$. 

\begin{exa}
In $\cat{Set}$ the monomorphisms are the injections, and the epimorphisms are
the surjections. A morphism $f : C \to D$ in $\cat{Set}$ that is both a
monomorphism and an epimorphism is an isomorphism. The same holds in the 
category $\cat{Mod} A$ of left modules over a ring $A$. 
\end{exa}

This example could be misleading, because the property of being an epimorphism 
is often not preserved by forgetful functors, as the next exercise shows. 

\begin{exer} \label{exer:1310}
Consider the category of rings $\cat{Rng}$.
Show that the forgetful functor 
$\cat{Rng} \to \cat{Set}$ respects monomorphisms, but it does not respect
epimorphisms. 
(Hint: show that the inclusion $\Z \to \mbb{Q}$ is an epimorphism in 
$\cat{Rng}$.)  
\end{exer}

By a {\em subobject} of an object $C \in \cat{C}$ we mean a monomorphism 
$f : C' \inj C$ in $\cat{C}$. We sometimes write $C' \subseteq C$ in this 
situation, but this is only notational (and does not mean inclusion of sets). 
We say that two subobjects $f_0 : C'_0 \inj C$ and $f_1 : C'_1 \inj C$ of $C$
are {\em isomorphic} if  there is an isomorphism $g : C'_0 \iso C'_1$ such that 
$f_1 \circ g = f_0$.

Likewise, by a {\em quotient} of $C$ we mean an 
epimorphism $g : C \surj C''$ in $\cat{C}$. There is an analogous notion of 
isomorphic quotients. 

\begin{exer} \label{exer:1676}
Let $\cat{C}$ be a category, and let $C$ be an object of $\cat{C}$. 
\begin{enumerate}
\item Suppose $f_0 : C'_0 \inj C$ and $f_1 : C'_1 \inj C$ are subobjects of
$C$. Show that there is at most one morphism $g : C'_0 \to C'_1$ such that 
$f_1 \circ g = f_0$; and if $g$ exists, then it is a monomorphism. 

\item Show that isomorphism is an equivalence relation on the set of subobjects 
of $C$. Show that the set of equivalence classes of subobjects of $C$ is 
partially ordered by ``inclusion''. (Ignore set-theoretical issues.) 

\item Formulate and prove the analogous statements for quotient objects. 
\end{enumerate}
\end{exer}

An {\em initial object} in a category $\cat{C}$ is an object 
$C_0 \in \cat{C}$,
such that for every object $C \in \cat{C}$ there is exactly one morphism 
$C_0 \to C$. Thus the set $\opn{Hom}_{\cat{C}}(C_0, C)$ is a singleton.
A {\em terminal object} in $\cat{C}$ is an object $C_{\infty} \in \cat{C}$,
such that for every object $C \in \cat{C}$ there is exactly one morphism 
$C \to C_{\infty}$.

\begin{dfn}
A {\em zero object} in a category $\cat{C}$ is an object which is both initial
and terminal. 
\end{dfn}

Initial, terminal and zero objects are unique up to unique isomorphisms (but
they need not exist).

\begin{exa}
In $\cat{Set}$, $\varnothing$ is an initial object, and every singleton is a
terminal object. There is no zero object.
\end{exa}

\begin{exa}
In $\cat{Mod} A$, every trivial module (with only the zero element) is a 
zero object, and we denote this module by $0$. This is allowed, since any other
zero module is uniquely isomorphic to it.  
\end{exa}

\mysubsection{Products and Coproducts}
Let $\cat{C}$ be a category. By a {\em collection of objects of $\cat{C}$} 
indexed by a (small) set $I$, we mean a function $I \to \opn{Ob}(\cat{C})$,
$i \mapsto C_i$. We usually denote this collection using the curly brackets 
notation: $\{ C_i \}_{i \in I}$.

Given a collection $\{ C_i \}_{i \in I}$ of objects of $\cat{C}$, 
its {\em product} is a pair $(C, \{ p_i \}_{i \in I})$
consisting of an object $C \in \cat{C}$, and a collection
$\{ p_i \}_{i \in I}$ of morphisms $p_i : C \to C_i$, called 
{\em projections}. The pair $(C, \{ p_i \}_{i \in I})$ 
must have this universal property: given an object $D$ and morphisms 
$f_i : D \to C_i$, there is a unique morphism $f : D \to C$  s.t.\ 
$f_i = p_i \circ f$. Of course if a product
$(C, \{ p_i \}_{i \in I})$ exists, then it is unique up to a
unique isomorphism; and we usually write
$\prod_{i \in I} C_i := C$, leaving the projection morphisms implicit. 

\begin{exa}
In $\cat{Set}$ and $\cat{Mod} A$ all products exist,
and they are the usual cartesian products. 
\end{exa}

For a collection $\{ C_i \}_{i \in I}$ of objects
of $\cat{C}$, its {\em coproduct} is a pair $(C, \{ e_i \}_{i \in I})$,
consisting of an object $C$ and a collection $\{ e_i \}_{i \in I}$
of morphisms $e_i : C_i \to C$, called {\em embeddings}. 
The pair $(C, \{ e_i \}_{i \in I})$ 
must have this universal property: given an object $D$ and morphisms 
$f_i : C_i \to D$, there is a unique morphism $f : C \to D$  s.t.\ 
$f_i = f \circ e_i$. If a coproduct
$(C, \{ e_i \}_{i \in I})$ exists, then it is unique up to a
unique isomorphism; and we write
$\coprod_{i \in I} C_i := C$,  leaving the embeddings implicit. 

\begin{exa} \label{exa:103}
In $\cat{Set}$ the coproduct is the disjoint union. 
In $\cat{Mod} A$ the coproduct is the direct sum. 
\end{exa}

Product and coproducts are very degenerate cases of {\em limits} and {\em 
colimits} respectively. In this book we will not need to use limits and 
colimits in their most general form. All we shall need is inverse limits and 
direct limits indexed by $\N$; and these will be recalled in Subsection 
\ref{subsec:inv-dir-lim} below. 

We do need to talk about {\em fibered products}. Let $\cat{C}$ be some 
category. Recall that a commutative diagram 
\begin{equation} \label{eqn:3625}
\UseTips \xymatrix @C=6ex @R=6ex {
E
\ar[d]_{g_1}
\ar[r]^{g_2}
&
D_2
\ar[d]^{f_2}
\\
D_1
\ar[r]^{f_1}
&
C
}
\end{equation}
in $\cat{C}$ is called {\em cartesian} if for every object $E' \in \cat{C}$, 
with morphisms 
$g_1' : E' \to D_1$ and $g_2' : E' \to D_2$ that satisfy 
$f_1 \circ g'_1 = f_2 \circ g'_2$, there exists a unique morphism 
$h : E' \to E$ such that 
$g'_i = g_i \circ h$.
See commutative diagram below. 
\[ \UseTips \xymatrix @C=6ex @R=6ex {
E'
\ar@{-->}@(d,ul)[ddr]_{g'_1}
\ar@{-->}@(r,ul)[drr]^{g'_2}
\ar@{-->}[dr]_{h}
\\
&
E
\ar[d]_{g_1}
\ar[r]^{g_2}
&
D_2
\ar[d]^{f_2}
\\
&
D_1
\ar[r]^{f_1}
&
C
} \]

A cartesian diagram is also called a {\em pullback diagram}, and the object $E$ 
is called the {\em fibered product} of $D_1$ and $D_2$ over $C$, with notation
$D_1 \times_C D_2 := E$.
This notation leaves the morphisms implicit. Of course if a fibered product 
exists, than it is unique up to a unique isomorphism that commutes with the 
given arrows. 

There is a dual notion: {\em fibered coproduct}. The input is morphisms 
$C \to D_1$ and $C \to D_2$ in $\cat{C}$, and the fibered coproduct 
$D_1 \sqcup_{C} D_1$ in $\cat{C}$ is just the fibered product in 
the opposite category $\cat{C}^{\mrm{op}}$.

\mysubsection{Equivalence of Categories}
Recall that a functor $F : \cat{C} \to \cat{D}$ is an {\em equivalence} if
there exist a functor $G : \cat{D} \to \cat{C}$, and isomorphisms of 
functors (i.e.\ natural isomorphisms)
$G \circ F \iso \opn{Id}_{\cat{C}}$ and 
$F \circ G \iso  \opn{Id}_{\cat{D}}$. Such a functor $G$ is called a {\em
quasi-inverse} of $F$, and it is unique up to isomorphism (if it exists).  

The functor $F : \cat{C} \to \cat{D}$ is {\em full} (resp.\  {\em faithful})
if for every $C_0, C_1 \in \cat{C}$ the function 
\[ F : \opn{Hom}_{\cat{C}}(C_0, C_1) \to 
\opn{Hom}_{\cat{D}} \bigl( F(C_0), F(C_1) \bigr) \]
is surjective (resp.\ injective). 

We know that $F : \cat{C} \to \cat{D}$ is an equivalence if and only if these 
two conditions hold:
\begin{enumerate}
\rmitem{i} $F$ is essentially surjective on objects. This means that for every
$D \in \cat{D}$  there is an isomorphism $F(C) \iso D$ for some $C \in \cat{C}$.

\rmitem{ii} $F$ is fully faithful (i.e.\ full and faithful). 
\end{enumerate}

\begin{exer}
If you are not sure about the last claim (characterization of equivalences), 
then prove it. (Hint: use the axiom of choice to construct a quasi-inverse of 
$F$.)  
\end{exer}

A functor $F : \cat{C} \to \cat{D}$ is called an {\em isomorphism of 
categories} if it is bijective on sets of objects and on sets of morphisms. It 
is clear that an isomorphism of categories is an equivalence. If $F$ is an 
isomorphism of categories, then it has an inverse isomorphism 
$F^{-1} : \cat{D} \to \cat{C}$, which is unique.

\mysubsection{Bifunctors} \label{subsec:bifunc}
Let $\cat{C}$ and $\cat{D}$ be categories. Their product is the category 
$\cat{C} \times \cat{D}$ defined as follows: the set of objects is 
\[ \opn{Ob}(\cat{C} \times \cat{D}) := 
\opn{Ob}(\cat{C}) \times \opn{Ob}(\cat{D}) . \] 
The sets of morphisms are 
\[ \opn{Hom}_{\cat{C} \times \cat{D}} \big( (C_0, D_0), (C_1, D_1) \big) :=
\opn{Hom}_{\cat{C} } (C_0, C_1) \times 
\opn{Hom}_{\cat{D} } (D_0, D_1)  . \]
The composition is 
\[ (f_1, g_1) \circ (f_0, g_0) := 
(f_1 \circ f_0, g_1 \circ g_0) , \]
and the identity morphism on an object $(C, D)$ is
$( \opn{id}_C,  \opn{id}_D)$. 

A {\em bifunctor} from $(\cat{C}, \cat{D})$ to $\cat{E}$ is a functor 
$F : \cat{C} \times \cat{D} \to \cat{E}$.
The extra information that is implicit when we call $F$ a bifunctor 
is that the source category $\cat{C} \times \cat{D}$ is a product.

\mysubsection{Representable Functors} \label{subsec:repfunc}

Let $\cat{C}$ be a category. An object $C \in \cat{C}$ gives rise 
to a functor 
\begin{equation} \label{eqn:3025}
\opn{Y}_{\cat{C}}(C) : \cat{C}^{\mrm{op}} \to \cat{Set} , \quad
\opn{Y}_{\cat{C}}(C) := \opn{Hom}_{\cat{C}}(-, C) .  
\end{equation}
Explicitly, the functor $\opn{Y}_{\cat{C}}(C)$ sends an object 
$D \in \cat{C}$ to the set 
$\opn{Y}_{\cat{C}}(C)(D) :=  \opn{Hom}_{\cat{C}}(D, C)$,
and a morphism
$\psi : D_0 \to D_1$ in $\cat{C}$ goes to the function 
\[ \opn{Y}_{\cat{C}}(C)(\psi) := \opn{Hom}_{\cat{C}}(\psi,  \opn{id}_C) : 
\opn{Hom}_{\cat{C}}(D_1, C) \to \opn{Hom}_{\cat{C}}(D_0, C) . \]

Now suppose we are given a morphism 
$\phi : C_0 \to C_1$ in $\cat{C}$. There is a morphism of functors (a natural 
transformation) 
\begin{equation} \label{eqn:3026}
\opn{Y}_{\cat{C}}(\phi) := \opn{Hom}_{\cat{C}}(-, \phi) : 
\opn{Y}_{\cat{C}}(C_0) \to \opn{Y}_{\cat{C}}(C_1) .
\end{equation}

Consider the category 
$\cat{Fun}(\cat{C}^{\mrm{op}}, \cat{Set})$,
whose objects are the functors $F : \cat{C}^{\mrm{op}} \to \cat{Set}$, and
whose morphisms are the morphisms of functors
$\eta : F_0 \to F_1$.  
There is a set-theoretic difficulty here: the sets of objects and morphisms 
of $\cat{Fun}(\cat{C}^{\mrm{op}}, \cat{Set})$
are too big (unless $\cat{C}$ is a small category), and this is not a 
$\cat{U}$-category. Hence we must enlarge the universe, as mentioned in 
Subsection \ref{subsec:set-theor}. 

\begin{dfn} \label{dfn:3025}
The {\em Yoneda functor} of the category $\cat{C}$ is the functor
\[ \opn{Y}_{\cat{C}} : \cat{C} \to \cat{Fun}(\cat{C}^{\mrm{op}}, \cat{Set}) \]
described by formulas (\ref{eqn:3025}) and (\ref{eqn:3026}).
\end{dfn}

\begin{thm}[Yoneda Lemma] \label{thm:3025}
The Yoneda functor $\opn{Y}_{\cat{C}}$ is fully faithful. 
\end{thm}

See \cite[Section III.2]{Mac2} or \cite[Section 1.4]{KaSc2} for a proof. The 
proof is not hard, but it is very confusing. 

A functor $F : \cat{C}^{\mrm{op}} \to \cat{Set}$ is called 
{\em representable} if there is an isomorphism of functors 
$\eta : F \iso \opn{Y}_{\cat{C}}(C)$ for some object $C \in \cat{C}$.
Such an object $C$ is said to {\em represent} the functor $F$. 
The Yoneda Lemma says that $\opn{Y}_{\cat{C}}$ is an equivalence
from $\cat{C}$ to the category of representable functors. 
Thus the pair $(C, \eta)$ is unique up to a unique isomorphism (if it exists).
Note that the isomorphism of sets
$\eta_C : F(C) \iso \opn{Y}_{\cat{C}}(C)(C)$ 
gives a special element $\til{\eta} \in F(C)$ such that 
$\eta_C(\til{\eta}) = \opn{id}_C$. 

Dually, any object $C \in \cat{C}$ gives rise to a functor
\begin{equation} \label{eqn:3027}
\opn{Y}_{\cat{C}}^{\vee}(C) : \cat{C} \to \cat{Set} , \quad
\opn{Y}_{\cat{C}}^{\vee}(C) := \opn{Hom}_{\cat{C}}(C, -) . 
\end{equation}
A morphism $\phi : C_0 \to C_1$ in $\cat{C}$ induces a morphism of 
functors
\begin{equation} \label{eqn:3028}
\opn{Y}_{\cat{C}}^{\vee}(\phi) := \opn{Hom}_{\cat{C}}(\phi, -) : 
\opn{Y}_{\cat{C}}^{\vee}(C_1) \to \opn{Y}_{\cat{C}}^{\vee}(C_0) .
\end{equation}

\begin{dfn} \label{dfn:3026}
The {\em dual Yoneda functor} of the category $\cat{C}$ is the functor
\[ \opn{Y}_{\cat{C}}^{\vee} : \cat{C}^{\mrm{op}} \to 
\cat{Fun}(\cat{C}, \cat{Set}) \]
described by formulas (\ref{eqn:3027}) and (\ref{eqn:3028}).
\end{dfn}

\begin{thm}[Dual Yoneda Lemma] \label{thm:3026}
The dual Yoneda functor $\opn{Y}_{\cat{C}}^{\vee}$ is fully faithful. 
\end{thm}

This is also proved in \cite[Section III.2]{Mac2} and 
\cite[Section 1.4]{KaSc2}. 

A functor $F : \cat{C} \to \cat{Set}$ is called {\em corepresentable} if there 
is an isomorphism of functors
$\eta : F \iso \opn{Y}_{\cat{C}}^{\vee}(C)$ for some object $C \in \cat{C}$.
The object $C$ is said to {\em corepresent} the functor $F$.
The dual Yoneda Lemma says that the functor $\opn{Y}_{\cat{C}}^{\vee}$ is 
an equivalence from $\cat{C}^{\mrm{op}}$ to the category of corepresentable 
functors. The identity automorphism $\opn{id}_C$ corresponds to 
a special element $\til{\eta} \in F(C)$.

\mysubsection{Inverse and Direct Limits} \label{subsec:inv-dir-lim}
 
We are only interested in direct and inverse limits indexed by the ordered set 
$\N$. For a general discussion see \cite{Mac1} or \cite{KaSc1}. 

Let $\cat{C}$ be a category. Recall that an 
{\em $\N$-indexed direct system} in $\cat{C}$ is data  
\[ \bigl( \{ C_k \}_{k \in \N}, \{ \mu_{k} \}_{k \in \N} \bigr) , \]
where $C_k$ are objects of $\cat{C}$, and 
$\mu_k : C_k \to C_{k + 1}$
are morphisms, that we call {\em transitions}.
A {\em direct limit} of this system is  data  
$\bigl( C , \{ \ep_{k} \}_{k \in \N} \bigr)$,
where $C \in \cat{C}$, and $\ep_k : C_k \to C$ are morphisms, that we call 
{\em abutments}, such that $\ep_{k + 1} \circ \mu_k = \ep_k$ for all $k$. The 
universal property required is this: if 
$\bigl( C' , \{ \ep'_{k} \}_{k \in \N} \bigr)$
is another pair such that $\ep'_{k + 1} \circ \mu_k = \ep'_k$, then
there is a unique morphism 
$\ep : C \to C'$ such that 
$\ep'_{k} = \ep \circ \ep_k$. 
If a direct limit $C$ exists, then of course it is unique, up to a unique 
isomorphism. We then write 
$\lim_{k \to} C_k := C$,
and call this the direct limit of the system $\{ C_k \}_{k \in \N}$, keeping 
the transitions and the abutments implicit. 
Sometimes we look at the morphisms 
\[ \mu_{k_0, k_1} := \mu_{k_1 - 1} \circ \cdots \circ \mu_{k_0} :
C_{k_0} \to C_{k_1} \]
for $k_0 < k_1$, and $\mu_{k, k} := \opn{id}_{C_k}$. 

By an {\em $\N$-indexed inverse system} in the category $\cat{C}$ we mean 
data 
\[ \bigl( \{ C_k \}_{k \in \N}, \{ \mu_k \}_{k \in \N} \bigr) , \]
where $\{ C_k \}_{k \in \N}$ is a collection of objects, and 
$\mu_k : C_{k + 1} \to C_{k}$ 
are morphisms, also called transitions. 
An {\em inverse limit} of this system is data  
$\bigl( C , \{ \ep_{k} \}_{k \in \N} \bigr)$,
where $C \in \cat{C}$, and $\ep_k : C \to C_k$ are morphisms, that we also call 
abutments, such that $\mu_k \circ \ep_{k + 1} = \ep_{k}$. These satisfy an 
analogous universal property. 
If an inverse limit $C$ exists, then it is unique, up to a unique 
isomorphism. We then write 
$\lim_{\lar k} C_k := C$,
and we call this the inverse limit of the system $\{ C_k \}_{k \in \N}$.
We define the morphisms 
\[ \mu_{k_0, k_1} := \mu_{k_0} \circ \cdots \circ \mu_{k_1 - 1} :
C_{k_1} \to C_{k_0} \]
for $k_0 < k_1$, and $\mu_{k, k} := \opn{id}_{C_k}$. 

\begin{exer} \label{exer:3010}
We can view the ordered set $\N$ as a category, with a single 
morphism $k \to l$ when $k \leq l$, and no morphisms otherwise. 
\begin{enumerate}
\item Interpret $\N$-indexed direct and inverse systems in $\cat{C}$ as 
functors $F : \N \to \cat{C}$ and $G : \N^{\mrm{op}} \to \cat{C}$ respectively. 

\item Let $\bar{\N} := \N \cup \{ \infty \}$. 
Interpret the direct and inverse limits of $F$ and $G$ respectively as 
functors $\bar{F} : \bar{\N} \to \cat{C}$ and 
$\bar{G} : \bar{\N}^{\mrm{op}} \to \cat{C}$,
extending $F$ and $G$, with suitable universal properties. 
\end{enumerate}
\end{exer}

\begin{exer} \label{exer:3011}
Prove that $\N$-indexed direct and inverse limits exist in the categories 
$\cat{Set}$ and $\cat{Mod} A$, for any ring $A$. Give explicit formulas.
\end{exer}

\begin{exa} \label{ea:1740}
Let $\cat{M}$ be the category of finite abelian groups. 
The inverse system $\{ M_k \}_{k \in \N}$, where
$M_k := \Z / (2^k)$, and the transition 
$\mu_k :  M_{k + 1} \to M_k$ is the canonical surjection, 
does not have an inverse limit in $\cat{M}$. We can also make 
$\{ M_k \}_{k \in \N}$ into a direct system, in which the transition
$\nu_k : M_k \to  M_{k + 1}$ 
is multiplication by $2$. The direct limit does not exist in $\cat{M}$.
\end{exa}

If  $\{ C_k \}_{k \in \N}$ is a direct system in $\cat{C}$, and 
$D \in \cat{C}$ is any object, then there is an induced inverse system 
$\bigl\{ \opn{Hom}_{\cat{C}}(C_k, D) \}_{k \in \N}$
in $\cat{Set}$, and it has a limit. If 
$C := \lim_{k \to} C_k$ exists, then the abutments $\ep_k : C_k \to C$ induce 
a morphism 
\begin{equation} \label{eqn:3010}
\opn{Hom}_{\cat{C}}(C, D) \to 
\lim_{\leftarrow k} \, \opn{Hom}_{\cat{C}}(C_k, D) 
\end{equation}
in $\cat{Set}$. 

Similarly, if $\{ C_k \}_{k \in \N}$ is an inverse system in $\cat{C}$, 
and $D \in \cat{C}$ is any object, then there is an induced inverse system 
$\bigl\{ \opn{Hom}_{\cat{C}}(D, C_k) \}_{k \in \N}$
in $\cat{Set}$, and it has a limit. If $C := \lim_{\lar k} C_k$ exists,
then the abutments $\ep_k : C \to C_k$ induce a morphism 
\begin{equation} \label{eqn:3011}
\opn{Hom}_{\cat{C}}(D, C) \to 
\lim_{\leftarrow k} \, \opn{Hom}_{\cat{C}}(D, C_k) , 
\end{equation}
in $\cat{Set}$.

\begin{prop} \label{prop:1535}
Let $\cat{C}$ be a category. 
\begin{enumerate}
\item Let  $\{ C_k \}_{k \in \N}$ be a direct system in $\cat{C}$,
and let $\bigl( C , \{ \ep_{k} \}_{k \in \N} \bigr)$ be data as in the 
definition of a direct limit. Then  
$C = \lim_{k \to} C_k$ if and only if for every object 
$D \in \cat{C}$, the function \tup{(\ref{eqn:3010})} is bijective. 

\item Let  $\{ C_k \}_{k \in \N}$ be an inverse system in $\cat{C}$,
and let $\bigl( C , \{ \ep_{k} \}_{k \in \N} \bigr)$ be data as in the 
definition of an inverse limit. Then  
$C = \lim_{\lar k} C_k$ if and only if for every object 
$D \in \cat{C}$, the function \tup{(\ref{eqn:3011})} is bijective. 
\end{enumerate}
\end{prop}

\begin{exer}  \label{exer:1535}
Prove Proposition \ref{prop:1535}.
\end{exer}

\begin{rem} \label{rem:4090}
In many cases inverse and direct limits do not exist in $\cat{C}$ because ``its 
objects are too small''. This happens in the category 
$\cat{Set}_{\mrm{fin}}$ of finite sets, and also in the category 
$\catt{Ab}_{\mrm{fin}}$ of finite abelian groups (see Example \ref{ea:1740} 
above).  

There is a very effective method to enlarge $\cat{C}$ just enough so that the 
bigger category will have the desired limits. This is done by means of the 
categories $\cat{Ind}_{}(\cat{C})$ and $\cat{Pro}_{}(\cat{C})$ of {\em 
ind-objects} and {\em pro-objects} of $\cat{C}$, respectively. 
See [KaSc1, Sec 1.11] and [KaSc2, Sec 6.1] for detailed discussions.

Here are two examples. For $\cat{C} = \cat{Set}_{\mrm{fin}}$, the 
category $\cat{Ind}_{}(\cat{Set}_{\mrm{fin}})$ is (canonically equivalent to) 
$\cat{Set}$. 
The category $\cat{Pro}_{}(\cat{Set}_{\mrm{fin}})$ is the category of 
totally disconnected compact Hausdorff topological spaces, 
whose morphisms are the continuous functions. 

For $\cat{C} = \catt{Ab}_{\mrm{fin}}$, 
the category $\cat{Ind}_{}(\cat{Ab}_{\mrm{fin}})$ of ind-objects is the 
category $\cat{Ab}$ of abelian groups. The category 
$\cat{Pro}_{}(\cat{Ab}_{\mrm{fin}})$ of pro-objects is the category of profinite
abelian groups, whose morphisms are the continuous group homomorphisms. 
\end{rem}

\cleardoublepage
\mysection{Abelian Categories and Additive Functors} \label{sec:abelian}

\AYcopyright

The concept of {\em abelian category} is an extremely useful abstraction of a
category of modules. It was introduced by A. Grothendieck in his foundational 
paper \cite{Gro} from 1957. 
References for this material are \cite{Mac1}, \cite{Mac2}, 
\cite{Rot}, \cite{HiSt}, \cite{We}, \cite{KaSc1} and \cite{KaSc2}.

\mysubsection{Linear Categories} \label{subsec:1}

\begin{dfn} \label{dfn:1083}
Let $\K$ be a commutative ring. 
A  {\em $\K$-linear category}%
\index{Linear category}%
\index{Linear category!K@$\K$-linear category}
is a category $\cat{M}$, endowed with 
a $\K$-module structure on each of the sets of morphisms  
$\opn{Hom}_{\cat{M}}(M_0, M_1)$. The condition is this:
\begin{itemize}
\item For all $M_0, M_1, M_2 \in \cat{M}$ the composition function 
\[ \begin{aligned}
& \opn{Hom}_{\cat{M}}(M_1, M_2) \times \opn{Hom}_{\cat{M}}(M_0, M_1)
\to \opn{Hom}_{\cat{M}}(M_0, M_2) \\
& (\phi_1, \phi_0) \mapsto \phi_1 \circ \phi_0 
\end{aligned} \]
is $\K$-bilinear.
\end{itemize}
If $\K = \Z$, we say that $\cat{M}$ is a {\em linear category}.
\end{dfn}

As already mentioned in Convention \ref{conv:2490}, in our book all linear 
categories are $\K$-linear. 

\begin{prop} \label{prop:1035}
Let $\cat{M}$ be a $\K$-linear category.
\begin{enumerate}
\item For every object $M \in \cat{M}$, the set 
$\opn{End}_{\cat{M}}(M) := \opn{Hom}_{\cat{M}}(M, M)$,
with its given addition operation, and with the operation of composition,
is a central $\K$-ring. 

\item For every two objects $M_0, M_1 \in \cat{M}$, the set 
$\opn{Hom}_{\cat{M}}(M_0, M_1)$,
with its given addition operation, and with the operations of composition,
is a left module over the ring $\opn{End}_{\cat{M}}(M_1)$, and 
a right module over the ring $\opn{End}_{\cat{M}}(M_0)$.
Furthermore, these left and right actions commute with each other. 
\end{enumerate}
\end{prop}

\begin{exer} \label{exer:3255}
Prove Proposition \ref{prop:1035}.
\end{exer}

This result can be reversed:

\begin{exa} \label{exa:104}
Let $A$ be a central $\K$-ring.  
Define a category $\cat{M}$ like this: there is a single object $M$, and its
set of morphisms is
$\opn{Hom}_{\cat{M}}(M, M) := A$.
Composition in $\cat{M}$ is the multiplication of $A$. Then $\cat{M}$ is a
$\K$-linear category.
\end{exa}

For a central $\K$-ring $A$, the opposite ring $A^{\mrm{op}}$ has the same 
$\K$-module structure as $A$, but the multiplication is reversed. 

\begin{exer} \label{exer:1030}
Let $A$ be a nonzero ring. Let $P, Q \in \cat{Mod} A$ 
be distinct free $A$-modules of rank $1$.
\begin{enumerate}
\item Prove that there is a ring isomorphism 
$\opn{End}_{\cat{Mod} A}(P) \cong A^{\mrm{op}}$.
Is this ring isomorphism canonical?

\item Let $\cat{M}$ be the full subcategory of $\cat{Mod} A$ on the set of 
objects $\{ P, Q \}$. Compare the linear category $\cat{M}$ to the ring of 
matrices $\opn{Mat}_2(A^{\mrm{op}})$.
\end{enumerate}
\end{exer}

\mysubsection{Additive Categories}

\begin{dfn} 
An  {\em additive category}%
\index{Additive category}
is a linear category $\cat{M}$
satisfying these two conditions:
\begin{enumerate}
\rmitem{i} $\cat{M}$ has a zero object $0$.

\rmitem{ii} $\cat{M}$ has finite coproducts.
\end{enumerate}
\end{dfn}

Observe that $\opn{Hom}_{\cat{M}}(M, N) \neq \varnothing$
for every $M, N \in \cat{M}$, since this is an abelian group. For the zero 
object $0 \in \cat{M}$ we have 
$\opn{Hom}_{\cat{M}}(M, 0) = \opn{Hom}_{\cat{M}}(0, M) = 0$, 
the zero abelian group.  We denote the unique arrows $0 \to M$ and 
$M \to 0$ also by $0$. So the numeral $0$ has a lot of meanings; but they are 
always (hopefully) clear from the context. 
The coproduct in a linear category $\cat{M}$ is usually denoted by $\bigoplus$, 
and is called the {\em direct sum}; cf.\ Example \ref{exa:103}.

\begin{exa}
Let $A$ be a central $\K$-ring. The category $\cat{Mod} A$ is a $\K$-linear 
additive category. 
The full subcategory $\cat{F} \subseteq \cat{Mod} A$ on the free modules is
also additive. 
\end{exa}

\begin{prop} \label{prop:1}
Let $\cat{M}$ be a linear category. Let $\{ M_i \}_{i \in I}$ be a  
finite collection of objects of $\cat{M}$, and assume the coproduct
$M = \bigoplus_{i \in I} M_i$ exists, with embeddings 
$e_i : M_i \to M$. 

\begin{enumerate}
\item For each $i$ let $p_i : M \to M_i$ be the unique morphism 
such that $p_i \circ e_i = \opn{id}_{M_i}$, and 
 $p_i \circ e_j = 0$ for $j \neq i$. 
Then $(M, \{ p_i \}_{i \in I})$ is a product of the collection 
$\{ M_i \}_{i \in I}$.

\item $\sum_{i \in I} e_i \circ p_i = \opn{id}_M$.
\end{enumerate}
\end{prop}

\begin{exer} \label{exer:2490}
Prove this proposition.
\end{exer}

Part (1) of Proposition \ref{prop:1} directly implies:

\begin{cor} \label{cor:1030}
An additive category has finite products. 
\end{cor}

\begin{dfn} \label{dfn:1025}
Let $\cat{M}$ be an additive category, and let $\cat{N}$ be a full subcategory
of $\cat{M}$. We say that $\cat{N}$ is a {\em full additive subcategory} of
$\cat{M}$ if $\cat{N}$ contains the zero object, and is closed under finite 
direct sums.
\end{dfn}

\begin{exer}
In the situation of Definition \ref{dfn:1025}, show that the category $\cat{N}$ 
is itself additive.
\end{exer}

\begin{exa} \label{exa:1030}
Consider the linear category $\cat{M}$ from Example \ref{exa:104},
built from a ring $A$. It does not have a zero object (unless the ring $A$ is 
the zero ring), so it is not additive. 

A more puzzling question is this: Does $\cat{M}$ have finite direct 
sums? This turns out to be equivalent to whether or not $A \cong A \oplus A$ as 
right $A$-modules. 
One can show that when $A$ is nonzero and commutative, or nonzero and 
noetherian, then $A \not\cong A \oplus A$ in $\cat{Mod} A^{\mrm{op}}$.
On the other hand, if we take a field $\K$, and a countable rank $\K$-module
$N$, then $A := \opn{End}_{\K}(N)$ will satisfy $A \cong A \oplus A$.
\end{exa}

\begin{prop} \label{prop:1030}
Let $\cat{M}$ be a linear category, and $N \in \cat{M}$. The following 
conditions are equivalent\tup{:}
\begin{enumerate}
\rmitem{i} The ring  $\opn{End}_{\cat{M}}(N)$ is trivial. 
\rmitem{ii} $N$ is a zero object of $\cat{M}$.
\end{enumerate}
\end{prop}

\begin{proof} \mbox{}

\smallskip \noindent
(ii) $\Rightarrow$ (i): Since the set $\opn{End}_{\cat{M}}(N)$ is a singleton, 
it must be the trivial ring. 

\medskip \noindent 
(i) $\Rightarrow$ (ii): If the ring $\opn{End}_{\cat{M}}(N)$ is 
trivial, then all left and right modules over it must be trivial. 
Now use Proposition \ref{prop:1035}(2). 
\end{proof}

\mysubsection{Abelian Categories} \label{subsec:ab-cats}

\begin{dfn}
Let $\cat{M}$ be an additive category, and let $f : M \to N$ be a morphism in
$\cat{M}$. A {\em kernel} of $f$ is a pair $(K, k)$, consisting of an object
$K \in \cat{M}$ and a morphism $k : K \to M$, with these two properties:
\begin{enumerate}
\rmitem{i}  $f \circ k = 0$. 
\rmitem{ii} If $k' : K' \to M$ is a morphism in $\cat{M}$ such that 
$f \circ k' = 0$, then there is a unique morphism $g : K' \to K$ 
such that $k' = k \circ g$.
\end{enumerate}
\end{dfn}

In other words, the object $K$ represents the functor 
$\cat{M}^{\mrm{op}} \to \cat{Ab}$, 
\[ K' \mapsto \{ k' \in \opn{Hom}_{\cat{M}}(K', M) \mid f \circ k' = 0 \} . \]
The kernel of $f$ is of course unique up to a unique isomorphism (if it
exists), and we denote if by $\opn{Ker}(f)$. Sometimes $\opn{Ker}(f)$ refers
only to the object $K$, and other times it refers only to the morphism
$k$; as usual, this should be clear from the context.

\begin{dfn}
Let $\cat{M}$ be an additive category, and let $f : M \to N$ be a morphism in
$\cat{M}$. A {\em cokernel} of $f$ is a pair $(C, c)$, consisting of an object
$C \in \cat{M}$ and a morphism $c : N \to C$, with these two properties:
\begin{enumerate}
\rmitem{i}  $c \circ f = 0$. 
\rmitem{ii} If $c' : N \to C'$ is a morphism in $\cat{M}$ such that 
$c' \circ f = 0$, then there is a unique morphism $g : C \to C'$ 
such that $c' = g \circ c$.
\end{enumerate}
\end{dfn}

In other words, the object $C$ corepresents the functor 
$\cat{M} \to \cat{Ab}$, 
\[ C' \mapsto \{ c' \in \opn{Hom}_{\cat{M}}(N, C') \mid c' \circ f = 0 \} . \]
The cokernel of $f$ is of course unique up to a unique isomorphism (if it
exists), and we denote if by $\opn{Coker}(f)$. Sometimes $\opn{Coker}(f)$ 
refers only to the object $C$, and other times it refers only to the morphism
$c$; as usual, this should be clear from the context.

\begin{exa}
In $\cat{Mod} A$ all kernels and cokernels exist. 
Given $f : M \to N$, the kernel is 
$k : K \to M$, where 
$K := \{ m \in M \mid f(m) = 0 \}$,
and the $k$ is the inclusion. The cokernel is $c : N \to C$, where 
$C := N / f(M)$, and $c$ is the canonical projection.
\end{exa}

\begin{prop}
Let  $\cat{M}$ be an additive category, and let $f : M \to N$ be a morphism 
in $\cat{M}$.
\begin{enumerate}
\item If $k : K \to M$ is a kernel of $f$, then $k$  is a monomorphism.
\item If $c : N \to C$ is a cokernel of $f$, then $c$ is an epimorphism.
\end{enumerate}
\end{prop}

\begin{exer} 
Prove the proposition.
\end{exer}

\begin{dfn}
Assume the additive category $\cat{M}$ has kernels and cokernels. 
Let $f : M \to N$ be a morphism in $\cat{M}$.
\begin{enumerate}
\item Define the {\em image} of $f$ to be 
$\opn{Im}(f) := \opn{Ker}(\opn{Coker}(f))$.

\item Define the {\em coimage} of $f$ to be 
$\opn{Coim}(f) := \opn{Coker}(\opn{Ker}(f))$.
\end{enumerate}
\end{dfn}

The image is familiar, but the coimage is probably not. The next diagram 
should help. 
We start with a morphism $f : M \to N$ in $\cat{M}$. 
The kernel and cokernel of $f$ fit into this diagram:
$K \xar{k} M \xar{f} N \xar{c} C$. 
Inserting $\al := \opn{Coker}(k) = \opn{Coim}(f)$ and 
$\be := \opn{Ker}(c) = \opn{Im}(f)$
we get the following commutative diagram (solid arrows):
\begin{equation} \label{eqn:1015}
\UseTips  \xymatrix @C=8ex @R=6ex {
K
\ar[r]^{k}
\ar[dr]_{0}
& 
M
\ar[r]^{f}
\ar[d]_{\al}
\ar@{-->}[dr]^{\ga}
&
N
\ar[r]^{c}
&
C
\\
&
M'
\ar@{-->}[r]_{f'}
&
N'
\ar[ur]_{0}
\ar[u]_{\be}
}
\end{equation}
Since $c \circ f = 0$ there is a unique morphism $\ga$ (the dashed arrow)
making the diagram commutative. Now 
$\be \circ \ga \circ k = f \circ k = 0$; and $\be$ is a monomorphism; so 
$\ga \circ k = 0$. Hence there is a unique morphism $f' : M' \to N'$ 
making the diagram commutative.
We conclude that $f : M \to N$ induces a morphism 
\begin{equation}  \label{eqn:2}
f' : \opn{Coim}(f) \to \opn{Im}(f) . 
\end{equation}

\begin{dfn} \label{dfn:1}
An {\em abelian category}%
\index{Abelian category}
is an additive category $\cat{M}$ with these two extra properties:
\begin{enumerate}
\rmitem{i} All morphisms in $\cat{M}$ admit kernels and cokernels.
\rmitem{ii} For every morphism $f : M \to N$ in $\cat{M}$, the induced morphism 
$f'$ in equation (\ref{eqn:2}) is an isomorphism.
\end{enumerate}
\end{dfn}

Here is a less precise but (maybe) easier to remember way to state property 
(ii). Because $M' = \opn{Coker}(\opn{Ker}(f))$ and 
$N' = \opn{Ker}(\opn{Coker}(f))$, 
we see that 
\begin{equation} \label{eqn:1019}
\opn{Coker}(\opn{Ker}(f)) = \opn{Ker}(\opn{Coker}(f))  .
\end{equation}

{}From now on we forget all about the coimage. 

\begin{exer}
For any ring $A$, prove that the category $\cat{Mod} A$ is abelian. 
\end{exer}

This includes the category $\cat{Ab} = \cat{Mod} \Z$, from which the name 
derives. 

\begin{dfn} \label{dfn:1015}
Let $\cat{M}$ be an abelian category, and let $\cat{N}$ be a full subcategory
of $\cat{M}$. We say that $\cat{N}$ is a 
{\em full abelian subcategory}%
\index{Abelian category! full  abelian subcategory of}
of $\cat{M}$ if the zero object belongs to $\cat{N}$, and $\cat{N}$ is closed  
in $\cat{M}$ under taking finite direct sums, kernels and cokernels.
\end{dfn}

\begin{exer}
In the situation of Definition \ref{dfn:1015}, show that the category $\cat{N}$ 
is itself abelian.
\end{exer}

\begin{exa}
Let $\cat{M}_1$ be the category of finitely generated abelian groups,
and let $\cat{M}_0$ be the category of finite abelian groups. 
Then $\cat{M}_1$ is a full abelian subcategory of $\cat{Ab}$, and 
$\cat{M}_0$ is a full abelian subcategory of $\cat{M}_1$. 
\end{exa}

\begin{exer}
Let $\cat{N}$ be the full subcategory of $\cat{Ab}$ whose objects are the
finitely generated free abelian groups. It is an additive subcategory of
$\cat{Ab}$  (since it is closed under direct sums).
\begin{enumerate}
\item Show that $\cat{N}$ is closed under kernels in $\cat{Ab}$.

\item Show that $\cat{N}$ is not closed under cokernels in $\cat{Ab}$, so it 
is not a full abelian subcategory of $\cat{Ab}$.

\item Show that $\cat{N}$ has cokernels (not the same as those of $\cat{Ab}$).
Still, it fails to be an abelian category.
\end{enumerate}
\end{exer}

\begin{exer}
The category $\cat{Grp}$ is not linear of course. Still, it does have a 
zero object (the trivial group). Show that $\cat{Grp}$ has kernels and 
cokernels, but condition (ii) of Definition \ref{dfn:1} fails. 
\end{exer}

\begin{exer}
Let $\cat{Hilb}$ be the category of Hilbert spaces over $\mbb{C}$. 
The morphisms are the continuous $\mbb{C}$-linear homomorphisms. Show that 
$\cat{Hilb}$ is a $\mbb{C}$-linear additive category with kernels and 
cokernels, 
but it is not an abelian category.
\end{exer}

\begin{exer}
Let $A$ be a ring. Show that $A$ is {\em left noetherian} iff the category 
$\cat{Mod}_{\mrm{f}} A$ of  finitely generated left modules is a full abelian 
subcategory of $\cat{Mod} A$. 
\end{exer}

\begin{exa} \label{exa:1}
Let $(X, \mcal{A})$ be a {\em ringed space}; namely $X$ is a topological space 
and $\mcal{A}$ is a sheaf of rings on $X$. 
There is a very detailed discussion of sheaves of modules in 
\cite[Section II.1]{Har} and \cite[Sections 2.1-2.2]{KaSc1}. 
Let us mention only one aspect of it. 

We denote by $\cat{Mod}^{\mrm{pr}} \mcal{A}$ the category of {\em presheaves} 
of left $\mcal{A}$-modules on $X$. This is an abelian category. 
Given a morphism $\phi : \mcal{M} \to \mcal{N}$ in 
$\cat{Mod}^{\mrm{pr}} \mcal{A}$, its kernel is the presheaf 
$\opn{Ker}^{\mrm{pr}}(\phi)$ defined by 
\[ \Gamma \bigl( U, \opn{Ker}^{\mrm{pr}}(\phi) \bigr) := 
\opn{Ker} \big( \phi : \Gamma(U, \mcal{M}) \to \Gamma(U, \mcal{N}) \big)  \]
on every open set $U \subseteq X$. 
The cokernel of $\phi$ is the  presheaf $\opn{Coker}^{\mrm{pr}}(\phi)$
defined by 
\[ \Gamma \bigl( U,  \opn{Coker}^{\mrm{pr}}(\phi) \bigr) : 
\opn{Coker} \big( \phi : \Gamma(U, \mcal{M}) \to \Gamma(U, \mcal{N}) \big) . \]

Now let $\cat{Mod} \mcal{A}$ be the full subcategory of 
$\cat{Mod}^{\mrm{pr}} \mcal{A}$ consisting of {\em sheaves}. 
It is a full additive subcategory of $\cat{Mod}^{\mrm{pr}} \mcal{A}$, closed 
under kernels. However, $\cat{Mod} \mcal{A}$ is not closed under 
cokernels inside $\cat{Mod}^{\mrm{pr}} \mcal{A}$,
and hence it is not a full abelian subcategory.

Nonetheless, $\cat{Mod} \mcal{A}$ is itself an abelian category, but
with different cokernels. Given a morphism 
$\phi : \mcal{M} \to \mcal{N}$ in $\cat{Mod} \mcal{A}$,
its cokernel in $\cat{Mod} \mcal{A}$ is the sheaf $\opn{Coker}(\phi)$ that's  
associated to the presheaf $\opn{Coker}^{\mrm{pr}}(\phi)$.
\end{exa}

\begin{prop} \label{prop:2490}
Let $\cat{M}$ be a linear category.
\begin{enumerate}
\item The opposite category $\cat{M}^{\mrm{op}}$ has a canonical structure 
of a linear category. 

\item If $\cat{M}$ is additive, then $\cat{M}^{\mrm{op}}$ is also additive.

\item If $\cat{M}$ is abelian, then $\cat{M}^{\mrm{op}}$ is also abelian.
\end{enumerate}
\end{prop}

\begin{proof}
(1) Since
$\opn{Hom}_{\cat{M}^{\mrm{op}}}(M, N) = \opn{Hom}_{\cat{M}}(N, M)$,
this is an abelian group. The bilinearity of the composition in 
$\cat{M}^{\mrm{op}}$ is clear.

\medskip \noindent 
(2) The zero objects in $\cat{M}$ and $\cat{M}^{\mrm{op}}$ 
are the same. Existence of finite coproducts in 
$\cat{M}^{\mrm{op}}$ is because of existence of finite products in $\cat{M}$;
see Proposition \ref{prop:1}(1). 

\medskip \noindent 
(3) $\cat{M}^{\mrm{op}}$ has kernels and cokernels, since 
$\opn{Ker}_{\cat{M}^{\mrm{op}}}(\opn{Op}(\phi)) = \opn{Coker}_{\cat{M}}(\phi)$ 
and vice versa. Also the symmetric condition (ii) of Definition \ref{dfn:1} 
holds. 
\end{proof}

\begin{prop}
Let $\phi : M \to N$ be a morphism in an abelian category $\cat{M}$. 
\begin{enumerate}
\item $\phi$ is a monomorphism if and only if $\opn{Ker}(\phi) = 0$. 
\item $\phi$ is an epimorphism  if and only if $\opn{Coker}(\phi) = 0$. 
\item $\phi$ is an isomorphism  if and only if it is both a monomorphism and an 
epimorphism.
\end{enumerate}
\end{prop}

\begin{exer}
Prove this proposition.
\end{exer}

Consider a diagram 
\begin{equation} \label{eqn:3665}
\bsym{S} = \bigl( \cdots  M_{-1} \xar{\phi_{-1}}
M_{0} \xar{\phi_0} M_{1} \xar{\phi_{1}} M_{2} \cdots \bigr)
\end{equation}
in an abelian category $\cat{M}$, extending finitely or infinitely to either 
side. Such a diagram is called a {\em sequence in $\cat{M}$}. 
An object $M_i$ appearing in $\bsym{S}$ is called {\em interior in $\bsym{S}$} 
if there is an object $M_{i - 1}$ appearing to the left of it, and an object 
$M_{i + 1}$ appearing to the right of it. 

\begin{dfn} \label{dfn:3635}
\index{Exact sequence}
Let $\bsym{S}$ be a sequence in the abelian category $\cat{M}$, with notation 
as in (\ref{eqn:3665}). 
\begin{enumerate}
\item Suppose $M_i$ is an interior object in $\bsym{S}$.
We say that the sequence $\bsym{S}$ is {\em exact at $M_i$} if 
$\opn{Im}(\phi_{i - 1}) =  \opn{Ker}(\phi_{i})$,
as subobjects of $M_i$. 

\item The sequence $\bsym{S}$ is said to be {\em exact} if it is exact at all 
of its interior objects.
\end{enumerate}
\end{dfn}

\begin{exa} \label{exa:3635}
A morphism $\phi : M \to N$ in an abelian category $\cat{M}$ is a monomorphism 
iff $0 \to M \xar{\phi} N$
is an exact sequence. The morphism $\phi$ is  an epimorphism iff the sequence
$M \xar{\phi} N \to 0$ is exact. 
\end{exa}

\begin{dfn} \label{dfn:4095}
A {\em short exact sequence} in an abelian category $\cat{M}$ is as exact 
sequence of the form 
\[ \bsym{S} = \bigl( 0 \to M_0 \xar{\phi_0} M_1 \xar{\phi_1} M_2 \to 0 \bigr) . 
\] 
\end{dfn}

\begin{prop} \label{prop:3637}
Let $\cat{M}$ be a $\K$-linear abelian category.
\begin{enumerate}
\item Let 
$0 \to M' \xar{\phi} M \xar{\psi} M''$
be an exact sequence in $\cat{M}$. Then for every $L \in \cat{M}$ the sequence 
\[ 0 \to \opn{Hom}_{\cat{M}}(L, M') \xar{\opn{Hom}(\opn{id}, \phi)} 
\opn{Hom}_{\cat{M}}(L, M) \xar{\opn{Hom}(\opn{id}, \psi)}
\opn{Hom}_{\cat{M}}(L, M'') \]
in $\cat{Mod} \K$ is exact. 

\item Let 
$M' \xar{\phi} M \xar{\psi} M'' \to 0$
be an exact sequence in $\cat{M}$. Then for every $N \in \cat{M}$ the sequence 
\[ 0 \to \opn{Hom}_{\cat{M}}(M'', N) \xar{\opn{Hom}(\psi, \opn{id})} 
\opn{Hom}_{\cat{M}}(M, N) \xar{\opn{Hom}(\phi, \opn{id})}
\opn{Hom}_{\cat{M}}(M', N) \]
in $\cat{Mod} \K$ is exact. 
\end{enumerate}
\end{prop}

\begin{exer} \label{exer:3635}
Prove Proposition \ref{prop:3637}. (Hint: use the definitions of kernel, 
cokernel and image.)
\end{exer}

\mysubsection{A Method for Producing Proofs in Abelian Categories} 
\label{subsec:sheaf-tricks}

A well-known difficulty in the theory of abelian categories is this: formulas 
that are easy to prove for a module category 
$\cat{M} = \cat{Mod} A$, using {\em elements}, are often 
very hard to prove in an abstract abelian category $\cat{M}$ (directly from the 
axioms). 

A neat solution to this difficulty was found by P. Freyd and B. Mitchell:

\begin{thm}[Freyd-Mitchell]
Let $\cat{M}$ be a small abelian category. Then $\cat{M}$ is equivalent to a
full abelian subcategory of $\cat{Mod} A$, for a suitable ring $A$. 
\end{thm}

\begin{rem} \label{rem:4091}
This is a deep and difficult result. See \cite{Fre}, a book that's basically 
devoted to proving this theorem. 

A modern proof can be found in
in \cite[Theorem 9.6.10]{KaSc2} -- but it too is very involved. Roughly 
speaking, they show that the abelian category $\cat{Ind}(\cat{M}^{\mrm{op}})$ 
of ind-objects of $\cat{M}$ has an injective cogenerator. This implies that the 
abelian category $\cat{Pro}(\cat{M})$ of pro-objects, that's equivalent 
to $\cat{Ind}(\cat{M}^{\mrm{op}})^{\mrm{op}}$,
has a projective generator, say $P$. Defining the ring 
$A := \opn{End}_{\cat{Pro}(\cat{M})}(P)^{\mrm{op}}$,
there is an equivalence of abelian categories
$\cat{Pro}(\cat{M}) \approx \cat{Mod} A$.
On the other hand, the Yoneda functor is a fully faithful embedding 
$\cat{M} \to \cat{Pro}(\cat{M})$.
\end{rem}

The Freyd-Mitchell Theorem implies that for 
purposes of finitary calculations in the abelian category $\cat{M}$ 
(e.g.\ checking whether a sequence is exact, 
see Definition \ref{dfn:3635}) we can assume that objects of $\cat{M}$ 
have elements. This often simplifies the work.

Since we do not give a proof of the Freyd-Mitchell Theorem in our book, 
we feel it is improper to use it. As a substitute, we provide the two 
``sheaf tricks'' below, namely Propositions \ref{prop:3635} and 
\ref{prop:3636}, with full proofs. Later in the book these sheaf tricks are 
used to give relatively easy proofs of several results on abstract abelian 
categories, most notably Theorem \ref{thm:3640} 
on the existence of the long exact 
cohomology sequence. The method of proof using these tricks is explained in  
Remark \ref{rem:3655}. It is not as slick as the method that the Freyd-Mitchell 
Theorem offers; but at least we have self-contained proofs. See Remark 
\ref{rem:3630} for some background (influence and history) on the sheaf tricks. 

First an important technical lemma. In an abelian category $\cat{M}$ we have 
finite products, and they are also coproducts (see Proposition 
\ref{prop:1}). 
 
\begin{lem} \label{lem:3630}
Let $\cat{M}$ be an abelian category. Consider a diagram 
\[ \tag{D}
\UseTips \xymatrix @C=6ex @R=6ex {
M_1 \times M_2
\ar[r]^(0.6){p_2}
\ar[d]_{p_1}
&
M_2
\ar[d]^{\phi_2}
\\
M_1
\ar[r]^{\phi_1}
&
N
} \]
in $\cat{M}$, where $p_i$ are the projections. 
Define the object 
\[ L := \opn{Ker}(\phi_1 \circ p_1 - \phi_2 \circ p_2) \sub M_1 \times M_2 , \]
with inclusion morphism 
$e : L \to M_1 \times M_2$.
Next define the morphisms 
$\psi_i := p_i \circ e : L \to M_i$. 
\begin{enumerate}
\item The diagram 
\[ \UseTips \xymatrix @C=6ex @R=6ex {
L
\ar[r]^{\psi_2}
\ar[d]_{\psi_1}
&
M_2
\ar[d]^{\phi_2}
\\
M_1
\ar[r]^{\phi_1}
&
N
} \]
is cartesian, and $L = M_1 \times_N M_2$. 

\item If $\phi_1$ is an epimorphism, then $\psi_2$ is an epimorphism.
\end{enumerate}
\end{lem}

Note that the diagram (D) is not assumed to be commutative. In case (D) does 
happen to be commutative, then 
$M_1 \times_N M_2 = M_1 \times M_2$. 

\begin{proof}
(1) The fact that $L$ (with the morphisms $\psi_i$) is the fibered product is 
immediate from the definitions of product and kernel. 

\medskip \noindent
(2) Here we follow [Mac2, Sec VIII.4]. 
Let $\rho$ be a morphism such that 
\[ \rho \circ (\phi_1 \circ p_1 - \phi_2 \circ p_2) = 0 . \]
Consider the embedding 
\[ e_1 : M_1 \to M_1 \oplus M_2 = M_1 \times M_2 . \]
Then 
\[ 0 = \rho \circ (\phi_1 \circ p_1 - \phi_2 \circ p_2) \circ e_1  
= \rho \circ \phi_1 . \]
Because $\phi_1$ is an epimorphism, it follows that $\rho = 0$. 
We conclude that 
\begin{equation} \label{eqn:3630}
\phi_1 \circ p_1 - \phi_2 \circ p_2 : M_1 \times M_2 \to N
\end{equation}
is an epimorphism. Thus (\ref{eqn:3630}) is the cokernel of $e$. 

Next let $\si : M_2 \to P$ be a morphism such that $\si \circ \psi_2 = 0$. 
Since $\psi_2 = p_2 \circ e$, we get 
$\si \circ p_2 \circ e = 0$. 
Hence $\si \circ p_2$ factors through $\opn{Coker}(e)$.
Namely there is a morphism $\si' : N \to P$ such that 
\begin{equation} \label{eqn:3626}
\si \circ p_2 = \si' \circ (\phi_1 \circ p_1 - \phi_2 \circ p_2) : 
M_1 \times M_2 \to P . 
\end{equation}
But $p_2 \circ e_1 = 0$, and therefore
\[ 0 = \si \circ p_2 \circ e_1 = 
\si' \circ (\phi_1 \circ p_1 - \phi_2 \circ p_2) \circ e_1 = 
\si' \circ \phi_1 . \]
As $\phi_1$ is an epimorphism, it follows that $\si' = 0$. 
Finally, using (\ref{eqn:3626}), we get 
\[ \si = \si \circ p_2 \circ e_2 = 
\si' \circ (\phi_1 \circ p_1 - \phi_2 \circ p_2) \circ e_2 = 
- \si' \circ \phi_2 = 0 . \]
We conclude that $\psi_2$ is an epimorphism. 
\end{proof}

Before giving the precise statements, here is a heuristic.%
\index{Sheaf Tricks (in abelian categories)|(}

\begin{rem} \label{rem:3655}
The sheaf tricks (Propositions \ref{prop:3635} and \ref{prop:3636}) work like 
this: we pretend that the objects of our $\K$-linear abelian 
category $\cat{M}$ are ``sheaves on an imaginary topological space $X$''; 
objects playing this role are denoted by letters $M, N, \ldots$. 
The objects of $\cat{M}$ are also ``open sets in the topological space $X$'';
and objects playing this role are denoted by letters $U, V, \ldots$.
Given a sheaf $M$ and an open set $U$, there is a $\K$-module 
$\Ga(U, M)$ of ``sections of $M$ over $U$''. 
These ``sections'' can be added and subtracted, being elements of a 
$\K$-module. Sometimes we require ``refinement'', i.e.\ replacing an 
``open set $U$'' by a ``covering $V \surj U$'', giving rise to an embedding 
$\K$-modules $\Ga(U, M) \inj \Ga(V, M)$. This allegory is made precise in 
Definition \ref{dfn:3621} below. 
\end{rem}

\begin{dfn}[Sheaf Metaphor] \label{dfn:3621}
Let $\cat{M}$ be a $\K$-linear abelian category. 
\begin{enumerate}
\item Objects of $\cat{M}$ are called sheaves or open sets, depending on the 
role they play in each context. 

\item For an open set $U \in \cat{M}$ and a sheaf $M \in \cat{M}$ we write 
$\Ga(U, M) := \lb \opn{Hom}_{\cat{M}}(U, M)$.
This is a $\K$-module, and we call it the {\em module of sections of 
the sheaf $M$ over the open set $U$}.

\item Given a morphism $\rho: V \to U$ of open sets in $\cat{M}$, and a sheaf 
$M \in \cat{M}$, we use this notation for the resulting $\K$-module 
homomorphism:
\[ \rho^* := \opn{Hom}_{\cat{M}}(\rho, \opn{id}_M) : 
\Ga(U, M) \to \Ga(V, M) . \]
We call $\rho^*$ the {\em pullback along $\rho$}.

\item Given a morphism $\phi : M \to N$ of sheaves in $\cat{M}$, and an open 
set $U \in \cat{M}$, we use this notation for the resulting $\K$-module 
homomorphism:
\[ \Ga(U, \phi) := \opn{Hom}_{\cat{M}}(\opn{id}_U, \phi) : 
\Ga(U, M) \to \Ga(U, N) . \]

\item If $\rho : V \to U$ is a morphism of open sets in $\cat{M}$ that is an  
epimorphism, then we call $\rho$ a {\em covering of $U$}. The homomorphism 
$\rho^*$ in this case is called the {\em restriction of $M$ from $U$ to $V$}.   
\end{enumerate}   
\end{dfn}

Note that given a covering $\rho : V \surj U$, the restriction 
homomorphism $\rho^*$ is injective (see Proposition \ref{prop:3637}(2)). By 
slight abuse of notation, when the covering $\rho$ is clear from the context, 
we will often identify the $\K$-module $\Ga(U, M)$ with its image in 
$\Ga(V, M)$ under $\rho^*$. 
Another convenient abuse of notation is writing
$\phi$ instead of $\Ga(U, \phi)$, for a morphism $\phi : M \to N$ in $\cat{M}$ 
and open set $U$. Combined, these notations give sense to such 
a commutative diagram
\begin{equation} \label{eqn:4935}
\UseTips \xymatrix @C=8ex @R=6ex {
\Ga(U, M)
\ar[r]^{\phi}
\ar[d]_{\sub}
&
\Ga(U, N)
\ar[d]^{\sub}
\\
\Ga(V, M)
\ar[r]^{\phi}
&
\Ga(V, N)
} 
\end{equation}
in $\cat{Mod} \K$. This will be used in item (3) of the next proposition. 

\begin{prop}[First Sheaf Trick] \label{prop:3635}
\index{Sheaf Tricks (in abelian categories)|)}
Let $\phi : M \to N$ be a morphism in an abelian category $\cat{M}$.
\begin{enumerate}
\item The morphism $\phi$ is zero if and only if
for every $U \in \cat{M}$ the homomorphism $\Ga(U, \phi)$ 
in $\cat{Mod} \K$ is zero.

\item  The morphism $\phi$ is a monomorphism if and only if for every $U \in 
\cat{M}$ the homomorphism $\Ga(U, \phi)$ 
in $\cat{Mod} \K$ is injective. 

\item  The morphism  $\phi$ is an epimorphism if and only if for every $U \in 
\cat{M}$ and every section $n \in \Ga(U, N)$, there exists a covering $V \surj 
U$ and a section $m \in \Ga(V, M)$, such that 
$\phi(m) = n$ in $\Ga(V, N)$.
\end{enumerate}
\end{prop}

The heuristic interpretation of item (3) of the proposition is this: 
$\phi : M \to N$ is an epimorphism in $\cat{M}$ if and only if it is ``locally 
surjective''.

\begin{proof}
(1) If $\phi = 0$ then for every $m \in \Ga(U, M)$ 
we have $\Ga(U, \phi)(m) = \phi \circ m = 0$ in $\Ga(U, N)$; thus 
$\Ga(U, \phi) = 0$. 
Conversely, if $\Ga(U, \phi) = 0$ for all $U$, then take $U := M$ and 
$m := \opn{id}_M \in \Ga(U, M)$. We obtain 
$\phi = \Ga(U, \phi)(m) = 0$. 

\medskip \noindent
(2) First assume $\phi$ is a monomorphism. Take an arbitrary object $U$ and 
a section $m \in \Ga(U, M)$. So $m : U \to M$ is a morphism in $\cat{M}$. 
 If $\phi \circ m = 0$, then, by the definition of a monomorphism, we must have 
$m = 0$. Thus $\Ga(U, \phi)$ is injective.

Conversely, assume that $\Ga(U, \phi)$ is injective for every $U$.
Take $U := \opn{Ker}(\phi)$ and let
$m : U \to M$ be the inclusion. Then 
$\Ga(U, \phi)(m) = \phi \circ m = 0$.
But then $m = 0$ and $U = 0$. Thus $\phi$ is a monomorphism. 

\medskip \noindent
(3) First assume $\phi$ is an epimorphism. Consider a section 
$n \in \Ga(U, N)$. So we have morphisms 
$\phi : M \surj N$ and $n : U \to N$. Let
$V := M \times_N U$, the fibered product. By Lemma \ref{lem:3630}(2)
the projection $\rho : V \to U$ is an epimorphism, i.e.\ a covering in our 
terminology. Now the other projection
$m : V \to M$ satisfies $\phi \circ m = n \circ \rho$. 
This means that $\phi(m) = n$ in $\Ga(V, N)$.   

Conversely, let us take $U := N$ and 
$n := \opn{id}_N \in \Ga(U, N)$. There exists an epimorphism 
$\rho : V \surj U$ and a morphism $m : V \to M$ such that 
$\phi \circ m = n \circ \rho = \rho$.
We see that $\phi \circ m$ is an epimorphism, and hence $\phi$ is an 
epimorphism.
\end{proof}

\begin{exa} \label{exa:4090}
Suppose 
\[ \bsym{S} = \bigl(  0 \to M' \xar{\phi} M \xar{\psi} M'' \to 0 \bigr) \]
is a sequence in $\cat{M}$. Here is how exactness of $\bsym{S}$ is tested using 
the first sheaf trick. 

By item (1) of the first sheaf trick, 
the condition $\psi \circ \phi = 0$ is same as 
$\Ga(U, \psi \circ \phi) = 0$ for every $U \in \cat{M}$. 

Exactness at $M'$ is the same as this condition: for every $U \in \cat{M}$ and 
every section $m' \in \Ga(U, M')$, if $\phi(m') = 0$ then $m' = 0$. 
This is by item (2). 

Exactness at $M''$ is the same as this condition: 
for every $U \in \cat{M}$ and every section 
$m'' \in \Ga(U, M'')$, there exists a covering $V \surj U$ and a section 
$m \in \Ga(V, M)$, such that 
$\psi(m) = m''$ in $\Ga(V, M'')$.

Finally, exactness at $M$ is the same as this condition (given that 
$\psi \circ \phi = 0$): for every $U \in \cat{M}$ and every section 
$m \in \Ga(U, M)$ such that $\psi(m) = 0$, there exists a covering 
$V \surj U$ and a section $m' \in \Ga(V, M')$, such that 
$\phi(m') = m$ in $\Ga(V, M'')$.
\end{exa}

\begin{prop}[Second Sheaf Trick] \label{prop:3636}
\index{Sheaf Tricks (in abelian categories)|)}
Let $M$ be an object in an abelian category $\cat{M}$. Suppose 
$\rho_1 : V_1 \to U$ and $\rho_2 : V_2 \to U$ are morphisms in $\cat{M}$, 
such that 
$(\rho_1, \rho_2) : V_1 \oplus V_2 \to U$
is an epimorphism. Let $W := V_1 \times_U V_2$, with projection morphisms 
$\si_1 : W \to V_1$ and $\si_2 : W \to V_2$. 
Then the sequence 
\[ 0 \to \Ga(U, M) \xar{(\rho_1^*, \rho_2^*)} \Ga(V_1, M) \times \Ga(V_2, M)
\xar{(\si_1^*, - \si_2^*)} \Ga(W, M) \]
in $\cat{Mod} \K$ is exact. 
\end{prop}

\begin{proof}
Note that $V_1 \times V_2 = V_1 \oplus V_2$. 
By Lemma \ref{lem:3630}(1) there is an exact sequence 
\[ 0 \to W \xar{(\si_1, -\si_2)} V_1 \oplus V_2 \xar{(\rho_1, \rho_2)} U 
\to 0 \]
in $\cat{M}$. By Proposition \ref{prop:3637}(2) we get an exact sequence
\[ 0 \to \Ga(U, M) \xar{(\rho_1^*, \rho_2^*)} \Ga(V_1 \oplus V_2, M)
\xar{(\si_1, -\si_2)^*} \Ga(W, M) \]
in $\cat{Ab}$. But 
\[ \Ga(V_1 \oplus V_2, M) \cong \Ga(V_1, M) \times \Ga(V_2, M) . \qedhere \]
\end{proof}

\begin{rem} \label{rem:3630}
As can be seen in the practical application of our sheaf tricks (e.g.\ 
in the proofs of Proposition \ref{prop:3620} and Theorem \ref{thm:3640})
these tricks reduce proofs about an abstract abelian category $\cat{M}$, to 
proofs that are just like in the concrete case of the  abelian category
$\cat{M} = \cat{Ab} X$ of sheaves of abelian groups on a topological 
space $X$. This is not as easy as the case 
$\cat{M} = \cat{Mod} A$ of modules over a ring $A$, that the Freyd-Mitchell 
Theorem permits, but -- at least for someone with experience in algebraic 
geometry -- our method is quite straightforward. 

Our sheaf tricks are inspired by \cite[tag 05PL]{SP}. There it is shown 
that the geometric allegory is genuine, except that instead of a topological 
space $X$, there is a {\em site} $X$ (in the sense of Grothendieck). The 
underlying category of the site $X$ is the abelian category $\cat{M}$ itself, 
and the coverings of $X$ are the epimorphisms in $\cat{M}$. It is proved in 
\cite{SP} that the category $\cat{M}$ embeds as a full abelian subcategory of 
$\cat{Ab} X$, the category of sheaves of abelian groups on $X$. 

Proposition \ref{prop:3635} is very similar to
\cite[Theorem VIII.3]{Mac2}, where what 
we call ``sections'' are called ``members''. But S. MacLane's method loses the 
abelian group structure (the members are not elements of abelian groups). 
An improvement of MacLane's method can be found in an unpublished note of 
G. Bergman \cite{Ber}. 
\end{rem}

\mysubsection{Additive Functors} \label{subsec:add-funcs}

\begin{dfn} \label{dfn:1031}
Let $\cat{M}$ and $\cat{N}$ be $\K$-linear categories. A functor 
$F : \cat{M} \to \cat{N}$ is called a {\em $\K$-linear functor}%
\index{Linear functor}%
\index{Linear functor!K@$\K$-linear functor}
if for every $M_0, M_1 \in \cat{M}$ the function 
\[ F : \opn{Hom}_{\cat{M}}(M_0, M_1) \to
\opn{Hom}_{\cat{N}} \bigl( F(M_0), F(M_1) \bigr) \]
is a $\K$-linear homomorphism. 

When the base ring $\K$ is implicit (cf.\ Convention \ref{conv:2490}), 
we sometimes say that $F$ is an {\em additive functor}, with the same meaning 
as $\K$-linear functor. 
\end{dfn}

Additive functors commute with finite direct sums. More precisely:

\begin{prop} \label{prop:2}
Let $F : \cat{M} \to \cat{N}$ be an additive functor between linear categories,
let $\{ M_i \}_{i \in I}$ be a finite collection of objects of $\cat{M}$, and
assume that  the direct sum $(M, \{ e_i \}_{i \in I})$ of the collection 
$\{ M_i \}_{i \in I}$ exists in $\cat{M}$. Then 
$\big( F(M), \{ F(e_i) \}_{i \in I} \big)$ is a direct sum of the collection 
$\{ F(M_i) \}_{i \in I}$ in $\cat{N}$.  
\end{prop}

\begin{exer} \label{exer:1825}
Prove Proposition \ref{prop:2}. (Hint: use Proposition \ref{prop:1}.)
\end{exer}

Note that the proposition above also talks about finite products, because of 
Proposition \ref{prop:1}.

\begin{prop} \label{prop:3600}
Suppose 
$F, G : \cat{K} \to \cat{L}$
are additive functors between linear categories, and 
$\eta : F \to G$ is a morphism of functors. 
Let $M, M' , N$ be objects of $\cat{K}$, and assume that 
$N \cong M \oplus M'$. Then the following two conditions are equivalent. 
\begin{enumerate}
\rmitem{i} $\eta_N : F(N) \to G(N)$ is an isomorphism. 

\rmitem{ii} $\eta_M : F(M) \to G(M)$ and 
$\eta_{M'} : F(M') \to G(M')$ are isomorphisms. 
\end{enumerate}
\end{prop}

\begin{exer} \label{exer:3601}
Prove Proposition \ref{prop:3600}.
\end{exer}

\begin{exa} \label{exa:4935}
Let $f : A \to B$ be a ring homomorphism. There are two additive functors 
associated to $f$~: the 
{\em restriction functor}%
\index{Restriction functor}
$\opn{Rest}_f : \cat{Mod} B \to  \cat{Mod} A$
and the {\em induction functor}\index{Induction functor}
$\opn{Ind}_f : \cat{Mod} A \to  \cat{Mod} B$.
Given a $B$-module $N$, the $A$-module $\opn{Rest}_f(N)$ has the 
same underlying $\K$-module, and $A$ acts on it through $f$.
For an $A$-module $M$, the induced $B$-module is 
$\opn{Ind}_f(M) := B \ot_A M$.
\end{exa}

\begin{prop} \label{prop:109}
Let $F : \cat{M} \to \cat{N}$ be an additive functor between linear
categories. Then\tup{:}
\begin{enumerate}
\item For every $M \in \cat{M}$ the function 
$F : \opn{End}_{\cat{M}}(M) \to \opn{End}_{\cat{N}}(F(M))$
is a ring homomorphism.

\item For every $M_0, M_1\in \cat{M}$ the function 
\[ F : \opn{Hom}_{\cat{M}}(M_0, M_1) \to 
\opn{Hom}_{\cat{N}} \bigl( F(M_0), F(M_1) \bigr) \]
is a homomorphism of left $\opn{End}_{\cat{M}}(M_1)$-modules, and of 
right $\opn{End}_{\cat{M}}(M_0)$-modules.

\item If $M$ is a zero object of $\cat{M}$, then $F(M)$ is a 
zero object of $\cat{N}$.
\end{enumerate}
\end{prop}

\begin{proof} \mbox{}

\smallskip \noindent
(1) By Definition \ref{dfn:1031} the function $F$ respects addition. By the 
definition of a functor, it respects multiplication and units. 

\medskip \noindent 
(2) Immediate from the definitions, like (1). 

\medskip \noindent 
(3) Combine part (1) with Proposition \ref{prop:1030}.
\end{proof}

\begin{dfn}
\index{Exact functor}
Let $F : \cat{M} \to \cat{N}$ be an additive functor between abelian
categories. 
\begin{enumerate}
\item $F$ is called {\em left exact} if it commutes with kernels. 
Namely for every morphism $\phi : M_0 \to M_1$ in $\cat{M}$, with 
kernel $k : K \to M_0$, the morphism $F(k) : F(K) \to F(M_0)$ is 
a kernel of $F(\phi) : F(M_0) \to F(M_1)$.

\item $F$ is called {\em right exact} if it commutes with cokernels. 
Namely for every morphism $\phi : M_0 \to M_1$ in $\cat{M}$, with 
cokernel $c : M_1 \to C$, the morphism $F(c) : F(M_1) \to F(C)$ is 
a cokernel of $F(\phi) : F(M_0) \to F(M_1)$.

\item $F$ is called {\em exact} if it is both left exact and right exact.
\end{enumerate}
\end{dfn}

Let us illustrate this. Suppose $\phi : M_0 \to M_1$
is a morphism in $\cat{M}$, with  kernel $K$ and cokernel $C$. Applying $F$ to
the sequence  
$K \xar{k} M_0 \xar{\phi} M_1 \xar{c} C$
in $\cat{M}$ 
we get the solid arrows in this diagram 
\[ \UseTips  \xymatrix @C=8ex @R=6ex {
F(K)
\ar[r]^{F(k)}
\ar@{-->}[dr]_{\psi}
& 
F(M_0)
\ar[r]^{F(\phi)}
&
F(M_1)
\ar@{-->}[d]
\ar[r]^{F(c)}
&
F(C)
\\
&
\opn{Ker}_{\cat{N}}(F(\phi))
\ar@{-->}[u]
&
\opn{Coker}_{\cat{N}}(F(\phi))
\ar@{-->}[ur]_{\chi}
} \]
in $\cat{N}$. 
Because $\cat{N}$ is abelian, we get the vertical dashed arrows: the kernel 
and cokernel of $F(\phi)$. 
The slanted dashed arrows exist and are unique because 
$F(\phi) \circ F(k) = 0$ and $F(c) \circ F(\phi) = 0$. 
Left exactness of $F$ requires $\psi$ to be an isomorphism, and 
right exactness requires $\chi$ to be an isomorphism.

Recall that a short exact sequence in $\cat{M}$ is an exact sequence of the 
form 
\begin{equation} \label{eqn:4096}
\bsym{S} = \bigl( 0 \to M_0 \xar{\phi_0} M_1 \xar{\phi_1} M_2 \to 0 \bigr) .
\end{equation}

\begin{prop}  \label{prop:1826}
Let $F : \cat{M} \to \cat{N}$ be an additive functor between abelian
categories. 
\begin{enumerate}
\item The functor $F$ is left exact if and only if 
for every short exact sequence $\bsym{S}$
in $\cat{M}$, with notation \tup{(\ref{eqn:4096})}, the sequence 
\[ 0 \to F(M_0) \xar{F(\phi_0)} F(M_1) \xar{F(\phi_0)} F(M_2)  \]
is exact in $\cat{N}$.

\item The functor $F$ is right exact if and only if 
for every short exact sequence $\bsym{S}$ in $\cat{M}$, with the notation 
with notation \tup{(\ref{eqn:4096})}, the sequence 
\[ F(M_0) \xar{F(\phi_0)} F(M_1) \xar{F(\phi_1)} F(M_2) \to 0 \]
is exact in $\cat{N}$.
\end{enumerate}
\end{prop}

\begin{exer} \label{exer:1826}
Prove Proposition \ref{prop:1826}. (Hint: $M_0 \cong \opn{Ker}(M_1 \to M_2)$ 
etc.)
\end{exer}

\begin{exa} \label{exa:4680}
Let $A$ be a commutative ring, and let $M$ be a fixed $A$-module. 
Define functors 
$F, G : \cat{Mod} A \to \cat{Mod} A$ like this:
$F(N) := M \ot_A N$ and $G(N) := \opn{Hom}_A(M, N)$.
Then $F$ is right exact and $G$ is left exact.
\end{exa}

\begin{exa} \label{exa:4681}
Let $f : A \to B$ be a ring homomorphism. 
In the notation of Example \ref{exa:4935}, the restriction functor 
$\opn{Rest}_f$ is exact, and the induction functor
$\opn{Ind}_f$ is right exact. 
\end{exa}

\begin{prop}
Let $F : \cat{M} \to \cat{N}$ be an additive functor between abelian
categories. If $F$ is an equivalence then it is exact.
\end{prop}

\begin{proof}
We will prove that $F$ respects kernels; the proof for cokernels is similar.
Take a morphism $\phi : M_0 \to M_1$ in $\cat{M}$, with kernel $K$. 
We have this diagram (solid arrows):
\[ \UseTips  \xymatrix @C=8ex @R=6ex {
M 
\ar@{-->}[d]_{\psi}
\ar@{-->}[dr]^{\theta}
\\
K
\ar[r]^{k}
& 
M_0
\ar[r]^{\phi}
&
M_1
} \]
Applying $F$ we obtain 
 this diagram (solid arrows):
\[ \UseTips  \xymatrix @C=8ex @R=6ex {
N = F(M)
\ar@{-->}[d]_{F(\psi)}
\ar@{-->}[dr]^{\bar{\theta}}
\\
F(K)
\ar[r]^{F(k)}
& 
F(M_0)
\ar[r]^{F(\phi)}
&
F(M_1)
} \]
in $\cat{N}$. 
Suppose $\bar{\theta} : N \to F(M_0)$ is a morphism in $\cat{N}$ s.t.\ 
$F(\phi) \circ \bar{\theta} = 0$.
Since $F$ is essentially surjective on objects, there is some $M \in \cat{M}$
with an isomorphism $\al : F(M) \iso N$. After replacing $N$ with $F(M)$ and 
$\bar{\theta}$ with $\bar{\theta} \circ \al$, we can assume that  
$N = F(M)$. 

Now since $F$ is fully faithful, there is a unique $\theta : M \to M_0$ s.t.\ 
$F(\theta) = \bar{\theta}$; and $\phi \circ \theta = 0$. 
So there is a unique $\psi : M \to K$ s.t.\  
$\theta = k \circ \psi$. It follows that 
$F(\psi) : F(M) \to F(K)$ is the unique morphism s.t.\ 
$\bar{\theta} = F(k) \circ F(\psi)$.
\end{proof}

Here is a result that could afford another proof of the previous proposition. 

\begin{prop} \label{prop:1827}
Let $F : \cat{M} \to \cat{N}$ be an additive functor between linear
categories. Assume $F$ is an equivalence, with quasi-inverse $G$. Then 
$G : \cat{N} \to \cat{M}$ is an additive functor.
\end{prop}

\begin{exer} \label{exer:1827}
Prove Proposition \ref{prop:1827}.
\end{exer}

\begin{dfn} \label{dfn:3620}
Consider abelian categories $\cat{M}$ and $\cat{N}$. 
Suppose we are given a sequence 
\[ \cdots F_{-1} \xar{\phi_{-1}} F_0 \xar{\phi_0} F_1 \xar{\phi_1} F_2 \cdots \]
(finite or infinite on either side), 
where each $F_i : \cat{M} \to \cat{N}$ is an additive functor, and each
$\phi_i : F_i \to F_{i + 1}$
is a morphism of functors. We say that this sequence is an {\em exact sequence 
of additive functors} if for every object $M \in \cat{M}$ the sequence 
\[ \cdots F_{-1}(M) \xar{\phi_{-1, M}} F_0(M) \xar{\phi_{0, M}} F_1(M) 
\xar{\phi_{1, M}} F_2(M) \cdots \]
in $\cat{N}$\ is exact.
\end{dfn}

\begin{prop} \label{prop:3620}
Let $\cat{M}$ and $\cat{N}$ be abelian categories, and let 
$F_0 \xar{\phi_0} F_1 \xar{\phi_1} F_2$
be a sequence of additive functors $\cat{M} \to \cat{N}$.
\begin{enumerate}
\item If  
$F_0 \xar{\phi_0} F_1 \xar{\phi_1} F_2 \to 0$
is an exact sequence of additive functors, and if the functors $F_0$ and $F_1$ 
are both right exact, then the functor $F_2$ is right exact. 

\item If  
$0 \to F_0 \xar{\phi_0} F_1 \xar{\phi_1} F_2$
is an exact sequence of additive functors, and if the functors $F_1$ and $F_2$ 
are both left exact, then the functor $F_0$ is left exact.
\end{enumerate}
\end{prop}

\begin{proof} \mbox{}

\smallskip \noindent
(1) Let 
$M' \xar{\si} M \xar{\tau} M'' \to 0$ 
be an exact sequence in $\cat{M}$. We must show that 
$F_2(M') \to F_2(M) \to F_2(M'') \to 0$ 
is an exact sequence in $\cat{N}$. Let us examine the commutative diagram 
\begin{equation} \label{eqn:3628}
\UseTips \xymatrix @C=8ex @R=6ex {
F_0(M')
\ar[r]^{\phi_{0, M'}}
\ar[d]^{F_0(\si)}
&
F_1(M')
\ar[d]^{F_1(\si)}
\ar[r]^{\phi_{1, M'}}
&
F_2(M')
\ar[r]
\ar[d]^{F_2(\si)}
&
0
\\
F_0(M)
\ar[r]^{\phi_{0, M}}
\ar[d]^{F_0(\tau)}
&
F_1(M)
\ar[r]^{\phi_{1, M}}
\ar[d]^{F_1(\tau)}
&
F_2(M)
\ar[r]
\ar[d]^{F_2(\tau)}
&
0
\\
F_0(M'')
\ar[r]^{\phi_{0, M''}}
\ar[d]
&
F_1(M'')
\ar[r]^{\phi_{1, M''}}
\ar[d]
&
F_2(M'')
\ar[r]
\ar[d]
&
0
\\
0
&
0
&
0
}
\end{equation}
in $\cat{N}$. It is known that the rows and the first two columns are exact. We 
must prove that the third column is exact. 

First let us prove that $F_2(\tau)$ is an epimorphism. This is easy: 
we know that $F_1(\tau)$ and $\phi_{1, M''}$ are epimorphisms; hence 
$\phi_{1, M''} \circ F_1(\tau) = F_2(\tau) \circ \phi_{1, M}$
is an epimorphism; and thus $F_2(\tau)$ is an epimorphism.

More challenging is the proof that the third column is exact at $F_2(M)$, 
namely that 
$\opn{Im}(F_2(\si)) \to \opn{Ker}(F_2(\tau))$
is an epimorphism. For this we shall use the first sheaf trick (Proposition 
\ref{prop:3635}) and a diagram chase. Consider a section 
$m_2 \in \Ga(U, \opn{Ker}(F_2(\tau)))$
on some ``open set'' $U$, i.e.\ for some object $U \in \cat{M}$. 
We shall prove that there exist a covering $V \surj U$ and a section 
$m'_2 \in \Ga(V, F_2(M'))$ such that 
$F_2(\si)(m'_2) = m_2$ 
in $\Ga(V, F_2(M))$. 

Because $\phi_{1, M}$ is an epimorphism, there is a covering 
$U_1 \surj U$ and a section 
$m_1 \in \Ga(U_1, F_1(M))$
such that 
$\phi_{1, M}(m_1) = m_2$ in $\Ga(U_1, F_2(M))$. 
Let
$m''_1 := F_1(\tau)(m_1)  \in \Ga(U_1, F_1(M''))$.
We have 
$\phi_{1,M''}(m''_1) = F_2(\tau)(m_2) = 0$. 
This means that 
$m''_1 \in \Ga(U_1, \opn{Ker}(\phi_{1,M''}))$.

The exactness of the third row says that for some covering 
$U_2 \surj U_1$ there is a section  
$m''_0 \in \Ga(U_2, F_0(M''))$
such that 
$\phi_{0, M''}(m''_0 ) = m''_1$ in $\Ga(U_2, F_1(M''))$.

The exactness of the first column implies that for some covering 
$U_3 \surj U_2$ there is a section  
$m_0 \in \Ga(U_3, F_0(M))$
such that 
$m''_0 = F_0(\tau)(m_0)$ in $\Ga(U_3, F_0(M''))$.
Define 
$\til{m}_1 := m_1 - \phi_{0, M}(m_0)$ in $\Ga(U_3, F_1(M))$.
Note that 
$\phi_{1, M}(\til{m}_1) = \phi_{1, M}(m_1) = m_2$ in $\Ga(U_3, F_2(M))$.
Also $F_1(\tau)(\til{m}_1) = 0$, i.e.\ 
$\til{m}_1 \in \Ga(U_3, \opn{Ker}(F_1(\tau)))$.

Due to the exactness of the second column, for some covering 
$V \surj U_3$ there is a section 
$m'_1 \in \Ga(V, F_1(M'))$
such that 
$\til{m}_1 = F_1(\si)(m'_1)$ in $\Ga(V, F_1(M))$.
Define
$m'_2 := \phi_{1, M'}(m'_1) \in \Ga(V, F_2(M'))$.
Then 
$F_2(\si)(m'_2) = \phi_{1, M}(\til{m}_1) = m_2$ in $\Ga(V, F_2(M))$.

\medskip \noindent 
(2) See next exercise.
\end{proof}

\begin{exer} \label{exer:3625}
Prove part (2) of Proposition \ref{dfn:3620}. (Hint: imitate the proof of part 
(1); but this is easier.)
\end{exer}

We end this subsection with a discussion of additive contravariant functors. 
Suppose $\cat{M}$ and $\cat{N}$ are linear categories. A contravariant functor 
$F : \cat{M} \to \cat{N}$ is said to be additive if it satisfies the condition 
in Definition \ref{dfn:1031}, with the obvious changes. 

\begin{prop} \label{prop:2491}
Let $\cat{M}$ and $\cat{N}$ be linear categories. Put on $\cat{M}^{\mrm{op}}$ 
the canonical linear structure \tup{(}see Proposition \tup{\ref{prop:2490})}.

\begin{enumerate}
\item The functor $\opn{Op} : \cat{M} \to \cat{M}^{\mrm{op}}$
is an additive contravariant functor. 

\item If $F : \cat{M} \to \cat{N}$ is an additive contravariant functor, then 
$F \circ \opn{Op} : \cat{M}^{\mrm{op}} \to \lb \cat{N}$
is an additive functor; and vice versa. 
\end{enumerate}
\end{prop}

\begin{exer} 
Prove Proposition \ref{prop:2491}.
\end{exer}

In view of Proposition \ref{prop:1826}, we can give an unambiguous definition 
of left and right exact contravariant functors. 

\begin{dfn} \label{dfn:4096}
Let $F : \cat{M} \to \cat{N}$ be an additive contravariant functor between 
abelian categories. 
\begin{enumerate}
\item $F$ a {\em left exact contravariant functor} if for 
every short exact sequence $\bsym{S}$ in $\cat{M}$, in the notation of 
(\ref{eqn:4096}), the sequence 
\[ 0 \to F(M_2) \xar{F(\phi_1)} F(M_1) \xar{F(\phi_0)} F(M_0) \]
in $\cat{N}$ is exact. 

\item $F$ is a {\em right exact contravariant functor} 
if for every short exact sequence $\bsym{S}$ in $\cat{M}$, the sequence 
\[ F(M_2) \xar{F(\phi_1)} F(M_1) \xar{F(\phi_0)} F(M_0) \to 0 \]
in $\cat{N}$ is exact.

\item $F$ is an {\em exact contravariant functor} if it sends every short 
exact sequence $\bsym{S}$ in $\cat{M}$ to a short exact sequence in $\cat{N}$. 
\end{enumerate}
\end{dfn}

\begin{prop} \label{prop:2492} 
Let $\cat{M}$ and $\cat{N}$ be abelian categories. 
Recall that $\cat{M}^{\mrm{op}}$ is also an abelian category.

\begin{enumerate}
\item The functor $\opn{Op} : \cat{M} \to \cat{M}^{\mrm{op}}$
is an exact contravariant functor. 

\item If $F : \cat{M} \to \cat{N}$ is an exact contravariant functor, then 
$F \circ \opn{Op} : \cat{M}^{\mrm{op}} \to  \cat{N}$
is an exact functor; and vice versa. Likewise for left exactness and right 
exactness. 
\end{enumerate}
\end{prop}

\begin{exer} 
Prove Proposition \ref{prop:2492}.
\end{exer}

\begin{exa} \label{exa:4682}
Let $A$ be a commutative ring, and let $M$ be a fixed $A$-module. 
Define the contravariant functor $F : \cat{Mod} A \to \cat{Mod} A$ to be 
$F(N) := \opn{Hom}_A(N, M)$.
Then $F$ is a left exact contravariant function.
\end{exa}

Sometimes $\cat{M}$ and $\cat{M}^{\mrm{op}}$ are equivalent as abelian 
categories, as the next exercise shows. For a counterexample see Remark 
\ref{rem:2510} below. 

\begin{exer} \label{exer:2515}
Let $\K$ be a field, and consider the category 
$\cat{M} := \cat{Mod}_{\mrm{f}} \K$
of finitely generated $\K$-modules (traditionally known as ``finite dimensional 
vector spaces over $\K$''). This is a $\K$-linear abelian category. 
Find a $\K$-linear equivalence
$F : \cat{M}^{\mrm{op}} \to \cat{M}$.
\end{exer}

\mysubsection{Projective Objects}
In this subsection $\cat{M}$ is an abelian category.

A {\em splitting} of an epimorphism 
$\psi : M \to M''$ in $\cat{M}$ is a morphism
$\al : M'' \to M$ s.t.\ $\psi \circ \al = \opn{id}_{M''}$.
A splitting of a monomorphism 
$\phi : M' \to M$ is  a morphism $\be : M \to M'$ s.t.\ 
$\be \circ \phi = \opn{id}_{M'}$.
A splitting of a short exact sequence 
\begin{equation} \label{eqn:1035}
0 \to M' \xar{\phi} M \xar{\psi} M'' \to 0
\end{equation}
is a splitting of the epimorphism $\psi$, or equivalently a splitting of the
mono\-morphism $\phi$. The short exact sequence is said to be {\em split} if it
has some splitting. 

\begin{exer} \label{exe:1035}
Show how to get from a splitting of $\phi$  to a splitting of $\psi$, and vice 
versa. Show how any of those gives rise to an isomorphism 
$M \cong M' \oplus M''$. 
\end{exer}

\begin{dfn}
An object $P \in \cat{M}$ is called a 
{\em projective object}%
\index{Abelian category! projective object in}
if for every morphism $\ga : P \to N$ and every {\em epimorphism}
$\psi : M \surj N$, there exists a morphism 
$\til{\ga} : P \to N$ such that 
$\psi \circ \til{\ga} = \ga$. 
\end{dfn}

This is described in the following commutative diagram in $\cat{M}$~: 
\[ \UseTips  \xymatrix @C=8ex @R=6ex {
&
P
\ar[d]^{\ga}
\ar@{-->}[dl]_{\til{\ga}}
\\
M 
\ar@{->>}[r]_{\psi}
&
N
} \]

\begin{prop} \label{prop:1670}
The following conditions are equivalent for $P \in \cat{M}$~\tup{:}
\begin{enumerate}
\rmitem{i} $P$ is projective.
\rmitem{ii} The additive functor 
$\opn{Hom}_{\cat{M}}(P, -) : \cat{M} \to \cat{Ab}$
is exact.
\rmitem{iii} Every short exact sequence \tup{(\ref{eqn:1035})} with 
$M'' = P$ is split.
\end{enumerate}
\end{prop}

\begin{proof}
Exercise.
\end{proof}

\begin{dfn}
We say $\cat{M}$ {\em has enough projectives} if every $M \in \cat{M}$ admits
an epimorphism $P \to M$ from a projective object $P$. 
\end{dfn}

\begin{exer} \label{exer:1670}
Let $A$ be a ring. 
\begin{enumerate}
\item Prove that an $A$-module $P$ is projective iff it is a 
direct summand of a free module; i.e.\ $P \oplus P' \cong Q$ for some module 
$P'$ and free module $Q$.

\item Prove that the category $\cat{Mod} A$ has enough projectives.
\end{enumerate} 
\end{exer}

\begin{exer} \label{exer:1675}
Let $\cat{M}$ be the category of finite abelian groups.
Prove that the only projective object in $\cat{M}$ is $0$. So $\cat{M}$ does 
not have enough projectives.
(Hint: use Proposition \ref{prop:1670}.)
\end{exer}

\begin{exa} \label{exa:2870}
Consider the scheme $X := \mbf{P}^1_{\K}$, the projective line over a field
$\K$, with structure sheaf $\mcal{O}_X$. 
The category $\cat{Coh} \mcal{O}_X$ of coherent $\mcal{O}_X$-modules is
abelian (it is a full abelian subcategory of $\cat{Mod} \mcal{O}_X$, cf.\ 
Example \ref{exa:1}). One can show that the only projective object of 
$\cat{Coh} \mcal{O}_X$ is $0$, but this is quite involved. 

Let us only indicate why the sheaf $\mcal{O}_X$ is not a projective
object of $\cat{Coh} \mcal{O}_X$. 
Denote by $t_0, t_1$ the homogeneous coordinates of $X$. These belong to 
$\Gamma(X, \mcal{O}_X(1))$, so each determines a homomorphism of sheaves
$t_j : \mcal{O}_X(i) \to \mcal{O}_X(i+1)$. We get a sequence 
\begin{equation} \label{eqn:5020}
0 \to \mcal{O}_X(-2) \xar{\sbmat{t_0 & t_1}} \mcal{O}_X(-1)^{\oplus 2} 
\xar{\sbmat{-t_1 \\ t_0}} \mcal{O}_X \to 0
\end{equation}
in $\cat{Coh} \mcal{O}_X$, which is known to be exact.
It is also known that 
$\Gamma(X, \mcal{O}_X) = \K$ and  
$\Gamma(X, \mcal{O}_X(-1)) = 0$. Therefore the sequence (\ref{eqn:5020}) is not 
split. 
\end{exa}

\mysubsection{Injective Objects}
In this subsection $\cat{M}$ is an abelian category.

\begin{dfn}
An object $I \in \cat{M}$ is called an 
{\em injective object}%
\index{Abelian category! injective object in}
if for every morphism $\ga : M \to I$ and every {\em monomorphism} 
$\psi : M \inj N$, there exists a morphism 
$\til{\ga} : N \to I$ such that 
$\til{\ga} \circ \psi = \ga$. 
\end{dfn}

This is depicted in the following commutative diagram in 
$\cat{M}$~:
\[ \UseTips  \xymatrix @C=8ex @R=6ex {
I
\\
M 
\ar[u]^{\ga}
\ar@{>->}[r]_{\psi}
&
N
\ar@{-->}[ul]_{\til{\ga}}
} \]

\begin{prop} \label{prop:2493}
The following conditions are equivalent for $I \in \cat{M}$\tup{:}
\begin{enumerate}
\rmitem{i} $I$ is injective.
\rmitem{ii} The additive functor 
$\opn{Hom}_{\cat{M}}(-, I) : \cat{M}^{\mrm{op}} \to \cat{Ab}$
is exact.

\rmitem{iii} Every short exact sequence \tup{(\ref{eqn:1035})} with 
$M' = I$ is split.
\end{enumerate}
\end{prop}

\begin{exer}
Prove  Proposition \ref{prop:2493}.
\end{exer}

Recall that $\opn{Op} : \cat{M} \to \cat{M}^{\mrm{op}}$ is an exact functor. 

\begin{prop} \label{prop:2494}
An object $J \in \cat{M}$ is injective if and only if the object 
$\opn{Op}(J) \in \cat{M}^{\mrm{op}}$ is projective. 
\end{prop}

\begin{exer}
Prove Proposition \ref{prop:2494}. 
\end{exer}

\begin{exa}
Let $A$ be a ring. Unlike projectives, the structure of injective objects in
$\cat{Mod} A$ is very complicated, and not much is known (except that they
exist). However if $A$ is a commutative noetherian ring then we know this:
every injective module $I$ is a direct sum of indecomposable injective 
modules; and the indecomposables are parameterized by $\opn{Spec}(A)$, the 
set of prime ideals of $A$. 
These facts are due to Matlis; see Subsection \ref{subsec:mr-inj-res} in the 
book. 
\end{exa}

\begin{dfn}
We say $\cat{M}$ {\em has enough injectives} if every $M \in \cat{M}$ admits
a monomorphism $M \to I$ to an injective object $I$. 
\end{dfn}

Here are a few results about injective objects. Recall that modules over a ring 
are always left modules by default. 

\begin{prop} \label{prop:1045}
Let $f : A \to B$ be a ring homomorphism, and let $I$ be an injective 
$A$-module. Then $J := \opn{Hom}_A(B, I)$ is an  injective $B$-module.
\end{prop}

\begin{proof}
Note that $B$ is a left $A$-module via $f$, and a right $B$-module. This makes
$J$ into a left $B$-module. In a formula: for $\phi \in J$ and $b, b' \in B$ we
have
$(b \cd \phi)(b') = \phi(b' \cd b)$. 

Now given any $N \in \cat{Mod} B$ there is an isomorphism 
\begin{equation} \label{eqn:4}
\opn{Hom}_B(N, J) =  \opn{Hom}_B(N, \opn{Hom}_A(B, I)) \cong 
\opn{Hom}_A(N, I) .
\end{equation}
 This is a natural isomorphism (of functors in $N$). So the functor 
$\opn{Hom}_B(-, J)$ is exact, and hence $J$ is injective. 
\end{proof}

\begin{thm}[Baer Criterion] \label{thm:151}
Let $A$ be a ring and $I$ an $A$-module. Assume that 
every $A$-module homomorphism $\a \to I$ from a left ideal $\a \subseteq A$
extends to a homomorphism $A \to I$. Then the module $I$ is injective.
\end{thm}

\begin{proof}
Consider an $A$-module $M$, a submodule $N \subseteq M$, and a homomorphism 
$\ga : N \to I$. We have to prove that $\ga$ extends to a homomorphism 
$M \to I$. Look at the pairs $(N', \ga')$ consisting of a submodule 
$N' \subseteq M$ that contains $N$, and a homomorphism 
$\ga' : N' \to I$ that extends $\ga$. The set of all such pairs is ordered by 
inclusion, and it satisfies the conditions of Zorn's Lemma. Therefore there 
exists a maximal pair $(N', \ga')$. We claim that $N' = M$. 

Otherwise, there is an element $m \in M$ that does not belong to $N'$. 
Define $N'' := N' + A \cd m$, so $N' \subsetneq N'' \sub M$. Let 
$\a := \{ a \in A \mid a \cd m \in N' \}$,
which is a left ideal of $A$. There is a short exact sequence 
\[ 0 \to \a \xar{\al} N' \oplus A \to N'' \to 0 \]
of $A$-modules, where 
$\al(a) := (a \cd m, -a)$. Let $\phi : \a \to I$ be the homomorphism
$\phi(a) := \ga'(a \cd m)$. 
By assumption, it extends to a homomorphism 
$\til{\phi} : A \to I$. We get a homomorphism 
$\ga' + \til{\phi}: N' \oplus A \to I$
that vanishes on the image of $\al$. 
Thus there is an induced homomorphism 
$\ga'' : N'' \to I$. This contradicts the maximality of $(N', \ga')$.
\end{proof}

\begin{lem} \label{lem:1045}
The $\Z$-module $\mbb{Q} / \Z$ is injective. 
\end{lem}

\begin{proof}
By the Baer criterion, it is enough to consider a homomorphism
$\ga : \a \to \mbb{Q} / \Z$ for an ideal $\a = n \cd \Z \subseteq \Z$. We may 
assume that $n \neq 0$. Say $\ga(n) = r + \Z$ with $r \in \mbb{Q}$. 
Then we can extend $\ga$ to $\til{\ga} : \Z \to \mbb{Q} / \Z$ with 
$\til{\ga}(1) := r / n + \Z$. 
\end{proof}

\begin{lem}  \label{lem:1046}
Let $\{ I_x \}_{x \in X}$ be a collection of injective objects of $\cat{M}$. If
the product $\prod_{x \in X} I_x$ exists in $\cat{M}$, then it is an 
injective object.
\end{lem}

\begin{proof}
Exercise.
\end{proof}

\begin{thm} \label{thm:1115}
Let $A$ be any ring. The category $\cat{Mod} A$ has enough injectives.
\end{thm}

\begin{proof} The proof is done in a few steps. 

\smallskip \noindent
Step 1. Here $A = \Z$. Take any nonzero $\Z$-module $M$ and any nonzero 
$m \in M$. Consider the cyclic submodule $M' := \Z \cd m \subseteq M$. 
There is a homomorphism $\ga' : M' \to \mbb{Q} / \Z$ s.t.\ 
$\ga'(m) \neq 0$. Indeed, if $M' \cong \Z$, then we take any 
$r \in \mbb{Q} - \Z$; and if $M' \cong \Z / (n)$ for some $n > 0$, then we 
take $r := 1 / n$. In either case, we define 
$\ga'(m) := r + \Z \in \mbb{Q} / \Z$.
Since $\mbb{Q} / \Z$ is an injective $\Z$-module, $\ga'$ extends 
to a homomorphism $\ga : M \to \mbb{Q} / \Z$.
By construction we have $\ga(m) \neq 0$. 

\medskip \noindent 
Step 2. Now $A$ is any ring, $M$ is any nonzero $A$-module, and $m \in M$ a 
nonzero element. Define the $A$-module 
$I := \opn{Hom}_{\Z}(A, \mbb{Q} / \Z)$, which, according to Lemma 
\ref{lem:1045} and Proposition \ref{prop:1045}, is an injective $A$-module. 
Let $\ga : M \to \mbb{Q} / \Z$ be a $\Z$-linear homomorphism such that 
$\ga(m) \neq 0$. Such $\ga$ exists by step 1. Let $\th : I \to \mbb{Q} / \Z$ 
be the $\Z$-linear homomorphism that sends an element $\chi \in I$ to 
$\chi(1) \in \mbb{Q} / \Z$.
The adjunction formula (\ref{eqn:4}) gives an $A$-module homomorphism 
$\psi : M \to I$ s.t.\ $\th \circ \psi = \ga$. We note that 
$(\th \circ \psi)(m) = \ga(m) \neq 0$, and hence $\psi(m) \neq 0$.  

\medskip \noindent 
Step 3. Here $A$ and $M$ are arbitrary. Let $I$ be as in step 2. For 
every nonzero  $m \in M$ there is an $A$-linear homomorphism $\psi_m : M \to I$ 
such that $\psi_m(m) \neq 0$.  For $m = 0$ let $\psi_0 : M \to I$ be an 
arbitrary homomorphism (e.g.\ $\psi_0 = 0$). 
Define the $A$-module $J := \prod_{m \in M} I$. There is a homomorphism 
$\psi := \prod_{m \in M} \psi_m : M \to J$,
and it is easy to check that $\psi$ 
is a monomorphism. By Lemma \ref{lem:1046}, $J$ is an injective $A$-module. 
\end{proof}

\begin{exer} \label{exer:1040}
At the price of getting a bigger injective module, we can make the 
construction of injective resolutions functorial. Let
$I := \opn{Hom}_{\Z}(A, \mbb{Q} / \Z)$ as above. Given an $A$-module $M$, 
consider the set 
$X(M) := \opn{Hom}_A(M, I) \cong \opn{Hom}_{\Z}(M, \mbb{Q} / \Z)$.
Let 
$J(M) := \prod_{\psi \in X(M)} I$. There is a ``tautological'' homomorphism 
$\phi_M : M \to J(M)$. Show that $\phi_M$ is a monomorphism, 
$J : M \mapsto J(M)$ is a functor, and $\phi : \opn{Id} \to J$ is a natural 
transformation. 

Is the functor $J : \cat{Mod} A \to \cat{Mod} A$ additive?
\end{exer}

\begin{exa}
Let $\cat{N}$ be the category of torsion abelian groups, and $\cat{M}$ the
category  of finite abelian groups. Then 
$\cat{N} \subseteq \cat{Ab}$
and $\cat{M} \subseteq \cat{N}$ are
full abelian subcategories. $\cat{M}$ has no projectives nor injectives except
$0$ (see Exercise \ref{exer:1675} regarding projectives). The only projective 
in $\cat{N}$ is $0$. However, it can be shown that $\cat{N}$ has enough 
injectives; see \cite[Lemma III.3.2]{Har} or \cite[Proposition 4.6]{Ye1}.
\end{exa}

\begin{prop} \label{prop:1046}
If $A$ is a left noetherian ring, then every direct sum of injective $A$-modules
is an injective module. 
\end{prop}

\begin{exer} \label{exer:1829}
Prove Proposition \ref{prop:1046}.
(Hint: use the Baer criterion.)
\end{exer}

\begin{exer} \label{exer:1045}
Here we study injectives in the category $\cat{Ab} = \cat{Mod} \Z$. 
By Lemma \ref{lem:1045}, the module $I := \Q / \Z$ is injective. For a 
(positive) prime number $p$, we denote by $\what{\Z}_p$ the ring of $p$-adic 
integers, and by $\what{\Q}_p$ its field of fractions (namely the $p$-adic 
completions of $\Z$ and $\Q$ respectively). Define the abelian group
$I_p := \what{\Q}_p / \what{\Z}_p$. 
\begin{enumerate}
\item Show that $I_p$ is an injective object of $\cat{Ab}$. 
\item Show that $I_p$ is indecomposable (i.e.\ it is not the direct sum of 
two nonzero objects). 
\item Show that $I \cong \bigoplus_p I_p$.
\item The theory (see Subsection \ref{subsec:mr-inj-res}) tells us that there 
is another indecomposable injective 
object in $\cat{Ab}$, besides the $I_p$. Try to identify it.
\end{enumerate}
\end{exer}

The abelian category $\cat{Mod} \mcal{A}$ associated to a ringed space $(X, 
\mcal{A})$ was introduced in Example \ref{exa:1}.

\begin{prop} \label{prop:2875}
Let $(X, \mcal{A})$ be a ringed space. The category $\cat{Mod} \mcal{A}$
has enough injectives. 
\end{prop}

\begin{proof}
Let $\mcal{M}$ be an $\mcal{A}$-module. Take a point $x \in X$. The stalk 
$\mcal{M}_x$ is a module over the ring $\mcal{A}_x$, and by Theorem 
\ref{thm:1115} we can
find an embedding $\phi_x : \mcal{M}_x \to I_x$ into an injective
$\mcal{A}_x$-module. 
Let $g_x : \{ x \} \to X$ be the inclusion, which we may view as a map of
ringed spaces from $( \{ x \}, \mcal{A}_x)$ to $(X, \mcal{A})$.
Define $\mcal{I}_x := {g_x}_* (I_x)$, which is an $\mcal{A}$-module (in fact it 
is a constant sheaf supported on the closed set $\ol{  \{ x \} } \subseteq X$). 
The adjunction formula gives rise to a sheaf homomorphism 
$\psi_x : \mcal{M} \to  \mcal{I}_x$. 
Since the functor $g_x^* : \cat{Mod} \mcal{A} \to \cat{Mod} \mcal{A}_x$
is exact, the adjunction formula shows that $\mcal{I}_x$ is an injective
object of $\cat{Mod} \mcal{A}$.

Finally let $\mcal{J} := \prod_{x \in X} \mcal{I}_x$. This is an injective
$\mcal{A}$-module. There is a homomorphism 
$\psi :=  \prod_{x \in X} \psi_x : \mcal{M} \to \mcal{J}$
in $\cat{Mod} \mcal{A}$. This is a monomorphism, since for every point $x$, 
letting $\mcal{J}_x$ be the stalk of the sheaf $\mcal{J}$ at $x$, the 
composition $\mcal{M}_x \xar{\psi_x} \mcal{J}_x \xar{p_x} \mcal{I}_x$
is the embedding $\phi_x : \mcal{M}_x \to I_x$.
\end{proof}

\begin{rem} \label{rem:2510}
Let $A$ be a nonzero ring, and consider the abelian category 
$\cat{M} := \cat{Mod} A$, the category of $A$-modules. A reasonable question to 
ask is this: Are the abelian categories $\cat{M}$ and $\cat{M}^{\mrm{op}}$ 
equivalent? The answer is negative. 
In fact, P. Freyd, in \cite[Exercise 5.B.3]{Fre}, shows that 
$(\cat{Mod} A)^{\mrm{op}}$ is not equivalent, as an abelian category, to 
$\cat{Mod} B$ for any ring $B$. The argument involves a delicate study of 
countable coproducts and products, and properties of Grothendieck abelian 
categories.  
\end{rem}

Here is a special case of the previous remark, that might shed some more light 
on the issue. 

\begin{exa} \label{exa:2990}
Consider the category $\cat{Mod} \Z = \cat{Ab}$ of abelian groups. 
Here is a proof that there does not exist an 
additive equivalence 
$F : \cat{Ab}^{\mrm{op}} \to \cat{Ab}$.
Suppose we had such an equivalence. Consider the object 
$P := \Z \in \cat{Ab}$, and let $I := F(P) \in \cat{Ab}$.
Because $P$ is an indecomposable projective object, and 
$F : \cat{Ab} \to \cat{Ab}$ is a contravariant equivalence,
the object $I$ has to be an indecomposable injective. 
The endomorphism rings are 
$\opn{End}_{\cat{Ab}}(I) \cong \opn{End}_{\cat{Ab}}(P)^{\mrm{op}} = 
\Z^{\mrm{op}} = \Z$.
However, the structure theorem for injective modules over commutative 
noetherian rings (Theorem \ref{thm:2026}) says that the only indecomposable 
injectives in $\cat{Ab}$ are $I = \what{\Q}_p / \what{\Z}_p$ and 
$I = \Q$; and their endomorphism rings are $\what{\Z}_p$ and $\Q$ respectively.
\end{exa}

\cleardoublepage
\mysection{Differential Graded Algebra} \label{sec:DG-algebra}

\AYcopyright

Recall that according to Convention \ref{conv:2490} 
there is a nonzero commutative base ring $\K$. By default all  
rings are $\K$-central, all linear categories are $\K$-linear,
all linear operations (such as ring homomorphisms and linear functors)
are $\K$-linear, and $\ot$ means $\ot_{\K}$. 

Throughout ``DG'' stands for ``differential graded''. 
There is some material about DG algebra in a few published references, such as
the book \cite{Mac1} and the papers \cite{Kel1}, \cite{To1}, \cite{Sao} and 
\cite{Avr}. However, for our 
purposes we need a much more detailed understanding of this theory, and this is 
what the present section provides.

\mysubsection{Graded Algebra} \label{subsec:gr-alg}

Before entering the DG world, it is good to understand the graded world. 

\begin{dfn} \label{dfn:4101}
A {\em cohomologically graded $\K$-module}%
\index{Cohomologically graded module}
is a $\K$-module $M$ equipped with a direct sum decomposition
$M = \bigoplus_{i \in \Z} M^i$
into $\K$-submodules. The $\K$-module $M^i$ is called the {\em homogeneous 
component of cohomological degree $i$} of $M$. The nonzero elements of $M^i$ 
are called {\em homogeneous elements of cohomological degree $i$}. 
\end{dfn}

From here until Section \ref{sec:alg-gra-rings} we are going to respect the 
following convention, which will simplify the discussion. 

\begin{conv} \label{conv:4125}
By ``graded $\K$-module'' we mean a {\em cohomologically graded $\K$-module}, 
as defined above.
\end{conv}

In Section \ref{sec:alg-gra-rings} we shall introduce {\em algebraically graded 
rings and modules}, and then we shall have to make a careful distinction 
between these notions. See Remark \ref{rem:2516} regarding commutativity in the 
two settings. 

Suppose $M$ and $N$ are graded $\K$-modules. For an integer $i$ let 
\[ (M \ot N)^i := \bigoplus_{j \in \Z} \, (M^j \ot N^{i - j}) . \]
Then 
\begin{equation} \label{eqn:1691}
M \ot N = \bigoplus_{i \in \Z} \, (M \ot N)^i  
\end{equation}
is a graded $\K$-module. 

A $\K$-linear homomorphism $\phi : M \to N$ is said to be {\em homogeneous of  
degree $i$} if $\phi(M^j) \sub N^{j + i}$ for all $j$. 
We denote by $\opn{Hom}_{\K}(M, N)^i$ the $\K$-module of degree $i$ 
homomorphisms $M \to N$. In other words 
\begin{equation} \label{eqn:4120}
\opn{Hom}_{\K}(M, N)^i := \prod_{j \in \Z} \, 
\opn{Hom}_{\K}(M^j, N^{j + i}) .
\end{equation}

\begin{dfn} \label{dfn:4120}
Let $M$ and $N$ be graded $\K$-modules.
\begin{enumerate}
\item The {\em module of graded $\K$-linear homomorphisms} from $M$ to $N$ is 
the graded $\K$-module
\[ \opn{Hom}_{\K}(M, N) := \bigoplus_{i \in \Z} \,
\opn{Hom}_{\K}(M, N)^i  . \]

\item A  degree $0$ homomorphism $\phi : M \to N$ is called a 
{\em strict homomorphism of graded $\K$-modules}.%
\index{Strict homomorphism!of graded modules} 
\end{enumerate}
\end{dfn}

If $M_0, M_1, M_2$ are graded $\K$-modules, and 
$\phi_k : M_{k} \to M_{k + 1}$ are $\K$-linear homomorphisms of degrees 
$i_k$, then $\phi_1 \circ \phi_0 : M_0 \to M_2$
is a $\K$-linear homomorphism of degree $i_0 + i_1$. 
The identity automorphism $\opn{id}_M : M \to M$ has degree $0$. 

\begin{dfn} \label{dfn:4100}
The {\em strict category of graded $\K$-modules}%
\index{Strict category! of graded modules} 
is the category $\dcat{G}_{\mrm{str}}(\K)$,
whose objects are the cohomologically graded $\K$-modules, and whose morphisms 
are the strict homomorphisms of cohomologically graded $\K$-mod\-ules.
\end{dfn}

It is easy to see that $\dcat{G}_{\mrm{str}}(\K)$ is a $\K$-linear abelian
category -- the kernels and cokernels are degreewise. 

\begin{rem} \label{rem:4120}
Let 
$\opn{Ungr} : \dcat{G}_{\mrm{str}}(\K) \to \dcat{M}(\K)$
be the functor that forgets the grading. It is faithful, but often not full. 
Namely the obvious homomorphism 
\[ \opn{Ungr} \bigl( \opn{Hom}_{\K}(M, N) \bigr) \to 
\opn{Hom}_{\K} \bigl( \opn{Ungr}(M), \opn{Ungr}(N) \bigr) \]
is injective but not bijective. 
See Remark \ref{rem:3795} for a discussion. 
\end{rem}

The tensor operation $(- \ot -)$ from (\ref{eqn:1691}) makes 
$\dcat{G}_{\mrm{str}}(\K)$ into a {\em monoidal $\K$-linear category},
with monoidal unit $\K$. This means that the bifunctor 
\[ (- \ot -) : \dcat{G}_{\mrm{str}}(\K) \times 
\dcat{G}_{\mrm{str}}(\K) \to \dcat{G}_{\mrm{str}}(\K) \]
satisfies the monoidal axioms (associativity up to a trifunctorial cocycle, 
etc., see \cite[Chapter XI]{Mac2}), and it is $\K$-bilinear. 
Moreover, given $M, N \in \dcat{G}_{\mrm{str}}(\K)$, let us define the 
{\em braiding isomorphism}%
\index{Braiding isomorphism}
\begin{equation} \label{eqn:4100}
\opn{br}_{M, N} : M \ot N \iso N \ot M , 
\end{equation}
\begin{equation} \label{eqn:4101}
\opn{br}_{M, N}(m \ot n) := (-1)^{i \cd j} \cd n \ot m 
\end{equation}
for homogeneous elements $m \in M^i$ and $n \in N^j$.
Because 
$\opn{br}_{N, M} \circ \opn{br}_{M, N} = \opn{id}_{M \ot N}$,
this makes $\dcat{G}_{\mrm{str}}(\K)$ into a {\em symmetric 
monoidal $\K$-linear category}. The braiding isomorphism (\ref{eqn:4101})
is often called the {\em Koszul sign rule}%
\index{Koszul sign rule}.
See Remark \ref{rem:2516} for a background discussion. 

\begin{exer} \label{exer:4100}
For an integer $l \geq 1$ let $\tau$ be a permutation of the set 
$\{ 1, \ldots, l \}$. Show that for every 
$M_1, \ldots, M_l \in \dcat{G}_{\mrm{str}}(\K)$ there is an isomorphism 
\[ \opn{br}_{\tau} : M_1 \ot \cdots \ot M_l \iso 
M_{\tau(1)} \ot \cdots \ot M_{\tau(l)} \]
in $\dcat{G}_{\mrm{str}}(\K)$, which is functorial in the sequence of objects
$\{ M_i \}_{i = 1, \ldots, l}$, functorial in the permutation $\tau$, monoidal
(for a partition $l = l_1 + l_2$), and 
$\opn{br}_{\tau} = \opn{br}_{M_1, M_2}$
when $l = 2$ and $\tau$ is the transposition.
See \cite[Theorem XI.1.1]{Mac2} or \cite{Berger} for the answer.
\end{exer}

\begin{dfn} \label{dfn:4103}
A {\em cohomologically graded central $\K$-ring}%
\index{Cohomologically graded ring}
is a central $\K$-ring $A$, equipped with a direct sum decomposition
$A =  \bigoplus_{i \in \Z} A^i$
into $\K$-sub\-modules, such that $1_A \in A^0$, and 
$A^i \cd A^j \sub A^{i + j}$.
\end{dfn}

\begin{dfn} \label{dfn:4104}
Let $A$ and $B$ be cohomologically graded central $\K$-rings. 
A {\em homomorphism of cohomologically graded central $\K$-rings} 
is a $\K$-ring homomorphism $f : A \to B$ that respects the gradings, 
namely $f(A^i) \sub B^i$.

The category of cohomologically graded central $\K$-rings is denoted by 
$\catt{GRng} \centover \K$.  
\end{dfn}

As always for ring homomorphisms, $f$ must preserve units, i.e.\ 
$f(1_A) = 1_B$. Note that $\K$ itself is a graded ring, concentrated in degree 
$0$; and it is the initial object of $\catt{GRng} \centover \K$.  

To simplify the discussion, we shall follow the next convention (until Section 
\ref{sec:alg-gra-rings}). 

\begin{conv} \label{conv:4126}
By ``graded ring'' we mean a {\em cohomologically graded central $\K$-ring}, 
as defined above. 
\end{conv}

Recall that by Convention \ref{conv:2490}, all ring homomorphisms, and that 
includes graded rings, are over $\K$.

\begin{exa} \label{exa:1705}
Let $M$ be a graded $\K$-module. Then the graded module
\[ \opn{End}_{\K}(M) := \opn{Hom}_{\K}(M, M) = 
\bigoplus_{i \in \Z} \, \opn{Hom}_{\K}(M, M)^i , \]
with the operation of composition, is a graded $\K$-ring. 
\end{exa}

The prototypical manifestation of the Koszul sign rule is the next 
definition. 

\begin{dfn} \label{dfn:4102}
Let $A = \bigoplus_{i \in \Z} A^i$ be a graded $\K$-ring.
We say that $A$ is a {\em weakly commutative graded ring}%
\index{Cohomologically graded ring! weakly commutative}
if 
$b \cd a = (-1)^{i \cd j} \cd a \cd b$
for all $a \in A^i$ and $b \in A^j$.
\end{dfn}

\begin{rem} \label{rem:4103}
Let us give a categorical explanation of Definition \ref{dfn:4102}, using the 
symmetric monoidal structure of $\dcat{G}_{\mrm{str}}(\K)$. 

The general categorical way to define a {\em ring object} in the symmetric
monoidal linear category
$\bigl( \dcat{G}_{\mrm{str}}(\K), \ot, \K, \opn{br} \big)$
is as an object $A \in \dcat{G}_{\mrm{str}}(\K)$ equipped with a 
multiplication morphism 
$\opn{m} : A \ot A \to A$
and a unit morphism $\opn{u} : \K \to A$, such that the triple 
$(A, \opn{m}, \opn{u})$ obeys the ring axioms. 
In our case it is precisely Definition \ref{dfn:4103}.

The {\em categorical commutativity condition} for the ring object 
$(A, \opn{m}, \opn{u})$  is that the diagram 
\[ \UseTips  \xymatrix @C=6ex @R=6ex {
A \ot A
\ar[dr]^{\opn{m}}
\ar[d]_{\opn{br}}
\\
A \ot A
\ar[r]^{\opn{m}}
&
A
} \]
in $\dcat{G}_{\mrm{str}}(\K)$ is commutative. This is what weak commutativity 
is, in Definition \ref{dfn:4102}.
\end{rem}

Here is a definition similar to Definition \ref{dfn:4102}.

\begin{dfn} \label{dfn:4105}
Let $A = \bigoplus_{i \in \Z} A^i$ be a graded ring.
\begin{enumerate}
\item Homogeneous elements $a \in A^i$ and $b \in A^j$ are said to {\em 
graded-commute} with each other if 
$b \cd a = (-1)^{i \cd j} \cd a \cd b$.

\item A homogeneous element $a \in A^i$ is called a {\em graded-central 
element} 
if it graded-commutes with all homogeneous elements of $A$. 

\item The {\em graded-center} of $A$ is the $\K$-submodule
$\opn{Cent}(A) \sub A$ generated by the homogeneous 
graded-central elements.
\end{enumerate}
\end{dfn}

\begin{exer} \label{exer:1706}
Let $A$ be a graded ring. Show that:
\begin{enumerate}
\item $\opn{Cent}(A)$ is a graded subring of $A$, it is weakly commutative, and 
it contains the image of the base ring $\K$. 

\item $A$ is weakly commutative iff $\opn{Cent}(A) = A$. 
\end{enumerate}
\end{exer}

Below are several sign formulas that are more subtle consequences of the Koszul 
sign rule. They can be traced -- with effort -- to the {\em biclosed monoidal 
structure} of $\dcat{G}_{\mrm{str}}(\K)$, namely to the interaction between 
the bifunctors $(- \ot -)$ and $\opn{Hom}_{\K}(-, -)$. 

Suppose that for $k = 0, 1$ we are given graded $\K$-module
homomorphisms $\phi_k : M_k \to N_k$ of degrees $i_k$. 
Then the homomorphism 
\[ \phi_0 \ot \phi_1 \in 
\opn{Hom}_{\K} \bigl( (M_0 \ot M_1, N_0 \ot N_1 \bigr)^{i_0 + i_1}  \]
acts on a tensor $m_0 \ot m_1 \in M_0 \ot M_1$, 
with $m_k \in M_k^{j_k}$, like this: 
\begin{equation} \label{eqn:4105}
(\phi_0 \ot \phi_1)(m_0 \ot m_1) := (-1)^{i_1 \cd j_0} \cd 
\phi_0(m_0) \ot \phi_1(m_1) \in N_0 \ot N_1 .
\end{equation}
The rule of thumb explaining this formula is that $\phi_1$ and $m_0$ were 
transposed. 

Suppose we are given graded $\K$-module
homomorphisms $\phi_0 : N_0 \to M_0$
and $\phi_1 : M_1 \to N_1$ of degrees $i_0$ and $i_1$. 
Then the homomorphism 
\[ \opn{Hom}(\phi_0, \phi_1) \in 
\opn{Hom}_{\K} \bigl( \opn{Hom}_{\K}(M_0, M_1), 
\opn{Hom}_{\K}(N_0, N_1) \bigr) ^{i_0 + i_1}  \]
acts on $\ga \in \opn{Hom}_{\K}(M_0, M_1)^j$
as follows: for an element $n_0 \in N^k_0$ we have
\begin{equation} \label{eqn:4106}
\opn{Hom}(\phi_0, \phi_1)(\ga)(n_0) := 
(-1)^{i_0 \cd (i_1 + j)}(\phi_1 \circ \ga \circ \phi_0)(n_0) \in 
N_1^{k + i_0 + i_1 + j} . 
\end{equation}
The sign is because $\phi_0$ jumped across $\phi_1$ and $\ga$. 

\begin{dfn} \label{dfn:4107} 
Let $A$ and $B$ be graded  rings. Then 
$A \ot B$ is a graded ring, with multiplication 
\[ (a_0 \ot b_0) \cd (a_1 \ot b_1) := 
(-1)^{i_1 \cd j_0} \cd (a_0 \cd a_1) \ot (b_0 \cd b_1)  \]
for elements $a_k \in A^{i_k}$ and  $b_k \in B^{j_k}$.
\end{dfn}

We shall require another notion of commutativity, that is not of categorical 
nature. 

\begin{dfn}  \label{dfn:3091} 
Let $A = \bigoplus_{i \in \Z} A^i$ be a graded ring.
\begin{enumerate}
\item The graded ring  $A$ is called 
{\em strongly commutative}%
\index{Cohomologically graded ring! strongly commutative}
if it is weakly commutative (Definition \ref{dfn:4102}), 
and also $a^2 = 0$ if $a \in A^i$ and $i$ is odd. 

\item The graded ring $A$ is called {\em nonpositive}%
\index{Cohomologically graded ring! nonpositive}
if $A^i = 0$ for all $i > 0$.

\item The graded ring $A$ is called a 
{\em commutative graded ring}%
\index{Cohomologically graded ring! commutative}
if it is nonpositive and strongly commutative. 
\end{enumerate}
\end{dfn}

This definition is taken from \cite{Ye9}. In \cite{YeZh3} the term 
``super-commutative'' was used instead of ``strongly commutative''. 
The name ``strongly commutative'' was suggested to us by J. Palmieri. 

\begin{exa} \label{exa:4625}
By a {\em graded set}%
\index{Graded set}
we mean a set $X$ partitioned as 
$X = \coprod_{i \in \Z} X^i$. 
The elements of $X^i$ are called variables of degree $i$. The 
{\em noncommutative graded polynomial ring}%
\index{Graded polynomial ring! noncommutative}
on the graded set $X$ is the graded ring $\K \bra{X}$, which is 
the free graded $\K$-module spanned by 
the monomials (i.e.\ words) $x_1 \cdots x_n$ in the elements of $X$, and its 
multiplication is by concatenation of monomials. 

The {\em strongly commutative graded polynomial ring}%
\index{Graded polynomial ring! commutative}
on the graded set $X$ is the graded ring 
$\K[X] := \K \bra{X} / I$,
where $I$ is the two-sided ideal of $\K \bra{X}$ generated by the elements 
$y \cd x - (-1)^{i \cd j} \cd x \cd y$ for all variables 
$x \in X^i$ and $y \in X^j$, together with the elements 
$z \cd z$ for all $z \in X^k$ and odd $k$. 
This is a strongly commutative graded ring. 
\end{exa}

\begin{exer} \label{exer:4105}
Show that if $A$ and $B$ are weakly (resp.\ strongly) commutative graded rings, 
then so is $A \ot B$. 
\end{exer}

\begin{rem} \label{rem:2516}
Weak commutativity is the obvious commutativity condition in the cohomologically
graded setting, when the Koszul sign rule%
\index{Koszul sign rule}
is imposed; see Remark 
\ref{rem:4103}. Of course there are many instances of commutative graded rings 
that do not involve the Koszul sign rule; see e.g.\ \cite{Eis} or \cite{Mats}. 
In our book we call them {\em algebraically graded rings}, and they are studied 
in Sections \ref{sec:alg-gra-rings}-\ref{sec:BDC}.

Strong commutativity has another reason. Its role is to guarantee that the 
strongly commutative graded polynomial ring $\K[X]$ from Example 
\ref{exa:4625} above will be a graded-free $\K$-module. Without this 
condition, the square of an odd variable $z$ would be a nonzero $2$-torsion 
element. 

Of course, if $2$ is invertible in $\K$ (e.g.\ if $\K$ contains $\Q$), then 
weak and strong commutativity of a graded central $\K$-ring coincide. Since 
most texts dealing with DG rings assume that $\Q \sub \K$, the subtle 
distinction we make is absent from them. 
\end{rem}

\begin{rem} \label{rem:2992}
Let $A = \bigoplus_{i \in \Z} A^i$ be a weakly commutative graded ring, and 
assume $A$ has nonzero odd elements. Then, after we forget the grading, the 
ring $A$ is no longer commutative (except in special cases, like in 
characteristic $2$). 

This phenomenon will reappear often in the study of DG algebra.
\end{rem}

\begin{rem} \label{rem:2991}
The origin of the Koszul sign rule%
\index{Koszul sign rule},
predating the work of J.L. Koszul, 
appears to be in the formula for the differential $\d$ of the tensor product of 
two complexes, that occurs in classical homological algebra (and is repeated 
here in Definition \ref{dfn:2992}). Without the sign we won't have 
$\d \circ \d = 0$. 

In his thesis \cite{Kosz2}, Koszul talks about the sign rule explicitly,
when discussing the Chevalley-Eilenberg complex. 
See \cite{Berger} for a detailed study of the Koszul sign rule and its 
group-theoretic interpretation.  

It should be noted that there are some sign inconsistencies in the book 
\cite{RD}. The Koszul sign rule could have prevented them. 

In our book we made a strenuous effort to have correct (namely consistent) 
signs. This did not always produce satisfactory outcomes -- for instance, in 
the context of the triangulated structure of the opposite homotopy category, as 
explained in Remark \ref{rem:2306}. 
\end{rem}

\begin{dfn} \label{dfn:4127}
Let $A$ be a graded ring. A 
{\em graded left $A$-module}%
\index{Graded module}
is a left $A$-module $M$, equipped with a $\K$-module decomposition 
$M =  \bigoplus_{i \in \Z} M^i$, such that 
$A^i \cd M^j \sub M^{i + j}$ for all $i, j$. 
\end{dfn}

We can also talk about graded right $A$-modules, and graded bimodules. 
But our default option (see Convention \ref{conv:2490}) is that modules are 
left modules.

\begin{exer} \label{exer:4110}
Let $M$ be a graded $\K$-module, $A$ a graded central $\K$-ring, and
$f : A \to \opn{End}_{\K}(M)$ a homomorphism in $\catt{GRng} \centover \K$.
\begin{enumerate}
\item Show that $M$ becomes a graded $A$-module, with action  
$a \cd m := f(a)(m)$.

\item Show that every graded $A$-module structure on $M$, that is consistent 
with the given graded $\K$-module structure, arises this way. 
\end{enumerate}  
\end{exer}

\begin{lem} \label{lem:1690}
Let $A$ be a graded ring, let $M$ be a graded right $A$-module, 
and let $N$ be a graded left $A$-module. Then the $\K$-module 
$M \ot_A N$ has a direct sum decomposition
\[ M \ot_A N = \bigoplus_{i \in \Z} \, (M \ot_A N)^i, \]
where $(M \ot_A N)^i$ is the $\K$-linear span of the tensors $m \ot n$, for 
all $j \in \Z$, $m \in M^j$ and  $n \in N^{i - j}$. 
\end{lem}

\begin{proof}
There is a canonical surjection of $\K$-modules 
$M \ot N \to M \ot_A N$.
Its kernel is the $\K$-submodule 
$L \subseteq M \ot N$ generated by the elements 
$(m \cd a) \ot n - m \ot (a \cd n)$,
for $m \in M^j$, $n \in N^k$ and $a \in A^l$.
So $L$ is a graded submodule of  $M \ot N$, and therefore so is the 
quotient. Finally, by formula (\ref{eqn:1691}) the $i$-th homogeneous component 
of $M \ot_A N$ is precisely  $(M \ot_A N)^i$. 
\end{proof}

\begin{dfn} \label{dfn:1690}
Let $A$ be a graded ring, and let $M, N$ be graded $A$-modules.
For each $i \in \Z$ define $\opn{Hom}_{A}(M, N)^i$ to be the subset of 
$\opn{Hom}_{\K}(M, N)^i$ consisting of the homomorphisms $\phi : M \to N$ such 
that
$\phi(a \cd m) = (-1)^{i \cd k} \cd a \cd \phi(m)$
for all $a \in A^k$. Next define the graded $\K$-module
\[ \opn{Hom}_{A}(M, N) := \bigoplus_{i \in \Z} \, \opn{Hom}_{A}(M, N)^i . \]
\end{dfn}

Suppose $\cat{C}$ is a $\K$-linear category (Definition \ref{dfn:1083}).
Since the composition of morphisms is $\K$-bilinear, for every  
triple of objects $M_0, M_1, M_2 \in \cat{C}$, composition can be expressed as
a $\K$-linear homomorphism 
\[ \opn{Hom}_{\cat{C}}(M_1, M_2) \ot 
\opn{Hom}_{\cat{C}}(M_0, M_1) \xar{} 
\opn{Hom}_{\cat{C}}(M_0, M_2) , \]
$\phi_1 \ot \phi_0 \mapsto \phi_1 \circ \phi_0$.
We refer to it as the composition homomorphism. It will be used in the 
following definition. 

\begin{dfn} \label{dfn:1692}
A {\em graded $\K$-linear category}%
\index{Graded category}
is a $\K$-linear category $\cat{C}$, 
endowed with a grading on each of the $\K$-modules 
$\opn{Hom}_{\cat{C}}(M_0, M_1)$. The conditions are these:
\begin{enumerate}
\rmitem{a} For every object $M$, the identity automorphism $\opn{id}_M$ has 
degree $0$. 

\rmitem{b} For every triple of objects $M_0, M_1, M_2 \in \cat{C}$, the 
composition homomorphism 
\[ \opn{Hom}_{\cat{C}}(M_1, M_2) \ot \opn{Hom}_{\cat{C}}(M_0, M_1)
\xar{} \opn{Hom}_{\cat{C}}(M_0, M_2) \]
is a strict homomorphism of graded $\K$-modules. 
\end{enumerate}
\end{dfn}

In item (b) we use the graded module structure on a tensor product from 
equation (\ref{eqn:1691}).
A morphism $\phi \in \opn{Hom}_{\cat{C}}(M_0, M_1)^i$ is called a 
morphism of degree $i$. 
 
\begin{dfn} \label{dfn:1693} 
Let $\cat{C}$ be a graded $\K$-linear category.
The {\em strict subcategory of $\cat{C}$}%
\index{Strict subcategory! of a graded category}
is the subcategory $\opn{Str}(\cat{C})$ on all objects of $\cat{C}$, whose 
morphisms are the degree $0$ morphisms of $\cat{C}$.
\end{dfn}

\begin{exa} \label{exa:1696}
Let $A$ be a graded ring. Define $\cat{GMod} A$ to be 
the category whose objects are the graded $A$-modules. For 
$M, N \in \cat{GMod} A$, the set of morphisms is the graded $\K$-module 
$\opn{Hom}_{\cat{GMod} A}(M, N) := \opn{Hom}_{A}(M, N)$ 
from Definition \ref{dfn:1690}. Then $\cat{GMod} A$ is a graded $\K$-linear 
category. The morphisms in the subcategory
$\cat{GMod}_{\mrm{str}} A := \opn{Str}(\cat{GMod} A)$ are the 
strict homomorphisms of graded $A$-modules.
We often write
$\dcat{G}(A) := \cat{GMod} A$ and 
$\dcat{G}_{\mrm{str}}(A) := \cat{GMod}_{\mrm{str}} A$.  
In the special case $A = \K$, the category $\dcat{G}_{\mrm{str}}(\K)$
already appeared in Definition \ref{dfn:4100}. 
\end{exa}

\begin{rem} \label{rem:1190}
The name ``strict homomorphism of graded modules'', and the corresponding 
notations $\cat{GMod}_{\mrm{str}} A = \dcat{G}_{\mrm{str}}(A)$, 
are new. We introduced them to distinguish the 
abelian category $\cat{GMod}_{\mrm{str}} A$ from the graded category 
$\cat{GMod} A$ that contains it.
See Definitions \ref{dfn:1063} and
\ref{dfn:1212} for the DG versions of these notions.
\end{rem}

\begin{dfn} \label{dfn:1694}
Let $\cat{C}$ and $\cat{D}$ be graded $\K$-linear categories. 
A functor $F : \cat{C} \to \cat{D}$
is called a 
{\em graded $\K$-linear functor}%
\index{Graded functor}
if it satisfies this condition: 
\begin{itemize}
\item[$\vartriangleright$] For every pair of objects $M_0, M_1 \in \cat{C}$, 
the function 
\[ F : \opn{Hom}_{\cat{C}}(M_0, M_1) \to  
\opn{Hom}_{\cat{D}} \bigl( F(M_0), F(M_1) \bigr) \] 
is a strict homomorphism of graded $\K$-modules. 
\end{itemize}
\end{dfn}

\begin{conv} \label{conv:4625}
To simplify the terminology, we shall often use the expressions ``graded 
category'' and ``graded functor'' as abbreviations for ``graded $\K$-linear
category'' and ``graded $\K$-linear functor'', respectively. 
\end{conv}
 
\begin{exa} \label{exa:1694}
Let $A$ be a graded ring. We can view $A$ as a category 
$\cat{A}$ with a single object, and it is a graded category.
If $f : A \to B$ is a homomorphism of graded rings, then 
passing to single-object categories we get a 
graded functor $F : \cat{A} \to \cat{B}$. 
\end{exa}

Recall that ``morphism of functors'' is synonymous with ``natural 
transformation''. 

\begin{dfn} \label{dfn:1695}
Let  $F, G : \cat{C} \to \cat{D}$
be graded functors between graded categories, and let $i \in \Z$. 
A {\em degree $i$ morphism of graded functors}%
\index{Graded functor! morphism of {\indash}s}
$\eta : F \to G$ is a collection 
$\eta = \{ \eta_M \}_{M \in \cat{C}}$ of morphisms 
$\eta_M \in \opn{Hom}_{\cat{D}} \bigl( F(M), G(M) \bigr)^i$,
such that for every morphism 
$\phi \in \opn{Hom}_{\cat{C}} (M_0, M_1)^j$ there is equality 
\[ G(\phi) \circ \eta_{M_0} = (-1)^{i \cd j} \cd  \eta_{M_1} \circ F(\phi) \]
inside 
$\opn{Hom}_{\cat{D}} \bigl( F(M_0), G(M_1) \bigr)^{i + j}$.
\end{dfn}

If $i$ is odd in the definition above, then after we forget the grading, 
$\eta : F \to G$ is usually no longer a morphism of functors; this is another 
instance of the phenomenon mentioned in Remark \ref{rem:2992}.

\begin{dfn}  \label{dfn:1700}
Let $\cat{M}$ be an abelian category. 
A {\em graded object in $\cat{M}$}%
\index{Abelian category! graded object in}
is a collection 
$\{ M^i \}_{i \in \Z}$ of objects $M^i \in \cat{M}$.
\end{dfn}

Because we did not assume that $\cat{M}$ has countable direct sums, the graded 
objects are ``external'' to $\cat{M}$; cf.\ Exercise \ref{exer:4112}.

Suppose $M = \{ M^i \}_{i \in \Z}$ and $N = \{ N^i \}_{i \in \Z}$
are graded objects in $\cat{M}$. For an integer $i$ we define the $\K$-module
\begin{equation} \label{eqn:2495}
\opn{Hom}_{\cat{M}}(M, N)^i := \prod_{j \in \Z} 
\opn{Hom}_{\cat{M}}(M^j, N^{j + i}) .
\end{equation}
We get a graded $\K$-module 
\begin{equation} \label{eqn:1700}
\opn{Hom}_{\cat{M}}(M, N) := \bigoplus_{i \in \Z} \,
\opn{Hom}_{\cat{M}}(M, N)^i . 
\end{equation}

\begin{dfn}  \label{dfn:1701}
Let $\cat{M}$ be an abelian category. 
The {\em category of graded objects in $\cat{M}$} is the graded linear 
category $\dcat{G}(\cat{M})$, whose objects are the graded objects 
in $\cat{M}$, and the morphism sets are the graded modules 
\[ \opn{Hom}_{\dcat{G}(\cat{M})}(M, N) := \opn{Hom}_{\cat{M}}(M, N) \]
from equation (\ref{eqn:1700}). 
The composition operation is the obvious one. 
\end{dfn}

\begin{exer}  \label{exer:4112}
Suppose $\cat{M} = \dcat{M}(A)$, the category of modules over a 
central $\K$-ring $A$. For every 
$M = \{ M^i \}_{i \in \Z} \in \dcat{G}(\cat{M})$ let 
$F(M) := \bigoplus_{i \in \Z} M^i$.
Then $F(M)$ is a graded $A$-module, as discussed earlier, so $F(M)$ is an 
object of the category $\dcat{G}(A)$ from Example \ref{exa:1696}.
Prove that 
$F : \dcat{G}(\cat{M}) \to \dcat{G}(A)$
is an isomorphism of graded categories. 
\end{exer}

In the next definition we combine graded rings and linear categories, to 
concoct a new hybrid.

\begin{dfn} \label{dfn:1702}
Let $\cat{M}$ be an abelian category, and let $A$ be a graded ring. A 
{\em graded $A$-module in  $\cat{M}$}%
\index{Abelian category! graded $A$-module in}
is an object 
$M \in \dcat{G}(\cat{M})$, together with a graded ring homomorphism 
$f : A \to \opn{End}_{\cat{M}}(M)$. 
The set of graded $A$-modules in $\cat{M}$ is denoted by 
$\dcat{G}(A, \cat{M})$. 
\end{dfn}

What the definition says is that an element 
$a \in A^i$ gives rise to a degree $i$ endomorphism 
$f(a)$ of the graded object $M  = \{ M^j \}_{j \in \Z}$. In turn, this means 
that for every $j$, 
$f(a) : M^j \to M^{j + i}$ is a morphism in $\cat{M}$.  
The operation $f$ satisfies $f(1_A) = \opn{id}_M$ and 
$f(a_1 \cd a_2) = f(a_1) \circ f(a_2)$. 

\begin{exa}  \label{exa:1701}
If $A = \K$, then $\dcat{G}(A, \cat{M}) = \dcat{G}(\cat{M})$; 
and if $\cat{M} = \cat{Mod} \K$, then 
$\dcat{G}(A, \cat{M}) = \dcat{G}(A)$.
\end{exa}

The next definition is a variant of Definition \ref{dfn:1690}.

\begin{dfn} \label{dfn:1703}
Let $\cat{M}$ be an abelian category, and let $A$ be a 
graded ring. For $M, N \in \dcat{G}(A, \cat{M})$ and $i \in \Z$ 
we define 
$\opn{Hom}_{A, \cat{M}}(M, N)^i$ to be the subset of 
$\opn{Hom}_{\cat{M}}(M, N)^i$ consisting of the  morphisms 
$\phi : M \to N$ such that
\[ \phi \circ f_M(a) = (-1)^{i \cd k} \cd f_N(a) \circ \phi \]
for all $a \in A^k$. Next let 
\[ \opn{Hom}_{A, \cat{M}}(M, N) := 
\bigoplus_{i \in \Z} \, \opn{Hom}_{A, \cat{M}}(M, N)^i . \]
This is a graded $\K$-module.
\end{dfn}

\begin{dfn} \label{dfn:1704}
Let $\cat{M}$ be an abelian category, and let $A$ be a graded   
ring. The {\em category of graded $A$-modules in $\cat{M}$}
is the graded category $\dcat{G}(A, \cat{M})$
whose objects are the graded $A$-modules in 
$\cat{M}$, and the morphism graded $\K$-modules are 
\[ \opn{Hom}_{\dcat{G}(A, \cat{M})}(M_0, M_1) :=
\opn{Hom}_{A, \cat{M}}(M_0, M_1)  \]
from Definition \ref{dfn:1703}. 
The compositions are those of $\dcat{G}(\cat{M})$. 
\end{dfn}

In the special case that $A$ is a ring (i.e.\ a graded ring 
concentrated in degree $0$), this is \cite[Definition 8.5.1]{KaSc2}. 

Notice that forgetting the action of $A$ is a faithful graded 
functor $\dcat{G}(A, \cat{M}) \to \dcat{G}(\cat{M})$.
As in every graded category, there is the subcategory
\begin{equation} \label{eqn:2360}
\dcat{G}_{\mrm{str}}(A, \cat{M}) := \opn{Str}(\dcat{G}(A, \cat{M}))
\sub \dcat{G}(A, \cat{M})
\end{equation}
of strict (i.e.\ degree $0$) morphisms. 

\begin{exer} \label{exer:1705}
Show that $\dcat{G}_{\mrm{str}}(A, \cat{M})$ is an abelian category. 
\end{exer}

\begin{rem} \label{rem:2517}
The reader may have noticed that we can talk about the graded category 
$\dcat{G}(\cat{M})$ for every linear category $\cat{M}$, regardless if it 
is abelian or not. We chose to restrict attention to the abelian case for a 
pedagogical reason: this will hopefully reduce confusion between the many sorts 
of graded (and later DG) categories that occur in our discussion. 
\end{rem}

\mysubsection{DG \texorpdfstring{$\K$}{K}-Modules}

Recall that Conventions \ref{conv:2490} and \ref{conv:4125} are in force. 

\begin{dfn} \label{dfn:1070}
A {\em DG $\K$-module}%
\index{Differential graded module}
 is a graded $\K$-module 
$M = \bigoplus_{i \in \Z} M^i$, together with a $\K$-linear homomorphism 
$\d_M : M \to M$ of degree $1$, called the {\em differential}, satisfying 
$\d_M \circ \d_M = 0$.
\end{dfn}

When there is no danger of confusion, we may write $\d$ instead of $\d_M$. 

\begin{dfn} \label{dfn:1084}
Let $M$ and $N$ be DG $\K$-modules. A 
{\em strict homomorphism of DG $\K$-modules}%
\index{Strict homomorphism! of DG modules}
is a $\K$-linear homomorphism $\phi : M \to N$ 
of degree $0$ that commutes with the differentials. 
The resulting category is denoted by $\cat{DGMod}_{\mrm{str}} \K$, or by the 
abbreviated notation $\dcat{C}_{\mrm{str}}(\K)$. 
\end{dfn}

It is easy to see that $\dcat{C}_{\mrm{str}}(\K)$ is a $\K$-linear abelian 
category. There is a forgetful functor 
$\dcat{C}_{\mrm{str}}(\K) \to \dcat{G}_{\mrm{str}}(\K)$, 
$(M, \d_M) \mapsto M$.

\begin{dfn} \label{dfn:2992}
Suppose $M$ and $N$ are DG $\K$-modules. 
\begin{enumerate}
\item 
\index{Differential graded module! tensor product of {\indash}s}
The graded $\K$-module structure on the tensor product $M \ot N$
was given in equation (\ref{eqn:1691}). 
We put on it the differential 
\[ \d(m \ot n) := \d_M(m) \ot n + (-1)^i \cd m \ot \d_N(n) \]
for $m \in M^i$ and $n \in N^j$. In this way $M \ot N$ becomes a DG 
$\K$-module. We sometimes write $\d_{M \ot N}$ for this differential.

\item 
\index{Differential graded module! Hom between {\indash}s}
The graded module $\opn{Hom}_{\K}(M, N)$
was introduced in Definition \ref{dfn:4120}. We give it this 
differential: 
\[ \d(\phi) := \d_N \circ \phi - (-1)^i \cd \phi \circ \d_M \]
for $\phi \in \opn{Hom}_{\K}(M, N)^i$. 
Thus $\opn{Hom}_{\K}(M, N)$ becomes a DG $\K$-module. 
We sometimes denote this differential by 
$\d_{\opn{Hom}_{\K}(M, N)}$.
\end{enumerate}
\end{dfn}

\begin{dfn} \label{dfn:2993}
Let $M$ be a DG $\K$-module, and let $i$ be an integer. 
\begin{enumerate}
\item The module of degree $i$ {\em cocycles}%
\index{Cocycles! module of}%
\index{1-Zi(M)@$\opn{Z}^i(M)$}
of $M$ is
\[ \opn{Z}^i(M) := \opn{Ker} \bigl( M^i \xar{\d_M} M^{i + 1} \bigr) . \]

\item The module of degree $i$ {\em coboundaries}%
\index{Coboundaries! module of}%
\index{1-Bi(M)@$\opn{B}^i(M)$}
of $M$ is 
\[ \opn{B}^i(M) := \opn{Im} \bigl( M^{i -1} \xar{\d_M} M^{i} \bigr) . \]

\item The module of degree $i$ {\em decocycles}%
\index{Decocycles! module of}%
\index{1-Yi(M)@$\opn{Y}^i(M)$}
of $M$ is 
\[ \opn{Y}^i(M) := \opn{Coker} \bigl( M^{i -1} \xar{\d_M} M^{i} \bigr) . \]

\item The $i$-th {\em cohomology}%
\index{Cohomology! module}%
\index{1-Hi(M)@$\opn{H}^i(M)$}
module of $M$ is 
\[ \opn{H}^i(M) := \opn{Z}^i(M) \, / \, \opn{B}^i(M)
\cong \opn{Coker} \bigl( M^{i - 1} \xar{\d_M} \opn{Z}^i(M) \bigr) . \]
\end{enumerate}
\end{dfn}

The definition above is standard, with the exception of item (3), which is 
a new contribution (both the notation $\opn{Y}^i(M)$ and the name 
``decocycles''). 
The fact that $\d_M \circ \d_M = 0$ implies that 
$\opn{B}^i(M) \sub \opn{Z}^i(M) \sub M^i$, so item (4) makes sense. 
On the other hand, since 
$\opn{Y}^i(M) = M^i / \opn{B}^i(M)$,
there is a canonical isomorphism 
\begin{equation} \label{eqn:2993}
\opn{H}^i(M) \cong \opn{Ker} \bigl( \opn{Y}^i(M) \xar{\d_M} M^{i + 1} \bigr) .
\end{equation}
The modules defined above are functorial in $M$; to be precise, these are 
$\K$-linear functors  
\begin{equation} \label{eqn:4127}
\opn{Z}^i, \opn{B}^i, \opn{Y}^i, \opn{H}^i  : 
\dcat{C}_{\mrm{str}}(\K) \to \dcat{M}(\K) . 
\end{equation}

We can rephrase Definition \ref{dfn:1084} using the notion of cocycles:  for 
DG $\K$-modules $M$ and $N$ there is equality
\begin{equation} \label{eqn:2520}
\opn{Hom}_{\dcat{C}_{\mrm{str}}(\K)}(M, N) = 
\opn{Z}^0 \bigl( \opn{Hom}_{\K}(M, N) \bigr)
\end{equation}
of submodules of $\opn{Hom}_{\K}(M, N)$.

\mysubsection{DG Rings and Modules} \label{subsec:DG-rng-mod}

Recall that Conventions \ref{conv:2490}, \ref{conv:4125} and \ref{conv:4126} 
are in place. 

\begin{dfn} \label{dfn:1071}
A {\em DG central $\K$-ring}
\index{Differential graded ring}
is a graded central $\K$-ring 
$A =  \bigoplus_{i \in \Z} A^i$ (see Definition \ref{dfn:4103}), together with 
a $\K$-linear homomorphism $\d_A : A \to A$ of degree $1$, called the {\em 
differential}, satisfying  the equation
$\d_A \circ \d_A = 0$, and the graded Leibniz rule
\index{Graded Leibniz rule}
\[ \d_A(a \cd b) =  \d_A(a) \cd b + (-1)^i \cd a \cd \d_A(b) \]
for all $a \in A^i$ and $b \in A^j$. 
\end{dfn}

We sometimes write $\d$ instead of $\d_A$.

\begin{dfn} \label{dfn:1074}
Let $A$ and $B$ be DG central $\K$-rings. A 
{\em homomorphism of DG central $\K$-rings} 
\index{Differential graded ring! homomorphism of {\indash}s}
$f : A \to B$ is a graded central 
$\K$-ring homomorphism that 
commutes with the differentials of $A$ and $B$.
The resulting category is denoted by $\catt{DGRng} \centover \K$. 
\end{dfn}

Central $\K$-rings are viewed as DG central $\K$-rings concentrated in 
degree $0$ (and with trivial differentials). Thus the category $\catt{Rng} 
\centover \K$ is a full subcategory of $\catt{DGRng} \centover \K$. 

\begin{conv} \label{conv:4626}
To simplify terminology, we usually write ``DG ring'' instead of 
``DG central $\K$-ring''. 
\end{conv}

\begin{dfn} \label{dfn:3090}
\index{Differential graded ring! weakly commutative}
\index{Differential graded ring! strongly commutative}
\index{Differential graded ring! nonpositive}
\index{Differential graded ring! commutative}
A DG ring $A$ is called {\em weakly commutative}, {\em strongly 
commutative}, {\em nonpositive} or {\em commutative} if it is so, respectively, 
as a graded ring (after forgetting the differential), in the sense of 
Definition \ref{dfn:3091}. 

The corresponding full subcategories of $\catt{DGRng} \centover \K$ are 
$\catt{DGRng}_{\mrm{wc}} / \K$,  $\catt{DGRng}_{\mrm{sc}} / \K$,
$\catt{DGRng}^{\leq 0} \centover \K$ and 
$\catt{DGRng}^{\leq 0}_{\mrm{sc}} / \K$.
\end{dfn}
 
When $\K = \Z$ we can write $\catt{DGRng} := \catt{DGRng} \centover \Z$, etc.

Here are few examples of DG rings. First a silly example.

\begin{exa} \label{exa:1095}
Let $A$ be a graded ring. Then $A$, with the zero differential, 
is a DG ring.
\end{exa}

\begin{exa} \label{exa:1096}
Let $X$ be a differentiable (i.e.\ of type $\mrm{C}^{\infty}$) manifold over 
$\K := \R$. The de Rham complex $A$ of $X$ is a DG $\R$-ring, with the wedge 
product and the exterior differential. See \cite[Section 2.9.7]{KaSc1} for 
details. This is a strongly commutative DG ring, in the sense of Definition 
\ref{dfn:3090}. However, it is almost never nonpositive. 
\end{exa}

The next example  is the algebraic analogue of the previous one.
 
\begin{exa} \label{exa:1097}
Let $C$ be a commutative $\K$-ring. Then the algebraic de Rham complex
$A := \Omega_{C / \K} = \bigoplus_{p \geq 0} \Omega^p_{C / \K}$ is a
DG $\K$-ring. It is a strongly commutative DG ring, and it is 
almost never nonpositive. 
See \cite[Exercise 16.15]{Eis} or \cite[Section 25]{Mats} for details. 
\end{exa}

\begin{exa} \label{exa:1098}
Let $M$ be a DG $\K$-module. Consider the DG $\K$-module \lb 
$\opn{End}_{\K}(M) := \opn{Hom}_{\K}(M, M)$
from Definition \ref{dfn:2992}(2). 
Composition of endomorphisms is an associative 
multiplication on $\opn{End}_{\K}(M)$ that respects the grading, and the graded 
Leibniz rule holds. We see that 
$\opn{End}_{\K}(M)$ is a DG $\K$-ring.
It is usually neither weakly commutative nor nonpositive. 
\end{exa}

\begin{exa} \label{exa:1099}
Let $C$ be a commutative ring and let $c \in C$ be an element. 
The {\em Koszul complex}
\index{Koszul complex}
 of $c$ is the DG $C$-module $\opn{K}(C; c)$ defined as 
follows. In degree $0$ we let $\opn{K}^0(C; c) := C$. 
In degree $-1$, $\opn{K}^{-1}(C; c)$ is a free $C$-module of rank $1$, with 
basis element $x$. All other homogeneous components are trivial. The 
differential $\d$ is determined by what it does to the basis element
$x \in \opn{K}^{-1}(C; c)$, and we let 
$\d(x) := c \in \opn{K}^{0}(C; c)$. 

We want to make $\opn{K}(C; c)$ into a strongly commutative DG ring
(in the sense of Definition \ref{dfn:3090}). Since $x$ is an odd element,
we must define the graded ring structure by 
$\opn{K}(C; c) := C \ot \K[x]$,
where $\K[x]$ is the strongly commutative graded polynomial ring from Example 
\ref{exa:4625}, on the graded set 
$X = X^{-1} := \{ x \}$; i.e.\ $x^2 = 0$. 
\end{exa}

\begin{exa} \label{exa:1100}
Let $A$ and $B$ be DG rings. The graded ring 
$A \ot B$ from Definition \ref{dfn:4107}, with the differential
from Definition \ref{dfn:2992}(1),
is a DG ring.
\end{exa}

\begin{exa} \label{exa:1101}
Let $C$ be a commutative ring and let 
$\bsym{c} = (c_1, \ldots, c_n)$ be a sequence of elements in $C$. By combining 
Examples \ref{exa:1099} and \ref{exa:1100} we obtain the Koszul complex 
\index{Koszul complex}
\[ \opn{K}(C; \bsym{c}) := \opn{K}(C; c_1) \ot_{C} \cdots \ot_{C}
\opn{K}(C; c_n) . \]
This is a strongly commutative DG $C$-ring.
\end{exa}

\begin{rem} \label{rem:4685}
In the book \cite[Section 17, Exercises]{Eis} there are a few variants of the 
Koszul complex. Note that in most of the classical literature (including 
\cite{Eis} and \cite{Mats}), the multiplicative structure of $\opn{K}(C; 
\bsym{c})$ has been ignored. 
\end{rem}

\begin{dfn} \label{dfn:1080}
Let $A$ be a DG ring. The {\em opposite DG ring}
\index{Differential graded ring! opposite}
$A^{\mrm{op}}$
is the same DG $\K$-module as $A$, but the multiplication $\cdot^{\mrm{op}}$ of 
$A^{\mrm{op}}$ is the multiplication $\cd$ of $A$, reversed and twisted by 
signs: 
$a \cdot^{\mrm{op}}  b := (-1)^{i \cd j} \cd b \cd a$
for $a \in A^i$ and $b \in A^j$. 
\end{dfn}

\begin{exer} \label{exer:2500}
Verify that $A^{\mrm{op}}$ is a DG ring. 
\end{exer}

Note that $A$ is weakly commutative iff $A = A^{\mrm{op}}$. 

\begin{rem} \label{rem:4130}
In algebraic topology and homotopy theory it is customary to use lower indices 
for DG rings (and to call them DG algebras). In other words, they use {\em 
homological grading}, as opposed to our cohomological grading. Thus our 
nonpositive DG ring $A =  \bigoplus_{i \leq 0} A^i$ becomes a nonnegative 
graded algebra $A =  \bigoplus_{i \geq 0} A_i$ in their language. 

Another notion that is common in homotopy theory is that of a 
{\em connective DG algebra}; this is a DG ring 
$A =  \bigoplus_{i \in \Z} A_i$
whose homology $\opn{H}(A) = \bigoplus_{i \in \Z} \opn{H}_i(A)$
is nonnegative. In our cohomological language, the analogue would be a 
{\em coconnective DG ring}, that is a DG ring $A$ whose cohomology 
$\opn{H}(A) = \bigoplus_{i \in \Z} \opn{H}^i(A)$
is a nonpositive graded ring. Cf.\ Exercise \ref{exer:1100} below. 

If $A$ is a coconnective DG ring, then its smart truncation  
$A' := \opn{smt}^{\leq 0}(A)$, in the sense of Definition \ref{dfn:2320}, 
is a nonpositive DG ring, and the inclusion $A' \to A$ is a quasi-isomorphism. 
\end{rem}

\begin{dfn} \label{dfn:1081}
Let $A$ be a DG ring. A {\em DG left $A$-module}
\index{Differential graded module}
is a graded left $A$-module $M = \bigoplus_{i \in \Z} M^i$, together with a 
$\K$-linear homomorphism $\d_M : M \to M$ of degree $1$, called the 
differential, satisfying $\d_M \circ \d_M = 0$ and 
\[ \d_M(a \cd m) =  \d_A(a) \cd m + (-1)^i \cd a \cd \d_M(m) \]
for $a \in A^i$ and $m \in M^j$.
\end{dfn}

DG right $A$-modules are defined likewise, but we won't deal with them 
much. This is because DG right $A$-modules are DG left modules over the 
opposite DG ring $A^{\mrm{op}}$. More precisely, if $M$ is a DG right 
$A$-module, then the formula 
$a \cd m := (-1)^{i \cd j} \cd m \cd a$,
for $a \in A^i$ and $m \in M^j$, makes $M$ into a DG left 
$A^{\mrm{op}}$-module.

As implied by Convention \ref{conv:2490}(5), all DG modules are by default 
DG left modules. 
 
\begin{prop} \label{prop:1130}
Let $A$ be a DG $\K$-ring, and let $M$ be a DG $\K$-module. 
\begin{enumerate}
\item Suppose $f : A \to \opn{End}_{\K}(M)$ is a DG $\K$-ring homomorphism. 
Then the formula $a \cd m := f(a)(m)$, for $a \in A^i$ and $m \in M^j$,
makes $M$ into a DG $A$-module. 

\item Conversely, every DG $A$-module structure on $M$, that's compatible with 
the DG $\K$-module structure, arises in this way from a 
DG $\K$-ring homomorphism $f : A \to \opn{End}_{\K}(M)$.
\end{enumerate}
\end{prop}

\begin{exer} \label{exer:2520}
Prove this proposition. 
\end{exer}

\begin{dfn} \label{dfn:1285}
Let $M$ and $N$ be DG $A$-modules. A
{\em strict homomorphism of DG $A$-modules}
\index{Strict homomorphism! of DG modules}
 $\phi : M \to N$ 
is a strict homomorphism of DG $\K$-modules (Definition \ref{dfn:1084})
that respects the action of $A$. The resulting category is denoted by 
$\cat{DGMod}_{\mrm{str}} A$, or by the short notation 
$\dcat{C}_{\mrm{str}}(A)$%
\index{1-Cstr(A)@$\dcat{C}_{\mrm{str}}(A)$}.
\end{dfn}

The category $\dcat{C}_{\mrm{str}}(A)$ is abelian. See Proposition 
\ref{prop:4130} for a more general statement. 

\begin{exer} \label{exer:1100}
Let $A$ be a DG ring. Show that the the module of cocycles 
$\opn{Z}(A) := \bigoplus_{i \in \Z} \opn{Z}^i(A)$
is a graded subring of $A$, and the module of coboundaries
$\opn{B}(A) := \bigoplus_{i \in \Z} \opn{B}^i(A)$
is a two-sided graded ideal of $\opn{Z}(A)$. Conclude that the cohomology module
$\opn{H}(A) := \bigoplus_{i \in \Z} \opn{H}^i(A)$
is a graded ring. 

Let $f : A \to B$ be a homomorphism of DG rings. Show that 
$\opn{H}(f) : \opn{H}(A) \to \opn{H}(B)$
is a graded ring homomorphism. 
\end{exer}

\begin{exer} \label{exer:2521}
Let $A$ be a DG ring. Given a DG $A$-module $M$, show that its cohomology 
$\opn{H}(M)$ is a graded $\opn{H}(A)$-module. If $\phi : M \to N$
is a homomorphism in $\dcat{C}_{\mrm{str}}(A)$, then 
$\opn{H}(\phi) : \opn{H}(M) \to \opn{H}(N)$
is a homomorphism in $\dcat{G}_{\mrm{str}}(\opn{H}(A))$.
Conclude that 
\begin{equation} \label{eqn:4130}
\opn{H} : \dcat{C}_{\mrm{str}}(A) \to \dcat{G}_{\mrm{str}}(\opn{H}(A)) 
\end{equation}
is a linear functor.  
\end{exer}

\begin{dfn} \label{dfn:1100}
\index{Differential graded module! tensor product of {\indash}s}
Let $A$ be a DG $\K$-ring, let $M$ be a right DG $A$-module, and 
let $N$ be a left DG $A$-module. By Lemma \ref{lem:1690}, $M \ot_A N$
is a graded $\K$-module. We make it into a DG $\K$-module with the differential 
from Definition \ref{dfn:2992}(1). 
\end{dfn}

\begin{dfn} \label{dfn:1101}
\index{Differential graded module! Hom between {\indash}s}
Let $A$ be a DG $\K$-ring, and let $M$ and $N$ be left DG $A$-modules.
The graded $\K$-module $\opn{Hom}_{A}(M, N)$ from Definition \ref{dfn:1690}
is made into a DG $\K$-module with the differential from
Definition \ref{dfn:2992}(2). 
\end{dfn}

Generalizing formula (\ref{eqn:2520}), for DG $A$-modules $M$ and $N$ 
there is equality
\[ \opn{Hom}_{\dcat{C}_{\mrm{str}}(A)}(M, N) = 
\opn{Z}^0 \bigl( \opn{Hom}_{A}(M, N) \bigr) . \]

\begin{prop} \label{prop:4131}
Let $A$ be a DG ring. 
\begin{enumerate}
\item The category $\dcat{C}_{\mrm{str}}(A)$ has products. Given a collection 
$\{ M_x \}_{x \in X}$ of DG $A$-modules, their product 
$M  = \prod_{x \in X} M_x$ in $\dcat{C}_{\mrm{str}}(A)$
is $M := \bigoplus_{i \in \Z}  M^i$, where 
$M^i := \prod_{x \in X} M_x^i$.

\item The category $\dcat{C}_{\mrm{str}}(A)$ has direct sums. Given a 
collection $\{ M_x \}_{x \in X}$ of DG $A$-modules, their direct sum 
$M  = \bigoplus_{x \in X} M_x$ in $\dcat{C}_{\mrm{str}}(A)$
is $M := \bigoplus_{i \in \Z}  M^i$, where 
$M^i := \bigoplus_{x \in X} M_x^i$.

\item The functor $\opn{H}$ from \tup{(\ref{eqn:4130})} commutes with all 
products and direct sums. 
\end{enumerate}
\end{prop}

\begin{exer} \label{exer:4130}
Prove this proposition. 
\end{exer}

\mysubsection{DG Categories} \label{susec:DGCats}

In Definition \ref{dfn:1692} we saw graded categories. Here is the DG version. 
 
\begin{dfn} \label{dfn:1063}
A {\em DG $\K$-linear category}
\index{Differential graded category}
is a $\K$-linear category $\cat{C}$, endowed 
with a DG $\K$-module structure on each of the morphism $\K$-modules 
$\opn{Hom}_{\cat{C}}(M_0, M_1)$. The conditions are these:
\begin{itemize}
\rmitem{a} For every object $M$, the identity automorphism $\opn{id}_M$ is a 
degree $0$ cocycle in $\opn{Hom}_{\cat{C}}(M, M)$. 

\rmitem{b} For every triple of objects $M_0, M_1, M_2 \in \cat{C}$, the 
composition homomorphism 
\[ \opn{Hom}_{\cat{C}}(M_1, M_2) \ot \opn{Hom}_{\cat{C}}(M_0, M_1)
\xar{} \opn{Hom}_{\cat{C}}(M_0, M_2) \]
is a strict homomorphism of DG $\K$-modules. 
\end{itemize}
\end{dfn}

The differential in a DG $\K$-linear category $\cat{C}$ is sometimes denoted by 
$\d_{\cat{C}}$; for instance 
\begin{equation} \label{eqn:4135}
\d_{\cat{C}} : \opn{Hom}_{\cat{C}}(M_0, M_1)^{i} \to 
\opn{Hom}_{\cat{C}}(M_0, M_1)^{i + 1 } . 
\end{equation}

The next convention extends Convention \ref{conv:4626}. 

\begin{conv} \label{conv:4627}
We shall often write ``DG category'' instead of the longer expression
``DG $\K$-linear category''. 
\end{conv}

If $\cat{C}'$ is a full subcategory of a DG category $\cat{C}$, then of course
$\cat{C}'$ is itself a DG category. 

\begin{dfn} \label{dfn:1286}
Let $\cat{C}$ be a DG category.
\begin{enumerate}
\item A morphism $\phi \in \opn{Hom}_{\cat{C}}(M, N)^i$
is called a {\em degree $i$ morphism}. 

\item A morphism $\phi \in \opn{Hom}_{\cat{C}}(M, N)$
is called a {\em cocycle} if $\d_{\cat{C}}(\phi) = 0$. 

\item A morphism $\phi : M \to N$ in $\cat{C}$ is 
called a {\em strict morphism}%
\index{Differential graded category! strict morphism in}
if it is a degree $0$ cocycle.
\end{enumerate}
\end{dfn}

\begin{lem} \label{lem:1211}
Let $\cat{C}$ be a DG category, and for $i = 0, 1, 2$ let 
$\phi_i : M_{i} \to M_{i + 1}$ be a morphism in $\cat{C}$ of degree $k_i$.  
\begin{enumerate}
\item The morphism $\phi_1 \circ \phi_0$ has degree $k_0 + k_1$, and 
\[ \d_{\cat{C}}(\phi_1 \circ \phi_0) = 
\d_{\cat{C}}(\phi_1) \circ \phi_0 +
(-1)^{k_1} \cd \phi_1 \circ \d_{\cat{C}}(\phi_0) . \]

\item If $\phi_0$ and $\phi_1$ are cocycles, then so is $\phi_1 \circ \phi_0$.

\item If  $\phi_1$ is a coboundary, and $\phi_0$ and $\phi_2$ are cocycles, 
then $\phi_2 \circ \phi_1 \circ \phi_0$ is a coboundary.
\end{enumerate}
\end{lem}

\begin{proof}
(1) This is just a rephrasing of item (b) in Definition \ref{dfn:1063}. 

\medskip \noindent 
(2) This is immediate from (1). 

\medskip \noindent 
(3) Say $\phi_1 = \d_{\cat{C}}(\psi_1)$ for some degree $k_1 - 1$ morphism 
$\psi_1 : M_{1} \to M_2$. Then 
\[ \phi_2 \circ \phi_1 \circ \phi_0 = 
\d_{\cat{C}} \bigl( (-1)^{k_2} \cd \phi_2 \circ \psi_1 \circ \phi_0 \bigr) .
\qedhere \]
\end{proof}

The previous lemma makes the next definition possible. 

\begin{dfn} \label{dfn:1212} \mbox{ }
Let $\cat{C}$ be a DG category.
\begin{enumerate}
\item The {\em strict subcategory}%
\index{Differential graded category! strict subcategory of}%
\index{Strict subcategory! of a DG category}
of $\cat{C}$ is the category 
$\opn{Str}(\cat{C})$, with the same objects as 
$\cat{C}$, but with strict morphisms only. Thus 
\[ \opn{Hom}_{\opn{Str}(\cat{C})}(M, N) = 
\opn{Z}^0 \bigl( \opn{Hom}_{\cat{C}}(M, N) \bigr) . \]

\item The {\em homotopy category}%
\index{Differential graded category! homotopy category of}
of $\cat{C}$ is the category 
$\opn{Ho}(\cat{C})$, with the same objects as $\cat{C}$, 
and with morphism sets 
\[  \opn{Hom}_{\opn{Ho}(\cat{C})}(M, N) := 
\opn{H}^0 \bigl( \opn{Hom}_{\cat{C}}(M, N) \bigr) . \]

\item We denote by 
$\opn{P} : \opn{Str}(\cat{C}) \to \opn{Ho}(\cat{C})$
the functor which is  the identity on objects, and sends a strict morphism to 
its homotopy class. 
\end{enumerate}
\end{dfn}

The categories $\opn{Str}(\cat{C})$ and $\opn{Ho}(\cat{C})$ are linear, and 
the inclusion functor $\opn{Str}(\cat{C}) \to \cat{C}$ and the projection 
functor $\opn{P} : \opn{Str}(\cat{C}) \to \opn{Ho}(\cat{C})$
are linear. The first is faithful, and the second is full.

\begin{exa} \label{exa:1286}
If $\cat{A}$ is a DG category, then for every object $x \in \cat{A}$, 
its set of endomorphisms $A := \opn{End}_{\cat{A}}(x)$ is a DG ring.
Conversely, every DG ring $A$ can be viewed as a DG category $\cat{A}$ 
with a single object.
\end{exa}

\begin{exa} \label{exa:1285}
Let $A$ be a DG ring. The set of DG $A$-modules forms a 
DG category $\cat{DGMod} A$, in which the morphism DG $\K$-modules are
\[ \opn{Hom}_{\cat{DGMod} A}(M, N) := \opn{Hom}_{A}(M, N) \]
from Definition \ref{dfn:1101}. 
We shall often write $\dcat{C}(A) := \cat{DGMod} A$%
\index{1-C(A)@$\dcat{C}(A)$}. 
The strict subcategory here is 
\[ \opn{Str}(\cat{DGMod} A) = \cat{DGMod}_{\mrm{str}} A 
= \dcat{C}_{\mrm{str}}(A) ; \]
cf.\  Definition \ref{dfn:1285}. 
\end{exa}

Here is a result that will be used later. 

\begin{prop} \label{prop:1220}
Let $\phi : M \to N$ be a degree $i$ isomorphism in a DG category $\cat{C}$. 
Assume $\phi$ is a cocycle, namely 
$\d_{\cat{C}}(\phi) = 0$. 
Then its inverse $\phi^{-1} : N \to M$ is also a cocycle.  
\end{prop}

\begin{proof}
According the Leibniz rule (Lemma \ref{lem:1211}(1)), and the fact that 
$\opn{id}_M$ is a cocycle, we have
\[ \begin{aligned}
& 0 = \d_{\cat{C}}(\opn{id}_M) = \d(\phi^{-1} \circ \phi)
\\ & \ 
= \d_{\cat{C}}(\phi^{-1}) \circ \phi + 
(-1)^{-i} \cd \phi^{-1} \circ \, \d_{\cat{C}}(\phi) 
= \d_{\cat{C}}(\phi^{-1}) \circ \phi . 
\end{aligned} \]
Because $\phi$ is an isomorphism, we conclude that 
$\d_{\cat{C}}(\phi^{-1}) = 0$. 
\end{proof}

\begin{rem} \label{rem:1270}
The fact that the concept ``DG category'' refers both to DG rings 
(Example \ref{exa:1286}) and to the categories of DG modules over them (Example 
\ref{exa:1285}) is a source of frequent confusion. See Remarks \ref{rem:2517} 
and \ref{rem:1245}. 
\end{rem}

\begin{rem} \label{rem:1175}
For other accounts of DG categories see the relatively old references 
\cite{Kelly}, \cite{Kel1}, \cite{BoKa}, or the more recent \cite{To1}. 
An internet search can give plenty more information,
including the relation to simplicial and infinity categories. 

In this book we shall be exclusively concerned with the DG categories 
$\dcat{C}(A, \cat{M})$, to be introduced in Subsection \ref{subsec:DGModinM},
which have a lot more structure than other DG categories.

See Remark \ref{rem:1245} regarding the DG category
$\dcat{C}(\cat{A})$ of DG modules over a DG category $\cat{A}$, 
in the sense of \cite{Kel1}; and Remark \ref{rem:4960} regarding 
$\mrm{A}_{\infty}$ categories. 
\end{rem}

\mysubsection{DG Functors} \label{subsec:DGFunc}

Here $\cat{C}$ and $\cat{D}$ are DG $\K$-linear categories (see Definition 
\ref{dfn:1063}). When we forget differentials,  $\cat{C}$ and $\cat{D}$ become 
graded $\K$-linear categories. So we can talk about graded $\K$-linear
functors $\cat{C} \to \cat{D}$, as in Definition \ref{dfn:1694}. 

Recall the meaning of a strict homomorphism of DG $\K$-modules:
it has degree $0$ and commutes with the differentials. 

\begin{dfn} \label{dfn:1065}
Let $\cat{C}$ and $\cat{D}$ be DG $\K$-linear categories.
A functor  $F : \cat{C} \to \cat{D}$
is called a {\em DG $\K$-linear functor}
\index{Differential graded functor} 
if it satisfies this condition: 
\begin{itemize}
\item[$\vartriangleright$] For every pair of objects $M_0, M_1 \in \cat{C}$, 
the function 
\[ F : \opn{Hom}_{\cat{C}}(M_0, M_1) \to  
\opn{Hom}_{\cat{D}} \bigl( F(M_0), F(M_1) \bigr) \] 
is a strict homomorphism of DG $\K$-modules. 
\end{itemize}
\end{dfn}

In other words, $F$ is a DG functor if it is a graded functor, and 
$\d_{\cat{D}} \circ F = F \circ \d_{\cat{C}}$
as degree $1$ homomorphisms 
\[ \opn{Hom}_{\cat{C}}(M_0, M_1) \to  
\opn{Hom}_{\cat{D}} \bigl( F(M_0), F(M_1) \bigr) . \]

To match Convention \ref{conv:4627}, here is: 

\begin{conv} \label{conv:4628}
We shall usually write ``DG functor'' instead of the longer expression
``DG $\K$-linear functor''. 
\end{conv}

\begin{exa} \label{exa:1250}
Let $f : A \to B$ be a homomorphism of DG rings, and let 
$\cat{A}$ and $\cat{B}$ be the corresponding single-object DG  
categories. Then $f$ gives rise to a DG functor 
$F : \cat{A} \to \cat{B}$.
\end{exa}

Other examples of DG functors, more relevant to our study, will be given 
in Subsection \ref{subsec:exa-DGfun}. 

\begin{dfn} \label{dfn:1195}
Let $F, G : \cat{C} \to \cat{D}$ be DG functors.
\begin{enumerate}
\item A {\em degree $i$ morphism of DG functors}%
\index{Differential graded functor! morphism of {\indash}s}
$\eta : F \to G$ is a 
degree $i$ morphism of graded functors, as in Definition \ref{dfn:1695}.

\item Let $\eta : F \to G$ be a degree $i$ morphism of DG functors.
For each object $M \in \cat{C}$ there is a degree $i + 1$ morphism 
$\d_{\cat{D}}(\eta_M) : F(M) \to G(M)$
in $\cat{D}$. 
We define the collection of morphisms  
$\d_{\cat{D}}(\eta) := 
\{ \d_{\cat{D}}(\eta_M) \}_{M \in \cat{C}}$.

\item A {\em strict morphism of DG functors}
\index{Differential graded functor! strict morphism of {\indash}s}
is a degree $0$ morphism of graded functors  
$\eta : F \to G$ such that $\d_{\cat{D}}(\eta) = 0$. 
\end{enumerate}
\end{dfn}

Explanation of item (2) of the definition: for an object $M \in \cat{C}$, we 
have
\[ \eta_M \in \opn{Hom}_{\cat{D}} \bigl( F(M), G(M) \bigr)^i . \]
Now $\opn{Hom}_{\cat{D}} \bigl( F(M), G(M) \bigr)$
is a DG $\K$-module with differential $\d_{\cat{D}}$; see 
formula (\ref{eqn:4135}). Thus 
\[ \d_{\cat{D}}(\eta_M) \in 
\opn{Hom}_{\cat{D}} \bigl( F(M), G(M) \bigr)^{i + 1} . \]

\begin{prop} \label{prop:1295}
In the situation of Definition \tup{\ref{dfn:1195}}, the collection of 
morphisms $\d_{\cat{D}}(\eta)$  is a degree $i + 1$ morphism of DG 
functors $F \to G$. 
\end{prop}

\begin{exer} \label{exer:2522}
Prove this proposition.
\end{exer}

The categories $\opn{Str}(\cat{C})$ and 
$\opn{Ho}(\cat{C})$ were introduced in 
Definition \ref{dfn:1212}.

\begin{prop} \label{prop:1150}
Let $F : \cat{C} \to \cat{D}$ be a DG functor. Then $F$ induces linear 
functors
$\opn{Str}(F) : \opn{Str}(\cat{C}) \to \opn{Str}(\cat{D})$
and
$\opn{Ho}(F) : \opn{Ho}(\cat{C}) \to  \opn{Ho}(\cat{D})$.
\end{prop}

\begin{proof}
Because $F$ is a DG functor, it sends $0$-cocycles in 
$\opn{Hom}_{\cat{C}}(M_0, M_1)$ to $0$-cocycles in 
$\opn{Hom}_{\cat{D}}(F(M_0), F(M_1))$. 
The same for $0$-coboundaries.  
\end{proof}

By abuse of notation, and when there is no danger for confusion, we will 
sometimes write $F$ instead of $\opn{Str}(F)$ or $\opn{Ho}(F)$.

\begin{exer} \label{exer:1250}
Let $\cat{A}$ and $\cat{C}$ be DG categories, and assume 
$\cat{A}$ is small. Define 
$\cat{DGFun}(\cat{A}, \cat{C})$ to be the set of DG functors 
$F : \cat{A} \to \cat{C}$.
Show that $\cat{DGFun}(\cat{A}, \cat{C})$ is a DG category, 
where the morphisms are from Definition \tup{\ref{dfn:1195}(1)}, and their 
differentials are from Definition \tup{\ref{dfn:1195}(2)}.
\end{exer}

\mysubsection{Complexes in Abelian Categories} \label{subsec:complexes}

Here we recall facts about complexes from the classical homological theory, and 
place them within our context. In this subsection $\cat{M}$ is a $\K$-linear 
abelian category. 

A {\em complex of objects of $\cat{M}$}, or a 
{\em complex in $\cat{M}$}%
\index{Abelian category! complex in},
is a diagram 
\begin{equation} \label{eqn:1065}
\bigl( \cdots \to M^{-1} \xar{\d_M^{-1}}  M^0 \xar{\d_M^0}  
M^1 \xar{\d_M^1} M^2 \to \cdots \bigr)
\end{equation}
of objects and morphisms in $\cat{M}$, such that $\d_M^{i+1} \circ \d_M^i = 0$.
The collection of objects $M := \{ M^i \}_{i \in \Z}$ 
is nothing but a graded object of $\cat{M}$, as defined in Subsection 
\ref{subsec:gr-alg}. The collection of morphisms 
$\d_M := \{ \d^i_M \}_{i \in \Z}$ 
is called a {\em differential}, or a {\em coboundary operator}.
Thus a complex is a pair $(M, \d_M)$ made up of a graded object $M$ and a 
differential $\d_M$ on it. 
We sometimes write $\d$ instead of $\d_M$ or $\d_M^i$. 
At other times we leave the differential implicit, and just refer to the 
complex as $M$. 

Let $N$ be another complex in $\cat{M}$. A {\em strict morphism of complexes} 
$\phi : M \to N$ is a collection 
$\phi = \{ \phi^i \}_{i \in \Z}$ of morphisms $\phi^i : M^i \to N^i$ in
$\cat{M}$, such that
\begin{equation} \label{eqn:1254}
\d^i_N \circ \phi^i = \phi^{i+1} \circ \d^i_M .
\end{equation}
Note that a strict morphism $\phi : M \to N$ can be viewed as a commutative 
diagram 
\[ \UseTips  \xymatrix @C=8ex @R=6ex {
\cdots 
\ar[r]
&
M^i
\ar[r]^{\d_M^i}
\ar[d]_{\phi^i}
&
M^{i+1}
\ar[r]
\ar[d]_{\phi^{i+1}}
&
\cdots
\\
\cdots 
\ar[r]
&
N^i
\ar[r]^{\d_N^i}
&
N^{i+1}
\ar[r]
&
\cdots
} \]
in $\cat{M}$. 
The identity automorphism $\opn{id}_M$ of the complex $M$ is a strict 
morphism.  

\begin{rem} \label{rem:1210}
In most textbooks, what we call ``strict morphism of complexes'' is simply 
called a ``morphism of complexes''. See Remark \ref{rem:1190} for an 
explanation. 
\end{rem}

Let us denote by 
$\dcat{C}_{\mrm{str}}(\cat{M})$%
\index{1-Cstr(M)@$\dcat{C}_{\mrm{str}}(\cat{M})$}
the category of complexes in 
$\cat{M}$, with strict morphisms. This is a $\K$-linear abelian category. 
Indeed, the direct sum of complexes is the degree-wise direct sum, i.e.\ 
$(M \oplus N)^i = M^i \oplus N^i$.
The same for kernels and cokernels.
If  $\cat{N}$ is a full abelian subcategory of $\cat{M}$, then  
$\dcat{C}_{\mrm{str}}(\cat{N})$ is a full abelian subcategory of 
$\dcat{C}_{\mrm{str}}(\cat{M})$.

A single object $M^0 \in \cat{M}$ can be viewed as a complex 
\[ M :=  \bigl( \cdots \to 0 \to M^0 \to 0 \to \cdots \bigr) , \] 
where $M^0$ is in degree $0$; the differential of this complex is of course 
zero. The assignment $M^0 \mapsto M$ is a fully faithful $\K$-linear functor 
$\cat{M} \to \dcat{C}_{\mrm{str}}(\cat{M})$.

Let $M$ and $N$ be complexes in $\cat{M}$. As in (\ref{eqn:1700}) there is a 
graded $\K$-module $\opn{Hom}_{\cat{M}}(M, N)$. It is a DG $\K$-module with 
differential $\d$ given by the formula 
\begin{equation} \label{eqn:1130}
\d(\phi) := \d_N \circ \phi - (-1)^i \cd \phi \circ \d_M 
\end{equation}
for $\phi \in \opn{Hom}_{\cat{M}}(M, N)^i$. 
It is easy to check that $\d \circ \d = 0$. 
We sometimes denote this differential by
$\d_{\opn{Hom}}$ or $\d_{\opn{Hom}_{\cat{M}}(M, N)}$.

Thus, an element $\phi \in \opn{Hom}_{\cat{M}}(M, N)^i$
is a collection $\phi = \{ \phi^j \}_{j \in \Z}$ of morphisms 
$\phi^j : M^j \to N^{j + i}$. In a diagram, for $i = 2$, it looks like this:
\[ \UseTips  \xymatrix @C=7ex @R=8ex {
\cdots 
\ar[r]
&
M^j
\ar[r]^{\d}
\ar[drr]^(0.6){\phi^j}
&
M^{j+1}
\ar[r]^{\d}
\ar[drr]^(0.6){\phi^{j+1}}
&
M^{j + 2}
\ar[r]^{\d}
&
M^{j + 3}
\ar[r]
&
\cdots
\\
\cdots 
\ar[r]
&
N^j
\ar[r]^{\d}
&
N^{j + 1}
\ar[r]^{\d}
&
N^{j + 2}
\ar[r]^{\d}
&
N^{j + 3}
\ar[r]
&
\cdots 
} \]
Warning: since $\phi$ does not have to commute with the differentials, this is 
usually not a commutative diagram! 

For a triple of complexes $M_0, M_1, M_2$ and degrees $i_0, i_1$ there are 
$\K$-linear homomorphisms
\[ \begin{aligned}
& 
\opn{Hom}_{\cat{M}}(M_1, M_2)^{i_1} \otimes_{} 
\opn{Hom}_{\cat{M}}(M_0, M_1)^{i_0} \xar{} 
\opn{Hom}_{\cat{M}}(M_0, M_2)^{i_0 + i_1} ,
\\
&
\phi_1 \ot \phi_0 \mapsto \phi_1 \circ \phi_0 . 
\end{aligned} \]

\begin{lem} \label{1213}
The composition homomorphism 
\[ \opn{Hom}_{\cat{M}}(M_1, M_2) \otimes_{} 
\opn{Hom}_{\cat{M}}(M_0, M_1) \xar{} 
\opn{Hom}_{\cat{M}}(M_0, M_2) \]
is a strict homomorphism of DG $\K$-modules. 
\end{lem}

\begin{exer}
Prove the lemma. 
\end{exer}

The lemma justifies the next definition.

\begin{dfn}
Let $\dcat{C}(\cat{M})$%
\index{1-C(M)@$\dcat{C}(\cat{M})$}
be the DG $\K$-linear category whose objects are the 
complexes in $\cat{M}$, and whose morphism DG $\K$-modules are 
$\opn{Hom}_{\cat{M}}(M, N)$, from formulas (\ref{eqn:1700}) and 
(\ref{eqn:1130}).  
\end{dfn}

It is clear, from comparing formulas (\ref{eqn:1130}) and (\ref{eqn:1254}), 
that the strict morphisms of complexes defined at the top of this subsection 
are the same as those from Definition \ref{dfn:1212}(1). In other words, 
$\opn{Str}(\dcat{C}(\cat{M})) = \dcat{C}_{\mrm{str}}(\cat{M})$.

\begin{rem}
A possible ambiguity could arise in the meaning of 
$\opn{Hom}_{\cat{M}}(M, N)$ 
if $M, N \in \cat{M}$. Does it mean the $\K$-module of morphisms in the 
category $\cat{M}$~? Or, if we view $M$ and $N$ as complexes by the canonical
embedding $\cat{M} \subseteq \dcat{C}(\cat{M})$, does 
$\opn{Hom}_{\cat{M}}(M, N)$
mean the complex of morphisms in the DG category $\dcat{C}(\cat{M})$~? 
It turns out that there is no actual difficulty: since the complex of 
$\K$-modules $\opn{Hom}_{\cat{M}}(M, N)$ is concentrated in degree $0$, 
we may view it as a single $\K$-module, and this is precisely the $\K$-module 
of 
morphisms in the category $\cat{M}$.
\end{rem}

When $\cat{M} = \cat{Mod} A$ for a ring $A$, there is no essential distinction 
between complexes and DG modules. The next proposition is the DG version of 
Exercise \ref{exer:4112}. 

\begin{prop} \label{prop:1065}
Let $A$ be a ring. Given a complex 
$M \in \dcat{C}(\cat{Mod} A)$, with notation as in \tup{(\ref{eqn:1065})}, 
define the DG $A$-module 
$F(M) := \bigoplus_{i \in \Z} M^i$,
with differential $\d := \sum_{i \in \Z} \d_M^i$. Then the functor 
$F : \dcat{C}(\cat{Mod} A) \to \cat{DGMod} A$
is an isomorphism of DG categories. 
\end{prop}

\begin{exer}
Prove this proposition. (Hint: choose good notation.)
\end{exer}

\mysubsection{The Long Exact Cohomology Sequence}
\label{subsec:long-cohom-seq}

As in the previous subsection, $\cat{M}$ is a $\K$-linear 
abelian category. Here we give a detailed proof of the existence and 
functoriality of the long exact cohomology sequence, using our {\em sheaf 
tricks} from Subsection \ref{subsec:sheaf-tricks}. This allows us to avoid the 
Freyd-Mitchell Theorem. 

Consider a short exact sequence
\begin{equation} \label{eqn:3640}
\bsym{E} = \bigl( 0 \to L \xar{\phi} M \xar{\psi} N \to 0 \bigr)
\end{equation}
in the abelian category $\dcat{C}_{\mrm{str}}(\cat{M})$. 
This means that $\phi$ and $\psi$ are morphisms in 
$\dcat{C}_{\mrm{str}}(\cat{M})$, and in each degree $i$ there is a short 
exact sequence 
\[ 0 \to L^i \xar{\phi^i} M^i \xar{\psi^i} N^i \to 0 \]
in $\cat{M}$. The next definitions and lemmas refer to the short exact 
sequence $\bsym{E}$. In them we shall make use the notation from 
Subsection \ref{subsec:sheaf-tricks}.

Let us denote by
$\pi_N : \opn{Z}^i(N) \surj \opn{H}^i(N)$
the canonical epimorphism in $\cat{M}$; and likewise for $L$ and $M$. 

\begin{dfn} \label{dfn:3645}
Let $\bar{n} \in \Ga(U, \opn{H}^i(N))$ for some $U \in \cat{M}$,
and let $V \surj U$ be a covering of $U$. 
By a {\em connecting triple} for $\bar{n}$ over $V$
we mean a triple $(n, m, l)$ of sections
$n \in \Ga(V, \opn{Z}^i(N))$, 
$m \in \Ga(V, M^i)$ and $l \in \Ga(V, \opn{Z}^{i + 1}(L))$,
such that
$\pi_N(n) = \bar{n} \in \Ga(V, \opn{H}^i(N))$, 
$\psi(m) = n \in \Ga(V, N^i)$
and
$\phi(l) = \d(m) \in  \Ga(V, M^{i + 1})$.
\end{dfn}

Here is the relevant diagram in $\cat{M}$~:
\[ \UseTips \xymatrix @C=6ex @R=6ex {
\opn{H}^{i + 1}(L)
&
\opn{Z}^{i + 1}(L)
\ar@{->>}[l]_{\pi_L}
\ar[r]^{\phi}
&
M^{i + 1}
&
\\
&
&
M^{i}
\ar[u]^{\d}
\ar[r]^{\psi}
&
\opn{Z}^{i}(N)
\ar@{->>}[r]^{\pi_N}
&
\opn{H}^{i}(N)
} \]

It is possible that some coverings do not admit a connecting triple. But some 
do:

\begin{lem} \label{lem:3650}
Let $\bar{n} \in \Ga(U, \opn{H}^i(N))$. 
There exists a covering $\rho : V \surj U$ that admits a connecting triple 
$(n, m, l)$ for $\bar{n}$.
\end{lem}

\begin{proof}
Because $\pi_N : \opn{Z}^{i}(N) \to \opn{H}^{i}(N)$ is an epimorphism, by 
the first sheaf trick there is a covering $V' \surj U$ with 
a section $n \in \Ga(V', \opn{Z}^{i}(N))$
such that $\pi_N(n) = \bar{n}$. 

Because $\psi : M^i \to N^i$ is an epimorphism, there is a covering 
$V'' \surj V'$ with 
a section $m \in \Ga(V'', M^i)$ such that $\psi(m) = n$. 

Consider the section $\d(m) \in \Ga(V'', M^{i + 1})$. 
We have 
$\psi(\d(m)) = \d(\psi(m)) = \d(n) = 0$.
By exactness at $M^{i + 1}$, there is a covering $V \surj V''$ with 
a section $l \in \Ga(V, L^{i + 1})$ such that $\phi(l) = \d(m)$
in $\Ga(V, M^{i + 1})$. 

Now 
$\phi(\d(l)) = \d(\phi(l)) = \d(\d(m))) = 0$.
Because $\phi$ is a monomorphism, it follows that $\d(l) = 0$. Thus 
$l \in \Ga(V, \opn{Z}^{i + 1}(L))$. 
\end{proof}

\begin{lem} \label{lem:3651} 
Let $\bar{n} \in \Ga(U, \opn{H}^i(N))$. Suppose $V \surj U$
is a covering that admits connecting triples 
$(n, m, l)$ and $(n', m', l')$ for $\bar{n}$.
Then 
$\pi_L(l') = \pi_L(l)$ in $\Ga(V, \opn{H}^{i + 1}(L))$.
\end{lem}

\begin{proof}
We know that $\pi_N(n) = \bar{n} = \pi_N(n')$ in $\Ga(V, \opn{H}^{i}(N))$. 
Thus 
$n - n' \in \opn{Ker}(\pi_N) \sub \Ga(V, \opn{Z}^{i}(N))$.
Looking at the  exact sequence 
\[ N^{i - 1} \xar{\d} \opn{Z}^{i}(N) \xar{\pi_N} \opn{H}^{i}(N) \to 0 \]
and invoking the first sheaf trick, we see that there's a covering 
$W \surj V$ with a section $\til{n} \in \Ga(W, N^{i - 1})$
such that $\d(\til{n}) = n - n'$. 

Again invoking the first sheaf trick, there's a covering 
$W' \surj W$ with a section $\til{m} \in \Ga(W', M^{i - 1})$
such that $\psi(\til{m}) = \til{n}$. 

Consider the section 
$(m - m') - \d(\til{m}) \in \Ga(W', M^{i})$.
It satisfies 
\[ \psi \bigl( (m - m') - \d(\til{m}) \bigr) = 
(n - n') - \d(\til{n}) = 0 . \]
By the first sheaf trick, there's a covering 
$W'' \surj W'$ with a section $\til{l} \in \Ga(W'', L^{i})$
such that 
$\phi(\til{l}) = (m - m') - \d(\til{m})$.
Now 
\[ \phi \bigl( \d(\til{l}) - (l - l') \bigr) = 
\d(m - m') - \phi(l - l') = 0 . \]
Because $\phi$ is a monomorphism, it follows that 
$l - l' = \d(\til{l})$. 
Therefore $\pi_L(l) = \pi_L(l')$ in $\Ga(W'', \opn{H}^{i + 1}(L))$. 
But $W'' \surj V$ is a covering, so 
$\pi_L(l) = \pi_L(l')$ in $\Ga(V, \opn{H}^{i + 1}(L))$. 
\end{proof}

Lemma \ref{lem:3651} justifies the next definition. 

\begin{dfn} \label{dfn:3650}
Let $\bar{n} \in \Ga(U, \opn{H}^i(N))$. Suppose $V \surj U$ is a covering 
that admits a connecting triple. We define 
$\de_V(\bar{n}) \in \Ga(V, \opn{H}^{i + 1}(L))$
to be the unique section such that 
$\de_V(\bar{n}) = \pi_L(l)$ for every 
connecting triple $(n, m, l)$ for $\bar{n}$ over $V$.
\end{dfn}

\begin{lem} \label{lem:3640}
Let $\bar{n} \in \Ga(U, \opn{H}^i(N))$, and let  
\[ \UseTips \xymatrix @C=6ex @R=6ex {
V'
\ar@{->>}[r]^{\rho'}
\ar[d]_{\tau}
&
U'
\ar[d]^{\si}
\\
V
\ar@{->>}[r]^{\rho}
&
U
} \]
be a commutative diagram in $\cat{M}$, in which the horizontal arrows are 
coverings. Suppose $V$ admits a connecting triple for $\bar{n}$.
Then $V'$ admits a connecting triple for 
$\bar{n}' := \si^*(\bar{n})\in \Ga(U', \opn{H}^i(N))$,
and there is equality 
$\de_{V'}(\bar{n}') = \tau^*(\de_V(\bar{n}))$
in $\Ga(V', \opn{H}^{i + 1}(L))$.
\end{lem}

\begin{proof}
If $(n, m, l)$ is a connecting triple for $\bar{n}$ over $V$, then 
$\bigl( \tau^*(n), \tau^*(m), \tau^*(l) \bigr)$
is a connecting triple for $\bar{n}'$ over $V'$. Thus
$\de_{V'}(\bar{n}') = \pi_L(\tau^*(l)) = \tau^*(\pi_L(l)) = 
\tau^*(\de_V(\bar{n}))$.
\end{proof}

\begin{lem} \label{lem:3653}
Let $\bar{n} \in \Ga(U, \opn{H}^i(N))$. Suppose $\rho : V \surj U$
is a covering that admits a connecting triple for $\bar{n}$.
Then $\de_V(\bar{n})$ lies in the subgroup 
$\Ga(U, \opn{H}^{i + 1}(L))$ of $\Ga(V, \opn{H}^{i + 1}(L))$.
\end{lem}

\begin{proof}
This is a ``descent'' argument, using the second sheaf trick.
Let $V_1, V_2 := V$ and  $\rho_1, \rho_2 := \rho$.
In the notation of Proposition \ref{prop:3636}, applied to the object 
$Q := \opn{H}^{i + 1}(L) \in \cat{M}$, there is an exact sequence 
\[ 0 \to \Ga(U, Q) \xar{(\rho_1^*, \rho_2^*)} \Ga(V_1, Q) \times \Ga(V_2, Q)
\xar{(\si_1^*, - \si_2^*)} \Ga(W, Q) .  \]
Here 
\[ \rho_1 \circ \si_1 = \rho_2 \circ \si_2 : 
W = V_1 \times_{U} V_2 \to U \]
is a covering too. 
According to Lemma \ref{lem:3640}, the section 
\[ \bigl( \de_{V_1}(\bar{n}), \de_{V_2}(\bar{n}) \bigr) \in 
\Ga(V_1, Q) \times \Ga(V_2, Q) \]
satisfies 
\[ (\si_1^*, - \si_2^*)\bigl( \de_{V_1}(\bar{n}), \de_{V_2}(\bar{n}) \bigr) =
\de_W(\bar{n}) - \de_W(\bar{n}) = 0 . \]
Hence $\de_V(\bar{n})$ is in the image of $\rho^*$. 
\end{proof}

\begin{lem} \label{lem:3641}
There is a unique morphism 
$\de_{\bsym{E}}^i : \opn{H}^{i}(N) \to \opn{H}^{i + 1}(L)$
in $\cat{M}$ with the following property\tup{:}
\begin{itemize}
\item[($\dag$)] For every $U \in \cat{M}$, every section 
$\bar{n} \in \Ga(U, \opn{H}^i(N))$, and every
covering $V \surj U$ that admits a connecting triple for $\bar{n}$, there is 
equality 
$\de_{\bsym{E}}^i(\bar{n}) = \de_V(\bar{n})$
in $\Ga(U, \opn{H}^{i + 1}(L))$.
\end{itemize}
\end{lem}

\begin{proof}
We look at the universal section of $\opn{H}^{i}(N)$, namely we take 
$U_0 := \opn{H}^{i}(N)$ and 
$\bar{n}_0 := \opn{id}_{U_0} \in \Ga(U_0, \opn{H}^i(N))$.
Let $\rho_0 : V_0 \surj U_0$ be a covering that admits a connecting triple for 
$\bar{n}_0$; see Lemma \ref{lem:3650}. According to Lemma \ref{lem:3653} we 
know that 
\[ \de_{V_0}(\bar{n}_0) \in \Ga(U_0, \opn{H}^{i + 1}(L)) =
\opn{Hom}_{\cat{M}} \bigl( \opn{H}^i(N), \opn{H}^{i + 1}(L) \bigr) . \]
Define 
$\de_{\bsym{E}}^i := \de_{V_0}(\bar{n}_0)$.

We need to prove that $\de_{\bsym{E}}^i$ has property ($\dag$). 
So let $U \in \cat{M}$, let $\bar{n} \in \Ga(U, \opn{H}^i(N))$, and let
$\rho : V \surj U$ a covering that admits a connecting triple for $\bar{n}$.
Define the morphism $\si := \bar{n} : U \to U_0$. 
We get a morphism 
$\si \circ \rho : V \to U_0$. 
Define $V' := V \times_{U_0} V_0$. By Lemma \ref{lem:3630}(2) the induced 
morphism $V' \to V$ is a covering. Therefore, by composing it with $\rho$  we 
get a covering $\rho' : V' \surj U$. 

Consider the commutative diagram 
\[ \UseTips \xymatrix @C=6ex @R=6ex {
V'
\ar@{->>}[r]^{\rho'}
\ar[d]_{\tau}
&
U
\ar[d]^{\si}
\\
V_0
\ar@{->>}[r]^{\rho_0}
&
U_0
} \]
in $\cat{M}$. The horizontal arrows are coverings. Now 
$\bar{n} = \si = \si \circ \opn{id}_{U_0} = \si^*(\bar{n}_0)$.
So, according to Lemma \ref{lem:3640}, there is equality
\[ \de_{V'}(\bar{n}) = \tau^*(\de_{V_0}(\bar{n}_0)) =
\tau^*(\de_{\bsym{E}}^i \circ \rho_0) = 
\de_{\bsym{E}}^i \circ \rho_0 \circ \tau = 
\de_{\bsym{E}}^i \circ \si \circ \rho' = \de_{\bsym{E}}^i(\bar{n}) \]
in $\Ga(V', \opn{H}^{i + 1}(L))$. 
But $\rho' : V' \surj U$ is a covering, and hence there is equality 
$\de_{V'}(\bar{n}) = \de_{\bsym{E}}^i(\bar{n})$
in $\Ga(U, \opn{H}^{i + 1}(L))$.
\end{proof}

\begin{dfn} \label{dfn:3640}
The morphism 
$\de_{\bsym{E}}^i : \opn{H}^{i}(N) \to \opn{H}^{i + 1}(L)$
is called the {\em $i$-th connecting morphism}
\index{Connecting morphism}
 of the short exact sequence
$\bsym{E}$ from (\ref{eqn:3640}). 
\end{dfn}

\begin{thm}[Long Exact Sequence in Cohomology] \label{thm:3640}
\index{Long exact sequence in cohomology}
Let $\cat{M}$ be an abelian category. 
Given a short exact sequence 
\[ \bsym{E} = \bigl( 0 \to L \xar{\phi} M \xar{\psi} N \to 0 \bigr) \]
in $\dcat{C}_{\mrm{str}}(\cat{M})$, 
the sequence 
\[ \cdots \to \opn{H}^{i}(L) \xar{\opn{H}^{i}(\phi)}
\opn{H}^{i}(M) \xar{\opn{H}^{i}(\psi)}
\opn{H}^{i}(N) \xar{\de_{\bsym{E}}^i} \opn{H}^{i + 1}(L) \to \cdots \]
in $\cat{M}$ is exact. 
\end{thm}

\begin{proof}
We need to check exactness at $\opn{H}^{i}(N)$, $\opn{H}^{i + 1}(L)$  and
$\opn{H}^{i}(M)$. This will be done using the first sheaf trick (Proposition 
\ref{prop:3635}) and the second sheaf trick (Proposition \ref{prop:3636}). 

\medskip \noindent 
$\triangleright$  Step 1. 
Exactness at $\opn{H}^{i}(N)$. This is done in two substeps. 
For simplification we shall write $\bar{\psi} := \opn{H}^{i}(\psi)$. 

\medskip \noindent 
$\triangleright$  $\triangleright$ Substep 1.a.
We start by proving that 
$\de_{\bsym{E}}^i \circ \bar{\psi} = 0$.
According to Proposition \ref{prop:3635}(1) it suffices 
to show that given an element 
$\bar{m} \in \Ga(U, \opn{H}^{i}(N))$, 
with image 
$\bar{n} := \bar{\psi}(\bar{m}) \in \Ga(U, \opn{H}^{i}(N))$, 
the element 
$\de_{\bsym{E}}^i(\bar{n})$ is zero. 

By Proposition \ref{prop:3635}(3) there is a covering $V \surj U$, and an 
element $m \in \Ga(V, \opn{Z}^{i}(M))$, 
such that 
$\pi_M(m) = \bar{m} \in  \Ga(V, \opn{H}^{i}(M))$. 
Letting $n := \psi(m) \in  \Ga(V, \opn{Z}^{i}(N))$, 
we see that $\pi_N(n) = \bar{n} \in  \Ga(V, \opn{H}^{i}(N))$.
Now $\d_M(m) = 0$, and hence, taking $l := 0$, the triple $(n, m, l)$
is a connecting triple for $\bar{n}$ over $V$. Therefore
$\de_{\bsym{E}}^i(\bar{n}) = \pi_L(l) = 0$.

\medskip \noindent 
$\triangleright$  $\triangleright$  Substep 1.b.
Now we prove that the morphism 
$\bar{\psi} : \opn{H}^{i}(M) \to \opn{Ker}(\de_{\bsym{E}}^i)$
is an epimorphism. 
We shall use Proposition \ref{prop:3635}(3). Let 
$\bar{n} \in \Ga(U, \opn{Ker}(\de_{\bsym{E}}^i))$; 
so $\bar{n} \in \Ga(U, \opn{H}^{i}(N))$ and
$\de_{\bsym{E}}^i(\bar{n}) = 0$. 
It suffices to find a covering $V' \surj U$ and a section 
$\bar{m} \in \Ga(V', \opn{H}^{i}(M))$ such that 
$\bar{\psi}(\bar{m}) = \bar{n}$ in $\Ga(V', \opn{H}^{i}(N))$. 

Take a covering $V \surj U$  that admits a connecting triple $(n, m, l)$ for 
$\bar{n}$. Then 
$\pi_L(l) = \de_{V}(\bar{n}) = \de_{\bsym{E}}^i(\bar{n}) = 0$
in $\Ga(V, \opn{H}^{i + 1}(L))$. 
By Proposition \ref{prop:3635}(3) there is some covering 
$V' \surj V$ and a section $l' \in \Ga(V', L^{i})$ s.t.\ 
$\d_L(l') = l$ in $\Ga(V', L^{i + 1})$. 
Define $m' := \phi(l') \in \Ga(V', M^{i})$. 
Now 
$\d_M(m - m') = \d_M(m) - \d_M(m') = \phi(l) - \phi(l) = 0$,
so $m - m'$ belongs to $\Ga(V', \opn{Z}^{i}(M))$. 
Because
$\psi(m - m') = \psi(m) = n$,
the cohomology class 
$\bar{m} := \pi_M(m - m')$ in $\Ga(V', \opn{H}^{i}(M))$
satisfies $\bar{\psi}(\bar{m}) = \bar{n}$
in $\Ga(V', \opn{H}^{i}(N))$.

\medskip \noindent 
$\triangleright$ Step 2. Exactness at $\opn{H}^{i + 1}(L)$. This is also done 
in two substeps. We will use the abbreviation 
$\bar{\phi} := \opn{H}^{i + 1}(\phi)$. 

\medskip \noindent 
$\triangleright$  $\triangleright$  Substep 2.a.
We start by proving that 
$\bar{\phi} \circ \de_{\bsym{E}}^i  = 0$.
Take $\bar{n} \in \Ga(U, \opn{H}^{i}(N))$, and let 
$\bar{l} := \de_{\bsym{E}}^i(\bar{n})$. We need to prove that 
$\bar{\phi}(\bar{l}) = 0$. 
Let $V \surj U$ be a covering that admits a connecting triple $(n, m, l)$ for 
$\bar{n}$. This says that $\bar{l} = \pi_L(l)$, and hence 
$\bar{\phi}(\bar{l}) = \pi_L(\phi(l)) = \pi_M(\d_M(m)) = 0$.

\medskip \noindent 
$\triangleright$  $\triangleright$  Substep 2.b.
Now we prove that 
$\de_{\bsym{E}}^i : \opn{H}^{i}(N) \to \opn{Ker}(\bar{\phi})$
is an epimorphism. 
Let $\bar{l} \in \Ga(U, \opn{H}^{i + 1}(L))$ be such that 
$\bar{\phi}(\bar{l}) = 0$. 
Take a covering $V \surj U$ for which there's a section 
$l \in \Ga(V, \opn{Z}^{i + 1}(L))$
s.t.\ $\pi_L(l) = \bar{l}$. The section 
$\phi(l) \in \Ga(V, \opn{Z}^{i + 1}(M))$
satisfies 
$\pi_M(\phi(l)) = \bar{\phi}(\bar{l}) = 0$
in $\Ga(V, \opn{H}^{i + 1}(M))$, and therefore there is some covering 
$V' \surj V$ and a section 
$m \in \Ga(V', M^{i})$ s.t.\ 
$\d_M(m) = \phi(l)$ in $\Ga(V', M^{i + 1})$.

Define $n := \psi(m) \in \Ga(V', N^{i})$. Then 
$\d_N(n) = \d_N(\psi(m)) = \psi(\d_M(m)) = \psi(\phi(l)) = 0$,
so in fact $n \in \Ga(V', \opn{Z}^{i}(N))$. 
Let $\bar{n} := \pi_N(n) \in \Ga(V', \opn{H}^{i}(N))$.
We see that $(n, m, l)$ is a connecting triple for $\bar{n}$ over $V'$. 
Therefore 
$\bar{l} = \pi_L(l) = \de_{V'}(\bar{n}) = \de_{\bsym{E}}^i(\bar{n})$
in $\Ga(V', \opn{H}^{i + 1}(L))$. 

\medskip \noindent 
$\triangleright$ Step 3. Exactness at $\opn{H}^{i}(M)$: 
let's use the abbreviations $\bar{\phi} := \opn{H}^{i}(\phi)$ and 
$\bar{\psi} := \opn{H}^{i}(\psi)$. 
It is clear that $\bar{\psi} \circ \bar{\phi} = 0$. 
It remains to prove that 
$\bar{\phi} : \opn{H}^{i}(L) \to \opn{Ker}(\bar{\psi})$
is an epimorphism. 

Let $\bar{m} \in \Ga(U, \opn{H}^{i}(M))$ be such that 
$\bar{\psi}(\bar{m}) = 0$. 
We are going to find a covering $W \to U$ and a section 
$\bar{l} \in \Ga(W, \opn{H}^{i}(L))$
such that $\bar{\phi}(\bar{l}) = \bar{m}$ in 
$\Ga(W, \opn{H}^{i}(M))$.

Take a covering $V \surj U$ for which there's a section 
$m \in \Ga(V, \opn{Z}^{i}(M))$
s.t.\ $\pi_M(m) = \bar{m}$. 
The section $n := \psi(m) \in \Ga(V, \opn{Z}^{i}(N))$
satisfies
$\pi_N(n) = \bar{\psi}(\bar{m}) = 0$
in $\Ga(V, \opn{H}^{i}(N))$, and therefore there is some covering 
$V' \surj V$ and a section 
$n' \in \Ga(V', N^{i - 1})$ s.t.\ 
$\d(n') = n$ in $\Ga(V', N^{i})$.

Take a covering $V'' \surj V'$ for which there is a section 
$m'' \in \Ga(V'', M^{i - 1})$ s.t.\ 
$\psi(m'') = n'$. The section 
$m - \d_M(m'') \in \Ga(V'', \opn{Z}^{i}(M))$
satisfies $\pi_M(m - \d(m'')) = \bar{m}$ in 
$\Ga(V'', \opn{H}^{i}(M)$, and also 
$\psi(m - \d_M(m'')) = \psi(m) - \psi(\d_M(m'')) = n - n = 0$.

There is a covering $W \surj V''$ with a section 
$l \in \Ga(W, L^i)$ s.t.\ $\phi(l) = m - \d_M(m'')$. 
Because $\phi$ is a monomorphism, and $m - \d_M(m'')$ is a cocycle, it follows 
that the section $l$ belongs to 
$\Ga(W, \opn{Z}^{i}(L))$. Its cohomology class 
$\bar{l} := \pi_L(l)$ in $\Ga(W, \opn{H}^{i}(L))$
satisfies 
$\bar{\phi}(\bar{l}) = \pi_M(\phi(l)) = \pi_M(m - \d_M(m'')) = \bar{m}$
in $\Ga(W, \opn{H}^{i}(M))$.  
\end{proof}

\begin{prop} \label{prop:3640}
Let $\chi : \bsym{E} \to \bsym{E}'$ be a morphism of short exact sequences in 
$\dcat{C}_{\mrm{str}}(\cat{M})$. Namely 
$\chi = (\chi_L, \chi_M, \chi_N)$ in the commutative diagram with exact rows
\[ \UseTips \xymatrix @C=8ex @R=6ex {
0
\ar[r]
&
L
\ar[r]^{\phi}
\ar[d]_{\chi_L}
&
M
\ar[r]^{\psi}
\ar[d]_{\chi_M}
&
N
\ar[r]
\ar[d]_{\chi_N}
&
0
\\
0
\ar[r]
&
L'
\ar[r]^{\phi'}
&
M'
\ar[r]^{\psi'}
&
N'
\ar[r]
&
0
} \]
in $\dcat{C}_{\mrm{str}}(\cat{M})$. Then for every $i$ the diagram 
\[ \UseTips \xymatrix @C=8ex @R=6ex {
\opn{H}^{i}(N) 
\ar[r]^{\de_{\bsym{E}}^i}
\ar[d]_{\opn{H}^{i}(\chi_N)}
&
\opn{H}^{i + 1}(L)
\ar[d]^{\opn{H}^{i + 1}(\chi_L)}
\\
\opn{H}^{i}(N') 
\ar[r]^{\de_{\bsym{E}'}^i}
&
\opn{H}^{i + 1}(L')
} \]
is commutative. 
\end{prop}

\begin{exer} \label{exer:3640}
Prove Proposition \ref{prop:3640}. (Hint: study the proof of the previous 
theorem.)
\end{exer}

\mysubsection{The DG Category 
\texorpdfstring{$\dcat{C}(A, \cat{M})$} {C(A,M)}} 
\label{subsec:DGModinM}

We now combine material from previous subsections. The concept introduced in 
the definition below is new. It is the DG version of Definition 
\ref{dfn:1702}. 

Conventions \ref{conv:2490}, \ref{conv:4626}, \ref{conv:4627} and 
\ref{conv:4628} are in place. Thus all linear structures and operations 
are $\K$-linear, and rings and bimodules are central over $\K$. 

Recall that given an  abelian category $\cat{M}$, the category of 
complexes $\dcat{C}(\cat{M})$ is a DG category. For a complex 
$M \in \dcat{C}(\cat{M})$ we have the DG ring 
$\opn{End}_{\cat{M}}(M) = \opn{Hom}_{\cat{M}}(M, M)$.
The multiplication in this ring is composition. 

\begin{dfn} \label{dfn:1130}
Let $\cat{M}$ be an abelian category, and let $A$ be a DG ring. A {\em DG 
$A$-module in  $\cat{M}$}
\index{Abelian category! DG $A$-module in}
is an object 
$M \in \dcat{C}(\cat{M})$, together with a DG ring homomorphism 
$f : A \to \opn{End}_{\cat{M}}(M)$. 
\end{dfn}

If $M$ is a DG $A$-module in $\cat{M}$, then after forgetting the 
differentials, $M$ becomes a graded $A$-module in $\cat{M}$.

\begin{dfn} \label{dfn:1131}
\index{Differential graded module! Hom between {\indash}s}
Let $\cat{M}$ be an abelian category, let $A$ be a DG ring, and let $M, N$ 
be DG $A$-modules in $\cat{M}$. In Definition \ref{dfn:1703} we introduced the 
graded $\K$-module $\opn{Hom}_{A, \cat{M}}(M, N)$. This is made into a DG 
$\K$-module with differential 
\[ \d(\phi) := \d_N \circ \phi - (-1)^i \cd \phi \circ \d_M \]
for $\phi \in \opn{Hom}_{A, \cat{M}}(M, N)^i$. 
\end{dfn}

When we have to be specific, we denote the differential of 
$\opn{Hom}_{A, \cat{M}}(M, N)$ by 
$\d_{\mrm{Hom}}$, $\d_{A, \cat{M}}$, or some similar expression. 

As we have seen before (in Lemmas \ref{1213} and \ref{lem:1211}), given 
morphisms
\[ \phi_k \in \opn{Hom}_{A, \cat{M}}(M_{k}, M_{k + 1})^{i_k} \]
for $k \in \{ 0,  1 \}$, we have  
\[ \phi_1 \circ \phi_0 \in 
\opn{Hom}_{A, \cat{M}}(M_0, M_2)^{i_0 + i_1} , \]
and 
\[ \d(\phi_1 \circ \phi_0) = 
\d(\phi_1) \circ \phi_0 + (-1)^{i_1} \cd \phi_1 \circ \d(\phi_0) . \]
Also the identity automorphism $\opn{id}_M$ belongs to 
$\opn{Hom}_{A, \cat{M}}(M, M)^{0}$, and $\d(\opn{id}_M) = 0$.
Therefore the next definition is legitimate. 

\begin{dfn} \label{dfn:1213}
\index{Differential graded category! of DG $A$-modules in $\cat{M}$}
Let $\cat{M}$ be an abelian category, and let $A$ be a DG ring. The 
DG category of DG $A$-modules in $\cat{M}$ is 
denoted by $\dcat{C}(A, \cat{M})$%
\index{1-C(A,M)@$\dcat{C}(A, \cat{M})$}.
The morphism DG $\K$-modules are
\[ \opn{Hom}_{\dcat{C}(A, \cat{M})}(M_0, M_1) :=
\opn{Hom}_{A, \cat{M}}(M_0, M_1) \]
from Definition \ref{dfn:1131}.
The composition is that of $\dcat{C}(\cat{M})$.
\end{dfn}

Notice that forgetting the action of $A$ is a faithful DG
functor $\dcat{C}(A, \cat{M}) \to \dcat{C}(\cat{M})$.

\begin{exa} \label{exa:1820}
If $A = \K$, then $\dcat{C}(A, \cat{M}) = \dcat{C}(\cat{M})$; 
and if $\cat{M} = \cat{Mod} \K$, then 
$\dcat{C}(A, \cat{M}) = \dcat{C}(A) = \cat{DGMod} A$.
\end{exa}

\begin{dfn} \label{dfn:1132}
In the situation of Definition \ref{dfn:1213}:
\begin{enumerate}
\item The strict category of $\dcat{C}(A, \cat{M})$%
\index{Strict category! of DG $A$-modules in $\cat{M}$}
(see Definition \ref{dfn:1212}(1)) is denoted by \lb 
$\dcat{C}_{\mrm{str}}(A, \cat{M})$%
\index{1-Cstr(A,M)@$\dcat{C}_{\mrm{str}}(A, \cat{M})$}.

\item The homotopy category of $\dcat{C}(A, \cat{M})$%
\index{Homotopy category! of DG $A$-modules in $\cat{M}$}%
\index{1-K(A,M)@$\dcat{K}(A, \cat{M})$}
(see Definition \ref{dfn:1212}(2)) is denoted by 
$\dcat{K}(A, \cat{M})$.
\end{enumerate}
\end{dfn}

To say this in words: a morphism $\phi : M \to N$ in 
$\dcat{C}(A, \cat{M})$ is strict if and only if it has degree $0$ and 
$\phi \circ \d_M = \d_N \circ \phi$. 
The morphisms in $\dcat{K}(A, \cat{M})$ are the homotopy classes of strict 
morphisms. 

An additive functor $F : \cat{M} \to \cat{N}$ between abelian categories is 
called {\em faithfully exact} if for every sequence $\bsym{E}$ in 
$\dcat{C}_{\mrm{str}}(A, \cat{M})$, the sequence $\bsym{E}$ is exact if and 
only if the sequence $F(\bsym{E})$ in $\cat{N}$ is exact.

\begin{prop} \label{prop:4130}
The category $\dcat{C}_{\mrm{str}}(A, \cat{M})$ is a $\K$-linear abelian 
category, and the forgetful functors 
\[ \dcat{C}_{\mrm{str}}(A, \cat{M}) \xar{} \dcat{C}_{\mrm{str}}(\cat{M}) 
\xar{\opn{Und}} \dcat{G}_{\mrm{str}}(\cat{M}) \]
are faithfully exact. 
\end{prop}

\begin{exer} \label{exer:4145}
Prove the proposition above. 
\end{exer}

Like Definition \ref{dfn:2993}, given a DG module 
$M \in \dcat{C}(A, \cat{M})$ and an 
integer $i$, we can consider the objects of degree $i$ cocycles 
$\opn{Z}^i(M)$,  decocycles $\opn{Y}^i(M)$, coboundaries $\opn{B}^i(M)$ and 
cohomology $\opn{H}^i(M)$. These are all objects of $\cat{M}$. As $M$ varies we 
get linear functors 
\begin{equation} \label{eqn:4145}
\opn{Z}, \opn{Y}, \opn{B}, \opn{H} : \dcat{C}_{\mrm{str}}(A, \cat{M}) \to 
\dcat{G}_{\mrm{str}}(\cat{M}) . 
\end{equation}
These objects are related by the following exact sequences in $\cat{M}$~:
\begin{equation} \label{eqn:4148}
0 \to \opn{Z}^i(M) \to M^i \xar{\d^i_M} M^{i + 1} , 
\end{equation}
\begin{equation} \label{eqn:4149}
M^{i - 1} \xar{\d^{i - 1}_M} M^{} \to \opn{Y}^i(M) \to 0 , 
\end{equation}
\begin{equation} \label{eqn:4146}
0 \to \opn{B}^i(M) \to \opn{Z}^i(M) \to \opn{H}^i(M) \to 0 , 
\end{equation}
\begin{equation} \label{eqn:4147}
0 \to \opn{H}^i(M) \to \opn{Y}^i(M) \to \opn{B}^{i + 1}(M) \to 0 . 
\end{equation}

\begin{rem} \label{rem:3980}
Actually, the cohomology functor $\opn{H}$ from equation (\ref{eqn:4145})
factors through the category $\dcat{G}_{\mrm{str}}(\opn{H}(A),  \cat{M})$. 
See formula (\ref{eqn:2360}). We will not need this fact. 
\end{rem}

\begin{prop} \label{prop:4145}
The functor $\opn{Z}$ in equation \tup{(\ref{eqn:4145})}
is left exact, and the functor
$\opn{Y}$ in equation \tup{(\ref{eqn:4145})}
is right exact.
\end{prop}

\begin{proof}
It is enough to consider each degree separately. And by Proposition 
\ref{prop:4130} we may ignore the DG ring $A$.
Let $F^i : \dcat{G}_{\mrm{str}}(\cat{M}) \to \cat{M}$
be the functor that sends a graded module $M$ to its degree $i$ component 
$M^i$. This is an exact functor. 

Formula (\ref{eqn:4148}) exhibits $\opn{Z}^i$ as the kernel of the homomorphism 
of functors $\d^i : F^i \to F^{i + 1}$. According to Proposition 
\ref{prop:3620}(2) the functor $\opn{Z}^i$ is left exact. 

Similarly, formula (\ref{eqn:4149}) exhibits $\opn{Y}^i$ as the cokernel of the 
homomorphism of functors 
$\d^{i - 1} : F^{i - 1} \to F^{i}$. According to Proposition 
\ref{prop:3620}(2) the functor $\opn{Y}^i$ is right exact. 
\end{proof}

\begin{prop} \label{prop:3140}
Let $\{ M_x \}_{x \in X}$ be a collection of objects of
$\dcat{C}(A, \cat{M})$, indexed by a set $X$. Assume that for every $i \in \Z$ 
the direct sum $M^i := \bigoplus_{x \in X} M^i_x$ 
exists in $\cat{M}$. Then the graded object 
$M := \{ M^i \}_{i \in \Z} \in \dcat{G}(\cat{M})$ 
has a canonical structure of DG $A$-module in $\cat{M}$, 
and $M$ is a direct sum of the collection $\{ M_x \}_{x \in X}$
in $\dcat{C}_{\mrm{str}}(A, \cat{M})$. 
\end{prop}

\begin{exer} \label{exer:3140}
Prove Proposition \ref{prop:3140}. 
(Hint: look at the proof of Theorem \ref{thm:2495}.)
\end{exer}

\begin{exa} \label{exa:2367}
Since $\dcat{M}(\K)$ has infinite direct sums, the proposition above shows that 
$\dcat{C}_{\mrm{str}}(A)$ has infinite direct sums. 
\end{exa}

\begin{prop} \label{prop:2365}
Given a DG ring $A$, let $A^{\natural}$ be the graded ring gotten by forgetting 
the differential. Consider the forgetful functor 
$\opn{Und} : \dcat{C}(A, \cat{M}) \to \dcat{G}(A^{\natural}, \cat{M})$
that forgets the differentials of DG modules, i.e.\ 
$\opn{Und}(M, \d_M) :=  M$.
\begin{enumerate}
\item $\opn{Und}$ is a fully faithful graded functor.

\item On the strict categories, 
$\opn{Str}(\opn{Und}) :  \dcat{C}_{\mrm{str}}(A, \cat{M}) \to 
\dcat{G}_{\mrm{str}}(A, \cat{M})$
is a faithfully exact additive functor.
\end{enumerate} 
\end{prop}

\begin{exer} \label{exer:2367}
Prove Proposition \ref{prop:2365}.
\end{exer}

A graded object $M = \{ M^i \}_{i \in \Z}$ in $\cat{M}$ is said to be {\em
bounded above} if the set $\{ i \mid M^i \neq 0 \}$ is bounded above. 
Likewise we define {\em bounded below} and {\em bounded} graded objects. 

\begin{dfn} \label{dfn:3220}
\index{1-C(A,M)@$\dcat{C}^{\star}(A, \cat{M})$}
\index{Boundedness! condition}
We define $\dcat{C}^-(A, \cat{M})$, $\dcat{C}^+(A, \cat{M})$ and 
$\dcat{C}^{\mrm{b}}(A, \cat{M})$ to be the full subcategories of 
$\dcat{C}(A, \cat{M})$ consisting of bounded above, bounded below and bounded
DG modules respectively. 
\end{dfn}

The symbols ``$-$'', ``$+$'' and ``$\mrm{b}$'' and 
``$\bra{\mrm{empty}}$'' (the empty symbol) are called {\em boundedness 
indicators}. We shall usually use the symbol ``$\star$'' to denote an 
unspecified boundedness indicator. Thus the expression 
$\dcat{C}^{\star}(A, \cat{M})$ can refer to any of the four full subcategories 
of $\dcat{C}(A, \cat{M})$ that are mentioned in the definition above.

\begin{rem} \label{rem:1245}
Here is a generalization of Definition \ref{dfn:1213}. Instead of a DG central
$\K$-ring $A$ we can take a small $\K$-linear DG category $\cat{A}$. 
We then define the $\K$-linear DG category 
$\dcat{C}(\cat{A}, \cat{M}) := \cat{DGFun}(\cat{A}, \dcat{C}(\cat{M}))$
as in Exercise \ref{exer:1250}.

This is indeed a generalization of Definition \ref{dfn:1213}: when $\cat{A}$ 
has a single object $x$, and we write $A := \opn{End}_{\cat{A}}(x)$, then 
the functor $M \mapsto M(x)$ is an isomorphism of DG categories
$\dcat{C}(\cat{A}, \cat{M}) \iso \dcat{C}(A, \cat{M})$. 

In the special case of $\cat{M} = \dcat{M}(\K) = \cat{Mod} \K$, the DG category 
$\dcat{C}(\cat{A}, \cat{M})$ is what Keller \cite{Kel1} calls the DG category 
of {\em DG $\cat{A}$-modules}. 

Almost everything we do in this book for $\dcat{C}(A, \cat{M})$
holds in the more general context of $\dcat{C}(\cat{A}, \cat{M})$.
However, in the more general context a lot of the intuition is lost, and some 
aspects become pretty cumbersome. This is the reason we have decided to stick 
with the less general context. 
\end{rem}

\begin{rem} \label{rem:4960}
The theory of DG categories has a generalization -- this is the theory of {\em 
$\mrm{A}_{\infty}$ categories}%
\index{Ainfty@$\mrm{A}_{\infty}$ category},
which will be outlined in this remark.
All set theoretical issues are going to be ignored.

The concept of an {\em $\mrm{A}_{\infty}$ algebra}, namely a single-object 
$\mrm{A}_{\infty}$ category, was introduced by J. Stasheff \cite{Sta1}, 
\cite{Sta2} in the 1960's, in relation with his work on {\em higher 
operations on H-spaces} in topology. 

$\mrm{A}_{\infty}$ categories are at the core of  {\em Homological Mirror 
Symmetry}, conceived by K. Fukaya and M. Kontsevich during the 1990's.
Indeed, a symplectic manifold $Y$ gives rise to the {\em Fukaya category} 
$\cat{F}(Y)$, which is an $\mrm{A}_{\infty}$ category, and is a categorification 
of {\em Floer cohomology}. See \cite{Fuk}, \cite{Kon}, \cite{KoSo}, \cite{Seid} 
and \cite{HMS}. 

Further developments on $\mrm{A}_{\infty}$ algebras and categories were made, 
among others, by T.V. Kadeishvili \cite{Kad}, K. Lef\`evre-Hasegawa \cite{Lef}, 
B. Keller \cite{Kel6}, \cite{Kel5}, and Canonaco-Ornaghi-Stellari \cite{COS}. 
For the analogous theory of {\em $\mrm{L}_{\infty}$ algebras},
generalizing DG Lie algebras, see \cite{Ye5.5}.

We shall only discuss strictly unital $\mrm{A}_{\infty}$ categories, giving 
very few details (full details can be found in \cite{Kel5} and \cite{COS}), and 
our notation deviates somewhat from the references
(because we want it to be compatible with the rest of the book). 
The base ring $\K$ is assumed to be a field. 
An $\mrm{A}_{\infty}$ category $\cat{A}$ over $\K$
has a set of objects $\opn{Ob}(\cat{A})$. For each pair of objects $x, y \in 
\opn{Ob}(\cat{A})$ there is a graded $\K$-module 
$\cat{A}(x, y) = \opn{Hom}_{\cat{A}}(x, y)$. 
For each $x \in \opn{Ob}(\cat{A})$ there is a degree 
$0$ element $\opn{id}_x \in \cat{A}(x, x)$.
For each $i \geq 1$ and each sequence 
$x_0, \ldots, x_i \in \opn{Ob}(\cat{A})$, there is a 
degree $2 - i$ homomorphism 
\begin{equation} \label{4960}
\opn{m}_i : \cat{A}(x_{i - 1}, x_{i}) \ot \cdots \ot  
\cat{A}(x_{0}, x_{1}) \to \cat{A}(x_{0}, x_{i}) 
\end{equation}
$\dcat{G}(\K)$. 
The operations $\opn{m}_i$ satisfy certain relations. 
Roughly speaking, the operation $\opn{m}_1$ is a differential, the operation 
$\opn{m}_2$ is a multiplication, the element $\opn{id}_x$ is a unit for the 
multiplication $\opn{m}_2$, the multiplication $\opn{m}_2$ is associative 
up to the hotomopy $\opn{m}_3$, etc.
The relations imply that the cohomology $\opn{H}(\cat{A})$, w.r.t.\ the 
differential $\opn{m}_1$, is a graded category, in the sense of 
Definition \ref{dfn:1692}.
If the higher (i.e.\ $i > 2$) operations $\opn{m}_i$ vanish, then $\cat{A}$ 
is a DG category, in the sense 
of Definition \ref{dfn:1063}.

Suppose $\cat{A}$ and $\cat{B}$ are $\mrm{A}_{\infty}$ categories. 
An {\em $\mrm{A}_{\infty}$ functor} $F : \cat{A} \to \cat{B}$
consists of a function
$F_{\mrm{ob}} : \opn{Ob}(\cat{A}) \to \opn{Ob}(\cat{B})$, 
and for each $i \geq 1$ and each sequence 
$x_0, \ldots, x_i \in \opn{Ob}(\cat{A})$, a degree $1 - i$ homomorphism
\begin{equation} \label{4961}
F_i : \cat{A}(x_{i - 1}, x_{i}) \ot \cdots \ot  \cat{A}(x_{0}, x_{1}) \to 
\cat{B} \bigl( F_{\mrm{ob}}(x_{0}), F_{\mrm{ob}}(x_{i}) \bigr) 
\end{equation}
in $\dcat{G}(\K)$. 
There are relations that the $F_i$ have to satisfy.
Given another $\mrm{A}_{\infty}$ functor $G : \cat{B} \to \cat{C}$, 
there is a composed $\mrm{A}_{\infty}$ functor
$G \circ F : \cat{A} \to \cat{C}$, which is hard to describe directly (it is 
usually defined in terms of {\em cocategories}, via the {\em bar 
construction}). 
We denote by $\cat{A_{\infty}Cat}$ the category whose objects are the 
$\mrm{A}_{\infty}$ categories and whose morphisms are the $\mrm{A}_{\infty}$ 
functors.  

For a pair $\cat{A}, \cat{B} \in \opn{Ob}(\cat{A_{\infty}Cat})$
we denote by  $\cat{A_{\infty}Fun}(\cat{A}, \cat{B})$ 
the set of $\mrm{A}_{\infty}$ functors from $\cat{A}$ to $\cat{B}$, so 
that \begin{equation} \label{eqn:4970}
\opn{Hom}_{\cat{A_{\infty}Cat}}(\cat{A}, \cat{B}) =
\cat{A_{\infty}Fun}(\cat{A}, \cat{B}) . 
\end{equation}
An important fact is that 
$\cat{A_{\infty}Fun}(\cat{A}, \cat{B})$
is itself the set of objects of an $\mrm{A}_{\infty}$ category, i.e.\ 
$\cat{A_{\infty}Fun}(\cat{A}, \cat{B}) \in \opn{Ob}(\cat{A_{\infty}Cat})$. 
If $\cat{B}$ is a DG category (i.e.\ its higher operations vanish), then 
$\cat{A_{\infty}Fun}(\cat{A}, \cat{B})$
is also a DG category.

We shall discuss more aspects of $\mrm{A}_{\infty}$ categories in Remarks 
\ref{rem:4970} and \ref{rem:4971}. 
\end{rem}

\mysubsection{Contravariant DG Functors} \label{subsec:contrvar-dg-func}

In this subsection we address, in a systematic fashion, the issue of reversing 
arrows in DG categories. As always, we work over a commutative base ring $\K$. 
 
\begin{dfn} \label{dfn:2495}
Let $\cat{C}$ and $\cat{D}$ be DG categories.
A {\em contravariant DG functor}%
\index{Differential graded functor! contravariant}
$F : \cat{C} \to \cat{D}$ consists of a function 
$F : \opn{Ob}(\cat{C}) \to \opn{Ob}(\cat{D})$,
and for each pair of objects $M_0, M_1 \in \opn{Ob}(\cat{C})$ a homomorphism 
\[ F : \opn{Hom}_{\cat{C}}(M_0, M_1) \to 
\opn{Hom}_{\cat{D}} \bigl( F(M_1), F(M_0) \bigr) \]
in $\dcat{C}_{\mrm{str}}(\K)$. The conditions are:
\begin{itemize}
\rmitem{a} Units: $F(\opn{id}_M) = \opn{id}_{F(M)}$. 

\rmitem{b} Graded reversed composition: suppose that for $k \in \{ 0, 1 \}$
we are given morphisms 
$\phi_k \in \opn{Hom}_{\cat{C}}(M_{k}, M_{k + 1})^{i_k}$.
Then there is equality 
\[ F(\phi_1 \circ \phi_0) = (-1)^{i_0 \cd i_1} \cd F(\phi_0) \circ F(\phi_1) \]
inside 
$\opn{Hom}_{\cat{D}} \bigl( F(M_2), F(M_0) \bigr)^{i_0 + i_1}$.  
\end{itemize}
\end{dfn}

Warning: a contravariant DG functor does not remain a contravariant functor
after the grading is forgotten; cf.\ Remark \ref{rem:2992}. 

Here is the categorical version of Definition \ref{dfn:1080}.

\begin{dfn} \label{dfn:1272}
Let $\cat{C}$ be a DG category. The 
{\em opposite DG category}%
\index{Differential graded category! opposite category of}
$\cat{C}^{\mrm{op}}$ has the same set of objects. The morphism DG $\K$-modules 
are
\[ \opn{Hom}_{\cat{C}^{\mrm{op}}} (M_0, M_1) := 
\opn{Hom}_{\cat{C}} (M_1, M_0) . \]
The composition $\circ^{\mrm{op}}$ of $\cat{C}^{\mrm{op}}$ is the composition 
$\circ$ of $\cat{C}$, reversed and multiplied by signs: 
\[ \phi_0 \circ^{\mrm{op}} \phi_1 := 
(-1)^{i_0 \cd i_1} \cd \phi_1 \circ \phi_0  \]
for morphisms 
$\phi_k \in \opn{Hom}_{\cat{C}} (M_k, M_{k + 1})^{i_k}$.
\end{dfn}

One needs to verify that this is indeed a DG category. 
This is basically the same verification as in Exercise \ref{exer:2500}.

As before, we define the operation 
$\opn{Op} : \cat{C} \to \cat{C}^{\mrm{op}}$
to be the identity on objects, and the identity on morphisms in reversed 
order, i.e.\ 
\[ \opn{Op} = \opn{id} : \opn{Hom}_{\cat{C}} (M_0, M_1) \iso 
\opn{Hom}_{\cat{C}^{\mrm{op}}} (M_1, M_0) . \]
Note that $(\cat{C}^{\mrm{op}})^{\mrm{op}} = \cat{C}$, 
and we denote the inverse operation $\cat{C}^{\mrm{op}} \to \cat{C}$
also by $\opn{Op}$. 

\begin{prop} \label{prop:2500}
Let $\cat{C}$, $\cat{D}$ and $\cat{E}$ be DG categories.
\begin{enumerate}
\item The operations $\opn{Op} : \cat{C} \to \cat{C}^{\mrm{op}}$
and $\opn{Op} : \cat{C}^{\mrm{op}} \to \cat{C}$
are contravariant DG functors. 

\item If $F : \cat{C} \to \cat{D}$ is a contravariant DG functor, 
then the composition $F \circ \opn{Op} : \cat{C}^{\mrm{op}} \to \cat{D}$
is a DG functor; and vice versa. 

\item If $F : \cat{C} \to \cat{D}$ and $G : \cat{D} \to \cat{E}$
are contravariant DG functors, 
then the composition $G \circ F : \cat{C} \to \cat{E}$
is a  DG functor. 
\end{enumerate}
\end{prop}

\begin{exer}  \label{exer:2495}
Prove the previous proposition. 
\end{exer}

Definitions \ref{dfn:1272} and \ref{dfn:2495} make sense for graded categories, 
by forgetting differentials. Thus for graded categories $\cat{C}$ and 
$\cat{D}$ we can talk about contravariant graded functors
$\cat{C} \to \cat{D}$, and about the graded category $\cat{C}^{\mrm{op}}$.

We already met $\dcat{G}(\cat{M})$, the category of graded objects in an
abelian category $\cat{M}$; see Definition \ref{dfn:1701}. It is a 
graded category. 
Its objects are collections $M = \{ M^i \}_{i \in \Z}$ of objects 
$M^i \in \cat{M}$. 

Let $\cat{M}$ and $\cat{N}$ be abelian categories, and let 
$F : \cat{M} \to \cat{N}$ be a contravariant linear functor. 
For a graded object 
$M = \{ M^i \}_{i \in \Z} \in \dcat{G}(\cat{M})$
let us define the graded object 
\begin{equation} \label{eqn:2496}
\dcat{G}(F)(M) := \{ N^i \}_{i \in \Z} \in \dcat{G}(\cat{N}) , \quad 
N^i := F(M^{-i}) \in \cat{N} .
\end{equation}

Next consider a pair of objects $M_0, M_1 \in \dcat{G}(\cat{M})$
and a degree $i$ morphism 
$\phi : M_0 \to M_1$ in $\dcat{G}(\cat{M})$. Thus 
\[ \phi = \{ \phi^{\lms j} \}_{j \in \Z} \in 
\opn{Hom}_{\dcat{G}(\cat{M})}(M_0, M_1)^i , \]
where, as in formula (\ref{eqn:2495}), the morphism $\phi^{\lms j}$ belongs to 
$\opn{Hom}_{\cat{M}}(M_0^j, M_1^{j + i})$.
We have objects 
$N_k := \dcat{G}(F)(M_k) \in \dcat{G}(\cat{N})$,
for $k \in \{ 0, 1 \}$, defined by (\ref{eqn:2496}). Explicitly, 
$N_k = \{ N_k^i \}_{i \in \Z}$ and $N_k^i = F(M_k^{-i})$. 
For each $j \in \Z$ define the morphism 
\begin{equation} \label{eqn:2505}
\psi^{\lms j} := (-1)^{i \cd j} \cd F(\phi^{-j - i}) \in
\opn{Hom}_{\cat{N}} \bigl( N_1^j, N_0^{j + i} \bigr) .
\end{equation}
Collecting them we obtain a morphism 
\begin{equation} \label{eqn:2497}
\dcat{G}(F)(\phi) := \{ \psi^{\lms j} \}_{j \in \Z}  \in 
\opn{Hom}_{\dcat{G}(\cat{N})}(N_1, N_0)^i . 
\end{equation}

\begin{lem} \label{lem:2495}
The assignments \tup{(\ref{eqn:2496})} and \tup{(\ref{eqn:2497})}
produce a contravariant graded functor 
$\dcat{G}(F) : \dcat{G}(\cat{M}) \to \dcat{G}(\cat{N})$.
\end{lem}

\begin{proof}
Since for morphisms of degree $0$ there is no sign twist, the identity 
automorphism 
$\opn{id}_M = \{ \opn{id}_{M^i} \}_{i \in \Z}$
of $M = \{ M^i \}_{i \in \Z}$ in $\dcat{G}(\cat{M})$
is sent to the identity automorphism of $\dcat{G}(F)(M)$ in 
$\dcat{G}(\cat{N})$. 

Next we look at morphisms 
$\phi_k = \{ \phi_k^{\lms j} \}_{j \in \Z}$
in 
$\opn{Hom}_{\dcat{G}(\cat{M})}(M_k, M_{k + 1})^{i_k}$ 
for $k = 0, 1$. 
The composition $\phi_1 \circ \phi_0$ has degree $i_0 + i_1$, and the $j$-th 
component of $\phi_1 \circ \phi_0$ is 
$\phi_1^{\lms j + i_0} \circ \phi_0^{\lms j}$.
Therefore the $j$-th component of $\dcat{G}(F)(\phi_1 \circ \phi_0)$ 
is 
\begin{equation} \label{eqn:2510}
\begin{aligned}
& \dcat{G}(F)(\phi_1 \circ \phi_0)^{\lms j} 
= (-1)^{\lms j \cd (i_0 + i_1)} \cd 
F \bigl( \phi_1^{-j - i_1} \circ \phi_0^{-j - (i_0 + i_1)} \bigr)
\\
& \quad = (-1)^{\lms j \cd (i_0 + i_1)} \cd 
F(\phi_0^{-j - (i_0 + i_1)}) \circ F(\phi_1^{-j - i_1}) .
\end{aligned}
\end{equation}
The  $j$-th component of $\dcat{G}(F)(\phi_k)$ is 
$\dcat{G}(F)(\phi_k)^{\lms j} = 
(-1)^{\lms j \cd i_k} \cd F(\phi_k^{-j - i_k})$. 
So the $j$-th component of 
$(-1)^{i_0 \cd i_1} \cd \dcat{G}(F)(\phi_0) \circ \dcat{G}(F)(\phi_1)$ 
is 
\begin{equation} \label{eqn:2511}
\begin{aligned}
& (-1)^{i_0 \cd i_1} \cd 
\bigl( \dcat{G}(F)(\phi_0) \circ \dcat{G}(F)(\phi_1) \bigr)^j 
\\
& \quad = (-1)^{i_0 \cd i_1} \cd 
(-1)^{(j + i_1) \cd i_0} \cd F(\phi_0^{-(j + i_1) - i_0})
\circ (-1)^{\lms j \cd i_1} \cd F(\phi_1^{-j - i_1}) . 
\end{aligned}
\end{equation}
We see that the morphisms (\ref{eqn:2510}) and (\ref{eqn:2511}) are equal. 
\end{proof}

Now we consider a complex $(M, \d_M) \in \dcat{C}(\cat{M})$. 
This is made up of a graded object 
$M  = \{ M^i \}_{i \in \Z} \in \dcat{G}(\cat{M})$ together 
with a differential $\d_M = \{ \d_M^i \}_{i \in \Z}$, 
where $\d_M^i : M^i \to M^{i + 1}$. 
We can view $\d_M$ as an element of 
$\opn{End}_{\dcat{G}(\cat{M})}(M)^1 = 
\opn{Hom}_{\dcat{G}(\cat{M})}(M, M) ^1$.
We specify a differential 
$\d_{\dcat{C}(F)(M)}$ on the graded object 
$\dcat{G}(F)(M) \in \dcat{G}(\cat{N})$
as follows: 
\begin{equation} \label{eqn:2506}
\d_{\dcat{C}(F)(M)} := - \dcat{G}(F)(\d_M) \in 
\opn{End}_{\dcat{G}(\cat{N})} \bigl( \dcat{G}(F)(M) \bigr)^1 . 
\end{equation}
To be explicit, the component 
\[ \d_{\dcat{C}(F)(M)}^i : \dcat{G}(F)(M)^i = F(M^{-i}) \to 
F(M^{-i - 1}) = \dcat{G}(F)(M)^{i + 1} \]
of $\d_{\dcat{C}(F)(M)}$ is, by (\ref{eqn:2505}), 
\begin{equation} \label{eqn:2365}
\d_{\dcat{C}(F)(M)}^i = (-1)^{i + 1} \cd F(\d_M^{-i - 1}) :
F(M^{-i}) \to F(M^{-i - 1}) . 
\end{equation}
This shows that our formula coincides with the one in 
\cite[Remark 1.8.11]{KaSc1}. 

The full subcategory $\dcat{C}^{\star}(\cat{N}) \sub \dcat{C}^{}(\cat{N})$, for 
a given boundedness condition $\star$, was introduced in Definition 
\ref{dfn:3220}. Each boundedness condition $\star$ has a {\em reversed 
boundedness condition} $-\star$, with the obvious meaning; e.g.\ if 
$\star$ is $+$ then $-\star$ is $-$. 

\begin{thm} \label{thm:3015} 
Let $\cat{M}$ and $\cat{N}$ be abelian categories, and let 
$F : \cat{M} \to \cat{N}$ be a contravariant linear functor. 
The assignments \tup{(\ref{eqn:2496})}, \tup{(\ref{eqn:2497})} and 
\tup{(\ref{eqn:2506})} produce a contravariant DG functor 
\[ \dcat{C}(F) : \dcat{C}(\cat{M}) \to \dcat{C}(\cat{N}) . \]
For every boundedness condition $\star$ we have 
$\dcat{C}(F) \bigl( \dcat{C}^{\star}(\cat{M}) \bigr) \sub
\dcat{C}^{- \star}(\cat{N})$,
where $-\star$ is the reversed boundedness condition. 
\end{thm}

\begin{proof}
We already know, by Lemma \ref{lem:2495}, that $\dcat{G}(F)$ is a graded 
functor. We need to prove that for a pair of DG modules 
$(M_0, \d_{M_0})$ and $(M_1, \d_{M_1})$ in $\dcat{C}(\cat{M})$ 
the strict homomorphism of graded $\K$-modules 
\[ \dcat{G}(F) : \opn{Hom}_{\dcat{G}(\cat{M})}(M_0, M_1) \to 
\opn{Hom}_{\dcat{G}(\cat{N})} \bigl( \dcat{G}(F)(M_1), \dcat{G}(F)(M_0) \bigr) 
\]
respects differentials. Take any 
$\phi \in \opn{Hom}_{\dcat{G}(\cat{M})}(M_0, M_1)^i$.
By definition we have 
$\d(\phi) = \d_{M_1} \circ \phi - (-1)^i \cd \phi \circ \d_{M_0}$.
Using the fact that $\dcat{G}(F)$ is a contravariant graded functor, we obtain 
these equalities: 
\[ \begin{aligned}
& \dcat{G}(F)(\d(\phi)) 
\\
& \ = (-1)^{i} \cd \dcat{G}(F)(\phi) \circ \dcat{G}(F)(\d_{M_1}) 
- (-1)^i \cd (-1)^i \cd \dcat{G}(F)(\d_{M_0}) \circ \dcat{G}(F)(\phi)
\\
& \ = \d_{\dcat{C}(F)(M_0)} \circ \dcat{G}(F)(\phi)
- (-1)^{i} \cd \dcat{G}(F)(\phi) \circ \d_{\dcat{C}(F)(M_1)} 
= \d(\dcat{G}(F)(\phi)) . 
\end{aligned} \]

The claim about the boundedness conditions is immediate from equation 
(\ref{eqn:2496}). 
\end{proof}

The sign appearing in formula (\ref{eqn:2506}) might seem arbitrary. 
Besides being the only sign for which Theorem \ref{thm:3015} holds, 
there is another explanation, which can be seen in the next exercise. 

\begin{exer} \label{exer:2505}
Take $\cat{M} = \cat{N} := \cat{Mod} \K$, and consider the contravariant 
additive functor $F := \opn{Hom}_{\K}(-, \K)$ from $\cat{M}$ to itself.
Let $M \in \dcat{C}(\cat{M})$; we can view $M$ as a complex of $\K$-modules or 
as a DG $\K$-module, as done in Proposition \ref{prop:1065}. 
Show that 
$\dcat{C}(F)(M) \cong \opn{Hom}_{\K}(M, \K)$
in $\dcat{C}_{\mrm{str}}(\K)$, where the second object is the graded module 
from formula (\ref{eqn:1700}), with the differential 
$\d$ from Definition \ref{dfn:2992}(2). 
\end{exer}

The next definition and theorem will help us later when studying contravariant 
triangulated functors.

\begin{dfn} \label{dfn:3015}
Let $A$ be a DG ring and let $\cat{M}$ be an abelian 
category. The {\em flipped category}
\index{Flipped category}
of $\dcat{C}(A, \cat{M})$ is the DG category 
$\dcat{C}(A, \cat{M})^{\mrm{flip}} := 
\dcat{C}(A^{\mrm{op}}, \cat{M}^{\mrm{op}})$. 
\end{dfn}

\begin{thm} \label{thm:2495}
Let $A$ be a DG ring and let $\cat{M}$ be an abelian category. Then\tup{:}
\begin{enumerate}
\item There is a canonical isomorphism of DG categories  
\[ \opn{Flip} : 
\dcat{C}(A, \cat{M})^{\mrm{op}} \iso 
\dcat{C}(A^{\mrm{op}}, \cat{M}^{\mrm{op}}) = \dcat{C}(A, \cat{M})^{\mrm{flip}} . 
\]

\item For every boundedness condition $\star$ we have 
\[ \opn{Flip} \bigl( \dcat{C}^{\star}(A, \cat{M})^{\mrm{op}} \bigr) = 
\dcat{C}^{-\star}(A^{\mrm{op}}, \cat{M}^{\mrm{op}}) , \]
where $-\star$ is the reversed boundedness condition. 

\item The induced isomorphism on the strict categories
\[ \opn{Str}(\opn{Flip}) : \dcat{C}_{\mrm{str}}(A, \cat{M})^{\mrm{op}} \iso 
\dcat{C}_{\mrm{str}}(A^{\mrm{op}}, \cat{M}^{\mrm{op}}) \]
is an exact functor, for the respective abelian category structures.

\item For every integer $i$, the functor $\opn{Str}(\opn{Flip})$ makes the 
diagram 
\[ \UseTips  \xymatrix @C=12ex @R=6ex {
\dcat{C}_{\mrm{str}}(A, \cat{M})^{\mrm{op}}
\ar[r]^{\opn{Str}(\opn{Flip})}_{\cong}
\ar[dr]_{\opn{H}^i}
&
\dcat{C}_{\mrm{str}}(A^{\mrm{op}}, \cat{M}^{\mrm{op}})
\ar[d]^{\opn{H}^{-i}}
\\
&
\cat{M}^{\mrm{op}}
} \]
commutative, up to an isomorphism of linear functors. 
\end{enumerate}
\end{thm}

\begin{proof} \mbox{}

\smallskip \noindent 
(1) According to Proposition \ref{prop:2500} there is a contravariant DG 
functor 
$\opn{Op} : \lb \dcat{C}(A, \cat{M})^{\mrm{op}} \to \dcat{C}(A, \cat{M})$.
It is bijective on objects and morphisms. We are going to construct a
contravariant DG functor
$E : \dcat{C}(A, \cat{M}) \to \dcat{C}(A^{\mrm{op}}, \cat{M}^{\mrm{op}})$
which is also bijective on objects and morphisms. The composed DG functor 
\[ \opn{Flip} := E \circ \opn{Op} : \dcat{C}(A, \cat{M})^{\mrm{op}} \to 
\dcat{C}(A^{\mrm{op}}, \cat{M}^{\mrm{op}}) \]
will have the desired properties. 

Let us construct $E$. We start with the contravariant additive functor 
$F := \opn{Op} : \cat{M} \to \cat{M}^{\mrm{op}}$.
Theorem \ref{thm:3015}   says that 
$\dcat{C}(F) : \dcat{C}(\cat{M}) \to \dcat{C}(\cat{M}^{\mrm{op}})$
is a contravariant DG functor. 
Recall that an object of $\dcat{C}(A, \cat{M})$ is a triple 
$(M, \d_M, f_M)$, where $M \in \dcat{G}(\cat{M})$; $\d_M$ is a differential on 
the graded object $M$; and 
$f_M : A \to \opn{End}_{\dcat{C}(\cat{M})}(M)$
is a DG ring homomorphism. See Definitions \ref{dfn:1703}, \ref{dfn:1130} and 
\ref{dfn:1213}.
Define 
$(N, \d_N) := \dcat{C}(F)(M, \d_M) \in \dcat{C}(\cat{M}^{\mrm{op}})$.
Since 
$\dcat{C}(F) : \opn{End}_{\dcat{C}(\cat{M})}(M, \d_M) \to 
\opn{End}_{\dcat{C}(\cat{M}^{\mrm{op}})}(N, \d_N)$
is a DG ring anti-homomorphism (by which we mean the single object version of a 
contravariant DG functor), and 
$\opn{Op} : A^{\mrm{op}} \to A$ 
is also such an anti-homomorphism, it follows that 
\[ f_N := \dcat{C}(F) \circ f_M \circ \opn{Op} : A^{\mrm{op}} \to 
\opn{End}_{\dcat{C}(\cat{M}^{\mrm{op}})}(N, \d_N) \]
is a DG ring homomorphism. Thus 
$E(M, \d_M, f_M) := (N, \d_N, f_N)$ 
is an object of $\dcat{C}(A^{\mrm{op}}, \cat{M}^{\mrm{op}})$.
In this way we have a function 
\[ E : \opn{Ob} \bigl( \dcat{C}(A, \cat{M}) \bigr) \to 
\opn{Ob} \bigl( \dcat{C}(A^{\mrm{op}}, \cat{M}^{\mrm{op}}) \bigr) , \]
and it is clearly bijective.  

The operation of $E$ on morphisms is of course that of $\dcat{C}(F)$. 
It remains to verify that the resulting morphisms in 
$\dcat{C}(\cat{M}^{\mrm{op}})$ respect the action of 
elements of $A^{\mrm{op}}$. Namely that the condition in Definition 
\ref{dfn:1703} is satisfied. Take any morphism 
\[ \phi \in \opn{Hom}_{\dcat{C}(A, \cat{M})} 
\bigl( (M_0, \d_{M_0}, f_{M_0}), (M_1, \d_{M_1}, f_{M_1}) \bigr)^i \]
and any element $a \in (A^{\mrm{op}})^j$; and write 
$(N_k, \d_{N_k}, f_{N_k}) := E(M_k, \d_{M_k}, f_{M_k})$
and
\[ \psi := \dcat{G}(F)(\phi) \in 
\opn{Hom}_{\dcat{C}(\cat{M}^{\mrm{op}})} 
\bigl( (N_1, \d_{N_1}), (N_0, \d_{N_0}) \bigr)^i . \]
We have to prove that 
\begin{equation} \label{eqn:2366}
\psi \circ f_{N_1}(a) = (-1)^{i \cd j} \cd f_{N_0}(a) \circ \psi . 
\end{equation}
This is done using Lemma \ref{lem:2495}, like in the proof of 
Theorem \ref{thm:3015}; and we leave this to the reader.

\medskip \noindent 
(2) Clear from Theorem \ref{thm:3015}.

\medskip \noindent 
(3) Exactness in the categories 
$\dcat{C}_{\mrm{str}}(A, \cat{M})^{\mrm{op}}$
and 
$\dcat{C}_{\mrm{str}}(A^{\mrm{op}}, \cat{M}^{\mrm{op}})$
is checked in each degree separately, and both are exactness in the abelian 
category $\cat{M}^{\mrm{op}}$. The functor
$\opn{Str}(\opn{Flip})$ is $\pm \opn{Id}_{\cat{M}^{\mrm{op}}}$
in each degree, so it is exact. 

\medskip \noindent 
(4) As complexes in $\cat{M}^{\mrm{op}}$, $M$ and $\opn{Flip}(M)$ are equal, up 
to the renumbering of degrees and the signs of the differentials in various 
degrees. So the cohomology objects satisfy 
$\opn{H}^i(M) \cong \opn{H}^{-i}(\opn{Flip}(M))$
in $\cat{M}^{\mrm{op}}$. 
\end{proof}

\begin{exer} \label{exer:2365}
Prove formula (\ref{eqn:2366}) above. 
\end{exer}

\begin{rem} \label{rem:2505}
Theorem \ref{thm:2495} will be used to introduce a triangulated structure on 
the category $\dcat{K}(A, \cat{M})^{\mrm{op}}$. This will be done in Subsection 
\ref{subsec:opp-hom-triang}.  

Combined with Proposition \ref{prop:2500}, Theorem \ref{thm:2495} allows us to 
replace a contravariant DG functor 
$F : \dcat{C}(A, \cat{M}) \to \cat{D}$
with a usual, covariant, DG functor 
$F \circ \opn{Flip}^{-1} : \dcat{C}(A^{\mrm{op}}, \cat{M}^{\mrm{op}})
\to \cat{D}$.
This replacement is going to be very useful when discussing formal properties, 
such as existence of derived functors etc. 

However, in practical terms (e.g.\ for producing resolutions of DG modules),
the category $\dcat{C}(A^{\mrm{op}}, \cat{M}^{\mrm{op}})$ is not very 
helpful. The reason is that the opposite abelian category 
$\cat{M}^{\mrm{op}}$ is almost always a synthetic construction (it does not 
``really exist in concrete terms''). See Remark \ref{rem:2510} and 
Example \ref{exa:2990}. 

We are going to maneuver between the two approaches for reversal of morphisms, 
each time choosing the more useful approach. 
\end{rem}

\cleardoublepage
\mysection{Translations and Standard Triangles} \label{sec:stand-tri}

\AYcopyright

As before, we fix a $\K$-linear abelian category $\cat{M}$,  and a
DG central $\K$-ring $A$. In this section we study the translation functor and 
the standard cone of a strict morphism, all in the context of the DG category 
$\dcat{C}(A, \cat{M})$. 

We then study properties of DG functors 
$F : \dcat{C}(A, \cat{M}) \to \dcat{C}(B, \cat{N})$
between such DG categories. In view of Theorem \ref{thm:2495} it suffices to 
look at covariant DG functors (and not to worry about contravariant DG 
functors).

\mysubsection{The Translation Functor} \label{subsec:translation}

The translation functor goes back to the beginnings of derived categories -- 
see Remark \ref{rem:1275}. The treatment in this subsection, with the
``little $\opn{t}$ operator'',  is taken from \cite[Section 1]{Ye9}. 
Here and later, the word ``operator'' is used as a 
synonym for ``morphism in a linear category''. 

\begin{dfn} \label{dfn:1260}
Let $M = \{ M^i \}_{i \in \Z}$ be a graded object in $\cat{M}$, i.e.\ an 
object of $\dcat{G}(\cat{M})$. The {\em translation}
\index{Translation! of a graded module}
of $M$ is the object 
\[ \opn{T}(M) = \bigl\{ \opn{T}(M)^i \bigr\}_{i \in \Z} \in \dcat{G}(\cat{M}) \]
defined as follows: the graded component of degree $i$ of $\opn{T}(M)$ is
$\opn{T}(M)^i := M^{i + 1}$.
\end{dfn}

\begin{dfn}[The little t operator] \label{dfn:1261}
Let $M = \{ M^i \}_{i \in \Z}$ be an object of $\dcat{G}(\cat{M})$. We define 
\[ \opn{t}_M : M \to \opn{T}(M) \]
\index{Little t operator}%
\index{1-TM@$\opn{t}_{M}$}
to be the degree $-1$ morphism of graded objects of $\cat{M}$, that for every 
degree $i$ is the identity morphism
${\opn{t}_M}|_{M^{i}} := \opn{id}_{M^{i}} :  M^{i} \iso \opn{T}(M)^{i - 1}$
of the object $M^{i}$ in $\cat{M}$.
\end{dfn}

Note that the morphism 
$\opn{t}_M : M \to \opn{T}(M)$ is a degree $-1$ isomorphism in the graded 
category $\dcat{G}(\cat{M})$. 

\begin{dfn} \label{dfn:2300}
For $M \in \dcat{G}(\cat{M})$ we define the morphism 
$\opn{t}_M^{-1} : \opn{T}(M) \to M$
in $\dcat{G}(\cat{M})$ to be the inverse of $\opn{t}_M$. 
\end{dfn}

Of course the morphism $\opn{t}_M^{-1}$ has degree $+1$. The reason for stating 
this definition is to avoid the potential confusion between the morphism 
$\opn{t}_M^{-1}$ in $\dcat{G}(\cat{M})$ and the degree $-1$ component of the 
morphism $\opn{t}_M$, which we denote by ${\opn{t}_M}|_{M^{-1}}$, as in 
Definition \ref{dfn:1261} above. 

\begin{dfn} \label{dfn:1170}
Let $M = \{ M^i \}_{i \in \Z}$ be a DG $A$-module in $\cat{M}$, i.e.\ an 
object of $\dcat{C}(A, \cat{M})$.
The {\em translation} of $M$ 
\index{Translation! of a DG module}
\index{1-T(M)@$\opn{T}(M)$}
is the object 
\[ \opn{T}(M) \in \dcat{C}(A, \cat{M}) \]
defined as follows. 
\begin{enumerate}
\item As a graded object of $\cat{M}$, it is as specified in 
Definition \ref{dfn:1260}. 

\item The differential $\d_{\opn{T}(M)}$
is defined by the formula
\[ \d_{\opn{T}(M)} := - \opn{t}_M \circ \, \d_M^{} \circ 
\opn{t}_M^{-1} . \]

\item Let $f_M : A \to \opn{End}_{\cat{M}}(M)$ be the DG ring homomorphism that 
determines the action of $A$ on $M$. Then 
\[ f_{\lms \opn{T}(M)} : A \to \opn{End}_{\cat{M}}(\opn{T}(M)) \]
is defined by
\[ f_{\lms \opn{T}(M)}(a) := (-1)^{\lms j} \cd 
\opn{t}_M \circ \, f_{M}(a) \circ \opn{t}_M^{-1} \]
for $a \in A^j$.
\end{enumerate}
\end{dfn}

Thus, the differential 
$\d_{\lms \opn{T}(M)}  = \{ \d_{\lms \opn{T}(M)}^i \}_{i \in \Z}$ makes this 
diagram in $\cat{M}$ commutative for every $i$~: 
\[ \UseTips  \xymatrix @C=8ex @R=6ex {
\opn{T}(M)^i 
\ar[r]^{\d_{\lms \opn{T}(M)}^i}
&
\opn{T}(M)^{i + 1}
\\
M^{i + 1}
\ar[u]^{\opn{t}_M}
\ar[r]^{- \d_M^{i + 1}}
&
M^{i + 2}
\ar[u]_{\opn{t}_M}
} \]
And the left $A$-module structure makes this diagram in $\cat{M}$ 
commutative for every $i$ and every $a \in A^j$~:
\[ \UseTips  \xymatrix @C=12ex @R=6ex {
\opn{T}(M)^i 
\ar[r]^{f_{\lms \opn{T}(M)}(a)}
&
\opn{T}(M)^{i + j}
\\
M^{i + 1}
\ar[u]^{\opn{t}_M}
\ar[r]^{ (-1)^{\lms j} \cd f_{\lms M}(a)}
&
M^{i + j + 1}
\ar[u]_{\opn{t}_M}
} \]

\begin{prop}  \label{prop:1176}
The morphisms $\opn{t}_M$ and $\opn{t}_M^{-1}$ are cocycles, in the DG 
$\K$-modules 
$\opn{Hom}_{A, \cat{M}} \bigl( M, \opn{T}(M) \bigr)$ and 
$\opn{Hom}_{A, \cat{M}} \bigl( \opn{T}(M), M \bigr)$
respectively.
\end{prop}

\begin{proof}
We use the notation $\d_{\mrm{Hom}}$ for the differential in the DG module \lb 
$\opn{Hom}_{A, \cat{M}} \bigl( M, \opn{T}(M) \bigr)$. 
Let us calculate. Because $\opn{t}_M$ has degree $-1$, we have
\[ \begin{aligned}
& \d_{\mrm{Hom}}(\opn{t}_M) = 
\d_{\opn{T}(M)} \circ  \opn{t}_M + \opn{t}_M \circ \, \d_{M}
\\ & \quad 
=  (- \opn{t}_M \circ \, \d_M^{} \circ \opn{t}_M^{-1}) 
\circ  \opn{t}_M + \opn{t}_M \circ \, \d_{M} = 0 . 
\end{aligned} \]

As for $\opn{t}_M^{-1}$~: this is done using the graded Leibniz rule, just like 
in the proof  Proposition \ref{prop:1220}.
\end{proof}

\begin{dfn} \label{dfn:1171}
\index{Translation! of a graded morphism}
Given a  morphism 
$\phi \in \opn{Hom}_{A, \cat{M}} (M, N)^{i}$
we define the morphism 
\[ \opn{T}(\phi)  \in 
\opn{Hom}_{A, \cat{M}} \bigl( \opn{T}(M), \opn{T}(N) \bigr)^i \]
to be
\[ \opn{T}(\phi) := (-1)^i \cd \opn{t}_N \circ \, \phi \circ 
\opn{t}_M^{-1}  .  \]
\end{dfn}

To clarify this definition, let us write 
$\phi = \{ \phi^{\lms j} \}_{j \in \Z}$, 
so that $\phi^{\lms j} : M^j \to N^{j + i}$
is a morphism in $\cat{M}$. Then 
$\opn{T}(\phi)^{\lms j} : \opn{T}(M)^{\lms j} \to \opn{T}(N)^{\lms j + i}$
is  
$\opn{T}(\phi)^{\lms j} = (-1)^i \cd \opn{t}_N \circ \, \phi^{\lms j + 1} 
\circ \opn{t}_M^{-1}$.
The corresponding commutative diagram in $\cat{M}$, for each $i, j$, is: 
\begin{equation} \label{eqn:1162}
 \UseTips  \xymatrix @C=12ex @R=6ex {
\opn{T}(M)^{\lms j}
\ar[r]^{\opn{T}(\phi)^{\lms j}}
&
\opn{T}(N)^{\lms j + i}
\\
M^{j + 1}
\ar[u]^{\opn{t}_M}
\ar[r]^{ (-1)^i \cd \phi^{\lms j + 1} }
&
N^{j + i + 1}
\ar[u]_{\opn{t}_N}
}
\end{equation}

\begin{thm} \label{thm:1260} 
Let $\cat{M}$ be an abelian category and let $A$ be a DG ring.
\begin{enumerate}
\item The assignments $M \mapsto \opn{T}(M)$ and 
$\phi \mapsto \opn{T}(\phi)$ are a  DG functor
\[ \opn{T} : \dcat{C}(A, \cat{M}) \to \dcat{C}(A, \cat{M}) . \]

\item The collection 
$\opn{t} := \{ \opn{t}_M \}_{M \in \dcat{C}(A, \cat{M})}$ 
is a degree $-1$ isomorphism 
\[ \opn{t} : \opn{Id} \to  \opn{T} \]
of DG functors from $\dcat{C}(A, \cat{M})$ to itself. 
\end{enumerate}
\end{thm}

\begin{proof} \mbox{}

\smallskip \noindent
(1) Take morphisms
$\phi_1 : M_0 \to M_1$ and $\phi_2 : M_1 \to M_2$, of degrees $i_1$ and $i_2$ 
respectively. Then 
\[ \begin{aligned}
& \opn{T}(\phi_2 \circ \phi_1) = 
(-1)^{i_1 + i_2} \cd \opn{t}_{M_2} \circ \, (\phi_2 \circ \phi_1) \circ 
\opn{t}_{M_0}^{-1}
\\ 
& \qquad  = (-1)^{i_1 + i_2} \cd \opn{t}_{M_2} \circ \, \phi_2 \circ 
(\opn{t}_{M_1}^{-1} \circ \opn{t}_{M_1}) \circ \, \phi_1 \circ  
\opn{t}_{M_0}^{-1} 
\\ 
& \qquad = \bigl( (-1)^{i_2} \cd \opn{t}_{M_2} \circ \, \phi_2 \circ 
\opn{t}_{M_1}^{-1}\bigr) \circ 
\bigl( (-1)^{i_1} \cd  \opn{t}_{M_1} \circ \, \phi_1 \circ \opn{t}_{M_0}^{-1} 
\bigr)
\\ 
& \qquad =
\opn{T}(\phi_2) \circ \opn{T}(\phi_1) . 
\end{aligned} \]
Clearly $\opn{T}(\opn{id}_{M}) = \opn{id}_{\opn{T}(M)}$, and 
$\opn{T}(\la \cd \phi + \psi) = \la \cd \opn{T}(\phi) + \opn{T}(\psi)$
for all $\la \in \K$ and 
$\phi, \psi \in \opn{Hom}_{A, \cat{M}} (M_0, M_1)^{i}$. 
So $\opn{T}$ is a $\K$-linear graded functor. 

By Proposition \ref{prop:1176} we know that 
$\d \circ \opn{t} = - \opn{t} \circ \, \d$
and 
$\d \circ \opn{t}^{-1} = - \opn{t}^{-1} \circ \, \d$.
This implies that for every morphism 
$\phi$ in $\dcat{C}(A, \cat{M})$, we have 
$\opn{T}(\d(\phi)) = \d(\opn{T}(\phi))$.
So $\opn{T}$ is a DG functor. 

\medskip \noindent
(2) Take some $\phi \in \opn{Hom}_{A, \cat{M}} (M_0, M_1)^{i}$. 
We have to prove that 
$\opn{t}_{M_1} \circ \, \phi = \lb (-1)^{i} \cd \opn{T}(\phi) \circ 
\opn{t}_{M_0}$ 
as elements of 
$\opn{Hom}_{A, \cat{M}} \bigl( M_0, \opn{T}(M_1) \bigr)^{i + 1}$.
But by Definition \ref{dfn:1171} we have 
\[ \opn{T}(\phi) \circ \opn{t}_{M_0} = 
\bigl( (-1)^{i} \cd \opn{t}_{M_1} \circ \, \phi \circ \opn{t}_{M_0}^{-1} \bigr) 
\circ \opn{t}_{M_0} =  (-1)^{i} \cd \opn{t}_{M_1} \circ \, \phi . \qedhere \]
\end{proof}

\begin{dfn} \label{dfn:1275}
\index{Translation functor! of the DG category $\dcat{C}(A, \cat{M})$}
We call $\opn{T}$ the {\em translation functor} of the DG category 
$\dcat{C}(A, \cat{M})$.  
\end{dfn}

\begin{cor} \label{cor:1300} \mbox{}
\begin{enumerate}
\item The functor $\opn{T}$ is an automorphism of the DG category 
$\dcat{C}(A, \cat{M})$.

\item For every $k, l \in \Z$ there is an equality of functors
$\opn{T}^{l} \circ \opn{T}^{k} = \opn{T}^{l + k}$.
\end{enumerate}
\end{cor}

\begin{proof}
(1) By the theorem, $\opn{T}$ is a DG functor. By definition, $\opn{T}$
is bijective on the set of objects of $\dcat{C}(A, \cat{M})$ and on the sets of 
morphisms. 

\medskip \noindent 
(2) By part (1) of this corollary, the inverse $\opn{T}^{-1}$ is a uniquely 
defined functor (not just up to an isomorphism of functors). 
Thus we can define $\opn{T}^k$, the $k$-th power of $\opn{T}$, and the equality 
stated holds. 
\end{proof}

\begin{prop} \label{prop:1935}
Consider any $M \in \dcat{C}(A, \cat{M})$.
\begin{enumerate}
\item There is equality 
$\opn{t}_{\opn{T}(M)} = - \opn{T}(\opn{t}_M)$
of degree $-1$ morphisms $\opn{T}(M) \to \opn{T}^2(M)$
in  $\dcat{C}(A, \cat{M})$.

\item There is equality 
$\opn{t}_{\opn{T}^{-1}(M)} = - \opn{T}^{-1}(\opn{t}_M)$
between these degree $-1$ morphisms  
$\opn{T}^{-1}(M) \to \opn{T}(\opn{T}^{-1}(M)) = M = \opn{T}^{-1}(\opn{T}(M))$
in  $\dcat{C}(A, \cat{M})$.
\end{enumerate}
\end{prop}

\begin{proof}
(1) This is an easy calculation, using Definition \ref{dfn:1171}:
$\opn{T}(\opn{t}_M) = \lb 
- \opn{t}_{\opn{T}(M)} \circ \, \opn{t}_{M} \circ \opn{t}_{M}^{-1} = 
- \opn{t}_{\opn{T}(M)}$.

\medskip \noindent 
(2) A similar calculation. 
\end{proof}

\begin{prop} \label{prop:4155}
Let $M \in \dcat{C}(A, \cat{M})$ and $i \in \Z$. There is equality 
$\d_{\opn{T}^i(M)} = \opn{T}^i(\d_M)$
in 
$\opn{Hom}_{A, \cat{M}} \bigl(\opn{T}^i(M) , \opn{T}^i(M) \bigr)^{1}$. 
\end{prop}

\begin{proof}
We start with $i = 1$. The differential $\d_M$ is an element of 
$\opn{Hom}_{\cat{M}}(M , M)^1$. 
By Definitions \ref{dfn:1170} and \ref{dfn:1171} we get 
$\d_{\opn{T}(M)} = - \opn{t}_M \circ \, \d_M \circ \opn{t}_M^{-1} 
=  \opn{T}(\d_{M})$.
Using induction on $i$ the assertion holds for all $i \geq 0$. 

For $i \leq 0$ we use descending induction on $i$. 
We assume that the assertion holds for $i$. 
Let us define $N := \opn{T}^{i - 1}(M)$, so 
$\opn{T}^{}(N) = \opn{T}^{i}(M)$. 
By the previous paragraph, with the DG module $N$, we know that 
$\d_{\opn{T}^{i}(M)} = \d_{\opn{T}^{}(N)} = \opn{T}(\d_{N}) = 
\opn{T}(\d_{\opn{T}^{i - 1}(M)})$.
Applying the functor $\opn{T}^{-1}$ to this equality we get 
$\opn{T}^{i - 1}(\d_{M}) = 
\opn{T}^{-1} \bigl( \opn{T}^{i}(\d_M) \bigr) = 
\opn{T}^{-1}(\d_{\opn{T}^{i}(M)})
= \d_{\opn{T}^{i - 1}(M)}$. 
\end{proof}

\begin{rem} \label{rem:1275}
There are several names in the literature for the translation functor
$\opn{T}$~: {\em twist}, {\em shift} and {\em suspension}. There are also 
several notations: $\opn{T}(M) = M[1] = \Sigma M$.
In the later part of this book we shall use the notation 
$M[k] := \opn{T}^{k}(M)$ for the $k$-th translation. 
\end{rem}

\mysubsection{The Standard Triangle of a Strict Morphism} \label{subsec:cone}

As before,  we fix an abelian category $\cat{M}$ and a DG ring $A$.
Here is the cone construction in $\dcat{C}(A, \cat{M})$, as it looks 
using the operator $\opn{t}$. 

\begin{dfn} \label{dfn:1172}
Let $\phi : M \to N$ be a strict morphism in $\dcat{C}(A, \cat{M})$.
The {\em standard cone of $\phi$}%
\index{Cone! standard {\indash} of strict morphism in $\dcat{C}(A, \cat{M})$}
is the object 
$\opn{Cone}(\phi) \in \dcat{C}(A, \cat{M})$
defined as follows. As a graded $A$-module in $\cat{M}$ we let
\[ \opn{Cone}(\phi) := N \oplus \opn{T}(M) . \]
The differential $\d_{\mrm{Cone}}$ is this: if we express 
the graded module as a column
\[ \opn{Cone}(\phi) = \bmat{N \\[0.1em] \opn{T}(M)} , \]
then $\d_{\mrm{Cone}}$ is left multiplication by the matrix
\[ \d_{\mrm{Cone}} :=
\bmat{\d_N & \phi \circ \opn{t}_M^{-1} \\[0.1em] 0 & \d_{\opn{T}(M)}} \]
of degree $1$ morphisms of graded $A$-modules in $\cat{M}$.
\end{dfn}

In other words, the morphism 
$\d_{\mrm{Cone}}^i : \opn{Cone}(\phi)^i \to \opn{Cone}(\phi)^{i + 1}$
is 
$\d_{\mrm{Cone}}^i = \d^i_N + \d^i_{\opn{T}(M)} + \phi^{i + 1} \circ 
\opn{t}_M^{-1}$,
where $\phi^{i + 1} \circ \opn{t}_M^{-1}$ is the composed morphism 
$\opn{T}(M)^i \xar{\opn{t}_M^{-1}} M^{i + 1} \xar{\phi^{i + i}}  N^{i + 1}$.

In the situation of Definition \ref{dfn:1172}, let us denote by 
\begin{equation} \label{eqn:1195}
e_{\phi} : N \to N \oplus \opn{T}(M) 
\end{equation}
the embedding, and by
\begin{equation} \label{eqn:1196}
p_{\phi} : N \oplus \opn{T}(M) \to \opn{T}(M)
\end{equation}
the projection. 
Thus, as matrices we have 
\[ e_{\phi} = \bmat{ \opn{id}_N  \\[0.1em] 0 } \quad \tup{and} \quad 
p_{\phi} = \bmat{ 0 &  \opn{id}_{\opn{T}(M)} } . \]
The standard cone of $\phi$ sits in the exact sequence 
\begin{equation} \label{eqn:1815}
0 \to N \xar{e_{\phi}} \opn{Cone}(\phi) \xar{p_{\phi}} \opn{T}(M) \to 0
\end{equation}
in the abelian category $\dcat{C}_{\mrm{str}}(A, \cat{M})$. 

\begin{dfn} \label{dfn:1140}
Let $\phi : M \to N$ be a morphism in $\dcat{C}_{\mrm{str}}(A, \cat{M})$.
The diagram 
\[ M \xar{ \, \phi \, } N \xar{ \, e_{\phi} \, } \opn{Cone}(\phi) 
\xar{ \, p_{\phi} \, } \opn{T}(M) \]
in $\dcat{C}_{\mrm{str}}(A, \cat{M})$ is called the {\em standard triangle} 
\index{Triangle! standard {\indash} of a strict morphism in 
$\dcat{C}(A, \cat{M})$}
associated to $\phi$. 
\end{dfn}

The cone construction is functorial, in the following sense. 

\begin{prop}\label{prop:1162}
Let 
\[ \UseTips  \xymatrix @C=8ex @R=6ex {
M_0 
\ar[r]^{\phi_0}
\ar[d]_{\psi}
& 
N_0 
\ar[d]^{\chi}
\\
M_1
\ar[r]^{\phi_1}
& 
N_1 
} \]
be a commutative diagram in $\dcat{C}_{\mrm{str}}(A, \cat{M})$. 
Then 
\begin{equation} \label{eqn:1113}
\bigl( \chi, \opn{T}(\psi) \bigr) : \opn{Cone}(\phi_0) \to \opn{Cone}(\phi_1)
\end{equation}
is a morphism in $\dcat{C}_{\mrm{str}}(A, \cat{M})$, and the diagram 
\[ \UseTips  \xymatrix @C=10ex @R=6ex {
M_0 
\ar[r]^{\phi_0}
\ar[d]_{\psi}
& 
N_0 
\ar[r]^(0.4){e_{\phi_0}}
\ar[d]_{\chi}
& 
\opn{Cone}(\phi_0)
\ar[r]^{p_{\phi_0}}
\ar[d]_{(\chi, \opn{T}(\psi))}
& 
\opn{T}(M_0)
\ar[d]_{\opn{T}(\psi)}
\\
M_1
\ar[r]^{\phi_1}
& 
N_1
\ar[r]^(0.4){e_{\phi_1}}
& 
\opn{Cone}(\phi_1)
\ar[r]^{p_{\phi_1}}
& 
\opn{T}(M_1)
} \]
in $\dcat{C}_{\mrm{str}}(A, \cat{M})$ is commutative. 
\end{prop}

\begin{proof}
This is a simple consequence of the definitions. 
\end{proof}

\mysubsection{The Gauge of a Graded Functor} \label{subsec:gauge}

Recall that we have an abelian category $\cat{M}$ and a DG ring $A$.
The next definition is new. 

\begin{dfn} \label{dfn:1214}
Let 
$F : \dcat{C}(A, \cat{M}) \to \dcat{C}(B, \cat{N})$
be a graded functor. For every object 
$M \in \dcat{C}(A, \cat{M})$ let 
\[ \ga_{F, M} := \d_{F(M)} - F(\d_M)  \in 
\opn{Hom}_{B, \cat{N}} \bigl( F(M), F(M) \bigr)^1 . \]
The collection of morphisms 
\[ \ga_F := \{ \ga_{F, M} \}_{M \in \dcat{C}(A, \cat{M})}  \]
is called the {\em gauge of $F$}.
\index{Gauge of a graded functor} 
\end{dfn}

The next theorem is due to R. Vyas. 

\begin{thm} \label{thm:1225}
The following two conditions are equivalent for a graded functor 
$F : \dcat{C}(A, \cat{M}) \to \dcat{C}(B, \cat{N})$.
\begin{enumerate}
\rmitem{i} $F$ is a DG functor.
 
\rmitem{ii} The gauge $\ga_F$ is a degree $1$ morphism of graded functors 
$\ga_F : F \to F$.  
\end{enumerate}
\end{thm}

\begin{proof}
Recall that $F$ is a DG functor (condition (i)) iff 
\begin{equation}  \label{eqn:1231}
(F \circ \d_{A, \cat{M}})(\phi) = (\d_{B, \cat{N}} \circ F)(\phi)
\end{equation}
for every 
$\phi \in \opn{Hom}_{A, \cat{M}}(M_{0}, M_{1})^{i}$.
And $\ga_F$ is a degree $1$ morphism of graded functors  (condition (ii))
iff 
\begin{equation}  \label{eqn:1232}
\ga_{F, M_1} \circ  F(\phi) = (-1)^i \cd F(\phi) \circ \ga_{F, M_0} 
\end{equation}
for every such $\phi$. 

Here is the calculation. 
Because $F$ is a graded functor, we get
\begin{equation}  \label{eqn:1225}
\begin{aligned}
& F \bigl( \d_{A, \cat{M}}(\phi) \bigr) = 
F \bigl( \d_{M_1} \circ \phi - (-1)^i \cd \phi \circ \d_{M_0} \bigr) 
\\
& \qquad = F(\d_{M_1}) \circ F(\phi) - (-1)^i \cd F(\phi) \circ F(\d_{M_0})
\end{aligned} 
\end{equation}
and 
\begin{equation}  \label{eqn:1226}
\d_{B, \cat{N}} \bigl( F(\phi)  \bigr) = 
\d_{F(M_1)} \circ F(\phi) - (-1)^i \cd F(\phi) \circ \d_{F(M_0)} .
\end{equation}
Using equations  (\ref{eqn:1225}) and (\ref{eqn:1226}),
and the definition of $\ga_F$, we obtain 
\begin{equation}  \label{eqn:1227}
\begin{aligned}
& (F \circ \d_{A, \cat{M}} - \d_{B, \cat{N}} \circ F)(\phi) = 
F \bigl( \d_{A, \cat{M}}(\phi) \bigr) - \d_{B, \cat{N}} \bigl( F(\phi)  \bigr)
\\
& \qquad = \bigl( F(\d_{M_1}) - \d_{F(M_1)} \bigr) \circ F(\phi) 
- (-1)^i \cd F(\phi) \circ \bigl( F(\d_{M_0}) - \d_{F(M_0)} \bigr)
\\
& \qquad = - \ga_{F, M_1} \circ  F(\phi) 
+ (-1)^i \cd F(\phi) \circ \ga_{F, M_0} . 
\end{aligned} 
\end{equation}

Finally, the vanishing of the first expression in (\ref{eqn:1227}) 
is the same as equality in (\ref{eqn:1231});  whereas the vanishing of the last 
expression in (\ref{eqn:1227}) is the same as equality in (\ref{eqn:1232}).
\end{proof}

\mysubsection{The Translation Isomorphism of a DG Functor}

Here we consider \lb belian categories $\cat{M}$ and $\cat{N}$, 
and DG rings $A$ and $B$.
The translation functor of the DG category $\dcat{C}(A, \cat{M})$ will be 
denoted here by $\opn{T}_{A, \cat{M}}$. 
For an object $M \in \dcat{C}(A, \cat{M})$, we have the little t operator 
\[ \opn{t}_M \in \opn{Hom}_{A, \cat{M}} \bigl( M, \opn{T}_{A, \cat{M}}(M) 
\bigr)^{-1} . \]
This is an isomorphism in the DG category $\dcat{C}(A, \cat{M})$. 
Likewise for the DG category $\dcat{C}(B, \cat{N})$.

\begin{dfn} \label{dfn:1150}
Let 
$F : \dcat{C}(A, \cat{M}) \to \dcat{C}(B, \cat{N})$
be a DG functor.
For an object $M \in \dcat{C}(A, \cat{M})$, let 
\[ \tau_{F, M} : F(\opn{T}_{A, \cat{M}}(M)) \to \opn{T}_{B, \cat{N}}(F(M)) \]
be the degree $0$ isomorphism
\[ \tau_{F, M} := \opn{t}_{F(M)} \circ \, F(\opn{t}_M)^{-1}  \]
in $\dcat{C}(B, \cat{N})$, called the {\em translation isomorphism}
\index{Translation isomorphism! of a DG functor}
of the functor $F$ at the object $M$. 
\end{dfn}

The isomorphism $\tau_{F, M}$ sits in the following commutative diagram
\[ \UseTips  \xymatrix @C=8ex @R=6ex {
F(\opn{T}_{A, \cat{M}}(M))
\ar[r]^{\tau_{F, M}}
&
\opn{T}_{B, \cat{N}}(F(M))
\\
F(M)
\ar[u]^{ F(\opn{t}_M) } 
\ar[ur]_{ \opn{t}_{F(M)} }
} \]
of isomorphisms in the category $\dcat{C}(B, \cat{N})$.

\begin{prop} \label{prop:1179}
$\tau_{F, M}$ is an isomorphism in $\dcat{C}_{\mrm{str}}(B, \cat{N})$.
\end{prop}

\begin{proof}
We know that $\tau_{F, M}$ is an isomorphism in $\dcat{C}(B, \cat{N})$.
It suffices to prove that both $\tau_{F, M}$ and its inverse  
$\tau_{F, M}^{-1}$ are 
strict morphisms. Now by Proposition \ref{prop:1176}, $\opn{t}_M$ and 
$\opn{t}_M^{-1}$ are cocycles. Therefore, $F(\opn{t}_M)$ and 
$F(\opn{t}_M)^{-1} = F(\opn{t}_M^{-1})$ are cocycles. For the same reason,
$\opn{t}_{F(M)}$ and $\opn{t}_{F(M)}^{-1}$  are cocycles. But 
$\tau_{F, M} = \opn{t}_{F(M)} \circ \, F(\opn{t}_M)^{-1}$, and 
$\tau_{F, M}^{-1} = F(\opn{t}_M) \, \circ \opn{t}_{F(M)}^{-1}$.
\end{proof}

\begin{thm} \label{thm:1150}
Let 
$F : \dcat{C}(A, \cat{M}) \to \dcat{C}(B, \cat{N})$
be a DG functor. Then the collection
\[ \tau_F := \{ \tau_{F, M} \}_{M \in \dcat{C}(A, \cat{M})} \] 
is an isomorphism 
\[ \tau_F : F \circ \opn{T}_{A, \cat{M}} \iso \opn{T}_{B, \cat{N}} \circ \, F \]
of functors 
$\dcat{C}_{\mrm{str}}(A, \cat{M}) \to \dcat{C}_{\mrm{str}}(B, \cat{N})$.
\end{thm}

The slogan summarizing this theorem is ``A DG functor commutes with 
translations''. 

\begin{proof}
In view of Proposition \ref{prop:1179}, all we need to prove is that 
$\tau_{F}$ is a morphism of functors (i.e.\ it is a natural transformation). 

Let $\phi : M_0 \to M_1$ be a morphism in 
$\dcat{C}_{\mrm{str}}(A, \cat{M})$. We must prove that the diagram 
\[ \UseTips  \xymatrix @C=8ex @R=6ex {
(F \circ \opn{T}_{A, \cat{M}})(M_0)
\ar[r]^{ \tau_{F, M_0} } 
\ar[d]_{ (F \circ \opn{T}_{A, \cat{M}})(\phi) }
&
(\opn{T}_{B, \cat{N}} \circ \, F)(M_0)
\ar[d]^{ (\opn{T}_{B, \cat{N}} \circ \, F)(\phi) }
\\
(F \circ \opn{T}_{A, \cat{M}})(M_1)
\ar[r]^{ \tau_{F, M_1} } 
&
(\opn{T}_{B, \cat{N}} \circ \, F)(M_1)
} \]
in $\dcat{C}_{\mrm{str}}(B, \cat{N})$ is commutative.
This will be true if the next diagram 
\[ \UseTips  \xymatrix @C=10ex @R=6ex {
(F \circ \opn{T}_{A, \cat{M}})(M_0)
\ar[d]_{ (F \circ \opn{T}_{A, \cat{M}})(\phi) }
&
F(M_0)
\ar[l]_(0.4){ F(\opn{t}_{M_0}) }
\ar[r]^(0.4){ \opn{t}_{F(M_0)} }
\ar[d]_{F(\phi)}
&
(\opn{T}_{B, \cat{N}} \circ \, F)(M_0)
\ar[d]^{ (\opn{T}_{B, \cat{N}} \circ \, F)(\phi) }
\\
(F \circ \opn{T}_{A, \cat{M}})(M_1)
&
F(M_1)
\ar[l]_(0.4){ F(\opn{t}_{M_1}) }
\ar[r]^(0.4){ \opn{t}_{F(M_1)} }
&
(\opn{T}_{B, \cat{N}} \circ \, F)(M_1)
} \]
in $\dcat{C}(B, \cat{N})$, whose horizontal arrows are isomorphisms, is 
commutative. For this to be true, it is enough to prove that both squares in 
this diagram are commutative. This is true by 
Theorem \ref{thm:1260}(2). 
\end{proof}

Recall that the translation $\opn{T}$ and all its powers are DG functors. 
To finish this subsection, we calculate their translation isomorphisms. 

\begin{prop} \label{prop:1300}
For every integer $k$, the translation isomorphism of the DG functor 
$\opn{T}^k$ is 
$\tau_{\opn{T}^k} = (-1)^k \cd \opn{id}_{\lms \opn{T}^{k + 1}}$,
where $\opn{id}_{\lms \opn{T}^{k + 1}}$ is the identity automorphism of the 
functor $\opn{T}^{k + 1}$. 
\end{prop}

\begin{proof}
By Definition \ref{dfn:1150} and Proposition \ref{prop:1935}(1), for $k = 1$ 
the formula is
$\tau_{\lms \opn{T}, M} = \opn{t}_{\opn{T}(M)} \circ \opn{T}(\opn{t}_M)^{-1} = 
- \opn{id}_{\lms \opn{T}^2(M)}$,
where $\opn{id}_{\lms \opn{T}^2(M)}$ is the identity automorphism of the DG 
module 
$\opn{T}^2(M)$. Hence $\tau_{\lms \opn{T}} =  - \opn{id}_{\opn{T}^2}$.
For other integers $k$ the calculation is similar.
\end{proof}

\mysubsection{Standard Triangles and DG Functors} \label{subsec:std-triangles}

In Subsection \ref{subsec:cone} we defined the standard triangle associated to 
a strict morphism; and in Subsection \ref{subsec:DGFunc} we defined DG 
functors. Now we show how these notions interact with each other. 
As before, we consider abelian categories $\cat{M}$ and $\cat{N}$, 
and DG rings $A$ and $B$.

Let 
$F : \dcat{C}(A, \cat{M}) \to \dcat{C}(B, \cat{N})$
be a DG functor. Given a morphism 
$\phi : M_0 \to M_1$ in $\dcat{C}_{\mrm{str}}(A, \cat{M})$,
we have a morphism 
$F(\phi) : F(M_0) \to F(M_1)$
in $\dcat{C}_{\mrm{str}}(B, \cat{N})$, 
and objects 
$F(\opn{Cone}_{A, \cat{M}}(\phi))$ and 
$\opn{Cone}_{B, \cat{N}}(F(\phi))$ 
in $\dcat{C}(B, \cat{N})$.
By definition, and using the fully faithful graded functor
$\opn{Und}$ from Proposition \ref{prop:2365},
there is a canonical isomorphism  
\begin{equation}  \label{eqn:1197}
\opn{Cone}_{A, \cat{M}}(\phi) \cong M_1 \oplus \opn{T}_{A, \cat{M}}(M_0)
\end{equation}
in $\dcat{G}_{\mrm{str}}(A, \cat{M})$.
Since $F$ is a linear functor, it commutes with finite direct sums,
and therefore there is a canonical isomorphism 
\begin{equation}  \label{eqn:1185}
F(\opn{Cone}_{A, \cat{M}}(\phi)) \cong F(M_1) \oplus 
F(\opn{T}_{A, \cat{M}}(M_0))
\end{equation}
in $\dcat{G}_{\mrm{str}}(A, \cat{M})$. 
And by definition there is a canonical isomorphism 
\begin{equation}  \label{eqn:1186}
\opn{Cone}_{B, \cat{N}}(F(\phi)) \cong F(M_1) \oplus 
\opn{T}_{B, \cat{N}}(F(M_0))
\end{equation}
in $\dcat{G}_{\mrm{str}}(A, \cat{M})$.
Warning: the isomorphisms (\ref{eqn:1197}), (\ref{eqn:1185}) and 
(\ref{eqn:1186}) might not commute with the differentials. 
The differentials on the right sides are  diagonal matrices, but on the left 
sides they are upper-triangular matrices (see Definition \ref{dfn:1172}). 

\begin{lem} \label{lem:1200}
Let 
$F , G : \dcat{G}(A, \cat{M}) \to \dcat{G}(B, \cat{N})$
be graded functors, and let $\eta : F \to G$ be a degree $j$ 
morphism of graded functors. Suppose 
$M \cong M_0 \oplus M_1$ in $\dcat{G}_{\mrm{str}}(A, \cat{M})$, with embeddings 
$e_i : M_i \to M$ and  projections $p_i : M \to M_i$. 
Then 
\[ \eta_M = \bigl( G(e_0), G(e_1) \bigr) \circ (\eta_{M_0}, \eta_{M_1}) \circ 
\bigl( F(p_0), F(p_1) \bigr) , \]
as degree $j$ morphisms $F(M) \to G(M)$ in 
$\dcat{G}(B, \cat{N})$. 
\end{lem}

The lemma says that the diagram 
\[ \UseTips  \xymatrix @C=14ex @R=6ex {
F(M) 
\ar[r]^(0.4){ (F(p_0), F(p_1)) }
\ar[d]_{ \eta_M }
& 
F(M_0) \oplus F(M_1)
\ar[d]^{ (\eta_{M_0}, \eta_{M_1}) }
\\
G(M)
& 
G(M_0) \oplus G(M_1)
\ar[l]_(0.6){ (G(e_0), G(e_1)) }
} \]
in $\dcat{G}(B, \cat{N})$ is commutative. 

\begin{proof}
It suffices to prove that the diagram below is commutative for $i = 0, 1$~:
\[ \UseTips  \xymatrix @C=10ex @R=6ex {
F(M_i)
\ar[r]^{ F(e_i) }
\ar[d]_{ \eta_{M_i} }
\ar@(u,u)[rr]^{ \opn{id} }
&
F(M) 
\ar[r]^{ F(p_i) }
\ar[d]_{ \eta_M }
& 
F(M_i) 
\ar[d]_{ \eta_{M_i} }
\\
G(M_i)
\ar[r]^{ G(e_i) }
\ar@(d,d)[rr]_{ \opn{id} }
&
G(M)
\ar[r]^{ G(p_i) }
& 
G(M_i) 
} \]
This is true because $\eta$ is a morphism of functors (a natural 
transformation). 
\end{proof}

\begin{thm} \label{thm:1185}
Let 
$F : \dcat{C}(A, \cat{M}) \to \dcat{C}(B, \cat{N})$
be a DG functor, and let 
$\phi : M_0 \to M_1$ be a morphism in $\dcat{C}_{\mrm{str}}(A, \cat{M})$.
Define the isomorphism 
\[ \opn{cone}(F, \phi) : F (\opn{Cone}_{A, \cat{M}}(\phi)) \iso 
\opn{Cone}_{B, \cat{N}} (F(\phi))  \]
in  $\dcat{G}_{\mrm{str}}(A, \cat{M})$ to be 
\[ \opn{cone}(F, \phi) := 
\bigl( \opn{id}_{F(M_1)}, \tau_{F, M_0} \bigr) . \]
Then\tup{:}
\begin{enumerate}
\item The isomorphism $\opn{cone}(F, \phi)$ commutes with the differentials, 
so it is an isomorphism in $\dcat{C}_{\mrm{str}}(A, \cat{M})$. 

\item The diagram 
\[ \UseTips  \xymatrix @C=7ex @R=6ex {
F(M_0) 
\ar[r]^{F(\phi)}
\ar[d]_{\opn{id}}
& 
F(M_1) 
\ar[r]^(0.4){F(e_{\phi})}
\ar[d]_{\opn{id}}
& 
F(\opn{Cone}_{A, \cat{M}}(\phi))
\ar[r]^{F(p_{\phi})}
\ar[d]_{\opn{cone}(F, \phi)}
& 
F(\opn{T}_{A, \cat{M}}(M_0))
\ar[d]_{ \tau_{F, M_0}}
\\
F(M_0)
\ar[r]^{F(\phi)}
& 
F(M_1)
\ar[r]^(0.4){e_{F(\phi)}}
& 
\opn{Cone}_{B, \cat{N}}(F(\phi))
\ar[r]^{p_{F(\phi)}}
& 
\opn{T}_{B, \cat{N}}(F(M_0))
} \]
in $\dcat{C}_{\mrm{str}}(B, \cat{N})$ is commutative. 
\end{enumerate}
\end{thm}

When defining $\opn{cone}(F, \phi)$ above, we are using the decompositions 
(\ref{eqn:1185}) and  (\ref{eqn:1186}) in the category 
$\dcat{G}_{\mrm{str}}(A, \cat{M})$, 
and the isomorphism $\tau_{F, M_0}$ from Definition 
\ref{dfn:1150}.

The slogan summarizing this theorem is ``A DG functor sends standard triangles 
to standard triangles''. 

\begin{proof} \mbox{}

\smallskip \noindent
(1) To save space let us write $\th := \opn{cone}(F, \phi)$. 
We have to prove that $\d_{B, \cat{N}}(\th) = 0$. 
Let's write 
$P := \opn{Cone}_{A, \cat{M}}(\phi)$
and 
$Q := \opn{Cone}_{B, \cat{N}}(F(\phi))$.
Recall that 
$\d_{B, \cat{N}}(\th) =  \d_{Q} \circ \th - \th \circ \d_{F(P)}$.
We have to prove that this is the zero element in 
$\opn{Hom}_{B, \cat{N}} \bigl( F(P), Q \bigr)^1$.
 
Writing the cones as column modules:
\[ P = \bmat{ M_1 \\[0.2em] \opn{T}_{A, \cat{M}}(M_0) }  
\quad \tup{and}  \quad 
Q = \bmat{ F(M_1) \\[0.2em] \opn{T}_{B, \cat{N}}(F(M_0)) } \, , \]
the matrices representing the morphisms in question are
\[ \th = 
\bmat{\opn{id}_{F(M_1)} & 0 \\[0.2em] 0 & \tau_{F, M_0} } 
\, , \quad 
\d_P = \bmat{\d_{M_1} & \phi \circ \opn{t}_{M_0}^{-1} 
\\[0.2em] 0 & \d_{\opn{T}_{A, \cat{M}}(M_0)} } \]
and \[  
\d_Q = \bmat{\d_{F(M_1)} & F(\phi) \circ \opn{t}_{F(M_0)}^{-1} 
\\[0.2em] 0 & \d_{\opn{T}_{B, \cat{N}}(F(M_0))} } \ . \] 

Let us write $\ga := \ga_F$ for simplicity. 
According to Theorem \ref{thm:1225}, the gauge 
$\ga : F \to F$ is a degree $1$ morphism of DG functors
$\dcat{C}(A, \cat{M}) \to \dcat{C}(B, \cat{N})$.
Because the decomposition (\ref{eqn:1197}) is in the 
category $\dcat{G}_{\mrm{str}}(A, \cat{M})$, Lemma \ref{lem:1200} tells us that 
$\ga_P$ decomposes too, i.e.\ 
\[ \ga_P =  \bmat{ \ga_{M_1} & 0 \\[0.2em] 0 & \ga_{\opn{T}_{A, \cat{M}}(M_0)} 
} . \]

By definition of $\ga_P$ we have 
$\d_{F(P)} = F(\d_P) + \ga_P$ in 
$\opn{Hom}_{B, \cat{N}} \bigl( F(P), F(P) \bigr)^1$.
It follows that 
\[ \begin{aligned}
& \d_{F(P)} = F(\d_P) + \ga_P  
\\[0.2em]
& \qquad 
= \bmat{ F(\d_{M_1})  & 
F(\phi \circ \opn{t}_{M_0}^{-1}) 
\\[0.2em] 0 & F(\d_{\opn{T}_{A, \cat{M}}(M_0)})  } +
\bmat{ \ga_{M_1} & 0 \\[0.2em] 0 & \ga_{\opn{T}_{A, \cat{M}}(M_0)} } 
\\[0.2em]
& \qquad 
= \bmat{ F(\d_{M_1}) +  \ga_{M_1} & 
F(\phi \circ \opn{t}_{M_0}^{-1}) 
\\[0.2em] 0 & 
F(\d_{\opn{T}_{A, \cat{M}}(M_0)}) + \ga_{\opn{T}_{A, \cat{M}}(M_0)} }
\\[0.2em]
& \qquad 
= \bmat{ \d_{F(M_1)} & 
F(\phi \circ \opn{t}_{M_0}^{-1}) 
\\[0.2em] 0 & \d_{F(\opn{T}_{A, \cat{M}}(M_0))} } .
\end{aligned}   \]

Finally we will check that 
$\th \circ \d_{F(P)}$ and $ \d_{Q} \circ \th$ are equal as matrices of 
morphisms. We do that in each matrix entry separately. The two left 
entries in the matrices $\th \circ \d_{F(P)}$ and $ \d_{Q} \circ \th$
agree trivially. The bottom right entries in these matrices
are $\tau_{F, M_0} \circ \d_{F(\opn{T}_{A, \cat{M}}(M_0))}$
and 
$\d_{\opn{T}_{B, \cat{N}}(F(M_0))} \circ \tau_{F, M_0}$
respectively; they are equal by Proposition \ref{prop:1179}. 
And in the top right  entries we have 
$F(\phi \, \circ \opn{t}_{M_0}^{-1})$
and 
$F(\phi) \, \circ \opn{t}_{F(M_0)}^{-1} \circ \, \tau_{F, M_0}$
respectively. Now
$F(\phi \, \circ \opn{t}_{M_0}^{-1}) = F(\phi) \circ F(\opn{t}_{M_0}^{-1})$;
so it suffices to prove that 
$F(\opn{t}_{M_0}^{-1}) = \opn{t}_{F(M_0)}^{-1} \circ \, \tau_{F, M_0}$.
This is immediate from the definition of $\tau_{F, M_0}$.

\medskip \noindent 
(2) By definition of $\th = \opn{cone}(F, \phi)$, the diagram is commutative in 
$\dcat{G}_{\mrm{str}}(A, \cat{M})$. But by part (1) we know that all morphisms 
in it commute with the differentials, so they lie in 
$\dcat{C}_{\mrm{str}}(B, \cat{N})$.
And the functor $\opn{Und}$ from Proposition \ref{prop:2365}
is faithful. 
\end{proof}

\begin{cor} \label{cor:1915}
In the situation of Theorem \tup{\ref{thm:1185}}, the diagram 
\[ \UseTips  \xymatrix @C=6.2ex @R=6ex {
F(M_0) 
\ar[r]^{F(\phi)}
\ar[d]_{\opn{id}}
& 
F(M_1) 
\ar[r]^(0.4){F(e_{\phi})}
\ar[d]_{\opn{id}}
& 
F(\opn{Cone}_{A, \cat{M}}(\phi))
\ar[rr]^{\tau_{F, M_0} \, \circ \, F(p_{\phi})}
\ar[d]_{\opn{cone}(F, \phi)}
& &
\opn{T}_{B, \cat{N}}(F(M_0))
\ar[d]_{\opn{id}}
\\
F(M_0)
\ar[r]^{F(\phi)}
& 
F(M_1)
\ar[r]^(0.4){e_{F(\phi)}}
& 
\opn{Cone}_{B, \cat{N}}(F(\phi))
\ar[rr]^{p_{F(\phi)}}
& &
\opn{T}_{B, \cat{N}}(F(M_0))
} \]
is an isomorphism of triangles in $\dcat{C}_{\mrm{str}}(B, \cat{N})$. 
\end{cor}

\begin{proof}
Just rearrange the diagram in item (2) of the theorem. 
\end{proof}

\mysubsection{Examples of DG Functors} \label{subsec:exa-DGfun}

Recall that $\cat{M}$ and $\cat{N}$ are abelian categories, and $A$ 
and $B$ are DG rings.
Here are four examples of DG functors, of various types. 
We work out in detail the translation isomorphism, the cone isomorphism and 
the gauge in each example. These examples 
should serve as templates for constructing other DG functors. 

\begin{exa} \label{exa:1117}
Here $A = B = \K$, so 
$\dcat{C}(A, \cat{M}) = \dcat{C}(\cat{M})$
and $\dcat{C}(B, \cat{N}) = \dcat{C}(\cat{N})$. 
Let $F : \cat{M} \to \cat{N}$ be a $\K$-linear functor. It extends to a functor 
$\dcat{C}(F) : \dcat{C}(\cat{M}) \to \dcat{C}(\cat{N})$
as follows: on objects, a complex 
\[ M = \bigl( \{ M^i \}_{i \in \Z}, \{ \d_M^i \}_{i \in \Z} \bigr) 
\in \dcat{C}(\cat{M}) \]
goes to the complex 
\[ \dcat{C}(F)(M) := \bigl( \{ F(M^i) \}, \{ F(\d_M^i) \} \bigr) \in 
\dcat{C}(\cat{N}) . \]
A morphism $\phi = \{ \phi^j \}_{j \in \Z}$ in $\dcat{C}(\cat{M})$ goes 
to the morphism 
$\dcat{C}(F)(\phi) : = \lb \{ F(\phi^j) \}_{j \in \Z}$
in $\dcat{C}(\cat{N})$. A slightly tedious calculation shows that $\dcat{C}(F)$ 
is a graded functor.

Given a complex $M \in \dcat{C}(\cat{M})$, let 
$N := \dcat{C}(F)(M) \in \dcat{C}(\cat{N})$.
Then the translation of $N$ is 
$\opn{T}_{\cat{N}}(N) = \dcat{C}(F) \bigl( \opn{T}_{\cat{M}}(M) \bigr)$;
and the little t operator of $N$ is 
$\opn{t}_{N} = \dcat{C}(F)(\opn{t}_{M})$. 
So the translation isomorphism 
\[ \tau_{\dcat{C}(F)} :  \dcat{C}(F) \circ \opn{T}_{\cat{M}} \iso 
\opn{T}_{\cat{N}} \circ \, \dcat{C}(F) \]
of functors 
$\dcat{C}_{\mrm{str}}(\cat{M}) \to \dcat{C}_{\mrm{str}}(\cat{N})$
is the identity automorphism of this functor. 

Let $\phi : M_0 \to M_1$ be a morphism in 
$\dcat{C}_{\mrm{str}}(\cat{M})$, 
whose image under $\dcat{C}(F)$ is the morphism 
$\psi : N_0 \to N_1$  in $\dcat{C}_{\mrm{str}}(\cat{N})$.
Then 
\[ \opn{Cone}(\psi) = N_1 \oplus \opn{T}_{\cat{N}}(N_0) =
\dcat{C}(F) \bigl( \opn{Cone}(\phi) \bigr)  \]
as graded objects in $\cat{N}$, with differential 
\[ \d_{\opn{Cone}(\psi)} = 
\bmat{\d_{N_1} & \psi \circ \opn{t}_{N_0}^{-1} \\[0.2em] 0 & \d_{\opn{T}(N_0)}} 
= \dcat{C}(F) \Bigl( 
\bmat{\d_{M_1} & \phi \circ \opn{t}_{M_0}^{-1} \\[0.2em] 0 & \d_{\opn{T}(M_0)}}
\Bigr) = \dcat{C}(F) \bigl( \d_{\opn{Cone}(\phi)} \bigr) . \]
We see that the cone isomorphism $\opn{cone}(F, \phi) $ is the identity 
automorphism of the DG module $\opn{Cone}(\psi)$. The gauge $\ga_{\dcat{C}(F)}$ 
of the graded functor $\dcat{C}(F)$ is zero. Therefore, by Theorem 
\ref{thm:1225}, $\dcat{C}(F)$ is a DG functor. 
\end{exa}

The next example is much more complicated, and we work out the full details 
(only once -- later on, such details will be left to the reader). 

\begin{exa} \label{exa:1118}
Let $A$ and $B$  be DG rings, and fix some 
$N \in \cat{DGMod} B \ot A^{\mrm{op}} = \dcat{C}(B \ot A^{\mrm{op}})$.
In other words, $N$ is a DG $B$-$A$-bimodule. 
For every $M \in \dcat{C}(A)$ we have a DG $\K$-module
$F(M) := N \ot_{A} M$,
as in Definition \ref{dfn:1100}. 
The differential of $F(M)$ is
\begin{equation} \label{eqn:1220}
\d_{F(M)} = \d_N \ot \opn{id}_M + \opn{id}_N \ot \, \d_M . 
\end{equation}
See formula (\ref{eqn:4105}) 
regarding the Koszul sign rule that's involved.
But $F(M)$ has the structure of a DG $B$-module: for every $b \in B$, 
$n \in N$ and $m \in M$, the action is 
$b \cd (n \ot m) := (b \cd n) \ot m$.
Clearly 
\[ F : \dcat{C}(A) = \cat{DGMod} A \to \dcat{C}(B) = \cat{DGMod} B \]
is a linear functor. We will show that it is actually a DG functor. 

Let $M_0, M_1 \in \dcat{C}(A)$, and consider the $\K$-linear homomorphism 
\begin{equation} \label{eqn:1120}
F : \opn{Hom}_A(M_0, M_1) \to 
\opn{Hom}_B(N \ot_{A} M_0, N \ot_{A} M_1) .
\end{equation}
Take any $\phi \in \opn{Hom}_A(M_0, M_1)^i$. Then 
\[  F(\phi) \in \opn{Hom}_B(N \ot_{A} M_0, N \ot_{A} M_1) \]
is the homomorphism that on a homogeneous tensor 
$n \ot m \in (N \ot_{A} M_0)^{k + j}$, 
with  $n \in N^k$ and $m \in M_0^j$, has the value
\[ F(\phi) (n \ot m) = (-1)^{i \cd k} \cd n \ot \phi(m) \in 
(N \ot_{A} M_1)^{k + j + i} .  \]
In other words, 
\begin{equation}  \label{eqn:1222}
F(\phi) = \opn{id}_N \ot \, \phi . 
\end{equation}
We see that the homomorphism $F(\phi)$ has degree $i$. So
$F$ is a graded functor. 

Let us calculate $\ga_F$, the gauge of $F$. 
From (\ref{eqn:1222}) and (\ref{eqn:1220}) we get 
$\ga_{F, M} = \d_N \ot \opn{id}_M$,
which is often a nonzero endomorphism of $F(M)$. Still, take any degree $i$ 
morphism $\phi : M_0 \to M_1$ 
in $\dcat{C}(A)$. Then 
\[ \begin{aligned}
& \ga_{M_1} \circ F(\phi) = 
(\d_N \ot \opn{id}_{M_1}) \circ (\opn{id}_N \ot \, \phi)
\\ & \qquad
= \d_N \ot \phi = 
(-1)^i \cd (\opn{id}_N  \ot \, \phi) \circ (\d_N \ot \opn{id}_{M_0})  
= (-1)^i \cd F(\phi) \circ \ga_{M_0} . 
\end{aligned} \]
We see that $\ga_F$ satisfies the condition of Definition \ref{dfn:1195}(1), 
which is really Definition \ref{dfn:1695}.
By Theorem \ref{thm:1225}, $F$ is a DG functor. (It is also possible to 
calculate directly that $F$ is a DG functor.) 

Finally let us figure out what is the translation isomorphism $\tau_F$ of the 
functor $F$. Take $M \in \dcat{C}(A)$. 
Then 
$\tau_{F, M} : F(\opn{T}_{A}(M)) \to \opn{T}_{B}(F(M))$
is an isomorphism in $\dcat{C}_{\mrm{str}}(B)$. 
By Definition \ref{dfn:1150} we have 
$\tau_{F, M} = \opn{t}_{F(M)} \circ \, F(\opn{t}_M)^{-1}$.
Take any $n \in N^k$ and $m \in M^{j + 1}$, so that 
$n \ot \opn{t}_M(m) \in (N \ot_A \opn{T}_{A}(M))^{k + j} 
= F(\opn{T}_{A}(M))^{k + j}$,
a degree $k + j$ element of $F(\opn{T}_{A}(M))$. But 
\[ n \ot \opn{t}_M(m) = (-1)^k \cd (\opn{id}_N \ot \opn{t}_M)(n \ot m) = 
(-1)^k \cd F(\opn{t}_M)(n \ot m) . \]
Therefore 
\[ \tau_{F, M}(n \ot \opn{t}_M(m)) = 
(-1)^k \cd \opn{t}_{F(M)}(n \ot m) \in \opn{T}_{B}(F(M))^{k + j} .  \]
Observe that when $N$ is concentrated in degree $0$, we are back in the 
situation of Example \ref{exa:1117}, in which there are no sign twists, and 
$\tau_{F, M}$ is the identity automorphism. 
\end{exa}

\begin{exa} \label{exa:1120}
Let $A$ and $B$  be DG rings, and fix some 
$N \in \cat{DGMod} A \ot B^{\mrm{op}} = \dcat{C}(B \ot A^{\mrm{op}})$.
For any $M \in \cat{DGMod} A$ we define 
$F(M) := \opn{Hom}_A(N, M)$.
This is a DG $B$-module: for every $b \in B^i$ and
$\phi \in \opn{Hom}_A(N, M)^j$, the homomorphism 
$b \cd \phi \in \opn{Hom}_A(N, M)^{i + j}$ has value
$(b \cd \phi)(n) := (-1)^{i \cd (j + k)} \cd \phi(n \cd b) \in M^{i + j + k}$
on $n \in N^k$. As in the previous example, 
\[ F : \dcat{C}(A) = \cat{DGMod} A \to \dcat{C}(B) = \cat{DGMod} B \]
is a $\K$-linear graded functor. 

The value of the gauge $\ga_F$ at $M \in \dcat{C}(A)$  is 
$\ga_{F, M} = \opn{Hom}(\d_N, \opn{id}_M)$.
See formula (\ref{eqn:4106})  
regarding this notation. Namely for 
$\psi \in F(M)^j = \opn{Hom}_A(N, M)^j$
we have 
$\ga_{F, M}(\psi) = (-1)^{\lms j} \cd \psi \circ \d_N$.
It is not too hard to check that $\ga_F$ is a degree $1$ morphism of graded 
functors. Hence, by Theorem \ref{thm:1225}, $F$ is a DG functor. 

The formula for the translation isomorphism $\tau_F$ is as follows. 
Take $M \in \dcat{C}(A)$. Then 
\[ \tau_{F, M} : F(\opn{T}_{A}(M)) = \opn{Hom}_A(N, \opn{T}_{A}(M))
\to \opn{T}_{B}(F(M)) = \opn{T}_{B}(\opn{Hom}_A(N, M)) \]
is, by definition, 
$\tau_{F, M} = \opn{t}_{F(M)} \circ \, F(\opn{t}_{M})^{-1}$.
Now
$F(\opn{t}_{M})^{-1} = \opn{Hom}(\opn{id}_N, \opn{t}_{M}^{-1})$.
So given any 
$\psi \in F(\opn{T}_{A}(M))^k$, we have 
$\tau_{F, M}(\psi) = \opn{t}_{F(M)}(\opn{t}_{M}^{-1} \circ \, \psi) \in 
\opn{T}_{B}(F(M))^k$.
\end{exa}

We end with a contravariant example. 

\begin{exa} \label{exa:2300}
Let $A$ be a commutative ring. Fix some complex $N \in \dcat{C}(A)$. 
For any $M \in \dcat{C}(A)$ let 
$F(M) := \opn{Hom}_A(M, N) \in \dcat{C}(A)$.
For every degree $i$ morphism 
$\phi : M_0 \to M_1$ in $\dcat{C}(A)$ let
$F(\phi) : F(M_1) \to F(M_0)$ be the degree $i$ morphism
$F(\phi) := \opn{Hom}_A(\phi, \opn{id}_N)$.

A direct calculation, that we leave to the reader, shows that 
\[ F : \opn{Hom}_{\dcat{C}(A)}(M_0, M_1) \to 
\opn{Hom}_{\dcat{C}(A)} \bigl( F(M_1), F(M_0) \bigr) \]
is a strict homomorphism of DG $\K$-modules. It also satisfies conditions (a) 
and (b) of Definition \ref{dfn:2495}. Thus we have a contravariant DG 
functor 
$F : \dcat{C}(A) \to \dcat{C}(A)$.
For contravariant DG functors we do not talk about translation isomorphisms or 
gauges. 
\end{exa}

We will return to the {\em DG bifunctors} $\opn{Hom}(-, -)$ and 
$(- \ot -)$ later in the book, in Subsection \ref{subsec:dg-bifun}.

%% file: block2_190413.tex

\renewcommand{\thisfile}{block2\_190328}  
  
\cleardoublepage
\mysection{Triangulated Categories and Functors} 
\label{sec:triangulated}  

\AYcopyright

In this section we introduce triangulated categories and triangulated 
functors. We prove that for a DG ring $A$ and an abelian category  
$\cat{M}$, the homotopy category $\dcat{K}(A, \cat{M})$
is triangulated. We also prove that a DG functor 
$F : \dcat{C}(A, \cat{M}) \to \dcat{C}(B, \cat{N})$
induces a triangulated functor 
$F : \dcat{K}(A, \cat{M}) \to \dcat{K}(B, \cat{N})$.

Recall that by Convention \ref{conv:2490}, there is a nonzero 
commutative base ring $\K$, which is implicit most of the time. All rings are 
central $\K$-rings, and all ring homomorphisms are over $\K$. All linear 
categories are $\K$-linear, and all 
linear (i.e.\ additive) functors between them are $\K$-linear.

\mysubsection{Triangulated Categories}

In this subsection we define triangulated categories, and make some remarks 
regarding them. 
  
Recall that a functor is called an isomorphism of categories if it is bijective 
on sets of objects and on sets of morphisms. 

\begin{dfn} \label{dfn:1155} 
\index{Translation functor! of an additive category}
\index{T-additive! category}
Let $\cat{K}$ be an additive category. A {\em translation functor} on $\cat{K}$ 
is an additive automorphism $\opn{T}$ of $\cat{K}$.
The pair $(\cat{K}, \opn{T})$ is called a {\em T-additive category}. 
\end{dfn}

\begin{rem}
Some texts give a more relaxed definition: the translation functor $\opn{T}$ is 
only required to be an 
additive auto-equivalence of $\cat{K}$. The resulting theory is more 
complicated (it is $2$-categorical, but most texts try to suppress this fact). 
\end{rem}

Later in the book (see Definition \ref{dfn:2120}) we will write 
$M[k] := \opn{T}^k(M)$, the $k$-th translation of an object $M$. 

\begin{dfn} \label{dfn:1159} 
\index{T-additive! functor}
\index{Translation isomorphism! of an additive functor}
Suppose $(\cat{K}, \opn{T}_{\cat{K}})$ and 
$(\cat{L}, \opn{T}_{\cat{L}})$ are T-additive categories.
A {\em T-additive functor} between them is a pair 
$(F, \tau)$, consisting of an additive functor $F : \cat{K} \to \cat{L}$, 
together with an isomorphism
\[ \tau :  F \circ \opn{T}_{\cat{K}}  \iso 
\opn{T}_{\cat{L}} \circ \, F   \]
of functors $\cat{K} \to \cat{L}$,
called a {\em translation isomorphism}. 
\end{dfn}

\begin{dfn} \label{dfn:1840}
Let $(\cat{K}_i, \opn{T}_i)$ be T-additive categories, for $i = 0, 1, 2$, and 
let 
\[ (F_i, \tau_i) : (\cat{K}_{i - 1}, \opn{T}_{i - 1}) \to 
(\cat{K}_i, \opn{T}_i) \]
be T-additive functors. The composition 
\[ (F, \tau) = (F_2, \tau_2) \circ (F_1, \tau_1) \]
is the T-additive functor 
$(\cat{K}_0, \opn{T}_0) \to (\cat{K}_2, \opn{T}_2)$
defined as follows: the functor is
$F := F_2 \circ F_1$, and the translation isomorphism 
$\tau : F \circ T_0 \iso T_2 \circ F$ 
is 
$\tau := \tau_2 \circ F_2(\tau_1)$. 
\end{dfn}

\begin{dfn} \label{dfn:1160} 
Suppose $(\cat{K}, \opn{T}_{\cat{K}})$ and 
$(\cat{L}, \opn{T}_{\cat{L}})$ are T-additive categories,
and 
\[ (F, \tau), (G, \nu) : (\cat{K}, \opn{T}_{\cat{K}}) \to 
(\cat{L}, \opn{T}_{\cat{L}}) \] 
are T-additive functors. A {\em morphism of T-additive functors}
\index{T-additive! morphism of {\indash} functors}
\[ \eta : (F, \tau) \to (G, \nu) \] 
is a morphism of functors $\eta : F \to G$, such that for every object 
$M \in \cat{K}$ this diagram in $\cat{L}$ is commutative:
\[ \UseTips  \xymatrix @C=8ex @R=6ex {
F(\opn{T}_{\cat{K}}(M))
\ar[r]^{\tau_M}
\ar[d]_{\eta_{\lms \opn{T}_{\cat{K}}(M)}}
&
\opn{T}_{\cat{L}}(F(M))
\ar[d]^{\opn{T}_{\cat{L}}(\eta_M)}
\\
G(\opn{T}_{\cat{K}}(M))
\ar[r]^{\nu_M}
&
\opn{T}_{\cat{L}}(G(M)) \ . 
} \]
\end{dfn}

\begin{dfn} \label{dfn:1156} 
Let $(\cat{K}, \opn{T})$ be a T-additive category. 
A {\em triangle} 
\index{Triangle! in a T-additive category}
in $(\cat{K}, \opn{T})$ is a diagram 
\[ L \xar{\al} M \xar{\be} N \xar{\ga} \opn{T}(L)  \]
in $\cat{K}$.
\end{dfn}

Sometimes, when the names of the morphisms in the triangle above are not 
important, we may use the more compact notation 
\begin{equation} \label{eqn:4910}
L \xar{} M \xar{} N \xar{ \, \triangle \, } .
\end{equation}

\begin{dfn} \label{dfn:1157} 
Let $(\cat{K}, \opn{T})$ be a T-additive category. Suppose 
\[  L \xar{\al} M \xar{\be} N \xar{\ga} \opn{T}(L) \]
and 
\[ L' \xar{\al'} M' \xar{\be'} N' \xar{\ga'} \opn{T}(L') \]
are triangles in  $(\cat{K}, \opn{T})$. 
A {\em morphism of triangles}
\index{Triangle! morphism of {\indash}s}
between them is a commutative diagram 
\[ \UseTips  \xymatrix @C=8ex @R=6ex { 
L 
\ar[r]^{\al}
\ar[d]_{\phi}
&
M
\ar[r]^{\be}
\ar[d]_{\psi}
& 
N
\ar[r]^{\ga}
\ar[d]_{\chi}
&
\opn{T}(L)
\ar[d]_{\opn{T}(\phi)}
\\
L' 
\ar[r]^{\al'}
&
M'
\ar[r]^{\be'}
& 
N'
\ar[r]^{\ga'}
&
\opn{T}(L') 
} \]
in $\cat{K}$. 

The morphism of triangles $(\phi, \psi, \chi)$ is called an isomorphism if 
$\phi, \psi$ and $\chi$ are all isomorphisms. 
\end{dfn}

\begin{rem} \label{rem:4690}
Why ``triangle''? This is because sometimes a triangle 
\[ L \xar{\al} M \xar{\be} N \xar{\ga} \opn{T}(L) \]
is written as a diagram
\[ \UseTips  \xymatrix @C=4ex @R=4ex {
& 
N 
\ar[dl]_{\ga}
\\
L
\ar[rr]_{\al}
& &
M
\ar[lu]_{\be}
} \]
See \cite[page 20]{RD}. 
\end{rem}

\begin{dfn} \label{dfn:1158} 
A {\em triangulated category}
\index{Triangulated category}
is a T-additive category
$(\cat{K}, \opn{T})$, \lb equipped with a set of triangles called  {\em 
distinguished triangles}. 
\index{Triangle! distinguished}
The following axioms have to be satisfied:
\begin{enumerate}[leftmargin=30pt]
\item[(TR1)] 
\begin{itemize}
\rmitem{a} Every triangle that is isomorphic to a distinguished
triangle is also a distinguished triangle. 

\rmitem{b} For every morphism $\al : L \to M$ in $\cat{K}$ there is a
distinguished triangle 
\[ L \xar{\al} M \xar{} N \xar{} \opn{T}(L) . \] 

\rmitem{c}  For every object $M$ the triangle 
\[ M \xar{\opn{id}_M} M \to 0 \to \opn{T}(M) \] 
is distinguished.
\end{itemize}

\item[(TR2)] A triangle 
\[ L \xar{\al} M \xar{\be} N \xar{\ga} \opn{T}(L) \]
is distinguished iff the triangle 
\[ M \xar{\be} N \xar{\ga} \opn{T}(L) \xar{- \opn{T}(\al)} \opn{T}(M) \] 
is distinguished.

\item[(TR3)]  Suppose 
\[ \UseTips  \xymatrix @C=8ex @R=6ex { 
L 
\ar[r]^{\al}
\ar[d]_{\phi}
&
M
\ar[r]^{\be}
\ar[d]_{\psi}
& 
N
\ar[r]^{\ga}
&
\opn{T}(L)
\\
L' 
\ar[r]^{\al'}
&
M'
\ar[r]^{\be'}
& 
N'
\ar[r]^{\ga'}
&
\opn{T}(L') 
} \]
is a commutative diagram in $\cat{K}$ in which the rows are distinguished 
triangles. 
Then there exists a morphism $\chi : N \to N'$ such that the diagram
\[ \UseTips  \xymatrix @C=8ex @R=6ex { 
L 
\ar[r]^{\al}
\ar[d]_{\phi}
&
M
\ar[r]^{\be}
\ar[d]_{\psi}
& 
N
\ar[r]^{\ga}
\ar[d]_{\chi}
&
\opn{T}(L)
\ar[d]_{\opn{T}(\phi)}
\\
L' 
\ar[r]^{\al'}
&
M'
\ar[r]^{\be'}
& 
N'
\ar[r]^{\ga'}
&
\opn{T}(L') \ .
} \]
is a morphism of triangles.

\item[(TR4)] Suppose we are given these three distinguished triangles:
\[ L \xar{\al} M \xar{\ga} P \xar{} \opn{T}(L) , \]
\[ M \xar{\be} N \xar{\ep} R \xar{} \opn{T}(M) , \]
\[ L \xar{\be \circ \al} N \xar{\de} Q \xar{} \opn{T}(L) . \]
Then there is a distinguished triangle 
\[  P \xar{\phi} Q \xar{\psi} R \xar{\rho} \opn{T}(P)  \]
making diagram (\ref{eqn:5022}) commutative. 
\begin{figure}[!htb]
\centering
\begin{equation} \label{eqn:5022}
\UseTips  \xymatrix @C=8ex @R=6ex { 
L 
\ar[r]^{\al}
\ar[d]_{\opn{id}}
&
M
\ar[r]^{\ga}
\ar[d]_{\be}
& 
P
\ar[r]^{}
\ar[d]_{\phi}
&
\opn{T}(L)
\ar[d]_{\opn{id}}
\\
L
\ar[r]^{\be \circ \al}
\ar[d]_{\al}
&
N
\ar[r]^{\de}
\ar[d]_{\opn{id}}
& 
Q
\ar[r]^{}
\ar[d]_{\psi}
&
\opn{T}(L)
\ar[d]_{\opn{T}(\al)}
\\
M
\ar[r]^{\be}
\ar[d]_{\ga}
&
N
\ar[r]^{\ep}
\ar[d]_{\de}
& 
R
\ar[r]^{}
\ar[d]_{\opn{id}}
&
\opn{T}(M)
\ar[d]_{\opn{T}(\ga)}
\\
P
\ar[r]^{\phi}
&
Q
\ar[r]^{\psi}
& 
R
\ar[r]^{\rho}
&
\opn{T}(P)
} 
\end{equation}
\end{figure}
\end{enumerate}
\end{dfn}

Here are a few remarks on triangulated categories. See Remark \ref{exer:3015} 
regarding triangulated categories arising from algebraic topology. 

\begin{rem} \label{rem:1280}
The numbering of the axioms we use is taken from \cite{RD}, and it agrees with 
the numbering in \cite{We}. The numbering in \cite{KaSc1} and \cite{KaSc2}  
is different.

In the situation that we care about, namely 
$\cat{K} = \dcat{K}(A, \cat{M})$,
the distinguished triangles will be those triangles that are isomorphic, in 
$\dcat{K}(A, \cat{M})$, to the standard triangles in $\dcat{C}(A, \cat{M})$
from Definition \ref{dfn:1140}. 
See Definition \ref{dfn:140} below for the precise statement. 

The object $N$ in item (b) of axiom (TR1) is referred to as a {\em cone} of the 
morphism $\al$. We should think of the cone as something combining 
``the cokernel'' and ``the kernel'' of $\al$. 
Axiom (TR3) guarantees a sort of weak functoriality of the cone. 
See Subsection \ref{subsec:props-tr-cats}, and in particular 
Remark \ref{rem:4160}, for more properties of the cone.

Axiom (TR2) says that if we ``turn'' a distinguished triangle (cf.\ Remark 
\ref{rem:4690}) we remain with a distinguished triangle.
\end{rem}

\begin{rem} \label{rem:1281}
The axiom (TR4) is called the {\em octahedral axiom}.
\index{Octahedral axiom}
It is supposed to replace the isomorphism  
\[ (N / L) / (M / L) \cong N / M \]
for objects $L \sub M \sub N$ in an abelian category $\cat{M}$. 
The octahedral axiom is needed for the theory of {\em t-structures}: it is 
used, in \cite{BBD}, to show that the heart of a t-structure is an abelian 
category (see Remark \ref{rem:3600}). This axiom is also needed to form 
Verdier quotients of triangulated categories (see Remark \ref{rem:4220}).
The book \cite{Ne1} has detailed discussions of the octahedral axiom. 

A T-additive category $(\cat{K}, \opn{T})$ that only satisfies axioms 
(TR1)-(TR3) is called a {\em pretriangulated category}. 
(The reader should not confuse ``pretriangulated category'', as used here, with 
the ``pretriangulated DG category'' from \cite{BoKa}; see Remark 
\ref{rem:1475}.) It is not known whether the octahedral axiom is a consequence 
of the other axioms; there was a recent paper by A. Maccioca 
(arxiv:1506.00887) claiming that, but it had a fatal error in it.

In our book the octahedral axiom does not play any role. For this reason we had 
excluded it from an earlier version of the book, in which we had discussed 
pretriangulated categories only. Our decision to include this axiom in the 
current version of the book, and thus to talk about triangulated categories 
(rather than about pretriangulated ones) is just to be more in line with the 
mainstream usage. With the exception of a longer proof of Theorem 
\ref{thm:1255} -- stating that $\dcat{K}(A, \cat{M})$ is a triangulated 
category -- there is virtually no change in the content of the book, and almost 
all definitions and results are valid for pretriangulated categories.  
\end{rem}

\begin{rem} \label{rem:2305}
The structure of triangulated categories is poorly understood, and there is no
classification of triangulated functors between them. This is true even 
for the derived category $\dcat{D}(A)$ of a ring $A$. For instance, a famous 
open question of J. Rickard is whether every triangulated auto-equivalence $F$ 
of $\dcat{D}(A)$ is isomorphic to $P \ot^{\mrm{L}}_{A} (-)$ for some tilting 
complex $P$. More on this issue in Section \ref{sec:perf-tilt-NC}. 

In this book the role of the triangulated structure on derived categories is 
secondary. It is mostly used for induction on the amplitude of the cohomology 
of complexes; cf.\ Section \ref{sec:adj-equ-cohdim}. 

Furthermore, one of the most important functors studied here -- the squaring 
operation $\opn{Sq}_{B / A}$ from Subsections \ref{subsec:squaring} and
\ref{subsec:RNCDC-uniq} -- is a functor from $\dcat{D}(B)$ to itself, which is 
not triangulated, nor even linear: it is a {\em quadratic functor}. 
\end{rem}

\mysubsection{Triangulated and Cohomological Functors} 
\label{subsec:tr-coh-funcs}

Suppose $\cat{K}$ and $\cat{L}$ are T-additive categories, with 
translation functors
$\opn{T}_{\cat{K}}$ and $\opn{T}_{\cat{L}}$ respectively. 
A T-additive functor $F : \cat{K} \to \cat{L}$ was defined in Definition 
\ref{dfn:1159}. A morphism of T-additive functors $\eta : F \to G$ was 
introduced in Definition \ref{dfn:1160}.

\begin{dfn} \label{dfn:1276}
Let  $\cat{K}$ and $\cat{L}$  be triangulated categories. 
\begin{enumerate}
\item A {\em triangulated functor} 
\index{Triangulated functor}
from $\cat{K}$ to $\cat{L}$
is a T-additive functor 
$(F, \tau) : \cat{K} \to \cat{L}$ 
that satisfies this condition: for every distinguished triangle 
\[ L \xar{\al} M \xar{\be} N \xar{\ga} \opn{T}_{\cat{K}}(L) \]
in $\cat{K}$, the triangle 
\[ F(L) \xar{F(\al)} F(M)  \xar{F(\be)} F(N)  
\xar{ \ \tau_L \circ \, F(\ga) \ } \opn{T}_{\cat{L}}(F(L)) \]
is a distinguished triangle in $\cat{L}$.

\item Suppose $(G, \nu) : \cat{K} \to \cat{L}$ is another triangulated functor.
A {\em morphism of triangulated functors}
\index{Triangulated functor! morphism of {\indash}s}
$\eta : (F, \tau) \to (G , \nu)$ is 
a morphism of T-additive functors, as in Definition \ref{dfn:1160}.
\end{enumerate}
\end{dfn}

Sometimes we keep the translation isomorphism $\tau$ implicit, and refer to $F$ 
as a triangulated functor.

\begin{dfn}  \label{dfn:1450}
Let $\cat{K}$ be a triangulated category. A 
{\em full triangulated subcategory}%
\index{Triangulated category! full triangulated subcategory of}
of $\cat{K}$ is a subcategory 
$\cat{L} \subseteq \cat{K}$ satisfying these conditions:
\begin{enumerate}
\rmitem{a} $\cat{L}$ is a full additive subcategory (see Definition 
\ref{dfn:1025}). 

\rmitem{b}
$\cat{L}$ is closed under translations, i.e.\ $L \in \cat{L}$ iff
$\opn{T}(L) \in \cat{L}$.

\rmitem{c} $\cat{L}$ is closed under distinguished triangles, i.e.\ if 
$L' \to L \to L'' \xar{\, \triangle\, }$ 
is a  distinguished triangle in $\cat{K}$ s.t.\ $L', L \in \cat{L}$, then 
also $L'' \in \cat{L}$.
\end{enumerate}
\end{dfn}

Observe that $\cat{L}$ itself is triangulated, and the inclusion 
$\cat{L} \to \cat{K}$ is a triangulated functor. 

\begin{prop} \label{prop:2305}
For $i = 0, 1, 2$ let $(\cat{K}_i, \opn{T}_i)$ be triangulated categories, and 
let 
$(F_i, \tau_i) : (\cat{K}_{i - 1}, \opn{T}_{i - 1}) \to 
(\cat{K}_i, \opn{T}_i)$ 
be triangulated functors. Define the T-additive functor 
$(F, \tau) := (F_2, \tau_2) \circ (F_1, \tau_1)$
as in Definition \tup{\ref{dfn:1840}}. Then 
$(F, \tau) : (\cat{K}_0, \opn{T}_0) \to (\cat{K}_2, \opn{T}_2)$
is a triangulated functor. 
\end{prop}

\begin{exer} \label{exer:2300}
Prove Proposition \ref{prop:2305}. 
\end{exer}

\begin{dfn} \label{dfn:1510}
Let $\cat{K}$ be a triangulated category, and let $\cat{M}$ be an 
abelian category. 
\begin{enumerate}
\item A {\em cohomological functor}%
\index{Cohomological functor}
$F : \cat{K} \to \cat{M}$ is an
additive functor, such that for every distinguished triangle 
\[ L \xar{\al} M \xar{\be} N \xar{\ga} \opn{T}(L) \] 
in $\cat{K}$, the sequence 
\[ F(L) \xar{F(\al)} F(M) \xar{F(\be)} F(N) \] 
is exact in $\cat{M}$.

\item A {\em contravariant cohomological functor}%
\index{Cohomological functor! contravariant}
$F : \cat{K} \to \cat{M}$ is 
a contravariant  additive functor, such that for every distinguished triangle 
\[ L \xar{\al} M \xar{\be} N \xar{\ga} \opn{T}(L) \] 
in $\cat{K}$, the sequence 
\[ F(N) \xar{F(\be)} F(M) \xar{F(\al)} F(L) \] 
is exact in $\cat{M}$.
\end{enumerate}
\end{dfn}

\begin{prop} \label{prop:1281}
Let  $F : \cat{K} \to \cat{M}$ be a cohomological functor, and let 
\[ L \xar{\al} M \xar{\be} N \xar{\ga} \opn{T}(L) \]
be a distinguished triangle in $\cat{K}$. Then the sequence 
\[ \begin{aligned}
&    \cdots \to F(\opn{T}^i(L)) \xar{F(\opn{T}^i(\al))} F(\opn{T}^i(M)) 
\xar{F(\opn{T}^i(\be))} F(\opn{T}^i(N)) 
\\[0.2em]
& \hspace{2ex}
\xar{F(\opn{T}^i(\ga))} F(\opn{T}^{i + 1}(L))
\xar{F(\opn{T}^{i + 1}(\al))} F(\opn{T}^{i + 1}(M)) \to \cdots 
\end{aligned} \]
in $\cat{M}$ is exact.
\end{prop}

\begin{proof}
By axiom (TR2) we have these distinguished triangles:
\[ \opn{T}^i(L) \xar{(-1)^{i} \cd \opn{T}^i(\al)} \opn{T}^i(M) 
\xar{(-1)^{i} \cd \opn{T}^i(\be)} \opn{T}^i(N) 
\xar{(-1)^{i} \cd \opn{T}^i(\ga)} \opn{T}^{i + 1}(L) , \]
\[ \opn{T}^i(M) \xar{(-1)^{i} \cd \opn{T}^i(\be)} \opn{T}^i(N) 
\xar{(-1)^{i} \cd \opn{T}^i(\ga)} \opn{T}^{i + 1}(L) 
\xar{(-1)^{i+1} \cd \opn{T}^{i + 1}(\al)} \opn{T}^{i + 1}(M) , \]
\[ \opn{T}^{i}(N) \xar{(-1)^{i} \cd \opn{T}^{i}(\ga)} 
\opn{T}^{i + 1}(L) \xar{(-1)^{i+1} \cd \opn{T}^{i + 1}(\al)} 
\opn{T}^{i + 1}(M) \xar{(-1)^{i+1} \cd 
\opn{T}^{i + 1}(\be)} \opn{T}^{i + 1}(N) . \]
Now use the definition of a cohomological functor, noting that multiplying 
morphisms in an exact sequence by $-1$ preserves exactness.
\end{proof}

\begin{prop} \label{prop:4}
Let $\cat{K}$ be a triangulated category.
For every $P \in \cat{K}$ the functor
\[ \opn{Hom}_{\cat{K}}(P, -) : \cat{K} \to \cat{Mod} \K \]
is a cohomological functor, and the functor 
\[ \opn{Hom}_{\cat{K}}(-, P) : \cat{K} \to \cat{Mod} \K \]
is a contravariant cohomological functor. 
\end{prop}

\begin{proof}
We will prove the covariant statement. The contravariant statement is 
proved similarly, and we leave this to the reader. 

Consider a distinguished triangle 
$L \xar{\al} M \xar{\be} N \xar{\ga} \opn{T}(L)$
in $\cat{K}$. We have to prove that the sequence 
\[ \opn{Hom}_{\cat{K}}(P, L) \xar{ \opn{Hom}(\opn{id}_P, \al) } 
\opn{Hom}_{\cat{K}}(P, M) \xar{ \opn{Hom}(\opn{id}_P, \be) } 
\opn{Hom}_{\cat{K}}(P, N)  \]
in $\cat{Mod} \K$ is exact. In view of Proposition \ref{prop:1280}, all we need 
to show is that for every $\psi : P \to M$ s.t.\ $\be \circ \psi = 0$, there is 
some  $\phi : P \to L$ s.t.\ $\psi = \al \circ \phi$.  In a picture, we must 
show that the diagram below (solid arrows) 
\[ \UseTips  \xymatrix @C=8ex @R=6ex { 
P
\ar[r]^{\opn{id}}
\ar@{-->}[d]_{\phi}
&
P
\ar[r]^{}
\ar[d]_{\psi}
& 
0
\ar[r]^{}
\ar[d]_{}
&
\opn{T}(P)
\ar@{-->}[d]_{\opn{T}(\phi)}
\\
L
\ar[r]^{\al}
&
M
\ar[r]^{\be}
& 
N
\ar[r]^{\ga}
&
\opn{T}(L) \ .
}  \]
can be completed (dashed arrow). 
This is true by (TR2) (= turning) and (TR3) (= extending). 
\end{proof}

\begin{exer} \label{exer:2301}
Prove the contravariant statement in the proposition above. 
\end{exer}

\begin{que} \label{eqn:1305}
Let  $\cat{K}$ and $\cat{L}$  be triangulated categories, and let 
$F : \cat{K} \to \cat{L}$ be an additive functor. Is it true that there is at 
most one isomorphism of functors 
$\tau :  F \circ \opn{T}_{\cat{K}}  \iso \opn{T}_{\cat{L}} \circ \, F$
such that the pair $(F, \tau)$ is a triangulated functor?
\end{que}

\mysubsection{Some Properties of Triangulated Categories} 
\label{subsec:props-tr-cats}

In this subsection we prove a few general results on triangulated categories. 
Recall that all triangulated categories, and all triangulated functors between 
them, are (implicitly) $\K$-linear.

\begin{prop} \label{prop:1280}
Let $\cat{K}$ be a triangulated category. If
$L \xar{\al} M \xar{\be} N \xar{\ga} \opn{T}(L)$
is a distinguished triangle
in $\cat{K}$, then $\be \circ \al = 0$. 
\end{prop}

\begin{proof}
By axioms (TR1) and (TR3) we have a commutative diagram 
\[ \UseTips  \xymatrix @C=8ex @R=6ex { 
L
\ar[r]^{\opn{id}_L}
\ar[d]_{\opn{id}_L}
&
L
\ar[r]^{}
\ar[d]_{\al}
& 
0
\ar[r]^{}
\ar[d]_{}
&
\opn{T}(L)
\ar[d]_{\opn{T}(\opn{id}_L)}
\\
L
\ar[r]^{\al}
&
M
\ar[r]^{\be}
& 
N
\ar[r]^{\ga}
&
\opn{T}(L) \ .
}  \]
We see that $\be \circ \al$ factors through $0$. 
\end{proof}

\begin{prop}  \label{prop:1283}
Let $\cat{K}$ be a triangulated category, and let 
\[ \UseTips  \xymatrix @C=8ex @R=6ex { 
L 
\ar[r]^{\al}
\ar[d]_{\phi}
&
M
\ar[r]^{\be}
\ar[d]_{\psi}
& 
N
\ar[r]^{\ga}
\ar[d]_{\chi}
&
\opn{T}(L)
\ar[d]_{\opn{T}(\phi)}
\\
L' 
\ar[r]^{\al'}
&
M'
\ar[r]^{\be'}
& 
N'
\ar[r]^{\ga'}
&
\opn{T}(L') \ .
}  \]
be a morphism of distinguished triangles. If $\phi$ and $\psi$ are isomorphisms,
then $\chi$ is also an isomorphism.
\end{prop}

\begin{proof}
Take an arbitrary $P \in \cat{K}$, and let 
$F := \opn{Hom}_{\cat{K}}(P, -)$. We get a 
commutative diagram 
\[ \UseTips  \xymatrix @C=7ex @R=6ex { 
F(L) 
\ar[r]^{F(\al)}
\ar[d]^{F(\phi)}
&
F(M)
\ar[r]^{F(\be)}
\ar[d]_{F(\psi)}
& 
F(N)
\ar[r]^{F(\ga)}
\ar[d]_{F(\chi)}
&
F(\opn{T}(L))
\ar[d]_{F(\opn{T}(\phi))}
\ar[r]^(0.45){F(\opn{T}(\al))}
&
F(\opn{T}(M)) 
\ar[d]_{F(\opn{T}(\psi))}
\\
F(L') 
\ar[r]^{F(\al')}
&
F(M')
\ar[r]^{F(\be')}
& 
F(N')
\ar[r]^{F(\ga')}
&
F(\opn{T}(L'))
\ar[r]^(0.47){F(\opn{T}(\al'))}
&
F(\opn{T}(M')) 
} \]
in $\cat{Mod} \K$. By Proposition \ref{prop:4} the rows in the diagram 
are exact sequences.
Since the other vertical arrows are isomorphisms, it follows that 
\[ F(\chi) : \opn{Hom}_{\cat{K}}(P, N) \to 
\opn{Hom}_{\cat{K}}(P, N') \]
is an isomorphism of abelian groups. By forgetting structure, we see that 
$F(\chi)$ is an isomorphism of sets. 

Consider the Yoneda functor
$\opn{Y}_{\cat{K}} : \cat{K} \to \cat{Fun}(\cat{K}^{\mrm{op}}, \cat{Set})$
from Subsection \ref{subsec:repfunc}. 
There is a morphism 
$ \opn{Y}_{\cat{K}}(\chi) : \opn{Y}_{\cat{K}}(N) \to  
\opn{Y}_{\cat{K}}(N')$
in 
$\cat{Fun}(\cat{K}^{\mrm{op}}, \cat{Set})$. 
For every object $P \in \cat{K}$, letting 
$F := \opn{Hom}_{\cat{K}}(P, -)$
as above, we have 
$\opn{Y}_{\cat{K}}(N)(P) = F(N)$ and
$\opn{Y}_{\cat{K}}(N')(P) = F(N')$.
The calculation above shows that 
\[ F(\chi) = \opn{Y}_{\cat{K}}(\chi)(P) : \opn{Y}_{\cat{K}}(N)(P) \to 
\opn{Y}_{\cat{K}}(N')(P) \]
is an isomorphism in $\cat{Set}$. 
Therefore $\opn{Y}_{\cat{K}}(\chi)$ is an isomorphism in 
$\cat{Fun}(\cat{K}^{\mrm{op}}, \cat{Set})$. 
According to Yoneda Lemma (Theorem \ref{thm:3025})
the morphism $\chi : N \to N'$ in $\cat{K}$ is an isomorphism. 
\end{proof}

\begin{exer} \label{exer:4690}
Try to prove the last proposition without the Yoneda Lemma. 
\end{exer}

Recall that if $\eta : F \to G$ is a morphism of functors
$\cat{K} \to \cat{L}$, then for an object $M \in \cat{K}$
we denote by $\eta_M : F(M) \to G(M)$ the corresponding morphism in 
$\cat{L}$. 

\begin{cor} \label{cor:3601}
Suppose 
$F, G : \cat{K} \to \cat{L}$ 
are triangulated functors between triangulated categories, and 
$\eta : F \to G$ is a morphism of triangulated functors. 
Let 
$L \xar{\al} M \xar{\be} N \xar{\ga} \opn{T}(L)$ 
be a distinguished triangle in $\cat{K}$. If $\eta_L$ and $\eta_M$
are isomorphisms, then $\eta_N$ is an isomorphism.
\end{cor}

\begin{proof}
After applying the functors $F$ and $G$, and the morphisms of functors 
$\tau_{F}$, $\tau_{G}$ and $\eta$, we obtain a 
commutative diagram 
\[  \UseTips  \xymatrix @C=8ex @R=6ex { 
F(L) 
\ar[r]^{F(\al)}
\ar[d]_{\eta_L}
&
F(M)
\ar[r]^{F(\be)}
\ar[d]_{\eta_M}
& 
F(N)
\ar[rr]^(0.43){\tau_{F, L} \circ \, F(\ga)}
\ar[d]_{\eta_N}
&
&
\opn{T}_{\cat{L}}(F(L))
\ar[d]_{\opn{T}_{\cat{L}}(\eta_L)}
\\
G(L) 
\ar[r]^{G(\al)}
&
G(M)
\ar[r]^{G(\be)}
& 
G(N)
\ar[rr]^(0.43){\tau_{G, L} \circ \, G(\ga)}
&
&
\opn{T}_{\cat{L}}(G(L)) 
} \]
in $\cat{L}$ in which the rows are distinguished triangles. 
According to Proposition \ref{prop:1283} the morphism $\eta_N$ is an 
isomorphism.
\end{proof}

\begin{cor} \label{cor:3605}   
Let $\cat{K}$ be a triangulated category, and let
$L \xar{\al} M \xar{} N \xar{\, \triangle\, }$
be a distinguished triangle in it.
The two conditions below are equivalent.  
\begin{enumerate}
\rmitem{i} $\al : L \to M$ is an isomorphism.
\rmitem{ii} $N \cong 0$. 
\end{enumerate}
\end{cor}

\begin{exer} \label{exer:3600}
Prove Corollary \ref{cor:3605}. (Hint: Use Proposition \ref{prop:1283} and 
axiom (TR1)(c).)
\end{exer}

\begin{dfn} \label{dfn:3600}
Let $\cat{K}$ be a triangulated category, and let 
$\al : L \to M$ be a morphism in $\cat{K}$. By axiom (TR1)(b) there exists a 
distinguished triangle 
\[ \tag{$\dag$} L \xar{\al} M \xar{\be} N \xar{\ga} \opn{T}(L) \]
in $\cat{K}$. The object $N$ is called a 
{\em cone of $\al$}%
\index{Cone! of morphism in triangulated category},
and the distinguished triangle ($\dag$) is called a 
{\em distinguished triangle built on $\al$}%
\index{Triangle! distinguished}. 
\end{dfn}

\begin{cor} \label{cor:3600}
In the situation of Definition \tup{\ref{dfn:3600}}, the object $N$ and the 
distinguished triangle \tup{($\dag$)} are unique up to isomorphism.
\end{cor}

\begin{proof}
Suppose 
\[ \tag{$\dag'$} L \xar{\al} M \xar{\be'} N' \xar{\ga'} \opn{T}(L) \]
is another distinguished triangle built on $\al$. By axiom (TR3) there is a 
commutative diagram 
\[ \UseTips  \xymatrix @C=8ex @R=6ex { 
L 
\ar[r]^{\al}
\ar[d]_{\opn{id}_L}
&
M
\ar[r]^{\be}
\ar[d]_{\opn{id}_M}
& 
N
\ar[r]^{\ga}
\ar[d]_{\chi}
&
\opn{T}(L)
\ar[d]_{\opn{T}(\opn{id}_L)}
\\
L 
\ar[r]^{\al}
&
M
\ar[r]^{\be'}
& 
N'
\ar[r]^{\ga'}
&
\opn{T}(L) \ .
}  \]
in $\cat{K}$. Proposition \ref{prop:1283} says that $\chi$ is an isomorphism. 
\end{proof}

\begin{rem} \label{rem:4160}
In general the isomorphism $\chi$ in the corollary above is not unique, and 
thus the cone of $\al$ is not functorial in the morphism $\al$. However, in 
some special cases $\chi$ is unique -- see Remark \ref{rem:3600} below. 

Note also that the standard cones and triangles in the DG category 
$\dcat{C}(A, \cat{M})$ are functorial, by Proposition \ref{prop:1162}.
This implies that in the triangulated categories 
$\dcat{K}(A, \cat{M})$ and $\dcat{D}(A, \cat{M})$
we can sometimes arrange to have functorial cones.
\end{rem}

\begin{prop} \label{prop:3605}
Let $\cat{K}$ be a triangulated category, and let
$L \xar{\al} M \xar{\be} N \xar{\ga} \opn{T}(L)$ 
be a distinguished triangle in it.
The two conditions below are equivalent.  
\begin{enumerate}
\rmitem{i} The morphism $\ga$ is zero.
\rmitem{ii} There exists a morphism 
$\tau : M \to L$ such that $\tau \circ \al = \opn{id}_L$. 
\end{enumerate}
\end{prop}

\begin{proof} \mbox{}

\smallskip \noindent
(i) $\Rightarrow$ (ii): By Propositions \ref{prop:4} the contravariant functor 
$\opn{Hom}_{\cat{K}}(-, L)$ is cohomological. Applying it to the given 
distinguished triangle, and using Proposition \ref{prop:1281}, we get  
an exact sequence of $\K$-modules 
\[ \opn{Hom}_{\cat{K}}(M, L) \xar{\opn{Hom}(\al, \opn{id}_L)}  
\opn{Hom}_{\cat{K}}(L, L) \xar{\opn{Hom}(\opn{T}^{-1}(\ga), \opn{id}_L)}
\opn{Hom}_{\cat{K}}(\opn{T}^{-1}(N), L) . \]
Since the homomorphism 
$\opn{Hom}_{\cat{K}}(\opn{T}^{-1}(\ga), \opn{id}_L)$ is zero, 
there exists some morphism $\tau \in \opn{Hom}_{\cat{K}}(M, L)$
such that 
$\opn{id}_L = \opn{Hom}_{\cat{K}}(\al, \opn{id}_L)(\tau) = \tau \circ \al$.

\medskip \noindent 
(ii) $\Rightarrow$ (i):
Let us examine the commutative diagram 
\[ \UseTips  \xymatrix @C=14ex @R=6ex { 
\opn{Hom}_{\cat{K}}(N, N)
\ar[d]_(0.45){\opn{Hom}(\opn{id}, \ga)}
\\
\opn{Hom}_{\cat{K}}(N, \opn{T}(L))
\ar[d]_{\opn{Hom}(\opn{id}, \opn{T}(\al))}
\ar@(r,ul)[dr]^{ \ \opn{Hom}(\opn{id}, \opn{id})}
\\
\opn{Hom}_{\cat{K}}(N, \opn{T}(M))
\ar[r]^{\opn{Hom}(\opn{id}, \opn{T}(\tau))}
&
\opn{Hom}_{\cat{K}}(N, \opn{T}(L))
}  \]
in $\cat{Mod} \K$. Because the homomorphism 
$\opn{Hom}_{\cat{K}}(\opn{id}, \opn{id})$
is bijective, it follows that 
$\opn{Hom}_{\cat{K}}(\opn{id}, \opn{T}(\al))$
is injective. But the column is an exact sequence (by Propositions \ref{prop:4} 
and \ref{prop:1281}), and therefore 
$\opn{Hom}_{\cat{K}}(\opn{id}_{N}, \ga) = 0$. We conclude that
$\ga = \opn{Hom}_{\cat{K}}(\opn{id}_{N}, \ga)(\opn{id}_{N}) = 0$.
\end{proof}

\begin{lem} \label{lem:3600}
Let $M, M'$ be objects in a triangulated category $\cat{K}$.
Consider the canonical morphisms 
\[ M  \xar{e} M \oplus M' \xar{p} M \quad \tup{and} \quad
M'  \xar{e'} M \oplus M' \xar{p'} M' . \]
Then 
\[ M \xar{e} M \oplus M' \xar{p'} M' \xar{0} \opn{T}(M) \]
is a distinguished triangle in $\cat{K}$.
\end{lem}

\begin{proof}
By axiom (TR1)(c) and axiom (TR2) there is a distinguished triangle 
$0 \xar{0} M' \xar{\opn{id}_{M'}} M' \xar{0} 0$.
By axiom (TR1)(b) there is a distinguished triangle
\begin{equation} \label{eqn:4160}
 M \xar{e} M \oplus M' \xar{\be} N \xar{\ga} \opn{T}(M) 
\end{equation}
in $\cat{K}$, for some object $N$. 
Because $p \circ e = \opn{id}_M$, Proposition \ref{prop:3605} 
says that $\ga = 0$. 

Next, since $p' \circ e = 0$, axiom (TR3) produces a morphism of triangles 
\begin{equation} \label{eqn:3600}
\UseTips  \xymatrix @C=8ex @R=6ex { 
M 
\ar[r]^(0.4){e}
\ar[d]_{0}
&
M \oplus M'
\ar[r]^(0.6){\be}
\ar[d]_{p'}
& 
N 
\ar[r]^(0.47){\ga = 0}
\ar@{-->}[d]_{\psi}
&
\opn{T}(M)
\ar[d]_{0}
\\
0
\ar[r]^{0}
&
M'
\ar[r]^{\opn{id}_{M'}}
&
M'
\ar[r]^{0}
& 
0
} 
\end{equation}
We claim that $\psi$ is an isomorphism. This is proved 
indirectly. For every object $L \in \cat{K}$ there is a commutative diagram 
\begin{equation} \label{eqn:3602}
\UseTips  \xymatrix @C=4ex @R=6ex { 
\opn{Hom}_{\cat{K}}(L, M)
\ar[r]
&
\opn{Hom}_{\cat{K}}(L, M \oplus M')
\ar[r]
\ar[d]_{\opn{Hom}(\opn{id}_L, p')}
& 
\opn{Hom}_{\cat{K}}(L, N)
\ar[r]
\ar[d]_{\opn{Hom}(\opn{id}_L, \psi)}
&
0
\\
0
\ar[r]
&
\opn{Hom}_{\cat{K}}(L, M')
\ar[r]^{\opn{id}}
& 
\opn{Hom}_{\cat{K}}(L, M')
\ar[r]
&
0
}
\end{equation}
in $\cat{Mod} \K$, that's gotten by applying the functor 
$\opn{Hom}_{\cat{K}}(L,- )$ to the left part of diagram (\ref{eqn:3600}).
The rows in (\ref{eqn:3602}) are exact sequences. 
Because 
\[ \opn{Hom}_{\cat{K}}(L, M \oplus M') \cong 
\opn{Hom}_{\cat{K}}(L, M) \oplus \opn{Hom}_{\cat{K}}(L, M') \]
we can replace (\ref{eqn:3602}) with the next commutative diagram
\begin{equation} \label{eqn:4161}
\UseTips  \xymatrix @C=4ex @R=6ex { 
0
\ar[r]
&
\opn{Hom}_{\cat{K}}(L, M')
\ar[r]
\ar[d]_{\opn{Hom}(\opn{id}_L, \opn{id}_{M'})}
& 
\opn{Hom}_{\cat{K}}(L, N)
\ar[r]
\ar[d]_{\opn{Hom}(\opn{id}_L, \psi)}
&
0
\\
0
\ar[r]
&
\opn{Hom}_{\cat{K}}(L, M')
\ar[r]^{\opn{id}}
& 
\opn{Hom}_{\cat{K}}(L, M')
\ar[r]
&
0
}
\end{equation}
with exact rows. We see that 
$\opn{Hom}(\opn{id}_L, \psi)$ is an isomorphism of $\K$-modules. Hence it is a  
bijection of sets. Using the Yoneda Lemma, like in the proof of Proposition 
\ref{prop:1283}, we conclude that $\psi$ is an isomorphism in $\cat{K}$. 

Finally, from the commutative diagram (\ref{eqn:3600}) we know that 
$\psi \circ \be = p'$. So we have this isomorphism of triangles: 
\begin{equation} \label{eqn:3605}
\UseTips  \xymatrix @C=7ex @R=6ex { 
M 
\ar[r]^(0.4){e}
\ar[d]_{\opn{id}}
&
M \oplus M'
\ar[r]^(0.64){\be}
\ar[d]_{\opn{id}_{}}
& 
N 
\ar[r]^(0.45){\ga = 0}
\ar[d]_{\psi}^{\cong}
&
\opn{T}(M)
\ar[d]_{\opn{id}}
\\
M
\ar[r]^(0.4){e}
&
M \oplus M'
\ar[r]^(0.64){p'}
&
M'
\ar[r]^(0.45){0}
& 
\opn{T}(M)
} 
\end{equation}
The first triangle is distinguished. By axiom (TR1)(a) the second triangle is 
also distinguished. 
\end{proof}

\begin{dfn} \label{dfn:3235}
Let $\cat{K}$ be a linear category, and let $M, N \in \cat{K}$.
\begin{enumerate}
\item The object $M$ is called a {\em retract}
\index{Retract}
of $N$ if there are morphisms 
$M \xar{e} N \xar{p} M$ such that $p \circ e = \opn{id}_M$. 
The morphism $e : M \to N$ is called an embedding. 

\item The object $M$ is called a {\em direct summand}
\index{Direct summand}
of $N$ if there is an 
object $M' \in \cat{K}$ and an isomorphism $M \oplus M' \cong N$. 
The corresponding morphism $e : M \to N$ is called an embedding. 
\end{enumerate}
\end{dfn}

Note that in both items above the embedding $e : M \to N$ is a monomorphism 
in $\cat{K}$. 

\begin{thm} \label{thm:4160}  
Let $\cat{K}$ be a triangulated category, and let $e : M \to N$
be a morphism in $\cat{K}$. 
The following two conditions are equivalent\tup{:}
\begin{enumerate}
\rmitem{i} There is a morphism $e' : M' \to N$ in $\cat{K}$ such that 
\[ (e, e') : M \oplus M' \to N \]
is an isomorphism. In other words, $M$ is a direct summand of $N$,
\index{Direct summand}
with embedding $e : M \to N$. 

\rmitem{ii} There is a morphism $p : N \to M$  in $\cat{K}$  such that 
$p \circ e = \opn{id}_M$. In other words, $M$ is a retract of $N$, 
\index{Retract}
with embedding $e : M \to N$. 
\end{enumerate}
\end{thm}

\begin{proof} \mbox{}

\smallskip \noindent 
(i) $\Rightarrow$ (ii): This is trivial (and only requires $\cat{K}$ to be a 
linear category). 

\medskip \noindent 
(ii) $\Rightarrow$ (i): Let 
$M \xar{e} N \xar{\be} M' \xar{\ga} \opn{T}(M)$
be a distinguished triangle in $\cat{K}$ built on $e$. By Proposition 
\ref{prop:3605} we know that $\ga = 0$. According to Lemma \ref{lem:3600},
$M \xar{\ep} M \oplus M' \xar{\pi'} M' \xar{0} \opn{T}(M)$
is also the distinguished triangle in $\cat{K}$;
here $\ep$ is the embedding and $\pi'$ is the projection. 
Turning these two distinguished triangles and using axiom (TR3) we get a 
morphism $\th$ such that the diagram 
\[ \UseTips  \xymatrix @C=8ex @R=6ex {
\opn{T}^{-1}(M') 
\ar[r]^(0.6){0}
\ar[d]_{\opn{id}}
&
M 
\ar[r]^(0.4){\ep}
\ar[d]_{\opn{id}}
&
M \oplus M'
\ar[r]^(0.58){\pi'}
\ar@{-->}[d]_{\th}
& 
M'
\ar[d]_{\opn{id}}
\\
\opn{T}^{-1}(M') 
\ar[r]^(0.6){0}
&
M
\ar[r]^(0.5){e}
& 
N 
\ar[r]^(0.45){\be}
&
M'
} \]
is commutative. By Proposition \ref{prop:1283}, $\th$ is an isomorphism. 
Finally we define $e' := \th \circ \ep'$, where 
$\ep' : M' \to M \oplus M'$ is the embedding. 
\end{proof}

\begin{dfn} \label{dfn:4915}
Let $\cat{K}$ be a triangulated category. A {\em saturated full triangulated 
subcategory} of $\cat{K}$ is a full triangulated subcategory 
$\cat{K}' \sub \cat{K}$ (Definition \ref{dfn:1450}) which is closed in $\cat{K}$
under taking direct summands and under isomorphisms. 
\end{dfn}

\begin{dfn} \label{dfn:4916}
Let $\cat{K}$ be a triangulated category and let $Z \sub \opn{Ob}(\cat{K})$ be 
a set of objects. The {\em saturated full triangulated subcategory of $\cat{K}$
generated by $Z$} is the smallest saturated full triangulated subcategory 
$\cat{K}' \sub \cat{K}$ such that $Z \sub \opn{Ob}(\cat{K}')$.
\end{dfn}

\begin{prop} \label{prop:4470}
Let $\cat{K}$ be a triangulated category, let $Z \sub \opn{Ob}(\cat{K})$ 
be a set of objects, and let $\cat{K}' \sub \cat{K}$ be the saturated full 
triangulated subcategory of $\cat{K}$ generated by $Z$. An object 
$M \in \cat{K}$ belongs to $\cat{K}'$ if 
and only if there is a sequence $M_0, \ldots, M_r$ of objects of $\cat{K}$, 
such that $M_r = M$, and for every $i \leq r$ at least one of the following 
conditions holds\tup{:}
\begin{itemize}
\rmitem{a} $M_i \in Z$. 

\rmitem{b} There is an isomorphism $M_i \cong \opn{T}^p(M_j)$ in $\cat{K}$
for some $j < i$ and $p \in \Z$. 

\rmitem{c} There is a distinguished triangle 
$M_j \to M_k \to M_i \xar{\, \triangle\, }$
in $\cat{K}$ for some $j, k < i$.

\rmitem{d} $M_i$ is a direct summand of $M_j$ in $\cat{K}$ for some $j < i$.
\end{itemize}
\end{prop}

\begin{prop} \label{prop:4603}
Let $\cat{K}$ and $\cat{L}$ be triangulated categories, 
let $F, G : \cat{K} \to \cat{L}$ be triangulated functors, and let 
$\ze : F \to G$ be a morphism of triangulated functors. 
Denote by $\cat{K}'$ the full subcategory of $\cat{K}$ on the objects
$M$ such that $\ze_M : F(M) \to G(M)$ is an isomorphism in $\cat{L}$.
Then $\cat{K}'$ is a saturated full triangulated subcategory of $\cat{K}$. 
\end{prop}

\begin{cor} \label{cor:4920}
Let $\cat{K}$ and $\cat{L}$ be triangulated categories, 
let $F : \cat{K} \to \cat{L}$ be a triangulated functor, and let
$\cat{K}'$ be the kernel of $F$, i.e.\ $\cat{K}' \sub \cat{K}$ is the 
full subcategory on the objects $M \in \cat{K}$ such that $F(M) = 0$. 
Then $\cat{K}'$ is a saturated full triangulated subcategory of $\cat{K}$.
\end{cor}

\begin{exer} \label{exer:4603}
Prove Proposition \ref{prop:4470}, Proposition \ref{prop:4603}  and Corollary 
\lb \ref{cor:4920}. 
\end{exer}

We end this subsection with two remarks. 

\begin{rem} \label{rem:4777}
The story of the concept ``saturated full triangulated 
subcategory'' is convoluted. In the text \cite{Ver0} from 1977, J.-L. Verdier 
introduced {\em \'epaisse full triangulated subcategories}
(a rather messy definition). However, J. Rickard proved (see  
\cite[Proposition 1.3]{Ric0}, published in 1989) that \'epaisse is equivalent 
to saturated, as defined above. Rickard did not give a name for this new 
and useful property. 

Verdier died in 1989. A few years later, in 1996, G. Maltsiniotis published 
Verdier's thesis \cite{Ver}. In this thesis the word ``\'epaisse'' is absent, 
and instead we find the definition of saturated categories, identical to 
Definition \ref{dfn:4915}. Maltsiniotis notes that Verdier was aware (although 
this is only implicit in his writing) of the equivalence between \'epaisse and 
saturated; indeed, these are precisely the categories that occur as kernels of 
triangulated functors (see Corollary \ref{cor:4920} above, and Remark 
\ref{rem:4220} on Verdier quotients and saturated denominator 
sets of morphisms). 
\end{rem}

\begin{rem} \label{rem:3600}
As mentioned above, in Remark \ref{rem:4160}, the cone of a morphism 
$\al : L \to M$ in a triangulated category $\cat{D}$ is usually not functorial. 
However, if $\cat{D}$ is endowed with a {\em t-structure},
\index{t-structure}
then some cones can be made functorial. Here is a quick explanation. 

The concept of t-structure was introduced by A.A. Beilinson, J. Bernstein and 
P. Deligne in the book \cite{BBD}. A t-structure on a triangulated category 
$\cat{D}$ consists of a pair $(\cat{D}^{\leq 0}, \cat{D}^{\geq 0})$ of full 
subcategories of $\cat{D}$, that satisfy a few axioms (see e.g.\ 
\cite[Chapter X]{KaSc1} or \cite[Definition 4.1]{YeZh2}). The prototypical 
example is the {\em standard t-structure} on $\cat{D} := \dcat{D}(\cat{M})$, 
the derived category of an abelian category $\cat{M}$
(see Section \ref{sec:der-cat} of the book). Here $\cat{D}^{\leq 0}$ 
is the full subcategory of $\dcat{D}(\cat{M})$ on the complexes $M$ with 
nonpositive cohomology, and 
$\cat{D}^{\geq 0}$ is the full subcategory of $\dcat{D}(\cat{M})$ 
on the complexes $M$ with nonnegative cohomology. Other t-structures on 
$\cat{D} = \dcat{D}(\cat{M})$ are often called {\em perverse t-structures}. 

The {\em heart} of the t-structure $(\cat{D}^{\leq 0}, \cat{D}^{\geq 0})$ is 
the category 
$\cat{D}^{0} := \cat{D}^{\leq 0} \cap \cat{D}^{\geq 0}$.
In the prototypical example above, the heart is the category of complexes $M$ 
with cohomology concentrated in degree $0$, so it is equivalent to $\cat{M}$. 
Something similar happens in general: the heart $\cat{D}^{0}$ is always an 
abelian category, and there is a cohomological functor 
$\mrm{H}^0 : \cat{D} \to \cat{D}^{0}$.
The short exact sequences 
$0 \to L \xar{\al} M \xar{\be} N \to 0$
in $\cat{D}^{0}$ are the distinguished triangles 
$L \xar{\al} M \xar{\be} N \xar{\, \triangle\, }$
in $\cat{D}$ such that $L, M, N \in \cat{D}^{0}$.

Here is how a t-structure $(\cat{D}^{\leq 0}, \cat{D}^{\geq 0})$ affects 
functoriality of cones. 
Axiom (TR3) guarantees that for a pair of distinguished triangles and 
morphism $\phi, \psi$ as in the solid diagram below, there is {\em at least 
one} morphism $\chi$ (the dashed arrow) making the whole diagram commutative.  
\[ \UseTips  \xymatrix @C=8ex @R=6ex { 
L 
\ar[r]^{\al}
\ar[d]_{\phi}
&
M
\ar[r]^{\be}
\ar[d]_{\psi}
& 
N
\ar[r]^{\ga}
\ar@{-->}[d]_{\chi}
&
\opn{T}(L)
\ar[d]_{\opn{T}(\phi)}
\\
L' 
\ar[r]^{\al'}
&
M'
\ar[r]^{\be'}
& 
N'
\ar[r]^{\ga'}
&
\opn{T}(L') 
} \]
According to \cite[Proposition 1.1.9]{BBD} (see also \cite[Lemma 4.10]{YeZh2}),
if $L \in \cat{D}^{\leq 0}$ and $\opn{T}(N') \in \cat{D}^{\geq 0}$, then this 
morphism $\chi$ is {\em unique}. 

The idea of a t-structure emerged in the study of algebraic topology,  
microlocal analysis and representations of Lie algebras. In that context the 
triangulated category $\cat{D}$ was 
the full subcategory of $\dcat{D}(\K_X)$ on complexes with suitable 
finiteness and constructibility conditions. The base ring $\K$ was either a 
field or the $l$-adic integers $\what{\Z}_l$. The space (or site) $X$ was 
either a topological space with a classical metric topology, or
a scheme with its \'etale topology. 
The objects of the heart where called 
{\em perverse sheaves}. This theory is explained in the original book 
\cite{BBD}, and also in \cite{Bor} and in \cite[Chapter X]{KaSc1}. 
In algebraic geometry one can also consider {\em perverse coherent sheaves}; 
see \cite{YeZh2} and \cite{AriBez}. 
\end{rem}

\mysubsection{The Homotopy Category is Triangulated} 
\label{subsec:K-is-triang}

In this subsection we consider a $\K$-linear abelian category $\cat{M}$ and a 
DG central $\K$-ring $A$, where $\K$ is the commutative base ring. These 
ingredients give rise to the $\K$-linear DG category $\dcat{C}(A, \cat{M})$ of 
DG $A$-modules in $\cat{M}$, as in Subsection \ref{subsec:DGModinM}. 

The strict category $\dcat{C}_{\mrm{str}}(A, \cat{M})$ and the 
homotopy category $\dcat{K}(A, \cat{M})$ were introduced in Definition
\ref{dfn:1132}. Recall that these $\K$-linear categories have the same objects 
as $\dcat{C}(A, \cat{M})$. The morphism $\K$-modules are 
\[ \opn{Hom}_{\dcat{C}_{\mrm{str}}(A, \cat{M})} (M_0, M_1) = 
\opn{Z}^0 \bigl( \opn{Hom}_{\dcat{C}(A, \cat{M})} (M_0, M_1) \bigr) \]
and 
\[ \opn{Hom}_{\dcat{K}(A, \cat{M})} (M_0, M_1) = 
\opn{H}^0 \bigl( \opn{Hom}_{\dcat{C}(A, \cat{M})} (M_0, M_1) \bigr) . \]
There is a full additive functor 
$\opn{P} : \dcat{C}_{\mrm{str}}(A, \cat{M}) \to \dcat{K}(A, \cat{M})$ 
which is the identity on objects, and on morphisms it sends a homomorphism to 
its homotopy class. 
Morphisms in $\dcat{K}(A, \cat{M})$ will often by decorated with a bar, like
$\bar{\phi}$. 

Consider the translation functor $\opn{T}$ from  Definition 
\ref{dfn:1275}. Since $\opn{T}$ is a DG functor from $\dcat{C}(A, \cat{M})$ to 
itself (see Corollary \ref{cor:1300}), it restricts to a linear functor
from $\dcat{C}_{\mrm{str}}(A, \cat{M})$ to itself, and it induces a linear 
functor $\bar{\opn{T}}$ from $\dcat{K}(A, \cat{M})$ to itself, such that 
$\opn{P} \circ \opn{T} = \bar{\opn{T}} \circ \opn{P}$.

\begin{prop} \label{prop:1416} 
\mbox{}
\begin{enumerate}
\item The category $\dcat{C}_{\mrm{str}}(A, \cat{M})$,
equipped with the translation functor $\opn{T}$, is a T-additive category. 

\item The category $\dcat{K}(A, \cat{M})$,
equipped with the translation functor $\bar{\opn{T}}$, is a T-additive 
category. 

\item Let $\tau : \opn{P} \circ \opn{T} \iso \bar{\opn{T}} \circ \opn{P}$
be the identity automorphism. Then the pair 
\[ (\opn{P}, \tau) : \dcat{C}_{\mrm{str}}(A, \cat{M}) \to \dcat{K}(A, \cat{M}) 
\]
is a T-additive functor.
\end{enumerate}
\end{prop}

\begin{proof}
(1) We need to prove that $\dcat{C}_{\mrm{str}}(A, \cat{M})$
is additive. Of course the zero complex is a zero object. 
Next we consider finite direct sums. 
Let $M_1, \ldots, M_r$ be a finite collection of objects in 
$\dcat{C}(A, \cat{M})$. Each $M_i$ is a DG $A$-module in $\cat{M}$, and we 
write it as $M_i = \{ M_i^j \}_{j \in \Z}$. 
In each degree $j$ the direct sum 
$M^j := \bigoplus_{i = 1}^r M^j_i$ 
exists in $\cat{M}$. Let $M := \{ M^j \}_{j \in \Z}$ be the resulting graded 
object in $\cat{M}$. The differential $\d_M : M^j \to M^{j + 1}$ exists 
by the universal property of direct sums; so we obtain a complex 
$M \in \dcat{C}(\cat{M})$. 
The DG $A$-module structure on $M$ is defined 
similarly: for $a \in A^k$, there is an induced degree $k$ morphism 
$f(a) : M \to M$ in $\dcat{C}(\cat{M})$. Thus $M$ becomes an object of
$\dcat{C}(A, \cat{M})$. But the embeddings $e_i : M_i \to M$ are strict 
morphisms, so $(M, \{ e_i \})$ is a coproduct of the collection $\{ M_i \}$ 
in $\dcat{C}_{\mrm{str}}(A, \cat{M})$.

\medskip \noindent 
(2) Now consider the category $\dcat{K}(A, \cat{M})$.
Because the functor 
$\opn{P} : \dcat{C}_{\mrm{str}}(A, \cat{M}) \to \dcat{K}(A, \cat{M})$
is additive, and is bijective on objects, part (1) above and Proposition 
\ref{prop:2} say that $\dcat{K}(A, \cat{M})$ is an additive category. 

\medskip \noindent 
(3) Clear. 
\end{proof}

From now on  we  denote by $\opn{T}$, instead of by $\bar{\opn{T}}$, the 
translation functor of $\dcat{K}(A, \cat{M})$.

\begin{dfn} \label{dfn:140}
A triangle 
\[ L \xar{\bar{\al}} M \xar{\bar{\be}} N \xar{\bar{\ga}} \opn{T}(L)  \]
in $\dcat{K}(A, \cat{M})$ is said to be a {\em distinguished triangle}
\index{Triangle! distinguished}
if there is a standard triangle 
\[ L' \xar{\al'} M' \xar{\be'} N' \xar{\ga'} \opn{T}(L')  \]
in $\dcat{C}_{\mrm{str}}(A, \cat{M})$, as in Definition \ref{dfn:1140},
and an isomorphism of triangles 
\[ \UseTips  \xymatrix @C=8ex @R=6ex { 
L' 
\ar[r]^{\opn{P}(\al')}
\ar[d]_{ \bar{\phi} } 
&
M'
\ar[r]^{\opn{P}(\be')}
\ar[d]_{ \bar{\psi} }
& 
N'
\ar[r]^{\opn{P}(\ga')}
\ar[d]_{ \bar{\chi} }
&
\opn{T}(L')
\ar[d]_{ \opn{T}(\bar{\phi}) }
\\
L 
\ar[r]^{\bar{\al}}
&
M
\ar[r]^{\bar{\be}}
& 
N
\ar[r]^{\bar{\ga}}
&
\opn{T}(L)
\ .
}  \]
in $\dcat{K}(A, \cat{M})$.
\end{dfn}

\begin{thm} \label{thm:1255}
The T-additive category $\dcat{K}(A, \cat{M})$, with the set of 
distinguished triangles defined above, is a triangulated category%
\index{Triangulated category! $\dcat{K}(A, \cat{M})$ is a}%
\index{1-K(A,M)@$\dcat{K}(A, \cat{M})$}.
\end{thm}

The proof is after three lemmas. 

\begin{lem}  \label{lem:2}
Let $M \in \dcat{C}(A, \cat{M})$, and consider the standard cone 
$N := \opn{Cone}(\opn{id}_M)$. Then the DG module $N$ is null-homotopic, i.e.\ 
$0 \to N$ is an isomorphism in $\dcat{K}(A, \cat{M})$.
\end{lem}

\begin{proof}
We shall exhibit a homotopy $\th$ from $0_N$ to $\opn{id}_N$. 
Recall from Subsection \ref{subsec:cone} that 
$N = \opn{Cone}(\opn{id}_M) = M \oplus \opn{T}(M) = 
\sbmat{ M \\[0.2em]\opn{T}(M) }$
as graded modules, with a differential whose matrix presentation is
$\d_{N} = \sbmat{\d_M & \opn{t}_M^{-1} \\[0.2em] 0 & \d_{\opn{T}(M)}}$.
And by the definition in Subsection \ref{subsec:translation} we have
$\d_{\opn{T}(M)} = - \opn{t}_M \circ \, \d_M \circ \opn{t}_M^{-1}$. 
Define 
$\th : N \to N$ to be the degree $-1$ morphism with  matrix presentation
$\th :=  \sbmat{ 0 & \, 0 \\[0.2em] \opn{t}_M & \, 0 }$.
Then, using the formulas above for $\d_N$ and $\d_{\opn{T}(M)}$, we get
\[ \d_N \circ \th + \th \circ \d_N =  
\bmat{ \opn{id}_M & 0 \\[0.1em] 0 & \opn{id}_{\opn{T}(M)} } = \opn{id}_N. 
\qedhere \]
\end{proof}

\begin{exer} \label{exer:1255}
Here is a generalization of Lemma \ref{lem:2}. Consider a morphism 
$\phi : M_0 \to M_1$ in $\dcat{C}_{\mrm{str}}(A, \cat{M})$. Show that the three 
conditions below are equivalent:
\begin{enumerate}
\rmitem{i} $\phi$ is a homotopy equivalence.
\rmitem{ii} $\bar{\phi}$ is an isomorphism in $\dcat{K}(A, \cat{M})$.
\rmitem{iii} The DG module $\opn{Cone}(\phi)$ is null-homotopic.
\end{enumerate}
Try to do this directly, not using Proposition \ref{prop:4}(2) and Theorem 
\ref{thm:1255}.
\end{exer}

The next lemma is based on \cite[Lemma 1.4.2]{KaSc1}.

\begin{lem} \label{lem:1}
Consider a morphism $\al : L \to M$ in 
$\dcat{C}_{\mrm{str}}(A, \cat{M})$, 
the standard  triangle  
$L \xar{\al} M \xar{\be} N \xar{\ga} \opn{T}(L)$
built on $\al$, and the standard  triangle
$M \xar{\be} N \xar{\phi} P \xar{\psi}  \opn{T}(M)$
built on $\be$, all in $\dcat{C}_{\mrm{str}}(A, \cat{M})$. 
So $N = \opn{Cone}(\al)$ and $P = \opn{Cone}(\be)$.
There is a morphism 
$\rho : \opn{T}(L) \to P$ in $\dcat{C}_{\mrm{str}}(A, \cat{M})$ s.t.\ 
$\bar{\rho}$ is an isomorphism in $\dcat{K}(A, \cat{M})$, and the diagram 
\[ \UseTips \xymatrix @C=8ex @R=6ex { 
M
\ar[r]^{\bar{\be}}
\ar[d]_{\opn{id}_M}
&
N
\ar[r]^{\bar{\ga}}
\ar[d]_{\opn{id}_N}
& 
\opn{T}(L)
\ar[r]^{- \opn{T}(\bar{\al})}
\ar[d]_{\bar{\rho}}
&
\opn{T}(M)
\ar[d]_{\opn{id}_{\opn{T}(M)}}
\\
M
\ar[r]^{\bar{\be}}
&
N
\ar[r]^{\bar{\phi}}
& 
P
\ar[r]^{\bar{\psi}}
&
\opn{T}(M) 
} \]
commutes in $\dcat{K}(A, \cat{M})$.
\end{lem}

\begin{proof}
Note that 
$N = M \oplus \opn{T}(L)$ and
$P = N \oplus \opn{T}(M) = M \oplus \opn{T}(L) \oplus \opn{T}(M)$
as graded modules. Thus $P$ and $\d_P$ have the following matrix presentations:
\[ P = \bmat{ M \\[0.2em] \opn{T}(L) \\[0.2em] \opn{T}(M) } \ , \quad
\d_P =  \bmat{ \d_M & \al \circ \opn{t}_L^{-1} & \opn{t}_M^{-1} 
\\[0.3em] 0 & \d_{\opn{T}(L)} & 0 
\\[0.3em] 0 &  0 & \d_{\opn{T}(M)} } \ . \]
Define morphisms $\rho: \opn{T}(L) \to P$ and $\chi : P \to \opn{T}(L)$ in 
$\dcat{C}_{\mrm{str}}(A, \cat{M})$ by the matrix presentations
\[ \rho := \bmat{ 0 \\[0.1em] \opn{id}_{\opn{T}(L)}  
\\[0.2em] - \opn{T}(\al) } \  , \quad
\chi := \bmat{ 0 & \opn{id}_{\opn{T}(L)} & 0 } \ .  \]
Direct calculations show that
$\chi \circ \rho = \opn{id}_{\opn{T}(L)}$, \,
$\rho \circ \ga = \rho \circ \chi \circ \phi$ and 
$\psi \circ \rho = \lb - \opn{T}(\al)$.

It remains to prove that $\rho \circ \chi$ is homotopic to $\opn{id}_P$. 
Define a degree $-1$ morphism $\th : P \to P$ by the matrix
\[ \th := \bmat{ 0 & 0 & 0 \\[0em] 0 & 0 & 0 
\\[0em] \opn{t}_{M} & 0 & 0} . \]
Then a direct calculation, using the equalities 
$\opn{t}_{M} \circ \, \d_M + \d_{\opn{T}(M)} \circ \opn{t}_{M} = 0$
and 
$\opn{T}(\al) =  \opn{t}_{M} \circ \, \al \circ \opn{t}_{L}^{-1}$
gives 
$\th \circ \d_P + \d_P \circ \th = \opn{id}_P - \rho  \circ \chi$.
\end{proof}

\begin{lem} \label{lem:1265} 
Consider a  standard triangle 
$L \xar{\al} M \xar{\be} N \xar{\ga} \opn{T}(L)$
in $\dcat{C}_{\mrm{str}}(A, \cat{M})$.
For every integer $k$, the triangle 
\[ \opn{T}^k(L) \xar{ \, \opn{T}^k(\al) \, } \opn{T}^k(M) 
\xar{ \, \opn{T}^k(\be) \, } \opn{T}^k(N) 
\xar{ \, (-1)^k \cd \opn{T}^k(\ga) \, } 
\opn{T}^{k + 1}(L)  \]
is isomorphic, in $\dcat{C}_{\mrm{str}}(A, \cat{M})$, to a standard triangle.
\end{lem}

\begin{proof}
Combine Corollary \ref{cor:1300}, Corollary \ref{cor:1915} with $F = \opn{T}$,  
and Proposition \ref{prop:1300}.
\end{proof}

\begin{proof}[Proof of Theorem \tup{\ref{thm:1255}}] 
We essentially follow the proof of \cite[Proposition 1.4.4]{KaSc1}, adding some
details.

\medskip \noindent
(TR1): By definition the set of distinguished triangles in 
$\dcat{K}(A, \cat{M})$ is closed under isomorphisms. This establishes item (a). 

As for item (b): consider any morphism $\bar{\al} : L \to M$ in
$\dcat{K}(A, \cat{M})$. It is represented by a morphism $\al : L \to M$ in 
$\dcat{C}_{\mrm{str}}(A, \cat{M})$. Take the standard triangle on $\al$ 
in $\dcat{C}_{\mrm{str}}(A, \cat{M})$. Its image in $\dcat{K}(A, \cat{M})$
has the desired property. 

Finally, Lemma \ref{lem:2} shows that the triangle 
$M \xar{\opn{id}_M} M \to 0 \to \opn{T}(M)$ 
is isomorphic in $\dcat{K}(A, \cat{M})$ to the triangle
$M \xar{\opn{id}_M} M \xar{\bar{e}} \opn{Cone}(\opn{id}_M) \xar{\bar{p}} 
\opn{T}(M)$.
The latter is the image of a standard triangle, and so it is distinguished.

\medskip \noindent
(TR2): Consider the triangles
\begin{equation} \label{eqn:1255}
L \xar{\bar{\al}} M \xar{\bar{\be}} N \xar{\bar{\ga}} \opn{T}(L)
\end{equation}
and
\begin{equation} \label{eqn:1256}
M \xar{\bar{\be}} N \xar{\bar{\ga}} \opn{T}(L)\xar{- \opn{T}(\bar{\al})} 
\opn{T}(M) 
\end{equation}
in $\dcat{K}(A, \cat{M})$. If (\ref{eqn:1255}) is distinguished, then by 
Lemma \ref{lem:1} so is (\ref{eqn:1256}). 

Conversely, if (\ref{eqn:1256}) is distinguished, then by turning it $5$ times, 
and using the previous step (namely by Lemma \ref{lem:1}), we see that the 
triangle
\[ \opn{T}^2(L) \xar{\opn{T}^2(\bar{\al})} \opn{T}^2(M) 
\xar{\opn{T}^2(\bar{\be})} \opn{T}^2(N) \xar{\opn{T}^2(\bar{\ga})} 
\opn{T}^3(L) \]
is distinguished. 
According to Lemma \ref{lem:1265} (with $k = -2$), the triangle 
gotten by applying $\opn{T}^{-2}$ to this is distinguished.
But this is just the triangle (\ref{eqn:1255}).

\medskip \noindent
(TR3): 
Consider a commutative diagram in $\dcat{K}(A, \cat{M})$~:
\begin{equation} \label{eqn:1915}
\UseTips \xymatrix @C=8ex @R=6ex { 
\bar{L} 
\ar[r]^{\bar{\al}}
\ar[d]_{\bar{\phi}}
&
\bar{M}
\ar[r]^{\bar{\be}}
\ar[d]_{\bar{\psi}}
& 
\bar{N}
\ar[r]^{\bar{\ga}}
&
\opn{T}(\bar{L})
\\
\bar{L}' 
\ar[r]^{\bar{\al}'}
&
\bar{M}'
\ar[r]^{\bar{\be}'}
& 
\bar{N}'
\ar[r]^{\bar{\ga}'}
&
\opn{T}(\bar{L}') 
}
\end{equation}
where the horizontal triangles are distinguished. 
By definition the rows in (\ref{eqn:1915}) are isomorphic in 
$\dcat{K}(A, \cat{M})$
to the images under the functor $\opn{P}$ of 
standard triangles in $\dcat{C}_{\mrm{str}}(A, \cat{M})$. These are the rows in 
diagram (\ref{eqn:1916}) below. 
The vertical morphisms in  (\ref{eqn:1915}) are also induced from morphisms in 
$\dcat{C}_{\mrm{str}}(A, \cat{M})$, i.e.\ 
$\bar{\phi} = \opn{P}(\phi)$ and  $\bar{\psi} = \opn{P}(\psi)$. 
Thus (\ref{eqn:1915}) is isomorphic to the image under $\opn{P}$ of the 
following diagram:
\begin{equation} \label{eqn:1916}
\UseTips \xymatrix @C=8ex @R=6ex { 
L 
\ar[r]^{\al}
\ar[d]_{\phi}
&
M
\ar[r]^{\be}
\ar[d]_{\psi}
& 
N
\ar[r]^{\ga}
&
\opn{T}(L)
\\
L' 
\ar[r]^{\al'}
&
M'
\ar[r]^{\be'}
& 
N'
\ar[r]^{\ga'}
&
\opn{T}(L') 
} 
\end{equation}
Warning: the diagram (\ref{eqn:1916}) is only  commutative up to homotopy in 
 $\dcat{C}_{\mrm{str}}(A, \cat{M})$.

Since the rows in (\ref{eqn:1916}) are standard triangles
(see Definition \ref{dfn:1140}), the objects $N$ and 
$N'$ are cones: $N = \opn{Cone}(\al)$ and $N' = \opn{Cone}(\al')$.
The commutativity up to homotopy of this diagram means 
that there is a degree $-1$ morphism 
$\th : L \to M'$ in $\dcat{C}(A, \cat{M})$ such that
$\al' \circ \phi = \psi \circ \al + \d(\th)$.

Define the morphism 
$\chi : N = \sbmat{ M \\[0.2em] \opn{T}(L) } \to 
N' = \sbmat{ M' \\[0.2em] \opn{T}(L') }$
by the matrix presentation 
\[ \chi :=  
\bmat{ \psi & \th \circ \opn{t}_{L}^{-1} \\[0.3em] 0 & \opn{T}(\phi) } . \]
An easy calculation shows that $\chi$ is a morphism in 
$\dcat{C}_{\mrm{str}}(A, \cat{M})$, and that there are equalities
$\opn{T}(\phi) \circ \ga = \ga' \circ \chi$
and 
$\chi \circ \be = \be' \circ \psi$. 
Therefore, when we apply the functor $\opn{P}$, and conjugate by the original 
isomorphism between (\ref{eqn:1915}) and the image of (\ref{eqn:1916}), 
we obtain a commutative diagram 
\[ \UseTips \xymatrix @C=8ex @R=6ex { 
\bar{L} 
\ar[r]^{\bar{\al}}
\ar[d]_{\bar{\phi}}
&
\bar{M}
\ar[r]^{\bar{\be}}
\ar[d]_{\bar{\psi}}
& 
\bar{N}
\ar[r]^{\bar{\ga}}
\ar[d]_{\bar{\chi}}
&
\opn{T}(\bar{L})
\ar[d]_{\opn{T}(\bar{\phi})}
\\
\bar{L}' 
\ar[r]^{\bar{\al}'}
&
\bar{M}'
\ar[r]^{\bar{\be}'}
& 
\bar{N}'
\ar[r]^{\bar{\ga}'}
&
\opn{T}(\bar{L}') 
} \]
in $\dcat{K}(A, \cat{M})$, where $\bar{\chi}$ is conjugate to
$\opn{P}(\chi)$. 

\medskip \noindent
(TR4): We may assume that the three given distinguished triangles are 
standard triangles in $\dcat{C}_{\mrm{str}}(A, \cat{M})$. 
Namely, we can assume that 
$\al : L \to M$ and $\be : M \to N$ are morphisms in 
$\dcat{C}_{\mrm{str}}(A, \cat{M})$; the DG modules $P, Q, R$ are 
$P = \opn{Cone}(\al)$, 
$Q = \opn{Cone}(\be \circ \al)$ and
$R = \opn{Cone}(\be)$; and the morphisms $\ga, \de, \ep$ in 
$\dcat{C}_{\mrm{str}}(A, \cat{M})$ are
$\ga = e_{\al}$, $\de = e_{\be \circ \al}$ and 
$\ep = e_{\be}$. All this in the notation of Subsection \ref{subsec:cone}. 

In matrix notation we have
\[ P = \bmat{M \\[0.1em] \opn{T}(L)} , \quad 
Q = \bmat{N \\[0.1em] \opn{T}(L)} , \quad 
R = \bmat{N \\[0.1em] \opn{T}(M)} . \]
We define the morphisms 
$\phi : P \to Q$ and $\psi : Q \to R$ in $\dcat{C}_{\mrm{str}}(A, \cat{M})$ 
by the matrix presentations 
\[ \phi := \bmat{\be & 0 \\[0.1em] 0 & \opn{id}_{\opn{T}(L)}} , \quad 
\psi := \bmat{\opn{id}_{N} & 0 \\[0.1em] 0 & \opn{T}(\al)} . \]
(We leave to the reader to verify that $\phi$ and $\psi$ commute with the 
differentials $\d_P$, $\d_Q$ and $\d_R$; this is just linear algebra, using 
the matrix presentations of the differentials of the cones from Definition 
\ref{dfn:1172}.)
Define the morphism 
$\rho : R \to \opn{T}(Q)$ in $\dcat{C}_{\mrm{str}}(A, \cat{M})$
to be the composition of the morphisms 
$R \to \opn{T}(M) \xar{\opn{T}(\ga)} \opn{T}(Q)$.
Then the diagram (\ref{eqn:5022}) in $\dcat{C}_{\mrm{str}}(A, \cat{M})$ is 
commutative.

It remains to prove that the triangle 
\begin{equation} \label{eqn:2485}
 P \xar{\bar{\phi}} Q \xar{\bar{\psi}} R \xar{\bar{\rho}} \opn{T}(P) 
\end{equation}
in $\dcat{K}(A, \cat{M})$ is distinguished. 
Let $C := \opn{Cone}(\phi)$; so we have a standard triangle
\begin{equation} \label{eqn:2486}
 P \xar{\phi} Q \xar{e_{\phi}} C \xar{p_{\phi}} \opn{T}(P)
\end{equation}
in $\dcat{C}_{\mrm{str}}(A, \cat{M})$. We are going to prove that the triangles 
(\ref{eqn:2485}) and (\ref{eqn:2486}) are isomorphic in 
$\dcat{K}(A, \cat{M})$, by producing an isomorphism 
$\bar{\chi} : C \iso R$ in $\dcat{K}(A, \cat{M})$ that makes the diagram 
\[ \UseTips  \xymatrix @C=8ex @R=6ex { 
P
\ar[r]^{\bar{\phi}}
\ar[d]^{\opn{id}}
&
Q
\ar[r]^{\bar{e}_{\phi}}
\ar[d]^{\opn{id}}
&
C
\ar[r]^{\bar{p}_{\phi}}
\ar[d]^{\bar{\chi}}
&
\opn{T}(P)
\ar[d]^{\opn{id}}
\\
P
\ar[r]^{\bar{\phi}}
&
Q
\ar[r]^{\bar{\psi}}
&
R
\ar[r]^{\bar{\rho}}
&
\opn{T}(P)
} \]
commutative. 

Here are the matrices for the object $C$, the morphism 
$\chi : C \to R$, and another morphism 
$\om: R \to C$, both in $\dcat{C}_{\mrm{str}}(A, \cat{M})$. 
\[ C = \bmat{N \\[0.1em] \opn{T}(L) \\[0.3em] \opn{T}(M) 
\\[0.1em] \opn{T}^2(L)} , \quad 
\chi := \bmat{\opn{id}_N & 0 & 0 & 0 
\\[0.1em] 0 & \opn{T}(\al) & \opn{id}_{\opn{T}(M)} & 0 } , \quad 
\om := \bmat{\opn{id}_N & 0 \\[0.1em] 0 & 0 
\\[0.1em] 0 & \opn{id}_{\opn{T}(M)} \\[0.1em] 0 & 0} . \]
Again, we leave it to the reader to check that $\chi$ and $\om$ commute with 
the differentials. It is easy to see that 
$\om \circ \psi = e_{\phi}$, 
$\rho \circ \chi = p_{\phi}$ and 
$\chi \circ \om = \opn{id}_R$. 

Finally we must find a homotopy between $\om \circ \chi$ and $\opn{id}_C$.
Consider the degree $-1$ endomorphism $\th$ of $C$~:
\[ \th := 
\bmat{0 & 0 & 0 & 0 
\\[0em] 0 & 0 & 0 & 0 
\\[0em] 0 & 0 & 0 & 0 
\\[0em] 0 & \opn{t}_{\opn{T}(L)} & 0 & 0} . \]
Then 
$\d_C \circ \th + \th \circ \d_C = \opn{id}_C - \, \om \circ \chi$.
\end{proof}

A full subcategory of a DG category is of course a DG category itself. 
Full additive (resp.\ triangulated) subcategories of additive
(resp.\ triangulated) categories were defined in 
Definition \ref{dfn:1025} (resp.\ \ref{dfn:1450}). 

\begin{cor} \label{cor:3015}
Let $\cat{C} \sub \dcat{C}(A, \cat{M})$ be a full subcategory satisfying these 
three conditions\tup{:}
\begin{itemize}
\rmitem{a} $\cat{C}_{\mrm{str}}$ is a full additive subcategory of 
$\dcat{C}_{\mrm{str}}(A, \cat{M})$. 

\rmitem{b} $\cat{C}$ is translation invariant, i.e.\
$M \in \cat{C}$ if and only $\opn{T}(M) \in \cat{C}$.

\rmitem{c} $\cat{C}$ is closed under standard cones, i.e.\ for every morphism 
$\phi$ in $\cat{C}_{\mrm{str}}$ the object $\opn{Cone}(\phi)$ belongs to 
$\cat{C}$.
\end{itemize}
Then $\cat{K} := \opn{Ho}(\cat{C})$
is a full triangulated subcategory of 
$\dcat{K}(A, \cat{M})$. 
\end{cor}

\begin{proof}
Each condition here implies the same numbered condition in Definition 
\ref{dfn:1450}. 
\end{proof}

Here is a partial converse to the corollary. 

\begin{prop} \label{prop:3015}
Suppose $\cat{K} \sub \dcat{K}(A, \cat{M})$ is a full triangulated subcategory. 
Let $\cat{C} \sub \dcat{C}(A, \cat{M})$ be the full subcategory on the set of 
objects $\opn{Ob}(\cat{K})$. Then the DG category $\cat{C}$ satisfies 
conditions \tup{(a)-(c)} of Corollary \tup{\ref{cor:3015}}, and 
$\cat{K} = \opn{Ho}(\cat{C})$.
\end{prop}

\begin{exer} \label{exer:3015}
Prove Proposition \ref{prop:3015}. 
\end{exer}

Here are two remarks about other sorts of triangulated categories.

\begin{rem} \label{rem:4911}
Let $\cat{E}$ be an exact category in the sense of D. Quillen \cite{Qu1}. 
Assume that $\cat{E}$ is a {\em Frobenius category}, i.e.\ $\cat{E}$ has 
enough injectives and enough projectives, and an object is injective if and only 
if it is projective. The 
{\em stable category}%
\index{Stable category}
of $\cat{E}$ is the 
quotient category $\ul{\cat{E}} := \cat{E} / \cat{E}_0$, where 
$\cat{E}_0 \sub \cat{E}$ is the two-sided ideal of morphisms 
factoring through a projective object. The set of objects of $\ul{\cat{E}}$ is 
the same as that of $\cat{E}$. D. Happel \cite{Hap0} has shown that 
$\ul{\cat{E}}$ is 
a triangulated category. 

By definition, an {\em algebraic triangulated category} is a 
triangulated category $\cat{K}$ which is equivalent to the stable category 
$\ul{\cat{E}}$ of some Frobenius category $\cat{E}$.
It turns out that all the derived categories studied in this book
(see Section \ref{sec:der-cat}) are in fact 
algebraic triangulated categories; see the paper \cite{Kel0} of B. Keller. 

Here are a couple of other examples of algebraic triangulated categories. 
Let $\cat{M}$ be an additive category, and let 
$\cat{E} := \dcat{C}_{\mrm{str}}(\cat{M})$, the strict category of complexes in 
$\cat{M}$. Put on $\cat{E}$ the exact structure in which the short exact 
sequences are the degreewise split sequences. Then $\cat{E}$ becomes a 
Frobenius category; the projective-injective objects in it are the contractible 
complexes. In this case the ideal $\cat{E}_0$ coincides with the null-homotopic 
morphisms in $\cat{E} = \dcat{C}_{\mrm{str}}(\cat{M})$, and hence 
$\ul{\cat{E}} = \dcat{K}(\cat{M})$, the homotopy category of $\cat{M}$. 

For the second example, suppose that $\K$ is a field and $A$ a noetherian 
$\K$-ring which is Gorenstein, i.e.\ $A$ has finite injective dimension 
as a left and as a right module over itself. The category of finitely generated 
$A$-modules is $\dcat{M}_{\mrm{f}}(A) = \cat{Mod}_{\mrm{f}} A$.
Let $\cat{E} \sub \dcat{M}_{\mrm{f}}(A)$ 
be the full subcategory on the modules $M$ such that 
$\opn{Ext}^i_A(M, A) = 0$ for all $i > 0$. These modules are sometimes called 
{\em maximal Cohen-Macaulay modules}. Put on $\cat{E}$ the exact structure in 
which the exact sequences are the exact sequences in $\dcat{M}_{\mrm{f}}(A)$
with all terms belonging to $\cat{E}$. Then $\cat{E}$ is a Frobenius category. 
The projective-injective objects of $\cat{E}$ are the modules in it that are
projective $A$-modules. Moreover, according to R.-O. Buchweitz \cite{Buch},
the triangulated category $\ul{\cat{E}}$ in this case is equivalent to the {\em 
stable derived category} of $A$, which is the Verdier quotient 
$\dcat{D}^{\mrm{b}}_{\mrm{f}}(A) / \dcat{D}(A)_{\mrm{perf}}$.
Here $\dcat{D}^{\mrm{b}}_{\mrm{f}}(A)$ is the derived category of bounded 
complexes with finitely generated cohomologies (see subsection 
\ref{subsec:thick-subcats}), $\dcat{D}(A)_{\mrm{perf}}$ is the full 
subcategory on the perfect complexes (see Section \ref{sec:perf-tilt-NC}), and 
Verdier quotients are explained in Remark \ref{rem:4220}. 
Note that the triangulated category 
$\dcat{D}^{\mrm{b}}_{\mrm{f}}(A) / \dcat{D}(A)_{\mrm{perf}}$
is called the {\em triangulated category of singularities}  by
D.O. Orlov \cite{Orl}, and this was later shortened to the {\em singularity 
category} of $A$.
\end{rem}

\begin{rem} \label{rem:4920}
In 1961, A. Dold and D. Puppe \cite{DoPu} tried to formalize the properties 
of the {\em stable homotopy category of topological spaces}. What they got is a 
theory almost identical to the theory of triangulated categories studied here 
(Definition \ref{dfn:1158}), with the exception of the octahedral axiom (TR4). 

In more modern terms we can consider the notion of a 
{\em stable model category}%
\index{Stable category}, 
that we now define. Let $\cat{C}$ be a pointed Quillen model category, meaning 
that it has a zero object. One defines a {\em suspension functor} $\Sigma$ and a 
{\em loop functor} $\Om$ on $\cat{C}$ using pushout and pullback diagrams 
involving the zero object. Let $\cat{D}$ be the {\em homotopy category} of 
$\cat{C}$. (In the terminology of our book, $\cat{D}$ is actually the {\em 
derived category} of $\cat{C}$, because it is gotten by inverting the weak 
equivalences in $\cat{C}$; see Sections 
\ref{sec:loc-cats}-\ref{sec:der-funcs}.) If the functors 
$\Sigma$ and $\Om$ induce auto-equivalences on $\cat{D}$, quasi-inverse to each 
other, then $\cat{C}$ is called a stable model category. In this case it can be 
shown that $\cat{D}$ is a triangulated category, with translation functor 
$\opn{T} := \Sigma$.
A more general notion is that of a {\em stable infinity category}
(see \cite{Lur2} or \cite{Hi2}). 

A {\em topological triangulated category} is, by definition, a triangulated 
category $\cat{K}$ which is equivalent to a full triangulated subcategory 
of the homotopy category $\cat{D}$ of a stable model category $\cat{C}$, as 
above. For more information on these topics see \cite{Ho1}, \cite{Sch}, 
\cite{ShSc}, \cite{Lur2} and \cite{CaSt}. 
\end{rem}

\newpage

\mysubsection{From DG Functors to Triangulated Functors} 
\label{subsec:DG-to-tri-fun}

We now add a second DG ring $B$, and a second abelian category $\cat{N}$.
DG functors were introduced in Subsection \ref{subsec:DGFunc}.

Consider a DG functor 
$F  : \dcat{C}(A, \cat{M}) \to \dcat{C}(B, \cat{N})$.
From Theorem \ref{thm:1150} we know that the translation isomorphism $\tau_F$ 
is an isomorphism of DG functors
$\tau_F : F \circ \opn{T}_{A, \cat{M}} \iso 
\opn{T}_{B, \cat{N}} \circ \, F$.
Therefore, when we pass to the homotopy categories, 
and writing $\bar{F} := \opn{Ho}(F)$
and $\bar{\tau}_F := \opn{Ho}(\tau_F)$,
we get a T-additive functor 
$(\bar{F}, \bar{\tau}_F): \dcat{K}(A, \cat{M}) \to \dcat{K}(B, \cat{N})$.

If
$G : \dcat{C}(A, \cat{M}) \to \dcat{C}(B, \cat{N})$
is another DG functor, then we have another T-additive functor 
$(\bar{G}, \bar{\tau}_G) :  \dcat{K}(A, \cat{M}) \to \dcat{K}(B, \cat{N})$.
And if $\eta : F \to G$ is a strict morphism of DG functors, then there is a 
morphism of additive functors 
$\bar{\eta} := \opn{Ho}(\eta) : \bar{F} \to \bar{G}$.
This notation will be used in the next theorem. 

\begin{thm} \label{thm:1265}
Let $A$ and $B$ be DG rings, let $\cat{M}$ and $\cat{N}$ be abelian categories, 
and let 
$F : \dcat{C}(A, \cat{M}) \to \dcat{C}(B, \cat{N})$
be a DG functor. 
\index{Differential graded functor}
\begin{enumerate}
\item The T-additive functor  
\[ (\bar{F}, \bar{\tau}_F) :  \dcat{K}(A, \cat{M}) \to \dcat{K}(B, \cat{N}) \]
is a triangulated functor. 
\index{Triangulated functor}

\item Suppose 
$G : \dcat{C}(A, \cat{M}) \to \dcat{C}(B, \cat{N})$
is another DG functor, and $\eta : F \to G$ is a strict morphism of DG functors.
\index{Differential graded functor! strict morphism of {\indash}s} 
Then 
\[ \bar{\eta} : (\bar{F}, \bar{\tau}_F) \to (\bar{G}, \bar{\tau}_G) \]
is a morphism of triangulated functors. 
\index{Triangulated functor! morphism of {\indash}s}
\end{enumerate}
\end{thm}

\begin{proof} \mbox{}

\smallskip \noindent 
(1) Take a distinguished triangle 
$L \xar{\bar{\al}} M \xar{\bar{\be}} N \xar{\bar{\ga}} \opn{T}(L)$
in $\dcat{K}(A, \cat{M})$. Since we are only interested in triangles up to 
isomorphism, we can assume that this is the image under the functor $\opn{P}$ 
of a standard triangle  
$L \xar{\al} M \xar{\be} N \xar{\ga} \opn{T}(L)$
in $\dcat{C}_{\mrm{str}}(A, \cat{M})$. According to Theorem \ref{thm:1185}
and Corollary \ref{cor:1915}, 
there is a  standard triangle  
$L' \xar{\al'} M' \xar{\be'} N' \xar{\ga'} \opn{T}(L')$
in $\dcat{C}_{\mrm{str}}(B, \cat{N})$, and a commutative diagram 
\[ \UseTips  \xymatrix @C=6ex @R=6ex { 
F(L) 
\ar[r]^{F(\al)}
\ar[d]_{\phi}
&
F(M) 
\ar[r]^{F(\be)}
\ar[d]_{\psi}
& 
F(N)
\ar[rr]^(0.44){\tau_{F, L} \, \circ \, F(\ga)}
\ar[d]_{\chi}
&
&
\opn{T}( F(L))
\ar[d]_{\opn{T}(\phi)}
\\
L' 
\ar[r]^{\al'}
&
M'
\ar[r]^{\be'}
& 
N'
\ar[rr]^{\ga'}
&
&
\opn{T}(L') 
} \]
in $\dcat{C}_{\mrm{str}}(B, \cat{N})$,
in which the vertical arrows are isomorphisms. (Actually, we can take 
$L' = F(L)$, $\phi = \opn{id}_{F(L)}$, etc.) 
After applying the functor $\opn{P}$ to this diagram, we see that the condition 
in Definition \ref{dfn:1276}(1) is satisfied. 

\medskip \noindent
(2) By Definition \ref{dfn:1276}(2), all we need to prove that $\bar{\eta}$ is 
a morphism of T-add\-itive functors. Let's use these abbreviations:
$\cat{K} := \dcat{K}(A, \cat{M})$ and 
$\cat{L} := \dcat{K}(B, \cat{N})$. 
Definition \ref{dfn:1160} requires that for every object 
$M \in \cat{K}$ the diagram 
\begin{equation} \label{eqn:3215}
\UseTips  \xymatrix @C=8ex @R=6ex {
\bar{F}(\opn{T}_{\cat{K}}(M))
\ar[r]^{\bar{\tau}_{F, M}}
\ar[d]_{\bar{\eta}_{\opn{T}_{\cat{K}}(M)}}
&
\opn{T}_{\cat{L}}(\bar{F}(M))
\ar[d]^{\opn{T}_{\cat{L}}(\bar{\eta}_M)}
\\
\bar{G}(\opn{T}_{\cat{K}}(M))
\ar[r]^{\bar{\tau}_{G, M}}
&
\opn{T}_{\cat{L}}(\bar{G}(M)) \ . 
} 
\end{equation}
in $\cat{L}$ is commutative.
Going back to Definitions \ref{dfn:1171} and \ref{dfn:1150}  we see that 
$\opn{T}_{\cat{L}}(\eta_M) = \opn{t}_{G(M)} \circ \, \eta_M \circ 
\opn{t}_{F(M)}^{-1}$,
$\tau_{F, M} = \opn{t}_{F(M)} \circ \, F(\opn{t}_M^{-1})$
and 
$\tau_{G, M} = \opn{t}_{G(M)} \circ G(\opn{t}_M^{-1})$, 
as morphisms in  $\dcat{C}(B, \cat{N})$.
Thus the path in (\ref{eqn:3215}) that starts by going right is represented by 
the morphism 
\[ \begin{aligned}
&
\opn{T}_{\cat{L}}(\eta_M) \circ \tau_{F, M} =
\opn{t}_{G(M)} \circ \, \eta_M \circ \opn{t}_{F(M)}^{-1} \circ 
\opn{t}_{F(M)} \circ \, F(\opn{t}_M^{-1})
\\
& \quad 
= \opn{t}_{G(M)} \circ \, \eta_M \circ F(\opn{t}_M^{-1})  
\end{aligned} \]
in $\dcat{C}_{\mrm{str}}(B, \cat{N})$, and other path is represented by 
$\opn{t}_{G(M)} \circ \, G(\opn{t}_M^{-1}) \circ 
\eta_{\lms \opn{T}_{\cat{K}}(M)}$.
Because $\opn{t}_{G(M)}$ is an isomorphism (of degree $-1$) in 
$\dcat{C}(B, \cat{N})$, it suffices to prove that 
$\eta_M \circ F(\opn{t}_M^{-1}) =  G(\opn{t}_M^{-1}) \circ
\eta_{\opn{T}_{\cat{K}}(M)}$
in $\dcat{C}(B, \cat{N})$; namely that the diagram 
\[ \UseTips  \xymatrix @C=8ex @R=6ex {
F(\opn{T}_{\cat{K}}(M))
\ar[r]^(0.55){F(\opn{t}_M^{-1})}
\ar[d]_{\eta_{\lms \opn{T}_{\cat{K}}(M)}}
&
F(M)
\ar[d]^{\eta_M}
\\
G(\opn{T}_{\cat{K}}(M))
\ar[r]^(0.55){G(\opn{t}_M^{-1})}
&
G(M)
} \]
is commutative. This is true because $\eta :F \to G$ is a strict morphism of DG 
functors, and here it acts on the degree $1$ morphism 
$\opn{t}_M^{-1} : \opn{T}_{\cat{K}}(M) \to M$
in $\dcat{C}(B, \cat{N})$.
\end{proof}

\begin{cor} \label{cor:1265}
For every integer $k$, the pair
$\bigl( \opn{T}^k, (-1)^k \cd \opn{id}_{\opn{T}^{k + 1}} \bigr)$
is a triangulated functor from $\dcat{K}(A, \cat{M})$ to itself. 
\end{cor}

\begin{proof}
Combine Theorem \ref{thm:1265} and Proposition \ref{prop:1300}. 
\end{proof}

\begin{rem} \label{rem:1475}
In \cite{BoKa}, A.I. Bondal and M.M. Kapranov introduced the concept of 
{\em pretriangulated DG category}%
\index{Differential graded category! pretriangulated}.
This is a DG category $\cat{C}$ for which the 
homotopy category $\opn{Ho}(\cat{C})$ is canonically triangulated (the details 
of the definition are too complicated to mention here). Our DG categories 
$\dcat{C}(A, \cat{M})$ are pretriangulated in the sense of \cite{BoKa}; but 
they have a lot more structure (e.g.\ the objects have cohomologies too). 

Suppose $\cat{C}$ and $\cat{C}'$ are pretriangulated DG categories. 
In \cite{BoKa} there is a (rather complicated) definition of {\em pre-exact DG 
functor} $F : \cat{C} \to \cat{C}'$. It is stated there that if $F$ is a 
pre-exact DG functor, then 
$\opn{Ho}(F) : \opn{Ho}(\cat{C}) \to \opn{Ho}(\cat{C}')$
is a triangulated functor. This is analogous to our Theorem \ref{thm:1265}.
Presumably, Theorems \ref{thm:1150} and \ref{thm:1185} imply that every DG 
functor 
$F : \dcat{C}(A, \cat{M}) \to \dcat{C}(A', \cat{M}')$
is pre-exact in the sense of \cite{BoKa}; but we did not verify this.  
\end{rem}

\mysubsection{The Opposite Homotopy Category is Triangulated}
\label{subsec:opp-hom-triang}

Here we introduce a canonical triangulated structure on the opposite homotopy 
category $\dcat{K}(A, \cat{M})^{\mrm{op}}$. This gives us a way to talk about 
contravariant triangulated functors whose source is a full subcategory of 
$\dcat{K}(A, \cat{M})$. Our solution is precise, but it is not totally 
satisfactory -- see Remark \ref{rem:2306}. 

We already gave a thorough treatment of contravariant DG functors in 
Subsection \ref{subsec:contrvar-dg-func}. In the previous subsection we 
explained precisely how to pass from DG functors to triangulated functors. 
In this subsection we treat the contravariant case.

As before, $A$ is a central DG $\K$-ring, and $\cat{M}$ is a $\K$-linear 
abelian category. The DG category of DG $A$-modules in $\cat{M}$ is 
$\dcat{C}(A, \cat{M})$. 
In Subsection \ref{subsec:contrvar-dg-func} we introduced the flipped DG 
category 
$\dcat{C}(A^{\mrm{op}}, \cat{M}^{\mrm{op}})$. 
In Theorem \ref{thm:2495} we had a canonical isomorphism of DG categories
\begin{equation} \label{eqn:3000}
\opn{Flip} : 
\dcat{C}(A, \cat{M})^{\mrm{op}} \iso 
\dcat{C}(A^{\mrm{op}}, \cat{M}^{\mrm{op}}) . 
\end{equation}

\begin{dfn} \label{dfn:2375}
The {\em flipped category}
\index{Homotopy category! flipped}
\index{1-K(A,M)flip@$\dcat{K}(A, \cat{M})^{\mrm{flip}}$}
of $\dcat{K}(A, \cat{M})$ is the triangulated category 
\[ \dcat{K}(A, \cat{M})^{\mrm{flip}} :=
\dcat{K}(A^{\mrm{op}}, \cat{M}^{\mrm{op}}) = 
\opn{Ho} \bigl( \dcat{C}(A^{\mrm{op}}, \cat{M}^{\mrm{op}}) \bigr) . \]
See Theorem \ref{thm:1255}. Its translation functor is denoted by 
$\opn{T}^{\mrm{flip}}$. 
\end{dfn}

Since 
$\opn{Ho} \bigl( \dcat{C}(A, \cat{M})^{\mrm{op}} \bigr) = 
\dcat{K}(A, \cat{M})^{\mrm{op}}$,
the isomorphism (\ref{eqn:3000}) 
induces a $\K$-linear isomorphism of categories 
\begin{equation} \label{eqn:2995}
\ol{\opn{Flip}} := \opn{Ho}(\opn{Flip}) :  
\dcat{K}(A, \cat{M})^{\mrm{op}} \iso \dcat{K}(A, \cat{M})^{\mrm{flip}} .
\end{equation}

\begin{dfn} \label{dfn:2995}
The category $\dcat{K}(A, \cat{M})^{\mrm{op}}$
\index{Homotopy category! opposite}
\index{1-K(A,M)op@$\dcat{K}(A, \cat{M})^{\mrm{op}}$}
is given the triangulated 
category structure induced from the flipped category 
$ \dcat{K}(A, \cat{M})^{\mrm{flip}}$
under the isomorphism (\ref{eqn:2995}). 
The translation functor of $\dcat{K}(A, \cat{M})^{\mrm{op}}$ is denoted by 
$\opn{T}^{\mrm{op}}$.
\end{dfn}

Thus 
\begin{equation} \label{eqn:3001}
\opn{T}^{\mrm{op}} = \ol{\opn{Flip}}^{\, -1} \circ 
\opn{T}^{\mrm{flip}} \circ \, \ol{\opn{Flip}} .
\end{equation}
The distinguished triangles of $\dcat{K}(A, \cat{M})^{\mrm{op}}$
are the images under $\ol{\opn{Flip}}^{\, -1}$ of the distinguished triangles 
of $ \dcat{K}(A, \cat{M})^{\mrm{flip}}$.
And tautologically, 
\begin{equation} \label{eqn:3002}
(\ol{\opn{Flip}}, \opn{id}) : 
\dcat{K}(A, \cat{M})^{\mrm{op}} \to  \dcat{K}(A, \cat{M})^{\mrm{flip}}
\end{equation}
is an isomorphism of triangulated categories. 

\begin{dfn} \label{dfn:3000}
Let $\cat{K} \sub \dcat{K}(A, \cat{M})$ be a full additive subcategory, and 
assume that $\cat{K}^{\mrm{op}}$ is a triangulated subcategory of 
$\dcat{K}(A, \cat{M})^{\mrm{op}}$.
Then we give $\cat{K}^{\mrm{op}}$ the triangulated structure induced from 
$\dcat{K}(A, \cat{M})^{\mrm{op}}$.
\index{1-Kop@$\cat{K}^{\mrm{op}}$} 
\end{dfn}

To say this a bit differently, the condition on $\cat{K}$ in the definition is 
that 
$\cat{K}^{\mrm{flip}} := \ol{\opn{Flip}}(\cat{K}^{\mrm{op}})
\sub  \dcat{K}(A, \cat{M})^{\mrm{flip}}$ 
is a triangulated subcategory. The triangulated structure we put on 
$\cat{K}^{\mrm{op}}$ is such that 
$(\ol{\opn{Flip}}, \opn{id}) : (\cat{K}^{\mrm{op}}, \opn{T}^{\mrm{op}})
\to (\cat{K}^{\mrm{flip}}, \opn{T}^{\mrm{flip}})$
is an isomorphism of triangulated categories. 

\begin{rem} \label{rem:3170}
It could happen (though we have no example of it) that 
$\cat{K} \sub \dcat{K}(A, \cat{M})$ is a full triangulated subcategory, 
yet $\cat{K}^{\mrm{op}} \sub \dcat{K}(A, \cat{M})^{\mrm{op}}$ 
is not a triangulated subcategory; or vice versa. 
More on this issue in Remark \ref{rem:2306}. 
\end{rem}

The remark above notwithstanding, in many important cases, such as in 
Propositions \ref{prop:3050} and
\ref{prop:3051}, both $\cat{K} \sub \dcat{K}(A, \cat{M})$
and $\cat{K}^{\mrm{op}} \sub \dcat{K}(A, \cat{M})^{\mrm{op}}$
are triangulated. 

\begin{dfn} \label{dfn:2996}
Let $\cat{L}$ be a triangulated category, and let 
$\cat{K} \sub \dcat{K}(A, \cat{M})$ be a full subcategory 
such that $\cat{K}^{\mrm{op}} \sub \dcat{K}(A, \cat{M})^{\mrm{op}}$ is 
triangulated. A {\em contravariant triangulated functor}
\index{Triangulated functor! contravariant}
from $\cat{K}$ to $\cat{L}$ is, by definition, a triangulated functor 
$(F, \tau) : \cat{K}^{\mrm{op}} \to \cat{L}$
in the sense of Definition \ref{dfn:1276}, where $\cat{K}^{\mrm{op}}$ has the 
triangulated structure from Definition \ref{dfn:3000}.
\end{dfn}

\begin{thm} \label{thm:2995}
Let $A$ and $B$ be DG rings, and let $\cat{M}$ and $\cat{N}$ be abelian 
categories. Let 
$\cat{C} \sub \dcat{C}(A, \cat{M})$ be a full subcategory, and let 
$F : \cat{C}^{\mrm{op}} \to \dcat{C}(B, \cat{N})$
be a DG functor%
\index{Differential graded functor! contravariant}.
Consider the homotopy category 
$\cat{K} := \opn{Ho}(\cat{C}) \sub 
\opn{Ho} \bigl( \dcat{C}(A, \cat{M}) \bigr) = \dcat{K}(A, \cat{M})$
and the induced additive functor 
$\bar{F} := \opn{Ho}(F) : \cat{K}^{\mrm{op}} \to \dcat{K}(B, \cat{N})$.
Suppose that 
$\cat{K}^{\mrm{op}} \sub \dcat{K}(A, \cat{M})^{\mrm{op}}$ is triangulated.
Then there is a canonical translation isomorphism 
$\bar{\tau}$ such that 
\[ (\bar{F}, \bar{\tau}) : \cat{K}^{\mrm{op}} \to \dcat{K}(B, \cat{N}) \]
is a triangulated functor. 
\index{Triangulated functor! contravariant}
\end{thm}

\begin{proof}
Define 
\[ \cat{C}^{\mrm{flip}} := \opn{Flip}(\cat{C}^{\mrm{op}}) \sub
\dcat{C}(A^{\mrm{op}}, \cat{M}^{\mrm{op}}) =
\dcat{C}(A, \cat{M})^{\mrm{flip}} . \]
This is a DG category, and 
$\opn{Flip} : \cat{C}^{\mrm{op}} \to \cat{C}^{\mrm{flip}}$
is an isomorphism of DG categories. Also
$\opn{Ho} (\cat{C}^{\mrm{flip}}) = \cat{K}^{\mrm{flip}}$.
There is a DG functor
\[ F^{\mrm{flip}} := F \circ \opn{Flip}^{-1} : \cat{C}^{\mrm{flip}} 
\to \dcat{C}(B, \cat{N}) , \]
and it induces an additive functor 
\[ \bar{F}^{\mrm{flip}} := \opn{Ho}(F^{\mrm{flip}}) : 
\cat{K}^{\mrm{flip}} \to \dcat{K}(B, \cat{N}) . \]
By Theorem \ref{thm:1265} there is a translation isomorphism 
$\bar{\tau}^{\mrm{flip}}$ such that 
\[ (\bar{F}^{\mrm{flip}}, \bar{\tau}^{\mrm{flip}}) : 
\cat{K}^{\mrm{flip}} \to \dcat{K}(B, \cat{N})  \]
is a triangulated functor.
Using formula (\ref{eqn:3001}) we obtain these equalities and isomorphisms:
\[ \begin{aligned}
& \bar{F} \circ \opn{T}^{\mrm{op}} = 
\bar{F} \circ \ol{\opn{Flip}}^{\, -1} \circ \opn{T}^{\mrm{flip}} \circ \, 
\ol{\opn{Flip}} 
= \bar{F}^{\mrm{flip}} \circ \opn{T}^{\mrm{flip}} \circ \, 
\ol{\opn{Flip}} 
\\
& \quad 
\xar{\bar{\tau}^{\mrm{flip}}}
\opn{T}_{\dcat{K}(B, \cat{N})} \circ \, \bar{F}^{\mrm{flip}} \, \circ  \, 
\ol{\opn{Flip}} 
= \opn{T}_{\dcat{K}(B, \cat{N})} \circ \, \bar{F} . 
\end{aligned} \]
We define $\bar{\tau}$ to be the composed isomorphism. By construction 
the pair $(\bar{F}, \bar{\tau})$ is a triangulated functor. 
\end{proof}

\begin{exer} \label{exer:3220}
Let $\cat{C} \sub \dcat{C}(A, \cat{M})$ be a full subcategory s.t.\ 
$\cat{K} := \opn{Ho}(\cat{C})$ is a triangulated subcategory of 
$\dcat{K}(A, \cat{M})$, and let 
$G : \cat{C} \to \dcat{C}(B, \cat{N})^{\mrm{op}}$
be a DG functor. Write $\bar{G} := \opn{Ho}(G)$. 
Show that there is a translation isomorphism $\bar{\nu}$ s.t.\ 
$(\bar{G}, \bar{\nu}) : \cat{K} \to \dcat{K}(B, \cat{N})^{\mrm{op}}$
is a triangulated functor. 
(Hint: study the proof of the theorem above.) 
\end{exer}

\begin{rem} \label{rem:2306}
Let $\cat{K}$ be a triangulated category. There is a way to introduce a 
triangulated structure on the opposite category $\cat{K}^{\mrm{op}}$. 
\index{Triangulated category! opposite}
The translation functor is 
$\opn{T}^{\mrm{op}} := \opn{Op} \circ \opn{T}^{-1} \circ \opn{Op}^{-1}$.
The distinguished triangles are defined to be 
\[ N \xar{\opn{Op}(\be)} M \xar{\opn{Op}(\al)} L 
\xar{\opn{Op}(- \opn{T}^{-1}(\ga))}  \opn{T}^{\mrm{op}}(N) , \]
where 
$L \xar{\al} M \xar{\be} N \xar{\ga} \opn{T}(L)$
goes over all distinguished triangles in $\cat{K}$.
See \cite[Remark 10.1.10(ii)]{KaSc2}. 
This approach has an advantage: it allows to talk about contravariant 
triangulated functors 
\index{Triangulated functor! contravariant}
$\cat{K} \to \cat{L}$ without any restriction on the 
category $\cat{K}$. This is in contrast to our quite restricted Definition 
\ref{dfn:2996}. The disadvantage of the approach presented in this paragraph   
is that there is no easy way to tie it with the DG theory (as opposed to our 
Theorems \ref{thm:1265} and \ref{thm:2995}). 

It would have been very pleasing if in the case 
$\cat{K} = \dcat{K}(A, \cat{M})$, the triangulated structure on 
$\dcat{K}(A, \cat{M})^{\mrm{op}}$ presented in the paragraph above would 
coincide with the triangulated structure from Definition \ref{dfn:2995}.
However, our calculations seem to indicate otherwise. 

Since in this book we are only interested in triangulated functors that are 
of DG origin (covariant, as in Theorem \ref{thm:1265}, or contravariant, as in 
Theorem \ref{thm:2995}), we decided to adopt Definition \ref{dfn:2996}.
It is a reliable, yet somewhat awkward approach. For instance, it requires us 
to perform particular calculations for the composition of contravariant 
triangulated functors, as in Lemma \ref{lem:2156}. 
\end{rem}

\cleardoublepage
\mysection{Localization of Categories} \label{sec:loc-cats}

\AYcopyright

This section is devoted to the general theory of {\em Ore 
localization of categories}. We are given a category $\cat{A}$ and a 
multiplicatively closed set of morphisms $\cat{S} \sub \cat{A}$.
The localized category $\cat{A}_{\cat{S}}$, gotten by formally inverting the 
morphisms in $\cat{S}$, always exists. The goal is to have a presentation of the
morphisms of $\cat{A}_{\cat{S}}$ as left or right fractions. 

In Section \ref{sec:der-cat} we shall apply the results of this section to 
triangulated categories.

\mysubsection{The Formalism of Localization}
We will start with a category $\cat{A}$, without even assuming it is 
linear. Still we use the notation $\cat{A}$, because it will be suggestive to 
think about a linear category $\cat{A}$ with a single object, which is just a 
ring $A$. The reason is that our localization procedure is the same as that in 
noncommutative ring theory -- even when the category is not linear, and it has 
multiple objects. In Subsection \ref{subsec:loc-lin} we treat linear 
categories, and thus we recover the ring theoretic localization as a special 
case. 

The emphasis will be on morphisms rather than on objects. Thus it will be 
convenient to write 
$\cat{A}(M, N) := \opn{Hom}_{\cat{A}}(M, N)$
for $M, N \in \opn{Ob}(\cat{A})$. 
We sometimes use the notation $a \in \cat{A}$ for a morphism 
$a \in \cat{A}(M, N)$, leaving the objects implicit. When we write $b \circ a$ 
for $a, b \in \cat{A}$, we implicitly mean that these morphisms are 
composable. 

For heuristic purposes, we can think of $\cat{A}$ as a linear category (e.g.\ 
living inside some category of modules), with objects $M, N, \ldots$.
For any given object $M$, we then have a genuine ring 
$\cat{A}(M) := \cat{A}(M, M)$. 

\begin{dfn} \label{dfn:1405}
Let $\cat{A}$ be a category. A 
{\em multiplicatively closed set of morphisms}%
\index{Multiplicatively closed set of morphisms}
in $\cat{A}$ is a subcategory $\cat{S} \subseteq \cat{A}$ such that 
$\opn{Ob}(\cat{S}) = \opn{Ob}(\cat{A})$.
\end{dfn}

In other words, for any pair of objects $M, N \in \cat{A}$ 
there is a subset $\cat{S}(M, N) \subseteq \cat{A}(M, N)$, 
such that $\opn{id}_M \in \cat{S}(M, M)$, and such that for any
$s \in \cat{S}(L, M)$ and $t \in \cat{S}(M, N)$, the composition
$t \circ s \in \cat{S}(L, N)$. 

Using our shorthand, we can write the definition like this: 
$\opn{id}_M \in \cat{S}$,
and $s, t \in \cat{S}$ implies $t \circ s \in \cat{S}$. 

If $A = \cat{A}$ is a single object linear category, namely a ring, then 
$S = \cat{S}$ is a multiplicatively closed set in the sense of ring theory. 

There are various notions of localization in the literature.  We restrict 
attention to two of them. Here is the first: 

\begin{dfn} \label{dfn:1406}
Let $\cat{S}$ be a multiplicatively closed set of morphisms in a category 
$\cat{A}$. A {\em localization} 
\index{Localization! of a category}
of $\cat{A}$ with respect to $\cat{S}$  
is a pair $(\cat{A}_{\cat{S}}, \opn{Q})$, consisting of a category 
$\cat{A}_{\cat{S}}$ and a functor 
$\opn{Q} : \cat{A} \to \cat{A}_{\cat{S}}$, called the localization functor,
having the following properties:
\begin{itemize}[leftmargin=32pt] 
\item[(Loc1)] There is equality 
$\opn{Ob}(\cat{A}_{\cat{S}}) = \opn{Ob}(\cat{A})$, and $\opn{Q}$ is the 
identity 
on objects.

\rmitem{Loc2}
For every $s \in \cat{S}$, the morphism $\opn{Q}(s) \in \cat{A}_{\cat{S}}$
is invertible (i.e.\ it is an isomorphism). 

\rmitem{Loc3} Suppose $\cat{B}$ is a category, and 
$F : \cat{A} \to \cat{B}$
is a functor such that $F(s)$ is invertible for every 
$s \in \cat{S}$. 
Then there is a unique functor 
$F_{\cat{S}} : \cat{A}_{\cat{S}} \to \cat{B}$
such that 
$F_{\cat{S}} \circ \opn{Q} = F$ 
as functors $\cat{A} \to \cat{B}$.
\end{itemize}
\end{dfn}

In a commutative diagram:
\[ \UseTips  \xymatrix @C=10ex @R=6ex  {
\cat{S}
\ar[r]^{\mrm{inc}}
&
\cat{A}
\ar@{-->}[r]^{F}="q" 
\ar[d]_{\opn{Q}}
& 
\cat{B}
\\
&
\cat{A}_{\cat{S}}
\ar@{-->}[ur]_{F_{\cat{S}}} ^(0.5){}="f"
} \]

\begin{exa} \label{exa:2310}
Suppose $\cat{A}$ has a single object, say $x$. Here is what the definition 
above says in this case. We have a ring 
$A := \cat{A}(x, x)$ and a multiplicatively closed set
$S := \cat{S}(x, x) \sub A$. The localization of $A$ with respect to $S$ is a 
ring $A_S := \cat{A}_{\cat{S}}(x, x)$, with a ring homomorphism 
$\opn{Q} : A \to A_S$, such that $\opn{Q}(S) \sub (A_S)^{\times}$,
the group of invertible elements of $A_S$. 
Condition (Loc3) says that any ring homomorphism
$F : A \to B$ such that $F(S) \sub B^{\times}$ factors uniquely through 
$A_S$. 
\end{exa}

\begin{thm} \label{thm:2310}
Let $\cat{S}$ be a multiplicatively closed set of morphisms in a category 
$\cat{A}$. A localization 
\index{Localization! of a category}
$(\cat{A}_{\cat{S}}, \opn{Q})$ of 
$\cat{A}$ with respect to $\cat{S}$, in the sense of Definition 
\tup{\ref{dfn:1406}}, exists, and it is unique up to a unique isomorphism. 
\end{thm}

\begin{proof}
Uniqueness: Suppose $(\cat{A}'_{\cat{S}}, \opn{Q}')$ is another localization.
By condition (Loc3), for both localizations, there are unique functors 
$G : \cat{A}_{\cat{S}} \to \cat{A}'_{\cat{S}}$
and 
$G' : \cat{A}'_{\cat{S}} \to \cat{A}_{\cat{S}}$,
that satisfy 
$\opn{Q}' = G \circ \opn{Q}$ and 
$\opn{Q} = G' \circ \opn{Q}'$. These functors must be the identity on 
objects, and are inverses to each other. 

Existence: Here we encounter the set theoretic issue alluded to in Subsection 
\ref{subsec:set-theor}. The category $\cat{A}_{\cat{S}}$ that we will construct 
might have to belong to a bigger universe than that in which $\cat{A}$ lives; 
this is because of its large morphism sets. 
We proceed as follows. Let $\cat{B}$ be the free category on the set of objects 
$\opn{Ob}(\cat{A})$, with set of generating morphisms 
$\cat{A} \sqcup \cat{S}^{\mrm{op}}$.
Thus, the morphisms in $\cat{B}$ are the finite
sequences of composable morphisms in $\cat{A} \sqcup \cat{S}^{\mrm{op}}$, and 
composition is concatenation. Inside $\cat{B}$ we have the congruence 
(two-sided 
ideal) $\cat{N} \sub \cat{B}$ generated by the relations of $\cat{A}$, 
the relations of $\cat{S}^{\mrm{op}}$, and the 
relations $\opn{Op}(s) \circ s = \opn{id}$ and 
$s \circ \opn{Op}(s) = \opn{id}$ for all $s \in \cat{S}$. 
The quotient category 
$\cat{A}_{\cat{S}} := \cat{B} / \cat{N}$ has the desired properties. 
Cf.\ \cite{GaZi} for more details. 
\end{proof}

\begin{prop} \label{prop:3030}
Let $\cat{S}$ be a multiplicatively closed set in a category $\cat{A}$, 
let $\cat{B}$ be a category, 
and let $\opn{Q} : \cat{A} \to \cat{B}$ be a functor. 
Consider the multiplicatively closed set 
$\cat{S}^{\mrm{op}} \sub \cat{A}^{\mrm{op}}$ and the functor 
$\opn{Q}^{\mrm{op}} : \cat{A}^{\mrm{op}} \to \cat{B}^{\mrm{op}}$. 
The two conditions below are equivalent\tup{:}
\begin{itemize}
\rmitem{i} The pair $(\cat{B}, \opn{Q})$ is a localization of 
$\cat{A}$ with respect to $\cat{S}$. 

\rmitem{ii} The pair $(\cat{B}^{\mrm{op}}, \opn{Q}^{\mrm{op}})$ is a  
localization of $\cat{A}^{\mrm{op}}$ with respect to $\cat{S}^{\mrm{op}}$. 
\end{itemize}
\end{prop}

\begin{exer} \label{exer:3031}
Prove Proposition \ref{prop:3030}. 
\end{exer}

Often the localization $\cat{A}_{\cat{S}}$ is of little value, because there is 
no practical way to describe the morphisms in it. This issue will be addressed 
in the next subsection.

\mysubsection{Ore Localization}
There is a better notion of localization. 
The references here are \cite{RD}, \cite{GaZi}, \cite{We}, \cite{KaSc1}, 
\cite{Ste} and \cite{Row}. The first four references talk about localization of 
categories; and the last two talk about noncommutative rings. 
It seems that historically, this {\em noncommutative calculus of fractions} was 
discovered by Ore and Asano in ring theory, around 1930. There was progress in 
the categorical side, notably by Gabriel around 1960. 

In this subsection we mostly follow the treatment of \cite{Ste}; but we 
sometimes use diagrams instead of formulas.

\begin{dfn}  \label{dfn:1407}
Let $\cat{S}$ be a multiplicatively closed set of morphisms in a category 
$\cat{A}$. A {\em right Ore localization}
\index{Localization! Ore {\indash} of a category}
of $\cat{A}$ with respect to 
$\cat{S}$ is a pair $(\cat{A}_{\cat{S}}, \opn{Q})$, consisting of a category 
$\cat{A}_{\cat{S}}$ and a functor 
$\opn{Q} : \cat{A} \to \cat{A}_{\cat{S}}$, having the following 
properties:
\begin{itemize}[leftmargin=30pt] 
\item[(RO1)] There is equality 
$\opn{Ob}(\cat{A}_{\cat{S}}) = \opn{Ob}(\cat{A})$, and $\opn{Q}$ is the 
identity on
objects.

\item[(RO2)]  For every $s \in \cat{S}$, the morphism $\opn{Q}(s) \in
\cat{A}_{\cat{S}}$
is an isomorphism. 

\item[(RO3)]  Every morphism $q \in \cat{A}_{\cat{S}}$ can be written as
$q = \opn{Q}(a) \circ \opn{Q}(s)^{-1}$ for some 
$a \in \cat{A}$ and $s \in \cat{S}$.

\item[(RO4)]  Suppose $a, b \in \cat{A}$ satisfy 
$\opn{Q}(a) = \opn{Q}(b)$. Then $a \circ s = b \circ s$ for some 
$s \in \cat{S}$. 
\end{itemize}
\end{dfn}

The letters ``RO'' stand for ``right Ore''. 
We refer to the expression 
$q = \opn{Q}(a) \circ \opn{Q}(s)^{-1}$
as a {\em right fraction representation of $q$}. 
Here is the left sided version of this definition: 

\begin{dfn}  \label{dfn:2310}  
Let $\cat{S}$ be a multiplicatively closed set of morphisms in a category 
$\cat{A}$. A {\em left Ore localization}
\index{Localization! Ore {\indash} of a category}
of $\cat{A}$ with respect to 
$\cat{S}$ is a pair $(\cat{A}_{\cat{S}}, \opn{Q})$, consisting of a category 
$\cat{A}_{\cat{S}}$ and a functor 
$\opn{Q} : \cat{A} \to \cat{A}_{\cat{S}}$,
having the following properties:
\begin{itemize}[leftmargin=30pt] 
\item[(LO1)] There is equality 
$\opn{Ob}(\cat{A}_{\cat{S}}) = \opn{Ob}(\cat{A})$, and $\opn{Q}$ is the 
identity on
objects.

\item[(LO2)]  For every $s \in \cat{S}$, the morphism $\opn{Q}(s) \in
\cat{A}_{\cat{S}}$
is an isomorphism. 

\item[(LO3)]  Every morphism $q \in \cat{A}_{\cat{S}}$ can be written as
$q = \opn{Q}(s)^{-1} \circ \opn{Q}(a)$ for some 
$a \in \cat{A}$ and $s \in \cat{S}$.

\item[(LO4)]  Suppose $a, b \in \cat{A}$ satisfy 
$\opn{Q}(a) = \opn{Q}(b)$. Then $s \circ a = s \circ b$
for some $s \in \cat{S}$. 
\end{itemize}
\end{dfn}

As in the right case, the letters ``LO'' stand for ``left Ore''; and we refer 
to the expression 
$q = \opn{Q}(s)^{-1} \circ \opn{Q}(a)$
as a {\em left fraction representation of $q$}. 

\begin{rem} \label{rem:1845}
The results to follow in this subsection will be stated for right Ore 
localizations only. They all have ``left'' versions, with identical proofs 
(just a matter of reversing some arrows or compositions),  and therefore they 
will be omitted. 
\end{rem}

To reinforce the last remark, we give: 

\begin{prop} \label{prop:1845}
Let $\cat{S}$ be a multiplicatively closed set in a category $\cat{A}$, 
let $\cat{B}$ be a category, 
and let $\opn{Q} : \cat{A} \to \cat{B}$ be a functor. 
Consider the multiplicatively closed set 
$\cat{S}^{\mrm{op}} \sub \cat{A}^{\mrm{op}}$ and the functor 
$\opn{Q}^{\mrm{op}} : \cat{A}^{\mrm{op}} \to \cat{B}^{\mrm{op}}$. 
The two conditions below are equivalent\tup{:}
\begin{itemize}
\rmitem{i} The pair $(\cat{B}, \opn{Q})$ is a left Ore localization of 
$\cat{A}$ with respect to $\cat{S}$. 

\rmitem{ii} The pair $(\cat{B}^{\mrm{op}}, \opn{Q}^{\mrm{op}})$ is a right Ore 
localization of $\cat{A}^{\mrm{op}}$ with respect to $\cat{S}^{\mrm{op}}$. 
\end{itemize}
\end{prop}

\begin{exer} \label{exer:1845}
Prove Proposition \ref{prop:1845}. 
\end{exer}

\begin{lem} \label{lem:3}
Let $(\cat{A}_{\cat{S}}, \opn{Q})$ be a right Ore localization, 
let $a_1, a_2 \in \cat{A}$ and $s_1, s_2 \in \cat{S}$.
The following conditions are equivalent\tup{:}
\begin{enumerate}
\rmitem{i} 
$\opn{Q}(a_1) \circ \opn{Q}(s_1)^{-1} =  
\opn{Q}(a_2) \circ \opn{Q}(s_2)^{-1}$
in $\cat{A}_{\cat{S}}$. 

\rmitem{ii} There are $b_1, b_2 \in \cat{A}$ s.t.\ 
$a_1 \circ b_1 = a_2 \circ b_2$,
and $s_1 \circ b_1  = s_2 \circ b_2 \in \cat{S}$.
\end{enumerate}
\end{lem}

\begin{proof} \mbox{}

\smallskip \noindent
(ii) $\Rightarrow$ (i):
Since $\opn{Q}(s_i)$ and $\opn{Q}(s_i \circ b_i)$ are invertible, it 
follows that 
$\opn{Q}(b_i)$ are invertible. So 
\[ \begin{aligned}
& \opn{Q}(a_1) \circ \opn{Q}(s_1)^{-1} = \opn{Q}(a_1) \circ 
\opn{Q}(b_1) \circ \opn{Q}(b_1)^{-1} \circ \opn{Q}(s_1)^{-1} 
\\
& \quad = 
\opn{Q}(a_2) \circ \opn{Q}(b_2) \circ \opn{Q}(b_2)^{-1} \circ 
\opn{Q}(s_2)^{-1} =  \opn{Q}(a_2) \circ \opn{Q}(s_2)^{-1} .
\end{aligned} \]

\medskip \noindent
(i) $\Rightarrow$ (ii):
By property (RO3) there are $c \in \cat{A}$ and $u \in \cat{S}$ s.t.\
\begin{equation} \label{eqn:1425}
\opn{Q}(s_2)^{-1} \circ \opn{Q}(s_1) = \opn{Q}(c) \circ \opn{Q}(u)^{-1} . 
\end{equation}
Rewriting this equation we get 
\begin{equation} \label{eqn:1426}
\opn{Q}(s_1 \circ u) = \opn{Q}(s_2 \circ c) . 
\end{equation}
It is given that 
$\opn{Q}(a_1) =  \opn{Q}(a_2) \circ \opn{Q}(s_2)^{-1} \circ \opn{Q}(s_1)$.
Plugging (\ref{eqn:1425}) into it we obtain 
$\opn{Q}(a_1) =  \opn{Q}(a_2) \circ \opn{Q}(c) \circ \opn{Q}(u)^{-1}$.
Rearranging this equation we get 
$\opn{Q}(a_1 \circ u) = \opn{Q}(a_2 \circ c)$. 
By property (RO4) there is 
$v \in \cat{S}$ s.t.\ 
$a _1 \circ u \circ v = a_2 \circ c \circ v$.
Likewise, from equation (\ref{eqn:1426}) and property (RO4), there is 
$v' \in \cat{S}$ s.t.\ 
$s_1 \circ u \circ v' = s_2 \circ c \circ v'$. 

Again using property (RO3), there are $d \in \cat{A}$ and $w \in \cat{S}$ s.t.\
$\opn{Q}(v)^{-1} \circ \opn{Q}(v') = \opn{Q}(d) \circ \opn{Q}(w)^{-1}$.
Rearranging we get 
$\opn{Q}(v' \circ w) = \opn{Q}(v \circ d)$.
By property (RO4) there is $w' \in \cat{S}$ s.t.\ 
$v' \circ w \circ w' = v \circ d \circ w'$.
Define 
$b_1 := u \circ v \circ  d \circ w'$ 
and
$b_2 := c \circ v \circ d \circ w'$.
Then 
\[ \begin{aligned}
& s_1 \circ b_1 = s_1 \circ u \circ v \circ  d \circ w' = 
s_1 \circ u \circ v' \circ w \circ w' 
\\
& \quad = s_2 \circ c \circ v' \circ w \circ w' = s_2 \circ b_2 , 
\end{aligned} \]
and it is in $\cat{S}$. 
Also 
\[ a_1 \circ b_1 = a_1 \circ u \circ v \circ  d \circ w' =
a_2 \circ c \circ v  \circ  d \circ w' = a_2 \circ b_2 . \qedhere \]
\end{proof}

\begin{prop} \label{prop:103}
A right Ore localization $(\cat{A}_{\cat{S}}, \opn{Q})$ is a localization in 
the sense of Definition \tup{\ref{dfn:1406}}. 
\end{prop}

\begin{proof}
Say  $\cat{B}$ is a category, and 
$F : \cat{A} \to \cat{B}$
is a functor such that $F(s)$ is an isomorphism for every 
$s \in \cat{S}$. 

The uniqueness of a functor
$F_{\cat{S}} : \cat{A}_{\cat{S}} \to \cat{B}$
satisfying $F_{\cat{S}} \circ \opn{Q} = F$  is clear from property (RO3).
We have to prove existence.

Define $F_{\cat{S}}$ to be $F$ on objects.
For a morphism $q$ in $\cat{A}_{\cat{S}}$, property (RO3) says that
there is a right fraction presentation 
$q = \opn{Q}(a_1) \circ \opn{Q}(s_1)^{-1}$,
with $a_1 \in \cat{A}$ and $s_1 \in \cat{S}$.
Let $F_{\cat{S}}(q) := F(a_1) \circ F(s_1)^{-1}$ in $\cat{B}$.
We have to prove that this is well-defined, i.e.\ it does not depend on the 
presentation. So suppose that 
$q = \opn{Q}(a_2) \circ \opn{Q}(s_2)^{-1}$
is another right fraction presentation of $q$.
Let $b_1, b_2 \in \cat{A}$ be as in Lemma \ref{lem:3}. 
Since $F(s_i)$ and $F(s_i \circ b_i)$ are invertible, then so is $F(b_i)$. 
We get
$F(a_2) = F(a_1) \circ F(b_1) \circ F(b_2)^{-1}$
and 
$F(s_2) = F(s_1) \circ F(b_1) \circ F(b_2)^{-1}$.
Hence 
$F(a_2) \circ F(s_2)^{-1}  = F(a_1) \circ F(s_1)^{-1}$, 
as required. 

It remains to prove that $F_{\cat{S}}$ is a functor. Since the identity 
$\opn{id}_M$ of the object $M$ in $\cat{A}_{\cat{S}}$ can be presented 
as 
$\opn{id}_M = \opn{Q}(\opn{id}_M) \circ \opn{Q}(\opn{id}_M)^{-1}$, 
we see that $F_{\cat{S}}(\opn{id}_M) = \opn{id}_{F(M)}$. 

Next let $q_1$ and $q_2$ be morphisms in $\cat{A}_{\cat{S}}$, such that the 
composition $q_2 \circ q_1$ exists (i.e.\ the target of $q_1$ is the source of 
$q_2$). We have to show that 
$F_{\cat{S}}(q_2 \circ q_1)$ equals
$F_{\cat{S}}(q_2) \circ F_{\cat{S}}(q_1)$.
Choose presentations 
$q_i = \opn{Q}(a_i) \circ \opn{Q}(s_i)^{-1}$, 
so that 
\begin{equation} \label{eqn:1421}
F_{\cat{S}}(q_2) \circ F_{\cat{S}}(q_1) = 
F(a_2) \circ F(s_2)^{-1} \circ F(a_1) \circ F(s_1)^{-1}  . 
\end{equation}
By property (RO3) there is a right fraction presentation 
\begin{equation} \label{eqn:1428}
\opn{Q}(s_2)^{-1} \circ \opn{Q}(a_1) = \opn{Q}(b) \circ \opn{Q}(t)^{-1} 
\end{equation}
for some $b \in \cat{A}$ and $t \in \cat{S}$. 
Because
$\opn{Q}(a_1 \circ t) = \opn{Q}(s_2 \circ b)$,
by (RO4) there is some $r \in \cat{S}$ such that 
$a_1 \circ t \circ r = s_2 \circ b \circ r$.
Therefore 
$F(a_1 \circ t \circ r) = F(s_2 \circ b \circ r)$,
which implies, by canceling the invertible morphism $F(r)$ and rearranging, that
\begin{equation} \label{eqn:1420}
F(s_2)^{-1} \circ F(a_1) = F(b) \circ F(t)^{-1}
\end{equation}
in $\cat{B}$. 

Let us continue. Using equation (\ref{eqn:1428}) we have  
\[ \begin{aligned}
& q_2 \circ q_1 = \opn{Q}(a_2) \circ \opn{Q}(s_2)^{-1} \circ 
\opn{Q}(a_1) \circ \opn{Q}(s_1)^{-1} 
\\ 
& \qquad 
= \opn{Q}(a_2) \circ \opn{Q}(b) \circ \opn{Q}(t)^{-1} \circ \opn{Q}(s_1)^{-1} =
\opn{Q}(a_2 \circ b) \circ \opn{Q}(s_1 \circ t)^{-1} . 
\end{aligned} \]
Using this presentation of $q_2 \circ q_1$, and the equality (\ref{eqn:1420}), 
we obtain 
\[ \begin{aligned}
& F_{\cat{S}}(q_2 \circ q_1) = 
F(a_2 \circ b) \circ F(s_1 \circ t)^{-1} 
= F(a_2) \circ F(b) \circ F(t)^{-1} \circ F(s_1)^{-1}
\\
& \qquad 
= F(a_2) \circ F(s_2)^{-1} \circ F(a_1) \circ F(s_1)^{-1} . 
\end{aligned} \]
This is the same as (\ref{eqn:1421}). 
\end{proof}

\begin{cor} \label{cor:1405}
Let $\cat{S}$ be a multiplicatively closed set of morphisms in a category 
$\cat{A}$. Assume that these two assertions hold\tup{:}
\begin{itemize}
\item The pair $(\cat{B}, \opn{Q})$ is either a right or a left Ore 
localization of $\cat{A}$ with respect to $\cat{S}$.

\item The pair $(\cat{B}', \opn{Q}')$ is either a right or a left 
Ore localization of $\cat{A}$ with respect to $\cat{S}$.
\end{itemize}
Then there is a unique isomorphism of localizations 
$(\cat{B}, \opn{Q}) \cong (\cat{B}', \opn{Q}')$.
\end{cor}

\begin{proof}
By Proposition \ref{prop:103} (in its right or left versions, as the case may 
be), both $(\cat{B}, \opn{Q})$ and $(\cat{B}', \opn{Q}')$ 
are localizations in the sense of Definition \ref{dfn:1406}. Hence, by 
Theorem \ref{thm:2310}, there is a unique isomorphism 
$(\cat{B}, \opn{Q}) \cong (\cat{B}, \opn{Q}')$.
\end{proof}

\begin{dfn} \label{dfn:21}
Let $\cat{S}$ be multiplicatively closed set of morphisms in a  category 
$\cat{A}$. We say that $\cat{S}$ is a {\em right denominator set}
\index{Denominator set}
if it satisfies these two conditions:
\begin{enumerate}[leftmargin=30pt] 
\item[(RD1)] (Right Ore condition) Given $a \in \cat{A}$ and $s \in \cat{S}$,
there exist $b \in \cat{A}$ and $t \in \cat{S}$ such that
$a \circ t = s \circ b$.

\item[(RD2)] (Right cancellation condition) 
Given $a, b \in \cat{A}$ and  $s \in \cat{S}$ such that
$s \circ a = s \circ b$, there exists $t \in \cat{S}$ such that $a \circ t = b 
\circ t$.
\end{enumerate}
\end{dfn}

In the definition above the implicit assumption is that the sources and targets 
of the morphisms match; e.g., in (RD1) the morphisms $a$ and $s$ have the same 
target. Here are the diagrams illustrating the definition: 
\begin{equation} \label{eqn:2310}
\UseTips  \xymatrix @C=5ex @R=4ex  {
& 
K
\ar@{-->}[dl]_{t}
\ar@{-->}[dr]^{b}
\\
M 
\ar[dr]_{a}
& &
N
\ar[dl]^{s}
\\
&
L 
} 
\qquad 
\UseTips  \xymatrix @C=7ex @R=6ex {
K
\ar@{-->}[r]^{t}
&
M
\ar@(dl,ul)[d]_{a}
\ar@(dr,ur)[d]^{b}
\\
&
N
\ar[r]_{s}
&
L
} 
\end{equation}

Now the left version of Definition \ref{dfn:21}:

\begin{dfn} \label{dfn:2311}
Let $\cat{S}$ be multiplicatively closed set of morphisms in a  category 
$\cat{A}$. We say that $\cat{S}$ is a {\em left denominator set}
\index{Denominator set}
if it satisfies these two conditions:
\begin{enumerate}[leftmargin=30pt] 
\item[(LD1)] (Left Ore condition) Given $a \in \cat{A}$ and $s \in \cat{S}$,
there exist $b \in \cat{A}$ and $t \in \cat{S}$ such that
$t \circ a = b \circ s$.

\item[(LD2)] (Left cancellation condition) 
Given $a, b \in \cat{A}$ and  $s \in \cat{S}$ such that
$a \circ s = b \circ s$, there exists $t \in \cat{S}$ such that $t \circ a = t 
\circ b$.
\end{enumerate}
\end{dfn}

\begin{prop} \label{prop:2310}
Let $\cat{S}$ be multiplicatively closed set of morphisms in a  category 
$\cat{A}$. The two conditions below are equivalent\tup{:}
\begin{itemize}
\rmitem{i} $\cat{S}$ is a left denominator set in $\cat{A}$. 

\rmitem{ii} $\cat{S}^{\mrm{op}}$ is a right denominator set in 
$\cat{A}^{\mrm{op}}$. 
\end{itemize}
\end{prop}

\begin{exer} \label{exer:2312}
Prove Proposition \ref{prop:2310}. 
\end{exer}

Here is the main theorem regarding Ore localization. 

\begin{thm} \label{thm:103}
The following conditions are equivalent for a category $\cat{A}$ 
and a multiplicatively closed set of morphisms $\cat{S} \subseteq \cat{A}$. 
\begin{enumerate}
\rmitem{i} The right Ore localization
\index{Localization! Ore {\indash} of a category}
$(\cat{A}_{\cat{S}}, \opn{Q})$ exists.

\rmitem{ii} $\cat{S}$ is a right denominator set.
\index{Denominator set}
\end{enumerate}
\end{thm}

The proof of Theorem \ref{thm:103} is after some preparation. 
The hard part is proving that (ii) $\Rightarrow$ (i). 
The general idea is the same as in 
commutative localization: we consider the set of pairs of morphisms
$\cat{A} \times \cat{S}$, and define a relation $\sim$ on it, with the hope 
that this is an equivalence relation, and that the quotient set 
$\cat{A}_{\cat{S}}$ will be a category, and it will have the desired 
properties. 

Let's assume that $\cat{S}$ is a right denominator set. 
For any $M, N \in \opn{Ob}(A)$ consider the set 
\[ (\cat{A} \times \cat{S})(M, N) :=  
\coprod_{L \in \opn{Ob}(A)} 
\cat{A}(L, N) \times \cat{S}(L, M) . \]

\begin{rem} \label{rem:1405}
The set $(\cat{A} \times \cat{S})(M, N)$ could be big, namely not an element 
of the initial universe $\cat{U}$. This would require the introduction of a 
larger universe, say $\cat{V}$, in which $\cat{U}$ is an element. 
And the resulting category $\cat{A}_{\cat{S}}$ will be a 
$\cat{V}$-category. This is the same issue we had in the proof of Theorem 
\ref{thm:2310}. 

We will ignore this issue. Moreover, in many cases of interest (derived 
categories where there are DG enhancements, such as the K-injective 
enhancement), there will be an alternative presentation of $\cat{A}_{\cat{S}}$
as a $\cat{U}$-category. We will refer to this when we get to it. 
\end{rem}

An element 
$(a, s) \in (\cat{A} \times \cat{S})(M, N)$ can be pictured as a diagram
\begin{equation} \label{eqn:1408}
\UseTips  \xymatrix @C=5ex @R=4ex  {
& L 
\ar[dl]_{s}
\ar[dr]^{a}
\\
M 
& &
N
}
\end{equation}
in $\cat{A}$.
This diagram will eventually represent the right fraction 
$\opn{Q}(a) \circ \opn{Q}(s)^{-1} : M \to N$.

\begin{dfn}  \label{dfn:1408}
We define a relation $\sim$ on the set $\cat{A} \times \cat{S}$
like this:
$(a_1, s_1) \sim (a_2, s_2)$ 
if there exist $b_1, b_2 \in \cat{A}$ s.t.\ 
$a_1 \circ b_1 = a_2 \circ b_2$ 
and
$s_1 \circ b_1 = s_2 \circ b_2 \in \cat{S}$.
\end{dfn}

Note that the relation $\sim$ imposes condition (ii) of Lemma \ref{lem:3}.
Here it is in a commutative diagram, in which we have made the objects explicit:
\begin{equation} \label{eqn:1409}
\UseTips  \xymatrix @C=6ex @R=6ex  {
&
K
\ar@{-->}[dl]_{b_1}
\ar@{-->}[dr]^{b_2}
\\
L_1
\ar[d]_{s_1}
\ar[drr]_(0.75){a_1}
& &
L_2
\ar[dll]_(0.75){s_2}
\ar[d]^{a_2}
\\
M 
& &
N
} 
\end{equation}
The paths (i.e.\ morphisms) ending at $M$ are in $\cat{S}$. 

\begin{lem} \label{lem:4}
If $\cat{S}$ is a right denominator set, then the relation $\sim$ is an 
equivalence.
\end{lem}

\begin{proof}
Reflexivity: take $K := L$ and $b_i := \opn{id}_L : L \to L$. 
Symmetry is trivial. 

Now to prove transitivity. 
Suppose we are given 
$(a_1, s_1) \sim (a_2, s_2)$ and $(a_2, s_2) \sim (a_3, s_3)$.
So we have the solid part of the first diagram in (\ref{eqn:4925}),
and it is commutative. The morphisms ending at $M$ are in $\cat{S}$.
\begin{equation} \label{eqn:4925}
\UseTips  \xymatrix @C=6.5ex @R=6ex  {
H
\ar@{-->}[r]^{u}
&
I
\ar@{-->}[dl]_{t}
\ar@{-->}[dr]^{d}
\\
K
\ar[d]_{b_1}
\ar[dr]^{b_2}
& & 
J
\ar[dl]_{c_2}
\ar[d]^{c_3}
\\
L_1
\ar[d]_{s_1}
\ar[drr]_(0.65){a_1}
&
L_2
\ar[dl]_(0.6){s_2}
\ar[dr]^(0.6){a_2}
&
L_3
\ar[dll]^(0.65){s_3}
\ar[d]^{a_3}
\\
*+++{M} 
& &
*+++{N}
}
\qquad 
\UseTips  \xymatrix @C=5ex @R=6ex  {
H
\ar[dd]_{b_1 \circ t \circ u}
\ar[ddrr]^{c_3 \circ d \circ u}
\\
\rule{0pt}{0.5em}
\\
L_1
\ar[d]_{s_1}
\ar[drr]_(0.7){a_1}
&
&
L_3
\ar[dll]^(0.7){s_3}
\ar[d]^{a_3}
\\
*+{M}
& &
*+{N}
} 
\end{equation}
 
By condition (RD1) applied to the morphisms $K \to M \leftarrow J$ there are 
$t \in \cat{S}$ and $d \in \cat{A}$ s.t.\ 
$(s_3 \circ c_3) \circ d = (s_1 \circ b_1) \circ t$.
Comparing the morphisms $I \to M$ in this diagram, we see that  
\[ s_2 \circ (b_2 \circ t) = s_1 \circ b_1 \circ t = s_3 \circ c_3 \circ d =
s_2 \circ (c_2 \circ d) . \]
By (RD2) there is $u \in \cat{S}$ s.t.\ 
$(b_2 \circ t) \circ u = (c_2 \circ d) \circ u$.
So all morphisms $H \to M$ are equal and belong to $\cat{S}$, and  all 
morphisms $H \to N$ are equal.
Now delete the object $L_2$ and the arrows going through it. 
Then delete the objects $I, J, K$,  but keep the paths going through them. 
We get the second diagram in (\ref{eqn:4925}). It is commutative, and all
morphisms ending at $M$ are in $\cat{S}$. This is evidence for 
$(a_1, s_1) \sim (a_3, s_3)$.
\end{proof}

\begin{proof}[Proof of Theorem \tup{\ref{thm:103}}] \mbox{}

\smallskip \noindent
Step 1. In this step we prove (i) $\Rightarrow$ (ii). 
Take $a \in \cat{A}$ and $s \in \cat{S}$ with the same target. 
Consider $q := \opn{Q}(s)^{-1} \circ \opn{Q}(a)$. By (RO3) there are 
$b \in \cat{A}$ and $t \in \cat{S}$ s.t.\ 
$q = \opn{Q}(b) \circ \opn{Q}(t)^{-1}$. 
So
$\opn{Q}(s \circ b) = \opn{Q}(a \circ t)$.
By (RO4) there is $u \in \cat{S}$ s.t.\ 
$(s \circ b) \circ u = (a \circ t) \circ u$.
We read this as 
$s \circ (b \circ u) = a \circ (t \circ u)$,
and note that $t \circ u \in \cat{S}$. So (RD1) holds.

Next $a, b \in \cat{A}$ and $s \in \cat{S}$ s.t.\ 
$s \circ a = s \circ b$.  Then 
$\opn{Q}(s \circ a) = \opn{Q}(s \circ b)$. But $\opn{Q}(s)$ is invertible, so 
$\opn{Q}(a) = \opn{Q}(b)$. By (RO4) there is $t \in \cat{S}$ s.t.\ 
$a \circ t = b \circ t$. We have proved (RD2).

\medskip \noindent
Step 2. Now we assume that condition (ii) holds, and we define 
the morphism sets $\cat{A}_{\cat{S}}(M, N)$, composition between them, and the 
identity morphisms. 

For any $M, N \in \opn{Ob}(\cat{A})$ let
\[ \cat{A}_{\cat{S}}(M, N) :=
(\cat{A} \times \cat{S})(M, N) \, / \, \sim , \]
the quotient set modulo the relation $\sim$ from 
Definition \ref{dfn:1408}, which is an
equivalence relation by Lemma \ref{lem:4}.

We define composition like this. Given 
$q_1 \in \cat{A}_{\cat{S}}(M_0, M_1)$
and 
$q_2 \in \cat{A}_{\cat{S}}(M_1, M_2)$,
choose representatives 
$(a_i, s_i) \in (\cat{A} \times \cat{S})(M_{i-1}, M_i)$.
We use the notation 
$q_i = \ol{ (a_i, s_i) }$ to indicate this. 
By (RD1) there are $c \in \cat{A}$ and $u \in \cat{S}$ s.t.\ 
$s_2 \circ c = a_1 \circ u$. 
The composition 
$q_2 \circ q_1 \in \cat{A}_{\cat{S}}(M_0, M_2)$
is defined to be
\[ q_2 \circ q_1 := \ol{ (a_2 \circ c, s_1 \circ u) } 
\in (\cat{A} \times \cat{S})(M_{0}, M_2)) . \]
The idea behind the formula can be seen in diagram (\ref{eqn:4926}). 
\begin{equation} \label{eqn:4926}
\UseTips  \xymatrix @C=3.5ex @R=3.5ex  {
& & 
K
\ar[dl]_{u}
\ar[dr]^{c}
\\
&
L_1
\ar[dl]_{s_1}
\ar[dr]_(0.5){a_1}
& & 
L_2
\ar[dl]_(0.5){s_2}
\ar[dr]^{a_2}
\\
M_0 
& &
M_1
& &
M_2
}
\end{equation}

We have to verify that this definition is independent of the representatives.
So suppose we take other representatives 
$q_i = \ol{ (a'_i, s'_i) }$, and we choose morphisms $u', c'$ to
construct the composition. This is the solid part of the diagram
(\ref{eqn:4927}), and it is a commutative diagram. We must prove that 
\[ \ol{ (a_2 \circ c, s_1 \circ u) }  = \ol{ (a'_2 \circ c', s'_1 \circ u') } . 
\]
\begin{equation} \label{eqn:4927}
\UseTips  \xymatrix @C=8ex @R=6ex  {
J_1
\ar@{-->}[d]_{b_1}
\ar@{-->}[dr]^(0.7){b'_1}
&
K
\ar[dl]^(0.7){u}
\ar[drr]^(0.7){c}
& & 
K'
\ar[dll]_(0.7){u'}
\ar[dr]_(0.7){c'}
& 
J_2
\ar@{-->}[dl]_(0.7){b_2}
\ar@{-->}[d]^{b'_2}
\\
L_1
\ar[d]_{s_1}
\ar[drr]_(0.65){a_1}
&
L'_1
\ar[dl]^(0.65){s'_1}
\ar[dr]^{a'_1}
& & 
L_2
\ar[dl]_{s_2}
\ar[dr]_(0.65){a_2}
& 
L'_2
\ar[dll]^(0.65){s'_2}
\ar[d]^{a'_2}
\\
M_0 
& &
*++{M_1}
& &
M_2
} 
\end{equation}

There are morphisms $b_i, b'_i$ that are evidence for 
$(a_i, s_i) \sim (a'_i, s'_i)$.
They are depicted as the dashed arrows in (\ref{eqn:4927}). 
That whole diagram is commutative. 
The morphisms $J_1 \to M_0$, $K \to M_0$, $K' \to M_0$ and 
$J_2 \to M_1$ are all in $\cat{S}$.

Choose $v_1 \in \cat{S}$ and $d_1 \in \cat{A}$ s.t.\ the first diagram in 
(\ref{eqn:4928}) is commutative. This can be done by (RD1). 
\begin{equation} \label{eqn:4928}
\UseTips \xymatrix @C=3.5ex @R=3.5ex {
\\
&
I_1
\ar[dl]_{v_1}
\ar[dr]^{d_1}
\\
J_1
\ar[dr]_{b_1}
& &
K
\ar[dl]^{u}
\\
&
L_1
} 
\quad  
\UseTips  \xymatrix @C=3.5ex @R=3.5ex {
I_2
\ar@{-->}[dr]^{\til{v}}
\\
& 
\til{I}_2
\ar@{-->}[dl]_{\til{v}_2}
\ar@{-->}[dr]^{\til{d}_2}
\\
K'
\ar[dr]_{c'}
& &
J_2
\ar[dl]^{b'_2}
\\
&
L'_2 
\ar[dr]_{s'_2}
\\
&
&
M_1
} 
\quad 
\xymatrix @C=3.5ex @R=3.5ex {
\\
&
I_2
\ar[dl]_{v_2}
\ar[dr]^{d_2}
\\
K'
\ar[dr]_{c'}
& &
J_2
\ar[dl]^{b'_2}
\\
&
L'_2
} 
\end{equation}

Consider the solid part of the middle diagram in (\ref{eqn:4928}).
Since $J_2 \to M_1$ is in $\cat{S}$, by (RD1) there are 
$\til{v}_2 \in \cat{S}$ and $\til{d}_2 \in \cat{A}$ s.t.\ the two paths
$\til{I}_2 \to M_1$ are equal.
By (RD2) there is $\til{v} \in \cat{S}$ s.t.\ the two paths 
$I_2 \to L'_2$ are equal. 
We get the commutative diagram in the middle of (\ref{eqn:4928}).
Next, defining 
$d_2 := \til{d}_2 \circ \til{v}$ and 
$v_2 := \til{v}_2 \circ \til{v} \in \cat{S}$,
we obtain the third commutative diagram in (\ref{eqn:4928}).

We now embed the first and third diagrams from (\ref{eqn:4928}) into the 
diagram (\ref{eqn:4927}). This gives us the solid diagram in (\ref{eqn:4929}),
and it is commutative. 
The morphisms $I_1 \to M_0$ belong to $\cat{S}$.  

\begin{equation} \label{eqn:4929}
\UseTips  \xymatrix @C=8ex @R=6ex  {
H''
\ar@{-->}[r]^{w''}
& 
H'
\ar@{-->}[r]^{w'}
& 
H
\ar@{-->}[dl]_{w}
\ar@{-->}[dr]^{e}
\\
&
I_1
\ar@{->}[dl]_{v_1}
\ar@{->}[d]^{d_1}
& & 
I_2
\ar@{->}[d]_{v_2}
\ar@{->}[dr]^{d_2}
\\
J_1
\ar[d]_{b_1}
\ar[dr]^(0.7){b'_1}
&
K
\ar[dl]^(0.7){u}
\ar[drr]^(0.7){c}
& & 
K'
\ar[dll]_(0.7){u'}
\ar[dr]_(0.7){c'}
& 
J_2
\ar[dl]_(0.7){b_2}
\ar[d]^{b'_2}
\\
L_1
\ar[d]_{s_1}
\ar[drr]_(0.65){a_1}
&
L'_1
\ar[dl]^(0.65){s'_1}
\ar[dr]^{a'_1}
& & 
L_2
\ar[dl]_{s_2}
\ar[dr]_(0.65){a_2}
& 
L'_2
\ar[dll]^(0.65){s'_2}
\ar[d]^{a'_2}
\\
M_0 
& &
*++{M_1}
& &
M_2
}
\end{equation}

Choose $w \in \cat{S}$ and $e \in \cat{A}$, starting at an object $H$, to fill 
the diagram $I_1 \to M_0 \leftarrow I_2$, using (RD1). The path 
$H \to I_1 \to M_0$ is in $\cat{S}$, and all the paths $H \to M_0$ are equal.   
But we could have failure of commutativity in the paths $H \to L'_1$ and 
$H \to L_2$. 

The two paths  $H \to L'_1$ in (\ref{eqn:4929}) satisfy 
\[ s'_1 \circ (b'_1 \circ v_1 \circ w) = s'_1 \circ (u' \circ v_2 \circ e) . \]
Therefore, by (RD2), there is $w' \in \cat{S}$ s.t.\ 
\[  (b'_1 \circ v_1 \circ w) \circ w' = (u' \circ v_2 \circ e) \circ w' . \]
Next, the two paths  $H' \to L_2$ satisfy 
\[ s_2 \circ (c \circ d_1 \circ w \circ w') = s_2 \circ (b_2 \circ d_2 \circ e 
\circ w') ; \]
this is because we can take a detour through $L'_1$. 
Therefore, again by (RD2), there is $w'' \in \cat{S}$ s.t.\ 
\[ (c \circ d_1 \circ w \circ w') \circ w'' = (b_2 \circ d_2 \circ e \circ w') 
\circ w'' . \]
Now all paths $H'' \to M_2$ in (\ref{eqn:4929}) are equal. All paths 
$H'' \to M_0$ are equal and are in $\cat{S}$. 

Erase the objects $M_1, J_1, J_2$ and all arrows touching them
from (\ref{eqn:4929}). 
Then erase $H, H'$, but keep the paths through them. 
We obtain the commutative diagram (\ref{eqn:4930}).
\begin{equation} \label{eqn:4930}
\UseTips  \xymatrix @C=8.5ex @R=6ex  {
&
I_1
\ar[d]^{d_1}
& 
H''
\ar[l]_{w \circ w' \circ w''}
\ar[r]^{e \circ w' \circ w''}
& 
I_2
\ar[d]_{v_2}
\\
&
K
\ar[dl]^(0.7){u}
\ar[drr]^(0.7){c}
& & 
K'
\ar[dll]_(0.7){u'}
\ar[dr]_(0.7){c'}
& 
\\
L_1
\ar[d]_{s_1}
&
L'_1
\ar[dl]^(0.65){s'_1}
& & 
L_2
\ar[dr]_(0.65){a_2}
& 
L'_2
\ar[d]^{a'_2}
\\
M_0 
& &
& &
M_2
}
\end{equation}
This is evidence for 
\[  (a_2 \circ c, s_1 \circ u) \sim (a'_2 \circ c', s'_1 \circ u') . \]
The proof that composition is well-defined is done.

The identity morphism $\opn{id}_M$ in $\cat{A}_{\cat{S}}$
of an object $M$ is $\ol{(\opn{id}_M, \opn{id}_M)}$.  

\medskip \noindent
Step 3. 
We have to verify the associativity and the identity properties of 
composition in $\cat{A}_{\cat{S}}$. Namely that $\cat{A}_{\cat{S}}$ is a
category. This seems to be not too hard, given Step
2, and we leave it as an exercise.

\medskip \noindent
Step 4.  
The functor $\opn{Q} : \cat{A} \to \cat{A}_{\cat{S}}$ is defined to be 
$\opn{Q}(M) := M$ on objects, and $\opn{Q}(a) := \ol{(a, \opn{id}_M)}$ for 
$a : M \to N$ in $\cat{A}$. We have to verify this is a functor... Again, an 
exercise.

\medskip \noindent
Step 5.  
Finally we verify properties (RO1)-(RO4). (RO1) is clear. The inverse of
$\opn{Q}(s)$ is $\ol{(\opn{id}, s)}$, so (RO2) holds. 

It is not hard to see that 
$\ol{(a, s)} = \ol{(a, \opn{id})} \circ \ol{(\opn{id}, s)}$;
this is (RO3). 

If $\opn{Q}(a_1) = \opn{Q}(a_2)$, then 
$(a_1, \opn{id}_M) \sim (a_2, \opn{id}_M)$; so there
are $b_1, b_2 \in \cat{A}$  s.t.\ 
$a_1 \circ b_1 = a_2 \circ b_2$ and 
$\opn{id} \circ \, b_1 = \opn{id} \circ \, b_2 \in \cat{S}$. 
Writing $s := b_1 \in \cat{S}$, we get $a_1 \circ s = a_2 \circ s$. This proves 
(RO4). 
\end{proof}

\begin{exer} \label{exer:3030}
Finish the details in the proof of Theorem \tup{\ref{thm:103}}. 
\end{exer}

\begin{prop}[Common Denominator] \label{prop:107}
Let $\cat{A}$ be a category, let $\cat{S}$ be a right
denominator set in $\cat{A}$, and let $(\cat{A}_{\cat{S}}, \opn{Q})$ be the 
right Ore
localization.
For every two morphisms $q_1, q_2 : M \to N$ in $\cat{A}_{\cat{S}}$ there
is a common right denominator. Namely we can write 
$q_i = \opn{Q}(a_i) \circ \opn{Q}(s)^{-1}$
for suitable $a_i \in \cat{A}$ and $s \in \cat{S}$. 
\end{prop}

\begin{proof}
Choose representatives
$q_i = \opn{Q}(a'_i) \circ \opn{Q}(s'_i)^{-1}$. 
By (RD1) applied to $L_1 \to M \leftarrow L_2$, 
there are $b \in \cat{A}$ and $t \in \cat{S}$ 
such that the part of diagram (\ref{eqn:4931}) that lies above $M$ is 
commutative:
\begin{equation} \label{eqn:4931}
\UseTips  \xymatrix @C=8ex @R=6ex  {
L_1
\ar[d]_{s'_1}
\ar[drr]_(0.75){a'_1}
& 
L
\ar[l]_{t}
\ar[r]^{b}
&
L_2
\ar[dll]^(0.75){s'_2}
\ar[d]^{a'_2}
\\
M 
& 
&
N
} 
\end{equation}
Write
$s := s'_1 \circ t = s'_2 \circ b$,  $a_1 := a'_1 \circ t$ and 
$a_2 := a'_2 \circ b$. 
By Lemma \ref{lem:3} we get 
$q_i = \opn{Q}(a_i) \circ \opn{Q}(s)^{-1}$.
\end{proof}

\begin{prop} \label{prop:3380} 
Let $\cat{A}$ be a category, let $\cat{S} \sub \cat{A}$ be a right denominator 
set, and let $\cat{A}' \sub \cat{A}$ be a full subcategory.  
Define $\cat{S}' := \cat{A}' \cap \, \cat{S}$. 
Assume these two conditions hold\tup{:}
\begin{itemize}
\rmitem{i} $\cat{S}'$ is a right denominator set in $\cat{A}'$.

\rmitem{ii} Let $M \in \opn{Ob}(\cat{A})$. If 
there exists a morphism $s : M \to L'$ in $\cat{S}$ with 
$L' \in \opn{Ob}(\cat{A}')$, 
there exists a morphism $t : K' \to M$ in $\cat{S}$ with 
$K' \in \opn{Ob}(\cat{A}')$.
\end{itemize}
Then the canonical functor 
$\cat{A}'_{\cat{S}'} \to \cat{A}_{\cat{S}}$ 
is fully faithful.
\end{prop}

\begin{proof}
Let's denote the inclusion functor by $F : \cat{A}' \to \cat{A}$. 
We want to prove that its localization 
$F_{\cat{S}'} : \cat{A}'_{\cat{S}'} \to \cat{A}_{\cat{S}}$ 
is fully faithful.

\medskip \noindent
Step 1. Let $L'_1, L'_2 \in \opn{Ob}(\cat{A}')$, and let 
$q : L'_1 \to L'_2$ be a
morphism in $\cat{A}_{\cat{S}}$. Choose a presentation 
$q = \opn{Q}(a) \circ \opn{Q}(s)^{-1}$ with
$s : M \to L'_1$ a morphism in $\cat{S}$
and $a : M \to L'_2$  a morphism in $\cat{A}$. 
This is possible because $\cat{S}$ is a right denominator set in $\cat{A}$. 
By condition (ii) we can find a morphism $t : K' \to M$ in $\cat{S}$ 
with $K' \in \opn{Ob}(\cat{A}')$. See next diagram. 
\[ \UseTips  \xymatrix @C=6ex @R=6ex  {
K'
\ar[r]^{t}
& 
M 
\ar[dl]_{s}
\ar[dr]^{a}
\\
L'_1
\ar@{-->}[rr]^{q}
& &
L'_2
} \]
Then  
$q = \opn{Q}(a \circ t) \circ \opn{Q}(s \circ t)^{-1}$.
But $s \circ t \in \cat{S}'$ and $a \circ t  \in \cat{A}'$, 
so $q$ is in the image of the functor $F_{\cat{S}'}$. We see that 
$F_{\cat{S}'}$ is full. 

\medskip \noindent
Step 2. Let $p', q' : L'_1 \to L'_2$ 
be morphisms in $\cat{A}'_{\cat{S}'}$
such that $F_{\cat{S}'}(p') = F_{\cat{S}'}(q')$.  
Let us denote the localization functor of $\cat{A}'$ by
$\opn{Q}' : \cat{A}' \to \cat{A}'_{\cat{S}'}$. 
Because $\cat{S}'$ is a right denominator set in $\cat{A}'$,
and using Proposition \ref{prop:107}, we can find presentations 
$p' = \opn{Q}'(a') \circ \opn{Q}'(s')^{-1}$
and
$q' = \opn{Q}'(b') \circ \opn{Q}'(s')^{-1}$
with morphisms 
$s' : N' \to L'_1$ in $\cat{S}'$
and 
$a', b' : N' \to L'_2$ in $\cat{A}'$.
Next, since $F_{\cat{S}'}(p') = F_{\cat{S}'}(q')$, Lemma \ref{lem:3} tells us 
that there are morphisms 
$u, v : M \to N'$ in $\cat{A}$ s.t.\ 
$a' \circ u = b' \circ v$, and 
$s' \circ u = s' \circ v \in \cat{S}$. 
Condition (ii), applied to the morphism 
$s' \circ u : M \to L'_1$, says that there is a morphism 
$t : K' \to M$ in $\cat{S}$ with source $K' \in \opn{Ob}(\cat{A}')$.
See diagram below. 
\[ \UseTips  \xymatrix @C=6ex @R=6ex  {
K'
\ar[r]^{t}
&
M 
\ar[r]^{u, v}
&
N'
\ar[dl]_{s'}
\ar[dr]^{a', b'}
\\
&
L'_1
\ar@{-->}[rr]^{p', q'}
& &
L'_2
} \]
Then we have 
\[ \begin{aligned}
&
p' = \opn{Q}'(a' \circ u \circ t) \circ \opn{Q}'(s' \circ u \circ t)^{-1} 
\\ & \quad 
= \opn{Q}'(b' \circ u \circ t) \circ \opn{Q}'(s' \circ u \circ t)^{-1} = q' .
\end{aligned} \] 
This proves that $F_{\cat{S}'}$ is faithful.
\end{proof}

In the next subsection we study Ore localization of linear categories, and give 
rings as examples (see Examples \ref{exa:4207} and \ref{exa:1405}). And in 
Section \ref{sec:der-cat} we study localization of triangulated categories. 
Here is a nonlinear example. 

\begin{exa} \label{exa:4205}
Fix a base commutative ring $\K$ and a topological space (or, more generally, a 
site) $X$. Let $\K_X$ be the constant sheaf of rings $\K$ on $X$.
Consider the category
$\catt{DGRng}^{\leq 0}_{\mrm{sc}} / \K_X$
of sheaves of commutative DG $\K_X$-rings on $X$ (see Definition 
\ref{dfn:3090} and Remark \ref{rem:4945}). For the sake of brevity let us 
just write 
$\catt{DGRng}_{X} := \catt{DGRng}^{\leq 0}_{\mrm{sc}} / \K_X$.
 
The localization of $\catt{DGRng}_X$ with respect to the quasi-isomorphisms in 
it is called the {\em derived category of commutative DG $\K_X$-rings}, and the 
notation is $\dcat{D}(\catt{DGRng}_X)$. 
The categorical localization functor is 
$\opn{Q} : \catt{DGRng}_X \to \dcat{D}(\catt{DGRng}_X)$,
and it is the identity on objects. 
(The unfortunate tradition in homotopy theory is to call such a localization 
the ``homotopy category'', and to denote it by 
``$\opn{Ho}(\catt{DGRng}_X)$''.)

As explained in the lecture notes \cite{Ye15}, every DG $\K_X$-ring $\AA$ 
admits a {\em semi-pseudo-free resolution} $\til{\AA} \to \AA$ in 
$\catt{DGRng}_X$. The semi-pseudo-free DG $\K_X$-rings have good lifting 
properties. Indeed, they behave like cofibrant objects in a Quillen model 
structure (even though there does not seems to be such a structure on 
$\catt{DGRng}_X$, except when $X$ is a discrete space).

There is a congruence on the category $\catt{DGRng}_X$, called the {\em 
quasi-homotopy relation}. (It turns out that this idea was already known to some 
homotopy theorists, see Remark \ref{rem:4205} below.) Let us denote by 
$\dcat{K}(\catt{DGRng}_X)$ the {\em homotopy category}, which is the 
quotient of $\catt{DGRng}_X$ modulo this congruence. So there is a functor 
$\opn{P} : \catt{DGRng}_X \to \dcat{K}(\catt{DGRng}_X)$,
which is the identity on objects, and it sends a DG ring homomorphism to its 
quasi-homotopy class. 

The functors $\opn{Q}$ and $\opn{P}$ fit into this commutative diagram:
\begin{equation} \label{eqn:4913}
\UseTips \xymatrix @C=8ex @R=6ex {
\catt{DGRng}_X
\ar[r]^{\opn{P}}
\ar@(u,u)[rr]^{\opn{Q}}
&
\dcat{K}(\catt{DGRng}_X)
\ar[r]^{\bar{\opn{Q}}}
&
\dcat{D}(\catt{DGRng}_X)
}
\end{equation}
The functor $\bar{\opn{Q}}$ is a {\em faithful right Ore localization}. 
The proof relies on the lifting properties of semi-pseudo-free DG $\K_X$-rings. 
See Remark \ref{rem:4945} for related material. 
\end{exa}

\begin{rem} \label{rem:4205}
Suppose $\cat{C}$ is a category, with a given multiplicatively closed set of 
morphisms $\cat{S} \sub \cat{C}$. If $\cat{C}$ admits a {\em Quillen model 
structure} for which $\cat{S}$ is the set of {\em weak equivalences}, then 
there is automatically a {\em calculus of fractions} for the localization
functor $\opn{Q} : \cat{C} \to \cat{C}_{\cat{S}}$.

For instance, suppose that all objects of $\cat{C}$ are {\em fibrant}. 
In this case $\cat{C}$ is called a {\em category of fibrant objects}. 
Then there is a congruence relation $\cat{H}$ on $\cat{C}$,
and a commutative diagram 
\[ \UseTips \xymatrix @C=8ex @R=6ex {
\catt{C}
\ar[r]^{\opn{P}}
\ar@(u,u)[rr]^{\opn{Q}}
&
\cat{C} / \cat{H}
\ar[r]^{\bar{\opn{Q}}}
&
\cat{C}_{\cat{S}}
} \]
in which, like in (\ref{eqn:4913}), the functor $\opn{P}$ is full, and the 
functor $\bar{\opn{Q}}$ is a right Ore localization. 
These facts are explained in \cite{nLab} (search for the emphasized text in 
this remark). 
\end{rem}

\begin{rem} \label{rem:4970}
This is a continuation of Remark \ref{rem:4960} on strictly unital 
$\mrm{A}_{\infty}$ categories%
\index{Ainfty@$\mrm{A}_{\infty}$ category}, 
and it is related to Example \ref{exa:4205} and Remark \ref{rem:4205}
above, and to Remark \ref{rem:4605} below. The references are \cite{Kel5} and 
\cite{COS}. Recall that the base ring $\K$ is a field. 

Let $\cat{A}$ and $\cat{B}$ be $\mrm{A}_{\infty}$ categories,
i.e.\ objects of the category $\cat{A_{\infty}Cat}$. 
An $\mrm{A}_{\infty}$ functor $F : \cat{A} \to \cat{B}$
is called a {\em quasi-equivalence} if 
$\opn{H}(F) : \opn{H}(\cat{A}) \to \opn{H}(\cat{B})$
is an equivalence of graded categories.
The {\em derived category} of $\cat{A_{\infty}Cat}$
is the localization $\dcat{D}(\cat{A_{\infty}Cat})$ of $\cat{A_{\infty}Cat}$ 
with respect to the quasi-equivalences. 
(As mentioned above, the tradition is to call 
$\dcat{D}(\cat{A_{\infty}Cat})$ the ``homotopy category'' of 
$\cat{A_{\infty}Cat}$, and to denote 
it by ``$\opn{Ho}(\cat{A_{\infty}Cat})$''.)
There is a localization functor 
$\opn{Q} : \cat{A_{\infty}Cat} \to \dcat{D}(\cat{A_{\infty}Cat})$, 
which is the identity on objects. 

Suppose $F, F' : \cat{A} \to \cat{B}$ are $\mrm{A}_{\infty}$ functors. So 
$F$ and $F'$ are objects of the $\mrm{A}_{\infty}$ category 
$\cat{A_{\infty}Fun}(\cat{A}, \cat{B})$.
We say that $F$ and $F'$ are {\em $\mrm{A}_{\infty}$ homotopic} if 
they are isomorphic in the category 
$\opn{H}^0 \bigl( \cat{A_{\infty}Fun}(\cat{A}, \cat{B}) \bigr)$, 
in which case we write $F \approx F'$. 
Define $\dcat{K}(\cat{A_{\infty}Cat})$
to be the quotient of $\cat{A_{\infty}Cat}$ modulo the relation $\approx$. 
There is a full functor 
$\opn{P} : \cat{A_{\infty}Cat} \to \dcat{K}(\cat{A_{\infty}Cat})$, 
which is the identity on objects. 
An $\mrm{A}_{\infty}$ functor $F : \cat{A} \to \cat{B}$
is called an {\em $\mrm{A}_{\infty}$ equivalence} if it is an isomorphism in 
$\dcat{K}(\cat{A_{\infty}Cat})$. This means that there is  
an $\mrm{A}_{\infty}$ functor $G : \cat{B} \to \cat{A}$, 
such that $G \circ F \approx \opn{Id}_{\cat{A}}$ 
and $F \circ G \approx \opn{Id}_{\cat{B}}$.
There is a commutative diagram 
\begin{equation} \label{eqn:4975}
\UseTips \xymatrix @C=8ex @R=6ex {
\cat{A_{\infty}Cat}
\ar[r]^{\opn{P}}
\ar@(u,u)[rr]^{\opn{Q}}
&
\dcat{K}(\cat{A_{\infty}Cat})
\ar[r]^{\bar{\opn{Q}}}
&
\dcat{D}(\cat{A_{\infty}Cat})
}
\end{equation}

An important theorem (see \cite[Theorem 9.2.0.4 ]{Lef}) says that 
an $\mrm{A}_{\infty}$ functor $F$ is a quasi-equivalence iff it is an
 $\mrm{A}_{\infty}$ equivalence. This implies that the functor 
$\bar{\opn{Q}}$ in diagram (\ref{eqn:4975}) is an isomorphism of categories. 

It is worthwhile to mention that the category $\cat{A_{\infty}Cat}$
does not have a Quillen model structure (for which the weak equivalences are 
the quasi-equivalences). 

Let $\cat{DGCat}$ be the category of $\K$-linear DG categories, where the 
morphisms are the DG functors (see Definition \ref{dfn:1065}). 
There is a fully faithful embedding 
$\cat{DGRng} \centover \K \to \cat{DGCat}$,
since the DG rings (see Definition \ref{dfn:1074}) are the 
single-object DG categories.
Because an $\mrm{A}_{\infty}$ category whose higher operations vanish is a DG 
category, there is a faithful (but not full) embedding
$\cat{DGCat} \to \cat{A_{\infty}Cat}$.
The {\em derived category} of $\cat{DGCat}$
is the localization $\dcat{D}(\cat{DGCat})$ of $\cat{DGCat}$ with 
respect to the quasi-equivalences in it. 
(Once more, the traditional notation for $\dcat{D}(\cat{DGCat})$
is ``$\opn{Ho}(\cat{DGCat})$''.)
The embedding of $\cat{DGCat}$ into $\cat{A_{\infty}Cat}$ induces a functor 
$\dcat{D}(\cat{DGCat}) \lb \to \dcat{D}(\cat{A_{\infty}Cat})$,
and this is an equivalence. Together with the 
equivalence $\bar{\opn{Q}}$ from (\ref{eqn:4975}), this implies that
every morphism $F : \cat{A} \to \cat{B}$ in 
$\dcat{D}(\cat{DGCat})$
can be represented by an $\mrm{A}_{\infty}$ functor 
$\til{F} : \cat{A} \to \cat{B}$, which is unique up to a unique isomorphism in 
$\dcat{K}(\cat{A_{\infty}Cat})$. 

Moreover, the category $\dcat{D}(\cat{DGCat})$ is 
symmetric monoidal, with tensor operation $\cat{A} \ot \cat{B}$. 
It turns out that the DG category 
$\cat{A_{\infty}Fun}(\cat{A}, \cat{B})$
is the {\em internal Hom object}, making  
$\dcat{D}(\cat{DGCat})$ into a {\em closed monoidal category}.
\end{rem}

\mysubsection{Localization of Linear Categories} \label{subsec:loc-lin}

Until now in this section we dealt with arbitrary categories. 
In this subsection our categories will be linear over the 
commutative base ring $\K$ (which will be implicit most of the time).
See Convention \ref{conv:2490}. 

For convenience we only talk about right denominator sets here. 
All the statements hold equally for left denominator sets; cf.\ Remark 
\ref{rem:1845} and Proposition \ref{prop:1845}.

\begin{thm} \label{thm:104}
\index{Localization! Ore {\indash} of a linear category}
Let $\cat{A}$ be a $\K$-linear category, let $\cat{S}$ be a right
denominator set in $\cat{A}$, and let $(\cat{A}_{\cat{S}}, \opn{Q})$ be the 
right Ore
localization.

\begin{enumerate}
\item  The category $\cat{A}_{\cat{S}}$ has a unique $\K$-linear structure,
such that $\opn{Q} : \cat{A} \to \cat{A}_{\cat{S}}$ is a $\K$-linear 
functor. 

\item Suppose $\cat{B}$ is another $\K$-linear category, and 
$F : \cat{A} \to \cat{B}$ is a $\K$-linear functor s.t.\
$F(s)$ is invertible for every $s \in \cat{S}$. Let 
$F_{\cat{S}} : \cat{A}_{\cat{S}} \to \cat{B}$ be the localization of $F$. Then 
$F_{\cat{S}}$ is a $\K$-linear functor. 

\item If $\cat{A}$ is an additive category, then so is $\cat{A}_{\cat{S}}$.
\end{enumerate}
\end{thm}

\begin{proof} \mbox{}

\smallskip \noindent
(1) Let $q_1, q_2 : M \to N$ be morphisms in $\cat{A}_{\cat{S}}$.
Choose common denominator presentations
$q_i = \opn{Q}(a_i) \circ \opn{Q}(s)^{-1}$. Since $\opn{Q}$ must be an additive 
functor, we have to define 
\begin{equation} \label{eqn:1430}
\opn{Q}(a_1) + \opn{Q}(a_2) := \opn{Q}(a_1 + a_2) .
\end{equation}
By the distributive law (bilinearity of composition) we must define 
$q_1 + q_2 := \opn{Q}(a_1 + a_2) \circ \opn{Q}(s)^{-1}$.
And for $\la \in \K$ we must define 
$\la \cdot q_i := \opn{Q}(\la \cd a_i) \circ \opn{Q}(s)^{-1}$.

The usual tricks are then used to prove independence of representatives. 
For instance, to prove that (\ref{eqn:1430}) is independent of choices, suppose 
that $\opn{Q}(a_1) = \opn{Q}(a'_1)$ and $\opn{Q}(a_2) = \opn{Q}(a'_2)$. 
Then, by (RO4), there are $t_1, t_2 \in \cat{S}$ such that 
$a_1 \circ t_1 = a'_1 \circ t_1$ and $a_2 \circ t_2 = a'_2 \circ t_2$. 
By (RD1) there exist $b \in \cat{A}$ and $v \in \cat{S}$ s.t.\ 
$t_1 \circ b = t_2 \circ v$. 
Let $t_3 := t_2 \circ v \in \cat{S}$.
Then 
$(a_1 + a_2) \circ t_3 = (a'_1 + a'_2) \circ t_3$,
and hence 
$\opn{Q}(a_1 + a_2) = \opn{Q}(a'_1 + a'_2)$. 

In this way $\cat{A}_{\cat{S}}$ is a $\K$-linear category, and $\opn{Q}$ is a 
$\K$-linear functor. 

\medskip \noindent (2)
The only option for $F_{\cat{S}}$ is 
$F_{\cat{S}}(q_i) := F(a_i) \circ F(s)^{-1}$. 
The usual tricks are used to prove independence of representatives.

\medskip \noindent 
(3) Clear from Propositions \ref{prop:2} and \ref{prop:109}. 
\end{proof}

\begin{exer} \label{exer:5025}
Finish the proofs of items (1) and (2) of the theorem. 
\end{exer}

\begin{rem} \label{rem:4202}
Let $A$ be a ring, which we can think of as a one-object linear category
$\cat{A}$. In this context, Theorem
\ref{thm:104} is one of the most important results in noncommutative ring 
theory. See the books \cite{McRo}, \cite{Row} or \cite{Ste}. 
\end{rem}

Here are somewhat concrete examples. 

\begin{exa} \label{exa:1405}
Suppose $A$ is a noncommutative ring (i.e.\ it is not necessarily commutative), 
and $S$ is a central  multiplicatively closed set in $A$ (i.e.\ 
$S \sub \opn{Cent}(A)$). Because the elements of $S$ commute with all elements 
of $A$, the denominator conditions hold automatically. 
We get a left and right Ore localization of rings 
$\opn{Q} : A \to A_S$.

Note that if $S$ contains a nilpotent element, then the ring $A_S$ is trivial. 
\end{exa}

The last observation should serve as a warning: localization can sometimes 
kill everything. Fortunately, the localization procedure that 
gives rise to the derived category does not cause any catastrophe, as we shall 
see in Proposition \ref{prop:1410}.

\begin{exa} \label{exa:4207}
Suppose $A$ is a noncommutative ring that is both noetherian and an integral 
domain (i.e.\ the only zero-divisor in $A$ is $0$). Let $S$ be the set of 
nonzero elements of $A$. Then $S$ is both a left and a right denominator set. 
See \cite[Theorem 2.1.15]{McRo}. 
\end{exa}
 
\begin{rem} \label{rem:1}
Suppose $A$ is a ring and $S$ is a right denominator set in it. 
Then the right Ore localization $A_S$ is {\em flat} as left $A$-module. 
See \cite[Theorem 3.1.20]{Row}. We have no idea if something like this is true 
for linear categories with more than one object. 
\end{rem}

\begin{prop} \label{prop:1840}
Let $(\cat{K}, \opn{T})$ be a T-additive category, let
$\cat{S}$ be a right denominator set in $\cat{K}$ such that 
$\opn{T}(\cat{S}) = \cat{S}$, and let 
$\opn{Q} : \cat{K} \to  \cat{K}_{\cat{S}}$ be the localization functor.
\begin{enumerate}
\item There is a unique additive automorphism
$\opn{T}_{\cat{S}}$ of the category $\cat{K}_{\cat{S}}$, such that 
$\opn{T}_{\cat{S}} \circ \opn{Q} = \opn{Q} \circ \opn{T}$
as functors $\cat{K} \to  \cat{K}_{\cat{S}}$.

\item Let $\tau$ be the identity automorphism of the functor 
$\opn{Q} \circ \opn{T}$. Then 
$(\opn{Q}, \tau) :  (\cat{K}, \opn{T}) \to 
(\cat{K}_{\cat{S}}, \opn{T}_{\cat{S}})$
is a T-additive functor. 
\end{enumerate}
\end{prop}

\begin{proof} \mbox{}

\smallskip \noindent
(1) By the assumption the functor 
$\opn{Q} \circ \opn{T} : \cat{K} \to  \cat{K}_{\cat{S}}$
sends the morphisms in $\cat{S}$ to isomorphisms. By the property (Loc3) of 
localization in Definition \ref{dfn:1406}, the functor 
$\opn{T}_{\cat{S}} : \cat{K}_{\cat{S}} \to  \cat{K}_{\cat{S}}$ 
satisfying
$\opn{T}_{\cat{S}} \circ \opn{Q} = \opn{Q} \circ \opn{T}$
exists and is unique. 
Similarly, there is a unique functor 
$\opn{T}^{-1}_{\cat{S}} : \cat{K}_{\cat{S}} \to  \cat{K}_{\cat{S}}$ 
satisfying
$\opn{T}^{-1}_{\cat{S}} \circ \opn{Q} = \opn{Q} \circ \opn{T}^{-1}$.
An easy calculation shows that 
$\opn{T}^{-1}_{\cat{S}} \circ \opn{T}_{\cat{S}} = \opn{Id} = 
\opn{T}_{\cat{S}} \circ \opn{T}^{-1}_{\cat{S}}$.
Hence $\opn{T}_{\cat{S}}$ is an automorphism of $\cat{K}_{\cat{S}}$.
By Theorem \ref{thm:104} the functor $\opn{T}_{\cat{S}}$ is additive. 

\medskip \noindent 
(2) This is clear. 
\end{proof}

The composition of T-additive functors was defined in Definition \ref{dfn:1840}.

\begin{prop} \label{prop:1841}
In the situation of Proposition \tup{\ref{prop:1840}}, suppose 
$(\cat{K}', \opn{T}')$ is another  T-additive category, and 
$(F, \nu) :  (\cat{K}, \opn{T}) \to (\cat{K}', \opn{T}')$
is a T-additive functor, such that $F(s)$ is invertible for every 
$s \in \cat{S}$. Let 
$F_{\cat{S}} : \cat{K}_{\cat{S}} \to \cat{K}'$
be the localized functor. 
Then there is a unique isomorphism 
$\nu_{\cat{S}} : F_{\cat{S}} \circ \opn{T}_{\cat{S}} \iso
\opn{T}' \circ \, F_{\cat{S}}$
of functors $\cat{K}_{\cat{S}} \to  \cat{K}'$, such that
$(F, \nu) = (F_{\cat{S}}, \nu_{\cat{S}}) \circ (\opn{Q}, \tau)$ 
as T-additive functors 
$(\cat{K}, \opn{T}) \to (\cat{K}', \opn{T}')$.
\end{prop}

\begin{exer}  \label{exer:1840}
Prove Proposition \ref{prop:1841}. 
\end{exer}

\cleardoublepage
\mysection{The Derived Category 
\texorpdfstring{$\dcat{D}(A, \cat{M})$}{D(A,M)}} 
\label{sec:der-cat}

\AYcopyright

In this section we introduce the main mathematical concept of the book: the 
{\em derived category} $\dcat{D}(A, \cat{M})$ of DG $A$-modules in $\cat{M}$. 
Here $A$ is a central DG $\K$-ring, and $\cat{M}$ is a $\K$-linear abelian 
category. The base ring $\K$ can be any nonzero commutative ring, and it will 
remain implicit most of the time.

\mysubsection{Localization of Triangulated Categories}

Let $\cat{K}$ be a triangulated category, with translation functor $\opn{T}$. 
As we did in Section \ref{sec:loc-cats}, we shall often write 
$a \in \cat{K}$ for a morphism 
$a \in \cat{K}(M, N) = \opn{Hom}_{\cat{K}}(M, N)$, leaving the objects 
implicit. 

\begin{prop} \label{prop:106}
Suppose 
$H : \cat{K} \to \cat{M}$ is a cohomological functor, where $\cat{M}$ is some
abelian category. Let 
\[ \cat{S} := \bigl\{ s \in \cat{K} \mid H( \opn{T}^i(s)) \tup{ is invertible 
for all } i \in \Z \bigr\} . \] 
Then $\cat{S}$ is a left and right denominator set in $\cat{K}$.
\end{prop}

\begin{proof}
It is clear that $\cat{S}$ is closed under composition and contains the
identity morphisms. So it is a multiplicatively closed set. 

Let's prove that condition (RD1) of Definition \ref{dfn:21} holds. Suppose we 
are given morphisms
$L \xar{a} N \xleftarrow{s} M$ with $s \in \cat{S}$. 
We need to find morphisms 
$L \xleftarrow{t} K \xar{b} M$ with $t \in \cat{S}$ and such that 
$a \circ t = s \circ b$.  
 
Consider the solid commutative diagram 
\[ \UseTips  \xymatrix @C=8ex @R=6ex { 
K
\ar[r]^{t}
\ar@{-->}[d]_{b}
&
L
\ar[r]^{c \circ a}
\ar[d]_{a}
& 
P
\ar[r]^{}
\ar[d]_{=}
& 
\opn{T}(K) 
\ar@{-->}[d]_{\opn{T}(b)}
\\
M
\ar[r]^{s}
&
N
\ar[r]^{c}
& 
P
\ar[r]^{}
&
\opn{T}(M) 
}  \]
where the bottom row is a distinguished triangle built on $M \xar{s} N$,
and and the top row is a distinguished triangle built on 
$L \xar{c \circ a} P$, then turned $120^\circ$ to the right. By axiom (TR3) 
there 
is a morphism $b$ making the diagram commutative. Thus 
$a \circ t = s \circ b$.  Since $H(\opn{T}^i(s))$ are
invertible for all $i \in \Z$, it follows that $H(\opn{T}^i(P)) = 0$. But 
then $H(\opn{T}^i(t))$ are invertible for all $i \in \Z$, so 
$t \in \cat{S}$. 

Next we prove condition (RD2) of Definition \ref{dfn:21}. Because we are in 
an additive category, this condition is simplified: given $a \in \cat{K}$ and 
$s \in \cat{S}$ satisfying 
$s \circ a = 0$, we have to find $t \in \cat{S}$ satisfying 
$a \circ t = 0$. 

Say the objects involved are $L \xar{a} M \xar{s} N$. 
Take a distinguished triangle built on $s$ and then turned:
$P \xar{b} M \xar{s} N \xar{} \opn{T}(P)$.
We get an exact sequence 
\[ \opn{Hom}_{\cat{K}}(L, P) \xar{b \, \circ \, (-)} 
\opn{Hom}_{\cat{K}}(L, M) \xar{s \, \circ \, (-)} 
\opn{Hom}_{\cat{K}}(L, N) . \]
Since $s \circ a = 0$, there is $c : L \to P$ s.t.\ $a = b \circ c$. 
Now look at a distinguished triangle  built on $c$, and then turned:
$K \xar{t} L \xar{c} P \xar{} \opn{T}(K)$.
We know that $c \circ t = 0$; hence 
$a \circ t = b \circ c \circ t = 0$. 
But ($s \in \cat{S}$) $\Rightarrow$ ($H(\opn{T}^i(P)) = 0$ for all $i$)
$\Rightarrow$ ($t \in \cat{S}$).

The left versions (LD1) and (LD2) are proved the same way. 
\end{proof}

\begin{dfn} \label{dfn:1495}
A {\em denominator set of cohomological origin}
\index{Denominator set! of cohomological origin}
in $\cat{K}$ is a denominator set $\cat{S} \sub \cat{K}$ 
that arises from a cohomological functor $H$, as in Proposition \ref{prop:106}.
The morphisms in $\cat{S}$ are called {\em quasi-isomorphisms relative to $H$}.
\end{dfn}

\begin{thm} \label{thm:105}
Let $(\cat{K}, \opn{T})$ be a triangulated category, let $\cat{S}$ be a 
denominator set of cohomological origin in $\cat{K}$, 
\index{Denominator set! of cohomological origin}
and let 
$(\opn{Q}, \tau) : (\cat{K}, \opn{T}) \to
(\cat{K}_{\cat{S}}, \opn{T}_{\cat{S}})$
be the T-additive functor from Proposition \tup{\ref{prop:1840}}. Then the 
T-additive category 
$(\cat{K}_{\cat{S}}, \opn{T}_{\cat{S}})$
has a unique triangulated structure 
\index{Triangulated category}
such that these two properties 
hold\tup{:}

\begin{enumerate}
\rmitem{i} The pair $(\opn{Q}, \tau)$ is a triangulated functor. 
\index{Triangulated functor}

\rmitem{ii} Suppose 
$(\cat{K}', \opn{T}')$ 
is another triangulated category, and 
$(F, \nu) : (\cat{K}, \opn{T}) \to (\cat{K}', \opn{T}')$
is a triangulated functor, such that
$F(s)$ is invertible for every $s \in \cat{S}$. Let 
$(F_{\cat{S}}, \nu_{\cat{S}}) : (\cat{K}_{\cat{S}}, \opn{T}_{\cat{S}}) \to 
(\cat{K}', \opn{T}')$
be the T-additive functor from Proposition \tup{\ref{prop:1841}}.
Then $(F_{\cat{S}}, \nu_{\cat{S}})$ is a triangulated functor. 
\end{enumerate}
\end{thm}

This result is stated as \cite[Proposition I.3.2]{RD} and 
as \cite[Proposition 1.6.9]{KaSc1}. Both sources do not give proofs (just 
hints). 

\begin{proof} 
Since $\cat{S}$ is of cohomological origin we have 
$\opn{T}(\cat{S}) = \cat{S}$. Recall that the translation isomorphism $\tau$ is 
the identity automorphism of the functor $\opn{Q} \circ \opn{T}$; see 
Proposition \ref{prop:1840}. So we will ignore it. 

\medskip \noindent
Step 1.
The distinguished triangles in $\cat{K}_{\cat{S}}$ are defined to be those
triangles that are isomorphic to the images under $\opn{Q}$ of distinguished 
triangles in $\cat{K}$. Let us verify the axioms of triangulated 
category. 

\medskip \noindent
(TR1). By definition every triangle that's isomorphic to a distinguished 
triangle is distinguished; and the triangle 
$M \xar{\opn{id}_M} M \to 0 \to \opn{T}(M)$ 
in $\cat{K}_{\cat{S}}$ is clearly distinguished.

Suppose we are given a morphism $\al : L \to M$ in $\cat{K}_{\cat{S}}$. We have 
to build a distinguished triangle on it. 
Choose a fraction presentation 
$\al = \opn{Q}(a) \circ \opn{Q}(s)^{-1}$.
Using condition (LD1) we can find $b \in \cat{K}$ and $t \in 
\cat{S}$ such that $t \circ a = b \circ s$.
These fit into the solid commutative diagram
\[ \UseTips  \xymatrix @C=8ex @R=4ex { 
L
\ar@(ur,ul)@{-->}[rr]^{\al}
\ar[dr]_{b}
&
K
\ar[l]_{s}
\ar[r]^{a}
&
M
\ar[dl]^{t}
\\
&
\til{K}
}  \]
in $\cat{K}$. (The dashed arrow $\al$ is in $\cat{K}_{\cat{S}}$.) 

Consider the solid commutative diagram below, where
the rows are  distinguished triangles  built on $a$  and $b$ respectively.
\begin{equation} \label{eqn:1435}
 \UseTips  \xymatrix @C=8ex @R=6ex { 
K
\ar[r]^{a}
\ar[d]_{s}
&
M
\ar[r]^{e}
\ar[d]_{t}
& 
N
\ar[r]^{c}
\ar@{-->}[d]_{u}
&
\opn{T}(K)
\ar[d]_{\opn{T}(s)}
\\
L
\ar[r]^{b}
&
\til{K}
\ar[r]^{}
& 
P
\ar[r]^{d}
&
\opn{T}(L)
} 
\end{equation}
By (TR3) there is a morphism $u$ that makes the whole diagram commutative.
Since $s, t \in \cat{S}$ and $H$ is a cohomological functor, it follows that 
$u \in \cat{S}$. 
Applying the functor $\opn{Q}$ to (\ref{eqn:1435}), and using the 
isomorphism 
$\opn{Q}(t) : M \to \til{K}$ 
to replace $\til{K}$ with $M$,  we get the commutative diagram
\[ \UseTips  \xymatrix @C=8ex @R=6ex { 
K
\ar[r]^{\opn{Q}(a)}
\ar[d]_{\opn{Q}(s)}
&
M
\ar[r]^{\opn{Q}(e)}
\ar[d]_{\opn{Q}(\opn{id}_M)}
& 
N
\ar[r]^{\opn{Q}(c)}
\ar[d]_{\opn{Q}(u)}
&
\opn{T}(K)
\ar[d]_{\opn{T}(\opn{Q}(s))}
\\
L
\ar[r]^{\al}
&
M
\ar[r]^{\opn{Q}(u \circ e)}
& 
P
\ar[r]^{\opn{Q}(d)}
&
\opn{T}(L)
}  \]
in $\cat{K}_{\cat{S}}$. The top row is a distinguished triangle, and the
vertical arrows are isomorphisms. So the bottom row is a distinguished 
triangle. This is the triangle we were looking for. 

\medskip \noindent 
(TR2). Turning: this is trivial. 

\medskip \noindent 
(TR3). We are given the solid commutative diagram in $\cat{K}_{\cat{S}}$,
where the rows are distinguished triangles:
\begin{equation} \label{eqn:104}
\UseTips \xymatrix @C=8ex @R=6ex { 
L 
\ar[r]^{\al}
\ar[d]_{\phi}
&
M
\ar[r]^{\be}
\ar[d]_{\psi}
& 
N
\ar[r]^{\ga}
\ar@{-->}[d]_{\chi}
&
\opn{T}(L)
\ar[d]_{\opn{T}(\phi)}
\\
L' 
\ar[r]^{\al'}
&
M'
\ar[r]^{\be'}
& 
N'
\ar[r]^{\ga'}
&
\opn{T}(L') 
} 
\end{equation}
and we have to find $\chi$ to complete the diagram. 

By replacing the rows with isomorphic triangles, we can assume they come from
$\cat{K}$.
Thus we can replace (\ref{eqn:104}) with this diagram:
\begin{equation} \label{eqn:110}
\UseTips \xymatrix @C=8ex @R=6ex { 
L 
\ar[r]^{\opn{Q}(\al)}
\ar[d]_{\phi}
&
M
\ar[r]^{\opn{Q}(\be)}
\ar[d]_{\psi}
& 
N
\ar[r]^{\opn{Q}(\ga)}
\ar@{-->}[d]_{\chi}
&
\opn{T}(L)
\ar[d]_{\opn{T}(\phi)}
\\
L' 
\ar[r]^{\opn{Q}(\al')}
&
M'
\ar[r]^{\opn{Q}(\be')}
& 
N'
\ar[r]^{\opn{Q}(\ga')}
&
\opn{T}(L') 
} 
\end{equation}
in which $\al, \be, \ga, \al', \be', \ga'$ are morphisms in $\cat{K}$. 
It is a commutative diagram. Let us choose fraction presentations 
$\phi = \opn{Q}(a) \circ \opn{Q}(s)^{-1}$ and $\psi = \opn{Q}(b) \circ 
\opn{Q}(t)^{-1}$. 
Then the solid diagram (\ref{eqn:110}) comes from applying $\opn{Q}$ to the 
diagram 
\begin{equation} \label{eqn:103}
\UseTips \xymatrix @C=8ex @R=5.5ex { 
L 
\ar[r]^{\al}
&
M
\ar[r]^{\be}
& 
N
\ar[r]^{\ga}
&
\opn{T}(L)
\\
\til{L} 
\ar[d]_{a}
\ar[u]^{s}
&
\til{M}
\ar[d]_{b}
\ar[u]^{t}
& 
&
\opn{T}(\til{L})
\ar[d]_{\opn{T}(a)}
\ar[u]^{\opn{T}(s)}
\\
L' 
\ar[r]^{\al'}
&
M'
\ar[r]^{\be'}
& 
N'
\ar[r]^{\ga'}
&
\opn{T}(L') 
} \end{equation}
in $\cat{K}$. Here the rows are distinguished triangles in $\cat{K}$; but the 
diagram might fail to be commutative. 

By axiom (RO3) we can find $c \in \cat{K}$ and $u \in \cat{S}$ s.t.\ 
$\opn{Q}(t)^{-1} \circ \opn{Q}(\al) \circ \opn{Q}(s) = \opn{Q}(c) \circ 
\opn{Q}(u)^{-1}$.
This is the solid diagram below:
\[ \UseTips \xymatrix @C=8ex @R5.5ex { 
& &
L 
\ar[r]^{\al}
&
M
\\
\til{L}''
\ar@{-->}[r]^{u'}
&
\til{L}'
\ar[r]^{u}
\ar@(dr,dl)[rr]_(0.2){c}
&
\til{L} 
\ar[d]_(0.3){a}
\ar[u]^{s}
&
\til{M}
\ar[d]_{b}
\ar[u]^{t}
\\
& &
L' 
\ar[r]^{\al'}
&
M'
} \]
Thus 
$\opn{Q}(\al \circ s \circ u) = \opn{Q}(t \circ c)$.
By (RO4) there is a morphism $u' \in \cat{S}$ s.t.\ 
$(\al \circ s \circ u) \circ u' = (t \circ c) \circ u'$.
We get
\[ \phi = \opn{Q}(a) \circ \opn{Q}(s)^{-1} = \opn{Q}(a \circ u \circ u') \circ 
\opn{Q}(s \circ 
u \circ u')^{-1} \]
in $\cat{K}_{\cat{S}}$. 
Thus, after substituting $\til{L} :=  \til{L}''$,
 $s := s \circ u \circ u'$,
$a := a \circ u \circ u'$
and $c := c \circ u'$, we get a new diagram
\begin{equation} \label{eqn:105}
\UseTips \xymatrix @C=8ex @R=5.5ex { 
L 
\ar[r]^{\al}
&
M
\ar[r]^{\be}
& 
N
\ar[r]^{\ga}
&
\opn{T}(L)
\\
\til{L} 
\ar[r]^{c}
\ar[d]_{a}
\ar[u]^{s}
&
\til{M}
\ar[d]_{b}
\ar[u]^{t}
& 
&
\opn{T}(\til{L})
\ar[d]_{\opn{T}(a)}
\ar[u]^{\opn{T}(s)}
\\
L' 
\ar[r]^{\al'}
&
M'
\ar[r]^{\be'}
& 
N'
\ar[r]^{\ga'}
&
\opn{T}(L') 
} 
\end{equation}
in $\cat{K}$ instead of (\ref{eqn:103}). In this new diagram the top left 
square is commutative; but maybe the bottom left square is not commutative. 

When we apply $\opn{Q}$ to the diagram (\ref{eqn:105}), the  whole diagram, 
including the bottom left square, becomes commutative, since (\ref{eqn:110}) is
commutative. Again using condition (RO4), there is 
$v \in \cat{S}$ s.t.\ 
$(\al' \circ a) \circ v = (b \circ c) \circ v$.
In a diagram:
\[ \UseTips \xymatrix @C=8ex @R=5.5ex { 
& 
L 
\ar[r]^{\al}
&
M
\\
\til{L}'
\ar[r]^{v}
&
\til{L} 
\ar[d]_{a}
\ar[u]^{s}
\ar[r]^{c}
&
\til{M}
\ar[d]_{b}
\ar[u]^{t}
\\
& 
L' 
\ar[r]^{\al'}
&
M'
} \]
Performing the replacements
$\til{L} :=  \til{L}'$,
$s := s \circ v$,
$c := c \circ v$ and 
$a := a \circ v$ we now have a commutative square also at the bottom
left of (\ref{eqn:105}). Since $\ga \circ \be = 0$ and 
$\ga' \circ \be' = 0$,  
in fact the whole diagram (\ref{eqn:105}) in $\cat{K}$ is now commutative.

Now by (TR1) we can embed the morphism $c$ in a distinguished
triangle. We get the solid diagram
\begin{equation} \label{eqn:106}
\UseTips \xymatrix @C=8ex @R=5.5ex { 
L 
\ar[r]^{\al}
&
M
\ar[r]^{\be}
& 
N
\ar[r]^{\ga}
&
\opn{T}(L)
\\
\til{L} 
\ar[r]^{c}
\ar[d]_{a}
\ar[u]^{s}
&
\til{M}
\ar[r]^{\til{\be}}
\ar[d]_{b}
\ar[u]^{t}
& 
\til{N}
\ar[r]^{\til{\ga}}
\ar@{-->}[d]_{d}
\ar@{-->}[u]^{w}
&
\opn{T}(\til{L})
\ar[d]_{\opn{T}(a)}
\ar[u]^{\opn{T}(s)}
\\
L' 
\ar[r]^{\al'}
&
M'
\ar[r]^{\be'}
& 
N'
\ar[r]^{\ga'}
&
\opn{T}(L') 
} 
\end{equation}
in $\cat{K}$. The rows are distinguished triangles. 
Since $\til{\ga} \circ \til{\be} = 0$, the solid diagram is commutative.  
By (TR3) there are morphisms $w$ and $d$ that make the whole diagram
commutative. 
Now by the long exact cohomology sequence (see Theorem \ref{thm:3640})
the morphism $w$ belongs to $\cat{S}$. The morphism 
$\chi := \opn{Q}(d) \circ \opn{Q}(w)^{-1} :  N \to N'$
solves the problem.

\medskip \noindent 
(TR4). We will not give any details, as this axiom is not an important feature 
for us (see Remark \ref{rem:1281}).  A proof can be found as part of the proof 
of \cite[Theorem 2.1.8, page 97]{Ne1}, where this axiom is called (TR4$'$). 

\medskip \noindent
Step 2. Suppose $(F, \nu)$ is a triangulated functor as in condition (ii).
By Proposition \tup{\ref{prop:1841}} this extends uniquely to a 
T-additive functor $(F_{\cat{S}}, \nu_{\cat{S}})$. 
The construction of the triangulated structure on 
$(\cat{K}_{\cat{S}}, \opn{T}_{\cat{S}})$ in the previous steps,
and the defining property of the translation isomorphism $\nu_{\cat{S}}$ in 
Proposition \tup{\ref{prop:1841}}, show that
$(F_{\cat{S}}, \nu_{\cat{S}})$ is a triangulated functor. 

\medskip \noindent
Step 3. At this point $(\cat{K}_{\cat{S}}, \opn{T}_{\cat{S}})$ 
is a triangulated category, and conditions (i)-(ii) of the theorem are 
satisfied. We need to prove the uniqueness of the triangulated structure on
$(\cat{K}_{\cat{S}}, \opn{T}_{\cat{S}})$. Condition (i) says that we can't have 
less distinguished triangles than those we declared. We can't have more
distinguished triangles, because of condition (ii). 
\end{proof}

\begin{prop} \label{prop:1410}
Consider the situation of Proposition \tup{\ref{prop:106}} and Theorem 
\tup{\ref{thm:105}}. 
\begin{enumerate}
\item The cohomological functor 
$H : \cat{K} \to \cat{M}$ factors into 
$H = H_{\cat{S}} \circ \opn{Q}$, where 
$H_{\cat{S}} : \cat{K}_{\cat{S}} \to \cat{M}$ is a cohomological functor.

\item Let $M$ be an object of $\cat{K}$. The object $\opn{Q}(M)$ is zero 
in $\cat{K}_{\cat{S}}$ iff the objects $H(\opn{T}^i((M))$ are zero in $\cat{M}$ 
for all $i$. 
\end{enumerate}
\end{prop}

\begin{proof} \mbox{}

\smallskip \noindent
(1) The existence and uniqueness of the functor $H_{\cat{S}}$ are by the 
universal property (Loc3) in Definition \ref{dfn:1406}. 
We leave it as an exercise to show that $H_{\cat{S}}$ is a cohomological 
functor. 

\medskip \noindent 
(2) Since $H_{\cat{S}}$ is an additive functor, 
if $\opn{Q}(M) = 0$, then so is 
$H(M) = H_{\cat{S}}(\opn{Q}(M))$. And of course $\opn{Q}(M) = 0$ iff
$\opn{Q}(\opn{T}^i(M)) = 0$ for all $i$. 

For the converse, let $\phi: 0 \to M$ be the zero morphism in $\cat{K}$. 
If $H(\opn{T}^i(M)) = 0$ for all $i$, then 
$H(\opn{T}^i(\phi)) : 0 \to H(\opn{T}^i(M))$
are isomorphisms for all $i$. Then $\phi \in \cat{S}$, and so 
$\opn{Q}(\phi) : 0 \to \opn{Q}(M)$
is an isomorphism in  $\cat{K}_{\cat{S}}$.
\end{proof}

\begin{prop} \label{prop:1480}
Let $\cat{K}$ be a triangulated category, let 
$\cat{S}$ be a denominator set of cohomological origin in $\cat{K}$, and let 
$\cat{K}'$ be a full triangulated subcategory of $\cat{K}$. 
Then $\cat{S}' := \cat{K}' \cap \, \cat{S}$ is a denominator set of 
cohomological origin in $\cat{K}'$, the Ore localization 
$\cat{K}'_{\cat{S}'}$ exists, and $\cat{K}'_{\cat{S}'}$ is a triangulated 
category. 
\end{prop}

\begin{proof}
Let $H :  \cat{K} \to \cat{M}$ be a cohomological functor that determines 
$\cat{S}$. The functor $H|_{\cat{K}'} : \cat{K}' \to \cat{M}$ is also 
cohomological, and the set of morphisms $\cat{S}'$ satisfies 
\[ \cat{S}' = \bigl\{ s \in  \cat{K}'
\mid  H|_{\cat{K}'}(\opn{T}^i(s))  \tup{ is an isomorphism for all } i \bigr\}  
. \]
Hence Proposition \ref{prop:106} and Theorem \ref{thm:105} apply.
\end{proof}

In the situation of the proposition, the localization functor is denoted by
$\opn{Q}' : \cat{K}' \to \cat{K}'_{\cat{S}'}$.

\begin{prop} \label{prop:1496}
In the situation of Proposition \tup{\ref{prop:1480}}, 
let $F : \cat{K}' \to \cat{E}$ be a triangulated functor into some 
triangulated category $\cat{E}$. Assume that for every 
$s \in \cat{S}'$, the morphism $F(s)$ is an isomorphism in $\cat{E}$. 
Then there is a unique triangulated functor 
$F_{\cat{S}'} : \cat{K}'_{\cat{S}'} \to \cat{E}$
that extends $F$\tup{;} Namely 
$F_{\cat{S}'} \circ \opn{Q}' = F$ as functors $\cat{K}' \to \cat{E}$.
\end{prop}

\begin{proof}
This is part of Theorem \ref{thm:105}. 
\end{proof}

In particular we can look at the functor  
$F : \cat{K}' \xar{\opn{inc}} \cat{K} \xar{\opn{Q}} \cat{K}_{\cat{S}}$, 
and its extension 
$F_{\cat{S}'} : \cat{K}'_{\cat{S}'} \to \cat{K}_{\cat{S}}$.
We are interested in sufficient conditions for the functor 
$F_{\cat{S}'}$ to be fully faithful.

\begin{prop} \label{prop:1320}  
Let $\cat{K}$ be a triangulated category, let 
$\cat{S}$ be a denominator set of cohomological origin in $\cat{K}$,
and let $\cat{K}' \subseteq \cat{K}$ be a full triangulated subcategory.  
Define 
$\cat{S}' := \cat{K}' \cap \, \cat{S}$. 
Assume either of these conditions holds\tup{:}
\begin{itemize}
\rmitem{r} Let $M \in \opn{Ob}(\cat{K})$. If 
there exists a morphism $s : M \to L$ in $\cat{S}$ with 
$L \in \opn{Ob}(\cat{K}')$, 
there exists a morphism $t : K \to M$ in $\cat{S}$ with 
$K \in \opn{Ob}(\cat{K}')$.

\rmitem{l} Let $M \in \opn{Ob}(\cat{K})$. If 
there exists a morphism $s : L \to M$ in $\cat{S}$ with 
$L \in \opn{Ob}(\cat{K}')$, 
there exists a morphism $t : M \to K$ in $\cat{S}$ with 
$K \in \opn{Ob}(\cat{K}')$.
\end{itemize}
Then the functor 
$F_{\cat{S}'} : \cat{K}'_{\cat{S}'} \to \cat{K}_{\cat{S}}$ 
is fully faithful.
\end{prop}

\begin{proof}
By  Proposition \ref{prop:1480}, $\cat{S}'$ is a multiplicatively closed 
set of cohomological origin. So by Proposition \ref{prop:106}
it is a right and left denominator set in $\cat{K}'$. 
According to Proposition \ref{prop:3380}, under condition (r) the functor 
$F_{\cat{S}'}$ is  fully faithful.

For condition (l) we pass to the opposite categories, using Proposition 
\ref{prop:2310}. 
\end{proof}

\begin{rem} \label{rem:4220}
Let $\cat{K}$ be a triangulated category, and let 
$\cat{N} \sub \cat{K}$ be a full triangulated subcategory. Define 
$\cat{S}$ to be the set of morphisms $s$ in $\cat{K}$ such that the 
cone of $s$ (see Definition \ref{dfn:3600}) is isomorphic to an object of 
$\cat{N}$.  It turns out that $\cat{S}$ is a left and right denominator set in 
$\cat{K}$. We call $\cat{S}$ the {\em denominator set associated to $\cat{N}$}. 
The localization $\cat{K}_{\cat{S}}$ is a triangulated 
category, and it is called the {\em Verdier quotient}%
\index{Verdier quotient}
of $\cat{K}$ by $\cat{N}$, with notation 
$\catt{K} / \catt{N} :=  \cat{K}_{\cat{S}}$. 
See \cite[Chapter 2]{Ver} or \cite[Section 1.6]{KaSc1}. 

We want to relate this to Remark \ref{rem:4777}. 
Let $\til{\cat{S}}$ be a left and right denominator set in $\cat{K}$.
A morphism $s$ in $\cat{K}$ is called a {\em divisor} of $\til{\cat{S}}$ 
if there exist morphisms $u, v \in \cat{K}$ such that 
$u \circ s$ and $s \circ v$ belong to $\til{\cat{S}}$.
We call $\til{\cat{S}}$ a {\em saturated denominator set} if 
it contains all its divisors.

Let $\cat{N} \sub \cat{K}$ be a full triangulated subcategory, and let 
$\cat{S}$ be the associated denominator set, as in the first 
paragraph. Let  $\til{\cat{N}}$ be the saturated 
full triangulated subcategory of $\cat{K}$ generated by $\cat{N}$, in the 
sense of Definition \ref{dfn:4916}. Let $\til{\cat{S}}$ be the set of morphisms
gotten by adjoining to $\cat{S}$ all its divisors. 
Then $\til{\cat{S}}$ is a saturated left and right denominator set of morphisms 
in $\cat{K}$, the canonical functor 
$\cat{K}_{\cat{S}} \to \cat{K}_{\til{\cat{S}}}$
is an isomorphism, the kernel of the localization functor 
$\opn{Q} : \cat{K} \to \cat{K}_{\cat{S}}$
(see Corollary \ref{cor:4920}) is $\til{\cat{N}}$,
and $\til{\cat{S}}$ is the denominator set associated to $\til{\cat{N}}$.
This is proved in \cite[Chapter 2]{Ver} and \cite[Proposition 7.1.20]{KaSc2}.

An important example of a Verdier quotient is the passage from 
the homotopy category $\dcat{K}(A, \cat{M})$ to the derived category 
$\dcat{D}(A, \cat{M})$, that we shall study in Subsection 
\ref{subsec:def-der-cat} below. Indeed, the set of quasi-isomorphisms 
$\dcat{S}(A, \cat{M})$ in $\dcat{K}(A, \cat{M})$
can be described as the denominator set associated to the 
full subcategory $\dcat{N}(A, \cat{M}) \sub \dcat{K}(A, \cat{M})$ on the 
acyclic complexes, which is a saturated full triangulated subcategory. 
\end{rem}

\mysubsection{Definition of the Derived Category} \label{subsec:def-der-cat}

We now specialize to the triangulated category 
$\cat{K} = \dcat{K}(A, \cat{M})$, for a central DG $\K$-ring $A$ and 
a $\K$-linear abelian category $\cat{M}$. 

Recall the cohomology functor 
$\opn{H} : \dcat{C}_{\mrm{str}}(A, \cat{M}) \to 
\dcat{G}_{\mrm{str}}(\cat{M})$
from Definition \ref{dfn:2993}. 

\begin{lem} \label{lem:2320}
Suppose $\phi, \psi : M \to N$ are morphisms in 
$\dcat{C}_{\mrm{str}}(A, \cat{M})$ such that 
$\psi - \phi$ is a coboundary in $\opn{Hom}_{A, \cat{M}}(M, N)$. 
Then $\opn{H}(\phi) = \opn{H}(\psi)$, as morphisms 
$\opn{H}(M) \to \opn{H}(N)$ in $\dcat{G}_{\mrm{str}}(\cat{M})$. 
\end{lem}

\begin{exer} \label{exer:2320}
Prove the lemma. 
\end{exer}

It follows that there is an induced functor 
\begin{equation} \label{eqn:4985}
\opn{H} : \dcat{K}(A, \cat{M}) \to \dcat{G}_{\mrm{str}}(\cat{M}) . 
\end{equation}
We shall be interested in its degree $0$ component 
$\mrm{H}^0$. Of course all other components can be recovered from $\mrm{H}^0$
by the translation functor: 
$\mrm{H}^i = \mrm{H}^0 \circ \opn{T}^i$.

\begin{prop} \label{prop:1475}
Let $\cat{M}$ be an abelian category and let $A$ be a DG ring. Then the functor 
$\mrm{H}^0 : \dcat{K}(A, \cat{M}) \to \cat{M}$
is cohomological. 
\end{prop}

\begin{proof}
Clearly $\mrm{H}^0$ is additive. Consider a distinguished triangle 
\begin{equation} \label{eqn:1437}
L \xar{\al} M \xar{\be} N \xar{\ga} \opn{T}(L)
\end{equation}
in $\dcat{K}(A, \cat{M})$. We must prove that 
\[ \opn{H}^{0}(L) \xar{\opn{H}^0(\al)} \opn{H}^{0}(M)
\xar{\opn{H}^0(\be)} \opn{H}^{0}(N) \]
is an exact sequence in $\cat{M}$. 

By Definition \ref{dfn:140} we can assume that the distinguished triangle 
(\ref{eqn:1437}) is the image of a standard triangle in  
$\dcat{C}(A, \cat{M})$. Namely that $N = \opn{Cone}(\al)$, the standard 
cone associated to $\al$, and the morphisms are 
$\be = e_{\al}$ and $\ga = p_{\al}$; see Definitions  \ref{dfn:1172} and 
\ref{dfn:1140}. To be explicit, in matrix notation we have 
$N = \sbmat{M \\[0.2em] \opn{T}(L)}$, 
$\d_N = \sbmat{\d_M & \al \circ \, \opn{t}_L^{-1} 
\\[0.2em] 0 & \d_{\opn{T}(L)}}$ 
and 
$\be = \sbmat{\opn{id}_M \\[0.2em] 0}$.

Let us use the abbreviations $\bar{\be} := \opn{H}^0(\be)$ and
$\bar{\al} := \opn{H}^0(\al)$. 
By Proposition \ref{prop:1280} we know that $\be \circ \al = 0$, and hence 
$\bar{\be} \circ \bar{\al} = 0$. 
It remains to prove that the morphism 
$\bar{\al}  : \opn{H}^0(L) \to \opn{Ker}(\bar{\be})$
in $\cat{M}$ is an epimorphism. 
This will be done using the first sheaf trick (Proposition \ref{prop:3635}(3)):
we will prove that given a section 
$\bar{m} \in \Ga(U, \opn{Ker}(\bar{\be}))$
on an open set $U \in \cat{M}$, there is a refinement
$V' \surj U$ and a section 
$\bar{l} \in \Ga(V', \opn{H}^0(L))$
such that $\bar{\al}(\bar{l}) = \bar{m}$
in $\Ga(V', \opn{Ker}(\bar{\be}))$. 
See Subsection \ref{subsec:sheaf-tricks} for the ``geometric'' terminology 
that's used in the sheaf tricks. 

Take some section 
$\bar{m} \in \Ga(U, \opn{Ker}(\bar{\be}))$
on an open set $U \in \cat{M}$.
This means that 
$\bar{m} \in \Ga(U, \opn{H}^0(M))$ and 
$\bar{\be}(\bar{m}) = 0$. 
There is a covering $V \surj U$ and a section 
$m \in \Ga(V, \opn{Z}^{0}(M))$ such that 
$\pi_{M}(m) = \bar{m}$. Now 
\[ \be(m) = \bmat{m \\[0.1em] 0} \in \Ga(V, \opn{Z}^{0}(N)) \sub 
\Ga(V, N^{0}) = 
\bmat{ \Ga(V, M^{0}) \\[0.1em] \Ga(V, L^{1})}  \]
in column notation. 
The vanishing of $\bar{\be}(\bar{m}) = \pi_N(\be(m))$ in 
$\Ga(V, \opn{H}^{0}(N))$ means that there is a covering 
$V' \surj V$ and a section 
$n \in \Ga(V', N^{-1})$
such that $\d_N(n) = \be(m)$. 
In column notation we have 
$n = \sbmat{m' \\[0.1em] l}$,
with $m' \in \Ga(V', M^{-1})$ and $l \in \Ga(V', L^{0})$. 
The equality $\be(m) = \d_N(n)$ becomes the matrix equality 
\[ \bmat{m \\[0.1em] 0} =
\bmat{\d_M & \al \\[0.1em] 0 & -\d_{L}} \cd 
\bmat{m' \\[0.1em] l} = 
\bmat{ \d_M(m') + \al(l) \\[0.1em] -\d_{L}(l) } . \]
We see that $l \in \Ga(V', \opn{Z}^{0}(L))$. 
Then the cohomology class 
$\bar{l} := \pi_L(l) \in \Ga(V', \opn{H}^{0}(L))$
satisfies 
\[ \bar{\al}(\bar{l}) = \pi_M(\al(l)) = \pi_M(\d_M(m') + \al(l)) = 
\pi_M(m) = \bar{m} \]
in $\Ga(V', \opn{H}^{0}(M))$, as required. 
\end{proof} 

\begin{dfn} \label{dfn:1475}
A morphism $\phi$ in  $\dcat{K}(A, \cat{M})$ is called a 
{\em quasi-isomorphism}
\index{Quasi-isomorphism}
if the morphisms $\opn{H}^i(\phi)$ in $\cat{M}$ are 
isomorphisms for all $i$. 
The set of quasi-isomorphisms in $\dcat{K}(A, \cat{M})$ is denoted by 
$\dcat{S}(A, \cat{M})$. 
\end{dfn}

By Proposition \ref{prop:1475} the functor $\opn{H}^0$ is cohomological.
According to Proposition \ref{prop:106}, $\dcat{S}(A, \cat{M})$ is a 
denominator set of cohomological origin; so Theorem \ref{thm:105} applies to 
it, 
and the next definition makes sense. 

\begin{dfn}
Let $\cat{M}$ be a $\K$-linear abelian category and $A$ a central DG $\K$-ring.
The {\em derived category of DG $A$-modules in $\cat{M}$}%
\index{Derived category! of DG $A$-modules in $\cat{M}$}%
\index{1-D(A,M)@$\dcat{D}(A, \cat{M})$}%
\index{1-D(A)@$\dcat{D}(A)$}%
\index{1-D(M)@$\dcat{D}(\cat{M})$}
is the $\K$-linear triangulated category 
\[ \dcat{D}(A, \cat{M}) := \dcat{K}(A, \cat{M})_{\dcat{S}(A, \cat{M})} . \]
The corresponding triangulated localization functor is 
\[ \opn{Q} : \dcat{K}(A, \cat{M}) \to \dcat{D}(A, \cat{M}) . \]
\end{dfn}

We also have the additive functor
$\opn{P} : \dcat{C}_{\mrm{str}}(A, \cat{M}) \to \dcat{K}(A, \cat{M})$
which sends a strict morphism of DG modules to its homotopy class.

\begin{dfn} \label{dfn:1480}
Let $\cat{M}$ be an abelian category and let $A$ be a DG ring.
Define the functor 
\[ \til{\opn{Q}} := \opn{Q} \circ \opn{P} : 
\dcat{C}_{\mrm{str}}(A, \cat{M}) \to \dcat{D}(A, \cat{M}) . \]
\end{dfn}

We thus get a commutative diagram of additive functors
\begin{equation} \label{eqn:5025}
\UseTips \xymatrix @C=8ex @R=6ex {
\dcat{C}_{\mrm{str}}(A, \cat{M}) 
\ar[r]^{\opn{P}}
\ar@(u,u)[rr]^{\til{\opn{Q}}}
&
\dcat{K}(A, \cat{M})
\ar[r]^{\opn{Q}}
&
\dcat{D}(A, \cat{M})
}
\end{equation}
that are all the identity on objects.

It is sometimes convenient to describe morphisms in $\dcat{D}(A, \cat{M})$
in terms of the functor $\til{\opn{Q}}$. A morphism 
$s \in \dcat{C}_{\mrm{str}}(A, \cat{M})$
is called a quasi-isomorphism if $\opn{P}(s)$ is a quasi-isomorphism
in $\dcat{K}(A, \cat{M})$; i.e.\ if $\opn{H}(s)$ is an isomorphism in
$\dcat{G}_{\mrm{str}}(\cat{M})$. 

\begin{prop} \label{prop:1481}
\mbox{}
\begin{enumerate}
\item Every morphism $\phi$ in $\dcat{D}(A, \cat{M})$ can be written as a right 
fraction 
$\phi = \til{\opn{Q}}(a) \circ \til{\opn{Q}}(s)^{-1}$
where $a, s \in \dcat{C}_{\mrm{str}}(A, \cat{M})$ and 
$s$ is a quasi-isomorphism. 

\item Let $a \in \dcat{C}_{\mrm{str}}(A, \cat{M})$. 
Then $\til{\opn{Q}}(a) = 0$ in $\dcat{D}(A, \cat{M})$ iff 
there exists a quasi-isomorphism $s$ in $\dcat{C}_{\mrm{str}}(A, \cat{M})$
such that $a \circ s$ is a coboundary in $\dcat{C}(A, \cat{M})$.
\end{enumerate}
\end{prop}

\begin{proof} \mbox{}

\smallskip \noindent
(1) This is because of property (RO3) of Definition \ref{dfn:1407} and the fact 
that $\opn{P}$ is full. 

\medskip \noindent 
(2) Let us write $\bar{a} := \opn{P}(a) \in \dcat{K}(A, \cat{M})$.
Since $\opn{Q}(\bar{a}) = 0$, by Property (RO4) of Definition \ref{dfn:1407}
there is a quasi-isomorphism 
$\bar{s} \in \dcat{K}(A, \cat{M})$
such that $\bar{a} \circ \bar{s} = 0$ in $\dcat{K}(A, \cat{M})$.
Choose a quasi-isomorphism $s \in \dcat{C}_{\mrm{str}}(A, \cat{M})$
such that $\bar{s} = \opn{P}(s)$. Then $\opn{P}(a \circ s) = 0$, and 
this means that $a \circ s$ is a coboundary. 
\end{proof}

Of course there is a left version of this proposition.
The next definition is about arbitrary (not necessarily linear) categories. 

\begin{dfn} \label{4980}
A functor $F : \cat{C} \to \cat{D}$ between categories is called  
{\em conservative} if for every morphism $\phi$ in $\cat{C}$, 
$\phi$ is an isomorphism if and only if $F(\phi)$ is an isomorphism. 
\end{dfn}

By Proposition \ref{prop:1410} the cohomology functor (\ref{eqn:4985}) extends 
uniquely to a functor 
\begin{equation} \label{eqn:4986}
\opn{H} : \dcat{D}(A, \cat{M}) \to \dcat{G}_{\mrm{str}}(\cat{M}) . 
\end{equation}

\begin{cor} \label{cor:2145}
The functor $\opn{H}$ in \tup{(\ref{eqn:4986})} is conservative.
\end{cor}

\begin{proof}
One implication is trivial. For the other implication, assume that 
$\phi$ is a morphism in $\dcat{D}(A, \cat{M})$ such that 
$\opn{H}(\phi)$ is an isomorphism. We can write $\phi$ as a right fraction: 
$\phi = \opn{Q}(a) \circ \opn{Q}(s)^{-1}$, 
where $a \in \dcat{K}(A, \cat{M})$ and 
$s \in \dcat{S}(A, \cat{M})$.
Then 
$\opn{H}(\phi) = \opn{H}(\opn{Q}(a)) \circ \opn{H}(\opn{Q}(s))^{-1}$,
and we see that $\opn{H}(\opn{Q}(a))$ is an isomorphism.
But $\opn{H}(a) = \opn{H}(\opn{Q}(a))$, so in fact 
$a \in \dcat{S}(A, \cat{M})$ too. Therefore 
$\opn{Q}(a)$ is an isomorphism in 
$\dcat{D}(A, \cat{M})$. It follows that $\phi$ is an isomorphism in 
$\dcat{D}(A, \cat{M})$.
\end{proof}

\begin{exer} \label{exer:1846}
Here $\cat{M} = \dcat{M}(\K)$, so 
$\dcat{K}(A, \cat{M}) = \dcat{K}(A)$. Show that the functor 
$\opn{H}^0 : \dcat{K}(A) \to \dcat{M}(\K)$ is corepresentable 
by the object $A \in \dcat{K}(A)$ (see Subsection 
\ref{subsec:repfunc}). 
\end{exer}

\begin{rem} \label{rem:4971}
This is a continuation of Remarks \ref{rem:4960} and \ref{rem:4970} on strictly 
unital 
$\mrm{A}_{\infty}$ categories%
\index{Ainfty@$\mrm{A}_{\infty}$ category}.
The base ring $\K$ is a field. 
The references are \cite{Kel5} and \cite{COS}.
Recall that if $\cat{A}$ is an  
$\mrm{A}_{\infty}$ category and $\cat{B}$ is a DG category, then 
$\cat{A_{\infty}Fun}(\cat{A}, \cat{B})$
is a DG category. 

Consider the special case is when $\cat{B} = \dcat{C}(\K)$, the DG category of 
complexes of $\K$-modules. An $\mrm{A}_{\infty}$ functor 
$M : \cat{A} \to \dcat{C}(\K)$ is called an {\em $\mrm{A}_{\infty}$ 
$\cat{A}$-module}. 
The $\mrm{A}_{\infty}$ $\cat{A}$-modules are the objects of the DG category 
$\dcat{C}_{\infty}(\cat{A}) := \cat{A_{\infty}Fun}(\cat{A}, \dcat{C}(\K))$. 
The strict subcategory $\opn{Str}(\dcat{C}_{\infty}(\cat{A}))$, in the sense of
Definition \ref{dfn:1212}(1), is what most texts call the {\em category 
of $\mrm{A}_{\infty}$ $\cat{A}$-modules}. 

A morphism $\phi : M \to N$ in $\opn{Str}(\dcat{C}_{\infty}(\cat{A}))$
is called a {\em quasi-isomorphism} if for every object 
$x \in \opn{Ob}(\cat{A})$ the morphism 
$\phi(x) : M(x) \to N(x)$ in 
$\opn{Str}(\dcat{C}(\K)) = \dcat{C}_{\mrm{str}}(\K)$ is a 
quasi-isomorphism. The localization of 
$\opn{Str}(\dcat{C}_{\infty}(\cat{A}))$ w.r.t.\ 
the quasi-isomorphisms is the {\em derived category of $\mrm{A}_{\infty}$
$\cat{A}$-modules}, denoted by $\dcat{D}_{\infty}(\cat{A})$.
(As usual, the tradition is to call $\dcat{D}_{\infty}(\cat{A})$ the ``homotopy 
category'' of $\dcat{C}_{\infty}(\cat{A})$, and to denote 
it by $\opn{Ho}(\dcat{C}_{\infty}(\cat{A}))$.)
There is a localization functor 
$\opn{Q} : \opn{Str}(\dcat{C}_{\infty}(\cat{A})) \to 
\dcat{D}_{\infty}(\cat{A})$.

We can also consider the category 
$\dcat{K}_{\infty}(\cat{A}) := \opn{H}^0(\dcat{C}_{\infty}(\cat{A}))$, see 
Definition \ref{dfn:1212}(2). There is a full functor 
$\opn{P} : \opn{Str}(\dcat{C}_{\infty}(\cat{A})) \to 
\dcat{K}_{\infty}(\cat{A})$,
and a commutative diagram 
\begin{equation} \label{eqn:4978}
\UseTips \xymatrix @C=10ex @R=6ex {
\opn{Str}(\dcat{C}_{\infty}(\cat{A}))
\ar[r]^{\opn{P}}
\ar@(u,u)[rr]^{\opn{Q}}
&
\dcat{K}_{\infty}(\cat{A})
\ar[r]^{\bar{\opn{Q}}}
&
\dcat{D}_{\infty}(\cat{A})
}
\end{equation}
All three functors are the identity on objects. 
The categories $\dcat{K}_{\infty}(\cat{A})$ and 
$\dcat{D}_{\infty}(\cat{A})$ are triangulated, and the functor $\bar{\opn{Q}}$ 
is an isomorphism of triangulated categories. 

Finally, let $\cat{A}$ be a DG category. Then the DG category 
$\dcat{C}_{}(\cat{A})$ of DG $\cat{A}$-modules embeds faithfully (but not 
fully) into the DG category $\dcat{C}_{\infty}(\cat{A})$.
The induced triangulated functor 
$\dcat{D}_{}(\cat{A}) \to \dcat{D}_{\infty}(\cat{A})$
is an equivalence. 
\end{rem}

\mysubsection{Boundedness Conditions in
\texorpdfstring{$\dcat{K}(A, \cat{M})$}{K(A,M)}}
\label{subsec:bd-in-homot}

We continue with a central DG $\K$-ring $A$ and a $\K$-linear abelian category 
$\cat{M}$. 
Boundedness conditions in $\dcat{C}(A, \cat{M})$ were introduced in Definition 
\ref{dfn:3220}. Similarly: 

\begin{dfn} \label{dfn:3221}
\index{1-K(A,M)@$\dcat{K}^{\star}(A, \cat{M})$}
\index{Boundedness! condition}
We define $\dcat{K}^-(A, \cat{M})$, $\dcat{K}^+(A, \cat{M})$ and 
$\dcat{K}^{\mrm{b}}(A, \cat{M})$ to be the full subcategories of 
$\dcat{K}(A, \cat{M})$ consisting of bounded above, bounded below and bounded
DG modules respectively. 
\end{dfn}

Of course 
$\dcat{K}^{\mrm{b}}(A, \cat{M}) = \dcat{K}^-(A, \cat{M})
\cap \dcat{K}^+(A, \cat{M})$.
The subcategories $\dcat{K}^{\star}(A, \cat{M})$, 
for $\star \in \{ -, +, \mrm{b} \}$, 
are full triangulated subcategories of $\dcat{K}(A, \cat{M})$; this is 
because the operations of translation 
and cone preserve the various boundedness conditions. 
Note that 
$\dcat{K}^{\star}(A, \cat{M}) = 
\opn{Ho} \bigl( \dcat{C}^{\star}(A, \cat{M}) \bigr)$.
As the next example shows, sometimes the categories 
$\dcat{C}^{\star}(A, \cat{M})$ and 
$\dcat{K}^{\star}(A, \cat{M})$ can be very degenerate. 

\begin{exa} \label{exa:1460}
Let $A$ be the DG ring $\K[t, t^{-1}]$, the ring of Laurent polynomials in the  
variable $t$ of degree $1$, with the zero differential. If 
$M = \{ M^i \}_{i \in \Z}$ is a nonzero object of 
$\dcat{C}(A, \cat{M})$, then $M^i \neq 0$ for all $i$. Therefore 
the categories $\dcat{C}^{\star}(A, \cat{M})$ and 
$\dcat{K}^{\star}(A, \cat{M})$ are zero for 
$\star \in \{ -, +, \mrm{b} \}$. 
\end{exa}

Let 
$\dcat{S}^{\star}(A, \cat{M}) := \dcat{K}^{\star}(A, \cat{M}) \cap
\dcat{S}(A, \cat{M})$,
the category of quasi-isomorphisms in $\dcat{K}^{\star}(A, \cat{M})$. 
As already mentioned, Theorem \ref{thm:105} applies here, so we can 
localize. 

\begin{dfn} \label{dfn:1460}
For $\star \in \{ -, +, \mrm{b} \}$ we define
\[  \dcat{D}^{\star}(A, \cat{M}) := 
\dcat{K}^{\star}(A, \cat{M})_{\dcat{S}^{\star}(A, \cat{M})} , \]
\index{1-D(A,M)@$\dcat{D}^{\star}(A, \cat{M})$}
\index{Boundedness! condition}
the Ore localization of $\dcat{K}^{\star}(A, \cat{M})$ with respect to 
$\dcat{S}^{\star}(A, \cat{M})$.
\end{dfn}

Here is another kind of boundedness condition. 

\begin{dfn} \label{dfn:1461}
For $\star \in \{ -, +, \mrm{b} \}$ we define
$\dcat{D}(A, \cat{M})^{\star}$
\index{1-D(A,M)@$\dcat{D}(A, \cat{M})^{\star}$}
\index{Boundedness! condition}
to be the full subcategory of $\dcat{D}(A, \cat{M})^{}$
on the complexes $M$ whose cohomology $\opn{H}(M)$ is of boundedness type 
$\star$. 
\end{dfn}

Of course $\dcat{D}(A, \cat{M})^{\star}$ is a full triangulated subcategory 
of $\dcat{D}(A, \cat{M})$. 
See Remark \ref{rem:1461} regarding Definitions \ref{dfn:1460} and 
\ref{dfn:1461}. 

Handling boundedness conditions requires the use of {\em truncations}. These 
will be introduced after the next proposition. 

\begin{prop} \label{prop:2165}
Let 
$0 \to L \xar{\phi} M \xar{\psi} N \to 0$ 
be a short exact sequence in $\dcat{C}_{\mrm{str}}(A, \cat{M})$. 
Then there is a morphism $\th : N \to \opn{T}(L)$ in $\dcat{D}(A, \cat{M})$ 
such that 
$L \xar{\opn{Q}(\phi)} M \xar{\opn{Q}(\psi)} N \xar{\th} \opn{T}(L)$
is a distinguished triangle in $\dcat{D}(A, \cat{M})$.
\end{prop}

\begin{proof}
We are following the proof of \cite[Proposition 1.7.5]{KaSc1}.
Let $\til{N}$ be the standard cone on $\phi$. In matrix notation, as 
in Definition \ref{dfn:1172}, we have
$\til{N} = \sbmat{M \\[0.2em] \opn{T}(L)}$
and 
$\d_{\til{N}} =
\sbmat{\d_M & \, \phi \, \circ \, \opn{t}^{-1} 
\\[0.2em] 0 & \, \d_{\opn{T}(L)}}$.
The object $\til{N}$ sits inside the standard triangle 
$L \xar{\phi} M \xar{\til{\psi}} \til{N} \xar{\til{\chi}} \opn{T}(L)$
in  $\dcat{C}_{\mrm{str}}(A, \cat{M})$, where 
$\til{\psi} := \sbmat{\opn{id} \\[0.2em] 0}$ 
and 
$\til{\chi} := \sbmat{0 & \, \opn{id}}$
in matrix notation. Define the morphism $\ga = \til{N} \to N$ to be the matrix
$\ga := \sbmat{\psi & \, 0}$. We get a commutative diagram
\[ \UseTips \xymatrix @C=8ex @R=5.5ex {
L
\ar[r]^{\phi}
&
M
\ar[r]^{\til{\psi}}
\ar[dr]_{\psi} 
&
\til{N}
\ar[r]^{\til{\chi}} 
\ar[d]^{\ga}
&
\opn{T}(L)
\\
&
&
N
} \]
in  $\dcat{C}_{\mrm{str}}(A, \cat{M})$.
We shall prove below that $\ga$ is a quasi-isomorphism. 
Then the morphism 
$\th := \opn{Q}(\til{\chi}) \circ \opn{Q}(\ga)^{-1} : N \to \opn{T}(L)$ 
will work. 

Let $\til{K}$ be the standard cone on $\opn{id}_L$, and let 
$\til{\be} : \til{K} \to \til{N}$ be the matrix morphism 
$\sbmat{\phi & \, 0 \\[0.1em] 0 & \, \opn{id}} : 
\sbmat{L \\[0.2em] \opn{T}(L)} \to
\sbmat{M \\[0.2em] \opn{T}(L)}$.
This fits into a short exact sequence 
$0 \to \til{K} \xar{\til{\be}} \til{N} \xar{\ga} N \to 0$
in $\dcat{C}_{\mrm{str}}(A, \cat{M})$.
But the DG module $\til{K}$ is acyclic, and therefore $\ga$ is a 
quasi-isomorphism, as claimed. 
\end{proof}

\begin{dfn} \label{dfn:2320}
Let $M \in \dcat{C}(\cat{M})$ and $i \in \Z$. The 
{\em smart truncation of $M$ below $i$}%
\index{Truncation of a complex! smart}%
\index{1-Smtleq)@$\opn{smt}^{\leq i}(M)$}
is the complex 
\[ \opn{smt}^{\leq i}(M) := 
\bigl( \cdots \to M^{i-2} \xar{\d} M^{i-1} \xar{\d} \opn{Z}^i(M) \to 0 \to
\cdots \bigr) .  \]
The {\em smart truncation of $M$ above $i$}%
\index{1-Smtgeq)@$\opn{smt}^{\geq i}(M)$}
is the complex 
\[ \opn{smt}^{\geq i}(M) := 
\bigl( \cdots \to 0 \to \opn{Y}^{i}(M) \xar{\d} M^{i+1} \xar{\d} M^{i+2}
\to \cdots \bigr) . \]
\end{dfn}

In the definition we use the objects of cocycles $\opn{Z}^i(M)$
and decocycles $\opn{Y}^i(M)$ from Definition \ref{dfn:2993}. Note that 
$\opn{smt}^{\leq i}(M)$ is a subcomplex of $M$, whereas 
$\opn{smt}^{\geq i}(M)$ a quotient complex of $M$. Denoting the inclusion and 
the projection by $e$ and $p$ respectively, we get a diagram 
\begin{equation} \label{eqn:3035}
\opn{smt}^{\leq i}(M) \xar{e} M \xar{p} \opn{smt}^{\geq i + 1}(M) 
\end{equation}
in $\dcat{C}_{\mrm{str}}(\cat{M})$. 
This diagram is functorial in $M \in \dcat{C}_{\mrm{str}}(\cat{M})$.
But diagram (\ref{eqn:3035}) is not an exact sequence in general  -- there is a 
defect in degree $i$. 

Recall that a DG ring $A$ is called nonpositive if $A^i = 0$ for all $i > 0$. 

\begin{prop} \label{prop:3170}
Assume $A$ is a nonpositive DG ring. 
\begin{enumerate}
\item The differential of every 
$M \in \dcat{C}_{\mrm{str}}(A, \cat{M})$ is $A^0$-linear.

\item The smart truncations from Definition \tup{\ref{dfn:2320}}
are functors from $\dcat{C}_{\mrm{str}}(A, \cat{M})$ to itself. 
\end{enumerate}
\end{prop}

\begin{exer} \label{exer:3170}
Prove this proposition. 
\end{exer}

\begin{prop} \label{prop:2320}
Assume $A$ is nonpositive. 
For every $M \in \dcat{C}(A, \cat{M})$ there is a distinguished triangle
\[ \opn{smt}^{\leq i}(M) \xar{e} M \xar{p} \opn{smt}^{\geq i + 1}(M) 
\xar{\th} \opn{T}(\opn{smt}^{\leq i}(M)) \]
in $\dcat{D}(A, \cat{M})$. Also, 
$\opn{H}^j(e) :  \opn{H}^j(\opn{smt}^{\leq i}(M))  \to \opn{H}^j(M)$
is an isomorphism in $\cat{M}$ for all $j \leq i$, and 
$\opn{H}^j(p) : \opn{H}^j(M) \to \opn{H}^j(\opn{smt}^{\geq i + 1}(M))$
is an isomorphism in $\cat{M}$ for all $j \geq i + 1$.
\end{prop}

\begin{proof}
The claims about $\opn{H}^j(e)$ and $\opn{H}^j(p)$ are trivial to verify. 
Now there is a short exact sequence 
\begin{equation} \label{eqn:3036}
0 \to \opn{smt}^{\leq i}(M) \xar{e} M \xar{p'} N \to 0 
\end{equation}
in $\dcat{C}_{\mrm{str}}(A, \cat{M})$, where 
\[ N := \bigl( \cdots \to 0 \to M^i / \opn{Z}^{i}(M) \xar{\d} M^{i+1} 
\xar{\d} M^{i+2} \to \cdots \bigr) . \]
According to Proposition \ref{prop:2165} we get a distinguished triangle 
\[ \opn{smt}^{\leq i}(M) \xar{e} M \xar{p'} N 
\xar{\th'} \opn{T}(\opn{smt}^{\leq i}(M)) \]
in $\dcat{D}(A, \cat{M})$. Next, there is an obvious quasi-isomorphism 
$\phi : N \to \opn{smt}^{\geq i + 1}(M)$
in $\dcat{C}_{\mrm{str}}(A, \cat{M})$
such that $p = \phi \circ p'$. We define the morphism 
\[ \th := \th' \circ \opn{Q}(\phi)^{-1} : \opn{smt}^{\geq i + 1}(M)
\to \opn{T}(\opn{smt}^{\leq i}(M)) \]
in $\dcat{D}(A, \cat{M})$. 
\end{proof}

\begin{prop} \label{prop:1463}
Assume $A$ is nonpositive.
For $\star \in \{ -, +, \mrm{b} \}$ the canonical functor
$\dcat{D}^{\star}(A, \cat{M}) \to \dcat{D}(A, \cat{M})^{\star}$ 
is an equivalence of triangulated categories. 
\end{prop}

\begin{proof} It is done in four steps. 

\smallskip \noindent
Step 1. Here we prove that the functor 
$F^- : \dcat{D}^-(A, \cat{M}) \to \dcat{D}(A, \cat{M})$ 
is fully faithful.
Let $s : M \to L$ be a quasi-isomorphism in 
$\dcat{K}(A, \cat{M})$ with 
$L \in \dcat{K}^-(A, \cat{M})$.
Say $L$ is concentrated in degrees $\leq i$.
Then $\opn{H}^j(M) = \opn{H}^j(L) = 0$ for all $j > i$. 
The smart truncation $\opn{smt}^{\leq i}(M)$ belongs to 
$\dcat{K}^-(A, \cat{M})$, and the inclusion
$t : \opn{smt}^{\leq i}(M) \to M$ is a quasi-isomorphism. 
According to Proposition \ref{prop:1320}, with 
$\cat{K} = \dcat{K}(A, \cat{M})$ and 
$\cat{K}' = \dcat{K}^-(A, \cat{M})$, and with condition (r), we see that 
$F^-$ is fully faithful.

\medskip \noindent 
Step 2. Here we prove that the functor 
$F^+ : \dcat{D}^+(\cat{M}) \to \dcat{D}(\cat{M})$ 
is fully faithful.
Let $s : L \to M$ be a quasi-isomorphism in 
$\dcat{K}(A, \cat{M})$ with $L \in \dcat{K}^+(A, \cat{M})$.
Say $L$ is concentrated in degrees $\geq i$. 
Then $\opn{H}^j(M) = \opn{H}^j(L) = 0$ for all $j < i$. 
The smart truncation 
$\opn{smt}^{\geq i}(M)$ belongs to $\dcat{K}^+(A, \cat{M})$, 
and the projection 
$t : M \to \opn{smt}^{\geq i}(M)$ is a quasi-isomorphism. 
According to Proposition \ref{prop:1320}, with condition (l), we see that 
$F^+$ is fully faithful.

\medskip \noindent 
Step 3. The arguments in step 1 we show that 
$\dcat{D}^{\mrm{b}}(A, \cat{M}) \to \dcat{D}^{+}(A, \cat{M})$
is fully faithful. And by step 2, 
$\dcat{D}^{+}(A, \cat{M}) \to \dcat{D}(A, \cat{M})$ is fully faithful. 
Therefore 
$\dcat{D}^{\mrm{b}}(A, \cat{M}) \to \dcat{D}(A, \cat{M})$
is fully faithful.

\medskip \noindent 
Step 4. Smart truncation shows that the functor
$\dcat{D}^{\star}(A, \cat{M}) \to \dcat{D}(A, \cat{M})^{\star}$ 
is essentially surjective on objects. 
\end{proof}

\begin{rem} \label{rem:1461}
\index{Boundedness! condition}
\index{1-D(M)@$\dcat{D}^{\star}(\cat{M})$}
\index{1-D(M)@$\dcat{D}(\cat{M})^{\star}$}
Definition \ref{dfn:1460} with $A = \K$, namely $\dcat{D}^{\star}(\cat{M})$,
is standard usage of this notation -- starting from \cite[Section I.4]{RD}, all 
the way to the recent book \cite[Definition 13.1.2]{KaSc2}. An exception is 
\cite[Definition tag=05RU]{SP}. 

On the other hand, definition \ref{dfn:1461} with $A = \K$, namely 
$\dcat{D}(\cat{M})^{\star}$, 
is almost never used (an exception is \cite{PSY3}). 

We are going to require both distinct concepts,  
$\dcat{D}^{\star}(A, \cat{M})$ and $\dcat{D}(A, \cat{M})^{\star}$,
until the end of Section \ref{sec:exist-resol}. Beginning with 
Section \ref{sec:adj-equ-cohdim}, in Definition \ref{dfn:2122}, the notation 
$\dcat{D}(A, \cat{M})^{\star}$ will be eliminated, and 
$\dcat{D}^{\star}(A, \cat{M})$ will be redefined to replace it. 
\end{rem}

\newpage
\mysubsection{Thick Subcategories of \texorpdfstring{$\cat{M}$}{M}}
\label{subsec:thick-subcats}

We begin by recalling a definition about abelian categories. 

\begin{dfn} \label{dfn:2323}
Let $\cat{M}$ be an abelian category. A 
{\em thick abelian subcategory}%
\index{Abelian category! thick abelian subcategory of}
of $\cat{M}$ is a full abelian subcategory $\cat{N}$ that is closed under
extensions. Namely if 
\[ 0 \to M' \to M \to M'' \to 0 \]
is a short exact sequence in $\cat{M}$ with $M', M'' \in \cat{N}$, then 
$M \in \cat{N}$ too. 
\end{dfn}

\begin{prop} \label{prop:3171} 
Let $\cat{M}$ be an abelian category, and let $\cat{M}' \sub \cat{M}$ be a 
thick abelian subcategory. Suppose 
$M_1 \to M_2 \to N \to M_3 \to M_4$ 
is an exact sequence in $\cat{M}$, and the objects $M_i$ belong $\cat{M}'$. 
Then 
$N \in \cat{M}'$ too. 
\end{prop}

\begin{exer} \label{exer:2881}
Prove this proposition.
\end{exer}

\begin{dfn} \label{dfn:2321}
Let $\cat{M}$ be an abelian category and $\cat{N} \sub \cat{M}$ a thick abelian 
subcategory. We denote by $\dcat{D}_{\cat{N}}(\cat{M})$ 
\index{1-DN(M)@$\dcat{D}_{\cat{N}}(\cat{M})$}
the full subcategory of $\dcat{D}(\cat{M})$
consisting of complexes $M$ such that 
$\opn{H}^i(M) \in \cat{N}$ for every $i$. 

Given a boundedness condition $\star$, 
\index{Boundedness! condition}
\index{1-DNstar(M)@$\dcat{D}^{\star}_{\cat{N}}(\cat{M})$}
we write 
$\dcat{D}_{\cat{N}}^{\star}(\cat{M}) := 
\dcat{D}_{\cat{N}}(\cat{M}) \cap \dcat{D}^{\star}(\cat{M})$
and
$\dcat{D}_{\cat{N}}(\cat{M})^{\star} := 
\dcat{D}_{\cat{N}}(\cat{M}) \cap \dcat{D}(\cat{M})^{\star}$.
\end{dfn}

\begin{prop} \label{prop:3052}
If $\cat{N}$ is a thick abelian subcategory of $\cat{M}$ then
$\dcat{D}_{\cat{N}}(\cat{M})$ is a full triangulated subcategory of 
$\dcat{D}(\cat{M})$.
\end{prop}

\begin{proof}
Clearly $\dcat{D}_{\cat{N}}(\cat{M})$ is closed under translations. 
Now suppose 
$M' \to M \to M'' \xar{\, \triangle\, }$ \, 
is a distinguished triangle in $\dcat{D}(\cat{M})$ such that 
$M', M \in \dcat{D}_{\cat{N}}(\cat{M})$; we have to show that $M''$
is also in $\dcat{D}_{\cat{N}}(\cat{M})$.
Consider the  exact sequence 
$\opn{H}^i(M') \to \opn{H}^i(M) \to
\opn{H}^i(M'') \to \opn{H}^{i+1}(M') \to \opn{H}^{i+1}(M)$.
The four outer objects belong to $\cat{N}$. 
Since $\cat{N}$ is a thick abelian subcategory of $\cat{M}$ it follows (using 
Proposition \ref{prop:3171}) that 
$\opn{H}^i(M'') \in \cat{N}$.
\end{proof}

\begin{exa}
Let $A$ be a noetherian commutative ring. 
The category $\cat{Mod}_{\mrm{f}} A$ of finitely generated
modules is a thick abelian subcategory of $\cat{Mod} A$.
\end{exa}

\begin{exa}
Consider $\cat{Mod} \Z = \cat{Ab}$. 
As above we have the thick abelian subcategory 
$\cat{Ab}_{\mrm{f}} = \cat{Mod}_{\mrm{f}} \Z$ of finitely generated
abelian groups.
There is also the thick abelian subcategory 
$\cat{Ab}_{\mrm{tors}}$ of torsion abelian groups (every element has a finite
order). The intersection of $\cat{Ab}_{\mrm{tors}}$ and 
$\cat{Ab}_{\mrm{f}}$ is the category $\cat{Ab}_{\mrm{fin}}$ of finite abelian
groups. This is also thick. 
\end{exa}

\begin{exa}
Let $X$ be a noetherian scheme (e.g.\ an algebraic variety over an 
algebraically closed field). Consider the abelian category 
$\cat{Mod} \mcal{O}_X$ of $\mcal{O}_X$-modules. In it there is the thick abelian
subcategory $\cat{QCoh} \mcal{O}_X$ of quasi-coherent sheaves, and in 
that there is the thick abelian
subcategory $\cat{Coh} \mcal{O}_X$ of coherent sheaves.

Passing to complexes, let us write 
$\dcat{C}(X) := \dcat{C}(\cat{Mod} \mcal{O}_X)$,
the DG category of unbounded complexes. Its strict subcategory is 
$\dcat{C}_{\cat{str}}(X)$, and for each boundedness condition $\star$ we have 
the full subcategory 
$\dcat{C}^{\star}_{\cat{str}}(X) \sub \dcat{C}_{\cat{str}}(X)$. 
For the derived category we write 
$\dcat{D}(X) := \dcat{D}(\cat{Mod} \mcal{O}_X)$.
Inside it there are the full triangulated subcategories
\[ \dcat{D}_{\mrm{c}}(X) := 
\dcat{D}_{\cat{Coh} \mcal{O}_X}(\cat{Mod} \mcal{O}_X) \sub
\dcat{D}_{\mrm{qc}}(X) := 
\dcat{D}_{\cat{QCoh} \mcal{O}_X}(\cat{Mod} \mcal{O}_X) . \]

One can show (see \cite[Corollary II.7.19]{RD}) that the canonical functor \lb
$\dcat{D}^{+}(\cat{QCoh} \mcal{O}_X) \to \dcat{D}^{+}_{\mrm{qc}}(X)$
is an equivalence. The proof is based on the deep fact 
(see \cite[Theorem II.7.18]{RD}) 
that $\cat{QCoh} \mcal{O}_X$ has enough injectives relative to 
$\cat{Mod} \mcal{O}_X$, in the sense of Definition \ref{dfn:2880}. 
Then Theorem \ref{thm:2880} says that every complex 
$\MM \in \dcat{D}^{+}_{\mrm{qc}}(X)$ admits a quasi-isomorphism 
$\MM \to \II$ in $\dcat{C}_{\cat{str}}(X)$, where $\II$ is a bounded below 
complex of injective quasi-coherent $\OO_X$-modules. 

Then, using smart truncations and the fact that every 
quasi-coherent sheaf on $X$ is the direct limit of its coherent subsheaves, one 
can show that the canonical functor 
$\dcat{D}^{\mrm{b}}(\cat{Coh} \mcal{O}_X) \to \dcat{D}^{\mrm{b}}_{\mrm{c}}(X)$
is an equivalence. See \cite[Proposition 3.5]{Huy}; when $X = 
\opn{Spec}(A)$ is affine see also Proposition \ref{prop:200} below. The 
category 
$\dcat{D}^{\mrm{b}}_{\mrm{c}}(X)$
is the focus of a significant part of the research in modern algebraic geometry 
(e.g.\ birational geometry, see \cite{Huy}). 

A totally different approach yields this result: if $X$ is separated and 
quasi-compact, then the canonical functor
$\dcat{D}^{}(\cat{QCoh} \mcal{O}_X) \to \dcat{D}^{}_{\mrm{qc}}(X)$
is an equivalence. See \cite[Corollary 5.5]{BoNe}.
\end{exa}

For a left noetherian ring $A$ we write 
$\dcat{D}_{\mrm{f}}(\cat{Mod} A) :=
\dcat{D}_{\cat{Mod}_{\mrm{f}} A}(\cat{Mod} A)$.

\begin{prop} \label{prop:200}
Let $A$ be a left noetherian ring and $\star \in \{ -, \mrm{b} \}$. Then the
canonical functor 
$\dcat{D}^{\star}(\cat{Mod}_{\mrm{f}} A) \to 
\dcat{D}_{\mrm{f}}(\cat{Mod} A)^{\star}$
is an equivalence of triangulated categories.
\end{prop}

\begin{proof}
Consider the functor 
$F : \dcat{D}^{-}(\cat{Mod}_{\mrm{f}} A) \to \dcat{D}(\cat{Mod} A)$.
Suppose $s : M \to L$ is a quasi-isomorphism in $\dcat{K}(\cat{Mod} A)$,
such that $L \in \dcat{K}^{-}(\cat{Mod}_{\mrm{f}} A)$. 
Then $M \in \dcat{D}_{\mrm{f}}(\cat{Mod} A)^{-}$.
A bit later (in Theorem \ref{thm:3340}) we will prove that $M$ 
admits a free resolution $P \to M$, where $P$ is a
bounded above complex of finitely generated free modules. 
Thus we get a quasi-isomorphism $t : P \to M$ with 
$P \in \dcat{K}^{-}(\cat{Mod}_{\mrm{f}} A)$.  By Proposition \ref{prop:1320} 
with condition (r) we conclude that $F$ is fully faithful. This also shows 
that the essential image of $F$ is 
$\dcat{D}_{\mrm{f}}(\cat{Mod} A)^{-}$.
 
Next consider the functor 
$G : \dcat{D}^{\mrm{b}}(\cat{Mod}_{\mrm{f}} A) \to 
\dcat{D}^{-}(\cat{Mod}_{\mrm{f}} A)$.
Suppose $s : L \to M$ is a quasi-isomorphism in 
$\dcat{K}^-(\cat{Mod}_{\mrm{f}} A)$
with $L \in \dcat{K}^{\mrm{b}}(\cat{Mod}_{\mrm{f}} A)$. 
Say $\opn{H}(L)$ is concentrated in the integer interval $[d_0, d_1]$. 
Then 
$t : M \to \opn{smt}^{\geq d_0}(M)$
is a quasi-isomorphism, and 
$\opn{smt}^{\geq d_0}(M) \in \dcat{K}^{\mrm{b}}(\cat{Mod}_{\mrm{f}} A)$. 
By Proposition \ref{prop:1320} with condition (l) we conclude that $G$ is 
fully faithful. Therefore the composition 
$F \circ G : \dcat{D}^{\mrm{b}}(\cat{Mod}_{\mrm{f}} A) \to 
\dcat{D}(\cat{Mod} A)$
is fully faithful. Suitable truncations ($\opn{smt}^{\geq d_0}$ and 
$\opn{smt}^{\leq d_1}$) show that the essential image of $F \circ G$ is 
$\dcat{D}_{\mrm{f}}(\cat{Mod} A)^{\mrm{b}}$.
\end{proof}

\mysubsection{The Embedding of \texorpdfstring{$\cat{M}$}{M} in 
\texorpdfstring{$\dcat{D}(\cat{M})$}{D(M)}}

Here again we only consider an abelian category $\cat{M}$. 

For $M, N \in \cat{M}$ there is no difference between the $\K$-modules
$\opn{Hom}_{\cat{M}}(M, N)$, $\opn{Hom}_{\dcat{C}(\cat{M})}(M, N)$ and 
$\opn{Hom}_{\dcat{K}(\cat{M})}(M, N)$. Thus the canonical functors 
$\cat{M} \to \dcat{C}(\cat{M})$ and 
$\cat{M} \to \dcat{K}(\cat{M})$ are fully faithful. 
The same is true for $\dcat{D}(\cat{M})$, but this requires a proof.

Let $\dcat{D}(\cat{M})^0$ be the full subcategory of 
$\dcat{D}(\cat{M})$ consisting of complexes whose cohomology is 
concentrated in degree $0$. This is an additive subcategory of 
$\dcat{D}(\cat{M})$.

\begin{prop} \label{prop:121}
The canonical functor $\cat{M} \to \dcat{D}(\cat{M})^0$ 
is an equivalence.
\end{prop}

\begin{proof}
Let's denote the canonical functor $\cat{M} \to \dcat{D}(\cat{M})^0$ 
by $F$. Under the fully faithful embedding 
$\cat{M} \subseteq \dcat{C}_{\mrm{str}}(\cat{M})$,
$F$ is just the restriction of the functor $\til{\opn{Q}}$ from Definition 
\ref{dfn:1480}.

The functor $\opn{H}^0 : \dcat{D}(\cat{M}) \to \cat{M}$
satisfies $\opn{H}^0 \circ \, F = \opn{Id}_{\cat{M}}$. This implies 
that $F$ is faithful.

Next we prove that $F$ is full. Take any objects $M, N \in \cat{M}$ and a 
morphism $q : M \to N$ in $\dcat{D}(\cat{M})$. 
By Proposition \ref{prop:1481} we know that  
$q = \til{\opn{Q}}(a) \circ \til{\opn{Q}}(s)^{-1}$ for some
morphisms $a : L \to N$ and $s : L \to M$ in 
$\dcat{C}_{\mrm{str}}(\cat{M})$, with $s$
a quasi-isomorphism. Let $L' := \opn{smt}^{\leq 0}(L)$, the smart truncation of 
$L$. Since $L \in \dcat{D}(\cat{M})^0$, the inclusion 
$u : L' \to L$ is a quasi-isomorphism  in $\dcat{C}_{\mrm{str}}(\cat{M})$. 
Writing $a' := a \circ u$ and $s' := s \circ u$, we see that $s'$ is a 
quasi-isomorphism, and
$q = \til{\opn{Q}}(a') \circ \til{\opn{Q}}(s')^{-1}$. 

Next let $L'' := \opn{smt}^{\geq 0}(L')$, the other smart truncation. 
The projection $v : L' \to L''$ is a surjective quasi-isomorphism in 
$\dcat{C}_{\mrm{str}}(\cat{M})$. 
Because $L''$ is a complex concentrated in degree $0$, we can view it as 
an object of $\cat{M}$. The morphisms $a'$ and $s'$ factor 
as $a' = a'' \circ v$ and $s' = s'' \circ v$, where 
$a'' : L'' \to N$ and $s'' : L'' \to M$ are morphisms 
in $\cat{M}$. But $s''$ is a quasi-isomorphism in 
$\dcat{C}_{\mrm{str}}(\cat{M})$, and so it is actually an isomorphism in 
$\cat{M}$. Therefore we have a morphism 
$a'' \circ (s'')^{-1} : M \to N$
in $\cat{M}$, and 
\[ \til{\opn{Q}}(a'' \circ (s'')^{-1}) = 
\til{\opn{Q}}(a'') \circ \til{\opn{Q}}(s'')^{-1} = 
\til{\opn{Q}}(a') \circ \til{\opn{Q}}(s')^{-1} =  q  . \]

Finally we have to prove that every $L \in \dcat{D}(\cat{M})^0$ is isomorphic, 
in $\dcat{D}(\cat{M})$, to a complex $L''$ that's concentrated in degree $0$. 
But we already showed it in the previous paragraphs.
\end{proof}

\begin{prop} \label{prop:1440}
Let $\cat{M}$ be an abelian category and let 
$0 \to L \xar{\phi}  M \xar{\psi} N \to 0$ 
be a diagram in $\cat{M}$. The following conditions are equivalent\tup{:}
\begin{enumerate}
\rmitem{i} The diagram is an exact sequence. 
\rmitem{ii}  There is a distinguished triangle 
$L \xar{ \, \til{\opn{Q}}(\phi) \, }  M \xar{ \, \til{\opn{Q}}(\psi) \, } N 
\xar{ \, \th \, } \opn{T}(L)$ 
in $\dcat{D}(\cat{M})$.
\end{enumerate}
\end{prop}

\begin{exer} \label{exer:1917}
Prove Proposition \ref{prop:1440}. 
\end{exer}

The last two propositions say that the abelian category $\cat{M}$, with its 
kernels and cokernels, can be recovered from the triangulated category 
$\dcat{D}(\cat{M})$. 

\begin{rem} \label{rem:2320}
Assume that the diagram in Proposition \ref{prop:1440} is an exact sequence.
By Proposition \ref{prop:3605}, this sequence is split if and only if the 
morphism $\th$ is zero. Furthermore, if $\cat{M}$ has enough injectives, then 
there is a canonical bijection between the set of 
isomorphism classes of extensions of $N$ by $L$, and the set 
$\opn{Hom}_{\dcat{D}(\cat{M})}(N, \opn{T}(L))$. See \cite[Section I.6]{RD}, or 
use Exercise \ref{exer:2125}, which says that 
$\opn{Hom}_{\dcat{D}(\cat{M})}(N, \opn{T}(L)) \cong 
\opn{Ext}^1_{\cat{M}}(N, L)$.
\end{rem}

\mysubsection{The Opposite Derived Category is Triangulated}
\label{subsec:opp-dercat-triang}

Here we deal with a DG $\K$-ring $A$ and a $\K$-linear 
abelian category $\cat{M}$. 
In Subsection \ref{subsec:opp-hom-triang} we put a canonical triangulated 
structure on the opposite homotopy category 
$\dcat{K}(A, \cat{M})^{\mrm{op}}$. This structure is such that the flip functor 
\[ \opn{Flip} : \dcat{K}(A, \cat{M})^{\mrm{op}} \to
\dcat{K}(A^{\mrm{op}}, \cat{M}^{\mrm{op}}) =  
\dcat{K}(A, \cat{M})^{\mrm{flip}} \]
is an isomorphism of triangulated categories. 
Here, unlike in Subsection \ref{subsec:opp-hom-triang},
we are omitting the overline decoration from $\opn{Flip}$. 
In the current subsection we push this triangulated structure to 
the opposite derived category $\dcat{D}(A, \cat{M})^{\mrm{op}}$.

But first let us present two types of full subcategories
$\cat{K} \sub \dcat{K}(A, \cat{M})$ that are triangulated, and also 
$\cat{K}^{\mrm{op}} \sub \dcat{K}(A, \cat{M})^{\mrm{op}}$ 
is triangulated. 

For any boundedness condition $\star$ (see Definition \ref{dfn:3221}), we know 
that $\dcat{K}^{\star}(A, \cat{M})$
is a full triangulated subcategory of $\dcat{K}^{}(A, \cat{M})$. 
The notation 
$\dcat{K}^{\star}(A, \cat{M})^{\mrm{op}}$
refers to the opposite category of 
$\dcat{K}^{\star}(A, \cat{M})$;
thus, $\dcat{K}^{\star}(A, \cat{M})^{\mrm{op}}$ is the full subcategory of 
$\dcat{K}^{}(A, \cat{M})^{\mrm{op}}$
on the DG modules satisfying the boundedness condition $\star$. 

\begin{prop} \label{prop:3050}
For every boundedness condition $\star$, the subcategory 
$\dcat{K}^{\star}(A, \cat{M})^{\mrm{op}}$
is triangulated in $\dcat{K}^{}(A, \cat{M})^{\mrm{op}}$. 
\end{prop}

\begin{proof}
Let's write $\cat{K} := \dcat{K}^{\star}(A, \cat{M})$ and
$\cat{K}^{\mrm{op}} := \dcat{K}^{\star}(A, \cat{M})^{\mrm{op}}$.
Define  
$\cat{K}^{\mrm{flip}} \lb := \opn{Flip}(\cat{K}^{\mrm{op}}) \sub 
\dcat{K}(A^{\mrm{op}}, \cat{M}^{\mrm{op}})$.
According to Theorem \ref{thm:2495}(2) we know that
$\cat{K}^{\mrm{flip}} = \dcat{K}^{-\star}(A^{\mrm{op}}, \cat{M}^{\mrm{op}})$,
where $-\star$ is the reversed boundedness condition. 
Thus $\cat{K}^{\mrm{flip}}$ is a full triangulated subcategory of
$\dcat{K}(A^{\mrm{op}}, \cat{M}^{\mrm{op}})$. 
By definition of the triangulated structure on 
$\dcat{K}^{}(A, \cat{M})^{\mrm{op}}$
we conclude that $\cat{K}^{\mrm{op}}$ is triangulated.
\end{proof}

Let $\cat{N} \sub \cat{M}$ be a thick abelian subcategory. Then 
$\dcat{K}_{\cat{N}}(\cat{M})$, see Definition \ref{dfn:2321},  
is a full triangulated subcategory of $\dcat{K}(\cat{M})$.
This is just like Proposition \ref{prop:3052}.
The opposite category 
$\dcat{K}_{\cat{N}}(\cat{M})^{\mrm{op}}$,
on the same set of objects, is a full additive subcategory of 
$\dcat{K}(\cat{M})^{\mrm{op}}$. 

\begin{prop} \label{prop:3051}
Let $\cat{N} \sub \cat{M}$ be a thick abelian subcategory.
Then 
$\dcat{K}_{\cat{N}}(\cat{M})^{\mrm{op}}$
is a full triangulated subcategory of $\dcat{K}(\cat{M})^{\mrm{op}}$. 
\end{prop}

\begin{proof}
Let's write 
$\cat{K} := \dcat{K}_{\cat{N}}(\cat{M})$. Then 
$\cat{K}^{\mrm{flip}} := \opn{Flip}(\cat{K}^{\mrm{op}}) \sub 
\dcat{K}(\cat{M}^{\mrm{op}})$
is, by Theorem \ref{thm:2495}(4), the category  
$\dcat{K}_{\cat{N}^{\mrm{op}}}(\cat{M}^{\mrm{op}})$. 
Thus $\cat{K}^{\mrm{flip}}$ is a triangulated subcategory. 
By definition of the triangulated structure on 
$\dcat{K}^{}(\cat{M})^{\mrm{op}}$
we conclude that $\cat{K}^{\mrm{op}}$ is triangulated.
\end{proof}

Later, in Section \ref{sec:resol},  we will see that several {\em resolving 
subcategories} are also full triangulated subcategories of 
$\dcat{K}(A, \cat{M})^{\mrm{op}}$.
These are relevant to Theorems \ref{thm:1470} and \ref{thm:1500} below. 

Taking intersections, we get more full triangulated subcategories of  
$\dcat{K}^{}(A, \cat{M})^{\mrm{op}}$.

Recall that 
$\dcat{S}(A, \cat{M})$ is the set of quasi-isomorphisms in 
$\dcat{K}(A, \cat{M})$. Let $\dcat{S}(A, \cat{M})^{\mrm{op}}$ be the 
set of quasi-isomorphisms in $\dcat{K}(A, \cat{M})^{\mrm{op}}$.
Note that being a quasi-isomorphism has nothing to do with the 
triangulated structure. Since a morphism $\psi$ in 
$\dcat{K}(A, \cat{M})$ is a quasi-isomorphism iff 
$\opn{Op}(\psi)$ is a quasi-isomorphism in 
$\dcat{K}(A, \cat{M})^{\mrm{op}}$, we have a bijection
$\opn{Op} : \dcat{S}(A, \cat{M}) \iso \dcat{S}(A, \cat{M})^{\mrm{op}}$.
We know that $\dcat{S}(A, \cat{M})$ is a left and right 
denominator set in 
$\dcat{K}(A, \cat{M})$, by Proposition \ref{prop:106}. 
Therefore, by Proposition \ref{prop:1845},
$\dcat{S}(A, \cat{M})^{\mrm{op}}$ is left and right denominator set in 
$\dcat{K}(A, \cat{M})^{\mrm{op}}$.
We need to know more: 

\begin{prop} \label{prop:2356}
Let $\cat{K}^{\mrm{op}}$ be a full triangulated subcategory of 
$\dcat{K}(A, \cat{M})^{\mrm{op}}$, and define
$\cat{S}^{\mrm{op}} := \cat{K}^{\mrm{op}} \cap \, 
\dcat{S}(A, \cat{M})^{\mrm{op}}$,
the set of quasi-isomorphisms in $\cat{K}^{\mrm{op}}$. 
Then\tup{:}
\begin{enumerate}
\item $\cat{S}^{\mrm{op}}$ is a denominator set of cohomological origin in 
$\cat{K}^{\mrm{op}}$.

\item The localized category 
$\cat{D}^{\mrm{op}} := (\cat{K}^{\mrm{op}})_{\cat{S}^{\mrm{op}}}$
has a unique triangulated structure s.t.\ the localization functor 
$\opn{Q}^{\mrm{op}} : \cat{K}^{\mrm{op}} \to \cat{D}^{\mrm{op}}$
is a triangulated functor. 
\end{enumerate}
\end{prop}

\begin{proof}
On the flip side, i.e.\ in 
$\dcat{K}(A^{\mrm{op}}, \cat{M}^{\mrm{op}})$,
we know that the set of quasi-iso\-morphisms
$\dcat{S}(A^{\mrm{op}}, \cat{M}^{\mrm{op}})$
is a denominator set of cohomological origin. But
by Theorem \ref{thm:2495}(4) we have 
$\dcat{S}(A, \cat{M})^{\mrm{op}} = 
\opn{Flip}^{-1} \bigl( \dcat{S}(A^{\mrm{op}}, \cat{M}^{\mrm{op}}) \bigr)$,
so this too is a denominator set of cohomological origin. 
Now we can use Theorem \ref{thm:105}. 
\end{proof}

Observe that 
$\cat{D}^{\mrm{op}} = (\cat{K}_{\cat{S}})^{\mrm{op}}$,
by Proposition \ref{prop:3030}. 
The situation is summarized by the following commutative diagram of functors:
\begin{equation} \label{eqn:3060}
\UseTips \xymatrix @C=10ex @R=6ex {
\cat{C}_{\mrm{str}}^{\mrm{op}}
\ar@{->>}[r]^{\opn{P}^{\mrm{op}}}
\ar[d]_{\opn{Flip}}^{\cong}
&
\cat{K}^{\mrm{op}}
\ar[r]^{\opn{Q}^{\mrm{op}}}
\ar[d]_{\opn{Flip}}^{\cong}
&
\cat{D}^{\mrm{op}}
\ar[d]_{\opn{Flip}}^{\cong}
\\
\cat{C}_{\mrm{str}}^{\mrm{flip}}
\ar@{->>}[r]^{\opn{P}^{\mrm{flip}}}
&
\cat{K}^{\mrm{flip}}
\ar[r]^{\opn{Q}^{\mrm{flip}}}
&
\cat{D}^{\mrm{flip}} 
} 
\end{equation}
Here $\cat{C}$ is the full subcategory of 
$\dcat{C}(A, \cat{M})$ on the objects of $\cat{K}$,
and 
$\cat{C}^{\mrm{flip}} := \opn{Flip}(\cat{C}) \sub 
\dcat{C}(A, \cat{M})^{\mrm{flip}} = 
\dcat{C}(A^{\mrm{op}}, \cat{M}^{\mrm{op}})$.
All these functors are bijective on objects. The vertical ones are bijective on 
morphisms too. The functors marked $\surj$ are surjective on morphisms (i.e.\ 
they are full). The first vertical arrow is an isomorphism of abelian 
categories, and the other two vertical arrows are isomorphisms of triangulated 
categories.

Warning: Given $\cat{K}^{\mrm{op}}$ like in Proposition \ref{prop:2356}, 
there is no reason for the canonical triangulated functor 
$\cat{D}^{\mrm{op}} \to \dcat{D}(A, \cat{M})^{\mrm{op}}$
to be fully faithful; see Proposition \ref{prop:1320} for sufficient 
conditions.

\begin{dfn} \label{dfn:3225}
Let $\cat{K} \sub \dcat{K}(A, \cat{M})$ be a full additive subcategory s.t.\ 
$\cat{K}^{\mrm{op}}$ is a  triangulated subcategory of 
$\dcat{K}(A, \cat{M})^{\mrm{op}}$, and let 
$\cat{S} := \cat{K} \cap \, \dcat{S}(A, \cat{M})$.
The category 
$\cat{D}^{\mrm{op}} := (\cat{K}^{\mrm{op}})_{\cat{S}^{\mrm{op}}}$
is given the triangulated structure from Proposition \ref{prop:2356}. 
\index{Derived category! opposite}
\end{dfn}

For $\cat{K} = \dcat{K}(A, \cat{M})$ we get a triangulated structure on 
$\cat{D}^{\mrm{op}} = \dcat{D}(A, \cat{M})^{\mrm{op}}$.

\begin{dfn} \label{dfn:3226} 
Let $\cat{K} \sub \dcat{K}(A, \cat{M})$ be a full additive subcategory s.t.\ 
$\cat{K}^{\mrm{op}}$ is a triangulated subcategory of 
$\dcat{K}(A, \cat{M})^{\mrm{op}}$.
Define
$\cat{D} := \cat{K}_{\cat{S}}$
and 
$\cat{D}^{\mrm{op}} := (\cat{K}^{\mrm{op}})_{\cat{S}^{\mrm{op}}}$, 
where $\cat{S} := \cat{K} \cap \, \dcat{S}(A, \cat{M})$.
Let $\cat{E}$ be some triangulated category. A 
{\em contravariant triangulated functor}%
\index{Triangulated functor! contravariant}
from $\cat{D}$ to $\cat{E}$ is, by definition, a 
triangulate functor 
$F : \cat{D}^{\mrm{op}}  \to \cat{E}$,
where $\cat{D}^{\mrm{op}}$ has the canonical triangulated structure from 
Definition \ref{dfn:3225}.
\end{dfn}

\cleardoublepage
\mysection{Derived Functors} \label{sec:der-funcs}

\AYcopyright

Suppose 
$F : \dcat{K}(A, \cat{M}) \to \cat{E}$
is a triangulated functor. Here $A$ is a DG ring, $\cat{M}$ is an abelian 
category, and $\cat{E}$ is a triangulated category. 
In this section we define the right and left derived functors
$\mrm{R} F, \, \mrm{L} F : \dcat{D}(A, \cat{M}) \to \cat{E}$ 
of $F$. These are also triangulated functors, satisfying certain universal 
properties. We shall prove the uniqueness of the derived functors, and their 
existence under suitable assumptions.
 
The universal properties of the derived functors are best stated in 
$2$-categorical language. This will be explained is the first two subsections. 
In Subsection \ref{subsec:abstr-der-funs} we define derived functors in the 
abstract setting (as opposed to the triangulated setup), and prove the main 
results for them. These results will then be specialized to various settings:
triangulated functors (Subsection \ref{subsec:tri-der-funs}), contravariant 
triangulated functors (Subsection \ref{subsec:cntr-tri-der-funs}), and 
triangulated bifunctors (Subsection \ref{subsec:tri-bifun}).

\mysubsection{\texorpdfstring{$2$}{2}-Categorical Notation} 
\label{subsec:2-cat-not}

In this section we are going to do a lot of work with morphisms of functors 
(i.e.\ natural transformations). The language and notation of ordinary category 
theory that we used so far is not adequate for this purpose. Therefore we will 
now introduce notation from the theory of {\em $2$-categories}.
(We will not give a definition of a $2$-category here; but it is basically the 
structure of $\cat{Cat}$ that is mentioned below.) 
For more details on $2$-categories the reader 
can look at \cite{Mac2} or \cite[Section 1]{Ye7}.

Consider the set $\cat{Cat}$ of all categories. The set-theoretical aspects 
are neglected, as explained in Subsection \ref{subsec:set-theor}.
(Briefly, the precise solution is this: $\cat{Cat}$ is the set of all 
$\cat{U}$-categories; so $\cat{Cat}$ is a subset of a bigger Grothendieck 
universe, say $\cat{V}$, and it is a $\cat{V}$-category.)

The set $\cat{Cat}$ is the set of objects of a $2$-category.
\index{Two@$2$-category} 
This means that in $\cat{Cat}$ there are two kinds of morphisms: 
{\em $1$-morphisms} between objects, and {\em $2$-morphisms} between 
$1$-morphisms. There are several kinds of compositions, and these have several 
properties. All this will be explained below. 

Suppose $\cat{C}_0$, $\cat{C}_1, \ldots$ are categories, namely 
objects of $\cat{Cat}$. The $1$-morph\-isms between them are the functors.
The notation is as usual: $F : \cat{C}_0 \to \cat{C}_1$ denotes a functor.

Suppose $F, G : \cat{C}_0 \to \cat{C}_1$ are functors. The $2$-morphisms from 
$F$ to $G$ are the morphisms of functors 
(i.e.\ the natural transformations), and the notation is 
$\eta : F \twoto G$. The double arrow is the distinguishing notation for  
$2$-morphisms. When specializing to an object
$M \in \cat{C}_0$ we revert to the single arrow notation, namely 
$\eta_M : F(M) \to G(M)$ is the corresponding morphism in $\cat{C}_1$. 
The diagram depicting this is
\[ \UseTips  \xymatrix @C=12ex @R=6ex {
\cat{C}_0
\ar@(ur,ul)[r]^(0.5){F} _(0.5){}="eta-s"
\ar@(dr,dl)[r]_(0.5){G}  ^(0.5){}="eta-t"
\ar@{=>}  "eta-s";"eta-t" _(0.5){\eta}
&
\cat{C}_1
} \]
We shall refer to such a diagram as a {\em $2$-diagram}. 
\index{Two@$2$-diagram} 

Each object (category) $\cat{C}$ has its identity $1$-morphism  (functor)
$\opn{Id}_{\cat{C}} : \cat{C} \to \cat{C}$.  
Each $1$-morphism $F$ has its identity $2$-morphism (natural transformation)
$\opn{id}_F : F \twoto F$.  

Now we consider compositions. For functors there is nothing new: given functors 
$F_1 : \cat{C}_{0} \to \cat{C}_{1}$ and 
$F_2 : \cat{C}_{1} \to \cat{C}_{2}$, their composition, that we now call 
{\em horizontal composition}, is the functor
$F_2 \circ F_1 : \cat{C}_{0} \to \cat{C}_{2}$. The diagram is 
\[ \UseTips  \xymatrix @C=12ex @R=6ex {
\cat{C}_0
\ar[r]^(0.5){F_1} 
\ar@(dr,dl)[rr]_(0.5){F_2 \, \circ \, F_1}  
&
\cat{C}_1
\ar[r]^(0.5){F_2} 
&
\cat{C}_2
} \]
This can be viewed as a commutative $1$-diagram, or as a shorthand for the 
$2$-diagram
\[ \UseTips  \xymatrix @C=12ex @R=6ex {
\cat{C}_0
\ar[r]^(0.5){F_1}
\ar @/_3em/ [rr]_(0.5){F_2 \, \circ \, F_1} _(0.5){}="t"
&
\cat{C}_1
\ar[r]^(0.5){F_2} 
\ar@{=>}  ; "t" _(0.5){\opn{id}}
&
\cat{C}_2
} \]
in which $\opn{id}$ is the identity $2$-morphism of $F_2 \circ F_1$. 

The complication begins with compositions of $2$-morphisms. 
Suppose we are given $1$-morphisms
$F_i, G_i : \cat{C}_{i - 1} \to \cat{C}_{i}$
and $2$-morphisms $\eta_i : F_i \twoto G_i$. In a diagram: 
\[ \UseTips  \xymatrix @C=12ex @R=6ex {
\cat{C}_0
\ar@(ur,ul)[r]^(0.5){F_1} _(0.5){}="eta1-s"
\ar@(dr,dl)[r]_(0.5){G_1}  ^(0.5){}="eta1-t"
\ar@{=>}  "eta1-s";"eta1-t" _(0.5){\eta_1}
&
\cat{C}_1
\ar@(ur,ul)[r]^(0.5){F_2} _(0.5){}="eta2-s"
\ar@(dr,dl)[r]_(0.5){G_2}  ^(0.5){}="eta2-t"
\ar@{=>}  "eta2-s";"eta2-t" _(0.5){\eta_2}
&
\cat{C}_2
} \]
The {\em horizontal composition} is the morphism of functors 
$\eta_2 \circ \eta_1 : F_2 \circ F_1  \twoto G_2 \circ G_1$.
The diagram is 
\[ \UseTips  \xymatrix @C=18ex @R=6ex {
\cat{C}_0
\ar@(ur,ul)[r]^(0.5){F_2 \, \circ \, F_1} _(0.5){}="eta-s"
\ar@(dr,dl)[r]_(0.5){G_2 \, \circ \, G_1}  ^(0.5){}="eta-t"
\ar@{=>}  "eta-s";"eta-t" _(0.5){\eta_2 \, \circ \, \eta_1}
&
\cat{C}_2
} \] 

\begin{exer} \label{exer:1314}
For an object $M \in \cat{C}_0$, give an explicit formula for the morphism 
\[ (\eta_2 \circ \eta_1)_M : (F_2 \circ F_1)(M) \to (G_2 \circ G_1)(M) \]
in the category $\cat{C}_2$.
\end{exer}

Suppose we are given $1$-morphisms 
$E, F, G : \cat{C}_{0} \to \cat{C}_{1}$, and $2$-morphisms \lb 
$\ze : E \twoto F$ and $\eta : F \twoto G$. 
The diagram depicting this is
\[ \UseTips  \xymatrix @C=16ex @R=6ex {
\cat{C}_0
\ar@(u,u)[r]^(0.5){E} _(0.5){}="zeta-s"
\ar[r]_(0.7){F}  ^(0.5){}="zeta-t" _(0.5){}="eta-s"
\ar@(d,d)[r]_(0.5){G}  ^(0.5){}="eta-t"
\ar@{=>}  "zeta-s";"zeta-t" _(0.5){\zeta}
\ar@{=>}  "eta-s";"eta-t" _(0.5){\eta}
&
\cat{C}_1
} \] 
The {\em vertical composition} of $\zeta$ and $\eta$ is the $2$-morphism  
$\eta * \zeta : E \to G$.
Notice the new symbol for this operation. 
The corresponding diagram is 
\[ \UseTips  \xymatrix @C=18ex @R=6ex {
\cat{C}_0
\ar@(ur,ul)[r]^(0.5){E} _(0.5){}="eta-s"
\ar@(dr,dl)[r]_(0.5){G}  ^(0.5){}="eta-t"
\ar@{=>}  "eta-s";"eta-t" _(0.5){\eta \, * \, \ze}
&
\cat{C}_1
} \] 

\begin{exer} \label{exer:1311}
For an object $M \in \cat{C}_0$, give an explicit formula for the morphism
$(\eta * \zeta)_M : E(M) \to G(M)$
in the category $\cat{C}_{1}$.
\end{exer}

Something intricate occurs in the situation shown in the next diagram. 
\[ \UseTips  \xymatrix @C=16ex @R=6ex {
\cat{C}_0
\ar@(u,u)[r]^(0.5){E_1} _(0.5){}="zeta1-s"
\ar[r]_(0.7){F_1}  ^(0.5){}="zeta1-t" _(0.5){}="eta1-s"
\ar@(d,d)[r]_(0.5){G_1}  ^(0.5){}="eta1-t"
\ar@{=>}  "zeta1-s";"zeta1-t" _(0.5){\zeta_1}
\ar@{=>}  "eta1-s";"eta1-t" _(0.5){\eta_1}
&
\cat{C}_1
\ar@(u,u)[r]^(0.5){E_2} _(0.5){}="zeta2-s"
\ar[r]_(0.7){F_2}  ^(0.5){}="zeta2-t" _(0.5){}="eta2-s"
\ar@(d,d)[r]_(0.5){G_2}  ^(0.5){}="eta2-t"
\ar@{=>}  "zeta2-s";"zeta2-t" _(0.5){\zeta_2}
\ar@{=>}  "eta2-s";"eta2-t" _(0.5){\eta_2}
&
\cat{C}_2
} \] 
It turns out that 
\[ (\eta_2 * \zeta_2) \circ (\eta_1 * \zeta_1) = 
(\eta_2 \circ \eta_1) * (\zeta_2 \circ \zeta_1) \]
as morphisms $E_2 \circ E_1 \twoto G_2 \circ G_1$. 
This is called the {\em exchange property}. 

\begin{exer} \label{exer:1312}
Prove the exchange property.
\end{exer}

Just like abstract categories, we can talk about 
triangulated categories. There is the $2$-category $\cat{TrCat}$ 
of all ($\K$-linear) triangulated categories. The objects here are 
the triangulated categories $(\cat{K}, \opn{T})$; the $1$-morphisms
are the triangulated functors $(F, \tau)$; and the $2$-morphisms
are the morphisms of triangulated functors $\eta$. This is what we are going 
to use later in this section.

\mysubsection{Functor Categories} \label{subsec:fun-cats}

In this subsection we isolate a part of $2$-category theory. This simplifies 
the discussion greatly. All set theoretical issues (sizes of sets) are 
neglected. As before, this can be treated by introducing a bigger universe. 

\begin{dfn} \label{dfn:2025}
Given categories $\cat{C}$ and $\cat{D}$, let 
$\cat{Fun}(\cat{C}, \cat{D})$
\index{1-Fun@$\cat{Fun}(\cat{C}, \cat{D})$}
be the category whose objects are the functors $F : \cat{C} \to \cat{D}$.
Given objects $F, G \in \cat{Fun}(\cat{C}, \cat{D})$,
the morphisms $\eta : F \twoto G$ in 
$\cat{Fun}(\cat{C}, \cat{D})$ are the morphisms of functors, i.e.\ the  natural 
transformations. 
\end{dfn}

For the case $\cat{D} = \cat{Set}$ this concept was already mentioned in 
Subsection \ref{subsec:repfunc}. 

In terms of the $2$-category $\cat{Cat}$ from the previous subsection,
$\cat{C}$ and $\cat{D}$ are objects of $\cat{Cat}$; the objects of 
$\cat{Fun}(\cat{C}, \cat{D})$ are $1$-morphisms in 
$\cat{Cat}$; and the morphisms  of 
$\cat{Fun}(\cat{C}, \cat{D})$ are $2$-morphisms in 
$\cat{Cat}$. Thus we made a ``reduction of order'', from $2$ to $1$, by 
passing from $\cat{Cat}$ to $\cat{Fun}(\cat{C}, \cat{D})$.

Suppose $G : \cat{C}' \to \cat{C}$ and 
$H : \cat{D} \to \cat{D}'$ are functors. There is an induced functor
\begin{equation} \label{eqn:2025}
\cat{Fun}(G, H) : \cat{Fun}(\cat{C}, \cat{D}) \to 
\cat{Fun}(\cat{C}', \cat{D}')
\end{equation}
defined by 
$\cat{F}(G, H)(F) := H \circ F \circ G$. 

\begin{prop} \label{prop:2025}
If $G$ and $H$ are equivalences, then the functor $\cat{Fun}(G, H)$ in 
\tup{(\ref{eqn:2025})} is an equivalence.
\end{prop}

\begin{exer} \label{exer:2035}
Prove Proposition \ref{prop:2025}.
\end{exer}

Recall that for a category $\cat{C}$ and a multiplicatively closed set of 
morphisms $\cat{S} \sub \cat{C}$ we denote by $\cat{C}_{\cat{S}}$ the 
localization. It comes with the localization functor 
$\opn{Q} : \cat{C} \to \cat{C}_{\cat{S}}$,
which is the identity on objects. See Definition \ref{dfn:1406}. 

For a category $\cat{E}$ let $\cat{E}^{\times} \sub \cat{E}$
be the category of isomorphisms; it has all the objects, but its  
morphisms are just the isomorphisms in $\cat{E}$. 

\begin{dfn} \label{dfn:2037}
Given categories $\cat{C}$ and $\cat{E}$, a multiplicatively closed set of 
morphisms $\cat{S} \sub \cat{C}$, and a functor $F : \cat{C} \to \cat{E}$,
we say that {\em $F$ is localizable to $\cat{S}$} if 
$F(\cat{S}) \sub  \cat{E}^{\times}$. 
We denote by
$\cat{Fun}_{\cat{S}}(\cat{C}, \cat{E})$%
\index{1-Fun@$\cat{Fun}_{\cat{S}}(\cat{C}, \cat{D})$}
the full subcategory of $\cat{Fun}(\cat{C}, \cat{E})$
on the functors that are localizable to $\cat{S}$.
\end{dfn}

Here is a useful formulation of the universal property of localization
of categories (see Definition \ref{dfn:1406}). 
Recall that a functor is an isomorphism of categories iff it is an equivalence 
that is bijective on sets of objects. 

\begin{prop} \label{prop:2026} 
Let $\cat{C}$ and $\cat{E}$ be categories, and let $\cat{S} \sub \cat{C}$ be a 
multiplicatively closed set of morphisms. Then the functor 
\[ \cat{Fun}(\opn{Q}, \opn{Id}_{\cat{E}}) : 
\cat{Fun}(\cat{C}_{\cat{S}}, \cat{E}) \to 
\cat{Fun}_{\cat{S}}(\cat{C}, \cat{E}) \]
is an isomorphism of categories.
\end{prop}

\begin{exer} \label{exer:2036}
Prove Proposition \ref{prop:2026}.
\end{exer}

By definition a bifunctor 
$F : \cat{C} \times \cat{D} \to \cat{E}$
is a functor from the product category $\cat{C} \times \cat{D}$.
See Subsection \ref{subsec:bifunc}.
It will be useful to retain both meanings; so we shall write 
\begin{equation} \label{eqn:2035}
\cat{BiFun}(\cat{C} \times \cat{D}, \cat{E}) := 
\cat{Fun}(\cat{C} \times \cat{D}, \cat{E}) ,
\end{equation}
where in the first expression we recall that $\cat{C} \times \cat{D}$ is a 
product. 

The next proposition describes bifunctors in a non-symmetric fashion. 

\begin{prop} \label{prop:2027} 
Let $\cat{C}$, $\cat{D}$ and $\cat{E}$ be categories. There is an isomorphism 
of categories 
\[ \Xi : \cat{Fun}(\cat{C} \times \cat{D}, \cat{E}) \to 
\cat{Fun}(\cat{C}, \cat{Fun}(\cat{D}, \cat{E})) \]
with the following formula\tup{:}
for a functor $F : \cat{C} \times \cat{D} \to \cat{E}$,
the functor
\[ \Xi(F) : \cat{C} \to \cat{Fun}(\cat{D}, \cat{E}) \]
is $\Xi(F)(C) := F(C, -)$.
\end{prop}

\begin{exer} \label{exer:2037}
Prove Proposition \ref{prop:2027}.
\end{exer}

\begin{prop}  \label{prop:2028} 
Let $\cat{C}$ and $\cat{D}$ be categories, and let 
$\cat{S} \sub \cat{C}$ and $\cat{T} \sub \cat{D}$
be multiplicatively closed sets of morphisms. Then
$\cat{S} \times \cat{T}$ is a multiplicatively closed set of morphisms
in $\cat{C} \times \cat{D}$, and the canonical functor 
\[ \Theta:  (\cat{C} \times \cat{D})_{\cat{S} \times \cat{T}} \to 
\cat{C}_{\cat{S}} \times \cat{D}_{\cat{T}} \]
is an isomorphism of categories. 
\end{prop}

\begin{proof}
It is clear that
$\cat{S} \times \cat{T} \sub \cat{C} \times \cat{D}$  
is closed under compositions. Let's spell out how the functor 
$\Theta$ arises. The functor
\[ \opn{Q}_{\cat{S} \times \cat{T}} : \cat{C} \times \cat{D} \to 
(\cat{C} \times \cat{D})_{\cat{S} \times \cat{T}} \]
is universal for functors 
$F : \cat{C} \times \cat{D} \to \cat{E}$
that invert $\cat{S} \times \cat{T}$; namely functors $F$ such that 
$F(s, t) \in \cat{E}^{\times}$ for all $s \in \cat{S}$ and $t \in \cat{T}$.
Since the functor 
\[ \opn{Q}_{\cat{S}} \times \opn{Q}_{\cat{T}} : \cat{C} \times \cat{D} \to 
\cat{C}_{\cat{S}} \times \cat{D}_{\cat{T}} \]
inverts $\cat{S} \times \cat{T}$, we get the functor $\Theta$. We will prove 
that the functor $\opn{Q}_{\cat{S}} \times \opn{Q}_{\cat{T}}$ 
has the same universal property as $\opn{Q}_{\cat{S} \times \cat{T}}$; this 
will imply that $\Theta$ is an isomorphism. 

Consider an arbitrary category $\cat{E}$. 
Invoking Propositions \ref{prop:2027} and \ref{prop:2026} we get a commutative 
diagram 

\begin{equation} \label{eqn:3383}
\UseTips  \xymatrix @C=5ex @R=8ex {
\cat{Fun}(\cat{C}_{\cat{S}} \times \cat{D}_{\cat{T}}, \cat{E})
\ar[d]^{\cat{Fun}
(\opn{Q}_{\cat{S}} \times \opn{Q}_{\cat{T}}, \opn{Id}_{\cat{E}})}
\ar[r]^(0.44){\cong}
&
\cat{Fun}(\cat{C}_{\cat{S}}, \cat{Fun}(\cat{D}_{\cat{T}}, \cat{E})) 
\ar[r]^{\cong}
&
\cat{Fun}_{\cat{S}}(\cat{C}, \cat{Fun}_{\cat{T}}(\cat{D}, \cat{E}))
\ar[d]^{\tup{f.f.\ emb.}}
\\
\cat{Fun}_{}(\cat{C} \times \cat{D}, \cat{E})
\ar[rr]^{\Xi}_{\cong}
&
&
\cat{Fun}_{}(\cat{C}, \cat{Fun}_{}(\cat{D}, \cat{E}))
}
\end{equation}

\noindent in which the right vertical arrow is a fully faithful embedding. 
Take any functor $F$ in the bottom left corner of (\ref{eqn:3383}), 
and let $F' := \Xi(F)$, which lives in the bottom right corner. 
Now $F'$ belongs to the top right corner of (\ref{eqn:3383})
iff $F'(s)(t) \in \cat{E}^{\times}$ for all $s \in \cat{S}$ and $t \in \cat{T}$.
But $F'(s)(t) = F(s, t)$, so this happens iff $F$ inverts 
$\cat{S} \times \cat{T}$.
\end{proof}

Denominator sets were introduced in Definition \ref{dfn:21}.

\begin{prop} \label{prop:2100}
In the situation of Proposition \tup{\ref{prop:2028}}, the following two 
conditions are equivalent\tup{:}
\begin{enumerate}
\rmitem{i} The multiplicatively closed sets  
$\cat{S} \sub \cat{C}$ and $\cat{T} \sub \cat{D}$ are left \tup{(}resp.\ 
right\tup{)} denominator sets. 

\rmitem{ii} The multiplicatively closed set
$\cat{S} \times \cat{T} \sub \cat{C} \times \cat{D}$ is a left 
\tup{(}resp.\ right\tup{)} denominator set. 
\end{enumerate}
\end{prop}

\begin{exer} \label{exer:2100}
Prove Proposition \ref{prop:2100}.
\end{exer}

\mysubsection{Abstract Derived Functors} \label{subsec:abstr-der-funs}

Here we deal with right and left derived functors in an abstract setup (as 
opposed to the triangulated setup). 

\begin{dfn}[Right Derived Functor] \label{dfn:2080}
Consider a category $\cat{K}$ and a multiplicatively closed set of morphisms 
$\cat{S} \sub \cat{K}$, with localization functor 
$\opn{Q} : \cat{K} \to \cat{K}_{\cat{S}}$. 
Let $F : \cat{K} \to \cat{E}$ be a functor. A
{\em right derived functor of $F$ with respect to $\cat{S}$}
\index{Derived functor! abstract right}
is a pair $(\mrm{R} F, \eta^{\mrm{R}})$, where 
\[ \mrm{R} F : \cat{K}_{\cat{S}} \to \cat{E} \]
is a functor, and 
\[ \eta^{\mrm{R}} : F \twoto \mrm{R} F \circ \opn{Q} \] 
is a morphism of functors, such that the following universal property holds:
\begin{itemize}
\item[(R)] Given any pair $(G, \th)$, consisting of a functor
$G : \cat{K}_{\cat{S}} \to \cat{E}$ 
and a morphism of functors 
$\th : F \twoto G \circ \opn{Q}$,
there is a unique morphism of functors 
$\mu : \mrm{R}^{} F \twoto G$
such that 
$\th = (\mu \circ \opn{id}_{\opn{Q}}) * \eta^{\mrm{R}}$.
\end{itemize}
\end{dfn}

Pictorially: there is a $2$-diagram 
\[ \UseTips  \xymatrix @C=12ex @R=8ex  {
\cat{K}
\ar[r]^{F}_(0.45){}="q" 
\ar[d]_{\opn{Q}}
& 
\cat{E} 
\\
\cat{K}_{\cat{S}}
\ar[ur]_{\mrm{R} F} ^(0.45){}="f"
\ar@{=>}  "q";"f" _{\eta^{\mrm{R}}}
} \]
For any other pair $(G, \th)$ there is a unique morphism $\mu$ that sits in 
this $2$-diagram 
\[ \UseTips  \xymatrix @C=14ex @R=10ex {
\cat{K}
\ar[r]^{F} _(0.45){}="q"  ^(0.55){}="q2" 
\ar[d]_{\opn{Q}}
& 
\cat{E} 
\\
\cat{K}_{\cat{S}}
\ar[ur] ^(0.45){}="f" ^(0.5){}="f1" 
\ar@{=>}  "q";"f" _{\eta^{\mrm{R}}}
\ar@(r,d)[ru]_{G} ^(0.55){}="g"  ^(0.48){}="g1"
\ar@{=>}  "q2";"g" ^(0.7){\th}
\ar@{=>}  "f";"g1" _(0.45){\mu}
} \]
such that the diagram of $2$-morphisms (with $*$ composition) 
\[ \UseTips \xymatrix @C=8ex @R=6ex { 
F
\ar@{=>}[d]_(0.45){\eta^{\mrm{R}}}
\ar@{=>}[dr]^{\th}
&
\\
\mrm{R} F \circ \opn{Q}
\ar@{=>}[r]_{\mu \, \circ \, \opn{id}_{\opn{Q}}}
&
G \circ \opn{Q}
} \]
is commutative. 

\begin{prop} \label{prop:1975}
If a right derived functor $(\mrm{R} F, \eta^{\mrm{R}})$ exists, then it is
unique, up to a unique isomorphism. Namely, if  $(G, \th)$ is another right 
derived functor of $F$, then there is a unique isomorphism of functors 
$\mu : \mrm{R} F \twoiso G$ such that 
$\th = (\mu \circ \opn{id}_{\opn{Q}}) * \eta^{\mrm{R}}$.
\end{prop}

\begin{proof}
Despite the apparent complication of the situation, the usual argument for 
uniqueness of universals applies. 
To be explicit, given another right derived functor $(G, \th)$, let 
$\mu : \mrm{R} F \twoto G$
be the unique morphism that's guaranteed by property (R) of the pair 
$(\mrm{R} F, \eta^{\mrm{R}})$. Then let 
$\nu : G \twoto \mrm{R} F$ be the unique morphism that's guaranteed by 
property (R) of the pair $(G, \th)$. Then the morphisms 
$\opn{id}_{\mrm{R} F}, \nu * \mu : \mrm{R} F \twoto \mrm{R} F$ 
satisfy 
\[ (\opn{id}_{\mrm{R} F} \circ \opn{id}_{\opn{Q}}) * \mu = \mu = 
((\nu * \mu) \circ \opn{id}_{\opn{Q}}) * \mu . \]
The uniqueness in property (R) of the pair $(\mrm{R} F, \eta^{\mrm{R}})$
implies that 
$\opn{id}_{\mrm{R} F} = \nu * \mu$. Likewise 
$\opn{id}_{G} = \mu * \nu$. 
Hence $\mu$ is an isomorphism of functors. 
\end{proof}

Here is a rather general existence result. 
Left denominator sets were introduced in Definition \ref{dfn:2311}. 

\begin{thm}[Existence of Right Derived Functor] \label{thm:2080}
In the situation of Definition \tup{\ref{dfn:2080}}, 
assume there is a full subcategory $\cat{J} \subseteq \cat{K}$
such the following three conditions hold\tup{:}
\begin{enumerate}
\rmitem{a} The multiplicatively closed set $\cat{S}$ 
is a left denominator set in $\cat{K}$. 

\rmitem{b} For every object $M \in \cat{K}$ there is a morphism 
$\rho : M \to I$ in $\cat{S}$, with target $I \in \cat{J}$. 

\rmitem{c} If $\psi$ is a morphism in $\cat{S} \cap \cat{J}$, then 
$F(\psi)$ is an isomorphism in $\cat{E}$.
\end{enumerate}
Then the right derived functor
\[ (\mrm{R} F, \eta^{\mrm{R}}) : \cat{K}_{\cat{S}} \to \cat{E} \]
exists. Moreover, for every object $I \in \cat{J}$ the morphism  
\[ \eta^{\mrm{R}}_{I} : F(I) \to \mrm{R} F(I) \]
in $\cat{E}$ is an isomorphism.
\end{thm}

Condition (b) says that {\em $\cat{K}$ has enough right $\cat{J}$-resolutions}. 
Condition (c) says that {\em $\cat{J}$ is an $F$-acyclic category}. 

Theorem \ref{thm:2080} is \cite[Proposition 7.3.2]{KaSc2}. However their 
notation is different: what we call ``left denominator set'', they call ``right 
multiplicative system''. 

We need a definition and a few lemmas before giving the proof of the theorem.
In them we assume the situation of the theorem. 

\begin{dfn} \label{dfn:2082}
In the situation of Theorem \ref{thm:2080}, by a {\em system of right 
$\cat{J}$-resolutions} we mean a pair $(I, \rho)$, where 
$I : \opn{Ob}(\cat{K}) \to \opn{Ob}(\cat{J})$
is a function, and 
$\rho = \{ \rho_M \}_{M \in \opn{Ob}(\cat{K})}$ 
is a collection of morphisms
$\rho_M : M \to I(M)$ in $\cat{S}$.
Moreover, if $M \in \opn{Ob}(\cat{J})$, then 
$I(M) = M$ and $\rho_M = \opn{id}_M$. 
\end{dfn}

Since here $\cat{K}$ has enough right $\cat{J}$-resolutions, it follows that
systems of right $\cat{J}$-resolutions $(I, \rho)$ exist.

Let us introduce some new notation that will make the proofs more readable:
\begin{equation} \label{eqn:2098}
\cat{K}' := \cat{J}, \quad \cat{S}' := \cat{J} \cap \cat{S}, 
\quad \cat{D} := \cat{K}_{\cat{S}} \quad \tup{and} 
\quad \cat{D}' := \cat{K}_{\cat{S}'}' . 
\end{equation}
The inclusion functor is $U : \cat{K}' \to \cat{K}$, and its localization is 
$V : \cat{D}' \to \cat{D}$.
These sit in a commutative diagram
\begin{equation} \label{eqn:2083}
\UseTips  \xymatrix @C=6ex @R=6ex {
\cat{K}'
\ar[r]^{U} _(0.5){}="s1"
\ar[d]_{\opn{Q}'}
&
\cat{K}
\ar[d]^{\opn{Q}}
\\
\cat{D}'
\ar[r]^{V} ^(0.5){}="t1"
&
\cat{D} 
}
\end{equation}

\begin{lem} \label{lem:2085}
The multiplicatively closed set $\cat{S}'$ is a left denominator set in 
$\cat{K}'$. 
\end{lem}

\begin{proof}
We need to verify conditions (LD1) and (LD2) in Definition \ref{dfn:2311}.

\medskip \noindent
(LD1): Given morphisms $a' : L' \to N'$ in $\cat{K}'$ and 
$s' : L' \to M'$ in $\cat{S}'$, we must find 
morphisms $b' : M' \to K'$ in $\cat{K}'$ and 
$t' : N' \to K'$ in $\cat{S}'$, such that $t' \circ a' = b' \circ s'$.
Because $\cat{S} \sub \cat{K}$ satisfies this condition, we can find 
morphisms $b : M' \to K$ in $\cat{K}$ and $t : N' \to K$ in $\cat{S}$
such that $t \circ a' = b \circ s'$. There is a morphism $\rho : K \to K'$ in 
$\cat{S}$ with target $K' \in \cat{K}'$. Then the morphisms $t' := \rho \circ 
t$ 
and $b' := \rho \circ b$ satisfy $t' \circ a' = b' \circ s'$, and 
$t' \in \cat{S}'$.

\medskip \noindent
(LD2): Given morphisms $a', b' : M' \to N'$ in $\cat{K}'$ and 
$s' : L' \to M'$ in $\cat{S}'$, that satisfy 
$a' \circ s' = b' \circ s'$, we must find 
a morphism $t' : N' \to K'$ in $\cat{S}'$ such that 
$t' \circ a' = t' \circ  b'$.
Because $\cat{S} \sub \cat{K}$ satisfies this condition, we can find 
a morphism $t : N' \to K$ in $\cat{S}$ such that 
$t \circ a' = t \circ  b'$.
There is a morphism $\rho : K \to K'$ in $\cat{S}$ with target 
$K' \in \cat{K}'$. Then the morphism $t' := \rho \circ t$ has the required 
property. 
\end{proof}

\begin{lem} \label{lem:2086}
The functor $V : \cat{D}' \to \cat{D}$ is an equivalence. 
\end{lem}

\begin{proof}
Condition (b) of the theorem implies that $V$ is essentially surjective on 
objects. We need to prove that $V$ is fully faithful. 
We shall use the left version of Proposition \ref{prop:3380}, namely $\cat{S}$ 
and $\cat{S}'$ are left denominator sets, and in condition (ii) the morphisms 
are $s : L' \to M$ and $t : M \to K'$. (See Remark \ref{rem:1845} and 
Propositions \ref{prop:3030} and \ref{prop:2310} regarding side changes.)
By Lemma \ref{lem:2085} the left version of condition (i) of 
Proposition \ref{prop:3380} holds. Condition (b) of the theorem implies 
the left version of condition (ii) of Proposition \ref{prop:3380}. 
Then the left version of Proposition \ref{prop:3380} says that 
$V$ is fully faithful. 
\end{proof}

\begin{lem} \label{lem:2081}
Suppose a system  of right $\cat{K}'$-resolutions 
$(I, \rho)$ has been chosen. Then the function 
$I : \opn{Ob}(\cat{K}) \to \opn{Ob}(\cat{K}')$ extends uniquely to a 
functor $I : \cat{D} \to \cat{D}'$, such that 
$I \circ V = \opn{Id}_{\cat{D}'}$, and 
$\opn{Q}(\rho) : \opn{Id}_{\cat{D}} \twoto V \circ I$
is an isomorphism of functors.
Therefore the functor $I$ is a a quasi-inverse of $V$. 
\end{lem}

The relevant $2$-diagram is this: 
\[ \UseTips  \xymatrix @C=10ex @R=8ex {
\cat{K}'
\ar[r]^{\opn{Q}'}
\ar[d]_{U}
&
\cat{D}'
\ar[dr]^{V} _(0.45){}="t" 
\ar[r]^{\opn{Id}}
&
\cat{D}'
\\
\cat{K}
\ar[r]^{\opn{Q}} _{}="q"
&
\cat{D} 
\ar[u]^{I}
\ar[r]_{\opn{Id}} ^(0.45){}="s"
\ar@{=>}  "s";"t" ^(0.45){\opn{Q}(\rho)}
&
\cat{D} 
\ar[u]_{I}
} \]
Recall that in a $2$-diagram, an empty polygon means it is commutative, namely 
it can be filled with 
$\xtwoto{\opn{id}}$. 

\begin{proof}
Consider a morphism $\psi : M \to N$ in $\cat{D}$.
Since $V : \cat{D}' \to \cat{D}$ is an equivalence, and since 
$V(I(M)) = I(M)$ and $V(I(N)) = I(N)$, there is a unique morphism 
$I(\psi) : I(M) \to I(N)$ 
in $\cat{D}'$ satisfying 
\begin{equation} \label{eqn:2085}
V(I(\psi)) := \opn{Q}(\rho_N) \circ \psi \circ \opn{Q}(\rho_M)^{-1} . 
\end{equation}
in $\cat{D}$. 

Let us check that $I : \cat{D} \to \cat{D}'$ is really a functor. 
Suppose $\phi : L \to M$ and $\psi : M \to N$ are morphisms in $\cat{D}$. 
Then 
\[ \begin{aligned}
& V \bigl( I(\psi) \circ I(\phi) \bigr) =  V(I(\psi)) \circ V(I(\phi))
\\ & \qquad 
 = \bigl( \opn{Q}(\rho_N) \circ \psi \circ \opn{Q}(\rho_M)^{-1} \bigr) \circ 
\bigl( \opn{Q}(\rho_M) \circ \phi \circ \opn{Q}(\rho_L)^{-1} \bigr)
\\ & \qquad 
= \opn{Q}(\rho_N) \circ (\psi \circ \phi) \circ \opn{Q}(\rho_L)^{-1} 
= V \bigl( I(\psi \circ \phi) \bigr) . 
\end{aligned} \]
Since $V$ is an equivalence, it follows that 
$I(\psi) \circ I(\phi) = I(\psi \circ \phi)$. 

Because $\rho_{M'} : M' \to I(M')$ is the identity for any object 
$M' \in \cat{K}'$, we see that there is equality 
$I \circ V = \opn{Id}_{\cat{D}'}$. 
By the defining formula (\ref{eqn:2085}) of $I(\psi)$ we have a commutative 
diagram 
\[ \UseTips  \xymatrix @C=10ex @R=6ex {
V(I(M))
\ar[r]^{V(I(\psi))}
&
V(I(N))
\\
M
\ar[r]^{\psi}
\ar[u]^{\opn{Q}(\rho_M)}
&
N
\ar[u]_{\opn{Q}(\rho_N)}
} \]
in $\cat{D}$. Hence $\opn{Q}(\rho) : \opn{Id}_{\cat{D}} \twoto V \circ I$
is an isomorphism of functors.
\end{proof}

Diagram (\ref{eqn:2083}) induces a commutative diagram of categories:
\begin{equation} \label{eqn:2084}
\UseTips \xymatrix @C=16ex @R=5.5ex {
\cat{Fun}(\cat{K}', \cat{E})
&
\cat{Fun}(\cat{K}, \cat{E})
\ar[l]_{\cat{Fun}(U, \opn{Id})}
\\
\cat{Fun}_{\cat{S}'}(\cat{K}', \cat{E})
\ar[u]_{\mrm{f.f. \, emb}}
&
\cat{Fun}_{\cat{S}}(\cat{K}, \cat{E})
\ar[u]^{\mrm{f.f. \, emb}}
\ar[l]^{\mrm{equiv}}_{\cat{Fun}(U, \opn{Id})}
\\
\cat{Fun}(\cat{D}', \cat{E})
\ar[u]^{\cat{Fun}(\opn{Q}', \opn{Id})}_{\mrm{isom}}
&
\cat{Fun}(\cat{D}, \cat{E})
\ar[l]_{\cat{Fun}(V, \opn{Id})}^{\mrm{equiv}}
\ar[u]_{\cat{Fun}(\opn{Q}, \opn{Id})}^{\mrm{isom}}
} 
\end{equation}
The vertical arrows marked ``f.f.\ emb'' are fully faithful embeddings by 
definition. According to Proposition \ref{prop:2026}  the vertical arrows 
marked ``isom'' are isomorphisms of categories. And by Lemma \ref{lem:2086} and 
Proposition \ref{prop:2025} the arrow $\cat{Fun}(V, \opn{Id})$  in the bottom 
row is an equivalence. As a consequence, the arrow $\cat{Fun}(U, \opn{Id})$ in 
the middle row is also an equivalence.

Let's introduce the notation 
$F' := F \circ U : \cat{K}' \to \cat{E}$.
This functor is an object of the category in the middle left of diagram 
(\ref{eqn:2084}), since, by condition (c) of Theorem \ref{thm:2080}, it inverts 
$\cat{S}'$.  

\begin{lem} \label{lem:3390}
Let $G : \cat{D} \to \cat{E}$ be a functor. Given a morphism 
$\th' : F' \twoto G \circ \opn{Q} \circ \, U$
of functors $\cat{K}' \to \cat{E}$, there is a unique morphism 
$\th : F \twoto G \circ \opn{Q}$
of functors $\cat{K} \to \cat{E}$ s.t.\ 
$\th' = \th \circ \opn{id}_{U}$. 
\end{lem}

Note that $G$ is an object in the category in the bottom right of diagram 
(\ref{eqn:2084}). The morphisms $\th'$ and $\th$ are in the middle left and top 
right respectively of this diagram. 

\begin{proof}
For every object $M \in \cat{K}$ there is a morphism 
$\rho : M \to I$ in $\cat{S}$ with target $I \in \cat{K}'$. 
A morphism of functors $\th$ satisfying $\th' = \th \circ \opn{id}_{U}$
will make this diagram commutative:
\begin{equation} \label{eqn:3390}
\UseTips  \xymatrix @C=8ex @R=6ex {
F(M)
\ar[r]^(0.4){\th_M}
\ar[d]_{F(\rho)}
&
(G \circ \opn{Q})(M)
\ar[d]^{(G \circ \opn{Q})(\rho)}_{\cong}
\\
F(I)
\ar[r]^(0.4){\th'_I}
&
(G \circ \opn{Q})(I)
} 
\end{equation}
We are using the facts that $I = U(I)$ and that $\opn{Q}(\rho)$ is an 
isomorphism. This proves the uniqueness of $\th$. 

For existence, let us choose a system  of right $\cat{K}'$-resolutions 
$(I, \rho)$, and define $\th = \{ \th_M \}$ using (\ref{eqn:3390}), namely 
\begin{equation} \label{eqn:3394}
\th_M := (G \circ \opn{Q})(\rho_M)^{-1} \circ \th'_{I(M)} \circ 
F(\rho_M) .
\end{equation}
We must prove that this is indeed a morphism of functors 
$\cat{K} \to \cat{E}$. Namely, for a given morphism 
$\phi : M \to N$ in $\cat{K}$, we have to prove that the diagram 
\begin{equation} \label{eqn:3391}
\UseTips  \xymatrix @C=8ex @R=6ex {
F(M)
\ar[r]^(0.4){\th_M}
\ar[d]_{F(\phi)}
&
(G \circ \opn{Q})(M)
\ar[d]^{(G \circ \opn{Q})(\phi)}
\\
F(N)
\ar[r]^(0.4){\th_N}
&
(G \circ \opn{Q})(N)
}
\end{equation}
in $\cat{E}$ is commutative. 

Lemma \ref{lem:2085} tells us that the morphism $I(\opn{Q}(\phi))$ in 
$\cat{D}'$ can be written as a left fraction
\begin{equation} \label{eqn:4715}
I(\opn{Q}(\phi)) = \opn{Q}'(\psi_1)^{-1} \circ \opn{Q}'(\psi_0)
\end{equation}
of morphisms $\psi_0 \in \cat{K}'$ and $\psi_1 \in \cat{S}'$. 
We get these diagrams: 
\begin{equation} \label{eqn:2072}
\UseTips \xymatrix @C=4ex @R=5.5ex {
M
\ar[rr]^{\phi}
\ar[d]_{\rho_M}
&
&
N
\ar[d]^{\rho_N}
\\
I(M)
\ar[dr]_{\psi_0}
&
&
I(N)
\ar[dl]^{\psi_1}
\\
&
J
} 
\quad \quad 
\UseTips \xymatrix @C=4,5ex @R=5.5ex {
M
\ar[rr]^{\opn{Q}(\phi)}
\ar[d]_{\opn{Q}(\rho_M)}
&
&
N
\ar[d]^{\opn{Q}(\rho_N)}
\\
I(M)
\ar[rr]^{ V(I(\opn{Q}(\phi))) }
\ar[dr]_{\opn{Q}(\psi_0)}
&
&
I(N)
\ar[dl]^{\opn{Q}(\psi_1)}
\\
&
J
} 
\end{equation}
The first diagram is in the category $\cat{K}$, and it might fail to be 
commutative. The second diagram is in the category $\cat{D}$, and it is 
commutative: the bottom triangle commutes by formula (\ref{eqn:4715}) and the 
equality $V \circ \opn{Q}' = \opn{Q}$; and the top square commutes by 
formula (\ref{eqn:2085}). By condition (LO4) of the left Ore localization 
$\opn{Q} : \cat{K} \to \cat{D}$, there is a morphism 
$\psi : J \to L$ in $\cat{S}$ such that 
$\psi \circ \psi_0 \circ \rho_M = \psi \circ \psi_1 \circ \rho_N \circ \phi$
in $\cat{K}$. There is the morphism 
$\rho_L : L \to I(L)$ in $\cat{S}$, whose target $I(L)$ belongs to $\cat{K}'$. 
Thus, after replacing the object $J$ with $I(L)$, the morphism 
$\psi_0$ by $\rho_L \circ \psi \circ \psi_0$, and the morphism 
$\psi_1$ by $\rho_L \circ \psi \circ \psi_1$, and noting that the latter is a 
morphism in $\cat{S}'$, we can now assume that the first diagram in 
(\ref{eqn:2072}) commutative too.

Now we embed (\ref{eqn:3391}) in the bigger diagram (\ref{eqn:3392})
in $\cat{E}$. 
Since $\rho_N \in \cat{S}$ and $\psi_1 \in \cat{S}'$, 
the morphisms $(G \circ \opn{Q})(\rho_N)$, $(G \circ \mrm{Q})(\psi_1)$ and 
$F(\psi_1)$ are isomorphisms. The top and bottom squares in (\ref{eqn:3392}) 
are commutative by the definition of $\th_M$ and $\th_N$, see formula 
(\ref{eqn:3394}). The left and right rounded shapes (those involving $J$) are 
commutative because the first diagram in (\ref{eqn:2072}) is commutative. 
\begin{equation} \label{eqn:3392}
\UseTips  \xymatrix @C=6ex @R=5.5ex {
&
F(I(M))
\ar[r]^(0.45){\th'_{I(M)}}
\ar @/_2.5em/ [ddl]_{F(\psi_0)}
&
(G \circ \opn{Q})(I(M))
\ar @/^3em/ [ddr]^{(G \circ \mrm{Q})(\psi_0)}
\\
&
F(M)
\ar[u]^{F(\rho_M)}
\ar[r]^(0.4){\th_M}
\ar[d]_{F(\phi)}
&
(G \circ \opn{Q})(M)
\ar[u]_{(G \circ \opn{Q})(\rho_M)}^{\cong}
\ar[d]^{(G \circ \opn{Q})(\phi)}
\\
F(J)
&
F(N)
\ar[r]^(0.4){\th_N}
\ar[d]_{F(\rho_N)}
&
(G \circ \opn{Q})(N)
\ar[d]^{(G \circ \opn{Q})(\rho_N)}_{\cong}
&
(G \circ \opn{Q})(J)
\\
&
F(I(N))
\ar@(l,d)[ul]^{F(\psi_1)}_{\cong}
\ar[r]^(0.45){\th'_{I(N)}}
&
(G \circ \opn{Q})(I(N))
\ar@(r,d)[ur]_{\ (G \circ \mrm{Q})(\psi_1)}^{\cong}
}
\end{equation}

\medskip 
Since $\th' : F' \twoto G \circ \opn{Q} \circ \, U$ is a morphism of functors, 
we have a commutative diagram 
\begin{equation} \label{eqn:4700}
\UseTips  \xymatrix @C=8ex @R=5.5ex {
F(I(M))
\ar[r]^(0.45){\th'_{I(M)}}
\ar[d]_{F(\psi_0)}
&
(G \circ \opn{Q})(I(M))
\ar[d]^{(G \circ \opn{Q})(\psi_0)}
\\
F(J)
\ar[r]^(0.45){\th'_{J}}
&
(G \circ \opn{Q})(J)
\\
F(I(N))
\ar[u]^{F(\psi_1)}
\ar[r]^(0.45){\th'_{I(N)}}
&
(G \circ \opn{Q})(I(N))
\ar[u]_{(G \circ \opn{Q})(\psi_1)}
}
\end{equation}

\medskip \noindent
When we erase the arrow $\th'_{J}$ from diagram (\ref{eqn:4700}), we obtain the 
outer boundary of diagram (\ref{eqn:3392}). Therefore diagram 
(\ref{eqn:3392}) is commutative. In particular, the middle square of diagram 
(\ref{eqn:3392}) is commutative, and it is precisely diagram (\ref{eqn:3391}).
\end{proof}

\begin{proof}[Proof of Theorem \tup{\ref{thm:2080}}] 
The proof is divided into four steps.  

\smallskip \noindent 
Step 1. Recall that the functor $F' = F \circ U$ lives in the middle left term 
in diagram (\ref{eqn:2084}). Because the arrow $\cat{Fun}(\opn{Q}', \opn{Id})$
is an isomorphism, there is a unique functor 
$\mrm{R} F'$ living in the bottom left term of diagram (\ref{eqn:2084}) that 
satisfies $\mrm{R} F' \circ \opn{Q}' = F'$. See next commutative diagram. 
\begin{equation} \label{eqn:2091}
\UseTips \xymatrix @C=10ex @R=6ex {
\cat{K}'
\ar[r]^{F'}
\ar[d]_{\opn{Q}'}
&
\cat{E}
\\
\cat{D}'
\ar[ur]_{\mrm{R} F'}
} 
\end{equation}

Let $\eta' := \opn{id}_{F'}$.
We claim that the pair $(\mrm{R} F', \eta')$ is a right derived functor of 
$F'$. Indeed, suppose we are given a pair $(G', \th')$, 
where $G'$ is a functor in the bottom left corner of diagram (\ref{eqn:2084}), 
and $\th' : F' \twoto G'  \circ \opn{Q}'$ 
is a morphism in the top left corner of that diagram. See the $2$-diagram
(\ref{eqn:2093}). 
Because the function 
\begin{equation} \label{eqn:2092} 
\opn{Hom}_{\cat{Fun}(\cat{D}', \cat{E})}
(\mrm{R} F', G') \to 
\opn{Hom}_{\cat{Fun}(\cat{K}', \cat{E})}
(F', G' \circ \opn{Q}') 
\end{equation}
is bijective -- this is the left edge of diagram (\ref{eqn:2084}) -- 
there is a unique morphism 
$\mu' : \mrm{R} F' \twoto G'$ that goes to $\th'$ under (\ref{eqn:2092}). 
\begin{equation} \label{eqn:2093} 
\UseTips  \xymatrix @C=14ex @R=10ex {
\cat{K}'
\ar[r]^{F'} _(0.45){}="q"  ^(0.6){}="q2" 
\ar[d]_{\opn{Q}'}
& 
\cat{E} 
\\
\cat{D}'
\ar[ur]^(0.3){\mrm{R} F'} ^(0.45){}="f" ^(0.5){}="f1" 
\ar@{=>}  "q";"f" _{\eta'}
\ar@(r,d)[ru]_{G'} ^(0.55){}="g"  ^(0.48){}="g1"
\ar@{=>}  "q2";"g" ^(0.7){\th'}
\ar@{=>}  "f";"g1" _(0.45){\mu'}
} 
\end{equation}

\medskip \noindent 
Step 2. Now we choose a system of right $\cat{K}'$-resolutions 
$(I, \rho)$, in the sense of Definition \ref{dfn:2082}.
By Lemma \ref{lem:2081}  we get an equivalence of categories
$I : \cat{D} \to \cat{D}'$, that is a quasi-inverse to $V$, and 
an isomorphism of functors 
$\opn{Q}(\rho) : \opn{Id}_{\cat{D}} \iso V \circ I$. 
See the following $2$-diagram (the solid arrows). 
\begin{equation} \label{eqn:2094} 
\UseTips  \xymatrix @C=10ex @R=8ex {
\cat{K}'
\ar[r]^{\opn{Q}'}
\ar[d]_{U}
\ar@(u,u)[rrr]^{F'}
&
\cat{D}'
\ar[dr]^{V} _(0.45){}="t" 
\ar[r]^{\opn{Id}}
&
\cat{D}'
\ar[r]^{\mrm{R} F'}
&
\cat{E}
\\
\cat{K}
\ar[r]_{\opn{Q}} _{}="q"
&
\cat{D} 
\ar[u]^{I}
\ar[r]_{\opn{Id}} ^(0.45){}="s"
\ar@{=>}  "s";"t" ^(0.45){\opn{Q}(\rho)}
&
\cat{D} 
\ar[u]_{I}
\ar@{-->}[ur]_{\mrm{R} F}
} 
\end{equation}

Define the functor 
\begin{equation} \label{eqn:2095}  
\mrm{R} F := \mrm{R} F' \circ I :  \cat{D} \to \cat{E} .
\end{equation}
It is the dashed arrow in diagram (\ref{eqn:2094}). 
So the functor $\mrm{R} F$ lives in the bottom right corner of 
(\ref{eqn:2084}), and $\mrm{R} F' = \mrm{R} F \circ V$.

\medskip \noindent 
Step 3. Consider Lemma \ref{lem:3390}, with the functor 
$G := \mrm{R} F$, and the morphism of functors
\[ \eta' = \opn{id}_{F'} : F' \twoto \mrm{R} F' \circ \opn{Q}' = 
\mrm{R} F \circ \opn{Q} \circ \, U . \]
The lemma says that there is a unique morphism of functors 
$\eta^{\mrm{R}} : F \twoto \mrm{R} F \circ \opn{Q}$
s.t.\ $\eta^{\mrm{R}} \circ \opn{id}_U = \eta'$. 

We can give an explicit formula for the morphism of functors $\eta^{\mrm{R}}$. 
Take an object $M \in \cat{K}$. Then the morphism 
$\eta^{\mrm{R}}_M : F(M) \to \mrm{R} F(M) = F(I(M))$
in $\cat{E}$ is nothing but 
\begin{equation} \label{eqn:2096} 
\eta^{\mrm{R}}_M := F(\rho_M) . 
\end{equation}
Here is the calculation. By formula (\ref{eqn:3394}) we have
\[ \eta^{\mrm{R}}_M = (\mrm{R} F \circ \opn{Q})(\rho_M)^{-1} \circ \eta'_{I(M)} 
\circ F(\rho_M) . \]
But 
$\eta'_{I(M)} = \opn{id}_{F(I(M))}$, 
\[ I(\opn{Q}(\rho_M)) = \opn{Q}(\rho_{I(M)}) \circ \opn{Q}(\rho_M) \circ 
\opn{Q}(\rho_M)^{-1} = \opn{id}_{I(M)} \]
by (\ref{eqn:2085}), and 
\[ (\mrm{R} F \circ \opn{Q})(\rho_M) = 
(\mrm{R} F' \circ I \circ \opn{Q})(\rho_M) = 
\mrm{R} F'(\opn{id}_{I(M)}) = \opn{id}_{F(I(M))} . \]
So everything else gets canceled, and we are left with (\ref{eqn:2096}). 

\medskip \noindent 
Step 4. It remains to prove that the pair $(\mrm{R} F, \eta^{\mrm{R}})$ is a 
right derived functor of $F$. Suppose $(G, \th)$ is a pair, where $G$ is a 
functor in the category in the bottom right corner of diagram (\ref{eqn:2084}), 
and $\th : F \twoto G  \circ \opn{Q}$ 
is a morphism in the top right corner of the diagram.
We are looking for a morphism 
$\mu : \mrm{R} F \twoto G$
in the bottom right category in diagram (\ref{eqn:2084}) for which 
$\th = (\mu \circ \opn{id}_{\opn{Q}}) * \eta^{\mrm{R}}$. 
Let $G' := G \circ V$, and let $\th' : F' \twoto G' \circ \opn{Q}'$
be the morphism in the top left corner of (\ref{eqn:2084}) corresponding to 
$\th$. Because of the equivalence $\cat{Fun}(V, \opn{Id})$, finding such $\mu$ 
is the same as finding a morphism
$\mu' : \mrm{R} F' \twoto G'$
in the bottom left category in diagram (\ref{eqn:2084}), satisfying
\begin{equation} \label{eqn:2075}
\th' = (\mu' \circ \opn{id}_{\opn{Q}'}) * \eta' .
\end{equation}

Finally, by step 1 the pair $(\mrm{R} F', \eta')$ is a right derived functor 
of $F'$. This says that there is a unique morphism 
$\mu'$ satisfying (\ref{eqn:2075}). 
\end{proof}

\begin{dfn} \label{dfn:1350}
The construction of the right derived functor $(\mrm{R} F, \eta^{\mrm{R}})$ in 
the proof of Theorem \ref{thm:2080}, and specifically formulas (\ref{eqn:2095}) 
and (\ref{eqn:2096}), is called a {\em presentation of 
$(\mrm{R} F, \eta^{\mrm{R}})$ by the system of right $\cat{J}$-resolutions 
$(I, \rho)$}. 
\end{dfn}

Of course any other right derived functor of $F$ (perhaps presented  
by another system of right $\cat{J}$-resolutions) is uniquely isomorphic to 
$(\mrm{R} F, \eta^{\mrm{R}})$. This is according to Proposition 
\ref{prop:1350}. In Sections \ref{sec:resol} and \ref{sec:exist-resol} we shall 
give several existence results for right resolutions by suitable acyclic 
objects. 

Now to left derived functors. 

\begin{dfn}[Left Derived Functor] \label{dfn:2090}
Consider a category $\cat{K}$ and a multiplicatively closed set of morphisms 
$\cat{S} \sub \cat{K}$, with localization functor 
$\opn{Q} : \cat{K} \to \cat{K}_{\cat{S}}$. 
Let $F : \cat{K} \to \cat{E}$ be a functor. A
{\em left derived functor of $F$ with respect to $\cat{S}$} 
\index{Derived functor! abstract left}
is a pair $(\mrm{L} F, \eta^{\mrm{L}})$, where 
\[ \mrm{L} F : \cat{K}_{\cat{S}} \to \cat{E} \]
is a functor, and 
\[ \eta^{\mrm{L}} : \mrm{L} F \circ \opn{Q} \twoto F \] 
is a morphism of functors, such that the following universal property holds:
\begin{itemize}
\item[(L)] Given any pair $(G, \th)$, consisting of a functor
$G : \cat{K}_{\cat{S}} \to \cat{E}$ 
and a morphism of functors 
$\th : G \circ \opn{Q} \twoto F$,
there is a unique morphism of functors 
$\mu : G  \twoto \mrm{L}^{} F$
such that 
$\th = \eta^{\mrm{L}} * (\mu \circ \opn{id}_{\opn{Q}})$.
\end{itemize}
\end{dfn}

Pictorially: there is a $2$-diagram 
\[ \UseTips  \xymatrix @C=12ex @R=8ex  {
\cat{K}
\ar[r]^{F} _(0.45){}="q" 
\ar[d]_{\opn{Q}}
& 
\cat{E} 
\\
\cat{K}_{\cat{S}}
\ar[ur]_{\mrm{L} F} ^(0.45){}="f"
\ar@{=>}  "f";"q" ^{\eta^{\mrm{L}}}
} \]
For any other pair $(G, \th)$ there is a unique morphism $\mu$ that sits in 
this $2$-diagram
\[ \UseTips  \xymatrix @C=14ex @R=10ex {
\cat{K}
\ar[r]^{F} _(0.45){}="q"  ^(0.55){}="q2" 
\ar[d]_{\opn{Q}}
& 
\cat{E} 
\\
\cat{K}_{\cat{S}}
\ar[ur] ^(0.45){}="f" ^(0.5){}="f1" 
\ar@{=>}  "f";"q" ^{\eta^{\mrm{L}}}
\ar@(r,d)[ru]_{G} ^(0.55){}="g"  ^(0.48){}="g1"
\ar@{=>}  "g";"q2" _(0.7){\th}
\ar@{=>}  "g1";"f" ^(0.45){\mu}
} \]
such that the diagram of $2$-morphisms (with $*$ composition) 
\[ \UseTips \xymatrix @C=8ex @R=6ex { 
F
&
\\
\mrm{L} F \circ \opn{Q}
\ar@{=>}[u]^{\eta^{\mrm{L}}}
&
G \circ \opn{Q}
\ar@{=>}[ul]_{\th}
\ar@{=>}[l]^(0.4){\mu \, \circ \, \opn{id}_{\opn{Q}}}
} \]
is commutative. 

We see that a left derived functor with target $\cat{E}$ amounts to a right 
derived functor with target $\cat{E}^{\mrm{op}}$. It means that new proofs are 
not needed. 

\begin{prop} \label{prop:2090}
If a left derived functor $(\mrm{L} F, \eta^{\mrm{L}})$ exists, then it is
unique, up to a unique isomorphism. Namely, if  $(G, \th)$ is another left 
derived functor of $F$, then there is a unique isomorphism of functors 
$\mu : G \twoiso \mrm{L} F$ such that 
$\th = \eta^{\mrm{L}} * (\mu \circ \opn{id}_{\opn{Q}})$.
\end{prop}

The proof is the same as that of Proposition \ref{prop:1975}, only some arrows 
have to be reversed. 

\begin{thm}[Existence of Left Derived Functor] \label{thm:2091}
In the situation of Definition \tup{\ref{dfn:2090}}, 
assume there is a full subcategory $\cat{P} \sub \cat{K}$
such the following three conditions hold\tup{:}
\begin{enumerate}
\rmitem{a} The multiplicatively closed set $\cat{S}$ 
is a right denominator set in $\cat{K}$. 

\rmitem{b} For every object $M \in \cat{K}$ there is a morphism 
$\rho : P \to M$ in $\cat{S}$, with source $P \in \cat{P}$. 

\rmitem{c} If $\psi$ is a morphism in $\cat{P} \cap \cat{S}$, then 
$F(\psi)$ is an isomorphism in $\cat{E}$.
\end{enumerate}
Then the left derived functor
\[ (\mrm{L} F, \eta^{\mrm{L}}) : \cat{K}_{\cat{S}} \to \cat{E} \]
exists. Moreover, for every object $P \in \cat{P}$ the morphism  
\[ \eta^{\mrm{L}}_{P} : \mrm{L} F(P) \to F(P) \]
in $\cat{E}$ is an isomorphism.
\end{thm}

Condition (b) says that {\em $\cat{K}$ has enough left $\cat{P}$-resolutions}. 
Condition (c) says that {\em $\cat{P}$ is an $F$-acyclic category}. 

The proof is the same as that of Theorem \ref{thm:2080}, only some arrows 
have to be reversed. We leave this as an exercise. 

\begin{exer} \label{exer:2330}
Prove Theorem \ref{thm:2091}, including the lemmas leading to it. 
(Hint: replace $\cat{E}$ with $\cat{E}^{\mrm{op}}$.)
\end{exer}

For reference we give the next definition.  

\begin{dfn} \label{dfn:2095}
In the situation of Theorem \ref{thm:2091}, by a {\em system of left 
$\cat{P}$-resolutions} we mean a pair $(P, \rho)$, where 
$P : \opn{Ob}(\cat{K}) \to \opn{Ob}(\cat{P})$
is a function, and 
$\rho = \{ \rho_M \}_{M \in \opn{Ob}(\cat{K})}$ 
is a collection of morphisms
$\rho_M : P(M) \to M$ in $\cat{S}$.
Moreover, if $M \in \opn{Ob}(\cat{P})$, then 
$P(M) = M$ and $\rho_M = \opn{id}_M$. 
\end{dfn}

When $\cat{K}$ has enough left $\cat{P}$-resolutions, there exists
a system of left $\cat{P}$-resolutions $(P, \rho)$ .

\begin{dfn} \label{dfn:2331}
The construction of the left derived functor $(\mrm{L} F, \eta^{\mrm{L}})$ in 
the proof of Theorem \ref{thm:2091} -- i.e.\ the left variant of the proof of 
Theorem \ref{thm:2080}, and specifically the left versions of formulas 
(\ref{eqn:2095}) and (\ref{eqn:2096}) -- 
is called a {\em presentation of $(\mrm{L} F, \eta^{\mrm{L}})$ by the system of 
left $\cat{P}$-resolutions $(P, \rho)$}. 
\end{dfn}

Of course any other left derived functor of $F$ (perhaps presented  
by another system of left $\cat{P}$-resolutions) is uniquely isomorphic to 
$(\mrm{L} F, \eta^{\mrm{L}})$. This is according to Proposition 
\ref{prop:2090}. In Sections \ref{sec:resol} and \ref{sec:exist-resol} we shall 
give several existence results for left resolutions by suitable acyclic 
objects. 

\begin{rem} \label{rem:4696}
The right derived functor $\mrm{R} F$ from Definition \ref{dfn:2080} is a {\em 
left Kan extension} of $F$ along $\opn{Q}$. Likewise, the left derived functor 
$\mrm{L} F$ from Definition \ref{dfn:2090} is a {\em right Kan extension} of 
$F$ along $\opn{Q}$. See \cite[Chapter X]{Mac2}.
\end{rem}

\mysubsection{Triangulated Derived Functors} \label{subsec:tri-der-funs}

In this subsection we specialize the definitions and results of the previous 
subsection to the case of triangulated functors between triangulated 
categories. 
There is a fixed nonzero commutative base ring $\K$, and all categories and 
functors here are $\K$-linear. 

Triangulated functors and morphisms between them were introduced in Definition
\ref{dfn:1276}. Recall that a triangulated functor is a pair $(F, \tau)$, 
consisting of a linear functor $F$ and a translation isomorphism $\tau$, 
that respect the distinguished triangles. 

Here is the triangulated version of Definition \ref{dfn:2025}. 

\begin{dfn} \label{dfn:2335}
Let $\cat{K}$ and $\cat{L}$ be triangulated categories. We define 
$\cat{TrFun}(\cat{K}, \cat{L})$
\index{1-TrFun@$\cat{TrFun}(\cat{K}, \cat{L})$}
to be the category whose objects are the triangulated functors 
$(F, \tau) : \cat{K} \to \cat{L}$. Given objects 
$(F, \tau)$ and $(G, \nu)$ of $\cat{TrFun}(\cat{K}, \cat{L})$,
the morphisms 
$\al : (F, \tau) \twoto (G, \nu)$ 
in  $\cat{TrFun}(\cat{K}, \cat{L})$
are the morphisms of triangulated functors. 
\end{dfn}

\begin{lem} \label{lem:2370}  
Let 
$(F, \tau) : (\cat{K}, \opn{T}_{\cat{K}}) \to (\cat{L}, \opn{T}_{\cat{L}})$
be a triangulated functor between 
triangulated categories. Assume $F$ is an equivalence \tup{(}of 
abstract categories\tup{)}, with quasi-inverse 
$G : \cat{L} \to \cat{K}$, and with adjunction isomorphisms 
$\al : G \circ F \twoiso \opn{Id}_{\cat{K}}$ and  
$\be : F \circ G \twoiso \opn{Id}_{\cat{L}}$. 
Then there is an isomorphism of functors 
$\nu : G \circ \opn{T}_{\cat{L}} \twoiso \opn{T}_{\cat{K}} \circ \, G$ 
such that 
$(G, \nu) : (\cat{L}, \opn{T}_{\cat{L}}) \to (\cat{K}, \opn{T}_{\cat{K}})$
is a triangulated functor, and $\al$ and $\be$ are isomorphisms of 
triangulated functors.
\end{lem}

\begin{proof}
It is well-known that $G$ is additive (or in our case, $\K$-linear); but 
since 
the proof is so easy, we shall reproduce it. Take any pair of objects 
$M, N \in \cat{L}$. We have to prove that the bijection
\[ G_{M, N} : \opn{Hom}_{\cat{L}}(M, N) \to 
\opn{Hom}_{\cat{K}} \bigl( G(M), G(N) \bigr) \]
is linear. But 
\[ G_{M, N} = (F_{G(M), G(N)})^{-1} \circ
\opn{Hom}_{\cat{L}}(\be_M, \be_N^{-1}) \]
as bijections (of sets) between these modules. Since 
$\opn{Hom}_{\cat{L}}(\be_M, \be_N^{-1})$
and $F^{-1}_{G(M), G(N)}$ are $\K$-linear, so is $G_{M, N}$.

We define the isomorphism of functors $\nu$ by the formula 
\[ \nu := (\al \circ \opn{id}_{\opn{T}_{\cat{K}} \circ \, G}) * 
(\opn{id}_G \circ \, \tau \circ \opn{id}_G)^{-1} * 
(\opn{id}_{G \circ \opn{T}_{\cat{L}}} \circ \, \be)^{-1} , \]
in terms of the $2$-categorical notation.
This gives rise to a commutative diagram of isomorphisms
\[ \UseTips \xymatrix @C=12ex @R=6ex { 
G \circ \opn{T}_{\cat{L}} \circ \, F \circ G
\ar@{=>}[d]_{\opn{id} \circ \, \be}
& 
G \circ F \circ \opn{T}_{\cat{K}} \circ \, G
\ar@{=>}[l]_{\opn{id} \circ \, \tau \, \circ \, \opn{id}}
\ar@{=>}[d]^{\al \, \circ \, \opn{id}}
\\
G \circ \opn{T}_{\cat{L}}
\ar@{=>}[r]^{\nu}
&
\opn{T}_{\cat{K}} \circ \, G
} \]
of additive functors $\cat{L} \to \cat{K}$. 
So the pair $(G, \nu)$ is a T-additive functor. 

The verification that $(G, \nu)$ preserves triangles (in the sense of 
Definition \ref{dfn:1276}(1)) is done like the proof of the additivity of 
$G$, but now using axiom (TR1)(a) from Definition \ref{dfn:1158}.
Also $\al$ and $\be$ are morphisms of triangulated functors 
(Definitions \ref{dfn:1276}(2) and \ref{dfn:1160}). 
We leave these verifications as an exercise. 
\end{proof}

\begin{exer}
Finish the proof above (the last assertions).  
\end{exer}

From here on in this subsection we shall usually keep the translation 
isomorphisms (such as $\tau$ and $\nu$ above) implicit. 

Suppose we are given triangulated functors 
$U : \cat{K}' \to \cat{K}$ and 
$V : \cat{L} \to \cat{L}'$. There is an induced functor 
\[ \cat{TrFun}(U, V) : 
\cat{TrFun}(\cat{K}, \cat{L}) \to \cat{TrFun}(\cat{K}', \cat{L}') ; \]
the formula is the same as in (\ref{eqn:2025}). 

\begin{lem} \label{lem:2335}
If $U$ and $V$ are equivalences, then the functor $\cat{TrFun}(U, V)$ is an 
equivalence.
\end{lem}

\begin{proof}
Use Proposition \ref{prop:2025} and  Lemma \ref{lem:2370}. 
\end{proof}

\begin{dfn} \label{dfn:2336}
Let $\cat{K}$ and $\cat{L}$ be triangulated categories,
and let $\cat{S} \sub \cat{K}$ be a denominator set of cohomological origin. 
We define 
$\cat{TrFun}_{\cat{S}}(\cat{K}, \cat{L})$
\index{1-TrFunS@$\cat{TrFun}_{\cat{S}}(\cat{K}, \cat{L})$}
to be the full subcategory of 
$\cat{TrFun}(\cat{K}, \cat{L})$
whose objects are the functors that are localizable to $\cat{S}$, in the 
sense of Definition \ref{dfn:2037}.
\end{dfn}

We know that the localization functor 
$\opn{Q} : \cat{K} \to \cat{K}_{\cat{S}}$ is a left and right Ore 
localization. 

\begin{lem} \label{lem:2336} 
Let $\cat{K}$ and $\cat{E}$ be triangulated categories,
and let $\cat{S} \sub \cat{K}$ be a denominator set of cohomological origin. 
Then the functor 
\[ \cat{TrFun}(\opn{Q}, \opn{Id}_{\cat{E}}) : 
\cat{TrFun}(\cat{K}_{\cat{S}}, \cat{E}) \to 
\cat{TrFun}_{\cat{S}}(\cat{K}, \cat{E}) \]
is an isomorphism of categories.
\end{lem}

\begin{proof}
Use Proposition \ref{prop:2026} and Theorem \ref{thm:105}. 
\end{proof}

Here is the triangulated version of Definition \ref{dfn:2080}.

\begin{dfn}  \label{dfn:131}
Let $F : \cat{K} \to \cat{E}$ be a triangulated functor between  
triangulated categories, and let 
$\cat{S} \sub \cat{K}$ be a denominator set of cohomological origin.    
A {\em triangulated right derived functor of $F$ with respect to $\cat{S}$} 
\index{Derived functor! triangulated right}
\index{Triangulated functor! right derived}
is a triangulated functor 
$\mrm{R}^{} F : \cat{K}_{\cat{S}} \to \cat{E}$,
together with a morphism 
$\eta^{\mrm{R}} : F \twoto \mrm{R}^{} F \circ \opn{Q}$ 
of triangulated functors $\cat{K} \to \cat{E}$. The pair 
$(\mrm{R} F, \eta^{\mrm{R}})$ must have this universal property:
\begin{itemize}
\item[(R)] Given any pair $(G, \th)$, consisting of a triangulated 
functor $G : \cat{K}_{\cat{S}} \to \cat{E}$ and a morphism of triangulated 
functors $\th : F \twoto G \circ \opn{Q}$, 
there is a unique morphism of triangulated functors 
$\mu : \mrm{R}^{} F \twoto G$
such that $\th = (\mu \circ \opn{id}_{\opn{Q}}) * \eta^{\mrm{R}}$. 
\end{itemize}
\end{dfn}

\begin{prop} \label{prop:1350}
If a triangulated right derived functor $(\mrm{R} F, \eta^{\mrm{R}})$ exists, 
then it is unique, up to a unique isomorphism. Namely, if  $(G, \th)$ is 
another triangulated right derived functor of $F$, then there is a unique 
isomorphism of 
triangulated functors $\mu : \mrm{R} F \twoiso G$ such that 
$\th = (\mu \circ \opn{id}_{\opn{Q}}) * \eta^{\mrm{R}}$. 
\end{prop}

The proof is the same as that of Proposition \ref{prop:1975}.
Now for the triangulated version of Theorem \ref{thm:2080}.

\begin{thm} \label{thm:1470}
In the situation of Definition \tup{\ref{dfn:131}}, 
assume there is a full triangulated subcategory 
$\cat{J} \subseteq \cat{K}$
with these two properties\tup{:}
\begin{enumerate}
\rmitem{a} Every object $M \in \cat{K}$ admits a morphism 
$\rho : M \to I$ in $\cat{S}$ with target $I \in \cat{J}$.

\rmitem{b} If $\psi$ is a morphism in $\cat{J} \cap \cat{S}$, 
then $F(\psi)$ is an isomorphism in $\cat{E}$.
\end{enumerate}
Then the triangulated right derived functor
$(\mrm{R} F, \eta^{\mrm{R}}) : \cat{K}_{\cat{S}} \to \cat{E}$ 
of $F$ with respect to $\cat{S}$ exists. 
\index{Derived functor! triangulated right}
\index{Triangulated functor! right derived}
Moreover, for every object $I \in \cat{J}$ the morphism  
$\eta^{\mrm{R}}_I : F(I) \to \mrm{R} F(I)$ 
in $\cat{E}$ is an isomorphism.
\end{thm}

\begin{proof}
It will be convenient to change notation to that used in the proof of Theorem 
\ref{thm:2080}. Let's define 
$\cat{K}' := \cat{J}$, $\cat{S}' := \cat{K}' \cap \cat{S}$ and
$\cat{D}' := \cat{K}'_{\cat{S}'}$. 
The localization functor of $\cat{K}'$ is 
$\opn{Q}' : \cat{K}' \to \cat{D}'$. 
The inclusion functor is 
$U : \cat{K}' \to \cat{K}$, and its localization is the functor
$V : \cat{D}' \to \cat{D}$.
We have this commutative diagram of triangulated functors between 
triangulated categories:
\begin{equation} \label{eqn:2340}
\UseTips \xymatrix @C=6ex @R=6ex {
\cat{K}'
\ar[r]^{U}
\ar[d]_{\opn{Q}'}
&
\cat{K}
\ar[d]^{\opn{Q}}
\\
\cat{D}'
\ar[r]^{V}
&
\cat{D}
} 
\end{equation}
The functor $U$ is fully faithful. 
By Lemma \ref{lem:2086} the functor $V$ is an equivalence.

Diagram (\ref{eqn:2340}) induces a commutative diagram of linear categories:
\begin{equation} \label{eqn:2346}
\UseTips \xymatrix @C=14ex @R=6ex {
\cat{TrFun}(\cat{K}', \cat{E})
&
\cat{TrFun}(\cat{K}, \cat{E})
\ar[l]_{\cat{Fun}(U, \opn{Id})}
\\
\cat{TrFun}_{\cat{S}'}(\cat{K}', \cat{E})
\ar[u]_{\mrm{f.f.\, emb}}
&
\cat{TrFun}_{\cat{S}}(\cat{K}, \cat{E})
\ar[u]^{\mrm{f.f.\, emb}}
\ar[l]^{\mrm{equiv}}_{\cat{Fun}(U, \opn{Id})}
\\
\cat{TrFun}(\cat{D}', \cat{E})
\ar[u]^{\cat{Fun}(\opn{Q}', \opn{Id})}_{\mrm{isom}}
&
\cat{TrFun}(\cat{D}, \cat{E})
\ar[l]_{\cat{Fun}(V, \opn{Id})}^{\mrm{equiv}}
\ar[u]_{\cat{Fun}(\opn{Q}, \opn{Id})}^{\mrm{isom}}
} 
\end{equation}

\noindent By definition the arrows marked ``f.f.\ emb'' are fully faithful 
embeddings. According to Lemma \ref{lem:2336}  the arrows 
$\cat{Fun}(\opn{Q}, \opn{Id})$ and $\cat{Fun}(\opn{Q}', \opn{Id})$
are isomorphisms of categories. By Lemma \ref{lem:2335} the 
arrow $\cat{Fun}(V, \opn{Id})$ is an equivalence. 
It follows that the arrow $\cat{Fun}(U, \opn{Id})$ in the middle row is an 
equivalence too. 

We know that $\cat{S} \sub \cat{K}$ is a left denominator set.
Condition (b) of the theorem says that $F$ sends morphisms in 
$\cat{S}'$ to isomorphisms in $\cat{E}$. 
Condition (a) there says that there are enough right 
$\cat{K}'$-resolutions in $\cat{K}$.
Thus we are in a position to use the abstract Theorem \ref{thm:2080}.
It says that there is an abstract right derived functor 
$(\mrm{R} F, \eta^{\mrm{R}}) : \cat{D} = \cat{K}_{\cat{S}} \to \cat{E}$ 
of $F$ with respect to $\cat{S}$. 
However, going over the proof of Theorem \ref{thm:2080},
we see that all constructions there can be made within the triangulated 
setting, namely in diagram (\ref{eqn:2346}) instead of in  diagram 
(\ref{eqn:2084}). 
Therefore $\mrm{R} F$ is an object of the category in the 
bottom right corner of (\ref{eqn:2346}), and the morphism 
$\eta^{\mrm{R}} : F \twoto \mrm{R} F \circ \opn{Q}$
is in the category in the top right corner of (\ref{eqn:2346}).
The triangulated variant of step 4 in the proof of Theorem \ref{thm:2080} 
shows that $(\mrm{R} F, \eta^{\mrm{R}})$ satisfies condition (R) of Definition 
\ref{dfn:131}.
\end{proof}

By slight abuse of notation, in the situation of Definition \ref{dfn:131} 
we sometimes refer to the triangulated right derived functor $\mrm{R} F$
just as ``the right derived functor of $F$''. As the next corollary shows, this 
does really cause a problem. 

\begin{cor} \label{cor:3225}
In the situation of Theorem \tup{\ref{thm:1470}}, the triangulated right 
derived functor $(\mrm{R} F, \eta^{\mrm{R}})$ is also an abstract right derived 
functor of $F$ with respect to $\cat{S}$, in the sense of Definition 
\tup{\ref{dfn:2080}}.
\end{cor}

\begin{proof}
This is seen in the proof of Theorem \tup{\ref{thm:1470}}.
\end{proof}

In the situation of Theorem \ref{thm:1470}, let $\cat{K}^{\star}$ be a full 
triangulated subcategory of $\cat{K}$. Define 
$\cat{S}^{\star} := \cat{K}^{\star} \cap \cat{S}$ and 
$\cat{J}^{\star} := \cat{K}^{\star} \cap \cat{J}$. 
Denote by $W : \cat{K}^{\star} \to \cat{K}$
the inclusion functor, and by 
$W_{\cat{S}^{\star}} : \cat{K}^{\star}_{\cat{S}^{\star}} \to \cat{K}_{\cat{S}}$
its localization. 
Warning: the functor $W_{\cat{S}^{\star}}$ is not necessarily fully faithful; 
cf.\ Proposition \ref{prop:1320}. 

\begin{prop} \label{prop:1355}
Assume that every object $M \in \cat{K}^{\star}$ admits a morphism 
$\rho : M \to I$ in $\cat{S}^{\star}$ with target $I \in \cat{J}^{\star}$. 
Then the pair 
$(\mrm{R}^{\star} F, \eta^{\star}) := 
(\mrm{R} F \circ W_{\cat{S}^{\star}}, \eta^{\mrm{R}} \circ \opn{id}_W)$
is a right derived functor of 
$F^{\star} := F \circ W : \cat{K}^{\star} \to \cat{E}$.
\end{prop}

\begin{exer} \label{exer:1355}
Prove the last proposition.  (Hint: Start by choosing a system of right 
$\cat{J}^{\star}$-resolutions in $\cat{K}^{\star}$. Then extend it 
to a system of right $\cat{J}^{}$-resolutions in $\cat{K}^{}$. 
Now follow the proofs of Theorems \ref{thm:1470} and \ref{thm:2080}.)
\end{exer}

Here is the specialization of Definition \ref{dfn:131} to categories of DG 
modules. First:  

\begin{setup} \label{set:3236}
We are given a DG ring $A$, an abelian category
$\cat{M}$, and a full triangulated subcategory
$\cat{K} \sub \dcat{K}(A, \cat{M})$. We write 
$\cat{S} := \cat{K} \cap \, \dcat{S}(A, \cat{M})$, 
the set of quasi-isomorphisms in $\cat{K}$, and 
$\cat{D} := \cat{K}_{\cat{S}}$, the derived category of $\cat{K}$. 
\end{setup}

\begin{dfn} \label{dfn:3236}
Under Setup \ref{set:3236}, let
$\cat{E}$ be a triangulated category, and let 
$F : \cat{K} \to \cat{E}$ 
be a triangulated functor.  
A {\em triangulated right derived functor} 
\index{Derived functor! triangulated right}
\index{Triangulated functor! right derived}
of $F$ is a triangulated right derived functor 
$(\mrm{R} F, \eta^{\mrm{R}}) : \cat{D} \to \cat{E}$
of $F$ with respect to $\cat{S}$, in the sense 
of Definition \ref{dfn:131}. 
\end{dfn}

\begin{exa} \label{exa:1305}
Suppose we start from an additive functor 
$F^0 : \cat{M} \to \cat{N}$
between abelian categories. We know how to extend it to a DG functor 
\[ \bcat{C}^{+}(F^0) : \bcat{C}^{+}(\cat{M}) \to \bcat{C}^{+}(\cat{N}) , \]
and this induces a triangulated functor 
\[ \bcat{K}^{+}(F^0) : \bcat{K}^{+}(\cat{M}) \to  \bcat{K}^{+}(\cat{N}) . \]
By composing with the localization functor 
$\opn{Q}^+_{\cat{N}} : \bcat{K}^{+}(\cat{N}) \to \bcat{D}^{+}(\cat{N})$
we get a triangulated functor 
\begin{equation} \label{eqn:5030}
F := \opn{Q}^+_{\cat{N}} \circ \, \bcat{K}^{+}(F^0) : \bcat{K}^{+}(\cat{M}) 
\to \bcat{D}^{+}(\cat{N}) .
\end{equation}
Define $\cat{K} := \bcat{K}^{+}(\cat{M})$,
$\cat{S} := \bcat{S}^{+}(\cat{M})$ and  
$\cat{E} := \bcat{D}^{+}(\cat{N})$, so in this notation we have a triangulated 
functor $F : \cat{K} \to \cat{E}$, and we are in the situation of Definition 
\ref{dfn:3236}. Note that the restriction of $F$ to the full subcategory 
$\cat{M} \sub \bcat{K}^{+}(\cat{M})$ coincides with the original functor $F^0$. 

Assume that the abelian category $\cat{M}$ has enough injectives; namely 
that every object $M \in \cat{M}$ admits a monomorphism to an injective object. 
Let $\cat{J}$ be the full subcategory of $\cat{K}$ on the bounded below 
complexes of injective objects. 
We will prove later (see Corollary \ref{cor:1922} and Theorem \ref{thm:2880}) 
that properties (a) and (b) of Theorem \ref{thm:1470} hold 
for this $\cat{J}$, regardless what $F$ is. 
Therefore the triangulated right derived functor 
\begin{equation} \label{eqn:5031}
(\mrm{R}^+ F, \eta^+) : \bcat{D}^{+}(\cat{M}) \to \bcat{D}^{+}(\cat{N})
\end{equation}
exists, and for every complex $I \in \cat{J}$ the morphism 
\begin{equation} \label{eqn:2347}
\eta^+_I : F(I) \to \mrm{R}^+ F(I) 
\end{equation}
in $\bcat{D}^{+}(\cat{N})$ is an isomorphism. 

Now assume the original functor $F^0 : \cat{M} \to \cat{N}$ is left exact.
For each $q \geq 0$ we have the classical right derived functor
$\mrm{R}^q F^0 : \cat{M} \to \cat{N}$, and $\mrm{R}^0 F^0 \cong F^0$. 
Formula (\ref{eqn:2347}) shows that for every $M \in \cat{M}$ and $q \geq 0$ 
there is an isomorphism
$\mrm{R}^q F^0(M) \cong \opn{H}^q(\mrm{R}^+ F(M))$  
in $\cat{N}$. As the object $M \in \cat{M}$ moves,
this becomes an isomorphism 
$\mrm{R}^q F^0 \cong \opn{H}^q \circ \, \mrm{R}^+ F$
of additive functors $\cat{M} \to \cat{N}$.

In this example we were careful to use distinct notations for the functor 
$F^0$ between abelian categories, and the induced triangulated functor 
$F = \opn{Q}^+_{\cat{N}} \circ \, \bcat{K}^{+}(F^0)$.
Later we will use the same notation for $F^0$ and $F$. 
\end{exa}

Our treatment of triangulated left derived functors will be 
brief: we will state the definitions and the main results, but won't give 
proofs, beyond a hint here and there on the passage from right to left derived 
functors. 

\begin{dfn}  \label{dfn:1360}
Let $F : \cat{K} \to \cat{E}$ be triangulated functor between  
triangulated categories, and let 
$\cat{S} \sub \cat{K}$ be a denominator set of cohomological origin.    
A {\em triangulated left derived functor of $F$ with respect to $\cat{S}$} 
\index{Derived functor! triangulated left}
\index{Triangulated functor! left derived}
is a triangulated functor 
$\mrm{L}^{} F : \cat{K}_{\cat{S}} \to \cat{E}$,
together with a morphism 
$\eta^{\mrm{L}} : \mrm{L}^{} F \circ \opn{Q} \twoto F$ 
of triangulated functors $\cat{K} \to \cat{E}$. The pair 
$(\mrm{L} F, \eta^{\mrm{L}})$ must have this universal property:
\begin{itemize}
\item[(L)] Given any pair $(G, \th)$, consisting of a triangulated 
functor $G : \cat{K}_{\cat{S}} \to \cat{E}$ and a morphism of triangulated 
functors $\th : G \circ \opn{Q} \twoto F$, 
there is a unique morphism of triangulated functors 
$\mu : G \twoto \mrm{L}^{} F$
such that $\th = \eta^{\mrm{L}} * (\mu \circ \opn{id}_{\opn{Q}})$.
\end{itemize}
\end{dfn}

\begin{prop} \label{prop:1360}
If a triangulated left derived functor $(\mrm{L} F, \eta^{\mrm{L}})$ exists, 
then it is unique, up to a unique isomorphism. Namely, if  $(G, \th)$ is 
another triangulated left derived functor of $F$, then there is a unique 
isomorphism of triangulated functors $\mu : G \twoiso \mrm{L} F$ such that 
$\th = \eta^{\mrm{L}} * (\mu \circ \opn{id}_{\opn{Q}})$. 
\end{prop}

The proof is the same as that of Proposition \ref{prop:1350}, with direction of 
arrows in $\cat{E}$ reversed. 

\begin{thm} \label{thm:1500}
In the situation of Definition \tup{\ref{dfn:1360}}, 
assume there is a full triangulated subcategory 
$\cat{P} \subseteq \cat{K}$
with these two properties\tup{:}
\begin{enumerate}
\rmitem{a} Every object $M \in \cat{K}$ admits a morphism 
$\rho : P \to M$ in $\cat{S}$ with source $P \in \cat{P}$.

\rmitem{b} If $\psi$ is a morphism in $\cat{P} \cap \cat{S}$, 
then $F(\psi)$ is an isomorphism in $\cat{E}$.
\end{enumerate}
Then the triangulated left derived functor
$(\mrm{L} F, \eta^{\mrm{L}}) : \cat{K}_{\cat{S}} \to \cat{E}$ 
\index{Derived functor! triangulated left}
\index{Triangulated functor! left derived}
of $F$ with respect to $\cat{S}$ exists. 
Moreover, for every object $P \in \cat{P}$ the morphism  
$\eta^{\mrm{L}}_P : \mrm{L} F(P) \to F(P)$
in $\cat{E}$ is an isomorphism.
\end{thm}

The proof is the same as that of Theorem \ref{thm:1470}, with direction of 
arrows in $\cat{E}$ reversed.

\begin{cor} \label{cor:3226}
In the situation of Theorem \tup{\ref{thm:1500}}, the triangulated left derived 
functor $(\mrm{L} F, \eta^{\mrm{L}})$ is also an abstract left derived functor 
of $F$ with respect to $\cat{S}$, in the sense of Definition 
\tup{\ref{dfn:2090}}.
\end{cor}

The proof is the same as that of Corollary \ref{cor:3225}, with direction of 
arrows in $\cat{E}$ reversed.

In the situation of Theorem \ref{thm:1500}, let $\cat{K}^{\star}$ be a full 
triangulated subcategory of $\cat{K}$. Define 
$\cat{S}^{\star} := \cat{K}^{\star} \cap \cat{S}$ and 
$\cat{P}^{\star} := \cat{K}^{\star} \cap \cat{P}$. 
Denote by $W : \cat{K}^{\star} \to \cat{K}$
the inclusion functor, and by 
$W_{\cat{S}^{\star}} : \cat{K}^{\star}_{\cat{S}^{\star}} \to \cat{K}_{\cat{S}}$
its localization. 

\begin{prop} \label{prop:1502}
Assume that every $M \in \cat{K}^{\star}$ admits a morphism 
$\rho : P \to M$ in $\cat{S}^{\star}$ with source $P \in \cat{P}^{\star}$. 
Then the pair 
$(\mrm{L}^{\star}, \eta^{\star}) := 
(\mrm{L} F \circ W_{\cat{S}^{\star}}, \eta^{\mrm{L}} \circ \opn{id}_W)$ 
is a triangulated left derived functor of 
$F^{\star} := F \circ W : \cat{K}^{\star} \to \cat{E}$.
\end{prop}

The proof is just like that of Proposition \ref{prop:1355} (which was an 
exercise).

Here is the specialization of Definition \ref{dfn:1360} to categories of DG 
modules.

\begin{dfn} \label{dfn:3237}
Under Setup \ref{set:3236}, let
$\cat{E}$ be a triangulated category, and let 
$F : \cat{K} \to \cat{E}$ be a triangulated functor.  
A {\em triangulated left derived functor}
\index{Derived functor! triangulated left}
\index{Triangulated functor! left derived}
of $F$ is a triangulated left derived functor 
$(\mrm{L} F, \eta^{\mrm{L}}) : \cat{D} \to \cat{E}$
of $F$ with respect to $\cat{S}$, in the sense 
of Definition \ref{dfn:1360}. 
\end{dfn}

\begin{exa} \label{exa:1501}
This is analogous to Example \ref{exa:1305}.
We start with an additive functor 
$F^0 : \cat{M} \to \cat{N}$ between abelian categories. We know how to extend 
it 
to a DG functor 
\[ \bcat{C}^{-}(F^0) : \bcat{C}^{-}(\cat{M}) \to \bcat{C}^{-}(\cat{N}) , \]
and this induces a triangulated functor 
\[ \bcat{K}^{-}(F^0) : \bcat{K}^{-}(\cat{M}) \to  \bcat{K}^{-}(\cat{N}) . \]
By composing with $\opn{Q}^-_{\cat{N}}$ we get a triangulated functor
\begin{equation} \label{eqn:5032}
F := \opn{Q}^-_{\cat{N}} \circ \, \bcat{K}^{-}(F^0) : \bcat{K}^{-}(\cat{M}) 
\to \bcat{D}^{-}(\cat{N}) . 
\end{equation}
Defining $\cat{K} := \bcat{K}^{-}(\cat{M})$,
$\cat{S} := \bcat{S}^{+}(\cat{M})$ and 
$\cat{E} := \bcat{D}^{-}(\cat{N})$,
we have a triangulated functor 
$F : \cat{K} \to \cat{E}$, and the situation is that of Definition 
\ref{dfn:3237}.

Assume that the abelian category $\cat{M}$ has enough projectives. 
Define $\cat{P}$ to be the full subcategory of $\cat{K}$ on the bounded above 
complexes of projective objects. 
We will prove later (in Corollary \ref{cor:1923}  and Theorem \ref{thm:1540}) 
that properties (a) and (b) of Theorem \ref{thm:1500} hold 
in this situation. Therefore we have a left derived functor
\begin{equation} \label{eqn:5033}
(\mrm{L}^- F, \eta^-) : \bcat{D}^{-}(\cat{M}) \to \bcat{D}^{-}(\cat{N}) . 
\end{equation}
For any $P \in \cat{P}$ the morphism 
\begin{equation} \label{eqn:2351}
\eta^-_P : \mrm{L}^- F(P) \to F(P) 
\end{equation}
in $\cat{E}$ is an isomorphism. 

Now assume the functor $F^0$ is right exact. Then $F^0$ has the classical 
left derived functors 
$\mrm{L}_q F^0 : \cat{M} \to \cat{N}$
for all $q \geq 0$, and $\mrm{L}_0 F^0 \cong F^0$. Formula (\ref{eqn:2351}) 
shows that for every $M \in \cat{M}$ there is an isomorphism 
$\mrm{L}_q F^0(M) \cong \opn{H}^{-q}(\mrm{L}^- F(M))$
in $\cat{N}$. As $M$ moves there is an isomorphism  
$\mrm{L}_q F^0 \cong \opn{H}^{-q} \circ \, \mrm{L}^- F$
of additive functors $\cat{M} \to \cat{N}$.

In this example, as in Example \ref{exer:1355}, we were careful to use 
distinct notations for the functor 
$F^0$ between abelian categories, and the induced triangulated functor 
$F = \opn{Q}^-_{\cat{N}} \circ \, \bcat{K}^{-}(F^0)$.
Later we will use the same notation for $F^0$ and $F$. 
\end{exa}

\begin{rem} \label{rem:4695}
The reader might ask the following question at this point: 
Suppose we are given additive functors
\begin{equation} \label{eqn:4705}
\cat{L} \xar{F} \cat{M} \xar{G} \cat{N}
\end{equation}
between abelian categories, and the categories 
$\cat{M}$ and $\cat{N}$ have enough injectives, so that the right derived 
functors 
\begin{equation} \label{eqn:4706}
\UseTips \xymatrix @C=8ex @R=8ex { 
\dcat{D}^+(\cat{L})
\ar[r]^{\mrm{R}^+ F}
\ar@(ur,ul)[rr]^{\mrm{R}^+ (G \circ F)}
&
\dcat{D}^+(\cat{M})
\ar[r]^{\mrm{R}^+ G}
&
\dcat{D}^+(\cat{N})
}
\end{equation}
exist (see Example \ref{exa:1305} above). Does this diagram commute (up to 
an isomorphism of triangulated functors)? Similar questions can be asked about 
left derived functors, etc. 

Here is an answer.  Assume that $F$ sends injectives in 
$\cat{L}$ to {\em right $G$-acyclic objects} of $\cat{M}$. See Definition
\ref{dfn:3887} (that talks about categories of graded modules, but the ideas 
are the same). By Lemma \ref{lem:3885} we know that for each bounded 
below complex of right $G$-acyclic objects $J$, the canonical morphism 
$G(J) \to \mrm{R}^+ G(J)$ in $\dcat{D}^+(\cat{M})$ is an isomorphism. 
Therefore, for every bounded below complex of injectives
$I \in \dcat{C}(\cat{L})$,
letting $J := F(I)$,  we have $J \cong \mrm{R}^+ F(I)$, and hence 
\[ \mrm{R}^+ (G \circ F)(I) \cong (G \circ F)(I) = G(F(I)) = G(J) \cong 
\mrm{R}^+ G(J) \cong \mrm{R}^+ G (\mrm{R}^+ F(I)) . \]
Thus there is an isomorphism of triangulated functors  
\begin{equation} \label{eqn:4707}
\mrm{R}^+ (G \circ F) \cong \mrm{R}^+ G \circ \mrm{R}^+  F : 
\dcat{D}^+(\cat{L}) \to \dcat{D}^+(\cat{N}) .
\end{equation}

The original application of the isomorphism (\ref{eqn:4707}),
in \cite[Proposition II.5.1]{RD}, was this: given morphisms of schemes 
$X \xar{f} Y \xar{g} Z$, the derived pushforward functors satisfy
\begin{equation} \label{eqn:5034}
\mrm{R}^+ (g_* \circ f_*) \cong \mrm{R}^+ g_* \circ \mrm{R}^+ f_*
\end{equation}
as functors $\dcat{D}^+(X) \to \dcat{D}^+(Z)$. 
This is true because the functor $f_*$ sends injective $\OO_X$-modules to 
flasque $\OO_Y$-modules. 

This pleasant phenomenon seems to fail in some important situations 
that we study in our book. For instance, take a noncommutative ring $A$ and 
$A$-bimodules $S, T$ that are tilting complexes but not invertible 
bimodules (see Subsection \ref{subsec:tilting-rings}).  
Consider the abelian categories 
$\cat{L} = \cat{M} = \cat{N} := \dcat{M}(A)$
and the functors $F := \opn{Hom}_A(T, -)$ and $G := \opn{Hom}_A(S, -)$. 
It is quite likely that the functors $\mrm{R}^+G \circ \mrm{R}^+F$ and 
$\mrm{R}^+(G \circ F)$ are not isomorphic. 
A concrete instance in which this can be tested is the ring $A$ 
from Example \ref{exa:4870}, and the bimodules $S = T := A^*$ there.
\end{rem}

\begin{rem} \label{rem:4705}
The composition of derived functors explained in the previous remark is also 
related to {\em spectral sequences}.
\index{Spectral sequence}
Given the input (\ref{eqn:4705}), in classical homological algebra it was 
impossible to ask for something like the commutative diagram (\ref{eqn:4706}).
The best approximation was to have a {\em convergent spectral sequence}
\begin{equation} \label{eqn:4695}
E^{i, j}_k = \mrm{R}^j G \circ \mrm{R}^i F \Rightarrow
\mrm{R}^n (G \circ F) . 
\end{equation}
Indeed, under the assumption in the previous remark (namely that
$F$ sends injective objects to right $G$-acyclic objects), such a convergent 
spectral sequence does exist; see   
\cite[Chapters XI-XII]{Mac1}, \cite[Section 11]{Rot} or 
\cite[Section 5.7]{We}. In the geometric example above, for which the 
isomorphism (\ref{eqn:5034}) holds, this is the {\em Leray spectral 
sequence}.  

There are several kinds of spectral sequences, but basically they are all of the 
same nature, i.e.\ approximations of the composition of derived functors. 

Most traditional applications of spectral sequences are performed more 
effectively by filtrations and truncations of complexes, as can be seen in 
Section \ref{sec:adj-equ-cohdim} of this book.
Still, there are a few instances in which spectral sequences 
cannot be left out, such as the {\em niveau spectral sequence} in algebraic 
geometry (see \cite[Chapter IV]{RD}), or its noncommutative relative, the {\em 
double Ext spectral sequence} for Auslander dualizing complexes 
(see Remark \ref{rem:3821} and the paper \cite{YeZh1}).  
\end{rem}

\mysubsection{Contravariant Triangulated Derived Functors} 
\label{subsec:cntr-tri-der-funs}

In this subsection there is a fixed nonzero commutative base ring $\K$, and 
all categories and functors here are $\K$-linear. The next setup is assumed. 

\begin{setup} \label{set:3225}
We are given a DG ring $A$, an abelian category
$\cat{M}$, and a full additive subcategory
$\cat{K} \sub \dcat{K}(A, \cat{M})$ s.t.\ 
$\cat{K}^{\mrm{op}} \sub \dcat{K}(A, \cat{M})^{\mrm{op}}$
is triangulated. We write 
$\cat{S} := \cat{K} \cap \, \dcat{S}(A, \cat{M})$, 
the set of quasi-isomorphisms in $\cat{K}$, and 
$\cat{D}^{\mrm{op}} := (\cat{K}^{\mrm{op}})_{\cat{S}^{\mrm{op}}}$,
the derived category of $\cat{K}^{\mrm{op}}$. 
\end{setup}

As explained in Subsections \ref{subsec:opp-hom-triang} and 
\ref{subsec:opp-dercat-triang}, the categories 
$\cat{K}^{\mrm{op}}$ and $\cat{D}^{\mrm{op}}$ 
have canonical triangulated structures, and the localization functor 
$\opn{Q}^{\mrm{op}} : \cat{K}^{\mrm{op}} \to \cat{D}^{\mrm{op}}$
is triangulated. 

\begin{dfn} \label{dfn:2356}
\index{Derived functor! contravariant triangulated right}
\index{Triangulated functor! contravariant right derived}
Let $\cat{E}$ be a triangulated category, and let 
$F : \cat{K}^{\mrm{op}} \to \cat{E}$
be a triangulated functor.  
A {\em triangulated right derived functor} of $F$ is a triangulated 
right derived functor 
$(\mrm{R} F, \eta^{\mrm{R}}) : \cat{D}^{\mrm{op}} \to \cat{E}$
of $F$ with respect to $\cat{S}^{\mrm{op}}$, in the sense 
of Definition \ref{dfn:131}. 
\end{dfn}

Likewise:

\begin{dfn} \label{dfn:2357} 
\index{Derived functor! contravariant triangulated left}
\index{Triangulated functor! contravariant left derived}
Let $\cat{E}$ be a triangulated category, and let 
$F : \cat{K}^{\mrm{op}} \to \cat{E}$ 
be a triangulated functor. A {\em triangulated left derived functor} of $F$ is 
a triangulated left derived functor 
$(\mrm{L} F, \eta^{\mrm{L}}) : \cat{D}^{\mrm{op}} \to \cat{E}$
of $F$ with respect to $\cat{S}^{\mrm{op}}$, in the sense 
of Definition \ref{dfn:1360}. 
\end{dfn}

The uniqueness results (Propositions \ref{prop:1350} and \ref{prop:1360})
and existence results (Theorems \ref{thm:1470} and \ref{thm:1500}) 
apply in the contravariant case as well. We shall return to the existence 
question the end of Section \ref{sec:resol}, when we talk about {\em 
resolving subcategories}; this will make matters more concrete. 
For now we give one example. 

\begin{exa} \label{exa:3055}
Let $\cat{M}$ and $\cat{N}$ be abelian categories, and let 
$F^0 : \cat{M} \to \cat{N}$
be a contravariant additive functor. This is the same as an additive functor
$F^0 \circ \opn{Op} : \cat{M}^{\mrm{op}} \to \cat{N}$.
As in Theorem \ref{thm:3015}, $F^0$ gives 
rise to a contravariant DG functor 
\[ \dcat{C}^-(F^0) : \dcat{C}^-(\cat{M}) \to \dcat{C}^+(\cat{N}) . \]
This means, according to Proposition \ref{prop:2500}, that we have a  DG 
functor 
\[ \dcat{C}^-(F^0) \circ \opn{Op} : \dcat{C}^-(\cat{M})^{\mrm{op}} \to 
\dcat{C}^+(\cat{N}) . \]
Passing to homotopy categories, and using Theorem \ref{thm:2995}, we obtain a 
triangulated functor 
\[ \opn{Ho} \bigl( \dcat{C}^-(F^0) \circ \opn{Op} \bigr) :
\dcat{K}^-(\cat{M})^{\mrm{op}} \to \dcat{K}^+(\cat{N}) \]
(with a translation isomorphism $\bar{\tau}$ that we ignore 
now). Composing this with the localization functor
$\opn{Q}_{\cat{N}}$ we have a triangulated functor 
\begin{equation} \label{eqn:5035}
F := \opn{Q}_{\cat{N}} \circ \, 
\opn{Ho} \bigl( \dcat{C}^-(F^0) \circ \opn{Op} \bigr) : 
\dcat{K}^-(\cat{M})^{\mrm{op}} \to 
\dcat{D}^+(\cat{N}) .
\end{equation}
Similarly there is a triangulated functor 
\begin{equation} \label{eqn:5036}
F^{\mrm{flip}} := \opn{Q}_{\cat{N}} \circ \, 
\opn{Ho} \bigl( \dcat{C}^+(F^0 \circ \opn{Op}) \bigr) :
\dcat{K}^+(\cat{M}^{\mrm{op}}) \to \dcat{D}^+(\cat{N}) .
\end{equation}
Observe that this triangulated functor is of the same sort as the one we had in 
(\ref{eqn:5030}). The functors (\ref{eqn:5035}) and (\ref{eqn:5036}) 
fit into the the following commutative diagram of triangulated functors 
\[ \UseTips \xymatrix @C=10ex @R=6ex {
\dcat{K}^-(\cat{M})^{\mrm{op}}
\ar[r]^{F}
\ar[d]_{\ol{\opn{Flip}}}^{\cong}
&
\dcat{D}^+(\cat{N})
\\
\dcat{K}^+(\cat{M}^{\mrm{op}}) 
\ar[ur]_{F^{\mrm{flip}}}
} \]
where $\ol{\opn{Flip}}$ is the isomorphism of triangulated categories from 
(\ref{eqn:2995}). 

The triangulated right derived functor 
\begin{equation} \label{eqn:5037}
\mrm{R}^- F : \dcat{D}^-(\cat{M})^{\mrm{op}} \to \dcat{D}^+(\cat{N})
\end{equation}
of the functor $F$ from (\ref{eqn:5035}) would sit inside this diagram of 
triangulated functors, which is 
commutative up to a unique isomorphism: 
\[ \UseTips \xymatrix @C=10ex @R=6ex {
\dcat{D}^-(\cat{M})^{\mrm{op}}
\ar[r]^{\mrm{R}^-F}
\ar[d]_{\ol{\opn{Flip}}}^{\cong}
&
\dcat{D}^+(\cat{N})
\\
\dcat{D}^+(\cat{M}^{\mrm{op}}) 
\ar[ur]_{\mrm{R}^+ F^{\mrm{flip}}}
} \]

Like in Example \ref{exa:1305}, to show existence of $\mrm{R}^+ F^{\mrm{flip}}$ 
using Theorem \ref{thm:1470},
and to calculate it, we would need a full triangulated subcategory 
$\cat{J} \sub \dcat{K}^+(\cat{M}^{\mrm{op}})$, 
such that every $M \in \dcat{K}^+(\cat{M}^{\mrm{op}})$
admits a quasi-isomorphism 
$M \to I$ with $I \in \cat{J}$; and 
every quasi-isomorphism $\psi$ in $\cat{J}$ goes to an isomorphism 
$F^{\mrm{flip}}(\psi)$ in $\dcat{D}^+(\cat{N})$.
This translates, by the flip, to a full triangulated subcategory
$\cat{P} \sub \dcat{K}^-(\cat{M})$ such that 
$\cat{J} = \ol{\opn{Flip}}(\cat{P})$. 
The conditions on $\cat{P}$ are these: 
each $M \in \dcat{K}^-(\cat{M})$
admits a quasi-isomorphism 
$P \to M$ with $P \in \cat{P}$; and each quasi-isomorphism $\phi$ 
in $\cat{P}$ goes to an isomorphism
$F(\phi)$ in $\dcat{D}^+(\cat{N})$.
For instance, if $\cat{M}$ has enough projectives, then according to Corollary 
\ref{cor:1900}, the category $\cat{P}$ of bounded above complexes of 
projectives satisfies these conditions. 

To calculate the derived functor $\mrm{R}^- F$ in equation (\ref{eqn:5037}) we 
choose a system of left $\cat{P}$-resolutions 
in $\dcat{K}^-(\cat{M})$, as in Definition \ref{dfn:2095}. 
Thus for every $M$ we have a chosen quasi-isomorphism 
$\rho_M : P(M) \to M$ in $\dcat{K}^{-}(\cat{M})$. We take
$\mrm{R}^- F(M) := F(P(M)) \in \dcat{D}^+(\cat{N})$.
The morphism
$\eta^{-}_M : F(M) \to \mrm{R}^- F(M)$
is 
$\eta^{-}_M := F(\rho_M) : F(M) \to F(P(M))$.

The story is very similar for the triangulated left derived functor 
\[ \mrm{L}^+ F : \dcat{D}^+(\cat{M})^{\mrm{op}} \to \dcat{D}^-(\cat{N})  \]
of 
\[ F := \opn{Q}_{\cat{N}} \circ \, 
\opn{Ho} \bigl( \dcat{C}^+(F^0) \circ \opn{Op} \bigr) : 
\dcat{K}^+(\cat{M})^{\mrm{op}} \to 
\dcat{D}^-(\cat{N}) , \]
so we won't tell it. 
\end{exa}

In later parts of the book we shall use the same symbol $F$ to denote also the 
functors $F^0$, $F^0 \circ \opn{Op}$ and  
$\dcat{C}^{\star}(F^0) \circ \opn{Op}$
from the example above. This will simplify our notation. 

Though indirectly defined (in Definition \ref{dfn:2995}), the  translation 
functor $\opn{T}^{\mrm{op}}$
of $\dcat{K}(A, \cat{M})^{\mrm{op}}$
is known to have some useful properties. Note that the objects of 
$\dcat{K}(A, \cat{M})^{\mrm{op}}$ are the same as those of 
$\dcat{K}(A, \cat{M})$. 

\begin{prop} \label{prop:2307}
For every object $M \in \dcat{K}(A, \cat{M})^{\mrm{op}}$ there are canonical 
isomorphisms 
$\opn{T}^{\mrm{op}}(M)^i \cong M^{i - 1}$
and 
$\opn{H}^i(\opn{T}^{\mrm{op}}(M)) \cong \opn{H}^{i - 1}(M)$
in $\cat{M}$. 
\end{prop}

\begin{proof}
This is immediate from Theorem \ref{thm:2495} and Definition \ref{dfn:2995}.
\end{proof}

Here are several useful properties of the triangulated structure on 
$\dcat{D}(A, \cat{M})^{\mrm{op}}$.

\begin{prop} \label{prop:3060}
Suppose 
\[ 0 \to L \xar{\phi} M \xar{\psi} N \to 0 \]
is an exact sequence in 
$\dcat{C}_{\mrm{str}}(A, \cat{M})^{\mrm{op}}$. 
Then there is a distinguished triangle 
\[ L \xar{\opn{Q}^{\mrm{op}}(\phi)} 
M \xar{\opn{Q}^{\mrm{op}}(\psi)} 
N \xar{\th} \opn{T}^{\mrm{op}}(L) \]
in $\dcat{D}(A, \cat{M})^{\mrm{op}}$.
\end{prop}

\begin{proof}
Apply the functor $\opn{Flip}$ to the given exact sequence. By Theorem 
\ref{thm:2495}(3) we obtain a short exact sequence 
\[ 0 \to \opn{Flip}(L) \xar{\opn{Flip}(\phi)} \opn{Flip}(M) 
\xar{\opn{Flip}(\psi)} \opn{Flip}(N) \to 0 \]
in 
$\dcat{C}_{\mrm{str}}(A^{\mrm{op}}, \cat{M}^{\mrm{op}})$. 
According to Proposition \ref{prop:2165} there is a  distinguished triangle
\[ \opn{Flip}(L) \xar{\opn{Q}^{\mrm{flip}}(\opn{Flip}(\phi))}
\opn{Flip}(M) \xar{\opn{Q}^{\mrm{flip}}(\opn{Flip}(\psi))}
\opn{Flip}(N) \xar{\th'}
\opn{T}^{\mrm{flip}}(\opn{Flip}(L)) \]
in 
$\dcat{D}(A^{\mrm{op}}, \cat{M}^{\mrm{op}})$.
Take the morphism 
\[ \th := \opn{Flip}^{-1}(\th') : N \to 
\opn{Flip}^{-1}(\opn{T}^{\mrm{flip}}(\opn{Flip}(L))) = 
\opn{T}^{\mrm{op}}(L) \]
in $\dcat{D}(A, \cat{M})^{\mrm{op}}$.
Since 
$\opn{Flip}^{-1}(\opn{Q}^{\mrm{flip}}(\opn{Flip}(\phi))) = 
\opn{Q}^{\mrm{op}}(\phi)$,
and likewise for $\psi$, we have the desired distinguished triangle. 
\end{proof}

Recall from Proposition \ref{prop:3170}
that if $A$ is a nonpositive DG ring, then for any 
$M \in \dcat{C}(A, \cat{M})$ and integer $i$ the smart truncations 
$\opn{smt}^{\geq i}(M), \, \opn{smt}^{\leq i}(M) \in \dcat{C}(A, \cat{M})$
exist. 

\begin{prop} \label{prop:3061}
Assume $A$ is a nonpositive DG ring. 
For every object $M \in \dcat{C}(A, \cat{M})$
there is a distinguished triangle 
\[ \opn{smt}^{\geq i + 1}(M) \to M \to 
\opn{smt}^{\leq i}(M) \to 
\opn{T}^{\mrm{op}}(\opn{smt}^{\geq i + 1}(M)) \]
in $\dcat{D}(A, \cat{M})^{\mrm{op}}$. The truncations are performed in 
$\dcat{C}(A, \cat{M})$ as above. 
\end{prop}

\begin{proof}
Let 
$N := \opn{Flip}(M) \in \dcat{D}(A^{\mrm{op}}, \cat{M}^{\mrm{op}})$. 
By Proposition \ref{prop:2320} there is a distinguished triangle
\[ \opn{smt}^{\leq i}(N) \xar{} N \xar{} \opn{smt}^{\geq i + 1}(N) 
\xar{} \opn{T}^{\mrm{flip}}(\opn{smt}^{\leq i}(N)) \]
in $\dcat{D}(A^{\mrm{op}}, \cat{M}^{\mrm{op}})$,
where now the truncations are done in 
$\dcat{C}(A^{\mrm{op}}, \cat{M}^{\mrm{op}})$. 
The formulas defining the flip functor in the proof of Theorem \ref{thm:2495}
show that \lb
$\opn{Flip}^{-1}(\opn{smt}^{\leq i}(N)) = \opn{smt}^{\geq i + 1}(M)$
and 
$\opn{Flip}^{-1}(\opn{smt}^{\geq i + 1}(N)) = \opn{smt}^{\leq i}(M)$.
\end{proof}

%% file: block3_190413.tex
\renewcommand{\thisfile}{block3\_190328}   

\cleardoublepage
\mysection{DG and Triangulated Bifunctors} \label{sec:DG-tri-bifun}

\AYcopyright 
 
In this section we extend the theory of triangulated derived functors, that was 
studied in Section \ref{sec:der-funcs}, to  
{\em triangulated bifunctors}. As stated in Convention \ref{conv:2490}, all 
linear structures and operations (such as categories and functors) are  
$\K$-linear, where $\K$ is a fixed commutative base 
ring. The symbol $\ot$ means $\ot_{\K}$.

\mysubsection{DG Bifunctors} \label{subsec:dg-bifun}
We already talked about bifunctors in Subsection \ref{subsec:bifunc}. That 
was for categories without further structure. Here we will consider 
DG categories, and matters become a bit more complicated. 

As a warmup, let us begin with linear bifunctors. 

\begin{dfn} \label{dfn:5040}
Let $\cat{C}_1$ and $\cat{C}_2$ be linear categories. 
We define the linear category $\cat{C}_1 \ot \cat{C}_2$ as follows. The 
set of objects is 
\[ \opn{Ob}(\cat{C}_1 \ot \cat{C}_2) :=
\opn{Ob}(\cat{C}_1) \times \opn{Ob}(\cat{C}_2) . \]
For each pair of objects 
$(M_1, M_2), (N_1, N_2)$ in $\cat{C}_1 \ot \cat{C}_2$,
i.e.\ $M_i, N_i \in \opn{Ob}(\cat{C}_i)$,
we let 
\[ \opn{Hom}_{\cat{C}_1 \ot \cat{C}_2}
\bigl( (M_1, M_2), (N_1, N_2) \bigr) := 
\opn{Hom}_{\cat{C}_1} (M_1, N_1) \ot 
\opn{Hom}_{\cat{C}_2} (M_2, N_2) . \]

Given morphisms 
$\phi_i \in \opn{Hom}_{\cat{C}_i}(L_i, M_i)$
and
$\psi_i \in \opn{Hom}_{\cat{C}_i}(M_i, N_i)$, 
for $i = 1, 2$, their tensors are morphisms 
\[ \phi_1 \ot \phi_2 \in \opn{Hom}_{\cat{C}_1 \ot \cat{C}_2}
\bigl( (L_1, L_2), (M_1, M_2) \bigr) \]
and 
\[ \psi_1 \ot \psi_2 \in 
\opn{Hom}_{\cat{C}_1 \ot \cat{C}_2}
\bigl( (M_1, M_2), (N_1, N_2) \bigr) . \]
Every morphism in $\cat{C}_1 \ot \cat{C}_2$ is a sum of such pure tensors. We 
define the composition to be 
\[ \begin{aligned}
& (\psi_1 \ot \psi_2) \circ (\phi_1 \ot \phi_2) := 
(\psi_1 \circ \phi_1) \ot (\psi_2 \circ \phi_2)
\\ & \quad 
\in \opn{Hom}_{\cat{C}_1 \ot \cat{C}_2}
\bigl( (L_1, L_2), (N_1, N_2) \bigr) .
\end{aligned} \]
\end{dfn}

\begin{dfn} \label{dfn:1965}
Let $\cat{C}_1$, $\cat{C}_2$ and $\cat{D}$ be linear categories. A 
{\em linear bifunctor}%
\index{Bifunctor! linear}
\[ F : \cat{C}_1 \times \cat{C}_2 \to \cat{D} \]
is, by definition, a linear functor 
$F : \cat{C}_1 \ot \cat{C}_2 \to \cat{D}$.
\end{dfn}

Now to the DG situation. 

\begin{dfn} \label{dfn:2010}
Let $\cat{C}_1$ and $\cat{C}_2$ be DG categories. 
We define the DG category $\cat{C}_1 \ot \cat{C}_2$ as follows: the 
set of objects is 
\[ \opn{Ob}(\cat{C}_1 \ot \cat{C}_2) :=
\opn{Ob}(\cat{C}_1) \times \opn{Ob}(\cat{C}_2) . \]
For each pair of objects 
$(M_1, M_2), (N_1, N_2)$ in $\cat{C}_1 \ot \cat{C}_2$
we let 
\[ \opn{Hom}_{\cat{C}_1 \ot \cat{C}_2}
\bigl( (M_1, M_2), (N_1, N_2) \bigr) := 
\opn{Hom}_{\cat{C}_1} (M_1, N_1) \ot 
\opn{Hom}_{\cat{C}_2} (M_2, N_2) . \]
The formula for the composition involves the Koszul sign rule. Given 
morphisms 
$\phi_i \in \opn{Hom}_{\cat{C}_i}(L_i, M_i)^{d_i}$
and
$\psi_i \in \opn{Hom}_{\cat{C}_i}(M_i, N_i)^{e_i}$, 
for $i = 1, 2$, their tensors are morphisms 
\[ \phi_1 \ot \phi_2 \in 
\opn{Hom}_{\cat{C}_1 \ot \cat{C}_2}
\bigl( (L_1, L_2), (M_1, M_2) \bigr)^{d_1 + d_2} \]
and 
\[ \psi_1 \ot \psi_2 \in 
\opn{Hom}_{\cat{C}_1 \ot \cat{C}_2}
\bigl( (M_1, M_2), (N_1, N_2) \bigr)^{e_1 + e_2} . \]
We define the composition to be 
\[ \begin{aligned}
& (\psi_1 \ot \psi_2) \circ (\phi_1 \ot \phi_2) := 
(-1)^{d_1 \cd e_2} \cd (\psi_1 \circ \phi_1) \ot
(\psi_2 \circ \phi_2)
\\ & \qquad 
\in \opn{Hom}_{\cat{C}_1 \ot \cat{C}_2}
\bigl( (L_1, L_2), (N_1, N_2) \bigr)^{d_1 + d_2 + e_1 + e_2} . 
\end{aligned} \]
\end{dfn}

\begin{exa}  \label{exa:2010}
Suppose $\cat{C}_1$ and $\cat{C}_2$ are single-object DG 
categories, say $\opn{Ob}(\cat{C}_i) = \{ x_i \}$. 
Then $\cat{C}_1 \ot \cat{C}_2$ also has a single-object, say 
$y := (x_1, x_2)$. 
The endomorphism DG rings satisfy
\[ \opn{End}_{\cat{C}_1 \ot \cat{C}_2}(y) = 
\opn{End}_{\cat{C}_1}(x_1) \ot \opn{End}_{\cat{C}_2}(x_2) . \]
See formula \ref{eqn:4106} and Example \ref{exa:1100}. 
\end{exa}

DG functors between DG categories were introduced in Definition \ref{dfn:1065}.

\begin{dfn} \label{dfn:1966}
Let $\cat{C}_1$, $\cat{C}_2$ and $\cat{D}$ be DG categories. A 
{\em DG bifunctor}%
\index{Bifunctor! differential graded}
\[ F : \cat{C}_1 \times \cat{C}_2 \to \cat{D} \]
is, by definition, a DG functor
$F : \cat{C}_1 \ot \cat{C}_2 \to \cat{D}$,
where $\cat{C}_1 \ot \cat{C}_2$ is the DG category from Definition
\ref{dfn:2010}. 
\end{dfn}

Warning: due to the signs that the composition of odd morphisms acquires, 
a DG bifunctor $F$ does not become, after we forget the DG structure,
a linear bifunctor in the sense of Definition \ref{dfn:1965}.

\begin{prop} \label{prop:5040}
Given a DG bifunctor
$F : \cat{C}_1 \times \cat{C}_2 \to \cat{D}$, there are induced linear 
bifunctors on the strict subcategories
\[ \opn{Str}(F) : \opn{Str}(\cat{C}_{1})  \times \opn{Str}(\cat{C}_{2}) 
\to  \opn{Str}(\cat{D}) \]
and on the homotopy categories
\[ \opn{Ho}(F) : \opn{Ho}(\cat{C}_{1})  \times \opn{Ho}(\cat{C}_{2}) 
\to  \opn{Ho}(\cat{D}) . \]
\end{prop}

\begin{exer} \label{exer:5040}
Prove this proposition. 
\end{exer}

Let $F$ be a DG bifunctor as in Definition \ref{dfn:1966}. If we fix an object 
$M_1 \in \cat{C}_{1}$, then 
$F(M_1, -) :  \cat{C}_2 \to \cat{D}$
is a DG functor. Similarly if we fix the second argument of $F$. 

For DG bifunctors there are several options for contravariance. 
In Subsection \ref{subsec:contrvar-dg-func} we talked about contravariant DG 
functors, and the opposite DG category $\cat{C}^{\mrm{op}}$ of a given DG 
category $\cat{C}$.  

\begin{dfn} \label{dfn:3125}
Let $\cat{C}_1$, $\cat{C}_2$ and $\cat{D}$ be DG categories. A 
{\em DG bifunctor that is contravariant in the first or the second 
argument}%
\index{Bifunctor! contravariant}
is, by definition, a DG bifunctor
\[ F : \cat{C}_1^{\diamondsuit_1} \times \cat{C}_2^{\diamondsuit_2} 
\to \cat{D} \]
as in Definition \ref{dfn:1966}, where the symbols $\diamondsuit_1$ and 
$\diamondsuit_2$ are either $\bra{\mrm{empty}}$ or $\opn{op}$, as the case may 
be. 
\end{dfn}

Here are the two main examples of DG bifunctors. We give each of them in 
the commutative version and the noncommutative version. 

\begin{exa} \label{exa:1966}
Consider a commutative ring $A$. The category of complexes of $A$-modules is 
the DG category $\dcat{C}(A)$, and we take 
$\cat{C}_1 = \cat{C}_2 = \cat{D} := \dcat{C}(A)$. 
For each pair of objects 
$M_1, M_2 \in \dcat{C}(A)$ there is an object 
$F(M_1, M_2) := M_1 \ot_{A} M_2 \in \dcat{C}(A)$.
This is the usual tensor product of complexes. 
We define the action of $F$ on morphisms as follows: given 
$\phi_i \in \opn{Hom}_{\dcat{C}(A)}(M_i, N_i)^{k_i}$, 
we let 
\[ F(\phi_1, \phi_2) := \phi_1 \ot \phi_2 \in
\opn{Hom}_{\dcat{C}(A)} \bigl( F(M_1, M_2), F(N_1, N_2) \bigr)^{k_1 + k_2} . \]
The result is a DG bifunctor
\[ F = (- \ot_A -) : \dcat{C}(A) \times \dcat{C}(A) \to \dcat{C}(A) . \]
\end{exa}

\begin{exa} \label{exa:1965}
Consider DG rings $A_0, A_1, A_2$ (possibly noncommutative, but 
$\K$-central). 
Define the DG categories 
$\cat{C}_i := \dcat{C}(A_{i - 1} \ot A_i^{\mrm{op}})$
and
$\cat{D} := \dcat{C}(A_{0} \ot A_{2}^{\mrm{op}})$.
For each pair of objects 
$M_1 \in \cat{C}_1$ and $M_2 \in \cat{C}_2$ there is a DG $\K$-module 
$F(M_1, M_2) := M_1 \ot_{A_1} M_2$;
see Definition \ref{dfn:1100}.
This has a canonical DG $(A_{0} \ot A_{2}^{\mrm{op}})$-module structure:
\[  (a_0 \ot a_2) \cd  (m_1 \ot m_2) := (-1)^{\lms j_2 \cd (k_1 + k_2)} \cd
(a_0 \cd m_1) \ot (m_2 \cd a_2) \]
for elements 
$a_i \in A_i^{j_i}$  and $m_i \in M_i^{k_i}$. 
In this way $F(M_1, M_2)$ becomes an object of $\cat{D}$. 
We define the action of $F$ on morphisms as follows: given 
$\phi_i \in \opn{Hom}_{\cat{C}_i}(M_i, N_i)^{k_i}$,
we let 
\[ F(\phi_1, \phi_2) := \phi_1 \ot \phi_2 \in 
\opn{Hom}_{\cat{D}} \bigl( F(M_1, M_2), F(N_1, N_2) \bigr)^{k_1 + k_2} . \]
The result is a DG bifunctor
\[ F = (- \ot_{A_1} - ) : 
\dcat{C}(A_{0} \ot A_1^{\mrm{op}}) \times 
\dcat{C}(A_{1} \ot A_2^{\mrm{op}}) \to 
\dcat{C}(A_{0} \ot A_2^{\mrm{op}}) . \]
Compare this example to the one-sided construction in Example \ref{exa:1118}. 
\end{exa}

\begin{exa} \label{exa:1967}
Again we take a commutative ring $A$, but now our bifunctor $F$ arises from 
Hom, and so there is contravariance in the first argument. 
We define the DG categories 
$\cat{C}_1 = \cat{C}_2 = \cat{D} := \dcat{C}(A)$,
the same as in Example \ref{exa:1966}.
For every pair of objects 
$M_1, M_2 \in \dcat{C}(A)$ there is an object 
$F(M_1, M_2) := \opn{Hom}_{A}(M_1, M_2) \in \dcat{C}(A)$.
This is the usual Hom complex. 
We define the action of $F$ on morphisms as follows: given 
\[ \phi_1 \in \opn{Hom}_{\cat{C}_1^{\mrm{op}}}(M_1, N_1)^{k_1} =
\opn{Hom}_{\dcat{C}(A)^{\mrm{op}}}(M_1, N_1)^{k_1} =
\opn{Hom}_{A}(N_1, M_1)^{k_1} \]
and
\[ \phi_2 \in \opn{Hom}_{\cat{C}_2}(M_2, N_2)^{k_2} =
\opn{Hom}_{\dcat{C}(A)}(M_2, N_2)^{k_2} =
\opn{Hom}_{A}(M_2, N_2)^{k_2} \]
we let 
\[ \begin{aligned}
& F(\phi_1, \phi_2) := \opn{Hom}(\phi_1, \phi_2) \in
\opn{Hom}_{A} \bigl( \opn{Hom}_{A}(M_1, M_2), \opn{Hom}_{A}(N_1, N_2)
\bigr)^{k_1 + k_2}
\\
& \qquad = 
\opn{Hom}_{\cat{D}} \bigl( F(M_1, M_2), F(N_1, N_2) \bigr)^{k_1 + k_2} .
\end{aligned} \]
The result is a DG bifunctor
\[ F = \opn{Hom}_A(-, -) : \dcat{C}(A)^{\mrm{op}} \times \dcat{C}(A) \to 
\dcat{C}(A) . \]
\end{exa}

\begin{exa} \label{exa:1968}
Consider DG rings $A_0, A_1, A_2$ (possibly noncommutative, but 
$\K$-central). There is a DG bifunctor 
\[ F := \opn{Hom}_{A_1}(-, -) :
\dcat{C}(A_1 \ot A_0^{\mrm{op}})^{\mrm{op}} \times 
\dcat{C}(A_1 \ot A_2^{\mrm{op}}) \to 
\dcat{C}(A_0 \ot A_2^{\mrm{op}}) . \]
The details here are so confusing that we just leave them out.
(We shall come back to this bifunctor in Sections 
\ref{sec:perf-tilt-NC}-\ref{sec:rigid-DC-NC}, 
and then we shall have to deal with the details.)
\end{exa}

\mysubsection{Triangulated Bifunctors} \label{subsec:tri-bifun}

T-additive categories and triangulated categories were defined in Section 
\ref{sec:triangulated}. 

Suppose $(\cat{K}_1, \opn{T}_1)$ and $(\cat{K}_2, \opn{T}_2)$ are 
T-additive categories (linear over the base ring $\K$). There are two induced 
translation automorphisms on the category 
$\cat{K}_1 \times \cat{K}_2$~: 
\begin{equation} \label{eqn:4050}
\opn{T}_1(M_1, M_2) := \bigl( \opn{T}_1(M_1), M_2 \bigr) 
\end{equation}
and 
\begin{equation} \label{eqn:4051}
\opn{T}_2(M_1, M_2) := \bigl( M_1, \opn{T}_2(M_2) \bigr) . 
\end{equation}
These two functors commute: 
$\opn{T}_2 \circ \opn{T}_1 = \opn{T}_1 \circ \opn{T}_2$.

\begin{dfn} \label{dfn:1970}
Let $(\cat{K}_1, \opn{T}_1)$, $(\cat{K}_2, \opn{T}_2)$ 
and $(\cat{L}, \opn{T})$ be T-additive categories. 
A {\em T-additive bifunctor} 
\[ (F, \tau_1, \tau_2) : (\cat{K}_1, \opn{T}_1) \times (\cat{K}_2, \opn{T}_2) 
\to (\cat{L}, \opn{T}) \]
is made up of an additive bifunctor 
$F : \cat{K}_1 \times \cat{K}_2 \to \cat{L}$,
as in Definition \ref{dfn:1965}, together with isomorphisms 
$\tau_i : F \circ \opn{T}_i \iso \opn{T} \circ \, F$ 
of bifunctors 
$\cat{K}_1 \times \cat{K}_2 \to \cat{L}$. 
The condition is that 
\begin{equation*} \tag{$\dag$}
\tau_1 \circ \tau_2 = - \tau_2 \circ \tau_1 
\end{equation*}
as isomorphisms of bifunctors 
\[ F \circ \opn{T}_2 \circ \opn{T}_1 = F \circ \opn{T}_1 \circ \opn{T}_2 
\twoiso \opn{T} \circ \opn{T} \circ \, F . \]
\end{dfn}

The reason for the sign appearing in condition ($\dag$) will become clear 
in the proof of Lemma \ref{lem:4225}. 

Let $(F, \tau_1, \tau_2)$ be a T-additive bifunctor as in Definition 
\ref{dfn:1970}. If we fix an object $M_1 \in \cat{K}_{1}$, then 
\[ \bigl( F(M_1, -), \tau_2 \bigr) :  (\cat{K}_2, \opn{T}_2) \to 
(\cat{L}, \opn{T}) \]
is a T-additive functor. Similarly if we fix the second argument of $F$. 

In the next exercises we let the reader establish several operations on 
T-additive bifunctors. 

\begin{exer} \label{exer:1975}
In the situation of Definition \ref{dfn:1970}, suppose 
$(G, \tau) : (\cat{L}, \opn{T}) \to (\cat{L}', \opn{T}')$
is a T-additive functor into a fourth T-additive category 
$(\cat{L}', \opn{T}')$. Write the explicit formula for the T-additive bifunctor 
\[ (G, \tau) \circ (F, \tau_1, \tau_2) : 
(\cat{K}_1, \opn{T}_1) \times (\cat{K}_2, \opn{T}_2) 
\to (\cat{L}', \opn{T}') . \]
This should be compared to Definition \ref{dfn:1840}. 
\end{exer}

\begin{exer} \label{exer:1976}
In the situation of Definition \ref{dfn:1970}, suppose 
\[ (F', \tau'_1, \tau'_2) : (\cat{K}_1, \opn{T}_1) \times 
(\cat{K}_2, \opn{T}_2) \to (\cat{L}, \opn{T}) \]
is another T-additive bifunctor. Write the definition of a morphism of 
T-additive bifunctors 
$\ze : (F, \tau_1, \tau_2) \to (F', \tau'_1, \tau'_2)$.
Use Definition \ref{dfn:1840} as a template. 
\end{exer}

\begin{exer} \label{exer:1971}
Give a definition of a T-additive trifunctor. Show that if $F$ and $G$ are 
T-additive bifunctors, then $G(-, F(-, -))$ and $G(F(-, -), -)$
are T-additive trifunctors (whenever these compositions makes sense). 
\end{exer}

We now move to triangulated categories. 

\begin{dfn} \label{dfn:1971}
Let $(\cat{K}_1, \opn{T}_1)$, $(\cat{K}_2, \opn{T}_2)$ 
and $(\cat{L}, \opn{T})$ be triangulated categories. 
A {\em triangulated bifunctor}%
\index{Bifunctor! triangulated}
\[ (F, \tau_1, \tau_2) : (\cat{K}_1, \opn{T}_1) \times (\cat{K}_2, \opn{T}_2) 
\to (\cat{L}, \opn{T}) \]
is a T-additive bifunctor that respects the triangulated structure in each 
argument. Namely, for every distinguished triangle 
\[ L_1 \xar{\al_1} M_1 \xar{\be_1} N_1 \xar{\ga_1} \opn{T}_1(L_1) \]
in $\cat{K}_1$, and every object $L_2 \in \cat{K}_2$,
the triangle 
\[ F(L_1, L_2) \xar{F(\al_1, \opn{id})} 
F(M_1, L_2) \xar{F(\be_1, \opn{id})} 
F(N_1, L_2) \xar{\tau_1 \circ F(\ga_1, \opn{id})} 
\opn{T}(F(L_1, L_2)) \]
in $\cat{L}$ is distinguished; and the same for distinguished triangles in the 
second argument. 
\end{dfn}

The operations on triangulated bifunctors are the same as those on T-additive 
bifunctors (see exercises above). 

We now discuss triangulated bifunctors in our favorite setup:
DG modules in abelian categories. This is done both in the covariant and in 
the contravariant direction, in either argument.
Recall that the opposite homotopy category
$\dcat{K}(A, \cat{M})^{\mrm{op}}$
was given a canonical triangulated structure in Subsection 
\ref{subsec:opp-hom-triang}.

\begin{setup} \label{set:3231}
We are given this data:
\begin{itemize}
\item DG rings $A_1$ and $A_2$, and abelian categories $\cat{M}_1$ and 
$\cat{M}_2$.

\item Direction indicators $\diamondsuit_1$ and $\diamondsuit_2$, that are 
either $\bra{\mrm{empty}}$ or $\mrm{op}$. 

\item Full subcategories 
$\cat{C}_i \sub \dcat{C}(A_i, \cat{M}_i)$, whose homotopy categories 
$\cat{K}_i := \lb \opn{Ho}(\cat{C}_i) \sub \dcat{K}(A_i, \cat{M}_i)$
are such that $\cat{K}_i^{\diamondsuit_i}$ is a triangulated subcategory of 
$\dcat{K}(A_i, \cat{M}_i)^{\diamondsuit_i}$.
\end{itemize} 
\end{setup}

The translation functor of the category 
$\dcat{K}(A_i, \cat{M}_i)$ is denoted by $\opn{T}_i$, and the  
translation functor of the category 
$\dcat{K}(A_i, \cat{M}_i)^{\diamondsuit_i}$ is denoted 
by $\opn{T}_i^{\diamondsuit_i}$. 

\begin{dfn} \label{dfn:2020}
Under Setup \ref{set:3231}, let $\cat{L}$ be another triangulated 
category. A
{\em triangulated bifunctor from $(\cat{K}_1, \cat{K}_2)$ to $\cat{L}$ that is 
contravariant in the first or the second argument}%
\index{Bifunctor! contravariant}%
\index{Bifunctor! triangulated}
is, by definition, a triangulated bifunctor
\[ (F, \tau_1, \tau_2) : 
(\cat{K}_1^{\diamondsuit_1}, \opn{T}_1^{\diamondsuit_1}) \times 
(\cat{K}_2^{\diamondsuit_2}, \opn{T}_2^{\diamondsuit_2}) 
\to (\cat{L}, \opn{T}) \]
as in Definition \ref{dfn:1971}, where the contravariance is according to the 
direction indicators $\diamondsuit_1$ and $\diamondsuit_2$.
\end{dfn}

Next we explain how to obtain triangulated bifunctors from DG bifunctors. 
This is done in the following setup:

\begin{setup} \label{set:3230}
In addition to the data from Setup \ref{set:3231}, we are given: 
\begin{itemize}
\item A DG ring $B$ and an abelian category $\cat{N}$.

\item A DG bifunctor
$F : \cat{C}_1^{\diamondsuit_1} \times \cat{C}_2^{\diamondsuit_2}
\to \dcat{C}(B, \cat{N})$.
\end{itemize}
\end{setup}

The translation functor of the category 
$\dcat{K}(B, \cat{N})$ is denoted by $\opn{T}$.  
For every object
$(M_1, M_2) \in \cat{C}_1 \times \cat{C}_2$ 
there are isomorphisms
\begin{equation} \label{eqn:2015}
\tau_{i, M_1, M_2} : 
F(\opn{T}_i(M_1, M_2))  \iso \opn{T}(F(M_1, M_2))
\end{equation}
in $\dcat{C}(B, \cat{N})$, arising from Definition \tup{\ref{dfn:1150}}.
Let us make this explicit (only for $i = 2$, since the case $i = 1$ is  
similar). Fixing the object $M_1$ we obtain a DG functor
$G :  \cat{C}_2 \to \dcat{C}(B, \cat{N})$,  
$G(M_2) := F(M_1, M_2)$. 
The isomorphism 
\[ \tau_{2, M_1, M_2} : G(\opn{T}_2(M_2)) \iso \opn{T}(G(M_2)) \]
is then 
\begin{equation} \label{eqn:5042}
\tau_{2, M_1, M_2} = \opn{t}_{G(M_2)} \circ \, G(\opn{t}_{M_2})^{-1} . 
\end{equation}
We are using the little t operator here. 

\begin{lem} \label{lem:2015}
Under Setups \tup{\ref{set:3231}} and \tup{\ref{set:3230}}, assume that 
the direction indicators $\diamondsuit_1$ and $\diamondsuit_2$ are both
$\bra{\mrm{empty}}$. Fix $i \in \{ 1, 2 \}$.
Letting the object $(M_1, M_2)$ in \tup{(\ref{eqn:2015})} vary, we get an 
isomorphism 
$\tau_i : F \circ \opn{T}_i \iso \opn{T} \circ \, F$
of additive bifunctors 
$\opn{Str}(\cat{C}_{1}) \times \opn{Str}(\cat{C}_{2}) \to
\dcat{C}_{\mrm{str}}(B, \cat{N})$.
\end{lem}

\begin{proof}
This is an almost immediate consequence of the fact that the little 
t operators are morphisms of functors (see Theorem \ref{thm:1260}(2)),
\end{proof}

According to Proposition \ref{prop:5040}, the DG bifunctor $F$ induces an 
additive bifunctor 
$F : \cat{K}_{1}^{\diamondsuit_1} \times \cat{K}_{2}^{\diamondsuit_2} \to 
\dcat{K}(B, \cat{N})$ 
on the homotopy categories. 

\begin{lem} \label{lem:4225}    
Under Setups \tup{\ref{set:3231}} and \tup{\ref{set:3230}}, assume that 
the direction indicators $\diamondsuit_1$ and $\diamondsuit_2$ are both
$\bra{\mrm{empty}}$. Then 
\[ (F, \tau_1, \tau_2) :  \cat{K}_1 \times \cat{K}_2 \to 
\dcat{K}(B, \cat{N}) \]
is a triangulated bifunctor. 
\end{lem}

\begin{proof}
The only challenge is to prove that $(F, \tau_1, \tau_2)$ is a T-additive 
bifunctor; and in that, all we have to prove is that condition 
($\dag$) in Definition \ref{eqn:4051} holds, namely that 
$\tau_1 \circ \tau_2 = - \tau_2 \circ \tau_1$.
The rest hinges on single-argument considerations, that were already handled in 
Theorems \ref{thm:1150} and \ref{thm:1265}. 

So let us prove the equality ($\dag$). Choose a pair of objects 
$(M_1, M_2)$. We have diagram (\ref{eqn:5043}) in the category 
$\dcat{C}(B, \cat{N})$. 
\begin{figure}
\begin{equation} \label{eqn:5043}
\UseTips  \xymatrix @C=5ex @R=6.5ex {
&
F(\opn{T}_1(M_1), \opn{T}_2(M_2))
\ar[dr]^{F(\opn{t}_{M_1}^{-1}, \opn{id})}
\ar[dl]_{F(\opn{id}, \opn{t}_{M_2}^{-1})}
\\
F(\opn{T}_1(M_1), M_2)
\ar[d]_{\opn{t}_{F(\opn{T}_1(M_1), M_2)}}
\ar[dr]^{F(\opn{t}_{M_1}^{-1}, \opn{id})}
&
&
F(M_1, \opn{T}_2(M_2))
\ar[d]^{\opn{t}_{F(M_1, \opn{T}_2(M_2))}}
\ar[dl]_{F(\opn{id}, \opn{t}_{M_2}^{-1})}
\\
\opn{T}(F(\opn{T}_1(M_1), M_2))
\ar[d]_{\opn{T}(F(\opn{t}_{M_1}^{-1}, \opn{id}))}
&
F(M_1, M_2)
\ar[dl]^{\ \ \opn{t}_{F(M_1, M_2)}}
\ar[dr]_{\opn{t}_{F(M_1, M_2)}\ }
&
\opn{T}(F(M_1, \opn{T}_2(M_2))
\ar[d]^{\opn{T}(F(\opn{id}, \opn{t}_{M_2}^{-1}))}
\\
\opn{T}(F(M_1, M_2))
\ar[dr]^{\opn{T}(\opn{t}_{F(M_1, M_2)})}
&
&
\opn{T}(F(M_1, M_2))
\ar[dl]_{\opn{T}(\opn{t}_{F(M_1, M_2)}) \ \ }
\\
&
\opn{T}(\opn{T}(F(M_1, M_2)))
} 
\end{equation}
\end{figure}
Going from top to bottom on the left edge is the 
morphism $\tau_1 \circ \tau_2$, and going on the right edge is the 
morphism $\tau_2 \circ \tau_1$.
The bottom diamond is trivially commutative. 
The two triangles, with common vertex at $F(M_1, M_2)$, are $(-1)$-commutative, 
because $\opn{t} : \opn{Id} \to \opn{T}$ is a degree $-1$ morphism of DG 
functors. Since they occur on both sides, these signs cancel each other. 
Finally, the top diamond is $(-1)$-commutative, because
\[ (\opn{t}_{M_1}^{-1}, \opn{id}) \circ (\opn{id}, \opn{t}_{M_2}^{-1}) = 
(\opn{t}_{M_1}^{-1}, \opn{t}_{M_2}^{-1}) = 
- (\opn{id}, \opn{t}_{M_2}^{-1}) \circ (\opn{t}_{M_1}^{-1}, \opn{id}) . 
\qedhere \]
\end{proof}

\begin{thm} \label{thm:4225} 
Under Setups \tup{\ref{set:3231}} and \tup{\ref{set:3230}},
there are canonical translation isomorphisms $\tau_1$ and $\tau_2$ such that 
\[ (F, \tau_1, \tau_2) : 
(\cat{K}_1^{\diamondsuit_1}, \opn{T}_1^{\diamondsuit_1}) \times 
(\cat{K}_2^{\diamondsuit_2}, \opn{T}_2^{\diamondsuit_2})
\to (\dcat{K}(B, \cat{N}), \opn{T})  \]
is a triangulated bifunctor. 
\end{thm}

\begin{proof}
For $i = 1, 2$ let us define the indicator 
\[ \heartsuit_i :=
\begin{cases}
\mrm{flip} & \tup{if} \hspace{1ex} \diamondsuit_i = \mrm{op}, 
\\
\bra{\mrm{empty}} & \tup{if} \hspace{1ex} \diamondsuit_i =  \bra{\mrm{empty}} .
\end{cases} \]
We get DG categories $\cat{C}_i^{\heartsuit_i}$: 
if $\diamondsuit_i =  \bra{\mrm{empty}}$ then
$\cat{C}_i^{\heartsuit_i} = \cat{C}_i \sub \dcat{C}(A_i, \cat{M}_i)$,
and if $\diamondsuit_i = \mrm{op}$ then 
\[ \cat{C}_i^{\heartsuit_i} = \cat{C}_i^{\mrm{flip}} 
= \opn{Flip}(\cat{C}_i^{\mrm{op}}) \sub 
\dcat{C}(A_i, \cat{M}_i)^{\mrm{flip}} = 
\dcat{C}(A_i^{\mrm{op}}, \cat{M}_i^{\mrm{op}}) , \]
as in the proof of Theorem \ref{thm:2995}. Likewise there are triangulated 
categories $\cat{K}_i^{\heartsuit_i}$. 

We start by composing $F$ with the isomorphism of DG categories 
$\opn{Flip}^{-1}$ from Theorem \ref{thm:2495}, in each coordinate $i$ for which 
$\diamondsuit_i = \opn{op}$. This gives a new DG bifunctor 
\[ F' : \cat{C}_1^{\heartsuit_1} \times \cat{C}_2^{\heartsuit_2} \to
\dcat{C}(B, \cat{N}) . \]
Now we can apply Lemma \ref{lem:4225}, to get a triangulated bifunctor 
\[ F' : \cat{K}_1^{\heartsuit_1} \times \cat{K}_2^{\heartsuit_2} \to
\dcat{K}(B, \cat{N}) . \]
Specifically, this means that we have a pair of translations 
$(\tau'_1, \tau'_2)$ such that \lb $(F', \tau'_1, \tau'_2)$ is 
a triangulated bifunctor. 

Finally we compose $F'$ with the isomorphism of triangulated categories 
$\opn{Flip}$ from formula (\ref{eqn:3002}), in each coordinate $i$ for which 
$\diamondsuit_i = \opn{op}$. This recovers the original bifunctor $F$ on
$\cat{K}_1^{\diamondsuit_1} \times \cat{K}_2^{\diamondsuit_2}$.
There are translations $\tau_i$, that are gotten from the 
translations $\tau'_i$ like in the proof of Theorem \ref{thm:2995}. 
\end{proof}

\mysubsection{Derived Bifunctors}  \label{subsec:der-bifun}

Here we explain what are right and left derived bifunctors of triangulated 
bifunctors. The definitions and results are very similar to the single-argument 
case. For the sake of simplicity, we shall mostly ignore the translation 
functors and the translation isomorphisms; enough was said about them in the 
previous subsection. 

The next setup will be used throughout this subsection.

\begin{setup} \label{set:1976}
The following are given: 
\begin{itemize}
\item Triangulated categories $\cat{K}_1$, $\cat{K}_2$ and $\cat{E}$.

\item A triangulated bifunctor 
$F : \cat{K}_1 \times \cat{K}_2 \to \cat{E}$.

\item Denominator sets of cohomological origin
$\cat{S}_1 \sub \cat{K}_1$ and $\cat{S}_2 \sub \cat{K}_2$.
\end{itemize}
We write $\cat{D}_i := (\cat{K}_i)_{\cat{S}_i}$, 
and $\opn{Q}_i : \cat{K}_i \to \cat{D}_i$ are the localization functors.
\end{setup}

The localized category 
$\cat{D}_i := (\cat{K}_i)_{\cat{S}_i}$
is triangulated, and the localization functor 
$\opn{Q}_i : \cat{K}_i \to \cat{D}_i$
is triangulated. On the product categories we get a functor 
\[ \opn{Q}_1 \times \opn{Q}_2 : 
\cat{K}_1 \times \cat{K}_2 \to \cat{D}_1 \times \cat{D}_2 . \]

Before embarking on the definitions of the derived bifunctors, we want to 
introduce the relevant categories of functors. These are the bifunctor variants 
of what we had in Subsection \ref{subsec:fun-cats}. 

\begin{dfn} \label{dfn:3130} 
We denote by 
$\cat{LinBiFun}(\cat{K}_1 \times \cat{K}_2, \cat{E})$
the category whose objects are the linear bifunctors 
$F : \cat{K}_1 \times \cat{K}_2 \to \cat{E}$.
The morphisms are the obvious ones. 
\end{dfn}

This is a linear category. Actually, this category is the same as the category 
$\cat{LinFun}(\cat{K}_1 \ot \cat{K}_2, \cat{E})$ 
of linear functors
$F : \cat{K}_1 \ot \cat{K}_2 \to \cat{E}$, 
only with the extra information that the source is a product category.
     
\begin{dfn} \label{dfn:2036}
We denote by 
$\cat{TrBiFun}(\cat{K}_1 \times \cat{K}_2, \cat{E})$
the category whose objects are the $\K$-linear triangulated bifunctors 
$(F, \tau_1, \tau_2) : \cat{K}_1 \times \cat{K}_2 \to \cat{E}$.
The morphisms are those of T-additive bifunctors. 
\end{dfn}

There is a functor 
\begin{equation} \label{eqn:3130}
\cat{TrBiFun}(\cat{K}_1 \times \cat{K}_2, \cat{E}) \to 
\cat{LinBiFun}(\cat{K}_1 \times \cat{K}_2, \cat{E}) 
\end{equation}
that forgets $(\tau_1, \tau_2)$. It is a faithful additive functor. 

Suppose $U_i : \cat{K}'_i \to \cat{K}_i$ are triangulated functors between 
triangulated categories. We get an induced additive functor
\begin{equation} \label{eqn:2036}
\cat{Fun}(U_1 \times U_2, \opn{Id}_{\cat{E}}) :
\cat{TrBiFun}(\cat{K}_1 \times \cat{K}_2, \cat{E}) \to 
\cat{TrBiFun}(\cat{K}'_1 \times \cat{K}'_2, \cat{E})
\end{equation}
The formula is 
$F \mapsto F \circ (U_1 \times U_2)$.

\begin{lem} \label{lem:2035}
If the functors $U_1$ and $U_2$ are equivalences, then the functor 
in \tup{(\ref{eqn:2036})} is an equivalence.
\end{lem}

\begin{proof}
This is basically the same as the proof of Lemma \ref{lem:2335}.
\end{proof}

As in Definition \ref{dfn:2037} we denote by 
\[ \cat{TrBiFun}_{\cat{S}_1 \times \cat{S}_2}
(\cat{K}_1 \times \cat{K}_2, \cat{E}) \sub 
\cat{TrBiFun}(\cat{K}_1 \times \cat{K}_2, \cat{E}) \]
the full subcategory on the triangulated bifunctors $F$ such that 
$F(\cat{S}_1 \times \cat{S}_2) \sub \cat{E}^{\times}$. 

\begin{lem} \label{lem:2036} 
The functor 
\[ \cat{Fun}(\opn{Q}_1 \times \opn{Q}_2, \opn{Id}_{\cat{E}}) :
\cat{TrBiFun}(\cat{D}_1 \times \cat{D}_2, \cat{E}) \to 
\cat{TrBiFun}_{\cat{S}_1 \times \cat{S}_2} (\cat{K}_1 \times \cat{K}_2, \cat{E})
\]
is an isomorphism of categories. 
\end{lem}

\begin{proof}
This is basically the same as the proof of Lemma \ref{lem:2336}, when 
combined with the isomorphism of triangulated categories 
$(\cat{K}_1 \times \cat{K}_2)_{\cat{S}_1 \times \cat{S}_2} \to 
\cat{D}_1 \times \cat{D}_2$ 
from Proposition \ref{prop:2028}. 
\end{proof}

We are ready to talk about derived bifunctors.  

\begin{dfn} \label{dfn:1975}
Under Setup \ref{set:1976}, 
a {\em triangulated right derived bifunctor of $F$ with respect to 
$\cat{S}_1 \times \cat{S}_2$}%
\index{Bifunctor! derived}
is a pair $(\mrm{R} F, \eta^{\mrm{R}})$, where 
\[ \mrm{R} F : \cat{D}_1 \times \cat{D}_2 \to \cat{E} \]
is a triangulated bifunctor, and 
\[ \eta^{\mrm{R}}: F \twoto \mrm{R} F \circ (\opn{Q}_1 \times \opn{Q}_2) \] 
is a morphism of triangulated bifunctors, 
such that the following universal property holds:
\begin{itemize}
\item[(R)] Given a pair $(G, \th)$, consisting of a triangulated 
bifunctor
$G : \cat{D}_1 \times \cat{D}_2 \to \cat{E}$ 
and a morphism of triangulated bifunctors 
$\th : F \twoto G \circ (\opn{Q}_1 \times \opn{Q}_2)$,
there is a unique morphism of triangulated bifunctors 
$\mu : \mrm{R}^{} F \twoto G$
such that 
$\th = (\mu \circ \opn{id}_{\opn{Q}_1 \times \opn{Q}_2}) * \eta^{\mrm{R}}$.
\end{itemize}
\end{dfn}

Here is a diagram showing property (R):
\begin{equation} \label{eqn:2076}
\UseTips  \xymatrix @C=14ex @R=10ex {
\cat{K}_1 \times \cat{K}_2
\ar[r]^{F} _(0.45){}="q"  ^(0.55){}="q2" 
\ar[d]_{\opn{Q}_1 \times \opn{Q}_2}
& 
\cat{E} 
\\
\cat{D}_1 \times \cat{D}_2
\ar[ur]^(0.32){\mrm{R} F} ^(0.45){}="f" ^(0.5){}="f1" 
\ar@{=>}  "q";"f" _{\eta^{\mrm{R}}}
\ar@{->}@(r,d)[ru]_{G} ^(0.55){}="g"  ^(0.48){}="g1"
\ar@{=>}  "q2";"g" ^(0.7){\th}
\ar@{=>}  "f";"g1" _(0.45){\mu}
} 
\end{equation}
 
\begin{prop} \label{prop:2092}
If a triangulated right derived bifunctor exists, then it is unique up to a 
unique isomorphism. 
\end{prop}

\begin{proof}
This is just like the proof of Proposition \ref{prop:1350}. We leave the small 
changes to the reader. 
\end{proof}

Existence in general is like Theorem \ref{thm:1470}, but a bit more complicated.

\begin{thm} \label{thm:2000}
Under Setup \tup{\ref{set:1976}}, 
assume there are full triangulated subcategories 
$\cat{J}_1 \subseteq \cat{K}_1$ and 
$\cat{J}_2 \subseteq \cat{K}_2$
with these two properties\tup{:}
\begin{enumerate}
\rmitem{a} 
If $\phi_1 : I_1 \to J_1$ is a morphism in $\cat{J}_1 \cap \cat{S}_1$
and $\phi_2 : I_2 \to J_2$ is a morphism in $\cat{J}_2 \cap \cat{S}_2$, 
then 
\[ F(\phi_1, \phi_2) : F(I_1, I_2) \to F(J_1, J_2) \]
is an isomorphism in $\cat{E}$.
 
\rmitem{b} For every $p \in \{1, 2\}$ and every object $M_p \in \cat{K}_p$,
there exists a morphism $\rho_p : M_p \to I_p$ in $\cat{S}_p$ with target 
$I_p \in \cat{J}_p$. 
\end{enumerate}
Then the triangulated right derived bifunctor%
\index{Bifunctor! derived}
\[ (\mrm{R} F, \eta^{\mrm{R}}) : \cat{D}_1 \times \cat{D}_2 \to \cat{E} \]
of $F$ with respect to $\cat{S}_1 \times \cat{S}_2$
exists. Moreover, for every pair of objects $I_1 \in \cat{J}_1$ and 
$I_2 \in \cat{J}_2$ the morphism  
\[ \eta^{\mrm{R}}_{I_1, I_2} : F(I_1, I_2) \to \mrm{R} F (I_1, I_2) \]
in $\cat{E}$ is an isomorphism.
\end{thm}

In applications we will see that either $\cat{J}_1 = \cat{K}_1$ or 
$\cat{J}_2 = \cat{K}_2$; namely we will only need to resolve in the second or 
in the first argument, respectively. 

\begin{proof}
It will be convenient to change notation.
For $p = 1, 2$ let's define 
$\cat{K}'_p := \cat{J}_p$, $\cat{S}'_p := \cat{K}'_p \cap \cat{S}_p$ and
$\cat{D}'_p := (\cat{K}'_p)_{\cat{S}'_p}$. 
The localization functors are 
$\opn{Q}'_p : \cat{K}'_p \to \cat{D}'_p$. 
The inclusions are 
$U_p : \cat{K}'_p \to \cat{K}_p$, and their localizations are the functors 
$V_p : \cat{D}'_p \to \cat{D}_p$.
By Lemma \ref{lem:2086} the functors $V_p$ are equivalences.

The situation is depicted in these diagrams. 
We have this commutative diagram of 
products of triangulated functors between products of triangulated 
categories:
\begin{equation} \label{eqn:2050}
\UseTips \xymatrix @C=10ex @R=6ex {
\cat{K}'_1 \times \cat{K}'_2
\ar[r]^{U_1 \times U_2}
\ar[d]_{\opn{Q}'_1 \times \opn{Q}'_2}
&
\cat{K}_1 \times \cat{K}_2
\ar[d]^{\opn{Q}_1 \times \opn{Q}_2}
\\
\cat{D}'_1 \times \cat{D}'_2
\ar[r]^{V_1 \times V_2}
&
\cat{D}_1 \times \cat{D}_2
} 
\end{equation}

\noindent 
The arrow $V_1 \times V_2$ is an equivalence. 
Diagram (\ref{eqn:2050}) induces a commutative diagram of linear categories:
\begin{equation} \label{eqn:2051}
\UseTips \xymatrix @C=16ex @R=5.5ex {
\cat{TrBiFun}(\cat{K}'_1 \times \cat{K}'_2, \cat{E})
&
\cat{TrBiFun}(\cat{K}_1 \times \cat{K}_2, \cat{E})
\ar[l]_{\cat{Fun}(U_1 \times U_2, \opn{Id})}
\\
\cat{TrBiFun}_{\cat{S}'_1 \times \cat{S}'_2}
(\cat{K}'_1 \times \cat{K}'_2, \cat{E})
\ar[u]_{\mrm{f.f.\, emb}}
&
\cat{TrBiFun}_{\cat{S}_1 \times \cat{S}_2}
(\cat{K}_1 \times \cat{K}_2, \cat{E})
\ar[u]^{\mrm{f.f.\, emb}}
\ar[l]^{\mrm{equiv}}_{\cat{Fun}(U_1 \times U_2, \opn{Id})}
\\
\cat{TrBiFun}(\cat{D}'_1 \times \cat{D}'_2, \cat{E})
\ar[u]^(0.47){\cat{Fun}(\opn{Q}'_1 \times \opn{Q}'_2, \opn{Id})}_{\mrm{isom}}
&
\cat{TrBiFun}(\cat{D}_1 \times \cat{D}_2, \cat{E})
\ar[l]_{\cat{Fun}(V_1 \times V_2, \opn{Id})}^{\mrm{equiv}}
\ar[u]_(0.47){\cat{Fun}(\opn{Q}_1 \times \opn{Q}_2, \opn{Id})}^{\mrm{isom}}
} 
\end{equation}
According to Lemmas \ref{lem:2035} and \ref{lem:2036}, the arrows in the 
diagram above that are marked ``isom'' or ``equiv'' are isomorphisms or 
equivalences, respectively. By definition the arrows marked ``f.f.\ emb'' are 
fully faithful embeddings. 

We know that $\cat{S}_p \sub \cat{K}_p$ are left denominator sets. Therefore 
(see Proposition \ref{prop:2100})
$\cat{S}_1 \times \cat{S}_2 \sub \cat{K}_1 \times \cat{K}_2$
is a  left denominator set.
Condition (a) of Theorem \ref{thm:2000} says that $F$ sends morphisms in 
$\cat{S}'_1 \times \cat{S}'_2$ to isomorphisms in $\cat{E}$. 
Condition (b) there says that there are enough right 
$(\cat{K}'_1 \times \cat{K}'_2)$-resolutions in 
$\cat{K}_1 \times \cat{K}_2$.

Thus we are in a position to use the abstract Theorem \ref{thm:2080}.
It says that there is an abstract right derived functor 
$(\mrm{R} F, \eta^{\mrm{R}}) :  \cat{D}_1 \times \cat{D}_2 \to \cat{E}$.
However, going over the proof of Theorem \ref{thm:2080},
we see that all constructions there can be made within the triangulated 
setting, namely in diagram (\ref{eqn:2051}) instead of in  diagram 
(\ref{eqn:2084}). Therefore $\mrm{R} F$ is an object of the category in the 
bottom right corner of (\ref{eqn:2051}), and the morphism 
$\eta^{\mrm{R}} : F \twoto \mrm{R} F \circ (\opn{Q}_1 \times \opn{Q}_2)$
is in the category in the top right corner of (\ref{eqn:2051}).
\end{proof}

We now treat left derived bifunctors. Due to the similarity to the right 
derived functors, we give only a few details here. 

\begin{dfn} \label{dfn:2105}
Under Setup \ref{set:1976}, 
a {\em triangulated left derived bifunctor of $F$ with respect to 
$\cat{S}_1 \times \cat{S}_2$}%
\index{Bifunctor! derived}
is a pair $(\mrm{L} F, \eta^{\mrm{L}})$, where 
\[ \mrm{L} F : \cat{D}_1 \times \cat{D}_2 \to \cat{E} \]
is a triangulated bifunctor, and 
\[ \eta^{\mrm{L}} :  \mrm{L} F \circ (\opn{Q}_1 \times \opn{Q}_2) \twoto F \] 
is a morphism of triangulated bifunctors, 
such that the following universal property holds:
\begin{itemize}
\item[(L)] Given a pair $(G, \th)$, consisting of a triangulated 
bifunctor
$G : \cat{D}_1 \times \cat{D}_2 \to \cat{E}$ 
and a morphism of triangulated bifunctors 
$\th :  G \circ (\opn{Q}_1 \times \opn{Q}_2) \twoto F$,
there is a unique morphism of triangulated bifunctors 
$\mu : G \twoto \mrm{L}^{} F$
such that 
$\th =  \eta^{\mrm{L}} * (\mu \circ \opn{id}_{\opn{Q}_1 \times \opn{Q}_2})$.
\end{itemize}
\end{dfn}

\begin{prop} \label{prop:2105}
If a triangulated left derived bifunctor exists, then it is unique up to a 
unique isomorphism. 
\end{prop}

\begin{proof}
This is the opposite of Proposition \ref{prop:2092}, and we leave it to the 
reader to make the adjustments.
\end{proof}

\begin{thm} \label{thm:2106}
In the situation of Definition \tup{\ref{dfn:2105}}, 
assume there are full triangulated subcategories 
$\cat{P}_1 \subseteq \cat{K}_1$ and 
$\cat{P}_2 \subseteq \cat{K}_2$
with these two properties\tup{:}
\begin{enumerate}
\rmitem{a} If $\phi_1 : P_1 \to Q_1$ is a
morphism in $\cat{P}_1 \cap \cat{S}_1$
and $\phi_2 : P_2 \to Q_2$ is a 
morphism in $\cat{P}_2 \cap \cat{S}_2$, then 
\[ F(\phi_1, \phi_2) : F(P_1, P_2) \to F(Q_1, Q_2) \]
is an isomorphism in $\cat{E}$.
 
\rmitem{b} 
For every $i \in \{1, 2\}$ and every object $M_i \in \cat{K}_i$,
there exists a morphism $\rho_i : P_i \to M_i$ in $\cat{S}_i$ with target 
$P_i \in \cat{P}_i$. 
\end{enumerate}
Then the triangulated left derived bifunctor%
\index{Bifunctor! derived}
\[ (\mrm{L} F, \eta^{\mrm{L}}) : \cat{D}_1 \times \cat{D}_2 \to \cat{E} \]
of $F$ with respect to $\cat{S}_1 \times \cat{S}_2$
exists. Moreover, for every pair of objects $P_1 \in \cat{P}_1$ and 
$P_2 \in \cat{P}_2$ the morphism  
\[ \eta^{\mrm{L}}_{P_1, P_2} :  \mrm{L} F(P_1, P_2) \to F(P_1, P_2) \]
in $\cat{E}$ is an isomorphism.
\end{thm}

\begin{exer} \label{exer:4635}
Prove Theorem \ref{thm:2106}.
(Hint: This is the opposite of Theorem \ref{thm:2000}. The proof is analogous.)
\end{exer}

In applications we will see that either $\cat{P}_1 = \cat{K}_1$ or 
$\cat{P}_2 = \cat{K}_2$; namely we will only need to resolve in the second or 
in the first argument, respectively. 

There are more specialized definitions of triangulated derived bifunctors 
whose sources are categories of DG modules, like Definitions \ref{dfn:3236} and 
\ref{dfn:3237}. We leave them to the reader. 

\begin{rem} \label{5050}
Deriving triangulated bifunctors that are contravariant in one or two arguments 
involves passage to opposite categories; see Definition \ref{dfn:2020}. 
Proposition \ref{prop:2356} tells us that when 
$\cat{K}_i^{\mrm{op}} \sub \dcat{K}(A_i, \cat{M}_i)^{\mrm{op}}$
are full triangulated subcategories, the sets 
$\cat{S}_i^{\mrm{op}}$
of quasi-isomorphisms in them are of cohomological origin. So we are still 
within the conditions of Setup \ref{set:1976}, and Theorems \ref{thm:2000} and 
\ref{thm:2106} can be applied. 
We leave it to the reader to work out the details in the various cases. 
However, the important case of the bifunctor $\opn{Hom}(-, -)$ 
will be given a full treatment in Subsection \ref{subsec:RHom}.
\end{rem}

\cleardoublepage
\mysection{Resolving Subcategories of 
\texorpdfstring{$\dcat{K}(A, \cat{M})$}{K(A,M)}} 
\label{sec:resol} 

\AYcopyright

In this section we are back in the more concrete setting: $A$ is a DG ring, 
and $\cat{M}$ is an abelian category, both over the commutative base ring 
$\K$; see Convention \ref{conv:2490}. We will define {\em K-projective} and 
{\em K-injective} DG modules in
$\bcat{K}(A, \cat{M})$. These DG modules form full triangulated 
subcategories of $\bcat{K}(A, \cat{M})$, and are concrete versions of the 
abstract categories $\cat{P}$ and $\cat{J}$, that played important roles in 
Subsection \ref{subsec:tri-der-funs}.
For $\bcat{K}(A)$ we also define {\em K-flat DG modules}.

\mysubsection{K-Injective DG Modules} \label{subsec:K-inj}

For every integer $p$ we have the $p$-th cohomology functor 
$\opn{H}^p : \bcat{C}_{\mrm{str}}(A, \cat{M}) \to \cat{M}$.
There is equality 
$\opn{H}^p = \opn{H}^0 \circ \opn{T}^p$,
where $\opn{T}$ is the translation functor of $\bcat{C}(A, \cat{M})$. 
The functors $\opn{H}^p$ pass to the homotopy category, and 
$\opn{H}^0 : \bcat{K}(A, \cat{M}) \to \cat{M}$
is a cohomological functor, in the sense of Definition \ref{dfn:1510}.

\begin{dfn} \label{dfn:1505}
A DG module $N \in \dcat{C}(A, \cat{M})$ is called {\em acyclic}%
\index{Differential graded module! acyclic}
if $\mrm{H}^p(N) = 0$ for all $p$. 
\end{dfn}

\begin{dfn} \label{dfn:1506}
A DG module $I \in \dcat{C}(A, \cat{M})$ is called 
{\em K-injective}%
\index{Differential graded module! K-injective}
if for every acyclic DG module $N \in \dcat{C}(A, \cat{M})$, the DG $\K$-module 
$\opn{Hom}_{A, \cat{M}}(N, I)$ is acyclic.
\end{dfn}

The definition above characterizes K-injectives as objects of 
$\dcat{C}(A, \cat{M})$. 
The next proposition shows that being K-injective is intrinsic to the 
triangulated category $\dcat{K}(A, \cat{M})$, with the cohomological functor 
$\opn{H}^0$ (that tells us which are the acyclic objects). 

\begin{prop} \label{prop:1515}
A DG module $I \in \dcat{K}^{}(A, \cat{M})$ is K-injective if and only if 
$\opn{Hom}_{\dcat{K}(A, \cat{M})}(N, I) = 0$ for every acyclic DG module 
$N \in \dcat{K}^{}(A, \cat{M})$.
\end{prop}

\begin{proof}
This is because for every integer $p$ we have 
\[ \opn{H}^p \bigl( \opn{Hom}_{A, \cat{M}}(N, I) \bigr) \cong
\opn{H}^0 \bigl( \opn{Hom}_{A, \cat{M}}(\opn{T}^{-p}(N), I) \bigr) 
= \opn{Hom}_{\dcat{K}(A, \cat{M})}(\opn{T}^{-p}(N), I) , \]
and $N$ is acyclic iff $\opn{T}^{-p}(N)$ is acyclic.
\end{proof}

\begin{dfn} \label{dfn:1507}
Let $M \in \dcat{C}(A, \cat{M})$. A 
{\em K-injective resolution}%
\index{Resolution! K-injective}
of $M$ is a quasi-isomorphism $\rho : M \to I$ in 
$\dcat{C}_{\mrm{str}}(A, \cat{M})$, where 
$I$ is a K-injective DG module. 
\end{dfn}

\begin{rem} \label{rem:5051}
Let $\rho : M \to I$ be a K-injective resolution in 
$\dcat{C}_{\mrm{str}}(A, \cat{M})$. 
By a slight abuse of notation, we shall often refer to the induced 
quasi-isomorphism 
$\opn{P}(\rho) : M \to I$ in $\dcat{K}(A, \cat{M})$
as a K-injective resolution of $M$. Cf.\ Definitions \ref{dfn:1508} and
\ref{dfn:1512} below. 
\end{rem}

In the next section we will prove existence of K-injective resolutions in 
several contexts. Here is an easy example of a K-injective complex. 

\begin{exer} \label{exer:1505}
Let $I \in \dcat{K}(\cat{M})$ be a complex of injective objects of 
$\cat{M}$, with zero differential. Prove that $I$ is K-injective. 
\end{exer}

\begin{dfn} \label{dfn:1509}
Let $\cat{K}$ be a full subcategory of $\dcat{K}^{}(A, \cat{M})$. 
The full subcategory of  $\cat{K}$ on the K-injective DG modules in it is 
denoted by $\cat{K}_{\mrm{inj}}$. In other words, 
$\cat{K}_{\mrm{inj}} = \dcat{K}^{}(A, \cat{M})_{\mrm{inj}} \cap \cat{K}$.
\end{dfn}

\begin{rem} \label{rem:5050}
Warning: the property of a DG module $I \in \cat{K}$ being K-injective is in 
general not intrinsic to the subcategory 
$\cat{K} \sub \dcat{K}(A, \cat{M})$. This is because the 
condition in Proposition \ref{prop:1515} (and in Definition \ref{dfn:1506}) 
has to be tested against all acyclic DG modules 
$N \in \dcat{K}(A, \cat{M})$.
\end{rem}

\begin{prop}  \label{prop:130}
If  $\cat{K}$ is a full triangulated subcategory of
$\dcat{K}^{}(A, \cat{M})$, then $\cat{K}_{\mrm{inj}}$ is a full triangulated 
subcategory of $\cat{K}$.
\end{prop}

\begin{proof}
It suffices to prove that 
$\dcat{K}^{}(A, \cat{M})_{\mrm{inj}}$ is a 
triangulated subcategory of $\dcat{K}^{}(A, \cat{M})$.
It is easy to see that $\dcat{K}^{}(A, \cat{M})_{\tup{inj}}$ is closed under
translations. Suppose 
$I \to J \to K \xar{\, \triangle\, }$
is a distinguished triangle in $\dcat{K}^{}(A, \cat{M})$ such that  
$I, J$ are K-injective DG modules. We have to show that $K$ is also
K-injective. Take any acyclic DG module $N \in \dcat{K}^{}(A, \cat{M})$. 
There is an exact sequence
\[ \opn{Hom}_{\dcat{K}(A, \cat{M})}(N, J) \to 
\opn{Hom}_{\dcat{K}(A, \cat{M})}(N, K)
\to \opn{Hom}_{\dcat{K}(A, \cat{M})}(N, \opn{T}(I))  \]
in $\cat{Mod} \K$. Because $J$ and $\opn{T}(I)$ are 
K-injectives, Proposition \ref{prop:1515} says that
\[ \opn{Hom}_{\dcat{K}(A, \cat{M})}(N, J) = 0 = 
\opn{Hom}_{\dcat{K}(A, \cat{M})}(N, \opn{T}(I)) . \]
Therefore 
$\opn{Hom}_{\dcat{K}(A, \cat{M})}(N, K) = 0$. 
But $N$ is an arbitrary acyclic DG module, so $K$ is K-injective.
\end{proof}

\begin{exa} \label{exa:1520}
Let $\star$ be some boundedness condition (namely $\mrm{b}$, $+$ or $-$). 
We know that $\dcat{K}^{\star}(A, \cat{M})$ 
is a full triangulated subcategory of $\dcat{K}^{}(A, \cat{M})$.
Hence $\dcat{K}^{\star}(A, \cat{M})_{\mrm{inj}}$
is a triangulated subcategory too.  
\end{exa}

\begin{dfn} \label{dfn:1508}
Let $\cat{K}$ be a full triangulated subcategory of 
$\dcat{K}(A, \cat{M})$. We say that $\cat{K}$ 
{\em has enough K-injectives}%
\index{Resolution! K-injective}
if every DG module $M \in \cat{K}$ admits a K-injective resolution inside 
$\cat{K}$. I.e.\ there is a quasi-isomorphism $\rho : M \to I$
where $I \in \cat{K}_{\mrm{inj}}$.
\end{dfn}

Here is the crucial fact regarding K-injectives. 

\begin{lem} \label{lem:1510}
Let  $s : I \to M$ be a quasi-isomorphism in $\dcat{K}(A, \cat{M})$,
and assume $I$ is K-injective. Then $s$ has a left inverse, namely 
there is a morphism $t : M \to I$ in $\dcat{K}(A, \cat{M})$ such that 
$t \circ s = \opn{id}_{I}$. 
Moreover, this morphism $t$ is unique. 
\end{lem}

\begin{proof}
Consider a distinguished triangle 
$I \xar{s} M \to N \xar{\, \triangle \, }$
in $\dcat{K}(A, \cat{M})$ that's built on $s$. The long exact 
cohomology sequence tells us that $N$ is an acyclic DG module. 
Therefore
$\opn{Hom}_{\dcat{K}^{}(A, \cat{M})}(\opn{T}^p(N), I) = 0$ 
for all $p$. 
The exact sequence 
\[ \begin{aligned}
& \opn{Hom}_{\dcat{K}(A, \cat{M})}(N, I) \to 
\opn{Hom}_{\dcat{K}(A, \cat{M})}(M, I)
\\
& \quad 
\to \opn{Hom}_{\dcat{K}(A, \cat{M})}(I, I) 
\to \opn{Hom}_{\dcat{K}(A, \cat{M})}(\opn{T}^{-1}(N), I) 
\end{aligned} \]
shows that the homomorphism
\[ (-) \circ s : \opn{Hom}_{\dcat{K}(A, \cat{M})}(M, I) \iso 
\opn{Hom}_{\dcat{K}(A, \cat{M})}(I, I) \]
is a bijection.
We take $t : M \to I$ to be the unique morphism in $\dcat{K}(A, \cat{M})$
such that $t \circ s = \opn{id}_I$. 
\end{proof}

\begin{thm} \label{thm:3135}
\index{Differential graded module! K-injective}
Let $A$ be a DG ring, let $\cat{M}$ be an abelian category, and let
$\cat{K}$ be a full triangulated subcategory of $\dcat{K}(A, \cat{M})$.
Write $\cat{S} := \cat{K} \cap \, \dcat{S}(A, \cat{M})$,
and let $\cat{D} := \cat{K}_{\cat{S}}$ be the localization, with localization 
functor $\opn{Q} : \cat{K} \to \cat{D}$. 
Then for every $M \in \cat{K}$ and $I \in \cat{K}_{\mrm{inj}}$
the homomorphism  
\[ \opn{Q}_{M, I} : \opn{Hom}_{\cat{K}}(M, I) \to \opn{Hom}_{\cat{D}}(M, I) \]
is bijective.  
\end{thm}

\begin{proof}
We know that $\opn{Q} : \cat{K} \to \cat{D}$ is a left Ore localization with 
respect to $\cat{S}$. 
Suppose $q : M \to I$ is a morphism in $\cat{D} = \cat{K}_{\cat{S}}$.
Let us present $q$ as a left fraction: 
$q = \opn{Q}(s)^{-1} \circ \opn{Q}(a)$, 
where $a : M \to N$ and $s : I \to N$ are morphisms in 
$\cat{K}$, and $s$ is a quasi-isomorphism. 
By Lemma \ref{lem:1510}, $s$ has a left inverse 
$t : M \to I$ in $\cat{K}$. We get a morphism 
$t \circ a : M \to I$ in $\cat{K}$, and an easy calculation shows that
$\opn{Q}(t \circ a) = q$ in $\cat{D}$.
This proves surjectivity of the function $\opn{Q}_{M, I}$.

Now let's prove injectivity of $\opn{Q}_{M, I}$.  
If $a : M \to I$ is a morphism in $\cat{K}$ such that 
$\opn{Q}_{M, I}(a) = 0$, then by axiom (LO4) of left Ore localization (see 
Definition \ref{dfn:2310}) there is a quasi-isomorphism $s : I \to L$ in 
$\cat{K}$ such that $s \circ a = 0$ in $\cat{K}$. Let $t$ be the left inverse 
of 
$s$. Then $a = t \circ s \circ a= 0$ in $\cat{K}$. 
\end{proof}

\begin{cor} \label{cor:3135}
In the situation of Theorem \tup{\ref{thm:3135}}, 
the localization functor 
$\opn{Q} : \cat{K}_{\mrm{inj}} \to \cat{D}$ 
is fully faithful.
\end{cor}

\begin{proof}
Since $\cat{K}_{\mrm{inj}}$ is a full subcategory of 
$\cat{K}$, for every $J, I \in \cat{K}_{\mrm{inj}}$
we have bijections
\[ \opn{Hom}_{\cat{K}_{\mrm{inj}}}(J, I) \iso \opn{Hom}_{\cat{K}}(J, I) 
\xar{\opn{Q}_{J, I}} \opn{Hom}_{\cat{D}}(J, I) , \]
where the bijection $\opn{Q}_{J, I}$ is by the theorem above. 
\end{proof}

\begin{cor} \label{cor:1515}
In the situation of Theorem \tup{\ref{thm:3135}}, if
$\cat{K}$ has enough K-injectives, then the localization functor
$\opn{Q} : \cat{K}_{\mrm{inj}} \to \cat{D}$
is an equivalence.
\end{cor}

\begin{proof}
By Corollary \ref{cor:3135} the functor $\opn{Q}$ is fully faithful. The 
extra condition guarantees that $\opn{Q}$ is essentially surjective on objects. 
\end{proof}

\begin{cor} \label{cor:1509}
Let $\star$ be any boundedness condition. If $\dcat{K}^{\star}(A, \cat{M})$ has 
enough K-injectives, then the triangulated functor
\[ \opn{Q} : \dcat{K}^{\star}(A, \cat{M})_{\mrm{inj}} \to
\dcat{D}^{\star}(A, \cat{M}) \]
is an equivalence.
\end{cor}

\begin{proof}
Since $\dcat{K}^{\star}(A, \cat{M})$ is a full 
triangulated subcategory of $\dcat{K}(A, \cat{M})$,
this is a special case of the Corollary \ref{cor:1515}.
\end{proof}

\begin{rem} \label{rem:1511}
These last results are of tremendous importance, both theoretically and 
practically. In the theory, Corollary \ref{cor:1509}
shows that the localized category 
$\dcat{D}^{\star}(A, \cat{M})$, 
which is too big to lie inside the original universe $\cat{U}$ (see Remark 
\ref{rem:1405}), is equivalent to a  $\cat{U}$-category. 
On the practical side, Corollary \ref{cor:3135} means that among K-injective 
objects we do not need fractions to represent morphisms. 
\end{rem}

\begin{cor} \label{cor:1511}
Let $\cat{K} \sub \cat{K}'$ be full triangulated subcategories of 
$\dcat{K}^{}(A, \cat{M})$.
Define 
$\cat{S} := \cat{K} \cap \, \dcat{S}(A, \cat{M})$,
$\cat{S}' := \cat{K}' \cap \, \dcat{S}(A, \cat{M})$,
$\cat{D} := \cat{K}_{\cat{S}}$
and
$\cat{D}' := \cat{K}'_{\cat{S}'}$.
If $\cat{K}$ and $\cat{K}'$ have enough K-injectives, then the canonical 
functor $\cat{D} \to \cat{D}'$ is fully faithful. 
\end{cor}

\begin{proof}
Combine Corollary \ref{cor:1515} with the fact that the canonical functor
$\cat{K} \to \cat{K}'$ is fully faithful.
\end{proof}

\begin{cor} \label{cor:1922}
Let $\phi : I \to J$ be a morphism in $\dcat{C}_{\mrm{str}}(A, \cat{M})$ 
between K-injective objects. Then $\phi$ is a homotopy equivalence if and only 
if it is a quasi-isomorphism. 
\end{cor}

\begin{proof}
One implication is trivial. For the reverse implication, if $\phi$ is a 
quasi-isomorphism then it is an isomorphism in 
$\dcat{D}(A, \cat{M})$, and by Corollary \ref{cor:3135} for 
$\cat{K} = \dcat{K}(A, \cat{M})$
we see that $\phi$ is an isomorphism in $\dcat{K}(A, \cat{M})$. 
\end{proof}
   
\begin{thm}  \label{thm:5050}   
\index{Derived functor! triangulated right}
In the situation of Theorem \tup{\ref{thm:3135}}, assume that $\cat{K}$ has 
enough K-injectives. Let $\cat{E}$ be a triangulated category, and let 
$F : \cat{K} \to \cat{E}$
be a triangulated functor. Then $F$ has a right derived functor 
$(\mrm{R} F, \eta^{\mrm{R}}) : \cat{D} \to \cat{E}$. 
Furthermore, for every $I \in \cat{K}_{\mrm{inj}}$ the morphism 
$\eta^{\mrm{R}}_I : F(I) \to \mrm{R} F(I)$
in $\cat{E}$ is an isomorphism. 
\end{thm}

\begin{proof}
We will use Theorem \ref{thm:1470}. In the notation of that theorem, let 
$\cat{J} := \cat{K}_{\mrm{inj}}$.
Condition (a) of that theorem holds (this is the ``enough K-injectives'' 
assumption). Next, Corollary \ref{cor:1922} implies that every 
quasi-isomorphism $\phi : I \to J$ in $\cat{K}_{\mrm{inj}}$ is actually an 
isomorphism. Therefore $F(\phi)$ is an isomorphism in $\cat{E}$, and this is 
condition (b) of  Theorem \ref{thm:1470}. 
\end{proof}

Here is another useful definition. It is a variant of Definition \ref{dfn:2082}.

\begin{dfn} \label{dfn:1512}
Let $\cat{K}$ be a full triangulated subcategory of 
$\dcat{K}(A, \cat{M})$, and assume $\cat{K}$ has enough K-injectives.
A  {\em system of K-injective resolutions} in $\cat{K}$ is a pair $(I, \rho)$, 
where $I : \opn{Ob}(\cat{K}) \to \opn{Ob}(\cat{K}_{\mrm{inj}})$
is a function, and 
$\rho = \{ \rho_M \}_{M \in \opn{Ob}(\cat{K})}$ 
is a collection of quasi-isomorphisms
$\rho_M : M \to I(M)$ in $\cat{K}$.
Moreover, if $M \in \opn{Ob}(\cat{K}_{\mrm{inj}})$, then 
$I(M) = M$ and $\rho_M = \opn{id}_M$. 
\end{dfn}

\begin{exa} \label{exa:1530}
Let $A$ be any DG ring. We will prove later that 
$\dcat{K}^{}(A)$ has enough K-injectives. Therefore, given a triangulated 
functor $F : \dcat{K}^{}(A) \to \cat{E}$ into a triangulated category 
$\cat{E}$, the right derived functor  
$(\mrm{R} F, \eta^{\mrm{R}}) : \dcat{D}(A) \to \cat{E}$ 
exists. Suppose we choose a system  of K-injective resolutions $(I, \rho)$ in 
$\dcat{K}(A)$. Then we get a presentation of $(\mrm{R} F, \eta^{\mrm{R}})$ as 
follows: $\mrm{R} F(M) := F(I(M))$ and $\eta^{\mrm{R}}_M := F(\rho_M)$. 
\end{exa}

The proposition below is a variant of Lemma \ref{lem:2081}.

\begin{prop} \label{prop:1510}
In the situation of Theorem \tup{\ref{thm:3135}}, suppose $\cat{K}$ has enough 
K-injectives, and let $(I, \rho)$ be a system of K-injective resolutions in 
$\cat{K}$. Then the function $I$ extends uniquely to a triangulated functor
$I : \cat{D} \to \cat{K}_{\mrm{inj}}$, such that 
$\opn{Id}_{\cat{K}_{\mrm{inj}}} = I \circ \opn{Q}|_{\cat{K}_{\mrm{inj}}}$, 
and 
$\rho : \opn{Id}_{\cat{D}} \twoto \opn{Q} \circ \, I$
is an isomorphism of triangulated functors.
\end{prop}

\begin{proof}
The proof is the same as that of Lemma \ref{lem:2081}, 
except that here we use Corollary \ref{cor:1515}. 
\end{proof}

The next corollary is a categorical interpretation of the last proposition.

\begin{cor}[Functorial K-Injective Resolutions] \label{cor:1920}
Let $\cat{K}$ be a full triangulated subcategory of 
$\dcat{K}(A, \cat{M})$, and assume $\cat{K}$ has enough K-injectives.
\begin{enumerate}
\item There is a triangulated functor 
$I : \cat{K} \to \cat{K}$ and a morphism of triangulated functors 
$\rho : \opn{Id}_{\cat{K}} \twoto I$, such that for every object $M \in \cat{K}$
the object $I(M)$ is K-injective, and the morphism 
$\rho_M : M \to I(M)$ is a quasi-isomorphism. 

\item If $(I', \rho')$ is another such pair, then there is a unique isomorphism 
of triangulated functors $\ze : I \twoiso I'$ such that 
$\rho' = \ze \circ \rho$. 
\end{enumerate}
\end{cor}

\begin{exer} \label{exer:1920}
Prove Corollary \ref{cor:1920}.
\end{exer}

Infinite direct sums in $\dcat{C}_{\mrm{str}}(A, \cat{M})$ were discussed in 
Proposition \ref{prop:3140}. Their existence depends only on the existence of 
infinite direct sums in $\cat{M}$. 

\begin{thm} \label{thm:3140}
\index{Direct sums! infinite {\indash} in $\dcat{D}(A, \cat{M})$}
Let $A$ be a DG ring, let $\cat{M}$ be an abelian category, and let 
$\{ M_x \}_{x \in X}$ be a collection of objects of 
$\dcat{C}(A, \cat{M})$. Assume that $\dcat{C}(A, \cat{M})$ has enough 
K-injectives, and that the direct sum 
$M := \bigoplus_{x \in X} M_x$
exists in $\dcat{C}_{\mrm{str}}(A, \cat{M})$. 
Then $M$ is the direct sum of the collection $\{ M_x \}_{x \in X}$ in 
$\dcat{D}(A, \cat{M})$. 
\end{thm}

\begin{proof}
Let's use these shorthands: 
$\cat{C}_{\mrm{str}} := \dcat{C}_{\mrm{str}}(A, \cat{M})$,
$\cat{K} := \dcat{K}(A, \cat{M})$ and 
$\cat{D} := \dcat{D}(A, \cat{M})$.
For each index $x$ there is a morphism $e_x : M_x \to M$ in 
$\cat{C}_{\mrm{str}}$,
that passes to a morphism 
$\opn{Q}(e_x) : M_x \to M$ in $\cat{D}$. We need to verify that the collection 
of morphisms $\{ \opn{Q}(e_x) \}_{x \in X}$ has the universal property of a 
coproduct. 

Take any object $N \in \cat{D}$, and choose a K-injective 
resolution $N \to J$ in $\cat{C}_{\mrm{str}}$.
There is an isomorphism of DG $\K$-modules 
\begin{equation} \label{eqn:3140}
\begin{aligned}
&
\opn{Hom}_{A, \cat{M}}(M, J) = 
\opn{Hom}_{A, \cat{M}} \Bigl( \bigoplus\nolimits_{x \in X} M_x, \, J \Bigr)
\\
& \quad 
\cong \prod\nolimits_{x \in X} \opn{Hom}_{A, \cat{M}}(M_x, J) .
\end{aligned}
\end{equation}
From it we obtain isomorphisms of $\K$-modules 
\begin{equation} \label{eqn:3145}
\begin{aligned}
&
\opn{Hom}_{\cat{D}}(M, N) \cong^{\dag} \opn{Hom}_{\cat{K}}(M, J)
= \opn{H}^0 \bigl( \opn{Hom}_{A, \cat{M}}(M, J) \bigr) 
\\
& \quad 
\cong^{\ddag} \opn{H}^0 \Bigl( 
\prod\nolimits_{x \in X} \opn{Hom}_{A, \cat{M}}(M_x, J) \Bigr) 
\\
& \quad 
\cong^{\sharp} \prod\nolimits_{x \in X}
\opn{H}^0 \bigl( \opn{Hom}_{A, \cat{M}}(M_x, J) \bigr) 
\\
& \quad 
= \prod\nolimits_{x \in X} \opn{Hom}_{\cat{K}}(M_x, J)
\cong^{\dag} \prod\nolimits_{x \in X} \opn{Hom}_{\cat{D}}(M_x, N) . 
\end{aligned} 
\end{equation}
The isomorphisms $\cong^{\dag}$ are due to Theorem \ref{thm:3135}; the 
isomorphism $\cong^{\ddag}$ is from formula (\ref{eqn:3140}); and the 
isomorphism $\cong^{\sharp}$ is because the functor $\opn{H}^0$ commutes 
with products. Since the composed isomorphism in (\ref{eqn:3145}) is induced by 
the collection of morphisms $\{ \opn{Q}(e_x) \}_{x \in X}$, we see that the 
universal property of a coproduct holds. 
\end{proof}

\begin{rem} \label{rem:5055}
The concept of a K-injective complex (i.e.\ a K-injective object of 
$\dcat{C}(\cat{M})$), as well as the similar concepts of K-projective and 
K-flat complexes that will be discussed below, were introduced by N. 
Spaltenstein \cite{Spa} in 1988. At about the same time other authors (B. Keller 
\cite{Kel1}, M. Bockstedt and A. Neeman \cite{BoNe}, J. Bernstein and V. 
Lunts \cite{BeLu}) discovered this concept 
independently, with other names (such as {\em homotopically injective complex}).
The texts \cite{BeLu} and \cite{Kel1} already talk about DG modules over DG 
rings. 
\end{rem}

\mysubsection{K-Projective DG Modules} \label{subsec:K-proj}

This subsection is dual to the previous one, and so we will be brief. 

\begin{dfn} \label{dfn:1520}
A DG module $P \in \dcat{C}(A, \cat{M})$ is called 
{\em K-projective}%
\index{Differential graded module! K-projective}
if for every acyclic DG module $N \in \dcat{C}(A, \cat{M})$, the DG $\K$-module 
$\opn{Hom}_{A, \cat{M}}(P, N)$ is acyclic.
\end{dfn}

\begin{dfn} \label{dfn:1521}
Let $M \in \dcat{C}(A, \cat{M})$. A 
{\em K-projective resolution}%
\index{Resolution! K-projective}
of $M$ is a quasi-isomorphism 
$\rho : P \to M$ in $\dcat{C}_{\mrm{str}}(A, \cat{M})$, where 
$P$ is a K-projective DG module. 
\end{dfn}

Remark \ref{rem:5051} on terminology applies here too. 

\begin{prop} \label{prop:1520}
A DG module $P \in \dcat{K}^{}(A, \cat{M})$ is K-projective if and only if 
$\opn{Hom}_{\dcat{K}(A, \cat{M})}(P, N) = 0$ for every acyclic DG module 
$N \in \dcat{K}^{}(A, \cat{M})$.
\end{prop}

The proof is like that of Proposition \ref{prop:1515}.

\begin{dfn} \label{dfn:1522}
Let $\cat{K}$ be a full subcategory of $\dcat{K}^{}(A, \cat{M})$. 
The full subcategory of  $\cat{K}$ on the K-projective DG modules in it is 
denoted by $\cat{K}_{\mrm{prj}}$. In other words, 
$\cat{K}_{\mrm{prj}} = \dcat{K}^{}(A, \cat{M})_{\mrm{prj}} \cap \cat{K}$.
\end{dfn}

The warning in Remark  \ref{rem:5050} applies here too. 

\begin{prop}  \label{prop:1521}
If  $\cat{K}$ is a full triangulated subcategory of $\dcat{K}^{}(A, 
\cat{M})$, then $\cat{K}_{\mrm{prj}}$ is a full triangulated subcategory of 
$\cat{K}$.
\end{prop}

The proof is like that of Proposition \ref{prop:130}.

\begin{exa} \label{exa:1521}
Let $\star$ be some boundedness condition (namely $\mrm{b}$, $+$ or $-$). 
Since $\dcat{K}^{\star}(A, \cat{M})$ 
is a full triangulated subcategory of $\dcat{K}^{}(A, \cat{M})$,
we see that $\dcat{K}^{\star}(A, \cat{M})_{\mrm{prj}}$
is a full triangulated subcategory too.  
\end{exa}

\begin{dfn} \label{dfn:1523}
Let $\cat{K}$ be a full triangulated subcategory of 
$\dcat{K}(A, \cat{M})$. We say that $\cat{K}$ 
{\em has enough K-projectives}%
\index{Resolution! K-projective}
if every DG module $M \in \cat{K}$ admits a K-projective resolution inside 
$\cat{K}$. I.e.\ there is a quasi-isomorphism $\rho : P \to M$
where $P \in \cat{K}_{\mrm{prj}}$.
\end{dfn}

\begin{lem} \label{lem:1521}
Let  $s : M \to P$ be a quasi-isomorphism in 
$\dcat{K}(A, \cat{M})$, and assume $P$ is K-projective. Then $s$ has 
a right inverse; namely there is a morphism $t : P \to M$ in
$\dcat{K}(A, \cat{M})$ such that $s \circ t = \opn{id}_P$.
Moreover, this morphism $t$ is unique.  
\end{lem}

The proof is almost the same as that of Lemma \ref{lem:1510}.

\begin{thm} \label{thm:3145}
\index{Differential graded module! K-projective}
Let $A$ be a DG ring, let $\cat{M}$ be an abelian category, and let
$\cat{K}$ be a full triangulated subcategory of $\dcat{K}(A, \cat{M})$.
Write $\cat{S} := \cat{K} \cap \, \dcat{S}(A, \cat{M})$,
and let $\cat{D} := \cat{K}_{\cat{S}}$ be the localization, with localization 
functor $\opn{Q} : \cat{K} \to \cat{D}$. 
Then for every $M \in \cat{K}$ and $P \in \cat{K}_{\mrm{prj}}$
the homomorphism  
\[ \opn{Q}_{P, M} : \opn{Hom}_{\cat{K}}(P, M) \to \opn{Hom}_{\cat{D}}(P, M) 
\]
is bijective.  
\end{thm}

The proof is like that of Theorem \ref{thm:3135}, only now we use the fact that
the functor $\opn{Q}$ is a right Ore localization. 

The results from here to the end of the subsection are also proved like their 
K-injective counterparts in Subsection \ref{subsec:K-inj}. 

\begin{cor} \label{cor:3145} 
In the situation of Theorem \tup{\ref{thm:3145}}, 
the localization functor 
$\opn{Q} : \cat{K}_{\mrm{prj}} \to \cat{D}$
is fully faithful.
\end{cor}
        
\begin{cor} \label{cor:1520}
In the situation of Theorem \tup{\ref{thm:3145}}, 
if $\cat{K}$ has enough K-pro\-jectives, then the localization functor
$\opn{Q} : \cat{K}_{\mrm{prj}} \to \cat{D}$
is an equivalence.
\end{cor}

\begin{cor} \label{cor:3146}
Let $\star$ be any boundedness condition. If $\dcat{K}^{\star}(A, \cat{M})$ has 
enough K-projectives, then the triangulated functor
\[ \opn{Q} : \dcat{K}^{\star}(A, \cat{M})_{\mrm{prj}} \to
\dcat{D}^{\star}(A, \cat{M}) \]
is an equivalence.
\end{cor}

\begin{cor} \label{cor:1521}
Let $\cat{K} \sub \cat{K}'$ be full triangulated subcategories of 
$\dcat{K}^{}(A, \cat{M})$.
Define 
$\cat{S} := \cat{K} \cap \, \dcat{S}(A, \cat{M})$,
$\cat{S}' := \cat{K}' \cap \, \dcat{S}(A, \cat{M})$,
$\cat{D} := \cat{K}_{\cat{S}}$
and
$\cat{D}' := \cat{K}'_{\cat{S}'}$.
If $\cat{K}$ and $\cat{K}'$ have enough K-projectives, then the canonical 
functor $\cat{D} \to \cat{D}'$ is fully faithful. 
\end{cor}

\begin{cor} \label{cor:1923} 
Let $\phi : P \to Q$ be a morphism in $\dcat{C}_{\mrm{str}}(A, \cat{M})$ 
between K-projective objects. Then $\phi$ is a homotopy equivalence if and only 
if it is a quasi-isomorphism. 
\end{cor}

\begin{thm}  \label{thm:5051}  
\index{Derived functor! triangulated left}
In the situation of Theorem \tup{\ref{thm:3145}}, assume that $\cat{K}$ has 
enough K-projectives. Let $\cat{E}$ be a triangulated category, and let 
$F : \cat{K} \to \cat{E}$
be a triangulated functor. Then $F$ has a left derived functor 
$(\mrm{L} F, \eta^{\mrm{L}}) : \cat{D} \to \cat{E}$.
Furthermore, for every $P \in \cat{K}_{\mrm{prj}}$ the morphism 
$\eta^{\mrm{L}}_P :\mrm{L} F(P) \to F(P)$
in $\cat{E}$ is an isomorphism. 
\end{thm}

\begin{dfn} \label{dfn:1524}
Let $\cat{K}$ be a full triangulated subcategory of 
$\dcat{K}(A, \cat{M})$, and assume $\cat{K}$ has enough K-projectives.
A  {\em system of K-projective resolutions} in $\cat{K}$ is a pair $(P, \rho)$, 
where $P : \opn{Ob}(\cat{K}) \to \opn{Ob}(\cat{K}_{\mrm{prj}})$
is a function, and 
$\rho = \{ \rho_M \}_{M \in \opn{Ob}(\cat{K})}$ 
is a collection of quasi-isomorphisms
$\rho_M : P(M) \to M$ in $\cat{K}$.
Moreover, if $M \in \opn{Ob}(\cat{K}_{\mrm{prj}})$, then 
$P(M) = M$ and $\rho_M = \opn{id}_M$. 
\end{dfn}

\begin{exa} \label{exa:1531}
Let $A$ be any DG ring. We will prove later that $\bcat{K}^{}(A)$
has enough K-projectives. Therefore, given a triangulated functor 
$F : \bcat{K}^{}(A) \to \cat{E}$ into a triangulated category $\cat{E}$, 
the left derived functor  
$(\mrm{L} F, \eta^{\mrm{L}}) : \bcat{D}^{}(A) \to \cat{E}$
exists. Suppose we choose a system  of K-projective resolutions $(P, \rho)$ in 
$\bcat{K}(A)$. Then we get a presentation of $(\mrm{L} F, \eta^{\mrm{R}})$ as 
follows: $\mrm{L} F(M) := F(P(M))$ and $\eta^{\mrm{L}}_M := F(\rho_M)$. 
\end{exa}

\begin{prop} \label{prop:1522}
In the situation of Theorem \tup{\ref{thm:3145}}, suppose $\cat{K}$ has enough 
K-projectives, and let $(P, \rho)$ be a system  of K-projective resolutions.
Then the function $P$ extends uniquely to a triangulated functor
$P : \cat{D} \to \cat{K}_{\mrm{prj}}$, such that 
$\opn{Id}_{\cat{K}_{\mrm{prj}}} = P \circ \opn{Q}|_{\cat{K}_{\mrm{prj}}}$, 
and 
$\rho :  \opn{Q} \circ \, P  \twoto \opn{Id}_{\cat{D}}$
is an isomorphism of triangulated functors.
\end{prop}

\begin{cor}[Functorial K-Projective Resolutions] \label{cor:1921}
Let $\cat{K}$ be a full triangulated subcategory of 
$\dcat{K}(A, \cat{M})$, and assume $\cat{K}$ has enough K-projectives.
\begin{enumerate}
\item There is a triangulated functor 
$P : \cat{K} \to \cat{K}$ and a morphism of triangulated functors 
$\rho : P \twoto \opn{Id}_{\cat{K}}$, such that for every object 
$M \in \cat{K}$
the object $P(M)$ is K-projective, and the morphism 
$\rho_M : P(M) \to M$ is a quasi-isomorphism. 

\item If $(P', \rho')$ is another such pair, then there is a unique isomorphism 
of triangulated functors $\ze : P' \twoiso P$ such that 
$\rho' =  \rho \circ \ze$. 
\end{enumerate}
\end{cor}

\mysubsection{K-Flat DG Modules} \label{subsec:K-flat}

Recall that for a DG ring $A$, its opposite DG ring is $A^{\mrm{op}}$.
The objects of $\dcat{C}(A^{\mrm{op}})$ are the right DG $A$-modules. 

\begin{dfn} \label{dfn:1525}
A DG module $P \in \dcat{C}(A)$ is called 
{\em K-flat}%
\index{Differential graded module! K-flat}
if for every acyclic DG module $N \in \dcat{C}(A^{\mrm{op}})$, the DG 
$\K$-module $N \ot_A P$ is acyclic.
\end{dfn}

\begin{dfn} \label{dfn:4720}
Let $M \in \dcat{C}(A)$. A 
{\em K-flat resolution}%
\index{Resolution! K-flat}
of $M$ is a quasi-iso\-morphism $\rho : P \to M$ in 
$\dcat{C}_{\mrm{str}}(A)$, where $P$ is a K-flat DG module. 
\end{dfn}

Remark \ref{rem:5051} on terminology applies here as well. 

\begin{dfn} \label{dfn:3145}
Let $\cat{K}$ be a full triangulated subcategory of $\bcat{K}(A)$.
We denote by $\cat{K}_{\mrm{flat}}$ the full subcategory of 
$\bcat{K}(A)$ on the K-flat complexes in it. 
Thus $\cat{K}_{\mrm{flat}} = \cat{K} \cap \bcat{K}(A)_{\mrm{flat}}$.
\end{dfn}

The warning in Remark  \ref{rem:5050} applies here too. 

\begin{prop} \label{prop:1525}
If $P \in \dcat{C}(A)$ is K-projective then it is K-flat. 
\end{prop}

\begin{proof}
Let $\K^*$ be an injective cogenerator of $\dcat{M}(\K) = \cat{Mod} \K$.
This means that $\K^*$ is an injective $\K$-module, such that every nonzero 
$\K$-module $W$ admits a nonzero homomorphism $W \to \K^*$. 
See Examples \ref{exa:1665} and \ref{exa:1665}. 
It is not hard to see that a DG $\K$-module $W$ is acyclic if and only if 
$\opn{Hom}_{\K}(W, \K^*)$ is acyclic. 

Take an acyclic complex $N \in \dcat{C}(A^{\mrm{op}})$. Then by Hom-tensor 
adjunction there is an isomorphism of DG $\K$-modules
\[ \opn{Hom}_{\K}(N \ot_A P, \K^*) \cong 
\opn{Hom}_{A} \bigl(P, \opn{Hom}_{\K}(N, \K^*) \bigr) . \]
The right side is acyclic by our assumptions. Hence so is the left side. It 
follows that $N \ot_A P$ is acyclic. 
\end{proof}

\begin{prop} \label{prop:1526}
A DG module $P \in \dcat{K}(A)$ is K-flat iff 
\[ \opn{Hom}_{\dcat{K}(A)} \bigl(P, \opn{Hom}_{\K}(N, J) \bigr) = 0 \]
for every acyclic $N \in \dcat{C}(A^{\mrm{op}})$ and every injective 
$J \in \cat{Mod} \K$.
\end{prop}

\begin{exer}
Prove Proposition \ref{prop:1526}. (Hint: look at the proof of Proposition 
\ref{prop:1525}.)
\end{exer}

\begin{prop} \label{prop:3145}
Let $\cat{K}$ be a full triangulated subcategory of $\bcat{K}(A)$.
Then $\cat{K}_{\mrm{flat}}$ is a full triangulated subcategory of
$\cat{K}$. 
\end{prop}

\begin{exer} \label{exer:3145}
Prove Proposition \ref{prop:3145}. (Hint: The challenge is to show that 
$\cat{K}_{\mrm{flat}}$ is closed under distinguished triangles. 
Use Proposition \ref{prop:1526}.)
\end{exer}

In Theorem \ref{thm:2107} on the left derived bifunctor 
$(- \ot^{\mrm{L}}_{A} -)$ of the bifunctor $(- \ot_A -)$ we will use K-flat 
resolutions. 

\begin{rem} \label{rem:1525} 
In view of Proposition \ref{prop:1525}, the reader might wonder why we bother 
with K-flat DG modules. The reason is that there are situations in which there 
are enough K-flat objects but not enough K-projective objects. One such 
situation is that of the category 
$\dcat{C}(\mcal{A}) = \dcat{C}(\cat{Mod} \mcal{A})$,
where $(X, \mcal{A})$ is a ringed space and $\cat{Mod} \mcal{A}$
is the category of sheaves of $\AA$-modules; see Example \ref{exa:4725}. 
Another situation is when K-flatness is needed is for DG bimodules over 
noncommutative DG rings, as in Sections \ref{sec:perf-tilt-NC} and later.
\end{rem}

\mysubsection{Opposite Resolving Subcategories}
\label{subsec:resol-op}

In this subsection we explain how K-injective and K-projective DG modules in 
$\dcat{K}(A, \cat{M})$ behave in the opposite category 
$\dcat{K}(A, \cat{M})^{\mrm{op}}$, and then use this 
to prove existence of contravariant triangulated derived functors.

Recall the flipping operation from Subsections \ref{subsec:contrvar-dg-func}, 
\ref{subsec:opp-hom-triang} and \ref{subsec:opp-dercat-triang}.
There is an isomorphism of DG categories 
\begin{equation} \label{eqn:3150}
\opn{Flip} : \dcat{C}(A, \cat{M})^{\mrm{op}} \iso 
\dcat{C}(A, \cat{M})^{\mrm{flip}} = 
\dcat{C}(A^{\mrm{op}}, \cat{M}^{\mrm{op}}) , 
\end{equation}
and the corresponding isomorphisms of the homotopy categories
\begin{equation} \label{eqn:3151}
\opn{Flip} : \dcat{K}(A, \cat{M})^{\mrm{op}} \iso 
\dcat{K}(A, \cat{M})^{\mrm{flip}} = 
\dcat{K}(A^{\mrm{op}}, \cat{M}^{\mrm{op}})
\end{equation}
and on the derived categories
\begin{equation} \label{eqn:3152}
\opn{Flip} : \dcat{D}(A, \cat{M})^{\mrm{op}} \iso 
\dcat{D}(A, \cat{M})^{\mrm{flip}} = 
\dcat{D}(A^{\mrm{op}}, \cat{M}^{\mrm{op}}) . 
\end{equation}
By definition of the triangulated structure on 
$\dcat{K}(A, \cat{M})^{\mrm{op}}$,
the isomorphisms (\ref{eqn:3151}) and  (\ref{eqn:3152})
are of triangulated categories. 

The objects of $\dcat{C}(A, \cat{M})$ and
$\dcat{C}(A, \cat{M})^{\mrm{op}}$ are the same, and the property of a DG module 
$M$ being acyclic is the same in both categories. 

\begin{lem} \label{lem:3150}
A DG module $M \in \dcat{C}(A, \cat{M})$ is acyclic iff the DG module \lb 
$\opn{Flip}(M) \in \dcat{C}(A^{\mrm{op}}, \cat{M}^{\mrm{op}})$
is acyclic. 
\end{lem}

\begin{proof}
This follows from Theorem \ref{thm:2495}(4). 
\end{proof}

\begin{lem} \label{lem:3151}
A DG module $P \in \dcat{K}(A, \cat{M})$ is K-projective iff
the DG module 
$I := \opn{Flip}(P) \in \dcat{K}(A^{\mrm{op}}, \cat{M}^{\mrm{op}})$
is K-injective. 
\end{lem}

\begin{proof} 
According to Lemma \ref{lem:3150}, the functor 
$\opn{Flip}$ gives a bijection between the set of acyclic DG modules 
$M \in \dcat{K}(A, \cat{M})$
and the set of acyclic DG modules 
$ N := \opn{Flip}(M) \in \dcat{K}(A^{\mrm{op}}, \cat{M}^{\mrm{op}})$.
Since the functor (\ref{eqn:3151}) is an equivalence, we 
have an equality and an isomorphism
\begin{equation} \label{eqn:5055}
\opn{Hom}_{\dcat{K}(A, \cat{M})} (P, M) = 
\opn{Hom}_{\dcat{K}(A, \cat{M})^{\mrm{op}}} (M, P) \cong 
\opn{Hom}_{\dcat{K}(A^{\mrm{op}}, \cat{M}^{\mrm{op}})} (N, I)  
\end{equation}
in $\dcat{M}(\K)$.

By Proposition \ref{prop:1515}, the DG 
module $I$ is K-injective in 
$\dcat{K}(A^{\mrm{op}}, \cat{M}^{\mrm{op}})$
iff for every acyclic DG module 
$N \in \dcat{K}(A^{\mrm{op}}, \cat{M}^{\mrm{op}})$
we have 
$\opn{Hom}_{\dcat{K}(A^{\mrm{op}}, \cat{M}^{\mrm{op}})}(N, I) = 0$.
Similarly, by Proposition \ref{prop:1520}, the DG 
module $P$ is K-projective in $\dcat{K}(A, \cat{M})$
iff for every acyclic DG module 
$M \in \dcat{K}(A, \cat{M})$ we have 
$\opn{Hom}_{\dcat{K}(A, \cat{M})}(P, M) = 0$.
In view of (\ref{eqn:5055}) we conclude that $I$ is K-injective iff $P$ is 
K-projective.
\end{proof}

Here are the contravariant modifications of Theorems \ref{thm:5050} and 
\ref{thm:5051}. 

\begin{setup} \label{set:3150}
We are given:
\begin{itemize}
\item A full additive subcategory
$\cat{K} \sub \dcat{K}(A, \cat{M})$, 
s.t.\ $\cat{K}^{\mrm{op}} \sub \dcat{K}(A, \cat{M})^{\mrm{op}}$
is triangulated. The localization of $\cat{K}^{\mrm{op}}$ at the 
quasi-isomorphisms in it is denoted by $\cat{D}^{\mrm{op}}$.  

\item A triangulated category $\cat{E}$ and  a triangulated functor
$F : \cat{K}^{\mrm{op}} \to \cat{E}$.
\end{itemize}
\end{setup}

\begin{thm} \label{thm:3150}
Under Setup \tup{\ref{set:3150}}, 
assume that $\cat{K}$ has enough K-projectives. 
Then the triangulated right derived functor 
$(\mrm{R} F, \eta^{\mrm{R}}) : \cat{D}^{\mrm{op}} \to \cat{E}$
of $F$ exists. Furthermore, for every K-projective DG module
$P \in \cat{K}$, the morphism 
$\eta^{\mrm{R}}_P : F(P) \to \mrm{R} F(P)$ 
in $\cat{E}$ is an isomorphism. 
\end{thm}

\begin{proof}
Let 
$\cat{K}^{\mrm{flip}} := \opn{Flip}(\cat{K}^{\mrm{op}}) \sub
\dcat{K}(A^{\mrm{op}}, \cat{M}^{\mrm{op}})$.
This is a triangulated category. Its localization w.r.t.\ the 
quasi-isomorphisms is denoted by $\cat{D}^{\mrm{flip}}$. We have isomorphisms of 
triangulated categories
$\opn{Flip} : \cat{K}^{\mrm{op}} \iso \cat{K}^{\mrm{flip}}$
and 
$\opn{Flip} : \cat{D}^{\mrm{op}} \iso \cat{D}^{\mrm{flip}}$.
Let 
$F^{\mrm{flip}} := F \circ \opn{Flip}^{-1} : \cat{K}^{\mrm{flip}} \to \cat{E}$.
This is a triangulated functor. 

Since $\opn{Flip}$ sends quasi-isomorphisms to 
quasi-isomorphisms, and by Lemma \ref{lem:3151}, each K-projective 
resolution
$\rho : P \to M$ of $M$ in $\cat{K}$ goes to a K-injective resolution 
$\opn{Flip}(\rho) : \opn{Flip}(M) \to \opn{Flip}(P)$
of $N := \opn{Flip}(M)$ in $\cat{K}^{\mrm{flip}}$. 
We are given that $\cat{K}$ has enough K-projectives. Therefore 
$\cat{K}^{\mrm{flip}}$ has enough K-injectives. 

Now we can apply Theorem \ref{thm:5050} to deduce that the right 
derived functor
$(\mrm{R} F^{\mrm{flip}}, \eta^{\mrm{flip}}) : \cat{D}^{\mrm{flip}} 
\to \cat{E}$
of $F^{\mrm{flip}}$ exists. Then 
\[ (\mrm{R} F, \eta^{\mrm{R}}) := 
(\mrm{R} F^{\mrm{flip}}, \eta^{\mrm{flip}}) \circ \opn{Flip} : 
\cat{D}^{\mrm{op}} \to \cat{E} \]
is a right derived functor of $F$.

The assertion about the morphism $\eta^{\mrm{R}}_P$ is clear from the 
discussion 
above. 
\end{proof}

\begin{thm} \label{thm:3151}
Under Setup \tup{\ref{set:3150}}, 
assume that $\cat{K}$ has enough K-injectives. 
Then the left derived functor 
$(\mrm{L} F, \eta^{\mrm{L}}) : \cat{D}^{\mrm{op}} \to \cat{E}$
of $F$ exists. Furthermore, for every K-injective DG module
$I \in \cat{K}$, the morphism 
$\eta^{\mrm{L}}_I : \mrm{L} F(I) \to F(I)$ 
in $\cat{E}$ is an isomorphism. 
\end{thm}

\begin{proof}
The same as that of the previous theorem, except that now we use Lemma
\ref{lem:3151} and Theorem \ref{thm:5051}. 
\end{proof}

\cleardoublepage
\mysection{Existence of Resolutions} \label{sec:exist-resol}     

\AYcopyright
 
In this section we continue in the more concrete setting: 
$\K$ is a nonzero commutative base ring, $A$ is a central DG $\K$-ring, 
and $\cat{M}$ is a $\K$-linear abelian category; see Convention \ref{conv:2490}.
The DG category $\dcat{C}^{\star}(A, \cat{M})$
was introduced in Section \ref{subsec:DGModinM}. 
We will prove existence of  K-projective and K-injective resolutions in 
$\dcat{C}^{\star}(A, \cat{M})$,
under suitable conditions.

\mysubsection{Direct and Inverse Limits of Complexes}
\label{subsec:lims-cplx}
We shall have to work with limits in this section. Limits in abstract abelian 
and DG categories (not to mention triangulated categories) are a very 
delicate issue. We will try to be as concrete as possible, in order to avoid 
pitfalls and confusion. The notation we use for direct and inverse limits
of systems indexed by $\N$ was presented in Subsection \ref{subsec:inv-dir-lim}.
 
\begin{prop} \label{prop:1660} \mbox{}
\begin{enumerate}
\item Let $\{ M_k \}_{k \in \N}$ be a direct system in 
$\dcat{C}_{\mrm{str}}(A, \cat{M})$. Assume that for every $i$ the direct limit 
$\lim_{k \to} M^i_k$ exists in $\cat{M}$. Then 
the direct limit 
$M = \lim_{k \to} M_k$ exists in $\dcat{C}_{\mrm{str}}(A, \cat{M})$, and in 
degree $i$ it is $M^i = \lim_{k \to} M^i_k$.

\item Let $\{ M_k \}_{k \in \N}$ be an inverse system in 
$\dcat{C}_{\mrm{str}}(A, \cat{M})$. Assume that for every $i$ the inverse limit 
$\lim_{\lar k} M^i_k$ exists in $\cat{M}$. Then 
the inverse limit 
$M = \lim_{\lar k} M_k$ exists in $\dcat{C}_{\mrm{str}}(A, \cat{M})$, and in 
degree $i$ it is $M^i = \lim_{\lar k} M^i_k$.
\end{enumerate}
\end{prop}

\begin{proof}
We will only prove item (1); the proof of item (2) is similar. 
For each integer $i$ define 
$M^i := \lim_{k \to} M^i_k \in \cat{M}$. 
By the universal property of the direct limit, the differentials 
$\d : M^i_k \to M^{i + 1}_k$ induce differentials 
$\d : M^i \to M^{i + 1}$, and in this way we obtain a complex 
$M := \{ M^i \}_{i \in \Z} \in \dcat{C}_{}(\cat{M})$.
Similarly, every element $a \in A^j$ induces morphisms 
$a : M^i \to M^{i + j}$ in $\cat{M}$, and thus $M$ becomes an object of 
$\dcat{C}(A, \cat{M})$. There are morphisms $M_k \to M$ in 
$\dcat{C}_{\mrm{str}}(A, \cat{M})$, and it is easy to see that these make $M$ 
into a direct limit of the system $\{ M_k \}_{k \in \N}$.
\end{proof}

Since limits exist in $\cat{M} = \cat{Mod} \K$, the proposition above says 
that they exist in $\bcat{C}_{\mrm{str}}(A)$. Likewise they exist in the 
category $\dcat{G}_{\mrm{str}}(\K)$ of graded $\K$-modules.

We say that a direct system $\{ M_k \}_{k \in \N}$ in $\cat{M}$ is {\em 
eventually stationary} if $\mu_k : M_k \to M_{k + 1}$ are isomorphisms for 
large $k$. Similarly we can talk about an eventually stationary inverse system. 
The limit of an eventually stationary system (direct or inverse) always exists: 
it is $M_k$ for large enough $k$. 

\begin{prop} \label{prop:1575} \mbox{}
\begin{enumerate}
\item Let $\{ M_k \}_{k \in \N}$ be a direct system in 
$\bcat{C}_{\mrm{str}}(A, \cat{M})$.  
Assume that for each $i$ the direct system $\{ M^i_k \}_{k \in \N}$ in 
$\cat{M}$ is eventually stationary. Then the direct limit 
$M = \lim_{k \to} M_k$ exists in $\dcat{C}_{\mrm{str}}(A, \cat{M})$, 
the  direct limit 
$\lim_{k \to} \opn{H}(M_k)$ exists in $\dcat{G}_{\mrm{str}}(\cat{M})$, 
and the canonical morphism 
$\lim_{k \to} \, \opn{H}(M_k) \to \opn{H}(M)$ 
in $\dcat{G}_{\mrm{str}}(\cat{M})$ is an isomorphism. 

\item Let $\{ M_k \}_{k \in \N}$ be an inverse system in 
$\bcat{C}_{\mrm{str}}(A, \cat{M})$. Assume that for each $i$ the inverse system 
$\{ M^i_k \}_{k \in \N}$ in $\cat{M}$ is eventually stationary. Then the 
inverse limit $M = \lim_{\leftarrow k} M_k$ exists in 
$\dcat{C}_{\mrm{str}}(A, \cat{M})$, the  inverse limit 
$\lim_{\leftarrow k} \opn{H}(M_k)$ exists in $\dcat{G}_{\mrm{str}}(\cat{M})$, 
and the canonical morphism 
$\opn{H}(M) \to \lim_{\leftarrow k} \, \opn{H}(M_k)$
in $\dcat{G}_{\mrm{str}}(\cat{M})$ is an isomorphism. 
\end{enumerate}
\end{prop}

\begin{proof}
Again, we only prove item (1). 
As mentioned above, for each $i$ the limit
$M^i = \lim_{k \to} M^i_k$ exists in $\cat{M}$.
By Proposition \ref{prop:1660} the limit 
$M = \lim_{k \to} M_k$ exists in $\dcat{C}_{\mrm{str}}(A, \cat{M})$.

Regarding the cohomology: fix an integer $i$. 
Take $k$ large enough such that 
$M^{i'}_k \to M^{i'}_{k'}$ are isomorphisms for all $k \leq k'$ and 
$i - 1 \leq  i' \leq i + 1$. Then 
$M^{i'}_{k'} \to M^{i'}$ 
are isomorphisms in this range, and therefore 
$\opn{H}^i(M_{k'}) \to \opn{H}^i(M)$ are isomorphisms for all $k \leq k'$. 
We see that the direct system 
$\{ \opn{H}^i(M_k) \}_{k \in \N}$
is eventually stationary, and its direct limit is 
$\opn{H}^i(M)$.
\end{proof}

When we drop the abstract abelian category $\cat{M}$, i.e.\ when we work 
with $\cat{M} = \cat{Mod} \K =  \dcat{M}(\K)$ and 
$\dcat{C}_{\mrm{str}}(A, \cat{M}) = \dcat{C}_{\mrm{str}}(A)$, 
there is no problem of 
existence of limits. The next proposition says that furthermore ``direct limits 
are exact'' in $\bcat{C}_{\mrm{str}}(A)$. 

\begin{prop} \label{prop:1585}
Let $\{ M_k \}_{k \in \N}$ be a direct system in $\bcat{C}_{\mrm{str}}(A)$,
with limit $M = \lim_{k \to} M_k$. 
Then the canonical homomorphism 
$\lim_{k \to} \opn{H}(M_k) \to \opn{H}(M)$ 
in $\dcat{G}_{\mrm{str}}(\K)$ is bijective. 
\end{prop}

\begin{exer}  \label{exer:1741}
Prove Proposition \ref{prop:1585}. (Hint: forget the action of $A$, and work 
with complexes of $\K$-modules.)
\end{exer}

Exactness of inverse limits tends to be much more complicated than that of 
direct limits, even for $\K$-modules. We always have to impose some condition 
on the inverse system to have exactness in the limit. 

Let 
$\bigl( \{ M_k \}_{k \in \N}, \{ \mu_k \}_{k \in \N} \bigr)$
be an inverse system in $\bcat{M}(\K)$.
For every $k \leq l$ let 
$M_{k, l} \sub M_k$ be the image of the homomorphism
$\mu_{k, l} : M_l \to M_k$.
Note that there are inclusions 
$M_{k, l + 1} \sub M_{k, l}$, so for fixed $k$ we have an inverse system 
$\{ M_{k, l} \}_{l \geq k}$
of submodules of $M_k$. 

\begin{dfn} \label{dfn:1590}
We say that an inverse system $\{ M_k \}_{k \in \N}$ in $\bcat{M}(\K)$
has the {\em Mittag-Leffler property}%
\index{Mittag-Leffler property}
if for every index $k$, the 
inverse system $\{ M_{k, l} \}_{l \geq k}$ 
is eventually stationary.
\end{dfn}

\begin{exa} \label{exa:1590}
If the inverse system 
$\bigl( \{ M_k \}_{k \in \N}, \{ \mu_k \}_{k \in \N} \bigr)$
in $\bcat{M}(\K)$ satisfies at least one of the following conditions, then it 
has the  Mittag-Leffler property:
\begin{itemize}
\rmitem{a} The system has surjective transitions. 

\rmitem{b} The system is eventually stationary. 

\rmitem{c} For every $k \in \N$ there exists some $l \geq k$ such that
$M_{k, l} = 0$. This is called the {\em trivial Mittag-Leffler property}, and 
one says that the system is {\em pro-zero}. 
\end{itemize}
\end{exa}

\begin{thm}[Mittag-Leffler Argument] \label{thm:1535}
\index{Mittag-Leffler property}
Let $\{ M_k \}_{k \in \N}$ be an inverse system in $\dcat{C}_{\mrm{str}}(A)$, 
with inverse limit 
$M = \lim_{\leftarrow k} M_k$. Assume the system satisfies these two 
conditions\tup{:}
\begin{enumerate}
\rmitem{a} For every $i \in \Z$ the inverse system
$\{ M^i_k \}_{k \in \N}$ in $\dcat{M}(\K)$
has the Mittag-Leffler property.

\rmitem{b} For every $i \in \Z$ the inverse system
$\{ \opn{H}^i(M_k) \}_{k \in \N}$ in $\dcat{M}(\K)$
has the Mittag-Leffler property.
\end{enumerate}
Then the canonical homomorphisms 
$\opn{H}^i(M) \to \lim_{\leftarrow k} \opn{H}^i(M_k)$ 
are bijective for all $i$.
\end{thm}

\begin{proof}
We can forget all about the graded $A$-module structure, and just view this as 
an inverse system in $\dcat{C}_{\mrm{str}}(\Z)$, i.e.\ an inverse system of 
complexes of abelian groups. Now this is a special case of  
\cite[Proposition 1.12.4]{KaSc1} or 
\cite[Chapitre $0_{\mrm{III}}$, Proposition 13.2.3]{{EGA-III}}.
\end{proof}

The most useful instance of the ML argument is this:

\begin{cor} \label{cor:1590}
Let $\{ M_k \}_{k \in \N}$ be an inverse system in $\dcat{C}_{\mrm{str}}(A)$, 
with inverse limit 
$M = \lim_{\leftarrow k} M_k$. Assume the system satisfies these two 
conditions\tup{:}
\begin{enumerate}
\rmitem{a} For every $i \in \Z$ the inverse system
$\{ M^i_k \}_{k \in \N}$ has surjective transitions. 

\rmitem{b} For every $k$ the DG module $M_k$ is acyclic.  
\end{enumerate}
Then $M$ is acyclic. 
\end{cor}

\begin{proof}
Conditions (a) and (b) here imply conditions (a) and (b) of Theorem 
\ref{thm:1535}, respectively. 
\end{proof}

\begin{exer} \label{exer:1904}
Prove Corollary \ref{cor:1590} directly, without resorting to Theorem 
\ref{thm:1535}. 
\end{exer}

By {\em generalized integers}%
\index{Generalized integer}
we mean elements of the ordered set 
$\Z \cup \{ \pm \infty \}$.  Recall that for a subset $S \subseteq \Z$, its 
infimum is $\opn{inf}(S) \in \Z \cup \{ \pm \infty \}$,
where $\opn{inf}(S) = + \infty$ iff $S = \varnothing$. 
Likewise the supremum is $\opn{sup}(S) \in \Z \cup \{ \pm \infty \}$,
where $\opn{sup}(S) = - \infty$ iff $S = \varnothing$.
For $i \in \Z \cup \{ \pm \infty \}$ 
the expression $-i$ has a value in $\Z \cup \{ \pm \infty \}$.
If $i_0, i_1 \in \Z \cup \{ \pm \infty \}$ satisfy 
$i_0 \leq i_1$, then $-i_1 \leq -i_0$. 
For $i_0, i_1 \in \Z \cup \{ \infty \}$, the expressions 
$i_0 + i_1$ and $-i_0 - i_1$ have values in 
$\Z \cup \{ \pm \infty \}$. 
A generalized integer $i$ is called {\em finite} if $i < \infty$,
i.e.\ if $i \in  \Z \cup \{ -\infty \}$. 

\begin{dfn} \label{dfn:4765}
Given $i_0 \in \Z \cup \{ -\infty \}$ and 
$i_1 \in \Z \cup \{ \infty \}$, with $i_0 \leq i_1$,
the {\em integer interval}%
\index{Integer interval}
with these endpoints is the set 
\[ [i_0, i_1] := \{ i \in \Z  \mid i_0 \leq i \leq i_1 \} \sub \Z . \]
There is also the empty integer interval $\varnothing$.
When the elements of an integer interval represent degrees, we 
sometimes call it a {\em degree interval}%
\index{Degree interval}.
\end{dfn}

We exclude the possibilities $i_0 = i_1 = \infty$ and $i_0 = i_1 = -\infty$
for endpoints of intervals.
In calculus notation we have 
$[i_0, \infty] = [i_0, \infty) \cap \Z$, etc.; 
but we avoid the open interval notation on purpose. 
Whenever we refer to an integer interval $[i_0, i_1]$, it is always 
assumed that the generalized integers $i_0$ and $i_1$ satisfy the conditions in 
Definition \ref{dfn:4765}.
Further properties of integer intervals will be studied in Subsection 
\ref{subsec:bounded-revis}. 

\begin{dfn} \label{dfn:3155}
Let $M = \{ M^i \}_{i \in \Z}$ be a graded object in $\cat{M}$,
and let 
$S := \{ i \in \Z \mid  M^i \neq 0 \} \sub \Z$. 
\begin{enumerate}
\item The {\em supremum} of $M$ is
$\sup(M) := \opn{sup}(S) \in \Z \cup \{ \pm \infty \}$. 

\item The {\em infimum} of $M$ is 
$\inf(M) := \opn{inf}(S) \in \Z \cup \{ \pm \infty \}$. 

\item The {\em amplitude} of $M$ is 
$\opn{amp}(M) := \sup (M) - \inf (M) \in \N \cup \{ \pm \infty \}$.
\end{enumerate}
\end{dfn}

Note that for $M = 0$ this reads $S = \varnothing$, 
$\inf (M) = \infty$, $\sup (M) = -\infty$ 
and $\opn{amp} (M) = -\infty$.
Thus $M$ is bounded (resp.\ bounded above, resp.\ bounded below) iff 
$\opn{amp} (M) < \infty$ (resp.\ $\opn{sup} (M) < \infty$, resp.\
$\opn{inf} (M) > -\infty$).

\mysubsection{Totalizations}
Consider a complex 
\begin{equation} \label{eqn:3300}
(M, \pa) = \bigl( \cdots \to M^{-1} \xar{\pa^{-1}}
M^{0} \xar{\pa^{0}} M^{1} \to \cdots \bigr) 
\end{equation}
with entries in the abelian category $\dcat{C}_{\mrm{str}}(A, \cat{M})$. 
This means that for each $i \in \Z$ there is a DG module 
$(M^{i}, \d_{M^i}) \in \dcat{C}(A, \cat{M})$, 
and a morphism 
$\pa^i : M^i \to M^{i + 1}$
in $\dcat{C}_{\mrm{str}}(A, \cat{M})$, 
such that $\pa^{i + 1} \circ \pa^i = 0$.
Recall that the strictness of $\pa^i$ says that it has degree $0$, and that
\begin{equation} \label{eqn:3365}
\pa^i \circ \d_{M^i}  = \d_{M^{i + 1}} \circ \pa^i .
\end{equation}
In terms of our symbolic notation, $(M, \pa)$ is an object of the DG category 
\lb $\dcat{C}(\dcat{C}_{\mrm{str}}(A, \cat{M}))$. 

Each DG module $(M^i, \d_{M^i})$ has its own internal structure:
\begin{equation} \label{eqn:3305}
(M^i, \d_{M^i}) = \bigl( \{ M^{i, j} \}_{j \in \Z}, 
\{ \d_{M^{i}}^{\lms j} \}_{j \in \Z} \bigr) . 
\end{equation}
Here $M^{i, j} \in \cat{M}$, 
$\d_{M^{i}}^{\lms j} : M^{i, j} \to M^{i, j + 1}$ 
is a morphism in $\cat{M}$, and the relation 
$\d_{M^{i}}^{\lms j + 1} \circ \d_{M^{i}}^{\lms j} = 0$ holds.
There is an action of the DG ring $A$ on $M^i$, prescribed by a DG $\K$-ring 
homomorphism 
$f_{M^i} : A \to \opn{End}_{\cat{M}}(M^i)$.

Sometimes we view $M$ as a double complex, and then we refer to its rows and 
columns: $M^i$ is the $i$-th column of $M$, and
$\{ M^{i, j} \}_{i \in \Z}$ is its $j$-th row.
See diagram (\ref{eqn:5057}). 
\begin{equation} \label{eqn:5057}
\UseTips \xymatrix @C=3ex @R=3ex {
&
\vdots
& &
\vdots
& &
\vdots
\\
\cdots
\ar[r]
& 
M^{-1, 1}
\ar[rr]^{\pa^{-1}}
\ar[u]
& &
M^{0, 1}
\ar[rr]^{\pa^{0}}
\ar[u]
& &
M^{1, 1}
\ar[r]
\ar[u]
& 
\cdots
\\ \\
\cdots
\ar[r]
& 
M^{-1, 0}
\ar[rr]^{\pa^{-1}}
\ar[uu]^{\d_{M^{-1}}^{0}}
& &
M^{0, 0}
\ar[rr]^{\pa^{0}}
\ar[uu]^{\d_{M^{0}}^{0}}
& &
M^{1, 0}
\ar[r]
\ar[uu]^{\d_{M^{1}}^{0}}
& 
\cdots
\\
&
*[u]{\vdots}
\ar[u]
& &
*[u]{\vdots}
\ar[u]
& &
*[u]{\vdots}
\ar[u]
} 
\end{equation}

As explained in Proposition \ref{prop:1660},
if the abelian category $\cat{M}$ has 
countable direct sums, then so does $\dcat{C}_{\mrm{str}}(A, \cat{M})$. 
Recall that $A^{\natural}$ is the underlying graded ring of the DG ring $A$. 

\begin{dfn} \label{dfn:3300}
Assume that $\cat{M}$ has countable direct sums. Given a complex 
$(M, \pa)$ with entries in $\dcat{C}_{\mrm{str}}(A, \cat{M})$, with notation as 
in \tup{(\ref{eqn:3300})}, its {\em direct sum totalization}%
\index{Totalization! direct sum}
is the object 
\[ \opn{Tot}^{\oplus}(M, \pa)  = 
\bigl( \opn{Tot}^{\oplus}(M), \d_{\mrm{Tot}} \bigr) \in  
\dcat{C}(A, \cat{M}) \]
defined as follows, in four stages. 
\begin{enumerate}
\item There is a graded object 
\[ \opn{Tot}^{\oplus}(M) := \bigoplus_{i \in \Z} \, 
\opn{T}^{-i}(M^i) \in \dcat{G}(A^{\natural}, \cat{M}) . \]

\item On each summand $\opn{T}^{-i}(M^i)$ of $\opn{Tot}^{\oplus}(M)$
there is a differential 
$\d_{\opn{T}^{-i}(M^i)}$, 
and we define the degree $1$ operator 
$\d_M := \bigoplus_{i \in \Z} \d_{\opn{T}^{-i}(M^i)}$ 
on $\opn{Tot}^{\oplus}(M)$.

\item For each $i$ let 
\[ \opn{tot}(\pa)^i := \opn{t}^{-1}_{\lms \opn{T}^{-(i + 1)}(M^{i + 1})} \circ 
\opn{T}^{-i}(\pa^i) : 
\opn{T}^{-i}(M^i) \to \opn{T}^{-(i + 1)}(M^{i + 1}) . \]
We define the degree $1$ operator 
$\opn{tot}(\pa) :=  \bigoplus_{i \in \Z} \, \opn{tot}(\pa)^i$
on the graded object $\opn{Tot}^{\oplus}(M)$.

\item The differential $\d_{\mrm{Tot}}$ on the graded object 
$\opn{Tot}^{\oplus}(M)$ is 
$\d_{\mrm{Tot}} := \d_M + \opn{tot}(\pa)$.
\end{enumerate}
\end{dfn}

For the definition to be valid we need to verify something: 

\begin{lem} \label{lem:3030}
The degree $1$ operator $\d_{\mrm{Tot}}$ makes the pair 
$\bigl( \opn{Tot}^{\oplus}(M), \d_{\mrm{Tot}} \bigr)$
into an object of $\dcat{C}(A, \cat{M})$. 
\end{lem}

\begin{proof}
We have to check that 
$\d_{\mrm{Tot}} \circ \d_{\mrm{Tot}} = 0$
and that $\opn{tot}(\pa)$ is $A$-linear. 
Below we shall verify the first condition; the second verification is similar, 
and is left to the reader. 

For a given $i \in \Z$ consider the composed morphism starting at the summand
$\opn{T}^{-i}(M^i)$ of the total DG module:
\[ \begin{aligned}
&
\opn{T}^{-i}(M^i) \xar{\d_{\mrm{Tot}}}
\opn{T}^{-i}(M^i) \oplus \opn{T}^{-(i + 1)}(M^{i + 1}) 
\\
& \quad 
\xar{\d_{\mrm{Tot}}}
\opn{T}^{-i}(M^i) \oplus \opn{T}^{-(i + 1)}(M^{i + 1})
\oplus \opn{T}^{-(i + 2)}(M^{i + 2}) . 
\end{aligned} \]
Let us view these direct sums as columns, with matrices of operators acting on 
them from the left. So we have to calculate the matrix product 
\[ 
\bmat{ \d_{\opn{T}^{-i}(M^i)} & 0
\\[0.3em]
\opn{t}^{-1} \circ \opn{T}^{-i}(\pa^i) & \d_{\opn{T}^{-(i + 1)}(M^{i + 1})}
\\[0.3em]
0 & \opn{t}^{-1} \circ \opn{T}^{-(i + 1)}(\pa^{i + 1}) }
\circ
\bmat{ \d_{\opn{T}^{-i}(M^i)} 
\\[0.3em]
\opn{t}^{-1} \circ \opn{T}^{-i}(\pa^i) } 
\]
The resulting column of operators has  
$\d_{\opn{T}^{-i}(M^i)} \circ \d_{\opn{T}^{-i}(M^i)} = 0$
in the first position. In the second position we have
\begin{equation} \label{eqn:3366}
\opn{t}^{-1} \circ \opn{T}^{-i}(\pa^i) \circ \d_{\opn{T}^{-i}(M^i)} +
\d_{\opn{T}^{-(i + 1)}(M^{i + 1})} \circ \opn{t}^{-1} 
\circ \opn{T}^{-i}(\pa^i) .
\end{equation}
Using Proposition \ref{prop:4155} we have 
\[ \d_{\opn{T}^{-(i + 1)}(M^{i + 1})} \circ  \opn{t}^{-1} = 
- \opn{t}^{-1} \circ \, \d_{\opn{T}^{-i}(M^{i + 1})} = 
- \opn{t}^{-1} \circ \opn{T}^{-i}(\d_{M^{i + 1}}) . \]
Therefore (\ref{eqn:3366}) equals 
\[ \opn{t}^{-1} \circ \opn{T}^{-i} \bigl( \pa^i \circ \d_{M^i} -
\d_{M^{i + 1}} \circ \pa^i \bigr) , \]
which is zero by (\ref{eqn:3365}). 
In the bottom of the column we have
\begin{equation} \label{eqn:3367}
\opn{t}^{-1} \circ \opn{T}^{-(i + 1)}(\pa^{i + 1}) \circ 
\opn{t}^{-1} \circ \opn{T}^{-i}(\pa^{i}) . 
\end{equation}
According to Definition \ref{dfn:1171} we know that 
\[ \opn{T}^{-(i + 1)}(\pa^{i + 1}) \circ \opn{t}^{-1} = 
\opn{t}^{-1} \circ \opn{T}(\opn{T}^{-(i + 1)}(\pa^{i + 1})) = 
\opn{t}^{-1} \circ \opn{T}^{-i}(\pa^{i + 1})) . \]
Therefore (\ref{eqn:3367}) equals 
$\opn{t}^{-1} \circ \opn{T}^{-i} \bigl( \pa^{i + 1} \circ \pa^{i} \bigr) = 0$. 
\end{proof}

\begin{exer} \label{exer:4750}
Finish the proof of the lemma (the $A$-linearity of $\opn{tot}(\pa)$). 
\end{exer}

\begin{exa} \label{exa:3300}
If $M^i = 0$ for $i \neq -1, 0$, then 
\[ \opn{Tot}^{\oplus}(M, \pa) = 
\opn{Cone}(\pa^{-1} : M^{-1} \to M^0) , \]
the standard cone from Definition \ref{dfn:1172}.
In Proposition \ref{prop:3305} there is  a more general statement. 
\end{exa}

In Definition \ref{dfn:2320} we talked about smart truncations of DG 
modules. Here is another sort of truncations. 

\begin{dfn} \label{dfn:3066}
Let $\cat{N}$ be an abelian category. 
For a complex $N \in \dcat{C}(\cat{N})$ its 
{\em stupid truncations}%
\index{Truncation of a complex! stupid}%
\index{1-Sttleq@$\opn{stt}^{\leq q}(N)$}
at an integer $q$ are 
\[ \opn{stt}^{\leq q}(N) := 
\bigl( \cdots \to N^{q - 1} \to N^{q} \to 0 \to 0 \to \cdots \bigr) \] 
and%
\index{1-Sttgeq@$\opn{stt}^{\geq q}(N)$}
\[ \opn{stt}^{\geq q}(N) := 
\bigl( \cdots \to 0 \to 0 \to N^{q} \to N^{q + 1} \to \cdots \bigr)  . \]
\end{dfn}

These truncations fit into a short exact sequence 
\begin{equation} \label{eqn:2147}
0 \to \opn{stt}^{\geq q}(N) \to N \to \opn{stt}^{\leq q - 1}(N) \to  0
\end{equation}
in $\dcat{C}_{\mrm{str}}(\cat{N})$. 

Warning: The stupid truncations do not apply to DG modules 
$N \in \dcat{C}(B, \cat{N})$, unless the DG ring $B$ is a ring.   
However, in what follows we will use stupid truncations in 
$\dcat{C}(\cat{N})$, where 
$\cat{N}$ is the abelian category $\dcat{C}_{\mrm{str}}(A, \cat{M})$. 
This could be a bit confusing. 

\begin{prop} \label{prop:3305}
Let 
\[ (M, \pa) =  \bigl( \cdots \to M^{-1} \xar{\pa^{-1}}
M^{0} \xar{\pa^{0}} M^{1} \to 0 \to \cdots \bigr) \]
be a complex with entries in the abelian category 
$\dcat{C}_{\mrm{str}}(A, \cat{M})$, let 
\[ (M', \pa) :=  \bigl( \cdots \to M^{-1} \xar{\pa^{-1}}
M^{0} \to 0 \to \cdots \bigr) , \]
namely the stupid truncation of $M$ below $0$, and let 
$\rho : \opn{Tot}^{\oplus}(M', \pa) \to M^1$
be the morphism in $\dcat{C}_{\mrm{str}}(A, \cat{M})$ induced by
$\pa^{0} : M^{0} \to M^1$. 
Then there is an isomorphism 
\[ \tag{\dag} \opn{T} \bigl( \opn{Tot}^{\oplus}(M, \pa) \bigr) \cong 
\opn{Cone} \bigl( - \rho : \opn{Tot}^{\oplus}(M', \pa) \to M^1 \bigr) \]
in $\dcat{C}_{\mrm{str}}(A, \cat{M})$. 
\end{prop}

\begin{proof}
By comparing translations we see that the objects on both sides of (\dag) match 
in each degree. Likewise for all the differentials that do not involve 
$\pa^0$. A calculation similar to that in the proof of Lemma 
\ref{lem:3030} shows that there is a sign change there, and this is why we need 
$-\rho$ on the right side of (\dag). 
\end{proof}

\begin{cor} \label{cor:3310}
If, in the situation of Proposition \tup{\ref{prop:3305}}, the DG module  
$\opn{Tot}^{\oplus}(M, \pa)$ is acyclic, then 
$\rho : \opn{Tot}^{\oplus}(M', \pa) \to M^1$ 
is a quasi-isomorphism. 
\end{cor}

\begin{proof}
If $\opn{Tot}^{\oplus}(M, \pa)$ is acyclic then so is 
$\opn{T} \bigl( \opn{Tot}^{\oplus}(M, \pa) \bigr)$. 
The proposition implies that $-\rho$ is a quasi-isomorphism. Therefore $\rho$ 
is a quasi-isomorphism. 
\end{proof}

\begin{dfn} \label{dfn:3335}
Let $M \in \dcat{C}(A, \cat{M})$.
\begin{enumerate}
\item A {\em filtration}%
\index{Filtration on a DG module}
on $M$ is a direct system
$\{ F_j(M) \}_{j \geq -1}$ in $\dcat{C}_{\mrm{str}}(A, \cat{M})$
consisting of subobjects $F_j(M) \sub M$, indexed by the interval 
$[-1, \infty] \sub \Z$. 

\item We say that $M = \lim_{j \to} F_j(M)$
if the direct limit exists, and the 
canonical morphism $\lim_{j \to} F_j(M) \to M$ in 
$\dcat{C}_{\mrm{str}}(A, \cat{M})$ is an isomorphism. 

\item The filtration $F$ induces, for each $j \geq 0$, the subquotient 
\[ \opn{Gr}^F_j(M) := F_{j}(M) / F_{j - 1}(M) \in \dcat{C}(A, \cat{M}) . \] 
\end{enumerate}
\end{dfn}

Our filtrations will always be ascending. But sometimes they will be indexed 
by a subinterval $[j_0, j_1]$ of $[-1, \infty]$; and in this case 
$\opn{Gr}^F_j(M)$ will be defined only for 
$j \in [j_0 + 1, j_1]$.
 
In the next proposition we take $\cat{M} = \dcat{M}(\K)$; but really 
all we need is that $\cat{M}$ has exact countable direct limits. 

Recall that for a DG $A$-module $M$ we denote by $\opn{Z}(M)$
the object of cocycles of $M$. It is an object of 
$\dcat{G}(\K)$.

\begin{prop} \label{prop:3300}
Let $(M, \pa)$ be a complex with entries in the abelian category 
$\dcat{C}_{\mrm{str}}(A)$, 
with internal structure \tup{(\ref{eqn:3305})}. Assume that the two conditions 
below hold\tup{:}
\begin{itemize}
\rmitem{i} For every $j \in \Z$ the complex  
\[ \bigl( \cdots \to M^{-1, j} \xar{\pa^{-1}}
M^{0, j} \xar{\pa^{0}} M^{1, j} \to \cdots \bigr) \]
with entries in $\dcat{M}(\K)$ is acyclic. 

\rmitem{ii} There is some $j_1 \in \Z$ such that for all 
$j \geq j_1$ the complex 
\[ \bigl( \cdots \to \opn{Z}^j(M^{-1}) \xar{\pa^{-1}}
\opn{Z}^j(M^{0}) \xar{\pa^{0}} \opn{Z}^j(M^{1}) \to \cdots \bigr) \]
with entries in $\dcat{M}(\K)$ is acyclic. 
\end{itemize}
Then the DG $A$-module $\opn{Tot}^{\oplus}(M, \pa)$ is acyclic. 
\end{prop}

\begin{proof} The proof is done in two steps. 

\smallskip \noindent
Step 1. Let us replace condition (ii) by the stronger condition: 
\begin{itemize}
\rmitem{ii$'$} There is some $j'_1 \in \Z$ such that 
$M^{i, j} = 0$ for all $j > j'_1$. 
\end{itemize}
For every $i$ we introduce a filtration 
$\{ F_q(M^i) \}_{q \geq -1}$
on the DG $\K$-module $M^i$ as follows: 
\[ F_q (M^i) := \opn{stt}^{\geq j'_1 - q}(M^i)
= \bigl( \cdots \to 0 \to M^{i, j_1' - q} \xar{\d_{M^i}} M^{i, j'_1 - q + 1} 
\to \cdots \bigr) , \]
the stupid truncation above $j'_1 - q$. 
These filtrations induces a filtration on the DG $\K$-module 
$\opn{Tot}^{\oplus}(M, \pa)$~: 
\[ F_q \bigl( \opn{Tot}^{\oplus}(M) \bigr) := 
\opn{Tot}^{\oplus} 
\bigl( \cdots \to F_q(M^{-1}) \xar{\pa^{-1}}
F_q(M^{0}) \xar{\pa^{0}} F_q(M^{1}) \to \cdots \bigr) . \]
We note that $F_{-1} \bigl( \opn{Tot}^{\oplus}(M) \bigr) = 0$, and 
\begin{equation} \label{eqn:3032}
\bigcup_{q \geq -1} \, F_q \bigl( \opn{Tot}^{\oplus}(M) \bigr) = 
\opn{Tot}^{\oplus}(M) . 
\end{equation}

For every $q \geq 0$ the DG $\K$-module 
$\opn{Gr}^F_q \bigl( \opn{Tot}^{\oplus}(M) \bigr)$
is isomorphic in $\dcat{C}_{\mrm{str}}(\K)$, up to signs of differentials, to 
the complex in condition (i) 
with index $j = j'_1 - q$, so it is acyclic. For $q = -1$ 
we trivially have an acyclic DG $\K$-module 
$F_{-1} \bigl( \opn{Tot}^{\oplus}(M) \bigr)$. 
For every $q \geq 0$ there is an exact sequence 
\[ 0 \to F_{q - 1} \bigl( \opn{Tot}^{\oplus}(M) \bigr) \to 
F_q \bigl( \opn{Tot}^{\oplus}(M) \bigr) \to 
\opn{Gr}^F_q \bigl( \opn{Tot}^{\oplus}(M) \bigr) \to 0 \]
in $\dcat{C}_{\mrm{str}}(\K)$. By induction on $q$ we conclude that 
$F_q \bigl( \opn{Tot}^{\oplus}(M) \bigr)$
is acyclic. Finally, by $(\ref{eqn:3032})$, and the fact that cohomology 
commutes with direct limits in $\dcat{C}_{\mrm{str}}(\K)$, we see that the DG 
module $\opn{Tot}^{\oplus}(M, \pa)$ is acyclic. 

\medskip \noindent 
Step 2. Now we assume conditions (i) and (ii). 
For each $i$ we introduce a filtration 
$\{ G_q(M^i) \}_{q \geq 0}$
on the DG $\K$-module $M^i$ as follows: 
\[ G_q(M^i) := \opn{smt}^{\leq j_1 + q}(M^i) = 
\bigl( \cdots \to M^{i, j_1 + q - 1} \xar{\d_{M^i}} 
\opn{Z}^{j_1 + q}(M^i) \to 0 \to  \cdots \bigr) , \]
the smart truncation below $j_1 + q$. 

For every value of $q$ the complex 
\[ G_q(M) := \bigl( \cdots \to G_q(M^{-1}) \xar{\pa^{-1}}
G_q(M^{0}) \xar{\pa^{0}} G_q(M^{1}) \to \cdots \bigr) \]
satisfies conditions (i) and (ii$'$), with 
$j'_1 := j_1 + q$. Indeed, for $j < j_1 + q$ the $j$-th row of $G_q(M)$ is 
\[ \bigl( \cdots \to M^{-1, j} \xar{\pa^{-1}}
M^{0, j} \xar{\pa^{0}} M^{1, j} \to \cdots \bigr) , \]
and it is exact by condition (i). 
For  $j = j_1 + q$ its $j$-th row is 
\[ \bigl( \cdots \to \opn{Z}^j(M^{-1}) \xar{\pa^{-1}}
\opn{Z}(M^{0}) \xar{\pa^{0}} \opn{Z}^j(M^{1}) \to \cdots \bigr) , \]
and it is exact by condition (ii). 
And for $j > j_1 + q$ the $j$-th row of $G_q(M)$ is zero. 
Therefore, by step 1, the DG module 
$G_q \bigl( \opn{Tot}^{\oplus}(M) \bigr) := 
\opn{Tot}^{\oplus} \bigl( G_q(M) \bigr)$ 
is acyclic. 

Finally, since $M = \lim_{q \to} G_q(M)$, 
it follows that 
\[ \opn{Tot}^{\oplus}(M) = 
\lim_{q \to} \, G_q \bigl( \opn{Tot}^{\oplus}(M) \bigr) . \]
So $\opn{Tot}^{\oplus}(M)$ is acyclic.
\end{proof}

The second half of this subsection deals with products instead of direct sums. 
As explained in Proposition \ref{prop:1660}, 
if the abelian category $\cat{M}$ has 
countable products, then so does $\dcat{C}_{\mrm{str}}(A, \cat{M})$. 

\begin{dfn} \label{dfn:3320}
Assume that $\cat{M}$ has countable products. Given a complex 
$(M, \pa)$ with entries in $\dcat{C}_{\mrm{str}}(A, \cat{M})$, with notation as 
in \tup{(\ref{eqn:3300})}, its 
{\em product totalization}%
\index{Totalization! product}
is the object 
\[ \opn{Tot}^{\Pi}(M, \pa)  = 
\bigl( \opn{Tot}^{\Pi}(M), \d_{\mrm{Tot}} \bigr) \in  
\dcat{C}(A, \cat{M}) \]
defined as follows, in four stages. 
\begin{enumerate}
\item There is a graded object 
\[ \opn{Tot}^{\Pi}(M) := \prod_{i \in \Z} \, 
\opn{T}^{-i}(M^i) \in \dcat{G}(A^{\natural}, \cat{M}) . \]

\item On each factor $\opn{T}^{-i}(M^i)$ of $\opn{Tot}^{\Pi}(M)$
there is a differential 
$\d_{\opn{T}^{-i}(M^i)}$, 
and we let 
$\d_M := \prod_{i \in \Z} \, \d_{\opn{T}^{-i}(M^i)}$,
which is a differential on $\opn{Tot}^{\Pi}(M)$. 

\item For each $i$ let 
\[ \opn{tot}(\pa)^i := \opn{t}^{-1}_{\lms \opn{T}^{-(i + 1)}(M^{i + 1})} \circ 
\opn{T}^{-i}(\pa^i) : 
\opn{T}^{-i}(M^i) \to \opn{T}^{-(i + 1)}(M^{i + 1}) . \]
We define the degree $1$ operator 
$\opn{tot}(\pa) := \prod_{i \in \Z} \, \opn{tot}(\pa)^i$
on the graded object $\opn{Tot}^{\Pi}(M)$.

\item The differential $\d_{\mrm{Tot}}$ on the graded object 
$\opn{Tot}^{\Pi}(M)$ is 
$\d_{\mrm{Tot}} := \d_M + \opn{tot}(\pa)$.
\end{enumerate}
\end{dfn}

Note that (if $\cat{M}$ also has countable direct sums) there is a canonical 
embedding 
\begin{equation} \label{eqn:3370}
\opn{Tot}^{\oplus}(M) \sub \opn{Tot}^{\Pi}(M) 
\end{equation}
in $\dcat{C}_{\mrm{str}}(A, \cat{M})$. 
In case $M^i = 0$ for $\abs{i} \gg 0$, then this is an equality.

For Definition \ref{dfn:3320} to be valid we need to verify something: 

\begin{lem} \label{lem:3031}
The degree $1$ operator $\d_{\mrm{Tot}}$ makes the pair 
$\bigl( \opn{Tot}^{\Pi}(M), \d_{\mrm{Tot}} \bigr)$
into an object of $\dcat{C}(A, \cat{M})$. 
\end{lem}

\begin{proof}
This is the same calculation as in the proof of Lemma \ref{lem:3030},
heuristically using the inclusion (\ref{eqn:3370}).
\end{proof}

\begin{prop} \label{prop:3320}
Let 
\[ (M, \pa) =  \bigl( \cdots \to 0 \to M^{-1} \xar{\pa^{-1}}
M^{0} \xar{\pa^{0}} M^{1} \to \cdots \bigr) \]
be a complex with entries in the abelian category 
$\dcat{C}_{\mrm{str}}(A, \cat{M})$, let 
\[ (M', \pa) :=  \bigl( \cdots \to 0 \to M^{0} \xar{\pa^{0}} M^{1} \to \cdots 
\bigr) , \]
namely the stupid truncation of $M$ above $0$, and let 
$\rho : M^{-1} \to \opn{Tot}^{\Pi}(M', \pa)$
be the morphism in $\dcat{C}_{\mrm{str}}(A, \cat{M})$ induced by
$\pa^{-1} : M^{-1} \to M^0$. Then there is an isomorphism 
\[ \opn{Tot}^{\Pi}(M, \pa) \cong 
\opn{Cone} \bigl( \rho : M^{-1} \to \opn{Tot}^{\Pi}(M', \pa) \bigr) \]
in $\dcat{C}_{\mrm{str}}(A, \cat{M})$. 
\end{prop}

\begin{proof}
Like the proof of Proposition \ref{prop:3305}, only now there are no sign 
issues. 
\end{proof}

This directly implies:

\begin{cor} \label{cor:3320}
If, in the situation of Proposition \tup{\ref{prop:3320}}, the DG module \lb 
$\opn{Tot}^{\Pi}(M, \pa)$ is acyclic, then 
$\rho : M^{-1} \to \opn{Tot}^{\Pi}(M', \pa)$
is a quasi-isomorphism. 
\end{cor}

In Subsection \ref{subsec:zero} 
we discussed quotients in categories. They are dual to subobjects. Similarly, 
we now introduce a notion dual to a filtration of an object. 

\begin{dfn} \label{dfn:3336}
Let $M \in \dcat{C}(A, \cat{M})$.
\begin{enumerate}
\item A {\em cofiltration}%
\index{Cofiltration on a DG module}
on $M$ is an inverse system
$\{ F_j(M) \}_{j \geq -1}$ in $\dcat{C}_{\mrm{str}}(A, \cat{M})$
consisting of quotients $M \surj F_j(M)$, indexed by the interval 
$[-1, \infty] \sub \Z$. 

\item We say that $M = \lim_{\lar j} F_j(M)$
if the inverse limit exists, and the 
canonical morphism $M \to \lim_{\lar j} F_j(M)$
in $\dcat{C}_{\mrm{str}}(A, \cat{M})$ is an isomorphism.

\item The cofiltration $F$ gives rise to the subquotients 
\[ \opn{Gr}_j^F(M) := \opn{Ker} \bigl( F_j(M) \to F_{j - 1}(M) \bigr) 
\in \dcat{C}(A, \cat{M}) . \]
\end{enumerate}
\end{dfn}

Recall that for a DG $A$-module $M$ we denote by 
$\opn{Y}(M) := M / \opn{B}(M)$; this is the object of decocycles of $M$, and 
it belongs to $\dcat{G}(\K)$.
See Definition \ref{dfn:2993} and equation \ref{eqn:4145}. 

\begin{prop} \label{prop:3321}
Let $(M, \pa)$ be a complex with entries in the abelian category 
$\dcat{C}_{\mrm{str}}(A)$, 
with internal structure \tup{(\ref{eqn:3305})}. Assume that the two conditions 
below hold\tup{:}
\begin{itemize}
\rmitem{i} For every $j \in \Z$ the complex  
\[ \bigl( \cdots \to M^{-1, j} \xar{\pa^{-1}}
M^{0, j} \xar{\pa^{0}} M^{1, j} \to \cdots \bigr) \]
with entries in $\dcat{M}(\K)$ is acyclic. 

\rmitem{ii} There is some $j_0 \in \Z$ such that for all 
$j \leq j_0$ the complex 
\[ \bigl( \cdots \to \opn{Y}^j(M^{-1}) \xar{\opn{Y}^j(\pa^{-1})}
\opn{Y}^j(M^{0}) \xar{\opn{Y}^j(\pa^{0})} \opn{Y}^j(M^{1}) \to \cdots \bigr) \]
with entries in $\dcat{M}(\K)$ is acyclic. 
\end{itemize}
Then the DG $A$-module $\opn{Tot}^{\Pi}(M, \pa)$ is acyclic. 
\end{prop}

\begin{proof} The proof is divided into two steps. 

\smallskip \noindent
Step 1. We replace condition (ii) by the stronger condition: 
\begin{itemize}
\rmitem{ii$'$} There is some $j'_{0} \in \Z$ such that 
$M^{i, j} = 0$ for all $j < j'_{0}$. 
\end{itemize}
For every $i$ we introduce a cofiltration 
$\{ F_q(M^i) \}_{q \geq -1}$
on the DG $\K$-module $M^i$ as follows: 
$F_q (M^i) := \opn{stt}^{\leq j'_{0} + q}(M^i)$, 
the stupid truncation below $j'_{0} + q$. 
These cofiltrations induce a cofiltration on the DG $\K$-module 
$\opn{Tot}^{\Pi}(M, \pa)$~: 
\[ F_q \bigl( \opn{Tot}^{\Pi}(M) \bigr) := 
\opn{Tot}^{\Pi} 
\bigl( \cdots \to F_q(M^{-1}) \xar{\pa^{-1}}
F_q(M^{0}) \xar{\pa^{0}} F_q(M^{1}) \to \cdots \bigr) . \]
We note that $F_{-1} \bigl( \opn{Tot}^{\Pi}(M) \bigr) = 0$, and 
\begin{equation} \label{eqn:3320}
\lim_{\lar q} \, F_q \bigl( \opn{Tot}^{\Pi}(M) \bigr) = 
\opn{Tot}^{\Pi}(M) . 
\end{equation}

For every $q \geq 0$ the DG $\K$-module 
$\opn{Gr}^F_q \bigl( \opn{Tot}^{\Pi}(M) \bigr)$
is isomorphic in $\dcat{C}_{\mrm{str}}(\K)$ to the complex in condition (i) 
with index $j = j'_{0} + q$ (up to a sign change of the differentials), so it 
is acyclic. For $q = -1$ we trivially have an acyclic DG $\K$-module 
$F_{-1} \bigl( \opn{Tot}^{\Pi}(M) \bigr)$. 
For every $q \geq 0$ there is an exact sequence 
\[ 0 \to \opn{Gr}^F_q \bigl( \opn{Tot}^{\Pi}(M) \bigr) \to 
F_q \bigl( \opn{Tot}^{\Pi}(M) \bigr) \to
F_{q - 1} \bigl( \opn{Tot}^{\Pi}(M) \bigr) \to 0 \]
in $\dcat{C}_{\mrm{str}}(\K)$. By induction on $q$ we conclude that 
$F_q \bigl( \opn{Tot}^{\Pi}(M) \bigr)$
is acyclic. Finally, because the inverse system of complexes of $\K$-modules
$(\ref{eqn:3320})$ has surjective transitions, the Mittag-Leffler argument 
(see Corollary \ref{cor:1590}) says that the DG module 
$\opn{Tot}^{\Pi}(M, \pa)$ is acyclic. 

\medskip \noindent 
Step 2. Here we assume conditions (i) and (ii). 
For each $i$ we introduce a cofiltration 
$\{ G_q(M^i) \}_{q \geq 0}$
on the DG $\K$-module $M^i$ as follows: 
\[ G_q(M^i) := \opn{smt}^{\geq j_0 - q}(M^i) = 
\bigl( \cdots \to 0 \to \opn{Y}^{j_0 - q}(M^i) 
\xar{\d_{M^i}} M^{i, j_0 - q + 1} \to \cdots \bigr) , \]
the smart truncation above $j_0 - q$. 

For each $q$ the complex 
\[ G_q(M) := \bigl( \cdots \to G_q(M^{-1}) \xar{G_q(\pa^{-1})}
G_q(M^{0}) \xar{G_q(\pa^{0})} G_q(M^{1}) \to \cdots \bigr) \]
satisfies conditions (i) and (ii$'$), with $j'_{0} := j_0 - q$.
Indeed, its $j$-th row is 
\[ \bigl( \cdots \to M^{-1, j} \xar{\pa^{-1}}
M^{0, j} \xar{\pa^{0}} M^{1, j} \to \cdots \bigr) \]
for $j > j_0 - q$; it is 
\[ \bigl( \cdots \to \opn{Y}^j(M^{-1}) \xar{\opn{Y}^j(\pa^{-1})}
\opn{Y}(M^{0}) \xar{\opn{Y}^j(\pa^{0})} \opn{Y}^j(M^{1}) \to \cdots \bigr) \]
for $j = j_0 - q$; and it is zero for $j < j_0 - q$. All these rows are exact. 
Therefore, by step 1, the DG module 
$G_q \bigl( \opn{Tot}^{\Pi}(M) \bigr) := \opn{Tot}^{\Pi} \bigl( G_q(M) \bigr)$
is acyclic. 

Finally, since $M = \lim_{\lar q} \, G_q(M)$, 
it follows that 
\[ \opn{Tot}^{\Pi}(M) = \lb \lim_{\lar q} \,
G_q \bigl( \opn{Tot}^{\Pi}(M) \bigr) . \]
This is an inverse system of complexes of $\K$-modules with surjective 
transitions. According to the Mittag-Leffler argument the limit
$\opn{Tot}^{\Pi}(M)$ is acyclic.
\end{proof}

\mysubsection{K-Projective Resolutions in 
\texorpdfstring{$\bcat{C}^-(\cat{M})$}{C-(M)}} 
\label{subsec:exis-K-prj}

Recall that $\cat{M}$ is a $\K$-linear abelian category, and 
$\dcat{C}(\cat{M})$ is the DG category of complexes in $\cat{M}$. The strict 
category $\dcat{C}_{\mrm{str}}(\cat{M})$ is abelian. 

The next definition is inspired by the work of B. Keller 
\cite[Section 3.1]{Kel1}. 

\begin{dfn} \label{dfn:1535}
Let $P$ be an object of $\dcat{C}(\cat{M})$. 
\begin{enumerate}
\item A {\em semi-projective filtration}%
\index{Filtration on a DG module! semi-projective}
on $P$ is a filtration 
$F = \{ F_j(P) \}_{j \geq -1}$
on $P$ as an object of $\dcat{C}_{\mrm{str}}(\cat{M})$, such that:
\begin{itemize}
\item $F_{-1}(P) = 0$. 
\item  Each $\opn{Gr}^F_j(P)$ 
is a complex of projective objects of $\cat{M}$ with zero differential.
\item $P = \lim_{j \to} F_j(P)$ in $\dcat{C}_{\mrm{str}}(\cat{M})$.
\end{itemize}

\item The complex $P$ is called a 
{\em semi-projective complex}%
\index{Complex in abelian category! semi-projective}
if it admits some semi-projective filtration. 
\end{enumerate}
\end{dfn}

K-projective objects in $\dcat{C}(\cat{M})$ were introduced in Definition 
\ref{dfn:1520}.

\begin{thm} \label{thm:1536}
\index{Complex in abelian category! semi-projective}
\index{Complex in abelian category! K-projective}
Let $\cat{M}$ be an abelian category, and let $P$ be a 
semi-projective complex in $\dcat{C}(\cat{M})$. Then $P$ is 
K-projective. 
\end{thm}

\begin{proof} The proof is in four steps. 

\medskip \noindent
Step 1. We start by proving that if $P = \opn{T}^k(Q)$, the translation of a 
projective object $Q \in \cat{M}$, then $P$ is K-projective. This is easy: 
given an acyclic complex $N \in \dcat{C}(\cat{M})$, we have 
\[ \opn{Hom}_{\cat{M}}(P, N) = 
\opn{Hom}_{\cat{M}} \bigl( \opn{T}^k(Q), N \bigr) \cong
\opn{T}^{-k} \bigl( \opn{Hom}_{\cat{M}}(Q, N) \bigr) \]
in $\dcat{C}_{\mrm{str}}(\K)$. But 
$\opn{Hom}_{\cat{M}}(Q, -)$ is an exact functor 
$\cat{M} \to \dcat{M}(\K)$, so $\opn{Hom}_{\cat{M}}(Q, N)$ is an acyclic 
complex. 

\medskip \noindent
Step 2. Now $P$ is a complex of projective objects of $\cat{M}$ with zero 
differential. This means that 
$P \cong \bigoplus_{k \in \Z} \, \opn{T}^k(Q_k)$ 
in $\dcat{C}_{\mrm{str}}(\cat{M})$, where each $Q_k$ is a projective object in 
$\cat{M}$. But then 
\[ \opn{Hom}_{\cat{M}}(P, N) \cong 
\prod_{k \in \Z} \, \opn{Hom}_{\cat{M}} \bigl( \opn{T}^k(Q_k), N \bigr)  \]
in $\dcat{C}_{\mrm{str}}(\K)$. 
This is an easy case of Proposition \ref{prop:1535}. By step 1 and the fact 
that a product of acyclic complexes in $\dcat{C}_{\mrm{str}}(\K)$ is acyclic 
(itself an easy case of the Mittag-Leffler argument), we conclude that 
$\opn{Hom}_{\cat{M}}(P, N)$ is acyclic. 

\medskip \noindent
Step 3. Fix a semi-projective filtration $F = \{ F_j(P) \}_{j \geq -1}$
on $P$. Here we prove that for every $j$ the complex  
$F_j(P)$ is K-projective. This is done by induction on $j \geq -1$. 
For $j = -1$ it is trivial. For $j \geq 0$ there is an exact sequence 
\begin{equation} \label{eqn:1555}
0 \to F_{j - 1}(P) \to F_j(P) \to \opn{Gr}^F_j(P) \to 0 
\end{equation}
in $\dcat{C}_{\mrm{str}}(\cat{M})$. 
In each degree $i \in \Z$ the exact sequence 
\[ 0 \to F_{j - 1}(P)^i \to F_j(P)^i \to \opn{Gr}^F_j(P)^i \to 0 \]
in $\cat{M}$ splits, because $\opn{Gr}^F_j(P)^i$ is a projective object.
Thus the exact sequence (\ref{eqn:1555}) is split exact in the abelian category 
$\bcat{G}_{\mrm{str}}(\cat{M})$ of graded objects in $\cat{M}$. 

Let $N \in \dcat{C}(\cat{M})$ be an acyclic complex. Applying the 
functor $\opn{Hom}_{\cat{M}}(-, N)$ 
to the sequence of complexes (\ref{eqn:1555}) we obtain a sequence 
\begin{equation} \label{eqn:1556}
0 \to \opn{Hom}_{\cat{M}} \bigl( \opn{Gr}^F_j(P), N \bigr) \to 
\opn{Hom}_{\cat{M}} \bigl( F_j(P), N \bigr) \to
\opn{Hom}_{\cat{M}} \bigl( F_{j - 1}(P), N \bigr) \to  0 
\end{equation}
in $\dcat{C}_{\mrm{str}}(\K)$. Because (\ref{eqn:1555}) is split exact in 
$\bcat{G}_{\mrm{str}}(\cat{M})$, the sequence (\ref{eqn:1556}) is split exact 
in $\bcat{G}_{\mrm{str}}(\K)$. Therefore (\ref{eqn:1556}) is exact in 
$\dcat{C}_{\mrm{str}}(\K)$.

By the induction hypothesis the complex 
$\opn{Hom}_{\cat{M}} \bigl( F_{j - 1}(P), N \bigr)$
is acyclic. By step 2 the complex
$\opn{Hom}_{\cat{M}} \bigl( \opn{Gr}^F_{j}(P), N \bigr)$
is acyclic. The long exact cohomology sequence associated to 
(\ref{eqn:1556}) shows that the complex 
$\opn{Hom}_{\cat{M}} \bigl( F_j(P), N \bigr)$
is acyclic too.

\medskip \noindent
Step 4. We keep the semi-projective filtration 
$F = \{ F_j(P) \}_{j \geq -1}$ from step 3. 
Take any acyclic complex $N \in \bcat{C}(\cat{M})$. By 
Proposition \ref{prop:1535} we know that 
\[ \opn{Hom}_{\cat{M}}(P, N) \cong 
\lim_{\leftarrow j} \, \opn{Hom}_{\cat{M}} \bigl( F_j(P), N \bigr) \]
in $\dcat{C}_{\mrm{str}}(\K)$. According to step 3 the complexes 
$\opn{Hom}_{\cat{M}} \bigl( F_j(P), N \bigr)$ are all acyclic. 
The exactness of the sequences (\ref{eqn:1556}) implies that 
the inverse system \lb 
$\bigl\{ \opn{Hom}_{\cat{M}} \bigl( F_j(P), N \bigr) \bigr\}_{j \geq -1}$
in $\dcat{C}_{\mrm{str}}(\K)$ has surjective transitions. Now the 
Mittag-Leffler 
argument (Corollary \ref{cor:1590}) says that the inverse limit complex
$\opn{Hom}_{\cat{M}}(P, N)$ is acyclic. 
\end{proof}

\begin{prop} \label{prop:1900}
Let $\cat{M}$ be an abelian category. If $P \in \dcat{C}(\cat{M})$ is a bounded 
above complex of projectives, then $P$ is a semi-projective complex.
\end{prop}

\begin{proof}
Say $P$ is nonzero and $\opn{sup}(P) = i_1 \in \Z$. 
For $j \geq -1$ define 
$F_j(P) := \opn{stt}^{\geq i_1 - j}(P)$,
the stupid truncation above $i_1 - j$ from Definition \ref{dfn:3066}.
Then $\{ F_j(P) \}_{j \geq -1}$ is a semi-projective filtration on $P$. 
\end{proof}

The next theorem is dual to \cite[Lemma I.4.6(1)]{RD}, in the sense of 
changing direction of arrows. We give a detailed proof. 

\begin{thm} \label{thm:1540}
Let $\cat{M}$ be an abelian category, and let $\cat{P} \sub \cat{M}$ be a
full subcategory such that every object $M \in \cat{M}$ admits an 
epimorphism $P \surj M$ from some object $P \in \cat{P}$. Then every complex 
$M \in \bcat{C}^{-}(\cat{M})$ admits a quasi-isomorphism 
$\rho : P \to  M$ in $\dcat{C}^-_{\mrm{str}}(\cat{M})$,
such that $\opn{sup}(P) = \opn{sup}(M)$, and 
each $P^i$ is an object of $\cat{P}$.
\end{thm}

\begin{proof}
This is divided into steps. 

\smallskip \noindent
Step 1. We can assume that $M \neq 0$. 
After translating $M$, we can assume that $\opn{sup}(M) = 0$. 

\medskip \noindent
Step 2. Say the differential of the complex $M$ is 
$\d^{i}_M : M^i \to M^{i + 1}$.
Let us choose an epimorphism $\rho^{0} : P^{0} \surj M^{0}$ in $\cat{M}$ 
from some object $P^{0} \in \cat{P}$. We get a morphism 
$\de^{0} :  M^{-1} \oplus P^0 \to M^0$
whose components are $\d^{-1}_M$ and $-\rho^{0}$. Thus 
$\opn{Ker}(\de^{0}) = M^{-1} \times_{M^{0}} P^0$,
the fibered product. 
Next we choose an epimorphism
$\psi^{-1} : P^{-1} \surj \opn{Ker}(\de^{0})$  
from some object $P^{-1} \in \cat{P}$.
So there is an exact sequence 
\begin{equation} \label{eqn:5060}
P^{-1} \xar{\psi^{-1}} M^{-1} \oplus P^0 \xar{\de^0} M^0 \to 0  .
\end{equation}
The components of $\psi^{-1}$ are denoted by 
$\rho^{-1} : P^{-1} \to  M^{-1}$ and 
$\d_P^{-1} : P^{-1} \to  P^{0}$.
We have this commutative diagram 
\begin{equation} \label{eqn:5058}
\UseTips \xymatrix @C=6ex @R=6ex {
&
P^{-1}
\ar[d]^{\rho^{-1}}
\ar[r]^{\d_P^{-1}}
&
P^0 
\ar[d]^{\rho^0}
\ar[r]
&
0
\\
M^{-2}
\ar[r]^{\d_M^{-2}}
&
M^{-1}
\ar[r]^{\d_M^{-1}}
&
M^0 
\ar[r]
&
0
} 
\end{equation}
in $\cat{M}$.

\medskip \noindent
Step 3. This is the inductive step. Here $i \leq -1$, and we already have 
objects $P^{i}, \ldots, P^0$ in $\cat{P}$, and morphisms 
$\rho^{i}, \ldots, \rho^0$ and $\d_P^{i}, \ldots, \d_P^{-1}$,
that fit into this commutative diagram
\begin{equation} \label{eqn:1745}
\UseTips \xymatrix @C=6ex @R=6ex {
&
P^i
\ar[d]^{\rho^i}
\ar[r]^{\d_P^i}
&
P^{i + 1}
\ar[d]^{\rho^{i + 1}}
\ar[r] 
&
\cdots 
\ar[r]^{\d_P^{-1}}
& 
P^0 
\ar[d]^{\rho^0}
\ar[r]
&
0
\\
M^{i - 1}
\ar[r]^{\d_M^{i - 1}}
&
M^i
\ar[r]^{\d_M^{i}}
&
M^{i + 1}
\ar[r] 
&
\cdots 
\ar[r]^{\d_M^{-1}}
& 
M^0 
\ar[r]
&
0
} 
\end{equation}
in $\cat{M}$. Also 
$\d_P^{j + 1} \circ \d_P^{j} = 0$ for all $j$ in the interval $[i, -1]$

Define the morphism 
$\de^{i}  :  M^{i - 1} \oplus P^{i} \to  M^{i} \oplus P^{i + 1}$
to be the one with components 
$\d_M^{i - 1}$, $- \rho^i$ and $- \d_P^{i}$.
Expressing direct sums of objects as columns, and letting matrices of morphisms 
act on them from the left, we have this matrix representation of $\de^i$~:
\begin{equation} \label{eqn:1905}
\de^i = \bmat{ \d_M^{i - 1} & - \rho^i \\[0.3em] 0 & - \d_P^i } .
\end{equation}
So 
$\opn{Ker}(\de^{i}) = M^{i - 1} \times_{(M^i \oplus P^{i + 1})} P^{i}$,
the fibered product. 
 
Let us choose an epimorphism
$\psi^{i - 1} : P^{i - 1} \surj  \opn{Ker}(\de^{i})$ 
from an object $P^{i - 1} \in \cat{P}$. 
We get an exact sequence 
\begin{equation} \label{eqn:1750}
P^{i - 1} \xar{\psi^{i - 1}} M^{i - 1} \oplus P^{i} \xar{\de^i}
M^{i} \oplus P^{i + 1} .
\end{equation}
The components of the morphism $\psi^{i - 1}$ are denoted by 
$\rho^{i - 1} : P^{i - 1} \to  M^{i - 1}$ and 
$\d_P^{i - 1} : P^{i - 1} \to  P^{i}$.
In matrix representation is 
\begin{equation} \label{eqn:5061}
\psi^{i - 1} = \bmat{ \rho^{i - 1} \\[0.3em] \d_P^{i - 1} } . 
\end{equation}
In this way we obtain the slightly bigger diagram than (\ref{eqn:1745}):
\begin{equation} \label{eqn:1751}
\UseTips \xymatrix @C=6ex @R=6ex {
P^{i - 1}
\ar[d]^{\rho^{i - 1}}
\ar[r]^{\d_P^{i - 1}}
&
P^i
\ar[d]^{\rho^i}
\ar[r]^{\d_P^i}
&
P^{i + 1}
\ar[d]^{\rho^{i + 1}}
\ar[r] 
&
\cdots 
\ar[r]
& 
P^0 
\ar[d]^{\rho^0}
\ar[r]
&
0
\\
M^{i - 1}
\ar[r]^{\d_M^{i - 1}}
&
M^i
\ar[r]^{\d_M^{i}}
&
M^{i + 1}
\ar[r] 
&
\cdots 
\ar[r]
& 
M^0 
\ar[r]
&
0
} 
\end{equation}
Because $\de^i \circ \psi^{i - 1} = 0$ in (\ref{eqn:1750}), it follows that
$\d_P^{i} \circ \d_P^{i - 1} = 0$, and also 
\begin{equation} \label{eqn:1753}
\rho^{i} \circ \d_P^{i - 1} = \d_M^{i - 1}  \circ \rho^{i - 1} ,
\end{equation}
showing that diagram (\ref{eqn:1751}) is commutative. 

\medskip \noindent
Step 4. We carry out the construction in step 3 inductively for all $i \leq -1$, 
thus obtaining a diagram like (\ref{eqn:1751}) that goes infinitely to the 
left. 

Letting $P^i := 0$ for positive $i$, the 
collection $P := \{ P^i \}_{i \in \Z}$ becomes a complex, with 
differential $\d_P := \{ \d_P^{i} \}_{i \in \Z}$.
By equation (\ref{eqn:1753}), the collection 
$\rho := \{ \rho^{i} \}_{i \in \Z}$ is a strict morphism 
of complexes $\rho : P \to M$. 

\medskip \noindent
Step 5. It remains to prove that $\rho$ is a quasi-isomorphism.

Take any $i \leq 0$. Let us examine this diagram: 
\begin{equation} \label{eqn:1752}
\UseTips \xymatrix @C=6ex @R=6ex {
P^{i - 1}
\ar[r]^{- \psi^{i - 1}}
\ar[d]^{(0, \opn{id})}
&
M^{i - 1} \oplus P^{i}
\ar[r]^{\de^{i}}
\ar[d]^{\opn{id}}
&
M^{i} \oplus P^{i + 1}
\ar[d]^{\opn{id}}
\\
M^{i - 2} \oplus P^{i - 1}
\ar[r]^{\de^{i - 1}}
&
M^{i - 1} \oplus P^{i}
\ar[r]^{\de^{i}}
&
M^i \oplus P^{i+ 1}
}
\end{equation}
Comparing the formula for $\de^{i - 1}$ in (\ref{eqn:1905}) to the 
formula for $- \psi^{i - 1}$ in (\ref{eqn:5061}), we see that this diagram is 
commutative. 
An easy calculation using (\ref{eqn:1753}) shows that 
$\de^{i} \circ \de^{i - 1} = 0$. 
The top row is exact, because (up to signs) for $i = 0$ it is part of the 
exact sequence (\ref{eqn:5060}), and  for $i \leq -1$ it is the exact sequence 
(\ref{eqn:1750}). 
It follows that the bottom row is also exact.

Let $N =  \{ N^i \}_{i \in \Z}$ be the complex with  components 
$N^i := M^{i - 1} \oplus P^{i}$ for all $i$. 
The differential
$\d_N= \{ \d_N^{i} \}_{i \in \Z}$
is $\d^i_N := -\de^{i} : N^i \to N^{i + 1}$.
We have $N^i = 0$ for $i \geq 2$. 
The exactness of the sequence (\ref{eqn:5060}) says that 
$\opn{H}^1(N) = 0$, and the exactness of the second row in 
(\ref{eqn:1752}) says that $\opn{H}^i(N) = 0$ for $i \leq 0$. 
Hence the complex $N$ is acyclic. 
On the other hand, by the definition of the morphisms 
$\de^i$ in (\ref{eqn:1905}), we see that $N$ is just the standard cone on the 
strict morphism of complexes
$\opn{T}^{-1}(\rho) : \opn{T}^{-1}(P) \to \opn{T}^{-1}(M)$.
See Definition \ref{dfn:1172}. 
Therefore $\opn{T}^{-1}(\rho)$ is a quasi-isomorphism, and so is $\rho$.
\end{proof}

\begin{exa} \label{exa:4725}
Let $(X, \AA)$ be a ringed space. Usually there aren't enough projective 
objects in the abelian category $\dcat{M}(\AA) = \cat{Mod} \AA$ of sheaves of 
$\AA$-modules on $X$; but there are always enough flat objects. A general method 
to produce flat sheaves is as follows. Given an open set $U \sub X$ with 
inclusion map $g : U \to X$,  
the extension by zero sheaf $g_{!}(\AA|_{U})$ is a flat $\AA$-module.  
Let us call $g_{!}(\AA|_{U})$ a pseudo-free $\AA$-module of pseudo-rank 
$1$ and pseudo-support $U$ . A {\em pseudo-free $\AA$-module} is, by 
definition, a sheaf $\PP$ which is a direct sum of rank $1$ pseudo-free 
$\AA$-modules (with varying pseudo-supports). A pseudo-free $\AA$-module is 
flat, and every $\AA$-module $\MM$ admits an epimorphism
$\PP \surj \MM$ from some pseudo-free $\AA$-module $\PP$. 
More on this in \cite[Proposition II.1.2]{RD} or 
\cite[Proposition 2.4.12]{KaSc1}. 
(Warning: in case $(X, \OO_X)$ is a scheme, the pseudo-free $\OO_X$-modules are 
usually not quasi-coherent.) 

Now let $\dcat{C}(\AA) =  \dcat{C}(\cat{Mod} \AA)$ be the category of complexes 
of $\AA$-modules. Because $\dcat{M}(\AA)$ has enough flat objects, 
Theorem \ref{thm:1540} tells us that every complex 
$\MM \in \dcat{C}^{-}(\AA)$ admits a quasi-isomorphism 
$\PP \to \MM$ from a bounded above complex of flat 
$\AA$-modules $\PP$. Similarly to Proposition \ref{prop:1900} and Theorem 
\ref{thm:1536}, one can show that such a complex $\PP$ is K-flat. 
The conclusion is that the subcategory 
$\dcat{C}^{-}(\AA)$ has enough K-flat objects. 

But in fact more is true: the category of unbounded complexes $\dcat{C}^{}(\AA)$ 
also has enough K-flat objects; see \cite{Spa}. 
For an even more general statement see Remark \ref{rem:4675} below. 
\end{exa}

\begin{dfn} \label{dfn:2895}
Let $\cat{M}$ be an abelian category, and let $\cat{M}' \sub \cat{M}$ be a 
full abelian subcategory. 
We say that {\em $\cat{M}'$ has enough projectives relative to $\cat{M}$} 
if every object $M \in \cat{M}'$ admits an epimorphism 
$P \surj M$, where $P$ is an object of $\cat{M}'$ that is projective in the 
bigger category $\cat{M}$.
\end{dfn}

Of course, if $\cat{M}'$ has enough projectives relative to $\cat{M}$
then $\cat{M}'$ itself has enough projectives. 

Thick abelian categories were defined in Definition \ref {dfn:2323}.
The next theorem is the opposite of \cite[Lemma I.4.6(3)]{RD}. 

\begin{thm} \label{thm:2895}
Let $\cat{M}$ be an abelian category, and let $\cat{M}' \sub \cat{M}$ be a 
thick abelian subcategory that has enough projectives relative to $\cat{M}$. 
Let $M \in \dcat{C}(\cat{M})$ be a complex with bounded above cohomology, 
such that $\opn{H}^i(M) \in \cat{M}'$ for all $i$.
Then there is a quasi-isomorphism $\rho : P \to M$ in 
$\dcat{C}_{\mrm{str}}(\cat{M})$, where 
$P \in \dcat{C}^-(\cat{M}')$, each $P^i$ is projective in $\cat{M}$,
and $\opn{sup}(P) = \opn{sup}(\opn{H}(M))$. 
\end{thm}

\begin{proof}
The proof of Theorem \ref{thm:2880}, reversed, works here. 
(Despite being located later in the book, Theorem \ref{thm:2880} is 
not logically dependent on the current theorem.)
To be explicit, let us take $\cat{N} := \cat{M}^{\mrm{op}}$ and 
$\cat{N}' := \cat{M}'^{\, \mrm{op}}$.
Since monomorphisms in $\cat{M}$ become epimorphisms in $\cat{N}$, 
and projective objects in $\cat{M}$ become injective objects in $\cat{N}$,
the full abelian subcategory $\cat{N}' \sub \cat{N}$ satisfies 
assumptions of Theorem \ref{thm:2880}. By Theorem \ref{thm:2495} we have a 
canonical isomorphism of categories 
$\bcat{C}^{+}_{\mrm{str}}(\cat{N}) \iso 
\dcat{C}^-_{\mrm{str}}(\cat{M})^{\mrm{op}}$
that respects the cohomology functor. 
Thus a quasi-isomorphism $N \to J$ in 
$\bcat{C}^{+}_{\mrm{str}}(\cat{N})$
gives rise to a quasi-isomorphism $P \to M$ in 
$\bcat{C}^{-}_{\mrm{str}}(\cat{M})$.
\end{proof}

\begin{cor} \label{cor:1900}
\index{Resolution! K-projective}
If $\cat{M}$ is an abelian category with enough projectives, then  
$\bcat{C}^{-}(\cat{M})$ has enough K-projectives. 
\end{cor}

\begin{proof}
According to either Theorem \ref{thm:1540} or Theorem \ref{thm:2895},
every $M \in \bcat{C}^{-}(\cat{M})$ admits a quasi-isomorphism $P \to M$
from a bounded above complex of projectives $P$. Now use Proposition 
\ref{prop:1900} and Theorem \ref{thm:1536}. 
\end{proof}

\begin{cor} \label{cor:1750}
Let $\cat{M}$ be an abelian category with enough projectives, 
and let $M \in \dcat{C}(\cat{M})$ be a complex with bounded above cohomology. 
Then $M$ has a K-projective resolution $P \to M$,
such that $\opn{sup}(P) = \opn{sup}(\opn{H}(M))$,
and every $P^i$ is an projective object of $\cat{M}$. 
\end{cor}

\begin{proof}
We may assume that $\opn{H}(M)$ is not zero. 
Let $i := \opn{sup}(\opn{H}(M)) \in \Z$, and take $N := \opn{smt}^{\leq i}(M)$,
the smart truncation from  Definition \ref{dfn:2320}. Then 
$N \to M$ is a quasi-isomorphism and $\opn{sup}(N) = i$.
According to either Theorem \ref{thm:1540} or Theorem \ref{thm:2895},
there is a quasi-isomorphism $P \to N$, where $P$ is a complex 
of projectives and $\opn{sup}(P) = i$. 
By Proposition \ref{prop:1900} and Theorem \ref{thm:1536} the complex $P$ is 
K-projective. The composed quasi-isomorphism $P \to M$ is what we are looking 
for.
\end{proof}

\begin{cor} \label{cor:3175}
Under the assumptions of Theorem \tup{\ref{thm:2895}}, the canonical functor 
$\dcat{D}^-(\cat{M}') \to \dcat{D}^-_{\cat{M}'}(\cat{M})$
is an equivalence. 
\end{cor}

\begin{proof}
Consider the commutative diagram 
\[ \UseTips \xymatrix @C=6ex @R=6ex {
\dcat{K}^-(\cat{M}')_{\catt{M}\tup{-}\mrm{prj}}
\ar@(l,u)[drr]^{F}
\ar[d]_{\opn{Q}}
&
\\
\dcat{D}^-(\cat{M}')
\ar[r]^{G}
& 
\dcat{D}^-_{\cat{M}'}(\cat{M})
\ar[r]^{\tup{f.f.}}
& 
\dcat{D}(\cat{M})
} \]
where $\dcat{K}^-(\cat{M}')_{\catt{M}\tup{-}\mrm{prj}}$ is the homotopy 
category of bounded above complexes of objects of $\cat{M}'$ that are 
projective in $\cat{M}$. Since these are K-projective complexes in 
$\dcat{K}(\cat{M})$, the functors $F$ and $\opn{Q}$ are fully faithful, by 
Corollary \ref{cor:3145}. The functor marked ``f.f.'' is the fully faithful 
inclusion of this full subcategory. We conclude that the functor $G$ is also 
fully faithful. The theorem tells us that the essential 
image of $G$ is $\dcat{D}^-_{\cat{M}'}(\cat{M})$. 
\end{proof}

Here is an important instance where Theorem \ref{thm:2895} and  
Corollary \ref{cor:3175} apply. 

\begin{exa} \label{exa:2900}
Let $A$ be a left noetherian ring. Consider the abelian category 
$\cat{M} := \dcat{M}(A) = \cat{Mod} A$, and the thick abelian subcategory 
$\cat{M}' := \dcat{M}_{\mrm{f}}(A) = \cat{Mod}_{\mrm{f}} A$ of finitely 
generated modules. Then $\cat{M}'$ has enough projectives relative to 
$\cat{M}$. 
Theorem \ref{thm:2895} tells us that if 
$M \in \dcat{C}(\cat{M})$ is a complex such that the modules 
$\opn{H}^i(M)$ are all finitely generated, and $\opn{H}^i(M) = 0$ for 
$i \gg 0$, then there is a resolution $P \to M$, where $P$ is a bounded above 
complex of finitely generated projective modules. 
Corollary \ref{cor:3175} says that the canonical functor 
\[ \dcat{D}^{-}(\dcat{M}_{\mrm{f}}(A)) = \dcat{D}^{-}(\cat{M}') \to 
\dcat{D}^{-}_{\cat{M}'}(\cat{M}) = \dcat{D}^{-}_{\mrm{f}}(\dcat{M}(A))
= \dcat{D}^{-}_{\mrm{f}}(A) \]
is an equivalence. 
See Example \ref{exa:1930} for another approach to this problem. 
\end{exa}

We end this subsection with a result on K-flat complexes over rings. 

\begin{prop} \label{prop:4535}
Let $A$ be a ring and $P$ a bounded above complex of flat $A$-modules. 
Then $P$ is a K-flat complex. 
\end{prop}

\begin{exer} \label{exer:4535}
Prove this proposition. (Hint: see the proofs of Proposition \ref{prop:1900} 
and Theorem \ref{thm:1536}.)
\end{exer}

\newpage
\mysubsection{K-Projective Resolutions in 
\texorpdfstring{$\bcat{C}(A)$}{C(A)}}
\label{subsec:exis-K-prj-dgmods}

In this subsection $A$ is a DG $\K$-ring (without any vanishing assumption). 

Let $A^{\natural}$ be the graded ring that we have after forgetting the 
differential of $A$. In Subsection \ref{subsec:DGModinM} we introduced the 
functor 
$\opn{Und} : \dcat{C}(A) \to \dcat{G}(A^{\natural})$
that forgets the differentials of DG modules. We shall now use the abbreviated 
form $M^{\natural}  := \opn{Und}(M)$. 

Recall that the translation $\opn{T}^{-i}(A)$ is a DG $A$-module in which the 
element $\opn{t}^{-i}(1_A)$ is in degree $i$. 
This element is a cocycle, and when we forget the differentials, the graded 
module $\opn{T}^{-i}(A)^{\natural}$ is free over the graded ring 
$A^{\natural}$, with basis $\opn{t}^{-i}(1_A)$. Therefore, for every DG 
$A$-module $M$ there is a canonical isomorphism
\begin{equation} \label{eqn:1755}
\opn{Hom}_{A} \bigl( \opn{T}^{-i}(A), M \bigr) \cong \opn{T}^{i}(M)
\end{equation}
in $\dcat{C}_{\mrm{str}}(\K)$, and canonical isomorphisms 
\begin{equation} \label{eqn:1600}
\opn{Hom}_{\dcat{C}_{\mrm{str}}(A)} \bigl( \opn{T}^{-i}(A), M \bigr) \cong
\opn{Z}^0  \bigl( \opn{Hom}_{A} ( \opn{T}^{-i}(A), M ) \bigr) \cong \opn{Z}^i(M)
\end{equation}
in $\dcat{M}(\K)$. 
(Actually, (\ref{eqn:1755}) is an isomorphism in $\dcat{C}_{\mrm{str}}(A)$,
but this uses the DG $A$-bimodule structure of $\opn{T}^{-i}(A)$.)

We begin with a definition that is very similar to Definition \ref{dfn:1535}. 

\begin{dfn} \label{dfn:1575}
Let $P$ be an object of $\dcat{C}(\cat{A})$. 
\begin{enumerate}
\item We say that $P$ is a {\em free DG $A$-module}%
\index{Differential graded module! free}
if there is an isomorphism
$P \cong \bigoplus_{s \in S} \opn{T}^{-i_s}(A)$
in $\dcat{C}_{\mrm{str}}(A)$, 
for some indexing set $S$ and some collection of integers 
$\{ i_s \}_{s \in S}$. 

\item A {\em semi-free filtration}%
\index{Filtration on a DG module! semi-free}
on $P$ is a filtration 
$F = \{ F_j(P) \}_{j \geq -1}$
of $P$ in $\dcat{C}_{\mrm{str}}(A)$, such that:
\begin{itemize}
\item $F_{-1}(P) = 0$. 
\item  Each $\opn{Gr}^F_j(P)$ is a free DG $A$-module.
\item $P = \lim_{j \to} F_j(P)$.
\end{itemize}

\item The DG module $P$ is called {\em semi-free}%
\index{Differential graded module! semi-free}
if it admits some semi-free filtration. 
\end{enumerate}
\end{dfn}

Note that the direct limit in item (2) is just the union of DG submodules.

There is an alternative description of semi-free filtrations, that is useful 
for understanding their structure. We will give it in Definition 
\ref{dfn:4732} below, after some combinatorial preliminaries. 

The next definition has already appeared within Example \ref{exa:4625}.

\begin{dfn} \label{dfn:4730}
A {\em graded set}%
\index{Graded set}
is a set $S$ that is partitioned into subsets: 
$S = \coprod_{i \in \Z} \, S^i$. 
The elements of $S^i$ are said to have degree $i$,
so that $S^i = \{ s \in S \mid \opn{deg}(s) = i \}$.
\end{dfn}

A free DG $A$-module $P$ as in Definition \ref{dfn:1575}(1) can be described as 
follows. The indexing set $S$ can be made into a graded set by defining
$S^k := \{ s \in S \mid i_s = k \}$.
In other words, $\opn{deg}(s) = i_s$. We denote by $A \cd s$ the free DG 
$A$-module with basis $s$, so that $\d(s) = 0$, and there is an isomorphism 
$A \cd s \iso \opn{T}^{-i_s}(A)$ in $\dcat{C}_{\mrm{str}}(A)$
sending $s \mapsto \opn{t}^{-i_s}(1_A)$. Writing 
\begin{equation} \label{eqn:5062}
A \cd S := \bigoplus_{s \in S} \ A \cd s
\end{equation}
we obtain $A \cd S \cong P$ as DG $A$-modules. 
Note that $A^{\natural} \cd S \cong P^{\natural}$
as graded $A^{\natural}$-modules. 

\begin{dfn} \label{dfn:4731}
Let $S$ be a graded set. A {\em filtration} on $S$ is a direct system
$F = \bigl\{ F_{j}(S) \bigr\}_{j \geq -1}$
of subsets of $S$, such that 
$F_{-1}(S) = \varnothing$ and $\bigcup_{j} F_{j}(S) = S$.  
The pair $(S, F)$ is called a {\em filtered graded set}%
\index{Graded set! filtered}.
\end{dfn}

Note that each $F_{j}(S)$ is itself a graded set, with degree $i$ component
$F_{j}(S)^i := F_{j}(S) \cap S^i$.

\begin{dfn} \label{dfn:4732}
Let $P$ be a DG $A$-module. A 
{\em semi-basis}%
\index{Semi-basis! of DG module}
of $P$ is a filtered graded set 
$(S, F)$, together with an isomorphism 
$P^{\natural}  \cong A^{\natural} \cd S$
of graded $A^{\natural}$-modules, such that under this isomorphism we have 
$\d_P(F_j(S)) \sub A^{\natural} \cd F_{j - 1}(S)$
for every $j \geq 0$. 
\end{dfn}

\begin{prop} \label{prop:4730}
Let $P$ be a DG $A$-module. 
\begin{enumerate}
\item If $P$ has a semi-basis $(S, F)$, then $P$ has an induced semi-free 
filtration $\{ F_j(P) \}_{j \geq -1}$ such that 
$F_j(P)^{\natural} = A^{\natural} \cd F_{j}(S)$
as graded $A^{\natural}$-modules for every $j \geq -1$.

\item Every semi-free filtration of $P$ is induced by some semi-basis. 
\end{enumerate}
\end{prop}

\begin{exer} \label{exer:4725}
Prove Proposition \ref{prop:4730}.
\end{exer}

\begin{exa} \label{exa:1575}
If $A$ is a ring, then a free DG $A$-module $P$ is a complex of free $A$-modules
with zero differential. A semi-free DG $A$-module $P$ is also a  complex of 
free $A$-modules, but there is a differential on it, and there is a subtle 
condition on $P$ imposed by the existence of a semi-free filtration. If the 
complex $P$ happens to be bounded above, then it is automatically semi-free, 
with a filtration like the one in the proof of Proposition \ref{prop:1900}. 
\end{exa}

\begin{exer} \label{exa:1577}
Find a ring $A$, and a complex $P$ of free $A$-modules, that is not semi-free. 
(Hint: Take the ring $A = \K[\ep]$ of dual numbers. Find a complex of free 
$A$-modules $P$ that is acyclic but not null-homotopic. Now use Theorem 
\ref{thm:1575} and Corollary \ref{cor:1923}  to deduce that $P$ is not 
semi-free.)
\end{exer}

\begin{prop} \label{prop:4636}
Let $A$ and $B$ be DG rings, and let $P \in \dcat{C}(A)$ and 
$Q \in \dcat{C}(B)$ be semi-free DG modules. Then 
$P \ot Q \in \dcat{C}(A \ot B)$
is a semi-free DG module.
\end{prop}

\begin{exer} \label{exer:4636}
Prove this proposition. (Hint: let $\bigl\{ F_j(P) \bigr\}_{j \geq -1}$
and \lb $\bigl\{ G_k(Q) \bigr\}_{k \geq -1}$ be semi-free filtrations. 
Define 
\[ E_{j}(P \ot Q) := \sum_{k = 0}^{j} \, 
\bigl( F_k(P) \ot G_{j - k}(Q) \bigr) \sub P \ot Q . \]
Show that this is a semi-free filtration.) 
\end{exer}

\begin{thm} \label{thm:1575}
\index{Differential graded module! K-projective}
\index{Differential graded module! semi-free}
Let $P$ be an object of $\dcat{C}(A)$. If $P$ is semi-free, then it is 
K-projective. 
\end{thm}

\begin{proof}
This is similar to the proof of Theorem \ref{thm:1536}. 

\smallskip \noindent
Step 1. We start by proving that if $P = \opn{T}^{-i}(A)$, a translation of 
$A$, 
then $P$ is K-projective. This is easy: 
given an acyclic $N \in \dcat{C}(A)$, we have 
\[ \opn{Hom}_{A}(P, N) = 
\opn{Hom}_{A} \bigl( \opn{T}^{-i}(A), N \bigr) \cong
\opn{T}^{i} \bigl( \opn{Hom}_{A}(A, N) \bigr)
\cong \opn{T}^{i}(N)\]
in $\dcat{C}_{\mrm{str}}(\K)$, and this is acyclic. 

\medskip \noindent
Step 2. Now $P$ is free, say 
$P \cong \bigoplus_{s \in S} \, \opn{T}^{-i_s}(A)$.
Then 
\[ \opn{Hom}_{A}(P, N) \cong 
\prod_{s \in S} \, \opn{Hom}_{A} \bigl( \opn{T}^{-i_s}(A), N \bigr) . \]
By step 1 and the fact that a product of acyclic complexes in 
$\dcat{C}_{\mrm{str}}(\K)$ is 
acyclic, we conclude that $\opn{Hom}_{\cat{M}}(P, N)$ is acyclic. 

\medskip \noindent
Step 3. Fix a semi-free filtration $F = \{ F_j(P) \}_{j \geq -1}$
of $P$. Here we prove that for every $j \geq -1$ the DG module  
$F_j(P)$ is K-projective. This is done by induction on $j \geq -1$. 
For $j = -1$ it is trivial. For $j \geq 0$ there is an exact sequence 
\begin{equation} \label{eqn:1580}
0 \to F_{j - 1}(P) \to F_j(P) \to \opn{Gr}^F_j(P) \to 0 
\end{equation}
in the abelian category $\dcat{C}_{\mrm{str}}(A)$. 
Because 
$\opn{Gr}^F_j(P)$
is a free DG module, it is a projective object in the abelian category 
$\bcat{G}_{\mrm{str}}(A^{\natural})$ of graded modules over the graded ring 
$A^{\natural}$. 
Therefore the sequence (\ref{eqn:1580}) is split exact in 
$\bcat{G}_{\mrm{str}}(A^{\natural})$.

Let $N \in \dcat{C}(A)$ be an acyclic DG module. Applying the 
functor $\opn{Hom}_{A}(-, N)$ 
to the sequence (\ref{eqn:1580}) we obtain a sequence 
\begin{equation} \label{eqn:1581}
0 \to \opn{Hom}_{A} \bigl( \opn{Gr}^F_j(P), N \bigr) \to 
\opn{Hom}_{A} \bigl( F_j(P), N \bigr) \to
\opn{Hom}_{A} \bigl( F_{j - 1}(P), N \bigr) \to  0 
\end{equation}
in $\dcat{C}_{\mrm{str}}(\K)$. If we forget differentials this is a sequence in 
$\bcat{G}_{\mrm{str}}(\K)$. Because (\ref{eqn:1580}) is split exact in 
$\bcat{G}_{\mrm{str}}(A^{\natural})$, it follows that  (\ref{eqn:1581}) is 
split 
exact in $\bcat{G}_{\mrm{str}}(\K)$. Therefore (\ref{eqn:1581}) is exact in 
$\dcat{C}_{\mrm{str}}(\K)$.

By the induction hypothesis the DG $\K$-module  
$\opn{Hom}_{A} \bigl( F_{j - 1}(P), N \bigr)$
is acyclic. By step 2 the DG module
$\opn{Hom}_{A} \bigl( \opn{Gr}^F_{j}(P), N \bigr)$
is acyclic. The long exact cohomology sequence associated to 
(\ref{eqn:1581}) shows that the DG module 
$\opn{Hom}_{A} \bigl( F_j(P), N \bigr)$
is acyclic too.

\medskip \noindent
Step 4. We keep the semi-free filtration 
$F = \{ F_j(P) \}_{j \geq -1}$ from step 3. 
Take any acyclic $N \in \bcat{C}(\cat{M})$. By 
Proposition \ref{prop:1535} we know that 
\[ \opn{Hom}_{A}(P, N) \cong 
\lim_{\leftarrow j} \, \opn{Hom}_{A} \bigl( F_j(P), N \bigr) \]
in $\dcat{C}_{\mrm{str}}(\K)$. According to step 3 the complexes 
$\opn{Hom}_{A} \bigl( F_j(P), N \bigr)$ are all acyclic. 
The exactness of the sequences (\ref{eqn:1581}) implies that 
the inverse system \lb 
$\bigl\{ \opn{Hom}_{A} \bigl( F_j(P), N \bigr) \bigr\}_{j \geq -1}$
in $\dcat{C}_{\mrm{str}}(\K)$ has surjective transitions. Now the 
Mittag-Leffler argument (Corollary \ref{cor:1590}) says that the inverse limit 
complex $\opn{Hom}_{A}(P, N)$ is acyclic. 
\end{proof}

Here is a result similar to Theorem \ref{thm:1540}. 

\begin{thm} \label{thm:3305}
\index{Differential graded module! semi-free}
Let $A$ be a DG ring. Every $M \in \bcat{C}(A)$ admits a quasi-iso\-morphism
$\rho : P \to M$ in $\dcat{C}_{\mrm{str}}(A)$ from a semi-free DG $A$-module 
$P$.
\end{thm}

The proof of this theorem was communicated to us by B. Keller.
First we need a construction and a lemma. 

Let us define the DG $A$-module 
\begin{equation} \label{eqn:3307}
C := \opn{Cone} \bigl( \opn{id} : \opn{T}^{-1}(A) \to \opn{T}^{-1}(A) \bigr) , 
\end{equation}
the standard cone of the identity automorphism of $\opn{T}^{-1}(A)$
in $\dcat{C}_{\mrm{str}}(A)$. Since
\[ C = \opn{T}^{-1}(A) \oplus \opn{T}^{}(\opn{T}^{-1}(A)) = 
\opn{T}^{-1}(A) \oplus A \]
as graded objects, we have the elements 
$e_0 := 1_A \in A \sub C^0$ and 
$e_1 := \opn{t}^{-1}(1_A) \in \opn{T}^{-1}(A)^1 \sub C^1$.
They satisfy $\d_C(e_0) = e_1$. 
Note that the DG module $C$ is semi-free, with semi-free filtration 
\begin{equation} \label{eqn:3332}
F_j(C) := 
\begin{cases}
0 & \tup{if} \ \ j = -1
\\
\opn{T}^{-1}(A) &  \tup{if} \ \ j = 0 
\\
C &  \tup{if} \ \ j \geq 1 .
\end{cases}
\end{equation}
As a graded $A^{\natural}$-module, $C^{\natural}$ is free, with basis 
$(e_0, e_1)$. Of course the DG module $C$ is contractible. 

\begin{lem} \label{lem:3306}
Let $M \in \dcat{C}(A)$.
\begin{enumerate}
\item There is a homomorphism $\psi : Q \to M$ 
in $\dcat{C}_{\mrm{str}}(A)$, such that 
$\opn{Z}(\psi) : \opn{Z}(Q) \to \opn{Z}(M)$ is surjective,
and $Q = \bigoplus_{s \in S} \opn{T}^{-i_s}(A)$ 
for a collection of integers $\{ i_s \}_{s \in S}$. 

\item There is a surjective homomorphism
$\psi' : Q' \to M$ in $\dcat{C}_{\mrm{str}}(A)$,
such that 
$Q' := \bigoplus_{s \in S'} \opn{T}^{-i_s}(C)$
for some collection of integers $\{ i_s \}_{s \in S'}$,
where $C$ is the DG $A$-module from  formula \tup{(\ref{eqn:3307})}.
\end{enumerate}
\end{lem}

To clarify the notation in the lemma, the collections of integers  
$\{ i_s \}_{s \in S}$ and $\{ i_s \}_{s \in S'}$ are distinct; namely the
function $i : s \mapsto i_s$ is $i : S \sqcup S' \to \Z$. 

\begin{proof} \mbox{}

\smallskip \noindent
(1) For every cocycle $m \in \opn{Z}^i(M)$
there a homomorphism $\psi_m : \opn{T}^{-i}(A) \to M$ 
that sends the element 
$\opn{t}^{-i}(1_A) \in \opn{Z}^i(\opn{T}^{-i}(A))$
to $m$. Thus, if $\{ m_s \}_{s \in S}$ is a collection of 
homogeneous cocycles that generates $\opn{Z}(M)$ as a $\K$-module, we get a 
homomorphism $\psi$ as claimed. 

\medskip \noindent
(2) For every element $m \in M^i$ there is a homomorphism 
$\psi'_m : \opn{T}^{-i}(C) \to M$ 
in $\dcat{C}_{\mrm{str}}(A)$ that sends 
$\opn{t}^{-i}(e_0) \mapsto m$ and 
$\opn{t}^{-i}(e_1) \mapsto (-1)^i \cd \d_M(m)$. 
Hence, by taking a collection $\{ m_s \}_{s \in S'}$ of homogeneous elements 
that  generates $M$ as a $\K$-module, we get a homomorphism $\psi'$ as claimed. 
\end{proof}

The DG $A$-module $Q$ in item (1) of the lemma is free. The DG $A$-module 
$Q'$ in item (2) is semi-free, with semi-free filtration induced by that of 
$C$, namely 
\begin{equation} \label{eqn:3311}
F_{j}(Q') := \bigoplus_{s \in S'} \, \opn{T}^{-i_s}(F_j(C))
\end{equation}
for $j \geq -1$, where $F_j(C)$ is from (\ref{eqn:3332}). 

\begin{proof}[Proof of Theorem \tup{\ref{thm:3305}}]
The proof is in several steps. 

\smallskip \noindent
Step 1. We are going to produce an exact sequence 
\begin{equation} \label{eqn:3350}
\cdots \to Q^{-2} \xar{\pa^{-2}} Q^{-1} \xar{\pa^{-1}}
Q^0 \xar{\eta} M \to 0
\end{equation}
in $\dcat{C}_{\mrm{str}}(A)$ with these properties: 
\begin{itemize}
\rmitem{a} The sequence 
\[ \cdots \to \opn{Z}(Q^{-2}) \xar{\opn{Z}(\pa^{-2})} \opn{Z}(Q^{-1}) 
\xar{\opn{Z}(\pa^{-1})} \opn{Z}(Q^0) \xar{\opn{Z}(\eta)} \opn{Z}(M) \to 0 \]
in $\dcat{G}_{\mrm{str}}(\K)$ is exact.

\rmitem{b} For every $i \leq 0$ the DG $A$-module $Q^i$ has the following 
decomposition 
in $\dcat{C}_{\mrm{str}}(A)$~: 
$Q^i = Q_i \oplus Q'_i$, where
$Q_i = \bigoplus_{s \in S_i} \, \opn{T}^{-j_s}(A)$
and 
$Q'_i = \lb \bigoplus_{s \in S'_i} \, \opn{T}^{-j_s}(C)$
for some collection of integers 
$\{ j_s \}_{s \in S_i \sqcup S'_i}$. 
\end{itemize}
This will be done recursively on $i \leq 0$. 

By Lemma \ref{lem:3306} we can find DG $A$-modules $Q_0, Q'_0$ 
and homomorphisms 
$\psi_0 : Q_0 \to M$, $\psi'_0 : Q'_0 \to M$, such that both 
$\opn{Z}(\psi_0)$ and $\psi'_0$ are surjective. 
The DG modules $Q_0, Q'_0$ are of the required 
form (property (b) with $i = 0$). We let 
$Q^0 := Q_0 \oplus Q'_0$ and $\eta := \psi_0 \oplus \psi'_0$. 

Let $N^0 := \opn{Ker}(\eta)$. Then the sequence 
\begin{equation} \label{eqn:3352}
0 \to N^0 \to Q^0 \xar{\eta} M \to 0 
\end{equation}
in $\dcat{C}_{\mrm{str}}(A)$ is exact. The sequence 
\[ 0 \to \opn{Z}(N^0) \to \opn{Z}(Q^0) \xar{\opn{Z}(\eta)} \opn{Z}(M) \to 0 \]
in $\dcat{G}_{\mrm{str}}(\K)$ is also exact. Indeed, the exactness at 
$\opn{Z}(N^0)$ and $\opn{Z}(Q^0)$ is because the functor 
$\opn{Z} : \dcat{C}_{\mrm{str}}(\K) \to \dcat{G}_{\mrm{str}}(\K)$
is left exact (see Proposition \ref{prop:4145}), and the sequence 
(\ref{eqn:3352}) is exact. 
The exactness at $\opn{Z}(M)$ is by the surjectivity of 
$\opn{Z}(\psi_0) : \opn{Z}(Q_0) \to \opn{Z}(M)$.

Next we repeat this procedure with $N^0$ instead of $M$, to obtain an exact 
sequence 
\begin{equation} \label{eqn:3351}
0 \to N^{-1} \to Q^{-1} \xar{\pa^{-1}} N^0 \to 0 ,
\end{equation}
with $Q^{-1}$ of the required form (property (b) with $i = -1$), and such that 
\[ 0 \to \opn{Z}(N^{-1}) \to \opn{Z}(Q^{-1}) \xar{\opn{Z}(\pa^{-1})} 
\opn{Z}(N^0) \to 0 \]
is also exact. We then splice (\ref{eqn:3352}) and (\ref{eqn:3351})
to get the sequence 
\[ 0 \to N^{-1} \to Q^{-1} \xar{\pa^{-1}} Q^0 \xar{\eta} M . \]

Continuing recursively we obtain an exact sequence (\ref{eqn:3350}) that has 
properties (a) and (b). 

\medskip \noindent 
Step 2. We view the exact sequence (\ref{eqn:3350}) as an acyclic complex with 
entries in $\dcat{C}_{\mrm{str}}(A)$, that has $Q^0$ in degree $0$. 
Define the total DG $A$-modules 
\[ P :=\opn{Tot}^{\oplus} \bigl( 
\cdots \to Q^{-2} \xar{\pa^{-2}} Q^{-1} \xar{\pa^{-1}} Q^0 \to 0 \to \cdots
\bigr) \]
and 
\[ P^{\mrm{aug}} := \opn{Tot}^{\oplus} \bigl(
\cdots \to Q^{-2} \xar{\pa^{-2}} Q^{-1} \xar{\pa^{-1}}
Q^0 \xar{\eta} M \to 0 \to \cdots \bigr) . \]

Because the sequence (\ref{eqn:3350}) and the sequence appearing in property 
(a) are exact, we can use Proposition \ref{prop:3300} (with any 
$j_1 \in \Z$). The conclusion is that the DG $A$-module $P^{\mrm{aug}}$ is 
acyclic. Then, by Corollary \ref{cor:3310}, the homomorphism 
$\rho : P \to M$ in $\dcat{C}_{\mrm{str}}(A)$
that is induced by $\eta$ is a quasi-isomorphism. 

\medskip \noindent
Step 3. It remains to produce a semi-free filtration 
$\{ F_j(P) \}_{j \geq -1}$ on the DG $A$-module $P$. Our formula is this: 
$F_{-1}(P) := 0$ of course. For $k \geq 0$ we let 
\[ F_{2 \cd k}(P) := 
\Bigr( \bigoplus\nolimits_{l = 0}^{-(k - 1)} \, \opn{T}^{-l} (Q^l) \Bigl) 
\, \oplus \,   
\opn{T}^{k} \bigl( Q_{-k} \oplus F_0(Q'_{-k}) \bigr) , \]
where $F_0(Q'_{-k})$ comes from (\ref{eqn:3311}), and 
\[ F_{2 {\cdot} k + 1}(P)  := 
\bigoplus\nolimits_{l = 0}^{-k} \, \opn{T}^{-l} (Q^l) . \]
We leave it to the reader to verify that this is a semi-free filtration. 
\end{proof}

\begin{exer} \label{exer:4752}
Finish step 3 of the proof. 
\end{exer}

\begin{cor} \label{cor:1580}
Let $A$ be a DG ring. The category $\bcat{C}(A)$ 
has enough K-pro\-jectives%
\index{Resolution! K-projective}. 
\end{cor}

\begin{proof}
Combine Theorems \ref{thm:1575} and \ref{thm:3305}. 
\end{proof}

Recall that a DG ring $A$ is nonpositive if $A^i = 0$ for all $i > 0$.
According to Proposition \ref{prop:3170}, when $A$ is nonpositive the smart 
truncation functors exist. 

\begin{cor} \label{cor:3340}
If the DG ring $A$ is nonpositive, then every $M \in \bcat{C}(A)$ admits a 
quasi-iso\-morphism $\rho : P \to M$ in $\dcat{C}_{\mrm{str}}(A)$ from a 
semi-free DG $A$-module $P$, such that 
$\opn{sup}(P) = \opn{sup}(\opn{H}(M))$. 
\end{cor}

\begin{proof}
Let $i_1 := \opn{sup}(\opn{H}(M))$. If $i_1 = \pm \infty$ then there is nothing 
to prove beyond Theorem \ref{thm:3305}; so let us assume that $i_1 \in \Z$. 
Define $M' := \opn{smt}^{\leq i_1}(M)$. Then $M' \to M$ is a quasi-isomorphism, 
and $\opn{sup}(M') = i_1$. By replacing $M$ with $M'$ we can now assume that 
$\opn{sup}(M) = i_1$.

Since $\opn{sup}(M) \leq i_1$, in step 1 of the proof of the theorem we can 
take the DG module $Q_0$ to be concentrated in degrees $\leq i_1$.
Another consequence of the fact that $\opn{sup}(M) \leq i_1$ is that 
$M^{i_1} = \opn{Z}^{i_1}(M)$. This means that the homomorphism 
$\psi_0 : Q^{i_1}_0 \to M^{i_1}$ is already surjective; so we don't need to 
cover $M^{i_1}$ by $Q'_0$, and thus we can take $Q'_0$ to be concentrated in 
degrees $\leq i_1$. In this way we can arrange to have $\opn{sup}(Q^0) \leq 
i_1$. But then we also have $\opn{sup}(N^0) \leq i_1$. Thus 
recursively we can arrange to have $\opn{sup}(Q^l) \leq i_1$ for all $l$. 
Therefore in step 2 of the proof we get $\opn{sup}(P) \leq i_1$. 

Finally, because $\opn{sup}(\opn{H}(M)) = i_1$ we must have 
$\opn{sup}(P) = i_1$.
\end{proof}

\begin{dfn} \label{dfn:3355}
Let $A$ be a DG ring. A DG $A$-module $P$ is called a {\em 
pseudo-finite free DG $A$-module} if there are numbers $i_1 \in \Z$ and 
$r_i \in \N$ such that 
$ P \cong \bigoplus_{i \leq i_1} \opn{T}^{-i}(A)^{\oplus r_i}$
in $\dcat{C}_{\mrm{str}}(A)$. 
\end{dfn}

\begin{dfn} \label{dfn:1925}
Let $A$ be a DG ring and let $P$ be a DG $A$-module. 
\begin{enumerate}
\item A {\em pseudo-finite semi-free filtration}%
\index{Filtration on a DG module! pseudo-finite semi-free}
on $P$ is a semi-free 
filtration $F = \lb \bigl\{ F_j(P) \bigr\}_{j \geq -1}$
satisfying this finiteness condition: there are numbers 
$i_1 \in \Z$ and $r_j \in \N$ such that 
$\opn{Gr}^F_j(P) \cong \opn{T}^{-i_1 + j}(A)^{\oplus r_j}$
in $\dcat{C}_{\mrm{str}}(A)$ for all $j \in \N$. 

\item We call $P$ a {\em pseudo-finite semi-free DG $A$-module}%
\index{Differential graded module! pseudo-finite semi-free}
if it admits some pseudo-finite semi-free filtration.
\end{enumerate}
\end{dfn}

Note that in item (1) of the definition above, for every $j \geq -1$ we have
\begin{equation} \label{eqn:4736}
F_j(P)^{\natural} \cong 
\bigoplus_{0 \leq k \leq j}  \opn{Gr}^F_j(P)^{\natural}
\cong \bigoplus_{0 \leq k \leq j}  
\opn{T}^{-i_1 + k}(A^{\natural})^{\oplus r_k}
\end{equation}
as graded $A^{\natural}$-modules. 

\begin{prop} \label{prop:4731}
Let $P$ be a DG $A$-module. The two conditions below are equivalent. 
\begin{enumerate}
\rmitem{i} $P$ is a pseudo-finite semi-free DG $A$-module.

\rmitem{ii} $P$ admits a semi-basis $(S, F)$ such that $S$ is bounded above, 
say by $i_1 \in \Z$\tup{;}
each of the sets $S^i$ is finite\tup{;} and 
$F_j(S) = \bigcup_{i \in [i_1 - j, i_1]}$
for every $j \geq -1$. 
\end{enumerate}
\end{prop}

\begin{exer} \label{exer:4731}
Prove Proposition \ref{prop:4731}. (Hint: the numbers $i_1$ in Definition 
\ref{dfn:1925}(1) and in condition (ii) are the same.)
\end{exer}
 
A graded $A^{\natural}$-module 
$P^{\natural}$ is called {\em pseudo-finite free} if there is an 
isomorphism of graded $A^{\natural}$-modules
$P^{\natural} \cong \bigoplus_{j \geq 0} 
\opn{T}^{-i_1 + j}(A^{\natural})^{\oplus r_j}$
for some numbers $i_1 \in \Z$ and $r_j \in \N$. 

\begin{prop} \label{prop:3175}
Let $A$ be a DG ring and let $P$ be a DG $A$-module. The following 
two conditions are equivalent\tup{:}
\begin{enumerate}
\rmitem{i} The DG $A$-module $P$ is pseudo-finite semi-free.

\rmitem{ii} The DG $A$-module $P$ admits a semi-free filtration 
$F = \{ F_j(P) \}_{j \geq -1}$ such that 
each $\opn{Gr}^F_j(P)$ is a pseudo-finite free DG $A$-module, and this 
asymptotic formula holds\tup{:} 
$\lim_{j \to \infty} \opn{sup} (\opn{Gr}^F_j(P)) = - \infty$.
\end{enumerate}
If $A$ is nonpositive, the conditions above are equivalent to\tup{:} 
\begin{enumerate}
\rmitem{iii} The graded $A^{\natural}$-module $P^{\natural}$ is pseudo-finite 
free.
\end{enumerate}
If $A$ is a ring \tup{(}i.e.\ $A^i = 0$ for all $i \neq 0$\tup{)}, the 
conditions above are equivalent to\tup{:} 
\begin{enumerate}
\rmitem{iv} $P$ is a bounded above complex of finitely generated 
free $A$-modules. 
\end{enumerate}
\end{prop}

\begin{exer} \label{exer:3175}
Prove the proposition above. (Hint: rephrase these conditions in terms of  
semi-bases.)
\end{exer}

\begin{rem} \label{rem:4635}
Condition (iv) of Proposition \ref{prop:3175} is a hint for the reason we chose 
the name ``pseudo-finite'' in Definition \ref{dfn:1925}, and also the name 
``pseudo-noetherian'' in Definitions \ref{dfn:4755} and \ref{dfn:3180} below 
(in conjunction with Theorem \ref{thm:3340}). Indeed, in \cite{SGA6} and 
\cite{SP}, a complex of modules over a ring $A$ is called {\em pseudo-coherent} 
if it is isomorphic, in $\dcat{D}(A)$, to a bounded above complex of finite 
free 
$A$-modules. 

For a generalization of the notion of pseudo-coherent complex to DG modules, 
see 
Definition \ref{dfn:4470}.
\end{rem}

\begin{dfn} \label{dfn:4755} 
A graded ring $A$ is called {\em left pseudo-noetherian} if it is
nonpositive, the ring $A^0$ is left noetherian, and each $A^i$ is a finitely 
generated (left) $A^0$-module. 
\end{dfn}

\begin{dfn} \label{dfn:3180} 
Let $A$ be a DG ring. 
\begin{enumerate}
\item We call $A$ a 
{\em cohomologically left pseudo-noetherian DG ring}%
\index{Differential graded ring! cohomologically pseudo-noetherian}
if its cohomology $\opn{H}(A)$ is a left pseudo-noetherian graded ring.

\item Suppose $A$ is cohomologically left pseudo-noetherian. We 
denote by $\dcat{D}_{\mrm{f}}(A)$ the full subcategory of 
$\dcat{D}(A)$ on the DG modules $M$ such that $\opn{H}^i(M)$ is a finitely 
generated $\opn{H}^0(A)$-module for every $i$. 
\end{enumerate}
\end{dfn}

Of course $\dcat{D}_{\mrm{f}}(A)$ is a full triangulated subcategory of 
$\dcat{D}(A)$. As usual we combine indicators: 
$\dcat{D}_{\mrm{f}}^{\star}(A) := \dcat{D}_{\mrm{f}}(A) \cap 
\dcat{D}^{\star}(A)$.

\begin{exa} \label{exa:3182}
Suppose $A^0$ is a nonzero noetherian commutative ring.
Let $X = \coprod_{i \leq 0} X^i$ be a 
nonpositive graded set, such that each graded piece $X^i$ is finite. 
Consider the strongly commutative graded polynomial ring
$A^0[X] := A^0 \ot \K [X]$
from Example \ref{exa:4625}.
The graded ring $A^0[X]$ is (left and right) pseudo-noetherian, but if $X$ is 
infinite then $A^0[X]$ is not noetherian (on either side).  

Likewise, let $Y = \coprod_{i < 0} Y^i$ be a negative graded set, such that 
each graded piece $Y^i$ is finite. Take the noncommutative graded 
polynomial ring
$A^0 \bra{Y} := A^0 \ot \K \bra{Y}$
from Example \ref{exa:4625} and Definition \ref{dfn:4306}. 
The graded ring $A^0\bra{Y}$ is (left and right) pseudo-noetherian, but if $Y$ 
has at least two elements, then $A^0\bra{Y}$ is not noetherian (on either side).
\end{exa}

\begin{rem} \label{rem:3205}
Let $A$ be a DG ring whose cohomology $\opn{H}(A)$ is a nonpositive graded 
ring. Define $A' := \opn{smt}^{\leq 0}(A)$, the smart truncation of 
$A$ below $0$. It is easy to see that $A'$ is 
a DG subring of $A$, and the inclusion $A' \to A$ is a DG ring 
quasi-isomorphism. By Theorem \ref{thm:2363} there is an equivalence of 
triangulated categories $\dcat{D}(A) \to \dcat{D}(A')$. 
This says that for many purposes we can assume that $A$ itself is a 
nonpositive DG ring. 
\end{rem}

\begin{thm} \label{thm:3340}
\index{Differential graded ring! cohomologically pseudo-noetherian}
\index{Differential graded module! pseudo-finite semi-free}
Let $A$ be a nonpositive cohomologically left pseudo-noether\-ian DG ring, and 
let $M \in \dcat{D}_{\mrm{f}}^{-}(A)$. 
Then there is a quasi-isomorphism 
$P \to M$ in $\dcat{C}_{\mrm{str}}(A)$, where $P$ is a pseudo-finite 
semi-free DG $A$-module, and $\opn{sup}(P) = \opn{sup}(\opn{H}(M))$. 
\end{thm}

\begin{proof}
The proof is divided into five parts. 

\smallskip \noindent 
Part 1. We may assume that $\opn{H}(M) \neq 0$, so that 
$i_1 := \opn{sup}(\opn{H}(M)) \in \Z$. 
We are going to construct a direct system 
\[ 0 = F_{-1}(P) \sub F_0(P) \sub F_1(P) \sub \cdots \]
of pseudo-finite semi-free DG $A$-modules, together with a direct system of DG 
$A$-module homomorphisms 
$F_j(\rho) : F_{j}(P) \to M$
for $j \geq 0$, having these properties:
\begin{enumerate}
\rmitem{a} Let $j \geq 0$. The homomorphism 
$\opn{H}^i(F_j(\rho)) : \opn{H}^i(F_{j}(P)) \to \opn{H}^i(M)$ 
is surjective for all $i$, and bijective for all $i > i_1 - j$. 

\rmitem{b} $F_0(P)$ is a pseudo-finite free DG $A$-module with 
$\opn{sup}(F_0(P)) = i_1$.

\rmitem{c} Let $j \geq 1$. There is an isomorphism 
$F_{j}(P) / F_{j - 1}(P) \cong \opn{T}^{-i_1 + j}(A)^{\oplus r_j}$
in $\dcat{C}_{\mrm{str}}(A)$, for some $r_j \in \N$.
\end{enumerate}
This will be done recursively on $j \geq 0$. . 

\medskip \noindent
Part 2. Here we deal with the cases $j = -1, 0$. 
For $j = -1$ we take $F_{-1}(P) := 0$ of course. 

Now we deal with $F_0(P)$. 
For every $i \leq i_1$ we choose finitely many generators of the 
$\opn{H}^0(A)$-module $\opn{H}^i(M)$, and then we lift them to elements of
$\opn{Z}^i(M)$. These choices give rise to a pseudo-finite free DG $A$-module 
$F_0(P)$ and a homomorphism 
$F_0(\rho) : F_0(P) \to M$ in 
$\dcat{C}_{\mrm{str}}(A)$, such that 
$\opn{sup}(F_0(P)) = i_1$, and the homomorphism 
$\opn{H}(F_0(\rho)) : \opn{H}(F_0(P)) \to \opn{H}(M)$ 
in $\dcat{G}_{\mrm{str}}(\K)$ is surjective.  
Condition (a) holds for $j = 0$, and so does condition (b). Condition (c) is 
not applicable here. 

\medskip \noindent
Part 3. 
Suppose that $j \geq 0$, and we have a direct system 
\[ 0 = F_{-1}(P) \sub F_0(P) \sub F_1(P) \sub \cdots \sub F_j(P) \]
of semi-free DG $A$-modules, together with a direct system of DG $A$-module 
homomorphisms
$F_{j'}(\rho) : F_{j'}(P) \to M$ 
for $j' \in [0, j]$, satisfying conditions (a)-(c) with index $j'$ 
instead of $j$. In this part we will construct the DG $A$-module 
$F_{j + 1}(P)$. 

Let 
\begin{equation} \label{eqn:5065}
\bar{N} := \opn{Ker} \bigl( \opn{H}^{i_1 - j}(F_j(\rho)) : 
\opn{H}^{i_1 - j}(F_{j}(P)) \to \opn{H}^{i_1 - j}(M) \bigr) .
\end{equation}
Because $F_{j}(P)$ is a pseudo-finite semi-free DG $A$-module, and $A$ is a 
cohomologically left pseudo-noetherian DG ring, it follows that $\bar{N}$ is a 
finitely generated $\opn{H}^0(A)$-module.
We choose elements 
$\bar{n}_1, \ldots, \bar{n}_{r_{j + 1}} \in \bar{N}$
that generate $\bar{N}$ as an $\opn{H}^0(A)$-module.
Then we lift them to elements
$n_1, \ldots, n_{r_{j + 1}} \in \opn{Z}^{i_1 - j}(F_{j}(P))$.

Consider the free DG $A$-module 
$Q' := \opn{T}^{-i_1 + j}(A)^{\oplus r_{j + 1}}$
with basis $(e'_1, \ldots, \lb e'_{r_{j + 1}})$ concentrated in degree 
$i_1 - j$. There is a homomorphism 
$\phi : Q' \to F_{j}(P)$ in $\dcat{C}_{\mrm{str}}(A)$
with formula $\phi(e'_k) := n_k$. Define 
\begin{equation} \label{eqn:4755}
F_{j + 1}(P) := \opn{Cone} \bigl( \phi : Q' \to F_{j}(P) \bigr) , 
\end{equation}
the standard cone of $\phi$. 
The DG $A$-module $F_{j + 1}(P)$ is pseudo-finite semi-free.
Condition (b) still holds, because $F_{0}(P)$ did not change. And 
$F_{j + 1}(P) / F_{j}(P) \cong Q$, where 
\begin{equation} \label{eqn:4760}
Q := \opn{T}(Q') \cong \opn{T}^{-i_1 + j + 1}(A)^{\oplus r_{j + 1}} , 
\end{equation}
so condition (c) holds also for $j + 1$. 

\medskip \noindent
Part 4. We continue part 3. Here we construct the 
homomorphism 
$F_{j + 1}(\rho) : F_{j + 1}(P) \to M$
in $\dcat{C}_{\mrm{str}}(A)$. 
For every $k \in [1, r_{j + 1}]$ let 
$e_k := \opn{t}_{Q'}(e'_k) \in Q^{i_1 - j - 1}$.
The sequence $(e_1, \ldots, e_{r_{j + 1}})$ is a basis of the free DG 
$A$-module $Q$. 
We have $\d_Q(e_k) = 0$ and $\opn{t}^{-1}_{Q}(e_k) = e'_k \in Q'$, 
so by the definition of the differential of the cone we have 
\[ \d_{F_{j + 1}(P)}(e_k) = \d_Q(e_k) + \phi(\opn{t}^{-1}_{Q}(e_k)) = 
\phi(e'_k) = n_k . \]
Now the cohomology class $\bar{n}_k$ of $n_k$ belongs to $\bar{N}$,
so 
$\opn{H}^{i_1 - j}(F_j(\rho))(\bar{n}_k) = 0$.
This means that the element 
$F_j(\rho)(n_k)$ is a coboundary in $M$, and we can lift it to an element 
$m_k \in M^{i_1 - j - 1}$.

By definition of the cone we have
$F_{j + 1}(P)^{\natural} = F_{j}(P)^{\natural} \oplus Q^{\natural}$ 
as graded $A^{\natural}$-modules.
We define the homomorphism 
$F_{j + 1}(\rho) : F_{j + 1}(P)^{\natural} \to M^{\natural}$ 
in $\dcat{G}_{\mrm{str}}(A^{\natural})$
to be the extension of $F_{j}(\rho)$ such that 
$F_{j + 1}(\rho)(e_k) := m_k$. Then 
$F_{j + 1}(\rho) : F_{j + 1}(P) \to M$ 
is actually a homomorphism in $\dcat{C}_{\mrm{str}}(A)$. 
         
There is equality $F_{j + 1}(P)^i = F_{j}(P)^i$
in all degrees $i \geq i_1 - j$.
Therefore  
$\opn{H}^i(F_{j + 1}(\rho)) = \opn{H}^i(F_j(\rho))$
for all $i > i_1 - j$, and these homomorphisms are bijective by assumption. 
The homomorphism 
$\opn{Z}^{i_1 - j}(F_{j}(P)) \to \opn{Z}^{i_1 - j}(F_{j + 1}(P))$
is surjective, and hence 
$\opn{H}^{i_1 - j}(F_{j}(P)) \to \opn{H}^{i_1 - j}(F_{j + 1}(P))$
is also surjective. Because 
the submodule $\bar{N} \sub \opn{H}^{i_1 - j}(F_j(P))$
goes to zero in $\opn{H}^{i_1 - j}(F_{j + 1}(P))$,
we have the following commutative diagram:
\[ \UseTips \xymatrix @C=3ex @R=6ex {
0
\ar[r]
&
\bar{N}
\ar[r]
\ar[d]
&
\opn{H}^{i_1 - j}(F_{j}(P))
\ar[rrr]^(0.53){\opn{H}^{i_1 - j}(F_{j}(\rho))}
\ar[d]^{\mrm{surj}}
& & &
\opn{H}^{i_1 - j}(M)
\ar[r]
\ar[d]^{\opn{id}}
&
0
\\
0
\ar[r]
&
0
\ar[r]
&
\opn{H}^{i_1 - j}(F_{j + 1}(P))
\ar[rrr]^(0.53){\opn{H}^{i_1 - j}(F_{j + 1}(\rho))}
& & &
\opn{H}^{i_1 - j}(M)
\ar[r]
&
0
} \]
The top row is exact, by (\ref{eqn:5065}). Therefore the bottom row is exact, 
i.e.\ \lb $\opn{H}^{i_1 - j}(F_{j + 1}(\rho))$ is bijective.
We conclude that condition (a) holds with index $j + 1$. 

\medskip \noindent
Part 5. Define 
$P := \lim_{j \to} F_{j}(P)$ and 
$\rho := \lim_{j \to} F_{j}(\rho)$
in $\dcat{C}_{\mrm{str}}(A)$. 
Condition (c) implies that for every $i$ the direct system $\{ F_j(P^i) \}_{j 
\geq -1}$ is eventually stationary. Hence, by condition (a), the homomorphism 
$\opn{H}^i(\rho)$ is bijective. So $\rho$ is a quasi-isomorphism. 

By condition (b) and (c) the filtration $\{ F_j(P) \}_{j \geq -1}$
on $P$ satisfies condition (ii) of Proposition \ref{prop:3175}. Therefore $P$ 
is a pseudo-finite semi-free DG $A$-module.
\end{proof}

\begin{exa} \label{exa:1930} 
A special yet very important case of Theorem \ref{thm:3340} is this: 
$A$ is a left noetherian ring, and $M$ is a complex of $A$-modules 
with bounded above cohomology, such that each $\opn{H}^i(M)$ is a finitely 
generated $A$-module. Then $M$ has a resolution 
$P \to M$, where $P$ is a complex of finitely generated free $A$-modules, and 
$\opn{sup}(P) = \opn{sup}(\opn{H}(M))$. 
Compare this to Example \ref{exa:2900}.
\end{exa}

\begin{rem} \label{rem:4675}
In Example \ref{exa:4725} we discussed the existence of K-flat resolutions of 
unbounded complexes of sheaves of modules. Here is a more general assertion. 

Consider a {\em DG ringed space} $(X, \AA)$, namely $\AA$ is a 
sheaf of DG rings on a topological space $X$. There are no commutativity or 
boundedness conditions on $\AA$. By combining the concept of a 
semi-free DG module from Definition \ref{dfn:1575} with the concept of a 
pseudo-free module from Example \ref{exa:4725}, we get something new: a {\em 
semi-pseudo-free DG $\AA$-module}. This is a K-flat DG $\AA$-module. The 
arguments in \cite{Ye15} imply that every DG $\AA$-module $\MM$ admits a 
quasi-isomorphism $\PP \to \MM$ from a semi-pseudo-free DG $\AA$-module $\PP$.
Therefore $\dcat{C}(\AA)$ has enough K-flat objects.  
A detailed proof of this assertion will be published in the future. 
\end{rem}

\mysubsection{K-Injective Resolutions in 
\texorpdfstring{$\bcat{C}^+(\cat{M})$}{C+(M)}} 
\label{subsec:exis-K-inj}

In this subsection $\cat{M}$ is an abelian category, and $\dcat{C}(\cat{M})$ 
is the category of complexes in $\cat{M}$. Cofiltrations and their limits were 
introduced in Definition \ref{dfn:3336}.

\begin{dfn} \label{dfn:1620}
Let $I$ be a complex in $\dcat{C}(\cat{M})$. 
\begin{enumerate}
\item A {\em semi-injective cofiltration}%
\index{Cofiltration on a DG module! semi-injective}
on $I$ is a cofiltration $G = \{ G_q(I) \}_{q \geq -1}$
in $\dcat{C}_{\mrm{str}}(\cat{M})$ such that:
\begin{itemize}
\item $G_{-1}(I) = 0$. 

\item  Each $\opn{Gr}_q^G(I)$ 
is a complex of injective objects of $\cat{M}$ with zero differential.

\item $I = \lim_{\lar q} \, G_q(I)$.
\end{itemize}

\item The complex $I$ is called a {\em semi-injective complex}%
\index{Complex in abelian category! semi-injective}
if it admits some semi-injective cofiltration. 
\end{enumerate}
\end{dfn}

\begin{thm} \label{thm:1620}
\index{Complex in abelian category! semi-injective}
\index{Complex in abelian category! K-injective}
Let $\cat{M}$ be an abelian category, and let $I$ be a 
semi-injective complex in $\dcat{C}(\cat{M})$. Then $I$ is 
K-injective. 
\end{thm}

\begin{proof}
The proof is very similar to that of Theorem \ref{thm:1536}. 
 
\medskip \noindent
Step 1. We start by proving that if $I = \opn{T}^p(J)$, the translation of an 
injective object $J \in \cat{M}$, then $I$ is K-injective. This is easy: 
given an acyclic complex $N \in \dcat{C}(\cat{M})$, we have 
\[ \opn{Hom}_{\cat{M}}(N, I) = 
\opn{Hom}_{\cat{M}} \bigl( N, \opn{T}^p(J) \bigr) \cong
\opn{T}^{p} \bigl( \opn{Hom}_{\cat{M}}(N, J) \bigr) \]
in $\dcat{C}_{\mrm{str}}(\K)$. But 
$\opn{Hom}_{\cat{M}}(-, J)$ is an exact functor 
$\cat{M} \to \dcat{M}(\K)$, so $\opn{Hom}_{\cat{M}}(N, J)$ is an acyclic 
complex. 

\medskip \noindent
Step 2. Now $I$ is a complex of injective objects of $\cat{M}$ with zero 
differential. This means that 
$I \cong \prod_{p \in \Z} \opn{T}^p(J_p)$
in $\dcat{C}_{\mrm{str}}(\cat{M})$, where each $J_p$ is an injective object in 
$\cat{M}$. But then 
\[ \opn{Hom}_{\cat{M}}(N, I) \cong 
\prod_{p \in \Z} \, \opn{Hom}_{\cat{M}} \bigl( N, \opn{T}^p(J_p) \bigr) . \]
This is an easy case of Proposition \ref{prop:1535}(2). By step 1 and the fact 
that a product of acyclic complexes in $\dcat{C}_{\mrm{str}}(\K)$ is acyclic 
(itself an easy case of the Mittag-Leffler argument), we conclude that 
$\opn{Hom}_{\cat{M}}(N, I)$ is acyclic. 

\medskip \noindent
Step 3. Fix a semi-injective cofiltration $G = \{ G_q(I) \}_{q \geq -1}$
of $I$. Here we prove that for every $q$ the complex  
$G_q(I)$ is K-injective. This is done by induction on $q$. 
For $q = -1$ it is trivial. For $q \geq 0$ there is an exact sequence of 
complexes 
\begin{equation} \label{eqn:1620}
0 \to \opn{Gr}_q^{G}(I) \to G_q(I) \to G_{q - 1}(I) \to 0 
\end{equation}
in $\dcat{C}_{\mrm{str}}(\cat{M})$. 
In each degree $p \in \Z$ the exact sequence 
\[ 0 \to \opn{Gr}_q^{G}(I)^p \to G_q(I)^p \to G_{q - 1}(I)^p \to 0 \]
in $\cat{M}$ splits, because $\opn{Gr}^G_{q}(I)^p$ is an injective 
object. Thus the exact sequence (\ref{eqn:1620}) is split in the category 
$\bcat{G}_{\mrm{str}}(\cat{M})$ of graded objects in $\cat{M}$. 

Let $N \in \dcat{C}(\cat{M})$ be an acyclic complex. Applying the 
functor $\opn{Hom}_{\cat{M}}(N, -)$ 
to the sequence of complexes (\ref{eqn:1620}) we obtain a sequence 
\begin{equation} \label{eqn:1621}
0 \to \opn{Hom}_{\cat{M}} \bigl( N, \opn{Gr}^G_{q}(I) \bigr) \to  
\opn{Hom}_{\cat{M}} \bigl( N, G_q(I) \bigr) \to 
\opn{Hom}_{\cat{M}} \bigl( N, G_{q - 1}(I) \bigr) \to  0
\end{equation}
in $\dcat{C}_{\mrm{str}}(\K)$. Because (\ref{eqn:1620}) is split exact in 
$\bcat{G}_{\mrm{str}}(\cat{M})$, the sequence (\ref{eqn:1621}) is split exact 
in $\bcat{G}_{\mrm{str}}(\K)$. Therefore (\ref{eqn:1621}) is exact in 
$\dcat{C}_{\mrm{str}}(\K)$.

By the induction hypothesis the complex 
$\opn{Hom}_{\cat{M}} \bigl( N, G_{q - 1}(I) \bigr)$
is acyclic. By step 2 the complex
$\opn{Hom}_{\cat{M}} \bigl( N, \opn{Gr}^G_{q}(I) \bigr)$
is acyclic. The long exact cohomology sequence associated to 
(\ref{eqn:1621}) shows that the complex 
$\opn{Hom}_{\cat{M}} \bigl( N, G_{q}(I) \bigr)$
is acyclic too.

\medskip \noindent
Step 4. We keep the semi-injective cofiltration 
$G = \{ G_q(I) \}_{q \geq -1}$ from step 3. 
Take any acyclic complex $N \in \bcat{C}(\cat{M})$. By 
Proposition \ref{prop:1535} we know that 
\[ \opn{Hom}_{\cat{M}}(N, I) \cong 
\lim_{\leftarrow q} \, \opn{Hom}_{\cat{M}} 
\bigl( N, G_q(I) \bigr) \]
in $\dcat{C}_{\mrm{str}}(\K)$. According to step 3 the complexes 
$\opn{Hom}_{\cat{M}} \bigl( N, G_q(I) \bigr)$ are all acyclic. 
The exactness of the sequences (\ref{eqn:1621}) implies that 
the inverse system \lb 
$\bigl\{ \opn{Hom}_{\cat{M}} \bigl( N, G_q(I) \bigr) \bigr\}_{q \geq -1}$
in $\dcat{C}_{\mrm{str}}(\K)$ has surjective transitions. Now the 
Mittag-Leffler argument (Corollary \ref{cor:1590}) says that the inverse limit 
complex $\opn{Hom}_{\cat{M}}(N, I)$ is acyclic. 
\end{proof}

\begin{prop} \label{prop:1635}
Let $\cat{M}$ be an abelian category. If $I$ is a bounded below complex of 
injectives, then $I$ is a semi-injective complex.
\end{prop}

\begin{proof}
We can assume that $I \neq 0$. Let 
$p_0 := \opn{inf}(I) \in \Z$. For $q \geq -1$ let 
$G_q(I) := \opn{stt}^{\leq p_0 + q}(I)$,
the stupid truncation from Definition \ref{dfn:3066}.
The cofiltration $G = \{ G_q(I) \}_{q \geq -1}$ is semi-injective. 
\end{proof}

The next theorem is \cite[Lemma I.4.6(1)]{RD}. See also 
\cite[Proposition 1.7.7(i)]{KaSc1}. 

\begin{thm} \label{thm:1630}
Let $\cat{M}$ be an abelian category, and let $\cat{J} \sub \cat{M}$ be a
full subcategory such that every object $M \in \cat{M}$ admits a 
monomorphism $M \inj I$ to some object $I \in \cat{J}$. Then every complex 
$M \in \bcat{C}^{+}(\cat{M})$ 
admits a quasi-isomorphism $\rho : M \to  I$ in 
$\dcat{C}^+_{\mrm{str}}(\cat{M})$,
such that $\opn{inf}(I) = \opn{inf}(M)$,
and each $I^p$ is an object of $\cat{J}$.
\end{thm}

\begin{proof}
The proof is the same as that of Theorem \ref{thm:1540}, except for a 
mechanical reversal of arrows. To be more explicit, let us
take $\cat{N} := \cat{M}^{\mrm{op}}$ and $\cat{Q} := \cat{J}$.
Since monomorphisms in $\cat{M}$ become epimorphisms in $\cat{N}$, the 
full subcategory $\cat{Q} \sub \cat{N}$ satisfies the 
assumptions of Theorem \ref{thm:1540}. By Theorem \ref{thm:2495} we have a 
canonical isomorphism of categories 
$\bcat{C}^{-}_{\mrm{str}}(\cat{N}) \iso 
\dcat{C}^+_{\mrm{str}}(\cat{M})^{\mrm{op}}$. 
Thus a quasi-isomorphism $Q \to N$ in 
$\bcat{C}^{-}_{\mrm{str}}(\cat{N})$
gives rise to a quasi-isomorphism $M \to I$ in 
$\bcat{C}^{+}_{\mrm{str}}(\cat{M})$.
\end{proof}

\begin{dfn} \label{dfn:2880}
Let $\cat{M}$ be an abelian category, and let $\cat{M}' \sub \cat{M}$ be a 
full abelian subcategory. 
We say that {\em $\cat{M}'$ has enough injectives relative to $\cat{M}$} 
if every object $M \in \cat{M}'$ admits a monomorphism 
$M \inj I$, where $I$ is an object of $\cat{M}'$ that is injective in 
the bigger category $\cat{M}$.
\end{dfn}

Of course, in this situation the category $\cat{M}'$ itself has enough 
injectives. 

Thick abelian categories were defined in Definition \ref{dfn:2323}. 
The next theorem is \cite[Lemma I.4.6(3)]{RD}. See also 
\cite[Proposition 1.7.11]{KaSc1}. 

\begin{thm} \label{thm:2880}
Let $\cat{M}$ be an abelian category, and let $\cat{M}' \sub \cat{M}$ be a 
thick abelian subcategory that has enough injectives relative to $\cat{M}$.
Let $M \in \dcat{C}(\cat{M})$ be a complex with bounded below cohomology, 
such that $\opn{H}^i(M) \in \cat{M}'$ for all $i$.
Then there is a quasi-isomorphism $\rho : M \to I$ in 
$\dcat{C}_{\mrm{str}}(\cat{M})$, such that
$I \in \dcat{C}^+(\cat{M}')$, each $I^p$ is an injective object in 
$\cat{M}$, and $\opn{inf}(I) = \opn{inf}(\opn{H}(M))$. 
\end{thm}

Before the proof we need some auxiliary material. 

Suppose we are given morphisms $\psi_1 : K \to L_1$ and 
$\psi_2 : K \to L_2$ in $\cat{M}$. The {\em fibered coproduct} is the object 
\begin{equation} \label{eqn:2885}
L_1 \oplus_K L_2 := \opn{Coker} 
\bigl( (\psi_1, -\psi_2) : K \to L_1 \oplus L_2 \bigr) 
\end{equation}
in $\cat{M}$. 
It has an obvious universal property. The commutative diagram 
\begin{equation} \label{eqn:5072}
\UseTips \xymatrix @C=8ex @R=6ex {
K
\ar[r]^{\psi_1}
\ar[d]_{\psi_2}
&
L_1
\ar[d]^{\ep_1}
\\
L_2
\ar[r]^(0.37){\ep_2}
& 
L_1 \oplus_K L_2
} 
\end{equation}
in which $\ep_i$ are the morphisms induced by the embeddings 
$L_i \to L_1 \oplus L_2$, is sometimes called a {\em pushout diagram}. 

\begin{lem} \label{lem:2880}
In the situation above\tup{:}
\begin{enumerate}
\item The sequence
\[ \opn{Ker}(\psi_1) \xar{\psi_2} L_2 \xar{\ep_2} 
L_1 \oplus_K L_2 \xar{\pi_1} \opn{Coker}(\psi_1) \to 0 , \]
in which $\pi_1$ is the morphism induced by the projection 
$L_1 \oplus L_2 \to L_1$, is exact. 

\item Suppose $L_1 \to L'_1$ be a monomorphism. Then the induced morphism 
$L_1 \oplus_K L_2 \to L'_1 \oplus_K L_2$ 
is a monomorphism.
\end{enumerate}
\end{lem}

\begin{exer} \label{exer:2880}
Prove Lemma \ref{lem:2880}. (Hint: use the first sheaf trick,
Proposition \ref{prop:3635}.)  
\end{exer}

\begin{proof}[Proof of Theorem \tup{\ref{thm:2880}}]
The proof is an adaptation of the proof of \cite[Proposition 1.7.11]{KaSc1}. 
In the proof we use the objects of cocycles $\opn{Z}^p(L)$, 
co\-boundaries $\opn{B}^p(L)$ and decocycles $\opn{Y}^p(L)$, that are 
associated to a complex $L$ and an integer $p$; see Definition \ref{dfn:2993}.

\medskip \noindent
Step 1. We can assume that $\opn{H}(M) \neq 0$. 
By translating $M$, we may further assume that $\opn{inf}(\opn{H}(M)) = 0$. 
Then, after replacing $M$ with its smart truncation $\opn{smt}^{\geq 0}(M)$, 
we can assume that $\opn{inf}(M) = 0$. 

We are going to construct an inverse system 
\begin{equation} \label{eqn:5105}
F_p(I) = \bigl( \cdots \to 0 \to I^0 \xar{\d_I^0} I^1 
\xar{\d_I^{1}} \cdots \xar{\d_I^{p - 1}} I^p \to 0 \to \cdots \bigr)
\end{equation}
in $\dcat{C}_{\mrm{str}}(\cat{M}')$, indexed by $p \geq -1$, such that
the objects $I^p \in \cat{M}'$ are injective in $\cat{M}$. 
The morphism $F_{p + 1}(I) \to F_{p}(I)$ will send the object 
$I^{p + 1} = F_{p + 1}(I)^{p + 1}$ to 
$0$, and will be the identity in all other degrees. 
Simultaneously we will 
construct an inverse system of morphisms 
\begin{equation} \label{eqn:5106}
F_p(\phi) : M \to F_p(I)
\end{equation}
in $\dcat{C}_{\mrm{str}}(\cat{M})$.
The construction will be inductive. 
The morphisms $F_p(\phi)$ will satisfy this condition:
\begin{itemize}
\rmitem{$\heartsuit$} The morphism 
$\opn{H}^q(F_p(\phi)) : \opn{H}^q(M) \to \opn{H}^q(F_p(I))$
is a monomorphism for all $q \leq p$ and an isomorphism for all
$q \leq p - 1$. 
\end{itemize}

Then the complex 
$I := \lim_{\lar p} F_{p}(I)$ in $\dcat{C}(\cat{M}')$
is going to be of the required kind, and the morphism 
$\phi := \lim_{\lar p} F_{p}(\phi)$
in $\dcat{C}_{\mrm{str}}(\cat{M})$ is going to be a quasi-isomorphism.
Note that these inverse limits are innocent, since in each degree $q$ the 
inverse systems $\{ F_{p}(I)^q \}_{p \geq -1}$ and 
$\{ F_{p}(\phi)^q \}_{p \geq -1}$ are eventually stationary. 

\medskip \noindent
Step 2. In this step we begin the construction. 
Here $p \geq -1$, and we have a complex $F_p(I)$ as in (\ref{eqn:5105})
of the kind described above,
and a morphism $F_p(\phi) : M \to F_p(I)$ in 
$\dcat{C}_{\mrm{str}}(\cat{M})$ satisfying condition ($\heartsuit$). 
If  $p = -1$ then $F_{-1}(I) = 0$ and $F_{-1}(\phi) = 0$ of course.

We claim that the objects 
$F_p(I)^q$, $\opn{H}^q(F_p(I))$, 
$\opn{B}^{q}(F_p(I))$, $\opn{Y}^{q}(F_p(I))$ and 
$\opn{Z}^{q}(F_p(I))$ belong to $\cat{M}'$ for all $q \in [0, p + 1]$. 
For $q = p + 1$ all these objects are zero. 
For $F_p(I)^q = I^q$ it is given. For 
$\opn{H}^p(F_p(I)) = \opn{Coker}(\d_I^{p - 1})$ it is because $\cat{M}'$ is 
a full abelian subcategory of $\cat{M}$. For $\opn{H}^q(F_p(I))$ when 
$q < p$ we use the isomorphisms $\opn{H}^q(F_p(\phi))$ and the fact that 
$\opn{H}^q(M) \in \cat{M}'$. As for the rest of the objects listed, this is 
shown by descending induction on $q$, starting from $q = p$, using these short 
exact sequences 
\[ 0 \to \opn{Z}^{q}(F_p(I)) \to F_p(I)^{q} \to 
\opn{B}^{q + 1}(F_p(I)) \to  0 , \]
\[ 0 \to \opn{B}^{q}(F_p(I)) \to F_p(I)^q \to \opn{Y}^{q}(F_p(I)) \to  0 , \]
\[ 0 \to \opn{B}^{q}(F_p(I)) \to \opn{Z}^{q}(F_p(I)) \to 
\opn{H}^{q}(F_p(I)) \to  0 , \]
and the fact that $\cat{M}'$ is a full abelian subcategory of $\cat{M}$. 

\medskip \noindent
Step 3. The construction continues. Let us denote the components of the 
morphism $F_p(\phi)$ by
$\phi^q := F_p(\phi)^q : M^q \to F_p(I)^q = I^q$
for $q \in [0, p]$; outside this degree range we have $F_p(\phi)^q = 0$. 

The differential $\d_M^{p} : M^{p} \to M^{p + 1}$
induces a morphism $\psi_1 : \opn{Y}^p(M) \to \opn{Z}^{p + 1}(M)$.
There is also the morphism 
$\psi_2 := \opn{Y}^p(\phi) : \opn{Y}^p(M) \to \opn{Y}^{p}(F_p(I))$.
We define $N$ to be the fibered coproduct
\begin{equation} \label{eqn:5107}
N := \opn{Z}^{p + 1}(M) \oplus_{\opn{Y}^p(M)} \opn{Y}^{p}(F_p(I))
\end{equation}
relative to the morphisms $\psi_1$ and $\psi_2$. 
In the notation of Lemma \ref{lem:2880}, with the objects  $K := \opn{Y}^p(M)$, 
$L_1 := \opn{Z}^{p + 1}(M)$ and $L_2 := \opn{Y}^{p}(F_p(I))$,
and with the morphism $\psi_1, \psi_2$ that we have just defined, there are 
morphisms 
$\ep_1 : \opn{Z}^{p + 1}(M) \to N$,
$\ep_2 : \opn{Y}^{p}(F_p(I)) \to N$ and
$\pi_1 : N \to \opn{Coker}(\psi_1) = \opn{H}^{p + 1}(M)$.
According to Lemma \ref{lem:2880}(1) there is an exact sequence 
\begin{equation} \label{eqn:5109}
\opn{H}^{p}(M) \xar{\psi_2} \opn{Y}^{p}(F_p(I)) \xar{\ep_2} 
N \xar{\pi_1} \opn{H}^{p + 1}(M) \to 0 
\end{equation}
in $\cat{M}$. 
However, letting 
$\be : \opn{H}^{p}(F_p(I)) \to \opn{Y}^{p}(F_p(I))$ 
be the canonical monomorphism, there is equality 
$\psi_2 = \be \circ \opn{H}^p(F_p(\phi))$. By condition ($\heartsuit$) with 
indices $q = p$ we know that $\opn{H}^p(F_p(\phi))$ is a monomorphism. 
Therefore $\psi_2$ is a monomorphism, and (\ref{eqn:5109}) extends to a 
slightly longer exact sequence 
\begin{equation} \label{eqn:5115}
0 \to \opn{H}^{p}(M) \xar{\psi_2} \opn{Y}^{p}(F_p(I)) \xar{\ep_2} 
N \xar{\pi_1} \opn{H}^{p + 1}(M) \to 0 
\end{equation}
in $\cat{M}$. 

Let $\ga : \opn{Z}^{p + 1}(M) \inj M^{p + 1}$ be the canonical monomorphism,
and define 
$\psi_1^+ := \ga \circ \psi_1 : \opn{Y}^{p}(M) \to M^{p + 1}$. 
Let 
$\psi_2^+ := \psi_2 :  \opn{Y}^p(M) \to \opn{Y}^{p}(F_p(I))$.
We have the corresponding fibered coproduct
\begin{equation} \label{eqn:5108}
N^+ := M^{p + 1} \oplus_{\opn{Y}^p(M)} \opn{Y}^{p}(F_p(I)) . 
\end{equation}
The monomorphism $\ga : \opn{Z}^{p + 1}(M) \inj M^{p + 1}$ induces a 
morphism $\om : N \to N^+$ between these fibered coproducts, and Lemma 
\ref{lem:2880}(2) says that $\om$ is also a monomorphism. 

\medskip \noindent
Step 4. We now produce the object $I^{p + 1} \in \cat{M}'$. 
In step 2 we showed that $\opn{Y}^{p}(F_p(I)) \in \cat{M}'$. It is given that 
$\opn{H}^{q}(M) \in \cat{M}'$ for all $q$. 
The exact sequence (\ref{eqn:5109}) and Proposition \ref{prop:3171} 
say that $N \in \cat{M}'$ too. 
By assumption there is a monomorphism
$\chi : N \inj I^{p + 1}$
into some object  $I^{p + 1} \in \cat{M}'$ that's injective in $\cat{M}$. 
Because $\om : N \to N^+$ is a monomorphism, there is a morphism 
$\chi^+ : N^+ \to I^{p + 1}$
such that $\chi^+ \circ \om = \chi$. 

We have a canonical epimorphism $\al_I : I^{p} \surj \opn{Y}^{p}(F_p(I))$. 
Define the morphism 
$\d_{I}^{p} : I^{p} \to I^{p + 1}$
by the formula 
$\d_{I}^{p} := \chi \circ \ep_2 \circ \al_I$. Since 
$\al_I \cd \d_{I}^{p - 1} = 0$, we see that 
$\d_{I}^{p} \circ \d_{I}^{p - 1} = 0$. 
This is depicted in the commutative diagram below. 
\[ \UseTips \xymatrix @C=6ex @R=6ex {
I^{p - 1}
\ar[r]_(0.6){\d_I^{p - 1}}
\ar@(ur,ul)[rr]^{0}
&
I^p
\ar[r]^(0.4){\al_I}
\ar@(dr,dl)[rrr]_{\d_I^p}
&
\opn{Y}^{p}(F_p(I))
\ar[r]^(0.6){\ep_2}
&
N
\ar[r]^{\chi}
&
I^{p + 1}
} \]
Therefore we obtain a new complex 
\begin{equation} \label{eqn:5117}
F_{p + 1}(I) := \bigl( \cdots \to 0 \to I^0 \xar{\d_I^0} I^1 
\xar{\d_I^{1}} \cdots \xar{\d_I^{p - 1}} I^p 
\xar{\d_I^{p}} I^{p + 1} \to 0 \to \cdots \bigr) 
\end{equation}
in $\dcat{C}(\cat{M}')$, and it projects onto $F_{p}(I)$. 

\medskip \noindent
Step 5. In this step we construct the morphism 
$\phi^{p + 1} : M^{p + 1} \to I^{p + 1}$. 
Let $\ep^+_1 : M^{p + 1} \to N^+$ be the morphism that's part of the fibered 
coproduct (\ref{eqn:5108}). 
We define the morphism 
$\phi^{p + 1} : M^{p + 1} \to I^{p + 1}$
by the formula 
$\phi^{p + 1} := \chi^+ \circ \ep_1^+$. 

Let's examine the next diagram, where
$\al_M : M^p \surj \opn{Y}^{p}(M)$ is the canonical epimorphism, and 
$\psi_1^+$, $\psi_2^+ = \psi_2$, $\ep_1^+$, $\ep_2^+$ and $\ep_2$ are the 
morphisms involved in  the two fibered coproducts. 
\begin{equation} \label{eqn:5110}
\UseTips \xymatrix @C=6ex @R=6ex {
M^p
\ar[r]^(0.4){\al_M}
\ar[d]^{\phi^p}
&
\opn{Y}^{p}(M)
\ar[r]^(0.55){\psi_1^+}
\ar[d]^{\psi_2^+}
&
M^{p + 1}
\ar[r]^{\ep_1^+}
&
N^+
\ar[r]^{\chi^+}
&
I^{p + 1}
\ar[d]^{\opn{id}}
\\
I^p
\ar[r]^(0.4){\al_I}
&
\opn{Y}^{p}(F_p(I))
\ar[rr]^{\ep_2}
\ar[urr]^{\ep_2^+}
&
&
N
\ar[r]^{\chi}
\ar[u]_{\om}
&
I^{p + 1}
} 
\end{equation}
The left square is commutative because $\al : \opn{Id} \to \opn{Y}$ is a 
morphism of functors. The triangle that has $\opn{Y}^{p}(M)$ as a vertex is 
commutative because $N^+$ is the fibered coproduct. The other triangle, the 
one that has $N$ as a vertex, is commutative because $\om$ is a morphism 
between fibered coproducts. The right square is commutative by the defining 
property of $\chi^+$. So the whole diagram is commutative. 

We now take part of the commutative diagram (\ref{eqn:5110}), and add to it 
three curved arrows. 
\begin{equation} \label{eqn:5111}
\UseTips \xymatrix @C=6ex @R=5.5ex {
M^p
\ar[r]^(0.4){\al_M}
\ar[d]_{\phi^p}
\ar@(ur,ul)[rr]^{\d_{M}^{p}}
&
\opn{Y}^{p}(M)
\ar[r]^(0.55){\psi_1^+}
&
M^{p + 1}
\ar[r]^{\ep_1^+}
\ar@(ur,ul)[rr]^{\phi^{p + 1}}
&
N^+
\ar[r]^{\chi^+}
&
I^{p + 1}
\ar[d]^{\opn{id}}
\\
I^p
\ar[r]^(0.4){\al_I}
\ar@(dr,dl)[rrrr]_{\d_{I}^{p}}
&
\opn{Y}^{p}(F_p(I))
\ar[rr]^{\ep_2}
&
&
N
\ar[r]^{\chi}
&
I^{p + 1}
} 
\end{equation}
The curved polygons are commutative, as can be seen from the definitions of the 
morphisms $\psi_1^+$, $\phi^{p + 1}$ and $\d_{I}^{p}$. We conclude that 
$\d_I^{p} \circ \phi^p = \phi^{p + 1} \circ \d_M^{p}$,
so there is a new morphism of complexes 
$F_{p + 1}(\phi) : M \to F_{p + 1}(I)$
whose degree $p + 1$ component is $\phi^{p + 1}$. 

\medskip \noindent
Step 6. Here we prove that $\opn{H}^{p + 1}(F_{p + 1}(\phi))$ is a 
monomorphism and $\opn{H}^{p}(F_{p + 1}(\phi))$ is an isomorphism.
This will verify condition ($\heartsuit$) with index $p + 1$.

Consider the fibered coproduct 
\begin{equation} \label{eqn:5116}
N^{\lozenge} := I^{p + 1} \oplus_{\opn{Y}^p(F_p(I))} \opn{Y}^{p}(F_p(I))
\end{equation} 
of the morphisms $\psi_1^{\lozenge}$ and $\psi_2^{\lozenge}$,
where $\psi_1^{\lozenge} : \opn{Y}^{p}(F_p(I)) \to I^{p + 1}$
is the morphism induced by $\d_I^p$, and 
$\psi_2^{\lozenge}$ is the identity of $\opn{Y}^{p}(F_p(I))$. 
Note that $\ep_1^{\lozenge} : I^{p + 1} \to N^{\lozenge}$ is an isomorphism, 
and 
$\opn{Coker}(\psi_1^{\lozenge}) = \opn{Y}^{p + 1}(F_{p + 1}(I)) = 
\opn{H}^{p + 1}(F_{p + 1}(I))$. 
The morphisms 
$\phi^{p + 1} : \opn{Z}^{p + 1}(M) \to I^{p + 1}$
and 
$\opn{Y}^{p}(F_{p}(\phi)) : \opn{Y}^{p}(M) \to \opn{Y}^{p}(F_{p}(I))$
induce a morphism 
$\om^{\lozenge} : N \to N^{\lozenge}$ between the fibered coproducts 
(\ref{eqn:5107}) 
and (\ref{eqn:5116}). According to Lemma \ref{lem:2880}(1) there is an 
induced commutative diagram
\begin{equation} \label{eqn:51168}
\UseTips \xymatrix @C=4ex @R=6ex {
0
\ar[r]
&
\opn{H}^{p}(M)
\ar[r]^(0.47){\psi_2}
\ar[d]_{\opn{H}^{p}(F_{p + 1}(\phi))}
&
\opn{Y}^{p}(F_{p}(I))
\ar[r]^(0.65){\ep_2}
\ar[d]_{\opn{id}}
&
N
\ar[r]^(0.4){\pi_2}
\ar[d]_{\om^{\lozenge}}
&
\opn{H}^{p + 1}(M)
\ar[r]
\ar[d]^{\opn{H}^{p + 1}(F_{p + 1}(\phi))}
&
0
\\
0
\ar[r]
&
\opn{H}^{p}(F_{p + 1}(I))
\ar[r]^{\psi_2^{\lozenge}}
&
\opn{Y}^{p}(F_{p}(I))
\ar[r]^(0.65){\ep_2^{\lozenge}}
&
N^{\lozenge}
\ar[r]^(0.3){\pi_2^{\lozenge}}
&
\opn{H}^{p + 1}(F_{p + 1}(I))
\ar[r]
&
0
} 
\end{equation}
The first row is an exact sequence, since it is (\ref{eqn:5115}).
The second row is also exact: exactness at the last three objects is by 
Lemma \ref{lem:2880}(1). As for exactness at $\opn{H}^{p}(F_{p + 1}(I))$,
we use the fact that the morphism 
$\opn{Y}^{p}(F_{p + 1}(I)) \to \opn{Y}^{p}(F_{p}(I))$
is an isomorphism, and its composition with  the canonical embedding 
$\opn{H}^{p}(F_{p + 1}(I)) \inj \opn{Y}^{p}(F_{p + 1}(I))$
is $\psi_2^{\lozenge}$, so the latter is a monomorphism.
A little calculation shows that $\om^{\lozenge} = \ep_1^{\lozenge} \circ \chi$, 
so it is a monomorphism. 

Finally, a diagram chase in (\ref{eqn:51168}), using the first sheaf trick, 
shows that $\opn{H}^{p + 1}(F_{p + 1}(\phi))$
is a monomorphism, and that $\opn{H}^{p}(F_{p + 1}(\phi))$
is an isomorphism. 
\end{proof}

\begin{cor} \label{cor:3185}
Under the assumptions of Theorem \tup{\ref{thm:2880}}, the canonical functor 
$\dcat{D}^+(\cat{M}') \to \dcat{D}^+_{\cat{M}'}(\cat{M})$
is an equivalence. 
\end{cor}

\begin{proof}
This is very similar to the proof of Corollary \ref{cor:3175}, and there is no 
need to repeat it. 
\end{proof}

Here is an important instance in which Theorem \ref{thm:2880} and Corollary 
\ref{cor:3185} apply. 

\begin{exa} \label{exa:2895}
Let $(X, \OO_X)$ be a noetherian scheme. Associated to it are these abelian 
categories: the category $\cat{M} := \cat{Mod} \OO_X$ of $\OO_X$-modules,
and the thick abelian subcategory $\cat{M}' := \cat{QCoh} \OO_X$
of quasi-coherent $\OO_X$-modules. 
According to \cite[Proposition II.7.6]{RD} the category $\cat{M}'$ has 
enough injectives relative to $\cat{M}$.
By Corollary \ref{cor:3185} the canonical functor
\[ \dcat{D}(\cat{QCoh} \OO_X) \to \dcat{D}^+_{\mrm{qc}}(\cat{Mod} \OO_X) \]
is an equivalence. 

For a more general statement, using other methods, see 
\cite[Corollary 5.5]{BoNe}. 
\end{exa}

\begin{cor} \label{cor:1630}
\index{Resolution! K-injective}
If $\cat{M}$ is an abelian category with enough injectives, then  
$\bcat{C}^{+}(\cat{M})$ has enough K-injectives. 
\end{cor}

\begin{proof}
According to either Theorem \ref{thm:1630} or Theorem \ref{thm:2880},
every $M \in \bcat{C}^{+}(\cat{M})$ admits a quasi-isomorphism $M \to I$
to bounded below complex of injectives $I$. Now use Proposition 
\ref{prop:1635} and Theorem \ref{thm:1620}. 
\end{proof}

\begin{cor} \label{cor:1902}
Let $\cat{M}$ be an abelian category with enough injectives, and let 
$M \in \dcat{C}(\cat{M})$ be a complex with bounded below cohomology. Then $M$ 
has a K-injective resolution $M \to I$,
such that $\opn{inf}(I) = \opn{inf}(\opn{H}(M))$, 
and every $I^p$ is an injective object of $\cat{M}$. 
\end{cor}

\begin{proof}
We may assume that $\opn{H}(M)$ is nonzero. Let 
$p := \opn{inf}(\opn{H}(M)) \in \Z$, and let 
$N := \opn{smt}^{\geq p}(M)$, the smart truncation from 
Definition \ref{dfn:2320}. So $M \to N$ is a quasi-isomorphism, and 
$\opn{inf}(N) = p$. 
According to either Theorem \ref{thm:1630} or Theorem \ref{thm:2880},
there is a quasi-isomorphism $N \to I$, where $I$ is a complex 
of injectives and $\opn{inf}(I) = p$. 
By Proposition \ref{prop:1635} and Theorem \ref{thm:1620} the complex $I$ is 
K-injective. The composed quasi-isomorphism $M \to I$ is what we are looking 
for. 
\end{proof}

\mysubsection{K-Injective Resolutions in 
\texorpdfstring{$\bcat{C}(A)$}{C(A)}} 
\label{subsec:exis-K-inj-dgmods}

Recall that we are working over a non\-zero commutative base ring $\K$, and 
$A$ is a central DG $\K$-ring. 

An {\em injective cogenerator} of the abelian category
$\dcat{M}(\K) = \cat{Mod} \K$ is an injective 
$\K$-module $\K^*$ with this property: if $M$ is a nonzero $\K$-module, then 
$\opn{Hom}_{\K}(M, \K^*)$ is nonzero. These always exist. Here are a few 
examples. 

\begin{exa} \label{exa:1665} 
For every $\K$ there is a canonical choice for an injective cogenerator: 
$\K^* := \opn{Hom}_{\Z}(\K, \Q / \Z)$.
See proof of Theorem \ref{thm:1115}.
Usually this is a very big module! 
\end{exa}

\begin{exa} \label{exa:1876}
Assume $\K$ is a complete noetherian local ring, with maximal ideal 
$\m$ and residue field $\mbb{k} = \K / \m$. In this case we would prefer to take
the smallest possible injective cogenerator $\K^*$, and this is the injective 
hull of $\mbb{k}$ as a $\K$-module. 

Here are some special cases. 
If $\K$ is a field, then $\K^* = \K = \mbb{k}$. 
If $\K = \what{\Z}_p$, the ring of $p$-adic integers, then 
$\mbb{k} = \mbb{F}_p$, and $\K^* \cong \what{\Q}_p / \what{\Z}_p$, which is the 
$p$-primary part of $\Q / \Z$. 
If $\K$ contains some field, then there exists a ring homomorphism 
$\mbb{k} \to \K$ that lifts the canonical surjection $\K \to \mbb{k}$.
(This is by the Cohen Structure Theorem, see \cite[Theorem 7.7]{Eis}.)
After choosing such a lifting, there is an isomorphism of $\K$-modules 
$\K^* \cong \opn{Hom}^{\mrm{cont}}_{\mbb{k}}(\K, \mbb{k})$,
where continuity is for the $\m$-adic topology on $\K$ and the discrete 
topology on $\mbb{k}$. 
\end{exa}

In this subsection we have the following standing convention.  

\begin{conv} \label{conv:5075}
There is a fixed injective cogenerator $\K^*$ of $\dcat{M}(\K)$.
\end{conv}

We view $\K^*$ as a DG $\K$-module concentrated in degree $0$ 
(with zero differential). 
For every $p \in \Z$ there is the DG $\K$-module $\opn{T}^{-p} (\K^*)$, 
which is concentrated in degree $p$.

It will be convenient to blur the distinction between DG $\K$-modules with 
zero differentials and graded $\K$-modules. Namely, if $N$ is a DG 
$\K$-module such that $\d_N = 0$, we will identify $N$ with the graded modules 
$N^{\natural}$, $\opn{H}(N)$, etc.

\begin{dfn} \label{dfn:5075}
For a DG $\K$-module $V$ we define the DG $\K$-module 
$V^* := \opn{Hom}_{\K}(V, \K^*)$. 
\end{dfn}

The operation $(-)^*$ is an exact contravariant functor from 
$\dcat{C}_{\mrm{str}}(\K)$ to itself, and also from 
$\dcat{G}_{\mrm{str}}(\K)$ to itself. Note that given 
$M \in \dcat{C}(A)$, its dual $M^*$ is an object of 
$\dcat{C}(A^{\mrm{op}})$, so that we have an exact functor 
\begin{equation} \label{eqn:5075}
(-)^* : \dcat{C}_{\mrm{str}}(A)^{\mrm{op}} \to 
\dcat{C}_{\mrm{str}}(A^{\mrm{op}}) . 
\end{equation}

\begin{dfn} \label{dfn:1805}
A DG $\K$-module $W$ is called {\em cofree} if there is an isomorphism
$W \cong V^*$ in $\dcat{C}_{\mrm{str}}(\K)$, for some free DG $\K$-module 
$V$. 
\end{dfn}

Since the differential of a free DG $\K$-module $V$ is zero, we see that 
a cofree DG module $W$ also has zero differential. 

\begin{lem} \label{lem:5075}
Let $V$ be a free DG $\K$-module, and let $W := V^*$. 
\begin{enumerate}
\item If $V \cong \boplus_{s \in S} \opn{T}^{p_s} (\K)$
for some indexing set $S$ and some collection of integers 
$\{ p_s \}_{s \in S}$, then 
$W \cong \prod_{s \in S} \opn{T}^{-p_s} (\K^*)$. 

\item As an object of the abelian category $\dcat{G}_{\mrm{str}}(\K)$, $W$ is 
injective. 
\end{enumerate}
\end{lem}

In item (1) the isomorphisms can be considered either in 
$\dcat{C}_{\mrm{str}}(\K)$ or in $\dcat{G}_{\mrm{str}}(\K)$; it doesn't matter, 
because $V$ and $W$ have zero differentials. However, item (2) is true only in 
the category $\dcat{G}_{\mrm{str}}(\K)$.

\begin{lem} \label{lem:3346}
Let $\phi : U \to V$ be a homomorphism in $\dcat{G}_{\mrm{str}}(\K)$. 
\begin{enumerate}
\item $\phi$ is injective iff $\phi^* : V^* \to U^*$ is surjective. 

\item $\phi$ is surjective iff $\phi^* : V^* \to U^*$ is injective.

\item The canonical homomorphism $U \to U^{**} = (U^*)^*$ in 
$\dcat{G}_{\mrm{str}}(\K)$ is injective. 

\item There exists an injective homomorphism 
$U \inj W$ in $\dcat{G}_{\mrm{str}}(\K)$
into some cofree DG $\K$-module $W$. 
\end{enumerate}
\end{lem}

\begin{exer}
Prove Lemmas \ref{lem:5075} and \ref{lem:3346}. 
\end{exer}

\begin{dfn} \label{dfn:3345}
Let $W$ be a cofree DG $\K$-module. The {\em cofree DG $A$-module coinduced 
from $W$} is the DG $A$-module
$I_W :=  \opn{Hom}_{\K}(A, W)$.
There is a homomorphism $\th_W : I_W \to W$
in $\dcat{C}_{\mrm{str}}(\K)$, with formula 
$\th_W(\chi) := \chi(1_A) \in W$. 
\end{dfn}

\begin{dfn} \label{dfn:3346}
A DG $A$-module $J$ is called {\em cofree}%
\index{Differential graded module! cofree}
if there is an isomorphism 
$J \cong I_W$ in $\dcat{C}_{\mrm{str}}(A)$ for some cofree DG $\K$-module 
$W$. 
\end{dfn}

\begin{lem} \label{lem:3345}
Consider the free DG $\K$-module 
$V := \boplus_{s \in S} \opn{T}^{p_s} (\K)$
and the cofree DG $\K$-module 
$W := V^*$. There are canonical isomorphisms 
\[ I_W \cong (A \ot V)^* \cong
\Bigl( \bigoplus\nolimits_{s \in S} \, \opn{T}^{p_s} (A) \Bigr)^*
\cong \prod\nolimits_{s \in S} \, \opn{T}^{-p_s} (A^*) \]
in $\dcat{C}_{\mrm{str}}(A)$.
\end{lem}

\begin{exer}
Prove the lemma above. 
\end{exer}

\begin{lem} \label{lem:3347}
Let $W$ be a cofree DG $\K$-module, and let $M$ be a DG $A$-module.
The homomorphism 
\[ \opn{Hom}(\opn{id}_M, \th_W) : \opn{Hom}_{A}(M, I_W) \to 
\opn{Hom}_{\K}(M, W) \]
in $\dcat{C}_{\mrm{str}}(\K)$ is an isomorphism. 
\end{lem}

\begin{proof}
The homomorphism
\[ \opn{Hom}(\opn{id}_M, \th_W) :
\opn{Hom}_{A} \bigl( M,  \opn{Hom}_{\K}(A, W) \bigr) \iso 
\opn{Hom}_{\K}(M, W) \]
is just adjunction for the DG ring homomorphism $\K \to A$, so it is bijective. 
\end{proof}

Recall that $\dcat{G}_{\mrm{str}}(A^{\natural})$ is the abelian category whose 
objects are the graded $A^{\natural}$-modules, and the morphisms are the 
$A$-linear homomorphisms of degree $0$. The forgetful functor 
$\opn{Und} : \dcat{C}_{\mrm{str}}(A) \to \dcat{G}_{\mrm{str}}(A^{\natural})$, 
$M \mapsto M^{\natural}$, is faithful. 

\begin{lem} \label{lem:1770}
Let $I$ be a cofree DG $A$-module. Then $I^{\natural}$ is an injective object 
of $\dcat{G}_{\mrm{str}}(A^{\natural})$.
\end{lem}

\begin{proof}
We can assume that $I = I_W$ for some cofree DG $\K$-module $W$. 
For every  $M \in \dcat{G}_{\mrm{str}}(A^{\natural})$ there are isomorphisms 
\[ \begin{aligned}
& \opn{Hom}_{\dcat{G}_{\mrm{str}}(A^{\natural})} \bigl( M, I_W^{\natural} 
\bigr) 
= \opn{Hom}_{A}(M, I_W)^0  
\\
& \qquad \cong^{\heartsuit} \opn{Hom}_{\K}(M, W)^0 =
\prod_{p \in \Z} \, \opn{Hom}_{\K}(M^{p}, W^p) 
\end{aligned} \]
in $\dcat{M}(\K)$.
The isomorphism $\cong^{\heartsuit}$ is by Lemma \ref{lem:3347}. 
For every $p$  the functor
$\dcat{G}_{\mrm{str}}(A^{\natural}) \to \dcat{M}(\K)$,  $M \mapsto M^p$, is 
exact. Because each $W^p$ is an injective object of $\dcat{M}(\K)$, 
the contravariant functor $\opn{Hom}_{\K}(-, W^p)$
from $\dcat{M}(\K)$ to itself is exact. 
And the product of exact functors into $\dcat{M}(\K)$ is exact. We conclude 
that 
the functor 
$\opn{Hom}_{\dcat{G}_{\mrm{str}}(A^{\natural})} \bigl( -, I_W^{\natural} \bigr)$
is exact. 
\end{proof}

Recall that for a DG $\K$-module $M$ we have the object of decocycles
\[ \opn{Y}(M) := \opn{Coker} \bigl( \d_M : \opn{T}^{-1}(M) \to M \bigr) =
M / \opn{B}(M) \in \dcat{G}(\K) . \]

\begin{lem} \label{lem:3348}
Let $W$ be a cofree DG $\K$-module, let $M$ be a DG $A$-module, and let 
$\chi : \opn{Y}(M) \to W$ be a homomorphism in $\dcat{G}_{\mrm{str}}(\K)$.
Then there is a unique homomorphism 
$\psi : M \to I_W$ in $\dcat{C}_{\mrm{str}}(A)$, such that the diagram 
\[ \UseTips \xymatrix @C=10ex @R=6ex {
\opn{Y}(M)
\ar[r]^{\opn{Y}(\psi)}
\ar@(u,u)[rr]^{\chi}
&
\opn{Y}(I_W)
\ar[r]^(0.45){\opn{Y}(\th_W)}
&
\opn{Y}(W) = W
} \]
in $\dcat{G}_{\mrm{str}}(\K)$ is commutative. 
\end{lem}

\begin{proof}
Since the differentials of $W$ and $\opn{Y}(M)$ are zero, the canonical 
homomorphism 
\[ \al : \opn{Hom}_{\K}(\opn{Y}(M), W)^0 = 
\opn{Z}^0 \bigl( \opn{Hom}_{\K}(\opn{Y}(M), W) \bigr) \to 
\opn{Z}^0 \bigl( \opn{Hom}_{\K}(M, W) \bigr) , \]
that's induced by the canonical surjection 
$M \surj \opn{Y}(M)$, is bijective. This gives us a unique homomorphism 
$\al(\chi) : M \to W$ in $\dcat{C}_{\mrm{str}}(\K)$. Next we use Lemma 
\ref{lem:3347} to obtain a unique 
$\psi : M \to I_W$ in $\dcat{C}_{\mrm{str}}(A)$ s.t.\ 
$\th_W \circ \psi = \al(\chi)$. This $\psi$ is what we are looking for. 
\end{proof}

The next definition is dual to Definition \ref{dfn:1575}. 

\begin{dfn} \label{dfn:1665}
Let $I$ be an object of $\dcat{C}(A)$.
\begin{enumerate}
\item A {\em semi-cofree cofiltration}%
\index{Cofiltration on a DG module! semi-cofree}
on $I$ is a cofiltration
$G = \{ G_q(I) \}_{q \geq -1}$ on $I$ in $\dcat{C}_{\mrm{str}}(A)$ such that:
\begin{itemize}
\item $G_{-1}(I) = 0$. 

\item Each $\opn{Gr}_q^G(I)$ is a cofree DG $A$-module.

\item $I = \lim_{\lar q} \, G_q(I)$.
\end{itemize}

\item The DG $A$-module $I$ is called a {\em semi-cofree}%
\index{Differential graded module! semi-cofree}
if it admits a semi-cofree cofiltration. 
\end{enumerate}
\end{dfn}

\begin{prop} \label{prop:3375}
If $I$ is a semi-cofree DG $A$-module, then $I^{\natural}$ is 
an injective object in the abelian category 
$\dcat{G}_{\mrm{str}}(A^{\natural})$. 
\end{prop}

\begin{prop} \label{prop:3185}
Assume $A$ is a ring. If $I$ is a semi-cofree DG $A$-module, then each $I^p$ is 
an injective $A$-module. 
\end{prop}

\begin{exer} \label{exer:3185}
Prove the two propositions above. 
\end{exer}

\begin{thm} \label{thm:1665}
\index{Differential graded module! semi-cofree}
\index{Differential graded module! K-injective}
Let $I$ be an object of $\dcat{C}(A)$. If $I$ is semi-cofree, then it is 
K-injective.
\end{thm}

\begin{proof}
The proof is very similar to those of Theorems \ref{thm:1536} and 
\ref{thm:1575}. But because the arguments involve limits, we shall give the 
full proof. 

\medskip \noindent
Step 1. Suppose $I$ is cofree; say 
$I \cong \prod_{s \in S} \, \opn{T}^{-p_s}(A^*)$. 
Lemma \ref{lem:3347} implies that for every DG $A$-module $N$ there is an 
isomorphism 
\[ \opn{Hom}_A(N, I) \cong 
\prod_{s \in S} \ \opn{Hom}_{\K} \bigl( \opn{T}^{p_s}(N), \K^* \bigr)  \]
of graded $\K$-modules.  It follows that if $N$ is acyclic, then so is 
$\opn{Hom}_A(N, I)$. 

\medskip \noindent
Step 2. Fix a semi-cofree cofiltration
$G = \{ G_q(I) \}_{q \geq -1}$ of $I$.
Here we prove that for every $q \geq -1$ the DG module  
$G_q(I)$ is K-injective. This is done by induction on $q \geq -1$. 
For $q= -1$ it is trivial. For $q \geq 0$ there is an exact sequence 
\begin{equation} \label{eqn:1766}
0 \to \opn{Gr}^G_q(I) \to G_q(I) \to G_{q - 1}(I) \to 0 
\end{equation}
in the category $\dcat{C}_{\mrm{str}}(A)$. 
Because $\opn{Gr}^G_q(I)$ is a cofree DG $A$-module, it is an injective object 
in the abelian category $\bcat{G}_{\mrm{str}}(A^{\natural})$; see Lemma 
\ref{lem:1770}. Therefore the sequence (\ref{eqn:1766}) is split exact in 
$\bcat{G}_{\mrm{str}}(A^{\natural})$.

Let $N \in \dcat{C}(A)$ be an acyclic DG module. Applying the 
functor $\opn{Hom}_{A}(N, -)$ 
to the sequence (\ref{eqn:1766}) we obtain a sequence 
\begin{equation} \label{eqn:1767}
0 \to \opn{Hom}_{A} \bigl( N, \opn{Gr}^G_q(I) \bigr) \to 
\opn{Hom}_{A} \bigl( N, G_q(I) \bigr) \to
\opn{Hom}_{A} \bigl( N, G_{q - 1}(I) \bigr) \to  0 
\end{equation}
in $\dcat{C}_{\mrm{str}}(\K)$. If we forget differentials this is a sequence in 
$\bcat{G}_{\mrm{str}}(\K)$. Because (\ref{eqn:1766}) is split exact in 
$\bcat{G}_{\mrm{str}}(A^{\natural})$, it follows that  (\ref{eqn:1767}) is 
split 
exact in $\bcat{G}_{\mrm{str}}(\K)$. Therefore (\ref{eqn:1767}) is exact in 
$\dcat{C}_{\mrm{str}}(\K)$.

By the induction hypothesis the DG $\K$-module  
$\opn{Hom}_{A} \bigl( N, G_{q - 1}(I) \bigr)$
is \lb acyclic. By step 1 the DG $\K$-module
$\opn{Hom}_{A} \bigl( N, \opn{Gr}^G_{q}(I) \bigr)$
is acyclic. The long exact cohomology sequence associated to 
(\ref{eqn:1767}) shows that the DG $\K$-module 
$\opn{Hom}_{A} \bigl( N, G_q(I) \bigr)$
is acyclic too.

\medskip \noindent
Step 3. We keep the semi-cofree cofiltration 
$G = \{ G_q(I) \}_{q \geq -1}$ from step 2. 
Take any acyclic $N \in \bcat{C}(A)$. By 
Proposition \ref{prop:1535} we know that 
\[ \opn{Hom}_{A}(N, I) \cong 
\lim_{\leftarrow j} \, \opn{Hom}_{A} \bigl( N, G_q(I) \bigr) \]
in $\dcat{C}_{\mrm{str}}(\K)$. According to step 2 the complexes 
$\opn{Hom}_{A} \bigl( N, G_q(I) \bigr)$ are all acyclic. 
The exactness of the sequences (\ref{eqn:1767}) implies that 
the inverse system \lb 
$\bigl\{ \opn{Hom}_{A} \bigl( N, G_q(I) \bigr) \bigr\}_{q \geq -1}$ 
in $\dcat{C}_{\mrm{str}}(\K)$ has surjective transitions. Now the 
Mittag-Leffler 
argument (Corollary \ref{cor:1590}) says that the inverse limit complex
$\opn{Hom}_{A}(N, I)$ is acyclic. 
\end{proof}

\begin{thm} \label{thm:3325}
\index{Differential graded module! semi-cofree}
Let $A$ be a DG ring. Every DG $A$-module $M$ admits a quasi-isomorphism
$\rho : M \to I$ in $\dcat{C}_{\mrm{str}}(A)$ to a semi-cofree DG $A$-module 
$I$.
\end{thm}

The proof of this theorem was communicated to us by B. Keller. We need a lemma 
first.

The semi-free DG $A$-module $C$ was defined in formula (\ref{eqn:3307}). 
Its dual DG $A$-module $C^* = \opn{Hom}_{\K}(C, \K^*)$  
is a semi-cofree DG $A$-module, with semi-cofree cofiltration 
$G_q(C^*) := F_q(C)^*$,
where $F_q(C)$ is from (\ref{eqn:3332}). 

\begin{lem} \label{lem:3325}
Let $M \in \dcat{C}(A)$.
\begin{enumerate}
\item There is a homomorphism $\psi : M \to J$ 
in $\dcat{C}_{\mrm{str}}(A)$, such that 
$\opn{Z}(\psi) : \opn{Z}(M) \to \opn{Z}(J)$ is injective,
and $J = \prod_{s \in S} \opn{T}^{-p_s}(A^*)$
for a collection of integers $\{ p_s \}_{s \in S}$. 

\item There is an injective homomorphism
$\psi' : M \to J'$ in $\dcat{C}_{\mrm{str}}(A)$,
such that 
$J' = \prod_{s \in S'} \, \opn{T}^{-p_s}(C^*)$
for some collection of integers $\{ p_s \}_{s \in S'}$.
\end{enumerate}
\end{lem}

To clarify the notation in the lemma: we have two distinct collections 
of integers, $\{ p_s \}_{s \in S}$ and 
$\{ p_s \}_{s \in S'}$. 

\begin{proof} \mbox{}

\smallskip \noindent
(1) By Lemma \ref{lem:3346} there is an injective homomorphism 
$\chi : \opn{Y}(M) \inj W$ 
for some cofree graded $\K$-module 
$W = \prod_{s \in S} \, \opn{T}^{-p_s}(\K^*)$.
Let $J := I_W$. 
According to Lemma \ref{lem:3348} there is a homomorphism 
$\psi : M \to J$ in $\dcat{C}_{\mrm{str}}(A)$ s.t.\ 
$\opn{Y}(\th_W) \circ \opn{Y}(\psi) = \chi$. Because $\chi$ is injective, so is
$\opn{Y}(\psi)$. 

\medskip \noindent
(2) Consider the DG $A^{\mrm{op}}$-module $M^*$. By Lemma \ref{lem:3306}
there is a collection of integers $\{ p_s \}_{s \in S'}$,
and a surjective homomorphism
$\phi : Q \to M^*$
in $\dcat{C}_{\mrm{str}}(A^{\mrm{op}})$,
where 
$Q := \bigoplus\nolimits_{s \in S'} \, \opn{T}^{p_s}(C)$. 
Dualizing we get a DG $A$-module 
$J' := Q^* \cong  \prod_{s \in S'} \opn{T}^{-p_s}(C^*)$, 
and an injective homomorphism
$\phi^* : M^{**} \to Q^* = J'$
in $\dcat{C}_{\mrm{str}}(A)$.
Composing $\phi^*$ with the canonical embedding 
$M \inj M^{**}$ gives us $\psi' : M \inj J'$. 
\end{proof}

The DG $A$-module $J$ in item (1) of the lemma is cofree. The DG $A$-module 
$J'$ in item (2) is semi-cofree, with semi-cofree cofiltration induced by that 
of $C^*$, namely 
\begin{equation} \label{eqn:3326}
G_{q}(J') :=  \prod_{s \in S'} \, \opn{T}^{-p_s}(G_q(C^*)) .
\end{equation}

\begin{proof}[Proof of Theorem \tup{\ref{thm:3325}}]
\mbox{}

\smallskip \noindent
Step 1. We are going to produce an exact sequence 
\begin{equation} \label{eqn:3338}
0 \to M \xar{\eta} J^0 \xar{\pa^0} J^1 \xar{\pa^1} J^2 \to \cdots  
\end{equation}
in $\dcat{C}_{\mrm{str}}(A)$ with these properties: 
\begin{itemize}
\rmitem{a} The sequence 
\[ 0 \to \opn{Y}(M) \xar{\opn{Y}(\eta)} \opn{Y}(J^0) \xar{\opn{Y}(\pa^0)} 
\opn{Y}(J^1) \xar{\opn{Y}(\pa^1)} \opn{Y}(J^2) \to \cdots \]
is exact in $\dcat{G}_{\mrm{str}}(\K)$.

\rmitem{b} For every $p \geq 0$ the DG $A$-module $J^p$ has the following 
decomposition 
in $\dcat{C}_{\mrm{str}}(A)$~: 
$J^p = J_p \oplus J'_p$, where
$J_p = \prod_{s \in S_p} \opn{T}^{-q_s}(A^*)$
and 
$J'_p = \lb \prod_{s \in S'_p} \opn{T}^{-q_s}(C^*)$ 
for some collection of integers 
$\{ q_s \}_{s \in S_p \sqcup S'_p}$. 
\end{itemize}
This will be done recursively on $p$. 

By Lemma \ref{lem:3325} we can find DG $A$-modules $J_0, J'_0$ 
and homomorphisms 
$\psi_0 : M \to J_0$, $\psi'_0 : M \to J'_0$, such that both $\opn{Y}(\psi_0)$
and $\psi'_0$ are injective. The DG modules $J_0, J'_0$ are of the required 
form (property (b) with $p = 0$). We let 
$J^0 := J_0 \oplus J'_0$ and $\eta := \psi_0 \oplus \psi'_0$. 

Let $N^0 := \opn{Coker}(\eta)$. Then the sequence 
\begin{equation} \label{eqn:3329}
0 \to M \xar{\eta} J^0 \to N^0 \to 0 
\end{equation}
in $\dcat{C}_{\mrm{str}}(A)$ is exact. The sequence 
\[ 0 \to \opn{Y}(M) \xar{\opn{Y}(\eta)} \opn{Y}(J^0) \to \opn{Y}(N^0) \to  0 \]
in $\dcat{G}_{\mrm{str}}(\K)$ is also exact. Indeed, the exactness at 
$\opn{Y}(J^0)$ and $\opn{Y}(N^0)$ is because the functor 
$\opn{Y} : \dcat{C}_{\mrm{str}}(\K) \to \dcat{G}_{\mrm{str}}(\K)$
is right exact (see Proposition \ref{prop:4145}), and the sequence 
(\ref{eqn:3329}) is exact. 
The exactness at $\opn{Y}(M)$ is by the injectivity of 
$\opn{Y}(\psi_0) : \opn{Y}(M) \to \opn{Y}(J_0)$.

Next we repeat this procedure with $N^0$ instead of $M$, to obtain an exact 
sequence 
\begin{equation} \label{eqn:3330}
0 \to N^0 \xar{\pa^0} J^1 \to N^1 \to 0 ,
\end{equation}
with $J^1$ of the required form (property (b) with $p = 1$), and such that 
\[ 0 \to \opn{Y}(N^0) \xar{\opn{Y}(\pa^0)} \opn{Y}(J^1) \to \opn{Y}(N^1) \to  
0 \]
is also exact. We then splice (\ref{eqn:3329}) and (\ref{eqn:3330})
to get the exact sequence 
\[ 0 \to M \xar{\eta} J^0 \xar{\pa^0} J^1 \to N^1 \to 0 . \]

Continuing recursively we obtain an exact sequence (\ref{eqn:3338}) that has 
properties (a) and (b). 

\medskip \noindent 
Step 2. We view the exact sequence (\ref{eqn:3338}) as an acyclic complex with 
entries in $\dcat{C}_{\mrm{str}}(A)$, that has $J^0$ in degree $0$. 
Define the total DG $A$-modules 
\[ I :=\opn{Tot}^{\Pi} \bigl( \cdots \to 0 \to
J^0 \xar{\pa^0} J^1 \xar{\pa^1} J^2 \to \cdots  \bigr) \]
and 
\[ I^{\mrm{aug}} := \opn{Tot}^{\Pi} \bigl( \cdots \to 0 \to M \xar{\eta} 
J^0 \xar{\pa^0} J^1 \xar{\pa^1} J^2 \to \cdots  \bigr) . \]

Because the sequence (\ref{eqn:3338}) and the sequence appearing in property 
(a) are exact, we can use Proposition \ref{prop:3321} (with any 
$j_0 \in \Z$). The conclusion is that the DG $A$-module $I^{\mrm{aug}}$ is 
acyclic. Then, by Corollary \ref{cor:3320}, the homomorphism 
$\rho : M \to I$ in $\dcat{C}_{\mrm{str}}(A)$ is a quasi-isomorphism. 

\medskip \noindent
Step 3. It remains to produce a semi-cofree cofiltration 
$\{ G_q(I) \}_{q \geq -1}$ on the DG $A$-module $I$. Our formula is this: 
$G_{-1}(I) := 0$ of course. For $r \geq 0$ we let 
\[ G_{2 \cd r}(I) := 
\Bigr( \bigoplus\nolimits_{p = 0}^{r - 1} \, \opn{T}^{-p} (J^p) \Bigl) 
\, \oplus \, \opn{T}^{-r} \bigl( J_r \oplus G_0(J'_r) \bigr) , \]
where $G_0(J'_r)$ comes from (\ref{eqn:3326}), and 
\[ G_{2 {\cdot} r + 1}(I)  := 
\bigoplus\nolimits_{p = 0}^{r} \, \opn{T}^{-p} (J^p) . \]
We leave it to the reader to verify that this is a semi-cofree cofiltration. 
\end{proof}

\begin{cor} \label{cor:1665}
Let $A$ be a DG ring. The category $\bcat{C}(A)$ has 
enough K-in\-jectives%
\index{Resolution! K-injective}.
\end{cor}

\begin{proof}
Combine Theorems \ref{thm:1665} and \ref{thm:3325}. 
\end{proof}

Recall that a DG ring $A$ is nonpositive if $A^i = 0$ for all $i > 0$. 
When $A$ is nonpositive, the DG module $A^*$ is concentrated in 
degrees $\geq 0$. 

\begin{cor} \label{cor:3345}
If the DG ring $A$ is nonpositive, then every $M \in \bcat{C}(A)$ admits a 
quasi-iso\-morphism $\rho : M \to I$ in $\dcat{C}_{\mrm{str}}(A)$ to a 
semi-cofree DG $A$-module $I$, such that 
$\opn{inf}(I) = \opn{inf}(\opn{H}(M))$. 
\end{cor}

\begin{proof}
Let $p_0 := \opn{inf}(\opn{H}(M))$. If $p_0 = \pm \infty$ then there is 
nothing to prove beyond Theorem \ref{thm:3325}; so let us assume 
that $p_0 \in \Z$. 
Define $M' := \opn{smt}^{\geq p_0}(M)$. Then $M \to M'$ is a quasi-isomorphism, 
and $\opn{inf}(M') = p_0$. By replacing $M$ with $M'$ we can now assume that 
$\opn{inf}(M) = p_0$.

Since $\opn{inf}(M) = p_0$ it follows that 
$M^{p_0} = \opn{Y}^{p_0}(M)$. This means that in step 1 of the proof of the 
theorem, the homomorphism 
$\psi_0 : M^{p_0} \to J^{p_0}_0$ is already injective. So, like in the proof of 
Corollary \ref{cor:3340}, we can arrange to have 
$\opn{inf}(J^0) \geq p_0$. Then, recursively, we can arrange to 
have $\opn{inf}(J^q) \geq p_0$ for all $q$. Therefore in step 2 of the proof we 
get $\opn{inf}(I) \geq p_0$. 

Finally, because $\opn{inf}(\opn{H}(M)) = p_0$ we must have 
$\opn{inf}(I) = p_0$.
\end{proof}

We end with the K-injective analogue of Remark \ref{rem:4675}.

\begin{rem} \label{rem:4676}
Suppose $(X, \AA)$ is a ringed space. We know that the abelian category 
$\dcat{M}(\AA) = \cat{Mod} \AA$ of sheaves of $\AA$-modules on $X$ has 
enough injective objects. Theorem \ref{thm:1630} tells us that the subcategory 
$\dcat{C}^{+}(\AA)$ has enough K-injective objects. But in fact more is true: 
the category of unbounded complexes $\dcat{C}^{}(\AA)$ also has enough 
K-injective objects; see \cite{Spa}. (There is another, more abstract proof of 
this statement: because $\dcat{M}(\AA)$ is a Grothendieck abelian category, it 
follows that $\dcat{C}^{}(\AA)$ has enough K-injectives. See 
\cite[Corollary 14.1.8]{KaSc2}.)

Now let's consider a DG ringed space $(X, \AA)$. Fix an injective cogenerator 
$\K^*$ of $\dcat{M}(\K)$. 
For a point $x \in X$, consider the cofree DG $\AA_x$-module 
$I_x := \opn{Hom}_{\K}(\AA_x, \K^*)$.
Pushing this module forward along the inclusion $\{ x \} \sub X$ we get a 
constant sheaf $\II_x$ with support the closed set 
$\ol{\{ x \}}$, that we refer to as the cofree DG $\AA$-module 
associated to $x$. By definition, a {\em cofree DG $\AA$-module} is a product 
of sheaves $\II_x$ (for varying points $x \in X$).  
Combining this notion with that of a semi-cofree DG module, we can talk about 
{\em semi-cofree DG  $\AA$-modules}. Like Theorem \ref{thm:1665}, one can show 
that a semi-cofree DG  $\AA$-module $\II$ is K-injective. 

Assume that the topological space $X$ is finite dimensional, in the sense that 
there is a natural number $d$ such that for every open set $U \sub X$, every 
sheaf of $\K$-module $\NN$ on $U$, and every $p > d$, the sheaf cohomology
$\opn{H}^p(U, \NN)$ vanishes. 
This condition is quite mild -- it is satisfied when $X$ is a finite 
dimensional noetherian scheme, and also when $X$ is a subspace of a finite 
dimensional topological manifold. 
Then a geometric variant of Theorem \ref{thm:3325} is 
true: every DG  $\AA$-module $\MM$ admits a quasi-isomorphism 
$\MM \to \II$ to a semi-cofree DG  $\AA$-module $\II$.
The finite dimensionality of $X$ is needed to invoke the ML argument. 
A detailed proof of this assertion will be published in the future. 
\end{rem}

\cleardoublepage
\mysection{Adjunctions, Equivalences and Cohomological Dimension} 
\label{sec:adj-equ-cohdim}

\AYcopyright

In this section we discuss the derived Hom and tensor bifunctors in several 
situations, and show how these bifunctors are related in adjunction 
formulas. Cohomological dimensions of functors and DG modules are introduced, 
and they are used to prove several theorems. 
Given a homomorphism between two DG rings, we study the derived restriction, 
induction and coinduction functors between the derived categories of these DG 
rings. We prove that for a quasi-isomorphism of DG rings the restriction 
functor is an equivalence. We also prove the existence of NC semi-free DG ring 
resolutions. 

As before, we work over a nonzero commutative base ring $\K$, and Convention 
\ref{conv:2490} is in place.

\mysubsection{Boundedness Conditions Revisited} 
\label{subsec:bounded-revis}

This subsection is about notation. We start by introducing some combinatorial 
definitions that will enable a precise treatment of cohomological dimensions of 
objects and functors (in Subsection \ref{subsec:way-out}). 

We also change our notation for boundedness conditions in derived 
categories (Definition \ref{dfn:2122}) and for the translation functor 
(Definition \ref{dfn:2120}). These changes will simplify and improve our 
discussion (while hopefully not creating confusion). See Remark \ref{rem:4765} 
for an explanation. 

{\em Integer intervals} were introduced in Definition \ref{dfn:4765}. 
{\em Boundedness conditions} were introduced in Subsections 
\ref{subsec:DGModinM} and  \ref{subsec:bd-in-homot}.
We need to refine these notions. 

\begin{dfn} \label{dfn:4780}
\index{Boundedness! condition}
Let $S$ be an integer interval.
\begin{enumerate}
\item We call $S$ a {\em bounded above integer interval} if either 
$S = [i_0, i_1]$ with $i_1 \in \Z$, or if $S = \varnothing$. 

\item We call $S$ a {\em bounded below integer interval} if either 
$S = [i_0, i_1]$ with $i_0 \in \Z$, or if $S = \varnothing$.

\item We call $S$ a {\em bounded integer interval} if it is both bounded above 
and bounded below; i.e.\ either $S = [i_0, i_1]$ with $i_0, i_1 \in \Z$, or 
$S = \varnothing$.
\end{enumerate}
\end{dfn}

The next definition was already mentioned in passing in 
Subsection \ref{subsec:DGModinM}.

\begin{dfn} \label{dfn:4781}
The boundedness conditions in Definition \ref{dfn:4780} have abbreviations, 
that we call {\em boundedness indicators}%
\index{Boundedness! indicator}.
\begin{enumerate}
\item The indicator ``$-$'' stands for ``bounded above''.

\item The indicator ``$+$'' stands for ``bounded below''.

\item The indicator ``$\mrm{b}$'' stands for ``bounded''.

\item The indicator ``$\bra{\mrm{empty}}$'' (the empty symbol)
stands for the empty condition (no condition at all). 
\end{enumerate}
\end{dfn}

The symbol ``$\star$'' is used to denote an unspecified boundedness 
indicator, i.e.\ it is a placeholder for one of the boundedness indicators in 
Definition \ref{dfn:4781}.

\begin{dfn} \label{dfn:4782}
Let $S$ be an integer interval.
\begin{enumerate}
\item If $S =  [i_0, i_1]$, then its {\em length}%
\index{Integer interval! length of}
is $i_1 - i_0 \in \N \cup \{ \infty \}$. If $S = \varnothing$, then then its 
length is $-\infty$.  

\item If $S =  [i_0, i_1]$, then the {\em reversed integer interval}%
\index{Integer interval! reversed}
is $-S := [-i_1, -i_0]$.
If $S = \varnothing$, then its reverse is also $\varnothing$.

\item If $S =  [i_0, i_1]$, and if $T = [j_0, j_1]$ is a second nonempty 
integer interval, then their sum%
\index{Integer interval! sum of {\indash}s}
is $S + T := [i_0 + j_0, i_1 + j_1]$.
If $S = \varnothing$, and it $T$ is a second integer interval, then the sum is
$S + T := \varnothing$. 
\end{enumerate}
\end{dfn}

Thus an integer interval $S$ (whether empty or not) has finite length if and 
only if it is bounded. 

\begin{dfn}  \label{dfn:2125}
Let $M = \{ M^i \}_{i \in \Z}$ be a graded object in an abelian category 
$\cat{M}$, i.e.\ $M \in \dcat{G}(\cat{M})$. 
\begin{enumerate}
\item We say that $M$ is {\em concentrated} in an integer interval 
$S$ if the set \lb $\{ i \in \Z \mid M^i \neq 0 \}$ is contained in $S$.

\item The {\em concentration}%
\index{Concentration! of graded module}
\index{1-Con(M)@$\opn{con}(M)$}
of $M$ is the smallest integer interval 
$\opn{con}(M)$ in which $M$ is concentrated.

\item We say that $M$ is of {\em boundedness type $\star$},
for some boundedness indicator $\star$
if the integer interval $\opn{con}(M)$ satisfies the boundedness condition 
$\star$.
\end{enumerate}
\end{dfn}

Observe that  $\opn{con}(M) = \varnothing$ if and only if $M = 0$. 
Definition \ref{dfn:2125}(3) agrees with the boundedness 
conditions on a graded module $M$ from Subsections 
\ref{subsec:DGModinM} and  \ref{subsec:bd-in-homot}.

In Definition \ref{dfn:3155} we introduced the notions $\inf(M)$ and $\sup(M)$ 
for a graded object $M$. Using them we can describe the 
concentration $\opn{con}(M) = [i_0, i_1]$ when $M$ is nonzero:
$i_0 = \opn{inf} (M)$ and $i_1 = \opn{sup} (M)$.
The amplitude $\opn{amp}(M)$ is the length of the interval $\opn{con}(M)$.

Recall that to a DG module $M \in \dcat{C}^{}(A, \cat{M})$ we associate
the graded object $M^{\natural} \in \dcat{G}(\cat{M})$ that 
is gotten by forgetting the differential. The cohomology $\opn{H}(M)$
is another  object of $\dcat{G}(\cat{M})$. 

\begin{dfn} \label{dfn:4785}
Let  $A$ be a DG ring, $\cat{M}$ an abelian category, and 
$M \in \dcat{C}^{}(A, \cat{M})$.
\begin{enumerate}
\item The {\em concentration} of $M$%
\index{Concentration! of DG module}
and the {\em boundedness type}%
\index{Boundedness! type of DG module}
of $M$ are the same as in Definition \ref{dfn:2125}, when $M$ is seen as an 
object of $\dcat{G}^{}(\cat{M})$.

\item The {\em cohomological concentration}%
\index{Concentration! cohomological {\indash} of DG module}
and the {\em cohomological boundedness type} of $M$%
\index{Boundedness! type of DG module}
are the concentration and boundedness type, respectively, of the 
graded object $\opn{H}(M) \in \dcat{G}^{}(\cat{M})$,
in the sense of Definition \ref{dfn:2125}.
\end{enumerate}
\end{dfn}

The next definition is in conflict with Definitions \ref{dfn:1460} and 
\ref{dfn:1461}; but we already warned (in Remark \ref{rem:1461}) that this 
change will take place. See Remark \ref{rem:4765} for further discussion. 

\begin{dfn} \label{dfn:2122}
Let $A$ be a DG ring, $\cat{M}$ an abelian category, and $\star$ a boundedness 
indicator. 
We denote by $\dcat{D}^{\star}(A, \cat{M})$%
\index{1-D(A,M)@$\dcat{D}^{\star}(A, \cat{M})$}%
\index{1-D(M)@$\dcat{D}^{\star}(\cat{M})$}%
\index{Boundedness! condition}%
\index{Concentration! cohomological {\indash} of DG module}
the  full subcategory of $\dcat{D}(A, \cat{M})$ 
on the DG modules $M$ of cohomological boundedness type $\star$, as in 
Definition \ref{dfn:4785} above. 
\end{dfn}

Thus, for example,  a DG module $M$ belongs to 
$\dcat{D}^{\mrm{b}}(A, \cat{M})$ if and only if 
$\opn{con}(\opn{H}(M))$ is a bounded integer interval. 
As explained in subsection \ref{subsec:bd-in-homot}, 
$\dcat{D}^{\star}(A, \cat{M})$ is a full triangulated subcategory of 
$\dcat{D}^{}(A, \cat{M})$.

The notations $\dcat{C}^{\star}(A, \cat{M})$ and 
$\dcat{K}^{\star}(A, \cat{M})$
remain the same, as in Definitions \ref{dfn:3220} and \ref{dfn:3221}. 
Namely, these are the full subcategories of $\dcat{C}^{}(A, \cat{M})$ 
and $\dcat{K}^{}(A, \cat{M})$ respectively on the DG modules $M$ 
of boundedness type $\star$.%
\index{1-C(A,M)@$\dcat{C}^{\star}(A, \cat{M})$}%
\index{1-K(A,M)@$\dcat{K}^{\star}(A, \cat{M})$}%
\index{1-C(M)@$\dcat{C}^{\star}(\cat{M})$}%
\index{1-K(M)@$\dcat{K}^{\star}(\cat{M})$}

\begin{dfn} \label{dfn:2120}
Let $A$ be a DG ring and $\cat{M}$ an abelian category. 
For a DG module $M \in \dcat{C}(A, \cat{M})$ and an integer $i$, we write 
$M[i] := \opn{T}^i(M)$,
the $i$-th translation of $M$, 
as in Definition \ref{dfn:1170}.%
\index{1-M[i]@$M[i]$}%
\index{1-Ti(M)@$\opn{T}^i(M)$}%
\index{Translation! of DG module}

This notation applies also to the homotopy category $\dcat{K}(A, \cat{M})$,
the derived category $\dcat{D}(A, \cat{M})$, and every other T-additive 
category. 
\end{dfn}

The notation $M[i]$ makes it difficult to use the little t operator, and to 
talk about translation isomorphisms, but hopefully we won't require them 
anymore.

\begin{rem} \label{rem:4765}
We should say a few explanatory words on the change in Definition 
\ref{dfn:2122}, that conflicts with earlier definitions. 

Recall that in Definitions \ref{dfn:1460} and \ref{dfn:1461} we made a 
distinction between 
\begin{equation} \label{4770}
\dcat{D}^{\star}(A, \cat{M}) = 
\dcat{K}^{\star}(A, \cat{M})_{\dcat{S}^{\star}(A, \cat{M})} ,
\end{equation}
the localization of $\dcat{K}^{\star}(A, \cat{M})$ with respect to the 
quasi-isomorphisms in it, and $\dcat{D}(A, \cat{M})^{\star}$, 
the full subcategory of $\dcat{D}(A, \cat{M})$ on the complexes of 
cohomological boundedness type $\star$. This distinction was necessary in the 
early part of the book, when we were studying foundations and resolutions.
Now that the foundations are behind us, we no longer have a need for the 
subcategory $\dcat{D}^{\star}(-)$ in the sense of formula (\ref{4770}), i.e.\ 
the definition based on boundedness type. 

Still, one might ask: Why make the notational switch from 
$\dcat{D}(-)^{\star}$ to $\dcat{D}^{\star}(-)$, for the definition based on 
cohomological boundedness type? 
There are two reasons for that. The first is typographic: 
the notation $\dcat{D}(-)^{\star}$ becomes cumbersome when 
coupled with other indicators, such as
$\dcat{D}_{\mrm{f}}(-)^{\star}$
or $(\dcat{D}(-)^{\star})^{\mrm{op}}$.
After Definition \ref{dfn:2122} we can refer to these categories by the more 
elegant expressions $\dcat{D}^{\star}_{\mrm{f}}(-)$
or $\dcat{D}^{\star}(-)^{\mrm{op}}$, respectively. 

Second, the expression $\dcat{D}(-)^{\star}$ is a novelty, and nobody else uses 
it, whereas the notation $\dcat{D}^{\star}(-)$ is common. According 
to Proposition \ref{prop:1463}, in most situations the two distinct meanings of
$\dcat{D}^{\star}(-)$, the old one and the new one, refer to nested 
subcategories of $\dcat{D}(-)$ whose inclusion is an equivalence. Furthermore, 
in a few recent texts on derived categories (such as 
\cite[Definition tag=05RU]{SP} and \cite[Definition 1.11]{Ye8}) we find the 
same usage of the expression $\dcat{D}^{\star}(-)$ as in Definition 
\ref{dfn:2122}.
\end{rem}

\mysubsection{The Bifunctor \texorpdfstring{$\opn{RHom}$}{RHom}}   
\label{subsec:RHom}

Consider a DG ring $A$ and an abelian category $\cat{M}$.
Like in Example \ref{exa:1967} we get a DG bifunctor 
\[ F := \opn{Hom}_{A, \cat{M}}(-, -) : \dcat{C}(A, \cat{M})^{\mrm{op}} \times
\dcat{C}(A, \cat{M}) \to \dcat{C}(\K) . \]
In Definition \ref{dfn:2995} we put a triangulated structure on the opposite 
homotopy category $\dcat{K}(A, \cat{M})^{\mrm{op}}$. 
According to Theorem \ref{thm:4225}, there is an induced  
triangulated bifunctor 
\[ F = \opn{Hom}_{A, \cat{M}}(-, -) : \dcat{K}(A, \cat{M})^{\mrm{op}} \times
\dcat{K}(A, \cat{M}) \to \dcat{K}(\K) . \]
Postcomposing with the localization functor 
$\opn{Q} : \dcat{K}(\K) \to \dcat{D}(\K)$, 
we obtain a triangulated bifunctor 
\[ F = \opn{Hom}_{A, \cat{M}}(-, -) : \dcat{K}(A, \cat{M})^{\mrm{op}} \times
\dcat{K}(A, \cat{M}) \to \dcat{D}(\K) . \]
Next we pick full additive subcategories 
$\cat{K}_1, \cat{K}_2 \sub \dcat{K}(A, \cat{M})$
s.t.\ $\cat{K}_1^{\mrm{op}} \sub \dcat{K}(A, \cat{M})^{\mrm{op}}$
and
$\cat{K}_2 \sub \dcat{K}(A, \cat{M})$
are triangulated. In practice this choice would be by some boundedness 
conditions; for instance
$\cat{K}_1 := \dcat{K}^{-}(\cat{M})$
or $\cat{K}_2 := \dcat{K}^{+}(\cat{M})$,  
cf.\ Corollaries \ref{cor:1900} and \ref{cor:1630} respectively. 
We want to construct the right derived 
bifunctor of the triangulated bifunctor 
\[ F = \opn{Hom}_{A, \cat{M}}(-, -) : \cat{K}_1^{\mrm{op}} \times
\cat{K}_2 \to \dcat{D}(\K) . \]
This is done in the next theorem.

\begin{thm} \label{thm:2005}
Let $A$ be a DG $\K$-ring, let $\cat{M}$ be a $\K$-linear abelian 
category, and let $\cat{K}_1, \cat{K}_2 \sub \dcat{K}(A, \cat{M})$ 
be full additive subcategories, such that $\cat{K}_1^{\mrm{op}}$ is a 
full triangulated subcategory of $\dcat{K}(A, \cat{M})^{\mrm{op}}$, 
and $\cat{K}_2$ is a full triangulated subcategory of 
$\dcat{K}(A, \cat{M})$.
Assume either that $\cat{K}_1$ has enough K-projectives, or that $\cat{K}_2$ 
has 
enough K-injectives.
Let $\cat{D}_i$ denote the localization of $\cat{K}_i$ with respect to the 
quasi-isomorphisms in it. 
Then the triangulated bifunctor 
\[  \opn{Hom}_{A, \cat{M}}(-, -) : \cat{K}_1^{\mrm{op}} \times
\cat{K}_2 \to \dcat{D}(\K) \]
has a right derived bifunctor%
\index{1-RHom(A,M)@$\opn{RHom}_{A, \cat{M}}$}
\[  \opn{RHom}_{A, \cat{M}}(-, -) : \cat{D}_1^{\mrm{op}} \times
\cat{D}_2 \to \dcat{D}(\K) . \]
Moreover, if $P_1 \in \cat{K}_1$ is K-projective, or if 
$I_2 \in \cat{K}_2$ is K-injective, then the morphism 
\[ \eta^{\mrm{R}}_{P_1, I_2} : \opn{Hom}_{A, \cat{M}}(P_1, I_2) \to 
\opn{RHom}_{A, \cat{M}}(P_1, I_2) \]
in $\dcat{D}(\K)$ is an isomorphism. 
\end{thm}

\begin{proof}
If $\cat{K}_2$ has enough K-injectives, then we can take 
$\cat{J}_2 := \cat{K}_{2, \mrm{inj}}$, the full subcategory on the K-injectives 
inside $\cat{K}_2$. And we take 
$\cat{J}_1 := \cat{K}_{1}$. 
We claim that the conditions of Theorem \ref{thm:2000} are satisfied. Condition 
(b) is simply the assumption that $\cat{K}_2$ has enough K-injectives. As for 
condition (a): this is Lemma 
\ref{lem:2000} below, with its condition (ii). 

On the other hand, if $\cat{K}_1$ has enough K-projectives, then we can take 
$\cat{J}_1 := \cat{K}_{1, \mrm{prj}}$, 
the full subcategory on the K-projectives inside $\cat{K}_1$. 
And we take $\cat{J}_2 := \cat{K}_{2}$.
We claim that the conditions of Theorem \ref{thm:2000}  are satisfied for 
$\cat{J}_1^{\mrm{op}} \sub \cat{K}_1^{\mrm{op}}$. 
Condition (b) is simply the assumption that $\cat{K}_1$ has enough 
K-projectives: a quasi-isomorphism 
$\rho : P \to M$ in $\cat{K}_1$ becomes a quasi-isomorphism 
$\mrm{Op}(\rho) : \mrm{Op}(M) \to \mrm{Op}(P)$ in 
$\cat{K}_1^{\mrm{op}}$, and 
$\mrm{Op}(P) \in \cat{J}_1^{\mrm{op}}$.
As for condition (a): this is Lemma \ref{lem:2000} below,  with its condition 
(i). 

The last assertion also follows from \ref{lem:2000}.
\end{proof}

\begin{lem}  \label{lem:2000}
Suppose $\phi_1 : Q_1 \to P_1$ and $\phi_2 : I_2 \to J_2$
are quasi-iso\-morphisms in $\dcat{C}(A, \cat{M})$, 
and either of the conditions below holds\tup{:}
\begin{enumerate}
\rmitem{i} $Q_1$ and $P_1$ are both K-projective. 

\rmitem{ii} $I_2$ and $J_2$ are both K-injective.
\end{enumerate}
Then the homomorphism
\[  \opn{Hom}_{A, \cat{M}}(\phi_1, \phi_2) : 
\opn{Hom}_{A, \cat{M}}(P_1, I_2) \to \opn{Hom}_{A, \cat{M}}(Q_1, J_2) \]
in $\dcat{C}(\K)$ is a quasi-isomorphism. 
\end{lem}

\begin{proof}
We will only prove the case where $Q_1, P_1$ are both K-projective; the other 
case is very similar. 

The homomorphism in question factors as follows:
\[ \opn{Hom}_{A, \cat{M}}(\phi_1, \phi_2) = 
\opn{Hom}_{A, \cat{M}}(\phi_1, \opn{id}_{J_2})
\circ 
\opn{Hom}_{A, \cat{M}}(\opn{id}_{P_1}, \phi_2) . \]
It suffices to prove that each of the factors is a quasi-isomorphism. This can 
be done by a messy direct calculation, but we will provide an indirect proof 
that relies on properties of the homotopy category 
$\cat{K} := \dcat{K}(A, \cat{M})$ that 
were already established. 

Let $K_2$ be the standard cone on the homomorphism 
$\phi_2 : I_2 \to J_2$. So $K_2$ is acyclic. Because $P_1$ is K-projective it 
follows that $\opn{Hom}_{A, \cat{M}}(P_1, K_2)$ is acyclic. 
Thus for every integer $l$ we have 
\begin{equation} \label{eqn:2046}
 \opn{Hom}_{\cat{K}}(P_1[-l], K_2) \cong 
\opn{H}^{l} \bigl( \opn{Hom}_{A, \cat{M}}(P_1, K_2) \bigr) = 0 . 
\end{equation}
Next, there is a distinguished triangle 
\begin{equation} \label{eqn:2045}
I_2 \xar{\phi_2} J_2 \xar{\be_2} K_2 \xar{\ga_2} I_2[1]
\end{equation}
in $\cat{K}$. Applying the cohomological functor
$\opn{Hom}_{\cat{K}}(P_1[-l], -)$
to the distinguished triangle (\ref{eqn:2045}) yields a long exact sequence, as 
explained in Subsection \ref{subsec:tr-coh-funcs}. From it we deduce that the 
homomorphisms 
\[ \opn{Hom}_{\cat{K}}(P_1[-l], I_2) \to 
\opn{Hom}_{\cat{K}}(P_1[-l], J_2) \]
are bijective for all $l$. Using the isomorphisms like (\ref{eqn:2046}) for 
$I_2$ and $J_2$ we see that 
\[ \opn{Hom}_{A, \cat{M}}(\opn{id}_{P_1}, \phi_2) : 
\opn{Hom}_{A, \cat{M}}(P_1, I_2) \to 
\opn{Hom}_{A, \cat{M}}(P_1, J_2) \]
is a quasi-isomorphism. 

According to Corollary \ref{cor:1923} the homomorphism 
$\phi_1 : Q_1 \to P_1$ is a homotopy equivalence; so it is an isomorphism in 
$\cat{K}$. Therefore for every integer $l$ the homomorphism 
\[ \opn{Hom}_{\cat{K}}(Q_1, J_2[l]) \to 
\opn{Hom}_{\cat{K}}(P_1, J_2[l]) \]
is bijective. As above we conclude that
\[ \opn{Hom}_{A, \cat{M}}(\phi_1, \opn{id}_{J_2}) : 
\opn{Hom}_{A, \cat{M}}(Q_1, J_2) \to 
\opn{Hom}_{A, \cat{M}}(P_1, J_2) \]
is a quasi-isomorphism. 
\end{proof}

\begin{rem} \label{rem:2040}
Theorem \ref{thm:2005} should be viewed as a template. It has 
commutative and noncommutative variants, that we shall study later. 
And there are geometric variants in which the source and target are categories 
of sheaves.
\end{rem}

We end this section with the connection between RHom and morphisms in the 
derived category. 

\begin{dfn} \label{dfn:2121}
Under the assumptions of Theorem \tup{\ref{thm:2005}},  
for DG modules $M_1 \in \cat{K}_1$ and $M_2 \in \cat{K}_2$,
and for an integer $i$, we write%
\index{1-ExtiA@$\opn{Ext}^i_{A, \cat{M}}(-, -)$}
\[ \opn{Ext}^i_{A, \cat{M}}(M_1, M_2) := 
\opn{H}^i \bigl( \opn{RHom}_{A, \cat{M}}(M_1, M_2) \bigr) \in \dcat{M}(\K) . \]
\end{dfn}

\begin{exer} \label{exer:2125}
Let $A$ be a ring. Prove that for modules 
$M_1, M_2 \in \dcat{M}(A)$ the $\K$-module
$\opn{Ext}^i_{A}(M_1, M_2)$ defined above is canonically isomorphic to the 
module in the classical definition. Moreover, this is true regardless of the 
choices of the subcategories $\cat{K}_1 $ and $\cat{K}_2$,
as long as $M_i \in \cat{K}_i$. 
\end{exer}

\begin{cor} \label{cor:2120}
Under the assumptions of Theorem \tup{\ref{thm:2005}}, there is an  
isomorphism 
\[  \opn{Ext}^0_{A, \cat{M}}(-, -) \iso 
\opn{Hom}_{\dcat{D}(A, \cat{M})}(-, -) \]
of additive bifunctors 
$\cat{D}_1^{\mrm{op}} \times \cat{D}_2 \to \dcat{M}(\K)$.
\end{cor}

\begin{exer} \label{exer:2140}
Prove Corollary \ref{cor:2120}. (Hint: use Theorems \ref{thm:3135} and 
\ref{thm:3145}.)
\end{exer}

\mysubsection{The Bifunctor \texorpdfstring{$\ot^{\mrm{L}}$}{LTors}}  
\label{subsec:LTors}

Consider a DG ring $A$. Like in Example \ref{exa:1965} we get a DG bifunctor 
\[ F := (- \ot_A -) : \dcat{C}(A^{\mrm{op}}) \times \dcat{C}(A) 
\to \dcat{C}(\K) . \]
Passing to homotopy categories, and postcomposing with 
$\opn{Q} : \dcat{K}(\K) \to \dcat{D}(\K)$, 
we obtain a triangulated bifunctor 
\[ F = (- \ot_A -) : \dcat{K}(A^{\mrm{op}}) \times
\dcat{K}(A) \to \dcat{D}(\K) . \]
Next we pick full triangulated subcategories 
$\cat{K}_1 \sub \dcat{K}(A^{\mrm{op}})$
and
$\cat{K}_2 \sub \dcat{K}(A)$. 
In practice this choice would be by some boundedness conditions; for instance
$\cat{K}_1 := \dcat{C}^{-}(A^{\mrm{op}})$
or $\cat{K}_2 := \dcat{C}^{-}(A)$, cf.\ Corollary \ref{cor:1900}.
We want to construct the left derived 
bifunctor of the triangulated bifunctor 
\[ F = (- \ot_A -) : \cat{K}_1 \times
\cat{K}_2 \to \dcat{D}(\K) . \]
This is done in the next theorem. 

\begin{thm} \label{thm:2107}
Let $A$ be a DG $\K$-ring, and let 
$\cat{K}_1 \sub \dcat{K}(A^{\mrm{op}})$
and $\cat{K}_2 \sub \dcat{K}(A)$
be full triangulated subcategories.
Assume that either $\cat{K}_1$ or $\cat{K}_2$ has enough K-flat objects.
Let $\cat{D}_i$ denote the localization of $\cat{K}_i$ with respect to the 
quasi-iso\-morphisms in it. 
Then the triangulated bifunctor 
\[  (- \ot_A -)  : \cat{K}_1 \times \cat{K}_2 \to \dcat{D}(\K) \]
has a left derived bifunctor%
\index{1-LTens@$(- \ot^{\mrm{L}}_A -)$} 
\[ (- \ot^{\mrm{L}}_A -)  : \cat{D}_1 \times \cat{D}_2 \to \dcat{D}(\K) . \]
Moreover, if either $P_1 \in \cat{K}_1$ or $P_2 \in \cat{K}_2$ is K-flat, 
then the morphism 
\[ \eta^{\mrm{L}}_{P_1, P_2} : P_1 \ot^{\mrm{L}}_A P_2 \to P_1 \ot_A P_2 \]
in $\dcat{D}(\K)$ is an isomorphism. 
\end{thm}

\begin{proof}
If $\cat{K}_2$ has enough K-flats, then we can take 
$\cat{P}_2 := \cat{K}_{2, \mrm{flat}}$, the full subcategory on the 
K-flats inside $\cat{K}_2$. And we take 
$\cat{P}_1 := \cat{K}_{1}$. We claim that the conditions of Theorem 
\ref{thm:2106} are satisfied. Condition 
(b) is simply the assumption that $\cat{K}_2$ has enough K-flats. As for 
condition (a): this is Lemma \ref{lem:2109} below,
with its condition (ii). 

The other case is proved the same way (but replacing sides).  
The last assertion also follows from \ref{lem:2109}.
\end{proof}

\begin{lem}  \label{lem:2109}
Suppose 
$\phi_1 : P_1 \to Q_1$ and $\phi_2 : P_2 \to Q_2$
are quasi-isomorphisms in $\dcat{C}(A^{\mrm{op}})$
and $\dcat{C}(A)$ respectively, and either of the conditions below holds\tup{:}
\begin{enumerate}
\rmitem{i} $Q_1$ and $P_1$ are both K-flat.
\rmitem{ii} $P_2$ and $Q_2$ are both K-flat.
\end{enumerate}
Then the homomorphism
\[ \phi_1 \ot \phi_2 : P_1 \ot_A P_2 \to Q_1 \ot_A Q_2 \]
in $\dcat{C}(\K)$ is a quasi-isomorphism. 
\end{lem}

\begin{proof}
We will only prove the lemma under condition (i); the other 
case is very similar. The homomorphism in question factors as follows:
\[ \phi_1 \ot \phi_2 = (\phi_1 \ot \opn{id}_{P_2}) \circ 
(\opn{id}_{P_1} \ot \, \phi_2 ) . \]
It suffices to prove that each of the factors is a quasi-isomorphism. This can 
be done by a messy direct calculation, but we will provide an indirect proof 
that relies on properties of the DG categories  
$\dcat{C}(A^{\mrm{op}})$ and $\dcat{C}(A)$ that 
were already established. 

First we shall prove that $\opn{id}_{P_1} \ot \, \phi_2$ is a quasi-isomorphism.
Let $R_2$ be the standard cone on the strict homomorphism 
$\phi_2 : P_2 \to Q_2$. So there is a standard triangle 
\begin{equation} \label{eqn:2110}
P_2 \xar{\phi_2} Q_2 \to R_2 \to P_2[1]
\end{equation}
in $\dcat{C}_{\mrm{str}}(A)$, and $R_2$ is acyclic. 
Applying the DG functor $P_1 \ot_A (-)$ to the triangle (\ref{eqn:2110}), and 
using Theorem \ref{thm:1185}, we see that there is a standard triangle 
\begin{equation} \label{eqn:2111}
 P_1 \ot_A P_2 \xar{\opn{id}_{P_1} \ot \, \phi_2} P_1 \ot_A Q_2 \to P_1 \ot_A 
R_2 \to (P_1 \ot_A P_2)[1]
\end{equation}
in $\dcat{C}(\K)$. This becomes a distinguished triangle in the triangulated 
category $\dcat{K}(\K)$. Thus there is a long exact sequence in cohomology 
associated to (\ref{eqn:2111}). Because $P_1$ is K-flat it follows that 
$P_1 \ot_A R_2$ is acyclic. We conclude that 
$ \opn{H}^i(\opn{id}_{P_1} \ot \, \phi_2)$
is bijective for all $i$. 

Now we shall prove that $\phi_1 \ot \opn{id}_{P_2}$ is a quasi-isomorphism.
Let $R_1 \in \dcat{C}(A^{\mrm{op}})$ be the standard cone on the homomorphism 
$\phi_1 : P_1 \to Q_1$. It is both acyclic and K-flat. 
Using standard triangles like (\ref{eqn:2110}) and (\ref{eqn:2111}) we reduce 
the problem to showing that $R_1 \ot_A P_2$ is acyclic. 
According to Corollary \ref{cor:1580} and Proposition \ref{prop:1525} there is 
a quasi-isomorphism 
$\til{P}_2 \to P_2$ in $\dcat{C}(A)$ from some K-flat DG module $\til{P}_2$. 
As already proved in the previous paragraph, since $R_1$ is K-flat, the 
homomorphism
$R_1 \ot_A \til{P}_2 \to R_1 \ot_A P_2$
is a quasi-isomorphism. But $R_1$ is acyclic and $\til{P}_2$ is K-flat, and 
therefore
$R_1 \ot_A \til{P}_2$ is acyclic. We conclude that 
$R_1 \ot_A P_2$ is acyclic, as required. 
\end{proof}

\begin{rem} \label{rem:2105} 
Theorem \ref{thm:2107} should be viewed as a template. It has commutative and 
noncommutative variants, that we will talk about later. 
And there are geometric variants where the source and target are categories of 
sheaves.
\end{rem}

\begin{dfn} \label{dfn:4295} 
Under the assumptions of Theorem \tup{\ref{thm:2107}},  
for DG modules $M_1 \in \cat{K}_1$ and $M_2 \in \cat{K}_2$,
and for an integer $i$, we write%
\index{1-ToriA@$\opn{Tor}_i^{A}(-, -)$}
\[ \opn{Tor}_i^{A}(M_1, M_2) := 
\opn{H}^{-i} ( M_1 \ot^{\mrm{L}}_{A} M_2 ) \in \dcat{M}(\K) . \]
\end{dfn}

\begin{exer} \label{exer:4295} 
Let $A$ be a ring. Prove that for modules
$M_1 \in \dcat{M}(A^{\mrm{op}})$ and $M_2 \in \dcat{M}(A)$ the $\K$-module
$\opn{Tor}_i^{A}(M_1, M_2)$ defined above is canonically isomorphic to the 
module in classical definition. Moreover, this is true regardless of the 
choices of the subcategories $\cat{K}_1 $ and $\cat{K}_2$,
as long as $M_i \in \cat{K}_i$. 
\end{exer}

Recall that the category $\dcat{D}(A)$ admits infinite direct sums. 

\begin{prop} \label{prop:4300}
For every $M \in \dcat{D}(A^{\mrm{op}})$ the functor 
\[ M \ot^{\mrm{L}}_{A} (-) : \dcat{D}(A) \to \dcat{D}(\K) \]
commutes with infinite direct sums. 
\end{prop}

\begin{proof}
Fix a K-flat resolution $P \to M$ in $\dcat{C}_{\mrm{str}}(A^{\mrm{op}})$, so 
we have an isomorphism of functors 
$M \ot^{\mrm{L}}_{A} (-) \cong P \ot_A (-)$.
The functor 
\[ P \ot_A (-) : \dcat{C}_{\mrm{str}}(A) \to 
\dcat{C}_{\mrm{str}}(\K) \]
commutes with infinite direct sums. But by Theorem \ref{thm:3140}
the direct sums in $\dcat{C}_{\mrm{str}}(-)$ and in
$\dcat{D}(-)$ are the same. 
\end{proof}

\mysubsection{Cohomological Dimensions of Functors and Objects} 
\label{subsec:way-out}

The material here is a refinement of the notion of ``way-out functors'' from 
\cite[Section II.7]{RD}. Most of it is taken from 
\cite{RD} and \cite{Ye8}. As always, there is a fixed base ring $\K$, and 
Convention \ref{conv:2490} is in force. 

Integer intervals and their properties were introduced in Subsection 
\ref{subsec:bounded-revis}.

\begin{dfn} \label{dfn:2126}
Let $A, B$ be DG rings, let $\cat{M}, \cat{N}$ be abelian categories, and let 
$\cat{C}_0 \sub \cat{C} \sub \dcat{D}^{}(A, \cat{M})$
be full subcategories. 
\begin{enumerate}
\item Let  
$F : \cat{C} \to \dcat{D}(B, \cat{N})$
be an additive functor, and let $S$ be an integer interval. We say that 
{\em $F$ has cohomological displacement at 
most $S$ relative to $\cat{C}_0$}%
\index{Cohomological displacement! of covariant functor}
if 
\[ \opn{con} \bigl( \opn{H}(F(M)) \bigr) \subseteq
\opn{con} \bigl( \opn{H}(M) \bigr) + S \]
for every $M \in \cat{C}_0$.

\item Let 
$F : \cat{C}^{\mrm{op}} \to \dcat{D}(B, \cat{N})$
be an additive functor, 
and let $S$ be an integer interval. We say that {\em $F$ has 
cohomological displacement at most 
$S$ relative to $\cat{C}_0$}%
\index{Cohomological displacement! of contravariant functor}
if 
\[ \opn{con} \bigl( \opn{H}(F(M)) \bigr) \subseteq
- \opn{con} \bigl( \opn{H}(M) \bigr) + S \]
for every $M \in \cat{C}_0$.

\item Let $F$ be as in item (1) or (2).
The {\em cohomological displacement of $F$ relative to $\cat{C}_0$} is the 
smallest integer interval $S$ for which $F$ has cohomological 
displacement at most $S$ relative to $\cat{C}_0$. 

\item Let $S$ be the cohomological displacement of $F$ relative to 
$\cat{C}_0$. The 
{\em cohomological dimension of $F$  relative to $\cat{C}_0$}%
\index{Cohomological dimension! of functor}
is defined to be the length of the integer interval $S$. 

\item In case $\cat{C}_0 = \cat{C}$, we omit the clause ``relative to 
$\cat{C}_0$'' in all items above. 
\end{enumerate}
\end{dfn}

To emphasize the most important case of item (4) of the definition: {\em The 
functor $F$ has finite cohomological dimension if its cohomological 
displacement 
is bounded}. 

\begin{exa} \label{exa:2140}
The functor $F$ is the zero functor iff it has cohomological displacement 
$\varnothing$ and cohomological dimension $-\infty$. 
\end{exa}

\begin{exa} \label{exa:4285}
Let $F : \dcat{D}(A, \cat{M}) \to \dcat{D}(B, \cat{N})$
be an additive functor of finite cohomological dimension.
Then for each boundedness indicator $\star$ we have 
$F \bigl( \dcat{D}^{\star}(A, \cat{M}) \bigr) \sub 
\dcat{D}^{\star}(B, \cat{N})$.
\end{exa}

\begin{exa} \label{exa:2125}
Consider a nonzero commutative ring $A = B$, and the abelian category
$\cat{M} = \cat{N} := \dcat{M}(\K)$. So 
$\dcat{D}(A, \cat{M}) = \dcat{D}(B, \cat{N}) = \dcat{D}(A)$.
Take $\cat{C} := \dcat{D}(A)$.
For the covariant case (item (1) in Definition \ref{dfn:2126}) take a nonzero 
projective module $P$, and let
\[ F := \opn{RHom}_A \bigl( P \oplus P[1], - \bigr) : \dcat{D}(A) \to 
\dcat{D}(A) . \]
Then $F$ has cohomological displacement $[0, 1]$. 
For the contravariant case (item (2)) take a nonzero 
injective module $I$, and let 
\[ F := \opn{RHom}_A \bigl( -, I \oplus I[1] \bigr) : \dcat{D}(A)^{\mrm{op}} 
\to \dcat{D}(A) . \]
Then $F$ has cohomological displacement $[-1, 0]$. 
In both cases the cohomological dimension of $F$ is $1$. 
\end{exa}

\begin{exa} \label{exa:2135}
Suppose $A$ and $B$ are rings and 
$F : \dcat{M}(A) \to \dcat{M}(B)$ is a left exact additive functor. 
We get a triangulated functor 
$\mrm{R} F : \dcat{D}(A) \to \dcat{D}(B)$,
and $\opn{H}^i(\mrm{R} F(M)) = \mrm{R}^i F(M)$ for all $M \in \dcat{M}(A)$. 
Taking $\cat{C} := \dcat{M}(A)$, with its canonical embedding into
$\dcat{D}(A)$, the cohomological dimension of $\mrm{R} F$ relative to 
$\dcat{M}(A)$ equals the usual right cohomological dimension of the functor 
$F$. 
\end{exa}

Next is an example about triangulated functors of finite 
cohomological dimensions. Another example, of a contravariant triangulated 
functor of finite cohomological dimension, is the duality functor $D$ 
associated to a dualizing complex; see Subsections \ref{subsec:du-cplxs},
\ref{subsec:graded-NCDC} and \ref{subsec:NC-DC}.

\begin{exa} \label{exa:4595}  
Let $A$ be a commutative ring, let $\bsym{a} = (a_1, \ldots, a_n)$ be a finite 
sequence of elements of $A$, and let $\a \sub A$ be the ideal generated by 
$\bsym{a}$. The {\em $\a$-adic completion functor} $\Lambda_{\a}$ and the 
{\em $\a$-torsion functor} $\Ga_{\a}$ have derived functors 
$\mrm{L} \Lambda_{\a}, \mrm{R} \Ga_{\a} : \dcat{D}(A) \to  \dcat{D}(A)$.

We call the sequence $\bsym{a}$ {\em weakly proregular} if  
the Koszul complexes associated to powers of $\bsym{a}$ satisfy a rather 
complicated asymptotic formula -- see \cite[Definition 4.21]{PSY}. If $A$ is 
noetherian then weak proregularity is 
automatic; but this condition is also satisfied in many important 
non-noetherian cases. 

When $\bsym{a}$ is weakly proregular, the functors 
$\mrm{L} \Lambda_{\a}$ and $\mrm{R} \Ga_{\a}$ have finite cohomological 
dimensions, bounded by the length of $\bsym{a}$.
Furthermore, the {\em commutative MGM Equivalence} holds. To explain it, we 
need to introduce two full triangulated subcategories of $\dcat{D}(A)$~: 
\begin{itemize}
\item The category $\dcat{D}(A)_{\mrm{com}}$ of {\em derived $\a$-adically 
complete complexes}. These are the complexes $M$ such that the canonical 
morphism $M \to \mrm{L} \Lambda_{\a}(M)$ is an isomorphism. 

\item The category $\dcat{D}(A)_{\mrm{tor}}$ of {\em derived $\a$-torsion
complexes}. These are the complexes $M$ such that the canonical morphism
$\mrm{R} \Ga_{\a}(M) \to M$ is an isomorphism. 
\end{itemize}
The MGM Equivalence says that the functor 
$\mrm{R} \Ga_{\a} : \dcat{D}(A)_{\mrm{com}} \to  \dcat{D}(A)_{\mrm{tor}}$ 
is an equivalence of triangulated categories, with quasi-inverse 
$\mrm{L} \Lambda_{\a}$.
The name ``MGM'' stands for ``Matlis-Greenlees-May''.
This duality was worked out in the paper \cite{PSY} by M. Porta, L. Shaul and 
A. Yekutieli; expanding earlier work by E. Matlis; A. Grothendieck; 
J. Greenlees and P. May;  L. Alonso, A. Jeremias and J. Lipman; and P. 
Schenzel. 

In Subsection \ref{subsec:NC-MGM} we give a noncommutative variant of the MGM 
Equivalence, adapted from the paper \cite{VyYe} by R. Vyas and Yekutieli. 
\end{exa}

\begin{rem} \label{rem:2126}
Assume that in Definition \ref{dfn:2126}(1) we take $A = B = \K$,  
$\cat{C} = \dcat{D}(\cat{M})$, and 
$F : \dcat{D}(\cat{M}) \to \dcat{D}(\cat{N})$
is a triangulated functor. The functor $F$ has bounded below (resp.\ above) 
cohomological displacement if and only if it is a {\em way-out right} (resp.\ 
{\em left}) functor, in the sense of \cite[Section I.7, Definition]{RD}.
For instance, if $F$ is a way-out right functor, with bounding integers 
$n_1$ and $n_2$ (as defined in \cite{RD}), then the cohomological displacement 
of $F$ is contained in the integer interval $[n_1 - n_2, \infty]$.
Conversely, if $F$ has cohomological displacement at most
$[i_0, \infty]$ for some integer $i_0$, then $F$ is way-out right, and for 
every $n_1 \in \Z$ the integer $n_2 :=  n_1 - i_0$ satisfies the condition in 
\cite{RD}. Likewise for way-out left functors.  
\end{rem}

\begin{dfn} \label{dfn:2128}
Let $\star_1, \star_2$ be boundedness indicators, and assume the right derived 
bifunctor 
\[ \opn{RHom}_{A, \cat{M}} : 
\dcat{D}^{\star_1}(A, \cat{M})^{\mrm{op}} \times 
\dcat{D}^{\star_2}(A, \cat{M}) \to \dcat{D}(\K) \]
exists. Let $S$ be an integer interval of length 
$i\in \N \cup \{\pm  \infty \}$.
\begin{enumerate}
\item Let $M \in \dcat{D}^{\star_1}(A, \cat{M})$, and let 
$\cat{C} \sub \dcat{D}^{\star_2}(A, \cat{M})$ be a full 
subcategory. We say that $M$ has {\em projective concentration} $S$
and {\em projective dimension}%
\index{Projective! concentration of DG module}%
\index{Projective! dimension of DG module}
$i$ relative to $\cat{C}$ if the functor
\[ \opn{RHom}_{A, \cat{M}}(M, -)|_{\cat{C}} : \cat{C} \to \dcat{D}(\K) \]
has cohomological displacement $-S$.

\item Let $M \in \dcat{D}^{\star_2}(A, \cat{M})$, and let 
$\cat{C} \sub \dcat{D}^{\star_1}(A, \cat{M})$ be a full  
subcategory. We say that $M$ has {\em injective concentration} $S$ 
and {\em injective dimension}%
\index{Injective! dimension of DG module}%
\index{Injective! concentration of DG module}
$i$ relative to $\cat{C}$ if the functor
\[ \opn{RHom}_{A, \cat{M}}(-, M)|_{\cat{C}^{\mrm{op}}} 
: \cat{C}^{\mrm{op}} \to \dcat{D}(\K) \]
has  cohomological displacement $S$.

\item If $\cat{C} = \dcat{D}^{}(A, \cat{M})$,
then we omit the clause ``relative to $\cat{C}$'' in items (1) and (2). 
\end{enumerate}
\end{dfn}

\begin{exa} \label{exa:2127}
Continuing with the setup of Example \ref{exa:2125}, the DG module
$P \oplus P[1]$ (resp.\ $I \oplus I[1]$) has projective (resp.\ injective)
concentration $[-1, 0]$.
\end{exa}

\begin{exa} \label{exa:2128}
Let $A$ be a DG ring such that $\opn{H}(A)$ is nonzero, and consider the free 
DG module $P := A \in \dcat{D}(A)$.
The functor 
$F := \opn{RHom}_{A}(P, -) : \dcat{D}(A) \to \dcat{D}(\K)$ 
is isomorphic to the forgetful functor, so it has cohomological displacement 
$[0, 0]$ and cohomological dimension $0$.
Thus the DG module $P$ has projective concentration $[0, 0]$ and projective 
dimension $0$.
Note however that the cohomology $\opn{H}(P)$ could be unbounded! 
\end{exa}

\begin{prop} \label{prop:4291} 
Let $\cat{M}$ be an abelian category with enough injectives. The following are 
equivalent for $M \in \cat{M}$~\tup{:}
\begin{enumerate} [itemsep=0.3ex]
\rmitem{i} $M$ is an injective object of $\cat{M}$.

\rmitem{ii} $\opn{Ext}^{p}_{\cat{M}}(N, M) = 0$ for every 
$N \in \cat{M}$ and every $p \geq 1$. 

\rmitem{iii}  $\opn{Ext}^1_{\cat{M}}(N, M) = 0$ for every 
$N \in \cat{M}$. 
\end{enumerate}
\end{prop}

Note that by Corollary \ref{cor:1902} and Definition \ref{dfn:2121}, the 
$\K$-modules $\opn{Ext}^p_{\cat{M}}(N, M)$ exist.

\begin{exer} \label{exer:4291}
Prove Proposition \ref{prop:4291}. (Hint: the proof is just like in the case 
$\cat{M} = \dcat{M}(A)$. It can be found in standard textbooks on homological 
algebra.)
\end{exer}

\begin{prop} \label{prop:2115}
Let $\cat{M}$ be an abelian category with enough injectives. The following are 
equivalent for $M \in \dcat{D}^+(\cat{M})$~\tup{:}
\begin{enumerate}
\rmitem{i} $M$ has finite injective dimension.

\rmitem{ii} $M$ has finite injective dimension relative to $\cat{M}$.

\rmitem{iii} There is a quasi-isomorphism $M \to I$ in 
$\dcat{C}_{\mrm{str}}(\cat{M})$ 
to a bounded complex of injective $A$-modules $I$.
\end{enumerate}
\end{prop}

Note that by Corollary \ref{cor:1902} we can apply Definition \ref{dfn:2128} 
with boundedness type $\star_2 = +$, so we can talk about the injective 
dimension of $M$. 

\begin{proof} \mbox{}

\smallskip \noindent
(i) $\Rightarrow$ (ii): This is trivial, since 
$\cat{M} \sub \dcat{D}(\cat{M})$. 

\medskip \noindent 
(ii)  $\Rightarrow$ (iii): We may assume that $\opn{H}(M)$ is nonzero. 
Let $S$ be the injective concentration of the complex $M$
relative to $\cat{M}$, as in Definition \ref{dfn:2128}; this is a bounded 
integer interval (possibly empty, a priory). 

Let $q$ be an integer such that $\opn{H}^q(M) \neq 0$.
Consider the object $N := \opn{Z}^q(M) \in \cat{M}$ and the canonical 
monomorphism $N \to M^q$. This can be viewed as a morphism 
$\phi : N[-q] \to M$ in $\dcat{C}_{\mrm{str}}(\cat{M})$, and then 
$\opn{H}^q(\phi) : N = \opn{H}^q(N[-q]) \to \opn{H}^q(M)$ 
is a nonzero morphism in $\cat{M}$. Therefore
$\opn{Q}(\phi) :  N[-q] \to M$
is a nonzero morphism in $\dcat{D}(\cat{M})$, so this is a nonzero element of 
$\opn{H}^q \bigl( \opn{RHom}_{\cat{M}}(N, M) \bigr)$. 
Looking at Definition \ref{dfn:2128} we see that $q \in S$.
Since $S$ is nonempty and bounded, it has to be 
$S = [q_0, q_1]$ for some integers
$q_0 = \opn{inf}(\opn{H}(M))$ and $q_1 = \opn{sup}(\opn{H}(M))$.

According to Corollary \ref{cor:1902} there is 
quasi-isomorphism $M \to J$, where $J$ is a complex of injective objects of 
$\cat{M}$, and $\opn{inf}(J) = q_0$. 
Take $I := \opn{smt}^{\leq q_1}(J)$, the smart truncation from Definition 
\ref{dfn:2320}. 
Then the canonical homomorphism 
$I \to J$ is a quasi-isomorphism. The complex $I$ is concentrated in the 
integer interval $[q_0, q_1]$, and
$I^q = J^q$ is an injective object for all $q < q_1$. 

Let us prove that $I^{q_1} = \opn{Z}^{q_1}(J)$ is also an injective object of 
$\cat{M}$. Classically we would use a cosyzygy argument. Here we use another 
trick. Define $I' := \opn{stt}^{\leq q_1 - 1}(I)$, 
the stupid truncation of $I$ from Definition \ref{dfn:3066}. so
\[ I' = \bigl( \cdots 0 \to I^{q_0} \to \cdots \to I^{q_1 - 1} \to 0 \to 
\cdots \bigr) . \]
This is a bounded complex of injective objects.
Consider the short exact sequence 
\[ 0 \to I^{q_1}[-q_1] \to I \to I' \to 0 \]
in $\dcat{C}^+_{\mrm{str}}(\cat{M})$.   
According to Proposition \ref{prop:2165} this gives a distinguished triangle 
\begin{equation} \label{eqn:2129}
I^{q_1}[-q_1] \to I \to I' \xar{ \ \triangle \ }  
\end{equation}
in $\dcat{D}^+(\cat{M})$. Take any object $L \in \cat{M}$. Applying 
$\opn{RHom}_{\cat{M}}(L, -)$ to the distinguished triangle (\ref{eqn:2129}) and 
then taking cohomologies, we get a long exact sequence
\[ \cdots \to \opn{Ext}_{\cat{M}}^{q + q_1 - 1}(L, I') \to 
\opn{Ext}_{\cat{M}}^{q}(L, I^{q_1}) \to 
\opn{Ext}_{\cat{M}}^{q + q_1}(L, I) \to \cdots \]
in $\dcat{M}(\K)$. For every $q > 0$ the $\K$-module 
$\opn{Ext}_{\cat{M}}^{q + q_1 - 1}(L, I')$ 
vanishes trivially. By the definition of the interval $[q_0, q_1]$, and the 
existence of an isomorphism $M \cong I$ in $\dcat{D}(\cat{M})$, for every 
$q > 0$ the $\K$-module $\opn{Ext}_{\cat{M}}^{q + q_1}(L, I)$ is zero. Hence 
$\opn{Ext}_{\cat{M}}^{q}(L, I^{q_1}) = 0$ for all $q > 0$. 
This proves that the object $I^{q_1}$ is injective (see Proposition 
\ref{prop:4291}). 

We have quasi-isomorphisms $M \to J$ and $I \to J$. Since $I$ is 
K-injective, there is a quasi-isomorphism $M \to I$. 

\medskip \noindent 
(iii) $\Rightarrow$ (i): This is also trivial.
\end{proof}

\begin{prop} \label{prop:4293} 
Let $\cat{M}$ be an abelian category with enough projectives. The following are 
equivalent for $M \in \cat{M}$~\tup{:}
\begin{enumerate} [itemsep=0.3ex]
\rmitem{i} $M$ is a projective object of $\cat{M}$.

\rmitem{ii} $\opn{Ext}^i_{\cat{M}}(M, N) = 0$ for every 
$N \in \cat{M}$ and every $i \geq 1$. 

\rmitem{iii}  $\opn{Ext}^1_{\cat{M}}(M, N) = 0$ for every 
$N \in \cat{M}$. 
\end{enumerate}
\end{prop}

Note that by Corollary \ref{cor:1750} and Definition \ref{dfn:2121}, the 
$\K$-modules $\opn{Ext}^i_{\cat{M}}(M, N)$ exist.

\begin{prop} \label{prop:4290}
Let $\cat{M}$ be an abelian category with enough projectives. The following are 
equivalent for $M \in \dcat{D}^-(\cat{M})$~\tup{:}
\begin{enumerate}
\rmitem{i} $M$ has finite projective dimension.

\rmitem{ii} $M$ has finite projective dimension relative to 
$\cat{M}$.

\rmitem{iii} There is a quasi-isomorphism $P \to M$ in 
$\dcat{C}_{\mrm{str}}(\cat{M})$ from a bounded complex of projective objects 
$P$.
\end{enumerate}
\end{prop}

Note that by Corollary \ref{cor:1750} we can apply Definition \ref{dfn:2128} 
with boundedness type $\star_1 = -$, so we can talk about the projective 
dimension of $M$. 

\begin{exer} \label{exer:2390}
Prove Propositions \ref{prop:4293} and \ref{prop:4290}. (Hint: compare to 
Propositions \ref{prop:4291} and \ref{prop:2115}.
Or look in any standard textbook on homological algebra.)
\end{exer}

\begin{dfn} \label{dfn:2390} 
Let $\star_1, \star_2$ be boundedness indicators, and assume the left 
derived bifunctor 
\[ (- \ot^{\mrm{L}}_A -) : 
\dcat{D}^{\star_1}(A^{\mrm{op}}) \times 
\dcat{D}^{\star_2}(A) \to \dcat{D}(\K) \]
exists. Let $S$ be an integer interval of length 
$i\in \N \cup \{\pm  \infty \}$.
\begin{enumerate}
\item Let $M \in \dcat{D}^{\star_2}(A)$, and let 
$\cat{C} \sub \dcat{D}^{\star_1}(A^{\mrm{op}})$ be a full  
subcategory. We say that $M$ has {\em flat concentration} $S$
and {\em flat dimension}%
\index{Flat! dimension of DG module}%
\index{Flat! concentration of DG module}
$i$ relative to $\cat{C}$ if the functor
$(- \ot^{\mrm{L}}_A M) |_{\cat{C}} : \cat{C} \to \dcat{D}(\K)$
has cohomological displacement $S$.

\item If $\cat{C} = \dcat{D}^{}(A^{\mrm{op}})$,
then we omit the clause ``relative to $\cat{C}$''. 
\end{enumerate}
\end{dfn}

\begin{prop} \label{prop:2390}
Let $A$ be a ring. The following are equivalent for $M \in \dcat{D}(A)$\tup{:}
\begin{enumerate}
\rmitem{i} $M$ has finite flat dimension.

\rmitem{ii} $M$ has finite flat dimension relative to 
$\dcat{M}(A^{\mrm{op}})$.

\rmitem{iii} There is an isomorphism $P \cong M$ in $\dcat{D}(A)$, 
where $P$ is a bounded complex of flat $A$-modules.
\end{enumerate}
\end{prop}

\begin{exer} \label{exer:2391}
Prove Proposition \ref{prop:2390}. (The proof is similar to that of Proposition 
\ref{prop:4290}. It can be found in many standard books on homological algebra.)
\end{exer}

\begin{rem} \label{rem:4785}
Propositions \ref{prop:2115}, \ref{prop:4290} and \ref{prop:2390} do not seem 
to have analogues for $M \in \dcat{D}(A)$, when $A$ is a DG ring that is not 
concentrated in degree $0$. 
\end{rem}

\mysubsection{Theorems on Functors Satisfying Finiteness Conditions} 
\label{subsec:thm-funcs-fin}

In this subsection we apply the definitions from the previous subsection 
(on cohomological dimensions of functors), combined with noetherian conditions,
to prove Theorems \ref{thm:2135}, \ref{thm:2160} and \ref{thm:4600}. 
These results (or their analogues) were originally proved in 
\cite[Section II.7]{RD}. 

As always, there is a fixed base ring $\K$, and 
Convention \ref{conv:2490} is in force.

\begin{dfn} \label{dfn:2135}
Suppose $A$ is a left noetherian ring. 
\begin{enumerate}
\item We denote by $\dcat{M}_{\mrm{f}}(A)$ the full subcategory of 
$\dcat{M}(A) = \cat{Mod}A$ on the finitely generated modules. 

\item  We denote by $\dcat{D}_{\mrm{f}}(A)$ the full subcategory of 
$\dcat{D}(A) = \dcat{D}(\cat{Mod} A)$ on the complexes with cohomology modules
in $\dcat{M}_{\mrm{f}}(A)$.
\end{enumerate}
\end{dfn}

Because $A$ is left noetherian, the category $\dcat{M}_{\mrm{f}}(A)$ is a thick
abelian subcategory of $\dcat{M}(A)$, and the category $\dcat{D}_{\mrm{f}}(A)$ 
is a full triangulated subcategory of $\dcat{D}(A)$.  
When viewed as a left module, 
$A \in \dcat{M}_{\mrm{f}}(A) \sub \dcat{D}_{\mrm{f}}^{\mrm{b}}(A)$.
Note the similarity to Definition \ref{dfn:3180}. 

\begin{thm} \label{thm:2135}
Let $A$ be a left noetherian ring, let $\cat{N}$ be an abelian category, 
let $\star$ be a boundedness indicator, let 
$F, G : \dcat{D}^{\star}_{\mrm{f}}(A) \to \dcat{D}(\cat{N})$ 
be triangulated functors, and let 
$\ze : F \to G$ be a morphism of triangulated functors. Assume that 
the morphism 
$\ze_A : F(A) \to G(A)$
in $\dcat{D}(\cat{N})$ is an isomorphism. 
\begin{enumerate}
\item If $\star = -$, and if $F$ and $G$ have bounded above cohomological 
displacements, then 
$\ze_M : F(M) \to G(M)$
is an isomorphism for every $M \in \dcat{D}^-_{\mrm{f}}(A)$.

\item If $\star = \bra{\mrm{empty}}$, and if $F$ and $G$ have bounded 
cohomological displacements, then $\ze_M$ is an isomorphism for every 
$M \in \dcat{D}_{\mrm{f}}(A)$.
\end{enumerate} 
\end{thm}

We shall require the next lemmas for the proof of the theorem. 

\begin{lem} \label{lem:2145}
Let $\cat{D}$ be a triangulated category, let 
$F, G : \cat{D} \to \dcat{D}(\cat{N})$ be triangulated functors, let 
$\ze : F \to G$ be a morphism of triangulated functors, and let  
$L \xar{\phi} M \to N \xar{ \ \triangle \ }$
be a distinguished triangle in $\cat{D}$.
\begin{enumerate}
\item If the morphisms $\ze_L$ and $\ze_M$ are both isomorphisms, then 
$\ze_N$ is an isomorphism. 

\item Let $j$ be an integer. If $\opn{H}^{j - 1}(F(N))$,
$\opn{H}^{j - 1}(G(N))$, $\opn{H}^{j}(F(N))$ and 
$\opn{H}^{j}(G(N))$ are all zero, and if  $\opn{H}^{j}(\ze_L)$
is an isomorphism, then $\opn{H}^{j}(\ze_M)$
is an isomorphism.
\end{enumerate} 
\end{lem}

\begin{proof}
(1) In $\dcat{D}(\cat{N})$ we get the commutative diagram 
\begin{equation} \label{eqn:2148}
\UseTips \xymatrix @C=8ex @R=6ex {
F(L)
\ar[r]
\ar[d]^{\ze_L}
&
F(M)
\ar[r]
\ar[d]^{\ze_M}
&
F(N)
\ar[r]
\ar[d]^{\ze_N}
&
F(L)[1]
\ar[d]^{\ze_L[1]}
\\
G(L)
\ar[r]
&
G(M)
\ar[r]
&
G(N)
\ar[r]
&
G(L)[1]
} 
\end{equation}
with horizontal distinguished triangles. According to 
Proposition \ref{prop:1283}, $\ze_N$ is an isomorphism. 

\medskip \noindent 
(2) Passing to cohomologies in (\ref{eqn:2148})
we have a commutative diagram with exact rows
\[ \UseTips \xymatrix @C=8ex @R=6ex {
\opn{H}^{j - 1}(F(N))
\ar[r]
\ar[d]^{\opn{H}^{j - 1}(\ze_{N})}
&
\opn{H}^{j}(F(L))
\ar[r]^(0.47){\opn{H}^{j}(F(\phi))}
\ar[d]^{\opn{H}^{j}(\ze_{L})}
&
\opn{H}^{j}(F(M))
\ar[r]
\ar[d]^{\opn{H}^{j}(\ze_{M})}
&
\opn{H}^{j}(F(N))
\ar[d]^{\opn{H}^{j}(\ze_{N})}
\\
\opn{H}^{j - 1}(G(N))
\ar[r]
&
\opn{H}^{j}(G(L))
\ar[r]^(0.47){\opn{H}^{j}(G(\phi))}
&
\opn{H}^{j}(G(M))
\ar[r]
&
\opn{H}^{j}(G(N)) 
} \]
in $\cat{N}$. The vanishing assumption implies that $\opn{H}^j(F(\phi))$ and
$\opn{H}^j(G(\phi))$ are isomorphisms. Hence 
$\opn{H}^{j}(\ze_{M})$ is an isomorphism.
\end{proof}

\begin{lem} \label{lem:2146}
Let $\cat{D}$ be a triangulated category, let 
$F, G : \cat{D} \to \dcat{D}(\cat{N})$ be triangulated functors, and let 
$\ze : F \to G$ be a a morphism of triangulated functors.
The following conditions are equivalent for $M \in \cat{D}$~\tup{:}
\begin{enumerate}
\rmitem{i} $\ze_M$ is an isomorphism. 

\rmitem{ii} $\ze_{M[i]}$ is an isomorphism for every integer $i$. 

\rmitem{iii} The morphism 
$\opn{H}^j(\ze_M) : \opn{H}^j(F(M)) \to \opn{H}^j(G(M))$
is an isomorphism for every integer $j$. 
\end{enumerate}
\end{lem}

\begin{proof}
The equivalence (i) $\Leftrightarrow$ (ii) is because both $F$ and $G$ are 
triangulated functors. 
The equivalence (i) $\Leftrightarrow$ (iii) is because
the functor 
$\opn{H} : \dcat{D}(\cat{N}) \to \dcat{G}_{\mrm{str}}(\cat{N})$ 
is conservative; see Corollary \ref{cor:2145}.
\end{proof}

\begin{proof}[Proof of Theorem \tup{\ref{thm:2135}}] \mbox{}

\smallskip \noindent
(1) First assume $P$ is a bounded complex of finitely generated free 
$A$-modules. Then $P$ is obtained from $A$ by finitely many standard cones and 
translations. By Lemmas \ref{lem:2145}(1) and \ref{lem:2146} it follows that 
$\ze_P$ is an isomorphism. 

Next let $P$ be a bounded above complex of finitely generated free $A$-modules.
Choose some integer $j$. 
Let $i_1$ be an integer such that the integer interval $[-\infty, i_1]$ 
contains the cohomological displacements of $F$ and $G$. 
Define $P' := \opn{stt}^{\leq j - i_1 - 2}(P)$, the stupid truncation of $P$ 
below $j - i_1 - 2$; and let $P'' := \opn{stt}^{\geq j - i_1 - 1}(P)$, the 
complementary  stupid truncation. See Definition \ref{dfn:3066}. 
According to Proposition \ref{prop:2165}, the short exact 
sequence (\ref{eqn:2147}) 
gives a distinguished triangle 
$P'' \to P \to P' \xar{\triangle}$ 
in $\dcat{D}^{-}_{\mrm{f}}(A)$. The complex $P''$ is a bounded complex of 
finitely 
generated free $A$-modules, so we already know that $\ze_{P''}$ is an 
isomorphism. Hence $\opn{H}^j(\ze_{P''})$ is an isomorphism.
On the other hand $\opn{H}(P')$ is concentrated in the degree interval 
$[-\infty, j - i_1 - 2]$. Therefore 
$\opn{H}^{k}(F(P')) = \opn{H}^k(G(P')) = 0$ for all 
$k \geq j - 1$. By Lemma \ref{lem:2145}(2), $\opn{H}^j(\ze_P)$ is an 
isomorphism. Because $j$ is arbitrary, Lemma \ref{lem:2146} says that 
$\ze_{P}$ is an isomorphism. 

Now take an arbitrary $M \in \dcat{D}^-_{\mrm{f}}(A)$.
By Theorem \ref{thm:3340} and Example \ref{exa:1930} there is a 
resolution $P \to M$, where $P$ is a bounded above complex of finitely 
generated free $A$-modules. Since $\ze_P$ is an isomorphism, so is $\ze_M$. 

\medskip \noindent
(2) Now we assume that the functors $F$ and $G$ have finite cohomological 
dimensions. Take any complex $M \in \dcat{D}_{\mrm{f}}(A)$. 
By Lemma \ref{lem:2146} it suffices to prove 
that $\opn{H}^j(\ze_{M})$ is an isomorphism for every integer $j$. 

Let $[i_0, i_1]$ be a bounded integer interval that contains the 
cohomological displacements of the functors $F$ and $G$. 
Define $M'' := \opn{smt}^{\leq j - i_0}(M)$, the smart truncation of $M$ 
below $j - i_0$; and let $M' := \opn{smt}^{\geq j - i_0 + 1}(M)$, the 
complementary smart truncation. See Definition \ref{dfn:2320}. 
According to Proposition \ref{prop:2320}
there is a distinguished triangle
$M'' \to M \to M' \xar{\triangle}$
in $\dcat{D}_{\mrm{f}}(A)$. 
The cohomologies of these complexes satisfy 
$\opn{H}^i(M'') = \opn{H}^i(M)$ and 
$\opn{H}^i(M') = 0$ for $i \leq j - i_0$; and 
$\opn{H}^i(M'') = 0$ and $\opn{H}^i(M') = \opn{H}^i(M)$ for 
$i \geq j - i_0 + 1$. 
Note that $M'' \in \dcat{D}^-_{\mrm{f}}(A)$. 

By part (1) we know that $\ze_{M''}$ is an 
isomorphism, and therefore also $\opn{H}^j(\ze_{M''})$ is an isomorphism. 
The cohomology $\opn{H}(M')$ is concentrated in the degree interval 
$[j - i_0 + 1, \infty]$, 
and therefore the cohomologies $\opn{H}(F(M'))$ and  $\opn{H}(G(M'))$
are concentrated in the interval $[j + 1, \infty]$. In particular the objects 
$\opn{H}^{j - 1}(F(M'))$, $\opn{H}^{j - 1}(G(M'))$, 
$\opn{H}^{j}(F(M'))$ and $\opn{H}^{j}(G(M'))$ are zero. 
By Lemma \ref{lem:2145}(2), $\opn{H}^j(\ze_M)$ is an isomorphism.
\end{proof}

The triangulated structure of $\dcat{D}(A)^{\mrm{op}}$ was introduced in 
Subsection \ref{subsec:opp-dercat-triang}. 
By our notational conventions, 
$\dcat{D}_{\mrm{f}}(A)^{\mrm{op}}$ is the full subcategory of 
$\dcat{D}(A)^{\mrm{op}}$ on the complexes with finitely generated cohomology 
modules, and it is triangulated, by Propositions \ref{prop:3051} and 
\ref{prop:2356}. 

\begin{exer} \label{exer:4325}
State and prove the contravariant modification of Theorem \lb \ref{thm:2135}.
(Hint: study the proofs of Theorems \ref{thm:2135} and \ref{thm:2160}.)
\end{exer}

\begin{thm} \label{thm:2160}
Let $A$ be a left noetherian ring, let $\cat{N}$ be an abelian category, let
$\cat{N}_0 \sub \cat{N}$ be a thick abelian subcategory, let $\star$ be a 
boundedness condition, and let 
$F : \dcat{D}^{\star}_{\mrm{f}}(A)^{\mrm{op}} \to \dcat{D}(\cat{N})$ 
be a triangulated functor. Assume that 
$F(A)$ belongs to $\dcat{D}_{\cat{N}_0}(\cat{N})$.
\begin{enumerate}
\item If $\star = -$, and if $F$ has bounded below cohomological displacement, 
then $F(M)$ belongs to $\dcat{D}_{\cat{N}_0}(\cat{N})$
for every $M \in \dcat{D}^-_{\mrm{f}}(A)$.

\item If $\star = \bra{\mrm{empty}}$, and if $F$ has bounded cohomological 
displacement, then $F(M)$ belongs to $\dcat{D}_{\cat{N}_0}(\cat{N})$
for every $M \in \dcat{D}_{\mrm{f}}(A)$.
\end{enumerate} 
\end{thm}

\begin{proof} \mbox{}

\smallskip \noindent
(1) First assume $P$ is a bounded complex of finitely generated free 
$A$-modules. The complex $P$ is obtained from the free module $A$ by finitely 
many extensions and translations in
$\dcat{C}_{\mrm{str}}(A)^{\mrm{op}}$. These come either from stupid 
truncations as in (\ref{eqn:2147}) or from breaking up finite direct sums. 
According to Proposition \ref{prop:3060}, $P$ is obtained from $A$ 
by finitely many standard cones and translations in 
$\dcat{D}^{\star}_{\mrm{f}}(A)^{\mrm{op}}$. 
Since $\dcat{D}_{\cat{N}_0}(\cat{N})$ is a full triangulated 
subcategory and $F$ is a triangulated functor, it follows that 
$F(P) \in \dcat{D}_{\cat{N}_0}(\cat{N})$.

Next let $P$ be a bounded above complex of finitely generated free $A$-modules.
Choose some integer $j$. We want to prove that $\opn{H}^j(F(P)) \in \cat{N}_0$. 
Let $i_0$ be an integer such that the integer interval $[i_0, \infty]$ 
contains the cohomological displacement of $F$. 
Define $P' := \opn{stt}^{\leq -j - 1 + i_0}(P)$, the stupid truncation of 
$P$ below $-j - 1 + i_0$; and let $P'' := \opn{stt}^{\geq j + i_0}(P)$, the 
complementary  stupid truncation. 
These truncations are done in the category 
$\dcat{C}_{\mrm{str}}(A)$. 
The short exact sequence (\ref{eqn:2147}) gives, upon applying $\opn{Op}$, 
a short exact sequence in 
$\dcat{C}_{\mrm{str}}(A)^{\mrm{op}}$.
According to Proposition \ref{prop:3060} 
there is a distinguished triangle 
$P' \to P \to P'' \xar{\triangle}$ 
in $\dcat{D}^-_{\mrm{f}}(A)^{\mrm{op}}$.
Since $F$ is a triangulated functor, there is a distinguished triangle 
$F(P') \to F(P) \to F(P'') \xar{\triangle}$ 
in $\dcat{D}(\cat{N})$.
The complex $P''$ is a bounded complex of finitely generated 
free $A$-modules, so we already know that
$F(P'') \in \dcat{D}_{\cat{N}_0}(\cat{N})$,
and in particular $\opn{H}^j(F(P'')) \in \cat{N}_0$. 
On the other hand $\opn{H}(P')$ is concentrated in the interval 
$[-\infty, -j - 1 + i_0]$. Therefore 
$\opn{H}(F(P'))$ is  concentrated in the interval
$[j + 1, \infty]$, and in particular 
$\opn{H}^{j - 1}(F(P')) = \opn{H}^j(F(P')) = 0$.
As we saw in the proof of Lemma \ref{lem:2145}(2), 
$\opn{H}^j(F(P'')) \to \opn{H}^j(F(P))$ is an isomorphism. 
The conclusion is that $\opn{H}^j(F(P)) \in \cat{N}_0$. 

Now take an arbitrary $M \in \dcat{D}^-_{\mrm{f}}(A)$.
There is a quasi-isomorphism $P \to M$, where $P$ is a bounded above complex of 
finitely generated free $A$-modules. So $F(M) \cong F(P)$, and thus 
$F(M) \in \dcat{D}_{\cat{N}_0}(\cat{N})$. 
 
\medskip \noindent
(2) Now we assume that the functor $F$ has finite cohomological dimension. Take 
any complex $M \in \dcat{D}_{\mrm{f}}(A)$. We want to prove 
that for every $j \in \Z$ the object $\opn{H}^j(F(M))$ lies in $\cat{N}_0$. 

Let $[i_0, i_1]$ be a bounded integer interval that contains the 
cohomological displacement of the functor $F$. 
Define $M'' := \opn{smt}^{\leq - j + 1 + i_1}(M)$, the smart truncation of $M$ 
below $-j + 1 + i_1$; and let $M' := \opn{smt}^{\geq -j + 2 + i_1}(M)$, the 
complementary smart truncation. These truncations are preformed in the category 
$\dcat{C}_{\mrm{str}}(A)$. 
By Proposition \ref{prop:3061} there is a distinguished triangle 
$M' \to M \to M'' \xar{\triangle}$ 
in $\dcat{D}_{\mrm{f}}(A)^{\mrm{op}}$.
Since $F$ is a triangulated functor, there is a distinguished triangle 
$F(M') \to F(M) \to F(M'') \xar{\triangle}$ 
in $\dcat{D}(\cat{N})$.
The cohomology of $M'$ is concentrated in the integer
interval $[- j + 2 + i_1, \infty]$, and therefore the cohomology of 
$F(M')$ is concentrated in the interval $[- \infty, j - 2]$. 
In particular the objects $\opn{H}^{j - 1}(F(M'))$ and 
$\opn{H}^{j}(F(M'))$ are zero.
By the proof of Lemma \ref{lem:2145}(2), the morphism 
$\opn{H}^j(F(M'')) \to \opn{H}^j(F(M))$ is an isomorphism.
But $M'' \in \dcat{D}^-_{\mrm{f}}(A)$, so as we proved in part (1), its 
cohomologies are inside $\cat{N}_0$. 
\end{proof}

Here is the covariant version of Theorem \ref{thm:2160}.

\begin{thm} \label{thm:4600} 
Let $A$ be a left noetherian ring, let $\cat{N}$ be an abelian category, let
$\cat{N}_0 \sub \cat{N}$ be a thick abelian subcategory, let $\star$ be a 
boundedness indicator, and let 
$F : \dcat{D}^{\star}_{\mrm{f}}(A) \to \dcat{D}(\cat{N})$
be a triangulated functor. Assume that 
$F(A)$ belongs to $\dcat{D}_{\cat{N}_0}(\cat{N})$.
\begin{enumerate}
\item If $\star = -$, and if $F$ has bounded above cohomological displacement, 
then $F(M)$ belongs to $\dcat{D}_{\cat{N}_0}(\cat{N})$
for every $M \in \dcat{D}^-_{\mrm{f}}(A)$.

\item If $\star = \bra{\mrm{empty}}$, and if $F$ has bounded cohomological 
displacement, then $F(M)$ belongs to $\dcat{D}_{\cat{N}_0}(\cat{N})$
for every $M \in \dcat{D}_{\mrm{f}}(A)$.
\end{enumerate} 
\end{thm}

\begin{exer} \label{exer:4600}
Prove Theorem \ref{thm:4600}.  (Hint: study the proofs of Theorems 
\ref{thm:2135} and \ref{thm:2160}.)
\end{exer}

\begin{rem} \label{rem:4600}
In Theorems \ref{thm:2135}, \ref{thm:2160} and \ref{thm:4600} the 
source category is $\dcat{D}^{\star}_{\mrm{f}}(A)$, for a left noetherian 
ring $A$. There are variants of these theorems in which the source category
is $\dcat{D}^{\star}_{\mrm{f}}(A)$, for a left pseudo-noetherian DG 
ring $A$; and also with source category 
$\dcat{D}^{\star}_{\cat{M}_0}(\cat{M})$
for suitable abelian categories $\cat{M}_0 \sub \cat{M}$. 
We leave it to the reader to investigate these variations. 
\end{rem}

\mysubsection{Derived Restriction, Induction and Coinduction Functors}
\label{subsec:adjs-NC}

As before, there is a nonzero commutative base ring $\K$. All DG rings are 
$\K$-central, and all homomorphisms between them are over $\K$.

\begin{dfn} \label{dfn:3086}
Let $u : A \to B$ be a DG ring homomorphism.
\begin{enumerate}
\item The {\em restriction functor}%
\index{Restriction functor}
is the $\K$-linear DG functor%
\index{1-Rest@$\opn{Rest}_u$}
$\opn{Rest}_u = \opn{Rest}_{B / A} : \dcat{C}(B) \to \dcat{C}(A)$,
that is the identity on the underlying DG 
$\K$-modules, and the $A$-action is via $u$.

\item Since the functor $\opn{Rest}_u$  is exact, it extends to a 
triangulated functor 
$\opn{Rest}_u : \dcat{D}(B) \to \dcat{D}(A)$.
\end{enumerate}
\end{dfn}

\begin{dfn} \label{dfn:3087}
Let $u : A \to B$ be a DG ring homomorphism.
\begin{enumerate}
\item The {\em induction functor}%
\index{Induction functor}%
\index{1-Ind@$\opn{Ind}_u$}
is the DG functor 
$\opn{Ind}_u  = \opn{Ind}_{B / A} : \dcat{C}(A) \to \dcat{C}(B)$
with formula $\opn{Ind}_u := B \ot_A (-)$.

\item The {\em derived induction functor}%
\index{Induction functor! derived}%
\index{1-LInd@$\opn{LInd}_u$}
$\opn{LInd}_u  = \opn{LInd}_{B / A} : \dcat{D}(A) \to \dcat{D}(B)$ 
is the triangulated left derived functor of $\opn{Ind}_u$, namely 
$\opn{LInd}_u := B \ot^{\mrm{L}}_A (-)$.
\end{enumerate}
\end{dfn}

For every $M \in \dcat{D}(A)$ there is the canonical morphism 
$\eta^{\mrm{L}}_{M} : \opn{LInd}_u(M) \to \opn{Ind}_u(M)$
in $\dcat{D}(B)$, which is part of the left derived functor. 

\begin{dfn} \label{dfn:3111}
Let $u : A \to B$ be a DG ring homomorphism, and let 
$M \in \dcat{D}(A)$ and $N \in \dcat{D}(B)$ be DG modules.
A morphism
$\la : M \to \opn{Rest}_u(N)$
in $\dcat{D}(A)$ is called a 
{\em forward morphism in $\dcat{D}(A)$  over $u$}%
\index{Forward morphism! in derived category}.

Similarly there are forward morphisms over $u$ in the categories 
$\dcat{C}(A)$, $\dcat{C}_{\mrm{str}}(A)$ and $\dcat{K}(A)$.
\end{dfn}

We often omit the functor $\opn{Rest}_u$, and just talk about a 
forward morphism $\la : M \to N$ in $\dcat{D}(A)$ over $u$, etc. 

For any DG module $M \in \dcat{C}(A)$ there is a canonical forward homomorphism
\begin{equation} \label{eqn:3194}
\opn{q}_{u, M} : M \to B \ot_A M = \opn{Ind}_u(M) \, , \quad 
\opn{q}_{u, M}(m) := 1_B \ot m 
\end{equation}
in $\dcat{C}_{\mrm{str}}(A)$ over $u$. 
Now let $N \in \dcat{C}(B)$. The usual change of ring adjunction formula gives 
a 
canonical isomorphism 
\begin{equation} \label{eqn:2401}
\opn{fadj}_{u, M, N} : 
\opn{Hom}_{\dcat{C}(A)}(M, N) \iso 
\opn{Hom}_{\dcat{C}(B)} \bigl( B \ot_A M, N  \bigr) 
\end{equation}
in $\dcat{C}_{\mrm{str}}(\K)$. It is characterized by the property that 
for every forward morphism $\la : M \to N$ there is equality 
$\opn{fadj}_{u, M, N}(\la) \circ \opn{q}_{u, M} = \la$.
We refer to the isomorphism $\opn{fadj}_{u, M, N}$ as 
{\em forward adjunction}. 
Since the iso\-morphism $\opn{fadj}_{u, M, N}$ is functorial in 
$M$ and $N$, we see that the functor $\opn{Ind}_u$ is a left adjoint of 
$\opn{Rest}_u$. 

The next theorem shows that this is also true on the derived level. 

\begin{thm} \label{thm:3085} 
Let $u : A \to B$ be a homomorphism of DG $\K$-rings.
\begin{enumerate}
\item For each $M \in \dcat{D}(A)$ there is a unique forward morphism%
\index{Forward morphism! canonical}
\index{1-QLuM@$\opn{q}_{u, M}^{\mrm{L}}$}
\[ \opn{q}_{u, M}^{\mrm{L}} : M \to B \ot^{\mrm{L}}_{A} M = 
\opn{LInd}_u(M) \]
in $\dcat{D}(A)$ over $u$, called the {\em canonical forward morphism},
which is functorial in $M$, and satisfies 
$\eta^{\mrm{L}}_M \circ \opn{q}_{u, M}^{\mrm{L}} = \opn{q}_{u, M}$.

\item For each $M \in \dcat{D}(A)$ and $N \in \dcat{D}(B)$
there is a unique $\K$-linear isomorphism%
\index{Forward adjunction! derived}%
\index{1-Fadj@$\opn{fadj}^{\mrm{L}}_{u, M, N}$}
\[ \opn{fadj}_{u, M, N}^{\mrm{L}} :  
\opn{Hom}_{\dcat{D}(A)} \bigl( M, \opn{Rest}_u(N) \bigr) \iso 
\opn{Hom}_{\dcat{D}(B)} \bigl( \opn{LInd}_u(M), N \bigr) , \]
called {\em derived forward adjunction}, such that 
$\opn{fadj}_{u, M, N}^{\mrm{L}}(\la) \circ \opn{q}_{u, M}^{\mrm{L}} = \la$
for all forward morphisms $\la : M \to N$ in $\dcat{D}(A)$ over $u$.
The isomorphism $\opn{fadj}_{u, M, N}^{\mrm{L}}$ is functorial in $M$ and $N$. 

\item The functor
$\opn{LInd}_u : \dcat{D}(A) \to \dcat{D}(B)$
is left adjoint to $\opn{Rest}_u$. 
\end{enumerate}
\end{thm}

Here are the commutative diagrams in $\dcat{D}(A)$ illustrating the theorem. 
\[ \UseTips \xymatrix @C=8ex @R=6ex {
M
\ar[d]_{\opn{q}_{u, M}^{\mrm{L}}}
\ar[dr]^{\opn{q}_{u, M}} 
\\
\opn{LInd}_u(M)
\ar[r]_{\eta^{\mrm{L}}_M}
&
*+{\opn{Ind}_u(M)}
} 
\qquad \quad 
\UseTips \xymatrix @C=14ex @R=6ex {
M
\ar[d]_{\opn{q}_{u, M}^{\mrm{L}}}
\ar[dr]^{\la} 
\\
\opn{LInd}_u(M)
\ar[r]_(0.55){\opn{fadj}_{u, M, N}^{\mrm{L}}(\la)}
&
*+++{N}
} 
\]

\begin{proof} \mbox{}

\smallskip \noindent 
(1) Let us choose a system of K-projective resolutions in $\dcat{C}(A)$, as in 
Definition \ref{dfn:1524}, and use it to present the derived functor 
$\opn{LInd}_u$. 
To construct an isomorphism $\opn{q}_{u, M}^{\mrm{L}}$
which is functorial in $M$, it suffices to consider 
a K-projective DG $A$-module $M = P$. But then $\opn{LInd}_u(P) = 
\opn{Ind}_u(P)$, 
$\eta^{\mrm{L}}_P = \opn{id}_P$,
and we have no choice but to define 
$\opn{q}_{u, P}^{\mrm{L}} := \opn{q}_{u, P}$. 

\medskip \noindent
(2) We keep the system of K-projective resolutions. 
As in item (1), it is enough to say what is $\opn{fadj}_{u, P, N}^{\mrm{L}}$ 
for a K-projective DG $A$-module $P$. Then there is a canonical isomorphism 
 \[ \opn{Q}_A : \opn{Hom}_{\dcat{K}(A)}(P, N) \iso 
\opn{Hom}_{\dcat{D}(A)}(P, N) . \]
Since $B \ot_A P = \opn{LInd}_u(P)$ is a K-projective DG $B$-module, we also 
have  a canonical isomorphism 
\[ \opn{Q}_B : \opn{Hom}_{\dcat{K}(B)}(B \ot_A P, N) \iso
\opn{Hom}_{\dcat{D}(B)}(B \ot_A P , N) . \]
But the usual adjunction formula (\ref{eqn:2401}) induces, by passing to 
$0$-th cohomology on both sides, an isomorphism 
\[ \opn{H}^0(\opn{fadj}_{u, P, N})  : \opn{Hom}_{\dcat{K}(A)}(P, N) \iso 
\opn{Hom}_{\dcat{K}(B)}(B \ot_A P, N) ,  \]
and it satisfies 
$\opn{H}^0(\opn{fadj}_{u, P, N}) (\la) \circ \opn{q}_{u, P} = \la$
for any forward morphism $\la : P \to N$ in $\dcat{K}(A)$. 
We define 
\[ \opn{fadj}_{u, P, N}^{\mrm{L}} := \opn{Q}_B \circ 
\opn{H}^0(\opn{fadj}_{u, P, N}) \circ \opn{Q}_A^{-1} . \]

\medskip \noindent
(3) Since $\opn{fadj}_{u, M, N}^{\mrm{L}}$ is a bifunctorial isomorphism, 
it determines an adjunction. 
\end{proof}

\begin{dfn} \label{dfn:3191}
Let $u : A \to B$ be a DG ring homomorphism, and let 
$M \in \dcat{D}(A)$ and $N \in \dcat{D}(B)$ be DG modules.
A forward morphism $\la : M \to N$ in $\dcat{D}(A)$ over $u$ is called a 
{\em nondegenerate forward morphism}%
\index{Forward morphism! in derived category}%
\index{Forward morphism! nondegenerate}%
\index{1-Fadj@$\opn{fadj}^{\mrm{L}}_{u, M, N}$}
if the corresponding morphism 
\[ \opn{fadj}_{u, M, N}^{\mrm{L}}(\la) : B \ot^{\mrm{L}}_A M = 
\opn{LInd}_u(M) \to N \] 
in $\dcat{D}(B)$ is an isomorphism.
\end{dfn}

\begin{exa} \label{exa:2404}
Given $M \in \dcat{D}(A)$, let $N := B \ot^{\mrm{L}}_A M \in \dcat{D}(B)$. 
The canonical derived forward morphism%
\index{1-Fadj@$\opn{fadj}^{\mrm{L}}_{u, M, N}$}%
\index{1-QLuM@$\opn{q}_{u, M}^{\mrm{L}}$} 
$\opn{q}_{u, M}^{\mrm{L}} : M \to N$ is a nondegenerate forward morphism
in $\dcat{D}(A)$ over $u$, because 
$\opn{fadj}_{u, M, N}^{\mrm{L}}(\opn{q}_{u, M}^{\mrm{L}}) = \opn{id}_N$.
\end{exa}

\begin{prop} \label{prop:3191}
Suppose $A \xar{u} B \xar{v} C$ are homomorphisms of DG rings,
$L \in \dcat{D}(A)$, $M \in \dcat{D}(B)$ and $N \in \dcat{D}(C)$ are DG 
modules, $\la : L \to M$ is a nondegenerate forward morphism
in $\dcat{D}(A)$ over $u$, and $\mu : M \to N$ is a nondegenerate forward 
morphism in $\dcat{D}(B)$ over $v$. Then 
$\mu \circ \la : L \to  N$
is a nondegenerate forward morphism in $\dcat{D}(A)$ over $v \circ u$.
\end{prop}

\begin{proof}
Since $\la$ is a nondegenerate forward morphism, it induces an isomorphism 
$N \cong B \ot^{\mrm{L}}_A M$. We can thus assume that  
$N = B \ot^{\mrm{L}}_A M$ and that 
$\la = \opn{q}_{u, M}^{\mrm{L}} : M \to B \ot^{\mrm{L}}_A M$.
Similarly, we can assume that 
$\mu = \opn{q}_{v, N}^{\mrm{L}} : N \to C \ot^{\mrm{L}}_B N$.
But now 
$\mu \circ \la  = \opn{q}_{v, N}^{\mrm{L}} \circ \opn{q}_{u, M}^{\mrm{L}} = 
\opn{q}_{v \circ u, L}^{\mrm{L}} : L \to N$,
and this is known to be a nondegenerate forward morphism (see Example 
\ref{exa:2404}). 
\end{proof}

\begin{dfn} \label{dfn:3095}
Let $u : A \to B$ be a homomorphism of DG $\K$-rings.
\begin{enumerate}
\item The {\em coinduction functor}%
\index{Coinduction functor}%
\index{1-CInd@$\opn{CInd}_u$}
is the $\K$-linear DG functor 
$\opn{CInd}_u = \opn{CInd}_{B / A} : \dcat{C}(A) \to \dcat{C}(B)$ 
with formula
$\opn{CInd}_u := \opn{Hom}_{A}(B, -)$.

\item The  {\em derived coinduction functor}%
\index{Coinduction functor! derived}%
\index{1-RCInd@$\opn{RCInd}_u$}
$\opn{RCInd}_u = \opn{RCInd}_{B / A} : \dcat{D}(A) \to \dcat{D}(B)$  
is the triangulated right derived functor of $\opn{CInd}_u$.
\end{enumerate}
\end{dfn}

Thus for each $M \in \dcat{D}(A)$ we have 
$\opn{RCInd}_u(M) = \opn{RHom}_{A}(B, M)$.
There is the canonical morphism 
$\eta^{\mrm{R}}_{M} : \opn{CInd}_u(M) \to \opn{RCInd}_u(M)$
in $\dcat{D}(B)$, which is part of the right derived functor,

\begin{dfn} \label{dfn:3110}
Let $u : A \to B$ be a DG ring homomorphism, and let 
$M \in \dcat{D}(A)$ and $N \in \dcat{D}(B)$ be DG modules.
A morphism
$\th : \opn{Rest}_u(N) \to M$
in $\dcat{D}(A)$ is called a {\em backward} (or {\em trace}) {\em morphism in 
$\dcat{D}(A)$ over $u$}%
\index{Backward morphism! in derived category}%
\index{Trace morphism! in derived category}.

Similarly there are backward morphisms over $u$ in the categories 
$\dcat{C}(A)$, $\dcat{C}_{\mrm{str}}(A)$ and $\dcat{K}(A)$.
\end{dfn}

We shall use the terms ``backward morphism'' and ``trace morphism'' 
synonymously. 
We often omit the functor $\opn{Rest}_u$, and just talk about a 
backward morphism $\th : N \to M$ in $\dcat{D}(A)$  over $u$. 

For any DG module $M \in \dcat{C}(A)$ there is a canonical backward homomorphism
$\opn{tr}_{u, M} : \opn{CInd}_u(M) = \opn{Hom}_A(B, M) \to M$,
$\opn{tr}_{u, M}(\phi) := \phi(1_B)$ 
in $\dcat{C}_{\mrm{str}}(A)$ over $u$. 
Now let $N \in \dcat{C}(B)$. The usual change of ring adjunction formula gives 
a canonical isomorphism 
\[ \opn{badj}_{u, M, N} : 
\opn{Hom}_{\dcat{C}(A)}(N, M) \iso 
\opn{Hom}_{\dcat{C}(B)} \bigl( N, \opn{Hom}_A(B, M) \bigr) \]
in $\dcat{C}_{\mrm{str}}(\K)$. It is characterized by the property that 
for every backward morphism $\th : N \to M$ there is equality 
$\opn{tr}_{u, M} \circ \opn{badj}_{u, M, N}(\th) = \th$.
We refer to the isomorphism  $\opn{badj}_{u, M, N}$ as 
{\em backward adjunction}. 
Since the isomorphism $\opn{badj}_{u, M, N}$ is 
functorial in $M$ and $N$, we see that the functor $\opn{CInd}_u$ is a right 
adjoint of $\opn{Rest}_u$. 

The next theorem shows that this is also true on the derived level. 

\begin{thm} \label{thm:3192} 
Let $u : A \to B$ be a homomorphism of DG $\K$-rings.
\begin{enumerate}
\item For each $M \in \dcat{D}(A)$ there is a unique backward morphism 
\[ \opn{tr}_{u, M}^{\mrm{R}} : \opn{RCInd}_u(M) = 
\opn{RHom}_{A}(B, M) \to  M \]
in $\dcat{D}(A)$ over $u$, called the {\em canonical backward morphism}%
\index{Backward morphism! canonical}%
\index{1-TruM@$\opn{tr}_{u, M}^{\mrm{R}}$},
which is functorial in $M$, and satisfies 
$\opn{tr}_{u, M}^{\mrm{R}}  \circ \, \eta^{\mrm{R}}_M  = \opn{tr}_{u, M}$.

\item For each $M \in \dcat{D}(A)$ and $N \in \dcat{D}(B)$
there is a unique $\K$-linear isomorphism 
\[ \opn{badj}_{u, M, N}^{\mrm{R}} :  
\opn{Hom}_{\dcat{D}(A)} \bigl( \opn{Rest}_u(N), M \bigr) \iso 
\opn{Hom}_{\dcat{D}(B)} \bigl( N, \opn{RCInd}_u(M) \bigr) , \]
called {\em derived backward adjunction}%
\index{Backward adjunction! derived}%
\index{1-Badj@$\opn{badj}^{\mrm{R}}_{u, M, N}$},
such that 
$\opn{tr}_{u, M}^{\mrm{R}} \circ \opn{badj}_{u, M, N}^{\mrm{R}}(\th) = \th$
for all backward morphisms $\th : N \to M$ in $\dcat{D}(A)$ over $u$.
The isomorphism $\opn{badj}_{u, M, N}^{\mrm{R}}$ is functorial in $M$ and $N$. 

\item The functor
$\opn{RCInd}_u : \dcat{D}(A) \to \dcat{D}(B)$
is right adjoint to $\opn{Rest}_u$. 
\end{enumerate}
\end{thm}

Here are the commutative diagrams in $\dcat{D}(A)$ illustrating the theorem. 
\[ \UseTips \xymatrix @C=8ex @R=6ex {
\opn{CInd}_u(M)
\ar[r]^{\eta^{\mrm{R}}_M}
\ar[dr]_{\opn{tr}_{u, M}} 
&
\opn{RCInd}_u(M)
\ar[d]^{\opn{tr}_{u, M}^{\mrm{R}}}
\\
&
M
} 
\qquad \quad 
\UseTips \xymatrix @C=14ex @R=6ex {
N
\ar[r]^(0.4){\opn{badj}_{u, M, N}^{\mrm{R}}(\th)}
\ar[dr]_{\th} 
& 
\opn{RCInd}_u(M)
\ar[d]^{\opn{tr}_{u, M}^{\mrm{R}}}
\\
&
M
} 
\]

\begin{proof}
The proof is very similar to that of Theorem \ref{thm:3085}. The only change is 
that now we choose a system of K-injective resolutions in $\dcat{C}(A)$. 
For a K-injective DG $A$-module $I$, the DG $B$-module
$\opn{Hom}_A(B, I)$ is K-injective, and it represents
$\opn{RCInd}_u(M)$. 
\end{proof}

\begin{dfn} \label{dfn:3190} 
Let $u : A \to B$ be a DG ring homomorphism, and let 
$M \in \dcat{D}(A)$ and $N \in \dcat{D}(B)$ be DG modules.
A backward morphism $\th : N \to M$ in $\dcat{D}(A)$ over $u$ is called  
{\em a nondegenerate backward morphism}%
\index{Backward morphism! nondegenerate}%
\index{Backward morphism! in derived category}%
\index{1-Badj@$\opn{badj}^{\mrm{R}}_{u, M, N}$}
if the corresponding morphism
\[ \opn{badj}_{u, M, N}^{\mrm{R}}(\th) : N \to \opn{RHom}_A(B, M) 
= \opn{RCInd}_u(M) \]
in $\dcat{D}(B)$ is an isomorphism.
\end{dfn}

\begin{exa} \label{exa:2402}
Given $M \in \dcat{D}(A)$, let 
$N := \opn{RHom}_{A}(B, M) \in \dcat{D}(B)$. 
The canonical backward morphism 
$\opn{tr}_{u, M}^{\mrm{R}} : N \to M$%
\index{1-Badj@$\opn{badj}^{\mrm{R}}_{u, M, N}$}%
\index{1-TruM@$\opn{tr}_{u, M}^{\mrm{R}}$} 
is a nondegenerate backward morphism
in $\dcat{D}(A)$ over $u$, because 
$\opn{badj}_{u, M, N}^{\mrm{R}}(\opn{tr}_{u, M}^{\mrm{R}}) = \opn{id}_N$.
\end{exa}

\begin{prop} \label{prop:3190}
Suppose $A \xar{u} B \xar{v} C$ are homomorphisms of DG rings,
$L \in \dcat{D}(A)$, $M \in \dcat{D}(B)$ and $N \in \dcat{D}(C)$ are DG 
modules, $\th : M \to L$ is a nondegenerate backward morphism
in $\dcat{D}(A)$ over $u$, and $\ze : N \to M$ is a nondegenerate backward 
morphism in $\dcat{D}(B)$ over $v$. Then 
$\th \circ \ze : N \to  L$
is a nondegenerate backward morphism in $\dcat{D}(A)$ over $v \circ u$.
\end{prop}

\begin{proof}
Since $\th$ is a nondegenerate backward morphism, it induces an isomorphism 
$M \cong \opn{RHom}_{A}(B, L)$. We can thus assume that  
$M = \opn{RHom}_{A}(B, L)$ and that 
$\th = \opn{tr}_{u, L}^{\mrm{R}}$.
Similarly, we can assume that 
$N = \opn{RHom}_{B}(C, M)$ and that 
$\ze = \opn{tr}_{v, M}^{\mrm{R}}$.
But now 
$\th \circ \ze  =
\opn{tr}_{u, L}^{\mrm{R}}  \circ \opn{tr}_{v, M}^{\mrm{R}} = 
\opn{tr}_{v \circ u, L}^{\mrm{R}} : N \to L$,
and this is known to be a nondegenerate backward morphism (see Example 
\ref{exa:2402}). 
\end{proof}

\begin{exa} \label{exa:2403}
If $A = B$ and $u = \opn{id}_A$, then 
backward and forward morphisms over $u$ are just morphisms in $\dcat{D}(A)$. 
Nondegenerate (backward or forward) morphisms are just isomorphisms in 
$\dcat{D}(A)$. 
\end{exa}

\mysubsection{DG Ring Quasi-Isomorphisms}
\label{subsec:DGR-quisoms}

Suppose $A$ and $B$ are DG $\K$-rings. A homomorphism of DG $\K$-rings 
$u : A \to B$ induces a homomorphism of graded $\K$-rings 
$\opn{H}(u) : \opn{H}(A) \to \opn{H}(B)$; cf.\ Example \ref{exer:1100}.

\begin{dfn} \label{dfn:3085}
A DG ring homomorphism $u : A \to B$ is called a 
{\em quasi-isomorphism of DG rings}%
\index{Differential graded ring! quasi-isomorphism of {\indash}s}
if $\opn{H}(u)$ is an isomorphism. 
\end{dfn}

Here is a fundamental result, due to B. Keller in 1994 \cite{Kel1}, 
and independently to V. Hinich in 1997 \cite{Hi1}. 
It is the justification behind the use of DG ring resolutions
(see Remarks \ref{rem:2285}, \ref{rem:1420} and \ref{rem:3691}). 

\begin{thm} \label{thm:2363}
Let $u : A \to B$ be a quasi-isomorphism of DG $\K$-rings%
\index{Differential graded ring! quasi-isomorphism of {\indash}s}.
Then\tup{:}
\begin{enumerate}
\item The restriction functor 
$\opn{Rest}_u : \dcat{D}(B) \to \dcat{D}(A)$
is an equivalence of $\K$-linear triangulated categories, with quasi-inverse
$\opn{LInd}_u$. 

\item For every $M \in \dcat{D}(B^{\mrm{op}})$ and 
$N \in \dcat{D}(B)$ there is an  isomorphism
$M \ot^{\mrm{L}}_A N \iso M \ot^{\mrm{L}}_B N$
in $\dcat{D}(\K)$. This isomorphism is functorial in $M$ and $N$.

\item For every $M, N \in \dcat{D}(B)$ there is an isomorphism 
$\opn{RHom}_{A}(M, N) \iso \lb \opn{RHom}_{B}(M, N)$
in $\dcat{D}(\K)$. This isomorphism is functorial in $M$ and $N$.
\end{enumerate}
\end{thm}

Notice that the restriction functor $\opn{Rest}_u$ is suppressed in parts (2) 
and (3) of the theorem. The proof hinges on the next lemma. 

\begin{lem} \label{lem:4770}
Let $u : A \to B$ be a quasi-isomorphism of DG $\K$-rings, let 
$N \in \dcat{D}(B)$, and let $\rho : P \to N$ be a K-projective resolution in 
$\dcat{C}_{\mrm{str}}(A)$. 
\begin{enumerate}
\item The homomorphism
$\opn{q}_{u, P} : P \to B \ot_A P$
in $\dcat{C}_{\mrm{str}}(A)$ from equation \lb \tup{(\ref{eqn:3194})}, and the 
homomorphism 
$\psi : B \ot_A P \to N$, $\psi(b \ot p) := b \cd \rho(p)$,  
in $\dcat{C}_{\mrm{str}}(B)$, are both quasi-isomorphisms.

\item $B \ot_A P$ is a K-projective DG $B$-module.
\end{enumerate}
\end{lem}

\begin{proof} \mbox{}

\smallskip \noindent
(1) Look at the commutative diagram 
\[ \UseTips \xymatrix @C=8ex @R=6ex {
&
B \ot_A P 
\ar[r]^(0.56){\psi}
&
N
\\
P
\ar[r]^(0.4){\cong}
\ar[ur]^{\opn{q}_{u, P}}
&
A \ot_A P
\ar[u]_{u \, \ot \, \opn{id}_P}
\ar[r]^{\opn{id}_A \ot \, \rho}
&
A \ot_A N 
\ar[u]_{\cong}
} \]
in $\dcat{C}_{\mrm{str}}(A)$, in which the arrows marked $\cong$ are the 
canonical isomorphisms. 
The homomorphism $u \ot \opn{id}_P$ is a quasi-isomorphism, because $u$ is a 
quasi-isomorphism and $P$ is K-flat over $A$. 
The homomorphism $\opn{id}_A \ot \, \rho$ is a quasi-isomorphism, because 
$\rho$ is a quasi-isomorphism and $A$ is K-flat over itself. Therefore
$\opn{q}_{u, P}$ and $\psi$ are quasi-isomorphisms.

\medskip \noindent 
(2) Take an acyclic DG $B$-module $N$. There is an isomorphism 
\[ \opn{Hom}_B(B \ot_A P, N) \cong \opn{Hom}_A(P, N) \]
in $\dcat{C}_{\mrm{str}}(\K)$, coming from adjunction for the DG ring 
homomorphism $u$. By assumption the DG module $\opn{Hom}_A(P, N)$ is acyclic. 
This shows that $B \ot_A P$ is K-projective over $B$. 
\end{proof}

\begin{proof}[Proof of Theorem \tup{\ref{thm:2363}}] \mbox{}

\smallskip \noindent
(1) Take any $N \in \dcat{D}(B)$. 
Choose a K-projective resolution $\rho : P \to N$ in 
$\dcat{C}_{\mrm{str}}(A)$, so that 
$(\opn{LInd}_u \circ \opn{Rest}_u)(N) \cong B \ot_A P$
in $\dcat{D}(B)$. By the lemma we have a quasi-isomorphism 
$\psi : B \ot_A P \to N$ in $\dcat{C}_{\mrm{str}}(B)$.
This means that we have an isomorphism
\[ \opn{Q}(\psi) : (\opn{LInd}_u \circ \opn{Rest}_u)(N) \iso N \]
in $\dcat{D}(B)$, and it is functorial in $N$. 

On the other hand, starting from a DG module $M \in \dcat{D}(A)$, 
and choosing a K-projective resolution $\rho : P \to M$
in $\dcat{C}_{\mrm{str}}(A)$, we have an isomorphism 
$\opn{LInd}_u(M) \cong B \ot_A P$
in $\dcat{D}(B)$. So there is an isomorphism 
$(\opn{Rest}_u \circ \opn{LInd}_u)(M) \cong \lb B \ot_A P$
in $\dcat{D}(A)$. According to the lemma there is a quasi-isomorphism 
$\opn{q}_{u, P} : P \to B \ot_A P$
in $\dcat{C}_{\mrm{str}}(A)$.
We obtain an isomorphism 
\[ \opn{Q}(\rho) \circ \opn{Q}(\opn{q}_{u, P})^{-1} : 
(\opn{Rest}_u \circ \opn{LInd}_u)(M) \iso M \]
in $\dcat{D}(A)$, and it is functorial in $M$. 

\medskip \noindent
(2) Choose a K-projective resolution $\rho : P \to M$ in 
$\dcat{C}(A^{\mrm{op}})$. 
This produces an isomorphism 
$\phi_1 : P \ot_A N \iso M \ot^{\mrm{L}}_A N$
in $\dcat{D}(\K)$. The DG module 
$P \ot_A B \in \dcat{C}(B^{\mrm{op}})$
is K-projective, and there is a quasi-isomorphism
$\psi : P \ot_A B \to M$ in 
$\dcat{C}_{\mrm{str}}(B^{\mrm{op}})$,
by the lemma (with $B^{\mrm{op}}$ instead of $B$). 
We get an isomorphism 
\[ \phi_2 := \opn{Q}(\psi) \ot^{\mrm{L}}_B \opn{id}_N : 
(P \ot_A B) \ot_B N \iso M \ot^{\mrm{L}}_B N \]
in $\dcat{D}(\K)$.
There is also the canonical isomorphism 
\[ \phi_3 : (P \ot_A B) \ot_B N \iso P \ot_A N \]
in  $\dcat{D}(\K)$. 
The functorial isomorphism we want is 
\[ \phi_2 \circ \phi_3^{-1} \circ \phi_1^{-1} : 
M \ot^{\mrm{L}}_A N \iso M \ot^{\mrm{L}}_B N \]
in $\dcat{D}(\K)$. 

\medskip \noindent
(3) Let us choose a K-projective resolution 
$\rho : P \to M$ in $\dcat{C}(A)$. 
This produces an isomorphism 
\[ \phi_1 : \opn{RHom}_{A}(M, N) \iso \opn{Hom}_{A}(P, N)  \]
in $\dcat{D}(\K)$. The DG module 
$B \ot_A P \in \dcat{C}(B)$ is K-projective, and the homomorphism 
$\psi : B \ot_A P \to M$ 
in $\dcat{C}_{\mrm{str}}(B)$ is a quasi-isomorphism (again using the lemma). 
In this way we have an isomorphism 
\[ \phi_2 : \opn{Hom}_{B}(B \ot_A P, N) \iso \opn{RHom}_{B}(M, N) \]
in $\dcat{D}(\K)$. And there is the canonical isomorphism of adjunction 
\[ \phi_3 : \opn{Hom}_{B}(B \ot_A P, N) \iso \opn{Hom}_{A}(P, N) \]
in $\dcat{D}(\K)$.
The functorial isomorphism we want is
\[ \phi_2 \circ \phi_3^{-1} \circ \phi_1 :
\opn{RHom}_{A}(M, N) \iso \opn{RHom}_{B}(M, N) . \qedhere \]
\end{proof}

\mysubsection{Existence of DG Ring Resolutions}
\label{subsec:DG-ring-res}

In Theorem \ref{thm:2363} we saw that a DG ring quasi-isomorphism 
$A \to B$ does not change the derived categories: the restriction functor 
$\dcat{D}(B) \to \dcat{D}(A)$ is an equivalence of triangulated categories. 

Sometimes it is advantageous to replace a given DG ring by a quasi-isomorphic 
one (see Remarks \ref{rem:4605}, \ref{rem:4945}, \ref{rem:2285}, 
\ref{rem:1420} and \ref{rem:3691}). 
In this subsection we study {\em semi-free DG rings}, and 
{\em semi-free DG ring resolutions}. 

As before, we fix a nonzero commutative base ring $\K$.
All DG rings are $\K$-central, namely we work 
in the category $\catt{DGRng} \centover \K$. We use the shorthand 
``NC'' for ``noncommutative''. 

Graded sets and filtered graded sets were introduced in Definitions 
\ref{dfn:4730} and \ref{dfn:4731} respectively. 
The next definition was mentioned briefly in Example \ref{exa:4625}.

\begin{dfn} \label{dfn:4306}
Let $X$ be a graded set. The {\em NC polynomial ring on $X$} is the graded 
$\K$-ring $\K \bra{X}$. 
\end{dfn}

Perhaps we should say a few words on the structure of the ring $\K \bra{X}$. 
Given $n \geq 0$ and elements 
$x_{1}, \ldots, x_{n} \in X$, their product
$x_{1} \cdots x_{n} \in \K \bra{X}$
is called a monomial (or a word), and its degree is 
\begin{equation} \label{eqn:4305}
\opn{deg}(x_{1} \cdots x_{n}) := \opn{deg}(x_{1}) + \cdots + 
\opn{deg}(x_{n}) \in \Z .
\end{equation}
For $n = 0$ the unique monomial is denoted by $1$. 
As a graded $\K$-module, $\K \bra{X}$ is free, with basis the set of 
monomials. Multiplication in $\K \bra{X}$ is by concatenation of monomials, 
extended $\K$-bilinearly. 

\begin{dfn} \label{dfn:4308}
Let $(X, F)$ be a filtered graded set. There is an induced filtration
$F = \bigl\{ F_{j} \bigl( \K \bra{X} \bigr) \bigr\}_{j \geq -1}$ on 
$\K \bra{X}$ by $\K$-submodules, as follows:
\[ F_j \bigl( \K \bra{X} \bigr) := 
\begin{cases}
0 & \tup{if} \quad j = -1 
\\
\K \bra{F_j(X)} & \tup{if} \quad j \geq 0  .   
\end{cases} \]
\end{dfn}

For every $j \geq 0$ the filtered piece 
$F_j \bigl( \K \bra{X} \bigr)$ is a graded subring of $\K \bra{X}$. And of 
course 
$\K \bra{X} = \bigcup_{j} F_j \bigl( \K \bra{X} \bigr)$.

Recall that given a DG ring 
$A$, we denote by $A^{\natural}$ the graded ring gotten by forgetting the 
differential $\d_A$. Here is the main definition of this subsection. 

\begin{dfn} \label{dfn:4309}
A {\em NC semi-free DG $\K$-ring}%
\index{Differential graded ring! NC semi-free}
is a DG $\K$-ring $A$ that 
admits an isomorphism of graded $\K$-rings
$A^{\natural} \cong \K \bra{X}$,
for some filtered graded set $(X, F)$, such that under this isomorphism there 
is an inclusion 
$\d_A \bigl( F_j(X) \bigr) \sub F_{j - 1} \bigl( \K \bra{X} \bigr)$ 
for every $j \geq 0$. Such a filtered graded set $(X, F)$ is called a 
{\em multiplicative semi-basis}%
\index{Semi-basis! multiplicative {\indash} of DG ring}
of $A$. 
\end{dfn}

In this book we are mostly interested in semi-free DG rings because of the 
proposition below. 

\begin{prop} \label{prop:4305}
If $A$ is a NC semi-free DG $\K$-ring, then as a DG $\K$-module $A$ is 
semi-free, and hence K-flat. 
\end{prop}

\begin{proof}
Suppose we are given an isomorphism 
$A^{\natural} \cong \K \bra{X}$
for some filtered graded set $(X, F)$, as in Definition \ref{dfn:4309}.
Let $Y$ be the set of monomials in the elements of $X$. The set $Y$ is 
graded, see (\ref{eqn:4305}), and it forms a basis of $A^{\natural}$ as a 
graded $\K$-module. However, since $X$ could have elements of positive degree, 
this is not enough to make $A$ into a semi-free DG $\K$-module; we need to put 
a suitable filtration on $Y$. Cf.\ Proposition \ref{prop:4730}.

For an element $x \in X$ let 
$\opn{ord}^F(x) := \opn{min} \{ j \mid x \in F_j(X) \} \in \N$.
Thus 
$F_j(X) = \{ x \in X \mid \opn{ord}^F(x) \leq j \}$.
We extend this order function to monomials by 
$\opn{ord}^F(x_{1} \cdots x_{n}) := \sum_{i} \opn{ord}^F(x_{i})$.
The element $1$ (the monomial of length $0$) gets order $0$. 
Then we define 
$G_j(Y) := \{ y \in Y \mid \opn{ord}^F(y) \leq j \}$.
This is a filtration 
$G = \{ G_{j}(Y) \}_{j \geq -1}$ of the graded set $Y$, and it induces a 
filtration $G = \{ G_{j}(A) \}_{j \geq -1}$ on $A$ 
by $\K$-submodules, where $G_{j}(A) \sub A$ is defined to be the $\K$-linear 
span of $G_{j}(Y)$. 

The graded Leibniz rule, together with the inclusion 
$\d_A (F_j(X)) \sub \lb F_{j - 1}(\K \bra{X})$
from Definition \ref{dfn:4309}, imply that 
$\d_A (G_{j}(A)) \sub G_{j - 1}(A)$
for every $j \geq 0$. And 
$\opn{Gr}^G_j(A)$ is a free DG $\K$-module, with basis 
$\opn{Gr}_j^G(Y) := \{ y \in Y \mid \opn{ord}^F(y) = j \}$.
Thus $G$ is a semi-free filtration on the DG $\K$-module $A$, and $A$ is a 
semi-free DG $\K$-module; see Definition \ref{dfn:1575}.

Finally, according to Theorem \ref{thm:1575} and Proposition \ref{prop:1525}, 
$A$ is a K-flat DG $\K$-module. 
\end{proof}

\begin{dfn} \label{dfn:4310}
Let $A$ be a DG $\K$-ring. A {\em NC semi-free DG ring resolution} of 
$A$ relative to $\K$ is a quasi-isomorphism $u : \til{A} \to A$ 
of DG $\K$-rings, where $\til{A}$ is a NC semi-free DG $\K$-ring.
We also call $u : \til{A} \to A$ a {\em NC semi-free resolution of $A$ in 
$\catt{DGRng} \centover \K$}. 
\end{dfn}

The important result of this subsection is the next existence theorem. 

\begin{thm} \label{thm:4310}
Let $A$ be a DG $\K$-ring. There exists a NC semi-free DG ring
resolution $u : \til{A} \to A$ of $A$ relative to $\K$. 
\end{thm}

The proof of the theorem comes after two lemmas. 
If $X$ and $Y$ are graded sets, a degree $i$ function $f : X \to Y$ is a 
function such that $f(X^j) \sub Y^{j + i}$ for all $j$. Likewise,
a degree $i$ function $f : X \to B$ to a graded ring $B$ is a function 
such that $f(X^j) \sub B^{j + i}$ for all $j$.

\begin{lem} \label{lem:4315}
Let $X$ be a graded set, and consider the graded $\K$-ring 
$B := \K \bra{X}$. Suppose $\d_X : X \to B$ is a degree $1$ function.
\begin{enumerate}
\item  The function $\d_X$ extends uniquely to a $\K$-linear degree $1$ 
derivation $\d_B : B \to B$. 

\item The derivation $\d_B$ is a differential on $B$, i.e.\ 
$\d_B \circ \d_B = 0$, if and only if $\d_B(\d_X((x)) = 0$ for every 
$x \in X$. 
\end{enumerate}
\end{lem}

\begin{proof} \mbox{}

\smallskip \noindent
(1) For elements $x_1, \ldots, x_n \in X$, with 
$k_j := \opn{deg}(x_j)$, define 
\begin{equation} \label{eqn:4315}
\d_B(x_1 \cdots x_n) := \sum_{j = 1}^n \, 
(-1)^{k_1 + \cdots + k_{j - 1}} \cd x_1 \cdots x_{j - 1} \cd 
\d_X(x_j) \cd x_{j + 1} \cdots x_{n} . 
\end{equation}
This extends uniquely to a $\K$-linear homomorphism $\d_B : B \to B$ of degree 
$1$. A calculation shows that the graded Leibniz rule holds, so this is  
degree $1$ derivation. 

\medskip \noindent
(2) It is enough to verify that $\d_B \circ \d_B = 0$ on monomials. 
This is an elementary calculation, using formula (\ref{eqn:4315}). 
\end{proof}

\begin{lem} \label{lem:4316}
Let $B, C \in \catt{DGRng} \centover \K$. Assume that $B$ is NC semi-free, with 
a multiplicative semi-basis $(X, F)$.  
Let $w : X \to C$ be a degree $0$ function. 
\begin{enumerate}
\item  The function $w$ extends uniquely to a graded $\K$-ring homomorphism 
$w : B^{\natural} \to C^{\natural}$. 

\item The graded ring homomorphism $w$ becomes a DG ring homomorphism 
$w : B \to C$ if and only if 
$\d_C(w(x)) = w(\d_B(x))$ for every $x \in X$. 
\end{enumerate}
\end{lem}

\begin{proof} 
Item (1) is clear. As for item (2), it is enough to verify that $\d_C \circ w = 
w \circ \d_B$ on monomials $x_1 \cdots x_n \in B$. This is done by induction on 
$n$, using the graded Leibniz rule. 
\end{proof}

\begin{proof}[Proof of Theorem \tup{\ref{thm:4310}}]
The strategy is this: we are going to construct an increasing sequence of 
semi-free DG rings 
\begin{equation} \label{eqn:4945}
F_0(\til{A}) \sub F_1(\til{A}) \sub F_2(\til{A}) \sub \cdots ,
\end{equation}
together with a compatible system of DG ring homomorphisms 
$F_j(u) : \lb F_j(\til{A}) \to A$.
At the same time we will construct an increasing sequence of graded sets
\begin{equation} \label{eqn:4946}
\varnothing = F_{-1}(X) \sub F_0(X) \sub F_1(X) \sub F_2(X) \sub \cdots ,
\end{equation}
such that there is equality of graded rings
$F_j(\til{A})^{\natural} = \K \bra{F_j(X)}$
for every $j \geq 0$. 

Then we will define the DG ring
\begin{equation} \label{eqn:4310}
\til{A} := \lim_{j \to} \,  F_j(\til{A}) , 
\end{equation}
the DG ring homomorphism
\begin{equation} \label{eqn:4311}
u := \lim_{j \to} \,  F_j(u) : \til{A} \to  A ,
\end{equation}
and the graded set 
\begin{equation} \label{eqn:4775}
X := \lim_{j \to} \,  F_j(X) . 
\end{equation}
The graded set $X$ is filtered by 
$F := \bigl\{ F_j(X) \bigr\}_{j \geq -1}$. By the construction 
$(X, F)$ will be a multiplicative semi-basis of $\til{A}$, so that the latter 
is a semi-free DG ring. The construction will also ensure that $u$ is a 
quasi-isomorphism. All this will be done in several steps. 

\medskip \noindent
Step 1. In this step we deal with $j = 0$. 
Let's choose a collection $\{ \bar{a}_x \}_{x \in F_0(X)}$
of nonzero homogeneous elements of $\opn{H}(A)$ that generates it as a 
$\K$-ring. We make $F_0(X)$ into a graded set by declaring
$\opn{deg}(x) := \opn{deg}(\bar{a}_x)$ for $x \in F_0(X)$.
Next, for each $x$ we choose a homogeneous cocycle 
$a_x \in \opn{Z}(A)$ that represents the cohomology class $\bar{a}_x$. 
In this way we obtain a degree $0$ function
$F_0(u) : F_0(X) \to A$, $F_0(u)(x) := a_x$. 

Define the DG ring 
$F_0(\til{A}) := \K \bra{F_0(X)}$, with zero differential. 
According to Lemma \ref{lem:4316} the function 
$F_0(u) : F_0(X) \to A$ extends uniquely to a DG ring homomorphism 
$F_0(u) : F_0(\til{A}) \to A$.
The choice of the collection $\{ \bar{a}_x \}_{x \in F_0(X)}$
guarantees that the graded $\K$-ring homomorphism 
$\opn{H}(F_0(u)) : \opn{H}(F_0(\til{A})) \to \opn{H}(A)$ 
is surjective.

\medskip \noindent
Step 2. In this step and the next one we deal with the induction. Let
$j \geq 0$. The assumption is that we have a 
semi-free DG ring $F_j(\til{A})$, with a multiplicative semi-basis
\begin{equation} \label{eqn:4776}
\varnothing = F_{-1}(X) \sub F_0(X) \sub \cdots \sub F_{j - 1}(X) \sub 
F_{j}(X) \sub F_{j}(X) \sub \cdots  
\end{equation}
(constant from level $j$ onward), and a DG ring homomorphism 
$F_j(u) : F_j(\til{A}) \to A$
such that the graded ring homomorphism 
$\opn{H}(F_j(u)) : \opn{H}(F_j(\til{A})) \to \opn{H}(A)$ 
is surjective.

Define the graded two-sided ideal 
$K_j := \opn{Ker} \bigl( \opn{H}(F_j(u)) \bigr) \sub \opn{H}(F_j(\til{A}))$.
So there is an exact sequence of graded $\opn{H}(F_j(\til{A}))$-bimodules 
\begin{equation} \label{eqn:4316}
0 \to K_j \to \opn{H}(F_j(\til{A})) \xar{\opn{H}(F_j(u))}
\opn{H}(A) \to  0 . 
\end{equation}
Choose a collection $\{ \bar{a}_y \}_{y \in Y_{j + 1}}$
of nonzero homogeneous elements of $K_j$ that generates it as an 
$\opn{H}(F_j(\til{A}))$-bimodule (i.e.\ as a two-sided ideal). 
We make $Y_{j + 1}$ into a graded set by declaring
$\opn{deg}(y) := \opn{deg}(\bar{a}_y) - 1$ 
for $y \in Y_{j + 1}$.
Next, for each element $y \in Y_{j + 1}$ we choose a homogeneous cocycle 
$a_y \in \opn{Z} \bigl( F_j(\til{A}) \bigr)$ that represents the cohomology 
class $\bar{a}_y$. In this way we obtain 
a collection $\{ a_y \}_{y \in Y_{j + 1}}$ of homogeneous cocycles in 
$F_j(\til{A})$. 

Define the graded set 
$F_{j + 1}(X) := F_{j}(X) \sqcup Y_{j + 1}$
and the graded ring \lb 
$F_{j + 1}(\til{A})^{\natural} := \K \bra{ F_{j + 1}(X) }$.
By definition $F_{j}(\til{A})^{\natural}$ is a graded subring of 
$F_{j + 1}(\til{A})^{\natural}$. 
Next, define the degree $1$ function 
$\d_{F_{j + 1}(X)} : F_{j+ 1}(X) \to F_{j + 1}(\til{A})$
by the formula  
\[ \d_{F_{j + 1}(X)}(x) := 
\begin{cases}
\d_{F_{j}(\til{A})}(x) & \tup{if} \ \ x \in F_{j}(X)
\\
a_x & \tup{if} \ \ x \in Y_{j + 1} .   
\end{cases} \]
In this formula we are using the collection of homogeneous cocycles
$\{ a_y \}_{y \in Y_{j + 1}}$ that has been chosen above. 
According to Lemma \ref{lem:4315} the function $\d_{F_{j + 1}(X)}$ extends 
uniquely to a differential on $F_{j + 1}(\til{A})$, making it into a DG ring. 
And of course $F_{j}(\til{A}) \sub F_{j + 1}(\til{A})$
is a DG subring. 

Let us denote by 
$v_j : F_{j}(\til{A}) \to F_{j + 1}(\til{A})$
the inclusion. 
Because the cocycles $a_y \in F_{j}(\til{A})$, for $y \in Y_{j + 1}$, become 
coboundaries in $F_{j + 1}(\til{A})$, 
we see that 
\begin{equation} \label{eqn:4317}
K_j \sub \opn{Ker} (\opn{H}(v_j)) \sub \opn{H}(F_{j}(\til{A})) .
\end{equation}

Observe that the DG ring $F_{j + 1}(\til{A})$ is semi-free, with 
multiplicative semi-basis
\[ \varnothing = F_{-1}(X) \sub F_0(X) \sub \cdots \sub F_{j}(X) \sub 
F_{j + 1}(X) \sub F_{j + 1}(X) \sub \cdots , \]
constant from level $j + 1$ onward.

\medskip \noindent
Step 3. In this step we continue step 2, to construct the DG ring homomorphism 
$F_{j + 1}(u) : F_{j + 1}(\til{A}) \to A$.
For each element $y \in Y_{j + 1}$ the cohomology class 
$\bar{a}_y \in \opn{H}(F_{j}(\til{A}))$ belongs to the kernel of 
$\opn{H}(F_j(u))$. This means that the cocycle 
$F_j(u)(a_y) \in A$ is a coboundary. So we can find a homogeneous element
$b_y \in A$ such that 
$\d_A(b_y) = F_j(u)(a_y)$. 
Define the degree $0$ function 
$F_{j + 1}(u) : F_{j + 1}(X) \to A$
by the formula 
\[ F_{j + 1}(u)(x) := 
\begin{cases}
F_{j}(u)(x) & \tup{if} \ \ x \in F_{j}(X) 
\\
b_{x} & \tup{if} \ \ x \in Y_{j + 1}  .   
\end{cases} \]
By Lemma \ref{lem:4316} this function extends uniquely to a DG ring 
homomorphism 
$F_{j + 1}(u) : F_{j + 1}(\til{A}) \to A$.
Note that the restriction of $F_{j + 1}(u)$ to $F_{j}(\til{A})$
coincides with $F_{j}(u)$. 

\medskip \noindent
Step 4. Applying steps 2-3 recursively, we obtain direct systems 
$\{ F_j(\til{A}) \}_{j \geq 0}$ and $\{ F_j(u) \}_{j \geq 0}$
of DG rings  and homomorphisms, and we can define the 
DG ring $\til{A}$ by formula (\ref{eqn:4310}), and the DG ring homomorphism 
$\til{u}$ by formula (\ref{eqn:4311}). 
The graded set $X$ from formula (\ref{eqn:4775}) is filtered, and by the 
construction in steps 1 and 2 the filtered graded set $(X, F)$ is a 
multiplicative semi-basis of $\til{A}$. 

It remains to prove that the graded ring homomorphism
$\opn{H}(u) :  \opn{H}(\til{A}) \to \opn{H}(A)$
is bijective. 
Recall that for each $j \geq0$ we denoted by $v_j$ the inclusion
$F_{j}(\til{A}) \sub F_{j + 1}(\til{A})$. 
Define the graded ring
$L^j := \opn{Im} (\opn{H}(v_j)) \sub \opn{H}(F_{j + 1}(\til{A}))$. 
We get a commutative diagram 
\[ \UseTips \xymatrix @C=3.5ex @R=5.5ex {
0 
\ar[r]
&
*++{K_j}
\ar@{>->}[rr]^(0.4){\mrm{inc}}
&
&
\opn{H}(F_{j}(\til{A}))
\ar@{->>}[rr]^{\opn{H}(F_{j}(u))}
\ar@{->>}[d]_{\al_j}
\ar@(l,l)[dd]_{\opn{H}(v_j)}
&
&
*++{\opn{H}(A)}
\ar[r]
&
0
\\
&
&
&
*++{L_j}  
\ar@{>->}[d]_{\mrm{inc}}
\ar@{->>}[urr]_{\be_j}
\\
&
&
&
\opn{H}(F_{j + 1}(\til{A}))
\ar@{->>}@(r,d)[uurr]_{\opn{H}(F_{j + 1}(u))}
} \]
in $\dcat{G}_{\mrm{str}}(\K)$. The top row is an exact sequence (it is 
(\ref{eqn:4316})). 
Because $\al_j$ is surjective, there is equality 
$\opn{Ker}(\be_j) = \al_j \bigl( \opn{Ker}(\opn{H}(F_{j}(u)) \bigr) 
= \al_j(K_j)$.
But by formula (\ref{eqn:4317}) we know that 
$\al_j(K_j) = 0$. The conclusion is that 
$\be_j : L_j \to \opn{H}(A)$ is bijective. 
Hence $\{ L_j \}_{j \geq 0}$ is a constant direct system, and 
$\lim_{j \to} \be_j : \lim_{j \to} L_j \to \opn{H}(A)$
is bijective.
 
Now the direct systems 
$\{ L_j \}_{j \geq 0}$ and 
$\bigl\{ \opn{H}(F_{j}(\til{A})) \bigr\}_{j \geq 0}$
are sandwiched. Therefore the second direct system has a limit too, and it is 
the same limit; i.e.\ 
\[ \lim_{j \to} \, \opn{H}(F_{j}(u)) : 
\lim_{j \to} \, \opn{H}(F_{j}(\til{A})) \to \opn{H}(A) \]
is bijective. Because cohomology commutes with direct limits, 
and by formula (\ref{eqn:4310}), this implies that $\opn{H}(u)$
is bijective.
\end{proof}

\begin{exer} \label{exer:4315}
Modify the proof of Theorem \ref{thm:4310} to make the quasi-iso\-morphism 
$u : \til{A} \to A$ surjective. 
\end{exer}

\begin{cor} \label{cor:4605}
If $\opn{H}(A)$ is nonpositive, then $A$ has a semi-free DG ring resolution 
$\til{A} \to A$  such that $\til{A}$ is nonpositive. 
\end{cor}

\begin{proof}
Under this assumption, the graded set $X$ in the proof of Theorem \ref{thm:4310} 
is nonpositive.
\end{proof}

\begin{rem} \label{rem:4605}
NC semi-free DG rings have important lifting properties within 
$\catt{DGRng} \centover \K$.
V. Hinich, in \cite{Hi1}, calls them {\em standard cofibrant objects} of the 
Quillen model structure on $\catt{DGRng} \centover \K$.

One can also consider the full subcategory $\catt{DGRng}^{\leq 0} \centover \K$
of nonpositive NC DG rings. NC nonpositive semi-free DG rings
are studied in detail in \cite{Ye9}. This is done in a slightly more general 
context: the base ring $\K$ is allowed to be a commutative DG ring (see next 
remark), and $\K \to A$ is a central homomorphism of DG 
rings (i.e.\ the image of $\K$ is in $\opn{Cent}(A)$). 
The lifting properties are the same. 

See Remark \ref{rem:4970} for the related $\mrm{A}_{\infty}$ theory. 
\end{rem}

\begin{rem} \label{rem:4945}
There are also {\em commutative semi-free DG rings}.
Recall that a commutative DG ring in this book means nonpositive and strongly 
commutative; see Definitions \ref{dfn:3091} and \ref{dfn:3090}. The
commutative DG rings form the 
category $\catt{DGRng}^{\leq 0}_{\mrm{sc}} / \K$,
which is a full subcategory of 
$\catt{DGRng} \centover \K$.

A nonpositive graded set $S$ gives rise to the commutative polynomial ring 
$\K[S]$, see Example \ref{exa:4625}. 
A commutative DG ring $\til{A}$ is called commutative semi-free if it admits    
a graded ring isomorphism $\til{A}^{\natural} \cong \K[S]$ for some 
nonpositive graded set $S$. The graded set $S$ is automatically filtered by 
degree, namely $F_j(S) := \bigcup_{i \geq -j} S^i$.
This implies (cf.\ Proposition \ref{prop:4305})
that $\til{A}$ is semi-free as a DG $\K$-module . 

In \cite{Ye9} it is proved that every 
$A \in \catt{DGRng}^{\leq 0}_{\mrm{sc}} / \K$
admits a surjective quasi-isomorphism $\til{A} \to A$ from a 
commutative semi-free DG ring $\til{A}$. 
It is also proved that the commutative semi-free DG rings have good lifting 
properties in $\catt{DGRng}^{\leq 0}_{\mrm{sc}} / \K$. 
The situation studied in \cite{Ye9} is a bit more general: the base ring $\K$ 
is allowed to be itself a commutative DG ring; so the theory becomes 
relative. 

If in Example \ref{exa:4205} we take the topological space $X$ to be a 
single point, then the category 
$\catt{DGRng}^{\leq 0}_{\mrm{sc}} / \K_X$
discussed there coincides with 
$\catt{DGRng}^{\leq 0}_{\mrm{sc}} / \K$,
and the semi-pseudo-free DG $\K_X$-rings are just the semi-free commutative DG 
rings mentioned here. 

Unlike the geometric setting, the category 
$\catt{DGRng}^{\leq 0}_{\mrm{sc}} / \K$
admits a Quillen model structure (at least when $\K$ contains $\Q$),
and the calculus of fractions for the derived category 
$\dcat{D}(\catt{DGRng}^{\leq 0}_{\mrm{sc}} / \K)$
can be deduced from it. 

The fully faithful embedding 
\[ \catt{DGRng}^{\leq 0}_{\mrm{sc}} / \K \to 
\catt{DGRng}^{\leq 0} \centover \K \]
induces a functor 
\[ \dcat{D}(\catt{DGRng}^{\leq 0}_{\mrm{sc}} / \K) \to
\dcat{D}(\catt{DGRng}^{\leq 0} \centover \K) \]
between the derived categories, which presumably is neither full nor 
faithful (but we did not explore this question seriously). 
\end{rem}

\mysubsection{The Derived Tensor-Evaluation Morphism}
\label{subsec:der-tens-eval}

In this subsection we pre\-sent an important theorem (Theorem 
\ref{thm:4320}), that will reappear -- with modifications -- several times in 
the book. Before presenting it, we need a few finiteness and boundedness 
conditions, some old and some new. 

Pseudo-finite semi-free DG modules were introduced in Definition 
\ref{dfn:1925}. The saturated full triangulated subcategory generated by a set 
of objects was defined in Definition \ref{dfn:4916}. 

\begin{dfn} \label{dfn:4470}
Let $A$ be a DG ring. A DG $A$-module $L$ is called 
{\em derived pseudo-finite}%
\index{Differential graded module! derived pseudo-finite}
if it belongs to the saturated full triangulated subcategory of $\dcat{D}(A)$ 
generated by the pseudo-finite semi-free DG $A$-modules. 
\end{dfn}

\begin{exa} \label{exa:4471}
Suppose $A$ is a nonpositive cohomologically left pseudo-noe\-ther\-ian DG 
ring. Then the derived pseudo-finite DG $A$-modules are precisely the 
objects of $\dcat{D}_{\mrm{f}}^{-}(A)$. One direction is by 
Theorem \ref{thm:3340}, and the other direction is trivial.  
\end{exa}

\begin{exa} \label{exa:4470}
If $L$ is an algebraically perfect DG $A$-module (to be defined later, 
see Definition \ref{dfn:3401}) then it is derived pseudo-finite.
\end{exa}

\begin{rem} \label{rem:4470}
Let $L$ be a DG $A$-module. Following \cite{SGA6}, let us say that 
$L$ is {\em pseudo-coherent}%
\index{Complex in abelian category! pseudo-coherent}%
\index{Differential graded module! pseudo-coherent}
if $L$ is isomorphic in $\dcat{D}(A)$ to a 
pseudo-finite semi-free DG $A$-module $P$.
Clearly pseudo-coherent implies derived pseudo-finite. According to 
\cite[Lemma tag=064X]{SP} the converse is true when $A$ is a ring; but we do 
not know if the converse is true in general. 

Anyhow, the attribute ``derived pseudo-finite'' is very useful, 
and it can replace ``pseudo-coherent'' for most applications. For instance, 
see  condition (i) in Theorem \ref{thm:4320} below. 
\end{rem}

In Definition \ref{dfn:2390} we saw the notion of bounded below flat 
concentration of a DG module. Here is a refinement of it. 

\begin{dfn} \label{dfn:4341}
Let $A$ be a DG ring. 
\begin{enumerate}
\item We say that a DG $A$-module $P$ has 
{\em tensor displacement}%
\index{Tensor displacement}
inside an integer interval $[i_0, i_1]$ if for every 
$N \in \dcat{C}(A^{\mrm{op}})$ this inclusion of integer intervals holds:
\[ \opn{con}(N \ot_A P) \sub \opn{con}(N) + [i_0, i_1] . \]

\item We say that a DG $A$-module $P$ has
{\em bounded below tensor displacement}%
\index{Tensor displacement! bounded below}
if it has tensor displacement inside 
the integer interval $[i_0, \infty]$ for some $i_0 \in \Z$. 

\item A DG $A$-module $M$ is said to have 
{\em derived bounded below tensor displacement}%
\index{Tensor displacement! derived bounded below}
if $M$ belongs to the saturated full triangulated subcategory of 
$\dcat{D}(A)$ generated by the K-flat DG $A$-modules $P$ that have bounded 
below tensor displacement.
\end{enumerate}
\end{dfn}

\begin{exa} \label{exa:4601}
The DG $A$-module $M := A$ has derived bounded below tensor displacement, 
regardless of the boundedness of $A$ or $\opn{H}(A)$.

More generally, if $M$ is an algebraically perfect DG $A$-module (see 
Definition \ref{dfn:3401}) then it has derived bounded below tensor 
displacement.
\end{exa}

\begin{exa}  \label{exa:4473}
If $A$ is a ring, and $M$ is a complex of $A$-modules of finite flat dimension, 
then $M$ has derived bounded below tensor displacement.
See Proposition \ref{prop:2390}. 
\end{exa}

\begin{rem} \label{rem:4471} 
Clearly, if $M$ has derived bounded below tensor displacement, then it has 
bounded below flat concentration.
We don't know if the converse is true in general. 
\end{rem}

\begin{rem} \label{rem:4602}
Definition \ref{dfn:4341} can be adapted to other kinds of boundedness 
conditions. It can also be adapted to Hom displacement, in the first or second 
argument, to give refined notions of projective and injective displacements.
\end{rem}
 
\begin{thm}[Tensor-Evaluation] \label{thm:4320}
Let $A$ and $B$ be DG rings, and let $L \in \dcat{D}(A)$, 
$M \in \dcat{D}(A \ot B^{\mrm{op}})$ and $N \in \dcat{D}(B)$
be DG modules. There is a morphism 
\[ \opn{ev}^{\mrm{R, L}}_{L, M, N} : \opn{RHom}_{A}(L, M) \ot^{\mrm{L}}_{B} N 
\to \opn{RHom}_{A}(L, M \ot^{\mrm{L}}_{B} N) \]
in $\dcat{D}(\K)$, called {\em derived tensor-evaluation}%
\index{Tensor-evaluation morphism! derived},
that enjoys the properties below. 
\begin{enumerate}
\item The morphism $\opn{ev}^{\mrm{R, L}}_{L, M, N}$ is functorial in the 
objects $L, M, N$.

\item If all conditions \tup{(i)-(iii)} below hold, then 
$\opn{ev}^{\mrm{R, L}}_{L, M, N}$ is an isomorphism. 
\begin{enumerate}
\rmitem{i} The DG $A$-module $L$ is derived pseudo-finite 
\tup{(}Definition \tup{\ref{dfn:4470})}.

\rmitem{ii} The DG $A$-$B$-bimodule $M$ is cohomologically bounded below
\tup{(}Definition \tup{\ref{dfn:4785}(2))}.

\rmitem{iii} The DG $B$-module $N$ has derived bounded below tensor 
displacement \tup{(}Definition \tup{\ref{dfn:4341}(3))}.
\end{enumerate}
\end{enumerate}
\end{thm}

We need a lemma first. 

\begin{lem} \label{lem:4735}
Suppose we are given\tup{:}
\begin{itemize}
\item A pseudo-finite semi-free DG $A$-module $P$, with pseudo-finite 
semi-free filtration $\{ F_j(P) \}_{j \geq -1}$.

\item A bounded below DG $A$-$B$-bimodule $M$.

\item A DG $B$-module $Q$ with bounded below tensor displacement.
\end{itemize}
Then for every degree $l$ there is an index $j$, depending on $l$, such that 
the canonical homomorphisms  
\[ \bigl( \opn{Hom}_{A}(P, M) \ot_{B} Q \bigr)^l \to
\bigl( \opn{Hom}_{A}(F_{j}(P), M) \ot_{B} Q \bigr)^l \]
and 
\[ \bigl( \opn{Hom}_{A}(P, M \ot_{B} Q) \bigr)^l \to
\bigl( \opn{Hom}_{A}(F_{j}(P), M \ot_{B} Q) \bigr)^l  \]
in $\dcat{M}(\K)$ are bijective. 
\end{lem}

\begin{proof}
In the notation of Definition \ref{dfn:1925}, there are numbers $i_1 \in \Z$ 
and $r_k \in \N$ such that
\begin{equation} \label{eqn:4735}
\bigl( P / F_j(P) \bigr)^{\natural} \cong 
\bigoplus_{k \geq j + 1} \opn{T}^{-i_1 + k}(A^{\natural})^{\oplus r_k}
\end{equation}
in $\dcat{G}_{\mrm{str}}(A^{\natural})$ for all $j \geq 0$. 
In terms of bases, the graded free $A^{\natural}$-module 
$\bigl( P / F_j(P) \bigr)^{\natural}$
has a basis concentrated in degrees $\leq i_1 - j - 1$.

Say $M$ is concentrated in the degree interval $[s_0, \infty]$,
and $Q$ has tensor displacement inside the degree interval $[t_0, \infty]$,
for some $s_0, t_0 \in \Z$. From equation (\ref{eqn:4735}) we see that for 
every $j$, both graded $\K$-modules
\[  \opn{Hom}_{A} \bigl(  P / F_j(P), M  \bigr) \ot_{B} Q 
\quad \tup{and} \quad 
\opn{Hom}_{A} \bigl(  P / F_j(P), M  \ot_{B} Q \bigr) \]
are concentrated in the degree interval 
\[ [-i_1 + j + 1  + s_0 + t_0, \infty] \sub \Z . \]
The split short  exact sequence 
\[ 0 \to F_j(P)^{\natural} \to P^{\natural} \to 
\bigl( P / F_j(P) \bigr)^{\natural} \to 0 \]
in $\dcat{G}_{\mrm{str}}(A)$ tells us that for every 
$j \geq l + i_1 - s_0 - t_0$
the homomorphisms in question are bijective. 
\end{proof}

\begin{proof}[Proof of Theorem \tup{\ref{thm:4320}}] \mbox{} 

\smallskip \noindent
(1) Choose a K-projective resolution 
$P \to L$ in $\dcat{C}_{\mrm{str}}(A)$ and a K-projective resolution 
$Q \to N$ in $\dcat{C}_{\mrm{str}}(B)$. These choices are unique up to 
homotopy. Then 
$\opn{ev}^{\mrm{R, L}}_{L, M, N}$ 
is represented by the obvious homomorphism 
\[ \opn{ev}_{P, M, Q} : \opn{Hom}_{A}(P, M) \ot_{B} Q \to 
\opn{Hom}_{A}(P, M \ot_{B} Q) \]
in $\dcat{C}_{\mrm{str}}(\K)$. 
The morphism $\opn{ev}^{\mrm{R, L}}_{L, M, N}$ in $\dcat{D}(\K)$ does not 
depend on the choices, and hence it is functorial. 

\medskip \noindent 
(2) This is done in several steps. 

\smallskip \noindent 
Step 1. Let $P$ be a finite semi-free DG $A$-module; namely $P$ admits a 
pseudo-finite semi-free filtration $F = \{ F_j(P) \}_{j \geq -1}$
(see Definition \ref{dfn:1925}) such that $P = F_{j_1}(P)$ for some 
$j_1 \in \N$. Then the graded $A^{\natural}$-module $P^{\natural}$ is finite 
free. This implies that 
\[ \opn{ev}_{P, M, Q} : \opn{Hom}_{A}(P, M) \ot_{B} Q \to 
\opn{Hom}_{A}(P, M \ot_{B} Q) \]
is an isomorphism in $\dcat{C}_{\mrm{str}}(\K)$ for all 
$M \in \dcat{C}(A \ot B^{\mrm{op}})$ and $Q \in \dcat{C}(B)$. 

\medskip \noindent 
Step 2. Assume that $M$ is bounded below. 
Let $Q$ be a K-flat DG $B$-module of bounded below tensor displacement, and let 
$P$ be a pseudo-finite semi-free DG $A$-module, with pseudo-finite semi-free 
filtration $F = \{ F_j(P) \}_{j \geq -1}$.
By Lemma \ref{lem:4735}, for every $l \in \Z$ there is an index $j$ such that  
\[ \bigl( \opn{Hom}_{A}(P, M) \ot_{B} Q \bigr)^l =
\bigl( \opn{Hom}_{A}(F_{j}(P), M) \ot_{B} Q \bigr)^l \]
and 
\[ \bigl( \opn{Hom}_{A}(P, M \ot_{B} Q) \bigr)^l  =
\bigl( \opn{Hom}_{A}(F_{j}(P), M \ot_{B} Q) \bigr)^l . \]
By step 1 the homomorphism 
\[ \opn{ev}_{F_{j}(P), M, Q} : 
\opn{Hom}_{A} \bigl( F_{j}(P), M \bigr) \ot_{B} Q \to 
\opn{Hom}_{A} \bigl( F_{j}(P), M \ot_{B} Q \bigr) \]
is bijective. We conclude that 
\[ \opn{ev}_{P, M, Q} : \bigl( \opn{Hom}_{A}(P, M) \ot_{B} Q \bigr)^l \to 
\opn{Hom}_{A}(P, M \ot_{B} Q)^l \]
is bijective. Since $l$ is arbitrary, it follows that 
$\opn{ev}_{P, M, Q}$ is an isomorphism in $\dcat{C}_{\mrm{str}}(\K)$.

\medskip \noindent 
Step 3. In step 2 we proved that the morphism 
\[ \opn{ev}^{\mrm{R, L}}_{L, M, N} : \opn{RHom}_{A}(L, M) \ot^{\mrm{L}}_{B} N 
\to \opn{RHom}_{A}(L, M \ot^{\mrm{L}}_{B} N) \]
in $\dcat{D}(\K)$ is an isomorphism if $M$ is bounded below, 
$L = P$ is a pseudo-finite semi-free DG 
$A$-module, and $N = Q$ is K-flat DG $B$-module of bounded below tensor 
displacement. If $M$ is only assumed to satisfy condition (ii), then it is 
quasi-isomorphic to a suitable smart truncation of it which is bounded below; 
therefore $\opn{ev}^{\mrm{R, L}}_{L, M, N}$ is an isomorphism under this 
assumption too. 

Fixing $M$ and $N$, the full subcategory of $\dcat{D}(A)$
on the objects $L$ for which $\opn{ev}^{\mrm{R, L}}_{L, M, N}$ is an 
isomorphism is a saturated full triangulated subcategory (by Proposition 
\ref{prop:4603}). Therefore  
$\opn{ev}^{\mrm{R, L}}_{L, M, N}$ is an isomorphism for every $L, M, N$ 
such that $L$ satisfies condition (i), $M$ satisfies condition (ii), and 
$N = Q$ is K-flat of bounded below tensor displacement.

Now we fix $L$ and $M$ that satisfy condition (i) and (ii) respectively. 
The full subcategory of $\dcat{D}(B)$ on the objects $N$ for which 
$\opn{ev}^{\mrm{R, L}}_{L, M, N}$ is an 
isomorphism is a saturated full triangulated subcategory. By the previous 
paragraph, this subcategory contains the K-flat complexes $Q$ of bounded below 
tensor displacement. Therefore it contains all the complexes $N$ satisfying 
condition (iii).
\end{proof}

\begin{rem} \label{rem:4677}
Let us explain the syntax of the expression 
$\opn{ev}^{\mrm{R, L}}_{L, M, N}$
in Theorem \ref{thm:4320}. The letters ``ev'' stand for ``evaluation''.
The letters ``R'' and ``L'' in the superscript refer to the 
left and right derived functors, whereas the letters ``L'', ``M'' and ``N'' 
in the subscript refer to the DG modules involved. Note the different fonts. 
This sort of syntax is used frequently in the book; cf.\ the expressions 
$\opn{hm}^{\mrm{R}}_{M}$ in Definition \ref{dfn:2394},
$\opn{ev}^{\mrm{R, R}}$ in Lemma \ref{lem:2156},
$\ttev^{\mrm{R, L}}_{\m, (-)}$ in Theorem \ref{thm:3750} and
$\opn{fadj}_{u, M, N}^{\mrm{L}}$ in Theorem \ref{thm:3085}.
\end{rem}

\mysubsection{Hom-Tensor Formulas for Weakly Commutative DG Rings}
\label{subsec:adjs-comm}

In this subsection we assume that all DG rings are weakly commutative  
(see Definition \ref{dfn:3090}).
For such a DG ring $A$, its $0$-th cohomology 
$\opn{H}^0(A)$ is a commutative ring.
 
\begin{prop} \label{prop:3107}
For a weakly commutative DG ring $A$, the category $\dcat{D}(A)$ is 
$\opn{H}^0(A)$-linear. 
\end{prop}

\begin{proof}
This is true already for the homotopy category. 
Indeed, for any $M, N \in \dcat{C}(A)$ the 
$\K$-module 
$\opn{H}^0 \bigl( \opn{Hom}_A(M, N) \bigr)$
is an $\opn{H}^0(A)$-module; and composition in $\dcat{K}(A)$ is 
$\opn{H}^0(A)$-bilinear. 
\end{proof}

Let $A$ be a DG ring. There is a canonical isomorphism 
 $u : A \iso A^{\mrm{op}}$ to its opposite. The formula is 
$u(a) := (-1)^i \cd a$
for $a \in A^i$. This implies that every left DG $A$-module can be made into a 
right DG $A$-module, and vice-versa. The formula relating the left and right 
actions is 
$m \cd a = (-1)^{i \cd j} \cd a \cd m$
for $a \in A^i$ and $m \in M^j$. 
On the level of categories we obtain an isomorphism of DG categories 
$\dcat{C}(A) \cong \dcat{C}(A^{\mrm{op}})$,
which is the identity functor on the underlying DG $\K$-modules. 

Weak commutativity makes the tensor and Hom functors $A$-bilinear 
(in the graded sense), and therefore their derived functors have more 
structure: they are $\opn{H}^0(A)$-bilinear triangulated bifunctors 
\begin{equation} \label{eqn:2471}
(- \ot^{\mrm{L}}_A -) : \dcat{D}(A) \times \dcat{D}(A) \to \dcat{D}(A) 
\end{equation}
and 
\begin{equation} \label{eqn:2472}
\opn{RHom}_{A}(-, -) : \dcat{D}(A)^{\mrm{op}} \times \dcat{D}(A) \to 
\dcat{D}(A) . 
\end{equation}

When the DG rings in Theorem \ref{thm:2363} are weakly commutative, 
this result can be amplified: 

\begin{prop} \label{prop:3150} 
In the situation of Theorem \tup{\ref{thm:2363}}, assume that $A$ and $B$ 
are weakly commutative. 
\begin{enumerate}
\item The restriction functor 
$\opn{Rest}_u : \dcat{D}(B) \to \dcat{D}(A)$
is an equivalence of $\opn{H}^0(A)$-linear triangulated categories, with 
quasi-inverse $\opn{LInd}_u$. 

\item For every $M, N \in \dcat{D}(B)$ there is an isomorphism
$M \ot^{\mrm{L}}_A N \iso M \ot^{\mrm{L}}_B N$
in $\dcat{D}(A)$. This isomorphism is functorial in $M$ and $N$.

\item For every $M, N \in \dcat{D}(B)$ there is an isomorphism 
$\opn{RHom}_{A}(M, N) \iso \lb \opn{RHom}_{B}(M, N)$
in $\dcat{D}(A)$. This isomorphism is functorial in $M$ and $N$.
\end{enumerate}
\end{prop}

\begin{proof} \mbox{}

\smallskip \noindent 
(1) The ring homomorphism 
$\opn{H}^0(u) : \opn{H}^0(A) \to \opn{H}^0(B)$
makes $\dcat{D}(B)$ into an $\opn{H}^0(A)$-linear category, and then 
$\opn{Rest}_u : \dcat{D}(B) \to \dcat{D}(A)$
becomes an $\opn{H}^0(A)$-linear functor. 

\medskip \noindent 
(2,3) Going over the steps in the proof of Theorem \ref{thm:2363}, we see that 
all the moves are $A$-linear (in the graded sense). 
\end{proof}

Suppose $A \to B$ is a DG ring homomorphism. For every 
$M \in \dcat{C}(A)$ and $N \in \dcat{C}(B)$, their tensor product 
$M \ot_A N$ is a DG $B$-module, via the action on $N$. 
Thus we get a DG bifunctor 
\begin{equation} \label{eqn:4955}
(- \ot_A -) : \dcat{C}(A) \times \dcat{C}(B) \to \dcat{C}(B) . 
\end{equation}
There is an isomorphism 
\begin{equation} \label{eqn:3101}
M \ot_A N \iso N \ot_A M 
\end{equation}
in $\dcat{C}_{\mrm{str}}(B)$, that uses the Koszul sign rule
$m \ot n \mapsto (-1)^{i \cd j} \cd n \ot m$
for homo\-geneous elements $m \in M^i$ and $n \in N^j$. 

\begin{prop} \label{prop:3103}
Let $A \to B$ be a homomorphism of weakly commutative DG rings. 
The DG bifunctor \tup{(\ref{eqn:4955})}
has a left derived bifunctor 
\[ (- \ot^{\mrm{L}}_{A} -) : \dcat{D}(A) \times \dcat{D}(B) \to \dcat{D}(B) . \]
This is an $\opn{H}^0(A)$-bilinear triangulated bifunctor.
\end{prop}

\begin{proof}
The left derived bifunctor exists by Theorem \ref{thm:2106}, using K-flat 
resolutions in the first argument (i.e.\ in $\dcat{C}_{\mrm{str}}(A)$). 
As for the $\opn{H}^0(A)$-bilinearity: this is already comes in the homotopy 
category level, namely for the bifunctor 
\[ (- \ot_{A} -) : \dcat{K}(A) \times \dcat{K}(B) \to \dcat{K}(B) . \qedhere \]
\end{proof}

Similarly we have an $\opn{H}^0(A)$-bilinear triangulated bifunctor
\[ (- \ot^{\mrm{L}}_{A} -) : \dcat{D}(B) \times \dcat{D}(A) \to \dcat{D}(B) . \]
When $A = B$ these operations coincide with (\ref{eqn:2471}). 

\begin{prop}[Symmetry] \label{prop:3102}
Let $A \to B$ be a homomorphism of weakly commutative DG rings. For every 
$M \in \dcat{D}(A)$ and $N \in \dcat{D}(B)$ there is an isomorphism 
$M \ot^{\mrm{L}}_{A} N \cong N \ot^{\mrm{L}}_{A} M$
in $\dcat{D}(B)$. This isomorphisms is functorial in $M$ and $N$.  
\end{prop}

\begin{proof}
Use the isomorphism (\ref{eqn:3101}). 
\end{proof}

\begin{prop}[Associativity] \label{prop:3100}
Let $A_1 \to A_2 \to A_3$ be homomorphisms between weakly commutative DG rings.
For each $i$ let $M_i \in \dcat{D}(A_i)$.
There is an isomorphism 
\[ (M_1 \ot^{\mrm{L}}_{A_1} M_2) \ot^{\mrm{L}}_{A_2} M_3 \cong 
M_1 \ot^{\mrm{L}}_{A_1} (M_2 \ot^{\mrm{L}}_{A_2} M_3) \]
in $\dcat{D}(A_3)$.
This isomorphisms is functorial in the $M_i$. 
\end{prop}

\begin{proof}
Choose K-flat resolutions $P_i \to M_i$ in $\dcat{C}(A_i)$. Then we have 
isomorphisms
\[ \begin{aligned}
& (M_1 \ot^{\mrm{L}}_{A_1} M_2) \ot^{\mrm{L}}_{A_2} M_3 
\cong (P_1 \ot_{A_1} P_2) \ot_{A_2} P_3
\\
& \qquad 
\cong^{\dag} P_1 \ot_{A_1} (P_2 \ot_{A_2} P_3)
\cong 
M_1 \ot^{\mrm{L}}_{A_1} (M_2 \ot^{\mrm{L}}_{A_2} M_3)  
\end{aligned} \]
in $\dcat{D}(A_3)$. The isomorphism $\cong^{\dag}$ is due to the usual 
associativity of $\ot$.
\end{proof}

Suppose $A \to B$ is a DG ring homomorphism. For every 
$M \in \dcat{C}(A)$ and $N \in \dcat{C}(B)$, the DG $\K$-module
$\opn{Hom}_{A}(N, M)$ has a DG $B$-module structure, via the action on $N$. 
Thus we get a DG bifunctor 
\begin{equation} \label{eqn:4956}
\opn{Hom}_{A}(-, -) : \dcat{C}(B)^{\mrm{op}} 
\times \dcat{C}(A) \to \dcat{C}(B) .
\end{equation}

\begin{prop} \label{prop:3105}
Let $A \to B$ be a homomorphism of weakly commutative DG rings. 
The DG bifunctor \tup{(\ref{eqn:4956})}
has a right derived bifunctor
\[ \opn{RHom}_{A}(-, -) : \dcat{D}(B)^{\mrm{op}} 
\times \dcat{D}(A) \to \dcat{D}(B) . \]
This is an $\opn{H}^0(A)$-bilinear triangulated bifunctor.
\end{prop}

\begin{proof}
The right derived bifunctor exists by Theorem \ref{thm:2000}, using K-injective
resolutions in the second argument (i.e.\ in $\dcat{C}_{\mrm{str}}(A)$). 
\end{proof}

\begin{prop}[Hom-Tensor Adjunction] \label{prop:3106}
Let $A_1 \to A_2 \to A_3$ be homomorphisms between weakly commutative DG rings.
For each $i$ let $M_i \in \dcat{D}(A_i)$.
There is an isomorphism 
\[ \opn{RHom}_{A_2} \bigl( M_3, \opn{RHom}_{A_1}(M_2, M_1) \bigr) \cong
\opn{RHom}_{A_1} \bigl( M_2 \ot^{\mrm{L}}_{A_2} M_3, M_1) \]
in $\dcat{D}(A_3)$.
This isomorphisms is functorial in the $M_i$. 
\end{prop}

\begin{proof}
Choose a K-flat resolution $P_2 \to M_2$ in $\dcat{C}(A_2)$,
and a K-injective resolution $M_1 \to I_1$ in $\dcat{C}(A_1)$. 
An easy adjunction calculation show that 
$\opn{Hom}_{A_1}(P_2, I_1)$ is K-injective in $\dcat{C}(A_2)$. 
Then we have isomorphisms
\[ \begin{aligned}
& \opn{RHom}_{A_2} \bigl( M_3, \opn{RHom}_{A_1}(M_2, M_1) \bigr)
\cong
\opn{Hom}_{A_2} \bigl( M_3, \opn{Hom}_{A_1}(P_2, I_1) \bigr)
\\
& \qquad 
\cong^{\dag} 
\opn{Hom}_{A_1} (P_2 \ot_{A_2} M_3, I_1) \cong 
\opn{RHom}_{A_1} \bigl( M_2 \ot^{\mrm{L}}_{A_2} M_3, M_1)
\end{aligned} \]
in $\dcat{D}(A_3)$. The isomorphism $\cong^{\dag}$ is due to the usual 
Hom-tensor adjunction. 
\end{proof}

\begin{rem} \label{rem:3100}
It is not hard to show that $\dcat{D}(A)$ is a closed symmetric monoidal 
category, with monoidal operation $(- \ot^{\mrm{L}}_{A} -)$ and monoidal unit 
object $A$. It needs just a bit more than Propositions \ref{prop:3100}, 
\ref{prop:3102} and \ref{prop:3106}. 

For a noncommutative version of this monoidal structure see 
Remark \ref{rem:3450} and Definition \ref{dfn:4530}.
\end{rem}

We end this subsection with a weakly-commutative variant of Theorem 
\ref{thm:4320}. Cohomologically pseudo-noether\-ian nonpositive DG rings were 
introduced in Definition \ref{dfn:3180}. 
DG modules with derived bounded below tensor displacement were defined in
Definition \ref{dfn:4341}.

\begin{thm} \label{thm:2700} 
Let $A \to B$ be a homomorphism between weakly commutative DG rings, and let
$L \in \dcat{D}(A)$ and $M, N \in \dcat{D}(B)$ be DG modules. 
There is a morphism 
\[ \opn{ev}^{\mrm{R, L}}_{L, M, N} : \opn{RHom}_{A}(L, M) \ot^{\mrm{L}}_B N \to 
\opn{RHom}_{A}(L, M \ot^{\mrm{L}}_B N) \]
in $\dcat{D}(B)$, called {\em derived tensor-evaluation}%
\index{Tensor-evaluation morphism! derived},
which is functorial in these DG modules. Moreover, if conditions 
\tup{(i)-(iii)} below hold, then $\opn{ev}^{\mrm{R, L}}_{L, M, N}$ is an 
isomorphism. 
\begin{enumerate}
\rmitem{i} The DG ring $A$ is nonpositive and cohomologically 
pseudo-noetherian, and $L \in \dcat{D}^{-}_{\mrm{f}}(A)$.

\rmitem{ii} The DG $B$-module $M$ is in $\dcat{D}^{+}(B)$.

\rmitem{iii}  The DG $B$-module $N$ has derived bounded below tensor 
displacement.
\end{enumerate}
\end{thm}

\begin{proof}
Choose a K-projective resolution 
$P \to L$ in $\dcat{C}_{\mrm{str}}(A)$ and a K-flat resolution 
$Q \to N$ in $\dcat{C}_{\mrm{str}}(B)$. 
Consider the obvious homomorphism 
\[ \opn{ev}_{P, M, Q} : \opn{Hom}_{A}(P, M) \ot_{B} Q \to 
\opn{Hom}_{A}(P, M \ot_{B} Q) \]
in $\dcat{C}_{\mrm{str}}(B)$. 
We take 
$\opn{ev}^{\mrm{R, L}}_{L, M, N} := \opn{Q}(\opn{ev}_{P, M, Q})$
in $\dcat{D}(B)$. The morphism $\opn{ev}^{\mrm{R, L}}_{L, M, N}$ does not 
depend on the choices, and hence it is functorial. 

Now assume that conditions (i)-(iii) hold. Because the restriction functor \lb 
$\dcat{D}(B) \to \dcat{D}(\K)$ is conservative, it is enough to show that 
$\opn{ev}^{\mrm{R, L}}_{L, M, N}$ is an isomorphism in $\dcat{D}(\K)$.
According to Theorem \ref{thm:3340}, $L$ is a derived pseudo-finite DG 
$A$-module (cf.\ Example \ref{exa:4471}). Now we can use Theorem 
\ref{thm:4320}.
\end{proof}

%% file: block4_190413.tex

\renewcommand{\thisfile}{block4\_190328}
 
\cleardoublepage 
\mysection{Dualizing Complexes over Commutative Rings} 
\label{sec:dual-cplx-comm-rng} 

\AYcopyright 

In this section we finally explain what was outlined, as a motivating 
discussion, in Subsection \ref{subsec:2160} of the Introduction. Dualizing 
complexes are perhaps the most compelling reason to study derived categories. In 
the commutative ring setting of the current section the technicalities are 
milder than in the geometric setting (see Remark \ref{rem:4192}) and the 
noncommutative ring setting (that is treated in Sections \ref{sec:BDC} and 
\ref{sec:rigid-DC-NC}). 

In the first and third subsections of this section we talk about {\em dualizing 
complexes} and {\em residue complexes}, respectively, over commutative rings. 
This material is based on the original treatment by A. Grothendieck in 
\cite{RD}, but with a much more detailed discussion. Sandwiched between them is 
a reminder on the  classification of injective modules, due to E. Matlis. 

The last two subsections are on {\em rigid dualizing complexes} in the sense of 
M. Van den Bergh. 

In this section all rings are commutative $\K$-rings by default (see 
Convention \ref{conv:2490}). However, the base ring $\K$ will remain implicit 
most of the time.

\mysubsection{Dualizing Complexes} \label{subsec:du-cplxs}

Let $A$ be a commutative ring. The category of $A$-modules is 
$\dcat{M}(A) = \cat{Mod} A$. Because $A$ is 
commutative, the Hom bifunctor has its target in $\dcat{M}(A)$~: 
\[ \opn{Hom}_A(-,-) : \dcat{M}(A)^{\mrm{op}} \times \dcat{M}(A)
\to \dcat{M}(A) . \]
Likewise for the right derived bifunctor: 
\[ \opn{RHom}_A(-,-) : \dcat{D}(A)^{\mrm{op}} \times \dcat{D}(A)
\to \dcat{D}(A) . \]
See Proposition \ref{prop:3105}.

Let $M \in \dcat{C}(A)$. The DG $A$-module 
$\opn{End}_A(M) = \opn{Hom}_A(M, M)$ 
is a noncommutative DG central $A$-ring; namely there is a central DG ring 
homomorphism 
\begin{equation} \label{eqn:2390}
\opn{hm}_{M} : A \to \opn{Hom}_A(M, M) 
\end{equation}
called {\em homothety}. 
When we forget the ring structure, $\opn{hm}_{M}$ becomes a homomorphism in
$\dcat{C}_{\mrm{str}}(A)$. 

\begin{dfn} \label{dfn:2392}
Given a complex $M \in \dcat{D}(A)$, the 
{\em derived homothety morphism}%
\index{Homothety morphism! derived}
\[ \opn{hm}^{\mrm{R}}_{M} : A \to \opn{RHom}_A(M, M) \]
is the morphism in $\dcat{D}(A)$ with this formula:
$\opn{hm}^{\mrm{R}}_{M} := \eta^{\mrm{R}}_{M,M} \circ \opn{Q}(\opn{hm}_{M})$.
\end{dfn}

In other words, the diagram 
\[ \UseTips \xymatrix @C=10ex @R=8ex {
A
\ar[r]^(0.35){\opn{Q}(\opn{hm}_{M})}
\ar@(u,u)[rr]^{\opn{hm}^{\mrm{R}}_{M}}
&
\opn{Hom}_A(M, M) 
\ar[r]^(0.45){\eta^{\mrm{R}}_{M,M}}
&
\opn{RHom}_A(M, M) 
} \]
in $\dcat{D}(A)$ is commutative.
The letter ``R'' in the superscript of the expression 
$\opn{hm}^{\mrm{R}}_{M}$ refers to ``right derived functor''. 

\begin{exer} \label{exer:2393}
Prove that if $\rho : M \to I$ is a K-injective resolution, then 
the diagram 
\[ \UseTips \xymatrix @C=10ex @R=6ex {
A
\ar[r]^(0.3){\opn{Q}(\opn{hm}_{I})}
\ar@(dr,l)[drr]_{\opn{hm}^{\mrm{R}}_{M}}
&
\opn{Hom}_A(I, I)
\ar[r]^{\eta^{\mrm{R}}_{I, I}}_{\cong}
&
\opn{RHom}_A(I, I)
\ar[d]^{\opn{RHom}(\opn{Q}(\rho), \opn{Q}(\rho)^{-1})}_{\cong}
\\
&
&
\opn{RHom}_A(M, M)
} \]
in $\dcat{D}(A)$ is commutative.
\end{exer}

\begin{exer} \label{exer:2165}
Formulate and prove a version of the previous exercise with a K-projective 
resolution of $M$. 
\end{exer}

\begin{dfn} \label{dfn:2394}
A complex $M \in \dcat{D}(A)$ is said to have the {\em derived Morita property}
\index{Derived Morita property} 
if the derived homothety morphism 
$\opn{hm}^{\mrm{R}}_{M} : A \to \opn{RHom}_A(M, M)$ 
in $\dcat{D}(A)$ is an isomorphism. 
\end{dfn}

\begin{prop} \label{prop:2335}
The following conditions are equivalent for a complex 
$M \in \dcat{D}(A)$~\tup{:}
\begin{itemize}
\rmitem{i} $M$ has the  derived Morita property.

\rmitem{ii} The homothety ring homomorphism 
$A \to \opn{End}_{\dcat{D}(A)}(M)$
is a bijective, and 
$\opn{Hom}_{\dcat{D}(A)} \bigl( M, M[p] \bigr) = 0$ 
for all $p \neq 0$. 

\rmitem{iii} The $A$-module 
$\opn{H}^0 \bigl( \opn{RHom}_{A}(M, M) \bigr)$ 
is free of rank $1$, with basis the element $\opn{id}_M$, and 
$\opn{H}^p \bigl( \opn{RHom}_{A}(M, M) \bigr) = 0$ 
for every $p \neq 0$.  
\end{itemize}
\end{prop}

\begin{exer}
Prove Proposition \ref{prop:2335}.
(Hint: see Corollary \ref{cor:2120} and the preceding material.) 
\end{exer}

\begin{rem} \label{rem:2390}
In some texts (e.g.\ in \cite{AIL}), a complex $M$ with the derived Morita 
property is called a {\em semi-dualizing complex}. This name is only partly 
justified, because this 
property occurs in the definition of a dualizing complex -- see Definition 
\ref{dfn:2155} below. However, there is a whole other class of 
complexes with the derived Morita property -- these are the {\em tilting 
complexes}. Often these two classes of complexes are disjoint. More on these 
notions, and their noncommutative variants, can be found in Sections 
\ref{sec:perf-tilt-NC} and \ref{sec:rigid-DC-NC} of the book. 
\end{rem}

From here on the commutative ring $A$ is assumed to be noetherian.
Recall that $\dcat{D}_{\mrm{f}}(A)$ is the full subcategory of 
$\dcat{D}(A)$ on the complexes with finitely generated cohomology modules. The 
subcategory $\dcat{D}_{\mrm{f}}(A)$ is triangulated. 

The next definition first appeared in \cite[Section V.2]{RD}. 
The injective dimension of a complex was defined in Definition \ref{dfn:2128}.

\begin{dfn} \label{dfn:2155}
Let $A$ be a noetherian commutative ring. A complex of $A$-modules $R$ is 
called a {\em dualizing complex} 
\index{Dualizing complex! commutative}
if it has the following three 
properties\tup{:} 
\begin{enumerate}
\rmitem{i} $R \in \dcat{D}^{\mrm{b}}_{\mrm{f}}(A)$. 

\rmitem{ii} $R$ has finite injective dimension. 

\rmitem{iii} $R$ has the derived Morita property.
\end{enumerate}
\end{dfn}

Recall that in the traditional literature (e.g.\ \cite{Mats}), a noetherian 
ring $A$ is called {\em regular} if all its local 
rings $A_{\p}$ are regular local rings. The {\em Krull dimension} of $A$ is the 
dimension of the scheme $\opn{Spec}(A)$; namely the supremum of the 
lengths of strictly ascending chains of prime ideals in $A$. 
In practice we never see regular rings that are not finite dimensional (there 
are only pretty exotic examples of them). 
The following definition will simplify matters for us:

\begin{dfn} \label{dfn:3075}
We shall say that a noetherian commutative ring $A$ is 
{\em regular}%
\index{Regular commutative ring} 
if it has 
finite Krull dimension, and all its local rings $A_{\p}$ are regular local 
rings, in the sense of \cite{Mats}. 
\end{dfn}

Every field $\K$, and the ring of integers $\Z$, are regular rings. 
If $A$ is regular, then so is the polynomial ring 
$A[t_1, \ldots, t_n]$ in $n < \infty$ variables, and also the 
localization of $A$ at every multiplicatively closed set. See 
\cite[Chapter 7]{Mats}. 

As proved by Serre (see \cite[Theorem 19.2]{Mats}) a regular ring $A$ has 
{\em finite global cohomological dimension}. This means that there is a number 
$d \in \N$ such that for all modules $M, N \in \dcat{M}(A)$ and all 
$q > d$, the modules $\opn{Ext}^q_A(M, N)$ vanish. This implies that every
$A$-module $M$ has injective, projective and flat dimensions $\leq d$.
It follows that every complex 
$M \in \dcat{D}^{\mrm{b}}(A)$ has finite injective, projective and flat 
dimensions (see Definitions \ref{dfn:2128} and \ref{dfn:2390}).  

\begin{exa} \label{exa:2171}
Let $A$ be a regular ring. Taking $R := A$ we see that $R$ satisfies condition 
(ii) of Definition \ref{dfn:2155}. The other two conditions hold regardless of 
the regularity of $A$. Thus the complex $R = A$ is a dualizing complex over the 
ring $A$. 

In  Subsection \ref{subsec:2160} of the Introduction we used this fact for 
the ring $A = \Z$. 
\end{exa}

\begin{dfn} \label{dfn:2156}
Given a dualizing complex $R \in \dcat{D}(A)$, the {\em duality functor}
\index{Duality functor} 
associated to it is the triangulated functor 
\[ D : \dcat{D}(A)^{\mrm{op}} \to \dcat{D}(A), \quad 
D := \opn{RHom}_A(-, R) .  \]
\end{dfn}

The notation ``$D$'' deliberately keeps the dualizing complex $R$ implicit. 
Note that upon applying the functor $D$ to the object $A \in \dcat{D}(A)$ we 
get $D(A) = R$. 

Let us choose a K-injective resolution $\rho : R \to I$. There is an 
isomorphism of triangulated functors
\begin{equation} \label{eqn:3235}
\opn{pres}_{I} : D \iso \opn{Hom}_A(-, I)
\end{equation}
from $\dcat{D}(A)^{\mrm{op}}$ to $\dcat{D}(A)$,
that we call a presentation of $D$. For every 
$M \in \dcat{D}(A)$ the diagram 
\begin{equation} \label{eqn:4195}
\UseTips \xymatrix @C=8ex @R=6ex {
D(M)
\ar[r]^(0.4){\opn{id}}
\ar[d]^{\cong}_{\opn{pres}_{I, M}}
&
\opn{RHom}_A(M, R)
\ar[d]^{\opn{RHom}(\opn{id}_M, \opn{Q}(\rho))}_{\cong}
\\
\opn{Hom}_A(M, I)
\ar[r]^{\eta^{\mrm{R}}_{M, I}}_{\cong}
&
\opn{RHom}_A(M, I)
}
\end{equation}
is commutative. (Note that we can choose $I$ to be a bounded complex of 
injectives, by Proposition \ref{prop:2115}.)

Let $M, I \in \dcat{C}(A)$. There is a homomorphism 
\begin{equation} \label{eqn:3236}
\opn{ev}_{M, I} : M \to \opn{Hom}_A \bigl( \opn{Hom}_A(M, I), I \bigr) 
\end{equation}
in $\dcat{C}_{\mrm{str}}(A)$, called {\em Hom-evaluation}, with formula 
\[ \opn{ev}_{M, N}(m)(\phi) := (-1)^{p \cd q} \cd \phi(m) \]
for $m \in M^p$ and $\phi \in \opn{Hom}_A(M, I)^q$. 

\begin{lem} \label{lem:2156} 
Let $R$ be a dualizing complex over $A$, with associated duality functor $D$. 
\begin{enumerate}
\item The functor 
$D \circ D : \dcat{D}(A) \to \dcat{D}(A)$
is triangulated.

\item There is a unique morphism  
$\opn{ev}^{\mrm{R, R}} : \opn{Id} \to D \circ D$ 
of triangulated functors from $\dcat{D}(A)$ to itself, called {\em derived 
Hom-eval\-uation}, such that for every K-injective resolution  
$\rho : R \to I$, and every complex $M \in \dcat{D}(A)$, the diagram 
\[ \UseTips \xymatrix @C=8ex @R=6ex {
M
\ar[r]^(0.4){\opn{ev}^{\mrm{R, R}}_M}
\ar[dr]_{\opn{Q}(\opn{ev}_{M, I}) \ \ }
&
D(D(M))
\ar[d]^{\cong}
\\
&
 \opn{Hom}_A \bigl( \opn{Hom}_A (M, I), I \bigr)
} \]
is commutative. Here the vertical isomorphism is a double application of the 
presentation $\opn{pres}_{I}$.

\item For $M = A$ there is equality 
$\opn{ev}^{\mrm{R, R}}_A = \opn{hm}^{\mrm{R}}_{R}$ 
of morphisms $A \to (D \circ D)(A)$ in $\dcat{D}(A)$.
\end{enumerate}
\end{lem}

\begin{proof}
(1) Choose a K-injective resolution $\rho : R \to I$. 
Let's write 
$F := \lb \opn{Hom}_A \bigl( \opn{Hom}_A (-, I), I \bigr)$.
The functor $F$ can be seen as a DG functor from $\dcat{C}(A)$ to itself, that 
is then made into a functor from $\dcat{D}(A)$ to itself. So by Theorem 
\ref{thm:1265}(1), $F$ is a triangulated functor (i.e.\ there is a 
translation isomorphism $\tau_F$ etc.). Using the presentation 
(\ref{eqn:3235}) twice we get an isomorphism
$D \circ D \cong  F$  of functors from $\dcat{D}(A)$ to itself. This makes 
$D \circ D$ triangulated. 

\medskip \noindent 
(2) With the K-injective resolution $\rho : R \to I$ above, we define 
the morphism $\opn{ev}$ so as to make the diagram commutative. According to
Theorem \ref{thm:1265}(2), $\opn{ev}$ is a morphism of triangulated functors. 
Because the resolution $\rho : R \to I$ is unique up to homotopy, 
we get the same morphism $\opn{ev}$ regardless of choice of resolution. 

\medskip \noindent 
(3) This is clear from Exercise \ref{exer:2393}. 
\end{proof}

\begin{exer} \label{exer:3245}
Suppose that $A$ is a regular ring, and we take the dualizing complex $R := A$. 
Let $P$ be a bounded complex of finitely generated projective $A$-modules. 
Show that:
\begin{enumerate}
\item There is a canonical isomorphism 
$D(P) \cong \opn{Hom}_{A}(P, A)$ 
in $\dcat{D}(A)$. 

\item There is a canonical isomorphism 
\[ (D \circ D)(P) \cong \opn{Hom}_{A} \bigl( \opn{Hom}_{A}(P, A), A \bigr) \]
in $\dcat{D}(A)$.

\item Under the isomorphism in item (2), the morphism 
$\opn{ev}^{\mrm{R, R}}_{P}$ goes to 
$\opn{Q}(\opn{ev}_{P, A})$. 

\item The homomorphism $\opn{ev}_P$ is a quasi-isomorphism. 
\end{enumerate}
\end{exer}

Here is the first important result regarding dualizing complexes. 

\begin{thm} \label{thm:2155}
Suppose $R$ is a dualizing complex over the noetherian commutative ring $A$, 
with associated duality functor $D$. Then for every complex 
$M \in \dcat{D}_{\mrm{f}}(A)$ the following hold\tup{:}
\begin{enumerate}
\item The complex $D(M)$ belongs to $\dcat{D}_{\mrm{f}}(A)$. 

\item The morphism 
\[ \opn{ev}^{\mrm{R, R}}_M : M \to D(D(M)) \]
in $\dcat{D}(A)$ is an isomorphism. 
\end{enumerate}
\end{thm}

\begin{proof}
(1) Condition (ii) of Definition \ref{dfn:2155} says that the functor $D$ has 
finite cohomological dimension. Condition (i) says that 
$D(A) \in \dcat{D}_{\mrm{f}}(A)$.
The assertion follows from Theorem \ref{thm:2160}, with 
$\cat{N}_0 := \dcat{M}_{\mrm{f}}(A)$. 

\medskip \noindent
(2) The composition $D \circ D$ is a triangulated functor with finite 
cohomological dimension (at most twice the injective dimension of $R$). The 
cohomological dimension of the identity functor
$\opn{Id}$ is $0$ (if $A \neq 0$). 
By condition (iii) of Definition \ref{dfn:2155} we know that 
$\opn{ev}_A$ is an isomorphism. 
Now we can use Theorem \ref{thm:2135}. 
\end{proof}

\begin{cor} \label{cor:2155}
Under the assumptions of Theorem \tup{\ref{thm:2155}}, let 
$\star$ be one of the boundedness conditions $\mrm{b}$, $+$, $-$ or 
$\bra{\mrm{empty}}$, and let  $-\star$ be the reversed boundedness condition.
Then the functor 
\[ D : \dcat{D}^{\star}_{\mrm{f}}(A)^{\mrm{op}} \to 
\dcat{D}^{-\star}_{\mrm{f}}(A) \]
is an equivalence of triangulated categories. 
\end{cor}

\begin{proof}
The previous theorem tells us that $D$ is its own quasi-inverse. The claim 
about the boundedness holds because $D$ has finite cohomological dimension.
\end{proof}

We saw that dualizing complexes exist over regular rings. This fact is used 
for the very general existence result Theorem \ref{thm:2170}. But first we need 
some preparation. 

Here are some adjunction properties related to homomorphisms between 
commutative rings. These are enhancements of the material from Subsection 
\ref{subsec:adjs-comm}.

Recall that a ring homomorphism $u : A \to B$ gives rise to a forgetful functor 
$\opn{Rest}_u : \dcat{D}(B) \to \dcat{D}(A)$.
This functor is going to be implicit in the discussion below. 

\begin{lem} \label{lem:2170}
Let $A \to B$ be a ring homomorphism. 
\begin{enumerate}
\item If $I \in \dcat{C}(A)$ is K-injective, then 
$J := \opn{Hom}_A(B, I) \in \dcat{C}(B)$ is K-injective. 

\item Given $M \in \dcat{D}(A)$, let us define 
$N := \opn{RHom}_A(B, M) \in \dcat{D}(B)$.
Then there is an isomorphism 
$\opn{RHom}_B(-, N) \cong \opn{RHom}_A(-, M)$
of triangulated functors
$\dcat{D}(B)^{\mrm{op}} \to \dcat{D}(B)$.
\end{enumerate}
\end{lem}

\begin{proof}
(1) This is an adjunction calculation. Suppose $L \in \dcat{C}(B)$ is acyclic. 
There are isomorphisms  
\begin{equation} \label{eqn:2182}
\opn{Hom}_B(L,J) \cong 
\opn{Hom}_B \bigl( L, \opn{Hom}_A(B, I) \bigr) \cong  
\opn{Hom}_A(L,I) 
\end{equation}
in $\dcat{C}(B)$. Since $I$ is K-injective over $A$, the complex
$\opn{Hom}_A(L,I)$ is acyclic. 

\medskip \noindent
(2) Choose a K-injective resolution $M \to I$ in $\dcat{C}(A)$. 
Let $J := \opn{Hom}_A(B, I)$. Then $N \to J$ is a K-injective 
resolution in 
$\dcat{C}(B)$. There are isomorphisms of triangulated functors 
\begin{equation} \label{eqn:2180}
\opn{RHom}_A(-, M) \cong \opn{Hom}_A(-, I) 
\end{equation}
and 
\begin{equation} \label{eqn:2181}
\opn{RHom}_B(-, N) \cong \opn{Hom}_B(-, J) ,
\end{equation}
where the functors (\ref{eqn:2180}) are contravariant functors from 
$\dcat{D}(A)$ to itself, and the functors (\ref{eqn:2181}) are contravariant 
functors from $\dcat{D}(B)$ to itself. But given $L \in \dcat{C}(B)$
we can view $\opn{Hom}_A(L, I)$ as a complex of $B$-modules, and in this way 
the functors (\ref{eqn:2180}) become contravariant triangulated functors
from $\dcat{D}(B)$ to itself. Formula (\ref{eqn:2182}) shows that the functors 
(\ref{eqn:2180}) and (\ref{eqn:2181}) are isomorphic.
\end{proof}

\begin{lem} \label{lem:2171}
Let $A \to B$ be a flat ring homomorphism, let 
$M \in \dcat{D}^{-}_{\mrm{f}}(A)$, and let 
$N \in \dcat{D}^{+}(A)$. Then there is an isomorphism 
\[ \opn{RHom}_A(M, N) \ot_A B \cong 
\opn{RHom}_B(B \ot_A M, B \ot_A N) \]
in $\dcat{D}(B)$. This isomorphism is functorial in $M$ and $N$. 
\end{lem}

\begin{proof}
Theorem \ref{thm:2700} tells us that 
\[  \opn{RHom}_A(M, N) \ot_A B \cong 
\opn{RHom}_A(M, B \ot_A N) \]
in $\dcat{D}(B)$. Then, by forward adjunction (Theorem \ref{thm:3085})
we get 
\[ \opn{RHom}_A(M, B \ot_A N) \cong \opn{RHom}_B(B \ot_A M, B \ot_A N) . \]
The functoriality is clear. 
\end{proof}

\begin{lem} \label{lem:2185} 
Let $I$ be an $A$-module. The following conditions are equivalent\tup{:}
\begin{enumerate}
\rmitem{i} $I$ is injective.

\rmitem{ii} For every finitely generated $A$-module $M$ the module 
$\opn{Ext}^1_A(M, I)$ is zero.  
\end{enumerate}
\end{lem}

\begin{exer} \label{exer:2185}
Prove Lemma \ref{lem:2185}. (Hint: use the Baer criterion, Theorem 
\ref{thm:151}.) 
\end{exer}

\begin{lem} \label{lem:2172} 
The injective dimension of a complex $N \in \dcat{D}(A)$ equals the 
cohomological dimension of the functor
\[ \opn{RHom}_{A}(-, N)|_{\dcat{M}_{\mrm{f}}(A)^{\mrm{op}}} : 
\dcat{M}_{\mrm{f}}(A)^{\mrm{op}} \to \dcat{D}(A) . \]
\end{lem}

\begin{proof}
By definition the injective dimension of $N$, say $d$,
is the cohomological dimension of the functor
$\opn{RHom}_{A}(-, N) : \dcat{D}(A)^{\mrm{op}} \to \dcat{D}(A)$.
Let $d'$ be  the cohomological dimension of the functor
$\opn{RHom}_{A}(-, N)|_{\dcat{M}_{\mrm{f}}(A)^{\mrm{op}}}$. 
Obviously the inequality $d \geq d'$ holds. 
For the reverse inequality we may assume that $\opn{H}(N)$ is nonzero
and $d' < \infty$. This implies that there are integers $q_1 = q_0 + d'$
such that for every $M \in \dcat{M}_{\mrm{f}}(A)$ there is an inclusion
\[ \opn{con} \bigl( \opn{H} \bigl( \opn{RHom}_{A}(M, N) \bigr) \bigr) \sub 
[q_0, q_1] . \]
In particular, for $M = A$, we get 
$\opn{con}(\opn{H}(N)) \sub [q_0, q_1]$. 
Let $N \to J$ be an injective resolution in 
$\dcat{C}(A)$ with $\opn{inf}(J) \geq q_0$.  
Take $I := \opn{smt}^{\leq q_1}(J)$, the smart truncation from
Definition \ref{dfn:2320}.
The proof of Proposition \ref{prop:2115}, plus Lemma \ref{lem:2185}, 
show that $N \to I$ is an injective resolution. But then 
$\opn{RHom}_{A}(-, N) \cong \opn{Hom}_{A}(-, I)$,
so this functor has cohomological displacement in the interval 
$[q_0, q_1]$, that has length $d'$. 
\end{proof}

Recall that a ring homomorphism $A \to B$ is called {\em finite} if it makes
$B$ into a finitely generated $A$-module.  

\begin{prop} \label{prop:2200}
Let $A \to B$ be a finite ring homomorphism, and let $R_A$ be a dualizing 
complex over $A$. Then the complex 
\[ R_B := \opn{RHom}_{A}(B, R_A) \in \dcat{D}(B) \]
is a dualizing complex over $B$. 
\end{prop}

\begin{proof}
Consider the duality functors
$D_{A} := \opn{RHom}_{A}(-, R_{A})$
and $D_{B} := \lb \opn{RHom}_B(-, R_B)$.
As explained in the proof of Lemma \ref{lem:2170}(2),
they are isomorphic as contravariant triangulated functors from $\dcat{D}(B)$ 
to itself. 
Since $R_{B} = D_A(B)$ and $B \in \dcat{D}^{\mrm{b}}_{\mrm{f}}(A)$,
by Corollary \ref{cor:2155} we have 
$R_B \in \dcat{D}^{\mrm{b}}_{\mrm{f}}(A)$. 
But then also $R_B \in \dcat{D}^{\mrm{b}}_{\mrm{f}}(B)$.
Next, because $D_B(L) \cong D_A(L)$ for all $L \in \dcat{D}(B)$, this implies 
that the cohomological dimension of $D_B$ is at most that of 
$D_A$, which is finite. We see that the injective dimension of the complex 
$R_B$ is finite. Lastly, there is an isomorphism
$D_B \circ D_B \cong D_A \circ D_A$
as functors from $\dcat{D}^{\mrm{b}}_{\mrm{f}}(B)$ to itself, 
and hence 
$\opn{ev}^{\mrm{R, R}} : \opn{Id} \to D_B \circ D_B$ is an isomorphism.
Applying this to the object 
$B \in \dcat{D}^{\mrm{b}}_{\mrm{f}}(B)$
we see that 
$\opn{hm}^{\mrm{R}}_{R_B} = \opn{ev}^{\mrm{R, R}}_{B} : B \to 
(D_B \circ D_{B})(B)$  
is an isomorphism. So $R_B$ has the derived Morita property.
The conclusion is that $R_B$ is a dualizing complex over 
$B$. 
\end{proof}

Recall that a ring homomorphism $A \to B$ is called a {\em localization} if $B$ 
is isomorphic, as an $A$-ring, to the localization $A_S$ of $A$ with respect to 
some multiplicatively closed subset $S$. 

\begin{prop} \label{prop:2201}
Let $A \to B$ be a localization ring homomorphism, and let $R_A$ be a 
dualizing complex over $A$. Then the complex 
\[ R_B := B \ot_A R_A \in \dcat{D}(B) \]
is a dualizing complex over $B$. 
\end{prop}

\begin{proof}
It is clear that 
$R_B \in \dcat{D}^{\mrm{b}}_{\mrm{f}}(B)$.
By Lemma \ref{lem:2172}, to compute the injective dimension of 
$R_B$ it is enough to look at $\opn{RHom}_B(M, R_B)$
for $M \in \dcat{M}_{\mrm{f}}(B)$. 
We can find a finitely generated $A$-submodule $M' \sub M$ such that 
$B \cd M' = M$; and then $M \cong B \ot_{A} M'$. Lemma \ref{lem:2171} 
tells us that 
\[ \opn{RHom}_B(M, R_B) \cong 
\opn{RHom}_{A}(M', R_{A}) \ot_{A} B . \]
We conclude that the injective dimension of $R_B$ is at most that of 
$R_{A}$, which  is finite. Lastly, by the same lemma we get an isomorphism 
\[ \opn{RHom}_B(R_B, R_B) \cong 
\opn{RHom}_{A}(R_{A}, R_{A}) \ot_{A} B  , \]
and it is compatible with the morphisms from $B$. Thus $R_B$ has the derived 
Morita property.  
\end{proof}

Recall that a ring homomorphism $A \to B$ is called {\em finite type} if 
$B$ is finitely generated as an $A$ ring (i.e.\ $B$ is a quotient of the 
polynomial ring $A[t_1, \ldots, t_n]$ for some $n < \infty$). 

\begin{dfn} \label{dfn:3211}
A ring homomorphism $u : A \to B$ is called {\em essentially finite type}
\index{Essentially finite type}
(EFT) if it can be factored into $u = u_{\mrm{loc}} \circ u_{\mrm{ft}}$, where 
$u_{\mrm{ft}} : A \to B_{\mrm{ft}}$ is a finite type ring homomorphism, and 
$u_{\mrm{loc}} : B_{\mrm{ft}} \to B$ is a localization homomorphism.
\end{dfn}

This is illustrated in the following commutative diagram. 
\[ \UseTips \xymatrix @C=14ex @R=6ex {
A
\ar[r]_{\tup{finite type}}^{u_{\mrm{ft}}}
\ar@(u,u)[rr]^{u}
&
B_{\mrm{ft}}
\ar[r]_{\tup{localization}}^{u_{\mrm{loc}}}
&
B 
} \]

\begin{prop} \label{prop:3255}
Let $u : A \to  B$ be an EFT ring homomorphism.
\begin{enumerate}
\item Let $v : B \to C$ be another EFT ring homomorphism. Then 
$v \circ u : A \to C$ is an EFT ring homomorphism.

\item Let $A \to A'$ be a ring homomorphism, and define 
$B' := A' \ot_A B$. Then the induced ring homomorphism 
$u' : A' \to B'$ is EFT. 

\item If the ring $A$ is noetherian, then $B$ is also noetherian. 
\end{enumerate}
\end{prop}

\begin{exer} \label{exer:3256}
Prove Proposition \ref{prop:3255}.
\end{exer}

\begin{exa} \label{exa:3255}
Let $\K$ be a noetherian ring, let $X$ be a finite type $\K$-scheme, and let 
$x \in X$ be a point. Then the local ring $\OO_{X, x}$ is an EFT $\K$-ring. 
\end{exa}

\begin{thm} \label{thm:2170}
Let $\K$ be a regular noetherian commutative ring, and let $A$ be an 
essentially finite type commutative $\K$-ring. Then $A$ has a dualizing 
complex%
\index{Dualizing complex! commutative}. 
\end{thm}

See Definition \ref{dfn:3075} regarding regular rings. 

\begin{proof}
The ring homomorphism $\K \to A$ can be factored as 
$\K \to A_{\mrm{pl}} \to A_{\mrm{ft}} \to A$, 
where $A_{\mrm{pl}} = \K[t_1, \ldots, t_n]$ is a polynomial ring,
$A_{\mrm{pl}} \to A_{\mrm{ft}}$ is surjective, and $A_{\mrm{ft}} \to A$ is a 
localization. (The subscripts stand for ``polynomial'' and ``finite type'' 
respectively.) According to \cite[Theorem 19.5]{Mats} the ring $A_{\mrm{pl}}$ 
is regular; so, as shown in Example \ref{exa:2171}, the complex
$R_{\mrm{pl}} := A_{\mrm{pl}}$ is a dualizing complex over $A_{\mrm{pl}}$. 

Define 
$R_{\mrm{ft}} := \opn{RHom}_{A_{\mrm{pl}}}(A_{\mrm{ft}}, R_{\mrm{pl}}) \in 
\dcat{D}(A_{\mrm{ft}})$. 
By Proposition \ref{prop:2200} this is a 
dualizing complex over $A_{\mrm{ft}}$. 
Finally define 
$R := A \ot_{A_{\mrm{ft}}} R_{\mrm{ft}} \in \dcat{D}^{\mrm{b}}_{\mrm{f}}(A)$. 
By Proposition \ref{prop:2201} this is dualizing complex over $A$. 
\end{proof}

The proof of Theorem \ref{thm:2170} might give the impression that $A$ could 
have a lot of nonisomorphic dualizing complexes. This is not quite true, as the 
next theorem demonstrates. 

\begin{thm} \label{thm:2175}
Let $A$ be a noetherian commutative ring with connected spectrum, and let $R$ 
and $R'$ be dualizing complexes over $A$. Then there is a rank $1$ projective 
$A$-module $L$ and an integer $d$, such that 
$R' \cong R \ot_A L[d]$ in $\dcat{D}(A)$. 
\end{thm}

We need some lemmas before proving the theorem. 

\begin{lem}[K\"unneth Trick] \label{lem:2177} 
Let $M, M' \in  \dcat{D}^{-}(A)$, and let $i, i' \in \Z$ be such that 
$\opn{sup}(\opn{H}(M)) \leq i$ and $\opn{sup}(\opn{H}(M')) \leq i'$. 
Then 
\[ \opn{H}^{i + i'}(M \ot^{\mrm{L}}_A M') \cong 
\opn{H}^i(M) \ot_A \opn{H}^{i'}(M')  \]
as $A$-modules. 
\end{lem}

\begin{exer} \label{exer:2178}
Prove Lemma \ref{lem:2177}.
\end{exer}

\begin{lem}[Projective Truncation Trick] \label{lem:2189} 
Let $M \in \dcat{D}^-(A)$, with $i_1 := \lb \opn{sup}(\opn{H}(M)) \in \Z$. 
Assume the $A$-module $P := \opn{H}^{i_1}(M)$ is projective. 
Then there is an isomorphism 
$M \cong \opn{smt}^{\leq i_1 - 1}(M) \oplus P[- i_1]$  
in $\dcat{D}(A)$. 
\end{lem}

\begin{exer} \label{exer:2189}
Prove Lemma \ref{lem:2189}. (Hint: first replace $M$ with its truncation 
$\opn{smt}^{\leq i_1}(M)$. Then prove that $P$ is a direct summand of 
$M^{i_1}$.) 
\end{exer}

By a {\em principal open set} in the affine scheme in $\opn{Spec}(A)$ we mean a 
set of the form $\opn{Spec}(A_s)$, where $A_s$ is the localization of $A$ at an 
element $s \in A$. Thus 
\[ \opn{Spec}(A_s) = \{ \p \in \opn{Spec}(A) \mid s \notin \p \} . \]

\begin{lem} \label{lem:2186} 
Let $M, M' \in \dcat{M}_{\mrm{f}}(A)$, and let $\p \sub A$ be a prime ideal. 
\begin{enumerate}
\item If $M_{\p} \neq 0$ and $M'_{\p} \neq 0$ then 
$M_{\p} \ot_{A_{\p}} M'_{\p} \neq 0$. 

\item If $M_{\p} \ot_{A_{\p}} M'_{\p} \cong A_{\p}$ then 
$M_{\p} \cong M'_{\p} \cong A_{\p}$.

\item If $M_{\p} \cong A_{\p}$, then there is a principal open neighborhood 
$\opn{Spec}(A_s)$ of $\p$ in $\opn{Spec}(A)$ such that 
$M_s \cong A_s$ as $A_s$-modules. 
\end{enumerate}
\end{lem}

\begin{exer} \label{exer:2187}
Prove Lemma \ref{lem:2186}. (Hint: use the Nakayama Lemma.) 
\end{exer}

Here is a pretty difficult technical lemma. 

\begin{lem} \label{lem:2175}
In the situation of the theorem, let 
$M, M' \in \dcat{D}^{-}_{\mrm{f}}(A)$
satisfy \lb 
$M \ot^{\mrm{L}}_A M' \cong A$ in $\dcat{D}(A)$. 
Then $M \cong L[d]$ in $\dcat{D}(A)$ for some rank $1$ projective $A$-module 
$L$ and an integer $d$.
\end{lem}

\begin{proof}
For each prime $\p \sub A$ let $M_{\p} := A_{\p} \ot_A M$, and define 
$e_{\p} := \opn{sup}(\opn{H}(M_{\p}))$ in $\Z \cup \{ -\infty \}$.
Define the generalized integer $e'_{\p}$ similarly. 

Fix one prime $\p$. 
The associativity and symmetry of the left derived tensor product imply the 
existence of these isomorphisms
\begin{equation} \label{2189}
\begin{aligned}
& M_{\p} \ot^{\mrm{L}}_{A_{\p}} M'_{\p} = 
(A_{\p} \ot^{\mrm{L}}_{A} M) \ot^{\mrm{L}}_{A_{\p}}
(A_{\p} \ot^{\mrm{L}}_{A} M')
\\
& \qquad 
\cong A_{\p} \ot^{\mrm{L}}_{A} (M \ot^{\mrm{L}}_A M') 
\cong A_{\p} \ot^{\mrm{L}}_{A} A \cong A_{\p} 
\end{aligned}
\end{equation}
in $\dcat{D}(A_{\p})$. Since $A_{\p} \neq 0$, it follows that
$\opn{H}(M_{\p}) \neq 0$ and $\opn{H}(M'_{\p}) \neq 0$.
So $e_{\p}, e'_{\p} \in \Z$, and 
$\opn{H}^{e_{\p}}(M_{\p})$, $\opn{H}^{e'_{\p}}(M'_{\p})$ are nonzero finitely 
generated $A_{\p}$-modules. By Lemma \ref{lem:2186}(1) we know that 
\[ \opn{H}^{e_{\p}}(M_{\p}) \ot_{A_{\p}} \opn{H}^{e'_{\p}}(M'_{\p}) \neq 0 . \]
According to Lemma \ref{lem:2177} we have 
\[ \opn{H}^{e_{\p}}(M_{\p}) \ot_{A_{\p}} \opn{H}^{e'_{\p}}(M'_{\p}) \cong 
\opn{H}^{(e_{\p} + e'_{\p})}(M_{\p} \ot^{\mrm{L}}_{A_{\p}} M'_{\p}) \cong 
\opn{H}^{(e_{\p} + e'_{\p})}(A_{\p}) . \]
But $A_{\p}$ is concentrated in degree $0$; this forces 
$e_{\p} + e'_{\p} = 0$ and 
\[ \opn{H}^{e_{\p}}(M_{\p}) \ot_{A_{\p}} \opn{H}^{e'_{\p}}(M'_{\p}) \cong 
A_{\p} \]
in $\dcat{D}(A_{\p})$. 
By Lemma \ref{lem:2186}(2) we now see that
$\opn{H}^{e_{\p}}(M_{\p}) \cong \opn{H}^{e'_{\p}}(M'_{\p}) \cong A_{\p}$. 
According to Lemma \ref{lem:2189} there are isomorphisms 
\begin{equation} \label{eqn:2196}
M_{\p} \cong A_{\p}[-e_{\p}] \oplus \opn{smt}^{\leq e_{\p} - 1}(M_{\p}) 
\end{equation}
and 
\[ M'_{\p} \cong A_{\p}[-e'_{\p}] \oplus \opn{smt}^{\leq e'_{\p} - 1}(M'_{\p})  
\]
in $\dcat{D}(A_{\p})$. 
These, with (\ref{2189}), give an isomorphism
\begin{equation} \label{eqn:2195}
\bigl( A_{\p}[-e_{\p}] \oplus \opn{smt}^{\leq e_{\p} - 1}(M_{\p}) \bigr) 
\ot^{\mrm{L}}_{A_{\p}}
\bigl( A_{\p}[-e'_{\p}] \oplus \opn{smt}^{\leq e'_{\p} - 1}(M'_{\p}) \bigr)
\cong A_{\p} .
\end{equation}
The left side of (\ref{eqn:2195}) is the direct sum of four objects. 
Passing to the cohomology of (\ref{eqn:2195}) we see that 
$N := \opn{H} \bigl( \opn{smt}^{\leq e_{\p} - 1}(M_{\p})[-e'_{\p}] \bigr)$ 
is a direct summand of $A_{\p}$.
But, since $e'_{\p} + e_{\p} = 0$, the graded module $N$ is concentrated in the 
degree interval $[\infty,  -1]$. It follows that $N = 0$. Therefore, by 
(\ref{eqn:2196}), we deduce that 
\begin{equation} \label{eqn:2197}
M_{\p} \cong A_{\p}[-e_{\p}] . 
\end{equation}

The calculation above works for every prime $\p$. 
From (\ref{eqn:2197}) we get 
\begin{equation} \label{eqn:2189}
A_{\p} \ot_A \opn{H}^i(M) \cong \opn{H}^i(M_{\p}) \cong 
\begin{cases}
A_{\p} & \tup{if} \ i = e_{\p} , 
\\
0 & \tup{otherwise} . 
\end{cases}
\end{equation}

We now use Lemma \ref{lem:2186}(3) to deduce that for every prime $\p$ there is 
an open neighborhood $U_{\p}$ of $\p$ in $\opn{Spec}(A)$ such that 
$\opn{H}^{e_{\p}}(M_{\q}) \cong A_{\q}$ for all $\q \in U_{\p}$.
This implies, by equation (\ref{eqn:2189}), that 
$e_{\q} = e_{\p}$.  Therefore $\p \mapsto e_{\p}$ is a locally constant function
$\opn{Spec}(A) \to \Z$. We assumed that $\opn{Spec}(A)$ is connected, and this 
implies that this is a constant function, say $e_{\p} = -d$ for some integer 
$d$. 

Define $L := \opn{H}^{-d}(M) \in \dcat{M}_{\mrm{f}}(A)$. 
Using truncation we see that $M \cong L[d]$ in 
$\dcat{D}(A)$. We know that $L_{\p} \cong A_{\p}$ for all primes $\p$. 
Finally, Lemma \ref{lem:2171} says that the $A$-module $L$ is projective. 
\end{proof}

\begin{rem} \label{rem:2185}
Lemma \ref{lem:2175} is actually true in much greater generality: 
the ring $A$ does not have to be noetherian, and we do not have to assume 
that the complexes $M$ and $M'$ have bounded above or finite cohomology. 
The proof is harder (see the proof of Theorem \ref{thm:3585}).
\end{rem}

\begin{proof}[Proof of Theorem \tup{\ref{thm:2175}}]
Define the duality functors
$D := \opn{RHom}_{A}(-, R)$ and \lb $D' := \opn{RHom}_{A}(-, R')$; these are 
finite 
dimensional contravariant triangulated functors from $\dcat{D}_{\mrm{f}}(A)$ to 
itself. And define $F := D' \circ D$ and $F' := D \circ D'$, that are finite 
dimensional (covariant) triangulated functors from $\dcat{D}_{\mrm{f}}(A)$ to 
itself. Let
\begin{equation} \label{eqn:2177}
M := F(A) = D'(D(A)) = \opn{RHom}_{A}(R, R')
\end{equation}
and
\[ M' := F'(A) = D(D'(A)) = \opn{RHom}_{A}(R', R) . \]
These are objects of $\dcat{D}^{\mrm{b}}_{\mrm{f}}(A)$. 

For every object $N \in \dcat{D}(A)$ there is a morphism 
\[ \psi_N : N \ot^{\mrm{L}}_A \opn{RHom}_{A}(R, R') \to 
\opn{RHom}_{A} \bigl( \opn{RHom}_{A} (N, R) , R' \bigr) \]
defined as follows: we choose a K-projective resolution $P \to N$
and a K-injective resolution $R' \to I'$. Then $\psi_N$ is represented by the 
obvious homomorphism of complexes 
\[ P \ot^{}_A \opn{Hom}_{A}(R, I') \to 
\opn{Hom}_{A} \bigl( \opn{Hom}_{A} (P, R) , I' \bigr) . \]
As $N$ changes, $\psi_N$ is a morphism of triangulated functors
$\psi : (-) \ot^{\mrm{L}}_A M \to D' \circ D  = F$.
For $N = A$ the morphism $\psi_A$ is an isomorphism, by equation 
(\ref{eqn:2177}). The functor $F$ has finite cohomological dimension, and the 
functor $(-) \ot^{\mrm{L}}_A M$ has bounded above cohomological displacement. 
According to Theorem \ref{thm:2135}, the morphism 
$\psi_N$ is an isomorphism for every $N \in \dcat{D}^-_{\mrm{f}}(A)$. 
In particular this is true for $N := M'$. So, using Theorem \ref{thm:2155}, 
we obtain 
\[ M' \ot^{\mrm{L}}_A M \cong (D' \circ D)(M') \cong 
(D' \circ D \circ D \circ D')(A) \cong A . \]
According to Lemma \ref{lem:2175} there is an isomorphism 
$M \cong L[d]$. Finally, using the isomorphism $\psi_R$, we get 
\[ R \ot_A L[d] \cong F(R) = D'(D(R)) \cong D'(A) = R' . \qedhere \]
\end{proof}

What if $\opn{Spec}(A)$ has more than one connected component?

\begin{dfn} \label{dfn:3210}
Let $A$ be a ring. 
\begin{enumerate}
\item We say that $A$ has a {\em finite connected component decomposition}
if 
\[ \opn{Spec}(A) = \coprod_{i = 1}^{m} \, \opn{Spec}(A_i) , \] 
a finite disjoint union of nonempty connected closed subschemes. 

\item If $A$ has a finite connected component decomposition as above, then 
the 
{\em connected component decomposition}%
\index{Connected component decomposition}
of $A$ is the ring isomorphism 
$A \cong \prod_{i = 1}^{m} \, A_i$.
\end{enumerate}
\end{dfn}

Of course, if $A$ has a finite connected component decomposition, then 
this decomposition is unique, up to renumbering. 

\begin{prop} \label{prop:5080}
If $A$ is noetherian, then it has a finite connected component decomposition.
\end{prop}

\begin{exer} \label{exer:5080} \mbox{}
\begin{enumerate}
\item Prove the proposition above. 

\item Find a ring $A$ that does not have a finite connected component 
decomposition.
\end{enumerate}
\end{exer}

\begin{cor} \label{cor:2175}
Let $R$ and $R'$ be dualizing complexes over $A$, and let 
$A = \prod_{i = 1}^m A_i$ be the connected component decomposition of $A$.
Then there is an isomorphism 
\[ R' \cong R \ot_A \bigl( L_1[d_1] \oplus \cdots \oplus L_m[d_m] \bigr) \]
in $\dcat{D}(A)$, where each $L_i$ is a rank $1$ projective 
$A_i$-module, and each $d_i$ is an integer. Furthermore, the modules $L_i$ are 
unique up to isomorphism, and the integers $d_i$ are unique. 
\end{cor}

\begin{exer} \label{exer:2176}
Prove Corollary \ref{cor:2175}.
\end{exer}

\begin{rem} \label{rem:2195}
A rank $1$ projective $A$-module $L$ is also called an {\em invertible 
$A$-module}. This is because $L$ is invertible for the tensor product. 
Recall that the group of isomorphism classes of invertible $A$-modules is 
the {\em commutative Picard group} $\opn{Pic}_A(A)$. 

The {\em commutative derived Picard group} $\opn{DPic}_A(A)$ is the abelian 
group $\opn{Pic}_A(A) \times \Z^m$ that classifies dualizing complexes over 
$A$, as in Corollary \ref{cor:2175}.

Now assume that $A$ is {\em noncommutative}, and flat central over a 
commutative ring $\K$. There are noncommutative versions of dualizing complexes 
and of ``invertible'' complexes, that are called {\em tilting complexes}. The 
latter form the nonabelian group $\opn{DPic}_{\K}(A)$, and it classifies 
noncommutative dualizing complexes. See
\cite{Ric1}, \cite{Ric2}, \cite{Kel1}, \cite{Ye4} and \cite{RoZi}. 
We shall study this material in Sections \ref{sec:perf-tilt-NC} 
and \ref{sec:rigid-DC-NC} of the book. 
\end{rem}

\begin{rem} \label{rem:2196}
The lack of uniqueness of dualizing complexes has always been a source of 
difficulty. A certain uniqueness or functoriality is needed, already for 
proving existence of dualizing complexes on schemes. 

In \cite{RD} Grothendieck utilized local and global duality in 
order to formulate a suitable uniqueness of dualizing 
complexes. This approach was very cumbersome (even without providing details!).

Since then there have been a few approaches in the literature to attack this 
difficulty. Generally speaking, these approaches came in two flavors:
\begin{itemize}
\item {\em Representability}. This started with P. Deligne's Appendix to 
\cite{RD}, and continued most notably in the work of A. Neeman, J. Lipman and 
their coauthors. See \cite{Ne2}, \cite{Li2} and their references. 

\item {\em Explicit Constructions}. Mostly in the early work of Lipman et 
al., including \cite{Li1} and \cite{LNS}, and in the work of A. Yekutieli 
\cite{Ye2}, and \cite{Ye3} and \cite{Ye5}. 
\end{itemize}

In Subsection \ref{subsec:rig-res-cplx-rng} of the book we will present {\em 
rigid residue complexes}, for which there is a built-in uniqueness, and even 
functoriality (see Remark \ref{rem:4191}). 
\end{rem}

\newpage

\mysubsection{Interlude: The Matlis Classification of Injective Modules}
\label{subsec:mr-inj-res}

We start with a few facts about injective modules over rings that are neither 
commutative nor noetherian. Sources for this material are \cite{Rot} and
\cite{Lam}. 

\begin{dfn} \label{dfn:2205} \mbox{}
\begin{enumerate}
\item Let $M$ be an $A$-module. A submodule $N \sub M$ is called an 
{\em essential submodule} 
\index{Essential! submodule}
if for every nonzero submodule $L \sub M$, the 
intersection $N \cap L$ is nonzero. In this case we also say that $M$ is an 
{\em essential extension} of $N$. 

\item An {\em essential monomorphism}
\index{Essential! monomorphism} 
is a monomorphism $\phi : N \inj M$
whose image is an essential submodule of $M$.   

\item Let $M$ be an $A$-module. An {\em injective hull}
\index{Injective hull}
(or {\em injective 
envelope}) of $M$ is an injective module $I$, together with an essential 
monomorphism $M \inj I$. 
\end{enumerate}
\end{dfn}

\begin{prop} \label{prop:2225}
Every $A$-module $M$ admits an injective hull. 
\end{prop}

\begin{proof}
See \cite[Theorem 3.30]{Rot} or \cite[Section 3.D]{Lam}. 
\end{proof}

There is a weak uniqueness result for injective hulls.

\begin{prop} \label{prop:2207}
Let $M$ be an $A$-module, and suppose $\phi : M \inj I$ and 
$\phi' : M \inj I'$ are monomorphisms into injective modules. 
\begin{enumerate}
\item If $\phi$ is essential, then there is a monomorphism 
$\psi : I \inj I'$ such that $\psi \circ \phi = \phi'$.

\item If $\phi'$ is also essential, then $\psi$ above is an isomorphism.
\end{enumerate}
\end{prop}

\begin{exer} \label{exer:2207}
Prove Proposition \ref{prop:2207}.
\end{exer}

In classical homological algebra we talk about the minimal injective resolution 
of a module $M$. Let us recall it. We start by taking the injective hull
$\phi : M \inj I^0$. This gives an exact sequence 
\[ 0 \to M \xar{\phi} I^0 \to M^1 \to 0 , \]
where $M^1$ is the cokernel of $\phi$. Then we take the injective hull
$M^1 \inj I^1$, and this gives a longer exact sequence 
\[ 0 \to M \to I^0 \to I^1 \to M^2 \to 0 , \]
and so on. We want to generalize this idea to complexes. 

\begin{dfn} \label{dfn:2210}
\mbox{}
\begin{enumerate}
\item A {\em minimal complex of injective $A$-modules}%
\index{Complex in abelian category! minimal injective}
is a bounded below complex of injective modules $I$, such that for every integer 
$q$ the submodule of cocycles $\opn{Z}^q(I) \sub I^q$ is essential. 

\item Let $M \in \dcat{D}^+(A)$. A 
{\em minimal injective resolution}%
\index{Resolution! minimal injective}
of $M$ is a quasi-iso\-morphism $M \to I$ into a minimal complex of injectives 
$I$.
\end{enumerate}
\end{dfn}

\begin{prop} \label{prop:2210}
Let $M \in \dcat{D}^+(A)$.
\begin{enumerate}
\item  There exists a minimal injective resolution $\phi : M \to I$.

\item If $\phi' : M \to I'$ is another minimal injective resolution, then there 
is an isomorphism $\psi : I \to I'$ in $\dcat{C}_{\mrm{str}}(A)$ such that 
$\phi' = \psi \circ \phi$.

\item If $M$ has finite injective dimension, then it has a bounded minimal 
injective resolution $I$. 
\end{enumerate}
\end{prop}

\begin{proof}
(1) We know that there is a quasi-isomorphism $M \to J$ where $J$ is a bounded 
below complex of injective modules. For every $q$ let $E^q$ be an injective 
hull of $\opn{Z}^q(J)$. By Proposition \ref{prop:2207}(1) we can assume that 
$E^q$ sits inside $J^q$ like this: $\opn{Z}^q(J) \sub E^q \sub J^q$.
Since $E^q$ is injective, we can decompose $J^q$ into a direct sum:
$J^q \cong E^q \oplus K^q$. The homomorphism 
$\d_J^q : K^q \to J^{q + 1}$ 
is a monomorphism since $K^q \cap \opn{Z}^q(J) = 0$.
And the image $\d_J^q(K^q)$ is contained in $E^{q + 1}$. 
Thus $\d_J^q(K^q)$ is a direct summand of $E^{q + 1}$, and this shows that the 
quotient
\[ I ^{q + 1} := E^{q + 1} / \d_J^q(K^q) \cong
J^{q + 1} / (K^{q + 1} \oplus \d_J^q(K^q)) \]
is an injective module. The canonical surjection of graded modules 
$\pi : J \to I$ is a homomorphism of complexes, with kernel the acyclic complex 
\[ \bigoplus_{q \in \Z} \ \Bigl( K^q[-q] \xar{\d_J^q} \d_J^q(K^q)[-q - 1] 
\Bigr) . \]
Therefore $\pi$ is a quasi-isomorphism. A short calculation shows that $I$ is a 
minimal complex of injectives, i.e.\ $\opn{Z}^q(I) \sub I^q$ is essential. 

\medskip \noindent
(2) See next exercise. (We will not need this fact.) 

\medskip \noindent
(3) According to Proposition \ref{prop:2115}, the complex $J$ that appears in 
item (1) can be chosen to be bounded. 
\end{proof}

\begin{exer} \label{exer:2215}
Prove Proposition \ref{prop:2210}(2). 
\end{exer}

\begin{rem} \label{rem:2210}
Important: the isomorphisms $\psi$ in Propositions \ref{prop:2207} and 
\ref{prop:2210} are not unique (see next exercise). We will see below (in 
Subsection \ref{subsec:rig-res-cplx-rng}) 
that a rigid residue complex is a minimal complex of injectives that has no 
nontrivial rigid automorphisms. 
\end{rem}

\begin{exer} \label{exer:2210}
Take $A := \K[[t]]$, the power series ring over a field $\K$. 
Let $M := A / (t)$, the trivial module (the residue field viewed as an 
$A$-module). 
\begin{enumerate}
\item Find a minimal injective resolution $M \to  I$. 

\item Find a nontrivial automorphism of the complex $I$ in 
$\dcat{C}_{\mrm{str}}(A)$ that fixes the submodule $M \sub I^0$. 
\end{enumerate}
\end{exer}

Now we add the noetherian condition. 

\begin{prop} \label{prop:2205}
Assume $A$ is a left noetherian ring. 
Let $\{ I_z \}_{z \in Z}$ be a collection of injective $A$-modules. Then 
$I := \bigoplus_{z \in Z} I_z$ is an injective $A$-module. 
\end{prop}

\begin{exer} \label{exer:2205}
Prove Proposition \ref{prop:2205}. (Hint: use the Baer criterion.)
\end{exer}

From here on in this subsection all rings are noetherian commutative. For 
them much more can be said.

Recall that a module $M$ is called {\em indecomposable} if it is not the 
direct sum of two nonzero modules.  

\begin{dfn} \label{dfn:2208}
Let $\a \sub A$ be an ideal.
\begin{enumerate}
\item Let $M$ be an $A$-module. The {\em 
$\a$-torsion submodule} of $M$ is the submodule $\Ga_{\a}(M)$ 
consisting of the elements that are annihilated by powers of $\a$. Thus 
\[ \Ga_{\a}(M) = \lim_{i \to} \, \opn{Hom}_A(A / \a^i, M) \sub M . \]

\item If $\Ga_{\a}(M) = M$ then $M$ is called an {\em $\a$-torsion module}. 

\item The functor 
$\Ga_{\a} : \dcat{M}(A) \to \dcat{M}(A)$
is called the {\em $\a$-torsion functor}. 
\end{enumerate}
\end{dfn}

Here are some important properties of the torsion functor.

\begin{prop} \label{prop:3085}
Let $\a$ be an ideal in $A$. 
\begin{enumerate}
\item The functor $\Ga_{\a}$ is left exact. 

\item The functor $\Ga_{\a}$ commutes with infinite direct sums. 
\end{enumerate}
\end{prop}

\begin{exer} \label{exer:3085}
Prove the proposition above. 
\end{exer}

Perhaps the most important theorem about injective modules over 
noetherian commutative rings is the following structural result due to E. 
Matlis \cite{Matl} from 1958. See also \cite[Section V.4]{Ste}, 
\cite[Sections 3.F and 3.I]{Lam}, \cite[Section 18]{Mats} and \cite{BrSh}.

For a prime ideal $\p \sub A$ we write 
$\kk(\p) := A_{\p} / \p_{\p}$, the residue field of $\p$. 

\begin{thm}[Matlis] \label{thm:2026}
Let $A$ be a noetherian commutative ring. 
\begin{enumerate}
\item Let $\p \sub A$ be a prime ideal, and let $J(\p)$ be the injective hull 
of the $A_{\p}$-module $\kk(\p)$. Then, as an $A$-module, $J(\p)$ is
injective, indecomposable and $\p$-torsion. 

\item Suppose $I$ is an indecomposable injective $A$-module. Then 
$I \cong J(\p)$ for a unique prime ideal $\p \sub A$. 

\item Every injective $A$ module $I$ is a direct sum of indecomposable 
injective $A$-modules. 
\end{enumerate}
\end{thm}

Theorem \ref{thm:2026} tells us that every injective $A$-module $I$ can be 
written as a direct sum 
\begin{equation} \label{eqn:2215}
I \cong \bigoplus_{\p \in \opn{Spec}(A)} \, J(\p)^{\oplus \mu_{\p}} 
\end{equation}
for a collection of cardinal numbers 
$\{ \mu_{\p} \}_{\p \in \opn{Spec}(A)}$, called the {\em Bass numbers}. 
General counting tricks can show that the multiplicity $\mu_{\p}$ is an 
invariant of $I$. But we can be more precise: 

\begin{prop} \label{prop:2208}
Let $I$ be an injective $A$-module, with direct sum decomposition 
\tup{(\ref{eqn:2215})}.
Then for every $\p$ there is equality
\[ \mu_{\p} = \opn{rank}_{\kk(\p)} \bigl( 
\opn{Hom}_{A_{\p}} \bigl( \kk(\p) , A_{\p} \ot_A I \bigr) \bigr) . \]
\end{prop}

\begin{proof}
Consider another prime $\q$. If $\q \nsubseteq \p$ then there is an element
$a \in \q - \p$, and then $a$ is both invertible and locally nilpotent on 
$A_{\p} \ot_A J(\q)$. This implies that $A_{\p} \ot_A J(\q) = 0$.
On the other hand, if $\q \sub \p$, then 
$A_{\p} \ot_A J(\q) \cong J(\q)$. Therefore 
$A_{\p} \ot_A I \cong  \bigoplus_{\q \sub \p} J(\p)^{\oplus \mu_{\p}}$.

Next, if $\q \varsubsetneq \p$, then there is an element $b \in \p - \q$, 
and it is both invertible and zero on the module 
$\opn{Hom}_{A_{\p}} \bigl( \kk(\p) , J(\q)  \bigr)$.
The implication is that this module is zero.

Finally, if $\q = \p$ then we have 
\[ \opn{Hom}_{A_{\p}} \bigl( \kk(\p) , J(\p)  \bigr) \cong 
\opn{Hom}_{A_{\p}} \bigl( \kk(\p) , \kk(\p)  \bigr) \cong 
\kk(\p) , \]
because the inclusion $\kk(\p) \sub J(\p)$ is essential. 

Since $\kk(\p)$ is a finitely generated $A_{\p}$-module, the functor 
$\opn{Hom}_{A_{\p}}(\kk(\p), -)$
commutes with infinite direct sums. Therefore, putting all these cases 
together, we see that 
\[ \opn{Hom}_{A_{\p}} \bigl( \kk(\p) , A_{\p} \ot_A I \bigr) \bigr) \cong 
\kk(\p)^{\oplus \mu_{\p}}  \]
as $\kk(\p)$-modules.
\end{proof}

\mysubsection{Residue Complexes} \label{subsec:res-cplx-ring}

In this subsection $A$ is a noetherian commutative ring. 
Here we introduce residue complexes (called residual complexes in \cite{RD}). 
Most of the material is taken from the original \cite{RD}. 
In Example \ref{exa:2245} we will see the relation between the geometry of 
$\opn{Spec}(A)$ and the structure of dualizing complexes over $A$ (continuing 
Example \ref{exa:2252} from the Introduction). 

\begin{lem} \label{lem:2200}
Let $R$ be a dualizing complex over $A$ and let $\p \sub A$ be a prime ideal. 
We write 
$R_{\p} := A_{\p} \ot_A R \in \dcat{D}(A_{\p})$.
There is an integer $d$ such that 
\[ \opn{Ext}^i_{A_{\p}} \bigl(\kk(\p), R_{\p} \bigr) \cong 
\begin{cases}
\kk(\p) & \tup{if} \ i = -d , 
\\
0 & \tup{otherwise} . 
\end{cases} \]
\end{lem}

\begin{proof}
By Proposition \ref{prop:2201}, $R_{\p}$ is a dualizing complex over the local 
ring $A_{\p}$. And by Proposition \ref{prop:2200},
$S := \opn{RHom}_{A_{\p}} \bigl( \kk(\p), R_{\p} \bigr)$ 
is a dualizing complex over the residue field $\kk(\p)$.  
Since $\kk(\p)$ is a field, it is a regular ring, and so it is a dualizing 
complex over itself. Theorem  \ref{thm:2175} tells us that 
$S \cong \kk(\p)[d]$ in $\dcat{D}(\kk(\p))$ 
for some integer $d$.  
\end{proof}

\begin{dfn} \label{dfn:2200}
The number $d$ in Lemma \ref{lem:2200} is called the {\em dimension of $\p$ 
relative to $R$}, and is denoted by $\opn{dim}_R(\p)$. 
In this way we obtain a function
\[  \opn{dim}_R : \opn{Spec}(A) \to  \Z , \]
called the {\em dimension function} 
\index{Dimension function}
associated to $R$. 
\end{dfn}

Let us recall a few notions regarding the combinatorics of prime ideals in a 
ring $A$. A prime ideal $\q$ is an {\em immediate specialization} of another 
prime $\p$ if $\p \subsetneqq \q$, and there is no other prime $\p'$ such that 
$\p \subsetneqq \p' \subsetneqq \q$. 
In other words, if the dimension of the local ring $A_{\q} / \p_{\q}$ is $1$. 

A {\em chain of prime ideals} in $A$ is a sequence 
$(\p_0, \ldots, \p_n)$ of primes such that 
$\p_i \subsetneqq \p_{i + 1}$ for all $i$. The number $n$ is the {\em length} 
of the chain. The chain is called {\em saturated} if for each $i$ the prime 
$\p_{i + 1}$ is an immediate specialization of $\p_i$. 

\begin{thm} \label{thm:2205}
Let $R$ be a dualizing complex over $A$ and let $\p, \q \sub A$ be prime 
ideals. Assume that $\q$ is an immediate specialization of $\p$. Then
\[ \opn{dim}_R(\q) = \opn{dim}_R(\p) - 1 . \]
\end{thm}

To prove this theorem we need a baby version of local cohomology: 
codimension $1$ only. 

Let $\a$ be an ideal in $A$. 
The torsion functor $\Ga_{\a}$ has a right derived functor
$\mrm{R} \Ga_{\a}$. For every complex $M \in \dcat{D}(A)$, the module 
$\opn{H}^p_{\a}(M) := \opn{H}^p(\mrm{R} \Ga_{\a}(M))$ 
is called the {\em $p$-th cohomology of $M$ with support in $\a$}. 
In case $A$ is a local ring and $\m$ is its maximal ideal, then 
$\opn{H}^p_{\m}(M)$ is also called the {\em local cohomology of $M$}. 

Now suppose $\a$ is a principal ideal in $A$, generated by an 
element $a$. Let $A_a  = A[a^{-1}]$ be the localized ring.
For any $A$-module $M$ we write $M_a = A_a \ot_A M$. 
There is a canonical exact sequence 
\begin{equation} \label{eqn:2206}
0 \to \Ga_{\a}(M) \to M \to M_a . 
\end{equation}

\begin{lem} \label{lem:2220}
Let $\a = (a)$ be a principal ideal in $A$. 
\begin{enumerate}
\item For every injective module $I$ the sequence 
$0 \to \Ga_{\a}(I) \to I \to I_a \to 0$
is exact. 

\item For every $M \in \dcat{D}^+(A)$ there is a long exact sequence of 
$A$-modules 
\[ \cdots \to \opn{H}^p_{\a}(M) \to \opn{H}^p(M) \to \opn{H}^p(M_a) \to 
\opn{H}^{p + 1}_{\a}(M) \to  \cdots . \]
\end{enumerate}
\end{lem}

\begin{proof}
(1) Let $J(\q)$ be an indecomposable injective $A$-module. According to Theorem 
\ref{thm:2026}(1), if $a \in \q$ then 
$\Ga_{\a}(J(\q)) = J(\q)$ and $J(\q)_a = 0$. But if $a \notin \q$
then $J(\q) = J(\q)_a$ and $\Ga_{\a}(J(\q)) = 0$.
By Theorem \ref{thm:2026} and Proposition \ref{prop:3085}(2)
we see that each injective module $I$ breaks up into 
a direct sum $I = \Ga_{\a}(I) \oplus I_a$, and this proves that the sequence 
is split exact.  

\medskip \noindent
(2) Choose a resolution $M \to I$ by a bounded below complex of injectives. 
We obtain an exact sequence of complexes as in item (1). The long 
exact sequence in cohomology
\[ \cdots \to \opn{H}^p(\Ga_{\a}(I)) \to \opn{H}^p(I) \to 
\opn{H}^p(I_a) \to \opn{H}^{p + 1}(\Ga_{\a}(I)) \to  \cdots  \]
is what we want. 
\end{proof}

\begin{lem} \label{lem:2235}
Suppose $A$ is an integral domain, with fraction field $K$, such that
$A \neq K$. Then $K$ is not a finitely generated $A$-module. 
\end{lem}

\begin{proof}
Let $a \in A$ be a nonzero element that is not invertible. Then 
\[ A \subsetneqq a^{-1} \cd A \subsetneqq a^{-2} \cd A \subsetneqq
\cdots \sub K \]
is an infinite ascending sequence of $A$-submodules of $K$. 
\end{proof}

\begin{lem} \label{lem:2221}
For every ideal $\a$ and every $M \in \dcat{D}(A)$ there is an isomorphism of 
$A$-modules 
$\opn{H}^p_{\a}(M) \cong \lim_{k \to} \opn{Ext}^p_A( A / \a^k, M)$.
\end{lem}

\begin{proof}
Choose a K-injective resolution $M \to I$. Then, using the fact that cohomology 
commutes with direct limits, we have
\[ \begin{aligned}
& \opn{H}^p_{\a}(M) \cong \opn{H}^p(\Ga_{\a}(I)) \cong 
\opn{H}^p \bigl( \lim\nolimits_{k \to} \, \opn{Hom}_A(A / \a^k, I) \bigr)
\\ & \quad \cong  
\lim_{k \to} \, \opn{H}^p \bigl( \opn{Hom}_A (A / \a^k, I) \bigr) 
\cong \lim_{k \to} \, \opn{Ext}^p_A(A / \a^k, M) .
\end{aligned} \]
\end{proof}

\begin{lem} \label{lem:2250}
Assume $A$ is local, with maximal ideal $\m$. Let $R$ be a dualizing complex 
over $A$, and let $d := \opn{dim}_R(\m)$. Then for every $i \neq -d$ 
the local cohomology $\opn{H}^{i}_{\m}(R)$ vanishes. 
\end{lem}

See Remark \ref{rem:2250} for more about $\opn{H}^{-d}_{\m}(R)$. 

\begin{proof}
We know that 
\[ \opn{Ext}^i_{A} \bigl( \kk(\m), R \bigr) \cong 
\begin{cases}
\kk(\m) & \tup{if} \ i = -d , 
\\
0 & \tup{otherwise} .
\end{cases} \]
Let $N$ be a finite length $A$-module. Since $N$ is gotten from the residue 
field $\kk(\m)$ by finitely many extensions, induction on the length of $N$ 
shows that 
$\opn{Ext}^i_{A}(N, R) = 0$ for all $i \neq -d$. 
This holds in particular for $N := A / \m^k$. 
Now use Lemma \ref{lem:2221}.
\end{proof}

\begin{proof}[Proof of Theorem \tup{\ref{thm:2205}}]
Define $d := \opn{dim}_R(\q)$ and $e := \opn{dim}_R(\p)$.
We need to prove that $e = d + 1$. 

By definition, $d$ is the unique integer s.t.\ 
$\opn{Ext}^{-d}_{A_{\q}} \bigl( \kk(\q), R_{\q} \bigr) \neq 0$.
Let's define $\bar{A} := A / \p$. 
We know that 
$\bar{R} := \opn{RHom}_{A}(\bar{A}, R)$
is a dualizing complex over $\bar{A}$. There are isomorphisms 
\[ \begin{aligned}
& \opn{RHom}_{A_{\q}} \bigl( \kk(\q), R_{\q} \bigr) \cong 
\opn{RHom}_{\bar{A}_{\q}} \bigl( \kk(\q), 
\opn{RHom}_{A_{\q}} \bigl( \bar{A}_{\q} , R_{\q} \bigr) \bigr)
\\
& \quad 
\cong \opn{RHom}_{\bar{A}_{\q}} \bigl( \kk(\q), \bar{R}_{\q} \bigr)
\end{aligned} \]
in $\dcat{D}(\bar{A}_{\q})$, coming from adjunction for the homomorphism 
$A_{\q} \to \bar{A}_{\q}$. 
There is also an isomorphism 
\[ \opn{RHom}_{A_{\p}} \bigl( \kk(\p), R_{\p} \bigr)  
\cong \opn{RHom}_{\bar{A}_{\p}} \bigl( \kk(\p), \bar{R}_{\p} \bigr) \]
in $\dcat{D}(\bar{A}_{\p})$.
Hence we can replace $A$ and $R$ with $\bar{A}_{\q}$ and $\bar{R}_{\q}$ 
respectively. 

Now we have $\p = 0$ and $A = A_{\q}$. Thus $A$ is a $1$-dimensional local 
integral domain, 
with only two primes ideals: $0 = \p$ and the maximal ideal $\q$. Take 
any nonzero element $a \in \q$. Then the localization $A_a$ is the field of 
fractions of $A$, i.e.\ $A_a = \kk(\p)$. 
On the other hand, letting $\a := (a) \sub A$, the quotient 
$A / \a$ is an artinian local ring. So $A / \a$ is a finite length 
$A$-module, the ideal $\a$ is $\q$-primary, and $\Ga_{\a} = \Ga_{\q}$. 
 
By Lemma \ref{lem:2220} there is an exact sequence of $A$-modules 
\[ \cdots \to \opn{H}^{-e}_{\a}(R) \to \opn{H}^{-e}(R) \xar{\phi} 
\opn{H}^{-e}(R_a) \to \opn{H}^{-e + 1}_{\a}(R) \to  \cdots . \]
Since $a \neq 0$ there are equalities $A_a = A_{\p} = \opn{Frac}(A) = \kk(\p)$. 
Then 
$\opn{H}^{-e}(R_a) \cong \kk(\p)$, and this is not a finitely generated 
$A$-module by Lemma \ref{lem:2235}. On the other hand the $A$-module 
$\opn{H}^{-e}(R)$ is finitely generated. 
We conclude that homomorphism $\phi$ is not surjective, and thus 
$\opn{H}^{-e + 1}_{\a}(R) \neq 0$. 
But 
$\opn{H}^{-e + 1}_{\a}(R) = \opn{H}^{-e + 1}_{\q}(R)$,
so according to Lemma \ref{lem:2250} we must have 
$-e + 1 = -d$. Thus $e = d + 1$ as claimed. 
\end{proof}

\begin{cor} \label{cor:2235}
If $A$ has a dualizing complex, then the Krull dimension of $A$ is finite. 
More precisely, if $R$ is a dualizing complex over $A$, then 
$\opn{dim}(A)$ is at most the injective dimension of $R$. 
\end{cor}

\begin{proof}
Let $[i_0, i_1]$ be the injective concentration of the complex $R$. See 
Definition \ref{dfn:2128}. This is a bounded interval. Since 
\[ \opn{Ext}^i_{A_{\p}}(\kk(\p), R_{\p}) \cong
\opn{Ext}^i_{A}(A / \p, R)_{\p} ,  \]
we see that 
$\opn{dim}_{R}(\p) \in -[i_0, i_1] = [-i_1, -i_0]$. 

Let $(\p_0, \ldots, \p_n)$ be a chain of prime ideals in $A$. Because $A$ is 
noetherian, we can squeeze more primes into this chain, until after finitely 
many steps it becomes saturated. According to Theorem \ref{thm:2205}, for a 
saturated chain we have 
$n = \opn{dim}_{R}(\p_0) - \opn{dim}_{R}(\p_n)$.
Therefore $n \leq i_1 - i_0$. 
\end{proof}

\begin{dfn} \label{dfn:2237}
The ring $A$ is called {\em catenary} if for every pair of primes 
$\p \sub \q$ there is a number $n_{\p, \q}$ such that for every saturated chain 
$(\p_0, \ldots, \p_n)$ with $\p_0 = \p$ and $\p_n = \q$, there is equality
$n = n_{\p, \q}$. 
\end{dfn}

\begin{cor} \label{cor:2236}
If $A$ has a dualizing complex, then it is catenary.
\end{cor}

\begin{proof}
Let $R$ be a dualizing complex over $A$. The proof of the previous corollary 
shows that the number 
$n_{\p, \q}= \opn{dim}_{R}(\p) - \opn{dim}_{R}(\q)$
has the desired property. 
\end{proof}

\begin{exa} \label{exa:2245}
This is a continuation of Example \ref{exa:2252} from the Introduction. 
Consider the ring 
$A = \mbb{R}[t_1, t_2, t_3] / (t_3 \cd t_1, \, t_3 \cd t_2)$.
The affine algebraic variety
$X = \opn{Spec}(A) \sub \mbf{A}^3_{\mbb{R}}$
is shown in figure \ref{fig:2250}. 
It is the union of a plane $Y$ and a line 
$Z$, meeting at the origin. 

\begin{figure}
\centering
\includegraphics[scale=0.3]{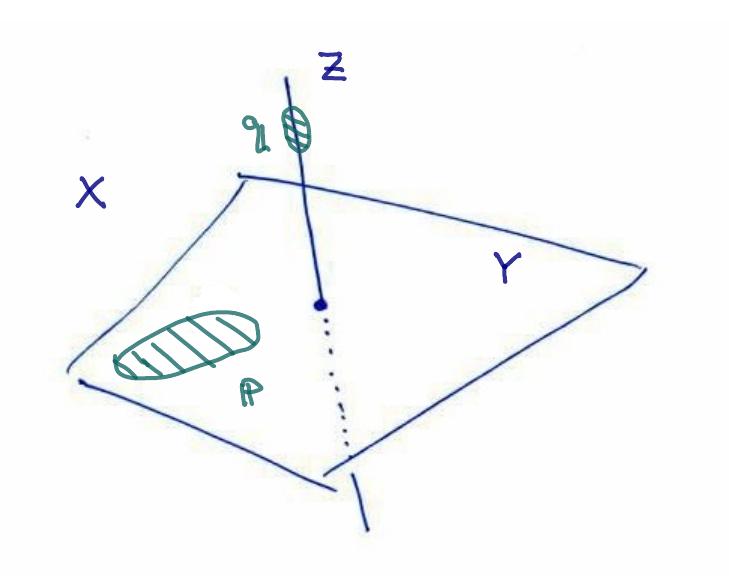}
\caption{An algebraic variety $X$ that is connected but not equidimensional: it 
has irreducible components $Y$ and $Z$ of dimensions $2$ and $1$ respectively.
The generic points $\p \in Y$ and $\q \in Z$ are shown.}
\label{fig:2250}
\end{figure}

Since the ring $A$ is finite type over the field $\R$, it has a dualizing 
complex $R$. We will now prove that there is some integer $i$ s.t.\ 
$\mrm{H}^i(R)$ and $\mrm{H}^{i+1}(R)$ are nonzero. 

Define the prime ideals 
$\m := (t_1, t_2, t_3)$, 
$\q := (t_1, t_2)$ and 
$\p := (t_3)$. Thus $\m$ is the origin, $\q$ is the generic point of the line
$Z = \opn{Spec}(A / \q)$, and $\p$ is the generic point of the plane
$Y = \opn{Spec}(A / \p)$. 
By translating $R$ as needed, we can assume that 
$\opn{dim}_R(\m) = 0$. Since $\m$ is an immediate specialization of $\q$, 
Theorem \ref{thm:2205} tells us that $\opn{dim}_R(\q) = 1$. 
Similarly, since every line in $Y$ passing through the origin gives 
rise to a saturated chain 
$(\p, \q', \m)$, we see that $\opn{dim}_R(\p) = 2$.

Since $\q$ is the generic point of $Z$, its local ring is the residue field: 
$A_{\q} = \kk(\q)$. We know that $\opn{dim}_R(\q) = 1$. Hence 
\[ \kk(\q) \cong \opn{Ext}^{-1}_{A_{\q}}(\kk(\q), R_{\q}) =
\opn{Ext}^{-1}_{A_{\q}}(A_{\q}, R_{\q}) \cong 
\opn{Ext}^{-1}_{A}(A, R)_{\q} \cong \opn{H}^{-1}(R)_{\q} . \]
Therefore $\opn{H}^{-1}(R) \neq 0$. 
A similar calculation involving $\p$ shows that 
$\opn{H}^{-2}(R) \neq 0$.
\end{exa}

\begin{exa} \label{2255}
Let $A$ be a local ring, with maximal ideal $\m$ and residue field  
$\kk(\m)$. Recall that $A$ is called {\em Gorenstein} if the free 
module $A$ has finite injective dimension. The ring $A$ is called called {\em 
Cohen-Macaulay} if its depth is equal to its dimension, where the depth of $A$ 
is the minimal integer $i$ such that 
$\opn{Ext}^i_A(\kk(\m), A) \neq 0$. It is known that Gorenstein implies 
Cohen-Macaulay. See \cite{Mats} for details. 

As is our usual practice (cf.\ Definition \ref{dfn:3075}) we shall say that 
a noetherian commutative ring $A$ is Cohen-Macaulay (resp.\ Gorenstein) if it 
has finite Krull dimension, and all its local rings $A_{\p}$ are Cohen-Macaulay 
(resp.\ Gorenstein) local rings, as defined above. 

Assume $A$ has a connected spectrum, and let $R$ be a dualizing complex over 
$A$. Grothendieck showed in \cite[Section V.9]{RD} that $A$ is a Cohen-Macaulay 
ring iff $R \cong L[d]$ for some finitely generated module $L$ and some integer 
$d$; the proof is not easy. It is however pretty easy to prove (using Theorem 
\ref{thm:2175}) that $A$ is a Gorenstein ring iff $R \cong L[d]$ for some 
invertible module $L$ and some integer $d$. 

There is a lot more to say about the relation between the CM (Cohen-Macaulay) 
property and duality. See Remark \ref{rem:2270}
\end{exa}

Recall that for each $\p \in \opn{Spec}(A)$ we denote by $J(\p)$ the 
corresponding indecomposable injective $A$-module. 

\begin{dfn} \label{dfn:2236}
A {\em residue complex} 
\index{Residue complex}
over $A$ is a complex of $A$-modules $\KK$ having these properties:
\begin{itemize}
\rmitem{i} $\KK$ is a dualizing complex. 

\rmitem{ii} For every integer $d$ there is an isomorphism of $A$-modules 
\[ \KK^{-d} \cong 
\bigoplus_{\substack{ \p \in \opn{Spec}(A) \\  
\opn{dim}_{\KK}(\p)  = d }}  J(\p) \ . \]
\end{itemize}
\end{dfn}

The reason we like residue complexes is this:

\begin{thm} \label{thm:2240}
Suppose $\KK$ and $\KK'$ are residue complexes over $A$ that have the same 
dimension function. Then the homomorphism 
\[ \opn{Q} : \opn{Hom}_{\dcat{C}_{\mrm{str}}(A)}(\KK, \KK') \to 
\opn{Hom}_{\dcat{D}(A)}(\KK, \KK') \]
is bijective. 
\end{thm}

In other words, for each morphism $\psi : \KK \to \KK'$ in 
$\dcat{D}(A)$ there is a unique strict homomorphism of complexes 
$\phi : \KK \to \KK'$ such that $\psi = \opn{Q}(\phi)$. 

\begin{proof}
Since the complex $\KK'$ is K-injective, by Theorem \ref{thm:3135}
we know that the homomorphism 
\[ \opn{Q} : \opn{Hom}_{\dcat{K}(A)}(\KK, \KK') \to 
\opn{Hom}_{\dcat{D}(A)}(\KK, \KK') \]
is bijective. And by definition the homomorphism 
\[ \opn{P} : \opn{Hom}_{\dcat{C}_{\mrm{str}}(A)}(\KK, \KK') \to 
\opn{Hom}_{\dcat{K}(A)}(\KK, \KK') \]
is surjective. It remains to prove that 
$\opn{Hom}_{A}(\KK, \KK')^{-1} = 0$,
i.e.\ here are no nonzero degree $-1$ homomorphisms $\ga : \KK \to \KK'$. 

The residue complexes $\KK$ and $\KK'$ decompose into indecomposable summands 
by the formula in property (ii) of Definition \ref{dfn:2236}. 
A homomorphism $\ga : \KK \to \KK'$ of degree $-1$ 
is nonzero iff at least one of its components 
$\ga_{\p, \q} : J(\p) \to J(\q)$ 
is nonzero, for some $J(\p) \sub \KK^{-i}$ and $J(\q) \sub \KK'^{\, -i - 1}$.
Denoting by $\opn{dim}$ the dimension function of both these dualizing 
complexes, we have $\opn{dim}(\p) = i$ and $\opn{dim}(\q) = i + 1$. 
But the lemma below says that $\q$ has to be a specialization of $\p$. 
Therefore, as in the proof of Corollary \ref{cor:2235}, there is an inequality 
in the oppose direction:  
$\opn{dim}(\p) \geq \opn{dim}(\q)$. 
We see that it is impossible to have a  nonzero degree $-1$ 
homomorphism $\ga : \KK \to \KK'$. 
\end{proof}

\begin{lem} \label{lem:2240}
Let $\p, \q$ be prime ideals. If there is a nonzero homomorphism 
$\ga : J(\p) \to J(\q)$, then 
$\q$ is a specialization of $\p$. 
\end{lem}

\begin{proof}
Assume $\q$ is not a specialization of $\p$; i.e.\ 
$\p \nsubseteq \q$. So there is an element $a \in \p - \q$. 
Let $\ga : J(\p) \to J(\q)$ be a homomorphism, and consider the module 
$N := \ga(J(\p)) \sub J(\q)$. Since $J(\p)$ is $\p$-torsion, the element $a$ 
acts on $N$ locally-nilpotently. On the other hand, $J(\q)$ is a module over 
$A_{\q}$, so $a$ acts invertibly on $J(\q)$, and hence it has zero annihilator 
in $N$. The conclusion is that $N = 0$. 
\end{proof}

Here is a general existence theorem. 
Minimal injective resolutions were defined in Definition \ref{dfn:2210}.
Their existence was proved in Proposition \ref{prop:2210}. 

\begin{thm} \label{thm:2245}
Suppose the ring $A$ has a dualizing complex $R$. 
Let $R \to \KK$ be a minimal injective resolution of $R$. Then $\KK$ is a 
residue complex over $A$. 
\end{thm}

The proof is after two lemmas. 

\begin{lem} \label{lem:2230}
Let $S \sub A$ be a multiplicatively closed set, with localization $A_S$. 
For each complex of $A$-modules $M$ we write $M_S := A_S \ot_A M$.
\begin{enumerate}
\item If $I$ is an injective $A$-module, then $I_S$ is an 
injective $A_S$-module.

\item If $I$ is an injective $A$-module and $M \sub I$ is an essential 
$A$-submodule, then $M_S \sub I_S$ is an essential $A_S$-submodule.

\item If $I$ is a minimal complex of injective $A$-modules, then 
$I_S$ is a minimal complex of injective $A_S$-modules.
\end{enumerate}
\end{lem}

\begin{proof}
(1) By Theorem \ref{thm:2026} there is a direct sum decomposition 
$I \cong I' \oplus I''$, where 
$I' \cong \bigoplus_{\p \cap S = \varnothing}  
J(\p)^{\oplus \mu_{\p}}$
and 
$I'' \cong \bigoplus_{\p \cap S \neq \varnothing} 
J(\p)^{\oplus \mu_{\p}}$.
If $\p \cap S = \varnothing$ then $J(\p) \cong J(\p)_S$ is an injective 
$A_S$-module; and if $\p \cap S \neq \varnothing$ then $J(\p)_S = 0$. 
We see that $I_S \cong I'$ is  an injective $A_S$-module.

\medskip \noindent 
(2) Denote by $\la : I \to I_S$ the canonical homomorphism. 
Under the decomposition $I \cong I' \oplus I''$ above, 
$\la|_{I'} : I' \to I_S$ is an isomorphism. 

Let $L$ be a nonzero $A_S$-submodule of $I_S$. 
Consider the $A$-submodule 
$L' := (\la|_{I'})^{-1}(L) \sub I'$;
so $\la|_{L'} : L' \to L$ is an isomorphism and $L' \neq 0$.  
Because $M \sub I$ is essential, the intersection $M \cap L'$ is 
nonzero. But then $\la(M \cap L')$ is a nonzero submodule of 
$M_S \cap L$, so $M_S \cap L$ is nonzero. 

\medskip \noindent 
(3) By part (1) the complex $I_S$ is a bounded below complex of injective 
$A_S$-modules. Exactness of localization shows that 
$\opn{Z}^n(I_S) = \opn{Z}^n(I)_S$ inside $I^n_S$; so by part (2) the inclusion 
$\opn{Z}^n(I_S) \inj I^n_S$ is essential.
\end{proof}

\begin{lem} \label{lem:2231}
Let $\a \sub A$ be an ideal, and define $B := A / \a$. 
\begin{enumerate}
\item If $I$ is an injective $A$-module, then 
$J := \opn{Hom}_A(B, I)$ is an injective $B$-module. 

\item Let $I$ and $J$ be as above. If $M \sub I$ is an essential $A$-submodule, 
then $N := \opn{Hom}_A(B, M)$ is an essential $B$-submodule of $J$.

\item If $I$ is a minimal complex of injective $A$-modules, then 
$J := \opn{Hom}_A(B, I)$ is a minimal complex of injective $B$-modules,
\end{enumerate}
\end{lem}

\begin{proof}
(1) This is immediate from adjunction. 

\medskip \noindent
(2) We identify $J$ and $N$ with the submodules of $I$ and $M$ respectively 
that are annihilated by $\a$. Let $L \sub J$ be a nonzero $B$-submodule.
Then $L$ is a nonzero $A$-submodule of $I$. 
Because $M$ is essential, the intersection $L \cap M$ is nonzero.
But $L \cap M$ is annihilated by $\a$, so it sits inside $N$, and in fact 
$L \cap M = L \cap N$. 

\medskip \noindent
(3) By part (1) the complex $J$ is a bounded below complex of injective 
$B$-modules. Left exactness of $\opn{Hom}_A(B, -)$ shows that 
$\opn{Z}^n(J) = \opn{Hom}_A(B, \opn{Z}^n(I))$ inside $J^n$; so by part (2) 
the inclusion $\opn{Z}^n(J) \inj J^n$ is essential.
\end{proof}

\begin{proof}[Proof of Theorem \tup{\ref{thm:2245}}]
Since $\KK \cong R$ in $\dcat{D}(A)$ it follows that $\KK$ is a dualizing 
complex. To show that $\KK$ has property (ii) of Definition \ref{dfn:2236} we 
have to count multiplicities. For every $\p$ and $d$ let $\mu_{\p, d}$ be the 
multiplicity of $J(\p)$ in $\KK^{-d}$, so that  
\[ \KK^{-d} \cong 
\bigoplus_{\p \in \opn{Spec}(A)}  J(\p)^{\oplus \mu_{\p, d}} . \]
We have to prove that 
\begin{equation} \label{eqn:2245}
\mu_{\p, d} = 
\begin{cases}
1 & \tup{if} \ \opn{dim}_{\KK}(\p) = d ,  
\\
0 & \tup{otherwise} . 
\end{cases}
\end{equation}
Now by Lemma \ref{lem:2230}(3) the complex $\KK_{\p} = A_{\p} \ot_A \KK$ is a  
minimal complex of injective  $A_{\p}$-modules. Because $\KK_{\p}$ is 
K-injective over $A_{\p}$ we get   
\[ \opn{Ext}^{-d}_{A_{\p}}(\kk(\p), R_{\p}) \cong
\opn{H}^{-d} \bigl( \opn{Hom}_{A_{\p}}(\kk(\p), \KK_{\p}) \bigr)  \]
as $\kk(\p)$-modules. By Lemma \ref{lem:2231}(3) the complex 
$\opn{Hom}_{A_{\p}}(\kk(\p), \KK_{\p})$
is a  minimal complex of injective  $\kk(\p)$-modules.
It is easy to see (and we leave this verification to the reader) that a 
minimal complex of injectives over a field must have zero differential. 
Therefore 
\[ \opn{H}^{-d} \bigl( \opn{Hom}_{A_{\p}}(\kk(\p), \KK_{\p}) \bigr) \cong 
\opn{Hom}_{A_{\p}}(\kk(\p), \KK_{\p}^{-d}) . \]
Now by arguments like in the proofs of Lemmas \ref{lem:2230}(1) 
and \ref{lem:2240} we know that 
\[ \opn{Hom}_{A_{\p}} \bigl( \kk(\p), J(\q)_{\p} \bigr) \cong 
\begin{cases}
\kk(\p) & \tup{if} \ \q = \p ,  
\\
0 & \tup{otherwise} . 
\end{cases} \]
It follows that 
\[ \opn{Hom}_{A_{\p}} \bigl( \kk(\p), \KK_{\p}^{-d} \bigr) \cong 
\kk(\p)^{\oplus \mu_{\p, d}} . \]
We see that 
\[ \opn{rank}_{\kk(\p)} \bigl( \opn{Ext}^{-d}_{A_{\p}}(\kk(\p), R_{\p}) \bigr) 
= \mu_{\p, d} . \]
But by Definition \ref{dfn:2200} this number satisfies (\ref{eqn:2245}). 
\end{proof}

\begin{cor} \label{cor:2245}
If $\KK$ is a residue complex over $A$, then it is a minimal complex 
of injective $A$-modules. 
\end{cor}

\begin{proof}
Let $\phi : \KK \to \KK'$ be a minimal injective resolution of $\KK$. According 
to Theorem \ref{thm:2245}, $\KK'$ is also a residue complex. Now 
$\opn{Q}(\phi) : \KK \to \KK'$ is an isomorphism in $\dcat{D}(A)$, so by Theorem
\ref{thm:2240} we know that $\phi : \KK \to \KK'$ is an isomorphism in 
$\dcat{C}_{\mrm{str}}(A)$. 
\end{proof}

\begin{exer} \label{exer:2245}
Find a direct proof of Corollary \ref{cor:2245}, without resorting to Theorems 
\ref{thm:2245} and \ref{thm:2240}. (Hint: look at the proof of Proposition 
\ref{prop:2210}.) 
\end{exer}

We end this subsection with three remarks. 

\begin{rem} \label{rem:2245}
Here is a brief explanation of {\em Matlis Duality}. For more details see 
\cite[Section V.5]{RD}, \cite[Theorem 18.6]{Mats} or \cite[Section 10.2]{BrSh}. 
Assume $A$ is a complete local ring with maximal ideal $\m$. As usual, the 
category of finitely generated $A$-modules is $\dcat{M}_{\mrm{f}}(A)$. There is 
also the category $\dcat{M}_{\mrm{a}}(A)$ of artinian $A$-modules. These are 
full abelian subcategories of $\dcat{M}(A)$. Note that these subcategories are 
characterized by dual properties: the objects of $\dcat{M}_{\mrm{f}}(A)$ are 
noetherian, i.e.\ they satisfy the ascending chain condition; and the objects 
of $\dcat{M}_{\mrm{a}}(A)$ satisfy the descending chain condition. 

Consider the indecomposable injective module $J(\m)$.  
The functor 
$D := \lb \opn{Hom}_A \bigl( -, J(\m) \bigr)$  
is exact of course. Matlis Duality asserts that
$D : \dcat{M}_{\mrm{f}}(A)^{\mrm{op}} \to \dcat{M}_{\mrm{a}}(A)$ 
is an equivalence, with quasi-inverse 
$D : \dcat{M}_{\mrm{a}}(A)^{\mrm{op}} \to \dcat{M}_{\mrm{f}}(A)$.

Later in this book we present a noncommutative graded version of Matlis 
Duality -- this is Theorem \ref{thm:4045}.
\end{rem}

\begin{rem} \label{rem:2250}
We now provide a brief discussion of {\em Local Duality}, based on 
\cite[Section V.6]{RD}. (There is a weaker variant of this result, 
for modules instead of complexes, that can be found in 
\cite[Theorem 11.2.6]{BrSh}.) Again $A$ is local, with maximal ideal $\m$. 
Let $R$ be a dualizing complex over $A$. By translating $R$ we can assume that 
$\opn{dim}_R(\m) = 0$. Lemma \ref{lem:2250} tells us that 
$\opn{H}^i_{\m}(R) = 0$ for all $i \neq 0$. A calculation, that relies on 
Matlis duality, shows that 
$\opn{H}^0_{\m}(R) \cong J(\m)$, the indecomposable injective corresponding to 
$\m$. 

Let us fix an isomorphism $\be : \opn{H}^0_{\m}(R) \iso J(\m)$. This induces 
a morphism 
\begin{equation} \label{eqn:2260}
\th_M : \mrm{R} \Ga_{\m}(M) \to 
\opn{Hom}_A \bigl( \opn{RHom}_A(M, R), J(\m) \bigr) ,
\end{equation}
functorial in $M \in \dcat{D}^+(A)$. The Local Duality Theorem
\cite[Theorem V.6.2]{RD} says that 
$\th_M$ is an isomorphism if $M \in \dcat{D}^+_{\mrm{f}}(A)$.

Here is a modern take on this theorem: we can construct the morphism $\th_M$ 
for all $M \in \dcat{D}(A)$. Let's replace $R$ by the residue complex $\KK$ 
(the minimal injective resolution of $R$). Then $\be$ 
is just an $A$-module isomorphism $\be : \KK^0 \iso J(\m)$. 
For each complex $M$ we choose a K-injective resolution $M \to I(M)$.
Then $\th_M$ is represented by the homomorphism 
\[ \til{\th}_M : \Ga_{\m}(I(M)) \to 
\opn{Hom}_A \bigl( \opn{Hom}_A \bigl( I(M), \KK \bigr), \KK^0 \bigr) \]
in $\dcat{C}_{\mrm{str}}(A)$
that sends an element $u \in \Ga_{\m} \bigl( I(M) \bigr)^p$ 
and a homomorphism 
$\phi \in \opn{Hom}_A \bigl( I(M), \KK \bigr)^{-p}$
to $(-1)^p \cd \phi(u) \in \KK^0$. 

We know that the functor $\mrm{R} \Ga_{\m}$ has finite cohomological 
dimension, \lb bounded by the number of generators of the ideal $\m$; see 
\cite{LC} or \cite{PSY}. The functor 
$\opn{RHom}_{A}(-, R)$
has finite cohomological dimension, which is the injective dimension of $R$. 
And the functor 
$\opn{Hom}_{A}(-, J(\m))$ has cohomological dimension $0$. 
Since $A \in \dcat{D}^+_{\mrm{f}}(A)$, the local duality theorem from \cite{RD}
tells us that $\th_A$ is an isomorphism. Now we can apply Theorem \ref{thm:2135}
to conclude that $\th_M$ is an isomorphism for every 
$M \in \dcat{D}_{\mrm{f}}(A)$.

Finally, let us mention that in Subsection \ref{subsec:bal-NCDC-LocDu} there is 
a noncommutative graded version of Local Duality. 
\end{rem}

\begin{rem} \label{rem:2270}
Here is more on the CM (Cohen-Macaulay) property and duality. 
Let $A$ be a noetherian ring with connected spectrum. Assume $A$ has a 
dualizing complex $R$, and corresponding dimension function $\opn{dim}_{R}$. 

Consider a complex $M \in  \dcat{D}_{\mrm{f}}^{\mrm{b}}(A)$.
In \cite{RD} Grothendieck defines $M$ to be a {\em CM complex with respect to 
$R$} if for every prime ideal $\p \sub A$ and every 
$i \neq - \opn{dim}_{R}(\p)$ the local cohomology satisfies 
$\opn{H}^i_{\p}(M_{\p}) = 0$.

It is proved in \cite{RD} that when $A$ is a regular ring, 
$R = A$, and $M$ is a finitely generated $A$-module, then $M$ is a CM module 
(in the conventional sense, see \cite{Mats}) iff it is a CM complex. 

Let $\dcat{D}_{\mrm{f}}^{0}(A)$ be the full subcategory of 
$\dcat{D}_{\mrm{f}}^{\mrm{b}}(A)$ on the complexes $M$ such that 
$\opn{H}^i(M) = 0$ for all $i \neq 0$.
We know that $\dcat{D}_{\mrm{f}}^{0}(A)$
is equivalent to 
$\dcat{M}_{\mrm{f}}(A) = \cat{Mod}_{\mrm{f}} A$.
In \cite{YeZh2} it was proved that the following
are equivalent for a complex $M \in  \dcat{D}_{\mrm{f}}^{\mrm{b}}(A)$~:
\begin{enumerate}
\rmitem{i} The complex $M$ is CM w.r.t.\ $R$. 

\rmitem{ii} The complex $\opn{RHom}_A(M, R)$ belongs to 
$\dcat{D}_{\mrm{f}}^{0}(A)$. 
\end{enumerate}

It follows that the CM complexes form an abelian subcategory of 
$\dcat{D}^{\mrm{b}}_{\mrm{f}}(A)$,
dual to $\dcat{M}_{\mrm{f}}(A)$. 
In fact, they are the heart of a perverse t-structure on 
$\dcat{D}^{\mrm{b}}_{\mrm{f}}(A)$,
and hence they deserve to be called {\em perverse finitely generated 
$A$-modules}. Geometrically, on the scheme $X := \opn{Spec}(A)$, 
the CM complexes inside $\dcat{D}^{\mrm{b}}_{\mrm{c}}(X)$
form a stack of abelian categories, and so they are 
{\em perverse coherent sheaves}. All this is explained in 
\cite[Section 6]{YeZh2}.
\end{rem}

\mysubsection{Van den Bergh Rigidity}
\label{subsec:squaring}

As we saw in Theorem \ref{thm:2175}, a dualizing complex $R$ 
over a noetherian commutative ring $A$ is not unique. This lack of uniqueness 
(not to mention any sort of functoriality!) was the source of major 
difficulties 
in \cite{RD}, first for gluing dualizing complexes on schemes, and then for 
producing the trace morphisms associated to maps of schemes. 

In 1997, M. Van den Bergh \cite{VdB} discovered the idea of {\em rigidity} for 
dualizing complexes. This was done in the context of noncommutative ring 
theory: $A$ is a noncommutative noetherian ring, central over a base field 
$\K$. The theory of noncommutative rigid dualizing complexes was developed 
further in several papers of J.J. Zhang and A. Yekutieli, among them 
\cite{YeZh1} and \cite{YeZh2}. We shall talk about this noncommutative theory in 
Section \ref{sec:rigid-DC-NC} of the book. 

Here we will deal with the commutative side only, which turns out to be 
extremely powerful. Before explaining it, let us first observe that this is 
one of the rare cases in which an idea originating from noncommutative 
algebra had significant impact in commutative algebra and algebraic geometry.  

In this subsection we study Van den Bergh rigidity in the following context: 

\begin{setup} \label{set:2285}
$A$ is a commutative ring, and $B$ is a flat commutative $A$-ring. 
\end{setup}

We introduce the notion of {\em rigid complex over $B$ relative to $A$},
and describe some of its properties. In Subsection 
\ref{subsec:rig-res-cplx-rng} we will discuss {\em rigid dualizing and residue
complexes}. This material is adapted from the papers \cite{YeZh3}, \cite{YeZh4} 
and \cite{Ye11}. 

The theory of rigid complexes does not really require the assumption
that $B$ is flat over $A$, but flatness makes the theory much easier. See 
Remark \ref{rem:2285}. 

Consider the enveloping ring $B \ot_A B$. It comes equipped with a few ring 
homomorphisms: 
\begin{equation} \label{eqn:2281}
B \xar{\eta_i} B \ot_{A} B \xar{\ep} B ,
\end{equation}
where $\eta_0(b) := b \ot 1$, $\eta_1(b) := 1 \ot b$, and  
$\ep(b_0 \ot b_1) := b_0 \cd b_1$. We view $B$ as a module over 
$B \ot_{A} B$ via $\ep$. Of course $\ep \circ \eta_i = \opn{id}_B$. 

Suppose we are given $B$-modules $M_0$ and $M_1$. Then the tensor product
$M_0 \ot_{A} M_1$ is a $(B \ot_{A} B)$-module. In this way we get an additive 
bifunctor
\[ (- \ot_A -) : \dcat{M}(B) \times \dcat{M}(B) \to \dcat{M}(B \ot_{A} B) . \]
Passing to complexes, and then to homotopy categories, we obtain a triangulated 
bifunctor
\begin{equation} \label{eqn:2276}
(- \ot_A -) : \dcat{K}(B) \times \dcat{K}(B) \to \dcat{K}(B \ot_{A} B) . 
\end{equation}

\begin{lem} \label{lem:2276}
The bifunctor \tup{(\ref{eqn:2276})} has a left derived bifunctor 
\[ (- \ot^{\mrm{L}}_A -) 
: \dcat{D}(B) \times \dcat{D}(B) \to \dcat{D}(B \ot_{A} B) . \]
If either $M_0$ or $M_1$ is a complex of $B$-modules that is K-flat over $A$, 
then the morphism 
$\eta^{\mrm{L}}_{M_0, M_1} : M_0 \ot^{\mrm{L}}_A M_1 \to M_0 \ot_A M_1$
in $\dcat{D}(B \ot_{A} B)$ is an isomorphism. 
\end{lem}

\begin{proof}
This is a variant of Theorem \ref{thm:2107}. 
We know by Corollary \ref{cor:1580} and Proposition \ref{prop:1525} 
that every complex $M \in \dcat{C}(B)$ admits a K-flat resolution 
$P \to M$. Because $B$ is flat over $A$, the complex $P$ is also K-flat over 
$A$. By Theorem \ref{thm:2106} the left derived bifunctor 
$(- \ot^{\mrm{L}}_A -)$ exists, and the condition on 
$\eta^{\mrm{L}}_{M_0, M_1}$ holds. 
\end{proof}

\begin{rem} \label{rem:2275}
The innocuous looking Lemma \ref{lem:2276} is actually of 
tremendous importance. Without the flatness of $A \to B$ we could do very 
little homological algebra of bimodules. Getting around the lack of flatness 
requires the use of flat DG ring resolutions, as explained in Remark 
\ref{rem:2285}. 
\end{rem}

Every module $L \in \dcat{M}(B)$ has an action by $B \ot_{A} B$ coming from the 
homomorphism $\ep$ in (\ref{eqn:2281}). Consider now a module 
$N \in \dcat{M}(B \ot_{A} B)$.
The $\K$-module $N$ has two possible $B$-module structures, coming from the 
homomorphisms $\eta_i$. Thus the $\K$-module
$\opn{Hom}_{B \ot_{A} B}(L, N)$
has three possible $B$-module structures: there is one action from the 
$B$-module structure on $L$, and there are two from the 
$B$-module structures on $N$. 

\begin{lem} \label{lem:2325}
The three $B$-module structures on $\opn{Hom}_{B \ot_{A} B}(L, N)$ coincide. 
\end{lem}

\begin{exer} \label{exer:2325}
Prove the lemma. 
\end{exer}

We are mostly interested in the $B$-module $L = B$. As the module $N$ changes,
we get an additive functor 
\[ \opn{Hom}_{B \ot_{A} B}(B, -) : \dcat{M}(B \ot_{A} B) \to \dcat{M}(B) . \]
Passing to complexes, and then to homotopy categories, we get a triangulated 
functor 
\[ \opn{Hom}_{B \ot_{A} B}(B, -) : \dcat{K}(B \ot_{A} B) \to \dcat{K}(B) . \]
This has a right derived functor
\begin{equation} \label{eqn:2279}
\opn{RHom}_{B \ot_{A} B}(B, -) :
\dcat{D}(B \ot_{A} B) \to \dcat{D}(B)  , 
\end{equation}
that is calculated by K-injective resolutions. Namely if 
$I \in \dcat{C}(B \ot_{A} B)$ is a K-injective complex, then the morphism 
\[ \eta^{\mrm{R}}_{B, I} : \opn{Hom}_{B \ot_{A} B}(B, I) \to 
\opn{RHom}_{B \ot_{A} B}(B, I) \]
in $\dcat{D}(B)$ is an isomorphism. 

By composing the bifunctor $(- \ot^{\mrm{L}}_A -)$ from Lemma \ref{lem:2276} 
and the functor \lb 
$\opn{RHom}_{B \ot_{A} B}(B, -)$ from (\ref{eqn:2279})
we obtain a triangulated bifunctor 
\begin{equation} \label{eqn:2280}
\opn{RHom}_{B \ot_{A} B}(B, - \ot^{\mrm{L}}_{A} -) :
\dcat{D}(B) \times \dcat{D}(B) \to \dcat{D}(B) .
\end{equation}

\begin{dfn} \label{dfn:2275} 
Under Setup \tup{\ref{set:2285}}, the {\em squaring operation} 
\index{Squaring operation}
\index{1-SqBA@$\opn{Sq}_{B / A}$}
is the functor 
\[ \opn{Sq}_{B / A} : \dcat{D}(B) \to \dcat{D}(B) \]
defined as follows:
\begin{enumerate}
\item For a complex $M \in \dcat{D}(B)$, its square is the complex 
\[ \opn{Sq}_{B / A}(M) := 
\opn{RHom}_{B \ot_{A} B}(B, M \ot^{\mrm{L}}_{A} M) \in \dcat{D}(B) . \]

\item For a morphism 
$\phi : M \to N$ in $\dcat{D}(B)$, its square is the morphism 
\[ \opn{Sq}_{B / A}(\phi) := 
\opn{RHom}_{B \ot_{A} B}(B, \phi \ot^{\mrm{L}}_{A} \phi) :
\opn{Sq}_{B / A}(M) \to \opn{Sq}_{B / A}(N) \]
in $\dcat{D}(B)$.
\end{enumerate}
\end{dfn}

It will be good to have an explicit formulation of the squaring operation. 
Let us first choose a K-projective resolution 
$\si : P \to M$ in $\dcat{C}(B)$. 
Note that $P$ is unique up to homotopy equivalence, and $\si$ is unique up to 
homotopy. Since $B$ is flat over $A$, the complex $P$ is K-flat over $A$. We 
get an isomorphism 
$\opn{pres}_{P} : P \ot_A P \iso M \ot^{\mrm{L}}_{A} M$ 
in $\dcat{D}(B \ot_{A} B)$, that we call a {\em presentation}. It is uniquely 
characterized by the commutativity of the diagram
\[ \UseTips \xymatrix @C=10ex @R=6ex {
M \ot^{\mrm{L}}_{A} M
\ar[r]^{\eta^{\mrm{L}}_{M, M}}
&
M \ot_A M
\\
P \ot^{\mrm{L}}_{A} P
\ar[r]_{\eta^{\mrm{L}}_{P, P}}^{\cong}
\ar[u]^{\opn{Q}(\si) \ot^{\mrm{L}}_A \opn{Q}(\si)}_{\cong}
&
P \ot_A P
\ar[ul]_{\opn{pres}_{P}}^(0.6){\cong}
\ar[u]_{\opn{Q}(\si \ot_A \si)}
} \]
in $\dcat{D}(B \ot_{A} B)$. 

Next we choose a K-injective resolution 
$\rho : P \ot_A P \to I$ in $\dcat{C}(B \ot_{A} B)$. 
The complex $I$ is unique up to homotopy equivalence,
and the homomorphism $\rho$ is unique up to homotopy. The resolution 
$\rho$ gives rise to an isomorphism
\[ \opn{pres}_{I} : \opn{Hom}_{B \ot_{A} B}(B, I) \iso 
\opn{RHom}_{B \ot_{A} B}(B, P \ot_{A} P) \]
in $\dcat{D}(B)$ such that the diagram 
\[ \UseTips \xymatrix @C=10ex @R=7ex {
\opn{Hom}_{B \ot_{A} B}(B, P \ot_{A} P)
\ar[r]^{\eta^{\mrm{R}}_{B, P \ot_{A} P}}
\ar[d]_{\opn{Q}(\opn{Hom}(\opn{id}, \rho))}
&
\opn{RHom}_{B \ot_{A} B}(B, P \ot_{A} P)
\ar[d]^{\opn{RHom}(\opn{id}, \opn{Q}(\rho))}_{\cong}
\\
\opn{Hom}_{B \ot_{A} B}(B, I)
\ar[r]_{\eta^{\mrm{R}}_{B, I}}^{\cong}
\ar[ur]^{\opn{pres}_{I}}_(0.6){\cong}
&
\opn{RHom}_{B \ot_{A} B}(B, I)
} \]
is commutative. 

The combination of the presentations $\opn{pres}_{P}$ and 
$\opn{pres}_{I}$ gives an isomorphism 
\[ \opn{pres}_{P, I} : 
\opn{Hom}_{B \ot_{A} B}(B, I) \iso \opn{Sq}_{B / A}(M) \]
in $\dcat{D}(B)$, that we also call a presentation. 

Let $\phi : M \to N$ be a morphism in $\dcat{D}(B)$. The morphism 
$\opn{Sq}_{B / A}(\phi)$ can also be made explicit using presentations. 
For that we need to choose a K-projective resolution 
$\si_N : Q \to N$ in $\dcat{C}(B)$, and a K-injective resolution 
$\rho_N : Q \ot_A Q \to J$ in $\dcat{C}(B \ot_{A} B)$. 
These provide us with a presentation 
\[ \opn{pres}_{Q, J} : \opn{Hom}_{B \ot_{A} B}(M, J) \iso \opn{Sq}_{B / A}(N) .
\]
There are homomorphisms
$\til{\phi} : P \to Q$ in $\dcat{C}_{\mrm{str}}(B)$,
and 
$\chi : I \to J$ in $\dcat{C}_{\mrm{str}}(B \ot_{A} B)$,
both unique up to homotopy, such that the diagrams 
\[ \UseTips \xymatrix @C=8ex @R=6ex {
P
\ar[r]^{\opn{Q}(\si)}_{\cong}
\ar[d]_{\opn{Q}(\til{\phi})}
&
M
\ar[d]^{\phi}
\\
Q
\ar[r]^{\opn{Q}(\si_N)}_{\cong}
&
N
} 
\qquad \quad 
\UseTips \xymatrix @C=9ex @R=6ex {
M \ot^{\mrm{L}}_A M
\ar[d]_{\phi \ot^{\mrm{L}}_A \phi}
&
P \ot_A P
\ar[l]_{\opn{pres}_{P}}^{\cong}
\ar[r]^(0.6){\opn{Q}(\rho)}_{\cong}
\ar[d]_{\opn{Q}(\til{\phi} \ot_A \til{\phi})}
&
I
\ar[d]^{\opn{Q}(\chi)}
\\
N \ot^{\mrm{L}}_A N
&
Q \ot_A Q
\ar[l]_{\opn{pres}_{Q}}^{\cong}
\ar[r]^(0.6){\opn{Q}(\rho_N)}_(0.6){\cong}
&
J
}
\]
in $\dcat{D}(B)$ and $\dcat{D}(B \ot_{A} B)$ respectively are commutative. 
See Subsections \ref{subsec:K-inj} and \ref{subsec:K-proj}. Then the diagram 
\begin{equation} \label{eqn:2329}
\UseTips \xymatrix @C=10ex @R=6ex {
\opn{Hom}_{B \ot_{A} B}(B, I)
\ar[r]^(0.55){\opn{pres}_{P, I}}_(0.55){\cong}
\ar[d]_{\opn{Q}(\opn{Hom}(\opn{id}, \chi))}
&
\opn{Sq}_{B / A}(M)
\ar[d]^{\opn{Sq}_{B / A}(\phi)}
\\
\opn{Hom}_{B \ot_{A} B}(B, J)
\ar[r]^(0.55){\opn{pres}_{Q, J}}_(0.55){\cong}
&
\opn{Sq}_{B / A}(N)
} 
\end{equation}
in $\dcat{D}(B)$ is commutative.

The squaring operation is not an additive functor. In fact, it is a 
{\em quadratic functor}:

\begin{thm} \label{thm:2325} 
Let $\phi : M \to N$ be a morphism in $\dcat{D}(B)$ and let $b \in B$.
Then 
\[ \opn{Sq}_{B / A}(b \cd \phi) = b^2 \cd \opn{Sq}_{B / A}(\phi) , \]
as morphisms 
$\opn{Sq}_{B / A}(M) \to \opn{Sq}_{B / A}(N)$ in $\dcat{D}(B)$. 
\end{thm}

\begin{proof}
We shall use presentations. Let 
$\til{\phi} : P \to Q$ be a homomorphism in $\dcat{C}_{\mrm{str}}(B)$
that represents $\phi$, as above. Then the homomorphism
$b \cd \til{\phi} : P \to Q$ 
$\dcat{C}_{\mrm{str}}(B)$ represents $b \cd \phi$.
Tensoring we get a homomorphism
\[ (b \cd \til{\phi}) \ot_A (b \cd \til{\phi}) : 
P \ot_A P \to Q \ot_A Q \]
$\dcat{C}_{\mrm{str}}(B \ot_A B)$.
But 
\[ (b \cd \til{\phi}) \ot_A (b \cd \til{\phi}) = 
(b \ot b) \cd (\til{\phi} \ot_A \til{\phi}) . \]
Hence on the K-injectives we get the homomorphism
$(b \ot b) \cd \chi : I \to  J$ 
$\dcat{C}_{\mrm{str}}(B \ot_A B)$.  We conclude that 
\[ \opn{Hom}_{B \ot_{A} B} \bigl( \opn{id}_B, (b \ot b) \cd \chi \bigr) : 
\opn{Hom}_{B \ot_{A} B}(B, I) \to 
\opn{Hom}_{B \ot_{A} B}(B, J) \]
represents 
$\opn{Sq}_{B / A}(b \cd \phi)$. Finally, by Lemma \ref{lem:2325} 
we know that 
\[ \begin{aligned}
& \opn{Hom}_{B \ot_{A} B} \bigl( \opn{id}_B, (b \ot b) \cd \chi \bigr) = 
\opn{Hom}_{B \ot_{A} B}(b^2 \cd \opn{id}_B, \chi) 
\\
& \quad = b^2 \cd  \opn{Hom}_{B \ot_{A} B}(\opn{id}_B, \chi) . 
\end{aligned} \]
\end{proof}

\begin{dfn} \label{dfn:2276} 
Let $M \in \dcat{D}(B)$. 
A {\em rigidifying isomorphism} 
\index{Rigidifying isomorphism}
for $M$ over $B$ 
relative to $A$ is an isomorphism 
$\rho : M \iso \opn{Sq}_{B / A}(M)$
in $\dcat{D}(B)$. 
\end{dfn}

\begin{dfn} \label{dfn:2295}
A {\em rigid complex} 
\index{Rigid complex}
over $B$ relative to $A$ is a pair 
$(M, \rho)$, consisting of a complex $M \in \dcat{D}(B)$ 
and a rigidifying isomorphism 
$\rho : M \iso \opn{Sq}_{B / A}(M)$
in $\dcat{D}(B)$. 
\end{dfn}

\begin{dfn} \label{dfn:2296}
Suppose $(M, \rho)$ and $(N, \si)$ are rigid complexes over $B$ relative 
to $A$. A {\em morphism of rigid complexes}
\index{Rigid complex! morphism of {\indash}s} 
$\phi : (M, \rho) \to (N, \si)$ 
is a morphism $\phi : M \to N$ in $\dcat{D}(B)$, such that the diagram 
\[ \UseTips \xymatrix @C=8ex @R=6ex {
M
\ar[r]^(0.35){\rho}_(0.35){\cong}
\ar[d]_{\phi}
&
\opn{Sq}_{B / A}(M)
\ar[d]^{\opn{Sq}_{B / A}(\phi)}
\\
N
\ar[r]^(0.35){\si}_(0.35){\cong}
&
\opn{Sq}_{B / A}(N)
} \]
in $\dcat{D}(B)$ is commutative. 

The category of rigid complexes over $B$ relative to $A$ is denoted by
\index{1-D(A)rigK@$\dcat{D}(A)_{\mrm{rig} / \K}$}
$\dcat{D}(B)_{\mrm{rig} / A}$. 
\end{dfn}

Recall (Definition \ref{dfn:2394}) that a complex $M \in \dcat{D}(B)$ 
has the derived Morita property 
if the derived homothety morphism 
$\opn{hm}^{\mrm{R}}_{M} : B \to \opn{RHom}_B(M, M)$
in $\dcat{D}(B)$ is an isomorphism. 

\begin{thm} \label{thm:2275}
Let $(M, \rho)$ be a rigid complex over $B$ relative to $A$. If $M$ has the 
derived Morita property, then the only automorphism of 
$(M, \rho)$ in $\dcat{D}(B)_{\mrm{rig} / A}$ is the identity. 
\end{thm}

\begin{proof}
Let 
$\phi : (M, \rho) \iso (M, \rho)$
be an automorphism in $\dcat{D}(B)_{\mrm{rig} / A}$.
By Proposition \ref{prop:2335}, there is a unique invertible element $b \in B$ 
such that $\phi = b \cd \opn{id}_M$, as morphisms $M \to M$ in 
$\dcat{D}(B)$. 

Next, according to Theorem \ref{thm:2325}, we have 
\[ \opn{Sq}_{B / A}(\phi) = \opn{Sq}_{B / A}(b \cd \opn{id}_M) = 
b^2 \cd \opn{Sq}_{B / A}(\opn{id}_M) . \]
Plugging this into the diagram in Definition \ref{dfn:2296} we get a 
commutative diagram 
\[ \UseTips \xymatrix @C=8ex @R=6ex {
M
\ar[r]^(0.35){\rho}_(0.35){\cong}
\ar[d]_{b \cd \opn{id}_M}
&
\opn{Sq}_{B / A}(M)
\ar[d]^{b^2 \cd \opn{id}_M}
\\
M
\ar[r]^(0.35){\rho}_(0.35){\cong}
&
\opn{Sq}_{B / A}(M)
} \]
in $\dcat{D}(B)$. 
Once more using Proposition \ref{prop:2335} we see that $b^2 = b$. 
Because $b$ is an invertible element, it follows that $b = 1$. Thus 
$\phi = \opn{id}_M$. 
\end{proof}

\begin{exa} \label{exa:2275}
Assume $B = A$, and take $M := B$. Then $B \ot_{A} B \cong B$, 
$M \ot^{\mrm{L}}_A M \cong M$, and there are canonical isomorphisms
\[ \opn{Sq}_{B / A}(M) = \opn{RHom}_{B \ot_{A} B}(B, M \ot^{\mrm{L}}_A M) \cong 
\opn{Hom}_{B}(B, M) \cong M . \]
Thus the pair $(M, \opn{id}_M)$ belongs to $\dcat{D}(B)_{\mrm{rig} / A}$.
Furthermore, the complex $M$ has the derived Morita property, so Theorem 
\ref{thm:2275} applies. 
\end{exa}

\begin{rem} \label{rem:3120}
To the reader who might object to Example \ref{exa:2275} for being silly, 
we have two things to say. First, the existence theorem for rigid dualizing 
complexes \ref{thm:3210}, whose proof is beyond the scope of this book, starts  
with the {\em tautological rigid structure} of $\K \in \dcat{D}(\K)$, where 
$\K$ is the regular base ring. The rigid structure is then propagated to all 
flat essentially finite type (FEFT) $\K$-rings using {\em induction} and {\em 
coinduction} of rigid structures. See the exercise and examples below for an 
indication of these procedures. 

The second fact we want to point out is that when $A$ is an FEFT $\K$-ring, 
then {\em $A$ has exactly one rigid complex} 
$(M, \rho) \in \dcat{D}(A)_{\mrm{rig} / \K}$ (up to a unique isomorphism) 
that is nonzero on each connected component of $\opn{Spec}(A)$; and it is the 
rigid dualizing complex of $A$ (that will be defined in the next subsection). 
This is proved in \cite{YeZh3}, \cite{AIL} and \cite{Sha3}. 
\end{rem}

\begin{exer} \label{exer:2395}
Assume $A \neq 0$, and take $B := A[t_1, \ldots, t_n]$, the polynomial ring in 
$n$ variables. 
\begin{enumerate}
\item Let $C := B \ot_A B$, and let $I$ be the kernel of the multiplication 
homomorphism $\ep : C \to B$. Show that $I$ is generated by the sequence 
$\bsym{c} = (c_1, \ldots, c_n)$, where 
$c_j := t_j \ot 1 - 1 \ot t_j$. 

\item Show that the Koszul complex 
$\opn{K}(C; \bsym{c})$ is a free resolution of $B$ over $C$. 
(Hint: for $1 \leq m \leq n$ define 
$B_m := A[t_1, \ldots, t_m]$ and 
$C_m := B_m \ot_A B_m$. 
By  direct calculation show that 
$\opn{K}(C_1; \bsym{c}_1) \to B_1$
is a quasi-isomorphism. Next, there is a ring isomorphism 
$B_m \ot_A B_1 \iso B_{m + 1}$, 
$1 \ot t_1 \mapsto t_{m + 1}$.
This induces a ring isomorphism 
$C_m \ot_A C_1 \cong C_{m + 1}$ 
and an isomorphism of complexes 
\[ \opn{K}(C_m; \bsym{c}_m) \ot_A \opn{K}(C_1; \bsym{c}_1)
\cong \opn{K}(C_{m + 1}; \bsym{c}_{m + 1}) . \]
Use this last formula and induction on $m$ to prove that 
$\opn{K}(C_m; \bsym{c}_m) \to B_m$ 
is a quasi-isomorphism.)

\item  Prove that 
\[ \opn{Ext}^i_{C}(B, C) \cong 
\begin{cases}
B & \tup{if} \quad i = n \\
0 & \tup{if} \quad i \neq n  . 
\end{cases} \]
(Hint: use the Koszul resolution from above.)

\item Conclude from item (3) that the complex $B[n] \in \dcat{D}(B)$ is rigid 
relative to $A$; namely that there is a rigidifying isomorphism 
$\rho : B[n] \iso \opn{Sq}_{B / A} \bigl( B[n] \bigr)$
in $\dcat{D}(B)$. 
\end{enumerate}
\end{exer}

\begin{exa} \label{exa:2576}
Let $A$ be a nonzero noetherian ring, and let
$B := \lb A[t_1, \ldots, t_n]$ as in the exercise above. 
Let $\Delta_{B / A} := \Om^n_{B / A}$, the module of degree $n$ 
differential forms of $B$ relative to $A$. 
It is the $n$-th exterior power of $\Om^1_{B / A}$, so it is a free 
$B$ module of rank $1$, with basis 
$\d(t_1) \wedge \cdots \wedge \d(t_n)$.
By the previous example, there exists some rigidifying isomorphism $\rho$ for 
the complex $\Delta_{B / A}[n]$. 

However, it can be shown (see \cite{YeZh3} or \cite{Ye11}) that the complex
$\Delta_{B / A}[n]$ has a {\em canonical rigidifying isomorphism relative to 
$A$}. I.e.\ there is a rigidifying isomorphism 
$\rho : \Delta_{B / A}[n] \iso 
\opn{Sq}_{B / A} \bigl( \Delta_{B / A}[n] \bigr)$
in $\dcat{D}(B)$ which is invariant under $A$-ring automorphisms of $B$. 

This rigidifying isomorphism $\rho$ is an incarnation of Grothendieck's 
Fundamental Local Isomorphism \cite[Proposition II.7.2]{RD} and the Residue 
Isomorphism
\cite[Theorem II.8.2]{RD}. 
\end{exa}

It is well-known that a finitely generated module $M$ over a noetherian ring 
$A$ is flat iff it is projective. See \cite[Corollary to Theorem 7.12]{Mats}. 

\begin{exa} \label{exa:2575}
Let $A$ be a noetherian ring, and let $A \to B$ be a finite flat ring 
homomorphism. So $B$ is a finitely generated projective $A$-module. Define 
$\Delta_{B / A} := \opn{Hom}_{A}(B, A) \in \dcat{M}(B)$.
It can be shown (see \cite{Ye11}) that the complex 
$\Delta_{B / A}$ has a canonical rigidifying isomorphism relative to $A$. I.e.\ 
there is a rigidifying isomorphism 
$\rho : \Delta_{B / A} \iso \opn{Sq}_{B / A}(\Delta_{B / A})$
in $\dcat{D}(B)$ which is invariant under $A$-ring automorphisms of $B$. 
\end{exa}

\begin{rem} \label{rem:4190}
The squaring operation $\opn{Sq}_{B / A}$ is functorial also in the ring $B$, 
in two directions. For this remark we assume that the ring $A$ is noetherian, 
and the rings $B$ and $C$ are flat essentially finite type $A$-rings. (See 
Remark \ref{rem:2285} regarding the flatness.)

Suppose $u : B \to C$ is a finite 
$A$-ring homomorphism, and $\th : N \to M$ is a backward (or trace) morphism 
in $\dcat{D}(B)$ over $u$, in the sense of Definition \ref{dfn:3110}. 
Then there is a backward morphism
\[ \opn{Sq}_{u / A}(\th) : \opn{Sq}_{C / A}(N) \to \opn{Sq}_{B / A}(M) \]
in $\dcat{D}(B)$ over $u$. The morphism $\opn{Sq}_{u / A}(\th)$ is functorial 
in $u$ and $\th$ in an obvious sense. Furthermore, under suitable finiteness 
conditions, if $\th$ is a nondegenerate backward morphism (see Definition 
\ref{dfn:3190}), then so is $\opn{Sq}_{u / A}(\th)$. We call this property the 
{\em trace functoriality} of the squaring operation. 

The trace functoriality of the squaring operation allows us to talk about {\em 
rigid trace morphisms}. With $u : B \to C$ as above, suppose that 
$(M, \rho) \in \dcat{D}(B)_{\mrm{rig} / A}$ and 
$(N, \si) \in \dcat{D}(C)_{\mrm{rig} / A}$. A {\em rigid trace morphism} 
$\th : (N, \si) \to (M, \rho)$
over $u$ relative to $A$ is a trace morphism $\th : N \to M$ over $u$ such 
that the diagram 
\[ \UseTips \xymatrix @C=8ex @R=6ex {
N
\ar[r]^(0.35){\si}_(0.35){\cong}
\ar[d]_{\th}
&
\opn{Sq}_{C / A}(N)
\ar[d]^{\opn{Sq}_{u / A}(\th)}
\\
M
\ar[r]^(0.35){\rho}_(0.35){\cong}
&
\opn{Sq}_{B / A}(M)
} \]
in $\dcat{D}(B)$ is commutative. 

Next, suppose $v : B \to C$ is an essentially \'etale homomorphism 
(i.e.\ formally \'etale, plus the EFT condition that was already assumed) of 
$A$-rings, and $\la : M \to N$ is a forward morphism over $v$, in 
the sense of Definition \ref{dfn:3111}. Then there is a forward morphism 
\[ \opn{q}_{v / A}(\la) : \opn{Sq}_{B / A}(M) \to \opn{Sq}_{C / A}(N) \]
in $\dcat{D}(B)$ over $u$. The morphism $\opn{Sq}_{v / A}(\la)$ is functorial 
in $v$ and $\la$ in an obvious sense. Furthermore, under suitable finiteness 
conditions, if $\la$ is a nondegenerate forward morphism (see Definition 
\ref{dfn:3191}), then so is $\opn{Sq}_{v / A}(\la)$. We call this property the 
{\em \'etale functoriality} of the squaring operation. 

The \'etale functoriality of the squaring operation allows us to talk about 
{\em rigid localization morphisms}. With $v : B \to C$ as above, suppose that 
$(M, \rho) \in \dcat{D}(B)_{\mrm{rig} / A}$, and 
$(N, \si) \in \dcat{D}(C)_{\mrm{rig} / A}$. A {\em rigid localization morphism} 
$\la : (M, \rho) \to (N, \si)$
over $v$ relative to $A$ is a forward morphism $\la : N \to M$ over $v$ such 
that the diagram 
\[ \UseTips \xymatrix @C=8ex @R=6ex {
M
\ar[r]^(0.35){\rho}_(0.35){\cong}
\ar[d]_{\la}
&
\opn{Sq}_{B / A}(M)
\ar[d]^{\opn{Sq}_{v / A}(\la)}
\\
N
\ar[r]^(0.35){\si}_(0.35){\cong}
&
\opn{Sq}_{C / A}(N)
} \]
in $\dcat{D}(B)$ is commutative. 

These functorialities of the squaring operation are hard to prove. An imprecise 
study of them was made in the papers \cite{ YeZh3} and \cite{YeZh4}. 
A correct treatment of the trace functoriality can be found in the paper 
\cite{Ye9}; and a correct treatment of the \'etale functoriality will 
appear in \cite{Ye11}. 
\end{rem}

\begin{rem} \label{rem:2347}
The squaring operation is related to {\em Hochschild cohomology}.
Assume for simplicity that $A$ is a field and $M$ is a $B$-module. 
Then for each $i$ the cohomology 
\[ \opn{H}^i (\opn{Sq}_{B / A}(M)) = 
\opn{Ext}^i_{B \ot_{A} B}(B, M \ot_A M) \]
is the $i$-th Hochschild cohomology with values in the $B$-bimodule 
$M \ot_A M$. 
For more on this material see \cite{AILN}, \cite{Sha1} and \cite{Sha2}. 
\end{rem}

\begin{rem} \label{rem:2285}
It is possible to avoid the assumption that $B$ is flat over $A$. 
This is done by choosing a commutative DG ring $\til{B}$ that is K-flat as a DG 
$A$-module, with a DG ring quasi-isomorphism $\til{B} \to B$ over $A$. 
Such resolutions always exist (see \cite{YeZh3} or \cite{Ye9}). Then we take 
\begin{equation} \label{eqn:2345}
\opn{Sq}_{B / A}(M) := 
\opn{RHom}_{\til{B} \ot_{A} \til{B}}(B, M \ot^{\mrm{L}}_A M) .
\end{equation}
This was the construction used by Zhang and Yekutieli in the paper 
\cite{YeZh3}.

Unfortunately there was a serious error in \cite{YeZh3}: 
we did not prove that formula (\ref{eqn:2345}) is independent of  
the resolution $\til{B}$. This error was discovered, and corrected,
by L. Avramov, S. Iyengar, J. Lipman and S. Nayak in their paper \cite{AILN}. 

There were ensuing errors in \cite{YeZh3} regarding the functoriality of the 
squaring operation in the ring $B$ (as described in Remark \ref{rem:4190}). 
The paper \cite{AILN} did not treat such functoriality at all, and the 
constructions and proofs of the trace functoriality were corrected only 
in our recent paper \cite{Ye9}. It is worthwhile to mention that the correct 
proofs (both in \cite{AILN} and \cite{Ye9}) rely on {\em noncommutative DG 
rings} and {\em DG bimodules} over them; and the squaring operation is 
replaced by the {\em rectangle operation}. 

Because the non-flat case is so much more complicated, we have decided not to 
reproduce it in the book. The interested reader can look up the research papers 
\cite{Ye9}, \cite{Ye11}, \cite{Ye12} and \cite{Ye13}, the survey article 
\cite{Ye5}, or the lecture notes \cite{Ye10}. 

A general treatment of derived categories of bimodules, based on K-flat DG ring 
resolutions, will be in the paper \cite{Ye14}.
\end{rem}

\mysubsection{Rigid Dualizing and Residue Complexes}
\label{subsec:rig-res-cplx-rng}

In this subsection we make these assumptions:

\begin{setup} \label{set:3210}
$\K$ is a nonzero regular noetherian commutative ring (Definition 
\ref{dfn:3075}). All other rings are in the 
category $\ringcfeft{\K}$ of FEFT commutative $\K$-rings, where 
``FEFT'' is short for ``flat essentially finite type''. 
\end{setup}

\begin{dfn} \label{dfn:2530}
A {\em rigid dualizing complex}%
\index{Dualizing complex! commutative rigid}
over $A$ relative to $\K$ is a rigid complex 
$(R, \rho)$ over $A$ relative to $\K$, as in Definition \ref{dfn:2295}, such 
that $R$ is a dualizing complex over $A$, in the sense of Definition 
\ref{dfn:2155}.
\end{dfn}

See Remark \ref{rem:2285} regarding the flatness condition. 

\begin{dfn} \label{dfn:4200}
For $A \in \ringcfeft{\K}$ we write  
$A^{\mrm{en}} := A \ot_{\K} A$ for the 
{\em enveloping ring}%
\index{Enveloping ring}%
\index{1-Aen@$A^{\mrm{en}}$}
of $A$ relative to $\K$.
\end{dfn}

With this notation the square of a complex $M \in \dcat{D}(A)$ relative to $\K$ 
is 
\[ \opn{Sq}_{A / \K}(M) = 
\opn{RHom}_{A^{\mrm{en}}}(A, M \otimes^{\mrm{L}}_{\K} M)
\in \dcat{D}(A)  . \]

Recall the category $\dcat{D}(A)_{\mrm{rig} / \K}$ of rigid complexes over 
$A$ relative to $\K$, from Definition \ref{dfn:2296}. 

\begin{thm}[Uniqueness] \label{thm:2280}
Let $A$ be a flat essentially finite type ring over the regular noetherian ring 
$\K$. If $A$ has a rigid dualizing complex $(R, \rho)$, then it is 
unique up to a unique isomorphism in $\dcat{D}(A)_{\mrm{rig} / \K}$. 
\end{thm}

\begin{proof}
Suppose $(R', \rho')$ is another rigid dualizing 
complex over $A$ relative to $\K$. Let 
$A = \prod_{i = 1}^r A_i$ be the connected component decomposition of $A$. 
Corollary \ref{cor:2175} says that 
$R' \cong R \ot^{\mrm{L}}_A P$,
where 
$P \cong \bigoplus_{i = 1}^r L_i[n_i]$
for integers $n_i$ and rank $1$ projective $A_i$-modules $L_i$. 
There is an isomorphism 
\begin{equation} \label{eqn:2545}
R' \otimes^{\mrm{L}}_{\K} R' = 
(R \ot^{\mrm{L}}_A P) \otimes^{\mrm{L}}_{\K} (R \ot^{\mrm{L}}_A P) \cong 
(R \otimes^{\mrm{L}}_{\K} R) \otimes^{\mrm{L}}_{A^{\mrm{en}}}
(P \otimes^{\mrm{L}}_{\K} P) 
\end{equation}
in $\dcat{D}(A^{\mrm{en}})$, and 
$P \otimes^{\mrm{L}}_{\K} P$
has finite flat dimension over $A^{\mrm{en}}$. 
So we have this sequence of isomorphisms in $\dcat{D}(A)$~:
\begin{equation} \label{eqn:2530}
\begin{aligned}
& R \ot^{\mrm{L}}_A P \cong R' \cong \opn{Sq}_{A / \K}(R')
=  \mrm{RHom}_{A^{\mrm{en}}}(A, R' \otimes^{\mrm{L}}_{\K} R') 
\\
& \quad 
\cong^{\diamondsuit} 
\mrm{RHom}_{A^{\mrm{en}}}
\bigl( A, (R \otimes^{\mrm{L}}_{\K} R) \otimes^{\mrm{L}}_{A^{\mrm{en}}}
(P \otimes^{\mrm{L}}_{\K} P) \bigr) 
\\
& \quad 
\cong^{\dag} \mrm{RHom}_{A^{\mrm{en}}}(A, R \otimes^{\mrm{L}}_{\K} R) 
\otimes^{\mrm{L}}_{A^{\mrm{en}}} (P \otimes^{\mrm{L}}_{\K} P)
\\
& \quad 
= \opn{Sq}_{A / \K}(R) 
\otimes^{\mrm{L}}_{A^{\mrm{en}}} (P \otimes^{\mrm{L}}_{\K} P)
\\
& \quad 
\cong R \otimes^{\mrm{L}}_{A^{\mrm{en}}} (P \otimes^{\mrm{L}}_{\K} P)
\cong R \ot^{\mrm{L}}_A P  \ot^{\mrm{L}}_A P .
\end{aligned}
\end{equation}
The isomorphism $\cong^{\diamondsuit}$ is by (\ref{eqn:2545}), and the 
isomorphism $\cong^{\dag}$ is by Theorem \ref{thm:2700}.
We also used the rigidifying isomorphisms of $R$ and $R'$. Now 
\[ \opn{RHom}_{A}(R, R \ot^{\mrm{L}}_A P) \cong 
\opn{RHom}_{A}(R, R) \ot^{\mrm{L}}_A P \cong P , \]
again using Theorem \ref{thm:2700}, and by the derived Morita property of 
$R$. Likewise there are isomorphisms
\[ \begin{aligned}
& 
\opn{RHom}_{A}(R, R \ot^{\mrm{L}}_A P \ot^{\mrm{L}}_A P) \cong 
\opn{RHom}_{A}(R, R \ot^{\mrm{L}}_A P) \ot^{\mrm{L}}_A P
\\
& \quad
\cong \opn{RHom}_{A}(R, R) \ot^{\mrm{L}}_A P \ot^{\mrm{L}}_A P
\cong P \ot^{\mrm{L}}_A P 
\end{aligned} \]
in $\dcat{D}(A)$.
Thus, together with (\ref{eqn:2530}), we deduce that 
$P \ot^{\mrm{L}}_A P \cong P$. 
But then on each connected component $A_i$ we have 
\[ L_i[n_i] \cong L_i[n_i] \ot_A L_i[n_i]  = 
(L_i \ot_A L_i)[2 \cd n_i] . \]
This implies that $L_i \cong A_i$ and $n_i = 0$. 
We see that actually $P \cong A$, so there is an isomorphism 
$\phi_0 : R \iso R'$ in $\dcat{D}(A)$. 

The isomorphism $\phi_0$ might not be rigid; but due to the derived 
Morita property of the complex $R'$, there is an invertible element $a \in A$ 
such that 
$\opn{Sq}_{A / \K}(\phi_0) \circ \rho = a \cd \rho' \circ \phi_0$ 
as isomorphisms 
$R \iso \opn{Sq}_{A / \K}(R')$. 
Define $\phi := a^{-1} \cd \phi_0$. 
Then, according to Theorem \ref{thm:2325}, we have 
\[ \opn{Sq}_{A / \K}(\phi) \circ \rho = 
a^{-2} \cd \opn{Sq}_{A / \K}(\phi_0) \circ \rho = 
a^{-2} \cd a \cd \rho' \circ \phi_0 = \rho' \circ \phi . \]
We see that
$\phi : (R, \rho) \iso (R', \rho')$
is a rigid isomorphism. Its uniqueness is already known by Theorem 
\ref{thm:2275}, since $R'$ has the derived Morita property.
\end{proof}

\begin{thm}[Existence] \label{thm:3210}
Let $A$ be a flat essentially finite type ring over the regular noetherian ring 
$\K$. Then $A$ has a rigid dualizing complex $(R_A, \rho_A)$. 
\end{thm}

The proof of this theorem is beyond the scope of the book, as it requires 
the study of {\em induced and coinduced rigid structures}. 
For a proof see \cite{Ye11}. (There is an incomplete proof in 
\cite[Theorem 1.1]{YeZh4}.)
The next two examples provide proofs in some cases. 

\begin{exa} \label{exa:4200}
Suppose $A = \K[t_1, \ldots, t_n]$, the polynomial ring in $n$ variables.
(More generally, we can take $A$ to be any essentially smooth $\K$-ring of 
relative dimension $n$). Define the complex 
$R_A := \Om^n_{A / \K}[n]$. According to Example \ref{exa:2576} there is a 
rigidifying isomorphism 
$\rho_A : R_A \iso \opn{Sq}_{A / \K}(R_A)$.
Because $A$ is a regular ring, it follows that $R_A$ is a dualizing complex. We 
see that $(R_A, \rho_A)$ is a rigid dualizing complex. 
\end{exa}

\begin{exa} \label{exa:4201}
Suppose $A$ is a finite type $\K$-ring. As explained in Example \ref{exa:4506},
$A$ has a noncommutative rigid dualizing complex, that is actually $A$-central. 
Thus it is a rigid dualizing complex over $A$ relative to $\K$ in the sense of 
Definition \ref{dfn:2530}. 
\end{exa}

The dimension function $\opn{dim}_R$ relative to a dualizing complex $R$
was introduced in Definition \ref{dfn:2200}. If $R' \cong R$, then of course 
the dimension functions satisfy 
$\opn{dim}_{R'} = \opn{dim}_R$. 
In view of the previous theorem, the next definition is valid. 

\begin{dfn} \label{dfn:2537}
Let $A \in \ringcfeft{\K}$. The {\em rigid dimension function 
relative to $\K$} is the function
\[   \opn{rig{.}dim}_{\K} : \opn{Spec}(A) \to \Z \]
given by the formula 
$\opn{rig{.}dim}_{\K}(\p) := \opn{dim}_R(\p)$,
where $R$ is any rigid dualizing complex over $A$ relative to $\K$. 
\end{dfn}

The next examples (taken from \cite{Ye11}) give an idea what the rigid 
dimension function looks like.  

\begin{exa} \label{exa:3235}
Let $\K$ be a field and $A$ a finite type $\K$-ring. 
Then for every $\p \in \opn{Spec}(A)$ there is equality
$\opn{rig{.}dim}_{\K}(\p) = \opn{dim}(A / \p)$, 
where the latter is the Krull dimension of the ring $A / \p$. 
\end{exa}

\begin{exa} \label{exa:3236}
Again $\K$ is a field, but now $L$ is a finitely generated field extension of 
$\K$; or in other words, $L$ is a field in $\ringcfeft{\K}$ . 
Then, writing $\p := (0) \sub L$, we have  
$\opn{rig{.}dim}_{\K}(\p) = \opn{tr{.}deg}_{\K}(L)$,
the transcendence degree of the field extension $\K \to L$. 
\end{exa}

\begin{exa} \label{exa:3237}
Take $\K = A = \Z$. For a maximal ideal $\m = (p) \sub \Z$ we 
have $\opn{rig{.}dim}_{\K}(\m) = -1$; and for the generic ideal 
$\p = (0) \sub \Z$ we have $\opn{rig{.}dim}_{\K}(\p) = 0$. 
\end{exa}

Residue complexes were introduced in Subsection \ref{subsec:res-cplx-ring}.

\begin{dfn} \label{dfn:2535}
A {\em rigid residue complex} 
\index{Residue complex! rigid}
over $A$ relative to $\K$ is a 
rigid complex $(\KK_A, \rho_A)$ over $A$ relative to $\K$, such that $\KK_A$ is 
a residue complex over $A$. 
\end{dfn}

Using the rigid dimension function relative to $\K$, we have this decomposition 
of the $A$-module $\KK^{-i}_A$ for each $i$~:
$\KK_A^{-i} \cong \bigoplus_{\opn{rig{.}dim}_{\K}(\p) = i} J(\p)$,
where $J(\p)$ is the indecomposable injective module corresponding to the prime 
ideal $\p$. 

In Definition \ref{dfn:2296} we introduced the category 
$\dcat{D}(A)_{\mrm{rig} / \K}$. Recall that the objects of 
$\dcat{D}(A)_{\mrm{rig} / \K}$
are rigid complexes $(M, \rho)$ over $A$ relative to $\K$; and the morphisms 
$\phi : (M, \rho) \to (N, \si)$
are the morphisms $\phi : M \to N$ in $\dcat{D}(A)$ for which there is equality 
$\si \circ \phi = \opn{Sq}_{A / \K}(\phi) \circ \rho$.
Rigid residue complexes live, or rather move, in another category. 

\begin{dfn} \label{dfn:2536}
The category $\dcat{C}(A)_{\mrm{rig} / \K}$
\index{1-C(A)rigK@$\dcat{C}(A)_{\mrm{rig} / \K}$}
is defined as follows. Its objects are the rigid complexes $(M, \rho)$ over $A$ 
relative to $\K$. Given two objects 
$(M, \rho)$ and $(N, \si)$, a morphism 
$\phi : (M, \rho) \to (N, \si)$
in $\dcat{C}(A)_{\mrm{rig} / \K}$ is a morphism 
$\phi : M \to  N$ in $\dcat{C}_{\mrm{str}}(A)$,
such that the diagram 
\[ \UseTips \xymatrix @C=8ex @R=6ex {
M
\ar[d]_{\opn{Q}(\phi)}
\ar[r]^(0.35){\rho}_(0.35){\cong}
&
\opn{Sq}_{A / \K}(M)
\ar[d]^{\opn{Sq}_{A / \K}(\opn{Q}(\phi))}
\\
N
\ar[r]^(0.35){\si}_(0.35){\cong}
&
\opn{Sq}_{A / \K}(N)
} \]
in $\dcat{D}(A)$ is commutative.
\end{dfn}

Let us emphasize the hybrid nature of the category 
$\dcat{C}(A)_{\mrm{rig} / \K}$~:
the morphisms are homomorphisms of complexes (literally degree $0$ 
homomorphisms of graded $A$-modules  $\phi : M \to N$ that commute with the 
differentials); but they must satisfy a compatibility condition (rigidity) in 
the derived category. 

\begin{thm} \label{thm:2535}
Let $A$ be an FEFT ring over the regular noetherian ring 
$\K$. The ring $A$ has a rigid residue complex $(\KK_{A}, \rho_{A})$ 
relative to $\K$, and it is unique, up to a unique isomorphism in 
$\dcat{C}(A)_{\mrm{rig} / \K}$.
\end{thm}

\begin{proof}
Existence: by Theorem \ref{thm:3210} there is a rigid dualizing complex 
$(R_{A}, \rho_{A})$ over $A$ relative to $\K$. Let 
$\KK_{A}$ be the minimal injective resolution of the complex $R_{A}$. 
According to Theorem \ref{thm:2245}, $\KK_{A}$ is a residue complex. It 
inherits the rigidifying isomorphism $\rho_{A}$ from $R_{A}$. So the 
pair $(\KK_{A}, \rho_{A})$ is a rigid residue complex over $A$ relative to 
$\K$. 

Uniqueness: suppose $(\KK', \rho')$ is another rigid residue complex over 
$A$ relative to $\K$.  Theorem \ref{thm:2280} tells us that there is a unique 
isomorphism 
$\phi : (\KK_{A}, \rho_{A}) \iso (\KK', \rho')$
in $\dcat{D}(A)_{\mrm{rig} / \K}$. 
Now the dimension functions of these two residue complexes are equal (both are 
$\opn{rig{.}dim}_{\K}$). So by Theorem \ref{thm:2240} the function 
\[ \opn{Q} : \opn{Hom}_{\dcat{C}(A)_{\mrm{rig} / \K}}
\bigl( (\KK_{A}, \rho_{A}), (\KK', \rho') \bigr) \to 
\opn{Hom}_{\dcat{D}(A)_{\mrm{rig} / \K}}
\bigl( (\KK_{A}, \rho_{A}), (\KK', \rho') \bigr) \]
is bijective. Thus $\phi$ lifts uniquely to an isomorphism in 
$\dcat{C}(A)_{\mrm{rig} / \K}$. 
\end{proof}

\begin{rem} \label{rem:4191}   
The rigid residue complex $\KK_{A}$ is functorial in the ring $A$ in two 
different ways, that we briefly explain here. As in Setup \ref{set:3210}, we 
work in the category $\ringcfeft{\K}$. The flatness condition can be removed 
(at a price, see Remark \ref{rem:2285}), but it seems that the assumption that 
the base ring $\K$ is regular, and the other rings are EFT over it, are 
necessary.  
 
Suppose $u : A \to B$ is a finite homomorphism in $\ringcfeft{\K}$. 
Then there is a homomorphism 
$\opn{tr}_{B / A} : \KK_B \to \KK_A$
in $\dcat{C}_{\mrm{str}}(A)$, called the {\em rigid trace homomorphism}. 
It is a nondegenerate rigid trace morphism over $u$ (see Remark \ref{rem:4190}),
namely it respects the rigidifying isomorphisms, and the induced homomorphism 
$\KK_B \to \opn{Hom}_{A}(B, \KK_A)$
in $\dcat{C}_{\mrm{str}}(B)$ is an isomorphism. Furthermore, generalizing 
Theorem \ref{thm:2535}, $\opn{tr}_{B / A}$ is the unique 
nondegenerate rigid trace morphism over $u$. If $B \to C$ is another finite 
homomorphism in $\ringcfeft{\K}$, then 
$\opn{tr}_{B / A} \circ \opn{tr}_{C / B}  = \opn{tr}_{C / A}$.

Suppose $v : A \to A'$ is an essentially \'etale homomorphism in 
$\ringcfeft{\K}$. Then there is a homomorphism 
$\opn{q}_{A' / A} : \KK_{A} \to \KK_{A'}$
in $\dcat{C}_{\mrm{str}}(A)$, called the {\em rigid localization homomorphism}. 
It is a nondegenerate rigid localization morphism over $v$ (see Remark 
\ref{rem:4190}), namely it respects the rigidifying isomorphisms, and the 
induced homomorphism 
$A' \ot_{A} \KK_{A} \to \KK_{A'}$
in $\dcat{C}_{\mrm{str}}(A')$ is an isomorphism. Furthermore,
$\opn{q}_{A' / A}$ is the unique nondegenerate rigid localization morphism over 
$v$. If $A' \to A''$ is another essentially \'etale homomorphism in 
$\ringcfeft{\K}$, then 
$\opn{q}_{A'' / A'} \circ \opn{q}_{A' / A} = \opn{q}_{A'' / A}$.

The two functorialities of the rigid residue complex $\KK_A$ commute with each 
other, in the following sense. Suppose we are given the first commutative 
diagram below in $\ringcfeft{\K}$, which is cocartesian (i.e.\ 
$B' \cong A' \ot_A B$), $A \to B$ is finite, and $A \to A'$ is essentially 
\'etale. Then the second diagram in $\dcat{C}_{\mrm{str}}(A)$ is commutative. 
\[ 
\UseTips \xymatrix @C=6ex @R=6ex {
A
\ar[r]^{}
\ar[d]_{}
&
B
\ar[d]
\\
A' 
\ar[r]
&
B'
}
\qquad \quad 
\UseTips \xymatrix @C=8ex @R=6ex {
\KK_{A}
\ar[d]_{\opn{q}_{A' / A}}
&
\KK_{B}
\ar[l]_{\opn{tr}_{B / A}}
\ar[d]^{\opn{q}_{B' / B}}
\\
\KK_{A'}
&
\KK_{B'}
\ar[l]_{\opn{tr}_{B' / A'}}
} 
\]

Imprecise proofs of these assertions can be found in the papers \cite{ YeZh3} 
and \cite{YeZh4}. A correct treatment will appear in \cite{Ye11}. 
\end{rem}

\begin{rem} \label{rem:4192}  
The functorialities of the rigid residue complexes that were outlined in Remark 
\ref{rem:4190} enable similar construction for schemes. In this remark we shall 
dispense with the flatness condition (invoking Remark \ref{rem:2285}); so  
$\K$ is a nonzero regular noetherian base ring (e.g.\ a field or the ring of 
integers $\Z$), and we consider the category $\catt{Sch} \eftover \K$
of EFT $\K$-schemes. 

Let $X$ be an EFT $\K$-scheme. A {\em rigid residue complex} on $X$ is a pair
$(\KK_X, \bsym{\rho}_X)$, consisting of a bounded complex of injective 
quasi-coherent $\OO_X$-modules $\KK_X$, and a {\em rigid structure} 
$\bsym{\rho}_X$ on $\KK_X$, to be explained below. 
For every affine open set $U \sub X$ which is 
{\em strictly EFT}, namely $A := \Ga(U, \OO_X)$ is an EFT $\K$-ring, the 
complex $\KK_A := \Ga(U, \KK_X)$
is equipped with a rigidifying isomorphism $\rho_U$, making it into a 
rigid residue complex over $A$ relative to $\K$.
If $U' \sub U$ is a smaller strictly EFT open set of $X$, 
then the restriction homomorphism 
$\KK_{A} \to \KK_{A'}$
must be the unique rigid localization homomorphism, for the given 
rigidifying isomorphisms $\rho_{U}$ and $\rho_{U'}$.
The collection of these rigidifying isomorphisms 
$\bsym{\rho}_X := \{ \rho_U \}$
is the rigid structure on $\KK_X$. 
The scheme $X$ admits a rigid residue complex $(\KK_X, \bsym{\rho}_X)$, and it 
is unique up to a unique rigid isomorphism. 

If $f : Y \to X$ is a finite map in $\catt{Sch} \eftover \K$, then 
there is the rigid trace homomorphism 
$\opn{tr}_f : f_{*}(\KK_Y) \to \KK_X$
in $\dcat{C}_{\mrm{str}}(X) := \dcat{C}_{\mrm{str}}(\cat{Mod} \OO_X)$; 
and it is nondegenerate, in the sense that the 
induced homomorphism 
$f_{*}(\KK_Y) \to \Hom_{X} \bigl( f_*(\OO_Y), \KK_X \bigr)$
in $\dcat{C}_{\mrm{str}}(X)$ is an isomorphism. 
If $g : X' \to X$ is an essentially \'etale map in 
$\catt{Sch} \eftover \K$, then 
there is the rigid localization homomorphism 
$\opn{q}_g : \KK_X \to g_{*}(\KK_{X'})$
in $\dcat{C}_{\mrm{str}}(X)$; and it is nondegenerate, in the sense that the 
induced homomorphism 
$g^*(\KK_X) \to \KK_{X'}$
is an isomorphism in $\dcat{C}_{\mrm{str}}(X')$. In this way the rigid residue 
complexes become a sheaf on the small \'etale site of $X$. 

Actually, an infinitesimal local construction gives, for every morphism 
$f : Y \to X$ $\catt{Sch} \eftover \K$, 
the {\em ind-rigid trace homomorphism}
$\opn{tr}_f : f_{*}(\KK_Y) \to \KK_X$
in $\dcat{G}_{\mrm{str}}(X)$; namely this is a homomorphism of graded 
quasi-coherent sheaves on $X$. The {\em Residue Theorem} says that when $f$ is 
proper, the homomorphism $\opn{tr}_f$ commutes with the differentials. 
And the {\em Duality Theorem} says that in the proper case, the ind-rigid trace 
$\opn{tr}_f$ induces an isomorphism 
\[ \mrm{R} f_{*} \bigl( \mrm{R} \Hom_{Y}(\NN, \KK_Y) \bigr) \to  
\mrm{R} \Hom_{X} \bigr( \mrm{R} f_{*}(\NN), \KK_X \bigr) \]
in $\dcat{D}(X)$ for every 
$\NN \in \dcat{D}_{\mrm{c}}^{\mrm{b}}(X)$.
This explicit rigid version of the original theorems of Grothendieck from 
\cite{RD} will appear in \cite{Ye12}. 
\end{rem}

\begin{rem} \label{rem:4193}
In this remark we consider {\em Deligne-Mumford stacks}, of finite type 
over the regular noetherian base ring $\K$. Due to the \'etale functoriality of 
the rigid residue complexes, as explained in Remark \ref{rem:4191}, every such 
stack $\mfrak{X}$ admits a rigid residue complex 
$(\KK_{\mfrak{X}}, \bsym{\rho}_{\mfrak{X}})$.
Here the rigid structure 
$\bsym{\rho}_{\mfrak{X}} = \{ \rho_{(U, g)} \}$
is indexed by the \'etale maps $g : U \to \mfrak{X}$ from strictly EFT affine 
$\K$-schemes $U$. 

To a map $f : \mfrak{Y} \to \mfrak{X}$ between DM stacks 
we associate the ind-rigid trace 
$\opn{tr}_{f} : f_{*}(\KK_{\mfrak{Y}}) \to \KK_{\mfrak{X}}$,
which is a homomorphism of graded quasi-coherent sheaves on $\mfrak{X}$. 
In order to construct the ind-rigid trace for stacks we need another property 
of rigid residue complexes over rings, that we call {\em \'etale codescent}. 
Suppose $v : A \to A'$ is a faithfully \'etale ring homomorphism of EFT 
$\K$-rings. Let $w_1, w_2 : A' \to A' \ot_{A} A'$ be the two inclusions. 
\'Etale codescent says that in every degree $i$ the sequence of $A$-module 
homomorphisms 
\[ \KK^i_{A' \ot_{A} A'} \xar{ \, \opn{tr}_{w_1} - \opn{tr}_{w_2} \, }
\KK^i_{A'} \xar{ \, \opn{tr}_{v} \, } \KK^i_A \to 0 \]
is exact. 

As for the Residue Theorem: we can only prove it for a proper map of DM stacks 
$f : \mfrak{Y} \to \mfrak{X}$ which is {\em coarsely schematic}.
Likewise, we can only prove the Duality Theorem for a proper map of DM stacks 
$f$ which is coarsely schematic and {\em tame}. These conditions, and a sketch 
of the proofs, can be found in the lecture notes \cite{Ye10}. 
These results shall be written in detail in the future paper \cite{Ye13}. 
Observe that duality is not expected to hold without the tameness condition; 
but the coarsely schematic condition seems to be a temporary technical hitch. 
\end{rem}

%% file: block5_190413.tex

\renewcommand{\thisfile}{block5\_190328}  

\cleardoublepage 
\mysection{Perfect and Tilting DG Modules over NC DG Rings} 
\label{sec:perf-tilt-NC}
 
\AYcopyright

Perfect and tilting complexes over noncommutative rings are 
among the important applications of derived categories to ring theory. 
In this section of the book we study the more general concepts of 
perfect DG modules and tilting DG bimodules over noncommutative DG rings. 

Throughout the section we adhere to the following convention, that extends 
Convention \ref{conv:2490}. 

\begin{conv} \label{conv:3485}
There is a fixed a nonzero commutative base ring $\K$ (e.g.\ a field, or 
the ring of integers $\Z$). By default all linear objects and operations are 
$\K$-linear. This means that the DG rings are DG central $\K$-rings (Definition 
\ref{dfn:1071}), the DG modules are $\K$-linear, the DG bimodules 
are central over $\K$, and all homomorphisms between them are $\K$-linear. 
Likewise all linear categories and functors are $\K$-linear.
The symbol $\ot$ stands for $\ot_{\K}$. 
\end{conv}

All the results in this section specialize to rings: the category 
$\catt{Rng} \centover \K$ of {\em central $\K$-rings} (also known as 
associative unital $\K$-algebras) is a full subcategory of the category
$\catt{DGRng} \centover \K$ of DG central $\K$-rings. 

From subsection \ref{subsec:flat-dg-rng} onward we shall impose a flatness 
condition on our DG rings.

\mysubsection{Algebraically Perfect DG Modules}
\label{subsec:perf-dgmods}
Here we define {\em algebraically perfect DG modules}, and prove several of 
their characterizations, in Theorems \ref{thm:3400}, \ref{thm:3415} and 
\ref{thm:3416}. See Remark \ref{rem:3491} regarding nomenclature. 

Semi-free filtrations of DG modules were introduced in Definition 
\ref{dfn:1575}. 

\begin{dfn} \label{dfn:3411}
Let $A$ be a DG ring and $P$ a DG $A$-module. 
\begin{enumerate}
\item Let $F = \bigl\{ F_j(P) \bigr\}_{j \geq -1}$ be a semi-free filtration of 
$P$. We say that $F$ has {\em finite extension length} if 
there is a number $j_1 \in \N$ such that $F_{j_1}(P) = P$.

\item The DG module $P$ is called {\em semi-free of finite extension length} 
if it admits some semi-free filtration $F$ that has finite extension length. 
\end{enumerate}
\end{dfn}

Recall that for a DG module $M$ and an integer $k$ we write  
$M[k] := \opn{T}^k(M)$, the $k$-th translation of $M$. See Definition 
\ref{dfn:2120}.

\begin{dfn} \label{dfn:3402}
Let $A$ be a DG ring and $P$ a DG $A$-module. 
\begin{enumerate}
\item We call $P$ a {\em finite free DG $A$-module} 
\index{Differential graded module! finite free}
if there is an isomorphism 
$P \cong \bigoplus_{i = 1}^r A[-k_i]$
in $\dcat{C}_{\mrm{str}}(A)$ for some $r \in \N$ and $k_i \in \Z$. 

\item A {\em finite semi-free filtration} 
\index{Filtration on a DG module! finite semi-free}
on $P$ is a semi-free filtration 
$F = \bigl\{ F_j(P) \bigr\}_{j \geq -1}$ that has these properties:
\begin{itemize}
\item Each $\opn{Gr}^F_j(P)$ is a finite free DG $A$-module.

\item The filtration $F$ has finite extension length (Definition 
\ref{dfn:3411}(1)). 
\end{itemize}

\item The DG module $P$ is called {\em finite semi-free} 
\index{Differential graded module! finite semi-free}
if it admits 
some finite semi-free filtration. 
\end{enumerate}
\end{dfn}

Saturated full triangulated subcategories, and those generated by a set of 
objects, were introduced in Definitions \ref{dfn:4915} and \ref{dfn:4916}.

\begin{dfn} \label{dfn:3401}
A DG $A$-module $L$ is called {\em algebraically perfect} 
\index{Differential graded module! algebraically perfect}
if $L$ belongs to the saturated full triangulated  subcategory of $\dcat{D}(A)$ 
generated by the DG module $A$. 
\end{dfn}

\begin{rem} \label{rem:3491}
The definition above is new. Earlier texts used the term {\em perfect}, and  
this referred, somewhat ambiguously, to several distinct properties of DG 
modules, or complexes, that are sometimes equivalent to each other. See Theorems 
\ref{thm:3400}, \ref{thm:3415} and \ref{thm:3416}. For commutative DG rings 
there is another notion -- that of a {\em geometrically perfect} DG module; see 
Remark \ref{rem:3425}. 
\end{rem}

The projective dimension of a DG module was introduced in Definition 
\ref{dfn:2128}.

\begin{prop} \label{prop:3415}
If $L \in \dcat{D}(A)$ is algebraically perfect, then it has finite projective 
dimension.
\end{prop}

\begin{proof}
The projective dimension of the DG module $A$ is zero. If $L$ has finite 
projective dimension then so do all its translates $L[k]$. 
If there's a distinguished triangle 
$L' \to L'' \to L \xar{\triangle}$
s.t.\ both $L'$ and $L''$ have finite projective dimension, then so does $L$. 
If $L$ is a direct summand of $L'$, and  $L'$ has finite projective 
dimension, then so does $L$.
Now use Proposition \ref{prop:4470}.
\end{proof}

Recall that the category $\dcat{D}(A)$ has arbitrary direct sums, and they are 
the direct sums in $\dcat{C}_{\mrm{str}}(A)$; see Theorem \ref{thm:3140} and 
Corollary \ref{cor:1665}. 

\begin{dfn} \label{dfn:3405}
A DG $A$-module $L$ is called a {\em compact object of $\dcat{D}(A)$}
\index{Differential graded module! compact}
if the functor  
\[ \opn{Hom}_{\dcat{D}(A)}(L, -) : \dcat{D}(A) \to \dcat{M}(\K) \]
commutes with infinite direct sums. 
\end{dfn}

Explanation: Let $\{ M_z \}_{z \in Z}$ be a collection of objects of 
$\dcat{D}(A)$, with direct sum 
$M := \bigoplus_{z \in Z} M_z$ and embeddings 
$e_z : M_z \to M$. The canonical homomorphism of $\K$-modules 
\begin{equation} \label{eqn:3480}
\begin{aligned}
& \Phi_L : 
\bigoplus_{z \in Z} \, \opn{Hom}_{\dcat{D}(A)}(L, M_z) \to
\opn{Hom}_{\dcat{D}(A)} (L, M) \ ,
\\[0.4em]
& \Phi_L := \bigl\{ \opn{Hom}(\opn{id}_L, e_z) \bigr\}_{z \in Z}
\end{aligned}
\end{equation}
is always injective. Compactness of $L$ says that $\Phi_L$ is bijective.  
Some texts use the name ``small'' instead of ``compact''; see Remark 
\ref{rem:3426}.  

\begin{prop} \label{prop:3545}
Let $A$ and $B$ be DG rings, $F : \dcat{D}(A) \to \dcat{D}(B)$ an 
equivalence of triangulated categories, and $L \in \dcat{D}(A)$. Then $L$ is a 
compact object of $\dcat{D}(A)$ if and only if $F(L)$ is a 
compact object of $\dcat{D}(B)$.  
\end{prop}

\begin{proof}
Suppose $\{ N_z \}_{z \in Z}$ is a collection of objects of $\dcat{D}(B)$. 
Because $F$ is an equivalence of categories, we can assume that 
$N_z = F(M_z)$ for some $M_z \in \dcat{D}(A)$.
The equivalence $F$ respects coproducts. Namely, if 
$M = \bigoplus_{z \in Z} M_z$ in $\dcat{D}(A)$, 
with embeddings $e_z : M_z \to M$, 
then the object $N := F(M) \in \dcat{D}(B)$, with the morphisms 
$F(e_z) : N_z \to N$,
is a coproduct of the collection of objects
$\{ N_z) \}_{z \in Z}$. 
There is a commutative diagram of $\K$-modules
\[  \UseTips  \xymatrix @C=10ex @R=6ex {
\bigoplus_{z \in Z} \opn{Hom}_{\dcat{D}(A)}(L, M_z) 
\ar[r]^(0.56){\Phi_L}
\ar[d]_{\bigoplus F}^{\cong}
&
\opn{Hom}_{\dcat{D}(A)}(L, M) 
\ar[d]^{F}_{\cong}
\\
\bigoplus_{z \in Z} \opn{Hom}_{\dcat{D}(B)} \bigl( F(L), N_z \bigr) 
\ar[r]^(0.56){\Phi_{F(L)}}
&
\opn{Hom}_{\dcat{D}(B)} \bigl( F(L), N \bigr) 
} \]
We see that if $\Phi_L$ is an isomorphism, then so is $\Phi_{F(L)}$.

For the reverse direction we do the same, but now for a quasi-inverse 
$G : \dcat{D}(B) \to \dcat{D}(A)$ of $F$, the object 
$F(L) \in \dcat{D}(B)$ and the object $L \cong G(F(L)) \in \dcat{D}(A)$.
\end{proof}

Here are a few lemmas about compact objects that will be needed for the proof 
of Theorem \ref{thm:3400}. 

\begin{lem} \label{lem:3405}
Let $L$ be a compact object of $\dcat{D}(A)$, let 
$\{ M_z \}_{z \in Z}$ be a collection of objects of $\dcat{D}(A)$, 
let $M := \bigoplus_{z \in Z} \, M_z$, 
and let $\phi : L \to M$
be a morphism in $\dcat{D}(A)$. Then there is a finite subset 
$Z_0 \sub Z$, such that writing 
$M_0 := \bigoplus_{z \in Z_0} \, M_z$, 
and denoting by $\ga_0 : M_0 \to M$ the inclusion, 
there exists a morphism 
$\phi_0 : L \to M_0$ in $\dcat{D}(A)$ satisfying 
$\phi = \opn{Q}(\ga_0) \circ \phi_0$. 
\end{lem}

\begin{proof}
We are given 
$\phi \in \opn{Hom}_{\dcat{D}(A)}(L, M)$.
Using the canonical isomorphism (\ref{eqn:3480})
that we get by the compactness of $L$, 
we can express $\phi$ as a finite sum, say  
$\phi = \Phi_L \bigl( \sum\nolimits_{z \in Z_0} \phi_z \bigr)$,
where $Z_0 \sub Z$ is a finite subset, and 
$\phi_z \in \opn{Hom}_{\dcat{D}(A)}(L, M_z)$.
The morphism 
\[ \phi_0 := \sum\nolimits_{z \in Z_0} \phi_z : L \to 
\bigoplus\nolimits_{z \in Z_0} \, M_z = M_0 \]
in $\dcat{D}(A)$ has the required property.
\end{proof}

\begin{rem} \label{rem:3405}
It might be helpful to say some words about the lemma. 
The inclusion $\ga_0 : M_0 \to M$ is a (split) monomorphism
in the abelian category $\dcat{C}_{\mrm{str}}(A)$. Suppose we choose a 
K-projective resolution $P \to L$. The morphisms $\phi : L \to M$
and $\phi_0 : L \to M_0$ in $\dcat{D}(A)$ are represented by 
homomorphisms 
$\til{\phi} : P \to M$ and $\til{\phi}_0 : P \to M_0$
in $\dcat{C}_{\mrm{str}}(A)$. 
There is no reason for the homomorphism $\til{\phi}$ to factor through $M_0$
in $\dcat{C}_{\mrm{str}}(A)$. 
All we can say is that {\em there is a homotopy}
$\til{\phi} \twoto \ga_0 \circ \til{\phi}_0$. 

Later, after proving Theorem \ref{thm:3400}, we will know that there is a 
finite semi-free DG $A$-module $Q$, with homomorphisms 
$P \xar{\al} Q \xar{\be} P$ in $\dcat{C}_{\mrm{str}}(A)$,
and with a homotopy 
$\be \circ \al \twoto \opn{id}_P$.  
The homomorphism 
$\til{\phi} \circ \be : Q \to M$
does factor through some finite direct sum $M_0$, and hence so does 
$\til{\phi} \circ \be \circ \al : P \to M$. 
The homotopy $\be \circ \al \twoto \opn{id}_P$
induces the homotopy $\til{\phi} \twoto \ga_0 \circ \til{\phi}_0$
that was mentioned above. 
\end{rem}

The next definition is taken from \cite{BoNe}. 

\begin{dfn} \label{dfn:3410}
Suppose 
$\bigl( \{ M_i \}_{i \in \N}, \{ \mu_i \}_{i \in \N} \bigr)$ 
is a direct system in $\dcat{D}(A)$. Let 
\[ \phi \, : \, \bigoplus\nolimits_{i \in \N} \, M_i \ \to \ 
\bigoplus\nolimits_{i \in \N} \, M_i \]
be the morphism in $\dcat{D}(A)$ defined by 
\[ \phi|_{M_i} := (\opn{id}, -\mu_i) : M_i \to M_i \oplus M_{i + 1} . \]
A {\em homotopy colimit} 
\index{Homotopy colimit}
of the direct system $\{ M_i \}_{i \in \N}$
is an object $M \in \dcat{D}(A)$ that is a cone on $\phi$, in the sense of 
Definition \ref{dfn:3600}. Namely $M$ sits in some distinguished triangle 
\[ \bigoplus\nolimits_{i \in \N} \, M_i \xar{\ \phi \ } 
\bigoplus\nolimits_{i \in \N} \, M_i 
\xar{\psi} M \xar{\ \triangle \ } \]
in $\dcat{D}(A)$.
\end{dfn}

Observe that a homotopy colimit exists, and it is unique up to a 
non\-unique isomorphism in $\dcat{D}(A)$. See Corollary \ref{cor:3600}.

\begin{lem} \label{lem:3610}
If 
$\bigl( \{ M_i \}_{i \in \N}, \{ \mu_i \}_{i \in \N} \bigr)$
is a direct system in $\dcat{C}_{\mrm{str}}(A)$, and we take 
$M := \lim_{i \to} M_i$,
the direct limit in $\dcat{C}_{\mrm{str}}(A)$, then $M$ is a 
homotopy colimit of the direct system 
$\bigl( \{ M_i \}_{i \in \N}, \{ \opn{Q}(\mu_i) \}_{i \in \N} \bigr)$ 
in $\dcat{D}(A)$
\end{lem}

\begin{proof}
Let 
\[ \til{\phi} \, : \, \bigoplus\nolimits_{i \in \N} \, M_i \ \to \ 
\bigoplus\nolimits_{i \in \N} \, M_i \]
be the morphism in $\dcat{C}_{\mrm{str}}(A)$ defined like in Definition 
\ref{dfn:3410} above; so the morphism $\phi$ in that definition equals 
$\opn{Q}(\til{\phi})$. 

An easy calculation shows that 
\[ 0 \to \bigoplus\nolimits_{i \in \N} \, M_i \xar{\ \til{\phi} \ } 
\bigoplus\nolimits_{i \in \N} \, M_i \xar{\ \til{\psi} \ } M \to 0 , \] 
where $\til{\psi}$ is the canonical homomorphism, is an exact sequence in 
$\dcat{C}_{\mrm{str}}(A)$. 
By Proposition \ref{prop:2165} there is a distinguished triangle 
\[ \bigoplus\nolimits_{i \in \N} \, M_i \xar{\ \opn{Q}(\til{\phi}) \ } 
\bigoplus\nolimits_{i \in \N} \, M_i \xar{\ \opn{Q}(\til{\psi}) \ } M 
\xar{ \, \triangle \, } \]
in $\dcat{D}(A)$. 
\end{proof}

\begin{lem} \label{lem:3410} 
Let $L$ be a compact object of $\dcat{D}(A)$, and 
let $\{ M_i \}_{i \in \N}$ be a direct system in $\dcat{D}(A)$, with homotopy 
colimit $M$. Then the morphism $\psi$ in Definition \tup{\ref{dfn:3410}}
induces an isomorphism
\[ \lim_{i \to} \, \opn{Hom}_{\dcat{D}(A)}(L, M_i) \iso 
\opn{Hom}_{\dcat{D}(A)}(L, M) \]
in $\dcat{M}(\K)$. 
\end{lem}

\begin{proof}
Applying the cohomological functor $\opn{Hom}_{\dcat{D}(A)}(L, -)$
to the distinguished triangle in Definition \ref{dfn:3410}, we get 
long exact sequence of $\K$-modules
\begin{equation}  \label{eqn:4790}
\begin{aligned}
&
\cdots \xar{\de^{k - 1}}
\opn{Hom}_{\dcat{D}(A)} \Bigl( L, 
\Bigl( \bigoplus\nolimits_{i \in \N} \, M_i \Bigr)[k] \Bigr)
\\
& \quad
\xar{ \opn{Hom}(L, \phi[k]) }
\opn{Hom}_{\dcat{D}(A)} \Bigl( L, 
\Bigl( \bigoplus\nolimits_{i \in \N} \, M_i \Bigr)[k] \Bigr)
\\
& \quad  
\xar{ \opn{Hom}(L, \psi[k]) }
\opn{Hom}_{\dcat{D}(A)} \bigl( L, M[k] \bigr) \xar{\de^{k}}  \cdots .
\end{aligned}
\end{equation}
Using the compactness of $L$ and the resulting isomorphism $\Phi_L$, we obtain 
from (\ref{eqn:4790}) the long exact sequence of $\K$-modules
\begin{equation}  \label{eqn:3411}
\begin{aligned}
&
\cdots \xar{\de^{k - 1}} 
\bigoplus\nolimits_{i \in \N} \, \opn{Hom}_{\dcat{D}(A)} 
\bigl( L, M_i[k] \bigr) 
\\
& \quad  
\xar{\phi^k} \bigoplus\nolimits_{i \in \N} \, 
\opn{Hom}_{\dcat{D}(A)} \bigl( L, M_i[k] \bigr) 
\\
& \quad  
\xar{\psi^k} \opn{Hom}_{\dcat{D}(A)} \bigl( L, M[k] \bigr)  \xar{\de^{k}}  
\cdots . 
\end{aligned}
\end{equation}
We know that for every $k$ the homomorphism $\phi^k$ above is 
injective. This implies that the connecting homomorphisms $\de^k$
in (\ref{eqn:3411}) are all zero. Hence, when we take $k = 0$, and we 
replace the two occurrences of  
``$\cdots$'' in (\ref{eqn:3411}) by ``$0$'', we get a short exact sequence.  
Finally, a calculation like in the proof of Lemma \ref{lem:3610} shows that 
\[ \opn{Coker}(\phi^0) \cong
\lim_{i \to} \, \opn{Hom}_{\dcat{D}(A)}(L, M_i) . \qedhere \]
\end{proof}

\begin{lem} \label{lem:3406} 
Let $L$ be a compact object of $\dcat{D}(A)$, let $M$ be a DG $A$-module, 
and let $\{ F_j(M) \}_{j \geq -1}$ be a filtration of $M$ in
$\dcat{C}_{\mrm{str}}(A)$ such that 
$\bigcup_j F_j(M) = M$. Given a morphism 
$\phi : L \to M$ in $\dcat{D}(A)$, there is an index $j_1 \in \N$, and a 
morphism $\phi_{j_1} : L \to F_{j_1}(M)$ in $\dcat{D}(A)$, 
such that 
$\phi = \opn{Q}(\ga_{j_1}) \circ \phi_{j_1}$,
where $\ga_{j_1} : F_{j_1}(M) \to M$ is the inclusion. 
\end{lem}

\begin{proof}
Here $M \cong \lim_{j \to} \, F_j(M)$ in $\dcat{C}_{\mrm{str}}(A)$. By Lemma 
\ref{lem:3610} we know that $M$ is a homotopy colimit of the direct system 
$\{ F_j(M) \}_{j \geq -1}$. According to Lemma \ref{lem:3410} we 
get a bijection 
\begin{equation} \label{eqn:3420}
\lim_{j \to} \, \opn{Hom}_{\dcat{D}(A)}(L, F_j(M)) \cong  
\opn{Hom}_{\dcat{D}(A)}(L, M) . 
\end{equation}
So there is some $j_1 \in \N$ and some 
$\phi_{j_1} \in \opn{Hom}_{\dcat{D}(A)}(L, F_{j_1}(M))$
that goes to $\phi$ under the bijection (\ref{eqn:3420}).
\end{proof}

The next lemma is adapted from \cite[Lemma tag=09R2]{SP}. 

\begin{lem} \label{lem:3412}
Let $L$ be a compact object of $\dcat{D}(A)$, let $P$ be a semi-free
DG $A$-module of finite extension length, and let $\phi : L \to P$ be a 
morphism 
in $\dcat{D}(A)$.
Then there is a finite semi-free DG $A$-module $P'$, with a morphism 
$\phi' : L \to P'$ in $\dcat{D}(A)$ and a monomorphism 
$\ga' : P' \to P$ in $\dcat{C}_{\mrm{str}}(A)$, such that 
$\phi = \opn{Q}(\ga') \circ \phi'$. 
\end{lem}

\begin{proof}
Let $\{ F_j(P) \}_{j \geq -1}$ be a semi-free filtration of $P$ of finite 
extension length; say $F_{j_1}(P) = P$. We will construct DG submodules
\[ P_0 \sub P_1 \sub \cdots \sub P_{j_1} \sub P_{j_1 + 1} \]
of $P$, such that properties (a)-(c) below hold for every 
$k \in [0, j_1 + 1]$. For each such $k$ and each $j \geq -1$ we define  
$F_j(P_k) := F_j(P) \cap P_k$.
\begin{itemize}
\item[(a)] The filtration $\{ F_j(P_k) \}_{j \geq -1}$ is a semi-free 
filtration on $P_k$.

\item[(b)] The free DG $A$-module $\opn{Gr}^F_j(P_k)$ is finite free for all 
$j \geq k$. 

\item[(c)] The morphism $\phi$ factors through $P_k$; namely there is a 
morphism $\phi_k : L \to P_k$ in $\dcat{D}(A)$ s.t.\ 
$\phi = \opn{Q}(\ga_k) \circ \phi_k$, where $\ga_k : P_k \to P$ is the 
inclusion. 
\end{itemize}

For $k = 0$, the DG module $P' := P_0$ is finite semi-free. The morphism 
$\phi' := \phi_0 : L \to P'$
in $\dcat{D}(A)$ and the inclusion 
$\ga' := \ga_0 : P' \to P$
in $\dcat{C}_{\mrm{str}}(A)$ 
will satisfy $\phi = \opn{Q}(\ga') \circ \phi'$. 

The construction of the DG modules $P_k$ and the morphisms $\phi_k$ is by 
descending induction on $k$, starting from $k = j + 1$. 
We start with $P_{j_1 +1 } := P$ and $\phi_{j_1 + 1} := \phi$. 
Properties (a)-(c) hold trivially here. 

Now to the inductive step. Suppose that $k \in [1, \ldots, j + 1]$, 
and we already found a DG submodule 
$P_k \sub P$ and a morphism $\phi_k : L \to P_k$ as required. 
We are going to produce a DG submodule 
$P_{k - 1} \sub P_k$ such that 
$\opn{Gr}^F_{j}(P_{k - 1}) = \opn{Gr}^F_{j}(P_k)$
for all $j \neq k - 1$, and 
$\opn{Gr}^F_{k - 1}(P_{k - 1})$ is a finite free DG submodule of 
$\opn{Gr}^F_{k - 1}(P_k)$. 
This means that we have to shrink $F_{k - 1}(P_k)$ suitably.

For every $j \geq 0$ the DG $A$-module $\opn{Gr}^F_{j}(P_k)$ is free, and 
we choose a basis for it, namely a collection 
$\{ \bar{p}_z \}_{z \in Z_j}$
of elements such that 
$\opn{Gr}^F_{j}(P_k) = \bigoplus_{z \in Z_j} A \cd \bar{p}_z$
in $\dcat{C}_{\mrm{str}}(A)$. 
Moreover, we choose these bases such that the indexing set 
$Z_j$ is finite for every $j \geq k$ (this is possible by property (b)), 
and $Z_j$ is empty for every $j \geq j_1 + 1$ (due to the fact that 
$F_{j_1}(P) = P$). 

By definition 
$\opn{Gr}^F_{j}(P_k) = F_{j}(P_k) / F_{j - 1}(P_k)$,
so we can lift the basis 
$\{ \bar{p}_z \}_{z \in Z_j}$ of $\opn{Gr}^F_{j}(P_k)$ to 
a collection $\{ p_z \}_{z \in Z_j}$ of homogeneous elements 
$p_z \in F_{j}(P_k)$. 
Let's define the indexing set $Z := \coprod_{j \geq k} Z_j$.
Then there are decompositions
\begin{equation} \label{eqn:5085}
P^{\natural}_k = F_{k - 2}(P_k)^{\natural} \oplus
\Bigl( \bigoplus\nolimits_{z \in Z_{k - 1}} \, A^{\natural} \cd p_z \Bigr) 
\oplus 
\Bigl( \bigoplus\nolimits_{z \in Z} \, A^{\natural} \cd p_z \Bigr) 
\end{equation}
and 
\[ F_{k - 1}(P_k)^{\natural} = F_{k - 2}(P_k)^{\natural} \oplus 
\Bigl( \bigoplus\nolimits_{z \in Z_{k - 1}} \, A^{\natural} \cd p_z \Bigr)   \]
in $\dcat{G}_{\mrm{str}}(A^{\natural})$, and each $A^{\natural} \cd p_z$ is a 
free graded $A^{\natural}$-module of rank $1$. 

Because $Z$ is a finite set, there is a finite subset $Y \sub Z_{k - 1}$ such 
that 
\[ \d(p_z) \in F_{k - 2}(P_k) \oplus
\Bigl( \bigoplus\nolimits_{y \in Y} \, A \cd p_y \Bigr) \oplus 
\Bigl( \bigoplus\nolimits_{z' \in Z} \, A \cd p_{z'} \Bigr) \]
for every $z \in Z$. 
For every $z \in Z_{k - 1}$ the element 
$\bar{p}_z \in \opn{Gr}^F_{k - 1}(P_k)$ satisfies 
$\d(\bar{p}_z) = 0$, and therefore the element 
$p_z \in F_{k - 1}(P_k)$ satisfies $\d(p_z) \in F_{k - 2}(P_k)$.
In particular this is true for $p_y$ with $y \in Y$. 
And $F_{k - 2}(P_k)$ is a DG submodule of $P_k$.
Combining all these observations we see that the graded submodule
\begin{equation} \label{eqn:5086}
Q^{\natural} := F_{k - 2}(P_k)^{\natural} \oplus
\Bigl( \bigoplus\nolimits_{y \in Y} \, A^{\natural} \cd p_y \Bigr) \oplus 
\Bigl( \bigoplus\nolimits_{z \in Z} \, A^{\natural} \cd p_z \Bigr) 
\end{equation}
of $P_k^{\natural}$ actually underlies a DG submodule $Q \sub P_k$. 
(Compare this to the decomposition (\ref{eqn:5085}) of $P_k^{\natural}$.)

Consider the quotient DG $A$-module $N := P_k / Q$,
with the canonical projection $\al: P_k \to N$, that's an epimorphism in 
$\dcat{C}_{\mrm{str}}(A)$. 
Define the set $X := Z_{k - 1} - Y$
and the elements 
$n_x := \al(p_x) \in N$. 
Then $N^{\natural}$ is a free graded $A^{\natural}$-module, with basis 
indexed by $X$, i.e.\ 
$N^{\natural} = \bigoplus\nolimits_{x \in X} A^{\natural} \cd n_x$.
Furthermore, because 
$\d(p_x) \in F_{k - 2}(P_k) \sub Q$
for all $x \in X$, we have $\d(n_x) = 0$;
so actually $N$ is a free DG $A$-module with basis $\{ n_x \}_{x \in 
X}$. 

Next consider the morphism
$\th : L \to N$ in $\dcat{D}(A)$ with formula 
$\th := \opn{Q}(\al) \circ \phi_k$.
By Lemma \ref{lem:3405} there is a finite subset $X_0 
\sub X$ such that $\th$ factors through the DG submodule
$N_0 := \bigoplus\nolimits_{x \in X_0} A \cd n_x$. 
I.e.\ there's a morphism $\th_0 : L \to N_0$ in $\dcat{D}(A)$
s.t.\ $\th = \opn{Q}(\be_0) \circ \th_0$, where 
$\be_0 : N_0 \to N$ is the inclusion (a monomorphism in 
$\dcat{C}_{\mrm{str}}(A)$). Define 
$\bar{N} := N / N_0 = \opn{Coker}(\be_0)$,
and let $\de : N \to \bar{N}$ be the projection
(an epimorphism in $\dcat{C}_{\mrm{str}}(A)$). 
Since $\de \circ \be_0 = 0$ in $\dcat{C}_{\mrm{str}}(A)$, 
it follows that 
\[ \opn{Q}(\de) \circ \th = \opn{Q}(\de) \circ \opn{Q}(\be_0) \circ \th_0 = 
\opn{Q}(\de \circ \be_0) \circ \th_0 = 0 \]
in $\dcat{D}(A)$. 

The subset $Y \cup X_0$ of $Z_{k - 1}$ is finite. 
We define 
\begin{equation} \label{eqn:3415}
P_{k - 1}^{\natural} := F_{k - 2}(P_k)^{\natural} \oplus
\Bigl( \bigoplus\nolimits_{z \in Y \cup X_0} \, A^{\natural}\cd p_z \Bigr) 
\oplus \Bigl( \bigoplus\nolimits_{z \in Z} \, A^{\natural} \cd p_z \Bigr) . 
\end{equation}
Then $P_{k - 1}^{\natural}$ is the underlying graded module of a DG submodule 
$P_{k - 1} \sub P_k$, slightly bigger than $Q$.
(Compare this formula to (\ref{eqn:5085}) and (\ref{eqn:5086}).)
By construction the DG module $P_{k - 1}$ has properties (a)-(b) 
with index $k - 1$. Indeed, 
\[ \opn{Gr}^F_j(P_{k - 1}) \cong 
\begin{cases}
\opn{Gr}^F_j(P_{k}) & \tup{if} \ j \geq k \ \tup{or} \ j \leq k - 2 
\\[0.4em]
\bigoplus\nolimits_{z \in Y \cup X_0} \, A \cd \bar{p}_z 
& \tup{if} \ j = k - 1 
\end{cases} \]

Let us denote by $\ep : P_{k - 1} \to P_{k}$ the inclusion. 
There is a short exact sequence
$0 \to P_{k - 1} \xar{\ep} P_k \xar{\de \circ \al} \bar{N} \to  0$
in $\dcat{C}_{\mrm{str}}(A)$, and this becomes a distinguished triangle 
$P_{k - 1} \xar{\ep} P_k \xar{\de \circ \al} \bar{N} \xar{\triangle}$
in $\dcat{D}(A)$. Applying the cohomological functor 
$\opn{Hom}_{\dcat{D}(A)}(L, -)$ to this distinguished triangle gives rise to a 
long exact sequence 
\[ \cdots \to \opn{Hom}_{\dcat{D}(A)}(L, P_{k - 1}) \to 
\opn{Hom}_{\dcat{D}(A)}(L, P_{k}) \to 
\opn{Hom}_{\dcat{D}(A)}(L, \bar{N}) \to \cdots . \]
The morphism $\phi_k$ belongs to the middle term here. 
We know that 
\[ \opn{Q}(\de \circ \al) \circ \phi_k = 
\opn{Q}(\de) \circ \opn{Q}(\al) \circ \phi_k = 
\opn{Q}(\de) \circ \th = 0 . \]
Thus there is a morphism 
$\phi_{k - 1} : L \to P_{k - 1}$
in $\dcat{D}(A)$ s.t.\ 
$\phi_k = \opn{Q}(\ep) \circ \phi_{k - 1}$. The morphism $\phi_{k - 1}$
has property (c) with index $k - 1$. 
\end{proof}

The derived tensor-evaluation morphism $\opn{ev}^{\mrm{R, L}}_{L, M, N}$ was 
introduced in Theorem \ref{thm:4320}. 
Finite semi-free, algebraically perfect and compact DG modules were defined 
in Definitions \ref{dfn:3402}(3), \ref{dfn:3401} and 
\ref{dfn:3405}, respectively. 

\begin{thm} \label{thm:3400}
Let $A$ be a DG ring and let $L$ be a DG $A$-module. The following four 
conditions are equivalent.
\begin{itemize}
\rmitem{i} $L$ is algebraically perfect. 

\rmitem{ii} $L$ is a direct summand in $\dcat{D}(A)$ of a finite 
semi-free DG 
$A$-module. 

\rmitem{iii} For every DG ring $B$, every DG module
$M \in \dcat{D}(A \ot_{} B^{\mrm{op}})$ and every DG module
$N \in \dcat{D}(B)$, the tensor-evaluation morphism 
\[ \opn{ev}^{\mrm{R, L}}_{L, M, N} : 
\opn{RHom}_{A}(L, M) \ot^{\mrm{L}}_{B} N \to 
\opn{RHom}_{A}(L, M \ot^{\mrm{L}}_{B} N) \]
in $\dcat{D}(\K)$ is an isomorphism. 

\rmitem{iv} $L$ is a compact object of $\dcat{D}(A)$.
\end{itemize}
\end{thm}

\begin{proof} \mbox{} 

\smallskip \noindent
(ii) $\Rightarrow$ (i): We use Proposition \ref{prop:4470}. A finite free DG 
$A$-module is algebraically perfect. A finite semi-free DG module is 
obtained by finitely many cones from finite free DG modules, so it is also 
algebraically perfect. Finally, our DG module $L$ is a direct summand of a 
finite semi-free DG module, so it is algebraically perfect.

\medskip \noindent 
(i) $\Rightarrow$ (iii): Fixing $M$ and $N$, we have two triangulated functors 
$F, G : \dcat{D}(A) \to \dcat{D}(\K)$,
namely 
$F := \opn{RHom}_{A}(-, M) \ot^{\mrm{L}}_{B} N$
and
$G := \opn{RHom}_{A}(-, M \ot^{\mrm{L}}_{B} N)$.
There is a morphism of triangulated functors 
$\ep_{} := \opn{ev}^{\mrm{R, L}}_{(-), M, N} : F \to  G$,
and we want to prove that $\ep_L : F(L) \to  G(L)$ is an isomorphism.
Trivially for $L := A$ the morphism $\ep_{L}$ is an isomorphism. 
According to Proposition \ref{prop:4603} the morphism $\ep_L$ is an isomorphism 
for every algebraically perfect DG module $L$. 

\medskip \noindent 
(iii) $\Rightarrow$ (iv): Let $\{ N_z \}_{z \in Z}$ be a collection of objects 
of $\dcat{D}(A)$. Define $N := \bigoplus_{z \in Z} N_z$.
We have to prove that the canonical homomorphism of $\K$-modules 
\[ \Phi_L : \bigoplus\nolimits_{z \in Z} \opn{Hom}_{\dcat{D}(A)}(L, N_z) \to 
\opn{Hom}_{\dcat{D}(A)} (L, N) \]
is bijective. Take $B := A$ and 
$M := A \in \dcat{D}(A \ot_{} A^{\mrm{op}})$. 
Consider the commutative diagram of $\K$-modules (\ref{eqn:4800}).
\begin{figure}
\begin{equation} \label{eqn:4800} 
\UseTips  \xymatrix @C=6ex @R=6ex {
\bigoplus_{z \in Z} \opn{H}^0 \bigl( \opn{RHom}_A(L, A) \ot^{\mrm{L}}_{A} N_z
\bigr) 
\ar[r]^(0.56){\th}_(0.56){\cong}
\ar[d]_{\bigoplus \opn{H}^0(\opn{ev}^{\mrm{R, L}}_{L, M, N_z})}^{\cong}
&
\opn{H}^0 \bigl( \opn{RHom}_A(L, A) \ot^{\mrm{L}}_{A} N \bigr) 
\ar[d]^{\opn{H}^0(\opn{ev}^{\mrm{R, L}}_{L, M, N})}_{\cong}
\\
\bigoplus_{z \in Z}  \opn{H}^0 \bigl( \opn{RHom}_A(L, A \ot^{\mrm{L}}_{A} N_z)
\bigr) 
\ar[r]
\ar[d]^{\cong}_{\bigoplus \ze_{N_z}}
&
\opn{H}^0 \bigl( \opn{RHom}_A(L, A \ot^{\mrm{L}}_{A} N) \bigr) 
\ar[d]_{\cong}^{\ze_{N}}
\\
\bigoplus_{z \in Z} \opn{Hom}_{\dcat{D}(A)}(L, N_z)
\ar[r]^(0.55){\Phi_L}
&
\opn{Hom}_{\dcat{D}(A)}(L, N)
}
\end{equation}
\end{figure}
The homomorphism $\th$ is bijective because both $(- \ot^{\mrm{L}}_{A} -)$ 
and $\opn{H}^0$ respect all direct sums 
(see Theorem \ref{thm:3140} and Propositions \ref{prop:1585} and 
\ref{prop:4300}).
The homomorphisms $\opn{H}^0(\opn{ev}^{\mrm{R, L}}_{L, M, N_z})$ and 
$\opn{H}^0(\opn{ev}^{\mrm{R, L}}_{L, M, N})$ are bijective by condition (iii). 
The homomorphisms $\ze_{N_z}$ and $\ze_{N}$ are bijective because 
$A \ot^{\mrm{L}}_{A} (-) \cong \opn{Id}$ and by Corollary \ref{cor:2120}. 
We conclude that $\Phi_L$ is bijective. 

\medskip \noindent 
(iv) $\Rightarrow$ (ii): 
Choose a semi-free resolution $\rho : P \to L$ in $\dcat{C}_{\mrm{str}}(A)$.
We have an isomorphism 
$\phi := \opn{Q}(\rho)^{-1} : L \to P$
in $\dcat{D}(A)$. Let $\{ F_j(P) \}_{j \geq -1}$ be a semi-free filtration on 
$P$. According to Lemma \ref{lem:3406} there is an index $j_1$ such that, 
letting $P' := F_{j_1}(P)$, and letting $\ga' : P' \to P$ be the inclusion, 
there exists a morphism $\phi' : L \to P'$ in $\dcat{D}(A)$ satisfying 
$\phi = \opn{Q}(\ga') \circ \phi'$. 

Now the DG module $P'$ is semi-free of finite extension length. Lemma 
\ref{lem:3412} says that there is a  finite semi-free DG $A$-module $P''$, with 
a morphism $\phi'' : L \to P''$ in $\dcat{D}(A)$, and a monomorphism 
$\ga'' : P'' \to P'$ in $\dcat{C}_{\mrm{str}}(A)$, such that 
$\phi' = \opn{Q}(\ga'') \circ \phi''$. 

Consider the morphism 
$\psi'' : P'' \to L$
in $\dcat{D}(A)$ defined by 
$\psi'' := \opn{Q}(\rho) \circ \opn{Q}(\ga') \circ \opn{Q}(\ga'')$.
We have equality 
\[ \psi'' \circ \phi'' = 
\opn{Q}(\rho) \circ \opn{Q}(\ga') \circ \opn{Q}(\ga'') \circ \phi'' = 
\opn{Q}(\rho) \circ \phi =  \opn{id}_L . \]
Thus $L$ is a retract in $\dcat{D}(A)$ of the finite semi-free DG module $P''$. 
But by Theorem \ref{thm:4160}, $L$ is then a direct summand of $P''$ in
$\dcat{D}(A)$.
\end{proof}

\begin{cor} \label{cor:3480}
Let $A$ and $B$ be DG rings, let 
$F : \dcat{D}(A) \to \dcat{D}(B)$ 
be an equivalence of triangulated categories, and let 
$L \in \dcat{D}(A)$. Then $L$ is an algebraically perfect DG $A$-module if and 
only if $F(L)$ is an algebraically perfect DG $B$-module.
\end{cor}

\begin{proof}
Combine Theorem \ref{thm:3400} and Proposition \ref{prop:3545}. 
\end{proof}

\begin{rem} \label{rem:3425}
Assume $A$ is a {\em commutative DG ring} (i.e.\ nonpositive and strongly 
commutative, see Definitions \ref{dfn:3091} and \ref{dfn:3090})  .
Then $\bar{A} := \opn{H}^0(A)$ is a commutative ring, and there is a DG ring 
homomorphism $A \to \bar{A}$. 
As explained in \cite{Ye10}, 
given a multiplicatively closed set $\bar{S} \sub \bar{A}$, letting 
$S \sub A^0$ be the preimage of $\bar{S}$, the localized commutative DG ring 
$A_S := A \ot_{A^0} A^0_S$ satisfies 
$\opn{H}^0(A_S) = \bar{A}_{\bar{S}}$. 
We are going to use the convenient notation $A_{\bar{S}} := A_S$.
In case $\bar{S} = \{ \bar{s}^k \}_{k \in \N}$
for  some element $\bar{s} \in \bar{A}$, we shall write 
$A_{\bar{s}} := A_{\bar{S}}$.

A sequence of elements $(\bar{s}_1, \ldots, \bar{s}_n)$ in $\bar{A}$ is called 
a {\em covering sequence} if 
$\opn{Spec}(\bar{A}) = \bigcup_{i} \, \opn{Spec}(\bar{A}_{\bar{s}_i})$.
This just means that 
$\bar{A} = \sum_{i} \bar{A} \cd \bar{s}_i$. 
For each $i$ there is a localized DG ring $A_{\bar{s}_i}$. 

A DG module $L \in \dcat{D}(A)$ is called {\em geometrically perfect} 
\index{Differential graded module! geometrically perfect}
if for some covering sequence $(\bar{s}_1, \ldots, \bar{s}_n)$ of $\bar{A}$
there are finite semi-free DG $A_{\bar{s}_i}$-modules $P_i$, and isomorphisms 
$P_i \cong A_{\bar{s}_i} \ot_A L$ in $\dcat{D}(A_{\bar{s}_i})$. 

In case $A$ is a commutative ring, for which $A = A^0 = \bar{A}$, the 
definition becomes simpler: a complex $L \in \dcat{D}(A)$ is geometrically 
perfect if there is a covering sequence $(s_1, \ldots, s_n)$ of 
$A$, bounded complexes of finite free $A_{s_i}$-modules $P_i$,
and isomorphisms $P_i \cong A_{s_i} \ot_A L$ in $\dcat{D}(A_{s_i})$. 
This is the classical definition of a perfect complex of $A$-modules, see 
\cite{SGA6}. 

According to \cite[Corollary 5.21]{Ye10}, a DG $A$-module $L$ is geometrically 
perfect iff it is algebraically perfect. And according to 
\cite[Theorem 5.11]{Ye10}, $L$ is geometrically perfect iff 
$L \in \dcat{D}^-(A)$, and its derived reduction 
$\bar{A} \ot^{\mrm{L}}_{A} L \in \dcat{D}(\bar{A})$ is geometrically perfect.
\end{rem}

In the special case of a ring we can say more. 

\begin{thm} \label{thm:3415}
Let $A$ be a ring and let $L$ be a DG $A$-module. The following two conditions 
are equivalent.
\begin{itemize}
\rmitem{i}  $L$ is an algebraically perfect DG $A$-module. 

\rmitem{ii} $L$ is isomorphic in $\dcat{D}(A)$ to a bounded complex of 
finitely generated projective $A$-modules.
\end{itemize}
\end{thm}

The proof is after this lemma. The concentration $\opn{con}(M)$ of a graded 
module was defined in Definition \ref{dfn:2125}. 

\begin{lem} \label{lem:3415}
In the situation of the theorem, if $L$ is algebraically perfect
and $\opn{H}(L) \neq 0$, then 
$\opn{con}(\opn{H}(L)) = [i_0, i_1]$ for some integers $i_0 \leq i_1$, and 
$\opn{H}^{i_1}(L)$ is a finitely presented $A$-module. 
\end{lem}

\begin{proof}
According to Theorem \ref{thm:3400}, $L$ is a direct summand in $\dcat{D}(A)$ 
of a finite semi-free DG $A$-module $P$. But $A$ is a ring, so $P$ is a bounded 
complex of (finitely generated free) $A$-modules. This implies that 
$\opn{H}(P)$ is bounded. But $\opn{H}(L)$ is a direct summand, in 
$\dcat{G}_{\mrm{str}}(A)$, of $\opn{H}(P)$; so it is also bounded. 
Because $\opn{H}(L)$ is nonzero, its concentration is a bounded nonempty integer
interval $[i_0, i_1]$. Note that 
$i_1 = \opn{sup}(\opn{H}(L))$. 

Corollary \ref{cor:3340}
says that there is a semi-free resolution $Q \to L$ with 
$\opn{sup}(Q) = i_1$. As in the proof of the implication 
(iv) $\Rightarrow$ (ii) in the proof of Theorem \ref{thm:3400}, 
there is a finite semi-free DG $A$-module $Q''$, with a monomorphism 
$Q'' \to Q$ in $\dcat{C}_{\mrm{str}}(A)$, 
such that $L$ is a direct summand of $Q''$ in $\dcat{D}(A)$. 
Now $Q''$ is a bounded complex of finitely generated free $A$-modules, with 
$\opn{sup}(Q'') = i_1$, and hence $\opn{H}^{i_1}(Q'')$ is a finitely 
presented $A$-module. But $\opn{H}^{i_1}(L)$ is a direct summand of 
$\opn{H}^{i_1}(Q'')$, so it too is a finitely presented $A$-module.
\end{proof}

\begin{proof}[Proof of the Theorem]
The implication (ii) $\Rightarrow$ (i) is easy; compare to the proof of the 
implication (ii) $\Rightarrow$ (i) of Theorem \ref{thm:3400}.

The implication (i) $\Rightarrow$ (ii) will be proved in two steps. 

\smallskip \noindent
Step 1. Here $L$ is an $A$-module, which is algebraically perfect as a DG 
$A$-module. According to Proposition \ref{prop:3415}, the DG module $L$ has 
finite projective dimension, say $d \in \N$. (The case $L \cong 0$ can be 
excluded.) This is also the projective dimension of the module $L$ in the 
classical sense (see Proposition \ref{prop:4290}). 
By Lemma \ref{lem:3415} the module $L$ is finitely generated. 
We will proceed by induction on $d$. 

If $d = 0$ then $P := L$ is already a finitely generated projective $A$-module. 
If $d \geq 1$, then we can find a finite free $A$-module $P^0$, with a 
surjection $\eta : P^0 \to L$ in $\dcat{M}(A)$. Let 
$L' := \opn{Ker}(\eta)$, so we have a short exact sequence 
$0 \to L' \to P^0 \xar{\eta} L \to 0$
in $\dcat{M}(A)$. 
According to Proposition \ref{prop:2165} this becomes a distinguished triangle 
\begin{equation} \label{eqn:3428}
L' \to P^0 \xar{\opn{Q}(\eta)} L \xar{\triangle}
\end{equation}
in $\dcat{D}(A)$.
Because both $L$ and $P^0$ are perfect, we see that the DG $A$-module $L'$ is 
also perfect. Thus $L'$ is a finitely generated $A$-module of finite 
projective dimension. However, by the usual syzygy argument 
(cf.\ \cite[Section 5.1.1]{Row}) -- or by examination of the long exact 
sequence gotten by applying 
$\opn{Hom}_{\dcat{D}(A)}(-, N)$, for an $A$-module $N$, 
to the distinguished triangle (\ref{eqn:3428}) -- we see that 
the projective dimension of $L'$ is $\leq d - 1$. Induction says 
that there is an isomorphism $L' \cong P'$ in $\dcat{D}(A)$ for some bounded 
complex of finitely generated projective $A$-modules $P'$.
Plugging this isomorphism into the distinguished triangle (\ref{eqn:3428})
we obtain a distinguished triangle 
\begin{equation} \label{eqn:3421}
P' \xar{\psi} P^0 \xar{\opn{Q}(\eta)} L \xar{\triangle} 
\end{equation}
in $\dcat{D}(A)$. Because $P'$ is K-projective, there is a homomorphism 
$\til{\psi} : P' \to P^0$ in $\dcat{C}_{\mrm{str}}(A)$ such that 
$\psi = \opn{Q}(\til{\psi})$. Let 
$P := \opn{Cone}(\til{\psi})$, the standard cone on $\til{\psi}$. 
Then $P$ is a bounded complex of finitely generated projective $A$-modules, and 
the distinguished triangle (\ref{eqn:3421}) gives rise to an isomorphism 
$P \cong L$ in $\dcat{D}(A)$.

\medskip \noindent
Step 2. Here $L$ is an arbitrary algebraically perfect DG $A$-module, nonzero 
in $\dcat{D}(A)$. By Lemma \ref{lem:3415} we have 
$\opn{con}(\opn{H}(L)) = [i_0, i_1]$ for integers $i_0 \leq i_1$, and 
$\opn{H}^{i_1}(L)$ is a finitely generated $A$-module. 
Let $n := \opn{amp}(\opn{H}(L)) = i_1 - i_0 \in \N$.
We proceed by induction on $n$. 

If $n = 0$ then $L$ is isomorphic to the translation of an $A$-module, and we 
are done by step 1. On the other hand, if $n \geq 1$, then choose a finite free 
$A$-module $P^{i_1}$, with a surjection 
$\bar{\eta} : P^{i_1} \to \opn{H}^{i_1}(L)$ in $\dcat{M}(A)$. 
This can be lifted to a homomorphism 
$\eta : P^{i_1}[-i_1] \to L$ in $\dcat{C}_{\mrm{str}}(A)$. 
Consider the DG $A$-module 
$L' := \opn{Cone}(\eta)$.
There is a distinguished triangle 
\begin{equation} \label{eqn:3423}
P^{i_1}[-i_1] \xar{\opn{Q}(\eta)} L \xar{} L' \xar{\triangle} 
\end{equation}
in $\dcat{D}(A)$. 
The DG module $L'$ is also algebraically perfect. There is a long exact 
cohomology sequence 
\[ \begin{aligned}
& \cdots \to 0 \to \opn{H}^{i_0}(L) \to \opn{H}^{i_0}(L') \to \cdots 
\to 0 \to \opn{H}^{i_1 - 2}(L) \to \opn{H}^{i_1 - 2}(L')
\\
& \quad 
\to 0 \to \opn{H}^{i_1 - 1}(L) \to \opn{H}^{i_1 - 1}(L') \to
P^{i_1} \xar{\bar{\eta}} \opn{H}^{i_1}(L) \to \opn{H}^{i_1}(L') \to 0 \to \cdots
\end{aligned} \]
Because $\bar{\eta}$ is surjective we know that $\opn{H}^{i_1}(L') = 0$.
Hence $\opn{con}(\opn{H}(L')) \sub [i_0, i_1 - 1]$
and $\opn{amp}(\opn{H}(L')) \leq n - 1$.

The induction assumption says that there is an isomorphism $L' \cong P'$
in $\dcat{D}(A)$ for some bounded complex of 
finitely generated projective $A$-modules $P'$.
Plugging this isomorphism into the triangle (\ref{eqn:3423}), and then turning 
it, we obtain a distinguished triangle 
\begin{equation} \label{eqn:3424} 
P'[-1] \xar{\psi} P^{i_1}[-i_1] \xar{\opn{Q}(\eta)} L \xar{\triangle} 
\end{equation}
in $\dcat{D}(A)$. Because $P'[-1]$ is K-projective, there is a homomorphism 
$\til{\psi} : P'[-1] \to P^{i_1}[-i_1]$ in $\dcat{C}_{\mrm{str}}(A)$ such that 
$\psi = \opn{Q}(\til{\psi})$. Let 
$P := \opn{Cone}(\til{\psi})$, the standard cone on $\til{\psi}$. 
Then $P$ is a bounded complex of finitely generated projective $A$-modules, and 
the distinguished triangle (\ref{eqn:3424}) gives rise to an isomorphism 
$P \cong L$ in $\dcat{D}(A)$.
\end{proof}

\begin{cor} \label{cor:4170}
Let $A$ be a ring and let $L$ be a complex of $A$-modules. The following two 
conditions are equivalent.
\begin{itemize}
\rmitem{i}  $L$ is a compact object of $\dcat{D}(A)$.

\rmitem{ii} $L$ is isomorphic in $\dcat{D}(A)$ to a bounded complex of 
finitely generated projective $A$-modules.
\end{itemize} 
\end{cor}

\begin{proof}
Combine Theorems \ref{thm:3415} and \ref{thm:3400}.
\end{proof}

Another special case is when the DG ring $A$ is nonpositive and cohomologically 
left pseudo-noetherian; see Definition \ref{dfn:3180}.

\begin{thm} \label{thm:3416}
Let $A$ be a cohomologically left pseudo-noetherian nonpositive DG ring,
and let $L$ be a DG $A$-module. The following two conditions are 
equivalent.
\begin{itemize}
\rmitem{i}  $L$ is an algebraically perfect DG $A$-module. 

\rmitem{ii} $L$ has finite projective dimension and belongs to 
$\dcat{D}^{-}_{\mrm{f}}(A)$.
\end{itemize}
\end{thm}

\begin{proof} \mbox{}

\nopagebreak
\smallskip \noindent
(i) $\twoto$ (ii): Since $A \in \dcat{D}^{-}_{\mrm{f}}(A)$, and 
$\dcat{D}^{-}_{\mrm{f}}(A)$ is a saturated full triangulated subcategory of 
$\dcat{D}(A)$, we see that $L \in \dcat{D}^{-}_{\mrm{f}}(A)$. By Proposition 
\ref{prop:3415}, $L$ has finite projective dimension. 

\medskip \noindent
(ii) $\twoto$ (i): We may assume that $\opn{H}(L) \neq 0$. 
Let $i_1 := \opn{sup}(\opn{H}(L)) \in \Z$, and let $n \in \N$ be the 
projective dimension of $L$.  By Theorem \ref{thm:3340}
there is a quasi-isomorphism $\rho : P \to L$ in $\dcat{C}_{\mrm{str}}(A)$ from 
a pseudo-finite semi-free DG $A$-module $P$ such that 
$\opn{sup}(P) = i_1$. Let $\{ F_j(P) \}_{j \geq -1}$ be a pseudo-finite 
semi-free filtration of $P$ as in Definition \ref{dfn:1925}(1),
with the same upper bound $i_1$.
Define $P' := F_{i_1 + n}(P)$ and $P'' := P / P'$. 
The DG $A$-module $P''$ is concentrated in the degree interval 
$[-\infty, i_1 -n - 1]$, and therefore 
\[ \opn{con} \bigl( \opn{H} ( \opn{RHom}_A(L, P'')) \bigr) \sub
[-\infty, - 1] . \]
This implies that 
\begin{equation} \label{eqn:3425}
 \opn{Hom}_{\dcat{D}(A)}(L, P'') \cong 
\opn{H}^0 \bigl( \opn{RHom}_A(L, P'') \bigr) = 0 .
\end{equation}

The DG $A$-module $P'$ is finite semi-free. 
Let $\ga' : P' \to P$ be the inclusion. 
From the distinguished triangle 
$P' \xar{\opn{Q}(\ga')} P \to P'' \xar{\triangle}$
in $\dcat{D}(A)$ we get an exact sequence of $\K$-modules
\[ \opn{Hom}_{\dcat{D}(A)}(L, P') \to 
\opn{Hom}_{\dcat{D}(A)}(L, P) \to 
\opn{Hom}_{\dcat{D}(A)}(L, P'') . \]
The isomorphism 
$\phi := \opn{Q}(\rho)^{-1} : L \to P$
belongs to the middle term above. But by (\ref{eqn:3425}) we see that $\phi$ 
comes from some morphism 
$\phi' : L \to P'$ in $\dcat{D}(A)$. 
Because $\opn{Q}(\rho \circ \ga') \circ \phi' = \opn{id}_L$, and using 
Theorem \ref{thm:4160},
we see that $L$ is a direct summand of $P'$ in $\dcat{D}(A)$. 
Theorem \ref{thm:3400} says that $L$ is algebraically perfect. 
\end{proof}

\begin{cor} \label{cor:3425}
Let $A$ be a left noetherian ring, and let $L$ be a DG $A$-module. The 
following 
two conditions are equivalent.
\begin{itemize}
\rmitem{i}  $L$ is an algebraically perfect DG $A$-module. 

\rmitem{ii} $L$ has finite projective dimension and belongs to 
$\dcat{D}^{\mrm{b}}_{\mrm{f}}(A)$.
\end{itemize}
\end{cor}

\begin{proof}
Combine Theorems \ref{thm:3415} and \ref{thm:3416}.
\end{proof}

\begin{rem} \label{rem:3426}
Here are some historical notes on perfect and compact objects.

Perfect complexes in algebraic geometry were introduced by A. Grothen\-dieck 
et.\ al.\ in \cite{SGA6}. For a scheme $X$, let $\dcat{D}_{\mrm{qc}}(X)$
be the derived category  of complexes of $\OO_X$-modules with 
quasi-coherent cohomology. The definition of a perfect complex given in 
\cite{SGA6} was this: $\LL$ is perfect if locally it is isomorphic to a bounded 
complex of finite rank free $\OO_X$-modules. This coincides with our notion of 
geometrically perfect complex, as in Remark \ref{rem:3425}, when $A$ is a 
commutative ring and $X = \opn{Spec}(A)$.

The idea of defining finiteness properties of an object $L$ in a category 
$\cat{D}$ via the commutation of the functor $\opn{Hom}_{\cat{D}}(L, -)$ with 
suitable coproducts goes back to the early days of category theory. Objects $L$ 
such that $\opn{Hom}_{\cat{D}}(L, -)$ commutes with coproducts were called 
{\em small}. Grothendieck \cite{Gro} used small objects in his construction of 
injective resolutions in what was eventually called a Grothendieck abelian 
category. P. Freyd \cite{Fre} used small objects in the proof of the 
Freyd-Mitchell Theorem. In algebraic topology, small objects were used by E. 
Brown \cite{Brown} in his celebrated Representability Theorem.

J. Rickard \cite{Ric1} proved a slightly weaker version of our Corollary 
\ref{cor:4170}. Shortly afterwards B. Keller \cite{Kel1} deduced from A. 
Neeman's work \cite{Ne0} a result that 
is essentially our Theorem \ref{thm:3400} -- he did not have condition (iii), 
but on the other hand he considered the more general derived category 
$\dcat{D}(\cat{A})$ of DG modules over a DG category $\cat{A}$, whereas we only 
consider $\dcat{D}(A)$ for a DG ring $A$, which is a single-object DG category. 

R.W. Thomason \cite{Thom} discovered that the perfect complexes in 
$\dcat{D}_{\mrm{qc}}(X)$, when $X$ is a quasi-projective scheme over a 
ring $\K$, are precisely the compact objects in $\dcat{D}_{\mrm{qc}}(X)$. 
Neeman \cite{Ne0} realized the connection between the work of Thomason and 
that of the topologists  A.K. Bousfield and D. Ravenel.
The terminology switch from ``small object'' to ``compact object'' seems to 
have taken place in \cite{Ne0}. 
Neeman's book \cite{Ne1} is  essentially devoted 
to the study of $\al$-compactly generated triangulated categories, for a 
regular cardinal number $\al$, and to generalizations of the Brown 
Representability Theorem. Note that what we called ``compact object'' in 
Definition \ref{dfn:3405} is, in the framework of \cite{Ne0}, an 
``$\aleph_0$-compact object''. 

The paper \cite{BonVdB} by A. Bondal and M. Van den Bergh was very influential 
in promoting the role of compact objects in algebraic geometry. 
In the last 20 years there has been a proliferation in the presence of compact 
(or perfect) objects in research in the areas of derived algebraic geometry,  
noncommutative algebraic geometry and mathematical physics -- much of this due 
to the influence of M. Kontsevich. For more details we recommend looking in the 
online reference \cite{nLab}. 
\end{rem}

\mysubsection{Derived Morita Theory} \label{subsec:der-morita}

Recall that Convention \ref{conv:3485} is in effect.
In particular all DG rings are $\K$-central, and 
$\ot = \ot_{\K}$.

\begin{dfn} \label{dfn:3525}
Let $A$ and $B$ be DG rings. The objects of the category 
$\dcat{C}(A \ot B^{\mrm{op}})$
are called 
\index{Differential graded bimodule}
{\em DG $A$-$B$-bimodules}. 
\end{dfn}

There is a commutative diagram (\ref{eqn:3535}) 
in $\catt{DGRng} \centover \K$.

\begin{equation} \label{eqn:3535}
\UseTips  \xymatrix @C=3ex @R=3ex {
&
\ A \ot B^{\mrm{op}}
\\
A
\ar[ur]
&
&
B^{\mrm{op}}
\ar[ul]
\\
&
\K 
\ar[ul]
\ar[ur]
\ar[uu]
}
\end{equation}

By restriction there is a commutative diagram of DG 
functors, and a commutative diagram of triangulated functors, see 
(\ref{eqn:3536}).

\begin{equation} \label{eqn:3536}
\UseTips  \xymatrix @C=3ex @R=3ex {
&
\ \dcat{C}(A \ot B^{\mrm{op}})
\ar[dl]_{\opn{Rest}_{A}}
\ar[dr]^{\opn{Rest}_{B^{\mrm{op}}}}
\ar[dd]^{\opn{Rest}_{\K}}
\\
\dcat{C}(A)
\ar[dr]_{\opn{Rest}_{\K}}
&
&
\dcat{C}(B^{\mrm{op}})
\ar[dl]^{\opn{Rest}_{\K}}
\\
&
\dcat{C}(\K) 
}
\quad 
\UseTips  \xymatrix @C=3ex @R=3ex {
&
\ \dcat{D}(A \ot B^{\mrm{op}})
\ar[dl]_{\opn{Rest}_{A}}
\ar[dr]^{\opn{Rest}_{B^{\mrm{op}}}}
\ar[dd]^{\opn{Rest}_{\K}}
\\
\dcat{D}(A)
\ar[dr]_{\opn{Rest}_{\K}}
&
&
\dcat{D}(B^{\mrm{op}})
\ar[dl]^{\opn{Rest}_{\K}}
\\
&
\dcat{D}(\K) 
}
\end{equation}

There are several manipulations of DG bimodules that we are going to use. 
First we note that $(A^{\mrm{op}})^{\mrm{op}} = A$.
Next, given DG $\K$-modules $M$ and $N$, let 
\begin{equation} \label{eqn:3550} 
\opn{br}_{M, N} : M \ot N \iso N \ot M
\end{equation}
be the isomorphism in $\dcat{C}(\K)$ defined by 
\begin{equation} \label{eqn:3553} 
\opn{br}_{M, N}(m \ot n) := (-1)^{i \cd j} \cd n \ot m
\end{equation}
for homogeneous elements $m \in M^i$ and $n \in N^j$. 
In Subsection \ref{subsec:gr-alg} the isomorphism $\opn{br}_{M, N}$ was called 
the braiding of the symmetric monoidal category $\dcat{C}(\K)$.
If $M \in \dcat{C}(A^{\mrm{op}})$, then we can view $M$ either as a 
left DG $A^{\mrm{op}}$-module, or as a right DG $A$-module.
The opposite holds for $N \in \dcat{C}(A)$. It is not hard to see that there 
is an isomorphism 
\begin{equation} \label{eqn:3551} 
M \ot_A N \iso N \ot_{A^{\mrm{op}}} M
\end{equation}
in $\dcat{C}(\K)$ with the same formula (\ref{eqn:3553}).
This works also for DG rings: there are isomorphisms 
\begin{equation} \label{eqn:3552} 
A \ot B^{\mrm{op}} \iso B^{\mrm{op}} \ot A
\end{equation}
and
\begin{equation} \label{eqn:3554} 
A \ot B^{\mrm{op}} \iso (B \ot A^{\mrm{op}})^{\mrm{op}}
\end{equation}
in $\catt{DGRng} \centover \K$, with a formulas like (\ref{eqn:3553}). 

There is a DG bifunctor
\begin{equation} \label{eqn:3537}
(- \ot^{}_{B} -) :  \dcat{C}(A \ot B^{\mrm{op}}) \times \dcat{C}(B)
\to \dcat{C}(A)  . 
\end{equation}
Here is a variation of Theorem \ref{thm:2107}.

\begin{prop} \label{prop:3535}
The bifunctor \tup{(\ref{eqn:3537})} has a triangulated left derived bifunctor
\[ (- \ot^{\mrm{L}}_{B} -) :
\dcat{D}(A \ot B^{\mrm{op}}) \times \dcat{D}(B)
\to \dcat{D}(A) . \]
If $P \in \dcat{D}(B)$ is a K-flat DG module, then for every 
$M \in \dcat{D}(A \ot B^{\mrm{op}})$
the morphism 
\[ \eta^{\mrm{L}}_{M, P} : M \ot^{\mrm{L}}_{B} P \to M \ot^{}_{B} P \]
in $\dcat{D}(A)$ is an isomorphism. 
\end{prop}

\begin{proof}
This is because $\dcat{C}(B)$ has enough K-flat objects. See
Theorem \ref{thm:2106} and Lemma \ref{lem:2109}.
\end{proof}

\begin{rem} \label{rem:3431}
There is a delicate issue here. Even though the category 
$\dcat{C}(A \ot B^{\mrm{op}})$ has enough K-flat objects, 
they can not be used to calculate $(- \ot^{\mrm{L}}_{B} -)$. 
The reason is this: given $N \in \dcat{C}(B)$,
the DG $(A \ot B^{\mrm{op}})$-modules that are acyclic 
for the DG functor $(-) \ot_B N$ are those that are 
{\em K-flat over $B^{\mrm{op}}$}. In general, a K-flat DG 
$(A \ot B^{\mrm{op}})$-module $M$ is not K-flat as a $B^{\mrm{op}}$-module. 

A sufficient condition for a K-flat DG $(A \ot B^{\mrm{op}})$-module 
$M$ to be K-flat over $B^{\mrm{op}}$ is that $A$ is K-flat as a DG 
$\K$-module. In Subsection \ref{subsec:flat-dg-rng} we will make heavy use of 
this fact. 
\end{rem}

There is another DG bifunctor we shall want to use:
\begin{equation} \label{eqn:3538}
\opn{Hom}_A(-, -) :  \dcat{C}(A \ot B^{\mrm{op}})^{\mrm{op}} \times 
\dcat{C}(A) \to \dcat{C}(B)  . 
\end{equation}
Here is a variation of Theorem \ref{thm:2005}. 

\begin{prop} \label{prop:3536}
The bifunctor \tup{(\ref{eqn:3538})} has a triangulated right derived 
bifunctor
\[ \opn{RHom}_A(-, -) :
\dcat{D}(A \ot B^{\mrm{op}})^{\mrm{op}} \times \dcat{D}(A)
\to \dcat{D}(B) . \]
If $I \in \dcat{D}(A)$ is a K-injective DG module, then for every 
$M \in \dcat{D}(A \ot B^{\mrm{op}})$
the morphism 
\[ \eta^{\mrm{R}}_{M, I} : \opn{Hom}_A(M, I) \to \opn{RHom}_A(M, I) \]
in $\dcat{D}(B)$ is an isomorphism. 
\end{prop}

\begin{proof}
This is because $\dcat{C}(A)$ has enough K-injective objects. See
Theorem \ref{thm:4225}, Theorem \ref{thm:2000} and 
Lemma \ref{lem:2000}.
\end{proof}

\begin{rem} \label{rem:3432}
Like in Remark \ref{rem:3431}, even though the category 
$\dcat{C}(A \ot  B^{\mrm{op}})$ has enough K-projective 
DG modules, they can not be used to calculate \lb $\opn{RHom}_A(-, N)$. 
This is because a K-projective DG $(A \ot  B^{\mrm{op}})$-module $P$ is 
not, in general, K-projective over $A$. 

A sufficient condition for a K-projective DG $(A \ot B^{\mrm{op}})$-module 
$P$ to be K-projective over $A$ is that $B$ is K-projective over $\K$. 
This is a bit stronger than the condition that $B$ is K-flat over $\K$.
\end{rem}

\begin{prop} \label{prop:3537}
Let $L \in \dcat{D}(A \ot B^{\mrm{op}})$, and consider the $\K$-linear 
triangulated functors
\[ F := \opn{RHom}_A(L, -) : \dcat{D}(A) \to \dcat{D}(B) \]
and 
\[ G := L \ot^{\mrm{L}}_{B} (-) : \dcat{D}(B) \to \dcat{D}(A) . \]
Then\tup{:}
\begin{enumerate}
\item The functor $F$ is a right adjoint of $G$.

\item The functor $F$ is an equivalence if and only if the functor $G$ is an 
equivalence. 
\end{enumerate}
\end{prop}

\begin{proof} \mbox{}

\smallskip \noindent
(1) Given $M \in \dcat{D}(B)$ and $N \in \dcat{D}(A)$ we choose 
a K-projective resolution $P \to M$ in $\dcat{C}_{\mrm{str}}(B)$ and a 
K-injective resolution $N \to I$ in $\dcat{C}_{\mrm{str}}(A)$.
We get isomorphisms 
$F(N) \cong F(I) \cong \opn{Hom}_A(L, I)$ 
and
$G(M) \cong G(P) \cong L \ot_B P$ 
in $\dcat{D}(B)$ and $\dcat{D}(A)$ respectively.
These give rise to isomorphisms of $\K$-modules
\[ \begin{aligned}
& \opn{Hom}_{\dcat{D}(B)} \bigl( M, F(N) \bigr) \cong^{\diamond}
\opn{H}^0 \bigl( \opn{RHom}_B \bigl( M, F(N) \bigr) \bigr)
\\
& \qquad \cong
\opn{H}^0 \bigl( \opn{Hom}_B \bigl( P, \opn{Hom}_A(L, I) \bigr) \bigr)
\cong^{\dag}
\opn{H}^0 \bigl( \opn{Hom}_A(L \ot_B P, I) \bigr)
\\
& \qquad \cong
\opn{H}^0 \bigl( \opn{RHom}_A \bigl( G(M), N \bigr) \bigr)
\cong^{\diamond} \opn{Hom}_{\dcat{D}(A)} \bigl( G(M), N \bigr) 
\end{aligned} . \]
The isomorphism $\cong^{\dag}$ is due to the noncommutative Hom-tensor 
adjunction (see \cite[Theorem 2.11]{Rot}), and the 
isomorphisms $\cong^{\diamond}$ are by Corollary \ref{cor:2120}.
The composed isomorphism 
\[ \opn{Hom}_{\dcat{D}(B)} \bigl( M, F(N) \bigr) \cong
\opn{Hom}_{\dcat{D}(A)} \bigl( G(M), N \bigr) \]
is functorial in $M$ and $N$. 

\medskip \noindent 
(2) Clear from (1). 
\end{proof}

\begin{dfn} \label{dfn:3430}
Let $A$ and $B$ be DG rings. A DG module
$L \in \dcat{D}(A \ot B^{\mrm{op}})$
is called a {\em pretilting DG $A$-$B$-bimodule}%
\index{Tilting! pre{\indash} DG bimodule}
if it satisfies the equivalent 
conditions in Proposition \ref{prop:3537}(2). 
\end{dfn}

\begin{rem} \label{rem:3490}
The literature has several different meanings for the word ``tilting''. 
See Remark \ref{rem:3436} for a historical survey. 
The commonly perceived meaning of a
tilting object, say as in the book \cite{AnHaKr}, is very close to what we call 
a pretilting DG bimodule (when $A$ and $B$ are rings). 

Notice that there is a lack of symmetry in the 
definition of a pretilting DG $A$-$B$-bimodule $L$.
This lack of symmetry will disappear when we talk about {\em tilting DG 
bimodules} in Subsection \ref{subsec:flat-dg-rng}. 
\end{rem}

\begin{dfn} \label{dfn:3446}
Let $A$ and $B$ be DG rings, and let 
$M \in \dcat{C}(A \ot B^{\mrm{op}})$. We say that $M$ is {\em K-flat} 
(resp.\ {\em K-injective}, resp.\ {\em K-projective}, resp.\ {\em semi-free},
resp.\ {\em algebraically perfect}) {\em over $A$}, or {\em on the $A$ side},
if $\opn{Rest}_A(M) \in \dcat{C}(A)$ is K-flat (resp.\  K-injective, resp.\ 
K-projective, resp.\ semi-free, resp.\ algebraically perfect).
Likewise we define the properties {\em on the $B^{\mrm{op}}$ side}. 
\end{dfn}

\begin{prop} \label{prop:3431}
Suppose $L \in \dcat{D}(A \ot B^{\mrm{op}})$ is a pretilting DG 
$A$-$B$-bimodule. Then $L$ is an algebraically perfect DG $A$-module.
\end{prop}

\begin{proof}
Under the equivalence $G : \dcat{D}(B) \to \dcat{D}(A)$ from Proposition 
\ref{prop:3537} we have $G(B) \cong L$. Of course $B \in \dcat{D}(B)$ is an 
algebraically perfect DG $B$-module. Now use Corollary \ref{cor:3480}. 
\end{proof}

Observe that an object $N \in \dcat{D}(A)$ is {\em nonzero} if it is 
not a zero object of the category $\dcat{D}(A)$, i.e.\ if $N \not\cong 0$ in
$\dcat{D}(A)$. This is equivalent to the condition that 
$\opn{H}(N) \neq 0$. The clarification is important for the next definition. 

\begin{dfn} \label{dfn:3431}
Let $\cat{E} \sub \dcat{D}(A)$ be a full triangulated subcategory that's closed 
under infinite direct sums. 
\begin{enumerate}
\item An object $L \in \dcat{D}(A)$ is called a 
{\em compact object relative to $\cat{E}$}%
\index{Differential graded module! compact}
if the functor 
\[ \opn{Hom}_{\dcat{D}(A)}(L, -) : \cat{E} \to \dcat{M}(\K) \]
commutes with infinite direct sums. 

\item An object $L \in \dcat{D}(A)$ is called a 
{\em generator relative to $\cat{E}$}%
\index{Differential graded module! generator}
if for every nonzero object $N \in \cat{E}$ there exists some 
$i \in \Z$ such that 
$\opn{Hom}_{\dcat{D}(A)}(L, N[i]) \neq 0$. 

\item An object $L \in \cat{E}$ is called a 
{\em compact generator of $\cat{E}$} 
if it is both a compact object relative to 
$\cat{E}$ and a generator 
relative to $\cat{E}$.
\end{enumerate}
\end{dfn}

\begin{exa} \label{exa:3430}
For every positive integer $r$ the free DG module 
$P := A^r$ is a compact generator of $\cat{E} := \dcat{D}(A)$.
\end{exa}

\begin{exa} \label{exa:3431}
Assume $A$ is a commutative noetherian ring and 
$\bsym{a} = \lb (a_1, \ldots, a_n)$ is a finite sequence of elements in it. 
(Or, more generally, $A$ is commutative and $\bsym{a}$ is a {\em weakly 
proregular sequence}, as in \cite{PSY}.)
Let $\a \sub A$ be the ideal generated by $\bsym{a}$. 
A complex $M \in \dcat{D}(A)$ is called {\em cohomologically $\a$-torsion}
if all its cohomology modules $\opn{H}^i(M)$ are $\a$-torsion. 
Consider the full subcategory 
$\cat{E} := \dcat{D}_{\a\tup{-tor}}(A)$ of $\dcat{D}(A)$ on the 
cohomologically $\a$-torsion complexes. 
This is full triangulated subcategory that's closed 
under infinite direct sums. 
The Koszul complex $\opn{K}(A; \bsym{a})$ is a compact generator of 
$\dcat{D}_{\a\tup{-tor}}(A)$. 
See \cite[Proposition 5.1]{PSY2}.
\end{exa}

\begin{prop} \label{prop:3547}
Let $A$ and $B$ be DG rings, $F : \dcat{D}(A) \to \dcat{D}(B)$ an 
equivalence of triangulated categories, 
$\cat{E} \sub \dcat{D}(A)$ a full triangulated subcategory that's closed 
under infinite direct sums, and $L \in \dcat{D}(A)$. 
\begin{enumerate}
\item $L$ is a compact object relative to $\cat{E}$ if and only if 
$F(L)$ is a compact object relative to $F(\cat{E})$. 

\item $L$ is a generator relative to $\cat{E}$ if and only if 
$F(L)$ is a generator relative to $F(\cat{E})$. 
\end{enumerate}
\end{prop}

\begin{proof}
For item (1) use the proof of Proposition \ref{prop:3545}. Item (2) is trivial. 
\end{proof}

\begin{dfn} \label{dfn:3465}
Let $A$ and $B$ be DG rings, let $\cat{E} \sub \dcat{D}(A)$ be a 
full triangulated subcategory that is closed under infinite direct sums, and 
let $L \in \dcat{D}(A \ot B^{\mrm{op}})$. 
\begin{enumerate}
\item We say that $L$ is a {\em compact object relative to $\cat{E}$ on the 
$A$ side} if $\opn{Rest}_A(L) \in \dcat{D}(A)$ is a compact object relative to 
$\cat{E}$, in the sense of Definition \lb \ref{dfn:3431}(1).

\item We say that $L$ is a {\em generator relative to $\cat{E}$ on the 
$A$ side} if $\opn{Rest}_A(L) \in \dcat{D}(A)$ is a generator relative to 
$\cat{E}$, in the sense of Definition \ref{dfn:3431}(2).

\item We say that $L$ is a {\em compact generator of $\cat{E}$ on the $A$ side} 
if $\opn{Rest}_A(L) \in \dcat{D}(A)$ is 
a compact generator of $\cat{E}$, in the sense of Definition \ref{dfn:3431}(3).
\end{enumerate}
\end{dfn}

\begin{lem} \label{lem:3436}
Let $\cat{E} \sub \dcat{D}(A)$ be a full triangulated subcategory that is 
closed under infinite direct sums, and let 
$L \in \dcat{D}(A \ot B^{\mrm{op}})$. The following 
two conditions are equivalent.
\begin{itemize}
\rmitem{i} The functor
$\opn{RHom}_A(L, -)|_{\cat{E}} : \cat{E} \to \dcat{D}(B)$
commutes with infinite direct sums.

\rmitem{ii} $L$ is a compact object relative to $\cat{E}$ on the $A$ side. 
\end{itemize}
\end{lem}

\begin{proof}
Let's write 
$F := \opn{RHom}_A(L, -)$. We know that for every $M \in \dcat{D}(A)$ and 
$j \in \Z$ there are isomorphisms
\[ \opn{Hom}_{\dcat{D}(A)}(L, M[j]) \cong 
\mrm{H}^0(F(M[j])) \cong \mrm{H}^j (F(M))  \]
in $\dcat{M}(\K)$, and they are functorial in $M$.
We see that $L$ is a compact object relative to $\cat{E}$ on the $A$ side iff
the functor
$\mrm{H} \circ F|_{\cat{E}} : \cat{E} \to \dcat{G}_{\mrm{str}}(\K)$
commutes with infinite direct sums. But the functor 
$ \mrm{H} : \dcat{D}(B) \to \dcat{G}_{\mrm{str}}(\K)$ 
commutes with infinite direct sums and is conservative (see Corollary 
\ref{cor:2145}). 
So $\mrm{H} \circ F|_{\cat{E}}$ commutes with infinite direct 
sums iff $F|_{\cat{E}}$ does. 
\end{proof}

\begin{lem} \label{lem:3435}
Let $A$ and $B$ be DG rings, let 
$F, G : \dcat{D}(A) \to \dcat{D}(B)$ be
triangulated functors that commute with infinite direct sums,
and let $\eta : F \to G$ be a morphism of triangulated functors. Assume that
$\eta_A: F(A) \to G(A)$ is an isomorphism. Then $\eta$ is an isomorphism. 
\end{lem}

\begin{proof}
Let us denote by $\cat{E}$ the full subcategory of $\dcat{D}(A)$ on the objects 
$M$ such that $\eta_M: F(M) \to G(M)$ is an isomorphism. We have to prove that 
$\cat{E} = \dcat{D}(A)$. 

By Proposition \ref{prop:4603}, 
the category $\cat{E}$ is a full triangulated subcategory of $\dcat{D}(A)$.
Since both functors $F, G$ commute with infinite direct sums,
we know that $\cat{E}$ is closed under infinite direct sums.

We are given that the free DG module $A$ belongs to $\cat{E}$.
Hence every free DG $A$-module $P$ is in $\cat{E}$.
Because $\cat{E}$ is triangulated, it follows that every semi-free DG 
$A$-module $P$ of finite extension length (Definition \ref{dfn:3411}) belongs 
to 
it. 

Next consider an arbitrary semi-free DG $A$-module $P$.
Let $\{ F_j(P) \}_{j \geq -1}$ be a semi-free filtration of $P$
(see Definition \ref{dfn:1575}). 
Then $P$ is a homotopy colimit of the direct system 
$\{ F_j(P) \}_{j \geq -1}$, and we have a distinguished triangle 
\begin{equation} \label{eqn:3440}
\bigoplus\nolimits_{j \in \N} F_j(P) \xar{\ \phi \ } 
\bigoplus\nolimits_{j \in \N} F_j(P) 
\xar{ \ \ } P \xar{\triangle}
\end{equation}
in $\dcat{D}(A)$; see Definition \ref{dfn:3410}. 
Each $F_j(P)$ is a semi-free DG $A$-module of finite extension length. It 
follows that $P \in \cat{E}$. 

Finally, every DG $A$-module $M$ admits a quasi-isomorphism 
$P \to M$ with $P$ semi-free. Therefore $M \in \cat{E}$.
\end{proof}

\begin{lem} \label{lem:3470}
Let $A$ and $B$ be DG rings, let $\cat{E}$ be a full triangulated 
subcategory of $\dcat{D}(A)$ which is closed under infinite direct sums
and isomorphisms, and let $G : \dcat{D}(B) \to \dcat{D}(A)$ be a
triangulated functor that commutes with infinite direct sums.
Assume that $G(B) \in \cat{E}$. Then the essential image of $G$ is contained 
in $\cat{E}$. 
\end{lem}

\begin{proof}
This is like the proof of Lemma \ref{lem:3435}.
\end{proof}

Let $A$ be a DG ring. Given a DG $A$-module $M$, the DG ring 
$\opn{End}_A(M)$ acts on $M$ from the left; and thus 
$B := \opn{End}_A(M)^{\mrm{op}}$ 
acts on $M$ from the right. Because the right action of $B$ on $M$ commutes 
with the left action of $A$ on it, we see that 
$M \in \dcat{C}(A \ot B^{\mrm{op}})$.

Recall that the category $\dcat{C}(A)$ has enough K-projective and K-injective 
objects. 

\begin{thm}[Derived Morita] \label{thm:3430}
\index{Derived Morita theory}
Let $A$ be a DG ring, let 
$\cat{E} \sub \dcat{D}(A)$ be a full triangulated subcategory which is closed 
under infinite direct sums, and let $L$ be a compact 
generator of $\cat{E}$. Choose an isomorphism 
$L \cong P$ in $\dcat{D}(A)$, where the DG $A$-module $P$ is either 
K-projective or K-injective, and define the DG ring
$B := \opn{End}_A(P)^{\mrm{op}}$.
Consider the triangulated functors 
\[ F := \opn{RHom}_A(P, -) : \dcat{D}(A) \to \dcat{D}(B) \] 
and 
\[ G := P \ot_B^{\mrm{L}} (-) : \dcat{D}(B) \to \dcat{D}(A) . \]
Then the functor 
\[ F|_{\cat{E}} : \cat{E} \to \dcat{D}(B) \]
is an equivalence, with quasi-inverse $G$. 
\end{thm}

This theorem is a variant of \cite[Theorem 8.2]{Kel1} by B. Keller.
In fact, Keller worked in the more general setting of a DG category $\cat{A}$, 
whereas we only treat a DG ring $A$. For the history of these ideas see Remark 
\ref{rem:3436}.
A comparison to classical Morita theory can be found in Example \ref{exa:3435}.

Note that if we choose $P$ to be a K-projective DG $A$-module, then 
$F \cong \lb \opn{Hom}_A(P, -)$. 
However, even then, $P$ will rarely be K-flat as a DG $B^{\mrm{op}}$-module, so 
$G$ is only the left derived functor of $P \ot_B (-)$. 

\begin{proof} 
The proof is in four steps. 

\smallskip \noindent 
Step 1. We can assume that the category $\cat{E}$ is closed under isomorphisms 
in $\dcat{D}(A)$. The functor $G$ commutes with infinite direct sums (see 
Proposition \ref{prop:4300}), 
and $G(B) \cong L \in \cat{E}$. According to Lemma \ref{lem:3470} the image of 
$G$ is contained in $\cat{E}$. 

\medskip \noindent 
Step 2. We already know that the functor $F$ is right adjoint to $G$, by 
Proposition \ref{prop:3537}. 
The corresponding morphisms of triangulated functors are denoted by
\begin{equation} \label{eqn:3471}
\th : \opn{Id}_{\dcat{D}(B)} \to F \circ G \quad \tup{and} \quad
\ze : G \circ F \to \opn{Id}_{\dcat{D}(A)} . 
\end{equation}
By step 1 we see that  
$F|_{\cat{E}} : \cat{E} \to \dcat{D}(B)$ 
is the right adjoint of 
$G : \dcat{D}(B) \to \cat{E}$.
Hence (\ref{eqn:3471}) restricts to morphisms of triangulated functors
\begin{equation} \label{eqn:3472}
\th : \opn{Id}_{\dcat{D}(B)} \to F|_{\cat{E}} \circ G \quad \tup{and} \quad
\ze : G \circ F|_{\cat{E}} \to \opn{Id}_{\cat{E}}  . 
\end{equation}
We will prove that $\th$ and $\ze$ in (\ref{eqn:3472}) are isomorphisms. 
By Lemma \ref{lem:3436} the functor $F|_{\cat{E}}$ commutes with infinite 
direct sums. Therefore both $F|_{\cat{E}} \circ G$ and 
$G \circ F|_{\cat{E}}$ commute with infinite direct sums.

\medskip \noindent 
Step 3. Now we will prove that $\th$ is an isomorphism of functors; i.e.\
for every $N \in \dcat{D}(B)$ the morphism
$\th_N :  N \to (F|_{\cat{E}} \circ G) (N)$ 
is an isomorphism. The functors $\opn{Id}_{\dcat{D}(B)}$ and 
$F|_{\cat{E}} \circ G$ both commute with infinite direct sums.
Therefore, by Lemma \ref{lem:3435}, 
it suffices to check that $\th_B$ is an isomorphism in $\dcat{D}(B)$. 
But $\th_B$ is represented (both when $P$ is K-projective and when it is  
K-injective) by the canonical homomorphism 
$\til{\th}_B : B \to \opn{Hom}_{A}(P, P \otimes_B B)$,
in $\dcat{C}_{\mrm{str}}(B)$, 
which is clearly bijective. 

\medskip \noindent 
Step 4. Finally we will prove that $\ze$ is an isomorphism of functors. 
Take any $M \in \cat{E}$, and 
consider the distinguished triangle
\begin{equation} \label{eqn:3473}
(G \circ F|_{\cat{E}})(M) \xar{\zeta_{M}} M \to M' \xar{\triangle}
\end{equation}
in $\cat{E}$, in which $M' \in \cat{E}$ is the cone of $\zeta_{M}$
(see Definition \ref{dfn:3600}).

Applying $F$ and using the functorial isomorphism $\th$ we get a 
distinguished triangle
\begin{equation} \label{eqn:3474}
F|_{\cat{E}}(M) \xar{\opn{id}_{F|_{\cat{E}}(M)}} F|_{\cat{E}}(M)
\to F|_{\cat{E}}(M') \xar{\triangle}
\end{equation}
in $\dcat{D}(B)$. This implies that 
$F|_{\cat{E}}(M') \cong 0$ in $\dcat{D}(B)$, and thus, after applying the 
forgetful functor, we get $F|_{\cat{E}}(M') \cong 0$ in $\dcat{D}(\K)$. But 
$F|_{\cat{E}}(M') \cong \opn{RHom}_{A}(L, M')$,
so 
$\opn{Hom}_{\dcat{D}(A)}(L, M'[i]) = 0$ 
for all $i$. Because $L$ is generator of $\cat{E}$, it follows that 
$M' \cong 0$ in $\cat{E}$. Going back to the distinguished 
triangle (\ref{eqn:3473}) we conclude that $\ze_M$ is an isomorphism.
\end{proof}

\begin{rem} \label{rem:3435}
Suppose that in the theorem we were to choose some other isomorphism 
$L \cong P'$ in $\dcat{D}(A)$ to a K-projective or K-injective DG $A$-module 
$P$. Then the DG ring $B' := \opn{End}_A(P')^{\mrm{op}}$
would be related to $B$ as follows. Without loss of generality we can assume 
that either $P$ is K-projective or $P'$ is K-injective. Then there is a 
quasi-isomorphism $\phi : P \to P'$ in $\dcat{C}_{\mrm{str}}(A)$, 
unique up to homotopy, that respects the given isomorphisms 
$L \cong P$ and $L \cong P'$ in $\dcat{D}(A)$. 
Define the DG bimodule 
$N := \opn{Hom}_A(P, P') \in \dcat{C}(B \ot B'^{\, \mrm{op}})$, 
and the matrix DG ring 
$C := \sbmat{B & \, N[-1] \\[0.2em] 0 & B'}$,
which has a differential depending on $\phi$. 
Then the obvious DG ring homomorphisms $C \to B$ and $C \to B'$ are 
quasi-isomorphisms. See \cite[Proposition 3.3]{PSY2}. 
\end{rem}

\begin{exa}  \label{exa:3432}
In case $\cat{E} = \dcat{D}(A)$, the theorem shows that the DG bimodule 
$P \in \dcat{D}(A \ot B^{\mrm{op}})$ 
is a pretilting DG $A$-$B$-bimodule. 
\end{exa}

\begin{exa}  \label{exa:3435}
Assume $A$ is a ring, and $P \in \dcat{M}(A) = \cat{Mod} A$ is a 
{\em progenerator}, i.e.\ $P$ is a finitely generated projective $A$-module, 
such that every nonzero $A$-module $N$ admits a nonzero homomorphism 
$P \to N$. Let $B := \opn{End}_A(P)^{\mrm{op}}$. 
Because $P$ is a compact generator of $\dcat{D}(A)$, the theorem 
says that 
\[ F := \opn{Hom}_A(P, -) : \dcat{D}(A) \to \dcat{D}(B) \]
is an equivalence of triangulated categories. 

Classical Morita Theory 
(see e.g.\ \cite[Section 4.1]{Row}) says that $F$ restricts to an equivalence 
of abelian categories 
\[ F = \opn{Hom}_A(P, -) : \dcat{M}(A) \to \dcat{M}(B) . \]
Furthermore, classical Morita theory tells us that there is a bijection between 
the set of isomorphism classes of such $A$-$B$-bimodules $P$, and the set of 
isomorphism classes of linear equivalences 
$F : \dcat{M}(A) \to \dcat{M}(B)$. 
This last assertion is an open problem in derived Morita theory; see Remark 
\ref{rem:3565}. 
\end{exa}

\begin{rem} \label{rem:3470}
Here is a geometric variant of Theorem \ref{thm:3430}. Let $(X, \OO_X)$ be a 
scheme. We denote by 
$\dcat{D}_{\mrm{qc}}(X) = \dcat{D}_{\mrm{qc}}(\cat{Mod} \OO_X)$ 
the derived category of 
(unbounded) complexes of $\OO_X$-modules with quasi-coherent cohomology 
sheaves. This is a triangulated category with infinite direct sums, and with 
enough K-injective resolutions (by \cite[Theorem 4.5]{Spa}). 

Suppose $\LL \in \dcat{D}_{\mrm{qc}}(X)$ is a compact generator. Choose a 
K-injective resolution $\LL \to \II$ in 
$\dcat{C}_{\mrm{str}}(X)$, and define the DG ring 
$B := \opn{End}_X(\II)^{\mrm{op}}$.
Then the functor 
\begin{equation} \label{eqn:5088}
\opn{RHom}_X(\II, -) : \dcat{D}_{\mrm{qc}}(X) \to \dcat{D}(B)
\end{equation}
is an equivalence, with quasi-inverse 
\[ \II \ot^{\mrm{L}}_{B} (-) : \dcat{D}(B) \to \dcat{D}_{\mrm{qc}}(X) . \]
This result is stated as \cite[Corollary 3.1.8]{BonVdB}, without such a precise 
formulation, and the proof is attributed to B. Keller. Presumably the proof of 
Theorem \ref{thm:3430} above, with slight changes, would work also in this 
geometric context. 

Having a compact generator of $\dcat{D}_{\mrm{qc}}(X)$ is quite a general 
feature -- by \cite[Theorem 3.1.1]{BonVdB} this is true if $X$ is quasi-compact 
and quasi-separated. 

The first such result is perhaps due to A. Beilinson \cite{Bei}, who showed 
that for $X = \mbf{P}^n_{\K}$, the $n$-dimensional projective space over a 
field $\K$, the sheaf 
\[ \LL := \bigoplus\nolimits_{i = 0}^n \, \OO_X(i) \]
is a compact generator of $\dcat{D}_{\mrm{qc}}(X)$. The DG ring $B$ appearing 
in the equivalence (\ref{eqn:5088}) is actually a finite $\K$-ring (the path 
ring of a Kronecker quiver modulo commutation relations).  
\end{rem}

\begin{rem} \label{rem:3436}
Tilting theory and derived Morita theory have their origins in the 
{\em representation theory of finite dimensional algebras} (in our terminology 
these are ``finite central $\K$-rings'', where $\K$ is a base field). Among the 
first examples of tilting are the {\em reflection functors} in the paper
\cite{BeGePo} by I.N. Bernstein, I.M. Gelfand and V.A. Ponomarev. Later these 
functors were understood to be of the form $\opn{Hom}_{A}(T, -)$ for a {\em 
tilting $A$-module} $T$. 

These ideas were generalized by M. Auslander, M. Platzeck and I. Reiten 
\cite{AuPlRe}, and later by S. Brenner and M.C.R. Butler \cite{BrnBu}, who 
coined the term {\em tilting functor}. 
These concepts were further clarified by  D. Happel, C.M. Ringel and 
K. Bongartz (see \cite{Hap0}, \cite{Hap} and \cite{HapRi}). In \cite{Hap0} 
Happel showed that for a tilting $A$-module $T$, 
with opposite endomorphism ring $B := \opn{End}_{A}(T)^{\mrm{op}}$, 
the functor $\opn{RHom}_{A}(T,-)$ is an equivalence of triangulated categories 
$\dcat{D}_{\mrm{f}}^{\mrm{b}}(A) \to \dcat{D}_{\mrm{f}}^{\mrm{b}}(B)$.  
This result was slightly generalized by E. Cline, B. Parshall and L. Scott
\cite{ClPaSc}. 

J. Rickard \cite{Ric1}, \cite{Ric2} was the first to talk about {\em tilting 
complexes} (as opposed to tilting modules). He introduced two-sided tilting 
complexes (that we will discuss in Subsection \ref{subsec:tilting}), and proved 
the celebrated Theorem \ref{thm:3565} (under some boundedness conditions). The 
generalization of derived Morita theory to unbounded derived categories, 
and from rings to DG categories, was done by B. Keller \cite{Kel1}.
In this seminal paper Keller also gave new construction of two-sided 
tilting complexes,  and characterized algebraic triangulated categories with a 
compact generator as those that are equivalent to derived categories of DG 
rings. 

For a thorough survey of tilting theory see the book \cite{AnHaKr}.
\end{rem}

\mysubsection{DG Bimodules over K-Flat DG Rings} \label{subsec:flat-dg-rng}

Recall that Convention \ref{conv:3485} is in effect.
In particular all DG rings are $\K$-central, and 
$\ot = \ot_{\K}$.

\begin{dfn} \label{dfn:3560}
Let $A$ be a DG $\K$-ring. We call $A$ a {\em K-flat DG $\K$-ring} 
\index{Differential graded ring! K-flat}
if $A$ is K-flat as a DG $\K$-module (see Definition \ref{dfn:1525}). 
The category of K-flat central DG $\K$-rings is denoted by 
$\catt{DGRng} \fcentover \K$. 
\end{dfn}

\begin{exa} \label{exa:3445}
If $\K$ is a {\em field}, then every DG  $\K$-ring $A$ is K-flat over 
$\K$.
\end{exa}

\begin{exa} \label{exa:3501}
If $A$ is a {\em semi-free} DG $\K$-ring (either commutative or 
noncommutative), then $A$ is a K-flat DG $\K$-ring. 
See Definition \ref{dfn:4309}, Proposition \ref{prop:4305} and Remark 
\ref{rem:4945}. 
\end{exa}

\begin{exa} \label{exa:3500}
If $A$ is a {\em nonpositive} DG $\K$-ring, and each 
$A^i$ is a {\em flat} $\K$-module, then $A$ is a K-flat DG $\K$-ring.
This includes the case of a flat $\K$-ring $A$. 
See Proposition \ref{prop:4535}. 
\end{exa}

\begin{exa} \label{exa:3525}
A very special case of Examples \ref{exa:3445} and \ref{exa:3500} is when 
$\K$ is a field and $A$ is a $\K$-ring (better known as a unital 
associative $\K$-algebra). 
This setting is the one commonly used in ring theory, see 
\cite{ATV}, \cite{ArZh}, \cite{Ye1}, \cite{VdB}, \cite{StaVdB}, \cite{YeZh1}. 
\end{exa}

\begin{dfn} \label{dfn:3595}
Let $A$ and $B$ be K-flat DG $\K$-rings. The 
{\em derived category of DG $B$-$A$-bimodules} 
\index{Derived category! of DG bimodules}
is the $\K$-linear triangulated category 
$\dcat{D}(B \ot A^{\mrm{op}})$. 
\end{dfn}

See Remark \ref{rem:1420} regarding the derived category of DG 
$B$-$A$-bimodules in the nonflat case.

\begin{dfn} \label{dfn:3450}
For a K-flat DG $\K$-ring $A$ we write 
$A^{\mrm{en}} := A \ot_{} A^{\mrm{op}}$, 
and call it the {\em enveloping DG ring} 
\index{Differential graded ring! enveloping}
of $A$ (relative to $\K$).
\end{dfn}

A explained in formula (\ref{eqn:3554}), the enveloping DG ring has a canonical 
isomorphism 
$A^{\mrm{en}} \iso (A^{\mrm{en}})^{\mrm{op}}$.

From here to the end of this section we assume the next convention, that 
strengthens Convention \ref{conv:3485}.  

\begin{conv} \label{conv:3445} 
In addition to the stipulations of Convention \ref{conv:3485}, we also 
assume by default that all DG rings are {\em K-flat central DG 
$\K$-rings}.
\end{conv}

Recall that a homomorphism of DG rings $f : A \to B$ induces a forgetful 
functor
$\opn{Rest}_f : \dcat{C}(B) \to \dcat{C}(A)$ 
that we call {\em restriction}. This is an exact DG functor, and thus 
it induces a triangulated functor  
\begin{equation} \label{eqn:3446}
\opn{Rest}_f : \dcat{D}(B) \to \dcat{D}(A) .
\end{equation}

\begin{prop} \label{prop:3445}
Let $f : A \to B$ be a DG ring homomorphism.
The triangulated functor $\opn{Rest}_f$ is conservative. Namely a morphism 
$\phi : M \to N$ in $\dcat{D}(B)$ is an isomorphism if and only if the morphism 
\[ \opn{Rest}_f(\phi) : \opn{Rest}_f(M) \to \opn{Rest}_f(N) \]
in  $\dcat{D}(A)$ is an isomorphism.
\end{prop}

\begin{proof}
We know by Corollary \ref{cor:2145} that the functor 
$\opn{H} : \dcat{D}(B) \to \dcat{G}(\K)$
is conservative. But 
$\opn{H} = \opn{H} \circ \opn{Rest}_f$.
\end{proof}

\begin{lem} \label{lem:3445}
Let $A$ and $B$ be DG rings.
\begin{enumerate}
\item If $P \in \dcat{C}(A \ot B^{\mrm{op}})$ is K-flat, then $P$ is K-flat 
over 
$A$.

\item  If $I \in \dcat{C}(A \ot B^{\mrm{op}})$ is K-injective, then $I$ is 
K-injective over $A$.

\item If $B$ is K-projective as a DG $\K$-module, and if 
$P \in \dcat{C}(A \ot_{\K} B^{\mrm{op}})$ is K-projective, 
then $P$ is K-projective over $A$.

\item If $B$ is a semi-free DG $\K$-ring, and if 
$P \in \dcat{C}(A \ot_{\K} B^{\mrm{op}})$ is semi-free, 
then $P$ is semi-free over $A$.
\end{enumerate}
\end{lem}

\begin{proof}
(1) and (2) are direct consequences of the canonical isomorphisms  
\[ M \ot_{A} P \cong (M \ot_{} B) \ot_{A \ot_{} B} P \]
and
\[ \opn{Hom}_A(N, I) \cong  \opn{Hom}_{A \ot_{} B}(N \ot_{} B, I)  \]
in $\dcat{C}_{\mrm{str}}(\K)$, 
for $M \in \dcat{C}(A^{\mrm{op}})$ and $N \in \dcat{C}(A)$,
together with the fact that $B$ is K-flat over $\K$. 
Items (3) and (4) are left as an exercise. 
\end{proof}

\begin{exer} \label{exer:3445}
Prove items (3) and (4) in this lemma. 
\end{exer}

\begin{prop} \label{prop:3450}
Let $A$, $B$ and $C$ be DG rings.
\begin{enumerate}
\item The DG bifunctor 
\[ (- \ot_{B} -) :  \dcat{C}(A \ot B^{\mrm{op}}) \times 
\dcat{C}(B \ot C^{\mrm{op}}) \to \dcat{C}(A \ot C^{\mrm{op}})  \]
has a triangulated left derived bifunctor 
\[ \bigl(  (- \ot^{\mrm{L}}_{B} -), \eta^{\mrm{L}} \bigr) :  
\dcat{D}(A \ot_{} B^{\mrm{op}}) \times 
\dcat{D}(B \ot_{} C^{\mrm{op}}) \to \dcat{D}(A \ot_{} C^{\mrm{op}}) . \]

\item Given $M \in \dcat{C}(A \ot_{} B^{\mrm{op}})$
and 
$N \in \dcat{C}(B \ot_{} C^{\mrm{op}})$,
such that $M$ is K-flat over $B^{\mrm{op}}$ or $N$ is K-flat over $B$, the 
morphism 
\[ \eta^{\mrm{L}}_{M, N} : M \ot_B^{\mrm{L}} N \to M \ot_{B} N \]
in $\dcat{D}(A \ot_{} C^{\mrm{op}})$ is an isomorphism. 

\item Suppose we are given DG ring homomorphisms 
$A' \to A$, $f : B' \to B$ and $C' \to C$, such that $f$ is a 
quasi-isomorphism. 
Then the diagram 
\[ \UseTips \xymatrix @C=10ex @R=6ex {
\dcat{D}(A \ot_{} B^{\mrm{op}}) \times \dcat{D}(B \ot_{} C^{\mrm{op}})
\ar[r]^(0.6){ (- \ot^{\mrm{L}}_{B} -) }
\ar[d]_{\opn{Rest} \times \opn{Rest}}
&
\dcat{D}(A \ot_{} C^{\mrm{op}})
\ar[d]^{\opn{Rest}}
\\
\dcat{D}(A' \ot_{} B'^{\, \mrm{op}}) \times \dcat{D}(B' \ot C'^{\, \mrm{op}})
\ar[r]^(0.6){ (- \ot^{\mrm{L}}_{B'} -) }
&
\dcat{D}(A' \ot_{} C'^{\, \mrm{op}})
} \]
is commutative up to an isomorphism of triangulated bifunctors. 

\item Suppose $D$ is another DG ring. Then there is an isomorphism  
\[ ((- \ot^{\mrm{L}}_{B} -) \ot^{\mrm{L}}_{C} -) \cong 
(- \ot^{\mrm{L}}_{B} (- \ot^{\mrm{L}}_{C} -)) \]
of triangulated trifunctors 
\[ \dcat{D}(A \ot_{} B^{\mrm{op}}) \times \dcat{D}(B \ot_{} C^{\mrm{op}}) 
\times 
\dcat{D}(C \ot_{} D^{\mrm{op}}) \to \dcat{D}(A \ot_{} D^{\mrm{op}}) . \]
\end{enumerate}
\end{prop}

\begin{proof} \mbox{}

\smallskip \noindent
(1) The DG bifunctor $(- \ot_{B} -)$ induces a triangulated bifunctor 
\[ (- \ot_{B} -) :  \dcat{K}(A \ot_{} B^{\mrm{op}}) \times 
\dcat{K}(B \ot_{} C^{\mrm{op}}) \to \dcat{K}(A \ot_{} C^{\mrm{op}})  \]
on the homotopy categories, in the obvious way. 
By Corollary \ref{cor:1580}, Proposition \ref{prop:1525} and Lemma 
\ref{lem:3445}(1), every DG bimodule $N \in \dcat{C}(B \ot_{} C^{\mrm{op}})$
admits a quasi-isomorphism 
$\rho_N : P \to N$, where 
$P \in \dcat{C}(B \ot_{} C^{\mrm{op}})$
is K-flat over $B$.
Given $M \in \dcat{D}(A \ot_{} B^{\mrm{op}})$, let us define the object
\[ M \ot^{\mrm{L}}_{B} N := M \ot_{B} P \in
\dcat{D}(A \ot_{} C^{\mrm{op}}) , \]
with the morphism 
\[ \eta^{\mrm{L}}_{M, N} := \opn{id}_M \ot_B \, \rho_N : 
M \ot_B^{\mrm{L}} N \to M \ot_{B} N . \]
The pair 
$\bigl( (- \ot^{\mrm{L}}_{B} -), \eta^{\mrm{L}} \bigr)$
is a left derived bifunctor of $(- \ot_{B} -)$. 

\medskip \noindent 
(2) Under either assumption the homomorphism 
$\opn{id}_M \ot_B \, \rho_N$
is a quasi-iso\-morphism; cf.\ Lemma \ref{lem:2109}.

\medskip \noindent 
(3) This is similar to the proof of item (2) of Theorem \ref{thm:2363}. 
Take $M \in \dcat{C}(A \ot B^{\mrm{op}})$ and 
$N \in \dcat{C}(B \ot C^{\mrm{op}})$. 
Let $\rho_N : P \to N$ be the resolution from item (1).
Choose a resolution 
$\rho'_P : P' \to P$ in $\dcat{C}_{\mrm{str}}(B' \ot C'^{\, \mrm{op}})$
where $P'$ is K-flat over $B'$. Then we have a canonical isomorphism 
\[ \opn{Rest}(M) \ot^{\mrm{L}}_{B'} \opn{Rest}(N) \cong M \ot_{B'} P' \]
in $\dcat{D}(A' \ot_{} C'^{\, \mrm{op}})$.
Now in the commutative diagram 
\[ \UseTips \xymatrix @C=10ex @R=6ex {
B' \ot_{B'} P'
\ar[r]^(0.6){\cong}
\ar[d]_{f \ot \, \opn{id}_{P'}}
&
P'
\ar[d]^{\rho'_P}
\\
B \ot_{B'} P'
\ar[r]^(0.6){\opn{id}_B \ot \, \rho'_P}
&
P
} \]
in $\dcat{C}_{\mrm{str}}(B')$ the vertical arrows are quasi-isomorphisms, 
because $P'$ is K-flat over $B'$. Therefore the bottom arrow 
$\opn{id}_B \ot \, \rho'_P$ is a quasi-isomorphism. 

We now look at this commutative diagram 
\[ \UseTips \xymatrix @C=6ex @R=6ex {
M \ot_{B'} P'
\ar@(ur,ul)[rr]^{\opn{id}_{M} \ot \, \rho'_P}
\ar[r]^(0.4){\cong}
&
M \ot_B B \ot_{B'} P'
\ar[r]^(0.6){\phi}
&
M \ot_B P
} \]
in $\dcat{C}_{\mrm{str}}(A' \ot_{} C'^{\, \mrm{op}})$, where 
$\phi := \opn{id}_M \ot \opn{id}_B \ot \, \rho'_P$. 
Because both $P$ and $B \ot_{B'} P'$ are K-flat over $B$, and 
$\opn{id}_B \ot \, \rho'_P$ is a quasi-isomorphism, 
Lemma \ref{lem:2109} tells us that $\phi$ is a quasi-isomorphism. 
Hence $\opn{id}_{M} \ot \, \rho'_P$ is a quasi-isomorphism. 
We get the desired canonical isomorphism 
\[ \opn{Q}(\opn{id}_{M} \ot \, \rho'_P) : 
\opn{Rest}(M) \ot^{\mrm{L}}_{B'} \opn{Rest}(N) \iso 
\opn{Rest}(M \ot^{\mrm{L}}_{B} N) \]
in $\dcat{D}(A' \ot_{} C'^{\, \mrm{op}})$.
 
\medskip \noindent 
(4) Given $M, N, P$ as above and 
$L \in \dcat{C}(C \ot_{} D^{\mrm{op}})$, 
we choose a quasi-iso\-morphism 
$\rho_L : Q \to L$, where 
$Q \in \dcat{C}(C \ot_{} D^{\mrm{op}})$
is K-flat over $C$. A small calculation shows that 
$P \ot_C Q$ is K-flat over $B$. The desired isomorphism 
\[ (M \ot^{\mrm{L}}_{B} N) \ot^{\mrm{L}}_{C} L \cong  
M \ot^{\mrm{L}}_{B} (N \ot^{\mrm{L}}_{C} L) \]
in $\dcat{D}(A \ot_{} D^{\mrm{op}})$
comes from the obvious isomorphism 
\[ (M \ot_{B} P) \ot_{C} Q \cong  
M \ot_{B} (P \ot_{C} Q) \]
in $\dcat{C}_{\mrm{str}}(A \ot_{} D^{\mrm{op}})$.
\end{proof}

Here is some notation: suppose we are given morphisms 
$\phi : M' \to M$ in 
$\dcat{D}(A \ot B^{\mrm{op}})$
and $\psi : N' \to N$ in $\dcat{D}(B \ot C^{\mrm{op}})$.
The result of applying the bifunctor $(- \ot_B^{\mrm{L}} -)$ is the morphism
\begin{equation} \label{eqn:3475}  
 \phi \ot^{\mrm{L}}_B \psi : M' \ot_B^{\mrm{L}} N' \to M \ot_B^{\mrm{L}} N
\end{equation}
in $\dcat{D}(A \ot C^{\mrm{op}})$. 

\begin{prop} \label{prop:1031}
Let $A$, $B$ and $C$ be DG rings. 
\begin{enumerate}
\item The DG bifunctor 
\[ \opn{Hom}_B(-, -) :  \dcat{C}(B \ot A^{\mrm{op}})^{\mrm{op}} \times 
\dcat{C}(B \ot C^{\mrm{op}}) \to \dcat{C}(A \ot C^{\mrm{op}})  \]
has a triangulated right derived bifunctor 
\[ \bigl( \opn{RHom}_B(-, -), \eta^{\mrm{R}} \bigr) :  
\dcat{D}(B \ot A^{\mrm{op}})^{\mrm{op}} \times 
\dcat{D}(B \ot C^{\mrm{op}}) \to \dcat{D}(A \ot C^{\mrm{op}}) . \]

\item Given $M \in \dcat{C}(B \ot A^{\mrm{op}})$
and $N \in \dcat{C}(B \ot C^{\mrm{op}})$,
such that either $M$ is K-projective over $B$ or $N$ is K-injective over $B$, 
the morphism 
\[ \eta^{\mrm{R}}_{M, N} : \opn{Hom}_B(M, N) \to \opn{RHom}_B(M, N) \]
in $\dcat{D}(A \ot C^{\mrm{op}})$ is an isomorphism. 

\item Suppose we are given DG ring homomorphisms 
$A' \to A$, $f : B' \to B$ and $C' \to C$, such that $f$ is a 
quasi-isomorphism. 
Then the diagram 
\[ \UseTips \xymatrix @C=14ex @R=6ex {
\dcat{D}(B \ot A^{\mrm{op}})^{\mrm{op}} \times \dcat{D}(B \ot C^{\mrm{op}})
\ar[r]^(0.6){ \opn{RHom}_B(-, -) }
\ar[d]_{\opn{Rest} \times \opn{Rest}}
&
\dcat{D}(A \ot C^{\mrm{op}})
\ar[d]^{\opn{Rest}}
\\
\dcat{D}(B' \ot A'^{\, \mrm{op}})^{\mrm{op}} \times 
\dcat{D}(B' \ot C'^{\, \mrm{op}})
\ar[r]^(0.6){ \opn{RHom}_{B'}(-, -) }
&
\dcat{D}(A' \ot C'^{\, \mrm{op}})
} \]
is commutative up to an isomorphism of triangulated bifunctors. 
\end{enumerate}
\end{prop}

\begin{exer} \label{exer:3560}
Prove Proposition \ref{prop:1031}. 
Hints: this is similar to the proof of items (1)-(3) of Proposition 
\ref{prop:3450}, but now we rely on Corollary \ref{cor:1665} and Lemma 
\ref{lem:3445}(2) for the existence of resolutions $\rho_N : N \to I$ 
in $\dcat{C}_{\mrm{str}}(B \ot C^{\mrm{op}})$ by DG bimodules $I$
that are K-injective over $B$. For item (3) use Lemma \ref{lem:2000}.
\end{exer}

Again there is related notation: suppose we are given morphisms
$\phi : M \to M'$ in $\dcat{D}(B \ot A^{\mrm{op}})$
and $\psi : N' \to N$ in $\dcat{D}(B \ot C^{\mrm{op}})$.
The result of applying the bifunctor $\opn{RHom}_B(-, -)$
is the morphism
\begin{equation} \label{eqn:3476}  
\opn{RHom}_B(\phi, \psi) : \opn{RHom}_B(M', N') \to \opn{RHom}_B(M, N) 
\end{equation}
in $\dcat{D}(A \ot C^{\mrm{op}})$.

\begin{prop} \label{prop:3550}
Let $A$, $B$, $C$ and $D$ be DG rings. For 
$M \in \dcat{D}(B \ot A^{\mrm{op}})$,
$N \in \dcat{D}(C \ot B^{\mrm{op}})$ and 
$L \in \dcat{D}(C \ot D^{\mrm{op}})$
there is an isomorphism 
\[ \opn{RHom}_B \bigl( M, \opn{RHom}_C(N, L) \bigr) \cong 
\opn{RHom}_C(N \ot_B^{\mrm{L}} M, L) \]
in $\dcat{D}(A \ot D^{\mrm{op}})$ called {\em derived Hom-tensor adjunction}. 
This isomorphism is functorial in the objects $M, N, L$.
\end{prop}

\begin{proof}
We choose a quasi-isomorphism $L \to J$ in 
$\dcat{C}_{\mrm{str}}(C \ot D^{\mrm{op}})$ into a DG module $J$ that is 
K-injective over $C$, and a quasi-isomorphism $Q \to N$ in 
$\dcat{C}_{\mrm{str}}(C \ot B^{\mrm{op}})$ from a DG module $Q$ that is K-flat 
over $B^{\mrm{op}}$. A calculation shows that the DG bimodule 
$\opn{Hom}_C(Q, J)$ is K-injective over $B$. 
The desired isomorphism comes from the obvious adjunction isomorphism 
\[ \opn{Hom}_B(M, \opn{Hom}_C(Q, J)) \cong 
\opn{Hom}_C(Q \ot_B M, J) \]
in $\dcat{C}_{\mrm{str}}(A \ot D^{\mrm{op}})$. 
\end{proof}

\begin{prop} \label{prop:3565}
Let $A$, $B$, $C$ and $D$ be DG rings. For 
$L \in \dcat{D}(A \ot C^{\mrm{op}})$,
$M \in \dcat{D}(A \ot B^{\mrm{op}})$ and 
$N \in \dcat{D}(B \ot D^{\mrm{op}})$
there is a morphism 
\[ \opn{ev}^{\mrm{R, L}}_{L, M, N} :
\opn{RHom}_A(L, M) \ot_B^{\mrm{L}} N \to 
\opn{RHom}_A(L, M \ot_B^{\mrm{L}} N) \]
in $\dcat{D}(C \ot D^{\mrm{op}})$, called {\em derived tensor-evaluation}. 
This morphism is functorial in the objects $L, M, N$.
Moreover, after applying the restriction functor  
$\dcat{D}(C \ot D^{\mrm{op}}) \to \dcat{D}(\K)$, 
the morphism $\opn{ev}^{\mrm{R, L}}_{L, M, N}$ coincides with the morphism from 
Theorem \tup{\ref{thm:4320}}. 
\end{prop}

The reason we can give stronger statement here, as compared to 
Theorem \ref{thm:4320}, is because here our DG rings are K-flat 
over $\K$. 

\begin{proof}
Choose a quasi-isomorphism $Q \to N$ in 
$\dcat{C}_{\mrm{str}}(B \ot D^{\mrm{op}})$ from a DG module $Q$ that is 
K-flat over $B$, a quasi-isomorphism $M \to I$ in 
$\dcat{C}_{\mrm{str}}(A \ot B^{\mrm{op}})$ into a DG module $I$ that is 
K-injective over $A$, and a quasi-isomorphism 
$\rho : I \ot_B Q \to J$ in $\dcat{C}_{\mrm{str}}(A \ot D^{\mrm{op}})$ into a 
DG module $J$ that is K-injective over $A$. Then 
$\opn{ev}^{\mrm{R, L}}_{L, M, N}$
is represented by the composed homomorphism 
\[  \opn{Hom}_A(L, I) \ot_B Q \xar{ \opn{ev}_{L, I, Q} }  
\opn{Hom}_A(L, I \ot_B Q) \xar{\til{\rho}} \opn{Hom}_A(L, J)  \]
in $\dcat{C}_{\mrm{str}}(C \ot D^{\mrm{op}})$,
where $\til{\rho} := \opn{Hom}(\opn{id}_L, \rho)$. 

When we forget $C$ and $D$, then we can choose a K-projective resolution 
$P \to L$ in $\dcat{C}_{\mrm{str}}(A)$. Then we have a commutative diagram 
\[ \UseTips \xymatrix @C=8ex @R=6ex {
\opn{Hom}_A(P, I) \ot_B Q
\ar[r]^{ \opn{ev}_{P, I, Q} }
&
\opn{Hom}_A(P, I \ot_B Q) 
\ar[r]^(0.55){\mrm{q.i.}}
&
\opn{Hom}_A(P, J)
\\
\opn{Hom}_A(L, I) \ot_B Q
\ar[r]^{ \opn{ev}_{L, I, Q} }
\ar[u]^{\mrm{q.i.}}
&
\opn{Hom}_A(L, I \ot_B Q) 
\ar[r]^(0.55){\til{\rho}}
\ar[u]
&
\opn{Hom}_A(L, J)
\ar[u]^{\mrm{q.i.}}
} \]
in $\dcat{C}_{\mrm{str}}(\K)$, in which the arrows marked ``q.i.'' are 
quasi-isomorphisms. By comparing to the proof of Theorem \ref{thm:4320}
this proves that last assertion.
\end{proof}

According to Propositions \ref{prop:3450}
and \ref{prop:1031} there are triangulated bifunctors 
\begin{equation} \label{eqn:3477}
(- \ot^{\mrm{L}}_{A} -) : \dcat{D}(A^{\mrm{en}}) \times 
\dcat{D}(A^{\mrm{en}}) \to \dcat{D}(A^{\mrm{en}}) ,
\end{equation}
\begin{equation} \label{eqn:3478}
(- \ot^{\mrm{L}}_{A} -) : \dcat{D}(A^{\mrm{en}}) \times 
\dcat{D}(A \ot B^{\mrm{op}}) \to \dcat{D}(A \ot B^{\mrm{op}}) 
\end{equation}
and
\begin{equation} \label{eqn:1479}
\opn{RHom}_A(-, -) : \dcat{D}(A^{\mrm{en}})^{\mrm{op}} \times 
\dcat{D}(A \ot B^{\mrm{op}}) \to \dcat{D}(A \ot B^{\mrm{op}}) .
\end{equation}

\begin{rem} \label{rem:3450}
The category $\dcat{D}(A^{\mrm{en}})$, with the tensor operation 
$(- \ot^{\mrm{L}}_{A} -)$, is a {\em monoidal category}.
Moreover, $\dcat{D}(A^{\mrm{en}})$ is a {\em biclosed monoidal category}: the 
two internal Hom operations are 
$\opn{RHom}_A(-, -)$ and $\opn{RHom}_{A^{\mrm{op}}}(-, -)$.
If $A$ is weakly commutative, then $\dcat{D}(A^{\mrm{en}})$ is a symmetric 
monoidal category. See \cite[Section 5]{VyYe}, \cite{Mac2} and \cite{nLab}. 
\end{rem}

\begin{rem} \label{rem:1420}
Here is an outline of our way to handle derived categories of DG bimodules in 
the absence of flatness. The idea is to choose K-flat resolutions 
$\til{A} \to A$ and $\til{B} \to B$ in the category 
$\catt{DGRng} \centover \K$
of DG central $\K$-rings. This can be done by 
Theorem \ref{thm:4310} and Proposition \ref{prop:4305}. 
Then the {\em derived category of DG $A$-$B$-bimodules} 
\index{Derived category! of DG bimodules}
is the triangulated 
category $\dcat{D}(\til{A} \ot_{} \til{B}^{\mrm{op}})$.
Note that the restriction functors 
$\dcat{D}(\til{A}) \to \dcat{D}(A)$ and 
$\dcat{D}(\til{B}^{\mrm{op}}) \to \dcat{D}(B^{\mrm{op}})$ 
are equivalences. This says that we have triangulated bifunctors 
\[ (- \ot^{\mrm{L}}_{\til{B}} -) : 
\dcat{D}(\til{A} \ot_{} \til{B}^{\mrm{op}}) \times 
\dcat{D}(B) \to \dcat{D}(A) \]
and
\[ \opn{RHom}_{\til{A}}(-, -) : 
\dcat{D}(\til{A} \ot_{} \til{B}^{\mrm{op}})^{\mrm{op}} \times 
\dcat{D}(A) \to \dcat{D}(B) . \]
We also have the ``tilde'' versions of the functors 
(\ref {eqn:3477}), (\ref{eqn:3478}) and (\ref{eqn:1479}). 
These are collectively called the {\em package of standard derived functors}. 

The triangulated category 
$\dcat{D}(\til{A} \ot_{} \til{B}^{\mrm{op}})$
is independent of the resolutions, up to a canonical equivalence.
The argument is this. 
Suppose we are given K-flat resolutions 
$\til{A}_i \to A$ and $\til{B}_i \to B$ in the category 
$\catt{DGRng} \centover \K$, indexed by $i \in I$ for some set $I$. 
For each $i, j \in I$, the 
the DG module $A \ot^{\mrm{L}}_{\K} B$
has an outside action by the DG ring 
$\til{A}_j \ot_{} \til{B}_j^{\mrm{op}}$
and an inside action the DG ring 
$\til{A}_i \ot_{} \til{B}_i^{\mrm{op}}$. 
(See Subsection \ref{subsec:RNCDC-uniq} regarding outside and inside actions.)
These give rise to a DG bimodule
\[ T_{i, j} := A \ot^{\mrm{L}}_{\K} B \in 
\dcat{D} \bigl( (\til{A}_j \ot_{} \til{B}_j^{\mrm{op}}) \ot_{}
(\til{A}_i \ot_{} \til{B}_i^{\mrm{op}})^{\mrm{op}} \bigr) . \]
It turns out that $T_{i, j}$ is a {\em tilting DG bimodule}, as defined in the 
next subsection; and hence the triangulated functor
\begin{equation} \label{eqn:5090}
G_{i, j} := T_{i, j} \ot^{\mrm{L}}_{\til{A}_i \ot \til{B}_i^{\mrm{op}}} (-) :
\dcat{D}(\til{A}_i \ot \til{B}_i^{\mrm{op}}) \to
\dcat{D}(\til{A}_j \ot \til{B}_j^{\mrm{op}}) 
\end{equation}
is an equivalence. For a third index $k \in I$ there is a canonical 
isomorphism
\[ T_{j, k} \ot^{\mrm{L}}_{\til{A}_j \ot_{} \til{B}_j^{\mrm{op}}} 
T_{i, j} \cong T_{i, k} 
\quad \tup{in} \quad 
\dcat{D} \bigl( (\til{A}_k \ot_{} \til{B}_k^{\mrm{op}}) \ot_{}
(\til{A}_i \ot_{} \til{B}_i^{\mrm{op}})^{\mrm{op}} \bigr) . \]
Thus there are isomorphisms of triangulated functors 
\begin{equation} \label{eqn:3501}
G_{j, k} \circ G_{i, j} \iso G_{i, k} ,
\end{equation}
and they satisfy the pentagon axiom. 
The equivalences $G_{i, j}$ respect the packages of standard derived 
functors.

We see that there is a $\K$-linear triangulated category that we will 
symbolically denote by 
$\dcat{D}(A \ot^{\mrm{L}}_{\K} B^{\mrm{op}})$,
and it is canonically equivalent to all the triangulated categories 
$\dcat{D}(\til{A}_i \ot_{} \til{B}_i^{\mrm{op}})$.
For more details see the lecture notes \cite{Ye16} 
or the paper \cite{Ye14}. 
\end{rem}

\begin{rem} \label{rem:1421}
Actually, we can say more. 
Recall the $2$-category $\catt{TrCat} / \K$ of $\K$-linear triangulated 
categories that was briefly mentioned at the end of Subsection 
\ref{subsec:2-cat-not}. 
There is a pseudofunctor 
\begin{equation} \label{eqn:3502}
\opn{DerCat} : \catt{DGRng} \centover \K \to \catt{TrCat} / \K
\end{equation}
that sends a DG ring $A$ to the triangulated category
$\opn{DerCat}(A) := \dcat{D}(A)$,
and a homomorphism $f : A \to B$ in $\catt{DGRng} \centover \K$
is sent to the triangulated functor 
\[ \opn{DerCat}(f) := \opn{LInd}_f = B \ot^{\mrm{L}}_{A} (-)
: \dcat{D}(A) \to \dcat{D}(B) . \]
Note that the category $\catt{DGRng} \centover \K$ is not linear (the morphisms 
sets have no linear structure), whereas the 
$2$-category $\catt{TrCat} / \K$ is linear (in the sense that the sets of 
$2$-morphisms are $\K$-modules). 

Define $\dcat{D}(\catt{DGRng} \centover \K)$
to be the categorical localization of the category  
$\catt{DGRng} \centover \K$ with respect to the quasi-isomorphisms in it. Then 
the pseudofunctor $\opn{DerCat}$ from (\ref{eqn:3502}) localizes, in the 
categorical sense, to a pseudofunctor
\[ \opn{DerCat} : \dcat{D}(\catt{DGRng} \centover \K) \to \catt{TrCat} / \K . \]

The (nonlinear) bifunctor 
\[ (- \ot_{\K} -) : 
(\catt{DGRng} \centover \K) \times (\catt{DGRng} \centover \K) \to
\catt{DGRng} \centover \K \]
has a left derived bifunctor 
\[ (- \ot^{\mrm{L}}_{\K} -) : 
\dcat{D}(\catt{DGRng} \centover \K) \times 
\dcat{D}(\catt{DGRng} \centover \K) \to 
\dcat{D}(\catt{DGRng} \centover \K) . \]
Given a pair of DG rings $A$ and $B$, we thus have DG ring 
$A \ot^{\mrm{L}}_{\K} B^{\mrm{op}}$, 
well-defined up to a unique isomorphism in 
$\dcat{D}(\catt{DGRng} \centover \K)$. 
The derived category of DG $A$-$B$-bimodules described at the end of Remark 
\ref{rem:1420} is canonically equivalent to 
$\opn{DerCat}(A \ot^{\mrm{L}}_{\K} B^{\mrm{op}})$
in $\catt{TrCat} / \K$. 
\end{rem}

\mysubsection{Tilting DG Bimodules} \label{subsec:tilting}

We continue with Convention \ref{conv:3445}. In particular all DG rings are 
K-flat central over the base ring $\K$, $\ot = \ot_{\K}$, 
and $A^{\mrm{en}} = A \ot A^{\mrm{op}}$ for a DG ring $A$. 

\begin{dfn} \label{dfn:3455}
Let $A$ and $B$ be DG rings.
An object $T \in \dcat{D}(B \ot_{} A^{\mrm{op}})$
is called a {\em tilting DG $B$-$A$-bimodule} 
\index{Tilting! DG bimodule}
if there exists some object $S \in \dcat{D}(A \ot_{} B^{\mrm{op}})$,
and isomorphisms 
$S \ot^{\mrm{L}}_B T \cong A$ in $\dcat{D}(A^{\mrm{en}})$
and
$T \ot^{\mrm{L}}_A S \cong B$ in $\dcat{D}(B^{\mrm{en}})$.

When $B = A$, so that $B \ot_{} A^{\mrm{op}} = A^{\mrm{en}}$,
we use the term {\em tilting DG $A$-bimodule} as an abbreviation for the term
tilting DG $A$-$A$-bimodule. 
\end{dfn}

It is clear from the symmetry of the definition that the object 
$S \in \dcat{D}(A \ot_{} B^{\mrm{op}})$ is a 
tilting DG $A$-$B$-bimodule. 

\begin{lem} \label{lem:3455}
Let $T \in \dcat{D}(B \ot A^{\mrm{op}})$ and 
$S, S' \in \dcat{D}(A \ot B^{\mrm{op}})$
satisfy 
$S \ot^{\mrm{L}}_B T \cong A$ in $\dcat{D}(A^{\mrm{en}})$ 
and
$T \ot^{\mrm{L}}_A S' \cong B$ in $\dcat{D}(B^{\mrm{en}})$.
Then 
$S' \cong S$ in $\dcat{D}(A \ot_{} B^{\mrm{op}})$.
Therefore $T$ is a tilting DG $B$-$A$-bimodule, and $S$ is a tilting 
DG $A$-$B$-bimodule. 
\end{lem}

\begin{proof}
Using the associativity of derived tensor products from Proposition \lb
\ref{prop:3450}(4)
we have isomorphisms 
\[ S \cong S \ot^{\mrm{L}}_B B \cong
S \ot^{\mrm{L}}_B (T \ot^{\mrm{L}}_A S') \cong 
(S \ot^{\mrm{L}}_B T) \ot^{\mrm{L}}_A S' \cong 
A \ot^{\mrm{L}}_A S' \cong S' \]
in $\dcat{D}(A \ot_{} B^{\mrm{op}})$. 
\end{proof}

The lemma implies that the tilting DG $A$-$B$-bimodule $S$ in Definition 
\ref{dfn:3455} is unique up to isomorphism.

\begin{dfn} \label{3460}
The tilting DG $A$-$B$-bimodule $S$ in Definition \ref{dfn:3455} is called 
the {\em quasi-inverse of $T$}. 
\index{Differential graded bimodule! quasi-inverse of tilting}
\end{dfn}

Here are some general properties of tilting DG bimodules. 

\begin{prop} \label{prop:3490}
Let $A, B, C$ be DG rings, let $T$ be a tilting DG $B$-$A$-bimodule, and let 
$S$ be a tilting DG $C$-$B$-bimodule. 
Then $S \ot^{\mrm{L}}_{B} T$ is a tilting $C$-$A$-bimodule. 
\end{prop}

\begin{proof}
Let $T^{\vee}$ be the quasi-inverse of $T$, and let 
$S^{\vee}$ be the quasi-inverse of $S$. Then
\[ \begin{aligned}
& (T^{\vee} \ot^{\mrm{L}}_{B} S^{\vee}) \ot^{\mrm{L}}_{C} 
(S \ot^{\mrm{L}}_{B} T) \cong 
T^{\vee} \ot^{\mrm{L}}_{B} (S^{\vee} \ot^{\mrm{L}}_{C} S )\ot^{\mrm{L}}_{B} T
\\ & \quad 
\cong T^{\vee} \ot^{\mrm{L}}_{B} B \ot^{\mrm{L}}_{B} T \cong A 
\end{aligned} \]
in $\dcat{D}(A^{\mrm{en}})$. Likewise 
\[ (S \ot^{\mrm{L}}_{B} T) \ot^{\mrm{L}}_{A} 
(T^{\vee} \ot^{\mrm{L}}_{B} S^{\vee}) \cong C \]
in $\dcat{D}(C^{\mrm{en}})$. 
\end{proof}

\begin{prop} \label{prop:3455}
Suppose $T$ is a tilting DG $B$-$A$-bimodule, and 
$S$ is a tilting DG $A$-$B$-bimodule. Then the triangulated functor 
\[ G_{T, S} : \dcat{D}(A^{\mrm{en}}) \to \dcat{D}(B^{\mrm{en}}) \]
with formula
$G_{T, S}(M) := T \ot^{\mrm{L}}_{A} M \ot^{\mrm{L}}_{A} S$ 
is an equivalence. 
\end{prop}

\begin{prop} \label{prop:3525}
Suppose $T$ is a tilting DG $B$-$A$-bimodule, with quasi-inverse 
$T^{\vee}$. Fix an isomorphism 
$T^{\vee} \ot^{\mrm{L}}_{B} T \iso A$
in $\dcat{D}(A^{\mrm{en}})$. 
Then for every $M_1, M_2 \in \dcat{D}(A^{\mrm{en}})$ there is an  
isomorphism 
\[ G_{T, T^{\vee}} (M_1 \ot^{\mrm{L}}_{A} M_2) \iso 
G_{T, T^{\vee}} (M_1) \ot^{\mrm{L}}_{B} G_{T, T^{\vee}} (M_2) \]
in $\dcat{D}(B^{\mrm{en}})$. This isomorphism is bifunctorial in
$(M_1, M_2)$. 
\end{prop}

\begin{exer} \label{exer:3590}
Prove Propositions \ref{prop:3455} and \ref{prop:3525}.
\end{exer}

\begin{dfn} \label{dfn:3457}
Let $A$ be a K-flat central DG $\K$-ring. The 
{\em noncommutative derived Picard group of $A$ relative to $\K$} 
\index{Derived Picard group}
\index{Picard group! noncommutative derived}
\index{1-DPic(A)@$\opn{DPic}_{\K}(A)$}
is the group $\opn{DPic}_{\K}(A)$, whose 
elements are the isomorphism classes in $\dcat{D}(A^{\mrm{en}})$ of the 
tilting DG $A$-bimodules. The multiplication in this group is induced from 
$(- \ot^{\mrm{L}}_{A} -)$, and the unit element is the class of the DG bimodule 
$A$. 
\end{dfn}

This definition is valid by Proposition \ref{prop:3490}. 
Propositions \ref{prop:3455} and \ref{prop:3525} say that: 

\begin{cor} \label{cor:3525}
A tilting DG $B$-$A$-bimodule $T$, with quasi-inverse $T^{\vee}$, 
and an isomorphism 
$T^{\vee} \ot^{\mrm{L}}_{B} T \cong A$ in $\dcat{D}(A^{\mrm{en}})$,
induce a group isomorphism 
\[ \opn{DPic}_{\K}(A) \iso \opn{DPic}_{\K}(B) . \]
\end{cor}

We don't know much about the structure of the derived Picard group 
$\opn{DPic}_{\K}(A)$ in general. However, when $A$ is either a ring
or a commutative DG ring, there are a few strong structural results -- see 
next subsection. 

\begin{rem} \label{rem:3595}
If $A$ is not K-flat over the base ring $\K$, then the derived Picard group 
should be defined using a K-flat resolution $\til{A} \to A$, as explained in 
Remark \ref{rem:1420}. Namely $\opn{DPic}_{\K}(A)$ is the group whose elements 
are the isomorphism classes in $\dcat{D}(\til{A}^{\mrm{en}})$ of the tilting DG 
$\til{A}$-bimodules, etc. This group is independent, up to a canonical 
group isomorphism, of the resolution $\til{A} \to A$. 
\end{rem}

The derived homothety morphism in the commutative setting,
and the related derived Morita property, were defined in subsection 
\ref{subsec:du-cplxs}. In the current noncommutative setting these notions 
become more involved, as we shall now see.

Given a DG bimodule $M \in \dcat{C}(B \ot_{\K} A^{\mrm{op}})$, 
there are DG ring homomorphisms 
\[  \opn{hm}_{M, A^{\mrm{op}}} : A^{\mrm{op}} \to \opn{End}_{B}(M) = 
\opn{Hom}_{B}(M, M) \]
and
\[  \opn{hm}_{M, B} : B \to 
\opn{End}_{A^{\mrm{op}}}(M) = \opn{Hom}_{A^{\mrm{op}}}(M, M) \]
that we call the {\em noncommutative homothety homomorphisms through 
$A^{\mrm{op}}$ and \lb through $B$} respectively. 
When we forget the ring structures, these become homomorphisms
\begin{equation} \label{eqn:4175}
\opn{hm}_{M, A^{\mrm{op}}} : A \to \opn{Hom}_{B}(M, M) 
\end{equation}
and
\begin{equation} \label{eqn:4176}
\opn{hm}_{M, B} : B \to \opn{Hom}_{A^{\mrm{op}}}(M, M)
\end{equation}
in $\dcat{C}_{\mrm{str}}(A^{\mrm{en}})$ and 
$\dcat{C}_{\mrm{str}}(B^{\mrm{en}})$ respectively. 

\begin{dfn} \label{dfn:3565}
Let $M \in \dcat{D}(B \ot_{} A^{\mrm{op}})$.
\index{Homothety morphism! noncommutative derived}
\begin{enumerate}
\item The {\em noncommutative derived homothety morphism of $M$ through
$A^{\mrm{op}}$} is the morphism 
\[ \opn{hm}^{\mrm{R}}_{M, A^{\mrm{op}}} := \eta^{\mrm{R}}_{M, M} \circ 
\opn{Q}(\opn{hm}_{M, A^{\mrm{op}}}) : A \to \opn{RHom}_{B}(M, M) \]
in $\dcat{D}(A^{\mrm{en}})$.

\item The {\em noncommutative derived homothety morphism of $M$ through $B$} is 
the morphism 
\[ \opn{hm}^{\mrm{R}}_{M, B} := \eta^{\mrm{R}}_{M, M} 
\circ \opn{Q}(\opn{hm}_{M, B}) : B \to \opn{RHom}_{A^{\mrm{op}}}(M, M)  \]
in $\dcat{D}(B^{\mrm{en}})$.
\end{enumerate}
\end{dfn}

Here is a commutative diagram in $\dcat{D}(A^{\mrm{en}})$ depicting the 
noncommutative derived homothety morphism through $A^{\mrm{op}}$.
\[ \UseTips \xymatrix @C=8ex @R=6ex {
A
\ar[d]_{\opn{Q}( \opn{hm}_{M, A^{\mrm{op}}} )}
\ar[dr]^{ \opn{hm}^{\mrm{R}}_{M, A^{\mrm{op}}} }
\\
\opn{Hom}_{B}(M, M)
\ar[r]^{\eta^{\mrm{R}}_{M, M}}
&
\opn{RHom}_{B}(M, M)
} \]
The morphism $\eta^{\mrm{R}}_{M, M}$ is part of the right derived bifunctor 
$\opn{RHom}_{B}(-, -)$; see Proposition \ref{prop:1031}. 

\begin{dfn} \label{3455}
Let $A$ and $B$ be DG rings, and let $M$ be an object of 
$\dcat{D}(B \ot_{} A^{\mrm{op}})$.
\index{Derived Morita property! noncommutative}
\begin{enumerate}
\item We say that $M$ has the {\em noncommutative derived Morita property on 
the $B$ side} if the derived homothety morphism
\[ \opn{hm}^{\mrm{R}}_{M, A^{\mrm{op}}} : A \to \opn{RHom}_{B}(M, M) \]
in $\dcat{D}(A^{\mrm{en}})$ is an isomorphism. 

\item We say that $M$ has the {\em noncommutative derived Morita property on 
the $A^{\mrm{op}}$ side} if the derived homothety morphism
\[ \opn{hm}^{\mrm{R}}_{M, B} : B \to 
\opn{RHom}_{A^{\mrm{op}}}(M, M)  \]
in $\dcat{D}(B^{\mrm{en}})$ is an isomorphism. 

\item We say that $M$ has the {\em noncommutative derived Morita property on 
both sides} if it has the noncommutative derived Morita property on 
the $B$ side and on the $A^{\mrm{op}}$ side.
\end{enumerate}
\end{dfn}

Recall the restriction functors from diagram (\ref{eqn:3536}). 
Compact generators were introduced in Definition \ref{dfn:3431}(3). 
The next definition resembles Definition \lb \ref{dfn:3465}.

\begin{dfn} \label{3505} 
Let $A$ and $B$ be DG rings, and let $T$ be an object of 
$\dcat{D}(B \ot_{} A^{\mrm{op}})$.
\index{Differential graded module! generator}
\index{Differential graded module! compact}
\begin{enumerate}
\item We say that $T$ is a {\em compact generator of $\dcat{D}(B)$},
or a  {\em compact generator on the $B$ side}, 
if $\opn{Rest}_B(T) \in \dcat{D}(B)$ is a compact generator. 

\item We say that $T$ is a {\em compact generator of 
$\dcat{D}(A^{\mrm{op}})$},
or a  {\em compact generator on the $A^{\mrm{op}}$ side}, if 
$\opn{Rest}_{A^{\mrm{op}}}(T) \in \dcat{D}(A^{\mrm{op}})$ is a compact 
generator. 

\item We say that $T$ is a {\em compact generator on both sides}
if it is a compact generator on the $B$ side and on the $A^{\mrm{op}}$ side.
\end{enumerate}
\end{dfn}

The next theorem is similar to results appearing in \cite{Ric2} and 
\cite{Kel2}. It is a derived version of the classical result for invertible 
bimodules over a ring (see Example \ref{exa:3435} or \cite[Section 4.1]{Row}). 

In Definition \ref{dfn:3430} we introduced the notion of pretilting DG 
bimodules, and in Definition \ref{dfn:3455} we introduced tilting DG bimodules. 

\begin{thm} \label{thm:3455}
Let $A$ and $B$ be K-flat DG central $\K$-rings. The following three 
conditions are equivalent for an object 
$T \in \dcat{D}(B \ot_{} A^{\mrm{op}})$. 
\begin{enumerate}
\rmitem{i} The DG $B$-$A$-bimodule $T$ is tilting.
\index{Differential graded bimodule! tilting}

\rmitem{ii} The DG $B$-$A$-bimodule $T$ is pretilting.
\index{Differential graded bimodule! pretilting}

\rmitem{iii} The DG $B$-$A$-bimodule $T$ is a compact generator on the 
$B$ side, and it has the noncommutative derived Morita property on the $B$ side.
\index{Differential graded module! compact}
\index{Differential graded module! generator}
\index{Derived Morita property! noncommutative}
\end{enumerate} 
\end{thm}

\begin{proof} \mbox{} 

\smallskip \noindent
(i) $\Rightarrow$ (ii): Let $S \in \dcat{D}(A \ot_{} B^{\mrm{op}})$ be the 
quasi-inverse of $T$. The functor 
\begin{equation} \label{eqn:3505}
G_T := T \ot^{\mrm{L}}_{A} (-) : \dcat{D}(A) \to \dcat{D}(B)
\end{equation}
has a quasi-inverse 
\begin{equation} \label{eqn:3513}
G_S := S \ot^{\mrm{L}}_{B} (-) : \dcat{D}(B) \to \dcat{D}(A) . 
\end{equation}
So $G_T$ is an equivalence, and the DG bimodule $T$ is pretilting. 

\medskip \noindent
(ii) $\Rightarrow$ (iii): By definition, the functor $G_T$ 
from formula (\ref{eqn:3505}) is an equivalence of triangulated categories. 
Also $G_T(A) \cong T$ in $\dcat{D}(B)$. 
Because $A$ is a compact generator of $\dcat{D}(A)$, Proposition 
\ref{prop:3545} says that $T$ is a compact generator of $\dcat{D}(B)$.

Next we shall prove that $T$ has the noncommutative derived Morita property 
on the $B$ side, namely that the derived homothety morphism 
\begin{equation} \label{eqn:3511}
\opn{hm}^{\mrm{R}}_{M, A^{\mrm{op}}} : A \to \opn{RHom}_B(T, T) 
\end{equation}
in $\dcat{D}(A^{\mrm{en}})$ is an isomorphism. 
According to Proposition \ref{prop:3445} the restriction functor
$\opn{Rest}_{A} : \dcat{D}(A^{\mrm{en}}) \to \dcat{D}(A)$ 
is conservative. Therefore, to prove that the morphism 
$\opn{hm}^{\mrm{R}}_{M, A^{\mrm{op}}}$ 
is an isomorphism in $\dcat{D}(A^{\mrm{en}})$, we can forget the 
$A^{\mrm{op}}$-module structures on the DG 
bimodules $A$ and $\opn{RHom}_B(T, T)$, and just prove that 
$\opn{Rest}_{A}(\opn{hm}^{\mrm{R}}_{M, A^{\mrm{op}}})$
is an isomorphism in $\dcat{D}(A)$. 

According to Proposition \ref{prop:3537} the functor 
\begin{equation} \label{eqn:3506}
F_T := \opn{RHom}_B(T, -) : \dcat{D}(B) \to \dcat{D}(A) , 
\end{equation}
which is the right adjoint of $G_T$, is an equivalence.
So there is an isomorphism of triangulated functors
$\ze : \opn{Id}_{\dcat{D}(A)} \iso F_T \circ G_T$
from $\dcat{D}(A)$ to itself. 
There is a diagram 
\[ \UseTips \xymatrix @C=6ex @R=6ex {
A
\ar[d]_{\opn{Rest}_{A}(\opn{hm}^{\mrm{R}}_{M, A^{\mrm{op}}}) }
\ar@(r,ul)[drr]^{\ze_A}_{\cong}
\\
\opn{RHom}_B(T, T)
\ar[r]^(0.45){\cong}
&
\opn{RHom}_B(T, T \ot^{\mrm{L}}_{A} A)
\ar[r]^(0.56){\cong}
&
(F_T \circ G_T)(A)
} \]
in $\dcat{D}(A)$, and a calculation (with elements, using a K-flat resolution 
$P \to T$ in $\dcat{C}_{\mrm{str}}(B \ot A^{\mrm{op}})$) shows that it is 
commutative. Because $\ze_A$ is an isomorphism, we conclude that 
$\opn{Rest}_{A}(\opn{hm}^{\mrm{R}}_{M, A^{\mrm{op}}})$ is an isomorphism. 

\medskip \noindent
(iii) $\Rightarrow$ (ii): 
Choose a resolution $T \to I$ in $\dcat{C}(B \ot_{} A^{\mrm{op}})$
such that $I$ is K-injective over $B$; this is possible by Lemma 
\ref{lem:3445}(2). (If $A$ is K-projective over $\K$ then we can also 
choose a resolution $P \to T$ in $\dcat{C}(B \ot_{} A^{\mrm{op}})$
such that $P$ is K-projective over $B$; and then the proof can proceed with $P$ 
instead of $I$.) Then the noncommutative derived homothety morphism 
$\opn{hm}^{\mrm{R}}_{T, A^{\mrm{op}}}$ in $\dcat{D}(A^{\mrm{en}})$ is 
represented by the canonical DG 
ring homomorphism $g : A \to A_I$, where 
$A_I := \opn{End}_B(I)^{\mrm{op}}$.
The noncommutative derived Morita property on the $B$ side says that 
$g$ is a quasi-isomorphism. Hence, by Theorem \ref{thm:2363}, 
the restriction functor 
$\opn{Rest}_{g} : \dcat{D}(A_I) \to \dcat{D}(A)$
is an equivalence. 

Because $I \cong T$ is a compact generator of $\dcat{D}(B)$, according to 
Theorem \ref{thm:3430} we know that the functor 
$F_I := \opn{RHom}_B(I, -) : \dcat{D}(B) \to \dcat{D}(A_I)$
is an equivalence. Taking $F_T$ to be the functor from (\ref{eqn:3506}),
we have a diagram of functors 
\[ \UseTips \xymatrix @C=8ex @R=6ex {
\dcat{D}(B)
\ar[d]_{F_I}
\ar[dr]^{F_T}
\\
\dcat{D}(A_I)
\ar[r]^{\opn{Rest}_{g}}
&
\dcat{D}(A)
} \]
that is commutative up to isomorphism (cf.\ Theorem \ref{thm:2363}(3)).
Therefore $F_T$ is an equivalence, and thus $T$ is pretilting. 

\medskip \noindent
(ii) $\Rightarrow$ (i): 
We assume that $T$ is a pretilting DG $B$-$A$-bimodule, i.e.\ 
the functors $G_T$ and $F_T$, from (\ref{eqn:3505}) and
(\ref{eqn:3506}) respectively, are equivalences.
Define the DG bimodule 
\begin{equation} \label{eqn:3507}
S := \opn{RHom}_B(T, B) \in \dcat{D}(A \ot_{} B^{\mrm{op}}) .
\end{equation}
Consider the tensor-evaluation morphism 
\begin{equation} \label{eqn:3508}
\opn{ev}^{\mrm{R, L}}_{T, B, T} : \opn{RHom}_B(T, B) 
\ot^{\mrm{L}}_{B} T \to 
\opn{RHom}_B(T, B \ot^{\mrm{L}}_{B} T)
\end{equation}
in $\dcat{D}(A^{\mrm{en}})$ from  Proposition \ref{prop:3565}.
Because the restriction functor
$\opn{Rest}_{\K} : \lb \dcat{D}(A^{\mrm{en}}) \to \dcat{D}(\K)$
is conservative, to prove that $\opn{ev}^{\mrm{R, L}}_{T, B, T}$ is an 
isomorphism in $\dcat{D}(A^{\mrm{en}})$, it suffices 
to prove it for 
$\opn{Rest}_{\K}(\opn{ev}^{\mrm{R, L}}_{T, B, T})$. 
But by Proposition \ref{prop:3431} the DG module $T$ is algebraically perfect 
on the $B$ side; so by Theorem \ref{thm:3400}, the morphism 
$\opn{Rest}_{\K}(\opn{ev}^{\mrm{R, L}}_{T, B, T})$ is an isomorphism.
Thus we get these isomorphisms 
\begin{equation} \label{eqn:3509}
\begin{aligned}
& S \ot^{\mrm{L}}_{B} T = \opn{RHom}_B(T, B) \ot^{\mrm{L}}_{B} T
\\
& \quad \xar{\opn{ev}^{\mrm{R, L}}_{T, B, T}}  
\opn{RHom}_B(T, B \ot^{\mrm{L}}_{B} T)
\cong \opn{RHom}_B(T, T) 
\end{aligned}
\end{equation}
in $\dcat{D}(A^{\mrm{en}})$. 
On the other hand, in the proof of ``(ii) $\Rightarrow$ (iii)'' above we 
already showed (formula (\ref{eqn:3511})) that 
$A \cong \opn{RHom}_B(T, T)$ in $\dcat{D}(A^{\mrm{en}})$. 
Combining this with (\ref{eqn:3509}) we deduce that 
\begin{equation} \label{eqn:3512}
S \ot^{\mrm{L}}_{B} T \cong A \ \ \tup{in} \ \ \dcat{D}(A^{\mrm{en}}) . 
\end{equation}

Next, consider the functor $G_S$ from formula (\ref{eqn:3513}), where now $S$ 
is the DG bimodule from (\ref{eqn:3507}). The same arguments as above -- those 
used for $\opn{ev}^{\mrm{R, L}}_{T, B, T}$ -- show that for every 
$N \in \dcat{D}(B)$ the morphism 
\[ \opn{ev}^{\mrm{R, L}}_{T, B, N} : \opn{RHom}_B(T, B) \ot^{\mrm{L}}_{B} N \to 
\opn{RHom}_B(T, B \ot^{\mrm{L}}_{B} N) \]
in $\dcat{D}(B)$ is an isomorphism. So, like (\ref{eqn:3509}), we get 
isomorphisms 
\begin{equation} \label{eqn:3514}
\begin{aligned}
& G_S(N) = S \ot^{\mrm{L}}_{B} N = \opn{RHom}_B(T, B) \ot^{\mrm{L}}_{B} N
\\
& \quad \xar{\opn{ev}^{\mrm{R, L}}_{T, B, N}}  
\opn{RHom}_B(T, B \ot^{\mrm{L}}_{B} N)
\cong \opn{RHom}_B(T, N) = F_T(N)
\end{aligned}
\end{equation}
in $\dcat{D}(A)$. Since these isomorphisms are functorial in $N$, we see that 
$G_S \cong F_T$ as functors. Therefore $G_S$ is an equivalence, and 
$S$ is a pretilting DG $A$-$B$-bimodule.

The same proof, but now for the pretilting DG $A$-$B$-bimodule $S$ instead of 
for $T$, shows that the DG $B$-$A$-bimodule 
\begin{equation} \label{eqn:3527}
T' := \opn{RHom}_A(S, A) \in \dcat{D}(B \ot_{} A^{\mrm{op}}) 
\end{equation}
satisfies the corresponding version of (\ref{eqn:3512}), i.e.\ 
$T' \ot^{\mrm{L}}_{A} S \cong B$
in $\dcat{D}(B^{\mrm{en}})$. 
By Lemma \ref{lem:3455} the DG $B$-$A$-bimodule $T \cong T'$ is tilting.
\end{proof}

\begin{cor} \label{cor:3455}
Let $T$ be a tilting DG $B$-$A$-bimodule. Then its quasi-inverse is the DG 
$A$-$B$-bimodule $S := \opn{RHom}_B(T, B)$.
\end{cor}

\begin{proof}
This was shown in the proof of the implication ``(ii) $\Rightarrow$ (i)''
above.
\end{proof}

\begin{cor} \label{cor:3505}
Let $T$ be a tilting DG $B$-$A$-bimodule. Then $T$ is a compact generator on 
both sides, and it has the noncommutative derived Morita property on both 
sides.
\end{cor}

\begin{proof}
By Theorem \ref{thm:3455}, $T$ is a compact generator on the $B$ side, 
and it has the noncommutative derived Morita property on the $B$ side. 

Using the DG ring isomorphism 
$A^{\mrm{op}} \ot B \iso B \ot A^{\mrm{op}}$
of (\ref{eqn:3552}), we can also view $T$ as an object of 
$\dcat{D}(A^{\mrm{op}} \ot B)$. And as such, $T$ is a tilting DG 
$A^{\mrm{op}}$-$B^{\mrm{op}}$-bimodule. Now Theorem \ref{thm:3455} tells us 
that $T$ is a compact generator on the $A^{\mrm{op}}$ side, 
and it has the noncommutative derived Morita property on the 
$A^{\mrm{op}}$ side. 
\end{proof}

The example of classical Morita Theory was already given; see Example 
\ref{exa:3435}. Next is an example that goes in another direction altogether.

\begin{exa} \label{eqn:3565}
Suppose $f : A \to B$ is a quasi-isomorphism in $\catt{DGRng} \fcentover \K$.
Then 
$T := B \in \dcat{D}(B \ot A^{\mrm{op}})$
is a tilting DG $B$-$A$-bimodule. 
The equivalence 
\[ F_T := \opn{RHom}_B(T, -) : \dcat{D}(B) \to \dcat{D}(A) \]
is just the restriction functor $\opn{Rest}_f$; and the equivalence 
\[ G_T := T \ot^{\mrm{L}}_{A} (-) :
\dcat{D}(A) \to \dcat{D}(B) \]
is just the derived induction functor $\opn{LInd}_f$. 
\end{exa}

\begin{thm}[Rickard-Keller] \label{thm:3565} 
Let $A$ and $B$ be K-flat DG central $\K$-rings. Assume that there exists a 
$\K$-linear equivalence of triangulated categories 
$F : \lb \dcat{D}(A) \to \dcat{D}(B)$,
and that $\opn{H}^i(A) = 0$ for all $i \neq 0$. 
Then there exists a tilting DG $B$-$A$-bimodule $T$. 
\end{thm}

This theorem is very similar to J. Rickard's \cite[Theorem 6.4]{Ric2}
and to B. Keller's \cite[Corollary 9.2]{Kel1}. See Remark \ref{rem:4189} for a 
brief discussion.

\begin{proof}
We begin with a warning: in the proof we will construct some new DG rings,
and they might fail to be K-flat over $\K$, thus violating Convention 
\ref{conv:3445}.

Consider the object $L := F(A) \in \dcat{D}(B)$.
It is a compact generator, by Propositions \ref{prop:3545} and \ref{prop:3547}.
Let us choose a K-projective resolution $P \to L$ in 
$\dcat{C}_{\mrm{str}}(B)$, and define the DG ring
$\til{A} := \opn{End}_B(P)^{\mrm{op}}$.
So $P \in \dcat{D}(B \ot_{\K} \til{A}^{\mrm{op}})$.

Because $F$ is an equivalence, we get ring isomorphisms 
\begin{equation} \label{eqn:3570}
\opn{H}^0(A) \cong \opn{End}_{\dcat{D}(A)}(A)^{\mrm{op}} \xar{F} 
\opn{End}_{\dcat{D}(B)}(L)^{\mrm{op}} \cong 
\opn{End}_{\dcat{D}(B)}(P)^{\mrm{op}} \cong \opn{H}^0(\til{A}) .
\end{equation}
For the same reason we have 
$\opn{H}^i(\til{A}) \cong \opn{H}^i(A)  = 0$ for all $i \neq 0$. 
Define the DG rings $A_1 := \opn{smt}^{\leq 0}(A)$ and 
$\til{A}_1 := \opn{smt}^{\leq 0}(\til{A})$,
and the rings $A_2 := \opn{H}^{0}(A)$ and
$\til{A}_2 := \opn{H}^{0}(\til{A})$. 
There are canonical quasi-isomorphisms 
$A_1 \to A$, $A_1 \to A_2$, 
$\til{A}_1 \to \til{A}$ and $\til{A}_1 \to \til{A}_2$ 
in $\catt{DGRng} \centover \K$. Equation (\ref{eqn:3570}) produces an 
isomorphism $A_2 \cong \til{A}_2$ in $\catt{Rng} \centover \K$.
Now the DG ring $B$ is K-flat over $\K$, so we have induced quasi-isomorphisms
\[ B \ot A^{\mrm{op}} \lar B \ot A_1^{\mrm{op}} \to B \ot A_2^{\mrm{op}}
\cong B \ot \til{A}_2^{\mrm{op}} \lar 
B \ot \til{A}_1^{\mrm{op}} \to B \ot \til{A}^{\mrm{op}} \]
in $\catt{DGRng} \centover \K$. According to Theorem \ref{thm:2363} we get 
$\K$-linear equivalences of  triangulated categories 
\[ \begin{aligned}
& \dcat{D}(B \ot A^{\mrm{op}}) \to \dcat{D}(B \ot A_1^{\mrm{op}}) \to 
\dcat{D}(B \ot A_2^{\mrm{op}}) 
\\
& \quad \to \dcat{D}(B \ot \til{A}_2^{\mrm{op}}) \to 
\dcat{D}(B \ot \til{A}_1^{\mrm{op}}) \to \dcat{D}(B \ot \til{A}^{\mrm{op}}) .
\end{aligned} \]
There is an object $T \in \dcat{D}(B \ot A^{\mrm{op}})$ that 
corresponds to $P \in \dcat{D}(B \ot \til{A}^{\mrm{op}})$
under this chain of equivalences. All these equivalence restrict to the 
identity functor of $\dcat{D}(B)$, and hence there is an isomorphism 
\begin{equation} \label{eqn:4181}
\opn{Rest}_B(T) \cong \opn{Rest}_B(P) \cong L 
\end{equation}
in $\dcat{D}(B)$. 

Since $L$ is a compact generator on the $B$ side, the same is true for $T$.
Finally, by equation (\ref{eqn:4181}), and using the equivalence $F$,  
we know that 
\begin{equation} \label{eqn:4182}
\opn{H}^i \bigl( \opn{RHom}_B(T, T) \bigr) \cong 
\opn{H}^i \bigl( \opn{RHom}_B(L, L) \bigr) \cong
\opn{H}^i \bigl( \opn{RHom}_A(A, A) \bigr) = 0
\end{equation}
for all $i \neq 0$. By equation (\ref{eqn:3570}) we know that 
\begin{equation} \label{eqn:4183}
\opn{H}^0 \bigl( \opn{RHom}_B(T, T) \bigr) \cong 
\opn{End}_{\dcat{D}(B)}(L) \cong \opn{End}_{\dcat{D}(A)}(A) \cong 
\opn{H}^0(A) 
\end{equation}
as $\opn{H}^0(A)$-bimodules. 
Once we trace the various morphisms, formulas (\ref{eqn:4182}) and 
(\ref{eqn:4183}) imply that the derived homothety morphism 
$\opn{hm}^{\mrm{R}}_{T, A^{\mrm{op}}} : A \to \lb \opn{RHom}_B(T, T)$
in $\cat{D}(A^{\mrm{en}})$ is an isomorphism. So $T$ has the derived Morita 
property on the $B$ side. 
By Theorem \ref{thm:3455}, $T$ is a tilting DG $B$-$A$-bimodule.
\end{proof}

\begin{rem} \label{rem:4189}
J. Rickard, in \cite[Theorem 6.4]{Ric1}, considers rings $A$ and $B$, and he 
does not assume that they are flat over a base ring $\K$. So instead of 
producing a tilting complex $T$, as we do in Theorem \ref{thm:3565} above, he 
actually produces a pretilting complex. Another difference is that Rickard 
considers an equivalence 
$F : \dcat{D}^{\mrm{b}}(A) \to \dcat{D}^{\mrm{b}}(B)$,
whereas we look at the unbounded derived categories. 

B. Keller, in \cite[Corollary 9.2]{Kel1}, works in much greater generality than 
we do: instead of DG rings $A$ and $B$, he looks at DG categories $\cat{A}$ and 
$\cat{B}$. Also he does not require $\cat{A}$ to be K-flat over $\K$. 

The apparent flatness limitation in our Theorem \ref{thm:3565} can be 
effectively overcome using K-flat resolutions 
$\til{A} \to A$ and $\til{B} \to B$ in $\catt{DGRng} \centover \K$,
as explained in Remark \ref{rem:1420}. Since the restriction functors 
$\dcat{D}(A) \to \dcat{D}(\til{A})$ and 
$\dcat{D}(B) \to \dcat{D}(\til{B})$ are equivalences,
we get an equivalence 
$\til{F} : \dcat{D}(\til{A}) \to \dcat{D}(\til{B})$.
Theorem \ref{thm:3565} produces a tilting DG bimodule 
$T \in \dcat{D}(\til{B} \ot \til{A}^{\mrm{op}})$. 
\end{rem}

\begin{rem} \label{rem:3565}
In the situation of Theorem \ref{thm:3565}, it is not known whether one can 
find a tilting DG $B$-$A$-bimodule $T$ such that 
$F \cong T \ot^{\mrm{L}}_{A} (-)$ as triangulated functors. 
This question was already raised by Rickard in \cite{Ric1}, and still remains 
open, except for a few special cases (see \cite{MiYe} for hereditary rings, and 
\cite{Ch} for triangular rings). 

If we drop the condition that the cohomology of $A$ is concentrated in  
degree $0$, then there is a counterexample to the assertion of Theorem 
\ref{thm:3565}, due to B. Shipley \cite[Section 5]{Shi}, based on her joint 
work with D. Dugger \cite{DuSh}. Here is its translation to our 
terminology. She considers the DG ring  
$A = \Z[\bar{x}] := \Z \bra{x} / (x^4)$,
where $x$ is a variable of degree $-1$, and $\d(\bar{x}) := 2$. 
(Note that $A$ is weakly commutative but not strongly commutative, since 
$\bar{x}^2 \neq 0$.) 
The second DG ring is $\opn{H}(A)$. A calculation shows that 
$\opn{H}(A) \cong  \mbb{F}_2[\bar{y}] := \mbb{F}_2 \bra{y} / (y^2)$,
where $y$ is a variable of degree $-2$, and $\d(\bar{y}) := 0$.
For the base ring we take $\K := \Z$. Because $\opn{H}(A)$ is not K-flat over 
$\K$, we choose a K-flat resolution $B \to \opn{H}(A)$ over $\K$. 

Shipley asserts that on the one hand there is an equivalence of triangulated 
categories $\dcat{D}(A) \to \dcat{D}(B)$, and on the other hand there does not 
exist a DG bimodule $T \in \dcat{D}(B \ot A^{\mrm{op}})$ 
that is a compact generator on the $B$ side, 
and has the noncommutative derived Morita property on the $B$ side
(condition (iii) of Theorem \ref{thm:3455}). Thus there does not exist a 
tilting DG $B$-$A$-bimodule. 
\end{rem}

\begin{rem} \label{rem:3596}
In algebraic geometry the role of tensoring with tilting bimodule complexes is 
played by {\em Fourier-Mukai transforms}. This theory is at the core of 
contemporary birational geometry. See the survey \cite{HiVdB} or the book 
\cite{Huy}. 
\end{rem}

\mysubsection{Tilting Bimodule Complexes over Rings} 
\label{subsec:tilting-rings}

We continue with Convention \ref{conv:3445}, but in this subsection 
we only look at rings. So our rings are flat central over the base ring $\K$, 
$\ot = \ot_{\K}$, and $A^{\mrm{en}} = A \ot A^{\mrm{op}}$ for a ring $A$. 
All bimodules are $\K$-central, and all additive functors are $\K$-linear.  
Rather than speaking about DG bimodules, here we speak about complexes of 
bimodules. By default our rings are nonzero. 

In classical Morita theory, an 
{\em invertible  $B$-$A$-bimodule}%
\index{Invertible! bimodule}
is a bimodule $P$ for which there exists an $A$-$B$-bimodule $Q$, such that  
$Q \ot_B P \cong A$ in $\dcat{M}(A^{\mrm{en}})$ 
and
$P \ot_A Q \cong B$ in $\dcat{M}(B^{\mrm{en}})$.
It is known that $P$ is a finitely generated projective module over $B$ and 
over $A^{\mrm{op}}$; and the reverse is true for $Q$. Furthermore, every 
$\K$-linear equivalence 
$F : \dcat{M}(A) \to \dcat{M}(B)$
is isomorphic, as a functor, to the functor $P \ot_A (-)$ for an invertible 
bimodule $P$, and this $P$ is unique up to a unique isomorphism. See 
\cite[Section 4.1]{Row} for proofs. 

We see that invertible bimodules are a very special kind of tilting bimodule 
complexes. Proposition \ref{prop:3575} clarifies matters. 
But first a lemma, which is a noncommutative variant of Lemma \ref{lem:2177}

\begin{lem}[K\"unneth Trick] \label{lem:3595} 
Let $M \in \dcat{D}^-(B \ot A^{\mrm{op}})$ and 
$N \in \dcat{D}^-(A \ot C^{\mrm{op}})$,  
and let $i_1, j_1 \in \Z$ be such that 
$\opn{sup}(\opn{H}(M)) \leq i_1$ and $\opn{sup}(\opn{H}(N)) \leq j_1$. 
Then 
\[ \opn{H}^{i_1 + j_1}(M \ot^{\mrm{L}}_A N) \cong 
\opn{H}^{i_1}(M) \ot_A \opn{H}^{j_1}(N)  \]
in $\dcat{M}(B \ot C^{\mrm{op}})$. 
\end{lem}

\begin{exer} \label{exer:3595}
Prove Lemma \ref{lem:3595}.
\end{exer}

\begin{prop} \label{prop:3575}
The following conditions are equivalent for a tilting bimodule complex 
$T \in \dcat{D}(B \ot A^{\mrm{op}})$, with quasi-inverse $S = T^{\vee}$. 
\begin{enumerate}
\rmitem{i} There is an isomorphism $T \cong P$ in 
$\dcat{D}(B \ot A^{\mrm{op}})$
for some invertible $B$-$A$-bimodule $P$. 

\rmitem{ii} $\opn{H}^0(T)$ is a projective $B$-module, and 
$\opn{H}^i(T) = 0$ for all $i \neq 0$. 

\rmitem{iii} $\opn{H}^i(T) = 0$ and $\opn{H}^i(S) = 0$ for all $i \neq 0$.  
\end{enumerate} 
\end{prop}

\begin{proof} \mbox{} 

\smallskip \noindent 
(i) $\Rightarrow$ (ii): this is trivial.

\medskip \noindent
(ii) $\Rightarrow$ (iii): 
Let $P := \opn{H}^0(T) \in \dcat{M}(B \ot A^{\mrm{op}})$,
so there is an isomorphism 
$T \cong P$ in $\dcat{D}(B \ot A^{\mrm{op}})$.
According to Corollary \ref{cor:3455} we have 
$S \cong \opn{RHom}_B(T, B) \cong \opn{Hom}_B(P, B)$.

\medskip \noindent
(iii) $\Rightarrow$ (i): 
Define $P := \opn{H}^{0}(T) \in \dcat{M}(B \ot A^{\mrm{op}})$
and $Q := \opn{H}^{0}(S) \in \dcat{M}(A \ot B^{\mrm{op}})$.
The K\"unneth Trick (Lemma \ref{lem:3595}) says that 
$P \ot_B Q \cong \opn{H}^{0}(T \ot^{\mrm{L}}_{B} S) \cong \opn{H}^{0}(A) = A$
and
$Q \ot_A P \cong \opn{H}^{0}(S \ot^{\mrm{L}}_{A} T) \cong \opn{H}^{0}(B) = B$
in $\dcat{M}(A^{\mrm{en}})$ and $\dcat{M}(B^{\mrm{en}})$ respectively. 
Thus $P$ and $Q$ are invertible bimodules. 
But $T \cong P$ in $\dcat{D}(B \ot A^{\mrm{op}})$. 
\end{proof}

The {\em Jacobson radical} of a ring $A$ is the intersection of all 
maximal left (or right) ideals of $A$. See \cite[Section 4.4.2]{Jac} or 
\cite[Section 2.5]{Row}. 

\begin{dfn}
A ring $A$, with Jacobson radical $\r$, is called {\em local} 
\index{Ring! local NC}
if $A / \r$ is a simple artinian ring. 
\end{dfn}

Note that this definition is wider than what is found in some books (e.g.\ 
\cite{Row}), where the condition is that $A / \r$ is a division ring. 

\begin{exa} \label{exa:3575}
If $C$ is a commutative local ring with maximal idea $\m$, and $n \geq 1$, 
then the matrix ring $A := \opn{Mat}_{n \times n}(C)$ is local, with Jacobson 
radical $\m \cd A = \opn{Mat}_{n \times n}(\m)$.  
\end{exa}

\begin{lem} \label{lem:3575}
Let $A$ be a local ring, let $M$ be a nonzero finitely generated right 
$A$-module, and let $N$ be a nonzero finitely generated left $A$-module. Then 
$M \ot_A N$ is nonzero. 
\end{lem}

\begin{proof}
Define $K := A / \r$, which is a simple artinian ring. By the noncommutative 
Nakayama Lemma (see \cite[Proposition 2.5.24]{Row}) the right $K$-module 
$M \ot_A K$ is nonzero, and also the left $K$-module $K \ot_A N$ is nonzero.
There is a canonical surjection 
\[ M \ot_A N \to (M \ot_A K) \ot_K (K \ot_A N) . \]
It suffices to prove that the target is nonzero. Thus we can assume that 
$A = K$ is a simple artinian ring. 

Every left module over $K$ is a  direct sum of simple ones; so we can assume 
that $N$ is simple. Likewise we can assume that $M$ is a simple right 
$K$-module. 

There is a ring isomorphism $K \cong \opn{Mat}_{n \times n}(D)$, 
where $D$ is a division ring and $n$ is a positive integer. The simple left 
$K$-module $N$ is isomorphic to $D^n$, seen as a column module; thus it is in 
fact a $K$-$D$-bimodule. Similarly $M \cong D^n$, seen as a row module; thus it 
is a $D$-$K$-bimodule. By an easy calculation (this is an elementary case of 
Morita equivalence) there is an isomorphism 
$M \ot_K N \cong D$ of $D$-$D$-bimodules. This is nonzero. 
\end{proof}

\begin{thm}[\cite{RoZi}, \cite{Ye4}] \label{thm:3575}
Let $A$ and $B$ be flat central $\K$-rings, with $A$ local, and let 
$T$ be a tilting $B$-$A$-bimodule complex%
\index{Tilting! complex of bimodules}%
\index{Invertible! bimodule}.
Then there is an iso\-morph\-ism $T \cong P[n]$
in $\dcat{D}(B \ot A^{\mrm{op}})$ for some 
invertible $B$-$A$-bimodule
$P$ and some integer $n$. 
\end{thm}

See Remark \ref{rem:3430} regarding the history of this theorem. 

\begin{proof}
Let $S \in \dcat{D}(A \ot B^{\mrm{op}})$
be the quasi-inverse of $T$. From Corollary \ref{cor:3505} and Theorem 
\ref{thm:3400} we know that $T$ and $S$ are algebraically perfect complexes on 
both sides. By Lemma \ref{lem:3415} we know that 
$\opn{con}(\opn{H}(T)) = [i_0, i_1]$ and 
$\opn{con}(\opn{H}(S)) = [j_0, j_1]$ for some integers $i_0 \leq i_1$ and 
$j_0 \leq j_1$; and also that $P := \opn{H}^{i_1}(T)$ and 
$Q := \opn{H}^{j_1}(S)$ are finitely presented modules on both 
sides. 

We now use the K\"unneth trick to obtain an isomorphism
\[ P \ot_A Q \cong 
\opn{H}^{i_1 + j_1}(T \ot^{\mrm{L}}_{A} S) \cong \opn{H}^{i_1 + j_1}(B) \]
in $\cat{M}(B^{\mrm{en}})$. By Lemma \ref{lem:3575} we know that 
$P \ot_A Q$ is nonzero. Therefore $i_1 + j_1 = 0$, and 
$P \ot_A Q \cong B$ in $\dcat{M}(B^{\mrm{en}})$. 
For similar reasons, on the reverse side we have 
$Q \ot_B P \cong \opn{H}^{0}(S \ot^{\mrm{L}}_{B} T) \cong A$
in $\dcat{M}(A^{\mrm{en}})$. Thus $P$ and $Q$ are invertible bimodules. 

We now restrict $T$ to $\dcat{D}(A^{\mrm{op}})$. Because its top cohomology 
module $\opn{H}^{i_1}(T) = P$ is projective, we can split it off (see Lemma 
\ref{lem:2189}, that holds also for a noncommutative ring), to get an 
isomorphism $T \cong T' \oplus P[-i_1]$ in $\dcat{D}(A^{\mrm{op}})$, 
with $\opn{sup}(\opn{H}(T')) \leq i_1 - 1$.
Similarly there's an isomorphism 
$S \cong S' \oplus Q[-j_1]$ in $\dcat{D}(A)$
with $\opn{sup}(\opn{H}(S')) \leq j_1 - 1$.
It follows that 
\[ \begin{aligned}
& B \cong T \ot^{\mrm{L}}_{A} S \cong 
\bigl( T' \oplus P[-i_1] \bigr) \ot^{\mrm{L}}_{A} 
\bigl( S' \oplus Q[-j_1] \bigr) 
\\
& \quad 
\cong (T' \ot^{\mrm{L}}_{A} S') \oplus
\bigl( P[-i_1] \ot_A S' \bigr) \oplus \bigl( T' \ot_A Q[-j_1] \bigr) 
\oplus B 
\end{aligned} \]
in $\dcat{D}(\K)$. The object $B \in \dcat{D}(\K)$ has cohomology concentrated 
in degree $0$. The direct summand $T' \ot_A Q[-j_1]$ has cohomology 
concentrated in degrees $\leq -1$, so its cohomology must be zero, and 
therefore $T' \ot_A Q[-j_1] = 0$ in $\dcat{D}(\K)$. 
But $Q$ is an invertible bimodule, and this forces $T' = 0$
in $\dcat{D}(A^{\mrm{op}})$. The conclusion is that $\opn{H}^{i}(T) = 0$ for 
all $i \neq i_1$.

Returning to the category $\dcat{D}(A \ot B^{\mrm{op}})$
we see that $T \cong P[-i_1]$. Finally we take $n := -i_1$. 
\end{proof}

\begin{cor} \label{cor:3580}
Let $A$ and $B$ be rings, with $A$ local.
If there is an equivalence of triangulated categories
$\dcat{D}(A) \to \dcat{D}(B)$, then there is an equivalence of 
abelian categories $\dcat{M}(A) \to \dcat{M}(B)$. In other words, $A$ and $B$ 
are Morita equivalent. 
\end{cor}

\begin{proof}
According to Theorem \ref{thm:3565} there exists a tilting $B$-$A$-bimodule 
complex $T$. Theorem \ref{thm:3575} says that $T \cong P[n]$ for an invertible 
$B$-$A$-bimodule $P$. Then 
$P \ot_A (-) : \dcat{M}(A) \to \dcat{M}(B)$
is an equivalence of abelian categories. 
\end{proof}

\begin{dfn} \label{dfn:3590}
For a central $\K$-ring $A$ we define the 
{\em noncommutative Picard group of $A$ relative to $\K$}
\index{Picard group! noncommutative}
to be the group $\opn{Pic}_{\K}(A)$, whose elements are 
the isomorphism classes in $\dcat{M}(A^{\mrm{en}})$ of the invertible 
bimodules. The operation is induced by $(- \ot_A -)$, and the unit element is 
the class of $A$. 
\end{dfn}

The derived Picard group was introduced in Definition \ref{dfn:3457}.

\begin{cor} \label{cor:3590}
If $A$ is a local ring, then 
\[ \opn{DPic}_{\K}(A) = \opn{Pic}_{\K}(A) \times \Z . \]
\end{cor}

\begin{proof}
By Theorem \ref{thm:3575} there is a surjective group homomorphism 
$\opn{Pic}_{\K}(A) \times \Z \to  \opn{DPic}_{\K}(A)$,
defined by $(P, n) \mapsto P[n]$.
It is injective because the functor 
$\dcat{M}(A^{\mrm{en}}) \to \dcat{D}(A^{\mrm{en}})$
is fully faithful. 
\end{proof}

The center of a ring $A$ is denoted by $\opn{Cent}(A)$. Of course 
$\opn{Cent}(A) = \opn{Cent}(A^{\mrm{op}})$. 
Given a complex $M \in  \dcat{D}(B \ot A^{\mrm{op}})$, 
there are ring homomorphisms 
\begin{equation} \label{eqn:3580}
\opn{chm}^{\dcat{D}}_{M, A^{\mrm{op}}} : \opn{Cent}(A) \to 
\opn{End}_{\dcat{D}(B \ot A^{\mrm{op}})}(M)
\end{equation}
and 
\begin{equation} \label{eqn:3581}
\opn{chm}^{\dcat{D}}_{M, B} : \opn{Cent}(B) \to 
\opn{End}_{\dcat{D}(B \ot A^{\mrm{op}})}(M) 
\end{equation}
that we call the {\em central homotheties} through $A^{\mrm{op}}$ and $B$
respectively. The formulas are the obvious ones: for an element 
$a \in \opn{Cent}(A)$ the action on $M$ is 
$\opn{chm}^{\dcat{D}}_{M, A^{\mrm{op}}}(a)(m) := m \cd a$ 
for $m \in M^i$. Likewise for $\opn{chm}^{\dcat{D}}_{M, B}$. 

\begin{lem} \label{lem:3580}
If $T$ is a tilting $B$-$A$-bimodule complex, then the ring homomorphisms 
$\opn{chm}^{\dcat{D}}_{M, A^{\mrm{op}}}$ and 
$\opn{chm}^{\dcat{D}}_{M, B}$ are both isomorphisms. 
\end{lem}

\begin{proof}
Let $S$ be the quasi-inverse of $T$. The functor 
\[ G := (-) \ot^{\mrm{L}}_{A} S : \dcat{D}(B \ot A^{\mrm{op}}) \to 
\dcat{D}(B \ot B^{\mrm{op}}) = \dcat{D}(B^{\mrm{en}}) \]
is an equivalence. It induces a ring isomorphism
\[ G : \opn{End}_{\dcat{D}(B \ot A^{\mrm{op}})}(T) \iso 
\opn{End}_{\dcat{D}(B^{\mrm{en}})}(B) , \]
and 
$G \circ \opn{chm}^{\dcat{D}}_{T, B} = \opn{chm}^{\dcat{D}}_{B, B}$.
But 
$\opn{End}_{\dcat{D}(B^{\mrm{en}})}(B) \cong 
\opn{End}_{\dcat{M}(B^{\mrm{en}})}(B)$,
and the ring homomorphism 
$\opn{Cent}(B) \to \opn{End}_{\dcat{M}(B^{\mrm{en}})}(B)$
is bijective. 

The proof for $\opn{chm}^{\dcat{D}}_{T, A^{\mrm{op}}}$ is similar.
\end{proof}

\begin{dfn} \label{dfn:5090}
Given a tilting $B$-$A$-bimodule complex $T$, we denote by 
$g_T : \opn{Cent}(A) \to \opn{Cent}(B)$
the ring isomorphism such that 
$\opn{chm}^{\dcat{D}}_{T, A^{\mrm{op}}}  = 
\opn{chm}^{\dcat{D}}_{T, B} \circ \, g_T$.
\end{dfn}

Even though the localization of a noncommutative ring $B$ with respect to a 
multiplicatively closed subset $Z \sub B$ is problematic in general, there is 
no difficulty at all if $Z \sub \opn{Cent}(B)$. In this case,
letting $C := \opn{Cent}(B)$, we get canonical $A$-ring isomorphisms
$B_{Z} \cong C_{Z} \ot_C B \cong B \ot_C C_Z$.

\begin{lem} \label{lem:3581}
Let $T$ be a tilting $B$-$A$-bimodule complex, and assume that the ring $A$ is 
commutative. Let $Z \sub A$ be a multiplicatively 
closed set, and define $A' := A_Z$ and $B' := B_{g_T(Z)}$. Then 
the complex 
$T' := B' \ot_B T \ot_A A'$ in 
$\dcat{D}(B' \ot A'^{\, \mrm{op}})$
is a tilting $B'$-$A'$-bimodule complex.
\end{lem}

\begin{proof}
The cohomology $\opn{H}(T)$ is a central graded $A$-bimodule, where the 
left action of $A$ on $\opn{H}(T)$ is via the isomorphism  $g_T$. (Warning: the 
complex of 
bimodules $T$ need not be central over $A$.) The flatness of $A \to A'$ and 
$B \to B'$ gives 
\[ \opn{H}(T') = \opn{H}(B' \ot_B T \ot_A A') \cong 
B' \ot_B \opn{H}(T) \ot_A A' \cong A' \ot_A \opn{H}(T) \ot_A A' . \]
Therefore 
$\opn{H}(T) \ot_A A' \to \opn{H}(T')$ 
is an isomorphism of graded modules, and this implies that 
\begin{equation} \label{eqn:3585}
T \ot_A A' \to T'
\end{equation}
is an isomorphism in 
$\dcat{D}(B \ot A'^{\, \mrm{op}})$. 

Let 
$S := \opn{RHom}_B(T, B)$ in $\dcat{D}(A \ot B^{\mrm{op}})$
be the quasi-inverse of $T$. 
In analogy to Definition \ref{dfn:5090}, there is a ring isomorphism 
$h_S : A \to \opn{Cent}(B)$
such that 
$\opn{chm}^{\dcat{D}}_{S, A}  = 
\opn{chm}^{\dcat{D}}_{S, B^{\mrm{op}}} \circ \, h_S$.
The cohomology $\opn{H}(S)$ is a central graded $A$-bimodule, 
where the right action of $A$ on $\opn{H}(S)$ is via $h_S$.

Let us consider the derived tensor product 
$S \ot^{\mrm{L}}_{B} T \in \dcat{D}(A^{\mrm{en}})$,
and the four ways that an element $a \in A$ acts on it as an 
object of $\dcat{D}(A^{\mrm{en}})$. 
The action of $a$ from the right side of $T$ the same as the action of 
$g_T(a) \in \opn{Cent}(B)$ from the left side of $T$. This is also the action 
of $g_T(a)$ from the right side of $S$, and also the action of 
$h_S^{-1}(g_T(a))$ from the left on $S$. 
But there are ring isomorphisms 
\[ \opn{End}_{\dcat{D}(A^{\mrm{en}})}(S \ot^{\mrm{L}}_{B} T) \cong 
\opn{End}_{\dcat{D}(A^{\mrm{en}})}(A) \cong A , \]
and $a$ acts on $A$ the same from both sides. The conclusion is that 
$h_S^{-1}(g_T(a)) = a$. 
Because $a$ is an arbitrary element, we see that $g_T = h_S$ as isomorphisms 
$A \iso \opn{Cent}(B)$. So we have a single ring isomorphism 
$A \cong \opn{Cent}(B)$, for which both $\opn{H}(T)$ and $\opn{H}(S)$ are 
central graded $A$-bimodules. 

Define 
$S' := A' \ot_A S \ot_B B' \in \dcat{D}(A' \ot B'^{\, \mrm{op}})$. 
The same arguments as for $T$ show that 
$A' \ot_A S \to S'$ is an isomorphism in 
$\dcat{D}(A' \ot B^{\mrm{op}})$.

We now calculate:
\[ \begin{aligned}
& S' \ot^{\mrm{L}}_{B'} T' = (A' \ot_A S \ot_B B') \ot^{\mrm{L}}_{B'}
(B' \ot_B T \ot_A A') 
\\
& \quad \cong (A' \ot_A S) \ot^{\mrm{L}}_{B} (B' \ot_B T \ot_A A') 
\cong^{\dag} (A' \ot_A S) \ot^{\mrm{L}}_{B} (T \ot_A A') 
\\
& \quad \cong A' \ot_A (S \ot^{\mrm{L}}_{B} T) \ot_A A') 
\cong A' \ot_A A \ot_A A' \cong A'
\end{aligned} \]
in $\dcat{D}(A'^{\, \mrm{en}})$. We used the associativity of 
$(- \ot^{\mrm{L}}_{(-)} -)$ several times, and also the ring isomorphism 
$A' \ot_A A' \cong A'$. The isomorphism $\cong^{\dag}$ is by formula 
(\ref{eqn:3585}). Similarly we show that 
$T' \ot^{\mrm{L}}_{A'} S' \cong B'$
in $\dcat{D}(B'^{\, \mrm{en}})$. 
\end{proof}

\begin{thm}[\cite{RoZi}, \cite{Ye4}] \label{thm:3585}
Let $A$ and $B$ be flat central $\K$-rings, and assume $A$ is commutative with 
connected spectrum. Let $T$ be a tilting 
$B$-$A$-bimodule complex%
\index{Differential graded bimodule! tilting}.
Then there is an isomorphism $T \cong P[n]$  
in $\dcat{D}(B \ot A^{\mrm{op}})$, for some 
invertible $B$-$A$-bimodule $P$ and some integer $n$. 
\end{thm}

See Remark \ref{rem:3430} regarding the history of this theorem. 

\begin{proof}
Because $A$ is commutative we have $A = A^{\mrm{op}}$.
The ring homomorphism $g_T$ makes $B$ into a central $A$-ring. 
We may assumed that $A \neq 0$, and hence $T \neq 0$.
The complex $T$ is algebraically perfect over $A$, 
so it has bounded cohomology, say with 
$\opn{sup}(\opn{H}(T)) = i_1 \in \Z$.
By Lemma \ref{lem:3415} the $A$-module $P := \opn{H}^{i_1}(T)$ is finitely 
presented.

For a prime $\p \in \opn{Spec}(A)$, with corresponding local ring
$A_{\p}$, we write $P_{\p} := P \otimes_A A_{\p}$.
Define $Y \sub \opn{Spec}(A)$ to be the support of $P$, i.e.\ 
\[ Y := \{ \p \in \opn{Spec}(A) \mid P_{\p} \neq 0 \} . \]
This is a nonempty set. 
Since $P$ is finitely generated it follows that $Y$ is a closed subset of 
$\opn{Spec}(A)$; see \cite[Proposition II.4.17]{Bou}.

Take any prime $\p \in Y$, and let $B_{\p} :=  B \otimes_A A_{\p}$.
Then, by Lemma \ref{lem:3581}, the complex
\[ T_{\p} := B_{\p} \otimes_B T \otimes_A A_{\p} \in 
\dcat{D}(B_{\p} \otimes A_{\p}^{\mrm{op}}) \]
is a tilting $B_{\p}$-$A_{\p}$-bimodule complex. Since
$\mrm{H}^{i_1}(T_{\p}) \cong P_{\p} \neq 0$,
Theorem \ref{thm:3575} implies that 
\begin{equation} \label{eqn:3587}
T_{\p} \cong P_{\p}[-i_1] \in \dcat{D}(B_{\p} \otimes A_{\p}^{\mrm{op}}) ,
\end{equation}
and that $P_{\p}$ is an invertible $B_{\p}$-$A_{\p}$-bimodule. 
In particular, $P_{\p}$ is a nonzero finitely generated projective 
$A_{\p}$-module. Thus $P_{\p}$ is a free $A_{\p}$-module, of rank 
$r_{\p} > 0$. Recall that $P$ is a finitely presented $A$-module. 
According to \cite[Section II.5.1,  Corollary]{Bou} there is an 
open neighborhood $U$ of $\p$ in $\opn{Spec}(A)$ on which $P$ is free of rank 
$r_{\p}$. In particular $P_{\q} \neq 0$ for all $\q \in U$.
Therefore $U \sub Y$. 

The conclusion is that $Y$ is also open in $\opn{Spec}(A)$.
Since $\opn{Spec}(A)$ is connected, it follows that 
$Y = \opn{Spec}(A)$. Another conclusion is that $P$ is projective as an 
$A$-module -- see \cite[Section II.5.2, Theorem 1]{Bou}.

Going back to equation (\ref{eqn:3587}) we see that 
$\mrm{H}^i(T)_{\p} \cong \mrm{H}^i(T_{\p}) = 0$
for all $i \neq i_1$. Therefore 
$\mrm{H}^i(T) = 0$ for $i \neq i_1$. By truncation we get an 
isomorphism $T \cong P[n]$ in 
$\dcat{D}(B \otimes A^{\mrm{op}})$,
where $n := -i_1$.
Finally, by Proposition \ref{prop:3575} the $B$-$A$-bimodule $P$ is invertible.
\end{proof}

Let $A$ be a commutative ring. An $A$-module $P$ can 
be viewed as a central $A$-bimodule. If $P$ is a rank $1$ projective 
$A$-module, then as a bimodule it is invertible. The usual Picard group of $A$
(see \cite[Section II.6]{Har})
is then the subgroup $\opn{Pic}_{A}(A)$ of $\opn{Pic}_{\K}(A)$, whose elements 
are the isomorphism classes of the central $A$-bimodules; and we refer to it 
here as the {\em commutative Picard group of $A$}. 

\begin{exer} \label{exer:3592}
Let $A$ be a commutative ring. We denote by $\opn{Aut}_{\K}(A)$ the group 
of ring automorphisms of $A$ (i.e.\ the Galois group). Show that 
\[ \opn{Pic}_{\K}(A) \cong \opn{Aut}_{\K}(A) \ltimes \opn{Pic}_{A}(A) . \]
(Hints: (1) Every invertible $A$-bimodule $P$ is isomorphic to 
$P' \ot_A P''$, where $P'$ is free of rank $1$ as a left $A$-module, and 
$P''$ is a central invertible bimodule.
This is proved like Theorem \ref{thm:3585}. 
(2) If $P'$ is an invertible $A$-bimodule that's free of rank $1$ as a left 
$A$-module, then it is also free as a right $A$-module, with the same basis 
element $e \in P'$. The proof is like those of Lemma \ref{lem:3575} and Theorem 
\ref{thm:3585}. Conclude that $P' \cong A(\ga)$, where 
$\ga \in \opn{Aut}_{\K}(A)$, and $A(\ga)$ is the twisted bimodule that is the
left $A$-module $A$, with right action 
$p \cd a := \ga(a) \cd p$ for $p \in A(\ga)$ and $a \in A$.)
\end{exer}

\begin{cor} \label{cor:3591} 
\index{1-DPic(A)@$\opn{DPic}_{\K}(A)$}
Let $A$ be a commutative ring with connected spectrum. Then 
\[ \opn{DPic}_{\K}(A) = \opn{Pic}_{\K}(A) \times \Z . \]
\end{cor}

\begin{proof}
The same as that of Corollary \ref{cor:3590}. 
\end{proof}

\begin{cor} \label{cor:3597} 
Let $A$ and $B$ be rings, and assume that $A$ is commutative ring with 
connected spectrum. 
If there is an equivalence of triangulated categories
$\dcat{D}(A) \to \dcat{D}(B)$, then there is an equivalence of 
abelian categories $\dcat{M}(A) \to \dcat{M}(B)$. In other words, $A$ and $B$ 
are Morita equivalent.
\end{cor}

\begin{proof}
According to Theorem \ref{thm:3565} there exists a tilting $B$-$A$-bimodule 
complex $T$. Theorem \ref{thm:3585} says that $T \cong P[n]$ for an 
invertible $B$-$A$-bimodule $P$. Then 
$P \ot_A (-) : \dcat{M}(A) \to \dcat{M}(B)$ 
is $\K$-linear equivalence of abelian categories. 
Moreover, $P$ is a projective $A$-module of rank $r$ for some positive integer 
$r$, and 
$B \cong \opn{End}_A(P)^{\mrm{op}}$. 
\end{proof}

\begin{rem} \label{rem:3610}
In the situation of Corollary \ref{cor:3597}, we can view $B$ as a central 
$A$-ring, using the homomorphism $g_T$ from Definition \ref{dfn:5090}. Then $B$ 
is an Azumaya $A$-ring; see \cite[Section 5.3]{Row}. Moreover, letting 
$X := \opn{Spec}(A)$, and letting $\BB$ be the sheafification of $B$ to $X$, 
then $\BB$ is a trivial Azumaya $\OO_X$-ring, in the sense of 
\cite[Section IV.2]{Mil}.
\end{rem}

Let $A$ be a commutative ring. If $\opn{Spec}(A)$ has a finite connected 
component decomposition, in the sense of Definition \ref{dfn:3210}(1),
then the connected component decomposition of $A$ is 
$A = \prod_{i = 1}^m \, A_i$, see Definition \ref{dfn:3210}(2).

\begin{cor} \label{cor:3585}
Let $A$ and $B$ be rings. Assume $A$ is commutative, and it 
has a finite connected component decomposition 
$A = \prod_{i = 1}^m \, A_i$.
Let $T$ be a tilting $B$-$A$-bimodule complex. We consider $B$ as a central 
$A$-ring via the isomorphism $g_T$ from Definition \tup{\ref{dfn:5090}}. Then 
$T \cong \bigoplus\nolimits_{i = 1}^m P_i[n_i]$
in $\dcat{D}(B \ot A^{\mrm{op}})$,
where $B_i := A_i \ot_A B$, $P_i$ is an invertible $B_i$-$A_i$-bimodule, and 
$n_i \in \Z$. 
\end{cor}

\begin{cor} \label{cor:3595} 
Let $A$ be a commutative ring, and assume it has 
a finite connected component decomposition $A = \prod_{i = 1}^m \, A_i$.
Then 
\[ \opn{DPic}_{\K}(A) = \opn{Pic}_{\K}(A) \times \Z^m . \] 
\end{cor}

\begin{cor} \label{cor:3598} 
Let $A$ and $B$ be rings. Assume that $A$ is commutative and it has a  
finite connected component decomposition. 
If there is an equivalence of triangulated categories
$\dcat{D}(A) \to \dcat{D}(B)$, then there is an equivalence of 
abelian categories $\dcat{M}(A) \to \dcat{M}(B)$. 
\end{cor}

\begin{exer} \label{exer:3585}
Prove Corollaries \ref{cor:3585}, \ref{cor:3595} and \ref{cor:3598}. Cf.\ 
Exercise \ref{exer:2176}. 
\end{exer}

The general case of these last three corollaries, namely when $A$ 
does not have a finite connected component decomposition, is dealt with in 
\cite[Theorem 6.13]{Ye8}. 

\begin{rem} \label{rem:4800}
The results in this subsection have variants in which the flatness condition
on the $\K$-rings $A$ and $B$ is dropped. 

For instance, Theorem \ref{thm:3575} can be modified as follows: If $A$ and $B$ 
are rings, with $A$ local, and if $T$ is a pretilting $B$-$A$-bimodule complex, 
then there is an isomorphism 
$T \cong P[n]$ in $\dcat{D}(A \ot B^{\mrm{op}})$ 
for some invertible $B$-$A$-bimodule $P$ and some integer $n$. 
Here is an outline of the proof. One chooses nonpositive K-flat resolutions
$\til{A} \to A$ and $\til{B} \to B$
in $\catt{DGRng} \centover \K$. 
Next a version of the K\"unneth trick (Lemma \ref{lem:3595}) for nonpositive DG 
rings has to be proved: For 
$M \in \dcat{D}^-(\til{B} \ot \til{A}^{\mrm{op}})$
and 
$N \in \dcat{D}^-(\til{A} \ot \til{B}^{\mrm{op}})$
satisfying 
$\opn{sup}(\opn{H}(M)) \leq i_1$ and $\opn{sup}(\opn{H}(N)) \leq j_1$,
there is an isomorphism
\[ \opn{H}^{i_1 + j_1}(M \ot^{\mrm{L}}_{\til{A}} N) \cong 
\opn{H}^{i_1}(M) \ot_A \opn{H}^{j_1}(N)  \]
in $\dcat{M}(B^{\mrm{en}})$. 
Then the proof of Theorem \ref{thm:3575} is 
repeated, for the tilting DG bimodule 
$T \in \dcat{D}(\til{A} \ot \til{B}^{\mrm{op}})$. 

The proofs of the other results are along the same lines. 
\end{rem}

\begin{rem} \label{rem:3590}
Assume $A$ is a commutative DG ring (Definitions \ref{dfn:3091} and 
\ref{dfn:3090}), with reduction $\bar{A} = \opn{H}^0(A)$. Here we do not need 
to assume flatness over the base ring $\K$. A DG $A$-module can 
be viewed as a central DG $A$-bimodule, so we have a derived tensor product 
\[ (- \ot^{\mrm{L}}_{A} -) : \dcat{D}(A) \times \dcat{D}(A) \to \dcat{D}(A) . 
\]
This operation is symmetric, i.e.\ 
$S \ot^{\mrm{L}}_{A} T \cong T \ot^{\mrm{L}}_{A} S$. 

A DG $A$-module $T$ is called a {\em
tilting DG $A$-module} if there is some $S \in \dcat{D}(A)$ such that 
$S \ot^{\mrm{L}}_{A} T \cong A$.
The {\em commutative derived Picard group} of $A$ is the group 
$\opn{DPic}_{A}(A)$, whose elements are the isomorphism classes in $\dcat{D}(A)$
of the tilting DG modules, the multiplication is induced by 
$(- \ot^{\mrm{L}}_{A} -)$, and the unit is the class of $A$. It is an abelian 
group. 

In \cite[Theorem 6.14]{Ye8} it was proved that the group homomorphism 
$\opn{DPic}_{A}(A) \to \opn{DPic}_{\bar{A}}(\bar{A})$,
which induced by the derived reduction functor 
$\bar{A} \ot^{\mrm{L}}_{A} (-)$,
is bijective. If the commutative ring $\bar{A}$
has a finite connected component decomposition 
$A = \prod_{i = 1}^m A_i$, then there is an isomorphism
\[ \opn{DPic}_{\bar{A}}(\bar{A}) \cong 
\opn{Pic}_{\bar{A}}(\bar{A}) \times \Z^m  . \]
This is according to Corollary \ref{cor:3595} with 
$\K = \bar{A}$. 
\end{rem}

As Corollaries \ref{cor:3590} and \ref{cor:3591} show, if the ring $A$ is 
either local, or commutative with connected spectrum, then the group 
$\opn{DPic}_{\K}(A)$ is not very interesting: it is 
\begin{equation} \label{4180}
\opn{DPic}_{\K}(A) = \opn{Pic}_{\K}(A) \times \bra{\si} , 
\end{equation}
where $\opn{Pic}_{\K}(A)$ is the classical contribution, and 
$\bra{\si} \cong \Z$ is the subgroup generated by the element $\si$, which is 
the class of the tilting complex $A[1]$.

The next example shows that matters are very different when $A$ is neither 
commutative nor local.

\begin{exa} \label{exa:4180}
Let $\K$ be an algebraically closed field and let $n$ be an integer $\geq 2$. 
Consider the ring $A$ of upper-triangular $n \times n$ matrices. For $n = 2$ it 
is $A = \sbmat{ \K & \K \\[0.1em] 0 & \K}$.
The group $\opn{DPic}_{\K}(A)$ contains a new element in this case. 
The classical NC Picard group $\opn{Pic}_{\K}(A)$ is trivial here. But the 
bimodule $A^* := \opn{Hom}_{\K}(A, \K)$
is tilting, and its class $\mu \in \opn{DPic}_{\K}(A)$ satisfies 
\begin{equation} \label{eqn:4185}
\mu^{n + 1} = \si^{n - 1} . 
\end{equation}
This was calculated by A. Yekutieli (for $n = 2$) and E. Kreines (for 
$n \geq 3$) in \cite{Ye4}. When $n = 2$ this says that 
$\si = \mu^3$, so the group $\opn{DPic}_{\K}(A)$ is larger than the classical 
part 
$\opn{Pic}_{\K}(A) \times \bra{\si} = \bra{\si} \cong \Z$.

Later, in \cite{MiYe}, J.-I. Miyachi and Yekutieli showed that the group 
$\opn{DPic}_{\K}(A)$ is abelian, it is generated by $\si$ and $\mu$, and the 
only relation is (\ref{eqn:4185}). Note that the ring $A$ is isomorphic to the 
path ring $\K[Q]$ of the Dynkin quiver $Q$ of type $\mrm{A}_n$ with all
arrows going in the same direction. 
The paper \cite{MiYe} contains calculations 
of the groups $\opn{DPic}_{\K}(A)$ for several other types of path rings of 
quivers. 

M. Kontsevich interprets the relation (\ref{eqn:4185}) as follows: 
the {\em fractional Calabi-Yau dimension} of the ``noncommutative smooth proper 
space'' $\dcat{D}^{\mrm{b}}_{\mrm{f}}(A)$ is $\frac{n - 1}{n + 1}$. 
More on this in Remark \ref{rem:4505} and Example \ref{exa:4870}. 
\end{exa}

\begin{rem} \label{rem:3430}
Here are a few historical notes on Theorems \ref{thm:3575} and  \ref{thm:3585}.
The concept of derived Picard group was discovered independently, around 
1997, by R. Rouquier and A. Zimmermann \cite{RoZi} (who had used the notation
$\opn{TrPic}(A)$ for this group) and Yekutieli \cite{Ye4}. The motivations of 
the two teams of authors were very different: Rouquier and Zimmermann were 
interested in invariants of finite dimensional algebras, whereas Yekutieli was 
trying to classify noncommutative dualizing complexes (See Section
\ref{sec:rigid-DC-NC}).
\end{rem}

\begin{rem} \label{rem:4650}
Suppose $\K$ is an algebraically closed field and $A$ is a finite $\K$-ring. 
It turns out that in this case the group $\opn{DPic}_{\K}(A)$ is a {\em locally 
algebraic group}. Namely there is a connected algebraic group 
$\opn{DPic}^0_{\K}(A)$, and $\opn{DPic}_{\K}(A)$ it is a (usually infinite) 
disjoint union of connected algebraic varieties, that are left and right cosets 
of  $\opn{DPic}^0_{\K}(A)$. This same sort of geometric structure can be found 
in the {\em Picard scheme} of a smooth projective algebraic variety $X$, in 
which the identity component is the {\em Jacobian variety} of $X$. But whereas 
the Picard scheme is an abelian group scheme, the group $\opn{DPic}_{\K}(A)$ is 
not 
abelian; and this means that there is a nontrivial geometric action of 
$\opn{DPic}^0_{\K}(A)$, by conjugation, on the connected components. 

The result above is in the paper \cite{Ye4.6}, and it is based on the papers 
\cite{HuSa} by B. Huisgen-Zimmermann and M. Saorin, and \cite{Rou} by Rouquier. 

Later B. Keller \cite{Kel3} used an abstraction of this idea -- essentially  
viewing the derived Picard group as a derived group stack (evaluating it on 
commutative artinian DG rings) -- to prove that the {\em Hochschild cohomology} 
of $A$ is the {\em Lie algebra}, in the derived sense, of $\opn{DPic}_{\K}(A)$. 
This allowed him to prove that the Hochschild cohomology of $A$, with its {\em 
Gerstenhaber Lie bracket}, is a derived invariant of $A$. 
\end{rem}

%% file: block6_190413.tex

\renewcommand{\thisfile}{block6\_190328}  

\cleardoublepage
\mysection{Algebraically Graded Noncommutative Rings}
\label{sec:alg-gra-rings}

\AYcopyright

In this section we study {\em algebraically graded rings}, and the categories 
of {\em algebraically graded modules} over them. These are the graded rings 
that appear in standard texts on commutative and noncommutative algebra, and 
are distinct from the {\em cohomologically graded rings} that underlie {\em DG 
rings}. 

In the subsequent sections of the book (Sections \ref{sec:der-tors-conn} and 
\ref{sec:BDC}) we shall concentrate on {\em connected graded rings}
(traditionally known as {\em connected graded algebras}), that behave 
like ``complete local rings'' within the algebraically graded context. 

The interest in algebraically graded rings, and especially in connected graded 
rings, stems from the prominent role they played in {\em noncommutative 
algebraic geometry}, as it was developed by M. Artin and his collaborators 
since around 1990. See the papers \cite{ATV}, \cite{ArZh} and \cite{StaVdB}.

\mysubsection{Categories of Algebraically Graded Modules} 
\label{subsec:alg-gr-mods}

Let $\K$ be a nonzero commutative base ring. An 
{\em algebraically graded $\K$-module}%
\index{Algebraically graded module}
is a $\K$-module $M$ with a direct sum decomposition 
$M = \bigoplus_{i \in \Z} M_i$
into $\K$-submodules. The submodule $M_i$ is called the homogeneous component 
of $M$ of {\em algebraic degree $i$}.
If $m \in M_i$ is a nonzero element, then we write 
$\opn{deg}(m) = i$. 

In Subsection \ref{subsec:gr-alg} we discussed {\em cohomologically graded rings 
and  modules}. There are three features that distinguish between 
cohomologically graded $\K$-mod\-ules and algebraically graded $\K$-modules.
The first distinguishing feature is the use of upper vs.\ lower 
indices. The second is that in an algebraically graded module 
there is never a differential involved; i.e.\ it is 
not the underlying graded module of a DG module. The third and most important 
change is in {\em signs of permutations}, and thus in {\em commutativity}. In 
the cohomologically graded setting the Koszul sign rule 
dictates that for graded $\K$-modules $M$ and $N$, the braiding isomorphism
$\opn{br}_{M, N} : M \ot N \to N \ot M$,
where $\ot := \ot_{\K}$, is 
$\opn{br}_{M, N}(m \ot n) = (-1)^{i \cd j} \cd n \ot m$
for homogeneous elements $m \in M^i$ and $n \in N^j$. 
But in the algebraically graded setting there are no signs: 
$\opn{br}_{M, N}(m \ot n) := n \ot m$
for $m \in M_i$ and $n \in N_j$. 
This is reflected in Definition \ref{dfn:3738} below. 

Let $M = \bigoplus_{i \in \Z} \, M_i$ be an algebraically graded $\K$-module.
A nonzero element $m \in M$ can be expressed uniquely as a sum 
\begin{equation} \label{eqn:4260}
m = m_{1} + \cdots + m_{r} , 
\end{equation}
such that $r \geq 1$; each $m_{i}$ is a nonzero homogeneous element; and 
$\opn{deg}(m_{i}) < \opn{deg}(m_{i + 1})$. 
This is called the {\em homogeneous component decomposition} of $m$. 
Of course $m$ is homogeneous if and only if $r = 1$. 

Suppose $M$ and $N$ are algebraically graded $\K$-modules. A $\K$-linear 
homomorphism $\phi : M \to N$
is said to be of {\em algebraic degree $i$} if
$\phi(M_j) \sub N_{j + i}$ for all $j \in \Z$. 
The $\K$-module of all $\K$-linear homomorphisms of degree $i$ is denoted by 
$\opn{Hom}_{\K}(M, N)_i$. Taking the direct sum we get an algebraically graded 
$\K$-module
\begin{equation}  \label{eqn:3701}
\opn{Hom}_{\K}(M, N) := \bigoplus_{i \in \Z} \, \opn{Hom}_{\K}(M, N)_i . 
\end{equation}

\begin{dfn}  \label{dfn:3700}
Let $M$ and $N$ be algebraically graded $\K$-modules. 
A {\em strict homomorphism of algebraically graded $\K$-modules}
\index{Algebraically graded module! strict homomorphism of {\indash}s}
$\phi : M \to N$ is $\K$-linear homomorphism of algebraic degree $0$.

The {\em category of algebraically graded $\K$-modules} is the category 
$\dcat{M}(\K, \mrm{gr})$, 
\index{1-M(K,gr)@$\dcat{M}(\K, \mrm{gr})$}
whose objects are the algebraically graded 
$\K$-modules, and whose morphisms are the strict homomorphisms.
\end{dfn}

Thus 
\[ \opn{Hom}_{\dcat{M}(\K, \mrm{gr})}(M, N) = \opn{Hom}_{\K}(M, N)_0 . \]
It is easy to see that $\dcat{M}(\K, \mrm{gr})$ is a $\K$-linear abelian 
category. The kernels and cokernels are degreewise. 

Suppose $M, N \in \dcat{M}(\K, \mrm{gr})$.
Their tensor product is also algebraically graded:
\[ (M \ot N)_i = \bigoplus_{j \in \Z} \, (M_j \ot N_{i - j}) \]
and 
\[ M \ot N = \bigoplus_{i \in \Z} \, (M \ot N)_i . \]

\begin{dfn}  \label{dfn:3701}
An {\em algebraically graded central $\K$-ring}
\index{Algebraically graded ring}
is a central $\K$-ring $A$, equipped with a direct sum decomposition 
$A = \bigoplus_{i \in \Z} A_i$ 
into $\K$-submodules, such that $1_A \in A_0$, and 
$A_i \cd A_j \sub A_{i + j}$ 
for all $i, j$.

A {\em homomorphism of algebraically graded central $\K$-rings}
\index{Algebraically graded ring! homomorphism of {\indash}s}
$f : A \to B$ is a $\K$-ring homomorphism such that 
$f(A_i) \sub B_i$ for all $i$. 

The algebraically graded central $\K$-rings form a category, that we
denote by $\catt{Rng}_{\mrm{gr}} \centover \K$.
\end{dfn}

The base ring $\K$ is an algebraically graded ring, concentrated in degree $0$.
It is the initial object of the category 
$\catt{Rng}_{\mrm{gr}} \centover \K$. 

In the traditional ring theory literature (e.g.\ \cite{McRo}, \cite{Row}, 
\cite{ATV}, \cite{ArZh}, \cite{Ye1} and \cite{YeZh0}) the algebraically graded 
central $\K$-rings were called ``associative unital graded $\K$-algebras''. 

\begin{exa} \label{exa:3701}
Let $M$ be an algebraically graded $\K$-module. Then 
$\opn{End}_{\K}(M) \lb  := \opn{Hom}_{\K}(M, M)$
is an algebraically graded central $\K$-ring.
The grading is according to (\ref{eqn:3701}), and multiplication is 
composition. 
\end{exa}

\begin{dfn}  \label{dfn:3702}
Let $A$ be an algebraically graded central $\K$-ring.
\begin{enumerate}
\item An {\em algebraically graded left $A$-module} 
\index{Algebraically graded module}
is a left $A$-module $M$, equipped with a direct sum decomposition 
$M = \bigoplus_{i \in \Z} \, M_i$
into $\K$-submodules, such that 
$A_i \cd M_j \sub M_{i + j}$ for all $i, j$.

\item Suppose $M$ and $N$ are algebraically graded left $A$-modules.
A $\K$-linear homomorphism $\phi : M \to N$ of degree $i$ is said to be a
{\em homomorphism of algebraically graded $A$-modules}
\index{Algebraically graded module! homomorphism of {\indash}s}
if $\phi(a \cd m) = a \cd \phi(m)$
for all $a \in A$ and $m \in M$. 

\item The $\K$-module of all $A$-linear homomorphisms of degree $i$ is denoted 
by $\opn{Hom}_{A}(M, N)_i$. Taking the direct sum we get an algebraically 
graded 
$\K$-module
\[ \opn{Hom}_{A}(M, N) := \bigoplus_{i \in \Z} \, \opn{Hom}_{A}(M, N)_i . \]
\end{enumerate}
\end{dfn}

Note that there is an inclusion 
$\opn{Hom}_{A}(M, N) \sub \opn{Hom}_{\K}(M, N)$
in $\dcat{M}(\K, \mrm{gr})$.

\begin{dfn}  \label{dfn:3704}
Let $A$ be an algebraically graded central $\K$-ring.
\begin{enumerate}
\item Suppose $M$ and $N$ are algebraically graded left $A$-modules.
The elements of $\opn{Hom}_{A}(M, N)_0$ are called 
{\em strict homomorphisms of algebraically graded $A$-modules}. 
\index{Algebraically graded module! strict homomorphism of {\indash}s}

\item The category of algebraically graded left $A$-modules, with strict 
homomorphisms, is denoted by $\dcat{M}(A, \mrm{gr})$. 
\index{1-M(A,gr)@$\dcat{M}(A, \mrm{gr})$}
\end{enumerate}
\end{dfn}

Thus 
\begin{equation}  \label{eqn:3736}
\opn{Hom}_{\dcat{M}(A, \mrm{gr})}(M, N) = \opn{Hom}_{A}(M, N)_0 . 
\end{equation}

Clearly, to put an algebraically graded left $A$-module structure on an
algebraically graded $\K$-module $M$ is the same as to give a
homomorphism $A \to \lb \opn{End}_{\K}(M)$ in 
$\catt{Rng}_{\mrm{gr}} \centover \K$.

Algebraically graded right $A$-modules, and algebraically graded 
$A$-$B$-bi\-modules, are defined similarly to Definition \ref{dfn:3702}. There 
are no signs anywhere. 
All homomorphisms are $\K$-linear, and all bimodules are central over $\K$. 

\begin{dfn}  \label{dfn:3735}
Let $A$ be an algebraically graded ring. Given an algebraically graded 
$A$-module $M$ and an integer $i$, the {\em $i$-th degree twist of $M$} 
\index{Algebraically graded module! degree twist of}
\index{1-M(i)@$M(i)$}
is the algebraically graded $A$-module $M(i)$ whose degree $j$ homogeneous 
component is $M(i)_j := M_{i + j}$.
The action of $A$ is not changed. 
\end{dfn}

The degree twist is an automorphism of the abelian category 
$\dcat{M}(A, \mrm{gr})$. And 
\[ \opn{Hom}_{\dcat{M}(A, \mrm{gr})}(M, N(i)) = \opn{Hom}_{A}(M, N)_i . \]
There is a similar twisting for graded right modules and for graded 
bimodules. 

There is a functor 
$\opn{Ungr} : \catt{Rng}_{\mrm{gr}} \centover \K \to 
\catt{Rng}_{} \centover \K$
that forgets the grading. For a fixed ring 
$A \in  \catt{Rng}_{\mrm{gr}} \centover \K$
there is a forgetful functor 
\begin{equation} \label{eqn:3715}
\opn{Ungr} : \dcat{M}(A, \mrm{gr}) \to \dcat{M}(\opn{Ungr}(A)) . 
\end{equation}
This is an exact faithful $\K$-linear functor. 

\begin{rem} \label{rem:3795}
The functor $\opn{Ungr}$ in (\ref{eqn:3715}) is usually not full, for two 
reasons. First, there could be nonzero homomorphisms of nonzero 
degree, so that 
\[ \opn{Hom}_{\dcat{M}(A, \mrm{gr})}(M, N) = 
\opn{Hom}_A(M, N)_0 \subsetneqq \opn{Hom}_A(M, N) . \]
Second, it is easy find an example where
\[ \opn{Ungr} \bigl( \opn{Hom}_A(M, N) \bigr) \subsetneqq
\opn{Hom}_{\opn{Ungr}(A)} \bigl( \opn{Ungr}(M), \opn{Ungr}(N) \bigr) . \]
Just take 
$M := \bigoplus_{i \in \Z} \, A(i)$ and $N := A$. Then
\[ \opn{Ungr} \bigl( \opn{Hom}_A(M, N) \bigr) \cong 
\bigoplus_{i \in \Z} \, \opn{Ungr}(A)  \]
and 
\[ \opn{Hom}_{\opn{Ungr}(A)} \bigl( \opn{Ungr}(M), \opn{Ungr}(N) \bigr) \cong 
\prod_{i \in \Z} \, \opn{Ungr}(A) \,  \]
in $\dcat{M}(\K)$.
\end{rem}

\begin{rem} \label{rem:3900}
Homomorphisms between algebraically graded modules and cohomologically graded 
modules are distinct, not only in the position of the degree indices. 
Compare Definition \ref{dfn:3702}(2) to Definition \ref{dfn:1690}. 
This distinction is due to the Koszul sign rule, that is present only in the 
cohomologically graded setting. 
\end{rem}

\begin{dfn}  \label{dfn:3738}
Let $A$ and $B$ be algebraically graded central $\K$-rings. 
\begin{enumerate}
\item The {\em opposite ring} of $A$ is the central $\K$-ring 
$A^{\mrm{op}}$, that is the same algebraically graded $\K$-module as $A$, with 
an isomorphism
$\opn{op} := \opn{id} : A \iso A^{\mrm{op}}$
in $\dcat{M}(\K, \mrm{gr})$. The multiplication $\cd^{\mrm{op}}$ of 
$A^{\mrm{op}}$ is 
\[ \opn{op}(a_1) \, \cd^{\mrm{op}} \, \opn{op}(a_2) := \opn{op}(a_2 \cd a_1) . 
\]

\item The ring $A$ is called a {\em commutative graded ring} 
\index{Algebraically graded ring! commutative}
if  $a_2 \cd a_1 = a_1 \cd a_2$ for all $a_1, a_2 \in A$. 

\item The tensor product $A \ot B$ is made into an algebraically graded ring 
with multiplication 
\[ (a_1 \ot b_1) \cd (a_2 \ot b_2) :=  
(a_1 \cd a_2) \ot (b_1 \ot b_2)  . \]

\item The {\em enveloping} ring of $A$ is the  algebraically graded ring 
$A^{\mrm{en}} := A \ot A^{\mrm{op}}$.
\end{enumerate}
\end{dfn}

In terms of opposite rings, an algebraically graded ring $A$ is commutative if 
and only if $A = A^{\mrm{op}}$. 
The category of algebraically graded right $A$-modules is identified with  
$\dcat{M}(A^{\mrm{op}}, \mrm{gr})$; 
and the category of $\K$-central algebraically graded $A$-$B$-bimodules is  
identified with $\dcat{M}(A \ot B^{\mrm{op}}, \mrm{gr})$. 

There are graded ring homomorphisms 
$A \to A \ot B \lar B$,
and corresponding restriction functors 
\begin{equation} \label{eqn:3740}
\dcat{M}(A, \mrm{gr}) \xlar{\opn{Rest}_A}
\dcat{M}(A \ot B, \mrm{gr}) \xar{\opn{Rest}_B}
\dcat{M}(B, \mrm{gr})  . 
\end{equation}
We shall usually omit any explicit reference to these restriction 
functors.

Given $M \in \dcat{M}(A, \mrm{gr})$ and 
$N \in \dcat{M}(B, \mrm{gr})$, the tensor product 
$M \ot N$ belongs to 
$\dcat{M}(A \ot B, \mrm{gr})$. 

Suppose $M \in \dcat{M}(A^{\mrm{op}}, \mrm{gr})$
and $N \in \dcat{M}(A, \mrm{gr})$.
Their tensor product over $A$ is also algebraically graded:
\begin{equation} \label{eqn:3720}
M \ot_A N = \bigoplus_{i \in \Z} \, (M \ot_A N)_i 
\in \dcat{M}(\K, \mrm{gr}) , 
\end{equation}
where $(M \ot_A N)_i$ is the image of $(M \ot N)_i$ under the canonical 
surjection 
$M \ot N \surj M \ot_A N$
in $\dcat{M}(\K, \mrm{gr})$. 
See Lemma \ref{lem:1690} for the corresponding statement in the cohomologically 
graded setting. 

There is an isomorphism of algebraically graded $\K$-rings
\begin{equation} \label{3721}
(A^{\mrm{en}})^{\mrm{op}} \iso A^{\mrm{en}} , \quad 
a_1 \ot a_2 \mapsto a_2 \ot a_1 . 
\end{equation}

From here on in this section, and in Sections \ref{sec:der-tors-conn}
and \ref{sec:BDC}, we assume the following convention. 

\begin{conv}  \label{conv:3700}
We fix a nonzero commutative base ring $\K$. The symbol $\ot$ means $\ot_{\K}$. 
All graded rings are assumed to be 
algebraically graded central $\K$-rings, and all graded ring homomorphisms are 
over $\K$. I.e.\  we work inside the category 
$\catt{Rng}_{\mrm{gr}} \centover \K$.

Let $A$ be a graded ring. Its enveloping ring is the graded ring 
$A^{\mrm{en}} := A \ot A^{\mrm{op}}$. 
All graded $A$-modules are assumed by default to be 
algebraically graded left $A$-modules.
For $M, N \in \dcat{M}(A, \mrm{gr})$ the expression 
$\opn{Hom}_A(M, N)$ refers to the graded $\K$-module in Definition 
\ref{dfn:3702}(3). 
All graded bimodules are assumed to be $\K$-central;  i.e.\ a graded
$A$-$B$-bimodule is an object of
$\dcat{M}(A \ot B^{\mrm{op}}, \mrm{gr})$.
\end{conv}

\begin{prop} \label{prop:3825}
Let $A$ be a graded ring. 
\begin{enumerate}
\item The category $\dcat{M}(A, \mrm{gr})$ is abelian. 

\item The category $\dcat{M}(A, \mrm{gr})$ has inverse limits. 

\item The category $\dcat{M}(A, \mrm{gr})$ has direct limits. Moreover, direct 
limits in $\dcat{M}(A, \mrm{gr})$ commute with the ungrading functor
\tup{(\ref{eqn:3715})}.
\end{enumerate}
\end{prop}

\begin{proof}
(1) The kernels and cokernels in $\dcat{M}(A, \mrm{gr})$ are degreewise, so the 
abelian property of $\dcat{M}(A, \mrm{gr})$ is inherited from the 
abelian property of $\dcat{M}(\K)$.

\medskip \noindent
(2) Inverse limits in  $\dcat{M}(A, \mrm{gr})$ are taken degreewise. 
Namely, given an inverse system  
$\{ M_x \}_{x \in X^{\mrm{op}}}$ of objects of $\dcat{M}(A, \mrm{gr})$,
indexed by a directed set $X$, let 
$M_i := \lim_{\lar x} (M_x)_i \in \dcat{M}(\K)$.
Then 
$M := \bigoplus_{i \in \Z} M_i \in \dcat{M}(A, \mrm{gr})$
is the inverse limit of $\{ M_x \}_{x \in X^{\mrm{op}}}$.

\medskip \noindent
(3) Direct limits in $\dcat{M}(A, \mrm{gr})$ are taken degreewise. 
Namely, given a direct system  
$\{ N_y \}_{y \in Y}$ of objects of $\dcat{M}(A, \mrm{gr})$,
indexed by a directed set $Y$, let 
$N_i := \lim_{y \to} (N_y)_i \in \dcat{M}(\K)$.
Then 
$N := \bigoplus_{i \in \Z} N_i \in \dcat{M}(A, \mrm{gr})$
is the direct limit of $\{ N_y \}_{y \in Y}$.

Because direct sums and direct limits commute with each other, it follows that 
$\lim_{y \to}$ commutes with $\opn{Ungr}$. 
\end{proof}

\begin{dfn} \label{dfn:4540}
A graded ring $A$  is called {\em nonnegative} 
\index{Algebraically graded ring! nonnegative}
if $A_i = 0$ for all $i < 0$; i.e.\ if $A = \bigoplus_{i \geq 0} A_i$.

If $A$ is a nonnegative graded ring $A$, then the two-sided graded ideal
$\m :=  \bigoplus_{i \geq 1} A_i$ is called the {\em augmentation ideal}. 
\index{Algebraically graded ring! augmentation ideal of}
It is the kernel of the surjective graded ring homomorphism $A \to A_0$, that we 
call the {\em augmentation homomorphism}. 

The augmentation ideal of the graded ring $A^{\mrm{op}}$ is denoted by 
$\m^{\mrm{op}}$. 
\end{dfn}

\begin{dfn} \label{dfn:4541}
Let $A$ be a graded ring $A$ and let 
$M = \bigoplus_{i \in \Z} M_i$ be a graded $A$-module.
\begin{enumerate}
\item We call $M$ a {\em bounded above} graded module if $M_i = 0$ for 
$i \gg 0$.

\item We call $M$ a {\em bounded below} graded module if $M_i = 0$ for 
$i \ll 0$.

\item We call $M$ a {\em bounded} graded module if it is both bounded above and 
bounded below.
\end{enumerate}
\end{dfn}

\begin{exa} \label{exa:4540}
If $A$ is a nonnegative graded ring, then as a graded module $A$ is bounded 
below. For such $A$, the bounded below graded modules behave like finitely 
generated modules over a noetherian local commutative ring -- see Proposition 
\ref {prop:4540} below.
\end{exa}

The next proposition is the first version of the graded Nakayama Lemma that we 
present; the second is Proposition \ref{prop:4590}. 

\begin{prop} \label{prop:4540}
Let $A$ be a nonnegative graded ring, with augmentation ideal $\m$. Let $M$ be 
a bounded below graded $A$-module, and let $M' \sub M$ be a graded 
$A$-submodule. If $M' + \m \cd M = M$ then $M' = M$. 
\end{prop}

\begin{proof}
Say $M = \bigoplus_{i \geq i_0} M_i$ for an integer $i_0$.
We will prove that $M'_i = M_i$
as $A_0$-modules for all $i$, by induction on $i \geq i_0 - 1$.
For $i =  i_0 - 1$ this is clear because both these $A_0$-modules are zero. 

Now take $i \geq i_0$, and assume that $M'_j = M_j$ for all  $j < i$. 
Suppose $m \in M_i$ is a nonzero element. 
Because $M = M' + \m \cd M$, we can express $m$ as follows:
$m = m' + \sum_{k = 1}^r a_k \cd m_k$,
where $m' \in M'_i$; $r \geq 0$; $a_k \in \m$ and $m_k \in M$ are 
nonzero homogeneous elements; 
and $\opn{deg}(a_k) + \opn{deg}(m_k) = i$.
But then $\opn{deg}(m_k) < i$, so by the induction hypothesis
we have $m_k \in M'$. We conclude that $m \in M'_i$.
\end{proof}

\begin{dfn}  \label{dfn:3706}
Let $A$ be a graded ring. 
\begin{enumerate}
\item A graded $A$-module $M$ is called a {\em graded-finite  
$A$-module} if $M$ can be generated by finitely many homogeneous elements. 

\item The full subcategory of $\dcat{M}(A, \mrm{gr})$ on the graded-finite
$A$-modules is denoted by $\dcat{M}_{\mrm{f}}(A, \mrm{gr})$.
\index{1-Mf(A,gr)@$\dcat{M}_{\mrm{f}}(A, \mrm{gr})$}

\item $A$ is called a {\em left graded-noetherian  ring} if every graded left 
ideal $\a \sub A$ is a graded-finite  $A$-module, in the sense of item (1) 
above. 

\item $A$ is called a {\em graded-noetherian  ring} if both $A$ and 
$A^{\mrm{op}}$ are left graded-noetherian  rings. 
\end{enumerate}
\end{dfn}

An obvious modification of item (1) above gives the definition of a 
graded-finite  right $A$-module. Of course the ring $A$ is called a right 
graded-noetherian  ring if $A^{\mrm{op}}$ is a left graded-noetherian  ring.

\begin{prop} \label{prop:4260}
If $A$ is left graded-noetherian ring, then $\dcat{M}_{\mrm{f}}(A, 
\mrm{gr})$ is a thick abelian subcategory of $\dcat{M}(A, \mrm{gr})$. 
\end{prop}

\begin{exer} \label{exer:4261}
Prove this proposition. 
\end{exer}

\begin{prop} \label{prop:3745}
Let $A$ be a graded ring and $M$ a graded $A$-module. 
The following two conditions are equivalent.
\begin{itemize}
\rmitem{i} $M$ is a graded-finite  $A$-module, in the sense of Definition 
\tup{\ref{dfn:3706}(1)}.

\rmitem{ii} $\opn{Ungr}(M)$ is a finitely generated $\opn{Ungr}(A)$-module. 
\end{itemize}
\end{prop}

\begin{exer} \label{exer:3745}
Prove this proposition. 
\end{exer}

\begin{thm} \label{thm:4245}
Let $A$ be a graded ring. The following two conditions are equivalent.
\begin{itemize}
\rmitem{i} The graded ring $A$ is left graded-noetherian, in the sense of 
Definition \tup{\ref{dfn:3706}(3)}.

\rmitem{ii} The ring $\opn{Ungr}(A)$ is left noetherian.  
\end{itemize}
\end{thm}

\begin{proof}
The case when $A$ is nonnegative is proved in many textbooks, and it is a 
nice exercise for the reader. 

The general case is much harder, and the few proofs in the 
literature are quite complicated. Here is a short proof, that is due to S. 
Snigerov. 

The implication (ii) $\Rightarrow$ (i) is trivial. So let us assume that $A$ is 
left graded-noetherian, and we shall prove that $\opn{Ungr}(A)$ 
is left noetherian. The proof proceeds in a few steps. 

\medskip \noindent
Step 1. Consider an $A_0$-submodule $V \sub A_j$ for some $j$, and let  
$\b := A \cd V \sub A$, which is a graded left ideal of $A$. 
By assumption, $\b$ is generated by finitely many homogeneous elements. 
An easy calculation shows that $V = \b \cap A_j$, and that 
$V$ is finitely generated as an $A_0$-module. The consequence is that  
the ring $A_0$ is left noetherian, and that every $A_j$ is a finitely generated 
$A_0$-module. 

\medskip \noindent
Step 2. Now we fix a nonzero left ideal $L \sub A$. We are going to find a 
finite number of generators for $L$ as an $A$-module. 

For a nonzero element 
$a \in L$ we consider its homogeneous component decomposition 
(\ref{eqn:4260}), and we define its bottom degree component
$\opn{bot}(a) := a_{1}$, and its top degree component 
$\opn{top}(a) := a_{r}$.
Let $\a_{\mrm{top}}$ be the left ideal in $A$ generated by the elements
$\opn{top}(a)$, where $a$ runs over all the nonzero elements of $L$. By 
assumption the graded ideal $\a_{\mrm{top}}$ is graded-finite; so we 
can choose a finite collection $\{ a_x \}_{x \in X_{\mrm{top}}}$ of nonzero 
elements of $L$, such that their top degree components $\opn{top}(a_x)$ 
generate the left ideal $\a_{\mrm{top}}$. Now let 
\[ d_1 := \opn{max} \, \bigl\{ \opn{deg}(\opn{top}(a_x)) \mid x \in 
X_{\mrm{top}} \bigr\} \in \Z . \]

\medskip \noindent
Step 3. Define $\a_{\mrm{bot}}$ to be the left ideal in $A$ generated by the 
elements $\opn{bot}(a)$, where $a$ runs over all the nonzero elements of 
$L \cap A_{\leq d_1}$. (It is possible that there are no such elements.)
By assumption the graded left ideal $\a_{\mrm{bot}}$ is graded-finite; 
so we can choose a finite collection $\{ a_x \}_{x \in X_{\mrm{bot}}}$ of 
nonzero elements of $L \cap A_{\leq d_1}$, such that their bottom degree 
components $\opn{bot}(a_x)$ generate the left ideal $\a_{\mrm{bot}}$. 
Now we let 
\[ d_0 := \opn{min} \, \bigl( \bigl\{ \opn{deg}(\opn{bot}(a_x)) \mid x \in 
X_{\mrm{bot}} \bigr\} \cup \{ 0 \} \bigr) \in \Z . \]

\medskip \noindent
Step 4. By step 1 the $A_0$-module
$V := L \cap \bigl( \bigoplus\nolimits_{i = d_0}^{d_1} A_i \bigr)$
is finitely generated. We choose a finite collection 
 $\{ a_x \}_{x \in X_{\mrm{mid}}}$ 
that generates $V$ as an $A_0$-module.

\medskip \noindent
Step 5. Define the finite indexing set  
$X := X_{\mrm{top}} \sqcup X_{\mrm{bot}} \sqcup X_{\mrm{mid}}$.
We claim that the collection $\{ a_x \}_{x \in X}$ 
generates the left ideal $L$. Take any element $a \in L$. 
By subtracting from $a$ a finite linear combination of elements 
from $\{ a_x \}_{x \in X_{\mrm{top}}}$, with coefficients taken from 
$A_{\geq 1}$, we can assume that 
$a \in L \cap A_{\leq d_1}$. 
Next, by subtracting from $a$ a finite linear combination of elements 
from $\{ a_x \}_{x \in X_{\mrm{bot}}}$, with coefficients taken from 
$A_{\leq -1}$, we can assume that $a \in V$, the $A_0$-submodule from step 4.
So $a$ is finite linear combination of elements 
from $\{ a_x \}_{x \in X_{\mrm{mid}}}$, with coefficients taken from $A_0$. 
\end{proof}

In view of Proposition \ref{prop:3745} and Theorem \ref{thm:4245}, we can use 
this more relaxed terminology without danger of ambiguity: 

\begin{dfn} \label{dfn:4262}
Let $A$ be a graded ring. 
\begin{enumerate}
\item We call $A$ a {\em left noetherian graded ring} 
\index{Algebraically graded ring! left noetherian}
if it is left graded-noetherian. Likewise for the attributes ``right noetherian 
graded ring'' and  ``noetherian graded ring''. 

\item Let $M$ be a graded $A$-module. We call $M$ a 
{\em finite graded $A$-module}
\index{Algebraically graded module! finite}
if it is graded-finite. Likewise for right modules. 
\end{enumerate}
\end{dfn}

We do not assume the rings here are noetherian. Partly this is because of the 
facts in the remark below. Instead we resort to {\em derived 
graded-pseudo-finite} modules and complexes, see Definition \ref{dfn:4532}. 

\begin{rem} \label{rem:3901}
Unlike the commutative theory, the Hilbert Basis Theorem does not hold in the 
noncommutative setting. Thus, a finitely generated $\K$-ring $A$ need 
not be noetherian; see Example \ref{exa:3755}. Likewise, if $A$ and $B$ are 
noetherian finitely generated $\K$-rings, the tensor product $A \ot B$ need not 
be noetherian; there is a counterexample in \cite{Rog}. 
See \cite{ArSmZh} for a deep discussion of permanence properties 
of NC noetherian rings. 
\end{rem}

We now pass to complexes of algebraically graded modules, retaining Convention 
\ref{conv:3700}. 
In Subsection \ref{subsec:complexes} we introduced the category 
$\dcat{C}(\cat{M})$ of complexes with entries in an abelian category $\cat{M}$.
This is a DG category. Its strict subcategory 
$\dcat{C}_{\mrm{str}}(\cat{M})$ is abelian. Here we take the abelian category 
$\cat{M} := \dcat{M}(A, \mrm{gr})$ for a graded ring $A$. 

\begin{dfn} \label{dfn:3730}
Let $A$ be a graded ring.
\begin{enumerate}
\item The category of complexes with entries in the abelian category 
$\dcat{M}(A, \mrm{gr})$ is the DG category 
$\dcat{C}(A, \mrm{gr}) := \dcat{C}(\dcat{M}(A, \mrm{gr}))$%
\index{1-C(A,gr)@$\dcat{C}(A, \mrm{gr})$}.

\item The strict subcategory of 
$\dcat{C}(A, \mrm{gr})$ is the abelian category  
$\dcat{C}_{\mrm{str}}(A, \mrm{gr}) := 
\dcat{C}_{\mrm{str}}(\dcat{M}(A, \mrm{gr}))$%
\index{1-Cstr(A,gr)@$\dcat{C}_{\mrm{str}}(A, \mrm{gr})$}.

\item Given a boundedness indicator $\star$, we denote by 
$\dcat{C}^{\star}(A, \mrm{gr})$ the full subcategory of 
$\dcat{C}(A, \mrm{gr})$ on the complexes $M$ having boundedness condition 
$\star$. 
\end{enumerate}
\end{dfn}

We now try to make Definition \ref{dfn:3730} more explicit. 
An object $M \in \dcat{C}(A, \mrm{gr})$ has these direct sum decompositions:
\begin{equation} \label{eqn:3730}
M = \bigoplus_{i \in \Z} \, M^i = \bigoplus_{i, j \in \Z} \, M^i_j . 
\end{equation}
Here each $M^i \in \dcat{M}(A, \mrm{gr})$, and each
$M^i_j \in \dcat{M}(\K)$. 
We call $M^i_j$ the homogeneous component of $M$ of bidegree 
$\sbmat{i \\[0em] j} \in \sbmat{\Z \\[0.1em] \Z} = \Z^2$.
The upper index $i$ is called the {\em cohomological degree}, and the lower 
index $j$ is called the {\em algebraic degree}. 
The differential $\d_{M}$ of the complex $M$ is an $A$-linear homomorphism 
$\d_M : M \to M$ of bidegree $\sbmat{1 \\[0.1em] 0}$
that satisfies $\d_M \circ \d_M = 0$. 

A morphism
$\phi : M \to N$ in the strict category 
$\dcat{C}_{\mrm{str}}(A, \mrm{gr})$
has bidegree $\sbmat{0 \\[0.1em] 0}$ and satisfies 
$\phi \circ \d_M = \d_N \circ \phi$.

Lastly, to clarify item (3) in the definition, a complex 
$M$ belongs to $\dcat{C}^{+}(A, \mrm{gr})$
if $M^i = 0$ for $i \ll 0$, etc. Indicators can be combined, so 
$\dcat{C}^{\star}_{\mrm{str}}(A, \mrm{gr})$
is a thick abelian subcategory of 
$\dcat{C}^{}_{\mrm{str}}(A, \mrm{gr})$.

Like for any abelian category, we have a homotopy category and a derived 
category. Here is the specific notation:

\begin{dfn} \label{dfn:3731} 
Let $A$ be a graded ring.
\begin{enumerate}
\item The homotopy category of complexes with entries in the abelian category 
$\dcat{M}(A, \mrm{gr})$ is the triangulated category%
\index{1-K(A,gr)@$\dcat{K}(A, \mrm{gr})$}
$\dcat{K}(A, \mrm{gr}) := \dcat{K}(\dcat{M}(A, \mrm{gr}))$.

\item The derived category of the abelian category $\dcat{M}(A, \mrm{gr})$ is 
the triangulated category%
\index{1-D(A,gr)@$\dcat{D}(A, \mrm{gr})$}
$\dcat{D}(A, \mrm{gr}) := \dcat{D}(\dcat{M}(A, \mrm{gr}))$.

\item Given a boundedness indicator $\star$, we denote by 
$\dcat{D}^{\star}(A, \mrm{gr})$ the full triangulated subcategory of 
$\dcat{D}(A, \mrm{gr})$ on the complexes $M$ whose cohomology $\opn{H}(M)$
has boundedness condition $\star$. 
\end{enumerate}
\end{dfn}

As always for derived categories, there are functors
\begin{equation} \label{eqn:3961}
\UseTips \xymatrix @C=6ex @R=6ex {
\dcat{C}^{}_{\mrm{str}}(A, \mrm{gr})
\ar[r]^{\opn{P}}
\ar@(ur,ul)[rr]^{\til{\opn{Q}}}
&
\dcat{K}(A, \mrm{gr})
\ar[r]^{\opn{Q}}
&
\dcat{D}(A, \mrm{gr}) \, . 
} 
\end{equation}
By abuse of notation we sometimes write $\opn{Q}$ rather than 
$\til{\opn{Q}}$; or, when there is no danger of confusion, we suppress the 
categorical localization functors $\opn{Q}$ and $\til{\opn{Q}}$ altogether. 

There is yet another category associated to $A$~: it is the category 
$\dcat{G}(A, \mrm{gr}) :=  \dcat{G}(\dcat{M}(A, \mrm{gr}))$
of cohomologically graded objects in $\dcat{M}(A, \mrm{gr})$. It has 
its strict subcategory $\dcat{G}_{\mrm{str}}(A, \mrm{gr})$. 
This is the target of the cohomology functor 
\[ \opn{H} : \dcat{C}_{\mrm{str}}(A, \mrm{gr})
\to \dcat{G}_{\mrm{str}}(A, \mrm{gr}) . \]

If $A \to B$ is a graded ring homomorphism that makes $B$ into a finite left 
$A$-module, and if $A$ is left noetherian, then of course $B$ is also left 
noetherian. Here is a partial converse: 

\begin{thm} \label{thm:4240}
Let $A$ be a nonnegative graded ring, and let $a \in A$ be a nonzero 
homogeneous central element of positive degree. Define the nonnegative graded 
ring $B := A / (a)$. If $B$ is left noetherian, then $A$ is also left 
noetherian. 
\end{thm}

This result is very similar to \cite[Lemma 8.2]{ATV}. There are 
three changes: we do not assume that the element $a$ is regular (i.e.\ not a 
zero-divisor), and we do not assume that $A$ is connected graded (only 
nonnegative). On the other hand, in \cite{ATV} the element $a$ is only 
required to be normalizing, not central. Our proof is essentially copied from 
\cite{ATV}. 

\begin{proof}
Assume, for the sake of contradiction, that the graded ring $A$ is not left 
noetherian. Then (taking Theorem \ref{thm:4245} into account) there is a graded 
left ideal $L \sub A$ that is not finitely generated. By Zorn's Lemma we can 
assume that $L$ is maximal with this property. 
Thus every graded left ideal $L^{\dag}$ such that 
$L \subneq L^{\dag} \sub A$ must be finitely generated.  

Consider the graded $A$-module $M := A / L$. Every graded $A$-module 
$M^{\dag} \sub M$ is of the form $M^{\dag} = L^{\dag} / L$, and therefore 
$M^{\dag}$ is a finitely generated $A$-module. 
In other words, $M$ is a noetherian object of $\dcat{M}(A, \mrm{gr})$.

Let $i := \opn{deg}(a)$. We examine the commutative diagram (\ref{4590}) in 
$\dcat{M}(A, \mrm{gr})$. 
Viewing the columns of this diagram as complexes, whose nonzero terms are in 
cohomological degrees $0$ and $1$, this is a short exact sequence in the 
abelian category $\dcat{C}_{\mrm{str}}(A, \mrm{gr})$. 
\begin{equation} \label{4590}
\UseTips \xymatrix @C=6ex @R=6ex {
0
\ar[r]
&
L
\ar[r]
&
A
\ar[r]
&
M
\ar[r]
&
0
\\
0
\ar[r]
&
L(-i)
\ar[u]^{a \cd (-)}
\ar[r]
&
A(-i)
\ar[u]^{a \cd (-)}
\ar[r]
&
M(-i)
\ar[u]^{a \cd (-)}
\ar[r]
&
0
}
\end{equation}

Define the graded $A$-modules $\bar{M}$ and $\bar{L}$ by these exact 
sequences:
\[ 0 \to \bar{M} \to M(-i) \xar{a \cd (-)} M \]
and 
\begin{equation} \label{eqn:4241}
L(-i) \xar{a \cd (-)} L \to \bar{L} \to 0 .
\end{equation}
The long exact cohomology sequence of (\ref{4590}) gives 
rise to this exact sequence in $\dcat{M}(A, \mrm{gr})$~:
\begin{equation} \label{eqn:4240}
\bar{M} \xar{\phi} \bar{L} \xar{\psi} B . 
\end{equation}
Because $M$ is a noetherian object of $\dcat{M}(A, \mrm{gr})$,
so is its degree twist $M(-i)$, and also the subobject $\bar{M} \sub M(-i)$. 
On the other hand, we are given that the ring  $B$ 
is noetherian; so it is a noetherian object of $\dcat{M}(A, \mrm{gr})$.
The exact sequence (\ref{eqn:4240}) shows that $\bar{L}$ is also a 
noetherian object of $\dcat{M}(A, \mrm{gr})$.
We see that $\bar{L}$ is a finitely generated $A$-module. 

Finally, let us choose finitely many homogeneous generators of the graded 
$A$-module $\bar{L}$, and then lift them to homogeneous elements of $L$. 
Let $L' \sub L$ be the $A$-submodule generated by this finite collection of 
elements. The exact sequence (\ref{eqn:4241}) shows that 
$L = (a \cd L) + L'$.
The graded $A$-module $L$ is bounded below, because $L \sub A$ and $A$ is 
nonnegative. The element $a$ belongs to the augmentation ideal $\m$ of $A$. 
Therefore, by Proposition \ref{prop:4540} (the graded Nakayama Lemma),
we have $L' = L$. This contradicts the assumption that $L$ is not finitely 
generated. 
\end{proof}

\mysubsection{Properties of Algebraically Graded Modules}
\label{subsec:props-ag-mods}

In this subsection we continue with Convention \ref{conv:3700}, to which we add:

\begin{conv}  \label{conv:4560}
The base ring $\K$ is a field. 
\end{conv}

See Remark \ref{rem:3740} regarding the requirement that $\K$ be a field. 
In a sense the category $\dcat{M}(\K, \mrm{gr})$ is boring under Convention
\ref{conv:4560}~:

\begin{prop} \label{prop:3851}
For an object $M \in \dcat{M}(\K, \mrm{gr})$ there is an isomorphism
$M \cong \bigoplus_{x \in X} \K(-i_x)$,
where $X$ is an indexing set and $\{ i_x \}_{x \in X}$ is a collection of 
integers. 
\end{prop}

\begin{proof}
This is a variation on the usual proof for $\K$-modules: by Zorn's Lemma there 
is a maximal linearly independent collection 
$\{ m_x \}_{x \in X}$ of homogeneous elements of $M$. 
This is a basis of $M$. The number $i_x$ is the degree of $m_x$.  
\end{proof}

\begin{dfn} \label{dfn:3795}
Let $A$ be a graded ring. 
We define the graded bimodule 
\[ A^* := \opn{Hom}_{\K}(A, \K) \in 
\dcat{M}(A^{\mrm{en}}, \mrm{gr}) , \]
see formula (\ref{eqn:3701}). 
\end{dfn}

\begin{dfn}  \label{dfn:3740}
Let $A$ be a graded ring.  
\begin{enumerate}
\item A graded $A$-module $P$ is called a {\em graded-free $A$-module}%
\index{Algebraically graded module! graded-free}
if
$P \cong \lb \bigoplus_{x \in X} A(-i_x)$
in $\dcat{M}(A, \mrm{gr})$
for some collection $\{ i_x \}_{x \in X}$ of integers.

\item A graded $A$-module $P$ is called a 
{\em graded-projective $A$-module}%
\index{Algebraically graded module! graded-projective}
if it is a projective object in the abelian category 
$\dcat{M}(A, \mrm{gr})$.

\item  A graded $A$-module $I$ is called a 
{\em graded-cofree $A$-module}%
\index{Algebraically graded module! graded-cofree}
if
$I \cong \lb \prod_{x \in X} A^*(-p_x)$
in $\dcat{M}(A, \mrm{gr})$ 
for some collection $\{ p_x \}_{x \in X}$ of integers.

\item A graded $A$-module $I$ is called a 
{\em graded-injective $A$-module}%
\index{Algebraically graded module! graded-injective}
if it is an injective object in the abelian category $\dcat{M}(A, \mrm{gr})$.

\item A graded $A$-module $P$ is called a {\em graded-flat $A$-module} 
\index{Algebraically graded module! graded-flat}
if the functor 
\[ (-) \ot_A P : \dcat{M}(A^{\mrm{op}}, \mrm{gr}) \to
\dcat{M}(\K, \mrm{gr}) \]
is exact. 
\end{enumerate}
\end{dfn}

Note that a graded-free $A$-module $P$ is of the form 
$P \cong A \ot_{} V$ for some $V \in \dcat{M}(\K, \mrm{gr})$;
and a graded-cofree $A$-module $I$ is of the form 
$I \cong \opn{Hom}_{\K}(A \ot W, \K)$ for some $W \in \dcat{M}(\K, \mrm{gr})$.

\begin{prop} \label{prop:3795}
Let $A$ and $B$ be graded rings.
\begin{enumerate}
\item There is an isomorphism 
\[ \opn{Hom}_{A}(-, A^*) \cong \opn{Hom}_{\K}(-, \K)  \]
of functors 
\[ \dcat{M}(A \ot B^{\mrm{op}}, \mrm{gr})^{\mrm{op}} \to 
\dcat{M}(B \ot A^{\mrm{op}}, \mrm{gr}) . \]

\item The graded bimodule $A^*$ is graded-injective over $A$, 
namely it is an injective object in the category
$\dcat{M}(A, \mrm{gr})$. 
\end{enumerate}
\end{prop}

\begin{proof}
This is by the adjunction isomorphism 
\[ \opn{Hom}_{A} \bigl( M, \opn{Hom}_{\K}(A, \K) \bigr) \cong
\opn{Hom}_{\K}(M, \K) \]
in $\dcat{M}(\K, \mrm{gr})$.
\end{proof}

The ungrading functor (\ref{eqn:3715}) preserves finiteness (i.e.\ 
being a finitely generated module). It also preserves freeness and 
projectivity, as the next proposition shows. 

\begin{prop} \label{prop:3715} 
Let $A$ be a graded ring, and let $P \in \dcat{M}(A, \mrm{gr})$. 
\begin{enumerate}
\item If $P$  is a graded-free $A$-module, then $\opn{Ungr}(P)$ is a 
free $\opn{Ungr}(A)$-module.

\item $P$ is a graded-projective $A$-module if and only if it is a direct 
summand, in $\dcat{M}(A, \mrm{gr})$, of a graded-free $A$-module.

\item If $P$ is a graded-projective $A$-module, then $\opn{Ungr}(P)$ is a 
projective $\opn{Ungr}(A)$-module.
\end{enumerate} 
\end{prop}

\begin{exer} \label{exer:3715} 
Prove Proposition \ref{prop:3715}.
\end{exer}

\begin{rem} \label{rem:4805}
In the situation of Proposition \ref{prop:3715}, one can show that $P$ is a 
graded-flat $A$-module if and only if  $\opn{Ungr}(P)$ is a flat 
$\opn{Ungr}(A)$-module. The proof is not easy. See 
\cite[Proposition 3.4.5]{NV}.
\end{rem}

For injectives things are much more complicated, as the next theorem and 
example show.

\begin{thm}[Van den Bergh] \label{thm:3715} 
Let $A$ be a left noetherian nonnegative graded ring, and let 
$I$ be graded-injective $A$-module. Then $\opn{Ungr}(I)$ has injective 
dimension at most $1$ as an $\opn{Ungr}(A)$-module.
\end{thm}

\begin{proof}
See \cite{Ye17}. A more general statement can be found in \cite{SoZa}.
\end{proof}

\begin{exa} \label{exa:3800}
Consider the graded ring $A := \K[t]$, the ring of polynomials in a variable 
$t$ of degree $1$. The graded $A$-module $I := \K[t, t^{-1}]$ is 
graded-injective, but $\opn{Ungr}(I)$ has injective dimension $1$ over 
$\opn{Ungr}(A)$. 
\end{exa}

\begin{prop} \label{prop:4905}  
\mbox{}
\begin{enumerate}
\item A graded-free $A$-module $P$ is graded-projective. 

\item A graded-projective $A$-module $P$ is graded-flat.

\item A graded-cofree $A$-module $I$ is graded-injective. 
\end{enumerate}
\end{prop}

\begin{exer} \label{exer:4045}
Prove this proposition. 
\end{exer}

\begin{prop} \label{prop:3900}
Let $A$ and $B$ be graded rings.
\begin{enumerate}
\item If $P \in \dcat{M}(A \ot B, \mrm{gr})$ is projective, then it is 
projective in $\dcat{M}(A, \mrm{gr})$.

\item If $P \in \dcat{M}(A \ot B, \mrm{gr})$ is flat, then it is 
flat in $\dcat{M}(A, \mrm{gr})$.

\item If $I \in \dcat{M}(A \ot B, \mrm{gr})$ is injective, then it is 
injective in $\dcat{M}(A, \mrm{gr})$.
\end{enumerate}
\end{prop}

\begin{proof}
(1) Take $N \in \dcat{C}(A, \mrm{gr})$. 
We have isomorphisms
\[ \begin{aligned}
& \opn{Hom}_{A}(P, N) \cong
\opn{Hom}_{A \ot B}(P, \opn{Hom}_A(A \ot B, N)) 
\\ & \quad  \cong \opn{Hom}_{A \ot B}(P, \opn{Hom}_{\K}(B, N)) 
\end{aligned} \]
in $\dcat{C}_{\mrm{str}}(\K, \mrm{gr})$. 
As an object of $\dcat{M}(\K, \mrm{gr})$, $B$ is projective.
Hence, if $N$ is an acyclic complex, then so is $\opn{Hom}_{A}(P, N)$.

\medskip \noindent
(2, 3) Like item (1), but using the fact that $B$ is graded-flat over $\K$, and 
using suitable adjunction isomorphisms.
\end{proof}

\begin{dfn} \label{dfn:4260}
Let $A$ be a graded ring.
\begin{enumerate}
\item Let $M$ be a graded $A$-module, and let $N \sub M$ 
be a graded submodule. We say that $N$ is a 
{\em graded-essential submodule}%
\index{Algebraically graded module! graded-essential submodule of}
of $M$ if for every nonzero graded submodule $M' \sub M$ the intersection 
$N \cap M'$ is nonzero. 

\item A homomorphism $\phi : N \to M$ in $\dcat{M}(A, \mrm{gr})$
is called a {\em graded-essential monomorphism}%
\index{Algebraically graded module! graded-essential monomorphism of {\indash}s}
if $\phi$ is an injective 
homomorphism, and $\phi(N)$ is a graded-essential submodule of $M$.
\end{enumerate}
\end{dfn}

Item (2) of the definition is a special case of the usual definition of an 
essential monomorphism in an abelian category. 

\begin{exer} \label{exer}
Let $N \sub M$ be graded $A$-modules. Show that $N$ is a graded-essential 
submodule of $M$ if and only if $\opn{Ungr}(N)$ is an essential 
$\opn{Ungr}(A)$-submodule of $\opn{Ungr}(M)$.
(Hint: prove that for every nonzero element $m \in M$ there is a homogeneous 
element $a \in A$ such that $a \cd m \in N$ and $a \cd m \neq 0$. Do this by 
induction on the number $r$ appearing in the decomposition (\ref{eqn:4260}).
See the proof of \cite[Lemma 3.3.13]{NV}, which is not totally correct, but 
can be fixed using our hint.)
\end{exer}

\begin{prop} \label{prop:3735}
Let $A$ be a graded ring. 
\begin{enumerate}
\item The category $\dcat{M}(A, \mrm{gr})$ has enough projectives.

\item The category $\dcat{M}(A, \mrm{gr})$ has enough injectives. Moreover, 
every $M \in \lb \dcat{M}(A, \mrm{gr})$ has an injective hull, namely there is 
a graded-essential mono\-morphism $M \inj I$ to a graded-injective module $I$. 
\end{enumerate}
\end{prop}

\begin{proof}
(1) Given a module $M \in \dcat{M}(A, \mrm{gr})$ and a homogeneous element 
$m \in M_i$, there is a homomorphism 
$A(-i) \to M$ that sends $1_A \mapsto m$. Therefore, by taking a sufficiently 
large direct sum of degree twists of $A$, there is a surjection 
$P \surj M$ from a graded-free module $P$. And graded-free modules are 
graded-projective, by Proposition \ref{prop:4905}. 

\medskip \noindent
(2) For every nonzero homogeneous element $m \in M_q$ there is a homomorphism 
$M \to \K(-q)$ in $\dcat{M}(\K, \mrm{gr})$
that is nonzero on $m$. By Proposition \ref{prop:3795}(1)
we get a homomorphism $M \to A^*(-q)$ in $\dcat{M}(A, \mrm{gr})$
that is nonzero on $m$. Thus, by taking a sufficiently large product 
$I :=  \prod_{x \in X} A^*(-q_x)$
in $\dcat{M}(A, \mrm{gr})$ we obtain a monomorphism $M \inj I$. And according 
Proposition \ref{prop:4905},  the graded-cofree object $I$ is injective in 
$\dcat{M}(A, \mrm{gr})$.

As for injective hulls: this is the same as in the ungraded case; see 
\cite[Theorem 3.30]{Rot} or \cite[Section 3.D]{Lam}.
\end{proof}

Now a structural result on injective objects in $\dcat{M}(A, \mrm{gr})$. 

\begin{thm} \label{thm:3795}
Let $A$ be a left noetherian graded ring. 
\begin{enumerate}
\item If $\{ I_x \}_{x \in X}$ is a collection of graded-injective $A$-modules, 
then $I := \bigoplus_{x \in X} I_x$ is a  graded-injective $A$-module. 

\item Every graded-injective $A$-module $I$ is a direct sum of indecomposable
graded-injective $A$-modules. 
\end{enumerate} 
\end{thm}

\begin{proof}
This is because $\dcat{M}(A, \mrm{gr})$ is a locally noetherian Grothendieck 
abelian category. See \cite[Section V.4]{Ste}.
Item (1) can be proved directly using the graded version of the Baer criterion 
(Theorem \ref{thm:151}). 
\end{proof}

\begin{dfn}  \label{dfn:3703}
A {\em connected graded $\K$-ring} 
\index{Algebraically graded ring! connected}
is a nonnegative graded $\K$-ring 
$A$ (see Conventions \ref{conv:3700} and 
\ref{conv:4560},  and Definition \ref{dfn:4540}), such that  the ring 
homomorphism $\K \to A_0$ is bijective, 
and $A_i$ is a finitely generated $\K$-module for all $i \geq 0$. 
\end{dfn}

Here are a few examples of connected graded rings. 

\begin{exa} \label{exa:3755}
Consider the noncommutative polynomial ring 
$A := \lb \K \bra{x_1, \ldots, x_n}$ 
in $n \geq 1$ variables, all of algebraic degree $1$. This is a connected 
graded $\K$-ring. If $n \geq 2$ then $A$ is not commutative and not noetherian 
(on either side). If $n = 1$ then 
$\K \bra{x_1} = \K[x_1]$, see next example. 
\end{exa}

\begin{exa} \label{exa:3757} 
Let
$A := \K [x_1, \ldots, x_n]$,
the commutative polynomial ring 
in $n \geq 1$ variables, all of algebraic degree $1$. This is a connected 
graded $\K$-ring, commutative and noetherian.  
\end{exa}

\begin{exa} \label{exa:3756} 
Let $\K[t]$ be the commutative polynomial ring in a variable $t$ of algebraic 
degree $1$. Next let 
\[ A := \K[t] \bra{x, y} / (y \cd x - x \cd y - t^2) , \]
where $x$ and $y$ are noncommuting variables of algebraic 
degree $1$ (that commute with $t$). This is the {\em homogeneous first Weyl 
algebra}. It is a connected graded $\K$-ring, noetherian, but not commutative. 
Indeed, if $\opn{char}(\K) = 0$, then the center of $A$ is $\K[t]$. 
\end{exa}

\begin{exa} \label{exa:4075} 
Let $\g$ be a finite Lie algebra over $\K$. Choose a basis 
$v_1, \ldots, v_n$ for $\g$ as a $\K$-module, and let 
$\ga_{i, j, k} \in \K$ be the constants describing the Lie bracket, i.e.\
\[ [v_i, v_j] = \sum_k \ga_{i, j, k} \cd v_k \in \g . \]
Let $\K[t]$ be the commutative polynomial ring in a variable $t$ of algebraic 
degree $1$. Let $x_1, \ldots, x_n$ be variables of algebraic degree $1$, and 
define the graded ring 
\[ A := \K[t] \bra{x_1, \ldots, x_n} \, / \, (I) , \]
where $I$ is the two-sided ideal generated by the quadratic elements 
\[ x_i \cd x_j - x_j \cd x_i  - \sum_{k = 1}^n \ga_{i, j, k} \cd x_k \cd t . \]
This is the {\em homogeneous universal enveloping ring} of $\g$. 
The quotient $A / (t - 1)$ is the  universal enveloping ring
$\opn{U}(\g)$. The ring $A$ is connected graded, noetherian, but often 
noncommutative (since the Lie algebra $\g$ embeds into the quotient ring 
$\opn{U}(\g)$ with its commutator Lie bracket). 
\end{exa}

\begin{exa} \label{exa:3775}
Suppose $A$ and $B$ are both connected graded $\K$-rings. Then 
$A \ot B$ is a connected graded $\K$-ring.
\end{exa}

\begin{dfn} \label{dfn:3724}
A graded $\K$-module $M = \bigoplus_{i \in \Z} \, M_i$ 
is called {\em degreewise finite} 
\index{Algebraically graded module! degreewise-finite}
if the $\K$-module $M_i$ is finitely generated for every $i \in \Z$. 

The full subcategory of $\dcat{M}(\K, \mrm{gr})$ on the 
degreewise finite modules is denoted by 
$\dcat{M}_{\mrm{dwf}}(\K, \mrm{gr})$. 
\end{dfn}

It is obvious that the category $\dcat{M}_{\mrm{dwf}}(\K, \mrm{gr})$ is a thick 
abelian subcategory of $\dcat{M}(\K, \mrm{gr})$, 
closed under subobjects and quotient objects.

\begin{dfn} \label{dfn:3725}
For $M \in \dcat{M}(\K, \mrm{gr})$ we let
\[ M^* := \opn{Hom}_{\K}(M, \K) \in \dcat{M}(\K, \mrm{gr}) . \]
This is the {\em graded $\K$-linear dual of $M$}. 
\index{Algebraically graded module! graded $\K$-linear dual}
\index{1-Mstar@$M^*$}
\end{dfn}

According to formula (\ref{eqn:3701}) we have
$M^* = \bigoplus_{i \in \Z} (M^*)_i$, where 
\[ (M^*)_i = \opn{Hom}_{\K}(M, \K)_i = \opn{Hom}_{\K}(M_{-i}, \K) = 
(M_{-i})^* . \]

\begin{lem} \label{lem:3725}
If $M \in \dcat{M}_{\mrm{dwf}}(\K, \mrm{gr})$ then 
$M^* \in \dcat{M}_{\mrm{dwf}}(\K, \mrm{gr})$, 
and the Hom-evaluation homomorphism
$\opn{ev}_M : M \to M^{**} = (M^*)^*$ 
in $\dcat{M}_{\mrm{dwf}}(\K, \mrm{gr})$ is bijective. 
Thus $(-)^*$ is a duality of the category 
$\dcat{M}_{\mrm{dwf}}(\K, \mrm{gr})$.
\end{lem}

The easy proof is omitted. 

\begin{dfn} \label{dfn:3723} 
Let $A$ be a graded ring. 
A module $M \in \dcat{M}(A, \mrm{gr})$ is called {\em degreewise $\K$-finite} 
\index{Algebraically graded module! degreewise-finite}
if it is degreewise finite as a graded $\K$-module. 

The full subcategory of $\dcat{M}(A, \mrm{gr})$ on the degreewise $\K$-finite 
modules is denoted by 
$\dcat{M}_{\mrm{dwf}}(A, \mrm{gr})$.
\index{1-Mdwf(A,gr)@$\dcat{M}_{\mrm{dwf}}(A, \mrm{gr})$}
\end{dfn}

Of course $\dcat{M}_{\mrm{dwf}}(A, \mrm{gr})$ is a thick abelian 
subcategory of $\dcat{M}(A, \mrm{gr})$,
closed under subobjects and quotients.

\begin{prop} \label{prop:4045} 
Let $A$ be a graded ring.  
If $M \in \dcat{M}_{\mrm{dwf}}(A, \mrm{gr})$ then 
$M^* \in \dcat{M}_{\mrm{dwf}}(A^{\mrm{op}}, \mrm{gr})$, 
and the Hom-evaluation homomorphism 
$\opn{ev}_M : M \to M^{**}$
in \lb $\dcat{M}(A, \mrm{gr})$ is bijective. Thus 
\[ (-)^* : \dcat{M}_{\mrm{dwf}}(A, \mrm{gr})^{\mrm{op}} \to 
\dcat{M}_{\mrm{dwf}}(A^{\mrm{op}}, \mrm{gr}) \]
is an equivalence of categories.
\end{prop}

\begin{exer} \label{exer:4075}
Prove this proposition. (Hint: use Lemma \ref{lem:3725}.)
\end{exer}

Noncommutative connected graded rings are very similar to complete commutative 
local rings. Here is the second variant of the graded Nakayama Lemma (the first 
was Proposition \ref{prop:4540}). 

\begin{prop} \label{prop:4590}    
Let $A$ be a connected graded ring, and let 
$M \in \dcat{M}(A, \mrm{gr})$ be a bounded below graded module. Define  
$V := \K \ot_A M \in \dcat{M}(\K, \mrm{gr})$.
\begin{enumerate}
\item There is a surjection 
$\pi : A \ot_{} V \surj M$
in $\dcat{M}(A, \mrm{gr})$, such that the diagram 
\[ \UseTips \xymatrix @C=8ex @R=6ex {
A \ot V
\ar[r]^(0.55){\pi}
\ar[d]_{\be}
&
M
\ar[d]^{\al}
\\
V
\ar[r]^{\opn{id}}
&
V
} \]
in which the morphisms $\al$ and $\be$ are the canonical projections, is 
commutative.  

\item If $V = 0$ then $M = 0$. 

\item If $M$ is a graded-projective $A$-module, then the surjection 
$\pi$ above is bijective. Hence $M$ is a graded-free $A$-module. 
\end{enumerate}
\end{prop}

Note that the projection $\be$ is induced by the augmentation homomorphism 
$A \to \K$, so $\opn{Ker}(\be) = \m \ot V$. Likewise, 
$\opn{Ker}(\al) = \m \cd M \sub M$.

\begin{proof}  \mbox{}

\smallskip \noindent
(1) Say $M = \bigoplus_{i \geq i_0} \, M_i$. 
Choose a section $\pi_0 : V \inj M$ in $ \dcat{M}(\K, \mrm{gr})$ of the 
canonical projection $\al : M \surj V$. 
Since $A \ot V$ is a graded-free $A$-module, $\pi_0$ extends to a homomorphism 
$\pi : A \ot_{} V \to M$ in $\dcat{M}(A, \mrm{gr})$.
See the commutative diagram below, in which $\ga(v) := 1_A \ot v$. 
\[ \UseTips \xymatrix @C=8ex @R=6ex {
V
\ar[dr]^{\pi_0}
\ar[d]_{\ga}
\\
A \ot V
\ar[r]^{\pi}
&
M
} \]
By Proposition \ref{prop:4540} the homomorphism $\pi$ is surjective. 

\medskip \noindent 
(2) Clear from the existence of the surjection $\pi$.

\medskip \noindent 
(3) Let $N := \opn{Ker}(\pi) \sub A \ot V$, so there is a short exact sequence 
\[ 0 \to N \to A \ot V \xar{\pi} M \to  0 \]
in $\dcat{M}(A, \mrm{gr})$, and $N$ is a bounded below graded $A$-module. 
If $M$ is graded-projective, then this sequence is split. 
It follows that $\K \ot_{A} N = 0$. By part (2) we see that $N = 0$. 
\end{proof}

\begin{dfn} \label{dfn:4265}
Let $A$ be a connected graded $\K$-ring, and let $M \in \dcat{M}(A, \mrm{gr})$. 
The {\em socle}%
\index{Algebraically graded module! socle of}%
\index{1-Soc(M)@$\opn{Soc}(M)$}
of $M$ is the graded $A$-submodule
$\opn{Soc}(M) := \opn{Hom}_A(\K, M) \sub M$.
\end{dfn}

In other words, 
$\opn{Soc}(M) = \{ m \in M \mid \m \cd m = 0 \}$.
The $A$-module structure of $\opn{Soc}(M)$ is through the augmentation 
homomorphism $A \to \K$; so Proposition \ref{prop:3851} applies to it. 

Here is a dual to the graded Nakayama Lemma. 

\begin{prop} \label{prop:4595}  
Let $A$ be a connected graded $\K$-ring, let 
$N \in \dcat{M}(A, \mrm{gr})$ be a bounded above graded module,
and let $W := \opn{Soc}(N)$. 
\begin{enumerate}
\item There is an injection
$\si : N \inj \opn{Hom}_{\K}(A, W)$
in $\dcat{M}(A, \mrm{gr})$,  such that the diagram 
\[ \UseTips \xymatrix @C=8ex @R=6ex {
W
\ar[r]^{\opn{id}}
\ar[d]_{\al}
&
W
\ar[d]^{\be}
\\
N
\ar[r]^(0.35){\si}
&
\opn{Hom}_{\K}(A, W)
} \]
in which $\al$ and $\be$ are the canonical embeddings, is commutative. 

\item If $W = 0$ then $N = 0$.

\item If $N$ is a graded-injective $A$-module, then the homomorphism $\si$ 
above is bijective. Hence $N$ is a graded-cofree $A$-module.
\end{enumerate}
\end{prop}

Note that $\be$ is induced from the augmentation homomorphism 
$A \to \K$, and 
$\opn{Im}(\be) = \opn{Soc} \bigl( \opn{Hom}_{\K}(A, W) \bigr)$.

\begin{proof} \mbox{} 

\smallskip \noindent 
(1) Say $N = \bigoplus_{i \leq i_1} \, N_i$ for some $i_1 \in \Z$.
Choose a splitting 
$\si_0 : N \surj W$ in $\dcat{M}(\K, \mrm{gr})$ of the inclusion 
$\al : W \inj N$. 
By adjunction there is an isomorphism 
\[ \opn{Hom}_{\K}(N, W) \cong 
\opn{Hom}_{A}(N,  \opn{Hom}_{\K}(A, W)) \]
in $\dcat{M}(\K, \mrm{gr})$, so there is a homomorphism 
$\si : N \to \opn{Hom}_{\K}(A, W)$ in $\dcat{M}(A, \mrm{gr})$ that lifts 
$\si_0$. See the commutative diagram below, where 
$\ga(\chi) := \chi(1_A)$.  
\[ \UseTips \xymatrix @C=8ex @R=6ex {
N
\ar[dr]_{\si_0}
\ar[r]^(0.35){\si}
&
\opn{Hom}_{\K}(A, W)
\ar[d]^{\ga}
\\
&
W
} \]
A descending induction on degree, starting from $i_1$, shows that 
$\si$ is an injection. 

\medskip \noindent 
(2) Clear from the existence of $\si$. 

\medskip \noindent 
(3) Let 
$M := \opn{Coker}(\si)$. 
If $N$ is graded-injective, then the exact sequence 
\[ 0 \to N \xar{\si} \opn{Hom}_{\K}(A, W) \to M \to 0 \]
in $\dcat{M}(A, \mrm{gr})$ is split. 
Applying the additive functor $\opn{Soc}$ to it, we see that 
$\opn{Soc}(M) = 0$. Since $M$ is a bounded above graded module, 
item (2) says that $M = 0$. 
\end{proof}

\begin{dfn} \label{dfn:3726} 
Let $A$ be a graded ring. 
A module $N \in \dcat{M}(A^{\mrm{op}}, \mrm{gr})$ is called a 
{\em cofinite graded $A^{\mrm{op}}$-module} 
\index{Algebraically graded module! cofinite}
if $N \cong M^*$ for some finite graded 
$A$-module $M$.

The full subcategory of $\dcat{M}(A^{\mrm{op}}, \mrm{gr})$ on the cofinite 
graded $A^{\mrm{op}}$-modules is denoted by 
\index{1-Mcof(A,gr)@$\dcat{M}_{\mrm{cof}}(A, \mrm{gr})$}
$\dcat{M}_{\mrm{cof}}(A^{\mrm{op}}, \mrm{gr})$.
\end{dfn}

Recall that an object $M$ in an abelian category $\cat{M}$ is called {\em 
noetherian} (resp.\ {\em artinian}) if it satisfies the ascending (resp.\ 
descending) chain condition on subobjects. 
Let us denote the corresponding full subcategories of $\cat{M}$ by
$\cat{M}_{\mrm{not}}$ and $\cat{M}_{\mrm{art}}$ respectively. 
Then $\cat{M}_{\mrm{not}}, \cat{M}_{\mrm{art}} \sub \cat{M}$
are thick abelian subcategories, closed under subobjects and 
quotient objects. 

Here is a variant of the Matlis Duality Theorem (see Remark \ref{rem:2245}):

\begin{thm}[NC Graded Matlis Duality] \label{thm:4045} 
\index{Matlis Duality Theorem! noncommutative graded} 
Assume $A$ is a left noetherian connected graded ring. 
\begin{enumerate}
\item The category $\dcat{M}_{\mrm{f}}(A, \mrm{gr})$ is a thick abelian 
subcategory of $\dcat{M}_{\mrm{dwf}}(A, \mrm{gr})$, 
closed under subobjects and quotient objects. 

\item The category $\dcat{M}_{\mrm{cof}}(A^{\mrm{op}}, \mrm{gr})$
is a thick abelian subcategory of 
$\dcat{M}_{\mrm{dwf}}(A^{\mrm{op}}, \mrm{gr})$, 
closed under subobjects and quotient objects. 

\item The functor 
\[ (-)^* : \dcat{M}_{\mrm{f}}(A, \mrm{gr})^{\mrm{op}} \to 
\dcat{M}_{\mrm{cof}}(A^{\mrm{op}}, \mrm{gr}) \]
is an equivalence. 

\item The modules in $\dcat{M}_{\mrm{f}}(A, \mrm{gr})$ are the noetherian 
objects in the abelian category $\dcat{M}(A, \mrm{gr})$.

\item The modules in $\dcat{M}_{\mrm{cof}}(A^{\mrm{op}}, \mrm{gr})$ are the 
artinian objects in the abelian category $\dcat{M}(A^{\mrm{op}}, \mrm{gr})$.
\end{enumerate}
\end{thm}

\begin{proof} \mbox{}

\smallskip \noindent
(1) Since $A$ is connected, as a left module it belongs to 
$\dcat{M}_{\mrm{dwf}}(A, \mrm{gr})$. Therefore every finite graded $A$-module 
belongs to $\dcat{M}_{\mrm{dwf}}(A, \mrm{gr})$, i.e.\
$\dcat{M}_{\mrm{f}}(A, \mrm{gr}) \sub \dcat{M}_{\mrm{dwf}}(A, \mrm{gr})$.

Because $A$ is left noetherian, it follows that 
$\dcat{M}_{\mrm{f}}(A, \mrm{gr})$
is a thick abelian subcategory of
$\dcat{M}(A, \mrm{gr})$, 
closed under subobjects and quotient objects.

\medskip \noindent 
(2, 3) By definition 
$\dcat{M}_{\mrm{cof}}(A^{\mrm{op}}, \mrm{gr})$
is the essential image under the functor $(-)^*$ of the category
$\dcat{M}_{\mrm{f}}(A, \mrm{gr})$.
Now use the equivalence of Proposition \ref{prop:4045} and part (1). 

\medskip \noindent 
(4) This is proved the same way as for ungraded left noetherian rings. 

\medskip \noindent 
(5) We know by (2, 3, 4) that the objects of 
$\dcat{M}_{\mrm{cof}}(A^{\mrm{op}}, \mrm{gr})$
are artinian objects in 
$\dcat{M}(A^{\mrm{op}}, \mrm{gr})$.

For the opposite direction, let 
$N \in \dcat{M}(A^{\mrm{op}}, \mrm{gr})$ be an artinian object. 
The descending chain condition forces $N$ to be bounded above, namely 
$N_i = 0$ for $i \gg 0$. 
Consider the socle
$W := \opn{Hom}_{A^{\mrm{op}}}(\K, N)$. 
This is a graded $A^{\mrm{op}}$-submodule of $N$, isomorphic to 
$\bigoplus_{x \in X}  \K(-i_x)$ for some indexing set $X$ and a collection of 
integers $\{ i_x \}_{x \in X}$. By the descending chain condition the set $X$ 
has to be finite. Hence 
\[ \bigoplus_{x \in X}  A^*(-i_x) \cong A^* \ot W \cong 
\opn{Hom}_{\K}(A, W) . \]
Proposition \ref{prop:4595} tells us that there is an injective 
homomorphism 
$N \inj \lb \bigoplus_{x \in X} A^*(-i_x)$.
We know that $A^*(-i_x) \in \dcat{M}_{\mrm{cof}}(A, \mrm{gr})$; 
and hence so is $N$. 
\end{proof}

\begin{rem} \label{rem:3903}
If $A$ is connected graded but not left noetherian, then there are cofinite 
graded $A^{\mrm{op}}$-modules that are not artinian. Take the ring 
$A := \K \bra{x_1, x_2}$ from Example \ref{exa:3755}. Then 
$A^* \in \dcat{M}(A^{\mrm{op}}, \mrm{gr})$
is cofinite but not artinian. 
\end{rem}

\mysubsection{Resolutions and Derived Functors} 
\label{gr-res-der-fun}

In this subsection we continue with Conventions \ref{conv:3700} and 
\ref{conv:4560}. Thus $\K$ is a base field and $\ot = \ot_{\K}$. 
In the next definitions we collect the algebraically graded versions of some 
definitions from Sections \ref{sec:resol} and \ref{sec:exist-resol}.

\begin{dfn} \label{dfn:4275}   
Let $A$ be a graded ring, and let $P \in \dcat{C}(A, \mrm{gr})$. 
\begin{enumerate}
\item  The complex $P$ is called a
{\em graded-free complex}%
\index{Complex of algebraically graded modules! graded-free}
if all the graded $A$-modules $P^i$ are 
graded-free, and the differential is zero. 

\item The complex $P$ is called a
{\em semi-graded-free complex}%
\index{Complex of algebraically graded modules! semi-graded-free}
if it admits a filtration 
$\{ F_j(P) \}_{j \geq -1}$ in $\dcat{C}_{\mrm{str}}(A, \mrm{gr})$, 
such that $F_{-1}(P) = 0$, 
$P = \bigcup_j F_j(P)$, and 
$\opn{Gr}^F_j(P)$ is a graded-free complex for every $j$. 

\item The complex $P$ is called a 
{\em K-graded-projective complex}%
\index{Complex of algebraically graded modules! K-graded-projective}
if for every acyclic complex $N \in \dcat{C}(A, \mrm{gr})$, the complex 
$\opn{Hom}_A(P, N)$ is acyclic.
\end{enumerate}
\end{dfn}

\begin{dfn} \label{dfn:4276}   
Let $A$ be a graded ring, and let $I \in \dcat{C}(A, \mrm{gr})$. 
\begin{enumerate}
\item  The complex $I$ is called a 
{\em graded-cofree complex}%
\index{Complex of algebraically graded modules! graded-cofree}
if all the $A$-modules $I^p$ are graded-cofree, and the differential is zero.

\item The complex $I$ is called a 
{\em semi-graded-cofree complex}%
\index{Complex of algebraically graded modules! semi-graded-cofree}
if it admits a cofiltration $\{ F_q(I) \}_{q \geq -1}$ in 
$\dcat{C}_{\mrm{str}}(A, \mrm{gr})$, such that $F_{-1}(I) = 0$, each 
$\opn{Gr}^F_q(I)$ is graded-cofree complex, and 
$I = \opn{lim}_{\lar q} F_q(I)$.

\item The complex $I$ is called a 
{\em K-graded-injective complex}%
\index{Complex of algebraically graded modules! K-graded-injective}
if for every acyclic complex $N \in \dcat{C}(A, \mrm{gr})$, the complex
$\opn{Hom}_A(N, I)$ is acyclic. 
\end{enumerate}
\end{dfn}

\begin{dfn} \label{dfn:4277} 
A complex $P \in \dcat{C}(A, \mrm{gr})$ is called a 
{\em K-graded-flat complex}%
\index{Complex of algebraically graded modules! K-graded-flat}
if for every acyclic complex 
$N \in \dcat{C}(A^{\mrm{op}}, \mrm{gr})$ the complex 
$N \ot_A P$ is acyclic. 
\end{dfn}

\begin{thm} \label{thm:3740}
Let $A$ be a graded ring, and let $M \in \dcat{C}(A, \mrm{gr})$.
\begin{enumerate}
\item There exists a quasi-isomorphism 
$\rho : P \to M$ in $\dcat{C}_{\mrm{str}}(A, \mrm{gr})$,
where $P$ is a semi-graded-free complex, 
and $\opn{sup}(P) = \opn{sup}(\opn{H}(M))$. 

\item  There exists a quasi-isomorphism 
$\rho : M \to I$ in $\dcat{C}_{\mrm{str}}(A, \mrm{gr})$,
where $I$ is a semi-graded-cofree complex, and
$\opn{inf}(I) = \opn{inf}(\opn{H}(M))$. 
\end{enumerate}
\end{thm}

\begin{proof} \mbox{}

\smallskip \noindent 
(1) This is a modification of Theorem \ref{thm:3305} and Corollary 
\ref{cor:3340}. The proof is the same, after a few obvious modifications. 

\medskip \noindent 
(2) This is a modification of Theorem \ref{thm:3325} and Corollary 
\ref{cor:3345}. The proof is the same, after a few obvious modifications. 
The injective object $\K^*$ we use here is $\K^* := \K$ of course. 
\end{proof}

The various kinds of complexes mentioned in the definitions above are related 
as follows. 

\begin{thm} \label{thm:4275}  
Let $A$ be a graded ring.
\begin{enumerate}
\item A complex $P \in \dcat{C}(A, \mrm{gr})$ is K-graded-projective if and 
only if it is 
isomorphic in $\dcat{K}(A, \mrm{gr})$ to a semi-graded-free complex $Q$.

\item If $P \in \dcat{C}(A, \mrm{gr})$ is K-graded-projective, then it is 
K-graded-flat. 

\item If $P \in \dcat{C}(A, \mrm{gr})$ is K-graded-projective, then 
$P^* \in \dcat{C}(A^{\mrm{op}}, \mrm{gr})$
is K-graded-injective. 

\item A complex $I \in \dcat{C}(A, \mrm{gr})$ is K-graded-injective iff it is 
isomorphic in $\dcat{K}(A, \mrm{gr})$ to a semi-graded-cofree complex $J$.
\end{enumerate}
\end{thm}

\begin{proof} \mbox{}

\smallskip \noindent
(1) By Theorem \ref{thm:1575}, that applies to this context too, a 
semi-graded-free complex $Q$ is K-graded-projective. And by Proposition 
\ref{prop:1520}, a complex $P \in \dcat{C}(A, \mrm{gr})$ is 
K-graded-projective iff $\opn{Hom}_{\dcat{K}(A, \mrm{gr})}(P, N) = 0$ for every 
acyclic complex $N \in \dcat{C}(A, \mrm{gr})$.

First let's assume that there is an isomorphism $P \cong Q$ in 
$\dcat{K}(A, \mrm{gr})$ to a semi-graded-free complex $Q$.
Then for every acyclic complex 
$N \in \dcat{C}(A, \mrm{gr})$ we get 
\[ \opn{Hom}_{\dcat{K}(A, \mrm{gr})}(P, N) \cong 
\opn{Hom}_{\dcat{K}(A, \mrm{gr})}(Q, N) = 0 . \]
Therefore $P$ is K-graded-projective.

Conversely, let's assume that $P$ is K-graded-projective. By Theorem 
\ref{thm:3740}(1) there is a quasi-isomorphism
$\rho : Q \to P$ in $\dcat{C}_{\mrm{str}}(A, \mrm{gr})$ from a 
semi-graded-free complex $Q$. By Corollary \ref{cor:1923}, 
$\opn{P}(\rho) : Q \to P$ is an 
isomorphism in $\dcat{K}(A, \mrm{gr})$.

\medskip \noindent
(2) Take $N \in \dcat{C}(A^{\mrm{op}}, \mrm{gr})$. 
We have the adjunction isomorphism
\[ \begin{aligned}
& (N \ot_A P)^* = \opn{Hom}_{\K}(N \ot_A P, \K) 
\\
& \quad \cong
\opn{Hom}_{A}(P, \opn{Hom}_{\K}(N, \K)) = \opn{Hom}_{A}(P, N^*) 
\end{aligned} \]
in $\dcat{C}_{\mrm{str}}(\K, \mrm{gr})$. 
If $N$ is acyclic then so are $N^*$, 
$\opn{Hom}_{A}(P, N^*)$ and $(N \ot_A P)^*$.
Since the functor $(-)^*$ is faithfully exact, we see that 
$N \ot_A P$ is acyclic too. Therefore $P$ is K-graded-flat.

\medskip \noindent
(3) Take $N \in \dcat{C}(A^{\mrm{op}}, \mrm{gr})$. 
We have the adjunction isomorphism
\[ \begin{aligned}
& \opn{Hom}_{A^{\mrm{op}}}(N, P^*) =
\opn{Hom}_{A^{\mrm{op}}}(N, \opn{Hom}_{\K}(P, \K))
\\
& \quad \cong \opn{Hom}_{\K}(N \ot_A P, \K) = (N \ot_A P)^* 
\end{aligned} \]
in $\dcat{C}_{\mrm{str}}(\K, \mrm{gr})$. 
By item (2) $P$ is K-graded-flat. If $N$ is acyclic then so are
$N \ot_A P$, $(N \ot_A P)^*$ and $\opn{Hom}_{A^{\mrm{op}}}(N, P^*)$. 
Therefore $P^*$ is K-graded-injective.

\medskip \noindent
(4) By Theorem \ref{thm:1665}, that applies to this context too, a 
semi-graded-cofree complex $J$ is K-graded-injective. And by Proposition 
\ref{prop:1515}, a complex 
$I \in \dcat{C}(A, \mrm{gr})$ is K-graded-injective iff
$\opn{Hom}_{\dcat{K}(A, \mrm{gr})}(N, I) = 0$ for every acyclic complex 
$N \in \dcat{C}(A, \mrm{gr})$.

First let's assume that there is an isomorphism $I \cong J$ in 
$\dcat{K}(A, \mrm{gr})$ to a semi-graded-cofree complex $J$.
Then for every acyclic complex 
$N \in \dcat{C}(A, \mrm{gr})$ we get 
\[ \opn{Hom}_{\dcat{K}(A, \mrm{gr})}(N, I) \cong 
\opn{Hom}_{\dcat{K}(A, \mrm{gr})}(N, J) = 0 . \]
Therefore $I$ is K-graded-injective.

Conversely, let's assume that $I$ is K-graded-injective. By Theorem 
\ref{thm:3740}(2) there is a quasi-isomorphism
$\rho : I \to J$ in $\dcat{C}_{\mrm{str}}(A, \mrm{gr})$ to a 
semi-graded-cofree complex $J$. By Corollary \ref{cor:1922}, 
$\opn{P}(\rho) : I \to J$ is an 
isomorphism in $\dcat{K}(A, \mrm{gr})$.
\end{proof}

The next proposition tells us that some resolving properties of complexes are 
preserved by the restriction functor
$\dcat{C}(A \ot B, \mrm{gr}) \to \dcat{C}(A, \mrm{gr})$.

\begin{prop} \label{prop:4906}  
Let $A$ and $B$ be graded rings.
\begin{enumerate}
\item If $P \in \dcat{C}(A \ot B, \mrm{gr})$ is K-graded-projective, then it is 
K-graded-projective in $\dcat{C}(A, \mrm{gr})$.

\item If $P \in \dcat{C}(A \ot B, \mrm{gr})$ is K-graded-flat, then it is 
K-graded-flat in $\dcat{C}(A, \mrm{gr})$.

\item If $I \in \dcat{C}(A \ot B, \mrm{gr})$ is K-graded-injective, then it is 
K-graded-injective in $\dcat{C}(A, \mrm{gr})$.
\end{enumerate}
\end{prop}

\begin{proof} \mbox{}

\smallskip \noindent
(1) Take $N \in \dcat{C}(A, \mrm{gr})$. 
We have isomorphisms
\[ \begin{aligned}
& \opn{Hom}_{A}(P, N) \cong
\opn{Hom}_{A \ot B}(P, \opn{Hom}_A(A \ot B, N)) 
\\ & \quad  \cong \opn{Hom}_{A \ot B}(P, \opn{Hom}_{\K}(B, N)) 
\end{aligned} \]
in $\dcat{C}_{\mrm{str}}(\K, \mrm{gr})$. 
As an object of $\dcat{C}(\K, \mrm{gr})$, $B$ is K-graded-projective.
Hence, if $N$ is acyclic, then so is $\opn{Hom}_{A}(P, N)$.

\medskip \noindent
(2, 3) Like item (1), but using the fact that $B$ is 
K-graded-flat over $\K$, and using other adjunction isomorphisms.
\end{proof}

\begin{cor} \label{cor:3830}
Let $A$ and $B$ be graded rings, and let 
$M \in \dcat{C}(A \ot B, \mrm{gr})$.
\begin{enumerate}
\item There exists a quasi-isomorphism 
$P \to M$ in $\dcat{C}_{\mrm{str}}(A \ot B, \mrm{gr})$,
where $P$ is a K-graded-projective complex over $A$, 
every $P^i$ is a graded-projective module over $A$, and
$\opn{sup}(P) = \opn{sup}(\opn{H}(M))$. 

\item There exists a quasi-isomorphism 
$P \to M$ in $\dcat{C}_{\mrm{str}}(A \ot B, \mrm{gr})$,
where $P$ is a K-graded-flat complex over $A$, 
every $P^i$ is a graded-flat module over $A$, and
$\opn{sup}(P) = \opn{sup}(\opn{H}(M))$. 

\item There exists a quasi-isomorphism 
$M \to I$ in $\dcat{C}_{\mrm{str}}(A \ot B, \mrm{gr})$,
where $I$ is a K-graded-injective complex over $A$, 
every $I^p$ is a graded-injective module over $A$, and
$\opn{inf}(I) = \opn{inf}(\opn{H}(M))$. 
\end{enumerate}
\end{cor}

\begin{proof} \mbox{}

\smallskip \noindent
(1) Take the resolution $P \to M$ from Theorem \ref{thm:3740}(1). 
In view of Theorem \ref{thm:4275}(1), Proposition \ref{prop:4906}(1) 
and Proposition \ref{prop:3715}(2),
the complex $P$ has the required properties. 

\medskip \noindent
(2) Use item (1), plus Theorem \ref{thm:4275}(2). 

\medskip \noindent
(3) Take the resolution $M \to I$ from Theorem \ref{thm:3740}(2).
This will do, by Theorem \ref{thm:4275}(4), Proposition \ref{prop:4906}(3), 
Proposition \ref{prop:4905}(3) and Proposition \ref{prop:3900}(3).
\end{proof}

\begin{cor} \label{cor:3740}
Let $A$, $B$ and $C$ be graded rings.
\begin{enumerate}
\item The triangulated right derived bifunctor 
\[ \qquad 
\opn{RHom}_A(-, -) : \dcat{D}(A  \ot B^{\mrm{op}}, \mrm{gr})^{\mrm{op}}
\times \dcat{D}(A  \ot C^{\mrm{op}}, \mrm{gr}) \to 
\dcat{D}(B \ot C^{\mrm{op}}, \mrm{gr}) \]
exists. 

\item Let $M \in \dcat{D}(A  \ot B^{\mrm{op}}, \mrm{gr})$ and
$N \in \dcat{D}(A  \ot C^{\mrm{op}}, \mrm{gr})$.
Assume that either $M$ is K-graded-projective over $A$, or $N$
is K-graded-injective over $A$. Then the morphism 
\[ \eta^{\mrm{R}}_{M, N} : \opn{Hom}_A(M, N) \to 
\opn{RHom}_A(M, N) \]
in $\dcat{D}(B \ot C^{\mrm{op}}, \mrm{gr})$
is an isomorphism. 

\item Let $B' \to B$ and $C' \to C$ be graded $\K$-ring homomorphisms. Then 
the diagram 
\[ \UseTips \xymatrix @C=13ex @R=6ex {
\dcat{D}(A  \ot B^{\mrm{op}}, \mrm{gr})^{\mrm{op}} \times 
\dcat{D}(A  \ot C^{\mrm{op}}, \mrm{gr})
\ar[r]^(0.6){\opn{RHom}_A(-, -)}
\ar[d]_{\opn{Rest} \times}^{\opn{Rest}}
&
\dcat{D}(B \ot C^{\mrm{op}}, \mrm{gr})
\ar[d]^{\opn{Rest}}
\\
\dcat{D}(A  \ot {B'}^{\, \mrm{op}}, \mrm{gr})^{\mrm{op}} \times 
\dcat{D}(A  \ot {C'}^{\, \mrm{op}}, \mrm{gr})
\ar[r]^(0.6){\opn{RHom}_A(-, -)}
&
\dcat{D}(B' \ot C'^{\, \mrm{op}}, \mrm{gr})
} \]
is commutative up to isomorphism. 
\end{enumerate}
\end{cor}

\begin{proof}
The same arguments used in the proof of Proposition \ref{prop:1031} work here, 
because there exist enough acyclic resolutions -- see Corollary \ref{cor:3830}.
\end{proof}

\begin{cor} \label{cor:3861}
The category $\dcat{D}(A, \mrm{gr})$ has infinite direct sums. On objects, the 
direct sum is degreewise, for the $\Z^2$ grading. 
\end{cor}

\begin{proof}
The proof of Theorem \ref{thm:3140} works here. 
\end{proof}

\begin{cor} \label{cor:3741} 
Let $A$, $B$ and $C$ be graded rings.
\begin{enumerate}
\item The triangulated left derived bifunctor 
\[ \qquad (- \ot^{\mrm{L}}_{A} -) : 
\dcat{D}(B \ot A^{\mrm{op}}, \mrm{gr})
\times \dcat{D}(A  \ot C^{\mrm{op}}, \mrm{gr}) \to 
\dcat{D}(B \ot C^{\mrm{op}}, \mrm{gr}) \]
exists. 

\item Let $M \in \dcat{D}(B \ot A^{\mrm{op}}, \mrm{gr})$ and
$N \in \dcat{D}(A  \ot C^{\mrm{op}}, \mrm{gr})$.
Assume that either $M$ is K-graded-flat over $A^{\mrm{op}}$, or $N$
is K-graded-flat over $A$. Then the morphism 
\[ \eta^{\mrm{L}}_{M, N} : M \ot^{\mrm{L}}_{A} N \to M \ot_{A} N \]
in $\dcat{D}(B \ot C^{\mrm{op}}, \mrm{gr})$
is an isomorphism. 

\item Let $B' \to B$ and $C' \to C$ be graded $\K$-ring homomorphisms. Then 
the diagram 
\[ \UseTips \xymatrix @C=11ex @R=6ex {
\dcat{D}(B \ot A^{\mrm{op}}, \mrm{gr}) \times 
\dcat{D}(A \ot C^{\mrm{op}}, \mrm{gr})
\ar[r]^(0.62){ (- \ot^{\mrm{L}}_{A} -) }
\ar[d]_{\opn{Rest} \times}^{\opn{Rest}}
&
\dcat{D}(B \ot C^{\mrm{op}}, \mrm{gr})
\ar[d]^{\opn{Rest}}
\\
\dcat{D}(B' \ot A^{\mrm{op}}, \mrm{gr}) \times 
\dcat{D}(A  \ot C'^{\, \mrm{op}}, \mrm{gr})
\ar[r]^(0.62){ (- \ot^{\mrm{L}}_{A} -) }
&
\dcat{D}(B' \ot C'^{\, \mrm{op}}, \mrm{gr})
} \]
is commutative up to isomorphism. 
\end{enumerate}
\end{cor}

\begin{proof}
The same arguments used in the proof of Proposition \ref{prop:3450} work here, 
because there exist enough acyclic resolutions -- see Corollary \ref{cor:3830}.
\end{proof}

\begin{cor} \label{cor:3835}
Let $A$, $B$ and $C$ be graded rings.
\begin{enumerate}
\item Given $M \in \dcat{M}(A  \ot B^{\mrm{op}}, \mrm{gr})$ and
$N \in \dcat{M}(A  \ot C^{\mrm{op}}, \mrm{gr})$, for every $p \in \N$ there is a 
canonical isomorphism 
\[ \opn{Ext}^p_A(M, N) \cong 
\opn{H}^p \bigl( \opn{RHom}_A(M, N) \bigr) \]
in $\dcat{M}(B \ot C^{\mrm{op}}, \mrm{gr})$.
Here $\opn{Ext}^p_A(-, -)$ is the classical $p$-th right derived bifunctor of 
$\opn{Hom}_A(-, -)$. 

\item Given $M \in \dcat{M}(B  \ot A^{\mrm{op}}, \mrm{gr})$ and
$N \in \dcat{M}(A  \ot C^{\mrm{op}}, \mrm{gr})$, for every $p \in \N$ there is a 
canonical isomorphism 
\[ \opn{Tor}_p^A(M, N) \cong 
\opn{H}^{-p} \bigl( M \ot^{\mrm{L}}_{A} N \bigr) \]
in $\dcat{M}(B \ot C^{\mrm{op}}, \mrm{gr})$.
Here $\opn{Tor}_p^A(-, -)$ is the classical $p$-th left derived bifunctor of 
$(- \ot_A -)$. 
\end{enumerate}
\end{cor}

\begin{proof}
This is clear from Corollary \ref{cor:3740}(2) and Corollary \ref{cor:3741}(2) 
respectively.
\end{proof}

\begin{cor} \label{cor:3742}
Let $A$, $B$, $C$ and $D$ be graded rings.
\begin{enumerate}
\item Let $L \in \dcat{D}(A \ot B^{\mrm{op}}, \mrm{gr})$, 
$M \in \dcat{D}(B \ot C^{\mrm{op}}, \mrm{gr})$ and 
$N \in \dcat{D}(C \ot D^{\mrm{op}}, \mrm{gr})$.
There is an isomorphism 
\[ L \ot^{\mrm{L}}_{B} (M \ot^{\mrm{L}}_{C} N) \cong 
(L \ot^{\mrm{L}}_{B} M) \ot^{\mrm{L}}_{C} N \]
in 
$\dcat{D}(A \ot D^{\mrm{op}}, \mrm{gr})$. 
This isomorphism is functorial in $L, M, N$.

\item Let $L \in \dcat{D}(B \ot A^{\mrm{op}}, \mrm{gr})$, 
$M \in \dcat{D}(C \ot B^{\mrm{op}}, \mrm{gr})$ and 
$N \in \dcat{D}(C \ot D^{\mrm{op}}, \mrm{gr})$.
There is an isomorphism 
\[ \opn{RHom}_{B} \bigl( L, \opn{RHom}_{C}(M, N) \bigr) \cong
\opn{RHom}_{C}(M \ot^{\mrm{L}}_{B} L, N) \]
in 
$\dcat{D}(A \ot D^{\mrm{op}}, \mrm{gr})$. 
This isomorphism is functorial in $L, M, N$.
\end{enumerate}
\end{cor}

\begin{proof}
Both assertions are immediate consequences of the associativity of the tensor 
product for modules, and the Hom-tensor adjunction for modules. 
\end{proof}

The triangulated category $\dcat{D}(A^{\mrm{en}}, \mrm{gr})$ has a biclosed 
monoidal structure on it, with monoidal operation 
$(- \ot^{\mrm{L}}_{A} -)$ and monoidal unit $A$. See Remark \ref{rem:3450} for 
the nongraded variant. In the next definition we 
present the pertinent operations related to this structure; these will be 
needed for several constructions. We do not need to know anything about 
the concept of monoidal structure beyond this terminology. 
The interested reader can look it up in \cite{Mac2}; the previous corollaries 
(almost) guarantee that the monoidal axioms hold in this context.

\begin{dfn}[Monoidal Operations] \label{dfn:4530}
\index{Monoidal category! structure on $\dcat{D}(A^{\mrm{en}}, \mrm{gr})$}
Let $A$ and $B$ be graded rings, and let 
$M \in \dcat{D}(A \ot B^{\mrm{op}}, \mrm{gr})$ and 
$N \in \dcat{D}(B \ot A^{\mrm{op}}, \mrm{gr})$ be complexes. 
\begin{enumerate}
\item The {\em left unitor}%
\index{Monoidal category! left unitor isomorphism}%
\index{1-Lu@$\opn{lu}$}
isomorphism for $M$ is the obvious isomorphism 
$\opn{lu} : A \ot_A^{\mrm{L}} M \iso M$
in $\dcat{D}(A \ot B^{\mrm{op}}, \mrm{gr})$. 

\item The {\em right unitor}%
\index{Monoidal category! right unitor isomorphism}%
\index{1-Ru@$\opn{ru}$}
isomorphism for $N$ is the obvious isomorphism 
$\opn{ru} : N \ot_A^{\mrm{L}} A \iso N$
in $\dcat{D}(B \ot A^{\mrm{op}}, \mrm{gr})$. 

\item The {\em left co-unitor}%
\index{Monoidal category! left co-unitor isomorphism}%
\index{1-Lcu@$\opn{lcu}$}
isomorphism for $M$ is the obvious isomorphism  
$\opn{lcu} : \opn{RHom}_{A}(A, M) \iso M$
in $\dcat{D}(A \ot B^{\mrm{op}}, \mrm{gr})$.
\end{enumerate}
\end{dfn}

Here is a notion that is dual to a minimal complex of injectives
(see Definition \ref{dfn:2210}).  
Recall that for a complex $P$, the module of degree $i$ coboundaries is 
$\opn{B}^i(P) = \opn{Im} \bigl( \d_P^{i - 1} : P^{i - 1} \xar{} P^i \bigr)
\sub P^i$.

\begin{dfn} \label{dfn:4531}
Let $A$ be a connected graded ring, with augmentation ideal $\m$. 
A {\em minimal complex of graded-free $A$-modules}%
\index{Complex of algebraically graded modules! minimal {\indash} graded-free
{\indash}}
is a bounded above complex $P$, such that for every $i$ the graded module $P^i$ 
is graded-free, and there is an inclusion 
$\opn{B}^i(P) \sub \m \cd P^i$.
\end{dfn}

\begin{prop} \label{prop:4530}
Let $A$ be a connected graded ring, and let $P$ be a bounded above complex
of graded-free $A$-modules.
The following conditions are equivalent\tup{:}
\begin{itemize}
\rmitem{i} $P$ is a minimal complex of graded-free $A$-modules.

\rmitem{ii} The complex $\K \ot_A P$ has zero differential.
\end{itemize}
\end{prop}

\begin{exer} \label{exer:4530}
Prove this proposition. 
\end{exer}

\begin{dfn} \label{dfn:4235}
Let $A$ be a connected graded ring, and let $M \in \dcat{M}(A, \mrm{gr})$. 
A {\em minimal graded-free resolution}%
\index{Resolution! minimal graded-free}
of $M$ is a quasi-isomorphism 
$\rho : P \to M$ in $\dcat{C}_{\mrm{str}}(A, \mrm{gr})$ from a minimal
complex of graded-free $A$-modules $P$.
\end{dfn}

\begin{prop} \label{prop:4265}   
Let $A$ be a connected graded ring, and let 
$M \in \dcat{M}(A, \mrm{gr})$ be a 
bounded below graded module. Then\tup{:}
\begin{enumerate}
\item $M$ has a minimal graded-free resolution $\eta : P \to M$.

\item For every $i$ there is an isomorphism 
$P^i \cong A \ot_{} \opn{Tor}^A_{-i}(\K, M)$
in $\dcat{M}(A, \mrm{gr})$.

\item The minimal graded-free resolution $\eta : P \to M$ in item \tup{(1)} 
above is unique, up to a non-unique isomorphism in 
$\dcat{C}_{\mrm{str}}(A,\mrm{gr})$.
\end{enumerate} 
\end{prop}

\begin{proof} \mbox{}

\smallskip \noindent
(1) Let $V_0 := \K \ot_A M$ and $P^0 := A \ot V_0$. By 
Proposition \ref{prop:4590} there is a surjection 
$\eta : P^0 \to M$. Let 
$M^{-1} := \opn{Ker}(\eta) \sub P^0$. 
By construction we have 
$M^{-1} \sub \m \ot V_0 = \m \cd P^0$. 
And of course $M^{-1}$ is a bounded below graded $A$-module. 

Next let $V_{1} := \K \ot_A M^{-1}$ and $P^{-1} := A \ot V_{1}$. 
By Proposition \ref{prop:4590} there is a surjection 
$\d^{-1} : P^{-1} \to M^{-1}$. The graded module  
$M^{-2} := \opn{Ker}(\d^{-1}) \sub P^{-1}$
satisfies $M^{-2} \sub \m \ot V_1 = \m \cd P^{-1}$.
And so on. 

\medskip \noindent
(2) Due to minimality, the differential of the complex $\K \ot_A P$ is zero; 
see Proposition \ref{prop:4530}. Therefore
\[ \opn{Tor}^A_{-i}(\K, M) \cong \opn{H}^i(\K \ot^{\mrm{L}}_{A} M) \cong 
\opn{H}^i(\K \ot_{A} P) \cong \K \ot_{A} P^{i} \cong V_{-i} . \]

\medskip \noindent
(3) Suppose $\eta' : P' \to M$ is some other minimal graded-free resolution.
Because $P$ and $P'$ are K-projective in $\dcat{C}(A, \mrm{gr})$, there is a 
homotopy equivalence $\phi : P \to P'$ in $\dcat{C}_{\mrm{str}}(A, \mrm{gr})$.
Therefore 
\begin{equation} \label{eqn:4530}
\opn{id} \ot \, \phi : \K \ot_A P \to \K \ot_A P'
\end{equation}
is a quasi-isomorphism. But by Proposition \ref{prop:4530}, these complexes 
have zero differentials, so (\ref{eqn:4530}) is an isomorphism in 
$\dcat{C}_{\mrm{str}}(\K, \mrm{gr})$. The graded Nakayama Lemma (Proposition 
\ref{prop:4590}(3)) implies that $\phi : P \to P'$ is an isomorphism in 
$\dcat{C}_{\mrm{str}}(A, \mrm{gr})$.
\end{proof}

\begin{exer} \label{exer:4531}
Try to generalize Proposition \ref{prop:4265} to a complex 
$M \in \lb \dcat{C}(A, \mrm{gr})$. The trick is to find the correct
boundedness conditions on $M$.
\end{exer}

\begin{cor} \label{cor:4237}
Assume $A$ is a left noetherian connected graded ring, and $M$ is a 
finite graded $A$-module. Then in the minimal graded-free 
resolution $\eta : P \to M$, the graded-free $A$-modules $P^i$ are all finite. 
\end{cor}

\begin{proof}
In the proof of Proposition \ref{prop:4265}(1), the graded $\K$-modules $V_i$ 
are all finite. 
\end{proof}

Here is a graded version of Definitions \ref{dfn:1925} and \ref{dfn:4470}. 

\begin{dfn} \label{dfn:4532} 
Let $A$ be a graded ring. 
\begin{enumerate}
\item A complex $P \in \dcat{C}(A, \mrm{gr})$
is called {\em pseudo-finite semi-graded-free}%
\index{Complex of algebraically graded modules! pseudo-finite semi-graded-free}
if it is a bounded above complex of finite graded-free $A$-modules.

\item A complex $L \in \dcat{D}(A, \mrm{gr})$
is called {\em derived graded-pseudo-finite}%
\index{Complex of algebraically graded modules! derived graded-pseudo-finite}
if it belongs to the saturated 
full triangulated subcategory of $\dcat{D}(A, \mrm{gr})$ generated by the 
pseudo-finite semi-graded-free complexes.
\end{enumerate}
\end{dfn}

\begin{exa} \label{exa:4390}
If the graded ring $A$ is connected and left noetherian, then the 
derived graded-pseudo-finite complexes over $A$ are precisely the objects 
of $\dcat{D}^{-}_{\mrm{f}}(A, \mrm{gr})$.
See Theorems \ref{thm:2895} or \ref{thm:3340}. 
\end{exa}

\begin{prop} \label{prop:4266}  
Assume $A$ is a left noetherian connected graded ring.
Consider $A$ as a graded $A^{\mrm{en}}$-module. Then in the minimal graded-free 
resolution $\rho : P \to A$ of $A$ over $A^{\mrm{en}}$, the graded-free 
$A^{\mrm{en}}$-modules $P^i$ are all finite. 
\end{prop}

The subtle point is that the graded ring $A^{\mrm{en}}$ might not be left 
noetherian -- see Remark \ref{rem:3901}. 

\begin{proof}
Let us write $P^{i} \cong A^{\mrm{en}} \ot V_{-i}$, where 
$V_{-i} \in \dcat{M}(\K, \mrm{gr})$. We shall prove that all the
$V_{-i}$ are finite graded $\K$-modules. 

When we restrict the quasi-isomorphism $\rho : P \to A$
to the category $\dcat{C}_{\mrm{str}}(A^{\mrm{op}}, \mrm{gr})$, 
it becomes a homotopy equivalence. Hence 
$\tau := \rho \ot \opn{id} : P \ot_{A} \K \to A \ot_{A} \K$
is a quasi-isomorphism in $\dcat{C}_{\mrm{str}}(A, \mrm{gr})$. 
Writing 
$Q := P \ot_{A} \K \in \dcat{C}(A, \mrm{gr})$,
we get a graded-free resolution 
$\tau : Q \to \K$ over $A$. Note that 
$Q^{i} \cong A \ot V_{-i}$ in $\dcat{C}_{\mrm{str}}(A, \mrm{gr})$.

Now 
\begin{equation} \label{eqn:5433}
\K \ot_A Q \cong \K \ot_A P \ot_A \K \cong \K \ot_{A^{\mrm{en}}} P
\end{equation}
in $\dcat{C}_{\mrm{str}}(\K, \mrm{gr})$. Because $P$ is a minimal 
complex of graded-free $A^{\mrm{en}}$-modules, 
the complex $\K \ot_{A^{\mrm{en}}} P$ has zero differential; see 
Proposition \ref{prop:4530}.
But then the complex $\K \ot_A Q$ has zero differential. We conclude that $Q$ 
is a minimal complex of graded-free $A$-modules, and hence 
$\tau : Q \to \K$ is a minimal graded-free resolution over $A$. Corollary 
\ref{cor:4237} says that the graded $A$-modules $Q^i$ are finite; and hence the 
graded $\K$-modules $V_{-i}$ are finite. 
\end{proof}

\begin{cor} \label{cor:4531}
If $A$ is a left noetherian connected graded ring, then $A$ is a 
derived graded-pseudo-finite complex over $A^{\mrm{en}}$. 
\end{cor}

\begin{proof}
This is an immediate consequence of Proposition \ref{prop:4266},
which says that the complex $P$ is pseudo-finite graded-semi-free over 
$A^{\mrm{en}}$.
\end{proof}

\begin{thm} \label{thm:4535}  
Let $A$, $B$, $C$ and $D$ be graded rings. For complexes 
$L \in \dcat{D}(A \ot C^{\mrm{op}}, \mrm{gr})$,
$M \in \dcat{D}(A \ot B^{\mrm{op}}, \mrm{gr})$ and 
$N \in \dcat{D}(B \ot D^{\mrm{op}}, \mrm{gr})$
there is a morphism 
\[ \opn{ev}^{\mrm{R, L}}_{L, M, N} :
\opn{RHom}_A(L, M) \ot_B^{\mrm{L}} N \to 
\opn{RHom}_A(L, M \ot_B^{\mrm{L}} N) \]
in $\dcat{D}(C \ot D^{\mrm{op}}, \mrm{gr})$, called 
{\em derived graded tensor-evaluation}.  
\index{Tensor-evaluation morphism! derived graded}
The morphism $\opn{ev}^{\mrm{R, L}}_{L, M, N}$ 
is functorial in the objects $L, M, N$.
Moreover, if all three conditions below hold, then 
$\opn{ev}^{\mrm{R, L}}_{L, M, N}$ is an isomorphism. 
\begin{enumerate}
\rmitem{i} The complex $L$ is derived graded-pseudo-finite over $A$.

\rmitem{ii} The complex $M$ has bounded below cohomology. 

\rmitem{iii} The complex $N$ has finite graded-flat dimension over $B$. 
\end{enumerate}
\end{thm}

Just to clarify, condition (ii) says that $\opn{H}^i(M) = 0$ for $i \ll 0$. 
We do not care about the boundedness (in algebraic degree) of the graded 
modules $\opn{H}^i(M)$. 

\begin{proof}
The proof of just like that of Theorem \ref{thm:4320}, but 
working in the graded derived categories. For condition (iii), notice 
that if $N$ has finite graded-flat dimension over $B$, then there is an 
isomorphism $P \cong N$ in $\dcat{D}(B, \mrm{gr})$, where $P$ is a bounded 
complex of graded-flat $B$-modules.
\end{proof}

\begin{dfn} \label{dfn:3732} 
Let $A$ be a graded ring. 
\begin{enumerate}
\item The full subcategory of $\dcat{D}(A, \mrm{gr})$ on the complexes $M$ 
whose cohomology modules $\opn{H}^p(M)$ belong to 
$\dcat{M}_{\mrm{dwf}}(A, \mrm{gr})$
is denoted by $\dcat{D}_{\mrm{dwf}}(A, \mrm{gr})$%
\index{1-Ddwf(A,gr)@$\dcat{D}_{\mrm{dwf}}(A, \mrm{gr})$}.

\item The full subcategory of $\dcat{D}(A, \mrm{gr})$ on the complexes 
whose cohomology modules belong to 
$\dcat{M}_{\mrm{f}}(A, \mrm{gr})$
is denoted by $\dcat{D}_{\mrm{f}}(A, \mrm{gr})$%
\index{1-Df(A,gr)@$\dcat{D}_{\mrm{f}}(A, \mrm{gr})$}.

\item The full subcategory of $\dcat{D}(A, \mrm{gr})$ on the complexes 
whose cohomology modules belong to 
$\dcat{M}_{\mrm{cof}}(A, \mrm{gr})$
is denoted by $\dcat{D}_{\mrm{cof}}(A, \mrm{gr})$%
\index{1-Dcof(A,gr)@$\dcat{D}_{\mrm{cof}}(A, \mrm{gr})$}.
\end{enumerate}
\end{dfn}

Because $\dcat{M}_{\mrm{dwf}}(A, \mrm{gr})$ is a a thick abelian 
subcategory of $\dcat{M}(A, \mrm{gr})$,
it follows that $\dcat{D}_{\mrm{dwf}}(A, \mrm{gr})$
is a full triangulated subcategory of 
$\dcat{D}(A, \mrm{gr})$. Likewise for 
$\dcat{D}_{\mrm{f}}(A, \mrm{gr})$
when $A$ is left noetherian, and for 
$\dcat{D}_{\mrm{cof}}(A, \mrm{gr})$
when $A$ is right noetherian connected.

Because $(-)^*$ is exact, it induces a triangulated functor 
\begin{equation} \label{eqn:3733}
(-)^* : \dcat{D}(A, \mrm{gr})^{\mrm{op}} \to 
\dcat{D}(A^{\mrm{op}}, \mrm{gr}) . 
\end{equation}

\begin{prop} \label{prop:4591}   
Let $A$ be a graded ring.
\begin{enumerate}
\item If $M \in \dcat{D}_{\mrm{dwf}}(A, \mrm{gr})$ then 
$M^* \in \dcat{D}_{\mrm{dwf}}(A^{\mrm{op}}, \mrm{gr})$.

\item The functor 
\[ (-)^* : \dcat{D}_{\mrm{dwf}}(A, \mrm{gr})^{\mrm{op}} \to 
\dcat{D}_{\mrm{dwf}}(A^{\mrm{op}}, \mrm{gr}) \]
is an equivalence of $\K$-linear triangulated categories. 

\item If $A$ is left noetherian connected, then the functor 
\[ (-)^* : \dcat{D}_{\mrm{f}}(A, \mrm{gr})^{\mrm{op}} \to 
\dcat{D}_{\mrm{cof}}(A^{\mrm{op}}, \mrm{gr}) \]
is an equivalence of $\K$-linear triangulated categories. 
\end{enumerate}
\end{prop}

\begin{exer} \label{exer:3830}
Prove this proposition. (Hint: use Proposition \ref{prop:4045} and 
Theorem \ref{thm:4045}.) 
\end{exer} 

The next result is from \cite{VdB}. 

\begin{thm} \label{thm:3800} 
Let $A$, $B$ and $C$  be graded $\K$-rings. 
\begin{enumerate}
\item For 
$M \in \dcat{D}(A \ot B^{\mrm{op}}, \mrm{gr})$ and
$N \in \dcat{D}(A \ot C^{\mrm{op}}, \mrm{gr})$
there is a morphism 
\[ \th_{M, N} :  \opn{RHom}_A(M, N) \to \opn{RHom}_{A^{\mrm{op}}}(N^*, M^*) \]
in $\dcat{D}(B \ot C^{\mrm{op}}, \mrm{gr})$, that is functorial in $M$ and $N$.

\item If
$M \in \dcat{D}_{\mrm{dwf}}(A \ot B^{\mrm{op}}, \mrm{gr})$ 
and
$N \in \dcat{D}_{\mrm{dwf}}(A \ot C^{\mrm{op}}, \mrm{gr})$
then $\th_{M, N}$ is an isomorphism. 
\end{enumerate}
\end{thm}

\begin{proof} \mbox{}

\smallskip \noindent
(1) Choose a K-graded-projective resolution $\rho : P \to M$ in 
$\dcat{C}_{\mrm{str}}(A \ot B^{\mrm{op}}, \mrm{gr})$. 
By Proposition \ref{prop:4906} the complex $P$ is 
K-graded-projective over $A$; so there is a canonical isomorphism 
\[ \phi_{\rho} : \opn{RHom}_A(M, N) \iso \opn{Hom}_A(P, N) \]
in $\dcat{D}(B \ot C^{\mrm{op}}, \mrm{gr})$. 
$\K$-linear duality gives a homomorphism
\[ \til{\th}_{P, N} : \opn{Hom}_A(P, N) \to 
\opn{Hom}_{A^{\mrm{op}}}(N^*, P^*) \]
in $\dcat{C}_{\mrm{str}}(B \ot C^{\mrm{op}}, \mrm{gr})$.

Now $\rho^* : M^* \to P^*$ is a quasi-isomorphism
in $\dcat{C}_{\mrm{str}}(B \ot A^{\mrm{op}}, \mrm{gr})$,
and according to Proposition \ref{prop:4906}(1) and 
Theorem \ref{thm:4275}(3) the the complex $P^*$ is 
K-graded-injective over $A^{\mrm{op}}$. Therefore there is a canonical 
isomorphism 
\[ \psi_{\rho^*} : \opn{Hom}_{A^{\mrm{op}}}(N^*, P^*) \iso 
\opn{RHom}_{A^{\mrm{op}}}(N^*, M^*) \]
in $\dcat{D}(B \ot C^{\mrm{op}}, \mrm{gr})$. 
The resulting morphism 
\[ \th_{M, N} := \psi_{\rho^*} \circ \opn{Q}(\til{\th}_{P, N}) \circ \phi_{\rho}
: \opn{RHom}_A(M, N) \to \opn{RHom}_{A^{\mrm{op}}}(M^*, N^*) \]
in $\dcat{D}(B \ot C^{\mrm{op}}, \mrm{gr})$ is functorial in $M$ and $N$.

\medskip \noindent
(2) Because the restriction functors between the derived categories are 
conservative (see Proposition \ref{prop:3445}), 
we can forget the rings $B$ and $C$, and just consider 
$\th_{M, N}$ as a morphism $\dcat{D}(\K, \mrm{gr})$. 

For every integer $i$ there is a canonical isomorphism
\begin{equation} \label{eqn:3805}
\opn{H}^i \bigl( \opn{RHom}_A(M, N) \bigr)  
\cong \opn{Hom}_{\dcat{D}_{\mrm{dwf}}( A, \mrm{gr})}
\bigl( M, N[i] \bigr) .
\end{equation}
Using the isomorphism 
$N^*[-i] \cong N[i]^*$ in 
$\dcat{D}(A^{\mrm{op}}, \mrm{gr})$ we also have  
\begin{equation} \label{eqn:3806}
\begin{aligned}
& \opn{H}^i \bigl( \opn{RHom}_{A^{\mrm{op}}}(N^*, M^*) \bigr) \cong 
\opn{Hom}_{\dcat{D}_{\mrm{dwf}}(A^{\mrm{op}}, \mrm{gr})}
\bigl( N^*, M^*[i] \bigr)
\\
& \quad \cong
\opn{Hom}_{\dcat{D}_{\mrm{dwf}}(A^{\mrm{op}}, \mrm{gr})}
\bigl( N[i]^*, M^* \bigr) . 
\end{aligned}
\end{equation}
The construction of $\th_{M, N}$ in the proof of item (1) above shows that the 
$\K$-linear homomorphism 
\[ \opn{H}^i(\th_{M, N}) :  \opn{H}^i \bigl( \opn{RHom}_A(M, N) \bigr) \to
\opn{H}^i \bigl( \opn{RHom}_{A^{\mrm{op}}}(N^*, M^*) \bigr) \]
becomes, after using the isomorphisms (\ref{eqn:3805}) and 
(\ref{eqn:3806}), the duality morphism 
\[ \opn{Hom}_{\dcat{D}_{\mrm{dwf}}(A, \mrm{gr})}
\bigl( M, N[i] \bigr)  \to  
\opn{Hom}_{\dcat{D}_{\mrm{dwf}}(A^{\mrm{op}}, \mrm{gr})}
\bigl( N[i]^*, M^* \bigr) . \]
According to Proposition \ref{prop:4591}(2) this is an isomorphism.
\end{proof}

\begin{rem} \label{rem:3740}
In this section, and the two sections following it, we work over a base field 
$\K$. This simplifies a lot of the constructions. It is however possible to 
relax this condition, and to work over a noetherian commutative base ring $\K$ 
that is not a field. But then there are (at least) two complications. 
The first complication is that we will have to fix an injective cogenerator 
$\K^*$ of $\dcat{M}(\K)$; cf.\ Subsection \ref{subsec:exis-K-inj-dgmods}. 
Matlis Duality will take the form 
$M^* := \opn{Hom}_{\K}(M, \K^*)$,
and $\opn{Hom}_{\K}(-, -)$ might also involve continuity in a delicate way. 

The second complication is that for the homological algebra to be effective, 
the graded $\K$-rings would have to be {\em flat}. 
The absence of flatness poses serious difficulties. This can be tackled using 
{\em K-flat DG ring resolutions}; and these DG rings would have to be bigraded 
(with upper and lower indices, as explained in the Subsection 
\ref{subsec:alg-gr-mods}). This sort 
of resolution has not been studied yet, as far as we know; but cf.\ 
\cite{Ye16} for some preliminary ideas. 
\end{rem}

\mysubsection{Artin-Schelter Regular Graded Rings} 
\label{AS-reg}

In this subsection we talk about an important class of graded rings: the 
{\em Artin-Schelter regular graded rings}. These are connected graded 
noncommutative rings, that homologically are similar to commutative polynomial 
rings. A great deal of noncommutative ring theory and homological algebra grew 
out of the study of this class of rings. See Remark \ref{rem:4653} for 
some historical notes.

We continue with Conventions \ref{conv:3700} and \ref{conv:4560}. 
Thus all graded rings are central over the base field $\K$. 

Suppose $A$ is a connected graded ring. Recall that the augmentation ideal of 
$A$ is denoted by $\m$. The opposite ring of $A$ is $A^{\mrm{op}}$, and the 
enveloping ring is $A^{\mrm{en}}$.
We view $A / \m \cong \K$ as a graded $A^{\mrm{en}}$-module.

\begin{dfn} \label{dfn:3920} 
Let $A$ be a noetherian connected graded ring. 
\begin{enumerate}
\item The graded ring $A$ is called a {\em regular graded ring}
\index{Algebraically graded ring! regular}
if it has finite graded global cohomological dimension; namely if there is a 
natural number $d$ such that 
$\opn{Ext}^i_A(M, N) = 0$ and
$\opn{Ext}^i_{A^{\mrm{op}}}(M', N') = 0$ 
for all $i > d$,  $M, N \in \dcat{M}(A, \mrm{gr})$ and 
$M', N' \in \dcat{M}(A^{\mrm{op}}, \mrm{gr})$.
The smallest such natural number $d$ is called the 
{\em graded global dimension}
\index{Algebraically graded ring! graded global dimension of}
of the ring $A$.

\item The graded ring $A$ is called a {\em Gorenstein graded ring}
\index{Algebraically graded ring! Gorenstein} 
if the graded bimodule $A$ has finite graded injective dimension 
over $A$ and over $A^{\mrm{op}}$; namely if there is a natural number $d$ 
such that $\opn{Ext}^i_A(M, A) = 0$ and
$\opn{Ext}^i_{A^{\mrm{op}}}(M', A) = 0$ 
for all $i > d$, $M \in \dcat{M}(A, \mrm{gr})$ and 
$M' \in \dcat{M}(A^{\mrm{op}}, \mrm{gr})$.
The smallest such natural number $d$ is called the 
{\em graded injective dimension} 
\index{Algebraically graded ring! graded injective dimension of}
of the ring $A$.
\end{enumerate}
\end{dfn}

Recall that the graded projective dimension of a graded $A$-module $M$ is 
smallest generalized integer $d \in \N \cup \{ \pm \infty \}$ 
such that $\opn{Ext}^i_A(M, N) = 0$ 
for all integers $i > d$ and all $N \in \dcat{M}(A, \mrm{gr})$.

The minimal graded-free resolution of a graded module was introduced in 
Definitions  \ref{dfn:4531} and \ref{dfn:4235}. 
 
\begin{thm} \label{thm:4548} 
Let $A$ be a noetherian connected graded ring, and let $d$ be a natural 
number. The following five conditions are equivalent. 
\begin{itemize}
\rmitem{i} The graded ring $A$ is regular, of graded global dimension $d$. 

\rmitem{ii} The graded projective dimension of the graded $A$-module $\K$ is 
$d$. 

\rmitem{iii} The graded $\K$-modules 
$\opn{Tor}_i^A(\K, \K)$ vanish for all $i > d$, but not for $i = d$. 

\rmitem{iv} The graded $\K$-modules 
$\opn{Ext}^i_A(\K, \K)$ vanish for all $i > d$, but not for $i = d$. 

\rmitem{v} Let $P \to A$ be the minimal graded-free resolution of the graded
$A^{\mrm{en}}$-module $A$. Then the graded $A^{\mrm{en}}$-modules $P^{-i}$  
vanish for all $i > d$, but not for $i = d$.  
\end{itemize}
\end{thm}

Observe that in this theorem, conditions (i), (iii) and (v) are 
op-symmetric, but not conditions (ii) and (iv). 

\begin{proof}
Consider the minimal graded-free resolution
$\rho : P \to A$ of $A$ over $A^{\mrm{en}}$. 
We express the graded-free $A^{\mrm{en}}$-module $P^{-i}$, for $i \geq 0$,  as 
\begin{equation} \label{eqn:4905}
P^{-i} =  A^{\mrm{en}} \ot V_{i} \cong A \ot V_i \ot A 
\end{equation}
with $V_{i} \in \dcat{M}_{\mrm{f}}(\K, \mrm{gr})$. 
The finiteness of the graded $\K$-modules $V_{i}$ is by Proposition 
\ref{prop:4266}. In formula (\ref{eqn:4905}) the ring $A^{\mrm{en}}$ acts on 
itself from the left, and it acts on $A \ot V_i \ot A$ by the outside action,
i.e.\ 
\[ (a_1 \ot \opn{op}(a_2)) \cd (b_1 \ot v \ot b_2) = 
(a_1 \cd b_1) \ot v \ot (b_2 \cd a_2) \]
for $a_1 \ot \opn{op}(a_2) \in A^{\mrm{en}}$ and 
$b_1 \ot v \ot b_2 \in A \ot V_i \ot A$.

As explained in the proof of Proposition \ref{prop:4266}, in the category \lb
$\dcat{C}_{\mrm{str}}(A^{\mrm{op}}, \mrm{gr})$
the homomorphism $\rho : P \to A$ is a homotopy equivalence. Hence, taking 
a module $M \in \dcat{M}(A, \mrm{gr})$, and writing 
$Q_M := P \ot_{A} M$, we get a quasi-isomorphism 
\[ \tau_M := \rho \ot \opn{id}_M : Q_M = P \ot_{A} M \to A \ot_{A} M \cong M  \]
in $\dcat{C}_{\mrm{str}}(A, \mrm{gr})$. 
Because $M$ is a graded-free $\K$-module, it follows that every
\[ Q_M^{-i} = P^{-i} \ot_{A} M \cong A \ot_{} V_{i} \ot_{} M \]
is also a graded-free $A$-module.
In this way we obtain a (functorial) graded-free resolution
$\tau_M : Q_M \to  M$ 
in $\dcat{C}_{\mrm{str}}(A, \mrm{gr})$. 

Assume now that condition (v) holds. This means that 
$V_{i} = 0$ for all $i > d$, and $V_{d} \neq 0$.
The graded-free resolution $Q_M$ has 
length $d$ for every nonzero $M \in \dcat{M}(A, \mrm{gr})$.
We see that  
\[  \opn{Ext}^i_A(M, N) \cong \opn{H}^i \bigl( \opn{Hom}_A(Q_M, N) \bigr)
= 0 \]
for all $N \in  \dcat{M}(A, \mrm{gr})$ and $i > d$. By the op-symmetry 
of condition (v), it also follows that 
$\opn{Ext}^i_{A^{\mrm{op}}}(M', N') = 0$
for all $M', N' \in  \dcat{M}(A^{\mrm{op}}, \mrm{gr})$ and $i > d$.
This says that the graded global dimension of $A$ is at most $d$. 
To get equality, we test for $M = N = \K$. 
In this case we know that $Q_{\K}$ is a minimal graded-free resolution of 
$\K$ over $A$ (see the proof of Proposition \ref{prop:4266}), so the complex 
$\K \ot_A Q_{\K}$ has trivial differential. Therefore
\begin{equation} \label{eqn:4906}
\begin{aligned}
&
\opn{Ext}^i_A(\K, \K) \cong 
\opn{H}^i \bigl( \opn{Hom}_A(Q_{\K}, \K) \bigr)
\\ & \quad 
\cong \opn{H}^i \bigl( \opn{Hom}_{\K}(\K \ot_A Q_{\K}, \K) \bigr) \cong 
\opn{Hom}_{\K}(V_{i}, \K) . 
\end{aligned}
\end{equation}
This is nonzero for $i = d$. 
Likewise over $A^{\mrm{op}}$. 
We have thus verified conditions (i), (ii) and (iv). 
Condition (iii) is verified similarly: the resolution $Q_{\K}$ shows that 
$\opn{Tor}_i^A(\K, \K) \cong V_{i}$
as graded $\K$-modules.

Let us prove the implication (i) $\Rightarrow$ (v):  since  
$\opn{Ext}^i_A(\K, \K) = 0$
for all $i > d$, equations (\ref{eqn:4906}) and (\ref{eqn:4905}) say that 
$P^{-i} = 0$ in this range. On the other hand, there is some pair 
$M, N \in \dcat{M}(A, \mrm{gr})$ such 
that $\opn{Ext}^d_A(M, N) \neq 0$, and this says that $P^{-d} \neq 0$.

The other implications are proved similarly, and we leave them as an exercise. 
\end{proof}

\begin{exer} \label{exer:4550}
Finish the proof of Theorem \ref{thm:4548}.
\end{exer}

The next condition first appeared in the paper \cite{ArSchl} by M. Artin and 
W. Schelter. Throughout, ``AS'' is an abbreviation for ``Artin-Schelter''.  

\begin{dfn} \label{dfn:4545}
Let $A$ be a noetherian connected graded ring. We say that $A$ satisfies the 
{\em AS condition}%
\index{Algebraically graded ring! AS condition on}
if there are integers $n$ and $l$ such that 
\[ \opn{Ext}^n_A(\K, A) \cong \opn{Ext}^n_{A^{\mrm{op}}}(\K, A) 
\cong \K(l) \]
as graded $\K$-modules, and  
\[ \opn{Ext}^i_A(\K, A) \cong \opn{Ext}^i_{A^{\mrm{op}}}(\K, A) \cong 0 \]
for all $i \neq n$.
The number $n$ is called the {\em AS dimension} of $A$, and the 
number $l$ is called the {\em AS index} of $A$.
\end{dfn}
 
In derived category terms, the AS condition says that 
\begin{equation} \label{eqn:4562}
\opn{RHom}_A(\K, A) \cong \opn{RHom}_{A^{\mrm{op}}}(\K, A) 
\cong \K(l)[-n]
\end{equation}
in $\dcat{D}(\K, \mrm{gr})$, 

\begin{dfn} \label{dfn:3782} 
Let $A$ be a noetherian connected graded ring. 
\begin{enumerate}
\item We say that $A$ is an {\em AS regular graded ring} 
\index{Algebraically graded ring! AS regular}
of dimension $n$ if it
is a regular graded ring (Definition \ref{dfn:3920}(1)), and it 
satisfies the {\em AS condition} (Definition \ref{dfn:4545}) with AS dimension 
$n$. 

\item We say that $A$ is an {\em AS Gorenstein graded ring} 
\index{Algebraically graded ring! AS Gorenstein}
of dimension $n$ if 
it is a Gorenstein graded ring (Definition \ref{dfn:3920}(2)), and it 
satisfies the {\em AS condition} (Definition \ref{dfn:4545}) with AS dimension 
$n$. 
\end{enumerate}
\end{dfn}

As we shall see later, in Subsection \ref{subsec:sym-der-tor}, the AS 
condition implies that the local cohomology of graded $A$-modules behaves as 
though $A$ is commutative. In the paper \cite{ArZh} by Artin and J.J. Zhang, 
subsequent to \cite{ArSchl}, this behavior was formalized as the {\em $\chi$ 
condition}; see Definition \ref{dfn:3713} below. 

Some texts do not make the assumption that the connected graded ring $A$ is 
noetherian, as we did in Definition \ref{dfn:3782}. However, without 
this assumption the theory is not as rich.  

It is natural to ask whether the AS dimension of the ring $A$ coincides 
with its graded global dimension (in the AS regular case), or with its 
graded injective dimension (in the AS Gorenstein case). The next 
propositions say that these numbers agree. 

\begin{prop} \label{prop:4560}
Let $A$ be an AS regular graded ring.
Then the AS dimension of $A$ equals its graded global dimension. 
\end{prop}

\begin{proof}
Let $n$ be the graded global dimension of $A$. 
Let $P \to A$ be the minimal graded-free resolution of the graded
$A^{\mrm{en}}$-module $A$. By Theorem \ref{thm:4548} we know that the 
graded $A^{\mrm{en}}$-modules $P^{-i}$  vanish 
for all $i > n$, but not for $i = n$.  
Let $Q_{\K} \to \K$ be the induced graded-free resolution of the $A$-module 
$\K$ from the proof of Theorem \ref{thm:4548}. Then the complex 
$Q_{\K}$ is concentrated in cohomological degrees $[-n, 0]$, and 
$Q^{-n}_{\K} \neq 0$. Also, as shown there, $Q_{\K} \to \K$ is a 
minimal graded-free resolution of $\K$ over $A$. 

Now let 
$T := \opn{Hom}_A(Q_{\K}, A) \in \dcat{C}(A^{\mrm{op}}, \mrm{gr})$,
so $\opn{RHom}_A(\K, A) \cong T$ 
in $\dcat{D}(\K, \mrm{gr})$. A quick calculation 
shows that $T$ is a minimal complex of graded-free
$A^{\mrm{op}}$-modules. This complex is concentrated in cohomological degrees 
$[0, n]$, and $T^n \neq 0$. Therefore 
$\opn{Ext}^n_A(\K, A) \cong \opn{H}^n (T) \neq 0$.
A similar calculation shows that 
$\opn{Ext}^n_{A^{\mrm{op}}}(\K, A) \neq 0$.
We see that the AS dimension of $A$ is $n$. 
\end{proof}

\begin{prop} \label{prop:4561}
Let $A$ be an AS Gorenstein graded ring.
Then the AS dimension of $A$ equals its graded injective dimension. 
\end{prop}

\begin{proof}
This is an immediate consequence of \cite[Theorem 4.5]{Jor}. 
\end{proof}

In Examples \ref{exa:3780} and \ref{exa:3990} we will present some important 
kinds of AS regular rings. 

\begin{dfn} \label{dfn:4805}
Let $A$ be a nonzero graded ring. 
A central element $a$ of $A$ is called a {\em regular central element}
\index{Regular central element! of a ring} 
if it is not a zero-divisor; i.e.\ $a \cd b = 0$ implies $b = 0$. 
\end{dfn}

If the element $a \in A$ is not central, then we have to distinguish between 
left and right regularity of $a$.

The next two theorems discuss the relation between the rings $A$ and $B$ from 
the following setup.

\begin{setup} \label{set:4565}
Let $A$ be a noetherian connected graded ring, and let $a \in A$ be a
homogeneous regular central element of positive degree. Define the connected 
graded ring $B := A / (a)$.
\end{setup}

\begin{thm} \label{thm:4546}
Assume Setup \tup{\ref{set:4565}}.
\begin{enumerate}
\item If  $B$ is a regular graded ring of graded global dimension $n - 1$, 
then $A$ is a regular graded ring of graded global dimension $n$. 

\item Moreover, if  $B$ is an AS regular graded ring of dimension $n - 1$,  
then $A$ is an AS regular graded of dimension $n$.
\end{enumerate}
\end{thm}

\begin{proof} \mbox{}

\smallskip \noindent
(1) By Theorem \ref{thm:4548}, to prove that $A$ is a regular graded ring
of graded global dimension $n$, 
it is enough to show that the cohomology of the complex
$\K \ot^{\mrm{L}}_{A} \K \in \dcat{D}(\K, \mrm{gr})$
has concentration $[-n, 0]$.
Because a graded-free resolution of $\K$ over $A$ is also a free 
resolution in the ungraded sense, and because we are only interested in the 
vanishing of cohomology, we can forget the algebraic grading, and just look at 
the complex
$\K \ot^{\mrm{L}}_{A} \K \in \dcat{D}(\K)$.
Similarly, because $B$ is a regular graded ring of graded global dimension 
$n - 1$, we know that the complex 
$\K \ot^{\mrm{L}}_{B} \K \in \dcat{D}(\K)$
satisfies 
\begin{equation} \label{eqn:4565}
\opn{con} \bigl( \opn{H}(\K \ot^{\mrm{L}}_{B} \K) \bigr) = [-n + 1, 0] . 
\end{equation}

Let $\til{B} := \opn{K}(A, a)$ be the Koszul complex on the regular central 
element $a \in A$; namely 
\[ \til{B}  = \bigl( \cdots \to 0 \to A \xar{a \cd (-)} A \to 0
\to \cdots \bigr) , \]
concentrated in cohomological degrees $-1$ and $0$. 
This is a DG ring: as a cohomologically graded ring we have 
$\til{B} := A \ot \K[x]$, where $\K[x]$ is the strongly commutative polynomial 
ring on the variable $x$ that has cohomological degree $-1$ 
(see Examples \ref{exa:4625} and \ref{exa:1099}). 
The differential is $\d(x) := a$. 
The canonical DG ring homomorphism $\til{B} \to B$ is a 
quasi-isomorphism. Also $\til{B}$ is semi-free DG module over $A$ and over 
$A^{\mrm{op}}$. According to Theorem \ref{thm:2363}(2) there is an isomorphism
\begin{equation} \label{eqn:4560}
\K \ot^{\mrm{L}}_{B} \K \cong \K \ot^{\mrm{L}}_{\til{B}} \K
\end{equation}
in $\dcat{D}(\K)$. 

By the associativity of the derived tensor product, and because $\K$ is a DG 
$A$-module via the DG ring homomorphism 
$A \to \til{B}$, there are isomorphisms
\begin{equation} \label{eqn:4550}
\K \ot^{\mrm{L}}_{A} \K \cong
\K \ot^{\mrm{L}}_{A} (\til{B} \ot^{\mrm{L}}_{\til{B}} \K) \cong
(\K \ot^{\mrm{L}}_{A} \til{B}) \ot^{\mrm{L}}_{\til{B}} \K
\end{equation}
in $\dcat{D}(\K)$. 
Now the image of the element $a$ in $\K$ is zero, so we obtain these 
isomorphisms
\[ \K \ot^{\mrm{L}}_{A} \til{B} \cong \K \ot_{A} \til{B} \cong 
\bigl( \K \xar{0 \cd (-)} \K \bigr) \cong \K[1] \oplus \K \]
in $\dcat{D}(\til{B}^{\mrm{op}})$.
Substituting this into (\ref{eqn:4550}) we get 
\[ \K \ot^{\mrm{L}}_{A} \K \cong 
\bigl( \K[1] \oplus \K \bigr) \ot^{\mrm{L}}_{\til{B}} \K
\cong (\K \ot^{\mrm{L}}_{\til{B}} \K)[1] \oplus 
(\K \ot^{\mrm{L}}_{\til{B}} \K) \]
in $\dcat{D}(\K)$. 
Hence 
\begin{equation} \label{eqn:4551}
\opn{H}^i \bigl( \K \ot^{\mrm{L}}_{A} \K \bigr) \cong 
\opn{H}^{i + 1} \bigl( \K \ot^{\mrm{L}}_{B} \K \bigr) \oplus 
\opn{H}^{i} \bigl( \K \ot^{\mrm{L}}_{B} \K \bigr) 
\end{equation}
as $\K$-modules. This, together with (\ref{eqn:4565}), imply that 
$\opn{con} \bigl( \opn{H}(\K \ot^{\mrm{L}}_{A} \K) \bigr) = [-n, 0]$,
as required. 

\medskip \noindent
(2) Now $B$ is an AS regular graded ring of dimension $n - 1$.
As explained in part (1) of the proof, we may neglect the algebraic grading, 
and we only need to prove that 
\begin{equation} \label{eqn:4570}
\opn{RHom}_{A}(\K, A) \cong \K[-n]
\end{equation}
in $\dcat{D}(\K)$.

Let $\til{B}$ be the DG ring from above. Because $\K$ is a DG $A$-module via 
the DG ring homomorphism $A \to \til{B}$, there is an adjunction isomorphism
\begin{equation} \label{eqn:4566} 
\opn{RHom}_{A}(\K, A) \cong 
\opn{RHom}_{\til{B}} \bigl( \K, \opn{RHom}_{A}(\til{B}, A) \bigr) 
\end{equation}
in $\dcat{D}(\K)$. 
Now $\til{B}$ is a semi-free DG module over $A$, so 
\[ \opn{RHom}_{A}(\til{B}, A) \cong 
\opn{Hom}_{A}(\til{B}, A) \cong 
\bigl( A \xar{(-) \cd a} A \bigr) \cong B[-1] \]
in $\dcat{D}(\til{B})$. Substituting this into (\ref{eqn:4566}), and using 
Theorem \ref{thm:2363}(2), we get 
\[ \opn{RHom}_{A}(\K, A) \cong 
\opn{RHom}_{\til{B}} \bigl( \K, B[-1] \bigr) 
\cong \opn{RHom}_{B}(\K, B)[-1] \]
in $\dcat{D}(\K)$.
The AS condition for $B$ says that 
$\opn{RHom}_{B}(\K, B) \cong \K[-n + 1]$ 
in $\dcat{D}(\K)$. Hence the isomorphism (\ref{eqn:4570}) holds. 
\end{proof}

\begin{rem} \label{rem:4560}
The proof of the theorem above can actually yield more, with some 
modifications. If we were to make the DG ring $\til{B}$ 
into a bigraded object, with both algebraic and cohomological grading, then we 
would obtain an isomorphism (\ref{eqn:4560}) in the derived category 
$\dcat{D}(\til{B}, \mrm{gr})$ of bigraded DG $\til{B}$-modules. 
The problem is that we did not develop this kind of theory in the book...

The outcome of this more refined construction would be that formula 
(\ref{eqn:4551}) is an isomorphism in $\dcat{M}(\K, \mrm{gr})$. 
This would give the ability to compare the AS indexes of $A$ and $B$. As can 
be guessed (cf.\ the easy case of Example \ref{exa:3990}, taking $a$ to be
one of the variables), if the AS index of $A$ is $l$, then the AS index of $B$ 
is $l - \opn{deg}(a)$. 
\end{rem}

\begin{thm} \label{thm:4562}
Assume Setup \tup{\ref{set:4565}}.
\begin{enumerate}
\item If $A$ is a Gorenstein graded ring of graded injective dimension 
$n$, then $B$ is a Gorenstein graded ring of graded injective dimension at most
$n - 1$. 

\item Moreover, if $A$ is an AS Gorenstein graded ring of dimension $n$, 
then $B$ is an AS Gorenstein graded ring of dimension $n - 1$.
\end{enumerate}
\end{thm}

\begin{exer} \label{exer:4560}
Prove Theorem \ref{thm:4562}. (Hint: Study the proof of Theorem 
\ref{thm:4546}.)
\end{exer}

\begin{rem} \label{rem:4653}
Here are some historical notes on the material in this subsection. 
Around 1985, M. Artin and W. Schelter \cite{ArSchl} began a systematic study of 
connected graded noncommutative rings (better known as graded algebras), with 
the hope of obtaining a classification in low dimensions, in analogy to the 
commutative projective geometry of curves and surfaces. This is how the concept 
of {\em Artin-Schelter regular graded rings} was born. 
Closely related concepts are {\em Sklyanin algebras} and {\em Koszul algebras}. 

The area of research initiated by Artin and Schelter is often called 
{\em NC algebraic geometry}. Indeed, in 1994 Artin and J.J. Zhang \cite{ArZh} 
introduced the concept of {\em NC projective scheme}. Given a noetherian 
connected graded ring $A$, let 
$\dcat{M}_{\mrm{tor}}(A, \mrm{gr}) \sub \dcat{M}(A, \mrm{gr})$
be the full subcategory on the $\m$-torsion modules. The quotient abelian 
category $\dcat{M}(A, \mrm{gr}) / \dcat{M}_{\mrm{tor}}(A, \mrm{gr})$
is viewed as the category of quasi-coherent sheaves on an imaginary scheme  
$\opn{Proj}(A)$. Note that if $A$ happens to be commutative, then, by the 
famous theorem of J.P. Serre, this quotient category is canonically equivalent 
to $\cat{QCoh} X$, where $X = \opn{Proj}(A)$ is the genuine projective scheme. 
The geometric metaphor of Artin and Zhang turned out to be very useful in 
predicting ``geometric properties'' of connected graded rings, that were later 
proved algebraically, working in $\dcat{M}(A, \mrm{gr})$ or its derived 
category $\dcat{D}(A, \mrm{gr})$. 

For AS regular graded rings of dimension $\leq 3$ generated in degree $1$, the 
classification was pretty much complete by 1991, in papers by Artin, jointly 
with J. Tate and M. Van den Bergh; see the papers \cite{ATV} and 
\cite{ATV2}, and the survey \cite{StaVdB}. Research on AS graded rings of 
higher dimensions, and with generators in various degrees, was later conducted 
by J.T. Stafford, S.P. Smith, Zhang, and their students and collaborators. 

In 2006, V. Ginzburg published the highly influential preprint \cite{Ginz} on 
CY (Calabi-Yau) algebras. Since then, the research on AS regular rings and 
their relatives has merged with the research on CY rings and categories. 
(It turns out that for connected graded rings, the twisted CY rings are 
precisely the AS regular rings.) See Remark \ref{rem:4505} for more on CY 
theory. Two of the latest papers on graded CY rings are 
\cite{ReRoZh} and \cite{ReRo}. 

In the subsequent sections of the book we shall deal with a few ideas that 
emerged as by-products of the Artin school of NC algebraic geometry: the $\chi$ 
condition, balanced dualizing complexes and rigid dualizing complexes.
\end{rem}

\cleardoublepage
\mysection{Derived Torsion over NC Graded Rings}
\label{sec:der-tors-conn}

\AYcopyright

As already mentioned, connected graded noncommutative rings are similar 
in several aspects to complete local commutative rings. 
In this section we concentrate on the derived $\m$-torsion functors 
associated to a connected graded ring $A$ with augmentation ideal $\m$. 

Some of the definitions and results here are quite old, going back to papers 
from the 1990's. Other definitions and results are adaptations to the 
connected graded setting of the work in the recent paper \cite{VyYe}. 

The {\em $\chi$ condition} on a noncommutative noetherian connected graded ring 
$A$ was introduced by M. Artin and J.J. Zhang \cite{ArZh} to guarantee that 
the {\em noncommutative projective scheme} $\opn{Proj}(A)$ will have good 
``geometric properties'', resembling the classical commutative case. In this 
section and the subsequent one, we will focus on another property that the 
$\chi$ condition implies: the {\em symmetry of derived torsion}; namely that 
from the point of view of derived torsion, the ring $A$ looks commutative. 
It is important to say that this phenomenon is {\em only seen in the derived 
category of graded bimodules} -- on the elementary level the ring $A$ is often 
terribly noncommutative (cf.\ Examples \ref{exa:4080} and \ref{exa:4083}), 
and the bimodules are far from being central!

In this section we follow Conventions \ref{conv:3700} and \ref{conv:4560}.
So $\K$ is a base field, and all rings are algebraically graded central 
$\K$-rings. The definitions and results from Section \ref{sec:alg-gra-rings} 
will be used here freely. See Remark \ref{rem:3740} regarding the possibility of 
replacing the base field $\K$ with a base commutative ring.

\mysubsection{Quasi-Compact Finite Dimensional Functors} 
\label{subsec:QC-FD-funs}

The content of this subsection is adapted from \cite[Section 1]{VyYe}. 
We state the results in the context of algebraically graded rings and modules, 
since here we need them for the graded $\m$-torsion functor 
$F := \Ga_{\m}$; yet these results hold in greater generality (see 
\cite{VyYe}). 

Recall from Subsection \ref{subsec:alg-gr-mods} that given a graded ring $A$, 
the category of graded $A$-modules is $\dcat{M}(A, \mrm{gr})$. 
This is a $\K$-linear abelian category. The category of complexes in 
$\dcat{M}(A, \mrm{gr})$ is $\dcat{C}(A, \mrm{gr})$, 
the homotopy category is $\dcat{K}(A, \mrm{gr})$, and the derived category is 
$\dcat{D}(A, \mrm{gr})$. 

In this subsection we consider the following setup:

\begin{setup} \label{set:4075}
$A$ and $B$ are graded rings, and 
$F : \dcat{M}(A, \mrm{gr}) \to \dcat{M}(B, \mrm{gr})$
is a left exact linear functor. 
\end{setup}

The functor $F$ extends to a triangulated functor
$F : \dcat{K}(A, \mrm{gr}) \to \dcat{K}(B, \mrm{gr})$
between the homotopy categories. Because there are enough K-injective 
resolutions, $F$ has a right derived functor 
$(\mrm{R} F, \eta^{\mrm{R}})$. Let us recall what this means: 
$\mrm{R} F : \dcat{D}(A, \mrm{gr}) \to \dcat{D}(B, \mrm{gr})$
is a triangulated functor, and 
$\eta^{\mrm{R}} : \opn{Q} \circ \, F \to \mrm{R} F \circ \opn{Q}$
is a morphism of triangulated functors 
$\dcat{K}(A, \mrm{gr}) \to \dcat{D}(B, \mrm{gr})$
that has a certain universal property (see Definition \ref{dfn:131}). 
The right derived functor $(\mrm{R} F, \eta^{\mrm{R}})$
can be constructed using a K-injective presentation: for each complex $M$ 
we choose a K-injective resolution 
$\rho_M : M \to I_M$
in $\dcat{C}_{\mrm{str}}(A, \mrm{gr})$, and then we take 
$\mrm{R} F(M) := F(I_M)$ and $\eta^{\mrm{R}}_M := F(\rho_M)$. 

Recall that the functor $F$ has the classical right derived functors  
$\mrm{R}^q F : \lb \dcat{M}(A, \mrm{gr}) \to \dcat{M}(B, \mrm{gr})$.
These can be expressed as follows: 
$\mrm{R}^q F(M) = \opn{H}^q(\mrm{R} F(M))$ 
for a module $M \in \dcat{M}(A, \mrm{gr})$. 
Since $F$ is left exact, the canonical homomorphism 
$F \to \mrm{R}^0 F$ is an isomorphism. 

Here is a classical definition:

\begin{dfn} \label{dfn:3885}
The {\em right cohomological dimension of $F$} 
\index{Cohomological dimension! right {\indash} of functor}
is 
\[ d := \opn{sup} \, \{ q \in \N \mid \mrm{R}^q F \neq 0 \} \in 
\N \cup \{ \pm \infty \} . \]
\end{dfn}

If $F \neq 0$ then $d \in \N \cup \{ \infty \}$, but for $F = 0$ the 
dimension is $d = -\infty$. The generalized integer $d$ is said to be finite 
if $d < \infty$; this terminology is designed so that the zero functor will 
have finite right cohomological dimension. 

Here is another standard definition. 

\begin{dfn} \label{dfn:3887}
A module $I \in \dcat{M}(A, \mrm{gr})$ is called a {\em right $F$-acyclic 
object} if $\mrm{R}^q F(I) = 0$ for all $q > 0$. 
\end{dfn}

Of course every injective object of $\dcat{M}(A, \mrm{gr})$ is right 
$F$-acyclic. But often there are many more right $F$-acyclic objects.

\begin{dfn} \label{dfn:3890}
A complex $I \in \dcat{C}(A, \mrm{gr})$ is called a {\em right $F$-acyclic 
complex} if the morphism $\eta^{\mrm{R}}_I : F(I) \to \mrm{R} F(I)$
in $\dcat{D}(B, \mrm{gr})$ is an isomorphism. 
\end{dfn}

Of course a K-injective complex is right $F$-acyclic, but often there 
are many more.
It is easy to see that a module 
$I \in \dcat{M}(A, \mrm{gr})$ is a right $F$-acyclic object
in the sense of Definition \ref{dfn:3887} if and only if it is right 
$F$-acyclic 
as a complex, in the sense of Definition \ref{dfn:3890}. 

\begin{lem} \label{lem:3885}
Let $I \in \dcat{C}(A, \mrm{gr})$ be a complex such that each of the modules  
$I^q$ is right $F$-acyclic. If $I$ is a bounded below complex, or if $F$ has 
finite right cohomological dimension, then $I$ is a right $F$-acyclic 
complex. 
\end{lem}

\begin{proof}
This is contained in the proof of \cite[Corollary I.5.3]{RD}, but for the sake 
of completeness we give a proof here. 

We can assume that neither $F$ nor $I$ are zero. 
First let us assume that $I$ is bounded below.
We can find a quasi-isomorphism $I \to J$ in 
$\dcat{C}_{\mrm{str}}(A, \mrm{gr})$, where $J$ is a bounded below complex 
of injective objects of $\dcat{M}(A, \mrm{gr})$, and thus it is a K-injective 
complex. We must show that 
$F(I) \to F(J)$ is a quasi-isomorphism. 
This amounts to showing that the complex $F(K)$ is acyclic, where $K$ is the 
standard cone on the quasi-isomorphism $I \to J$. Thus we reduce the problem to 
proving that for an acyclic bounded below complex $I$ made up of right 
$F$-acyclic objects, the complex $F(I)$ is acyclic. 

Say $\opn{inf}(I) = d_0 \in \Z$. For any $q$ let 
$\opn{Z}^q(I) := \opn{Ker}(I^q \to I^{q + 1})$, the object of degree 
$q$ cocycles. The acyclicity of the complex $I$ says that we have 
short exact sequences
\[ 0 \to \opn{Z}^{q}(I) \to I^{q} \to \opn{Z}^{q + 1}(I) \to 0 . \]
By induction on $q \geq d_0 - 1$, using the long exact cohomology sequence, we 
prove that each $\opn{Z}^{q}(I)$ is a right $F$-acyclic object.
Therefore the sequences 
\[ 0 \to F(\opn{Z}^{q}(I)) \to F(I^{q}) \to F(\opn{Z}^{q + 1}(I)) \to 0  \]
are exact for all $q$. Splicing together three sequences like this, for  
$q - 1$, $q$ and $q + 1$, shows that the sequence 
$F(I^{q - 1}) \to F(I^q) \to F(I^{q + 1})$
is exact. Thus the complex $F(I)$ is acyclic. 

Now $I$ is no longer bounded below, but the  right cohomological dimension
of $F$ is finite, say $d \in \N$. 
Let $I \to J$ be a quasi-isomorphism in 
$\dcat{C}_{\mrm{str}}(A, \mrm{gr})$, where 
$J$ is a K-injective complex made up of injective objects
(see Corollary \ref{cor:3830}(3)). We must show that 
$F(I) \to F(J)$ is a quasi-isomorphism. As above, by replacing $I$ with the 
standard cone on $I \to J$, we reduce the problem to proving 
that for an acyclic complex $I$ made up of right $F$-acyclic objects, the 
complex $F(I)$ is acyclic. 

Fix an an integer $q$. Let $M := \opn{Z}^q(I)$, and let 
\[ J := (\cdots \to 0 \to I^q \to  I^{q + 1} \to \cdots) , \]
the complex with $I^q$ placed in degree $0$ (so $J$ is a shift of a stupid 
truncation of $I$). By the previous part of the proof 
we know that $\mrm{R} F(J) \cong F(J)$ in $\dcat{D}(B, \mrm{gr})$. 
We have a quasi-isomorphism $M \to J$. Hence, for every $p > d$, 
by our assumption on $F$ we have
\[ 0 = \opn{H}^p(\mrm{R} F(M)) \cong \opn{H}^p(\mrm{R} F(J)) \cong
\opn{H}^p(F(J)) \cong \opn{H}^{p + q}(F(I))  . \]
Since $q$ was chosen arbitrarily, we conclude that $F(I)$ is acyclic.
\end{proof}

Recall the 
{\em cohomological dimension}%
\index{Cohomological dimension! of triangulated functor}%
\index{Additive functor! finite dimensional}
of the triangulated functor 
$\mrm{R} F : \dcat{D}(A, \mrm{gr}) \to \dcat{D}(B, \mrm{gr})$
from Definition \ref{dfn:2126}.
In what follows we use the various notions introduced in Subsection
\ref{subsec:way-out}, 
adapted to the graded setting, such as the concentration interval 
$\opn{con}(M) \sub \Z$ of an 
object $M \in \dcat{G}(A, \mrm{gr})$. 

\begin{lem} \label{lem:3891}
The cohomological dimension of the right derived functor $\mrm{R} F$
equals the right cohomological dimension of $F$. 
\end{lem}

\begin{proof}
As explained in Example \ref{exa:2135}, the cohomological dimension of 
$(\mrm{R} F)|_{\dcat{M}(A)}$ equals the right cohomological dimension of the 
functor $F$. We need to prove the reverse inequality. 

We can assume that $F \neq 0$, and it has finite right cohomological 
dimension, say $d \in \N$. We shall prove that the cohomological displacement 
of the derived functor $\mrm{R} F$ is contained in the interval
$[d_0, d_1] := [0, d]$. 
To be explicit, we shall prove that for every complex
$M \in \dcat{D}(A, \mrm{gr})$
the integer interval 
$\opn{con} \bigl( \opn{H}(\mrm{R} F(M)) \bigr)$ is contained in the interval 
$\opn{con} \bigl( \opn{H}(M) \bigr) + [0, d]$. This is done by cases, and we 
can assume that $M \neq 0$. 

\medskip \noindent 
Case 1. If $M$ is not in $\dcat{D}^-(A, \mrm{gr})$, 
i.e.\ $\opn{con}(\opn{H}(M)) = [d_0, \infty]$ for some 
$d_0 \in \Z \cup \{ -\infty \}$, then we can take a K-injective resolution 
$M \to I$ such that $\opn{inf}(I) = d_0$. Then 
$\opn{H}(\mrm{R} F(M)) = \opn{H}(F(I))$
is concentrated in $[d_0, \infty]$. 

\medskip \noindent 
Case 2. Here $M \in \dcat{D}^-(A, \mrm{gr})$, so that 
$\opn{con}(\opn{H}(M)) = [d_0, d_1]$ for some 
$d_0 \in \Z \cup \{ -\infty \}$ and $d_1 \in \Z$. 
We now take a  K-injective resolution $M \to I$ such that $\opn{inf}(I) = d_0$
and $I$ is made up of injective objects. We must prove that 
$\opn{H}^q(F(I)) = 0$ for all $q > d_1 + d$. This is done like in the 
proof of Lemma \ref{lem:3885}. 
Let $N := \opn{Z}^{d_1}(I)$, and let 
\[ J := (\cdots \to 0 \to I^{d_1} \to  I^{d_1 + 1} \to \cdots) , \]
the complex with $I^{d_1}$ placed in degree $0$. 
So there is a quasi-isomorphism $N \to J$, and it is an injective resolution of 
$N$. Hence for every $p > d$ we have
\[ 0 = \opn{H}^p(\mrm{R} F(N)) \cong \opn{H}^p(F(J)) \cong 
\opn{H}^{d_1 + p}(F(I))  . \qedhere \]
\end{proof}

\begin{dfn} \label{dfn:3888}
Let $\cat{M}$ and $\cat{N}$ be linear categories with arbitrarily direct sums. 
A linear functor 
$G : \cat{M} \to \cat{N}$ is called {\em quasi-compact} 
\index{Additive functor! quasi-compact}
if it commutes with infinite direct sums. Namely, for every collection 
$\{ M_x \}_{x \in X}$ of objects of $\cat{M}$, the canonical morphism
\[ \bigoplus\nolimits_{x \in X} G(M_x) \to 
G \left( \bigoplus\nolimits_{x \in X} M_x \right) \]
in $\cat{N}$ is an isomorphism. 
\end{dfn}

The name ``quasi-compact functor'' is inspired by the property of pushforward 
of quasi-coherent sheaves along a quasi-compact map of schemes. 

\begin{lem} \label{lem:3890}
Assume that the functors $\mrm{R}^q F$ are quasi-compact, for all $q > 0$. Let 
$\{ I_x \}_{x \in X}$ be a collection of right $F$-acyclic objects of 
$\dcat{M}(A, \mrm{gr})$.
Then the object 
$I := \bigoplus_{x \in X} I_x \in \dcat{M}(A, \mrm{gr})$
is right $F$-acyclic.
\end{lem}

\begin{proof}
Take some $q > 0$. Because $\mrm{R}^qF$ is quasi-compact, the canonical 
homomorphism 
$\bigoplus_{x \in X} \mrm{R}^q F(I_x) \to \mrm{R}^q F(I)$
in $\cat{N}$ is an isomorphism. But by assumption, $\mrm{R}^q F(I_x) = 0$ 
for all $x$. 
\end{proof}

\begin{lem} \label{lem:1058}
Assume the functors $\mrm{R}^q F$ are quasi-compact for all $q \geq 0$. Let 
$\{ I_x \}_{x \in X}$ be a collection of right $F$-acyclic 
complexes in $\dcat{C}(A, \mrm{gr})$, and define 
$I := \bigoplus_{x \in X} I_x \in \dcat{C}(A, \mrm{gr})$.
If $I$ is a bounded below complex, or if $F$ has finite right cohomological 
dimension, then $I$ is a right $F$-acyclic complex.
\end{lem}

\begin{proof}
For every index $x$ let us choose a quasi-isomorphism 
$\phi_x : I_x \to J_x$, where $J_x$ is a K-injective complex consisting of 
injective objects, and 
$\opn{inf}(J_x) \geq \opn{inf}(I_x)$. Define 
$J := \bigoplus_{x \in X} J_x$.
We get a quasi-isomorphism 
$\phi := \bigoplus_{x \in X} \phi_x : I \to J$
in $\dcat{C}_{\mrm{str}}(A, \mrm{gr})$. 
By construction, if $I$ is a bounded below complex, then so is $J$. 

For every $x \in X$ there is a commutative diagram 
\[ \UseTips \xymatrix @C=9ex @R=6ex {
F(I_x)
\ar[r]^{F(\phi_x)}
\ar[d]_{\eta^{\mrm{R}}_{I_x}}
& 
F(J_x)
\ar[d]^{\eta^{\mrm{R}}_{J_x}}
\\
\mrm{R} F(I_x)
\ar[r]^(0.5){\mrm{R} F(\phi_x)}
& 
\mrm{R} F(J_x)
} \]
in $\dcat{D}(B, \mrm{gr})$. The vertical arrows are isomorphisms because both 
$I_x$ and $J_x$ are right $F$-acyclic complexes. The morphism $\mrm{R} 
F(\phi_x)$ is 
also an isomorphism. It follows that $F(\phi_x)$ is an isomorphism in 
$\dcat{D}(B, \mrm{gr})$; and therefore it is a quasi-isomorphism in 
$\dcat{C}_{\mrm{str}}(B, \mrm{gr})$. 

Next consider this commutative diagram in 
$\dcat{C}_{\mrm{str}}(B, \mrm{gr})$~:
\[ \UseTips \xymatrix @C=14ex @R=6ex {
\bigoplus_{x \in X} \, F(I_x)
\ar[r]^{ \bigoplus_{x \in X} \, F(\phi_x) }
\ar[d]_{}
& 
\bigoplus_{x \in X} \, F(J_x)
\ar[d]
\\
F(I)
\ar[r]^{ F(\phi) }
& 
F(J)
} \]
Because each $F(\phi_x)$ is a quasi-isomorphism, it follows that 
the top horizontal arrow is a quasi-isomorphism.
Because the functor $F = \mrm{R}^0 F$ is quasi-compact, the vertical arrows 
are isomorphisms. We conclude that $F(\phi)$ is a quasi-isomorphism in 
$\dcat{C}_{\mrm{str}}(B, \mrm{gr})$.

Finally we look at this  commutative diagram in $\dcat{D}(B, \mrm{gr})$~:
\[ \UseTips \xymatrix @C=8ex @R=6ex {
F(I)
\ar[r]^{F(\phi)}
\ar[d]_{\eta^{\mrm{R}}_I}
& 
F(J)
\ar[d]^{\eta^{\mrm{R}}_J}
\\
\mrm{R} F(I)
\ar[r]^{ \mrm{R} F(\phi) }
& 
\mrm{R} F(J)
} \]
We know that the morphisms $F(\phi)$ and $\mrm{R} F(\phi)$ are isomorphisms. 
For every $p$, the object  
$J^p = \bigoplus_{x \in X} \, J_x^p$ 
in $\dcat{M}(A, \mrm{gr})$ is a direct sum of injective objects. So according 
to 
Lemma \ref{lem:3890}, $J^p$ is a right $F$-acyclic object.
By Lemma \ref{lem:3885}, in either of the two cases, the complex $J$ is a 
right $F$-acyclic complex.
This says that the morphism $\eta^{\mrm{R}}_J$ is an isomorphism.
We conclude that the morphism $\eta^{\mrm{R}}_I$ is an isomorphism too, and 
this says that $I$ is a right $F$-acyclic complex.
\end{proof}

\begin{thm} \label{thm:1001}
Under Setup \tup{\ref{set:4075}}, assume the functor $F$
has finite right cohomological dimension, and the functors $\mrm{R}^q F$, for 
all $q \geq 0$, are quasi-compact%
\index{Additive functor! quasi-compact}%
\index{Additive functor! finite dimensional}.
Then the triangulated functor 
$\mrm{R} F : \dcat{D}(A, \mrm{gr}) \to \dcat{D}(B, \mrm{gr})$
has finite cohomological dimensional and is quasi-compact.
\end{thm}

\begin{proof}
By Lemma \ref{lem:3891} the functor $\mrm{R} F$ is finite 
dimensional. We need to  prove that $\mrm{R} F$ commutes with 
infinite direct sums. Namely, consider a collection \lb $\{ M_x \}_{x \in X}$ 
of complexes in $\dcat{C}(A, \mrm{gr})$, and let 
$M := \bigoplus_{x \in X} M_x$. 
We have to prove that the canonical morphism 
$\bigoplus_{x \in X} \mrm{R} F(M_x) \to  \mrm{R} F(M)$
in $\dcat{D}(B, \mrm{gr})$ is an isomorphism. 

For every index $x$ we choose a quasi-isomorphism 
$\phi_x : M_x \to I_x$ in $\dcat{C}_{\mrm{str}}(A, \mrm{gr})$,
where $I_x$ is a K-injective complex. 
Define $I := \bigoplus_{x \in X} I_x$, so there is a quasi-isomorphism 
$\phi : M \to I$ in $\dcat{C}_{\mrm{str}}(A, \mrm{gr})$. We get a commutative 
diagram 
\[ \UseTips \xymatrix @C=14ex @R=6ex {
\bigoplus_{x \in X} \mrm{R} F(M_x)
\ar[r]^{ \bigoplus \mrm{R} F(\phi_x) } 
\ar[d]_{}
& 
\bigoplus_{x \in X} \mrm{R} F(I_x)
\ar[d]
\\
\mrm{R} F(M)
\ar[r]^{ \mrm{R} F(\phi) }
& 
\mrm{R} F(I)
} \]
in $\dcat{D}(B, \mrm{gr})$, in which the horizontal arrows are isomorphisms. 
Therefore it suffices to prove that the canonical morphism 
\begin{equation} \label{eqn:1060}
\bigoplus\nolimits_{x \in X} \mrm{R} F(I_x) \to \mrm{R} F(I) 
\end{equation}
in $\dcat{D}(B, \mrm{gr})$ is an isomorphism. 

Consider this commutative diagram in $\dcat{D}(B, \mrm{gr})$~:
\[ \UseTips \xymatrix @C=12ex @R=6ex {
\bigoplus_{x \in X} \, F(I_x)
\ar[r]^{ \bigoplus \eta^{\mrm{R}}_{I_x} }
\ar[d]_{}
& 
\bigoplus_{x \in X} \, \mrm{R} F(I_x)
\ar[d]
\\
F(I)
\ar[r]^{ \eta^{\mrm{R}}_I }
& 
\mrm{R} F(I)
} \]
For each $x$ the morphism $\eta^{\mrm{R}}_{I_x}$ is an isomorphism; and hence 
the top horizontal arrow is an isomorphism. 
Because $F = \mrm{R}^0 F$ is quasi-compact, the left vertical arrow is an 
isomorphism (in $\dcat{C}_{\mrm{str}}(B, \mrm{gr})$, and so also in 
$\dcat{D}(B, \mrm{gr})$). By Lemma \ref{lem:1058} the complex $I$ is right 
$F$-acyclic, and therefore $\eta^{\mrm{R}}_I$ is 
an isomorphism. Hence the remaining arrow is an isomorphism; but this is the 
morphism (\ref{eqn:1060}).
\end{proof}

\mysubsection{Weakly Stable and Idempotent Copointed Functors}
\label{subsec:idem-cop-fun}

The content of this subsection is adapted from \cite[Section 2]{VyYe}. 
Recall that we follow Conventions \ref{conv:3700} and \ref{conv:4560}. Like 
Subsection \ref{subsec:QC-FD-funs}, the results of this subsection apply in much 
greater generality. We state the results in the special context of 
algebraically graded rings and modules, since we need them for the graded 
$\m$-torsion functor $F := \Ga_{\m}$.

\begin{dfn} \label{dfn:3895}
Let $A$ be a graded ring and 
$F : \dcat{M}(A, \mrm{gr}) \to \dcat{M}(A, \mrm{gr})$
a left exact linear functor. 
\begin{enumerate}
\item The functor $F$ is called {\em stable} 
\index{Additive functor! stable}
if for every injective object $I \in \dcat{M}(A, \mrm{gr})$, the object $F(I)$ 
is injective.  

\item The functor $F$ is called {\em weakly stable} 
\index{Additive functor! weakly stable}
if for every injective object $I \in \dcat{M}(A, \mrm{gr})$, the object $F(I)$ 
is right $F$-acyclic, in the sense of Definition \ref{dfn:3887}.
\end{enumerate}
\end{dfn}

Clearly stable implies weakly stable. 

\begin{rem} \label{rem:3895}
The name ``stable'' comes from the 
theory of {\em torsion classes}; see \cite[Section VI.7]{Ste}.
Indeed, a torsion class is called stable if the corresponding torsion functor 
$F$ is a stable functor, as defined above. 

The name ``weakly stable'' was coined in \cite{VyYe}, where it was shown that 
for a finitely generated ideal $\a$ in a commutative ring $A$, the torsion 
functor $\Ga_{\a}$ is weakly stable iff the ideal $\a$ is {\em weakly 
proregular}. We should mention the well-known fact that if the commutative 
ring $A$ is noetherian, then for every ideal $\a \sub A$ the torsion functor 
$\Ga_{\a}$ is stable (cf.\ \cite[Lemma III.3.2]{Har}; or use the Matlis 
classification of injective $A$-modules in Subsection \ref{subsec:mr-inj-res}).
\end{rem}

\begin{dfn} \label{dfn:3896} 
Let $\cat{M}$ be a linear category (e.g.\ $\dcat{M}(A, \mrm{gr})$ or
$\dcat{D}(A, \mrm{gr})$ for a graded ring $A$), with identity functor 
$\opn{Id}_{\cat{M}}$.
\begin{enumerate}
\item A {\em copointed linear functor}%
\index{Additive functor! copointed}
on $\cat{M}$ is a pair 
$(F, \si)$, consisting of a linear functor 
$F : \cat{M} \to \cat{M}$ and a morphism of functors 
$\si : F \to \opn{Id}_{\cat{M}}$.

\item The copointed linear functor $(F, \si)$ is called 
{\em idempotent}%
\index{Additive functor! idempotent}
if the morphisms 
\[ \si_{F(M)}, \, F(\si_M) : F(F(M)) \to F(M) \]
are isomorphisms for all objects $M \in \cat{M}$. 

\item If $\cat{M}$ is a triangulated category, $F$ is a triangulated functor, 
and $\si$ is a morphism of triangulated functors, then we call 
$(F, \si)$ a {\em copointed triangulated functor}. 
\end{enumerate}
\end{dfn}

The name ``copointed'' is explained in Remark \ref{rem:3896}. 

Weak stability and idempotence together have the following effect. 

\begin{lem} \label{lem:3895}
Let $(F, \si)$ be an idempotent copointed linear functor on 
$\dcat{M}(A, \mrm{gr})$, and assume that $F$ is left exact and weakly stable. 
If $I$ is a right $F$-acyclic object of $\dcat{M}(A, \mrm{gr})$, then $F(I)$ is 
also a right $F$-acyclic object of $\dcat{M}(A, \mrm{gr})$.
\end{lem}

\begin{proof}
Choose an injective resolution $\rho : I \to J$; i.e.\ $J$ is a complex of 
injectives concentrated in nonnegative degrees, and $\rho$ is a 
quasi-isomorphism in $\dcat{C}_{\mrm{str}}(A, \mrm{gr})$. 
Since $\opn{H}^q(F(J)) \cong \mrm{R}^q F(I)$, and since $I$ is a right 
$F$-acyclic object, we see that the homomorphism of complexes
$F(\rho) : F(I) \to F(J)$ is a quasi-isomorphism. 
Therefore both $F(\rho)$ and $\mrm{R} F(F(\rho))$ are isomorphisms in 
$\dcat{D}(A, \mrm{gr})$. 

The weak stability of $F$ implies that $F(J)$ is a bounded below complex of 
right $F$-acyclic objects. According to Lemma \ref{lem:3885} the 
complex $F(J)$ is right $F$-acyclic. This means that the morphism 
$\eta^{\mrm{R}}_{F(J)} : F(J) \to \mrm{R} F (J)$
is an isomorphism. 

The idempotence of the copointed functor $F$ says that the morphisms 
$\si_{F(I)}$ and $\si_{F(J)}$ are isomorphisms in 
$\dcat{C}_{\mrm{str}}(A, \mrm{gr})$. 
We get a commutative diagram in $\dcat{D}(A, \mrm{gr})$~: 
\[ \UseTips \xymatrix @C=10ex @R=6ex {
F(I)
\ar[d]_{F(\rho)}^{\cong}
&
F(F(I))
\ar[l]_{\si_{F(I)}}^{\cong}
\ar[r]^{\eta^{\mrm{R}}_{F(I)}}
\ar[d]_{F(F(\rho))}
&
\mrm{R} F(F(I))
\ar[d]_{\mrm{R} F(F(\rho))}^{\cong}
\\
F(J)
&
F(F(J))
\ar[l]_{\si_{F(J)}}^{\cong}
\ar[r]^{\eta^{\mrm{R}}_{F(J)}}_{\cong}
&
\mrm{R} F(F(J))
} \]
We conclude that $\eta^{\mrm{R}}_{F(I)}$ is an isomorphism; and this means that 
the object $F(I)$ is right $F$-acyclic. 
\end{proof}

The next lemma is a generalization of \cite[Proposition 3.10]{PSY}.

\begin{lem} \label{lem:1063}
Suppose we are given a copointed linear functor $(F, \si)$ on \lb
$\dcat{M}(A, \mrm{gr})$. 
Then there is a unique morphism 
$\si^{\mrm{R}} : \mrm{R} F \to \opn{Id}_{\dcat{D}(A, \mrm{gr})}$
of triangulated functors from $\dcat{D}(A, \mrm{gr})$ to itself, satisfying 
this condition\tup{:} for every complex 
$M \in \dcat{D}(A, \mrm{gr})$ there is equality 
$\si^{\mrm{R}}_M \circ \eta^{\mrm{R}}_M = \si_M$
of morphisms $F(M) \to M$ in $\dcat{D}(A, \mrm{gr})$. 
\end{lem}

In a commutative diagram:
\begin{equation} \label{eqn:1455}
\UseTips \xymatrix @C=6ex @R=6ex {
F(M)
\ar[r]^{ \eta^{\mrm{R}}_M }
\ar[dr]_{ \si_M }
&
\mrm{R} F(M)
\ar[d]^{ \si^{\mrm{R}}_M }
\\
&
M
} 
\end{equation}

\begin{proof}
The existence of the morphism $\si^{\mrm{R}}$ comes for free from the universal 
property of the right derived functor. Still, for later reference, we give the 
construction. 

For a K-injective complex $I$ the morphism 
$\eta^{\mrm{R}}_I : F(I) \to \mrm{R} F(I)$ in $\dcat{D}(A, \mrm{gr})$ is an 
isomorphism, and we define 
$\si^{\mrm{R}}_I : \mrm{R} F(I) \to I$
to be 
$\si^{\mrm{R}}_I := \si_I \circ (\eta^{\mrm{R}}_I)^{-1}$.

For an arbitrary complex $M$ we choose a quasi-isomorphism 
$\rho : M \to I$ into a K-injective complex, and then we let 
$\si^{\mrm{R}}_M := \rho^{-1} \circ \si^{\mrm{R}}_I \circ \mrm{R} F(\rho)$
in $\dcat{D}(A, \mrm{gr})$. The corresponding commutative diagram in 
$\dcat{D}(A, \mrm{gr})$ is (\ref{eqn:4830}). 
\begin{equation} \label{eqn:4830}
\UseTips \xymatrix @C=10ex @R=6ex {
F(M)
\ar[r]^{\eta^{\mrm{R}}_M}
\ar[d]_{F(\rho)}
\ar@(ur,ul)[rr]^{\si_M}
&
\mrm{R} F(M)
\ar[d]_{\mrm{R} F(\rho)}^{\cong}
\ar[r]^{\si^{\mrm{R}}_M}
&
M
\ar[d]_{\rho}^{\cong}
\\
F(I)
\ar[r]^{\eta^{\mrm{R}}_I}_{\cong}
\ar@(dr,dl)[rr]_{\si_I}
&
\mrm{R} F(I)
\ar[r]^{\si^{\mrm{R}}_I}
&
I
} 
\end{equation}
It is easy to see that the collection of morphisms 
$\{ \si^{\mrm{R}}_M \}_{M \in \dcat{D}(A, \mrm{gr})}$ 
has the desired properties. 
\end{proof}

In this way we obtain a copointed triangulated functor 
$(\mrm{R} F, \si^{\mrm{R}})$ on $\dcat{D}(A, \mrm{gr})$.

\begin{thm} \label{thm:1004}
Let $(F, \si)$ be an  idempotent copointed linear functor on \lb
$\dcat{M}(A, \mrm{gr})$, and 
assume that $F$ is left exact and weakly stable.
\index{Additive functor! idempotent}
\index{Additive functor! copointed}
\index{Additive functor! weakly stable}
\begin{enumerate}
\item The copointed triangulated 
functor $(\mrm{R} F, \si^{\mrm{R}})$ on 
$\dcat{D}^+(A, \mrm{gr})$
is idempotent. 

\item If $F$ has finite right cohomological dimension, then the copointed 
triangulated functor $(\mrm{R} F, \si^{\mrm{R}})$ on 
$\dcat{D}(A, \mrm{gr})$
is idempotent. 
\end{enumerate}
\end{thm}

\begin{proof}
Let $M \in \dcat{D}(A, \mrm{gr})$. Choose a K-injective resolution 
$M \to I$ in $\dcat{C}_{\mrm{str}}(A, \mrm{gr})$, such that 
$I$ is a complex consisting of injective objects of $\dcat{M}(A, \mrm{gr})$, 
and $\opn{inf}(I) = \opn{inf}(\opn{H}(M))$. 
It suffices to prove that the morphisms 
\[ \si^{\mrm{R}}_{\mrm{R} F(I)}, \, \mrm{R} F(\si^{\mrm{R}}_I) : 
\mrm{R} F(\mrm{R} F(I)) \to \mrm{R} F(I) \]
in $\dcat{D}(A, \mrm{gr})$ are isomorphisms.

Note that by our choice, if $M \in \dcat{D}^+(A, \mrm{gr})$ 
then $I$ is a bounded below complex. 
Since each $I^q$ is an injective object, it is right $F$-acyclic. 
Because $F$ is a weakly stable functor, each of the objects $F(I^q)$ 
is right $F$-acyclic too. If the functor $F$ has finite right cohomological 
dimension, or if $I$ is bounded below, Lemma \ref{lem:3885} says that the 
complexes $I$ and $F(I)$ are both right $F$-acyclic complexes.  

Consider the diagram 
\begin{equation} \label{eqn:3897}
\UseTips \xymatrix @C=10ex @R=6ex {
F(F(I))
\ar[d]^{ F(\si_I) }
\ar[r]^{ \eta^{\mrm{R}}_{F(I)} }
&
\mrm{R}F(F(I))
\ar[d]^(0.55){ \mrm{R}F(\si_I) }
\ar[r]^{ \mrm{R}F(\eta^{\mrm{R}}_I) }
&
\mrm{R}F(\mrm{R}F(I))
\ar[d]^(0.55){ \mrm{R}F(\si^{\mrm{R}}_I) }
\\
F(I)
\ar[r]^{ \eta^{\mrm{R}}_{I} }
& 
\mrm{R}F(I)
\ar[r]^{ \opn{id} }
&
\mrm{R}F(I)
}
\end{equation}
in $\dcat{D}(A, \mrm{gr})$. 
The left square is commutative: it is gotten from the vertical
morphism $\si_I : F(I) \to I$, to which we apply in the horizontal direction 
the morphism of functors $\eta^{\mrm{R}} : F \to \mrm{R} F$. The right square 
is also commutative: it comes from applying the functor $\mrm{R}F$ to the 
commutative diagram 
\[ \UseTips \xymatrix @C=8ex @R=6ex {
F(I)
\ar[r]^{\eta^{\mrm{R}}_I}
\ar[d]_{\si_I}
&
\mrm{R}F(I)
\ar[d]^{\si^{\mrm{R}}_I}
\\
I
\ar[r]^{ \opn{id} }
&
I
} \]
that characterizes $\si^{\mrm{R}}_I$. 
Because $I$ and $F(I)$ are right $F$-acyclic complexes, the morphisms 
$\eta^{\mrm{R}}_I$ and $\eta^{\mrm{R}}_{F(I)}$ are isomorphisms. Hence 
$\mrm{R}F(\eta^{\mrm{R}}_I)$ is an 
isomorphism. So the horizontal morphisms in the diagram \ref{eqn:3897} are all 
isomorphisms. We are given that $F$ is idempotent, and thus $F(\si_{I})$ is an
isomorphism.  The conclusion of this discussion is that 
$\mrm{R}F(\si^{\mrm{R}}_I)$ is an isomorphism.

Next, let $\phi : \mrm{R}F(I) \to J$ be an isomorphism in 
$\dcat{D}(A, \mrm{gr})$ 
to a K-injective complex $J$ such that 
$\opn{inf}(J) = \opn{inf}(\opn{H}(\mrm{R}F(I)))$.
We know that 
$\eta^{\mrm{R}}_I : F(I) \to \mrm{R}F(I)$ is an isomorphism, and therefore the 
composed morphism 
$\phi \circ \eta^{\mrm{R}}_I : F(I) \to J$ 
is an isomorphism in $\dcat{D}(A, \mrm{gr})$. 
Let $\til{\psi} : F(I) \to J$ be a 
quasi-isomorphism in $\dcat{C}_{\mrm{str}}(A, \mrm{gr})$ representing 
$\phi \circ \eta^{\mrm{R}}_I$. Since $F(I)$ and $J$ are both right $F$-acyclic 
complexes, it follows that 
$F(\til{\psi}) : F(F(I)) \to F(J)$ 
is a  quasi-isomorphism. Consider the following commutative diagram 
\[ \UseTips \xymatrix @C=8ex @R=6ex {
F(F(I))
\ar[r]^{F(\til{\psi})}
\ar[d]^{ \si_{F(I)} } 
&
F(J)
\ar[r]^{\eta^{\mrm{R}}_{J}}
\ar[d]^{\si_{J}}
& 
\mrm{R}F(J)
\ar[d]^(0.55){ \si^{\mrm{R}}_J }
&
\mrm{R}F(\mrm{R}F(I))
\ar[d]^(0.55){ \si^{\mrm{R}}_{\mrm{R}F(I)} }
\ar[l]_{ \mrm{R}F(\phi) }
\\
F(I)
\ar[r]^{\til{\psi}}
& 
J
\ar[r]^{\opn{id}}
& 
J
& 
\mrm{R}F(I)
\ar[l]_{\phi} 
} \]
in $\dcat{D}(A, \mrm{gr})$. All 
horizontal arrows here are isomorphisms. 
We are given that $F$ is idempotent, and thus $\si_{F(I)}$ is an
isomorphism. The conclusion is that 
$\si^{\mrm{R}}_{\mrm{R}F(I)}$ is an isomorphism.
\end{proof}

We end this subsection with a notion dual to ``copointed functor''. 

\begin{dfn} \label{dfn:3897} 
Let $\cat{M}$ be a linear category (e.g.\ $\dcat{M}(A, \mrm{gr})$ or 
$\dcat{D}(A, \mrm{gr})$), 
with identity functor $\opn{Id}_{\cat{M}}$.
\begin{enumerate}
\item A {\em pointed linear functor} 
\index{Additive functor! pointed}
on $\cat{M}$ is a pair 
$(G, \tau)$, consisting of a linear functor 
$G : \cat{M} \to \cat{M}$ and a morphism of functors 
$\tau : \opn{Id}_{\cat{M}} \to G$.

\item The pointed linear functor $(G, \tau)$ is called {\em idempotent} 
\index{Additive functor! idempotent}
if the morphisms 
\[ \tau_{G(M)}, \, G(\tau_M) : G(M) \to G(G(M)) \]
are isomorphisms for all objects $M \in \cat{M}$. 

\item If $\cat{M}$ is a triangulated category, $G$ is a triangulated functor, 
and $\tau$ is a morphism of triangulated functors, then we call 
$(G, \tau)$ a {\em pointed triangulated functor}. 
\end{enumerate}
\end{dfn}

\begin{rem} \label{rem:3896}
Idempotent copointed functors already appeared in the literature under another 
name: {\em idempotent comonads}. Another name for (nearly) the same notion is 
a (Bousfield) colocalization functor, see  e.g.\ \cite{Kr2}. 
Dually, idempotent pointed functors are the same 
thing as {\em idempotent monads}. See \cite{nLab} for a discussion of these 
concepts. 

In \cite[Section 4.1]{KaSc2}, what we call an idempotent pointed functor is 
called a {\em projector}. It is proved there that for an idempotent 
pointed functor $(G, \tau)$, and an object $M \in \cat{M}$,  
the isomorphisms $\tau_{G(M)}$ and $G(\tau_M)$ are equal. 
The same proof (with arrows reversed) shows that for an idempotent 
copointed functor $(F, \si)$, and an object $M$,  
the morphisms $\si_{F(M)}$ and $F(\si_M)$ are equal. We shall not require these 
facts. 
\end{rem}

\mysubsection{Graded Torsion: Weak Stability and Idempotence} 
\label{subsec:der-tors-graded}

In this subsection we begin the study of derived torsion over a connected 
graded ring. We adhere to Conventions \ref{conv:3700} and \ref{conv:4560}.
Thus $\K$ is a field, and all graded $\K$-rings are 
algebraically graded central $\K$-rings. The symbol $\ot$ means $\ot_{\K}$. 
The enveloping ring of a graded $\K$-ring $A$ is 
the graded $\K$-ring 
$A^{\mrm{en}} = A \ot A^{\mrm{op}}$. 
The category of left algebraically graded $A$-modules is 
$\dcat{M}(A, \mrm{gr})$, and it is a $\K$-linear abelian category.
Its derived category is $\dcat{D}(A, \mrm{gr})$.

Recall the notion of left noetherian connected graded $\K$-ring from 
Definitions \ref{dfn:4262} and \ref{dfn:3703}. In this subsection we will 
assume the following setup:

\begin{setup} \label{setup:3755}
We are given a left noetherian connected graded $\K$-ring $A$, with 
augmentation ideal 
$\m = \bigoplus_{i \geq 1} A_i \sub A$.
We are also given graded $\K$-rings $B$ and $C$.
\end{setup}

The graded rings $B$ and $C$ are auxiliary -- they are just placeholders, 
allowing us to write more general formulas. For instance, 
$\dcat{M}(A \ot B^{\mrm{op}}, \mrm{gr})$
is $\dcat{M}(A^{\mrm{en}}, \mrm{gr})$ when we take $B := A$, but it is 
$\dcat{M}(A, \mrm{gr})$ when $B := \K$. 

Since $A$ is left noetherian, the augmentation ideal $\m$ is a
finite graded left $A$-module. 

We shall often make implicit use of the canonical isomorphisms 
\begin{equation} \label{eqn:3850}
A \ot_{A} M \iso M \quad \tup{and} \quad \opn{Hom}_A(A, M) \iso M
\end{equation}
in $\dcat{M}(A \ot B^{\mrm{op}}, \mrm{gr})$. There are similar canonical 
isomorphisms
\begin{equation} \label{eqn:3790}
A \ot^{\mrm{L}}_{A} M \iso M \quad \tup{and} \quad \opn{RHom}_A(A, M) \iso M
\end{equation}
in $\dcat{D}(A \ot B^{\mrm{op}}, \mrm{gr})$. Also we view
$\K \cong A / \m$ as a graded $A^{\mrm{en}}$-module via the augmentation ring 
homomorphism $A \to \K$. 

We are interested in graded $\m$-torsion. This makes sense for bimodules, and 
hence the slightly complicated definition below. 

\begin{dfn} \label{dfn:3745}
Under Setup \ref{setup:3755}:
\begin{enumerate}
\item Let $M \in \dcat{M}(A \ot B^{\mrm{op}}, \mrm{gr})$. 
An element $m \in M$ is called an {\em $\m$-torsion element} 
\index{Torsion! element}
if $\m^j \cd m = 0$ for $j \gg 0$. 

\item The set of $\m$-torsion elements of $M$ is denoted by 
$\Ga_{\m}(M)$, 
\index{1-Gamma-m@$\Ga_{\m}$}
and it is called the {\em $\m$-torsion submodule} of $M$.

\item We call $M$ an {\em $\m$-torsion module} 
\index{Torsion! module}
if $\Ga_{\m}(M) = M$.

\item We denote by 
$\dcat{M}_{\mrm{tor}}(A \ot B^{\mrm{op}}, \mrm{gr})$
the full subcategory of 
$\dcat{M}(A \ot B^{\mrm{op}}, \mrm{gr})$
on the $\m$-torsion modules. 
\end{enumerate}
\end{dfn}

The fact that $\Ga_{\m}(M)$ is a graded $(A \ot B^{\mrm{op}})$-submodule of $M$ 
is easy to see. Indeed, we can express the torsion submodule as follows: 
\begin{equation} \label{eqn:3772}
\Ga_{\m}(M) = \lim_{j \to} \, 
\opn{Hom}_A( A / \m^{j}, M) \sub \opn{Hom}_A(A, M) , 
\end{equation}
using the identification (\ref{eqn:3850}), 
where $j \geq 1$, $\m^j$ is the $j$-fold product 
$\m \cdots \m \sub A$, 
$A \to A / \m^{j + 1} \to A / \m^{j}$ are the canonical graded ring 
surjections, and the operation $\lim_{j \to}$ is in the category 
$\dcat{M}(A \ot B^{\mrm{op}}, \mrm{gr})$. 

It is clear that 
\[ \Ga_{\m} : \dcat{M}(A \ot B^{\mrm{op}}, \mrm{gr}) \to 
\dcat{M}_{\mrm{tor}}(A \ot B^{\mrm{op}}, \mrm{gr}) \]
is a linear functor%
\index{Torsion! functor},
and the inclusions 
$\Ga_{\m}(M) \sub M$ assemble into a monomorphism 
$\si_{} : \Ga_{\m} \inj \opn{Id}$
of functors from $\dcat{M}(A \ot B^{\mrm{op}}, \mrm{gr})$ to itself. 

The ring $B$ plays an auxiliary role in Definition \ref{dfn:3745}. 
Because there is no mention of torsion for $B^{\mrm{op}}$-modules, the notation 
$\dcat{M}_{\mrm{tor}}(A \ot B^{\mrm{op}}, \mrm{gr})$
is not ambiguous. However, in Subsection \ref{subsec:sym-der-tor}, where 
torsion for $B^{\mrm{op}}$ modules does become an option, we will switch to 
the richer notation 
$\dcat{M}_{(\mrm{tor}, ..)}(A \ot B^{\mrm{op}}, \mrm{gr})$.
See Definition \ref{dfn:3815}. 

\begin{lem} \label{lem:3835}    
Under Setup \tup{\ref{setup:3755}} the following hold.
\begin{enumerate}
\item The functor $\Ga_{\m}$ commutes with graded ring homomorphisms 
$C \to B$, in the sense that the diagram 
\[ \UseTips \xymatrix @C=6ex @R=6ex {
\dcat{M}(A  \ot B^{\mrm{op}}, \mrm{gr})
\ar[r]^(0.45){\Ga_{\m}}
\ar[d]_{\opn{Rest}}
&
\dcat{M}_{\mrm{tor}}(A \ot B^{\mrm{op}}, \mrm{gr})
\ar[d]^{\opn{Rest}}
\\
\dcat{M}(A  \ot C^{\mrm{op}}, \mrm{gr})
\ar[r]^(0.45){\Ga_{\m}}
&
\dcat{M}_{\mrm{tor}}(A \ot C^{\mrm{op}}, \mrm{gr})
} \]
is commutative. 

\item The functor $\Ga_{\m}$ is idempotent, in the sense that for every 
$M \in \lb \dcat{M}(A \ot B^{\mrm{op}}, \mrm{gr})$ 
the homomorphisms 
\[ \si_{\Ga_{\m}(M)} \, , \,   \Ga_{\m}(\si_{M}) \, : \,
\Ga_{\m}(\Ga_{\m}(M)) \to \Ga_{\m}(M) \]
are bijective and are equal to each other.

\item Let $\phi : M \to N$ be a homomorphism in 
$\dcat{M}(A \ot B^{\mrm{op}}, \mrm{gr})$,
such that $M$ is $\m$-torsion. Then 
there is a unique homomorphism 
$\phi' : M \to \Ga_{\m}(N)$ satisfying $\phi = \si_{N} \circ \phi'$.

\item The functor $\Ga_{\m}$ is left exact.

\item The functor $\Ga_{\m}$ commutes with direct limits. 

\item The category $\dcat{M}_{\mrm{tor}}(A \ot B^{\mrm{op}}, \mrm{gr})$
is a thick abelian subcategory of \lb
$\dcat{M}(A \ot B^{\mrm{op}}, \mrm{gr})$,
closed under taking subobjects and quotients.
\end{enumerate}
\end{lem}

\begin{exer} \label{exer:3810}
Prove Lemma \ref{lem:3835}. (Hint: for item (6) you will need the fact that $A$ 
is left noetherian.)
\end{exer}

Recall that for a  module $M \in \dcat{M}(A, \mrm{gr})$, 
its socle is $\opn{Soc}(M) = \opn{Hom}_A(\K, M)$.
The socle is a functor, and there is a monomorphism 
$\opn{Soc} \inj \Ga_{\m}$ 
of functors from $\dcat{M}(A, \mrm{gr})$ to itself.

\begin{lem} \label{lem:5095}
Let $M \in \dcat{M}(A, \mrm{gr})$, and let 
$W := \opn{Soc}(M)$. If $M$ is an $\m$-torsion module, then $W$ is a 
graded-essential submodule of $M$.  
\end{lem}

\begin{proof}
Take any nonzero graded $A$-submodule 
$M' \sub M$. Let $m \in M'$ be some nonzero homogeneous element.
Because $M$ is a torsion module, yet $m \neq 0$, 
there is a unique $j \in \N$ such that 
$\m^{j + 1} \cd m  = 0$ but $\m^{j} \cd m \neq 0$. 
So there is a homogeneous element $a \in \m^{j}$ such that $a \cd m \neq 0$. 
The element $a \cd m$ is then a nonzero homogeneous element in 
$M' \cap W$. 
\end{proof}

The $\K$-linear dual (see Definition \ref{dfn:3725}) of the ring $A$ is the 
graded bimodule 
$A^* = \opn{Hom}_{\K}(A, \K) \in \dcat{M}(A^{\mrm{en}}, \mrm{gr})$.
The augmentation ring homomorphism $A \to \K$ induces, by dualizing, a 
canonical monomorphism 
\begin{equation} \label{eqn:4810}
\K \inj A^* 
\end{equation}
in $\dcat{M}(A^{\mrm{en}}, \mrm{gr})$.
For every $W \in \dcat{M}(\K, \mrm{gr})$ there is an induced monomorphism 
\begin{equation} \label{eqn:5095}
W \cong W \ot \K \inj W \ot A^*  
\end{equation}
in $\dcat{M}(A^{\mrm{en}}, \mrm{gr})$.

\begin{lem} \label{lem:5100}   
Under Setup \tup{\ref{setup:3755}}, let $W \in \dcat{M}(\K, \mrm{gr})$, and 
define 
$I := A^* \ot W \in \dcat{M}(A, \mrm{gr})$.
Consider the canonical monomorphism $W \inj I$ 
in $\dcat{M}(A, \mrm{gr})$ from \tup{(\ref{eqn:5095})}. 
Then\tup{:}
\begin{enumerate}
\item $I$ is a graded-injective $A$-module.

\item $I$ is an $\m$-torsion $A$-module.

\item $W = \opn{Soc}(I)$.

\item $W$ is a graded-essential submodule of $I$. 
\end{enumerate}
\end{lem}

\begin{proof} \mbox{}

\smallskip \noindent
(1) Say $W = \bigoplus_{x \in X} \K(-i_x)$. Then 
$I \cong \bigoplus_{x \in X} A^*(-i_x)$ in $\dcat{M}(A, \mrm{gr})$.
According to Proposition \ref{prop:3795}(2) and Theorem \ref{thm:3795}(1), the 
graded $A$-module $I$ is graded-injective. 

\medskip \noindent
(2) The $A$-module $A^*$ is $\m$-torsion, and a direct sum of torsion modules 
is torsion (see Lemma \ref{lem:3835}(5)). 

\medskip \noindent
(3) This is because $\K = \opn{Soc}(A^*)$, via the canonical monomorphism 
(\ref{eqn:4810}). 

\medskip \noindent
(4) Combine items (2) and (3) with Lemma \ref{lem:5095}.
\end{proof}

\begin{prop} \label{prop:4811}  
Under Setup \tup{\ref{setup:3755}}, let 
$M \in \dcat{M}(A, \mrm{gr})$, and let 
$W := \opn{Soc}(M)$. The three conditions below are equivalent. 
\begin{enumerate}
\rmitem{i} $M$ is an $\m$-torsion graded $A$-module. 

\rmitem{ii} $W$ is a graded-essential $A$-submodule of $M$. 

\rmitem{iii} There is a graded-essential monomorphism $\ep : M \inj A^* \ot W$
in $\dcat{M}(A, \mrm{gr})$, that restricts to the identity on $W$. 
\end{enumerate}
\end{prop}

\begin{proof} \mbox{}

\smallskip \noindent 
(i) $\Rightarrow$ (ii): This is Lemma \ref{lem:5095}.

\medskip \noindent 
(ii) $\Rightarrow$ (iii): By Lemma \ref{lem:5100}
the graded $A$-module $A^* \ot W$ is 
graded-injective. 
Hence there is a homomorphism 
$\ep : M \to A^* \ot W$ in $\dcat{M}(A, \mrm{gr})$
that makes the diagram 
\[ \UseTips \xymatrix @C=8ex @R=6ex {
*++{W}
\ar@{>->}[r]
\ar@{>->}[d]
&
M
\ar@{-->}[dl]^{\ep}
\\
*+{A^* \ot W}
} \]
in $\dcat{M}(A, \mrm{gr})$ commutative.
Let $M' := \opn{Ker}(\ep)$. 
If $M'$ were nonzero, then, because $W \sub M$ is a graded-essential
submodule, we would have 
$W \cap M' \neq 0$. But $W \to A^* \ot W$ is a monomorphism.
We conclude that $\ep$ is a monomorphism too. 

Finally, by Lemma \ref{lem:5100}(4) we know that $W$ is a 
graded-essential submodule of $A^* \ot W$. It follows that the monomorphism 
$\ep$ is graded-essential.

\medskip \noindent 
(iii) $\Rightarrow$ (i): The graded $A$-module $A^*$ is torsion.
We know that the subcategory 
$\dcat{M}_{\mrm{tor}}(A, \mrm{gr})$ of $\dcat{M}(A, \mrm{gr})$
is closed under taking arbitrary direct sums and subobjects. 
Hence $M$ is $\m$-torsion. 
\end{proof}

Cofinite graded modules were defined in Definition \ref{dfn:3726}. 

\begin{cor} \label{cor:4270}
Assume $A$ is a noetherian connected graded ring. Let 
$M \in \dcat{M}_{\mrm{tor}}(A, \mrm{gr})$. 
The two conditions below are equivalent.
\begin{enumerate}
\rmitem{i} $M$ is a cofinite graded $A$-module.  

\rmitem{ii} $W := \opn{Soc}(M)$ is a finite graded $\K$-module.
\end{enumerate}
\end{cor}

\begin{proof}
Graded Matlis Duality (Theorem \ref{thm:4045}) tells us that 
the cofinite graded $A$-modules are the artinian objects in the abelian 
category $\dcat{M}(A, \mrm{gr})$. 

\medskip \noindent 
(i) $\Rightarrow$ (ii): We are given that $M$ is a cofinite graded $A$-module. 
Hence so is its submodule $W$. Since $W$ is an 
$A$-module via $\K$, it must be an artinian graded $\K$-module. So 
$W$ is finite over $\K$. 

\medskip \noindent 
(ii) $\Rightarrow$ (i): Since $A^*$ is a cofinite graded $A$-module, and 
$A^* \ot W$ is a finite direct sum of degree twists of $A^*$,
we see that $A^* \ot W$ is a cofinite graded $A$-module. By Proposition 
\ref{prop:4811} we know that $M$ is isomorphic to a graded submodule of 
$A^* \ot W$; so $M$ is also a cofinite graded $A$-module.
\end{proof}

\begin{rem} \label{rem:4270}
We do not know if the corollary remains true if $A$ is just left noetherian. 
\end{rem}

Now we pass to derived categories and functors. 

\begin{dfn} \label{dfn:3755}
Under Setup \ref{setup:3755}, we denote by 
$\dcat{D}_{\mrm{tor}}(A \ot B^{\mrm{op}}, \mrm{gr})$
\index{1-Dtor(A-B,gr)@$\dcat{D}_{\mrm{tor}}(A \ot B^{\mrm{op}}, \mrm{gr})$}
the full subcategory of 
$\dcat{D}(A \ot B^{\mrm{op}}, \mrm{gr})$
on the complexes $M$ whose cohomology modules $\opn{H}^p(M)$
belong to $\dcat{M}_{\mrm{tor}}(A \ot B^{\mrm{op}}, \mrm{gr})$
for all $p$.
\end{dfn}

\begin{prop} \label{prop:3850}   
Under Setup \tup{\ref{setup:3755}}  the following hold.
\begin{enumerate}
\item The functor $\Ga_{\m}$ has a triangulated right derived functor 
\[ \mrm{R} \Ga_{\m} : \dcat{D}(A \ot B^{\mrm{op}}, \mrm{gr}) \to
\dcat{D}(A \ot B^{\mrm{op}}, \mrm{gr}) . \]
\index{1-RGamma-m@$\mrm{R} \Ga_{\m}$}
\index{Torsion! derived {\indash} functor}
If $I \in \dcat{D}(A \ot B^{\mrm{op}}, \mrm{gr})$
is K-graded-injective over $A$, then the morphism 
$\eta^{\mrm{R}}_{I} : \Ga_{\m}(I) \to \mrm{R} \Ga_{\m}(I)$
is an isomorphism.

\item There is a unique morphism  
$\si^{\mrm{R}}_{} : \mrm{R} \Ga_{\m} \to \opn{Id}$%
\index{1-SiR@$\si^{\mrm{R}}_{}$}
of triangulated functors from 
$\dcat{D}(A \ot B^{\mrm{op}}, \mrm{gr})$
to itself, such that for every 
$M \in \dcat{D}(A \ot B^{\mrm{op}}, \mrm{gr})$
there is equality 
$\si^{\mrm{R}}_{M} \circ \eta^{\mrm{R}}_{M} = \opn{Q}(\si_M)$
of morphisms $\Ga_{\m}(M) \to M$ in $\dcat{D}(A \ot B^{\mrm{op}}, \mrm{gr})$. 

\item The category $\dcat{D}_{\mrm{tor}}(A \ot B^{\mrm{op}}, \mrm{gr})$
is a full triangulated subcategory of \nl
$\dcat{D}(A \ot B^{\mrm{op}}, \mrm{gr})$, and it contains the image of the 
functor 
$\mrm{R} \Ga_{\m}$. 

\item The functor $\mrm{R} \Ga_{\m}$ commutes with graded $\K$-ring 
homomorphisms $C \to B$, in the sense that the diagram 
\[ \UseTips \xymatrix @C=6ex @R=6ex {
\dcat{D}(A  \ot B^{\mrm{op}}, \mrm{gr})
\ar[r]^(0.45){ \mrm{R}\Ga_{\m}}
\ar[d]_{\opn{Rest}}
&
\dcat{D}_{\mrm{tor}}(A \ot B^{\mrm{op}}, \mrm{gr})
\ar[d]^{\opn{Rest}}
\\
\dcat{D}(A  \ot C^{\mrm{op}}, \mrm{gr})
\ar[r]^(0.45){\mrm{R} \Ga_{\m}}
&
\dcat{D}_{\mrm{tor}}(A \ot C^{\mrm{op}}, \mrm{gr})
} \]
is commutative up to isomorphism. 
Likewise for the morphism of triangulated functors $\si^{\mrm{R}}_{}$.

\item For every $M \in \dcat{D}(A \ot B^{\mrm{op}}, \mrm{gr})$ and $p \in \Z$ 
there is an isomorphism 
\[ \opn{H}^p(\mrm{R} \Ga_{\m}(M)) \cong 
\lim_{j \to} \, \opn{Ext}^p_A( A / \m^{j}, M) \]
in $\dcat{M}(A \ot B^{\mrm{op}}, \mrm{gr})$. This isomorphism is functorial in 
$M$. 
\end{enumerate}
\end{prop}

\begin{proof}
(1) By Corollary \ref{cor:3740}(3) every 
$M \in \dcat{C}_{}(A \ot B^{\mrm{op}}, \mrm{gr})$
has a resolution $\rho : M \to I$ 
in $\dcat{C}_{\mrm{str}}(A \ot B^{\mrm{op}}, \mrm{gr})$
by a complex $I$ that is K-graded-injective over $A$. If 
$\psi : I \to I'$ is a quasi-isomorphism in 
$\dcat{C}_{\mrm{str}}(A \ot B^{\mrm{op}}, \mrm{gr})$ 
between such complexes, 
then it is a homotopy equivalence in 
$\dcat{C}_{\mrm{str}}(A, \mrm{gr})$, 
and hence $\Ga_{\m}(\psi)$ is a quasi-isomorphism. 
We see that the full subcategory 
$\cat{J} \sub \dcat{K}(A \ot B^{\mrm{op}}, \mrm{gr})$
on these complexes $I$ satisfies the conditions of Theorem 
\ref{thm:1470}. Therefore the right derived functor 
$\mrm{R} \Ga_{\m}$ exists, and $\eta^{\mrm{R}}_{I}$ is an isomorphism for 
all $I \in \cat{J}$. 

\medskip \noindent
(2) This is a special case of Lemma \ref{lem:1063}. 

\medskip \noindent
(3) Clear from Lemma \ref{lem:3835}(6). 

\medskip \noindent
(4) This is immediate from items (1) and (2). 

\medskip \noindent
(5) Use formula (\ref{eqn:3772}), and the fact that cohomology commutes with 
direct limits.  
\end{proof}

\begin{dfn} \label{dfn:4045}
Under Setup \tup{\ref{setup:3755}}, the {\em $i$-th local cohomology functor} 
is the functor 
\[ \opn{H}^i_{\m} := \mrm{R}^i \Ga_{\m} : 
\dcat{M}(A \ot B^{\mrm{op}}, \mrm{gr}) \to 
\dcat{M}(A \ot B^{\mrm{op}}, \mrm{gr}) , \]
\index{1-Him@$\opn{H}^i_{\m}$}
the $i$-th right derived functor of $\Ga_{\m}$. 
\end{dfn}

\begin{cor} \label{cor:3836}
For every $i \in \N$ there is an isomorphism 
$\opn{H}^i_{\m} \cong \opn{H}^i \circ \, \mrm{R} \Ga_{\m}$
of functors from $\dcat{M}(A  \ot B^{\mrm{op}}, \mrm{gr})$ to itself. 
\end{cor}

\begin{proof}
This is immediate from Proposition \ref{prop:3850}(1). 
\end{proof}

\begin{prop} \label{prop:4050}   
For every $i \in \N$ the functor 
\[ \opn{H}^i_{\m} : \dcat{M}(A \ot B^{\mrm{op}}, \mrm{gr}) \to 
\dcat{M}(A \ot B^{\mrm{op}}, \mrm{gr}) \]
commutes with direct limits. 
\end{prop}

\begin{proof}
In view of Proposition \ref{prop:3850}(4), with $C := \K$, we can assume that 
$B = \K$, so $A \ot B^{\mrm{op}} = A$. 
Since direct limits commute with each other, by Proposition \ref{prop:3850}(5) 
it suffices to prove that for every $i, j$ the functor 
$\opn{Ext}^i_A( A / \m^{j}, -)$
commutes with direct limits. Now $A$ is left noetherian, so we can find a 
resolution $P \to A / \m^{j}$ by a nonpositive complex, where each $P^i$ is a 
finite graded-free $A$-module. And the functor 
\[ \opn{Ext}^i_A( A / \m^{j}, -) \cong 
\opn{H}^i \circ  \opn{Hom}_A( P, -) : \dcat{M}(A, \mrm{gr}) \to
\dcat{M}(A, \mrm{gr}) \]
commutes with infinite direct limits. 
\end{proof}

Definition \ref{dfn:3885}, when used in our context, says that 
the right cohomological dimension of the functor $\Ga_{\m}$ is 
\begin{equation} \label{eqn:3962}
d := \opn{sup} \, \{ p \in \N \mid \mrm{H}^p_{\m} \neq 0 \} \in 
\N \cup \{ \pm \infty \} .
\end{equation}

The cohomological dimension of a triangulated functor was introduced in 
Definition \ref{dfn:2126}. A triangulated functor is called 
{\em quasi-compact} if it commutes with infinite direct sums (Definition 
\ref{dfn:3888}). 

\begin{thm} \label{thm:4050}  
Under Setup \tup{\ref{setup:3755}}, assume the torsion functor
$\Ga_{\m} : \lb \dcat{M}(A, \mrm{gr}) \to \dcat{M}(A, \mrm{gr})$ 
has finite right cohomological dimension. Then the derived torsion functor 
\[ \mrm{R} \Ga_{\m} : \dcat{D}(A \ot B^{\mrm{op}}, \mrm{gr}) \to 
\dcat{D}(A \ot B^{\mrm{op}}, \mrm{gr}) \]
is quasi-compact and has finite cohomological dimension. 
\index{Additive functor! quasi-compact}
\index{Additive functor! finite dimensional}
\end{thm}

\begin{proof}
According to Proposition \ref{prop:3850}(4), the functor 
\[ \Ga_{\m} : \dcat{M}(A \ot B^{\mrm{op}}, \mrm{gr}) \to 
\dcat{M}(A \ot B^{\mrm{op}}, \mrm{gr}) \]
also has finite right cohomological dimension.
The right derived functors
\[ \mrm{R}^q \Ga_{\m} = \opn{H}^q_{\m} : 
\dcat{M}(A \ot B^{\mrm{op}}, \mrm{gr}) \to 
\dcat{M}(A \ot B^{\mrm{op}}, \mrm{gr}) \]
are quasi-compact by Proposition \ref{prop:4050}. So Theorem \ref{thm:1001} 
applies. 
\end{proof}

\begin{lem} \label{lem:3838}
Let $I \in \dcat{M}(A, \mrm{gr})$, and let 
$W := \opn{Soc}(I)$. The following three conditions are equivalent.
\begin{enumerate}
\rmitem{i} $I$ is graded-injective and $\m$-torsion. 

\rmitem{ii} $I$ is graded-injective and $W$ is an essential submodule of $I$. 

\rmitem{iii} There is an isomorphism 
$I \cong A^* \ot W$ in $\dcat{C}(A, \mrm{gr})$ which is the identity on $W$. 
\end{enumerate}
\end{lem}

In condition (iii) we view $W$ as a submodule of $A^* \ot W$ using the 
canonical embedding (\ref{eqn:5095}).

\begin{proof} \mbox{}

\smallskip \noindent
(i) $\Rightarrow$ (ii): This is Lemma \ref{lem:5095}.

\medskip \noindent
(ii) $\Rightarrow$ (iii): 
By Proposition \ref{prop:4811} there is an essential monomorphism 
$\ep : I \to A^* \ot W$ in $\dcat{M}(A, \mrm{gr})$
that extends the identity of $W$. Because $I$ is graded-injective, the 
monomorphism $\ep$ is split, by some epimorphism $\si : A^* \ot W \to I$. 
But $W$ is a graded-essential submodule of $A^* \ot W$, and therefore $\si$ is 
also a monomorphism. We conclude that both $\si$ and $\ep$ are 
isomorphisms. 

\medskip \noindent
(iii) $\Rightarrow$ (i): 
According to Lemma \ref{lem:5100}, the graded $A$-module $A^* \ot W$ is 
graded-injective and $\m$-torsion.
\end{proof}

\begin{lem} \label{lem:3900}
Let $M \in \dcat{M}(A, \mrm{gr})$ be an $\m$-torsion module, and let $I$ be an 
injective hull of $M$ in $\dcat{M}(A, \mrm{gr})$. Then 
$I \cong A^* \ot W$, where $W$ is the socle of $M$. 
In particular, $I$ is an  $\m$-torsion module.
\end{lem}

\begin{proof}
By Proposition \ref{prop:4811} there is a graded-essential monomorphism 
$\ep : M \inj A^* \ot W$ in $\dcat{M}(A, \mrm{gr})$.
From Lemma \ref{lem:5100}(1) we know that $A^* \ot W$ is a 
graded-injective module over $A$. So by the uniqueness (up to isomorphism) 
of graded-injective hulls, we get an isomorphism 
$I \cong A^* \ot W$. Finally, by Lemma \ref{lem:5100}(2),
$A^* \ot W$ is $\m$-torsion. 
\end{proof}

\begin{lem} \label{lem:3839}
Let $I$ be a graded-injective $A$-module. Then there is an isomorphism
$\Ga_{\m}(I) \cong A^* \ot W$ in $\dcat{M}(A, \mrm{gr})$,
where $W$ is the socle of $I$. 
\end{lem}

\begin{proof}
Consider the essential monomorphism $\ep : W \inj A^* \ot W$ 
from Lemma \ref{lem:5100}(4), and 
the monomorphism $\tau : W \inj I$.
Because $I$ is graded-injective, there is a homomorphism 
$\psi : A^* \ot W \to I$ in $\dcat{M}(A, \mrm{gr})$
such that $\psi \circ \ep = \tau$. Since $\ep$ is a graded-essential 
monomorphism, the usual calculation (see the proof of the implication 
(ii) $\Rightarrow$ (iii) in Proposition \ref{prop:4811}) shows that 
$\psi$ is a monomorphism. Let $J \sub I$ be the image of $\psi$,
so $\psi : A^* \ot W \iso J$. We know (from Lemma \ref{lem:5100}) 
that $J$ is graded-injective and $\m$-torsion. Hence there is a direct sum 
decomposition $I = J \oplus J'$
in $\dcat{M}(A, \mrm{gr})$, and $\Ga_{\m}(J) = J$.

By Proposition \ref{prop:4811} we know that 
$W \cap J' = \opn{Soc}(J')$ is a graded-essential submodule of $\Ga_{\m}(J')$; 
but $W \sub J$, so $W \cap J' = 0$, and therefore $\Ga_{\m}(J') = 0$.
The conclusion is that $\Ga_{\m}(I) = J$. 
\end{proof}

\begin{thm} \label{thm:3990}   
Under Setup \tup{\ref{setup:3755}}, the torsion functor
\[ \Ga_{\m} : \dcat{M}(A, \mrm{gr}) \to \dcat{M}(A, \mrm{gr}) \]
is stable, namely for every graded-injective $A$-module $I$
the $\m$-torsion submodule $\Ga_{\m}(I)$ is graded-injective. 
\end{thm}

\begin{proof}
This is immediate from Lemma \ref{lem:3839} and
Lemma \ref{lem:5100}(1). 
\end{proof}

Here is the graded NC version of Definition \ref{dfn:2210}.

\begin{dfn} \label{dfn:3837} \mbox{}
\begin{enumerate}
\item A {\em minimal complex of graded-injective $A$-modules}%
\index{Complex of algebraically graded modules! minimal {\indash}
graded-injective {\indash}}
is a bounded below complex $I$ of graded-injective $A$-modules, such 
that for every $p$ the submodule $\opn{Z}^p(I) \sub I^p$ is graded-essential.

\item Let $M \in \dcat{D}^+(A, \mrm{gr})$. A 
{\em minimal graded-injective resolution}%
\index{Resolution! minimal graded-injective}
of $M$ is a  quasi-isomorphism $\rho : M \to I$ 
in $\dcat{C}_{\mrm{str}}(A, \mrm{gr})$, where $I$ is a minimal complex of 
graded-injective $A$-modules.
\end{enumerate}
\end{dfn}

\begin{lem} \label{lem:3837}
Let $M \in \dcat{D}^+(A, \mrm{gr})$. Then $M$ has a minimal 
graded-injective resolution.
\end{lem}

\begin{proof}
The proof of Proposition \ref{prop:2210} works here too. 
\end{proof}

There is uniqueness for minimal graded-injective resolutions, but we won't need 
it. 

For the next proposition, and for later results, it will be convenient to use 
this notation:

\begin{notation} \label{not:3965}
If the torsion functor $\Ga_{\m}$ has finite right cohomological dimension, 
then we let $\star := \bra{\tup{empty}}$, the empty boundedness indicator; 
otherwise we let $\star := +$. Thus $\dcat{D}^{\star}$ is either 
$\dcat{D}$ or $\dcat{D}^+$.
\end{notation}

\begin{prop} \label{prop:3965}
With Notation \tup{\ref{not:3965}}, let 
$M \in \dcat{D}^{\star}(A \ot B^{\mrm{op}}, \mrm{gr})$. The following 
conditions are equivalent\tup{:}
\begin{enumerate}
\rmitem{i} $M \in \dcat{D}^{\star}_{\mrm{tor}}(A \ot B^{\mrm{op}}, \mrm{gr})$.

\rmitem{ii} The morphism 
$\si^{\mrm{R}}_{M} : \mrm{R} \Ga_{\m}(M) \to M$
in $\dcat{D}(A \ot B^{\mrm{op}}, \mrm{gr})$
is an isomorphism.
\end{enumerate}
\end{prop}

\begin{proof}
The implication (ii) $\Rightarrow$ (i) is due to Proposition 
\ref{prop:3850}(3).  

For the other implication, we can forget the ring $B$; it is not relevant. 
First let us assume that 
$M \in \dcat{D}^{+}_{\mrm{tor}}(A, \mrm{gr})$.
By Lemma \ref{lem:3900} the injective hull in $\dcat{M}(A, \mrm{gr})$ of 
a torsion module is torsion. 
Then, according to Theorem \ref{thm:2880}, there is a quasi-isomorphism 
$M \to I$ in $\dcat{C}^+_{\mrm{str}}(A, \mrm{gr})$,
where $I$ is a complex of torsion graded-injective $A$-modules.
It follows that homomorphism 
$\si_I : \Ga_{\m}(I) \to I$
in $\dcat{C}^+_{\mrm{str}}(A, \mrm{gr})$, that represents $\si^{\mrm{R}}_{M}$,
is an isomorphism. 

When $\opn{H}(M)$ is not bounded below, but $\Ga_{\m}$ has finite right 
cohomological dimension, we can use smart truncation, as in the proof of 
Theorem \ref{thm:2135}(2).
\end{proof}

\begin{prop} \label{prop:3990}   
With Notation \tup{\ref{not:3965}},
Let $M \in \dcat{D}^{\star}_{\mrm{tor}}(A \ot B^{\mrm{op}}, \mrm{gr})$ and 
$N \in \dcat{D}^+(A \ot C^{\mrm{op}}, \mrm{gr})$.
Then there is an isomorphism 
\[ \opn{RHom}_{A}(M, N) \iso 
\opn{RHom}_{A} \bigl( M, \mrm{R} \Ga_{\m}(N) \bigr) \]
in $\dcat{D}(B \ot C^{\mrm{op}}, \mrm{gr})$, 
and it is functorial in $M$ and $N$. 
\end{prop}

\begin{proof}
Choose a resolution $M \to I_M$ in 
$\dcat{C}^{\star}_{\mrm{str}}(A \ot B^{\mrm{op}}, \mrm{gr})$
where $I_M$ is K-graded-injective over $A$, each $I_M^p$ is 
graded-injective over $A$, and 
$\opn{inf}(I_M) = \opn{inf}(\opn{H}((M))$; this is possible by Corollary 
\ref{cor:3830}(3). Choose the same sort of resolution 
$N \to I_N$ in $\dcat{C}^+_{\mrm{str}}(A \ot C^{\mrm{op}}, \mrm{gr})$.
Then 
\begin{equation} \label{eqn:3990}
\opn{RHom}_{A}(M, N) \cong \opn{Hom}_{A}(I_M, I_N) 
\end{equation}
in $\dcat{D}(B \ot C^{\mrm{op}}, \mrm{gr})$,
and $\mrm{R} \Ga_{\m}(N) \cong \Ga_{\m}(I_N)$ in 
$\dcat{D}(A \ot C^{\mrm{op}}, \mrm{gr})$. By Theorem \ref{thm:3990} each 
$\Ga_{\m}(I^p_N)$ is graded-injective as an $A$-module, so 
$\Ga_{\m}(I_N)$ is a bounded below complex of graded-injective $A$-modules, and 
therefore it is a K-graded-injective complex over $A$. 
We see that  
\begin{equation} \label{eqn:3991}
\opn{RHom}_{A} \bigl( M, \mrm{R} \Ga_{\m}(N) \bigr) \cong
\opn{Hom}_{A} \bigl( I_M, \Ga_{\m}(I_N) \bigr)
\end{equation}
in $\dcat{D}(B \ot C^{\mrm{op}}, \mrm{gr})$. 

By Proposition \ref{prop:3965} (in both cases of Notation \ref{not:3965}) 
the homomorphism $\si_{I_M} : \Ga_{\m}(I_M) \to I_M$ is a quasi-isomorphism. 
Hence in the commutative diagram 
\[ \UseTips \xymatrix @C=12ex @R=6ex {
\opn{Hom}_{A}(I_M, I_N) 
\ar[r]^(0.45){ \opn{Hom}(\si_{I_M}, \opn{id}) }
&
\opn{Hom}_{A} \bigl( \Ga_{\m}(I_M), I_N \bigr)
\\
\opn{Hom}_{A} \bigl( I_M, \Ga_{\m}(I_N) \bigr)
\ar[r]^(0.45){ \opn{Hom}(\si_{I_M}, \opn{id}) }
\ar[u]^{ \opn{Hom}(\opn{id}, \si_{I_N}) }
&
\opn{Hom}_{A} \bigl( \Ga_{\m}(I_M), \Ga_{\m}(I_N) \bigr) 
\ar[u]_{ \opn{Hom}(\opn{id}, \si_{I_N}) }
} \]
in $\dcat{C}_{\mrm{str}}(B \ot C^{\mrm{op}}, \mrm{gr})$
the two horizontal arrows are quasi-isomorphisms.
The right vertical arrow is an isomorphism 
in $\dcat{C}_{\mrm{str}}(B \ot C^{\mrm{op}}, \mrm{gr})$ by Lemma 
\ref{lem:3835}(3). 
It follows that the left vertical arrow is a quasi-isomorphism.
Combining this with the isomorphisms (\ref{eqn:3990}) and 
(\ref{eqn:3991}) we obtain the desired isomorphism above. 
\end{proof}

The next definition is a specialization of Definitions \ref{dfn:3887}, 
\ref{dfn:3890} and \ref{dfn:3895}.

\begin{dfn} \label{dfn:3748} 
Under Setup \tup{\ref{setup:3755}}:
\begin{enumerate}
\item A module $N \in \dcat{M}(A \ot B^{\mrm{op}}, \mrm{gr})$
is called {\em graded-$\m$-flasque}%
\index{Algebraically graded module! graded-$\m$-flasque}
if it is a right $\Ga_{\m}$-acyclic object, i.e.\ if 
$\mrm{H}^q_{\m}(N) = 0$ for all $q > 0$.

\item The functor 
\[ \Ga_{\m} : \dcat{M}(A \ot B^{\mrm{op}}, \mrm{gr}) \to
\dcat{M}(A \ot B^{\mrm{op}}, \mrm{gr}) \]
is called {\em weakly stable}%
\index{Torsion! weakly stable {\indash} functor}
if for every graded-injective module
$I \in \lb \dcat{M}(A \ot B^{\mrm{op}}, \mrm{gr})$, 
the module  
$\Ga_{\m}(I) \in \dcat{M}(A \ot B^{\mrm{op}}, \mrm{gr})$
is graded-$\m$-flasque. 

\item A complex  $N \in \dcat{D}(A \ot B^{\mrm{op}}, \mrm{gr})$
is called a {\em K-graded-$\m$-flasque complex}%
\index{Complex of algebraically graded modules! K-graded-$\m$-flasque}
if it is a right $\Ga_{\m}$-acyclic complex, i.e.\ if the morphism 
$\eta^{\mrm{R}}_{N} : \Ga_{\m}(N) \to \mrm{R} \Ga_{\m}(N)$
in $\dcat{D}(A \ot B^{\mrm{op}}, \mrm{gr})$ is an isomorphism. 
\end{enumerate}
\end{dfn}

The term ``flasque'' in this meaning seems to have been first used in 
\cite{YeZh5}.

\begin{exa} \label{exa:3745}
Suppose $M \in \dcat{M}(A, \mrm{gr})$ is an $\m$-torsion module. Consider its 
minimal graded-injective resolution 
$0 \to M \to I^0 \to I^1 \to  \cdots$. 
By Lemma \ref{lem:3900} and induction on $q$, we see that all the graded 
modules $I^q$ are $\m$-torsion; in fact 
$I^q \cong A^* \ot W^q$ for some $W^q \in \dcat{M}(\K, \mrm{gr})$. 
Thus
$\mrm{R} \Ga_{\m}(M) \cong \Ga_{\m}(I) = I \cong M$,
so $\mrm{R}^q \Ga_{\m}(M) = 0$ for all $q > 0$, and $M$ is a 
graded-$\m$-flasque $A$-module. 
\end{exa}

\begin{thm} \label{thm:3900}
Under Setup \tup{\ref{setup:3755}}, the torsion functor
\[ \Ga_{\m} : \dcat{M}(A \ot B^{\mrm{op}}, \mrm{gr}) \to 
\dcat{M}(A \ot B^{\mrm{op}}, \mrm{gr}) \]
is weakly stable.
\end{thm}

\begin{proof}
Take an injective object $I \in \dcat{M}(A \ot B^{\mrm{op}}, \mrm{gr})$. 
According to Proposition \ref{prop:3900}(3), the module $I$ is 
graded-injective over $A$. By Lemma \ref{lem:3839} we know that 
there is an isomorphism 
$\Ga_{\m}(I) \cong A^* \ot W$
in $\dcat{M}(A, \mrm{gr})$ for some graded $\K$-module $W$. This is a 
graded-injective $A$-module, and 
therefore $\Ga_{\m}(I)$ is graded-$\m$-flasque as an object of 
$\dcat{M}(A, \mrm{gr})$. 
Finally, by Proposition \ref{prop:3850}(4), the property of being 
graded-$\m$-flasque is insensitive to the restriction functor 
$\dcat{M}(A \ot B^{\mrm{op}}, \mrm{gr}) \to \dcat{M}(A, \mrm{gr})$,
so $\Ga_{\m}(I)$ is graded-$\m$-flasque as an object of 
$\dcat{M}(A \ot B^{\mrm{op}}, \mrm{gr})$. 
\end{proof}

In Proposition \ref{prop:3850}(2) we presented a morphism 
$\si^{\mrm{R}}_{} : \mrm{R} \Ga_{\m} \to \opn{Id}$
of triangulated functors from 
$\dcat{D}^{\star}(A \ot B^{\mrm{op}}, \mrm{gr})$
to itself. Here we are using Notation \ref{not:3965}. 
According to Definition \ref{dfn:3896} the pair 
$(\mrm{R} \Ga_{\m}, \si^{\mrm{R}})$ 
is a copointed triangulated functor. Recall that 
$(\mrm{R} \Ga_{\m}, \si^{\mrm{R}})$ is an 
idempotent copointed triangulated functor 
if for every complex
$M \in \dcat{D}^{\star}(A \ot B^{\mrm{op}}, \mrm{gr})$
the morphisms 
\[ \si^{\mrm{R}}_{\mrm{R} \Ga_{\m}(M)} \, , \,   
\mrm{R} \Ga_{\m}(\si^{\mrm{R}}_{M}) \, : \,
\mrm{R} \Ga_{\m}(\mrm{R} \Ga_{\m}(M)) \to \mrm{R} \Ga_{\m}(M) \]
in $\dcat{D}(A \ot B^{\mrm{op}}, \mrm{gr})$
are isomorphisms. 

\begin{cor}[\cite{VyYe}] \label{cor:3966}
With Notation \tup{\ref{not:3965}}, the copointed triangulated functor 
$(\mrm{R} \Ga_{\m}, \si^{\mrm{R}})$
on $\dcat{D}^{\star}(A \ot B^{\mrm{op}}, \mrm{gr})$ is idempotent.
\end{cor}

\begin{proof}
By Theorem \ref{thm:3900} the functor $\Ga_{\m}$ is weakly stable, and 
by Lemma \ref{lem:3835}(2) the copointed functor $(\Ga_{\m}, \si)$ is 
idempotent. Thus Theorem \ref{thm:1004} applies. 
\end{proof}

We end the subsection with an example showing a torsion functor as in 
Theorem \ref{thm:3900} which is weakly stable but not stable.
The example is a slight modification of an example in \cite{VyYe}. 

\begin{exa} \label{exa:4077}
Let $A := \K[t]$, the polynomial ring in a variable $t$ of degree $1$. It is 
connected noetherian, and its augmentation ideal $\m$ is generated by $t$. 
Let $B := \K[s_0, s_1, \ldots]$, the commutative polynomial ring in countably 
many variables $s_0, s_1, \ldots$ of degree $1$. So $B$ is a graded ring, 
and it is not noetherian. Theorem \ref{thm:3900} tells us that the functor 
\[ \Ga_{\m} : \dcat{M}(A \ot B, \mrm{gr}) \to \dcat{M}(A \ot B, \mrm{gr}) \]
is weakly stable. (Since $B$ is commutative, it doesn't matter whether we write 
$B$ or $B^{\mrm{op}}$). 
We will prove that $\Ga_{\m}$ is not stable. 

Consider a countable collection $\{ I_p \}_{p \in \N}$ of graded-injective 
$B$-modules. Then 
$I := \prod_{p \in \N} I_p(-p)$
is a graded-injective $B$-module, and 
therefore $J := \opn{Hom}_B(A \ot B, I)$ is a graded-injective 
$(A \ot B)$-module. If the functor $\Ga_{\m}$ on 
$\dcat{M}(A \ot B, \mrm{gr})$ were stable, this would imply 
that $\Ga_{\m}(J)$ is also a graded-injective $(A \ot B)$-module. Since 
$B \to A \ot B$ is flat, $\Ga_{\m}(J)$ is then a graded-injective $B$-module
(see Proposition \ref{prop:3900}). 

As graded $B$-modules there are isomorphisms
\[ A \ot B \cong \bigoplus\nolimits_{p \in \N} B \cd t^p 
\cong \bigoplus\nolimits_{p \in \N} B(-p) . \]
Therefore 
\[ J = \opn{Hom}_B(A \ot B, I) \cong \prod\nolimits_{p \in \N} I(p)  \]
as graded $B$-modules, and 
\[ \Ga_{\m}(J) \cong \lim_{p \to} \, 
\Bigl( \prod\nolimits_{q = 0}^p I(q) \Bigr) \cong 
\bigoplus\nolimits_{p \in \N} I(p) . \]
Now for each $p$ the graded-module $I_p$ is a direct summand of $I(p)$, and 
hence $\bigoplus_{p \in \N} I_p$ is a direct summand of $\Ga_{\m}(J)$, as 
graded $B$-modules. Therefore $\bigoplus_{p \in \N} I_p$ is a graded-injective 
$B$-module. 

The conclusion is that a countable direct sum of graded-injective $B$-modules 
is graded-injective. According to the Bass-Papp Theorem 
\cite[Theorem 3.46]{Lam}, that holds also in the graded sense (cf.\ Proposition 
\ref{prop:3745}), the graded ring $B$ is 
noetherian. This is a contradiction. 
\end{exa}

\mysubsection{Representability of Derived Torsion} 
\label{subsec:repr-der-tors}

The results of this subsection are an adaptation of results from 
\cite[Section 7]{VyYe} to the connected graded setting. 

In Section \ref{sec:BDC} we are going to use the representability of derived 
torsion to give an alternative proof of M. Van den Bergh's theorem on 
existence of balanced dualizing complexes (see Theorem \ref{thm:3713}). 

Recall that Conventions \ref{conv:3700} and \ref{conv:4560}
are in force. In this subsection we assume the following strengthening of Setup 
\ref{setup:3755}.

\begin{setup} \label{setup:3900}
We are given a left noetherian connected graded $\K$-ring $A$, with 
augmentation ideal 
$\m = \bigoplus_{i \geq 1} A_i \sub A$,
such that the torsion functor $\Ga_{\m}$ has finite right cohomological 
dimension. We are also given a graded $\K$-ring $B$.
\end{setup}

As before, the graded ring $B$ plays an auxiliary role, a placeholder for $A$ 
or $\K$. 

\begin{dfn} \label{dfn:3750}
Under Setup \ref{setup:3900}, the 
{\em dedualizing complex of $A$}%
\index{Dedualizing complex}
is the complex 
$P_A := \mrm{R} \Ga_{\m}(A) \in \dcat{D}(A^{\mrm{en}}, \mrm{gr})$.
\end{dfn}

The name ``dedualizing'' is due to L. Positselsky \cite{Pos} . Regarding the  
lack of left-right symmetry in this definition, see Corollary \ref{cor:3875} 
and Proposition \ref{prop:3712}.

As a special case of Proposition \ref{prop:3850}, there is a morphism 
$\si^{\mrm{R}}_{A} : P_A = \mrm{R} \Ga_{\m}(A) \lb \to A$
in $\dcat{D}(A^{\mrm{en}}, \mrm{gr})$. 

\begin{thm}[Representability of Derived Torsion, \cite{VyYe}] \label{thm:3750} 
\index{Torsion! representability of derived {\indash}}
Under Setup \tup{\ref{setup:3900}} there is a unique isomorphism 
\[ \ttev^{\mrm{R, L}}_{\m, (-)} \, : 
\, P_A \ot^{\mrm{L}}_{A} (-) \, \iso \, \mrm{R} \Ga_{\m}(-) \]
of triangulated functors from 
$\dcat{D}(A \ot B^{\mrm{op}}, \mrm{gr})$
to itself, called {\em derived $\m$-torsion tensor-evaluation},
\index{Tensor-evaluation morphism! derived $\m$-torsion}
\index{1-EvRLm@$\ttev^{\mrm{R, L}}_{\m, (-)}$}
such that for every complex 
$M \in \dcat{D}(A \ot B^{\mrm{op}}, \mrm{gr})$
the diagram 
\[ \tag{$\heartsuit$}
\UseTips \xymatrix @C=10ex @R=6ex {
P_A \ot^{\mrm{L}}_{A} M
\ar[r]^{ \ttev^{\mrm{R, L}}_{\m, M} }_{\cong}
\ar[d]_{\si^{\mrm{R}}_{A} \, \ot^{\mrm{L}}_{A} \, \opn{id}_M}
& 
\mrm{R} \Ga_{\m}(M)
\ar[d]^{\si^{\mrm{R}}_{M}}
\\
A \ot^{\mrm{L}}_{A} M
\ar[r]^{\opn{lu}}_{\cong}
& 
M
} \]
is commutative.
\end{thm}

The isomorphism marked ``$\opn{lu}$'' in diagram ($\heartsuit$) is the left 
unitor isomorphism of the closed monoidal structure, see equation 
(\ref{eqn:3790}). 

The proof of the theorem requires the next lemma.

\begin{lem} \label{lem:3905}
Let 
$F, G : \dcat{D}(A, \mrm{gr}) \to \dcat{D}(B, \mrm{gr})$
be quasi-compact triangulated functors, and let $\eta : F \to G$
be a morphism of triangulated functors. Assume that
$\eta_A : F (A) \to G (A)$
is an isomorphism. Then $\eta$ is an isomorphism. 
\end{lem}

\begin{proof}
We are given that $\eta_A$ is an isomorphism. Because the algebraic degree 
twist $M \mapsto M(i)$ is an automorphism of $\dcat{D}(A, \mrm{gr})$,
it follows that $\eta_{A(i)}$ is an isomorphism for every integer $i$.
Likewise for the cohomological degree shift (i.e.\ the translation 
$M \mapsto M[i]$ on complexes), so $\eta_{A(i)[j]}$ is an isomorphism for every 
integer $j$. Both functors are quasi-compact, and therefore 
$\eta_P$ is an isomorphism for every graded-free complex of $A$-modules 
$P \cong \bigoplus_{x \in X} A(i_x)[j_x]$. 

Suppose we are given a distinguished triangle
$M' \to M \to M'' \xar{\triangle}$
in $\dcat{D}(A, \mrm{gr})$, such that two of the three morphisms 
$\eta_{M'}$, $\eta_{M}$ and $\eta_{M''}$ are isomorphisms. Because $\eta$ is a 
morphism of triangulated functors, it follows that third morphism is
also an isomorphism. 

Next consider a semi-graded-free DG module $P$. Take a filtration \lb 
$\{ F_j(P) \}_{j \geq -1}$ of $P$ as in Definition \ref{dfn:4275}(2).
For every $j$ we have a distinguished triangle
\[ F_{j - 1}(P) \xar{\th_j} F_j(P) \to \opn{Gr}^F_j(P) \xar{\triangle} \]
in $\dcat{D}(A, \mrm{gr})$, where 
$\th_j : F_{j - 1}(P) \to F_j(P)$ is the inclusion. 
Since $\opn{Gr}^F_j(P)$ is a graded-free complex, by
induction on $j$ we conclude that $\eta_{F_j(P)}$ is an isomorphism for every 
$j \geq 0$. The homotopy colimit construction
(see Definition \ref{dfn:3410})
gives a distinguished triangle
\[ \bigoplus_{j \in \N} F_j(P) \xar{\ \Theta \ } \bigoplus_{j \in \N} F_j(P) 
\to P \xar{\triangle} \]
in $\dcat{D}(A, \mrm{gr})$, where 
\[ \Theta|_{F_{j - 1}(P)} := 
(\opn{id}_{}, - \th_j) : F_{j - 1}(P) \to F_{j - 1}(P) \oplus F_{j}(P) . \]
By quasi-compactness we know that 
$\eta_{\bigoplus F_j(P)}$ is an isomorphism. Therefore $\eta_P$ 
is an isomorphism. 

Finally, every $M \in \dcat{D}(A, \mrm{gr})$ admits an  
isomorphism $M \cong P$ with $P$ semi-graded-free. Therefore $\eta_M$ is an 
isomorphism.
\end{proof}

\begin{proof}[Proof of Theorem \tup{\ref{thm:3750}}] \mbox{}

\smallskip \noindent 
Step 1. We begin by constructing the morphism of triangulated functors 
$\ttev^{\mrm{R, L}}_{\m, (-)}$.
For each complex $M \in \dcat{C}(A \ot B^{\mrm{op}}, \mrm{gr})$
we choose a K-graded-injective resolution 
$\rho_M : M \to I_M$ and a K-graded-flat resolution 
$\th_M : Q_M \to M$,
both in $\dcat{C}_{\mrm{str}}(A \ot B^{\mrm{op}}, \mrm{gr})$.
Note that $I_M$ is K-graded-$\m$-flasque over $A$, and $Q_M$ is K-graded-flat 
over $A$. We use these choices for presentations of the right derived functor 
\[ (\mrm{R} \Ga_{\m}, \eta^{\mrm{R}}) : 
\dcat{D}(A \ot B^{\mrm{op}}, \mrm{gr}) \to 
\dcat{D}(A \ot B^{\mrm{op}}, \mrm{gr}) \] 
and the left derived bifunctor 
\[ (- \ot_{A}^{\mrm{L}} -, \eta^{\mrm{L}}) : 
\dcat{D}(A^{\mrm{en}}, \mrm{gr}) \times 
\dcat{D}(A \ot B^{\mrm{op}}, \mrm{gr}) \to 
\dcat{D}(A \ot B^{\mrm{op}}, \mrm{gr}) . \]
Let us also choose a K-graded-injective resolution 
$\psi : A \to J$ in $\dcat{C}_{\mrm{str}}(A^{\mrm{en}}, \mrm{gr})$. 
With these choices we have the following presentations:
$P_A = \Ga_{\m}(J)$,
\begin{equation} \label{eqn:3910}
\mrm{R} \Ga_{\m}(M) = \Ga_{\m}(I_M)
\end{equation}
and
\begin{equation} \label{eqn:3911}
P_A \ot_A^{\mrm{L}} M = \Ga_{\m}(J) \ot_{A} Q_M .
\end{equation}

Suppose 
$N \in \dcat{C}(A \ot B^{\mrm{op}}, \mrm{gr})$
is some complex. Given homogeneous elements 
$x \in \Ga_{\m}(J)$ and $n \in N$, the tensor 
$x \ot n$ belongs to $\Ga_{\m}(J \ot_{A} N)$.
In this way we obtain a homomorphism 
\begin{equation} \label{eqn:3905} 
\ttev_{\m, J, N} : \Ga_{\m}(J) \ot_{A} N \to \Ga_{\m}(J 
\ot_{A} N) 
\end{equation}
in $\dcat{C}_{\mrm{str}}(A \ot B^{\mrm{op}}, \mrm{gr})$, 
which is functorial in $N$. 

Consider a complex $M \in \dcat{C}(A \ot B^{\mrm{op}}, \mrm{gr})$. The choices 
we made give rise to this solid diagram 
\begin{equation} \label{eqn:3906}  
\UseTips \xymatrix @C=8ex @R=6ex {
A \ot_{A} Q_M
\ar[r]^(0.6){\opn{lu}}
\ar[d]_{ \psi \, \ot_A \, \opn{id} }
&
Q_M
\ar[r]^{\th_{M}}
&
M
\ar[d]^{ \rho_M }
\\
J \ot_{A} Q_M 
\ar@{-->}[rr]^(0.6){\chi_M}
&
&
I_M
} 
\end{equation}
in $\dcat{C}_{\mrm{str}}(A \ot B^{\mrm{op}}, \mrm{gr})$.
The homomorphisms $\psi \, \ot_A \, \opn{id}$, $\th_{M} \circ \opn{lu}$  and
$\rho_M$ are quasi-isomorphisms.
Because $I_M$ is K-injective in $\dcat{C}(A \ot B^{\mrm{op}}, \mrm{gr})$, 
there is a quasi-iso\-morphism
$\chi_M : J \ot_{A} Q_M \to I_M$ 
that makes this diagram commutative up to homotopy. 

We now form diagram (\ref{eqn:3908})
in $\dcat{C}_{\mrm{str}}(A \ot B^{\mrm{op}}, \mrm{gr})$. 
\begin{figure}[!htb]
\begin{equation} \label{eqn:3908}  
 \UseTips \xymatrix @C=12ex @R=6ex {
\Ga_{\m}(J) \ot_{A} Q_M
\ar[r]^{ \ttev_{\m, J, Q_M} }
\ar[d]_{ \si_{J} \, \ot_A \, \opn{id} }
&
\Ga_{\m}(J \ot_{A} Q_M)
\ar[r]^(0.55){ \Ga_{\m}(\chi_M) } 
\ar[d]_{ \si_{(J \ot_{A} Q_M)} }
&
\Ga_{\m}(I_M)
\ar[d]_{ \si_{I_M} }
\\
J \ot_{A} Q_M
\ar[r]^{ \opn{id} }
&
J \ot_{A} Q_M
\ar[r]^(0.55){ \chi_M }
&
I_M
\\
A \ot_{A} Q_M 
\ar[u]^{ \psi \, \ot_A \, \opn{id} }
\ar[r]^{ \opn{id} }
&
A \ot_{A} Q_M 
\ar[u]^{ \psi \, \ot_A \, \opn{id} }
\ar[r]^(0.55){\th_{M} \circ \, \opn{lu}}
&
M
\ar[u]^{ \rho_M } 
} 
\end{equation}
\end{figure}
Here $\ttev_{\m, J, Q_M}$ is the homomorphism from (\ref{eqn:3905}). The 
diagram (\ref{eqn:3908}) is commutative up to homotopy. (Actually all small 
squares, except the bottom right one, are commutative in the strict sense.)
The vertical arrows $\psi \, \ot_A \, \opn{id}$ and $\rho_M$ are 
quasi-isomorphisms. Passing to $\dcat{D}(A \ot B^{\mrm{op}}, \mrm{gr})$
we get a commutative diagram, with vertical isomorphisms between the second and 
third rows. The diagram with the four extreme objects only is the one we are 
looking for. By construction it is a commutative diagram
in $\dcat{D}(A \ot B^{\mrm{op}}, \mrm{gr})$, 
and it is functorial in $M$. The morphism 
\[ \ttev^{\mrm{R, L}}_{\m, M} : P_A \ot_A^{\mrm{L}} M \to 
\mrm{R} \Ga_{\m}(M) \] 
that we want is represented by 
\[ \Ga_{\m}(\chi_M) \circ \ttev_{\m, J, Q_M}  : 
\Ga_{\m}(J) \ot_{A} Q_M \to \Ga_{\m}(I_M) . \]
See (\ref{eqn:3910}) and (\ref{eqn:3911}). 

\medskip \noindent
Step 2. In this step we prove that $\ttev^{\mrm{R, L}}_{\m, M}$
is an isomorphism for every complex 
$M \in \dcat{D}(A \ot B^{\mrm{op}}, \mrm{gr})$.
Because the functor $\opn{Rest}$ is conservative, it suffices to prove that the 
morphism 
\[ \opn{Rest}(\ttev^{\mrm{R, L}}_{\m, M}) : 
\opn{Rest}(P_A \ot_A^{\mrm{L}} M) \to \opn{Rest}(\mrm{R} \Ga_{\m}(M)) \]
in $\dcat{D}(A, \mrm{gr})$ is an isomorphism. Going over all the details of the 
construction above, and noting that 
$\rho_M : M \to I_M$ and $\th_M : Q_M \to M$ 
are K-flat and K-injective resolutions, respectively, also in the category
$\dcat{C}_{\mrm{str}}(A, \mrm{gr})$, we might as well forget about the ring 
$B$. 

So now we are in the case $B = \K$, $A \ot B^{\mrm{op}} = A$, and we want to 
prove that $\ttev^{\mrm{R, L}}_{\m, M}$ is an isomorphism for every 
$M \in \dcat{D}(A, \mrm{gr})$.
By Theorem \ref{thm:4050} the functor $\mrm{R} \Ga_{\m}$
on $\dcat{D}(A, \mrm{gr})$ is quasi-compact. The functor 
$P \ot_{A}^{\mrm{L}} (-)$ is also quasi-compact. This means that we can use 
Lemma \ref{lem:3905}, and it tells us 
that it suffices to prove that $\ttev^{\mrm{R, L}}_{\m, M}$ is an 
isomorphism for $M = A$. 

Let us examine the morphism $\ttev^{\mrm{R, L}}_{\m, A}$, 
i.e.\ $\ttev^{\mrm{R, L}}_{\m, M}$ for $M = A$. 
We can choose the K-flat resolution
$\th_A : Q_A \to A$ in $\dcat{C}_{\mrm{str}}(A, \mrm{gr})$ to be the identity 
of 
$A$. Also, we can choose the K-injective resolution 
$\rho_A : A \to I_A$ in $\dcat{C}_{\mrm{str}}(A, \mrm{gr})$ to be the 
restriction of $\psi : A \to J$. Then the homomorphism 
$\chi_A : J \ot_A Q_A \to I_A$
in diagram (\ref{eqn:3906}) can be chosen to be 
$\chi_A = \opn{id} \ot_A \opn{id}$. 
We get a commutative diagram 
\[ \UseTips \xymatrix @C=12ex @R=6ex {
\Ga_{\m}(J) \ot_{A} Q_A
\ar[r]^{ \ttev_{\m, J, Q_A} }
\ar[d]_{ \opn{id} \ot_A \th_A }
&
\Ga_{\m}(J \ot_{A} Q_A)
\ar[r]^(0.6){ \Ga_{\m}(\chi_A) } 
\ar[d]_{  \Ga_{\m}(\opn{id} \ot_A \th_A) }
&
\Ga_{\m}(I_A)
\ar[d]_{\opn{id}}
\\
\Ga_{\m}(J) \ot_{A} A
\ar[r]^{ \ttev_{\m, J, A} }
&
\Ga_{\m}(J \ot_{A} A)
\ar[r]^(0.6){\Ga_{\m}(\opn{ru})} 
&
\Ga_{\m}(J)
} \]
in $\dcat{C}_{\mrm{str}}(A, \mrm{gr})$. 
Here $\opn{ru} : J \ot_A A \iso J$ is the canonical isomorphism (the right 
unitor). The horizontal arrows in the second row, and the vertical arrows, are 
all bijective. We conclude that 
$\Ga_{\m}(\chi_A) \circ \ttev_{\m, J, Q_A}$
is bijective. Hence $\ttev^{\mrm{R, L}}_{\m, A}$ 
is an isomorphism in $\dcat{D}(A, \mrm{gr})$.

\medskip \noindent
Step 3. It remains to prove the uniqueness of 
$\ttev^{\mrm{R, L}}_{\m, (-)}$. By applying the functor 
$\mrm{R} \Ga_{\m}$ to the diagram ($\heartsuit$), and then making use of the 
morphism of functors
$\si^{\mrm{R}} : \mrm{R} \Ga_{\m} \to \opn{Id}$, 
we obtain the commutative diagram (\ref{eqn:4820})
in $\dcat{D}(A \ot B^{\mrm{op}}, \mrm{gr})$. 
\begin{figure}[!htb]
\begin{equation} \label{eqn:4820}
\UseTips \xymatrix @C=15ex @R=6ex {
P_A \ot^{\mrm{L}}_{A} M
\ar[r]^{ \ttev^{\mrm{R, L}}_{\m, M} }_{\cong}
& 
\mrm{R} \Ga_{\m}(M)
\\
\mrm{R} \Ga_{\m}(P_A \ot^{\mrm{L}}_{A} M)
\ar[u]^{\si^{\mrm{R}}_{(P_A \ot^{\mrm{L}}_{A} M)}}
\ar[r]^{ \mrm{R} \Ga_{\m}(\ttev^{\mrm{R, L}}_{\m, M}) }_{\cong}
\ar[d]_{\mrm{R} \Ga_{\m}(\si^{\mrm{R}}_{A} \, \ot^{\mrm{L}}_{A} \, \opn{id}_M)}
& 
\mrm{R} \Ga_{\m}(\mrm{R} \Ga_{\m}(M))
\ar[u]_{\si^{\mrm{R}}_{\mrm{R} \Ga_{\m}(M)}}
\ar[d]^{\mrm{R} \Ga_{\m}(\si^{\mrm{R}}_{M})}
\\
\mrm{R} \Ga_{\m}(A \ot^{\mrm{L}}_{A} M)
\ar[r]^{\mrm{R} \Ga_{\m}(\opn{lu})}_{\cong}
& 
\mrm{R} \Ga_{\m}(M)
}
\end{equation}
\end{figure}
The vertical morphisms 
$\si^{\mrm{R}}_{\mrm{R} \Ga_{\m}(M)}$
and 
$\mrm{R} \Ga_{\m}(\si^{\mrm{R}}_{M})$
are isomorphisms because of the idempotence
(Corollary \ref{cor:3966}).
Therefore the other two vertical arrows are isomorphisms. We see that the 
isomorphism $\ttev^{\mrm{R, L}}_{\m, M}$ can be expressed as the 
composition of other isomorphisms; so it is unique. 
\end{proof}

The next result is due to Van den Bergh \cite{VdB}. We give another 
proof, based on Theorem \ref{thm:3750}. 
Recall the $\K$-linear duality functor 
\[ (-)^* : \dcat{D}(A \ot B^{\mrm{op}}, \mrm{gr})^{\mrm{op}} \to 
\dcat{D}(B \ot A^{\mrm{op}}, \mrm{gr}) . \]
Its formula is  
\[  (-)^* := \opn{Hom}_{\K}(-, \K) \cong \opn{Hom}_{A}(-, A^*) . \]

\begin{cor}[Van den Bergh's Local Duality, \cite{VdB}] \label{cor:3901}
\index{Local Duality Theorem}
Under Setup \tup{\ref{setup:3900}}, for every complex 
$M \in \dcat{D}(A \ot B^{\mrm{op}},  \mrm{gr})$ there is an 
isomorphism 
\[ \mrm{R} \Ga_{\m}(M)^* \cong 
\opn{RHom}_A \bigl( M, (P_A)^* \bigr) \]
in $\dcat{D}(B \ot A^{\mrm{op}}, \mrm{gr})$.
This isomorphism is functorial in $M$. 
\end{cor}

\begin{proof}
We have the following functorial isomorphisms
\[ \begin{aligned}
& \mrm{R} \Ga_{\m}(M)^* \cong^{\dag} (P_A \ot^{\mrm{L}}_{A} M)^* 
= \opn{Hom}_{\K}(P_A \ot^{\mrm{L}}_{A} M, \K )
\\ & 
\quad \cong^{} \opn{RHom}_A(P_A \ot^{\mrm{L}}_{A} M, A^* )
\\ & 
\quad \cong^{\heartsuit} \opn{RHom}_A \bigl( M, 
\opn{RHom}_A(P_A , A^* ) \bigr) 
\\ & 
\quad \cong^{} \opn{RHom}_A \bigl( M, (P_A)^* \bigr) 
\end{aligned} \]
in $\dcat{D}(B \ot A^{\mrm{op}}, \mrm{gr})$.
The isomorphism $\cong^{\dag}$ is according to Theorem \ref{thm:3750}.
The isomorphism $\cong^{\heartsuit}$ comes from derived Hom-tensor adjunction.
\end{proof}

\mysubsection{Symmetry of Derived Torsion} 
\label{subsec:sym-der-tor}

The results of this subsection are an adaptation of results from 
\cite[Section 8]{VyYe} to the connected graded setting. 
Among other things, we prove that the $\chi$ condition of Artin-Zhang implies 
symmetry of derived torsion. 

In Section \ref{sec:BDC} below we are going to use the symmetry of derived 
torsion to give an alternative proof of Van den Bergh's theorem on existence of 
balanced dualizing complexes (see Theorem \ref{thm:3713}). 

We continue with Conventions \ref{conv:3700}  and \ref{conv:4560}.
Thus $\K$ is a field, and all graded rings are algebraically graded central 
$\K$-rings. The following setup is assumed in this subsection: 

\begin{setup} \label{setup:3855} 
We are given a noetherian connected graded $\K$-ring $A$, with augmentation 
ideal $\m$. The augmentation ideal of the opposite ring $A^{\mrm{op}}$ is 
$\m^{\mrm{op}}$.
\end{setup}

The assumption that $A$ is noetherian means that   
both $A$ and $A^{\mrm{op}}$ are left noetherian. 
The enveloping ring of $A$ is 
$A^{\mrm{en}} = A \ot A^{\mrm{op}}$, and it is often not noetherian.

As explained in the previous subsection, there is a torsion functor 
$\Ga_{\m} : \dcat{M}(A, \mrm{gr}) \lb \to \dcat{M}(A, \mrm{gr})$.
But by the same token, there is also a torsion functor 
$\Ga_{\m^{\mrm{op}}} : \lb \dcat{M}(A^{\mrm{op}}, \mrm{gr}) \to 
\dcat{M}(A^{\mrm{op}}, \mrm{gr})$.
On the category $\dcat{M}(A^{\mrm{en}}, \mrm{gr})$ we have two 
torsion functors:
\[ \Ga_{\m}, \Ga_{\m^{\mrm{op}}} : \dcat{M}(A^{\mrm{en}}, \mrm{gr}) \to 
\dcat{M}(A^{\mrm{en}}, \mrm{gr}) . \]
These functors commute with each other. Indeed, for any 
$M \in \dcat{M}(A^{\mrm{en}}, \mrm{gr})$ there is equality 
$\Ga_{\m}(\Ga_{\m^{\mrm{op}}}(M)) = \Ga_{\m^{\mrm{op}}}(\Ga_{\m}(M))$
of graded $A^{\mrm{en}}$-submodules of $M$. 

\begin{dfn} \label{dfn:3766} 
Under Setup \tup{\ref{setup:3855}}, we say that $A$ has 
{\em finite local cohomological dimension}%
\index{Algebraically graded ring! connected {\indash} of finite local 
cohomological dimension}
if the torsion functors $\Ga_{\m}$ and $\Ga_{\m^{\mrm{op}}}$,
on $\dcat{M}(A, \mrm{gr})$ and 
$\dcat{M}(A^{\mrm{op}}, \mrm{gr})$ respectively, have finite right 
cohomological dimensions, 
in the sense of Definition \ref{dfn:3885}. Namely if there is a natural number 
$d$ such that $\mrm{H}^p_{\m} = 0$ and $\mrm{H}^p_{\m^{\mrm{op}}} = 0$ 
for all $p > d$. See also formula (\ref{eqn:3962}). 
\end{dfn}

Note that the vanishing of the derived functors  
$\mrm{H}^p_{\m} = \mrm{R}^p \Ga_{\m}$, and thus the cohomological dimension of 
the functor $\Ga_{\m}$, is the same on bimodules and on modules, by Proposition 
\ref{prop:3850}(4). Likewise for $\mrm{H}^p_{\m^{\mrm{op}}}$.

As explained in Proposition \ref{prop:3850}(1), there are derived torsion 
functors 
\begin{equation} \label{eqn:3817}
\mrm{R} \Ga_{\m}, \mrm{R} \Ga_{\m^{\mrm{op}}} : 
\dcat{D}(A^{\mrm{en}}, \mrm{gr}) \to \dcat{D}(A^{\mrm{en}}, \mrm{gr}) .
\end{equation}
There is no reason for these derived functors to commute with each other in 
general. 

\begin{dfn} \label{dfn:3717} 
Under Setup \tup{\ref{setup:3855}}, the morphisms from Proposition \lb 
\ref{prop:3850}(2) associated to the derived functors 
in formula (\ref{eqn:3817}) are denoted by%
\index{1-RGamma-m@$\mrm{R} \Ga_{\m}$}%
\index{1-SiRm@$\si^{\mrm{R}}_{\m}$}
$\si^{\mrm{R}}_{\m} : \mrm{R} \Ga_{\m} \to \opn{Id}$
and 
$\si^{\mrm{R}}_{\m^{\mrm{op}}} : \mrm{R} \Ga_{\m^{\mrm{op}}} \to \opn{Id}$.
\end{dfn}

\begin{dfn} \label{dfn:3815} 
Under Setup \tup{\ref{setup:3855}}: 
\begin{enumerate}
\item We denote by 
$\dcat{M}_{(\mrm{tor}, ..)}(A^{\mrm{en}})$ the 
full subcategory of $\dcat{M}(A^{\mrm{en}})$ on the bimodules $M$ that are 
$\m$-torsion as $A$-modules; i.e.\ 
such that $\opn{Rest}_A(M) \in \dcat{M}_{\mrm{tor}}(A)$.

\item We denote by 
$\dcat{M}_{(.., \mrm{tor})}(A^{\mrm{en}})$ the 
full subcategory of $\dcat{M}(A^{\mrm{en}})$ on the bimodules $M$ that are 
$\m^{\mrm{op}}$-torsion as $A^{\mrm{op}}$-modules.

\item We let 
\[ \dcat{M}_{(\mrm{tor}, \mrm{tor})}(A^{\mrm{en}}) := 
\dcat{M}_{(\mrm{tor}, ..)}(A^{\mrm{en}}) \cap 
\dcat{M}_{(.., \mrm{tor})}(A^{\mrm{en}}) . \]

\item Let $\star, \diamond$ be torsion indicators, i.e.\ ``$\mrm{tor}$'' or 
``$..$''. We denote by 
$\dcat{D}_{(\star, \diamond)}(A^{\mrm{en}})$
the full subcategory of 
$\dcat{D}(A^{\mrm{en}})$ on the complexes $M$ such that 
$\opn{H}^p(M) \in \dcat{M}_{(\star, \diamond)}(A^{\mrm{en}})$ for 
every integer $p$.
\end{enumerate}
\end{dfn}

Trivially, for $M \in \dcat{D}(A^{\mrm{en}}, \mrm{gr})$ 
the $A$-modules $\opn{H}^p(\mrm{R} \Ga_{\m}(M))$ are $\m$-torsion, 
and the $A^{\mrm{op}}$-modules 
$\opn{H}^p(\mrm{R} \Ga_{\m^{\mrm{op}}}(M))$ are $\m^{\mrm{op}}$-torsion. 
Sometimes more happens:

\begin{dfn} \label{dfn:3720}
Under Setup \tup{\ref{setup:3855}}, a
complex $M \in \dcat{D}(A^{\mrm{en}}, 
\mrm{gr})$ is said to have 
{\em weakly symmetric derived $\m$-torsion}%
\index{Complex of algebraically graded modules! with weakly symmetric derived 
$\m$-torsion}
if these two conditions hold:
\begin{itemize}
\item For every $p \in \Z$ the bimodule  $\opn{H}^p(\mrm{R} \Ga_{\m}(M))$ 
is $\m^{\mrm{op}}$-torsion.

\item For every $p \in \Z$ the bimodule 
$\opn{H}^p(\mrm{R} \Ga_{\m^{\mrm{op}}}(M))$ is $\m$-torsion.
\end{itemize} 
\end{dfn}

In terms of Definition \ref{dfn:3815}, a complex 
$M \in \dcat{D}(A^{\mrm{en}}, \mrm{gr})$
has weakly symmetric derived $\m$-torsion if 
\index{1-Dtortor(A)@$\dcat{D}_{(\mrm{tor}, \mrm{tor})}(A^{\mrm{en}}, \mrm{gr})$}
\[ \mrm{R} \Ga_{\m}(M) \, , \, \mrm{R} \Ga_{\m^{\mrm{op}}}(M) \, \in \,
\dcat{D}_{(\mrm{tor}, \mrm{tor})}(A^{\mrm{en}}, \mrm{gr}) . \]

\begin{dfn} \label{dfn:3779}
Under Setup \tup{\ref{setup:3855}}, a complex 
$M \in \dcat{D}(A^{\mrm{en}}, \mrm{gr})$ is said to have 
{\em symmetric derived $\m$-torsion}%
\index{Complex of algebraically graded modules! with symmetric derived 
$\m$-torsion}
if there is an isomorphism 
\[ \ep_{M} : \mrm{R} \Ga_{\m}(M) \iso 
\mrm{R} \Ga_{\m^{\mrm{op}}}(M) \]
in $\dcat{D}(A^{\mrm{en}}, \mrm{gr})$, 
such that the diagram
\[ \UseTips \xymatrix @C=8ex @R=6ex {
\mrm{R} \Ga_{\m}(M)
\ar[r]^{\ep_{M}}_{\cong}
\ar[dr]_{\si^{\mrm{R}}_{\m, M}}
&
\mrm{R} \Ga_{\m^{\mrm{op}}}(M)
\ar[d]^{\si^{\mrm{R}}_{\m^{\mrm{op}}, M}}
\\
&
M
} \]
in $\dcat{D}(A^{\mrm{en}}, \mrm{gr})$ is commutative.
Such an isomorphism $\ep_M$ is called a 
{\em symmetry isomorphism}%
\index{Complex of algebraically graded modules! symmetry isomorphism of}.
\end{dfn}

Of course if $M$ has symmetric derived $\m$-torsion, then it has 
weakly symmetric derived $\m$-torsion.

\begin{thm}[Symmetric Derived Torsion, \cite{VyYe}] 
\label{thm:3720}
\index{Algebraically graded ring! connected {\indash} of finite local 
cohomological dimension}%
\index{Complex of algebraically graded modules! with weakly symmetric derived 
$\m$-torsion}%
\index{Complex of algebraically graded modules! with symmetric derived 
$\m$-torsion}
Let $A$ be a noetherian connected graded $\K$-ring of finite local
cohomological dimension. 
If $M \in \dcat{D}(A^{\mrm{en}}, \mrm{gr})$
has weakly symmetric derived $\m$-torsion, then $M$ has 
symmetric derived $\m$-torsion.

Moreover, the symmetry isomorphism $\ep_{M}$ is unique, and it is functorial in 
such complexes $M$. 
\end{thm}

\begin{proof}
In Theorem \ref{thm:3750}, with $B := A$,
we have an isomorphism 
\[ \ttev^{\mrm{R, L}}_{\m, (-)} \, : 
\, P_A \ot^{\mrm{L}}_{A} (-) \, \iso \, \mrm{R} \Ga_{\m}(-) \]
of triangulated functors from 
$\dcat{D}(A^{\mrm{en}}, \mrm{gr})$ to itself. The same theorem, but with the 
roles of $A$ and $A^{\mrm{op}}$ exchanged, gives an isomorphism 
\[ \ttev^{\mrm{R, L}}_{\m^{\mrm{op}}, (-)}  
\, : \, (-) \ot^{\mrm{L}}_{A} P_{A^{\mrm{op}}} \, 
\iso \, \mrm{R} \Ga_{\m^{\mrm{op}}}(-) \]
of triangulated functors from 
$\dcat{D}(A^{\mrm{en}}, \mrm{gr})$ to itself. Here 
$P_{A^{\mrm{op}}} := \mrm{R} \Ga_{\m^{\mrm{op}}}(A) \in 
\dcat{D}(A^{\mrm{en}}, \mrm{gr})$.

Consider diagram (\ref{eqn:3871})  in $\dcat{D}(A^{\mrm{en}}, \mrm{gr})$.
In this diagram, $\al$ is the associativity isomorphism of the derived tensor 
product, and $\opn{ru}, \opn{lu}$ are the monoidal unitor isomorphisms -- all 
belonging to the monoidal structure of $\dcat{D}(A^{\mrm{en}}, \mrm{gr})$.
A quick check, using K-flat resolutions and elements in them, shows that this 
is a commutative diagram. 

\begin{figure}[!htb]
\begin{equation} \label{eqn:3871} 
\UseTips \xymatrix @C=8ex @R=6ex {
(P_A \ot^{\mrm{L}}_{A} M) \ot^{\mrm{L}}_{A} P_{A^{\mrm{op}}}
\ar[d]_{ \opn{id}  \ot^{\mrm{L}}_{A} \si^{\mrm{R}}_{\m^{\mrm{op}}, A} }
\ar[r]^{\al}_{\cong}
&
P_A \ot^{\mrm{L}}_{A} (M \ot^{\mrm{L}}_{A} P_{A^{\mrm{op}}})
\ar[d]^{ \si^{\mrm{R}}_{\m, A}  \ot^{\mrm{L}}_{A} \opn{id} }
\\
(P_A \ot^{\mrm{L}}_{A} M) \ot^{\mrm{L}}_{A} A
\ar[d]_{\opn{ru}}
&
A \ot^{\mrm{L}}_{A} (M \ot^{\mrm{L}}_{A} P_{A^{\mrm{op}}})
\ar[d]^{\opn{lu}}
\\
P_A \ot^{\mrm{L}}_{A} M 
\ar[dr]_{ \si^{\mrm{R}}_{\m, A}  \ot^{\mrm{L}}_{A} \opn{id} \ }
&
M \ot^{\mrm{L}}_{A} P_{A^{\mrm{op}}}
\ar[d]^{ \opn{id} \ot^{\mrm{L}}_{A} \si^{\mrm{R}}_{\m^{\mrm{op}}, A} }
\\
&
M
} 
\end{equation}
\end{figure}

By Theorem \ref{thm:3750}, diagram (\ref{eqn:3871}) gives rise to 
diagram (\ref{eqn:3872}), with its solid arrows only. 
This new diagram in $\dcat{D}(A^{\mrm{en}}, \mrm{gr})$ (solid arrows 
only) is isomorphic to the previous one (without its second row), via the 
various isomorphisms 
$\ttev^{\mrm{R, L}}_{\m, (-)}$ and 
$\ttev^{\mrm{R, L}}_{\m^{\mrm{op}}, (-)}$.

\begin{figure}[!htb]
\begin{equation} \label{eqn:3872}
\UseTips \xymatrix @C=8ex @R=6ex {
\mrm{R} \Ga_{\m^{\mrm{op}}}(\mrm{R} \Ga_{\m}(M))
\ar[r]^{\al'}_{\cong}
\ar[d]_{ \si^{\mrm{R}}_{\m^{\mrm{op}}, \mrm{R} \Ga_{\m}(M)} }
&
\mrm{R} \Ga_{\m}(\mrm{R} \Ga_{\m^{\mrm{op}}}(M))
\ar[d]^{ \si^{\mrm{R}}_{\m, \mrm{R} \Ga_{\m^{\mrm{op}}}(M)} }
\\
\mrm{R} \Ga_{\m}(M)
\ar[dr]_{ \si^{\mrm{R}}_{\m, M} }
\ar@{-->}[r]^{\ep_M}
&
\mrm{R} \Ga_{\m^{\mrm{op}}}(M)
\ar[d]^{ \si^{\mrm{R}}_{\m^{\mrm{op}}, M} }
\\
&
M
} 
\end{equation}
\end{figure}

Because 
$\mrm{R} \Ga_{\m}(M) \in 
\dcat{D}_{(\mrm{tor}, \mrm{tor})}(A^{\mrm{en}}, \mrm{gr})$,
Proposition \ref{prop:3965}, applied with $A^{\mrm{op}}$ instead 
of with $A$, says that 
$\si^{\mrm{R}}_{\m^{\mrm{op}}, \mrm{R} \Ga_{\m}(M)}$
is an isomorphism. Likewise, because 
$\mrm{R} \Ga_{\m^{\mrm{op}}}(M) \in 
\dcat{D}_{(\mrm{tor}, \mrm{tor})}(A^{\mrm{en}}, \mrm{gr})$,
Proposition \ref{prop:3965} says that 
$\si^{\mrm{R}}_{\m, \mrm{R} \Ga_{\m^{\mrm{op}}}(M)}$
is an isomorphism. We define $\ep_M$ to be the unique isomorphism (the dashed 
arrow) that makes diagram (\ref{eqn:3872}) commutative.

Diagram (\ref{eqn:3872}) shows that $\ep_M$ is unique: it has to be 
\[ \ep_M = \si^{\mrm{R}}_{\m, \mrm{R} \Ga_{\m^{\mrm{op}}}(M)} \circ 
\al' \circ (\si^{\mrm{R}}_{\m^{\mrm{op}}, \mrm{R} \Ga_{\m}(M)})^{-1} . \]
The functoriality of $\ep_M$ is a consequence of the functoriality of diagram 
(\ref{eqn:3872}). 
\end{proof}

\begin{cor} \label{cor:3875}
Under the assumptions of Theorem \tup{\ref{thm:3720}}, if the bimodule $A$ has 
has weakly symmetric derived $\m$-torsion, then there is a unique isomorphism 
$\ep_A : P_A \iso P_{A^{\mrm{op}}}$
in $\dcat{D}(A^{\mrm{en}}, \mrm{gr})$
such that 
$\si^{\mrm{R}}_{\m, A} = \si^{\mrm{R}}_{\m^{\mrm{op}}, A} \circ \ep_A$
as morphisms $P_A \to A$. 
\end{cor}

\begin{proof}
Take $M = A$ in the theorem. 
\end{proof}

The socle functor from Definition \ref{dfn:4265} extends to complexes. 

\begin{lem} \label{lem:3840}
\index{Complex of algebraically graded modules! minimal {\indash}
graded-injective {\indash}}
Let $I$ be a minimal complex of graded-injective $A$-modules, and let 
$W := \opn{Soc}(I)$. 
Then the differential of the complex $W$ is zero. 
\end{lem}

\begin{proof}
Take some nonzero homogeneous element $w \in W^p$, say $w \in W^p_q$. 
Let $U := \K \cd w \sub W^p_q$, so $U$ is a nonzero graded 
$A$-submodule of $I^p$. Because $\opn{Z}^p(I)$ is an essential graded 
$A$-submodule of $I^p$ we must have 
$\opn{Z}^p(I) \cap U \neq 0$. Thus $w \in \opn{Z}^p(I)$, so 
$\d(w) = 0$.
\end{proof}

Cofinite graded $A$-modules were defined in Definition \ref{dfn:3726}. 
By graded Matlis duality (Theorem \ref{thm:4045}) the cofinite graded 
$A$-modules are the artinian objects of $\dcat{M}(A, \mrm{gr})$. 
The full subcategory $\dcat{M}_{\mrm{cof}}(A, \mrm{gr})$ of cofinite modules 
is a thick abelian subcategory of $\dcat{M}(A, \mrm{gr})$, closed under 
subobjects and quotients.

\begin{lem} \label{lem:3856}
Consider a complex $M \in \dcat{D}^+(A, \mrm{gr})$ and an integer $p$. 
\begin{enumerate}
\item If $\opn{Ext}^p_A(\K, M)$ is a finite graded $\K$-module, then 
$\mrm{H}^p_{\m}(M)$ is a cofinite graded $A$-module. 

\item If $\mrm{H}^{q}_{\m}(M)$ is a cofinite graded $A$-module for every 
$q \leq p$, then $\opn{Ext}^{q}_A(\K, M)$ is a finite graded $\K$-module
for every $q \leq p$.
\end{enumerate}
\end{lem}

\begin{proof} \mbox{}

\smallskip \noindent
(1) Let $M \to I$ be a minimal graded-injective resolution, and let 
$W := \opn{Soc}(I)$. Then $\opn{Ext}^p_A(\K, M) \cong \opn{H}^p(W)$. 
But by Lemma \ref{lem:3840} we know that the differential of $W$ is zero, so in 
fact $\opn{Ext}^p_A(\K, M) \cong W^p$. Therefore $W$ is a finite graded 
$\K$-module. 

Next, let $J := \Ga_{\m}(I)$. This complex has the property that 
$\mrm{H}^p_{\m}(M) \cong \mrm{H}^p(J)$. According to Lemma \ref{lem:3839} we 
have an isomorphism $J^p \cong A^* \ot W^p$. 
Since $A^*$ is a cofinite graded $A$-module, and $J^p$ is a finite direct sum 
of degree twists of $A^*$, it is also a cofinite graded $A$-module. But 
$\mrm{H}^p_{\m}(M)$ is a subquotient of $J^p$, so it too is a cofinite graded 
$A$-module. 

\medskip \noindent 
(2) We continue with the minimal graded-injective resolution $M \to I$. 
We may assume that $\opn{H}(M) \neq 0$. 
Let $p_0 := \opn{inf}(\opn{H}(M))$; so 
$p_0 = \opn{inf}(I)$. Take the subcomplexes $W \sub J \sub I$ as above. We know 
that $J^q \cong A^* \ot W^q$ for every $q$. We will prove that 
$W^q$ is a finite graded $\K$-module for every $q \leq p$, by induction on $q$, 
starting with $q = p_0$. 

So take an integer $q$ in the range $p_0 \leq q \leq p$, and assume that 
$W^{q - 1}$ is a finite graded $\K$-module. 
There is an exact sequence 
\[ J^{q - 1} \xar{\d} \opn{Z}^q(J) \to \mrm{H}^q_{\m}(M) \to  0  \]
in $\dcat{M}(A, \mrm{gr})$. 
Since $J^{q - 1} \cong A^* \ot W^{q - 1}$, this is a cofinite graded 
$A$-module, as we already noticed. We also know that $\mrm{H}^q_{\m}(M)$ is 
cofinite. It follows that $\opn{Z}^q(J)$ is a cofinite graded $A$-module.
But by Lemma \ref{lem:3840} we have 
$W^q \sub \opn{Z}^q(J)$, so $W^q$ is a cofinite graded $A$-module annihilated 
by $\m$. Therefore $W^q$ is a finite graded $\K$-module.
\end{proof}

\begin{dfn}[Artin-Zhang, \cite{ArZh}] \label{dfn:3713}
Assume Setup \tup{\ref{setup:3855}}.
\begin{enumerate}
\item We say that the ring $A$ satisfies the {\em left $\chi$ condition} 
if for every $M \in \dcat{M}_{\mrm{f}}(A, \mrm{gr})$
and every integer $p$, the graded $\K$-module 
$\opn{Ext}^p_{A}(\K, M)$ is finite. 

\item The ring $A$ is said to satisfy the 
{\em $\chi$ condition}%
\index{Algebraically graded ring! chi@$\chi$ condition on connected}
if both graded rings $A$ and $A^{\mrm{op}}$ satisfy the left $\chi$ condition.
\end{enumerate}
\end{dfn}

Another way to state the left $\chi$ condition is this: 
$M \in \dcat{M}_{\mrm{f}}(A, \mrm{gr})$ implies 
$\opn{RHom}_{A}(\K, M) \in \dcat{D}_{\mrm{f}}(\K, \mrm{gr})$.

\begin{rem} \label{rem:3855}
What we call the ``left $\chi$ condition'' was actually called the ``left 
$\chi^{\circ}$ condition'' in \cite[Definition 3.2]{ArZh}. 
Their ``left $\chi$ condition'' 
(in \cite[Definition 3.7]{ArZh}) is much more complicated to state. 
However, for a left noetherian connected graded ring $A$ (as we have here) 
these two conditions are equivalent, by \cite[Proposition 3.11(2)]{ArZh}. 
\end{rem}

\begin{dfn} \label{dfn:3711} 
Assume Setup \tup{\ref{setup:3855} }.
\begin{enumerate}
\item We say that the ring $A$ satisfies the {\em left special $\chi$ 
condition} if for every integer $p$, the graded $\K$-module 
$\opn{Ext}^p_{A}(\K, A)$ is finite. 

\item The ring $A$ is said to satisfy the 
{\em special $\chi$ condition}%
\index{Algebraically graded ring! special $\chi$ condition on connected}
if both graded rings $A$ and $A^{\mrm{op}}$ satisfy the left special $\chi$ 
condition.
\end{enumerate}
\end{dfn}

\begin{prop} \label{prop:3855}
The following two conditions are equivalent.
\begin{enumerate}
\rmitem{i} $A$ satisfies the left $\chi$ condition. 

\rmitem{ii} For every $M \in \dcat{M}_{\mrm{f}}(A, \mrm{gr})$ and every $p$ 
the graded $A$-module $\mrm{H}^p_{\m}(M)$ is cofinite. 
\end{enumerate}
\end{prop}

\begin{proof}
This is immediate from Lemma \ref{lem:3856}. 
\end{proof}

\begin{dfn} \label{dfn:4050} 
We denote by 
$\dcat{D}^{\mrm{b}}_{(\mrm{f}, \mrm{f})}(A^{\mrm{en}}, \mrm{gr})$ the 
full subcategory of \lb 
$\dcat{D}^{\mrm{b}}(A^{\mrm{en}}, \mrm{gr})$ on the 
complexes of bimodules $M$ whose cohomology bimodules $\opn{H}^p(M)$ are 
finite graded modules over $A$ and over $A^{\mrm{op}}$. 
\end{dfn}

\begin{prop} \label{prop:3712}
\mbox{}
\begin{enumerate}
\item If $A$ satisfies the $\chi$ condition, then it satisfies 
the special $\chi$ condition. 

\item If $A$ satisfies the $\chi$ condition, then every
$M \in \dcat{D}^{\mrm{b}}_{(\mrm{f}, \mrm{f})}(A^{\mrm{en}}, \mrm{gr})$
has weakly symmetric derived $\m$-torsion. 

\item If $A$ satisfies the special $\chi$ condition, 
then the bimodule $A$ has weakly symmetric derived $\m$-torsion. 
\end{enumerate}
\end{prop}

\begin{proof} \mbox{}

\smallskip \noindent 
(1) This is trivial.

\medskip \noindent
(2) Take some 
$M \in \dcat{D}^{\mrm{b}}_{(\mrm{f}, \mrm{f})}(A^{\mrm{en}}, \mrm{gr})$.
Using smart truncation and induction on the amplitude of cohomology, 
we can assume that 
$M \in \dcat{M}_{(\mrm{f}, \mrm{f})}(A^{\mrm{en}}, \mrm{gr})$.

The $\chi$ condition tells us that $\opn{Ext}^p_{A}(\K, M)$ and 
$\opn{Ext}^p_{A^{\mrm{op}}}(\K, M)$ 
are finite as graded $\K$-modules, for all $p$. For every $q \geq 1$ the graded 
bimodule $A / \m^q$ is gotten from $\K$ by finitely many degree twists and 
extensions in $\dcat{M}(A^{\mrm{en}}, \mrm{gr})$. Therefore 
$\opn{Ext}^p_{A}(A / \m^q, M)$ and 
$\opn{Ext}^p_{A^{\mrm{op}}}(A / \m^q, M)$
are finite as graded $\K$-modules, for all $p$ and all $q \geq 1$. 
Thus as left and as right graded $A$-modules, 
$\opn{Ext}^p_{A}(A / \m^q, M)$ and 
$\opn{Ext}^p_{A^{\mrm{op}}}(A / \m^q, M)$
are of finite length, and thus torsion on both sides. Passing to the direct 
limit using Proposition \ref{prop:3850}(5), we see that $\mrm{H}^p_{\m}(M)$ and 
$\mrm{H}^p_{\m^{\mrm{op}}}(M)$ are torsion $A$-modules on both sides. 

\medskip \noindent
(3) This is just a special case of the proof of item (2). 
\end{proof}

\begin{lem} \label{lem:3870}
Assume the functor $\Ga_{\m}$ has finite cohomological dimension, and  $A$ has 
the left special $\chi$ condition. Then for every
$M \in \dcat{D}_{\mrm{f}}(A, \mrm{gr})$
and every integer $i$, the  graded $A$-module $\opn{H}^i_{\m}(M)$ is cofinite. 
\end{lem}

\begin{proof}
Fix an integer $i$. We will prove that the graded $A$-module 
$\opn{H}^i_{\m}(M)$ is cofinite. Say $\Ga_{\m}$ has cohomological dimension 
$\leq d$ for some natural number $d$. 

\medskip \noindent
Step 1. Consider 
$M' := \opn{smt}^{\leq i + 1}(M)$,
the smart truncation of the complex $M$ below $i + 1$. 
There is a distinguished triangle 
$M' \to M \to M'' \xar{{\,} \triangle {\,}}$
in $\dcat{D}(A, \mrm{gr})$, 
and hence there is a distinguished triangle 
\[ \mrm{R} \Ga_{\m}(M') \to \mrm{R} \Ga_{\m}(M) \to 
\mrm{R} \Ga_{\m}(M'') \xar{ \, \triangle \, } . \]
Because $\mrm{H}^j_{\m}(M'') = 0$ for $j \leq i + 1$, we see that 
$\mrm{H}^i_{\m}(M) \cong \mrm{H}^i_{\m}(M')$. 
Therefore, by replacing $M$ with $M'$, we can assume that 
$M \in \dcat{D}^-_{\mrm{f}}(A, \mrm{gr})$.

\medskip \noindent
Step 2. Now $M \in \dcat{D}^-_{\mrm{f}}(A, \mrm{gr})$.
Choose a resolution $P \to M$, where $P$ is a bounded above complex of finite 
graded-free $A$-modules. Let 
$P' := \opn{stt}^{\geq i - d -2}(P)$, the stupid truncation of $P$ above 
$i - d -2$. Then 
$\mrm{H}^i_{\m}(M) \cong \opn{H}^i(P')$.
But $P'$ is gotten from $A$ by finitely many cones, finite direct sums, 
translations and degree twists. Since $A$ has the left special $\chi$ 
condition, and using Lemma \ref{lem:3856}(1), the graded $A$-module 
$\mrm{H}^j_{\m}(A)$ is cofinite for every $j$. Therefore  
$\mrm{H}^j_{\m}(P')$ is cofinite for every $j$; including $j = i$. 
\end{proof}

Here is a converse to the trivial implication in Proposition 
\ref{prop:3712}(1), under an extra finiteness condition. 

\begin{prop} \label{prop:3870}
If $A$ satisfies the special $\chi$ condition and has finite local 
cohomological 
dimension, then $A$ satisfies the $\chi$ condition. 
\end{prop}

\begin{proof}
This is an immediate consequence of Lemma \ref{lem:3870} and Proposition 
\ref{prop:3855}, applied to both $A$ and $A^{\mrm{op}}$. 
\end{proof}

AS regular graded rings were introduced in Definition \ref{dfn:3782}. 

\begin{cor} \label{cor:3870}
Let $A$ be an AS regular graded ring. 
Then $A$ satisfies the $\chi$ condition and has finite local 
cohomological dimension.
\end{cor}

\begin{proof}
The special $\chi$ condition is an immediate consequence of the AS condition
(Definition \ref{dfn:4545}). 
The graded global cohomological dimension of $A$ bounds its local cohomological 
dimension, so the latter is finite.
By Proposition \ref{prop:3870}, $A$ satisfies the $\chi$ condition.
\end{proof}

Actually this is true also for AS Gorenstein rings, but the proof is harder. 
See Corollary \ref{cor:3780} and Corollary \ref{cor:4655}.

The $\chi$ condition, and the finiteness of local cohomological dimension, 
occur in many important examples of noncommutative graded rings. 
The main results that ensure it involve a lifting argument -- 
Theorems \ref{thm:4546} and \ref{thm:4250}.
Here are several examples. 

\begin{exa} \label{exa:4081}
Let $A := \K [x_1, \ldots, x_n]$, the commutative polynomial ring in 
$n \geq 1$ variables from Example \ref{exa:3757}. 
This is an AS regular graded ring of dimension $n$.
By Corollary \ref{cor:3870} the graded ring $A$ has the $\chi$ condition and 
finite local cohomological dimension.
\end{exa}

\begin{exa} \label{exa:4082}
If $B$ is a commutative noetherian connected graded $\K$-ring, then there is a 
finite homomorphism $A \to B$ from a polynomial ring $A$ as in the previous 
example. By Corollaries \ref{cor:4030} and \ref{cor:4031} the ring $B$ 
satisfies the $\chi$ condition and has finite local cohomological dimension.
\end{exa}

\begin{exa} \label{exa:4080}
Let $A$ be the homogeneous Weyl algebra from Example \ref{exa:3756}.
The element $t \in A$ is central regular of degree $1$, and 
$A / (t) \cong \K[x, y]$,
the commutative polynomial ring in two variables of degree $1$. 
The ring $\K[x, y]$ is AS regular of dimension $2$, so by Theorem 
\ref{thm:4546} the ring $A$ is AS regular of dimension $3$.
By Corollary \ref{cor:3870} the graded ring $A$ has the $\chi$ condition and 
finite local cohomological dimension.

As noted before, in characteristic $0$ the ring $A$ is very noncommutative. 
\end{exa}

\begin{exa} \label{exa:4083}
Let $A$ be the homogeneous universal enveloping algebra from Example 
\ref{exa:4075}. The element $t \in A$ is central regular of degree $1$, and 
$A / (t) \cong \K[x_1, \ldots, x_n]$,
the commutative polynomial ring in $n$ variables. 
The polynomial ring is AS regular of dimension $n$, so 
by Theorem \ref{thm:4546} the ring $A$ is AS regular of dimension $n + 1$.
By Corollary \ref{cor:3870} the graded ring $A$ has the $\chi$ condition and 
finite local cohomological dimension.

As noted in Example \ref{exa:4075}, the ring $A$ could be very 
noncommutative, if $\g$ is very nonabelian (e.g.\ semisimple). 
\end{exa}

We end this subsection with a result that ties the $\chi$ condition to symmetry 
of derived torsion. Recall that 
$\dcat{D}_{(\mrm{f}, \mrm{f})}(A^{\mrm{en}}, \mrm{gr})$
is the full subcategory of 
$\dcat{D}(A^{\mrm{en}}, \mrm{gr})$
on the complexes $M$ such that for every $p$ the graded bimodule $\opn{H}^p(M)$ 
is finite over $A$ and $A^{\mrm{op}}$. 

\begin{thm} \label{thm:4051}
\index{Algebraically graded ring! connected {\indash} of finite local 
cohomological dimension}
\index{Algebraically graded ring! chi@$\chi$ condition on connected}
\index{Complex of algebraically graded modules! with symmetric derived 
$\m$-torsion}
Under Setup \tup{\ref{setup:3855}}, assume $A$ satisfies the $\chi$ condition 
and has finite local cohomological dimension. Then there is a unique 
isomorphism 
\[ \ep : \mrm{R} \Ga_{\m} \iso \mrm{R} \Ga_{\m^{\mrm{op}}} \]
of triangulated functors 
$\dcat{D}^{\mrm{b}}_{(\mrm{f}, \mrm{f})}(A^{\mrm{en}}, \mrm{gr}) \to 
\dcat{D}(A^{\mrm{en}}, \mrm{gr})$, 
such that the diagram 
\[ \UseTips \xymatrix @C=8ex @R=6ex {
\mrm{R} \Ga_{\m}
\ar[r]^{\ep}_{\cong}
\ar[dr]_{\si^{\mrm{R}}_{\m}}
&
\mrm{R} \Ga_{\m^{\mrm{op}}}
\ar[d]^{\si^{\mrm{R}}_{\m^{\mrm{op}}}}
\\
&
\opn{Id}
} \]
is commutative.
\end{thm}

\begin{proof}
According to Proposition \ref{prop:3712}(2) every complex 
$M \in \dcat{D}^{\mrm{b}}_{(\mrm{f}, \mrm{f})}(A^{\mrm{en}}, \mrm{gr})$
has weakly symmetric derived $\m$-torsion. By Theorem \ref{thm:3720} such a 
complex $M$ has symmetric derive torsion, and the symmetry isomorphism $\ep_M$ 
is unique and functorial. 
\end{proof}

\mysubsection{NC MGM Equivalence}
\label{subsec:NC-MGM}

The {\em commutative MGM Equivalence} was \lb explained in Example 
\ref{exa:4595}.
A noncommutative variant of this theory was recently developed by 
R. Vyas and A. Yekutieli \cite{VyYe}, and in this subsection we give an 
adaptation of it to the case of a NC connected graded ring $A$. 
In Section \ref{sec:BDC} below we shall use the NC MGM Equivalence to prove the
existence of balanced trace morphisms (Theorem \ref{thm:4030}), and in Section 
\ref{sec:rigid-DC-NC} the NC MGM Equivalence will help us prove that a balanced 
dualizing complex is graded rigid (Theorem \ref{thm:4400}). 

We continue with Conventions \ref{conv:3700} and \ref{conv:4560}, so $\K$ is a 
base field, and $A$ and $B$ are graded $\K$-rings.

Recall (from Remark \ref{rem:3450}) that the derived category 
$\dcat{D}(A^{\mrm{en}}, \mrm{gr})$ has a biclosed monoidal structure on it, 
with monoidal operation $(- \ot^{\mrm{L}}_{A} -)$ and monoidal unit $A$.
\index{Monoidal category} 
For complexes 
$M \in \dcat{D}(A \ot B^{\mrm{op}}, \mrm{gr})$
and 
$N \in \dcat{D}(B \ot A^{\mrm{op}}, \mrm{gr})$
the left and right unitor isomorphisms are%
\index{Monoidal category! left unitor isomorphism}%
\index{1-Lu@$\opn{lu}$}%
\index{Monoidal category! right unitor isomorphism}%
\index{1-Ru@$\opn{ru}$} 
$\opn{lu} : A \ot^{\mrm{L}}_{A} M \iso M$ 
and
$\opn{ru} : M \ot^{\mrm{L}}_{A} A \iso M$
respectively. The left co-unitor isomorphism%
\index{Monoidal category! left co-unitor isomorphism}%
\index{1-Lcu@$\opn{lcu}$} 
is $\opn{lcu} : \opn{RHom}_{A}(A, M) \iso M$.

\begin{dfn}[\cite{VyYe}] \label{dfn:3995} \mbox{}
\begin{enumerate}
\item A {\em copointed object}%
\index{Monoidal category! copointed object in}
in the monoidal category 
$\dcat{D}(A^{\mrm{en}}, \mrm{gr})$ 
is a pair $(P, \si)$, consisting of an object 
$P \in \dcat{D}(A^{\mrm{en}}, \mrm{gr})$ and a morphism 
$\si : P \to A$ in $\dcat{D}(A^{\mrm{en}}, \mrm{gr})$.

\item The copointed object $(P, \si)$ is called 
{\em idempotent}%
\index{Monoidal category! idempotent copointed object in}
if the morphisms 
\[  \opn{lu} \circ \, (\si \ot^{\mrm{L}}_A \opn{id}) , \ 
\opn{ru} \circ \, (\opn{id}  \ot^{\mrm{L}}_A \, \si) : 
\ P \ot_A^{\mrm{L}} P \to P \]
in $\dcat{D}(A^{\mrm{en}}, \mrm{gr})$ are both isomorphisms.
\end{enumerate}
\end{dfn}

\begin{rem} \label{rem:3995}
Here is an explanation of the name ``copointed''. In a monoidal category, with 
unit object $A$, a {\em point} of an object $P$ is a morphism $A \to P$. 
It thus makes sense to refer to a morphism in the dual direction,  
say $\si : P \to A$, as a {\em copoint} of $P$. 
\end{rem}

\begin{dfn} \label{dfn:3996}
Let  $(P, \si)$ be a copointed object in $\dcat{D}(A^{\mrm{en}}, \mrm{gr})$. 
\begin{enumerate}
\item Define the triangulated functors 
\[ F , G : \dcat{D}(A \ot B^{\mrm{op}}, \mrm{gr}) \to 
\dcat{D}(A \ot B^{\mrm{op}}, \mrm{gr}) \]
to be
$F := P \ot_{A}^{\mrm{L}} (-)$
and
$G := \opn{RHom}_{A}(P, - )$.
 
\item Let $\si : F \to \opn{Id}$ 
and $\tau : \opn{Id} \to G$
be the morphisms of triangulated functors from 
$\dcat{D}(A \ot B^{\mrm{op}}, \mrm{gr})$ to itself 
that are induced by the morphism $\si : P \to A$. Namely
\[ \si_M  := \opn{lu} \circ \, (\si \ot_{A}^{\mrm{L}} \opn{id}_M) : 
F(M) = P \ot_{A}^{\mrm{L}} M \to M  \]
and 
\[ \tau_M := \opn{RHom}_A(\si, \opn{id}_M)  \circ \opn{lcu}^{-1} : 
M \to \opn{RHom}_{A}(P, M) = G(M)  . \]
\end{enumerate}
We refer to $(F, \si)$ and $(G, \tau)$ as the 
{\em copointed and pointed triangulated functors induced by the copointed 
object $(P, \si)$}%
\index{Additive functor! copointed}%
\index{Additive functor! pointed}. 
\end{dfn}

Item (2) of the definition is shown in the 
commutative diagrams below 
in the category $\dcat{D}(A \ot B^{\mrm{op}}, \mrm{gr})$. 

\[ \UseTips \xymatrix @C=8ex @R=6ex {
P \ot_{A}^{\mrm{L}} M
\ar[d]_{\si \ot_{A}^{\mrm{L}} \opn{id}_M}
\ar[dr]^{\si_M}
\\
A \ot_{A}^{\mrm{L}} M
\ar[r]_(0.6){\opn{lu}}
&
M
}
\qquad
\UseTips \xymatrix @C=8ex @R=8ex {
M
\ar[dr]_{\tau_M}
&
\opn{RHom}_{A}(A, M)
\ar[d]^{\opn{RHom}_A(\si, \opn{id}_M)}
\ar[l]_(0.6){\opn{lcu}}
\\
&
\opn{RHom}_{A}(P, M)
} \]

\begin{dfn} \label{dfn:3997}
Let $(F, \si)$ and $(G, \tau)$ be the copointed and pointed triangulated 
functors on $\dcat{D}(A \ot B^{\mrm{op}}, \mrm{gr})$ from Definition 
\tup{\ref{dfn:3996}}. 
\begin{enumerate}
\item Let $\dcat{D}(A \ot B^{\mrm{op}}, \mrm{gr})_F$ be the full 
triangulated subcategory of $\dcat{D}(A \ot B^{\mrm{op}}, \mrm{gr})$ 
on the set of objects  
\[ \bigr\{ M \mid \si_M : F(M) \to M \tup{ is an isomorphism} \, \bigl\} . \]

\item  Let $\dcat{D}(A \ot B^{\mrm{op}}, \mrm{gr})_G$ be the full 
triangulated subcategory of $\dcat{D}(A \ot B^{\mrm{op}}, \mrm{gr})$ 
on the set of objects 
\[ \bigr\{ M \mid \tau_M : M \to G(M) \tup{ is an isomorphism} \, \bigl\} . \]
\end{enumerate}
\end{dfn}

Idempotent (co)pointed triangulated functors were introduced in Definitions 
\ref{dfn:3896} and \ref{dfn:3897}.

\begin{lem} \label{lem:3995}
If the copointed object  $(P, \si)$ is idempotent, then the copointed 
triangulated functor $(F, \si)$ and the 
pointed triangulated functor $(G, \tau)$ on 
$\dcat{D}(A \ot B^{\mrm{op}}, \mrm{gr})$ are idempotent.
\end{lem}

\begin{proof}
For $M \in \dcat{D}(A \ot B^{\mrm{op}}, \mrm{gr})$ there are equalities
(up to the associativity isomorphism of $(- \ot_A^{\mrm{L}} -)$, 
that should be inserted in the locations marked by ``$\dag$''): 
\begin{equation} \label{eqn:4006}
F(\si_M) = \opn{id}_P \ot_A^{\mrm{L}} \, 
\bigl( \opn{lu} \circ \, (\si \ot_{A}^{\mrm{L}} \opn{id}_M) \bigr) 
=^{\dag} \, 
\bigl( \opn{ru} \circ \, (\opn{id}_P \ot_A^{\mrm{L}} \, \si) \bigr)
\ot_{A}^{\mrm{L}} \opn{id}_M 
\end{equation}
and 
\begin{equation} \label{eqn:3996}
\si_{F(M)} = \opn{lu} \circ \, \bigl( \si \ot_A^{\mrm{L}} 
(\opn{id}_P \ot_A^{\mrm{L}} \opn{id}_M) \bigr) 
=^{\dag} \, 
\bigl( \opn{lu} \circ (\si \ot_A^{\mrm{L}} \opn{id}_P) \bigr) \ot_A^{\mrm{L}} 
\opn{id}_M 
\end{equation}
of morphisms
\[ F(F(M)) = P \ot_A^{\mrm{L}} P \ot_A^{\mrm{L}} M \to 
F(M) =  P \ot_A^{\mrm{L}} M  \]
in $\dcat{D}(A \ot B^{\mrm{op}}, \mrm{gr})$.
Because both $\opn{lu} \circ (\si \ot_A^{\mrm{L}} \opn{id}_P)$ and 
$\opn{ru} \circ (\opn{id}_P \ot_A^{\mrm{L}} \, \si)$ 
are isomorphisms in $\cat{D}(A^{\mrm{en}})$, it follows that the morphisms
$F(\si_M)$ and $\si_{F(M)}$ are isomorphisms in 
$\cat{D}(A \ot B^{\mrm{op}})$. 

There are also  equalities (up to the associativity and adjunction isomorphisms 
of $(- \ot_A^{\mrm{L}} -)$ and $\opn{RHom}_A(-,-)$, 
that should be inserted in the locations marked by by ``$\ddag$''): 
\begin{equation} \label{eqn:3997}
\begin{aligned}
& G(\tau_M) = 
\opn{RHom}_A \bigl( \opn{id}_P,  \opn{RHom}_A(\si, \opn{id}_M) \bigr) 
\circ 
\opn{RHom}_A(\opn{id}_P, \opn{lcu}) \\
& \quad =^{\ddag} \, 
\opn{RHom}_A(\si \ot_A^{\mrm{L}} \opn{id}_P, \opn{id}_M) \circ  
\opn{RHom}_A(\opn{lu}, \opn{id}_M) 
\end{aligned}
\end{equation}
and 
\begin{equation} \label{eqn:3998}
\begin{aligned}
& \tau_{G(M)} = \opn{RHom}_A(\si, \opn{id}_{\opn{RHom}_A(P, M)}) \circ 
\opn{lcu} \\
& \quad =^{\ddag} \, 
\opn{RHom}_A(\opn{id}_P \ot_A^{\mrm{L}} \, \si, \opn{id}_M) 
\circ  \opn{RHom}_A(\opn{ru}, \opn{id}_M) 
\end{aligned}
\end{equation}
of morphisms
\[ \begin{aligned}
& G(M) = \opn{RHom}_A(P, M) \to G(G(M))
\\
&
\quad = \opn{RHom}_A \bigl( P, \opn{RHom}_A(P, M) \bigr) \cong
\opn{RHom}_A(P \ot_A^{\mrm{L}} P, M) 
\end{aligned}   \]
in $\dcat{D}(A \ot B^{\mrm{op}}, \mrm{gr})$. 
Because both 
$\opn{lu} \circ \, (\si \ot_A^{\mrm{L}} \opn{id}_P)$ and 
$\opn{ru} \circ \, (\opn{id}_P \ot_A^{\mrm{L}} \, \si)$ 
are isomorphisms  in $\dcat{D}(A^{\mrm{en}}, \mrm{gr})$, it follows that the 
morphisms $G(\tau_M)$ and $\tau_{G(M)}$ are isomorphisms in 
$\dcat{D}(A \ot B^{\mrm{op}}, \mrm{gr})$. 
\end{proof}

\begin{lem} \label{lem:3996}
Consider the functors $F$ and $G$ on $\dcat{D}(A \ot B^{\mrm{op}}, \mrm{gr})$ 
from Definition \tup{\ref{dfn:3996}}. For every 
$M, N \in \dcat{D}(A \ot B^{\mrm{op}}, \mrm{gr})$ there is a bijection 
\[ \opn{Hom}_{\dcat{D}(A \ot B^{\mrm{op}}, \mrm{gr})} \bigl( F(M), N \bigr) 
\cong 
\opn{Hom}_{\dcat{D}(A \ot B^{\mrm{op}}, \mrm{gr})} \bigl( M, G(N) \bigr) , \]
and it is functorial in $M$ and $N$.  
\end{lem}

\begin{proof}
Choose a K-injective resolution $N \to J$ in 
$\dcat{C}(A \ot B^{\mrm{op}}, \mrm{gr})$,
and a K-flat resolution $\til{P} \to P$ in $\dcat{C}(A^{\mrm{en}}, \mrm{gr})$.
The usual Hom-tensor adjunction gives rise to an isomorphism 
\begin{equation} \label{eqn:4000}
\opn{Hom}_{A \ot B^{\mrm{op}}}(\til{P} \ot_A M, J) \cong 
\opn{Hom}_{A \ot B^{\mrm{op}}} \bigl( M, \opn{Hom}_{A}(\til{P}, J) \bigr)
\end{equation}
in $\dcat{C}_{\mrm{str}}(\K, \mrm{gr})$. From this we deduce that 
$\opn{Hom}_{A}(\til{P}, J)$
is K-injective in \lb $\dcat{C}(A \ot B^{\mrm{op}}, \mrm{gr})$.
We see that the isomorphism (\ref{eqn:4000}) represents an isomorphism 
\begin{equation} \label{eqn:4001}
\opn{RHom}_{A \ot B^{\mrm{op}}}(P \ot_A^{\mrm{L}} M, N) \cong 
\opn{RHom}_{A \ot B^{\mrm{op}}} \bigl( M, \opn{RHom}_{A}(P, N) \bigr)
\end{equation}
in $\dcat{D}(\K, \mrm{gr})$. Taking $\opn{H}^0$ in (\ref{eqn:4001}) 
gives us the isomorphism 
\[ \opn{Hom}_{\dcat{D}(A \ot B^{\mrm{op}}, \mrm{gr})}(P \ot_A^{\mrm{L}} M, N) 
\cong 
\opn{Hom}_{\dcat{D}(A \ot B^{\mrm{op}}, \mrm{gr})}
\bigl( M, \opn{RHom}_{A}(P, N) \bigr) \]
in $\dcat{M}(\K, \mrm{gr})$. This is what we want. 
\end{proof}

\begin{lem} \label{lem:4000}
Consider the functors $F$ and $G$ on $\dcat{D}(A \ot B^{\mrm{op}}, \mrm{gr})$ 
from Definition \tup{\ref{dfn:3996}}.  
Assume that the copointed object $(P, \rho)$ is idempotent. 
Then the kernel of $F$ equals the kernel 
of $G$. Namely for every $M \in \dcat{D}(A \ot B^{\mrm{op}}, \mrm{gr})$ we have 
$F(M) = 0$ if and only if $G(M) = 0$. 
\end{lem}

\begin{proof}
We shall use the adjunction formula from Lemma \ref{lem:3996}, with $N = M$. 

First assume $F(M) = 0$. Then 
$\opn{Hom}_{\dcat{D}(A \ot B^{\mrm{op}}, \mrm{gr})} \bigl( F(M), M \bigr)$
is zero, and by Lemma \ref{lem:3996} we see that 
$\opn{Hom}_{\dcat{D}(A \ot B^{\mrm{op}}, \mrm{gr})} \bigl( M, G(M) \bigr)$ 
is zero too. This implies that the morphism $\tau_M : M \to G(M)$ is zero.
Applying $G$ to it we deduce that the morphism 
$G(\tau_M) : G(M) \to G(G(M))$
is zero. But by Lemma \ref{lem:3995} the pointed functor 
$(G, \tau)$ is idempotent, and this means that $G(\tau_M)$ is an isomorphism. 
Therefore $G(M) = 0$. 

Now assume that $G(M) = 0$. Again using Lemma \ref{lem:3996}, but now in 
the reverse direction, we see that the morphism
$\si_M : F(M) \to M$ is zero. Therefore the morphism
$F(\si_M) : F(F(M)) \to F(M)$
is zero. But by Lemma \ref{lem:3995} the copointed functor 
$(F, \si)$ is idempotent, and this means that $F(\si_M)$ is an isomorphism. 
Therefore $F(M) = 0$. 
\end{proof}

Recall that Convention \ref{conv:3700} is in effect. 

\begin{thm}[Abstract Equivalence, \cite{VyYe}] \label{thm:4000}
\index{Monoidal category! idempotent copointed object in}
\index{Additive functor! copointed}
\index{Additive functor! pointed}
\index{Additive functor! idempotent}
Let $A$ and $B$ be graded $\K$-rings, and let $(P, \si)$ be an 
idempotent copointed object in $\dcat{D}(A^{\mrm{en}}, \mrm{gr})$.
Consider the triangulated functors 
\[ F, G : \dcat{D}(A \ot B^{\mrm{op}}, \mrm{gr}) \to 
\dcat{D}(A \ot B^{\mrm{op}}, \mrm{gr}) \]
and the categories 
$\dcat{D}(A \ot B^{\mrm{op}}, \mrm{gr})_F$ and 
$\dcat{D}(A \ot B^{\mrm{op}}, \mrm{gr})_G$ 
from Definitions \tup{\ref{dfn:3996}}
and \tup{\ref{dfn:3997}}. The following hold\tup{:}
\begin{enumerate}
\item The functor $G$ is a right adjoint to $F$. 

\item The copointed triangulated functor $(F, \si)$ and the 
pointed triangulated functor $(G, \tau)$ are idempotent.

\item The categories  $\dcat{D}(A \ot B^{\mrm{op}}, \mrm{gr})_F$ and 
$\dcat{D}(A \ot B^{\mrm{op}}, \mrm{gr})_G$ are the essential 
images of the functors $F$ and $G$ respectively. 

\item The functor 
\[ F : \dcat{D}(A \ot B^{\mrm{op}}, \mrm{gr})_G \to 
\dcat{D}(A \ot B^{\mrm{op}}, \mrm{gr})_F \]
is an equivalence of triangulated categories, with quasi-inverse $G$. 
\end{enumerate}
\end{thm}

\begin{proof}
(1) This is Lemma \ref{lem:3996}. 

\medskip \noindent
(2) This is Lemma \ref{lem:3995}. 

\medskip \noindent
(3) Take an object $M \in \dcat{D}(A \ot B^{\mrm{op}}, \mrm{gr})_G$. Then 
$M \cong G(M)$, so that $M$ is in the essential image of $G$. Conversely, 
suppose there is an isomorphism $\phi : M \iso G(N)$ for some $
N \in \dcat{D}(A \ot B^{\mrm{op}}, \mrm{gr})$. 
We have to prove that $\tau_M$ is an isomorphism. There is a commutative 
diagram 
\[ \UseTips \xymatrix @C=8ex @R=6ex {
M
\ar[r]^{ \phi }
\ar[d]_{ \tau_M }
& 
G(N)
\ar[d]^{ \tau_{G(N)} }
\\
G(M)
\ar[r]^(0.44){ G(\phi) }
&
G(G(N))
} \]
in $\dcat{D}(A \ot B^{\mrm{op}}, \mrm{gr})$ with horizontal isomorphisms.
The idempotence of $(G, \tau)$ says that the morphism $\tau_{G(N)}$ is an 
isomorphism. Therefore $\tau_M$ is an isomorphism. 

A similar argument (with reversed arrows, and using the idempotence of 
$(F, \si)$) tells us that the essential image of 
$F$ is $\dcat{D}(A \ot B^{\mrm{op}}, \mrm{gr})_F$.

\medskip \noindent 
(4) The morphism $\si : P \to A$  sits inside a 
distinguished triangle 
\begin{equation} \label{eqn:4005}
P \xar{\si} A \to N \xar{ \, \triangle \, }  
\end{equation}
in $\dcat{D}(A^{\mrm{en}}, \mrm{gr})$. Let us apply the functor 
$P \ot_A^{\mrm{L}} (-)$ to (\ref{eqn:4005}). We get a distinguished triangle 
\[ P \ot_A^{\mrm{L}} P 
\xar{\opn{ru} \circ \, (\opn{id} \ot_A^{\mrm{L}} \, \si)}
P \to P \ot_A^{\mrm{L}} N \xar{ \, \triangle \, } \]
in $\dcat{D}(A^{\mrm{en}}, \mrm{gr})$.
By the idempotence of $(P, \si)$, the first morphism above is an isomorphism; 
and hence $P \ot_A^{\mrm{L}} N = 0$. Therefore for every 
$M \in \dcat{D}(A \ot B^{\mrm{op}}, \mrm{gr})$ the complex
$F(N \ot_A^{\mrm{L}} M) = P \ot_A^{\mrm{L}} N \ot_A^{\mrm{L}} M$
is zero. Lemma \ref{lem:4000} tells us that 
\begin{equation} \label{eqn:4821}
G(N \ot_A^{\mrm{L}} M) = \opn{RHom}_{A}(P, N \ot_A^{\mrm{L}} M) 
\end{equation}
is zero. 

Now we go back to the distinguished triangle (\ref{eqn:4005}) and we apply to 
it the functor $(-) \ot_A^{\mrm{L}} M$, and then the functor
$\opn{RHom}_{A}(P, - )$. The result is the distinguished triangle 
\[ \opn{RHom}_{A}(P, P \ot_A^{\mrm{L}} M) \xar{\al_M} 
\opn{RHom}_{A}(P, M) \to 
\opn{RHom}_{A}(P, N \ot_A^{\mrm{L}} M)  \xar{ \, \triangle \, } \]
in $\dcat{D}(A \ot B^{\mrm{op}}, \mrm{gr})$. 
By (\ref{eqn:4821}) the third object in this triangle is zero, and it follows 
that $\al_M : G(F(M)) \to G(M)$ is an isomorphism. 
If moreover $M \in \dcat{D}(A \ot B^{\mrm{op}}, \mrm{gr})_G$, then $\tau_M$ is 
an isomorphism too, and thus we have an isomorphism 
\begin{equation} \label{eqn:4822}
\tau_M^{-1} \circ \al_M : G(F(M)) \to M
\end{equation}
in $\dcat{D}(A \ot B^{\mrm{op}}, \mrm{gr})_G$ that's functorial in $M$. 

Similarly, if we apply the functor 
$(-) \ot_A^{\mrm{L}} P$ to (\ref{eqn:4005}), we get a distinguished triangle 
\[ P \ot_A^{\mrm{L}} P 
\xar{\opn{lu} \circ \, (\si \, \ot_A^{\mrm{L}} \opn{id})}
P \to N \ot_A^{\mrm{L}} P \xar{\vartriangle} \]
in $\dcat{D}(A^{\mrm{en}}, \mrm{gr})$.
By the idempotence condition, the first morphism above is an isomorphism; and 
hence $N \ot_A^{\mrm{L}} P = 0$. Therefore for every 
$M \in \dcat{D}(A \ot B^{\mrm{op}}, \mrm{gr})$ the complex
\[ G \bigl( \opn{RHom}_{A}(N, M) \bigr) = 
\opn{RHom}_{A} \bigl( P, \opn{RHom}_{A}(N, M) \bigr) \cong 
\opn{RHom}_{A}(N \ot_A^{\mrm{L}} P, M) \]
is zero. Lemma \ref{lem:4000} tells us that 
\begin{equation} \label{eqn:4823}
F \bigl( \opn{RHom}_{A}(N, M) \bigr) = 
P \ot_A^{\mrm{L}} \opn{RHom}_{A}(N, M) 
\end{equation}
is zero. 

Next we apply the functor $\opn{RHom}_{A}(-, M)$, 
and then the functor $P \ot_A^{\mrm{L}} (-)$, to the distinguished triangle 
(\ref{eqn:4005}). We obtain distinguished triangle 
\[ P \ot_A^{\mrm{L}} \opn{RHom}_{A}(N, M)  \to 
P \ot_A^{\mrm{L}} M \xar{\be_M}
P \ot_A^{\mrm{L}} \opn{RHom}_{A}(P, M)  \xar{ \, \triangle \, } \]
in $\dcat{D}(A \ot B^{\mrm{op}}, \mrm{gr})$. By (\ref{eqn:4823}) the first 
object in this triangle is zero, and so 
$\be_M : F(M) \to F(G(M))$ is an isomorphism. 
If moreover $M \in \dcat{D}(A \ot B^{\mrm{op}}, \mrm{gr})_F$, then $\si_M$ 
is an isomorphism too, and thus we have an isomorphism 
\begin{equation} \label{eqn:4824}
\be_M \circ \si_M^{-1} : M \to F(G(M))
\end{equation}
in $\dcat{D}(A \ot B^{\mrm{op}}, \mrm{gr})_F$ that's functorial in $M$.

The isomorphisms (\ref{eqn:4822}) and (\ref{eqn:4824}) tell us that 
$F$ and $G$ are equivalences between 
$\dcat{D}(A \ot B^{\mrm{op}}, \mrm{gr})_G$ and
$\dcat{D}(A \ot B^{\mrm{op}}, \mrm{gr})_F$, quasi-inverse to each other. 
\end{proof}

We now leave the abstract setting, and return to $\m$-torsion. 
So we are in the situation of Setup \ref{setup:3900}. 
The dedualizing complex 
$P_A := \Ga_{\m}(A) \in \dcat{D}(A^{\mrm{en}}, \mrm{gr})$
from Definition \ref{dfn:3750} is equipped with the morphism 
$\si^{\mrm{R}}_{A} : P_A = \mrm{R} \Ga_{\m}(A) \to A$
from Proposition \ref{prop:3850}. 

\begin{dfn} \label{dfn:4825}
The pair $(P_A, \si^{\mrm{R}}_{A})$ is called the 
{\em dedualizing copointed object}%
\index{Monoidal category! dedualizing copointed object of the {\indash} 
$\dcat{D}(A^{\mrm{en}}, \mrm{gr})$}
of the monoidal category $\dcat{D}(A^{\mrm{en}}, \mrm{gr})$. 
\end{dfn}

\begin{thm}[\cite{VyYe}] \label{thm:4010}
Under Setup \tup{\ref{setup:3900}}, the dedualizing copointed object 
$(P_A, \si^{\mrm{R}}_{A})$ in the monoidal category 
$\dcat{D}(A^{\mrm{en}}, \mrm{gr})$ is idempotent. 
\index{Monoidal category! idempotent copointed object in}
\end{thm}

\begin{proof}
Let's write $\si := \si^{\mrm{R}}_{A}$ and $P := P_A$. 
We shall start by proving that 
\begin{equation}  \label{eqn:4010}
\opn{lu} \circ \, (\si \, \ot_A^{\mrm{L}} \, \opn{id}) : 
P \ot_{A}^{\mrm{L}} P \to P 
\end{equation}
is an isomorphism in $\dcat{D}(A^{\mrm{en}}, \mrm{gr})$. 
Because the forgetful functor 
$\opn{Rest} : \dcat{D}(A^{\mrm{en}}, \mrm{gr}) \lb \to \dcat{D}(A, \mrm{gr})$
is conservative, it is enough if we prove that 
$\opn{Rest}(\opn{lu} \circ (\si \, \ot_A^{\mrm{L}} \, \opn{id}))$
is an isomorphism. Let us introduce the temporary notation 
$P' := \opn{Rest}(P) \in \dcat{D}(A, \mrm{gr})$. 
With this notation, what we have to show is that 
\begin{equation}  \label{eqn:1372}
\opn{lu} \circ \, (\si \, \ot_A^{\mrm{L}} \, \opn{id}) : 
P \ot_{A}^{\mrm{L}} P' \to P' 
\end{equation}
is an isomorphism in $\dcat{D}(A, \mrm{gr})$.

Consider Theorem \ref{thm:3750} with $B = \K$ and 
$M = P' \in \dcat{D}(A, \mrm{gr})$. There is a commutative diagram 
\[ \UseTips \xymatrix @C=12ex @R=6ex {
P \ot_{A}^{\mrm{L}} P'
\ar[r]^{ \ttev^{\mrm{R, L}}_{\m, P'} }
\ar[d]_{ \si \, \ot_A^{\mrm{L}} \, \opn{id} }
&
\mrm{R} \Ga_{\m}(P')
\ar[d]^{ \si^{\mrm{R}}_{P'} }
\\
A \ot^{\mrm{L}}_{A} P'
\ar[r]^{ \opn{lu} }
&
P'
} \]
in $\dcat{D}(A, \mrm{gr})$, and the horizontal arrows are isomorphisms. 
It suffices to prove that 
$\si^{\mrm{R}}_{P'} : \mrm{R} \Ga_{\m}(P') \to P'$
is an isomorphism in $\dcat{D}(A, \mrm{gr})$. 
But there is an isomorphism 
$P' \cong \mrm{R} \Ga_{\m}(A')$,
where $A' := \opn{Rest}(A) \in \dcat{D}(A, \mrm{gr})$. 
So what we need to prove is that 
\[ \si^{\mrm{R}}_{\mrm{R} \Ga_{\m}(A')} : 
\mrm{R} \Ga_{\m}(\mrm{R} \Ga_{\m}(A')) \to \mrm{R} \Ga_{\m}(A') 
\]
is an isomorphism in $\dcat{D}(A, \mrm{gr})$. This is true because the 
copointed triangulated functor 
$(\mrm{R} \Ga_{\m}, \si^{\mrm{R}})$ on $\dcat{D}(A, \mrm{gr})$
is idempotent; see Corollary \ref{cor:3966} with $B := \K$. 

Now we are going to prove that 
\begin{equation} \label{eqn:4011}
\opn{ru} \circ \, (\opn{id}  \ot_A^{\mrm{L}} \, \si)  : 
P \ot_{A}^{\mrm{L}} P \to P 
\end{equation}
is an isomorphism in $\dcat{D}(A^{\mrm{en}}, \mrm{gr})$. 

We have this commutative (up to a canonical isomorphism) diagram in 
$\dcat{D}(A^{\mrm{en}}, \mrm{gr})$
\begin{equation} \label{eqn:4012}
\UseTips \xymatrix @C=12ex @R=6ex {
P \ot_{A}^{\mrm{L}} P
\ar[d]_{ \opn{id}  \ot_A^{\mrm{L}} \, \si }
&
P \ot_{A}^{\mrm{L}} P \ot_{A}^{\mrm{L}} A
\ar[l]_{ \opn{id}  \ot_A^{\mrm{L}} \opn{ru} }
\ar[d]^{ \opn{id}  \ot_A^{\mrm{L}} \, \si \, \ot_A^{\mrm{L}} \opn{id} }
\\
P \ot_{A}^{\mrm{L}} A
&
P \ot_{A}^{\mrm{L}} A \ot^{\mrm{L}}_{A} A
\ar[l]_{ \opn{id}  \ot_A^{\mrm{L}} \opn{ru} }
} 
\end{equation}

Using Theorem \ref{thm:3750} with $B = A$ and 
$M = A \in \dcat{D}(A^{\mrm{en}}, \mrm{gr})$, 
we have this commutative diagram in $\dcat{D}(A^{\mrm{en}}, \mrm{gr})$~:
\begin{equation} \label{eqn:4013}
\UseTips \xymatrix @C=12ex @R=6ex {
P_A \ot_{A}^{\mrm{L}} A
\ar[r]^{ \ttev^{\mrm{R, L}}_{\m, A} }
\ar[d]_{ \si \, \ot_A^{\mrm{L}} \, \opn{id} }
&
\mrm{R} \Ga_{\m}(A)
\ar[d]^{ \si^{\mrm{R}}_{A} }
\\
A \ot^{\mrm{L}}_{A} A
\ar[r]^{ \opn{lu} }
&
A
} 
\end{equation}
Applying the functor $P \ot_A^{\mrm{L}} (-)$ to this diagram, we obtain this 
commutative diagram 
\begin{equation} \label{eqn:4014}
\UseTips \xymatrix @C=15ex @R=6ex {
P \ot_{A}^{\mrm{L}} P \ot_{A}^{\mrm{L}} A
\ar[d]_{ \opn{id}  \ot_A^{\mrm{L}} \, \si \, \ot_A^{\mrm{L}} \, \opn{id} }
\ar[r]^{ \opn{id} \ot_{A}^{\mrm{L}} \, \ttev^{\mrm{R, L}}_{\m, A} }
&
P \ot_{A}^{\mrm{L}} \mrm{R} \Ga_{\m}(A)
\ar[d]_{ \opn{id}  \ot_A^{\mrm{L}} \, \si^{\mrm{R}}_A }
\\
P \ot_{A}^{\mrm{L}} A \ot^{\mrm{L}}_{A} A
\ar[r]^{ \opn{id} \ot_{A}^{\mrm{L}} \opn{lu} }
&
P \ot_{A}^{\mrm{L}} A 
}
\end{equation}
in $\dcat{D}(A^{\mrm{en}}, \mrm{gr})$.
The last move is using the fact that 
$\ttev^{\mrm{R, L}}_{\m, (-)}$ is an isomorphism of functors; this 
yields the next commutative diagram: 
\begin{equation} \label{eqn:4015}
\UseTips \xymatrix @C=15ex @R=6ex {
P \ot_{A}^{\mrm{L}} \mrm{R} \Ga_{\m}(A)
\ar[d]_{ \opn{id}  \ot_A^{\mrm{L}} \, \si^{\mrm{R}}_A }
\ar[r]^{ \ttev^{\mrm{R, L}}_{\m, \mrm{R} \Ga_{\m}(A)}} 
&
\mrm{R} \Ga_{\m}
(\mrm{R} \Ga_{\m}(A))
\ar[d]_{\mrm{R} \Ga_{\m}(\si^{\mrm{R}}_A)}
\\
P \ot_{A}^{\mrm{L}} A 
\ar[r]^{ \ttev^{\mrm{R, L}}_{\m, A} }
&
\mrm{R} \Ga_{\m}(A)
}
\end{equation}
All horizontal arrows in diagrams (\ref{eqn:4012}), (\ref{eqn:4013}), 
(\ref{eqn:4014}) and (\ref{eqn:4015}) are isomorphisms. 
Since $\opn{ru}$ is an isomorphism, to prove that (\ref{eqn:4011}) is an 
isomorphism, it is enough to prove that 
the morphism $\opn{id}  \ot_A^{\mrm{L}} \, \si$,
which is the left vertical arrow in diagram (\ref{eqn:4012}), is an 
isomorphism. Now this last morphism coincides with the 
left vertical arrow in diagram (\ref{eqn:4015}). 
Therefore it is enough to prove that the right vertical arrow in 
diagram (\ref{eqn:4015}) is an isomorphism. This is 
\begin{equation} \label{eqn:4016}
\mrm{R} \Ga_{\m}(\si^{\mrm{R}}_A) : 
\mrm{R} \Ga_{\m} (\mrm{R} \Ga_{\m}(A)) \to
\mrm{R} \Ga_{\m}(A) . 
\end{equation}
Because the functor $\opn{Rest}$ is conservative, it suffices to prove that 
\[ \mrm{R} \Ga_{\m}(\si^{\mrm{R}}_{A'}) : 
\mrm{R} \Ga_{\m}(\mrm{R} \Ga_{\m}(A')) \to
\mrm{R} \Ga_{\m}(A') \]
is an isomorphism in $\dcat{D}(A, \mrm{gr})$, where, as before, we write  
$A' := \opn{Rest}(A) \in \dcat{D}(A, \mrm{gr})$.
This is true because the copointed triangulated functor 
$(\mrm{R} \Ga_{\m}, \si^{\mrm{R}})$ on $\dcat{D}(A, \mrm{gr})$
is idempotent. 
\end{proof}

\begin{dfn} \label{dfn:3790}
Under Setup \tup{\ref{setup:3900}}, consider the dedualizing copointed object 
$(P_A, \si^{\mrm{R}}_{A})$ in the monoidal category 
$\dcat{D}(A^{\mrm{en}}, \mrm{gr})$; see Definition \ref{dfn:4825}.
\begin{enumerate}
\item The {\em abstract derived $\m$-adic completion functor}%
\index{Completion! abstract derived {\indash} functor}%
\index{1-ADCm@$\opn{ADC}_{\m}$}
is the triangulated \lb pointed functor $(\opn{ADC}_{\m}, \tau)$ 
on $\dcat{D}(A \ot B^{\mrm{op}}, \mrm{gr})$ 
that's induced by the copointed object 
$(P_A, \si^{\mrm{R}}_{A})$, in the sense of Definition \ref{dfn:3996}

\item The full subcategory of $\dcat{D}(A \ot B^{\mrm{op}}, \mrm{gr})$ on the 
complexes $M$ such that 
$\tau_{M} : M \to \opn{ADC}_{\m}(M)$
is an isomorphism is denoted by 
$\dcat{D}(A \ot B^{\mrm{op}}, \mrm{gr})_{\mrm{com}}$.
\end{enumerate}
\end{dfn}

To put item (1) in explicit terms,  
$\opn{ADC}_{\m}(M) = \opn{RHom}_A(P_A, M)$,
and the morphism 
$\tau_M : M \to \opn{ADC}_{\m}(M)$ is 
$\tau_M = \opn{RHom}_A(\si^{\mrm{R}}_{A}, \opn{id}_M)$. 

We already have the copointed triangulated functor 
$(\mrm{R} \Ga_{\m}, \si^{\mrm{R}})$ on the category 
$\dcat{D}(A \ot B^{\mrm{op}}, \mrm{gr})$. 
In analogy with Definition \ref{dfn:3790}(2) we make the next definition. 

\begin{dfn} \label{dfn:4020}
Under Setup \tup{\ref{setup:3900}}, we denote by 
$\dcat{D}(A \ot B^{\mrm{op}}, \mrm{gr})_{\mrm{tor}}$
the full subcategory of $\dcat{D}(A \ot B^{\mrm{op}}, \mrm{gr})$ on the 
complexes $M$ such that
$\si^{\mrm{R}}_{M} : \mrm{R} \Ga_{\m}(M) \to M$ 
is an isomorphism.
\end{dfn}

Proposition \ref{prop:3965} says that 
\[ \dcat{D}(A \ot B^{\mrm{op}}, \mrm{gr})_{\mrm{tor}} = 
\dcat{D}_{\mrm{tor}}(A \ot B^{\mrm{op}}, \mrm{gr}) , \]
where the latter is the full subcategory of 
$\dcat{D}(A \ot B^{\mrm{op}}, \mrm{gr})$ on the complexes $M$ whose cohomology 
modules are $\m$-torsion. 

\begin{thm}[NC MGM Equivalence, \cite{VyYe}] \label{thm:3792} 
\index{Noncommutative MGM Equivalence}
Under Setup \tup{\ref{setup:3900}} the following hold.
\begin{enumerate}
\item The categories  
$\dcat{D}(A \ot B^{\mrm{op}}, \mrm{gr})_{\mrm{com}}$ and 
$\dcat{D}(A \ot B^{\mrm{op}}, \mrm{gr})_{\mrm{tor}}$
are the essential images of the functors $\opn{ADC}_{\m}$ and 
$\mrm{R} \Ga_{\m}$ respectively.  

\item The pointed triangulated functor
$(\opn{ADC}_{\m}, \tau)$ and the copointed triangulated functor
$(\mrm{R} \Ga_{\m}, \si^{\mrm{R}})$ are idempotent. 

\item The functor 
\[ \mrm{R} \Ga_{\m} : \dcat{D}(A \ot B^{\mrm{op}}, \mrm{gr})_{\mrm{com}} \to 
\dcat{D}(A \ot B^{\mrm{op}}, \mrm{gr})_{\mrm{tor}} \]
is an equivalence of triangulated categories, with quasi-inverse 
$\opn{ADC}_{\m}$.
\end{enumerate}
\end{thm}

\begin{proof}
By Definition \ref{dfn:3790}, the pointed triangulated functor 
$(\opn{ADC}_{\m}, \tau)$
is the one induced by the copointed object $(P_A, \si^{\mrm{R}}_{A})$.
And by Theorem \ref{thm:3750}, the copointed triangulated functor 
$(\mrm{R} \Ga_{\m}, \si^{\mrm{R}})$
is the one induced by the copointed object $(P_A, \si^{\mrm{R}}_{A})$.
According to Theorem \ref{thm:4010}, the copointed object 
$(P_A, \si^{\mrm{R}}_{A})$ is idempotent.
This means that we can use Theorem \ref{thm:4000}  on abstract equivalence. 
Item (1) here is then a special case of item (3) of Theorem \ref{thm:4000};
item (2) here is a special case of item (2) of Theorem \ref{thm:4000};
and item (3) here is a special case of item (4) of Theorem \ref{thm:4000}.
\end{proof}

\begin{exa} \label{exa:4575}
In Example \ref{exa:4595} we presented the commutative MGM Equivalence. 
Recall the setting: a commutative ring $A$, and an ideal $\a \sub A$ 
generated by a weakly proregular sequence $\bsym{a}$. The module category there 
was $\dcat{M}(A)$. 
If we want to transfer this to the connected graded setting 
of the current section, where the module category is $\dcat{M}(A, \mrm{gr})$,
then we should take $A$ to be a commutative connected 
graded $\K$-ring, and $\bsym{a}$ is a sequence of homogeneous elements that 
generates the augmentation ideal $\m$. Then $A = \K[\bsym{a}]$ is a finitely 
generated commutative $\K$-ring, and thus noetherian. 
Because $A$ is commutative, left $A$-modules are the same as central 
$A$-bimodules. This is the reason we work with the category
$\dcat{M}(A, \mrm{gr})$ and not with $\dcat{M}(A^{\mrm{en}}, \mrm{gr})$.

Regardless of which setting we choose (graded or ungraded), the 
dedualizing complex $P_A$ over the ring $A$ (as an object of 
$\dcat{D}(A, \mrm{gr})$ or $\dcat{D}(A)$ respectively) can be made explicit:
$P_A \cong \opn{K}^{\vee}_{\infty}(A; \bsym{a}) \cong 
\opn{Tel}(A; \bsym{a})$, 
where $\opn{K}^{\vee}_{\infty}(A; \bsym{a})$ is the {\em infinite dual Koszul 
complex}, and $\opn{Tel}(A; \bsym{a})$ is the {\em telescope complex}. 
Since the telescope complex is a bounded complex of free $A$-modules (of 
countable rank), the derived torsion and completion functors have  
particularly nice presentations:
$\mrm{R} \Ga_{\m} = \opn{Tel}(A; \bsym{a}) \ot_A (-)$
and 
$\opn{ADC}_{\m} = \opn{Hom}_A \bigl( \opn{Tel}(A; \bsym{a}), - \bigr)$.
Also the abstract derived completion functor $\opn{ADC}_{\m}$ is canonically 
isomorphism to $\mrm{L} \Lambda_{\m}$, the left derived functor of the 
$\m$-adic completion functor $\Lambda_{\m}$. 
See \cite{PSY} for details. 
\end{exa}

\cleardoublepage
\mysection{Balanced Dualizing Complexes over NC Graded Rings}
\label{sec:BDC}

\AYcopyright

Let $A = \bigoplus_{i \geq 0} A_i$ be a noncommutative noetherian connected 
graded  $\K$-ring, with augmentation ideal $\m = \bigoplus_{i \geq 1} A_i$. 
A {\em balanced dualizing complex} over $A$ is a dualizing complex $R$ over $A$,
in the noncommutative graded sense, that satisfies the 
{\em Noncommutative Graded Local Duality Theorem} (Theorem \ref{thm:4616} 
below). Balanced dualizing complexes were introduced in 1992 by A. Yekutieli 
\cite{Ye1}.

In this section we define balanced dualizing complexes, and then we prove their 
uniqueness and existence, the local duality theorem, and the trace 
functoriality. One of the main results of this section is Corollary 
\ref{cor:4585}, which says that the following two properties are equivalent:
\begin{itemize}
\item[(i)] The ring $A$ satisfies the $\chi$ condition, and it has 
finite local cohomological dimension.

\item[(ii)] The ring $A$ has a balanced dualizing complex.
\end{itemize}
This result is the product of the combined efforts of 
Yekutieli, J.J. Zhang and M. Van den Bergh.
The $\chi$ condition was already discussed in Subsection 
\ref{subsec:sym-der-tor}. M. Artin and Zhang had introduced the $\chi$ 
condition and the finiteness of local cohomological dimension in their 1994 
paper \cite{ArZh}, to ensure that the {\em noncommutative projective scheme} 
$\opn{Proj}(A)$ has ``good geometric properties''. On the other hand, as 
mentioned above, balanced dualizing complexes were designed to satisfy the 
Noncommutative Local Duality Theorem. It is 
remarkable that in the end, these two properties turned out to be equivalent. 

Throughout this section we adhere to Conventions \ref{conv:3700} and 
\ref{conv:4560}. Thus $\K$ is a base field, and all graded rings are 
algebraically graded central $\K$-rings, 
as defined in Section \ref{sec:alg-gra-rings}.
See Remark \ref{rem:4082} regarding the complete case ($A$ is adically complete 
instead of graded) and the arithmetic case (the base ring $\K$ is not a field).

\mysubsection{Graded NC Dualizing Complexes}
\label{subsec:graded-NCDC}

In this subsection we introduce {\em graded noncommutative dualizing 
complexes}. This is a generalization of Grothendieck's original commutative 
definition from \cite{RD} (see Section \ref{sec:dual-cplx-comm-rng}). 

Consider a graded $\K$-ring $A$. 
Complexes of graded $A$-modules were studied in Subsections 
\ref{subsec:alg-gr-mods} and \ref{gr-res-der-fun}, and here we use these 
constructions and results. 
Let us recall that $\dcat{M}(A, \mrm{gr})$ is the abelian category of graded 
$A$-modules, and $\dcat{D}(A, \mrm{gr})$ is its derived category. 

Suppose $B$ is a second graded ring. Given complexes
$M \in \dcat{D}(A \ot B^{\mrm{op}}, \mrm{gr})$ and 
$N \in \dcat{D}(A^{\mrm{en}}, \mrm{gr})$, 
there is the 
{\em noncommutative derived Hom-evaluation morphism}%
\index{Hom-evaluation morphism! NC derived}%
\index{1-EvRRMN@$\opn{ev}^{\mrm{R, R}}_{M, N}$}
\begin{equation} \label{eqn:3915}
\opn{ev}^{\mrm{R, R}}_{M, N} : M \to 
\opn{RHom}_{A^{\mrm{op}}} \bigl( \opn{RHom}_{A}(M, N), N \bigr) 
\end{equation}
in $\dcat{D}(A \ot B^{\mrm{op}}, \mrm{gr})$.
For a choice of a K-injective resolution 
$N \to I$ in $\dcat{C}_{\mrm{str}}(A^{\mrm{en}}, \mrm{gr})$,
the morphism $\opn{ev}^{\mrm{R, R}}_{M, N}$ 
is represented by the homomorphism 
\[ \opn{ev}_{M, I} : 
M \to \opn{Hom}_{A^{\mrm{op}}} \bigl( \opn{Hom}_{A}(M, I), I \bigr) \]
in $\dcat{C}_{\mrm{str}}(A \ot B^{\mrm{op}}, \mrm{gr})$, 
whose formula is 
$\opn{ev}_{M, I}(m)(\phi) := (-1)^{p \cd q} \cd \phi(m)$
for $m \in M^p$ and $\phi \in \opn{Hom}_{A}(M, I)^q$. 

In the special case $B = A$ and 
$M = A \in \dcat{D}(A^{\mrm{en}}, \mrm{gr})$
we have the canonical isomorphisms 
$N \cong \opn{RHom}_{A}(A, N) \cong \opn{RHom}_{A^{\mrm{op}}}(A, N)$,
i.e.\ the left and right co-unitor isomorphisms of the biclosed monoidal 
structure. Then the derived NC Hom-evaluation morphism 
$\opn{ev}^{\mrm{R, R}}_{M, N}$ from (\ref{eqn:3915}) specializes to the 
{\em NC derived homothety morphism through $A$}%
\index{Homothety morphism! noncommutative derived}%
\index{1-HmRNA@$\opn{hm}^{\mrm{R}}_{N, A}$}~:
\begin{equation} \label{eqn:3916}
\opn{hm}^{\mrm{R}}_{N, A}
= \opn{ev}^{\mrm{R, R}}_{A, N} : A \to \opn{RHom}_{A^{\mrm{op}}}(N, N) 
\end{equation}
in $\dcat{D}(A^{\mrm{en}}, \mrm{gr})$.
Similarly, after exchanging $A$ and $A^{\mrm{op}}$, we recover 
the {\em NC derived homothety morphism through $A^{\mrm{op}}$}~:
\begin{equation} \label{eqn:3917}
\opn{hm}^{\mrm{R}}_{N, A^{\mrm{op}}} =
\opn{ev}^{\mrm{R, R}}_{A, N} : A \to \opn{RHom}_{A}(N, N) ,
\end{equation}
also in $\dcat{D}(A^{\mrm{en}}, \mrm{gr})$.
This should be compared to Definition \ref{dfn:3565}, which refers to the 
ungraded NC setting. 

\begin{dfn}[\cite{Ye1}] \label{dfn:3716} 
Let $A$ be a noetherian graded ring. A 
{\em graded NC dualizing complex}%
\index{Dualizing complex! graded NC}
over $A$ is a complex 
$R \in \dcat{D}^{\mrm{b}}(A^{\mrm{en}}, \mrm{gr})$
with the following three properties\tup{:}
\begin{enumerate}
\rmitem{i} {\em Finiteness of cohomology}: for every integer $p$ the 
$A$-bimodule 
$\opn{H}^p(R)$ is a finite graded module over $A$ and over 
$A^{\mrm{op}}$. 

\rmitem{ii} {\em Finite injective dimension}: the complex $R$ has finite 
graded-injective dimension over $A$ and over $A^{\mrm{op}}$. 

\rmitem{iii} {\em NC Derived Morita property}%
\index{Derived Morita property! noncommutative}:
the noncommutative derived homothety morphisms 
\[ \opn{hm}^{\mrm{R}}_{R, A} : A \to \opn{RHom}_{A^{\mrm{op}}}(R, R) \]
and
\[ \opn{hm}^{\mrm{R}}_{R, A^{\mrm{op}}} : A \to \opn{RHom}_A(R, R) \]
in $\dcat{D}(A^{\mrm{en}}, \mrm{gr})$
are both isomorphisms.
\end{enumerate}
\end{dfn}

Condition (i) can be restated as 
$R \in \dcat{D}_{(\mrm{f}, \mrm{f})}(A^{\mrm{en}}, \mrm{gr})$.
It can be shown (using Corollary \ref{cor:3830}(3) and smart truncation, as in 
the proof of Proposition \ref{prop:2115}) 
that condition (ii) is equivalent to the existence of an 
isomorphism $R \cong I$ in $\dcat{D}(A^{\mrm{en}}, \mrm{gr})$,
where $I$ is a bounded complex, and every bimodule $I^p$ is 
graded-injective over $A$ and over $A^{\mrm{op}}$.

\begin{dfn} \label{dfn:3820}
In the situation of Definition \ref{dfn:3716}, 
given another graded ring $B$, the 
{\em duality functors}%
\index{Duality functor}
associated to the dualizing complex $R$ are the triangulated functors 
\[ D_A : \dcat{D}(A \ot B^{\mrm{op}}, \mrm{gr})^{\mrm{op}} \to 
\dcat{D}(B \ot A^{\mrm{op}}, \mrm{gr}), \quad 
D_A := \opn{RHom}_A(-, R)   \]
and 
\[ D_{A^{\mrm{op}}} : 
\dcat{D}(B \ot A^{\mrm{op}}, \mrm{gr})^{\mrm{op}} 
\to \dcat{D}(A \ot B^{\mrm{op}}, \mrm{gr}), \quad 
D_{A^{\mrm{op}}} := \opn{RHom}_{A^{\mrm{op}}}(-, R) .  \]
\end{dfn}

Since the dualizing complex $R$ is usually clear from the context, we can omit 
it from the notation of the duality functors 
$D_A$ and $D_{A^{\mrm{op}}}$. We can also omit $R$ from the notation of 
the related derived Hom-evaluation morphisms, and just write 
$\opn{ev}^{\mrm{R, R}}_{M} : M \to D_{A^{\mrm{op}}}(D_{A}(M))$
and
$\opn{ev}^{\mrm{R, R}}_{M'} : M' \to D_{A}(D_{A^{\mrm{op}}}(M'))$
for complexes
$M \in  \dcat{D}(A \ot B^{\mrm{op}}, \mrm{gr})$
and 
$M' \in  \dcat{D}(B \ot A^{\mrm{op}}, \mrm{gr})$. 

\begin{thm} \label{thm:3820}
Let $A$ and $B$ be graded rings, with $A$ noetherian.
Let $R \in \dcat{D}(A^{\mrm{en}}, \mrm{gr})$ be a graded NC 
dualizing complex over $A$, with associated duality functors $D_A$ and 
$D_{A^{\mrm{op}}}$. Let $\star$ be a boundedness indicator, and let 
$M \in \lb \dcat{D}^{\star}_{(\mrm{f}, ..)}(A \ot B^{\mrm{op}}, \mrm{gr})$.
Then the following hold\tup{:}
\begin{enumerate}
\item The complex $D_A(M)$ belongs to 
$\dcat{D}^{-\star}_{(.., \mrm{f})}(B \ot A^{\mrm{op}}, \mrm{gr})$, 
where $-\star$ is the reversed boundedness indicator.

\item The derived Hom-evaluation morphism 
$\opn{ev}^{\mrm{R, R}}_{M} : M \to D_{A^{\mrm{op}}}(D_{A}(M))$
in $\dcat{D}(A \ot B^{\mrm{op}}, \mrm{gr})$ is an isomorphism. 

\item The functor
\[ D_{A} : 
\dcat{D}^{\star}_{(\mrm{f}, ..)}(A \ot B^{\mrm{op}}, \mrm{gr})^{\mrm{op}} \to
\dcat{D}^{-\star}_{(.., \mrm{f})}(B \ot A^{\mrm{op}}, 
\mrm{gr}) \]
is an equivalence of triangulated categories, with quasi-inverse 
$D_{A^{\mrm{op}}}$. 
\end{enumerate}
\end{thm}

\begin{proof} \mbox{}

\smallskip \noindent
(1) We can forget the ring $B$. Since the functor 
$D_A : \dcat{D}(A, \mrm{gr})^{\mrm{op}} \to \dcat{D}(A^{\mrm{op}}, \mrm{gr})$
has finite cohomological dimension, and since
$D_A(A) = R \in \dcat{M}_{\mrm{f}}(A^{\mrm{op}}, \mrm{gr})$,
this is a consequence of Theorem \ref{thm:2160}(2), slightly 
modified to handle the algebraically graded situation. 

\medskip \noindent 
(2) Because the restriction functor 
$\dcat{D}(A \ot B^{\mrm{op}}, \mrm{gr}) \to \dcat{D}(A, \mrm{gr})$
is conservative, we can forget about the ring $B$. The derived Morita 
property says that 
$\opn{ev}^{\mrm{R, R}}_{M} : M \to D_{A^{\mrm{op}}}(D_{A}(M))$
is an isomorphism for $M = A$. The functors $\opn{Id}$ and 
$D_{A^{\mrm{op}}} \circ D_{A}$
have finite cohomological dimensions. The assertion is then a consequence of 
Theorem \ref{thm:2135}(2), slightly modified to handle 
the algebraically graded situation. 

\medskip \noindent 
(3) This is clear from items (2) and (3).
\end{proof}

We shall require a notion of twisting in the NC graded setting. 

\begin{dfn} \label{dfn:3780}
Let $\phi$ be an automorphism of $A$ in the category 
$\catt{Rng}_{\mrm{gr}} \centover \K$, and let $i$ be an integer. 
The {\em $(\phi, i)$-twist of the bimodule $A$}%
\index{Twisted bimodule}%
\index{1-Afi@$A(\phi, i)$}
is the object
$A(\phi, i) \in \dcat{M}(A^{\mrm{en}}, \mrm{gr})$
defined as follows: as a left graded $A$-module it is the graded-free 
$A$-module $A(i)$ with basis element $e$ in degree $-i$. The right $A$-module 
action is given by the formula 
$e \cd^{\phi} a := \phi(a) \cd e$ for $a \in A$. 
\end{dfn}

To be more explicit, an element $m \in A(\phi, i)$ can be expressed uniquely as 
$m = b \cd e$ with $b \in A$. Then 
$m \cd^{\phi} a = (b \cd e) \cd^{\phi} a = (b \cd \phi(a)) \cd e 
\in A(\phi, i)$.

\begin{dfn} \label{dfn:3925}
An {\em invertible graded $A$-bimodule} is a bimodule 
$L \in \lb \dcat{M}(A^{\mrm{en}}, \mrm{gr})$
such that there exists some 
$L^{\vee} \in \dcat{M}(A^{\mrm{en}}, \mrm{gr})$
satisfying 
$L \ot_A L^{\vee} \cong L^{\vee} \ot_A L \cong A$
in $\dcat{M}(A^{\mrm{en}}, \mrm{gr})$.
The bimodule $L^{\vee}$ is called a quasi-inverse of $L$. 
\end{dfn}

\begin{prop} \label{prop:9325}  
Assume $A$ is connected graded. 
Let $L$ be an invertible graded $A$-bimodule.
Then $L \cong A(\phi, i)$, for a unique automorphism $\phi$ and a unique 
integer $i$. 
\end{prop}

\begin{proof}
Take a quasi-inverse $L^{\vee}$ of $L$. 
Since the tensor product commutes with the ungrading functor, we see that 
$L \ot_A L^{\vee} \cong L^{\vee} \ot_A L \cong A$ 
in $\dcat{M}(\opn{Ungr}(A)^{\mrm{en}})$. Classical Morita theory (see 
e.g.\ \cite[Chapter 4]{Row}) says that $L$ 
and $L^{\vee}$ are finite projective $\opn{Ungr}(A)$-modules on both sides. 
Hence they are flat over $A$ and $A^{\mrm{op}}$. 

Consider the obvious surjection 
\[ \K \cong \K \ot_A L^{\vee} \ot_A L \surj 
(\K \ot_A L^{\vee} \ot_{A} \K) \ot_{} (\K \ot_A L)  \]
in $\dcat{M}(A^{\mrm{en}}, \mrm{gr})$.
Since $L$ is a nonzero finite graded $A$-module and $A^{\mrm{op}}$-module, by 
the graded Nakayama Lemma (Proposition \ref{prop:4590}) the $\K$-modules  
$\K \ot_A L^{\vee} \ot_{A} \K$ and $\K \ot_A L$ are 
nonzero. Therefore $\K \ot_A L$ is a graded-free $\K$-module of rank $1$, so  
$\K \ot_A L \cong \K(i)$ in $\dcat{M}(\K, \mrm{gr})$ for some 
integer $i$. But $L$ is a flat $A$-module, and therefore, as in the classical 
commutative proof, we get an isomorphism
$L \cong A(i)$ in $\dcat{M}(A, \mrm{gr})$. 

By symmetry we also have an isomorphism 
$L \cong A(i)$ in $\dcat{M}(A^{\mrm{op}}, \mrm{gr})$.
Thus there is an element $e \in L_{-i}$ which is a basis of $L$ on both sides. 
For $a \in A$ let $\phi(a) \in A$ be the unique element such that 
$e \cd a = \phi(a) \cd e$. Then $\phi$ is a $\K$-ring automorphism of $A$, and 
$L \cong  A(\phi, i)$ in $\dcat{M}(A^{\mrm{en}}, \mrm{gr})$.

Because the only invertible elements of $A$ are the nonzero elements of 
$A_0 \cong \K$, and these are central elements, $A$ does not have nontrivial 
inner automorphisms. Thus the automorphism $\phi$ is unique. 
\end{proof}

The next theorem is a variant of Theorem \ref{thm:2175}. 

\begin{thm}[Uniqueness of NC Graded DC, \cite{Ye1}] \label{thm:3780} 
Let $A$ be a noetherian connected graded ring, and let  $R$ and $R'$ be graded 
NC dualizing complexes over $A$. Then there is an isomorphism
\[ R' \cong R \ot_A A(\phi, i)[j] \]
in $\dcat{D}(A^{\mrm{en}}, \mrm{gr})$, for a unique automorphism $\phi$ and 
unique integers $i, j$. 
\end{thm}

For the proof we shall need a few lemmas. 
Recall that a module $N \in \dcat{M}(A, \mrm{gr})$ is called bounded below 
if $N_i = 0$ for $i \ll 0$. 

\begin{lem} \label{lem:3916}
Let 
$T, T' \in \dcat{D}^-(A^{\mrm{en}}, \mrm{gr})$ 
satisfy these conditions\tup{:}
\begin{itemize}
\item For every $k$ the graded bimodules $\opn{H}^k(T)$ and $\opn{H}^k(T')$
are bounded below. 

\item There are isomorphisms
$T \ot^{\mrm{L}}_{A} T' \cong T' \ot^{\mrm{L}}_{A} T \cong A$
in $\dcat{D}(A^{\mrm{en}}, \mrm{gr})$.
\end{itemize}
Then $T \cong A(\phi, i)[j]$, for a unique automorphism $\phi$ and unique 
integers $i, j$. 
\end{lem}

\begin{proof}
This is similar to the proof of Lemma \ref{lem:2175}. Define 
$j_1 := \opn{sup}(\opn{H}(T))$ and 
$j'_1 := \opn{sup}(\opn{H}(T'))$. 
By smart truncation we can assume that 
$j_1 = \opn{sup}(T)$ and $j'_1 := \opn{sup}(T')$.
The K\"unneth trick (see Lemma \ref{lem:2177}; it works also in the 
graded setting) we have 
\[ \opn{H}^{j_1}(T) \ot_A \opn{H}^{j_1'}(T') \cong 
\opn{H}^{j_1 + j'_1}(T \ot^{\mrm{L}}_{A} T') \cong \opn{H}^{j_1 + j'_1}(A) \]
in $\dcat{M}(A^{\mrm{en}}, \mrm{gr})$. 
The graded Nakayama Lemma (Proposition \ref{prop:4590}) implies that
$\opn{H}^{j_1}(T) \ot_A \opn{H}^{j'_1}(T') \neq 0$.
Hence we must have $j_1 + j_1' = 0$ and 
$\opn{H}^{j_1}(T) \ot_A \opn{H}^{j'_1}(T') \cong A$.
By symmetry, we also have 
$\opn{H}^{j'_1}(T') \ot_A \opn{H}^{j_1}(T) \cong A$.
According to Proposition \ref{prop:9325} we know that 
$\opn{H}^{j_1}(T) \cong L$, where $L := A(\phi, i)$ for some $\phi$ and $i$. 
Therefore $\opn{H}^{j'_1}(T') \cong L'$,
where $L' := A(\phi^{-1}, -i)$. 

We now forget the left $A$-module structure of $T$, and the right $A$-module 
structure of $T'$. By the projective truncation trick (Lemma \ref{lem:2189})
we get an isomorphism $T \cong L[-j_1] \oplus N$ in 
$\dcat{D}(A^{\mrm{op}}, \mrm{gr})$, where 
$\opn{sup}(N) \leq j_1 - 1$. Similarly there is a an isomorphism 
$T' \cong L'[-j'_1] \oplus N'$ in 
$\dcat{D}(A, \mrm{gr})$, where $\opn{sup}(N') \leq j'_1 - 1$.
Then there is an isomorphism 
\begin{equation} \label{eqn:5125}
\begin{aligned}
& A \cong \opn{H}(T \ot^{\mrm{L}}_{A} T') \cong (L \ot_{A} L') 
\\ 
& \quad 
\oplus \bigl( L[-j_1] \ot_{A} \opn{H}(N') \bigr)
\oplus \bigl( \opn{H}(N) \ot_{A} L'[-j'_1] \bigr)  \oplus 
\opn{H}(N \ot^{\mrm{L}}_{A} N')
\end{aligned} 
\end{equation}
in $\dcat{G}_{\mrm{str}}(\K, \mrm{gr})$. But $A$ is concentrated in 
cohomological degree $0$, and this forces the three summands 
in the second line of (\ref{eqn:5125}) to be 
zero. Since $L'$ is graded-free of rank $1$ over $A$, it follows that 
$\opn{H}(N) = 0$. Therefore $\opn{H}^{k}(T) = 0$ for all $k \neq j_1$.
Letting $j := -j_1$ we get 
$T \cong \opn{H}^{-j}(T)[j] \cong L[j] = A(\phi, i)[j]$
in $\dcat{D}(A^{\mrm{en}}, \mrm{gr})$. 
\end{proof}

\begin{lem} \label{lem:3917}
Let $R, R' \in \dcat{D}(A^{\mrm{en}}, \mrm{gr})$ 
be graded NC dualizing complexes, and let 
$M \in \dcat{D}_{}(A^{\mrm{en}},  \mrm{gr})$.
Then there a morphism 
\[ \psi_M : M \ot^{\mrm{L}}_{A} \opn{RHom}_A(R, R') \iso 
\opn{RHom}_A \bigl( \opn{RHom}_{A^{\mrm{op}}}(M, R), R' \bigr) \]
in $\dcat{D}(A^{\mrm{en}}, \mrm{gr})$, that is functorial in $M$. 
If $M \in \dcat{D}^-_{(.., \mrm{f})}(A^{\mrm{en}},  \mrm{gr})$
then $\psi_M$ is an isomorphism. 
\end{lem}

\begin{proof}
This is the NC graded version of the isomorphism that was used in the proof of 
Theorem \ref{thm:2175}. To define $\psi_M$ we choose a K-projective 
resolution $P \to M$, and K-injective resolutions 
$R \to I$ and $R' \to I'$, all in 
$\dcat{C}_{\mrm{str}}(A^{\mrm{en}}, \mrm{gr})$.
Then $\psi_M$ is represented by the obvious homomorphism 
\[ \til{\psi}_{P, I, I'} : P \ot_{A} \opn{Hom}_A(I, I') \to 
\opn{Hom}_A \bigl( \opn{Hom}_{A^{\mrm{op}}}(P, I), I' \bigr) \]
in $\dcat{C}_{\mrm{str}}(A^{\mrm{en}}, \mrm{gr})$.

To show that $\psi_M$ is an isomorphism when 
$M \in \dcat{D}^-_{(.., \mrm{f})}(A^{\mrm{en}},  \mrm{gr})$,
we can forget the $A$-module structure on $M$, and view $\psi_M$ as a morphism 
in $\dcat{D}(\K, \mrm{gr})$. So we have a morphism $\psi_M$ between 
triangulated functors 
$\dcat{D}^-_{\mrm{f}}(A^{\mrm{op}}, \mrm{gr}) \to \dcat{D}(\K, \mrm{gr})$.
It is clear that $\psi_A$ is an isomorphism. 
By Theorem \ref{thm:2135}(1), 
which is valid also in the algebraically graded 
context, $\psi_M$ is an isomorphism for every 
$M \in \dcat{D}^-_{\mrm{f}}(A^{\mrm{op}}, \mrm{gr})$. 
\end{proof}

\begin{proof}[Proof of Theorem \tup{\ref{thm:3780}}]
Define the duality functors 
\[ D_A := \opn{RHom}_A(-, R) , \quad 
D_{A^{\mrm{op}}} := \opn{RHom}_{A^{\mrm{op}}}(-, R) , \]
\[ D'_A := \opn{RHom}_A(-, R') , \quad 
D'_{A^{\mrm{op}}} := \opn{RHom}_{A^{\mrm{op}}}(-, R') . \]
These are contravariant triangulated functors from 
$\dcat{D}(A^{\mrm{en}}, \mrm{gr})$ to itself. 
Then define the complexes 
$T := (D'_A \circ D_{A^{\mrm{op}}})(A) = \opn{RHom}_{A}(R, R')$
and 
$T' := \lb (D_A \circ D'_{A^{\mrm{op}}})(A) = \opn{RHom}_{A}(R', R)$
in $\dcat{D}(A^{\mrm{en}}, \mrm{gr})$. 
Note that by Theorem \ref{thm:3820} we have 
$T, T' \in \dcat{D}^{\mrm{b}}_{(\mrm{f}, ..)}(A^{\mrm{en}}, \mrm{gr})$.
Hence the graded bimodules $\opn{H}^k(T)$ and $\opn{H}^k(T')$ are all bounded 
below. 

According to Lemma \ref{lem:3917} and Theorem \ref{thm:3820} (applied twice) 
there are isomorphisms 
\[ \begin{aligned}
& T' \ot^{\mrm{L}}_{A} T \cong (D'_A \circ D_{A^{\mrm{op}}})(T') = 
(D'_A \circ D_{A^{\mrm{op}}} \circ D_A \circ D'_{A^{\mrm{op}}})(A) 
\\ & \quad 
\cong (D'_A \circ D'_{A^{\mrm{op}}})(A) \cong A 
\end{aligned} \]
in $\dcat{D}(A^{\mrm{en}}, \mrm{gr})$. 
By symmetry, there is also an isomorphism 
$T \ot^{\mrm{L}}_{A} T' \cong A$ 
in $\dcat{D}(A^{\mrm{en}}, \mrm{gr})$. 
By Lemma \ref{lem:3916} there is an isomorphism
$T \cong A(\phi, i)[j]$
in $\dcat{D}(A^{\mrm{en}}, \mrm{gr})$, for a unique automorphism $\phi$ and 
unique integers $i, j$. 

Finally we obtain these isomorphisms 
\[ \begin{aligned}
& R \ot_A A(\phi, i)[j] \cong 
R \ot^{\mrm{L}}_{A} T \cong^{\dag}
(D'_A \circ D_{A^{\mrm{op}}})(R) 
\\ & \quad 
\cong (D'_A \circ D_{A^{\mrm{op}}} \circ D_A)(A) \cong D'_A(A) = R'
\end{aligned} \]
in $\dcat{D}(A^{\mrm{en}}, \mrm{gr})$.
The isomorphism $\cong^{\dag}$ is from Lemma \ref{lem:3917}. 
\end{proof}

\mysubsection{Balanced DC: Definition, Uniqueness and Local Duality}
\label{subsec:bal-NCDC-LocDu}

Our goal in this subsection is to relate graded NC dualizing complexes with 
derived torsion. We continue with Conventions \ref{conv:3700} and 
\ref{conv:4560}, and we assume the following setup:

\begin{setup} \label{setup:3925}
$\K$ is a base field, and $A$ is a noetherian connected graded 
$\K$-ring, with augmentation ideal $\m$. 
The augmentation ideal of the opposite ring $A^{\mrm{op}}$ is 
$\m^{\mrm{op}}$. 
\end{setup}

Recall the graded bimodule $A^* = \opn{Hom}_{\K}(A, \K)$. 

\begin{dfn}[\cite{Ye1}] \label{dfn:3714} 
Under Setup \tup{\ref{setup:3925}}, a {\em balanced dualizing complex} 
\index{Dualizing complex! balanced}
over $A$ is a pair $(R, \be)$, where:
\begin{itemize}
\item[(B1)] $R \in \dcat{D}(A^{\mrm{en}}, \mrm{gr})$ is a graded NC 
dualizing 
complex over $A$ (Definition \ref{dfn:3716}), with symmetric derived 
$\m$-torsion (Definition \ref{dfn:3779}). 

\item[(B2)] $\be : \mrm{R} \Ga_{\m}(R) \iso A^*$
is an isomorphism in $\dcat{D}(A^{\mrm{en}}, \mrm{gr})$, called a 
{\em balancing isomorphism}.
\end{itemize}
\end{dfn}

The lack of left-right symmetry in item (B2) of this 
definition will be removed in Corollary \ref{cor:3762} below. 

\begin{rem} \label{rem:4030}
Definition \ref{dfn:3714} is a bit more sophisticated than the 
original definition in \cite{Ye1}. There the condition was that 
$\mrm{R} \Ga_{\m}(R) \cong \mrm{R} \Ga_{\m^{\mrm{op}}}(R) \cong A^*$
in $\dcat{D}(A^{\mrm{en}}, \mrm{gr})$, but a balancing isomorphism 
$\be$ was not specified. The current improved definition is influenced by later 
research, especially \cite{YeZh1} and \cite{VyYe}.
\end{rem}

The next theorem is implicit in \cite{Ye1} and \cite{YeZh1}. 

\begin{thm}[Uniqueness of BDC] \label{thm:3710}
Under Setup \tup{\ref{setup:3925}}, suppose that $(R, \be)$ and $(R', \be')$ 
are balanced dualizing complexes over $A$. 
Then there is a unique isomorphism 
$\psi : R' \iso R$
in $\dcat{D}(A^{\mrm{en}}, \mrm{gr})$, 
such that 
$\be \circ \mrm{R} \Ga_{\m}(\psi) = \be'$
as isomorphisms 
$\mrm{R} \Ga_{\m}(R') \iso A^*$
in $\dcat{D}(A^{\mrm{en}}, \mrm{gr})$.
\end{thm}

\begin{proof}
By Theorem \ref{thm:3780} there is an isomorphism 
$\psi^{\dag} : R'\iso R \ot_A L$
in $\dcat{D}(A^{\mrm{en}}, \mrm{gr})$, where 
$L = A(\phi, i)[j]$ is a translate of a twisted bimodule. 
Applying the functor $\mrm{R} \Ga_{\m}$ to $\psi^{\dag}$ we get a diagram of 
isomorphisms
\begin{equation} \label{eqn:3920}
\UseTips \xymatrix @C=10ex @R=6ex {
\mrm{R} \Ga_{\m}(R')
\ar[r]^(0.40){ \mrm{R} \Ga_{\m}(\psi^{\dag}) }_(0.40){\cong}
\ar[d]_{\be'}^{\cong}
&
\mrm{R} \Ga_{\m} \bigl( R \ot_A L \bigr)
\ar[r]^{\ga}_{\cong}
&
\mrm{R} \Ga_{\m}(R) \ot_A L
\ar[d]_{\be \, \ot \, \opn{id}}^{\cong}
\\
A^*
&
&
A^* \ot_A L
}
\end{equation}
in $\dcat{D}(A^{\mrm{en}}, \mrm{gr})$.
The isomorphism $\ga$ is the obvious one. 
Therefore $A^* \cong A^* \ot_A L$ in $\dcat{M}(A^{\mrm{en}}, \mrm{gr})$.
Due to the NC Graded Matlis Duality
(Theorem \ref{thm:4045}) we have $L \cong A$ in 
$\dcat{M}(A^{\mrm{en}}, \mrm{gr})$. 

Let us now rewrite diagram (\ref{eqn:3920}) with $L = A$. This is the solid 
diagram
\begin{equation} \label{eqn:3921}
\UseTips \xymatrix @C=10ex @R=6ex {
\mrm{R} \Ga_{\m}(R')
\ar[r]^(0.50){ \mrm{R} \Ga_{\m}(\psi^{\dag}) }_(0.50){\cong}
\ar[d]_{\be'}^{\cong}
&
\mrm{R} \Ga_{\m}(R)
\ar[d]_{\be}^{\cong}
\\
A^*
\ar@{-->}[r]^{c \cd (-)}_{\cong}
&
A^*
}
\end{equation}
in $\dcat{D}(A^{\mrm{en}}, \mrm{gr})$.
Because the automorphisms of $A^*$ in $\dcat{D}(A^{\mrm{en}}, \mrm{gr})$ are 
multiplication by nonzero elements of $\K$, there is a unique 
$c \in \K^{\times}$ for which the diagram (\ref{eqn:3921}) is commutative.
Hence $\psi := c^{-1} \cd \psi^{\dag}$
is the unique isomorphism $R' \iso R$ satisfying 
$\be \circ \mrm{R} \Ga_{\m}(\psi) = \be'$. 
\end{proof}

The next theorem is a noncommutative version of Grothendieck's Local Duality 
Theorem \cite[Theorem V.6.2]{RD}. Balanced dualizing complexes were 
invented in order to make this theorem hold. 
Theorem \ref{thm:4616} first appeared in \cite{Ye1}, in a slightly different 
formulation, and with a complicated proof. The current formulation of the 
theorem, as well as the proof, are taken from 
\cite[Proposition 3.4]{CWZ}. We thank R. Vyas for pointing out this proof to 
us. 

\begin{thm}[Local Duality] \label{thm:4616}
\index{Local Duality Theorem}  
Under Setup \tup{\ref{setup:3925}}, let 
$(R, \be)$ be a balanced dualizing complex over $A$.
There is a morphism
\[ \xi : \opn{RHom}_A(-, R) \to \bigl( \mrm{R} \Ga_{\m} (-) \bigr)^* \] 
of triangulated functors 
$\dcat{D}(A, \mrm{gr})^{\mrm{op}} \to \dcat{D}(A^{\mrm{op}}, \mrm{gr})$,
such that for every 
$M \in \lb \dcat{D}^{+}_{\mrm{f}}(A, \mrm{gr})$ 
the morphism 
\[ \xi_M : \opn{RHom}_A(M, R) \to \bigl( \mrm{R} \Ga_{\m} (M) \bigr)^* \]
in $\dcat{D}(A^{\mrm{op}}, \mrm{gr})$
is an isomorphism. 
\end{thm}

Note that this result resembles Van den Bergh's Local Duality, Corollary 
\ref{cor:3901} -- but this is misleading: the assumptions are not the same. 

\begin{proof} 
The proof is divided into three steps. 

\smallskip \noindent 
Step 1. Choose K-injective resolutions $R \to J$ and
$\rho : \Ga_{\m}(J) \to K$ in $\dcat{C}_{\mrm{str}}(A^{\mrm{en}}, \mrm{gr})$.
For each $M \in \dcat{D}(A, \mrm{gr})$ choose a K-injective resolution 
$M \to I_M$ in $\dcat{C}_{\mrm{str}}(A, \mrm{gr})$. 
There is an obvious homomorphism 
\[ \phi_M : \opn{Hom}_A(I_M, J) \to 
\opn{Hom}_A \bigl( \Ga_{\m}(I_M), \Ga_{\m}(J) \bigr)  \] 
in $\dcat{C}_{\mrm{str}}(A^{\mrm{op}}, \mrm{gr})$. 
Next we have the homomorphism 
\[ \opn{Hom}_A(\opn{id}, \rho) :  
\opn{Hom}_A \bigl( \Ga_{\m}(I_M), \Ga_{\m}(J) \bigr) \ \to \ 
\opn{Hom}_A \bigl( \Ga_{\m}(I_M), K \bigr) \] 
in $\dcat{C}_{\mrm{str}}(A^{\mrm{op}}, \mrm{gr})$. 
Composing these homomorphisms, and then going to the derived category, we 
obtain a morphism 
\begin{equation} \label{eqn:4835}
\begin{aligned}
& 
\th_M := \opn{Q} \bigl( \opn{Hom}_A(\opn{id}, \rho) \circ \phi_M \bigr) :
\\ & \quad 
\opn{RHom}_A(M, R) \to \opn{RHom}_A \bigl( \mrm{R} \Ga_{\m}(M), 
\mrm{R} \Ga_{\m}(R) \bigr)
\end{aligned}
\end{equation}
in $\dcat{D}(A^{\mrm{op}}, \mrm{gr})$. 
We have the balancing isomorphism 
\begin{equation} \label{eqn:4616}
\be : \mrm{R} \Ga_{\m}(R) \iso A^*
\end{equation}
in $\dcat{D}(A^{\mrm{en}}, \mrm{gr})$.
According to Proposition \ref{prop:3795} there is an isomorphism 
\[ \al_M : 
\opn{RHom}_A \bigl( \mrm{R} \Ga_{\m}(R), A^* \bigr) \iso 
\bigl( \mrm{R} \Ga_{\m}(R) \bigr)^* \]
in $\dcat{D}(A^{\mrm{op}}, \mrm{gr})$, and it is functorial in $M$.  
Let us define the morphism 
\begin{equation} \label{eqn:4615}
\xi_M :=  \al_M \circ \opn{RHom}(\opn{id}, \be) \circ \th_M :
\opn{RHom}_A(M, R) \to \bigl( \mrm{R} \Ga_{\m} (M) \bigr)^*
\end{equation}
in $\dcat{D}(A^{\mrm{op}}, \mrm{gr})$.
By construction the morphism $\xi_M$ is functorial in $M$. 

\medskip \noindent
Step 2. Take the complex $M := R$. 
By definition $R$ has the NC derived Morita property on the $A$ side.
This implies that the element 
$\opn{id}_R \in \opn{H}^0 \bigl( \opn{RHom}_A(R, R) \bigr)$
is a basis of this rank $1$ graded-free $A^{\mrm{op}}$-module. 
Because of the isomorphism (\ref{eqn:4616}) and by graded Matlis Duality,
we know that $\mrm{R} \Ga_{\m}(R)$ has the NC derived Morita property on the 
$A$ side. Therefore the element 
\[ \opn{id}_{\mrm{R} \Ga_{\m}(R)} \in 
\opn{H}^0 \bigl( \opn{RHom}_A \bigl( \mrm{R} \Ga_{\m}(R),
\mrm{R} \Ga_{\m}(R) \bigr) \bigr) \]
is a basis of this rank $1$ graded-free $A^{\mrm{op}}$-module. 
The construction of the morphism $\th_R$ in (\ref{eqn:4835}) shows that 
$\opn{H}^0(\th_R)(\opn{id}_R) = \opn{id}_{\mrm{R} \Ga_{\m}(R)}$.
We conclude that $\opn{H}^0(\th_R)$ is an isomorphism. Since all other 
cohomologies of the complexes below vanish, it follows that 
\begin{equation} \label{eqn:4838}
\begin{aligned}
& 
\th_R  : \opn{RHom}_A(R, R) \to 
\opn{RHom}_A \bigl( \mrm{R} \Ga_{\m}(R), \mrm{R} \Ga_{\m}(R) \bigr)
\end{aligned}
\end{equation}
is an isomorphism in 
$\dcat{D}(A^{\mrm{op}}, \mrm{gr})$.
Therefore, in view of (\ref{eqn:4615}), $\xi_R$ is also an isomorphism. 

\medskip \noindent
Step 3.
Define the functor 
\[ F := \bigl( \mrm{R} \Ga_{\m} (-) \bigr)^* : 
\dcat{D}^{}(A, \mrm{gr})^{\mrm{op}} \to 
\dcat{D}^{}(A^{\mrm{op}}, \mrm{gr}) . \]
Thus $\xi : D_{A} \to F$
is a morphism of triangulated functors, and we must prove that 
$\xi_M : D_{A}(M) \to F(M)$
is an isomorphism for every 
$M \in \dcat{D}^{+}_{\mrm{f}}(A, \mrm{gr})$.
We know that 
$D_{A^{\mrm{op}}} : \dcat{D}^{-}_{\mrm{f}}(A^{\mrm{op}}, \mrm{gr})^{\mrm{op}} 
\to \dcat{D}^{+}_{\mrm{f}}(A, \mrm{gr})$
is an equivalence. Thus it suffices to prove that the morphism  
\begin{equation} \label{eqn:4839}
\ze := \xi \circ \opn{id}_{D_{A^{\mrm{op}}}} : 
D_{A} \circ D_{A^{\mrm{op}}} \to F \circ D_{A^{\mrm{op}}} 
\end{equation}
of functors
$\dcat{D}^{-}_{\mrm{f}}(A^{\mrm{op}}, \mrm{gr}) \to 
\dcat{D}(A^{\mrm{op}}, \mrm{gr})$
is an isomorphism. 

Consider the object 
$A \in \dcat{D}^{-}_{\mrm{f}}(A^{\mrm{op}}, \mrm{gr})$.
By formula (\ref{eqn:4839}) the diagram 
\[ \UseTips \xymatrix @C=8ex @R=6ex {
(D_A \circ D_{A^{\mrm{op}}})(A)
\ar[d]_{\ze_A}
\ar[r]^{\opn{id}}
&
\opn{RHom}_A(R, R) 
\ar[d]^{\xi_R}
\\
(F \circ D_{A^{\mrm{op}}})(A)
\ar[r]^{\opn{id}}
&
\bigl( \mrm{R} \Ga_{\m} (R) \bigr)^* 
} \]
in $\dcat{D}(A^{\mrm{op}}, \mrm{gr})$ is commutative. 
In step 2 we proved that $\xi_R$ is an isomorphism. Hence $\ze_A$ is an 
isomorphism. 

The functors $D_{A}$ and $D_{A^{\mrm{op}}}$ have bounded cohomological 
displacements, and hence so does their composition 
$D_{A} \circ D_{A^{\mrm{op}}}$. (In fact the composition is isomorphic to the 
identity functor, so it has cohomological displacement $[0, 0]$.) 
The functor $\mrm{R} \Ga_{\m}$ has bounded below 
cohomological displacement, and the contravariant functor $(-)^*$ has 
cohomological displacement $[0, 0]$, so their composition $F$ has 
bounded above cohomological displacement.
We conclude that the functors 
$D_{A} \circ D_{A^{\mrm{op}}}$ and $F \circ D_{A^{\mrm{op}}}$ have 
bounded above cohomological displacements.
We have already established that $\ze_A$ is an isomorphism. 
Theorem \ref{thm:2135}(1), which is valid also in the algebraically graded 
setting, says that $\ze_N$ is an isomorphism for every 
$N \in \dcat{D}^{-}_{\mrm{f}}(A^{\mrm{op}}, \mrm{gr})$.
\end{proof}

\begin{cor}[\cite{YeZh0}] \label{cor:4655}  
Under Setup \tup{\ref{setup:3925}}, assume $A$ has a balanced dualizing 
complex. Then $A$ satisfies the $\chi$ condition, and it has finite local
cohomological dimension. 
\end{cor}

\begin{proof}
Given $M \in \dcat{M}_{\mrm{f}}(A, \mrm{gr})$, Theorem \ref{thm:4616}
provides an isomorphism 
\[ \xi_M : \opn{RHom}_A(M, R) \iso \bigl( \mrm{R} \Ga_{\m}(M) \bigr)^* \] 
in $\dcat{D}(A^{\mrm{op}}, \mrm{gr})$. 
This gives an isomorphism 
\begin{equation} \label{eqn:4655}
\opn{H}^{-p}(\xi_M) : \opn{Ext}^{-p}_A(M, R) \iso \bigl( \mrm{H}^p_{\m}(M) 
\bigr)^* 
\end{equation}
in $\dcat{M}(A^{\mrm{op}}, \mrm{gr})$
for every $p$. By Theorem \ref{thm:3820}(1) we know that 
$\opn{Ext}^{-p}_A(M, R) \in \dcat{M}_{\mrm{f}}(A^{\mrm{op}}, \mrm{gr})$.
Using NC graded Matlis Duality (Theorem \ref{thm:4045})  we get
\begin{equation} \label{eqn:3840}
\mrm{H}^p_{\m}(M) \cong \bigl( \opn{Ext}^{-p}_A(M, R) \bigr)^* \in 
\dcat{M}_{\mrm{cof}}(A, \mrm{gr}) . 
\end{equation}
By Proposition \ref{prop:3855} we conclude that $A$ satisfies the left $\chi$ 
condition. The op-symmetry of the situation implies that $A^{\mrm{op}}$ also 
satisfies the left $\chi$ condition. 

Because $R \in \dcat{D}^{\mrm{b}}(A^{\mrm{en}}, \mrm{gr})$
and is nonzero, we have 
$q_0 := \opn{inf}(\opn{H}(R)) \in \Z$. 
So $\opn{Ext}^{q}_A(M, R) = 0$
for all $M \in \dcat{M}(A, \mrm{gr})$ and $q < q_0$.
Taking $p := -q$ and using (\ref{eqn:3840}) we get 
$\mrm{H}^p(\mrm{R} \Ga_{\m}(M)) = \mrm{H}^p_{\m}(M) = 0$
for all $p > -q_0$. So the functor $\mrm{R} \Ga_{\m}$ has finite 
cohomological dimension. The op-symmetry of the situation implies that the 
functor $\mrm{R} \Ga_{\m^{\mrm{op}}}$ also has finite cohomological 
dimension.
\end{proof}

\begin{cor} \label{cor:3762}
Let $(R, \be)$ be a balanced dualizing complex over $A$. Then 
there is a unique symmetry isomorphism
$\ep_{R} : \mrm{R} \Ga_{\m}(R) \iso \mrm{R} \Ga_{\m^{\mrm{op}}}(R)$
in $\dcat{D}(A^{\mrm{en}}, \mrm{gr})$, in the sense of Definition 
\tup{\ref{dfn:3779}}. Hence there is a unique isomorphism
$\be^{\mrm{op}} : \mrm{R} \Ga_{\m^{\mrm{op}}}(R) \iso A^*$ 
$\dcat{D}(A^{\mrm{en}}, \mrm{gr})$ such that 
$\be^{\mrm{op}} \circ \ep_{R} = \be$.
\end{cor}

\begin{proof}
By Corollary \ref{cor:4655} the ring $A$ has finite local
cohomological dimension. So Theorem \ref{thm:3720} applies. 
\end{proof}

\mysubsection{Balanced DC: Existence}
\label{subsec:bal-NCDC-exis}

In this subsection we are going to prove the existence of balanced dualizing 
complexes in two important situations. We shall adhere to Conventions 
\ref{conv:3700} and \ref{conv:4560}, and to Setup \ref{setup:3925}.

The next lemma is \cite[Proposition 4.4]{Ye1}. 

\begin{lem} \label{lem:3855}
Let $R$ be a graded NC dualizing complex over $A$, with duality functors 
$D_A, D_{A^{\mrm{op}}} : \dcat{D}(A^{\mrm{en}}, \mrm{gr}) \to 
\dcat{D}(A^{\mrm{en}}, \mrm{gr})$.
The following four conditions are equivalent.  
\begin{enumerate}
\rmitem{i} There is an isomorphism 
$D_A(\K) \cong \K$ in $\dcat{D}(A^{\mrm{en}}, \mrm{gr})$. 

\rmitem{i$'$} There is an isomorphism 
$D_{A^{\mrm{op}}}(\K) \cong \K$ in $\dcat{D}(A^{\mrm{en}}, \mrm{gr})$. 

\rmitem{ii} There is an isomorphism 
$\mrm{R} \Ga_{\m}(R) \cong A^* \ot_A L$
in $\dcat{D}(A^{\mrm{en}}, \mrm{gr})$, for some 
graded invertible $A$-bimodule $L$ generated in algebraic degree $0$. 

\rmitem{ii$'$} There is an isomorphism 
$\mrm{R} \Ga_{\m^{\mrm{op}}}(R) \cong L' \ot_A A^*$
in $\dcat{D}(A^{\mrm{en}}, \mrm{gr})$, for some 
graded invertible $A$-bimodule $L'$ generated in algebraic degree $0$. 
\end{enumerate}
\end{lem}

Note that by Proposition \ref{prop:9325}, a graded invertible $A$-bimodule $L$ 
generated in algebraic degree $0$ is isomorphic to $A(\phi, 0)$ for some 
automorphism $\phi$. 

\begin{proof} \mbox{}

\smallskip \noindent
(i) $\Rightarrow$ (i$'$): 
It is given that $D_A(\K) \cong \K$. 
By Theorem  \ref{thm:3820} we know that
$\K \cong D_{A^{\mrm{op}}}(D_A(\K))$. 
Together we get $D_{A^{\mrm{op}}}(\K) \cong \K$. 

\medskip \noindent 
(i$'$) $\Rightarrow$ (i): The same, by op-symmetry (replacing $A$ with 
$A^{\mrm{op}}$).

\medskip \noindent 
$\bigl($(i) and (i$'$)$\bigr)$ $\Rightarrow$ (ii):
Recall that a graded invertible bimodule $L$ is free of rank $1$ on both sides. 
So it suffices to prove separately that 
\begin{equation} \label{eqn:3922}
\mrm{R} \Ga_{\m}(R) \cong A^* \ \  \tup{in} \ \dcat{D}(A, \mrm{gr}) 
\end{equation}
and
\begin{equation} \label{eqn:3923}
\mrm{R} \Ga_{\m}(R) \cong A^* \ \ \tup{in} \ 
\dcat{D}(A^{\mrm{op}}, \mrm{gr}) .
\end{equation}

Let $R \to I$ be a minimal graded-injective resolution of $R$ over $A$, so that 
$\mrm{R} \Ga_{\m}(R) \cong \Ga_{\m}(I)$
in $\dcat{D}(A, \mrm{gr})$. 
By Lemma \ref{lem:3840} the subcomplex 
$W := \opn{Soc}(I) = \opn{Hom}_{A}(\K, I) \sub I$
has zero differential. But  
$D_A(\K) = \opn{RHom}_A(\K, R) \lb \cong \opn{Hom}_{A}(\K, I) = W$ 
in $\dcat{D}(A, \mrm{gr})$.
From condition (i) we see that $W \cong \K$ in 
$\dcat{D}(A, \mrm{gr})$; and therefore also 
$W \cong \K$ in $\dcat{C}_{\mrm{str}}(A, \mrm{gr})$.
According to Lemma \ref{lem:3839} and Proposition \ref{prop:4811}, 
$W^p \sub \Ga_{\m}(I^p)$ is an essential submodule.
Again using Lemma \ref{lem:3839} we conclude that 
$\Ga_{\m}(I) \cong A^*$ in $\dcat{C}_{\mrm{str}}(A^{\mrm{en}}, \mrm{gr})$. 
Hence the isomorphism (\ref{eqn:3922}) holds. 

For the other isomorphism, let 
$\dcat{M}_{\mrm{f}/ \K}(A, \mrm{gr}) \sub \dcat{M}(A, \mrm{gr})$
be the full subcategory on the modules that are finite over $\K$, or in other 
words, the finite length $A$-modules. Conditions (i) and (i$'$) imply, using 
induction on length of modules, that the functor 
$\opn{H}^0 \circ \, D_A : \dcat{M}_{\mrm{f}/ \K}(A, \mrm{gr})^{\mrm{op}} \to
\dcat{M}_{\mrm{f}/ \K}(A^{\mrm{op}}, \mrm{gr})$
is an equivalence, with quasi-inverse $\opn{H}^0 \circ \, D_{A^{\mrm{op}}}$. 
For every $j \geq 1$ we have 
$A / \m^j \in \dcat{M}_{\mrm{f}/ \K}(A, \mrm{gr})$. 
The duality provides an isomorphism 
\[ \opn{Hom}_{A^{\mrm{op}}} 
\bigl( \K, \opn{H}^0(D_{A}(A / \m^j)) \bigr)
\cong \opn{Hom}_{A} (A / \m^j, \K) \cong \K . \]
This we rewrite as 
\[ \opn{Hom}_{A^{\mrm{op}}} 
\bigl( \K, \opn{Ext}^0_A(A / \m^j, R) \bigr) \cong \K . \]
By functoriality we can pass to the limit. Using Proposition \ref{prop:3850}(5)
we obtain an isomorphism 
$\opn{Hom}_{A^{\mrm{op}}} \bigl( \K, \opn{H}^0_{\m}(R) \bigr) \cong \K$.
I.e.\ the socle of the graded $A^{\mrm{op}}$-module 
$\opn{H}^0_{\m}(R)$ is $\K$. 
By Proposition \ref{prop:4811} there is an essential 
monomorphism 
$\psi : \opn{H}^0_{\m}(R) \inj A^*$
in $\dcat{M}(A^{\mrm{op}}, \mrm{gr})$. 
But from the previous calculation we already know that 
$\opn{H}^0_{\m}(R) \cong A^*$ in $\dcat{M}(A, \mrm{gr})$, and hence also in 
$\dcat{M}(\K, \mrm{gr})$. As these are degreewise finite graded $\K$-modules,
$\psi$ must be bijective. 
We conclude that the isomorphism (\ref{eqn:3923}) holds.

\medskip \noindent 
$\bigl($(i) and (i$'$)$\bigr)$ $\Rightarrow$ (ii$'$): 
The same, by op-symmetry.

\medskip \noindent 
(ii) $\Rightarrow$ (i): 
By Proposition \ref{prop:3990} there is an isomorphism 
\[ \opn{RHom}_A(\K, R) \cong
\opn{RHom}_A \bigl( \K, \mrm{R} \Ga_{\m}(R) \bigr)  \]
in $\dcat{D}(A^{\mrm{en}}, \mrm{gr})$. 
Now condition (ii) implies that 
$\mrm{R} \Ga_{\m}(R) \cong A^*$
in $\dcat{D}(A, \mrm{gr})$. So we get 
\[ \opn{RHom}_A(\K, R) \cong \opn{RHom}_A(\K, A^*) 
\cong \opn{Hom}_A(\K, A^*) \cong \K \]
in $\dcat{D}(\K, \mrm{gr})$.
This implies that  
$\opn{RHom}_A(\K, R) \cong \K$
in $\dcat{D}(A^{\mrm{en}}, \mrm{gr})$.

\medskip \noindent 
(ii$'$) $\Rightarrow$ (i$'$): The same, by op-symmetry. 
\end{proof}

\begin{dfn}[\cite{Ye1}] \label{dfn:3781} 
Under Setup \tup{\ref{setup:3925}}, a {\em prebalanced dualizing complex} over 
$A$ is a graded NC dualizing complex $R$ over $A$
that satisfies the equivalent conditions of Lemma \ref{lem:3855}.
\end{dfn}

Here is the first of two existence theorems for balanced dualizing complexes. 
The second is Theorem \ref{thm:3713}. 

\begin{thm}[Existence of BDC, \cite{Ye1}] \label{thm:3781}
Under Setup \tup{\ref{setup:3925}}, assume that $A$ has a prebalanced 
dualizing complex. Then $A$ has a balanced dualizing complex.
\index{Dualizing complex! balanced}
\end{thm}

The proof comes after the next lemma. Let 
$\dcat{M}_{\mrm{f}/ \K}(A^{\mrm{en}}, \mrm{gr}) \sub 
\dcat{M}(A^{\mrm{en}}, \mrm{gr})$
be the full subcategory on the bimodules that are finite over $\K$.

\begin{lem} \label{lem:3925}
Suppose $R$ is a prebalanced dualizing complex over $A$ such that 
$\mrm{R} \Ga_{\m}(R) \cong A^*$ in 
$\dcat{D}(A^{\mrm{en}}, \mrm{gr})$. 
Then there is an isomorphism 
$D_A \cong \opn{Hom}_{\K}(-, \K)$
of functors 
$\dcat{M}_{\mrm{f}/ \K}(A^{\mrm{en}}, \mrm{gr})^{\mrm{op}} \to
\dcat{M}_{\mrm{f}/ \K}(A^{\mrm{en}}, \mrm{gr})$.
\end{lem}

\begin{proof}
Let us choose an isomorphism 
$\be : \mrm{R} \Ga_{\m}(R) \iso A^*$
in $\dcat{D}(A^{\mrm{en}}, \mrm{gr})$. 
Take $M \in \dcat{M}_{\mrm{f}/ \K}(A^{\mrm{en}}, \mrm{gr})$. 
There is a sequence of functorial isomorphisms 
\[ \begin{aligned}
& D_A(M) = \opn{RHom}_A(M, R)
\cong^{\dag} \opn{RHom}_A \bigl( M, \mrm{R} \Ga_{\m}(R) \bigr)
\\ & \quad 
\cong^{\ddag} \opn{RHom}_A(M, A^*) 
\cong \opn{Hom}_{\K}(M, \K) 
\end{aligned} \]
in $\dcat{D}(A^{\mrm{en}}, \mrm{gr})$. 
The isomorphism $\cong^{\dag}$ is by Proposition \ref{prop:3990}, and the 
isomorphism $\cong^{\ddag}$ comes from $\be$. 
\end{proof}

\begin{proof}[Proof of Theorem \tup{\ref{thm:3781}}]
The proof is in two steps. 

\smallskip \noindent 
Step 1. Let $R'$ be a prebalanced dualizing complex over $A$. 
This means that there is an isomorphism 
\begin{equation} \label{eqn:4575}
\mrm{R} \Ga_{\m}(R') \cong A^* \ot_A L 
\end{equation}
in $\dcat{D}(A^{\mrm{en}}, \mrm{gr})$, for some 
graded invertible $A$-bimodule $L$ generated in algebraic degree $0$. Define 
$L^{\vee} := \opn{Hom}_A(L, A)$. Thus, if $L \cong A(\phi, i)$, then 
$L^{\vee} \cong A(\phi^{-1}, -i)$ in $\dcat{M}(A^{\mrm{en}}, \mrm{gr})$.

Next let 
\begin{equation} \label{eqn:3929}
R := R'  \ot_A L^{\vee} \in \dcat{D}(A^{\mrm{en}}, \mrm{gr}) ,
\end{equation}
which is a graded dualizing complex. 
It is clear from (\ref{eqn:4575}) that there is an isomorphism 
\begin{equation} \label{eqn:3927}
\be : \mrm{R} \Ga_{\m}(R) \iso A^* \ \ \tup{in} \ \  
\dcat{D}(A^{\mrm{en}}, \mrm{gr}) ,
\end{equation}
and this is the balancing isomorphism that we take. 
To complete the proof that $(R, \be)$ is a balanced dualizing complex, it 
remains to prove that 
\begin{equation} \label{eqn:3925}
\mrm{R} \Ga_{\m^{\mrm{op}}}(R) \cong A^*
\ \ \tup{in} \ \ \dcat{D}(A^{\mrm{en}}, \mrm{gr}) .
\end{equation}
This will be done in the second step. 
(We can't use Theorem \ref{thm:3720} on symmetric derived torsion, since  
that would involve circular reasoning.)

\medskip \noindent 
Step 2. The graded dualizing complex $R$ is prebalanced, by (\ref{eqn:3927}). 
So according to Lemma \ref{lem:3855} there is an isomorphism 
\begin{equation} \label{eqn:3928}
\mrm{R} \Ga_{\m^{\mrm{op}}}(R) \cong L^{\diamond} \ot_A A^* 
\end{equation}
in $\dcat{D}(A^{\mrm{en}}, \mrm{gr})$, for some 
invertible graded $A$-bimodule $L^{\diamond}$ generated in algebraic degree 
$0$. 

For every $j \geq 1$ the bimodule 
$A / \m^j$ belongs to $\dcat{M}_{\mrm{f}/ \K}(A^{\mrm{en}}, \mrm{gr})$,
and it is $\m$-torsion and $\m^{\mrm{op}}$-torsion. 
Therefore we get isomorphisms 
\begin{equation} \label{eqn:3992}
\begin{aligned}
& D_{A^{\mrm{op}}}(A / \m^j) = \opn{RHom}_{A^{\mrm{op}}}(A / \m^j, R)
\\ & \quad 
\cong^{\mrm{(a)}} 
\opn{RHom}_{A^{\mrm{op}}} 
\bigl( A / \m^j, \mrm{R} \Ga_{\m^{\mrm{op}}}(R) \bigr)
\\ & \quad 
\cong^{\mrm{(b)}} 
\opn{RHom}_{A^{\mrm{op}}}(A / \m^j, L^{\diamond} \ot_A  A^*) 
\\ & \quad 
\cong L^{\diamond} \ot_A \opn{RHom}_{A^{\mrm{op}}}(A / \m^j, A^*) 
\\ & \quad 
\cong L^{\diamond} \ot_A  \opn{Hom}_{\K}(A / \m^j, \K) 
= L^{\diamond} \ot_A (A / \m^j)^* 
\end{aligned} 
\end{equation}
in $\dcat{D}(A^{\mrm{en}}, \mrm{gr})$.
Explanation: the isomorphism $\cong^{\mrm{(a)}}$ is by Proposition 
\ref{prop:3990} (transcribed to the ring $A^{\mrm{op}}$); and the isomorphism 
$\cong^{\mrm{(b)}}$ is from equation (\ref{eqn:3928}).

There are also isomorphisms 
\begin{equation} \label{eqn:3993}
\begin{aligned}
&
A / \m^j \cong^1 D_A \bigl( D_{A^{\mrm{op}}}(A / \m^j) \bigr)
\cong^2 D_A \bigl( L^{\diamond} \ot_A (A / \m^j)^* \bigr)
\\ & \quad 
\cong^3 \bigl( L^{\diamond} \ot_A (A / \m^j)^* \bigr)^*
\cong^4 (A / \m^j) \ot_A (L^{\diamond})^{\vee}
\end{aligned} 
\end{equation}
in $\dcat{M}(A^{\mrm{en}}, \mrm{gr})$, where
$(L^{\diamond})^{\vee} := \opn{Hom}_A(L^{\diamond}, A)$.
Here is how they arise: the isomorphism $\cong^1$ is by Theorem \ref{thm:3820};
the isomorphism $\cong^2$ is from equation (\ref{eqn:3992}); 
and the isomorphism $\cong^3$ is due to Lemma \ref{lem:3925}.
For the isomorphism $\cong^4$, say 
$L^{\diamond} \cong A(\nu, k)$; then 
$(L^{\diamond})^{\vee} \cong A(\nu^{-1}, -k)$ and 
$(L^{\diamond})^{*} \cong A^*(\nu^{-1}, -k) = A^* \ot_{A} A(\nu^{-1}, -k)$.

Because the isomorphisms (\ref{eqn:3993}) hold for every $j \geq 1$, it follows 
that $(L^{\diamond})^{\vee} \cong A$ in 
$\dcat{M}(A^{\mrm{en}}, \mrm{gr})$,
and hence 
$L^{\diamond} \cong A$ in 
$\dcat{M}(A^{\mrm{en}}, \mrm{gr})$.
Then the isomorphism (\ref{eqn:3928}) becomes (\ref{eqn:3925}). 
\end{proof}

Recall the AS Gorenstein graded rings from Definition 
\ref{dfn:3782}. For an AS Gorenstein graded ring $A$ there are isomorphisms 
\[ \opn{RHom}_{A}(\K, A) \cong  \opn{RHom}_{A^{\mrm{op}}}(\K, A) \cong
\K(l)[-n] \]
in $\dcat{D}(\K, \mrm{gr})$. The numbers $n$ and $l$ are called the dimension 
and AS index of $A$, respectively. 

\begin{cor}[\cite{Ye1}] \label{cor:3780}
Let $A$ be an AS Gorenstein graded ring, of dimension $n$ and 
with AS index $l$. Then there is a unique automorphism $\phi$ of $A$ 
in $\catt{Rng}_{\mrm{gr}} \centover \K$ such that the complex 
$R_A := A(\phi, -l)[n]  \in \dcat{D}(A^{\mrm{en}}, \mrm{gr})$
is a balanced dualizing complex over $A$. 
\end{cor}

\begin{proof}
The complex 
$R' := A(-l)[n] \in \dcat{D}(A^{\mrm{en}}, \mrm{gr})$
is a prebalanced dualizing complex. Theorem \ref{thm:3781} asserts that $A$ has 
a balanced dualizing complex $R_A$. The construction of $R_A$ in formula 
(\ref{eqn:3929}) in the proof of  Theorem \ref{thm:3781} shows that 
$R_A = R' \ot_A A(\phi, 0) \cong A(\phi, -l)[n]$ 
for some automorphism $\phi$ of $A$. 

The uniqueness of $n, l$ and $\phi$ is a consequence of Theorems \ref{thm:3710} 
and \ref{thm:3780}. 
\end{proof}

Here are two important examples of AS regular graded rings and their balanced 
dualizing complexes. 

\begin{exa} \label{exa:3780}
Suppose $A$ is an {\em elliptic $3$-dimensional Artin-Schelter regular graded 
ring}; see \cite{ATV}. There is a normalizing regular element $g \in A_3$, and 
it defines an automorphism $\phi_g$ of $A$ by the formula 
$\phi_g(a) \cd g = g \cd a$ for $a \in A$. There is a constant 
$\la \in \K^{\times}$
arising from the elliptic curve associated to $A$, and it  defines 
an automorphism $\phi_{\la}$ of $A$ by the formula 
$\phi_{\la}(a) := \la^i \cd a$ for $a \in A_i$. According to 
\cite[Theorem 7.18]{Ye1} the balanced dualizing complex of $A$ is 
$R_A = A(\phi_g \circ \phi_{\la}, -3)[3]$.
\end{exa}

\begin{exa} \label{exa:3990}
Let $A := \K[t_1, \ldots, t_n]$, the commutative polynomial ring in $n$ 
variables, of algebraic degrees $\opn{deg}(t_i) = l_i \geq 1$. This is an 
Artin-Schelter regular graded ring, and a calculation with a Koszul complex 
shows that $\opn{RHom}_A(\K, A) \cong \K(l)[-n]$
in $\dcat{D}(A^{\mrm{en}}, \mrm{gr})$, where 
$l := \sum_{i = 1}^n l_i$. 
Because 
$\mrm{R} \Ga_{\m}(A) \cong A^*(l)[-n]$
is a central $A$-bimodule, it follows that the balanced dualizing dualizing 
complex of $A$ is $R_A := A(-l)[n]$. 

If we want to tie this with the commutative theory from Section
\ref{sec:dual-cplx-comm-rng}, and to algebraic geometry,
then we should take 
$R_A := \Om^n_{A / \K}[n]$,
where $\Om^n_{A / \K}$ is the graded-free $A$-module generated by the 
differential form 
$\d(t_1) \wedge \cdots \wedge \d(t_n)$ 
that has algebraic degree $l$. 
\end{exa}

Now an example that relies on Example \ref{exa:3990}. 

\begin{exa} \label{exa:3991}
Suppose $B$ is a commutative noetherian connected graded $\K$-ring. 
We can find a finite homomorphism $A \to B$ in 
$\catt{Rng}_{\mrm{gr}} \centover \K$ from a commutative polynomial ring $A$. 
Let $R_A := A(-l)[n]$ be the balanced dualizing complex of $A$, as in the 
previous example. Define 
$R_B := \opn{RHom}_{A}(B, R_A) \in \dcat{D}(B^{\mrm{en}}, \mrm{gr})$,
where we use a commutative graded-injective resolution $R_A \to I$, so that 
$R_B$ is a complex of central graded $B$-modules. 
By the arguments of Proposition \ref{prop:2200}, $R_B$ is graded NC 
dualizing complex over $B$. A calculation that we won't perform (but that is an 
easy variant of the proof of Theorem \ref{thm:4030} on the balanced 
trace morphism) shows that $R_B$ is balanced. 

The moral is that the commutative algebraically graded duality theory embeds 
within the noncommutative algebraically graded duality theory.  
\end{exa}

The next theorem is an important converse to Corollary \ref{cor:4655}.

Recall that for a graded NC dualizing complex $R$, the derived homothety 
morphism through $A^{\mrm{op}}$, namely the morphism 
$\opn{hm}^{\mrm{R}}_{R, A^{\mrm{op}}} : A \to \opn{RHom}_{A}(R, R)$ 
in \lb $\dcat{D}(A^{\mrm{en}}, \mrm{gr})$, is an isomorphism. See condition 
(iii) of Definition  \ref{dfn:3716}. 
By Van den Bergh's Local Duality (Corollary \ref{cor:3901}), for every complex 
$M \in \dcat{D}(A^{\mrm{en}},  \mrm{gr})$ there is an isomorphism 
\begin{equation} \label{eqn:4580}
\de_M : \opn{RHom}_A \bigl( M, (P_A)^* \bigr) \iso \mrm{R} \Ga_{\m}(M)^* 
\end{equation}
in $\dcat{D}(A^{\mrm{en}}, \mrm{gr})$.
Here 
$P_A = \mrm{R} \Ga_{\m}(A) \in \dcat{M}(A^{\mrm{en}}, \mrm{gr})$
is the dedualizing complex of $A$. 

\begin{thm}[Existence of BDC, Van den Bergh \cite{VdB}] \label{thm:3713} 
Under Setup \tup{\ref{setup:3925}}, assume the ring $A$ satisfies the special 
$\chi$ condition, and it has finite local cohomological dimension. 
Then $A$ has a balanced dualizing complex $(R_A, \be_A)$. 
\index{Dualizing complex! balanced}

More explicitly, the complex 
$R_A := (P_A)^* \in \dcat{D}(A^{\mrm{en}}, \mrm{gr})$
is a graded dualizing complex over $A$ with symmetric derived $\m$-torsion.
It has a balancing isomorphism
$\be_A : \mrm{R} \Ga_{\m}(R_A) \iso A^*$, 
which is the unique isomorphism 
in $\dcat{D}(A^{\mrm{en}}, \mrm{gr})$
such that the diagram 
\[ \UseTips \xymatrix @C=10ex @R=8ex {
A
\ar[r]_(0.3){\opn{hm}^{\mrm{R}}_{R, A^{\mrm{op}}}}^(0.3){\cong}
\ar@(u,u)[rr]^{(\be_A)^*}
&
\opn{RHom}_{A}(R_A, R_A)
\ar[r]_(0.57){\de_{R_A}}^(0.6){\cong}
&
\mrm{R} \Ga_{\m}(R_A)^*
} \]
in $\dcat{D}(A^{\mrm{en}}, \mrm{gr})$ is commutative.
\end{thm}

\begin{proof} 
The proof proceeds in 5 steps. 

\smallskip \noindent 
Step 1. We are given that $A$ has finite local cohomological dimension and it 
satisfies the special $\chi$ condition. By Proposition \ref{prop:3712}, the 
bimodule $A$ has weakly symmetric derived $\m$-torsion. Theorem \ref{thm:3720} 
says that the bimodule $A$ has symmetric derived $\m$-torsion; namely that there 
is an isomorphism 
\begin{equation} \label{eqn:3980}
\ep_A : P_A = \mrm{R} \Ga_{\m}(A) \iso \mrm{R} \Ga_{\m^{\mrm{op}}}(A) 
\end{equation}
in $\dcat{D}^{\mrm{b}}(A^{\mrm{en}}, \mrm{gr})$. 

The left $\chi$ condition for $A$, with Proposition \ref{prop:3855}, tell us 
that the complex $\mrm{R} \Ga_{\m}(A)$ has cofinite cohomology modules over 
$A$. By the same token, the left $\chi$ condition for $A^{\mrm{op}}$
tells us that the complex $\mrm{R} \Ga_{\m^{\mrm{op}}}(A)$ has cofinite 
cohomology modules over $A^{\mrm{op}}$. Using the isomorphism (\ref{eqn:3980}) 
we conclude that 
$P_A \in \dcat{D}^{\mrm{b}}_{(\mrm{cof}, \mrm{cof})}(A^{\mrm{en}}, \mrm{gr})$.
Now Graded Matlis Duality (Theorem \ref{thm:4045}) says that the complex $R_A$ 
satisfies 
\begin{equation} \label{eqn:4582}
R_A = (P_A)^* \in 
\dcat{D}^{\mrm{b}}_{(\mrm{f}, \mrm{f})}(A^{\mrm{en}}, \mrm{gr}) .
\end{equation}
This is condition (i) of Definition \ref{dfn:3716}. 

\medskip \noindent Step 2. 
Because the functor $\mrm{R} \Ga_{\m}$ has finite cohomological dimension, 
Van den Bergh's Local Duality (Corollary \ref{cor:3901})
says that the functor 
$\opn{RHom}_A \bigl( -, R \bigr) \lb \cong \mrm{R} \Ga_{\m}(-)^*$
also has finite cohomological dimension. This means that the complex $R_A$ has 
finite graded-injective dimension over $A$. Similarly, because 
$\mrm{R} \Ga_{\m^{\mrm{op}}}$ has finite cohomological dimension, the complex 
$R_A$ has finite graded-injective dimension over $A^{\mrm{op}}$. 
This is condition (ii) of Definition \ref{dfn:3716}.

\medskip \noindent Step 3.  
Now we shall prove that $R_A$ has the derived Morita property on the $A$ side. 
We have these isomorphisms in the category 
$\dcat{D}(A^{\mrm{en}}, \mrm{gr})$~:
\begin{equation} \label{eqn:3994}
\begin{aligned}
& \opn{RHom}^{}_A(R_A, R_A) = 
\opn{RHom}^{}_A \bigl( \mrm{R} \Ga_{\m}(A)^*, \mrm{R} \Ga_{\m}(A)^*
\bigr)  \\
& \quad \cong^{1}
\opn{RHom}^{}_{A^{\mrm{op}}} \bigl( \mrm{R} \Ga_{\m}(A), \mrm{R} 
\Ga_{\m}(A) \bigr) \\
& \quad \cong^2 
\opn{RHom}^{}_{A^{\mrm{op}}} 
\bigl( \mrm{R} \Ga_{\m^{\mrm{op}}}(A), \mrm{R} \Ga_{\m^{\mrm{op}}}(A) \bigr) 
\\ & \quad 
\cong^3 \opn{RHom}^{}_{A^{\mrm{op}}} 
\bigl( \mrm{R} \Ga_{\m^{\mrm{op}}}(A), A \bigr) 
\\ & \quad 
\cong^4 \opn{RHom}^{}_{A} 
\bigl( A^*, \mrm{R} \Ga_{\m^{\mrm{op}}}(A)^* \bigr) 
\\ & \quad 
\cong^5 \opn{RHom}^{}_{A} 
\bigl( A^*, \mrm{R} \Ga_{\m}(A)^* \bigr) 
\\ & \quad 
\cong^6 (A^{*})^* \cong^7 A \, .
\end{aligned}
\end{equation}
The isomorphisms $\cong^{1}$, $\cong^{4}$ and $\cong^{7}$ are from graded 
Matlis duality (Theorem \ref{thm:4045}); 
the isomorphisms $\cong^{2}$ and $\cong^{5}$ are by the symmetry isomorphism 
$\ep_A$; the isomorphism $\cong^{3}$ is due to Proposition \ref{prop:3990}; and 
the isomorphism $\cong^{6}$ is due to Van den Bergh's Local Duality (Corollary 
\ref{cor:3901}). 
These isomorphisms respect the derived homothety morphisms 
through $A^{\mrm{op}}$ (in lines 1, 5 and 6), the derived homothety morphisms
through $A$ (in lines 2, 3 and 4); and the canonical 
isomorphisms $A \to (A^*)^*$ and $A \to A$ in line 7. 
(Compare to the proof of Theorem \ref{thm:3681}, and the use of Lemma 
\ref{lem:4361} in it.)
The conclusion is that the derived homothety morphism
$\opn{hm}^{\mrm{R}}_{R, A^{\mrm{op}}} : A \to \opn{RHom}_A(R_A, R_A)$
is an isomorphism. 
This establishes the derived Morita property on the $A$ side for $R_A$. 

\medskip \noindent Step 4.
Replacing $A$ with $A^{\mrm{op}}$, using the symmetry isomorphism $\ep_A$, the 
same calculation as in step 3 shows that $R_A$ has the derived Morita property 
on the $A^{\mrm{op}}$ side. 
Thus condition (iii) of Definition \ref{dfn:3716} is verified. We conclude that 
$R_A$ is a graded NC dualizing complex. 

\medskip \noindent Step 5.  
In this step we prove that $R_A$ is balanced. 
By formula (\ref{eqn:4582}), the derived Morita property on the $A$ side, and 
Van den Bergh's Local Duality (Corollary \ref{cor:3901}), there are isomorphisms
\begin{equation} \label{eqn:4581}
A \underset{\cong}{\xar{\opn{hm}^{\mrm{R}}_{R_A, A^{\mrm{op}}}}} 
\opn{RHom}_{A}(R_A, R_A) \underset{\cong}{\xar{\de_{R_A}}}
\bigl( \mrm{R} \Ga_{\m}(R_A) \bigr)^* 
\end{equation}
in $\dcat{D}(A^{\mrm{en}}, \mrm{gr})$.
See the isomorphism (\ref{eqn:4580}) with $M := R_A = (P_A)^*$.
After we apply the functor $(-)^*$ to (\ref{eqn:4581}), Graded Matlis Duality 
(Theorem \ref{thm:4045}) says that the morphism 
\[ \be := 
\bigl( \de_{R_A} \circ \opn{hm}^{\mrm{R}}_{R_A, A^{\mrm{op}}} \bigr)^* : 
\Ga_{\m}(R_A) \to A^* \]
in $\dcat{D}(A^{\mrm{en}}, \mrm{gr})$ is an isomorphism. This is a 
balancing isomorphism for $R_A$, so item (B2) of Definition 
\ref{dfn:3714} is satisfied.

Similarly, the Derived Morita property on the $A^{\mrm{op}}$ side, 
Van den Bergh's Local Duality for $A^{\mrm{op}}$, and the symmetry isomorphism 
$\ep_A$ give us these isomorphisms 
\[ A \underset{\cong}{\xar{\opn{hm}^{\mrm{R}}_{R_A, A}}} 
\opn{RHom}_{A^{\mrm{op}}}(R_A, R_A) \underset{\cong}{\xar{\de_{R_A}}}
\bigl( \mrm{R} \Ga_{\m^{\mrm{op}}}(R_A) \bigr)^* \]
in $\dcat{D}(A^{\mrm{en}}, \mrm{gr})$.
By Graded Matlis Duality, this yields an isomorphism 
$\mrm{R} \Ga_{\m^{\mrm{op}}}(R_A) \lb \cong A^*$. We see that 
$R_A$ has weakly symmetric derived torsion; and hence, by Theorem 
\ref{thm:3720}, $R_A$ has symmetric derived torsion. This is condition (B1) of 
Definition \ref{dfn:3714}.
\end{proof}

\begin{cor}[Two Equivalent Properties] \label{cor:4585} 
Under Setup \tup{\ref{setup:3925}}, the two properties below are 
equivalent\tup{:}
\begin{enumerate}
\rmitem{i}  The graded ring $A$ satisfies the $\chi$ condition, and it has 
finite local cohomological dimension. 
\index{Algebraically graded ring! chi@$\chi$ condition on connected}
\index{Algebraically graded ring! connected {\indash} of finite local 
cohomological dimension}

\rmitem{ii} The graded ring $A$ has a balanced dualizing complex.
\index{Dualizing complex! balanced}
\end{enumerate}
\end{cor}

\begin{proof}
The implication (i) $\Rightarrow$ (ii) is the Van den Bergh Existence 
Theorem \ref{thm:3713}. 
The reverse implication is Corollary \ref{cor:4655}.
\end{proof}

Here are several other consequences of Theorem \ref{thm:3713}. 
Recall that a dualizing complex $R_A$ over $A$ gives rise to the 
associated duality functors 
\begin{equation} \label{eqn:5100}
D_A, D_{A^{\mrm{op}}} : \dcat{D}(A^{\mrm{en}}, \mrm{gr})^{\mrm{op}} \to
\dcat{D}(A^{\mrm{en}}, \mrm{gr}) ,
\end{equation}
namely 
$D_A = \opn{RHom}_{A}(-, R_A)$ and 
$D_{A^{\mrm{op}}} = \opn{RHom}_{A^{\mrm{op}}}(-, R_A)$.

\begin{thm} \label{thm:4055} 
Under Setup \tup{\ref{setup:3925}}, let $R_A$ be a balanced dualizing complex 
over $A$, with associated duality 
functors $D_A$ and $D_{A^{\mrm{op}}}$.
\begin{enumerate}
\item There is an isomorphism 
$D_A \cong D_{A^{\mrm{op}}}$
between the triangulated functors in \tup{(\ref{eqn:5100})}.
 
\item If 
$M \in \dcat{D}^{\mrm{b}}_{(\mrm{f}, \mrm{f})}(A^{\mrm{en}}, \mrm{gr})$
then $D_A(M)$ belongs to  
$\dcat{D}^{\mrm{b}}_{(\mrm{f}, \mrm{f})}(A^{\mrm{en}}, \mrm{gr})$.

\item The functor 
\[ D_A : \dcat{D}^{\mrm{b}}_{(\mrm{f}, \mrm{f})}(A^{\mrm{en}}, 
\mrm{gr})^{\mrm{op}} \to
\dcat{D}^{\mrm{b}}_{(\mrm{f}, \mrm{f})}(A^{\mrm{en}}, \mrm{gr}) \]
is an equivalence of triangulated categories, and it is its own quasi-inverse. 
\end{enumerate}
\end{thm}

\begin{proof} \mbox{}

\smallskip \noindent
(1) By Theorems \ref{thm:3713} and \ref{thm:3710} there is a canonical 
isomorphism 
$R_A \cong (P_A)^*$ in $\dcat{D}(A^{\mrm{en}}, \mrm{gr})$. 

Take a complex 
$M \in \dcat{D}^{\mrm{b}}_{(\mrm{f}, \mrm{f})}(A^{\mrm{en}}, \mrm{gr})$.
We know that $A$ has the $\chi$ condition and finite local cohomological 
dimension. By Proposition \ref{prop:3712} the complex $M$ has symmetric derived 
torsion. This means that there's an isomorphism 
$\ep_M : \mrm{R} \Ga_{\m}(M) \iso \mrm{R} \Ga_{\m^{\mrm{op}}}(M)$
in $\dcat{D}(A^{\mrm{en}}, \mrm{gr})$. 
By Theorem \ref{thm:3750}  and its opposite version (namely the theorem applied 
to the ring $A^{\mrm{op}}$ instead of $A$), this translates to an isomorphism 
\begin{equation} \label{eqn:5101}
P_A \ot^{\mrm{L}}_{A} M \cong M \ot^{\mrm{L}}_{A} P_A
\end{equation}
in $\dcat{D}(A^{\mrm{en}}, \mrm{gr})$.
Applying $(-)^*$ to  (\ref{eqn:5101}) we obtain an isomorphism 
\begin{equation} \label{eqn:4080}
(P_A \ot^{\mrm{L}}_{A} M)^* \cong (M \ot^{\mrm{L}}_{A} P_A)^*
\end{equation}
in $\dcat{D}(A^{\mrm{en}}, \mrm{gr})$.
Now by Hom-tensor adjunction we get isomorphisms
\begin{equation} \label{eqn:4081}
\begin{aligned}
& (P_A \ot^{\mrm{L}}_{A} M)^* = 
\opn{Hom}_{\K}(P_A \ot^{\mrm{L}}_{A} M , \K) 
\\
& \quad 
\cong 
\opn{RHom}_{A} \bigl( M, \opn{RHom}_{\K}(P_A, \K) \bigr) 
\\
& \quad 
\cong \opn{RHom}_{A} \bigl( M, (P_A)^* \bigr) 
\cong \opn{RHom}_{A}(M, R_A) = D_A(M) .
\end{aligned} 
\end{equation}
Similarly there are isomorphisms
\begin{equation} \label{eqn:4082}
\begin{aligned}
& (M \ot^{\mrm{L}}_{A} P_A)^* = 
\opn{Hom}_{\K}(M \ot^{\mrm{L}}_{A} P_A , \K) 
\\
& \quad 
\cong 
\opn{RHom}_{A^{\mrm{op}}} \bigl( M, \opn{RHom}_{\K}(P_A, \K) \bigr) 
\\
& \quad 
\cong \opn{RHom}_{A^{\mrm{op}}} \bigl( M, (P_A)^* \bigr) 
\cong \opn{RHom}_{A^{\mrm{op}}}(M, R_A) = D_{A^{\mrm{op}}}(M) .
\end{aligned}
\end{equation}
All these isomorphisms are functorial in the argument $M$.
By combining the isomorphisms (\ref{eqn:4080}), (\ref{eqn:4081}) and 
(\ref{eqn:4082}) we get an isomorphism of functors 
$D_A \cong D_{A^{\mrm{op}}}$. 

\medskip \noindent
(2)  By Theorem \ref{thm:3820} we know that 
$D_A(M)$ belongs to 
$\dcat{D}^{\mrm{b}}_{(.., \mrm{f})}(A^{\mrm{en}}, \mrm{gr})$.
The same theorem, applied to the ring $A^{\mrm{op}}$, tells us that 
$D_{A^{\mrm{op}}}(M)$ belongs to 
$\dcat{D}^{\mrm{b}}_{(\mrm{f}, ..)}(A^{\mrm{en}}, \mrm{gr})$.
Thus 
$D_A(M) \cong D_{A^{\mrm{op}}}(M) \in 
\dcat{D}^{\mrm{b}}_{(\mrm{f}, \mrm{f})}(A^{\mrm{en}}, \mrm{gr})$.

\medskip \noindent
(3) We already know, by Theorem \ref{thm:3820}, that 
$D_{A^{\mrm{op}}} \circ D_A \cong \opn{Id}$
as triangulated functors from 
$\dcat{D}^{\mrm{b}}_{(\mrm{f}, \mrm{f})}(A^{\mrm{en}}, \mrm{gr})$
to itself. By item (2) the essential image of 
$D_A$ is 
$\dcat{D}^{\mrm{b}}_{(\mrm{f}, \mrm{f})}(A^{\mrm{en}}, \mrm{gr})$. 
And by item (1) the functors $D_A$ and $D_{A^{\mrm{op}}}$ are isomorphic 
on this triangulated category. 
\end{proof}

\begin{lem} \label{lem:4060}
Under Setup \tup{\ref{setup:3925}}, assume $A$ satisfies the $\chi$ condition 
and has finite local cohomological dimension. 
Then there are inclusions 
$\dcat{D}^{}_{\mrm{f}}(A, \mrm{gr}) \sub
\dcat{D}^{}(A, \mrm{gr})_{\mrm{com}}$
and
$\dcat{D}^{}_{\mrm{cof}}(A, \mrm{gr}) \sub
\dcat{D}^{}(A, \mrm{gr})_{\mrm{tor}}$.
\end{lem}

\begin{proof}
The inclusion 
$\dcat{D}_{\mrm{cof}}^{}(A, \mrm{gr}) \sub
\dcat{D}^{}_{\mrm{tor}}(A, \mrm{gr})$
is trivial. And Proposition \ref{prop:3965} says that 
$\dcat{D}^{}_{\mrm{tor}}(A, \mrm{gr}) = 
\dcat{D}^{}(A, \mrm{gr})_{\mrm{tor}}$.

Now to the complete complexes. 
Take a complex 
$M \in \dcat{D}^{}_{\mrm{f}}(A, \mrm{gr})$.
Then there are isomorphisms 
\begin{equation} \label{eqn:4025}
\begin{aligned}
& \opn{ADC}_{\m}(M) = \opn{RHom}_{A}(P_A, M) 
\\ & \quad
\cong^1 \opn{RHom}_{A^{\mrm{op}}} \bigl( M^*, (P_A)^* \bigr) 
= \opn{RHom}_{A^{\mrm{op}}}(M^*, R_A) 
\\ & \quad
\cong^2 \opn{RHom}_{A^{\mrm{op}}} \bigl( M^*, \mrm{R} \Ga_{\m^{\mrm{op}}}(R_A) 
\bigr) 
\\ & \quad
\cong^3 \opn{RHom}_{A^{\mrm{op}}}(M^*, A^*) 
\cong (M^*)^* \cong^4 M
\end{aligned} 
\end{equation}
in $\dcat{D}(A, \mrm{gr})$. 
The isomorphism $\cong^1$ is due to Theorem \ref{thm:3800}, that applies 
because $M, P_A \in \dcat{C}_{\mrm{dwf}}(A, \mrm{gr})$. 
The isomorphism $\cong^2$ comes from Proposition \ref{prop:3990}, and it 
applies since 
$M^* \in \dcat{D}^{}_{\mrm{cof}}(A^{\mrm{op}}, \mrm{gr}) \sub 
\dcat{D}^{}_{\mrm{tor}}(A^{\mrm{op}}, \mrm{gr})$
and 
$R_A \in \dcat{D}^{\mrm{b}}(A^{\mrm{op}}, \mrm{gr})$. 
The isomorphism $\cong^3$ is because $R_A$ is a balanced dualizing complex.  
The isomorphism $\cong^4$ is by Proposition \ref{prop:4591}, that applies 
because $M \in \dcat{C}_{\mrm{dwf}}(A, \mrm{gr})$. 
We see that $M$ is in the essential image of the functor $\opn{ADC}_{\m}$. 
According to Theorem \ref{thm:3792}(1) it follows that 
$M \in \dcat{D}^{\mrm{b}}(A, \mrm{gr})_{\mrm{com}}$. 
\end{proof}

The next theorem is a refinement of the NC graded MGM Equivalence (Theorem 
\ref{thm:3792}). Given another graded ring $B$, there are triangulated functors 
\[ \mrm{R} \Ga_{\m}, \, \opn{ADC}_{\m} : 
\dcat{D}(A \ot B^{\mrm{op}}, \mrm{gr}) \to 
\dcat{D}(A \ot B^{\mrm{op}}, \mrm{gr}) . \]
They are the derived $\m$-torsion and the abstract derived $\m$-adic 
completion, respectively. 

\begin{thm} \label{thm:4070}
Under Setup \tup{\ref{setup:3925}}, assume $A$ satisfies the $\chi$ condition 
and has finite local cohomological dimension. Let $B$ be some graded ring. 
\begin{enumerate}
\item If 
$M \in \dcat{D}^{\mrm{b}}_{(\mrm{f}, ..)}(A \ot B^{\mrm{op}}, \mrm{gr})$
then 
$\mrm{R} \Ga_{\m}(M) \in 
\dcat{D}^{\mrm{b}}_{(\mrm{cof}, ..)}(A \ot B^{\mrm{op}}, \mrm{gr})$.

\item For 
$N \in \dcat{D}^{\mrm{b}}_{(\mrm{cof}, ..)}(A \ot B^{\mrm{op}}, \mrm{gr})$
there is an isomorphism 
$\opn{ADC}_{\m}(N) \cong \lb D_{A^{\mrm{op}}}(N^*)$ 
in $\dcat{D}(A \ot B^{\mrm{op}}, \mrm{gr})$,
and it is functorial in $N$. 

\item If 
$N \in \dcat{D}^{\mrm{b}}_{(\mrm{cof}, ..)}(A \ot B^{\mrm{op}}, \mrm{gr})$
then 
$\opn{ADC}_{\m}(N) \in 
\dcat{D}^{\mrm{b}}_{(\mrm{f}, ..)}(A \ot B^{\mrm{op}}, \mrm{gr})$.

\item The functor 
\[ \mrm{R} \Ga_{\m} : 
\dcat{D}^{\mrm{b}}_{(\mrm{f}, ..)}(A \ot B^{\mrm{op}}, \mrm{gr}) \to 
\dcat{D}^{\mrm{b}}_{(\mrm{cof}, ..)}(A \ot B^{\mrm{op}}, \mrm{gr}) \]
is an equivalence of triangulated categories, with quasi-inverse 
$\opn{ADC}_{\m}$. 
\end{enumerate}
\end{thm}

\begin{proof} \mbox{}

\smallskip \noindent
(1) The ring $B$ is irrelevant, so we may omit it. 
The $\chi$ condition and finiteness of local cohomology imply that 
$\mrm{R} \Ga_{\m}(M) \in \dcat{D}^{\mrm{b}}_{\mrm{cof}}(A, \mrm{gr})$.

\medskip \noindent 
(2) Take a complex 
$N \in \dcat{D}^{\mrm{b}}_{\mrm{cof}}(A \ot B^{\mrm{op}}, \mrm{gr})$.
We have these isomorphisms 
\[ \begin{aligned}
& \opn{ADC}_{\m}(N) = \opn{RHom}_{A}(P_A, N) 
\cong^{\dag} \opn{RHom}_{A^{\mrm{op}}} \bigl( N^*, (P_A)^* \bigr) 
\\ & \quad
= \opn{RHom}_{A^{\mrm{op}}}(N^*, R_A) 
= D_{A^{\mrm{op}}}(N^*) 
\end{aligned} \]
in $\dcat{D}(A\ot B^{\mrm{op}}, \mrm{gr})$.
The isomorphism $\cong^{\dag}$ is due to Theorem \ref{thm:3800}, that applies 
because 
$N \in \dcat{C}_{\mrm{dwf}}(A\ot B^{\mrm{op}}, \mrm{gr})$ and 
$P_A \in \dcat{C}_{\mrm{dwf}}(A^{\mrm{en}}, \mrm{gr})$. 

\medskip \noindent 
(3) Again we can forget about $B$. 
Take $N \in  \dcat{D}^{\mrm{b}}_{\mrm{cof}}(A, \mrm{gr})$.
Then  
$N^* \in \dcat{D}_{\mrm{f}}(A^{\mrm{op}}, \mrm{gr})$,
so by Theorem \ref{thm:3820} we have 
$D_{A^{\mrm{op}}}(N^*) \in \dcat{D}^{\mrm{b}}_{\mrm{f}}(A, \mrm{gr})$.
But by item (2) we know that 
$\opn{ADC}_{\m}(N) \cong D_{A^{\mrm{op}}}(N^*)$. 
            
\medskip \noindent 
(4) Items (1) and (3) tell us that the functors 
$\mrm{R} \Ga_{\m}$ and $\opn{ADC}_{\m}$ send the categories 
$\dcat{D}^{\mrm{b}}_{(\mrm{f}, ..)}(A \ot B^{\mrm{op}}, \mrm{gr})$
and 
$\dcat{D}^{\mrm{b}}_{(\mrm{cof}, ..)}(A \ot B^{\mrm{op}}, \mrm{gr})$
to each other, respectively. 
Lemma \ref{lem:4060} says that the NC MGM Equivalence (Theorem \ref{thm:3792}) 
applies here, namely that there is an isomorphism of triangulated functors 
$\opn{ADC}_{\m} \circ \, \mrm{R} \Ga_{\m} \cong \opn{Id}$
from 
$\dcat{D}^{\mrm{b}}_{(\mrm{f}, ..)}(A \ot B^{\mrm{op}}, \mrm{gr})$
to itself, and an isomorphism of triangulated functors 
$\mrm{R} \Ga_{\m} \circ \opn{ADC}_{\m} \cong \opn{Id}$
from
$\dcat{D}^{\mrm{b}}_{(\mrm{cof}, ..)}(A \ot B^{\mrm{op}}, \mrm{gr})$
to itself.
\end{proof}

In Theorem \ref{thm:4051} we saw that there is an isomorphism 
$\ep : \mrm{R} \Ga_{\m} \iso \mrm{R} \Ga_{\m^{\mrm{op}}}$
of triangulated functors 
$\dcat{D}^{\mrm{b}}_{(\mrm{f}, \mrm{f})}(A^{\mrm{en}}) \to 
\dcat{D}(A^{\mrm{en}}, \mrm{gr})$.
The next corollary produces a similar isomorphism for the abstract derived 
completion functors 
\begin{equation} \label{eqn:4058}
\opn{ADC}_{\m}, \opn{ADC}_{\m^{\mrm{op}}} : \dcat{D}(A^{\mrm{en}}, \mrm{gr}) 
\to 
\dcat{D}(A^{\mrm{en}}, \mrm{gr}) , 
\end{equation}
where we define
$\opn{ADC}_{\m^{\mrm{op}}} := \opn{RHom}_{A^{\mrm{op}}}(P_A, -)$. 

\begin{cor} \label{cor:4080}   
Under Setup \tup{\ref{setup:3925}}, assume $A$ satisfies the $\chi$ condition 
and has finite local cohomological dimension. Then there is an isomorphism 
$\ga : \lb \opn{ADC}_{\m} \iso \opn{ADC}_{\m^{\mrm{op}}}$
of triangulated functors 
$\dcat{D}^{\mrm{b}}_{(\mrm{cof}, \mrm{cof})}(A^{\mrm{en}}, \mrm{gr}) \to 
\dcat{D}(A^{\mrm{en}}, \mrm{gr})$.
\end{cor}

\begin{proof}
Take a complex 
$N \in \dcat{D}^{\mrm{b}}_{(\mrm{cof}, \mrm{cof})}(A^{\mrm{en}}, \mrm{gr})$.
By Theorem \ref{thm:4070}(2), there is an isomorphism 
$\opn{ADC}_{\m}(N) \cong D_{A^{\mrm{op}}}(N^*)$
in $\dcat{D}^{\mrm{b}}(A^{\mrm{en}}, \mrm{gr})$.
Similarly (replacing $A$ with $A^{\mrm{op}}$) there is an isomorphism 
$\opn{ADC}_{\m^{\mrm{op}}}(N) \cong D_{A}(N^*)$
in $\dcat{D}^{\mrm{b}}(A^{\mrm{en}}, \mrm{gr})$.
Now 
$N^* \in \dcat{D}^{\mrm{b}}_{(\mrm{f}, \mrm{f})}(A^{\mrm{en}}, \mrm{gr})$,
so by Theorem \ref{thm:4055} we have 
$D_{A^{\mrm{op}}}(N^*) \cong D_{A}(N^*)$.
Since we relied on functorial isomorphisms, this is also a functorial 
isomorphism.   
\end{proof}

\begin{cor} \label{cor:4070}
Assume $A$ satisfies the $\chi$ condition and has finite local cohomological 
dimension. Then the functor 
\[ \mrm{R} \Ga_{\m} : 
\dcat{D}^{\mrm{b}}_{(\mrm{f}, \mrm{f})}(A^{\mrm{en}}, \mrm{gr}) \to
\dcat{D}^{\mrm{b}}_{(\mrm{cof}, \mrm{cof})}(A^{\mrm{en}}, \mrm{gr}) \]
is an equivalence of triangulated categories, with quasi-inverse
$\opn{ADC}_{\m}$. 
\end{cor}

\begin{proof}
Taking $B := A$ in Theorem \ref{thm:4070}(4) we get an equivalence 
\begin{equation} \label{eqn:4845}
\mrm{R} \Ga_{\m} : 
\dcat{D}^{\mrm{b}}_{(\mrm{f}, ..)}(A^{\mrm{en}}, \mrm{gr}) \to
\dcat{D}^{\mrm{b}}_{(\mrm{cof}, ..)}(A^{\mrm{en}}, \mrm{gr}) ,
\end{equation}
with with quasi-inverse $\opn{ADC}_{\m}$. 
From the same theorem, but switching $A$ and $A^{\mrm{op}}$, 
we get an equivalence 
\[ \mrm{R} \Ga_{\m^{\mrm{op}}} : 
\dcat{D}^{\mrm{b}}_{(.., \mrm{f})}(A^{\mrm{en}}, \mrm{gr}) \to
\dcat{D}^{\mrm{b}}_{(.., \mrm{cof})}(A^{\mrm{en}}, \mrm{gr}) , \]
with with quasi-inverse $\opn{ADC}_{\m^{\mrm{op}}}$. 
According to Theorem \ref{thm:4051} there is an isomorphism  
$\mrm{R} \Ga_{\m}(M) \cong \mrm{R} \Ga_{\m^{\mrm{op}}}(M)$
for 
$M \in \dcat{D}^{\mrm{b}}_{(\mrm{f}, \mrm{f})}(A^{\mrm{en}}, \mrm{gr})$.
Therefore, for such a complex $M$ we have 
$\mrm{R} \Ga_{\m}(M) \in 
\dcat{D}^{\mrm{b}}_{(\mrm{cof}, \mrm{cof})}(A^{\mrm{en}}, \mrm{gr})$.
Corollary \ref{cor:4080} says that 
$\opn{ADC}_{\m}(N) \cong \opn{ADC}_{\m^{\mrm{op}}}(N)$
for a complex 
$N \in \dcat{D}^{\mrm{b}}_{(\mrm{cof}, \mrm{cof})}(A^{\mrm{en}}, \mrm{gr})$; 
so for such a complex $N$ we have 
$\opn{ADC}_{\m}(N) \in 
\dcat{D}^{\mrm{b}}_{(\mrm{f}, \mrm{f})}(A^{\mrm{en}}, \mrm{gr})$.
The conclusion is that the equivalence (\ref{eqn:4845}) restricts to 
the equivalence stated in the corollary. 
\end{proof}

\begin{que} \label{que:3990}
A consequence of Theorem \ref{thm:3750} is that the complex $P_A$ has finite 
flat dimension over $A$. Is it true that the complex $P_A$ has finite 
projective dimension over $A$~? I.e., does the functor 
\[ \opn{ADC}_{\m} = \opn{RHom}_{A}(P_A, -) : \dcat{D}(A, \mrm{gr}) \to 
\dcat{D}(A, \mrm{gr}) \]
have finite cohomological dimension?

In the commutative weakly proregular case this is true -- see Examples 
\ref{exa:4595} and \ref{exa:4575}. Indeed, if $\a$ is a WPR ideal in a 
commutative ring $A$, then the cohomological dimension of the functor
$\opn{ADC}_{\a} \cong \mrm{L} \Lambda_{\a}$ 
is at most the length of a finite generating sequence of $\a$. 
\end{que}

\begin{rem} \label{rem:3857}
The proof of Theorem \ref{thm:3713} provided here is not the original proof 
from \cite{VdB}. Our proof, which relies on the NC MGM Equivalence,
seems to allow a generalization to the case of a base ring $\K$ that is not a 
field (as long as the ring $A$ is flat over $\K$). 
See \cite{VyYe} for some results in this direction. 
\end{rem}

\begin{rem} \label{rem:4082}
Q.S. Wu and J.J. Zhang \cite{WuZh}, \cite{WuZh2} studied prebalanced dualizing 
complexes over complete semilocal noncommutative rings that contain a field. 
Some of their technical results are included in Subsection 
\ref{subsec:bal-trace} below. 

An ongoing project of R. Vyas and A. Yekutieli (a continuation of \cite{VyYe}) 
aims to study balanced dualizing complexes over complete semilocal 
noncommutative rings in the arithmetic setting, namely without the presence of 
a base field. A prototypical example is the ring 
$A = \K[[G]] / I$, where $\K = \what{\Z}_p$, $G$ is a compact $p$-adic Lie 
group, $\K[[G]]$ is the noncommutative Iwasawa algebra, and 
$I \sub \K[[G]]$ is some two-sided ideal.  Because the ring $A$ could fail to 
be flat over $\K$, the balanced dualizing complex of $A$ will be an object of 
the {\em derived category of bimodules} 
$\dcat{D}(\til{A}^{\mrm{en}})$, where 
$\til{A} \to A$ is a K-flat DG ring resolution of $A$ and 
$\til{A}^{\mrm{en}} := \til{A} \ot_{\K} \til{A}^{\mrm{op}}$;
see \cite{Ye16}. 
\end{rem}

\begin{rem} \label{rem:4081}
Suppose $A$ is a commutative noetherian complete local ring, with maximal ideal 
$\m$. In the paper \cite{AlJeLi} the authors discuss two kinds of dualizing 
complexes over $A$~: the usual dualizing complexes (see Subsection 
\ref{subsec:du-cplxs}), 
that they call {\em c-dualizing complexes}, which have finite (and thus 
complete) cohomology modules; and the {\em t-dualizing complexes}, 
that have cofinite (and thus torsion) cohomology modules.
The torsion injective module $I = A^*$ (the injective hull of the residue 
field) is a t-dualizing complex in their sense. One of the results of 
\cite{AlJeLi} is that the MGM Equivalence (that they, in that early paper, 
called {\em GM Duality}) exchanges t-dualizing complexes and c-dualizing 
complexes. In particular, by applying the derived $\m$-adic completion functor 
$\mrm{L} \Lambda_{\m}$ one obtains a c-dualizing complex 
$R := \mrm{L} \Lambda_{\m}(A^*)$. 

It was observed by Vyas that the Van den Bergh Existence Theorem (Theorem 
\ref{thm:3713}) can be understood as a noncommutative variant of that very same 
result. In our noncommutative graded setting we must replace the true derived 
$\m$-adic completion functor $\mrm{L} \Lambda_{\m}$ with its abstract 
avatar $\opn{ADC}_{\m}$. The graded $A$-bimodule $A^*$ is a t-dualizing complex 
(if we adjust the definition from \cite{AlJeLi} to the NC connected graded 
setting), and indeed, as Theorem \ref{thm:3713} says, the complex
\[ R_A := \opn{ADC}_{\m}(A^*) = \opn{RHom}_{A}(P_A, A^*) = (P_A)^* \]
is a graded NC dualizing complex (as defined in Definition \ref{dfn:3716}).
\end{rem}

\mysubsection{Balanced Trace Morphisms} 
\label{subsec:bal-trace}

In this subsection we adhere to Conventions \ref{conv:3700} and \ref{conv:4560}.

\begin{dfn} \label{dfn:3765}
A graded ring homomorphism $f : A \to B$ is called {\em finite} if it makes $B$ 
into a finite graded $A$-module on both sides. 
\end{dfn}

Throughout this subsection we assume the next setup: 

\begin{setup} \label{setup:3765} 
We are given noetherian connected graded $\K$-rings $A$ and $B$, with 
augmentation ideals $\m$ and $\n$ respectively, and a finite graded
$\K$-ring homomorphism $f : A \to B$. 
We are also given a homomorphism of graded $\K$-rings  $g : C \to D$. 
\end{setup}

The rings $C$ and $D$ will play auxiliary roles, to allow us to handle 
bimodules. In practice $C$ will be 
either $A$ or $\K$, and $D$ will be either $B$ or $\K$.
The graded ring homomorphism 
$f \ot g^{\mrm{op}} : A \ot C^{\mrm{op}} \to B \ot D^{\mrm{op}}$
induces a  restriction functor
\[ \opn{Rest} = \opn{Rest}_{f \ot g^{\mrm{op}}}
 : \dcat{D}(B \ot D^{\mrm{op}}, \mrm{gr}) \to 
\dcat{D}(A \ot C^{\mrm{op}}, \mrm{gr}) . \]
This restriction functor will often remain implicit. 

\begin{thm} \label{thm:3783}
Assume Setup \tup{\ref{setup:3765}}.
\begin{enumerate}
\item There is a canonical morphism
\[ \opn{Rest} \circ \, \mrm{R} \Ga_{\n} \to 
\mrm{R} \Ga_{\m} \circ \opn{Rest} \]
of triangulated functors 
$\dcat{D}(B \ot D^{\mrm{op}}, \mrm{gr}) \to 
\dcat{D}(A \ot C^{\mrm{op}}, \mrm{gr})$.

\item Let $N \in \dcat{D}(B \ot D^{\mrm{op}}, \mrm{gr})$, and assume either 
that 
$N \in \dcat{D}^+(B \ot D^{\mrm{op}}, \mrm{gr})$ 
or that the functor $\Ga_{\m}$ has finite right cohomological dimension. 
Then the canonical morphism
$\mrm{R} \Ga_{\n}(N) \to \mrm{R} \Ga_{\m}(N)$
in $\dcat{D}(A \ot C^{\mrm{op}}, \mrm{gr})$ is an isomorphism. 
\end{enumerate}
\end{thm}

For the proof we need a few lemmas.

\begin{lem} \label{lem:3785}
For every $N \in \dcat{M}(B \ot D^{\mrm{op}}, \mrm{gr})$
there is equality 
$\Ga_{\m}(N) = \Ga_{\n}(N)$
of $(B \ot D^{\mrm{op}})$-submodules of $N$. 
\end{lem}

\begin{proof}
Take any integer $q \geq 1$,
and let 
$\m^q := \m \cdots \m$, the $q$-fold power of the ideal $\m$. 
Because $B$ is a finite $A$-module, it follows that 
$B / (\m^q \cd B) \cong (A / \m^q) \ot_A B$
is a finite module over $A / \m^q$, and hence it is a finite $\K$-module. 
This implies that the graded module $B / (\m^q \cd B)$ is concentrated in a 
finite algebraic degree interval, say $[0, i_1]$. 
So for every $i > i_1$ there is equality $(\m^q \cd B)_i = B_i$. 
Taking $q' := \opn{max}(q, i_1 + 1)$ we get 
$\n^{q'} \sub \m^q \cd B$. 

The argument above shows that for every $q \geq 1$ there exists some 
$q' \geq q$ for which $\n^{q'} \sub \m^q \cd B$. 
Trivially there is an inclusion $\m^q \cd B \sub \n^q$.
These inclusions yield 
\[ \opn{Hom}_{B} ( B / \n^{q'} , N) \sub
\opn{Hom}_{A} ( A / \m^q , N) \sub 
\opn{Hom}_{B} ( B / \n^q , N) \sub N . \]
Combined with formula (\ref{eqn:3772}) the assertion is proved. 
\end{proof}

Recall the ML condition from Definition \ref{dfn:1590}.
An inverse system $\{ N_q \}_{q \in \N}$ 
in $\dcat{M}(\K, \mrm{gr})$ is said to have the {\em trivial ML 
property}, or is called {\em prozero}, if for every $q$ there is some 
$q' \geq q$ such that the homomorphism $N_{q'} \to N_q$ is zero. See 
\cite[Definition 3.5.6]{We}. 

\begin{lem} \label{lem:3768}
Let $\{ N_q \}_{q \in \N}$ be an inverse system in 
$\dcat{M}(\K, \mrm{gr})$ that has the ML property. The following are 
equivalent\tup{:}
\begin{enumerate}
\rmitem{i} The inverse system has the trivial ML property.

\rmitem{ii} $\opn{lim}_{\lar q} N_q = 0$. 
\end{enumerate}
\end{lem}

\begin{lem} \label{lem:3820} 
Let $\{ N_q \}_{q \in \N}$ be an inverse system in 
$\dcat{M}(A, \mrm{gr})$ that has the trivial ML property, and let 
$M \in \dcat{M}(A, \mrm{gr})$. Then 
$\lim_{q \to} \opn{Hom}_{A}(N_q, M) = 0$.
\end{lem}

\begin{exer} \label{exer:3765}
Prove Lemmas \ref{lem:3768} and \ref{lem:3820}.
\end{exer}

For a finite graded $A$-module $M$, and for numbers $p, q \in \N$, we 
define \begin{equation} \label{eqn:3765}
F_{p, q}(M) := \opn{Tor}^A_p(A / \m^{q + 1}, M) = 
\opn{H}^{-p} \bigl(  (A / \m^{q + 1})  \ot^{\mrm{L}}_{A} M \bigr) 
\in \dcat{M}(A, \mrm{gr}) . 
\end{equation}
Then, fixing $p$, we define 
\begin{equation} \label{eqn:3766}
F_{p}(M) := \opn{lim}_{\lar q} F_{p, q}(M) \in \dcat{M}(A, \mrm{gr}) . 
\end{equation}
In this way we obtain functors 
\begin{equation} \label{eqn:5202}
F_{p, q}, F_{p} : \dcat{M}_{\mrm{f}}(A, \mrm{gr}) \to \dcat{M}(A, \mrm{gr}) .
\end{equation}

Lemmas \ref{lem:3765} and \ref{lem:3766} below are essentially
\cite[Lemma 6.4 and Proposition 6.5]{WuZh}. 

\begin{lem} \label{lem:3765}
Let $M$ be a finite graded $A$-module and let $p \geq 1$. Then the 
inverse system $\{ F_{p, q}(M) \}_{q \in \N}$ 
has the trivial ML property. 
\end{lem}

\begin{proof} \mbox{}

\smallskip \noindent
Step 1. Here we show that for fixed $p$, the inverse 
system $\{ F_{p, q}(M) \}_{q \in \N}$ 
satisfies the ML condition. 

Because $A$ is left noetherian and $M$ is a finite graded $A$-module, we can 
find a resolution 
$\cdots \to L^{-1} \to L^0 \to M \to 0$
in $\dcat{M}_{}(A, \mrm{gr})$,
such that each $L^{-p}$ is a finite graded-free $A$-module. 
Then 
\[ F_{p, q}(M) = 
\opn{H}^{-p} \bigl( (A / \m^{q + 1}) \ot^{\mrm{L}}_{A} M \bigr) 
\cong \opn{H}^{-p} \bigl( (A / \m^{q + 1}) \ot_{A} L \bigr) \]
is a finite $(A / \m^{q + 1})$-module. So as a graded $A$-module it is 
artinian. 
For fixed $q$, and for all $q' \geq q$, the descending chain condition on the 
submodules
$\opn{Im} \bigl( F_{p, q'}(M) \to F_{p, q}(M) \bigr) \sub F_{p, q}(M)$
tells us that this eventually stabilizes. This is the ML property. 

\medskip \noindent
Step 2. Let 
$0 \to M' \to M \to M'' \to 0$
be a short exact sequence in $\dcat{M}_{\mrm{f}}(A, \mrm{gr})$.
For each $q$ we have a long exact Tor sequence 
\begin{equation} \label{eqn:3767}
\begin{aligned}
& \cdots \to F_{1, q}(M') \to F_{1, q}(M) \to F_{1, q}(M'') 
\\
& \quad \to F_{0, q}(M') \to F_{0, q}(M) \to F_{0, q}(M'') \to 0  
\end{aligned}
\end{equation}
As $q$ varies, these long exact sequences form an inverse system. By step 1 we 
know that in every position the ML property is satisfied. Therefore by Theorem
\ref{thm:1535}, applied in each degree separately, i.e.\ to the inverse systems 
in $\dcat{M}(\K)$, the inverse limit 
\begin{equation} \label{eqn:3768}
\begin{aligned}
& \cdots \to F_{1}(M') \to F_{1}(M) \to F_{1}(M'') 
\\
& \quad \to F_{0}(M') \to F_{0}(M) \to F_{0}(M'') \to 0  
\end{aligned}
\end{equation}
is an exact sequence. 

\medskip \noindent
Step 3. For a fixed degree $k$, we have $(\m^{q + 1} \cd M)_k = 0$
when $q \gg 0$. This is because $M$ is a finite graded $A$-module, and thus it 
is bounded below in algebraic degree. It follows that the canonical homomorphism
$M_k \to (M / \m^{q + 1} \cd M)_k$ 
is bijective for $q \gg 0$. Going to the limit in $q$ we see that the canonical 
homomorphism
\[ M \to F_0(M) =  \opn{lim}_{\lar q} F_{0, q}(M) = 
\opn{lim}_{\lar q} (M / \m^{q + 1} \cd M) \]
is bijective. Therefore there is an isomorphism of functors 
$F_0 \cong \opn{Id}$. We deduce that $F_0$ is an exact functor on 
$\dcat{M}_{\mrm{f}}(A, \mrm{gr})$. 

\medskip \noindent
Step 4. We now prove, by induction on $p$, that the functor $F_p$ is zero for 
every $p \geq 1$. 

Given a module $M \in \dcat{M}_{\mrm{f}}(A, \mrm{gr})$, 
we choose a surjection $L \surj M$ from a finite graded-free module $L$. 
We get a short exact sequence
$0 \to N \to L \to M \to 0$
in $\dcat{M}_{\mrm{f}}(A, \mrm{gr})$.
Hence, by step 2, there is a long exact sequence 
\begin{equation} \label{eqn:3769}
\begin{aligned}
& \cdots \to F_{2}(N) \to F_{2}(L) \to F_{2}(M) 
\\
& \quad \xar{\de^1} F_{1}(N) \to F_{1}(L) \to F_{1}(M) 
\\
& \quad \xar{\de^0} F_{0}(N) \to F_{0}(L) \to F_{0}(M) \to 0  
\end{aligned}
\end{equation}
in $\dcat{M}(A, \mrm{gr})$;  cf.\ the exact sequence (\ref{eqn:3768}).
Now $F_{p, q}(L) = 0$ for all $q \geq 0$ and all $p > 0$, because $L$ is flat; 
so in the limit we get $F_{p}(L) = 0$ for all $p > 0$. 

Since the functor $F_0$ is exact (by step 3), we see that $\de^0 = 0$. 
Also we have $F_{1}(L) = 0$. Therefore $F_{1}(M) = 0$. 

But $M$ was an arbitrary object of 
$\dcat{M}_{\mrm{f}}(A, \mrm{gr})$, so in fact the functor $F_1$ is zero. 
This implies that in the exact sequence (\ref{eqn:3769}) 
we have $F_{1}(N) = 0$. We also have  $F_{2}(L) = 0$. 
This implies that $F_2(M) = 0$. And so on. This finishes the proof that 
$F_p = 0$ for all $p \geq 1$. 

\medskip \noindent
Step 5. Take some $M \in \dcat{M}_{\mrm{f}}(A, \mrm{gr})$.
We know that $F_p(M) = 0$ for every $p \geq 1$. By step 1 we also 
know that for every $p$ inverse system $\{ F_{p, q}(M) \}_{q \in \N}$ 
has the ML property. Lemma \ref{lem:3768} says that 
for every $p \geq 1$ the inverse system $\{ F_{p, q}(M) \}_{q \in \N}$ 
has the trivial ML property.
\end{proof}

\begin{lem} \label{lem:3766}
If $J$ is a graded-injective $B$-module, then as an $A$-module, $J$ is 
graded-$\m$-flasque. 
\end{lem}

\begin{proof}
Take some $p > 0$. We have to show that 
$\mrm{R}^p \Ga_{\m}(J) = 0$. 
Now 
\[ \mrm{R}^p \Ga_{\m}(J) 
\cong \lim_{q \to} \, \opn{Ext}^p_A \bigl( A / \m^{q + 1}, J) . \]
But 
\[ \opn{Ext}^p_A \bigl( A / \m^{q + 1}, J) = 
\opn{H}^p \bigl( \opn{RHom}_A ( A / \m^{q + 1}, J) \bigr) . \]
By derived Hom-tensor adjunction we have a canonical isomorphism 
\[ \opn{RHom}_A ( A / \m^{q + 1}, J) \cong 
\opn{RHom}_B \bigl( B \ot^{\mrm{L}}_{A} (A / \m^{q + 1}), J \bigr) . \]
in $\dcat{D}(A, \mrm{gr})$. Because $J$ is graded-injective, we have 
\[ \opn{RHom}_B \bigl( B \ot^{\mrm{L}}_{A} (A / \m^{q + 1}), J \bigr) \cong 
\opn{Hom}_B \bigl( B \ot^{\mrm{L}}_{A} (A / \m^{q + 1}), J \bigr) , \]
and also  
\[ \begin{aligned}
& \opn{H}^p \bigl( 
\opn{Hom}_B \bigl( B \ot^{\mrm{L}}_{A} (A / \m^{q + 1}), J \bigr) \cong
\opn{Hom}_B \bigl( \opn{H}^{-p} \bigl(
B \ot^{\mrm{L}}_{A} (A / \m^{q + 1}) \bigr), J \bigr) 
\\
& \quad \cong 
\opn{Hom}_B \bigl( \opn{Tor}_p^A (B, A / \m^{q + 1}), J \bigr) . 
\end{aligned} \]
Putting it all together we get 
\[ \mrm{R}^p \Ga_{\m}(J) \cong 
\lim_{q \to} \, \opn{Hom}_B \bigl( F_{p, q}(B) , J \bigr) , \]
where 
$F_{p, q}(B) := \opn{Tor}_p^A (B, A / \m^{q + 1})$, 
like in formula (\ref{eqn:3765}).
Now $B$ is a finite $A^{\mrm{op}}$-module. 
According to Lemma \ref{lem:3765}, applied to the ring $A^{\mrm{op}}$ instead 
of 
$A$, the inverse system 
$\{ F_{p, q}(M) \}_{q \in \N}$ 
has the trivial ML property. By Lemma \ref{lem:3820} we see that 
$\lim_{q \to} \opn{Hom}_B \bigl( F_{p, q}(B) , J \bigr) = 0$.
\end{proof}

\begin{proof}[Proof of Theorem \tup{\ref{thm:3783}}] \mbox{}

\smallskip \noindent
(1) For a complex 
$N \in \dcat{D}(B \ot D^{\mrm{op}}, \mrm{gr})$
we choose a K-injective resolution $\rho : N \to J$ in 
$\dcat{C}_{\mrm{str}}(B \ot D^{\mrm{op}}, \mrm{gr})$.
Then we choose a  K-injective resolution $\si : J \to I$ in 
$\dcat{C}_{\mrm{str}}(A \ot C^{\mrm{op}}, \mrm{gr})$. (We are hiding the 
restriction functor $\opn{Rest}$ here.) So 
$\si \circ \rho : N \to I$ is a K-injective resolution 
in $\dcat{C}_{\mrm{str}}(A \ot C^{\mrm{op}}, \mrm{gr})$.
We get presentations 
$\mrm{R} \Ga_{\n}(N) = \Ga_{\n}(J)$ and 
$\mrm{R} \Ga_{\m}(N) = \Ga_{\m}(I)$. 
By Lemma \ref{lem:3785} there is equality 
$\Ga_{\m}(J) = \Ga_{\n}(J)$. 
The morphism $\mrm{R} \Ga_{\n}(N) \to \mrm{R} \Ga_{\m}(N)$
in $\dcat{D}(A \ot C^{\mrm{op}}), \mrm{gr})$
we are looking for is the one represented by the homomorphism 
$\Ga_{\m}(\si) : \Ga_{\n}(J) = \Ga_{\m}(J) \to \Ga_{\m}(I)$
in $\dcat{C}_{\mrm{str}}(A \ot C^{\mrm{op}}, \mrm{gr})$. 
It does not depend on the choices made. 

\medskip \noindent
(2) Now we choose the K-injective resolution $\rho : N \to J$ in 
$\dcat{C}_{\mrm{str}}(B \ot D^{\mrm{op}}, \mrm{gr})$ more carefully: 
we take the semi-graded-cofree resolution from Theorem \ref{thm:3740}(2).
Then the complex $J$ is K-injective in $\dcat{C}(B \ot D^{\mrm{op}}, \mrm{gr})$, 
each $J^p$ is injective in $\dcat{M}(B \ot D^{\mrm{op}}, \mrm{gr})$, and 
$\opn{inf}(J) = \opn{inf}(\opn{H}(N))$. 
According to Lemma \ref{lem:3766}, $J$ is a complex of graded-$\m$-flasque 
modules. By Lemma \ref{lem:3885}, under either of the two conditions, 
the complex $J$ is graded-$\m$-flasque. Therefore the homomorphism 
$\Ga_{\m}(\si) : \Ga_{\m}(J) \to \Ga_{\m}(I)$
in $\dcat{C}_{\mrm{str}}(A \ot C^{\mrm{op}}, \mrm{gr})$ is a quasi-isomorphism. 
\end{proof}

\begin{cor} \label{cor:4030}
If the functor $\Ga_{\m}$ has finite right cohomological dimension, then 
the functor $\Ga_{\n}$ has finite right cohomological dimension.
\end{cor}

\begin{proof}
Part 2 of Theorem \ref{thm:3783} implies that the right cohomological dimension 
of $\Ga_{\n}$ is at most the right cohomological dimension of $\Ga_{\m}$.  
\end{proof}

\begin{cor} \label{cor:4031}
If the ring $A$ satisfies the left $\chi$ condition, then so does the ring $B$. 
\end{cor}

\begin{proof}
A graded $B$-module $N$ is finite iff it is finite as a graded $A$-module, and 
the same for cofiniteness. Let $N \in \dcat{M}_{\mrm{f}}(B, \mrm{gr})$. 
By Proposition \ref{prop:3855} the graded $A$-modules 
$\opn{H}^p_{\m}(N)$ are all cofinite.
But by part 2 of Theorem \ref{thm:3783} we know that 
$\opn{H}^p_{\n}(N) \cong \opn{H}^p_{\m}(N)$ as graded $A$-modules. 
So the graded $B$-modules $\opn{H}^p_{\n}(N)$ are all cofinite.
Again using Proposition \ref{prop:3855} we see that $B$ satisfies the left 
$\chi$ condition.
\end{proof}

\begin{cor} \label{cor:4032}
If the ring $A$ has a balanced dualizing complex $(R_A, \be_A)$, then the ring 
$B$ has a balanced dualizing complex $(R_B, \be_B)$.
\end{cor}

\begin{proof}
According to Corollary \ref{cor:4585}, $A$ satisfies the $\chi$ condition, and 
it has finite local cohomological dimension. Corollaries \ref{cor:4030} and 
\ref{cor:4031}, applied to the finite graded ring homomorphisms
$A \to B$ and $A^{\mrm{op}} \to B^{\mrm{op}}$, tell us that 
$B$ satisfies the $\chi$ condition, and it has finite local cohomological 
dimension. Again using Corollary \ref{cor:4585}, now in the reverse direction, 
we conclude that $B$ has a balanced dualizing complex.
\end{proof}

The graded ring homomorphism $f : A \to B$ induces, by applying the functor 
$(-)^* = \opn{Hom}_{\K}(-, \K)$, a homomorphism
$f^* : B^* \to A^*$ in $\dcat{M}(A^{\mrm{en}}, \mrm{gr})$. 

\begin{dfn} \label{dfn:4030}
Under Setup \tup{\ref{setup:3765}}, assume $A$ and $B$ have balanced dualizing 
complexes $(R_A, \be_A)$ and $(R_B, \be_B)$ respectively. A 
{\em balanced trace morphism} 
\index{Dualizing complex! balanced trace morphism of balanced {\indash}es}
\index{1-Trf@$\opn{tr}_f$}
is a morphism 
\[ \opn{tr}_f = \opn{tr}_{B / A} : R_B \to R_A \]
in $\dcat{D}(A^{\mrm{en}}, \mrm{gr})$,
such that the diagram 
\[ \UseTips \xymatrix @C=10ex @R=6ex {
\mrm{R} \Ga_{\m}(R_B)
\ar[r]^{ \mrm{R} \Ga_{\m}(\opn{tr}_f) }
\ar[d]_{\be_B}^{\cong}
&
\mrm{R} \Ga_{\m}(R_A)
\ar[d]^{\be_A}_{\cong}
\\
B^*
\ar[r]^{f^*}
&
A^*
} \]
in $\dcat{D}(A^{\mrm{en}}, \mrm{gr})$ is commutative. 
Here we use the canonical isomorphism 
$\mrm{R} \Ga_{\n}(R_B) \lb \cong \mrm{R} \Ga_{\m}(R_B)$
from Theorem \tup{\ref{thm:3783}(2)}.
\end{dfn}

In Subsection \ref{subsec:adjs-NC} we studied backward morphisms for objects in 
derived categories. This makes sense also for complexes of graded bimodules. 
Namely, given complexes
$M \in \dcat{D}(A \ot C^{\mrm{op}}, \mrm{gr})$
and
$N \in \dcat{D}(B \ot D^{\mrm{op}}, \mrm{gr})$, 
we can talk about a {\em backward morphism} 
$\th : N \to M$ in $\dcat{D}(A \ot C^{\mrm{op}}, \mrm{gr})$. 
The backward morphism $\th$ induces two morphisms by adjunction. 
First there is the morphism 
$\opn{badj}^{\mrm{R}}_{A}(\th) : N \to \opn{RHom}_{A}(B, M)$ 
in $\dcat{D}(B \ot C^{\mrm{op}}, \mrm{gr})$; this is the {\em derived backward 
adjunction morphism on the $A$ side}. The second is 
$\opn{badj}^{\mrm{R}}_{C^{\mrm{op}}}(\th) : N \to 
\opn{RHom}_{C^{\mrm{op}}}(D, M)$
in $\dcat{D}(A \ot D^{\mrm{op}}, \mrm{gr})$; this is the {\em derived backward 
adjunction morphism on the $C^{\mrm{op}}$ side}.
The constructions are just like those in Theorem \ref{thm:3192}(2), 
but now using the graded bimodule resolutions from Subsection 
\ref{gr-res-der-fun}.

\begin{dfn} \label{dfn:4031}
Consider a backward morphism $\th : N \to M$ 
in \lb $\dcat{D}(A \ot C^{\mrm{op}}, \mrm{gr})$. 
\begin{enumerate}
\item $\th$ is called a {\em nondegenerate backward morphism on the $A$ side} 
if the morphism $\opn{badj}^{\mrm{R}}_{A}(\th)$ is an isomorphism. 

\item $\th$ is called a {\em nondegenerate backward morphism on the 
$C^{\mrm{op}}$ side} if the morphism 
$\opn{badj}^{\mrm{R}}_{C^{\mrm{op}}}(\th)$ 
is an isomorphism. 

\item  $\th$ is said to be a {\em nondegenerate backward morphism on both 
sides} if it is nondegenerate both on the $A$ side and on the 
$C^{\mrm{op}}$ side. 
\end{enumerate}
\end{dfn}

\begin{exa} \label{exa:4585}
The backward morphism $f^* : B^* \to A^*$ 
in $\dcat{D}(A^{\mrm{en}}, \mrm{gr})$ is a nondegenerate backward morphism 
on both sides. 
\end{exa}

\begin{thm} \label{thm:4030}
The following hold in the situation of Definition \tup{\ref{dfn:4030}}.
\begin{enumerate}
\item There exists a unique 
balanced trace morphism 
$\opn{tr}_f : R_B \to R_A$.

\item The balanced trace morphism $\opn{tr}_f$ is a nondegenerate backward 
morphism  on both sides.
\end{enumerate}
\end{thm}

We shall need several lemmas for the proof. In analogy to Definition 
\ref{dfn:3750}, the dedualizing complex of $B$ is 
$P_B := \mrm{R} \Ga_{\n}(B) \cong \mrm{R} \Ga_{\n^{\mrm{op}}}(B)$ 
in $\dcat{D}(B^{\mrm{en}}, \mrm{gr})$.

\begin{lem} \label{lem:4035}
There are canonical isomorphisms 
$P_B \cong B \ot^{\mrm{L}}_{A} P_A$ and
$P_B \cong P_A \ot^{\mrm{L}}_{A} B$
in $\dcat{D}(B \ot A^{\mrm{op}}, \mrm{gr})$ and 
$\dcat{D}(A \ot B^{\mrm{op}}, \mrm{gr})$ respectively.
\end{lem}

\begin{proof}
According to Theorems \ref{thm:3783}(2) and \ref{thm:3750} there are 
isomorphisms
$P_B = \mrm{R} \Ga_{\n}(B) \cong \mrm{R} \Ga_{\m}(B) \cong 
P_A \ot^{\mrm{L}}_{A} B$
in $\dcat{D}(A \ot B^{\mrm{op}}, \mrm{gr})$. The same theorems, but applied to 
the opposite rings, give us isomorphisms
$P_B \cong \mrm{R} \Ga_{\n^{\mrm{op}}}(B) \cong 
\mrm{R} \Ga_{\m^{\mrm{op}}}(B) \cong B \ot^{\mrm{L}}_{A} P_A$
in $\dcat{D}(B \ot A^{\mrm{op}}, \mrm{gr})$.
\end{proof}

Like Definition \ref{dfn:3790}, the abstract derived $\n$-adic completion 
functor is 
\[ \opn{ADC}_{\n} := \opn{RHom}_{B}(P_B, -) : 
\dcat{D}(B \ot D^{\mrm{op}}, \mrm{gr}) \to 
\dcat{D}(B \ot D^{\mrm{op}}, \mrm{gr}) . \]

\begin{lem} \label{lem:4036}
For $N \in \dcat{D}(B \ot D, \mrm{gr})$
there is an isomorphism 
$\opn{ADC}_{\n}(N) \cong \opn{ADC}_{\m}(N)$
in $\dcat{D}(A \ot D^{\mrm{op}}, \mrm{gr})$.
It is functorial in $N$. 
\end{lem}

\begin{proof}
We have isomorphisms 
\[ \begin{aligned}
& \opn{ADC}_{\n}(N) = \opn{RHom}_{B}(P_B, N) \cong^1 
\opn{RHom}_{B}(B \ot^{\mrm{L}}_{A} P_A, N) 
\\
& \qquad 
\cong^2 \opn{RHom}_{A}(P_A, N) = \opn{ADC}_{\m}(N)
\end{aligned} \]
in $\dcat{D}(A \ot D^{\mrm{op}}, \mrm{gr})$. 
The isomorphism $\cong^1$ comes from Lemma \ref{lem:4035}, and the 
isomorphism $\cong^2$ is by adjunction. 
\end{proof}

In the next lemma we take $C = D := \K$. 

\begin{lem} \label{lem:4065}
Let $N \in \dcat{D}^{+}_{\mrm{f}}(B, \mrm{gr})$,
let $M \in \dcat{D}^{+}_{\mrm{f}}(A, \mrm{gr})$,
and let 
$\th : N \to M$ be a backward morphism in 
$\dcat{D}(A, \mrm{gr})$. The following conditions are 
equivalent\tup{:}
\begin{enumerate}
\rmitem{i} $\th$ is a nondegenerate backward morphism on the $A$ side.  

\rmitem{ii} The morphism 
$\mrm{R} \Ga_{\m}(\th) : \mrm{R} \Ga_{\m}(N) \to \mrm{R} \Ga_{\m}(M)$
in $\dcat{D}(A, \mrm{gr})$ is a nondegenerate backward morphism on the $A$ 
side. 
\end{enumerate}
\end{lem}

\begin{proof}
Let us choose bounded below injective resolutions 
$N \to J$ and $M \to I$ in 
$\dcat{C}_{\mrm{str}}(B, \mrm{gr})$ and 
$\dcat{C}_{\mrm{str}}(A, \mrm{gr})$ respectively. Then the backward morphism 
$\th$ is represented by a homomorphism 
$\til{\th} : J \to I$ 
in $\dcat{C}_{\mrm{str}}(A, \mrm{gr})$. 

The morphism 
$\opn{badj}^{\mrm{R}}_{A}(\th)  : N \to \opn{RHom}_{A}(B, M)$
in $\dcat{D}(B, \mrm{gr})$ is represented by the homomorphism 
$\opn{badj}_{A}(\til{\th}) : J \to \opn{Hom}_{A}(B, I)$
in $\dcat{C}_{\mrm{str}}(B, \mrm{gr})$.
The formula for it it this: 
$\opn{badj}_{A}(\til{\th})(j)(b) = \til{\th}(b \cd j)$ for $j \in J$ 
and $b \in B$. 
Now $\opn{Hom}_{A}(B, I)$ is a K-graded-injective complex over $B$, by 
adjunction. By Lemma \ref{lem:3766} the complexes $J$ and 
$\opn{Hom}_{A}(B, I)$ are graded-$\m$-flasque. So the morphism 
\[ \mrm{R} \Ga_{\m} \bigl( \opn{badj}_{A}^{\mrm{R}}(\th) \bigr) : 
\mrm{R} \Ga_{\m}(N) \to \mrm{R} \Ga_{\m} \bigl( \opn{RHom}_{A}(B, M) \bigr) \]
is represented by 
\begin{equation} \label{eqn:4850}
\Ga_{\m}(\opn{badj}_{A}(\til{\th})) : \Ga_{\m}(J) \to  
\Ga_{\m} \bigl( \opn{Hom}_{A}(B, I) \bigr) . 
\end{equation}

The morphism $\mrm{R} \Ga_{\m}(\th)$ is represented by the homomorphism
$\Ga_{\m}(\til{\th}) : \Ga_{\m}(J) \to \Ga_{\m}(I)$.
Since the complex $\Ga_{\m}(I)$ is K-graded-injective over $A$ (by 
stability, see Theorem \ref{thm:3990}), it follows that 
$\opn{badj}_A^{\mrm{R}} \bigl( \mrm{R} \Ga_{\m}(\th) \bigr)$
is represented by 
\begin{equation} \label{eqn:4851}
\opn{badj}_A(\Ga_{\m}(\til{\th})) : \Ga_{\m}(J) \to
\opn{Hom}_A(B, \Ga_{\m}(I)) . 
\end{equation}

A calculation, using the fact that $B$ is a finite $A$-module, shows that 
\begin{equation} \label{eqn:4853}
 \Ga_{\m} \bigl( \opn{Hom}_{A}(B, I) \bigr) = 
\opn{Hom}_{A} \bigl( B, \Ga_{\m}(I) \bigr) 
\end{equation}
as subcomplexes of $\opn{Hom}_{A}(B, I)$. 
Let us consider the diagram 
\begin{equation} \label{eqn:4854}
\UseTips \xymatrix @C=11ex @R=6ex {
\Ga_{\m}(J)
\ar[r]^(0.38){ \opn{badj}_A(\Ga_{\m}(\til{\th})) }
\ar[dr]_(0.45){ \Ga_{\m}(\opn{badj}_A(\til{\th})) \ \ }
&
\opn{Hom}_{A} \bigl( B, \Ga_{\m}(I) \bigr)
\ar[d]^{=}
\ar[dr]^{\sub} 
\\
&
\Ga_{\m} \bigl( \opn{Hom}_{A}(B, I) \bigr)
\ar[r]^{\sub} 
&
\opn{Hom}_{A}(B, I)
}
\end{equation}
in $\dcat{C}_{\mrm{str}}(A, \mrm{gr})$. It is commutative: the triangle on the 
right is just the equality (\ref{eqn:4853}). And an element 
$j \in \Ga_{\m}(J)$ goes by all paths to the homomorphism 
$b \mapsto \til{\th}(b \cd j)$ in $\opn{Hom}_{A}(B, I)$. 

Passing to the derived category, (\ref{eqn:4854}) gives the 
commutative diagram
\begin{equation} \label{eqn:4855}
\UseTips \xymatrix @C=14ex @R=6ex {
\mrm{R} \Ga_{\m}(N)
\ar[r]^(0.4){ \opn{badj}_A^{\mrm{R}} (\mrm{R} \Ga_{\m}(\th)) }
\ar[dr]_(0.45){ \mrm{R} \Ga_{\m} (\opn{badj}_A^{\mrm{R}}(\th)) \ \ }
&
\opn{RHom}_{A} \bigl( B, \mrm{R} \Ga_{\m}(M) \bigr)
\ar[d]^{\cong}
\\
&
\mrm{R} \Ga_{\m} \bigl( \opn{RHom}_{A}(B, M) \bigr) 
}
\end{equation}
in $\dcat{D}(A, \mrm{gr})$.

Since both $N$ and $\opn{RHom}_{A}(B, M)$ belong to 
$\dcat{D}^+_{\mrm{f}}(A, \mrm{gr})$, by Lemma \ref{lem:4060} they belong to 
$\dcat{D}^{+}(A, \mrm{gr})_{\mrm{com}}$.
By Theorem \ref{thm:3792} (the MGM Equivalence) we know that the morphism
$\opn{badj}_A^{\mrm{R}}(\th)$ is an isomorphism if and only if 
the morphism \lb $\mrm{R} \Ga_{\m}(\opn{badj}_A^{\mrm{R}}(\th))$ is an 
isomorphism. Diagram (\ref{eqn:4855}) says that this happens if and only if 
$\opn{badj}_A^{\mrm{R}} \bigl( \mrm{R} \Ga_{\m}(\th) \bigr)$
is an isomorphism. Going back to the definition, we conclude that $\th$ is 
nondegenerate if and only if $\mrm{R} \Ga_{\m}(\th)$ is nondegenerate. 
\end{proof}

\begin{proof}[Proof of Theorem \tup{\ref{thm:4030}}] \mbox{}

\smallskip \noindent
(1) According to Theorems \ref{thm:3713} and \ref{thm:3710} there is an 
isomorphism 
$R_A \cong (P_A)^* \cong \opn{ADC}_{\m}(A^*)$
in $\dcat{D}(A^{\mrm{en}}, \mrm{gr})$. 
The same theorems, together with Lemma \ref{lem:4036}
(taking $C := A$ and $D := B$), say that 
$R_B \cong (P_B)^* \cong \opn{ADC}_{\n}(B^*) \cong \opn{ADC}_{\m}(B^*)$
in $\dcat{D}(A^{\mrm{en}}, \mrm{gr})$. 

The equivalence in Corollary \ref{cor:4070} says that there are 
isomorphisms 
\begin{equation} \label{eqn:4067}
\mrm{R} \Ga_{\m}(R_B) \cong (\mrm{R} \Ga_{\m} \circ \opn{ADC}_{\m})(B^*)
\cong B^*
\end{equation}
and 
\begin{equation} \label{eqn:4068}
\mrm{R} \Ga_{\m}(R_A) \cong (\mrm{R} \Ga_{\m} \circ \opn{ADC}_{\m})(A^*)
\cong A^*
\end{equation}
in $\dcat{D}(A^{\mrm{en}}, \mrm{gr})$.
And that there is a unique morphism 
$\th : R_B \to R_A$
in $\dcat{D}(A^{\mrm{en}}, \mrm{gr})$ such that the diagram 
\begin{equation} \label{eqn:4066}
\UseTips \xymatrix @C=8ex @R=6ex {
\mrm{R} \Ga_{\m}(R_B)
\ar[r]^{\mrm{R} \Ga_{\m}(\th)}
\ar[d]_{\cong}
&
\mrm{R} \Ga_{\m}(R_A)
\ar[d]^{\cong}
\\
B^*
\ar[r]^(0.5){f^*}
&
A^*
} 
\end{equation}
in $\dcat{D}(A^{\mrm{en}}, \mrm{gr})$, in which the vertical isomorphisms are 
(\ref{eqn:4067}) and (\ref{eqn:4068}), is commutative. 

But to get a balanced trace morphism $\opn{tr}_f$ we need a commutative diagram 
like (\ref{eqn:4066}) in which the vertical isomorphisms are $\be_B$ and 
$\be_A$ respectively. Because the automorphisms of $A^*$ and $B^*$ in 
$\dcat{D}(A^{\mrm{en}}, \mrm{gr})$
are multiplication by nonzero elements of $\K$, there is a unique 
$c \in \K^{\times}$ such that 
$\opn{tr}_f := c \cd \th : R_B \to R_A$
is a balanced trace morphism.

\medskip \noindent 
(2) To see that $\opn{tr}_f$ is a nondegenerate backward morphism on the 
$A$ side we can forget the $A^{\mrm{op}}$-structure.
So we can view $\opn{tr}_f$ as a backward morphism
$\opn{tr}_f : R_B \to R_A$ in $\dcat{D}(A, \mrm{gr})$.
Now
$\mrm{R} \Ga_{\m}(\opn{tr}_f) = \be_{A}^{-1} \circ f^* \circ \be_B$.
In Example \ref{exa:4585} we saw that $f^* : B^* \to A^*$ is a nondegenerate 
backward morphism in $\dcat{D}(A, \mrm{gr})$. The morphisms $\be_A$ and $\be_B$ 
are isomorphisms. Therefore 
\[ \mrm{R} \Ga_{\m}(\opn{tr}_f) : \mrm{R} \Ga_{\m}(R_B) \to 
\mrm{R} \Ga_{\m}(R_A) \]
is a nondegenerate backward morphism in $\dcat{D}(A, \mrm{gr})$. 
According to Lemma \ref{lem:4065} the morphism $\opn{tr}_f$
is a nondegenerate backward morphism in $\dcat{D}(A, \mrm{gr})$. 

As for the $A^{\mrm{op}}$ side: the symmetry from Corollary \ref{cor:4080} 
gives isomorphisms 
$R_A \cong \opn{ADC}_{\m^{\mrm{op}}}(A^*)$
and
$R_B \cong \opn{ADC}_{\m^{\mrm{op}}}(B^*)$
in $\dcat{D}(A^{\mrm{en}}, \mrm{gr})$. 
Going over the constructions in the proof of item (1) above, we see that
there is equality 
$\mrm{R} \Ga_{\m^{\mrm{op}}}(\opn{tr}_f) = 
d \cd \be_{A}^{-1} \circ f^* \circ \be_B$
for some $d \in \K^{\times}$.
Therefore $\mrm{R} \Ga_{\m^{\mrm{op}}}(\opn{tr}_f)$ 
is a nondegenerate backward morphism 
$\dcat{D}(A^{\mrm{op}}, \mrm{gr})$.
Now Lemma \ref{lem:4065}, transcribed to the ring $A^{\mrm{op}}$, 
says that $\opn{tr}_f$ is a 
nondegenerate backward morphism in
$\dcat{D}(A^{\mrm{op}}, \mrm{gr})$.
\end{proof}

\begin{cor} \label{cor:4065}
Let $A \xar{f} B \xar{g} C$ be finite homomorphisms in 
$\catt{Rng}_{\mrm{gr}} \centover \K$, 
and assume these rings have balanced dualizing 
complexes $(R_A, \be_A)$, $(R_B, \be_B)$ and $(R_C, \be_C)$
respectively. 
Then 
$\opn{tr}_{g \circ f} = \opn{tr}_f \circ \opn{tr}_g$
as morphisms $R_C \to R_A$ in $\dcat{D}(A^{\mrm{en}}, \mrm{gr})$. 
\end{cor}

\begin{proof}
The follows from the uniqueness in Theorem \ref{thm:4030} and the fact that 
$(g \circ f)^* = f^* \circ g^*$ as homomorphisms $C^* \to A^*$ in 
$\dcat{M}(A, \mrm{gr})$. 
\end{proof}

We end this subsection with a partial converse to Corollaries 
\ref{cor:4030} and \ref{cor:4031}.
Recall that a central element $a \in A$ is called {\em regular} if it is not a 
zero divisor, i.e.\ the only element $b \in A$ such that $a \cd b = 0$ is 
$b = 0$. 

\begin{thm} \label{thm:4250}
\index{Algebraically graded ring! chi@$\chi$ condition on connected}
\index{Algebraically graded ring! connected {\indash} of finite local 
cohomological dimension}
Let $A$ be a noetherian connected graded ring, and let $a \in A$ be a regular
homogeneous central element of positive degree. Define the connected graded 
ring $B := A / (a)$. If $B$ satisfies the $\chi$ condition and it has finite 
local cohomological dimension, then $A$ also  satisfies the $\chi$ condition 
and it has finite local cohomological dimension.
\end{thm}

This is \cite[Theorem 8.8]{ArZh}. The original proof does not use 
derived categories, and is quite different from the proof below. 

\begin{proof}
By op-symmetry it suffices to prove that $A$ satisfies the left $\chi$ 
condition and that the functor $\mrm{R} \Ga_{\m}$ has finite cohomological 
dimension. 

Let $d \in \N$ be the local cohomological dimension of $B$, and let 
$i \geq 1$ be the degree of the element $a$. 
In view of Propositions \ref{prop:3855} and \ref{prop:4050}, all we need to 
prove is that the following two conditions hold for every 
$M \in \dcat{M}_{\mrm{f}}(A, \mrm{gr})$~:
\begin{itemize}
\rmitem{i} $\opn{H}^p_{\m}(M)$ is a cofinite $A$-module for every $p$.

\rmitem{ii} $\opn{H}^p_{\m}(M) = 0$ for every $p > d + 1$. 
\end{itemize}

Let us denote the augmentation ideal of $B$ by $\n$. 
The short exact sequence 
$0 \to A \xar{a \cd (-)} A(i) \to B(i) \to 0$
in $\dcat{M}(A^{\mrm{en}}, \mrm{gr})$ is viewed as a distinguished triangle
$A \xar{a \cd (-)} A(i) \to B(i) \xar{ \, \triangle \, }$
in $\dcat{D}(A^{\mrm{en}}, \mrm{gr})$. 
Applying the functor $(-) \ot^{\mrm{L}}_{A} M$ to it, we obtain a 
distinguished triangle 
\begin{equation} \label{4251}
M \xar{a \cd (-)} M(i) \to N \xar{ \, \triangle \, }
\end{equation}
in $\dcat{D}(A, \mrm{gr})$, where by definition 
$N := B(i) \ot^{\mrm{L}}_{A} M \in \dcat{D}(B, \mrm{gr})$.
Note that 
$\opn{H}^{p}(N) \in \dcat{M}_{\mrm{f}}(B, \mrm{gr})$
for all $p$, and $\opn{H}^{p}(N) = 0$ unless $-1 \leq p \leq 0$. 

Next we apply the functor $\mrm{R} \Ga_{\m}$ to (\ref{4251}), and use Theorem 
\ref{thm:3783}, to obtain this distinguished triangle 
\begin{equation} \label{4252}
\mrm{R} \Ga_{\m}(M) \xar{a \cd (-)} \mrm{R} \Ga_{\m}(M(i)) \to 
\mrm{R} \Ga_{\n}(N) \xar{ \, \triangle \, }
\end{equation}
in $\dcat{D}(A, \mrm{gr})$.
Taking cohomologies in (\ref{4252}) we obtain this exact sequence 
\begin{equation} \label{4253}
\opn{H}^{p - 1}_{\n}(N) \to 
\opn{H}^{p}_{\m}(M) \xar{a \cd (-)}
\opn{H}^{p}_{\m}(M(i)) \to  
\opn{H}^{p}_{\n}(N) 
\end{equation}
in $\dcat{M}(A, \mrm{gr})$. 
Let use write 
\[ K^p := \opn{Ker} \bigl( \opn{H}^{p}_{\m}(M) \xar{a \cd (-)}
\opn{H}^{p}_{\m}(M(i)) \bigr) \sub \opn{H}^{p}_{\m}(M) . \]
Since $a \in \m$, there is an inclusion 
\begin{equation} \label{4254}
\opn{Soc} \bigl( \opn{H}^{p}_{\m}(M) \bigr) =
\opn{Hom}_{A} \bigl( \K, \opn{H}^{p}_{\m}(M) \bigr) \sub K^p . 
\end{equation}

For $p > d + 1$ we have 
$\opn{H}^{p - 1}_{\n}(N) = 0$ and 
$\opn{H}^{p}_{\n}(N) = 0$. 
The exact sequence (\ref{4253}) says that the homomorphism 
$a \cd (-) : \opn{H}^{p}_{\m}(M) \to \opn{H}^{p}_{\m}(M(i))$
is bijective. So $K^p = 0$, and by (\ref{4254}) we see that 
$\opn{Soc} \bigl( \opn{H}^{p}_{\m}(M) \bigr) = 0$. 
But $\opn{H}^{p}_{\m}(M)$ is an $\m$-torsion module, so according to 
Proposition \ref{prop:4811} we get $\opn{H}^{p}_{\m}(M) = 0$.
This establishes condition (ii). 

Finally take any $p$.
By the $\chi$ condition for $B$ and by Proposition \ref{prop:3855} we know 
that $\opn{H}^{p - 1}_{\n}(N)$ is a cofinite graded $B$-module; and thus it is 
a cofinite graded $A$-module. According to the NC Graded Matlis Duality 
(Theorem 
\ref{thm:4045}) the cofinite graded $A$-modules are the artinian objects in the 
abelian category $\dcat{M}(A, \mrm{gr})$.
From the exact sequence (\ref{4253}) we obtain a surjection 
$\opn{H}^{p - 1}_{\n}(N) \to K^p$.
Therefore $K^p$ is a cofinite graded $A$-module. 
The inclusion (\ref{4254}) says that 
$\opn{Soc} \bigl( \opn{H}^{p}_{\m}(M) \bigr)$ 
is a cofinite graded $A$-module; and thus it is a finite graded $\K$-module. 
Hence, by Corollary \ref{cor:4270}, 
$\opn{H}^{p}_{\m}(M)$ is a cofinite graded $A$-module.
This is condition (i). 
\end{proof}

\cleardoublepage
\mysection{Rigid Noncommutative Dualizing Complexes} 
\label{sec:rigid-DC-NC} 

\AYcopyright 

In this section of the book we are going to work with noncommutative (namely 
not necessarily commutative) rings. We shall often use the 
abbreviation ``NC'' for ``noncommutative''. 
The goal is to introduce {\em rigid NC dualizing complexes} in the sense of M. 
Van den Bergh \cite{VdB}, and to prove their existence and uniqueness under 
certain conditions.

\mysubsection{Noncommutative Dualizing Complexes}
\label{subsec:NC-DC}

Here we will define NC dualizing complexes over a NC ring $A$. 
The definition is almost identical to that in the NC graded setting (see 
Subsection \ref{subsec:graded-NCDC}). We will prove 
that the NC dualizing complexes over $A$ are 
parameterized by the derived Picard group of $A$. The content of this 
subsection is adapted from the paper \cite{Ye4}.

Let $\K$ be a commutative ring. 
Recall (from Definition \ref{dfn:4125}) that a {\em central $\K$-ring} is a 
ring $A$ equipped with a ring homomorphism $\K \to \opn{Cent}(A)$, where 
$\opn{Cent}(A)$ is the center of $A$. 
In more traditional texts, a central $\K$-ring is called an {\em associative 
unital $\K$-algebra}. The category of central $\K$-rings, with $\K$-ring 
homomorphisms, is denoted by $\catt{Rng} \centover \K$. 

Here is the convention that is in force throughout this section:

\begin{conv} \label{conv:3670} 
There is a  base field $\K$. All rings are by default central 
$\K$-rings, and all homomorphisms between $\K$-rings are over $\K$; namely we 
work within the category $\catt{Rng} \centover \K$.
All bimodules are $\K$-central, and all additive functors are $\K$-linear.
We use the notation $\ot$ for $\ot_{\K}$. For a ring $A$ we write 
$A^{\mrm{en}} := A \ot A^{\mrm{op}}$, the enveloping ring of $A$. 
\end{conv}

See Remark \ref{rem:3691} regarding the assumption that the base ring 
$\K$ is a field.    

Recall that a NC ring $A$ is called noetherian if it
is both left noetherian and right noetherian.
See Remark \ref{rem:3901} regarding the failure of the noetherian property to 
be preserved under finitely generated ring extensions. 

Let $A$ be a ring. As before, $\dcat{M}(A)$ is the category of left 
$A$-modules. The category of right $A$-modules is $\dcat{M}(A^{\mrm{op}})$,
and the category of $A$-bimodules is $\dcat{M}(A^{\mrm{en}})$.
Given another ring $B$, the category of $A$-$B$-bimodules is 
$\dcat{M}(A \ot B^{\mrm{op}})$.
These are $\K$-linear abelian categories. 

The derived category of $\dcat{M}(A)$ is 
$\dcat{D}(A) := \dcat{D}(\dcat{M}(A))$; its objects are the 
complexes of (left) $A$-modules. 
The other derived categories are  
$\dcat{D}(A^{\mrm{op}}) := \dcat{D}(\dcat{M}(A^{\mrm{op}}))$,
$\dcat{D}(A^{\mrm{en}}) := \dcat{D}(\dcat{M}(A^{\mrm{en}}))$
and 
$\dcat{D}(A \ot B^{\mrm{op}}) := \dcat{D}(\dcat{M}(A \ot B^{\mrm{op}}))$.

Consider the canonical ring anti-automorphism 
$\opn{op} : B \to B^{\mrm{op}}$, that is the identity on the underlying 
$\K$-module. There is a canonical ring isomorphism 
\begin{equation} \label{eqn:4358}
A \ot B^{\mrm{op}} \iso B^{\mrm{op}} \ot A ,  \quad  
a \ot \opn{op}(b) \mapsto \opn{op}(b) \ot a .
\end{equation}
The induced isomorphism on the bimodule categories 
$\dcat{M}(A \ot B^{\mrm{op}}) \iso \lb \dcat{M}(B^{\mrm{op}} \ot A)$
is the identity on the underlying $\K$-modules. Indeed, for 
$M \in \dcat{M}(A \ot B^{\mrm{op}})$, $m \in M$ and 
$a \ot \opn{op}(b) \in A \ot B^{\mrm{op}}$, we have 
$(\opn{op}(b) \ot a) \cd m = a \cd m \cd b = (a \ot \opn{op}(b)) \cd m$.
Now take $B = A$ in (\ref{eqn:4358}). Since 
$(A^{\mrm{op}})^{\mrm{op}} = A$, there is a canonical ring isomorphism 
\begin{equation} \label{eqn:4450}
A^{\mrm{en}} = A \ot A^{\mrm{op}} \iso 
A^{\mrm{op}} \ot A = (A^{\mrm{op}})^{\mrm{en}} 
\end{equation}
and an induced isomorphism on the bimodule categories 
\begin{equation} \label{eqn:4355}
\dcat{M}(A^{\mrm{en}}) \iso \dcat{M}((A^{\mrm{op}})^{\mrm{en}}) .
\end{equation}
All these categorical relations pass to derived categories. 

The ring $A^{\mrm{en}}$ is isomorphic to its opposite ring. 
There is a ring isomorphism 
\begin{equation} \label{eqn:4875}
\opn{trop} : A^{\mrm{en}} \iso (A^{\mrm{en}})^{\mrm{op}}
\end{equation}
with this formula: for $a_1, a_2 \in A$ it sends the element 
$a_1 \ot \opn{op}(a_2) \in A^{\mrm{en}}$ to 
\begin{equation} \label{eqn:4876}
\opn{trop}(a_1 \ot \opn{op}(a_2)) := 
\opn{op}(a_2 \ot \opn{op}(a_1)) \in (A^{\mrm{en}})^{\mrm{op}} .
\end{equation}
On the underlying $\K$-module $A \ot A$ this is just the transposition 
$a_1 \ot a_2 \mapsto a_2 \ot a_1$, and hence the name ``transpose opposite'' 
and the acronym ``trop''. 

Usually an $A$-$A$-bimodule $M$ is seen as a left 
$A^{\mrm{en}}$-module, with this action of an element 
$a_1 \ot \opn{op}(a_2) \in A^{\mrm{en}}$ on an element $m \in M$~:
\begin{equation} \label{eqn:4887}
(a_1 \ot \opn{op}(a_2)) \cd m = 
a_1 \cd m \cd a_2 .
\end{equation}
But using the ring isomorphism $\opn{trop}$ from formula (\ref{eqn:4875}) we 
can also consider $M$ as a left $(A^{\mrm{en}})^{\mrm{op}}$-module, and thus as 
a right $A^{\mrm{en}}$-module. Explicitly, the right action of an element 
$a_2 \ot \opn{op}(a_1) \in A^{\mrm{en}}$ on an element $m \in M$ is:
\begin{equation} \label{eqn:4878}
\begin{aligned}
&
m \cd (a_2 \ot \opn{op}(a_1)) =^{\heartsuit}
\opn{op}(a_2 \ot \opn{op}(a_1)) \cd m
\\ & \quad 
=^{\dag} (a_1 \ot \opn{op}(a_2)) \cd m =^{\diamondsuit} a_1 \cd m \cd a_2 . 
\end{aligned}
\end{equation}
In equality $=^{\heartsuit}$ we move from a right action by $A^{\mrm{en}}$
to the left action by $(A^{\mrm{en}})^{\mrm{op}}$, and in 
equality $=^{\dag}$ we made use of formula (\ref{eqn:4876})
to move from $(A^{\mrm{en}})^{\mrm{op}}$ to $A^{\mrm{en}}$.
Finally, the equality $=^{\diamondsuit}$ is formula (\ref{eqn:4887}). 
This game of musical chairs will be played in Subsections 
\ref{subsec:Fil-exit-R-NCDC} and \ref{subsec:RNCDC-uniq}.

Let $A$ and $B$ be rings. The canonical ring homomorphisms 
\[ A \to A \ot B^{\mrm{op}} \lar B^{\mrm{op}} \]
induce restriction functors 
\begin{equation} \label{eqn:3680}
\dcat{D}(A) \xleftarrow{\opn{Rest}_A} \dcat{D}(A \ot B^{\mrm{op}}) 
\xar{\opn{Rest}_{B^{\mrm{op}}}} \dcat{D}(B^{\mrm{op}}) . 
\end{equation}
See the commutative diagram \ref{eqn:3536}. We shall usually keep 
these restriction functors implicit, and instead use terminology like in 
Definition \ref{dfn:3446}.

In Definition \ref{dfn:2135} we introduced this notation: for a ring $A$ 
we denote by $\dcat{M}_{\mrm{f}}(A)$ the full subcategory of 
$\dcat{M}(A)$ on the finite (i.e.\ finitely generated) $A$-modules.
And $\dcat{D}_{\mrm{f}}(A)$ is the full subcategory of 
$\dcat{D}(A)$ on the complexes of $A$-modules $M$ such that 
$\opn{H}^i(M) \in \dcat{M}_{\mrm{f}}(A)$ for every $i$. 

If $A$ is left noetherian, then $\dcat{M}_{\mrm{f}}(A)$ is a thick abelian 
subcategory of $\dcat{M}(A)$; and hence $\dcat{D}_{\mrm{f}}(A)$ is a full 
triangulated subcategory of $\dcat{D}(A)$. 
As usual we can combine indicators: 
$\dcat{D}_{\mrm{f}}^{\star}(A) := \dcat{D}_{\mrm{f}}(A) \cap 
\dcat{D}^{\star}(A)$,
where $\star$ is some boundedness indicator ($+$, $-$, $\mrm{b}$ or
$\bra{\mrm{empty}}$). 

Now to finiteness of $A$-$B$-bimodules. The expression  
$\dcat{M}_{(\mrm{f}, ..)}(A \ot B^{\mrm{op}})$
denotes the full subcategory of $\dcat{M}(A \ot B^{\mrm{op}})$
on the bimodules that are finite over $A$; and the expression 
$\dcat{M}_{(.., \mrm{f})}(A \ot B^{\mrm{op}})$
denotes the full subcategory on the bimodules that are finite over 
$B^{\mrm{op}}$. Of course we define 
\[ \dcat{M}_{(\mrm{f}, \mrm{f})}(A \ot B^{\mrm{op}}) :=
\dcat{M}_{(\mrm{f}, ..)}(A \ot B^{\mrm{op}}) \cap 
\dcat{M}_{(.., \mrm{f})}(A \ot B^{\mrm{op}}) . \]
The same notation shall apply to derived categories. 

To a complex of bimodules 
$M \in \dcat{D}(A \ot B^{\mrm{op}})$ 
we associated two derived homothety morphisms in Subsection 
\ref{subsec:tilting}. Let us recall them. 
The derived homothety morphism through $B^{\mrm{op}}$ is the morphism 
\begin{equation} \label{equ:4336}
\opn{hm}^{\mrm{R}}_{M, B^{\mrm{op}}} : B \to \opn{RHom}_{A}(M, M)  
\end{equation}
in $\dcat{D}(B^{\mrm{en}})$; and the derived homothety morphism through $A$ 
is the morphism 
\begin{equation} \label{equ:4337}
\opn{hm}^{\mrm{R}}_{M, A} : A \to \opn{RHom}_{B^{\mrm{op}}}(M, M) 
\end{equation}
in $\dcat{D}(A^{\mrm{en}})$. 
According to Definition \ref{3455}, the complex $M$ has the  
noncommutative derived Morita property on the $A$ side (resp.\ on the 
$B^{\mrm{op}}$ side) if the morphism 
$\opn{hm}^{\mrm{R}}_{M, B^{\mrm{op}}}$
(resp.\ $\opn{hm}^{\mrm{R}}_{M, A}$) is an isomorphism. 

\begin{dfn}[\cite{Ye4}] \label{dfn:3675} 
Let $A$ be a noetherian ring. A 
{\em noncommutative dualizing complex over $A$} 
\index{Dualizing complex! noncommutative}
is an object 
$R \in \dcat{D}^{\mrm{b}}(A^{\mrm{en}})$ satisfying the 
following three conditions:
\begin{itemize}
\rmitem{i} For every $p$ the 
$A$-bimodule $\opn{H}^p(R)$ is finite over $A$ and over $A^{\mrm{op}}$. 

\rmitem{ii} The complex $R$ has finite injective 
dimension over $A$ and over $A^{\mrm{op}}$. 

\rmitem{iii} The complex of bimodules $R$ has the noncommutative derived Morita 
property on both sides. 
\end{itemize}
\end{dfn}

Note that conditions (i) and (ii) are one-sided, namely they refer separately 
to $\opn{Rest}_{A}(R) \in \dcat{D}(A)$ and 
to  $\opn{Rest}_{A^{\mrm{op}}}(R) \in \dcat{D}(A^{\mrm{op}})$; 
whereas condition (iii) is more complicated, and we will study it further in 
Lemma \ref{lem:4361} below. Also note that condition (i) can be stated as ``$R$ 
belongs to $\dcat{D}_{(\mrm{f}, \mrm{f})}(A^{\mrm{en}})$''. 

Before going on, we need to say something about the op-symmetry of  
Definition \ref{dfn:3675}. Given a complex 
$R \in \dcat{D}(A^{\mrm{en}})$, by formula (\ref{eqn:4355}) the complex $R$ is 
also an object of $\dcat{D}((A^{\mrm{op}})^{\mrm{en}})$, and we could ask 
whether $R$ is a dualizing complex over the ring $A^{\mrm{op}}$. The answer is 
positive of course, because all three conditions are op-symmetric. We record 
this fact as the next proposition. 

\begin{prop} \label{prop:4355}
If $R$ is a NC dualizing complex over $A$, then $R$ is also a 
NC dualizing complex over $A^{\mrm{op}}$. 
\end{prop}

\begin{exa} \label{exa:3695} 
A NC ring $A$ is called {\em regular} if it has finite global cohomological 
dimension on both sides. 
This is very similar to the graded definition (see Definition \ref{dfn:3920}).
Let us recall what it means: there is some 
natural number $d$ such that $\opn{Ext}^i_A(M, N) = 0$ for all $i > d$ and 
all $M, N \in \dcat{M}(A)$; and also  
$\opn{Ext}^i_{A^{\mrm{op}}}(M', N') = 0$ for all $i > d$ and 
all $M', N' \in \dcat{M}(A^{\mrm{op}})$. 
The smallest such number $d$ is the global cohomological dimension of $A$. 

Such rings are easy to find; for instance, if $\K$ is a field, and $\g$ is 
finite Lie algebra over $\K$, then the {\em universal enveloping ring}
$A := \opn{U}(\g)$ is regular; the global cohomological dimension of $A$ is 
$d := \opn{rank}_{\K}(\g)$.  

A weaker condition is this: $A$ is called {\em Gorenstein} if $A$ has finite 
injective dimension both as a left module and as a right module over itself. 
Namely there is some natural number $d$ such that 
$\opn{Ext}^i_A(M, A) = 0$ for all $i > d$ and all $M \in \dcat{M}(A)$; and 
also $\opn{Ext}^i_{A^{\mrm{op}}}(M', A) = 0$ for all $i > d$ and 
all $M' \in \dcat{M}(A^{\mrm{op}})$. 

Assume $A$ is Gorenstein. Then the complex of bimodules 
$R := A$ satisfies condition (ii) of Definition \ref{dfn:3675}. 
Conditions (i) and (iii) of this definition are automatically true. We conclude 
that $R := A$ is a NC dualizing complex over itself. 
\end{exa}

\begin{dfn} \label{dfn:3680}
Let $A$ be a noetherian ring, let $R$ be a NC dualizing complex over $A$, and 
let $B$ be a second ring. 
The {\em duality functors} 
\index{Duality functor}
associated to $R$ are the triangulated functors 
\[ D_A : \dcat{D}(A \ot B^{\mrm{op}})^{\mrm{op}} \to 
\dcat{D}(B \ot A^{\mrm{op}}), \quad 
D_A := \opn{RHom}_A(-, R)   \]
and 
\[ D_{A^{\mrm{op}}} : \dcat{D}(B \ot A^{\mrm{op}})^{\mrm{op}} 
\to \dcat{D}(A \ot B^{\mrm{op}}), \quad 
D_{A^{\mrm{op}}} := \opn{RHom}_{A^{\mrm{op}}}(-, R) .  \]
\end{dfn}

Notice that the expressions $D_A$ and $D_{A^{\mrm{op}}}$ leave $R$ and $B$ 
implicit; this is because the dualizing complex $R$ is taken to be fixed, and 
the second ring $B$ is less important. 

The NC derived Hom-evaluation morphisms were already discussed in Subsection 
\ref{subsec:graded-NCDC}, in the algebraically graded context. The nongraded 
definition is very similar. Given a complex 
$M \in \dcat{D}(A \ot B^{\mrm{op}})$, there is a morphism 
\begin{equation} \label{equ:4335}
\opn{ev}^{\mrm{R, R}}_{M} : M \to D_{A^{\mrm{op}}}(D_{A}(M)) = 
\opn{RHom}_{A^{\mrm{op}}} \bigl( \opn{RHom}_{A}(M, R), R \bigr)
\end{equation}
in $\dcat{D}(A \ot B^{\mrm{op}})$, which is functorial in $M$. If we choose a 
K-injective resolution $R \to I$ in $\dcat{D}(A^{\mrm{en}})$, then we have an 
isomorphism 
\[ D_{A^{\mrm{op}}}(D_{A}(M)) \cong 
\opn{Hom}_{A^\mrm{op}} \bigl( \opn{Hom}_{A}(M, I), I \bigr) \]
in $\dcat{D}(A \ot B^{\mrm{op}})$, and the morphism 
$\opn{ev}^{\mrm{R, R}}_{M}$ is represented by the Hom-evaluation homomorphism
\[ \opn{ev}_{M, I} : M \to 
\opn{Hom}_{A^\mrm{op}} \bigl( \opn{Hom}_{A}(M, I), I \bigr)  \]
in $\dcat{C}_{\mrm{str}}(A \ot B^{\mrm{op}})$. 
For this we use the results of Subsection \ref{subsec:flat-dg-rng}, that apply 
because $A$ and $B$ are flat over $\K$. 

In the special case when $B = A$ and $M = A$, the morphism 
$\opn{ev}^{\mrm{R, R}}_{A}$ recovers the derived 
homothety morphism through $A$, see (\ref{equ:4337}); i.e.\ 
\begin{equation} \label{eqn:4356}
\opn{hm}^{\mrm{R}}_{R, A} = \opn{ev}^{\mrm{R, R}}_{A} : 
A \to \opn{RHom}_{A^{\mrm{op}}}(R, R)
\end{equation}
in $\dcat{D}(A^{\mrm{en}})$. 
This formula makes implicit use of the left co-unitor isomorphism
$\opn{lcu} : \opn{RHom}_{A}(A, R) \iso R$ 
in $\dcat{D}(A^{\mrm{en}})$. 
Likewise for the homothety morphism $\opn{hm}^{\mrm{R}}_{R, A^{\mrm{op}}}$.

\begin{thm} \label{thm:3680}
Let $A$ be a noetherian ring, let $R$ be a NC dualizing complex over $A$, let 
$B$ be a second ring, and let $D_A$ and $D_{A^{\mrm{op}}}$ be the  
duality functors associated to $R$. Let $\star$ be a boundedness indicator, and 
let $M \in \dcat{D}^{\star}_{(\mrm{f}, ..)}(A \ot B^{\mrm{op}})$.
Then the following hold\tup{:}
\begin{enumerate}
\item The complex $D_A(M)$ belongs to 
$\dcat{D}^{-\star}_{(.., \mrm{f})}(B \ot A^{\mrm{op}})$, where $-\star$ 
is the reversed boundedness indicator. 

\item The derived Hom-evaluation morphism 
$\opn{ev}^{\mrm{R, R}}_{M} : M \to D_{A^{\mrm{op}}}(D_{A}(M))$
in $\dcat{D}(A \ot B^{\mrm{op}})$ is an isomorphism. 
\end{enumerate}
\end{thm}

\begin{proof}
The proof is the same as that of Theorem \ref{thm:3820}, with the obvious 
modifications (i.e.\ neglecting the algebraic grading). 
It all boils down to Theorems \ref{thm:2160} and \ref{thm:2135}.
\end{proof}

\begin{cor} \label{cor:3680}
Under the assumptions of Theorem \tup{\ref{thm:3680}}, 
the functor
\[ D_{A} : \dcat{D}^{\star}_{(\mrm{f}, ..)}(A \ot B^{\mrm{op}})^{\mrm{op}} \to 
\dcat{D}^{-\star}_{(.., \mrm{f})}(B \ot A^{\mrm{op}}) \]
is an equivalence of triangulated categories, with quasi-inverse
$D_{A^{\mrm{op}}}$. 
\end{cor}

\begin{proof}
First we note that the op-symmetry (Proposition \ref{prop:4355}) implies that  
Theorem \ref{thm:3680} is true after exchanging $A$ with $A^{\mrm{op}}$. 
Thus 
$\opn{ev}^{\mrm{R, R}}_{N} : N \to D_{A}(D_{A^{\mrm{op}}}(N))$
is an isomorphism for every 
$N \in \dcat{D}^{-\star}_{(.., \mrm{f})}(B \ot A^{\mrm{op}})$. 
Now the assertion is clear. 
\end{proof}

We need a better understanding of the derived homothety morphisms. 
Suppose $M \in \dcat{C}(A \ot B^{\mrm{op}})$. The right action of $B$ on $M$ is 
the homothety ring homomorphism through $B^{\mrm{op}}$,
$\opn{hm}_{M, B^{\mrm{op}}} : B^{\mrm{op}} \to
\opn{End}_{\dcat{C}_{\mrm{str}}(A)}(M)$.
Namely for elements $b \in B$ and $m \in M^i$ we have 
$\opn{hm}_{M, B^{\mrm{op}}}(\opn{op}(b))(m) = m \cd b \in M^i$.
By composing $\opn{hm}_{M, B^{\mrm{op}}}$ with the localization functor 
$\opn{Q} : \dcat{C}_{\mrm{str}}(A) \to \dcat{D}(A)$
we obtain the ring homomorphism 
\begin{equation} \label{eqn:4357}
\opn{hm}_{M, B^{\mrm{op}}}^{\dcat{D}} := 
\opn{Q} \circ \opn{hm}_{M, B^{\mrm{op}}} : 
B^{\mrm{op}} \to \opn{End}_{\dcat{D}(A)}(M) .
\end{equation}

\begin{lem} \label{lem:4361}
Let $M \in \dcat{D}(A \ot B^{\mrm{op}})$.
The three conditions below are equivalent. 
\begin{itemize}
\rmitem{i} $M$ has the derived Morita property on the $A$ side, i.e.\ 
\tup{(\ref{equ:4336})} is an isomorphism.

\rmitem{ii} For every $p \neq 0$ the module 
$\opn{Hom}_{\dcat{D}(A)}(M, M[p])$
is zero, and the ring homomorphism 
$\opn{hm}_{M, B^{\mrm{op}}}^{\dcat{D}} : B^{\mrm{op}} \to
\opn{End}_{\dcat{D}(A)}(M)$
is bijective. 

\rmitem{iii} For every $p \neq 0$ the module
$\opn{H}^p \bigl( \opn{RHom}_{A}(M, M) \bigr)$
is zero, and the $B^{\mrm{op}}$-module 
$\opn{H}^0 \bigl( \opn{RHom}_{A}(M, M) \bigr)$
is free, with basis the element 
$\opn{id}_M$. 
\end{itemize}
\end{lem}

\begin{proof}
This is clear from the canonical isomorphisms 
\[ \opn{Hom}_{\dcat{D}(A)} \bigl( M, M[p] \bigr) \cong 
\opn{H}^p \bigl( \opn{RHom}_{A}(M, M) \bigr) \]
of Corollary \ref{cor:2120}. 
\end{proof}

\begin{lem} \label{lem:4362}
Let $T \in \dcat{D}(A^{\mrm{en}})$, and consider the triangulated functor 
\[ G := T \ot^{\mrm{L}}_{A} (-) : 
\dcat{D}(A \ot B^{\mrm{op}}) \to \dcat{D}(A \ot B^{\mrm{op}}) . \]
For every $M \in \dcat{D}(A \ot B^{\mrm{op}})$  
the diagram 
\[ \UseTips \xymatrix @C=12ex @R=6ex {
B^{\mrm{op}}
\ar[r]^(0.4){ \opn{hm}_{M, B^{\mrm{op}}}^{\dcat{D}} }
\ar[dr]_(0.45){ \opn{hm}_{G(M), B^{\mrm{op}}}^{\dcat{D}} }
&
\opn{End}_{\dcat{D}(A)}(M)
\ar[d]^{G}
\\
&
\opn{End}_{\dcat{D}(A)}(G(M))
} \]
of ring homomorphisms is commutative. 
\end{lem}

\begin{proof}
This is because the action of $B^{\mrm{op}}$ on $M$ is categorical, i.e.\ by
endomorphisms in $\dcat{D}(A)$. 
\end{proof}

\begin{lem} \label{lem:4363}
Let $R \in \dcat{D}(A^{\mrm{en}})$, and consider the triangulated functor 
\[ D := \opn{RHom}_{A}(-, R) : 
\dcat{D}(A \ot B^{\mrm{op}})^{\mrm{op}} \to \dcat{D}(B \ot A^{\mrm{op}}) . \]
For every $M \in \dcat{D}(A \ot B^{\mrm{op}})$  
the diagram 
\[ \UseTips \xymatrix @C=12ex @R=6ex {
B^{\mrm{op}}
\ar[r]^(0.4){ \opn{hm}_{M, B^{\mrm{op}}}^{\dcat{D}} }
\ar[dr]_(0.45){ \opn{Op}(\opn{hm}_{D(M), B}^{\dcat{D}}) \quad }
&
\opn{End}_{\dcat{D}(A)}(M)
\ar[d]^{D}
\\
&
\opn{End}_{\dcat{D}(A^{\mrm{op}})}(D(M))^{\mrm{op}}
} \]
of ring homomorphisms is commutative. 
\end{lem}

\begin{proof}
Same as the previous lemma, only here $D$ is contravariant. 
\end{proof}

Tilting DG bimodules were defined in Definition \ref{dfn:3455}. Here we call 
them {\em tilting complexes}. 
Given a tilting complex $T \in \dcat{D}(A^{\mrm{en}})$, we know that its 
quasi-inverse $T^{\vee}$ satisfies
$T^{\vee} \cong \opn{RHom}_{A}(T, A) \cong  \opn{RHom}_{A^{\mrm{op}}}(T, A)$ 
in $\dcat{D}(A^{\mrm{en}})$. See Corollary \ref{cor:3455} with $B = A$,
and use the ring isomorphism (\ref{eqn:4450}). 

\begin{lem} \label{lem:4350}
Let $T$ be a tilting complex complex over $A$, with 
quasi-inverse $T^{\vee}$. Then there is an isomorphism 
$T \ot^{\mrm{L}}_{A} (-) \cong \opn{RHom}_{A}(T^{\vee}, -)$
of triangulated functors from 
$\dcat{D}(A \ot B^{\mrm{op}})$ to itself. 
\end{lem}

\begin{proof}
Let's write $S := T^{\vee}$. Then 
$T \cong S^{\vee} \cong \opn{RHom}_{A}(S, A)$. 
We know that $S$ is perfect on the $A$ side (see Corollary \ref{cor:3505}),
so by Theorem \ref{thm:3400} we have 
\begin{align*}
&
T \ot^{\mrm{L}}_{A} M \cong \opn{RHom}_{A}(S, A) \ot^{\mrm{L}}_{A} M 
\\ & \quad 
\cong \opn{RHom}_{A}(S, A \ot^{\mrm{L}}_{A} M) \cong 
\opn{RHom}_{A}(T^{\vee}, M)  . \qedhere
\end{align*} 
\end{proof}

\begin{thm}[\cite{Ye4}] \label{thm:3681}
Let $A$ be a noetherian ring, and let $R$ be a NC dualizing complex over $A$. 
\index{Dualizing complex! noncommutative}
\index{Tilting complex!}
\begin{enumerate}
\item Given a tilting complex $T$ over $A$, the complex
$R' := T \ot^{\mrm{L}}_{A} R$ 
is also a NC dualizing complex over $A$. 

\item Given a second NC dualizing complex $R'$ over $A$, the complex 
$T := \lb \opn{RHom}_A(R, R')$ 
is a tilting complex over $A$, and 
$R' \cong T \ot^{\mrm{L}}_{A} R$
in $\dcat{D}(A^{\mrm{en}})$. 

\item If $T$ is a tilting complex over $A$, and if
$T \ot^{\mrm{L}}_{A} R \cong R$ 
in $\dcat{D}(A^{\mrm{en}})$,  
then $T \cong A$ in $\dcat{D}(A^{\mrm{en}})$.
\end{enumerate}
\end{thm}

\begin{proof} \mbox{}

\smallskip \noindent
(1) We need to verify that $R'$ satisfies conditions (i)-(iii) of Definition 
\ref{dfn:3675}. 

Consider the equivalence of triangulated categories 
$G := T \ot^{\mrm{L}}_{A} (-) : \dcat{D}(A) \to \dcat{D}(A)$.
The functor $G$ has finite cohomological dimension (see Theorem \ref{thm:3400}, 
Corollary \ref{cor:3505} and Proposition \ref{prop:3415}). Now 
$G(A) = T$, and  
$T \in \dcat{D}^{\mrm{b}}_{\mrm{f}}(A)$
by Theorem \ref{thm:3415}, so Theorem \ref{thm:4600}
tells us that $G(M) \in \dcat{D}^{\mrm{b}}_{\mrm{f}}(A)$
for every $M \in \dcat{D}^{\mrm{b}}_{\mrm{f}}(A)$. 
Taking $M := R$ we see that 
$R' = G(R) \in \dcat{D}^{\mrm{b}}_{\mrm{f}}(A)$. 

Let $T^{\vee}$ be the quasi-inverse of $T$. By Lemma \ref{lem:4350}, with 
$B = A$, we see that 
\begin{equation} \label{eqn:4365}
R' = T \ot^{\mrm{L}}_{A} R \cong \opn{RHom}_{A}(T^{\vee}, R) = 
D_A(T^{\vee})
\end{equation}
in $\dcat{D}(A^{\mrm{en}})$. The tilting complex $T^{\vee}$ belongs to 
$\dcat{D}^{\mrm{b}}_{\mrm{f}}(A)$, so according to Theorem \ref{thm:3680}(1) 
with $B = \K$, we see that 
$D_A(T^{\vee}) \in \dcat{D}^{\mrm{b}}_{\mrm{f}}(A^{\mrm{op}})$. 
Thus $R' \in \dcat{D}^{\mrm{b}}_{\mrm{f}}(A^{\mrm{op}})$.
We have now verified condition (i) of Definition \ref{dfn:3675} for $R'$. 

To prove that $R'$ has finite injective dimension over $A$, we need to find a 
uniform bound on the cohomology of $\opn{RHom}_{A}(M, R')$, for all 
$M \in \dcat{M}(A)$. This is the same as finding a natural number $d$ such 
that $\opn{Hom}_{\dcat{D}(A)}(M, R'[i]) = 0$ for all 
$M \in \dcat{M}(A)$ and $\abs{i} > d$. 
We shall use the equivalence of triangulated categories 
\[ G^{\vee} := T^{\vee} \ot^{\mrm{L}}_{A} (-) \cong 
\opn{RHom}_{A}(T, -) : \dcat{D}(A) \to \dcat{D}(A) . \]
Since $R \cong G^{\vee}(R')$, we obtain a $\K$-module isomorphism
\[ G^{\vee} : \opn{Hom}_{\dcat{D}(A)} \bigl( M, R'[i] \bigr) \iso
\opn{Hom}_{\dcat{D}(A)} \bigl( G^{\vee}(M), R[i] \bigr)  \]
for every $i$. The functor $G^{\vee}$ has finite cohomological dimension, and 
$R$ has finite injective dimension; their sum $d$ serves as the required bound. 

Now let us prove that $R'$ has finite injective dimension over $A^{\mrm{op}}$.
Take a module $M \in \dcat{M}(A^{\mrm{op}})$. 
We have these isomorphisms for every integer $i$~:
\[ \begin{aligned}
& 
\opn{H}^i \bigl( \opn{RHom}_{A^{\mrm{op}}}(M, R') \bigr) \cong^{1}
\opn{Hom}_{\dcat{D}(A^{\mrm{op}})} \bigl( M, R'[i] \bigr)
\\ & \qquad 
\cong^{2} 
\opn{Hom}_{\dcat{D}(A^{\mrm{op}})} \bigl( M, D_A(T^{\vee})[i] \bigr)
\\ & \qquad 
\cong^{3} 
\opn{Hom}_{\dcat{D}(A)} \bigl( T^{\vee}, D_{A^{\mrm{op}}}(M)[i] \bigr)
\\ & \qquad 
\cong^{4} \opn{H}^i  \bigl(
\opn{RHom}_{A} \bigl( T^{\vee}, \opn{RHom}_{A^{\mrm{op}}}(M , R) \bigr)
\bigr)
\\ & \qquad
\cong^{5} \opn{H}^i  \bigl( 
T \ot^{\mrm{L}}_{A} \opn{RHom}_{A^{\mrm{op}}}(M, R) \bigr) 
\\ & \qquad
= \opn{H}^i  \bigl( (G \circ D_{A^{\mrm{op}}})(M) \bigr) . 
\end{aligned} \]
The justifications are as follows: 
\begin{itemize}
\item[$\cong^{1}$~:] This is by Corollary \ref{cor:2120}.

\item[$\cong^{2}$~:] This is by formula (\ref{eqn:4365}).

\item[$\cong^{3}$~:] This relies on Corollary \ref{cor:3680}.

\item[$\cong^{4}$~:] Again we use Corollary \ref{cor:2120}.

\item[$\cong^{5}$~:] It is due to Lemma \ref{lem:4350}.
\end{itemize}
The functors $G$ and $D_{A^{\mrm{op}}}$ have finite cohomological dimensions, 
so we have a uniform bound on the cohomology of 
$(G \circ D_{A^{\mrm{op}}})(M)$.
This completes the proof that $R'$ satisfies condition (ii) of Definition 
\ref{dfn:3675}. 

To prove that $R'$ has the derived Morita property on the $A$ side we 
use the fact that $R$ has the derived Morita property on the $A$ side, and we
invoke Lemma \ref{lem:4361} with $B := A$ and $M := R$. Since $R' = G(R)$, for 
every $i \in \Z$ we have an isomorphism of $\K$-modules 
\[ G : \opn{Hom}_{\dcat{D}(A)} \bigl( R, R[i] \bigr) \iso
\opn{Hom}_{\dcat{D}(A)} \bigl( R', R'[i] \bigr) , \]
and this is zero when $i \neq 0$. For $i = 0$ we use Lemma \ref{lem:4362}:
there is a commutative diagram of rings
\[ \UseTips \xymatrix @C=12ex @R=6ex {
A^{\mrm{op}}
\ar[r]^(0.4){\opn{hm}_{A^{\mrm{op}}, R}^{\dcat{D}}}_{\cong}
\ar[dr]_(0.45){\opn{hm}_{A^{\mrm{op}}, R'}^{\dcat{D}}}
&
\opn{End}_{\dcat{D}(A)}(R)
\ar[d]^{G}_{\cong}
\\
&
\opn{End}_{\dcat{D}(A)}(R')
} \]
and the vertical arrow $G$ is an isomorphism. 
Therefore $\opn{hm}_{A^{\mrm{op}}, R'}^{\dcat{D}}$ is a ring isomorphism. 

Finally, to prove that $R'$ has the derived Morita property on the 
$A^{\mrm{op}}$ side we use the fact that the tilting complex 
$T^{\vee}$ has the derived Morita property on the $A$ side 
(by Corollary \ref{cor:3505}), together with Lemma \ref{lem:4361}, applied to 
the complexes 
$T^{\vee} \in \dcat{M}(A \ot A^{\mrm{op}})$
and $R' \in \dcat{M}(A^{\mrm{op}} \ot A)$.
Recall from (\ref{eqn:4365}) that 
$R' \cong D_A(T^{\vee})$ in $\dcat{D}(A^{\mrm{en}})$.
By Corollary \ref{cor:3680}, for every $i \in \Z$ we 
have an isomorphism of $\K$-modules 
\[ D_A : \opn{Hom}_{\dcat{D}(A)} \bigl( T^{\vee}, T^{\vee}[i] \bigr) \iso
\opn{Hom}_{\dcat{D}(A^{\mrm{op}})} \bigl( R', R'[i] \bigr) , \]
and this is zero when $i \neq 0$. For $i = 0$ we use Lemma \ref{lem:4363}:
there is a commutative diagram of rings
\[ \UseTips \xymatrix @C=12ex @R=6ex {
A^{\mrm{op}}
\ar[r]^(0.4){\opn{hm}_{A^{\mrm{op}}, T^{\vee}}^{\dcat{D}}}_{\cong}
\ar[dr]_(0.45){ \opn{Op}(\opn{hm}_{A, R'}^{\dcat{D}}) \ }
&
\opn{End}_{\dcat{D}(A)}(T^{\vee})
\ar[d]^{D_A}_{\cong}
\\
&
\opn{End}_{\dcat{D}(A^{\mrm{op}})}(R')^{\mrm{op}}
} \]
It follows that 
$\opn{hm}_{A, R'}^{\dcat{D}} : A \to \opn{End}_{\dcat{D}(A^{\mrm{op}})}(R')$
is a ring isomorphism. 

\medskip \noindent 
(2) The proof of this item is a minor modification of the proof of Theorem 
\ref{thm:3780}. Indeed, an ungraded version of Lemma \ref{lem:3917} holds here, 
and we use Theorem \ref{thm:3680} instead of Theorem \ref{thm:3820}. 

\medskip \noindent 
(3) Assume that $T \ot^{\mrm{L}}_{A} R \cong R$ in $\dcat{D}(A^{\mrm{en}})$.
Of course $R = D_A(A)$. Formula (\ref{eqn:4365}) says that 
$D_A(A) \cong D_A(T^{\vee})$ in $\dcat{D}(A^{\mrm{en}})$.
By Corollary \ref{cor:3680} we conclude that 
$A \cong T^{\vee}$ in $\dcat{D}(A^{\mrm{en}})$.
But then $A \cong T$  in $\dcat{D}(A^{\mrm{en}})$.
\end{proof}

The op-symmetry that gave rise to Proposition \ref{prop:4355} implies that 
several variations of Theorem \ref{thm:3681} are true. We wish to write down 
only one of them, that will be needed later. 

\begin{cor} \label{cor:4365}
Let $A$ be a noetherian ring, and let $R$ and $R'$ be NC dualizing complexes 
over $A$. Then there is a tilting complex $T'$ such that 
$R' \cong R \ot^{\mrm{L}}_{A} T'$
in $\dcat{D}(A^{\mrm{en}})$. 
\end{cor}

\begin{proof}
This is a trivial (yet confusing) consequence of item (2) of  
Theorem \ref{thm:3681}, when viewing $R$ and $R'$ as NC dualizing complexes 
over the ring $A^{\mrm{op}}$. 
\end{proof}

Recall that a left action of a group $G$ on a nonempty set $X$ is called {\em 
simply transitive} if there is one $G$-orbit in $X$, and the stabilizer of 
every point $x \in X$ is trivial. I.e.\ for every (or equivalently, for some) 
point $x \in X$ the action function 
$G \times \{ x \} \to X$, $(g, x) \mapsto g(x)$, is bijective. 

\begin{cor} \label{cor:3690}
Let $A$ be a ring, and assume it has at least one NC dualizing complex. 
\index{1-DPic(A)@$\opn{DPic}_{\K}(A)$}
\index{Picard group! noncommutative derived}
Then the left action of the group $\opn{DPic}_{\K}(A)$ on the set of 
isomorphism classes of NC dualizing complexes, induced by the bifunctor
$(T, R) \mapsto T \ot^{\mrm{L}}_{A} R$,
is simply transitive.  
\end{cor}

\begin{proof}
By Theorem \ref{thm:3681}(1) this is a well-defined action:
$T \ot^{\mrm{L}}_{A} R$ is a NC dualizing complex. 
The action is transitive by Theorem \ref{thm:3681}(2), and it has trivial 
stabilizers by Theorem \ref{thm:3681}(3).
\end{proof}

\begin{rem} \label{rem:3696}
Suppose $A$ is a noetherian commutative $\K$-ring.
A noncommutative dualizing complex over $A$ is not the same as a commutative
dualizing complex over $A$, in the sense of Definition \ref{dfn:2155}. Indeed, 
there is a $\K$-ring 
homomorphism 
$A^{\mrm{en}} \to A$, and this induces a restriction functor 
$\dcat{D}(A) \to \dcat{D}(A^{\mrm{en}})$. 
If $R \in \dcat{D}(A)$ is a commutative dualizing complex over $A$, then
its image in $\dcat{D}(A^{\mrm{en}})$ is a noncommutative dualizing complex 
over $A$. However, according to Theorem \ref{thm:3681}(1), 
given any tilting complex $T$ over $A$, the complex 
$R' := T \ot^{\mrm{L}}_{A} R$
is also a NC dualizing complex over $A$. And $R'$ can very easily fail to be in 
the essential image of $\dcat{D}(A)$; e.g.\ by taking 
$T := A(\psi)$, the twist by an automorphism $\psi$ of the trivial bimodule 
$A$. 
\end{rem}

\begin{rem} \label{rem:3691}
In case the base ring $\K$ is not a field, but $A$ is a {\em flat noetherian 
$\K$-ring}, then the definitions and results in this subsection are 
still valid. 

In case the ring $A$ is noetherian, but is not flat over the base ring $\K$, we 
can still define NC dualizing complexes over it. For that we choose a 
K-flat DG ring resolution $\til{A} \to A$ over $\K$; see 
Theorem \ref{thm:4310} and Proposition \ref{prop:4305}. 
A NC dualizing complex over $A$ is then a complex
$R \in \dcat{D}(\til{A}^{\mrm{en}})$ 
that satisfies the three conditions of Definition \ref{dfn:3675}.
As explained in Remark \ref{rem:1420}, the triangulated category 
$\dcat{D}(\til{A}^{\mrm{en}})$ is independent of the resolution 
$\til{A}$, up to a canonical equivalence. This means that the dualizing complex 
$R$ is also independent of the resolution.

It is expected that Theorems \ref{thm:3680} and \ref{thm:3681} will still hold 
in this general situation; but we did not verify this assertion. 

However, flatness alone is not sufficient for many of the results 
in the following subsections (mainly in Subsection 
\ref{subsec:Fil-exit-R-NCDC}, in which we rely on results from Section 
\ref{sec:BDC}, that assumed Convention \ref{conv:3700}.) Therefore we decided 
to assume $\K$ is a field in the current subsection as well.
\end{rem}

\mysubsection{Rigid NC DC: Definition and Uniqueness}
\label{subsec:RNCDC-uniq}

In this subsection we introduce rigid NC dualizing complexes, and prove their 
uniqueness. The material is mostly from Van den Bergh's paper \cite{VdB}, with 
some improvements coming from \cite{Ye4}. 
We continue with Convention \ref{conv:3670}; in particular, $\K$ is a base 
field, and all rings are central over $\K$.

Let $A$ be a ring, and let 
$M, N \in \dcat{M}(A^{\mrm{en}})$.
The $\K$-module $M \ot N$ has two $A$-module structures and two 
$A^{\mrm{op}}$-module structures, all commuting with each other. 
We organize them into an {\em outside action} and an {\em inside 
action} of the ring $A^{\mrm{en}}$. Here are the formulas: given elements 
$m \in M$, $n \in N$ and $a_1, a_2 \in A$, we let
\begin{equation} \label{eqn:4215}
(a_1 \ot \opn{op}(a_2)) \, \cd^{\mrm{out}} \, (m \ot n) := 
(a_1 \cd m) \ot (n \cd a_2) 
\end{equation}
and
\begin{equation} \label{eqn:4216}
(a_1 \ot \opn{op}(a_2)) \, \cd^{\mrm{in}} \, (m \ot n) := 
(m \cd a_2) \ot (a_1 \cd n) . 
\end{equation}
In these formulas we are using the canonical anti-automorphism
$\opn{op} : A \to A^{\mrm{op}}$. 
We shall denote the two copies 
of $A^{\mrm{en}}$ that act on $M \ot N$ by
\begin{equation} \label{eqn:4460}
A^{\mrm{en, out}} = A^{\mrm{out}} \ot A^{\mrm{op, out}}
\end{equation}
and
\begin{equation} \label{eqn:4461}
A^{\mrm{en, in}} = A^{\mrm{in}} \ot A^{\mrm{op, in}} . 
\end{equation}
Namely the ring $A^{\mrm{en, out}}$ acts on $M \ot N$ by the outside action
(\ref{eqn:4215}), 
and the ring $A^{\mrm{en, in}}$ acts by the inside action (\ref{eqn:4216}). 

Let us define the ring 
\begin{equation} \label{eqn:4411}
A^{\mrm{four}} := A^{\mrm{en, out}} \ot A^{\mrm{en, in}}  .
\end{equation}
There are ring homomorphisms 
\begin{equation} \label{eqn:4500}
A^{\mrm{en, out}} \to A^{\mrm{four}} \lar \, A^{\mrm{en, in}}  .
\end{equation}
With this notation we have 
$M \ot N \in \dcat{M}(A^{\mrm{four}})$.
The action of the ring $A^{\mrm{four}}$ on the module $M \ot N$ is this, 
explicitly:  an element 
\[ \bigl( a_1 \ot \opn{op}(a_2) \bigr) \ot 
\bigl( a_3 \ot \opn{op}(a_4) \bigr) \in A^{\mrm{four}}\]
acts on an element 
$m \ot n \in M \ot N$ by
\[ \bigl( (a_1 \ot \opn{op}(a_2)) \ot (a_3 \ot \opn{op}(a_4)) \bigr) 
\cd (m \ot n) = 
(a_1 \cd m \cd a_4) \ot (a_3 \cd n \cd a_2) . \]

In some situations we prefer to view the inside action of $A^{\mrm{en}}$ on
$M \ot N$ as a right action. This is possible by the isomorphism 
$\opn{trop}$ from formula (\ref{eqn:4875}). 
To be explicit, for elements $m \ot n \in M \ot N$ and 
$a_3 \ot \opn{op}(a_4) \in A^{\mrm{en, in}}$
the action is as follows, using formula (\ref{eqn:4878}):
\begin{equation} \label{eqn:4877}
\begin{aligned}
&
(m \ot n) \, \cd^{\mrm{in}} \, (a_3 \ot \opn{op}(a_4)) = 
(a_4 \ot \opn{op}(a_3))\, \cd^{\mrm{in}} \, (m \ot n)
\\ & \quad 
= (m \cd a_3) \ot (a_4 \cd n) . 
\end{aligned}
\end{equation}

\begin{exa} \label{exa:5125}
Suppose we take $M = N := A \in \dcat{M}(A^{\mrm{en}})$.
Then $P := M \ot N = A \ot A$ can also be viewed as 
$P \cong A^{\mrm{en}}$, and then in it an 
$A^{\mrm{en}}$-bimodule. The left $A^{\mrm{en}}$-module structure on 
$P \cong A^{\mrm{en}}$ is this: for a module element $b_1 \ot b_2 \in P$ and a 
ring element $a_1 \ot \opn{op}(a_2) \in A^{\mrm{en}}$ we have 
\[ \bigl( a_1 \ot \opn{op}(a_2) \bigr) \cd (b_1 \ot b_2) = 
(a_1 \cd b_1) \ot (b_2 \cd a_2) ; \]
so this is the outside action on $P = A \ot A$. The right $A^{\mrm{en}}$-module 
structure on $P \cong A^{\mrm{en}}$ is this: for a ring element 
$a_3 \ot \opn{op}(a_4) \in A^{\mrm{en}}$ we have 
\[ (b_1 \ot b_2) \cd \bigl( a_3 \ot \opn{op}(a_4) \bigr) = 
(b_1 \cd a_3) \ot (a_4 \cd b_2) ; \]
so this is the inside action on $P = A \ot A$, cf.\ (\ref{eqn:4877}). 
\end{exa}

Given bimodules $L, M, N \in \dcat{M}(A^{\mrm{en}})$,
we denote by
$\opn{Hom}_{A^{\mrm{en, out}}}(L, M \ot N)$
the $\K$-module of homomorphisms $\phi : L \to M \ot N$ 
that are $A^{\mrm{en}}$-linear for the outside action of $A^{\mrm{en}}$ on 
$M \ot N$. This is an 
$A^{\mrm{en}}$-module by the inside action on $M \ot N$. 
Namely 
\[ \opn{Hom}_{A^{\mrm{en, out}}}(L, M \ot N)
\in \dcat{M}(A^{\mrm{en, in}}) = \dcat{M}(A^{\mrm{en}}) . \]

All these four-module operations extend in an obvious way to complexes, and 
can be derived. 

\begin{dfn} \label{dfn:4215}
Let $A$ be a ring and $M \in \dcat{D}(A^{\mrm{en}})$. 
The {\em NC square of $M$} 
\index{Squaring operation! noncommutative}
\index{1-SqA@$\opn{Sq}_{A}(M)$}
is the complex
\[ \opn{Sq}_{A}(M) := 
\opn{RHom}_{A^{\mrm{en, out}}}(A, M \ot M) \in 
\dcat{D}(A^{\mrm{en, in}}) = \dcat{D}(A^{\mrm{en}}) . \]
\end{dfn}

Of course the square is relative to the base field $\K$, which is implicit in 
the formulas. 

The square of the complex $M$ is calculated as follows. First note that $M$ is 
K-flat over $\K$, since $\K$ is a field; so 
$M \ot M = M \ot^{\mrm{L}}_{\K} M \in \dcat{D}(A^{\mrm{four}})$.
We choose a K-injective resolution 
$M \ot M \to I$ in $\dcat{C}_{\mrm{str}}(A^{\mrm{four}})$. 
The complex $I$ is unique up to a homotopy equivalence, that itself is unique 
up to homotopy (in other words, $I$ is unique up to a unique isomorphism in 
$\dcat{K}(A^{\mrm{four}})$). Because of the flatness of the ring homomorphisms 
(\ref{eqn:4500}), the complex 
$I$ is K-injective over $A^{\mrm{en, out}}$, and 
\begin{equation} \label{eqn:4370}
\opn{Sq}_{A}(M) = \opn{Hom}_{A^{\mrm{en, out}}}(A, I) \in 
\dcat{D}(A^{\mrm{en, in}}) = \dcat{D}(A^{\mrm{en}}) .
\end{equation}

\begin{dfn} \label{dfn:4370}
Let $\phi : M \to N$ be a morphism in $\dcat{D}(A^{\mrm{en}})$. 
The {\em NC square of $\phi$} 
\index{1-SqA@$\opn{Sq}_{A}(\phi)$}
is the morphism
\[ \opn{Sq}_{A}(\phi) := 
\opn{RHom}_{A^{\mrm{en, out}}}(\opn{id}_A, \phi \ot \phi) :
\opn{Sq}_{A}(M) \to \opn{Sq}_{A}(N) \]
in $\dcat{D}(A^{\mrm{en}})$.
\end{dfn}

The square of the morphism $\phi$ is calculated as follows. Choose a 
K-projective resolution $P \to M$ in $\dcat{C}_{\mrm{str}}(A^{\mrm{en}})$. 
Then $\phi$ is represented by a homomorphism 
$\til{\phi} : P \to N$ in 
$\dcat{C}_{\mrm{str}}(A^{\mrm{en}})$. The morphism 
$\phi \ot \phi : M \ot M \to N \ot N$ 
in $\dcat{D}(A^{\mrm{four}})$ is represented by the homomorphism 
$\til{\phi}  \ot \til{\phi} : P \ot P \to N \ot N$ 
in $\dcat{C}_{\mrm{str}}(A^{\mrm{four}})$.
Next we choose K-injective resolutions 
$P \ot P \to I$ and $N \ot N \to J$ in 
$\dcat{C}_{\mrm{str}}(A^{\mrm{four}})$. 
The homomorphism $\til{\phi}  \ot \til{\phi}$ induces a homomorphism 
$\til{\psi} : I \to J$ 
in $\dcat{C}_{\mrm{str}}(A^{\mrm{four}})$. Note that $\til{\psi}$ is unique up 
to homotopy. Finally
\begin{align*}
& \opn{Sq}_{A}(\phi) = 
\opn{Q} \bigl( \opn{Hom}_{A^{\mrm{en, out}}}(\opn{id}_A, \til{\psi}) \bigr) :
\\ & \quad 
\opn{Sq}_{A}(M) = \opn{Hom}_{A^{\mrm{en, out}}}(A, I)
\to \opn{Sq}_{A}(N) = \opn{Hom}_{A^{\mrm{en, out}}}(A, J) 
\end{align*}
in $\dcat{D}(A^{\mrm{en, in}}) = \dcat{D}(A^{\mrm{en}}) $.

\begin{dfn} \label{dfn:4216}
Let $A$ be a ring. A {\em NC rigid complex over $A$} 
\index{Rigid complex! noncommutative}
is a pair $(M, \rho)$, where $M \in \dcat{D}(A^{\mrm{en}})$, and 
$\rho : M \iso \opn{Sq}_{A}(M)$
is an isomorphism in $\dcat{D}(A^{\mrm{en}})$, called a {\em NC rigidifying 
isomorphism}. 
\end{dfn}

\begin{dfn} \label{dfn:4372} 
Let $(M, \rho)$ and $(N, \si)$ be NC rigid complexes over $A$. A
{\em NC rigid morphism over $A$} 
\[ \phi : (M, \rho) \to (N, \si) \]
is a morphism $\phi : M \to N$ 
in $\dcat{D}(A^{\mrm{en}})$, such that the diagram 
\[ \UseTips \xymatrix @C=8ex @R=6ex {
M
\ar[r]^(0.35){\rho}
\ar[d]_{\phi}
&
\opn{Sq}_{A}(M)
\ar[d]^{\opn{Sq}_{A}(\phi)}
\\
N
\ar[r]^(0.35){\si}
&
\opn{Sq}_{A}(N)
} \]
in $\dcat{D}(A^{\mrm{en}})$ is commutative. 
\end{dfn}

\begin{dfn} \label{dfn:4380}
Let $A$ be a ring. The category of {\em NC rigid complexes over $A$} 
is the category whose objects are the NC rigid complexes from Definition 
\ref{dfn:4216}, and whose morphisms are the NC rigid morphisms  Definition 
\ref{dfn:4372}.
\end{dfn}

\begin{dfn}[Van den Bergh \cite{VdB}] \label{dfn:4375}
Let $A$ be a noetherian ring. A {\em rigid NC dualizing complex over $A$} 
\index{Rigid complex! noncommutative {\indash} dualizing {\indash}}
\index{Dualizing complex! noncommutative rigid}
is a NC rigid complex $(R, \rho)$, such that $R$ is a NC dualizing complex over 
$A$, in the sense of Definition \ref{dfn:3675}.
\end{dfn}

Again, we remind that the notion of rigidity is relative to the base field 
$\K$. 

Recall that the center of the ring $A$ is denoted by $\opn{Cent}(A)$. 
Of course $\opn{Cent}(A) = \opn{Cent}(A^{\mrm{op}})$. 
Given a bimodule $M \in \dcat{M}(A^{\mrm{en}})$, there are two ring 
homomorphisms
\[ \opn{chm}_{M, A}, \, \opn{chm}_{M, A^{\mrm{op}}} : 
\opn{Cent}(A) \to \opn{End}_{\dcat{M}(A^{\mrm{en}})}(M) , \]
that we call the 
{\em central homotheties}%
\index{Homothety morphism! central}
through $A$ and through 
$A^{\mrm{op}}$ respectively.  In terms of elements they are 
$\opn{chm}_{M, A}(a)(m) := a \cd m$ and
$\opn{chm}_{M, A^{\mrm{op}}}(a)(m) := m \cd a$ 
for $a \in \opn{Cent}(A)$ and $m \in M$. 

\begin{exa} \label{exa:4375}
For the bimodule $M = A$, both ring homomorphisms
\[ \opn{chm}_{A, A}, \, \opn{chm}_{A, A^{\mrm{op}}} : 
\opn{Cent}(A) \to \opn{End}_{\dcat{M}(A^{\mrm{en}})}(A) \]
are isomorphisms, and they are equal. 
\end{exa}

The ring homomorphisms $\opn{chm}_{M, A}$ and $\opn{chm}_{M, A^{\mrm{op}}}$ 
extend to complexes, namely for  $M \in \dcat{C}(A^{\mrm{en}})$ 
there are ring homomorphisms 
\[ \opn{chm}_{M, A}, \, \opn{chm}_{M, A^{\mrm{op}}} : 
\opn{Cent}(A) \to \opn{End}_{\dcat{C}_{\mrm{str}}(A^{\mrm{en}})}(M) . \]
By postcomposing them with the categorical localization functor $\opn{Q}$ we 
obtain ring homomorphisms 
\[ \opn{chm}^{\dcat{D}}_{M, A}, \, \opn{chm}^{\dcat{D}}_{M,  A^{\mrm{op}}} : 
\opn{Cent}(A) \to \opn{End}_{\dcat{D}(A^{\mrm{en}})}(M) . \]
These central homotheties were already used in Subsection 
\ref{subsec:tilting-rings}.

\begin{exa} \label{exa:4376}
For the complex of bimodules $M = A$, both ring homomorphisms
\[ \opn{chm}^{\dcat{D}}_{M, A}, \, \opn{chm}^{\dcat{D}}_{M,  A^{\mrm{op}}} : 
\opn{Cent}(A) \to \opn{End}_{\dcat{D}(A^{\mrm{en}})}(A) \]
are isomorphisms, and they are equal. This is because of the fully faithful 
embedding 
$\dcat{M}(A^{\mrm{en}}) \to \dcat{D}(A^{\mrm{en}})$ and the previous example. 
\end{exa}

\begin{lem} \label{lem:4380}
Given $R \in \dcat{D}(A^{\mrm{en}})$, consider the functor 
\[ D := \opn{Hom}_A(-, R) : \dcat{D}(A^{\mrm{en}})^{\mrm{op}} \to
\dcat{D}(A^{\mrm{en}}) . \]
Then for every $M \in \dcat{D}(A^{\mrm{en}})$ the diagram of rings 
\[ \UseTips \xymatrix @C=12ex @R=6ex {
\opn{Cent}(A)
\ar[r]^(0.4){ \opn{chm}^{\dcat{D}}_{M,  A^{\mrm{op}}}}_{}
\ar[dr]_(0.4){\opn{chm}^{\dcat{D}}_{D(M), A}}
&
\opn{End}_{\dcat{D}(A^{\mrm{en}})}(M)
\ar[d]^{D}_{}
\\
&
\opn{End}_{\dcat{D}(A^{\mrm{en}})} \bigl( D(M) \bigr)^{\mrm{op}}
} \]
is commutative. 
\end{lem}

\begin{proof}
Choose a K-injective resolution $R \to I$ in 
$\dcat{C}_{\mrm{str}}(A^{\mrm{en}})$, 
and consider the functor 
$\til{D} := \opn{Hom}_A(-, I)$. 
So $D \cong \til{D}$ as functors 
$\dcat{D}(A^{\mrm{en}})^{\mrm{op}} \to \dcat{D}(A^{\mrm{en}})$. 
Because of the contravariance of $\til{D}$, the diagram of rings 
\[ \UseTips \xymatrix @C=12ex @R=6ex {
\opn{Cent}(A)
\ar[r]^(0.4){ \opn{chm}_{M,  A^{\mrm{op}}}}_{}
\ar[dr]_(0.4){\opn{chm}_{\til{D}(M), A}}
&
\opn{End}_{\dcat{C}_{\mrm{str}}(A^{\mrm{en}})}(M)
\ar[d]^{\til{D}}_{}
\\
&
\opn{End}_{\dcat{C}_{\mrm{str}}(A^{\mrm{en}})} 
\bigl( \til{D}(M) \bigr)^{\mrm{op}}
} \]
is commutative. After applying the localization functor $\opn{Q}$ to the 
vertical arrow $\til{D}$ in the diagram above, we obtain the commutative 
diagram 
on the derived level. 
\end{proof}

\begin{lem} \label{lem:4381}
Given $M, N \in \dcat{D}(A^{\mrm{en}})$, consider the complex  
\[ L := \opn{RHom}_{A^{\mrm{en, out}}}(A, M \ot N) \in 
\dcat{D}(A^{\mrm{en, in}}) = \dcat{D}(A^{\mrm{en}}) . \]
Let $a, a', a'' \in \opn{Cent}(A)$ be elements such that 
$\opn{chm}^{\dcat{D}}_{N, A}(a) = \opn{chm}^{\dcat{D}}_{N, A^{\mrm{op}}}(a')$ 
in \lb $\opn{End}_{\dcat{D}(A^{\mrm{en}})}(N)$ 
and 
$\opn{chm}^{\dcat{D}}_{M, A}(a') = 
\opn{chm}^{\dcat{D}}_{M, A^{\mrm{op}}}(a'')$
in $\opn{End}_{\dcat{D}(A^{\mrm{en}})}(M)$.
Then \lb 
$\opn{chm}_{L, A}^{\dcat{D}}(a) = \opn{chm}^{\dcat{D}}_{L, A^{\mrm{op}}}(a'')$ 
in $\opn{End}_{\dcat{D}(A^{\mrm{en}})}(L)$.
\end{lem}

\begin{proof}
The proof resembles that of Lemma \ref{lem:3581}.
Consider the object 
$M \ot N\in \dcat{D}(A^{\mrm{four}})$. 
There are four ring homomorphisms 
\[ \opn{chm}^{\dcat{D}}_{M \ot N, B} : \opn{Cent}(A) \to 
\opn{End}_{\dcat{D}(A^{\mrm{four}})}(M \ot N) , \]
corresponding to these four options for the ring $B$,
in terms of formulas (\ref{eqn:4460}) and (\ref{eqn:4461}):
$B := A^{\mrm{out}}, A^{\mrm{op, out}}, A^{\mrm{in}} , A^{\mrm{op, in}}$. 
We know that 
\begin{equation} \label{eqn:4455}
\begin{aligned}
& \opn{chm}^{\dcat{D}}_{M \ot N, A^{\mrm{in}}}(a) = 
\opn{chm}^{\dcat{D}}_{M \ot N, A^{\mrm{op, out}}}(a') 
\\ & 
\opn{chm}^{\dcat{D}}_{M \ot N, A^{\mrm{out}}}(a') = 
\opn{chm}^{\dcat{D}}_{M \ot N, A^{\mrm{op, in}}}(a'') . 
\end{aligned}
\end{equation}

Let $M \ot N \to J$ be a K-injective resolution in 
$\dcat{C}_{\mrm{str}}(A^{\mrm{four}})$, 
so 
\[ \opn{End}_{\dcat{D}(A^{\mrm{four}})}(M \ot N) \cong 
\opn{H}^0 \bigl( \opn{End}_{\dcat{C}_{\mrm{str}}(A^{\mrm{four}})}(J) \bigr) \]
as rings. Formulas (\ref{eqn:4455}) imply that 
\begin{equation} \label{eqn:4456}
\begin{aligned}
& \opn{chm}_{J, A^{\mrm{in}}}(a) \sim  
\opn{chm}_{J, A^{\mrm{op, out}}}(a') , 
\\ & 
\opn{chm}_{J, A^{\mrm{out}}}(a') \sim
\opn{chm}_{J, A^{\mrm{op, in}}}(a'') , 
\end{aligned}
\end{equation}
where ``$\sim$'' means the homotopy relation on degree $0$ cocycles in the 
DG ring $\opn{End}_{A^{\mrm{four}}}(J)$. 

Now let us write 
$\til{L} = \opn{Hom}_{A^{\mrm{en, out}}}(A, J) \in 
\dcat{C}(A^{\mrm{en, in}})$.
The complex $\til{L}$ is K-injective over $A^{\mrm{en, in}}$, 
again by the flatness of the ring homomorphism (\ref{eqn:4500}), 
and $\til{L} \cong L$ in $\dcat{D}(A^{\mrm{en, in}})$. 
Hence there is a ring isomorphism 
\begin{equation} \label{eqn:4459}
\opn{H}^0 \bigl( \opn{End}_{A^{\mrm{en, in}}}(\til{L}) \bigr) \cong
\opn{End}_{\dcat{D}(A^{\mrm{en, in}})}(L) . 
\end{equation}
Now $A$ is a central bimodule over its center.
This implies that the two actions of $\opn{Cent}(A)$ on $\til{L}$ that 
come from $\opn{chm}_{J, A^{\mrm{out}}}$ and 
$\opn{chm}_{J, A^{\mrm{op, out}}}$ are equal. 
The two actions of $\opn{Cent}(A)$ on $\til{L}$ that 
come from $\opn{chm}_{J, A^{\mrm{in}}}$ and 
$\opn{chm}_{J, A^{\mrm{op, in}}}$
are the actions  $\opn{chm}_{\til{L}, A^{\mrm{in}}}$ and 
$\opn{chm}_{\til{L}, A^{\mrm{op, in}}}$
respectively. From formulas (\ref{eqn:4456}) we conclude that 
$\opn{chm}_{\til{L}, A^{\mrm{in}}}(a) \sim
\opn{chm}_{\til{L}, A^{\mrm{op, in}}}(a'')$. 
By formula (\ref{eqn:4459}) we see that 
$\opn{chm}_{L, A^{\mrm{in}}}^{\dcat{D}}(a) = 
\opn{chm}^{\dcat{D}}_{L, A^{\mrm{op, in}}}(a'')$
as claimed. 
\end{proof}

\begin{thm}[\cite{Ye4}] \label{thm:4375}
Let $A$ be a noetherian ring, and let $R$ be a NC dualizing complex over $A$.
\begin{enumerate}
\item The two ring homomorphisms
\[ \opn{chm}^{\dcat{D}}_{R, A}, \, \opn{chm}^{\dcat{D}}_{R, A^{\mrm{op}}} : 
\opn{Cent}(A) \to \opn{End}_{\dcat{D}(A^{\mrm{en}})}(R) \]
are bijective.

\item If $R$ is a rigid NC dualizing complex, then 
$\opn{chm}^{\dcat{D}}_{R, A} = \opn{chm}^{\dcat{D}}_{R, A^{\mrm{op}}}$.
\end{enumerate}
\end{thm}

\begin{proof} \mbox{}

\smallskip \noindent 
(1) Recall the duality functor $D_A$ from Definition \ref{dfn:3680}, with 
$B = A$. By Corollary \ref{cor:3680} it gives rise to an equivalence 
$D_{A} : \dcat{D}^{\mrm{b}}_{(\mrm{f}, ..)}(A^{\mrm{en}})^{\mrm{op}} 
\to \dcat{D}^{\mrm{b}}_{(.., \mrm{f})}(A^{\mrm{en}})$.
By Lemma \ref{lem:4380}, for every 
$M \in \dcat{D}^{\mrm{b}}_{(\mrm{f}, ..)}(A^{\mrm{en}})$
we get a commutative diagram of rings 
\[ \UseTips \xymatrix @C=12ex @R=6ex {
\opn{Cent}(A)
\ar[r]^(0.4){\opn{chm}^{\dcat{D}}_{M, A^{\mrm{op}}}}_{}
\ar[dr]_(0.4){\opn{chm}^{\dcat{D}}_{D_A(M), A}}
&
\opn{End}_{\dcat{D}(A^{\mrm{en}})}(M)
\ar[d]^{D_A}_{\cong}
\\
&
\opn{End}_{\dcat{D}(A^{\mrm{en}})} \bigl( D_A(M) \bigr)^{\mrm{op}}
} \]
Note that the vertical arrow $D_A$ is an isomorphism. 
Taking $M = R$ we have $D_A(R) \cong A$ in $\dcat{D}(A^{\mrm{en}})$. As we saw 
in Example \ref{exa:4376}, the homomorphism 
$\opn{chm}^{\dcat{D}}_{A, A}$ 
is bijective. This proves that 
$\opn{chm}^{\dcat{D}}_{R, A^{\mrm{op}}}$ is bijective. 

Likewise for $\opn{chm}^{\dcat{D}}_{R, A}$, but this time we use the 
duality functor $D_{A^{\mrm{op}}}$.

\medskip \noindent
(2) Now there is a rigidifying isomorphism 
$\rho : R \iso \opn{Sq}_{A}(R)$ in $\dcat{D}(A^{\mrm{en}})$. 
Let $f : \opn{Cent}(A) \to \opn{Cent}(A)$
be the ring automorphism 
$f := (\opn{chm}^{\dcat{D}}_{R, A^{\mrm{op}}})^{-1} \circ \, 
\opn{chm}^{\dcat{D}}_{R, A}$.
We need to prove that $f = \opn{id}$.

For this we shall use Lemma \ref{lem:4381}, with 
$M = N := R$, so $L \cong \opn{Sq}_A(R)$ in $\dcat{D}(A^{\mrm{en}})$. 
Take an element $a \in \opn{Cent}(A)$, and let 
$a' := f(a)$ and $a'' := f(a')$ in $\opn{Cent}(A)$. 
The lemma says that 
$\opn{chm}_{L, A}^{\dcat{D}}(a) = \opn{chm}^{\dcat{D}}_{L, A^{\mrm{op}}}(a'')$. 
But by the definition of $f$ and $a'$ we have  
$\opn{chm}_{R, A}^{\dcat{D}}(a) = \opn{chm}^{\dcat{D}}_{R, A^{\mrm{op}}}(a')$. 
Since $R \cong L$ in $\dcat{D}(A^{\mrm{en}})$, 
and since $\opn{chm}^{\dcat{D}}_{R, A^{\mrm{op}}}$
is bijective, we conclude that $a' = a''$. This means that $f(a) = f(f(a))$, 
so $a = f(a)$. Finally, the element $a$ was arbitrary, and hence
$f = \opn{id}$. 
\end{proof}

\begin{lem} \label{lem:4384}
Let $(R, \rho)$ and $(R', \rho')$ be rigid dualizing complexes, 
let $\phi : R \to R'$ be a morphism in $\dcat{D}(A^{\mrm{en}})$, and let
$a \in \opn{Cent}(A)$. Then 
\[ \opn{Sq}_{A} \bigl( \phi \circ \opn{chm}^{\dcat{D}}_{R, A}(a) \bigr) = 
\opn{Sq}_{A}(\phi) \circ \opn{chm}^{\dcat{D}}_{\opn{Sq}_{A}(R), A}(a^2) \]
as morphisms $\opn{Sq}_{A}(R) \to \opn{Sq}_{A}(R')$ in 
$\dcat{D}(A^{\mrm{en}})$.
\end{lem}

\begin{proof}
Because $\opn{Sq}_{A}$ is a functor, there is equality 
\[ \opn{Sq}_{A} \bigl( \phi \circ \opn{chm}^{\dcat{D}}_{R, A}(a) \bigr) =
\opn{Sq}_{A}(\phi) \circ 
\opn{Sq}_{A} \bigl( \opn{chm}^{\dcat{D}}_{R, A}(a) \bigr) .  \]
So we can assume that $R' = R$ and 
$\phi = \opn{id}_R$. We need to prove that the equality 
\begin{equation} \label{eqn:4501}
\opn{Sq}_{A} \bigl( \opn{chm}^{\dcat{D}}_{R, A}(a) \bigr) = 
\opn{chm}^{\dcat{D}}_{\opn{Sq}_{A}(R), A}(a^2) 
\end{equation}
holds for these endomorphisms of $\opn{Sq}_{A}(R)$ in $\dcat{D}(A^{\mrm{en}})$.

From Definition \ref{dfn:4370}, applied to the morphism 
$\opn{chm}^{\dcat{D}}_{R, A}(a) : R \to R$ 
in $\dcat{D}(A^{\mrm{en}})$, we see that 
\begin{equation} \label{eqn:4462}
\opn{Sq}_{A} \bigl( \opn{chm}^{\dcat{D}}_{R, A}(a) \bigr) = 
\opn{chm}^{\dcat{D}}_{\opn{Sq}_{A}(R), A^{\mrm{out}}}(a) \circ 
\opn{chm}^{\dcat{D}}_{\opn{Sq}_{A}(R), A^{\mrm{in}}}(a) 
\end{equation}
as endomorphisms of $\opn{Sq}_{A}(R)$ in 
$\dcat{D}(A^{\mrm{en}}) = \dcat{D}(A^{\mrm{en, in}})$.
But because $R$ is a rigid dualizing complex,  Theorem \ref{dfn:4370}(2) 
says that 
$\opn{chm}^{\dcat{D}}_{R, A}(a) = \opn{chm}^{\dcat{D}}_{R, A^{\mrm{op}}}(a)$.
Therefore, in formula (\ref{eqn:4462}) we can replace the action of $a$ 
through $A^{\mrm{out}}$ with its action through $A^{\mrm{op, in}}$. This gives 
\begin{equation} \label{eqn:4464}
\opn{Sq}_{A} \bigl( \opn{chm}^{\dcat{D}}_{R, A}(a) \bigr) = 
\opn{chm}^{\dcat{D}}_{\opn{Sq}_{A}(R), A^{\mrm{op, in}}}(a) \circ 
\opn{chm}^{\dcat{D}}_{\opn{Sq}_{A}(R), A^{\mrm{in}}}(a) . 
\end{equation}
Next, using Lemma \ref{lem:4381}, with $M = N := R$ and $a' = a'' := a$, 
we obtain 
\[ \opn{chm}^{\dcat{D}}_{\opn{Sq}_{A}(R), A^{\mrm{op, in}}}(a) =
\opn{chm}^{\dcat{D}}_{\opn{Sq}_{A}(R), A^{\mrm{in}}}(a) . \]
Plugging this into (\ref{eqn:4464}) we get 
\[ \opn{Sq}_{A} \bigl( \opn{chm}^{\dcat{D}}_{R, A}(a) \bigr) = 
\opn{chm}^{\dcat{D}}_{\opn{Sq}_{A}(R), A^{\mrm{in}}}(a)^2 . \]
Finally, since 
$\opn{chm}^{\dcat{D}}_{\opn{Sq}_{A}(R), A^{\mrm{in}}}$ 
is a ring homomorphism, we end up with (\ref{eqn:4501}). 
\end{proof}

In Definition \ref{dfn:4470} we introduced derived pseudo-finite DG modules
over a DG ring $A$. Recall that a DG module $M \in \dcat{D}(A)$ is called 
derived pseudo-finite if it belongs to the saturated full triangulated 
subcategory of $\dcat{D}(A)$ generated by the pseudo-finite semi-free DG 
modules. When $A$ is a ring, the pseudo-finite semi-free DG $A$-modules are 
precisely the bounded above complexes of finite free $A$-modules. 

\begin{thm}[Uniqueness, \cite{VdB}, \cite{Ye4}] \label{thm:4380} 
Let $A$ be a noetherian ring, and assume $A$ is a derived pseudo-finite complex 
over $A^{\mrm{en}}$. Suppose $(R, \rho)$ is a NC rigid dualizing complex over 
$A$. Then $(R, \rho)$ is unique, up to a unique NC rigid isomorphism.
\end{thm}

\begin{proof}
Suppose $(R', \rho')$ is another NC rigid dualizing complex over $A$.
According to Theorem \ref{thm:3681} and Corollary \ref{cor:4365}, there are 
tilting complexes $T$ and $T'$ over $A$ such that 
\begin{equation} \label{eqn:4380}
R' \cong T \ot^{\mrm{L}}_{A} R \cong R \ot^{\mrm{L}}_{A} T'
\end{equation}
in $\dcat{D}(A^{\mrm{en}})$. 
We have the following sequence of isomorphisms in $\dcat{D}(A^{\mrm{en}})$~: 
\begin{equation} \label{eqn:4381}
\begin{aligned}
&
T \ot^{\mrm{L}}_{A} R \cong^{1} R' 
\cong^2 \opn{RHom}_{A^{\mrm{en, out}}}(A, R' \ot R')
\\ & \qquad 
\cong^{1} \opn{RHom}_{A^{\mrm{en, out}}} \bigl( A, (R \ot^{\mrm{L}}_{A} T') \ot 
(T \ot^{\mrm{L}}_{A} R) \bigr)
\\ & \qquad 
\cong^3 \opn{RHom}_{A^{\mrm{en, out}}} \bigl( A, (R \ot R) 
\ot^{\mrm{L}}_{A^{\mrm{en}}} (T' \ot T) \bigr) 
\\ & \qquad 
\cong^4 \opn{RHom}_{A^{\mrm{en, out}}} ( A, R \ot R) 
\ot^{\mrm{L}}_{A^{\mrm{en}}} (T' \ot T)  
\\ & \qquad 
\cong^5 R \ot^{\mrm{L}}_{A^{\mrm{en}}} (T' \ot T) 
\cong^6 T \ot^{\mrm{L}}_{A} R \ot^{\mrm{L}}_{A} T' .
\end{aligned}
\end{equation}
Here are the explanations for the various isomorphisms:
\begin{itemize}
\item[$\cong^{1}$~:] These are due to (\ref{eqn:4380}). 

\item[$\cong^{2}$~:] This is the rigidifying isomorphism $\rho'$.

\item[$\cong^{3}$~:] This isomorphism is by rearranging the tensor factors, and 
in it the inside action of $A^{\mrm{en}}$ on $R \ot R$ is viewed as a right 
action, and this matches the left outside action of $A^{\mrm{en}}$ on 
$T' \ot T$. Cf.\ formulas (\ref{eqn:4877}) and (\ref{eqn:4215}).

\item[$\cong^{4}$~:] It is an application of the tensor-evaluation 
isomorphism (see Theorem \ref{thm:4320}). This is where we need $A$ to be a 
derived pseudo-finite complex over $A^{\mrm{en}}$ -- so that 
condition (i) of Theorem \ref{thm:4320} will hold. 
Condition (ii) that theorem holds since the complex $R \ot R$ 
has bounded cohomology. As for condition (iii) of the theorem: 
the complex $T' \ot T$ is a tilting complex over the ring $A^{\mrm{en}}$,
so it has finite flat dimension, and thus it certainly has derived bounded 
below tensor displacement. 

\item[$\cong^{5}$~:] This uses the rigidifying isomorphism $\rho$.

\item[$\cong^{6}$~:] It is a rearranging the tensor factors. Recall 
that the left action of $A^{\mrm{en}}$ on $T' \ot T$ is the outside action, 
cf.\ formula(\ref{eqn:4215}).
\end{itemize}

Let $T^{\vee}$ be the quasi-inverse of $T$. 
After applying $T^{\vee} \ot^{\mrm{L}}_{A} (-)$
to the isomorphisms (\ref{eqn:4381}) we get 
$R \cong R \ot^{\mrm{L}}_{A} T'$
in $\dcat{D}(A^{\mrm{en}})$. 
Equation (\ref{eqn:4380}) tells us that there is an isomorphism 
$\phi^{\dag} : R \iso R'$ in $\dcat{D}(A^{\mrm{en}})$. 

The isomorphism $\phi^{\dag}$ above need not be rigid. 
What we do know is that both morphisms 
$\rho' \circ \phi^{\dag}, \, \opn{Sq}_{A}(\phi^{\dag}) \circ \rho : 
R \to \opn{Sq}_{A}(R')$ 
are isomorphisms in $\dcat{D}(A^{\mrm{en}})$.
By Theorem \ref{thm:4375} the automorphisms of $R$ in $\dcat{D}(A^{\mrm{en}})$ 
are all of the form $\opn{chm}^{\dcat{D}}_{R, A}(a)$, 
for elements $a \in \opn{Cent}(A)^{\times}$.
Thus there is a unique invertible central element $a \in A$ such that  
\begin{equation} \label{eqn:4469}
\opn{Sq}_{A}(\phi^{\dag}) \circ \rho = 
\rho' \circ \phi^{\dag} \circ \opn{chm}^{\dcat{D}}_{R, A}(a) . 
\end{equation}
Define the isomorphism 
\[ \phi := \phi^{\dag} \circ \opn{chm}^{\dcat{D}}_{R, A}(a^{-1})
: R \iso R' \]
in $\dcat{D}(A^{\mrm{en}})$. 
Then we have the following equalities:
\[ \begin{aligned}
& 
\opn{Sq}_{A}(\phi) \circ \rho = 
\opn{Sq}_{A} \bigl( \phi^{\dag} \circ \opn{chm}^{\dcat{D}}_{R, A}(a^{-1}) 
\bigr) \circ \rho
\\ & \quad 
=^{\mrm{(i)}} \opn{Sq}_{A}(\phi^{\dag}) \circ 
\opn{chm}^{\dcat{D}}_{\opn{Sq}_{A}(R), A}(a^{-2}) \circ \rho
\\ & \quad 
=^{\mrm{(ii)}} \opn{Sq}_{A}(\phi^{\dag}) \circ \rho \circ 
\opn{chm}^{\dcat{D}}_{R, A}(a^{-2})
\\ & \quad 
=^{\mrm{(iii)}} \rho' \circ \phi^{\dag} \circ \opn{chm}^{\dcat{D}}_{R, A}(a)
\circ \opn{chm}^{\dcat{D}}_{R, A}(a^{-2})
\\ & \quad 
=^{\mrm{(iv)}} \rho' \circ \phi^{\dag} \circ \opn{chm}^{\dcat{D}}_{R, A}(a^{-1})
= \rho' \circ \phi . 
\end{aligned} \]
The equality $=^{\mrm{(i)}}$ is due to Lemma \ref{lem:4384}. In 
the equality $=^{\mrm{(ii)}}$ we have used the fact that $\rho$
is a $\opn{Cent}(A)$-linear morphism for the action through $A$, so $\rho$ 
commutes with $\opn{chm}^{\dcat{D}}_{R, A}(a^{-2})$.
Equality $=^{\mrm{(iii)}}$ is from (\ref{eqn:4469}). 
And equality $=^{\mrm{(iv)}}$ is because $\opn{chm}^{\dcat{D}}_{R, A}$ is a 
ring homomorphism. 
We see that $\phi : R \iso R'$ is a rigid isomorphism. The uniqueness of the 
element $a$ implies the uniqueness of the rigid isomorphism $\phi$ (by the 
same sort of calculation). 
\end{proof}

\begin{rem} \label{rem:4455}
As mentioned in Remark \ref{rem:3691}, for the purposes of Subsection 
\ref{subsec:NC-DC}, the only requirement on the base ring $\K$ is that the 
central $\K$-ring $A$ is flat over it. 

In this subsection that is not enough. Suppose that $A$ is flat over $\K$.
In order to prove the uniqueness of a 
rigid dualizing complex $R \in \dcat{D}^{\mrm{b}}(A^{\mrm{en}})$,
utilizing Theorem \ref{thm:4320}, we need that the complex 
$R \ot^{\mrm{L}}_{\K} R \in \dcat{D}(A^{\mrm{en}})$ 
will have bounded below cohomology. The only way we know to guarantee this is 
if the base ring $\K$ is {\em regular}. Notice that the regularity requirement 
is also needed in the commutative theory (see Setup \ref{set:3210} and Remark 
\ref{rem:4192}). 

Presumably the regularity of the base ring $\K$ is the only obstruction.
The lack of flatness of $A$ can be handled using a suitable DG ring resolution 
$\til{A} \to A$, as explained in Remark \ref{rem:3691}.
Thus we predict that the definitions and results of this subsection will remain 
valid if $\K$ is a regular nonzero commutative base ring, and $A$ is a 
noetherian central $\K$-ring. 
\end{rem}

\mysubsection{Interlude: Graded Rings of Laurent Type} 
\label{subsec:str-gr-rings}

Algebraically graded rings were introduced in Section 
\ref{sec:alg-gra-rings}, 
and we now return to them, and refer to them simply as graded rings. 
In this section we specialize to a particular kind of graded ring, that is 
needed for the next two subsections. 
The two main results here are Theorems \ref{thm:4510} and \ref{thm:4511}. 

We continue with Convention \ref{conv:3670}; in particular, $\K$ is a base 
field, and all rings are central over $\K$.

\begin{dfn} \label{dfn:4510}
A graded ring $\til{B}$ is called a {\em graded ring of Laurent type} 
\index{Algebraically graded ring! of Laurent type}
if there is an invertible central element $\til{c}$ in $\til{B}$ of degree $1$. 
Such an element $\til{c}$ is called a {\em uniformizer} of $\til{B}$. 
\end{dfn}

Of course if $\til{B}$ is a graded ring of Laurent type, then 
$\til{B}_i \cd \til{B}_j = \til{B}_{i + j}$
for all $i, j \in \Z$. A graded ring with this property is called a 
{\em strongly graded ring}. 

The next lemma describes the structure of graded rings of Laurent type. 
Let $\K[t, t^{-1}]$ be the ring of Laurent polynomials in the degree $1$ 
variable $t$. Given a ring $B$ (not graded), we define the graded ring of 
Laurent type 
$B[t, t^{-1}] := B \ot_{\K} \K[t, t^{-1}]$.
The uniformizer is $t$, and the degree $0$ component is $B$.  

\begin{lem} \label{lem:4510}
Let $\til{B}$ be a graded ring of Laurent type, with uniformizer $\til{c}$, and 
define $B := \til{B}_0$, the degree $0$ subring.
Then there is a unique isomorphism of graded rings 
$B[t, t^{-1}] \iso \til{B}$ 
that is the identity on $B$, and  sends $t \mapsto \til{c}$.
\end{lem}

\begin{exer} \label{exer:4511}
Prove Lemma \ref{lem:4510}.
\end{exer}

Note that for $\til{B}$ as in the lemma, there is also a surjective ring 
homomorphism $\til{B} \to B$ that sends $\til{c} \mapsto 1$. This is not a 
graded ring homomorphism. 

Recall that given a graded ring $\til{B}$, the category of graded 
$\til{B}$-modules is $\dcat{M}(\til{B}, \mrm{gr})$.
For a graded module $\til{M} \in \dcat{M}(\til{B}, \mrm{gr})$,
let us write 
$\opn{Deg}_0(\til{M}) := \til{M}_0$,
the homogeneous component of degree $0$ of $\til{M}$. 

\begin{lem} \label{lem:4390}
Assume $\til{B}$ is a graded ring of Laurent type, with $B := \til{B}_0$. 
\begin{enumerate}
\item The functor 
$\opn{Deg}_0 : \dcat{M}(\til{B}, \mrm{gr}) \to \dcat{M}(B)$
is an equivalence of abelian categories, with quasi-inverse 
$N \mapsto \til{B} \ot_{B} N \cong \K[t, t^{-1}] \ot N$.

\item There is a bifunctorial isomorphism 
\[ \opn{Deg}_0 \bigl( \opn{Hom}_{\til{B}}(\til{M} , \til{N}) \bigr) \iso 
\opn{Hom}_{B} \bigl( \opn{Deg}_0(\til{M}), \opn{Deg}_0(\til{N}) \bigr) \]
for $\til{M} , \til{N} \in \dcat{M}(\til{B}, \mrm{gr})$. It sends 
$\til{\phi} \mapsto \til{\phi}|_{\opn{Deg}_0(\til{M})}$.

\item If the graded ring $\til{B}$ is left noetherian, then the ring 
$B$ is left noetherian, and the functor 
$\opn{Deg}_0 : \dcat{M}_{\mrm{f}}(\til{B}, \mrm{gr}) \to 
\dcat{M}_{\mrm{f}}(B)$
is an equivalence  of abelian categories. 
\end{enumerate}
\end{lem}

Since $\opn{Deg}_0$ is exact, it extends to derived categories. 

\begin{lem} \label{lem:4400}
Assume $\til{B}$ is a graded ring of Laurent type, with $B := \til{B}_0$. 
\begin{enumerate}
\item For every boundedness condition $\star$, the functor 
$\opn{Deg}_0 : \dcat{D}^{\star}(\til{B}, \mrm{gr}) \to \dcat{D}^{\star}(B)$
is an equivalence of triangulated categories. 

\item If the graded ring $\til{B}$ is left noetherian, then 
$\opn{Deg}_0 : \dcat{D}^{\star}_{\mrm{f}}(\til{B}, \mrm{gr}) \to 
\dcat{D}^{\star}_{\mrm{f}}(B)$
is an equivalence. 
\end{enumerate}
\end{lem}

\begin{exer} \label{exer:4445}
Prove Lemmas \ref{lem:4390} and \ref{lem:4400}.
\end{exer}

\begin{lem} \label{lem:4401}
Let $\til{B}$ and $\til{C}$ be graded rings, with 
$B := \til{B}_0$ and $C := \til{C}_0$. Assume $\til{B}$ is 
a graded ring of Laurent type.
For $\til{M}, \til{N} \in \dcat{D}(\til{B} \ot \til{C}^{\mrm{op}}, \mrm{gr})$ 
there is an isomorphism 
\[ \opn{Deg}_0 \bigl( \opn{RHom}_{\til{B}}(\til{M}, \til{N}) \bigr) \iso  
\opn{RHom}_{B} \bigl( \opn{Deg}_0(\til{M}), 
\opn{Deg}_0(\til{N}) \bigr) \]
in $\dcat{D}(C^{\mrm{en}})$, which is functorial in these complexes.
\end{lem}

\begin{proof}
Choose a semi-graded-free resolution $\til{P} \to \til{M}$ over 
$\til{B} \ot \til{C}^{\mrm{op}}$. 
Because $\til{C}$ is a graded-free $\K$-module, it follows that the complex 
$\til{P}$ is semi-graded-free over $\til{B}$. 
Since $B = \opn{Deg}_0(\til{B})$, the complex $\opn{Deg}_0(\til{P})$ is 
semi-free over the ring $B$; so 
$\opn{Deg}_0(\til{P}) \to \opn{Deg}_0(\til{M})$
is a semi-free resolution over $B$. 
We get these isomorphisms: 
\[ \begin{aligned}
&
\opn{Deg}_0 \bigl( \opn{RHom}_{\til{B}}(\til{M}, \til{N}) \bigr) 
\cong 
\opn{Deg}_0 \bigl( \opn{Hom}_{\til{B}}(\til{P}, \til{N}) \bigr) 
\\ & \quad 
\cong^{\dag} \opn{Hom}_{B} \bigl( \opn{Deg}_0(\til{P}), \opn{Deg}_0(\til{N}) 
\bigr)
\cong \opn{RHom}_{B} \bigl( \opn{Deg}_0(\til{M}), \opn{Deg}_0(\til{N}) \bigr) 
\end{aligned} \]
in $\dcat{D}(C^{\mrm{en}})$. 
The isomorphism $\cong^{\dag}$ comes from Lemma \ref{lem:4390}(2). 
\end{proof}

\begin{lem} \label{lem:4505}
Let $\til{B}$ be a graded ring, and let $\til{c}$ be a central homogeneous 
element of $\til{B}$ of positive degree. 
Then the central element $\til{c} - 1$ is regular in $\til{B}$.
\end{lem}

\begin{proof}
Take a nonzero element $\til{b} \in \til{B}$,
with homogeneous component decomposition 
$\til{b} = \til{b}_{1} + \cdots + \til{b}_{r}$, 
see (\ref{eqn:4260}). So either $r = 1$ and $\til{b} = \til{b}_{1}$
is homogeneous, or $r \geq 2$, and then
$\opn{deg}(\til{b}_{1}) < \opn{deg}(\til{b}_{i})$ for all 
$2 \leq i \leq r$. Then the homogeneous component decomposition
of $\til{b} \cd (\til{c} - 1)$ is 
\[ \til{b} \cd (\til{c} - 1) = - \til{b}_{1} + 
(\tup{higher degree terms}) . \]
The element $-\til{b}_{1}$ is nonzero, and therefore 
$\til{b} \cd (\til{c} - 1) \neq 0$. 
\end{proof}

Derived pseudo-finite complexes were introduced in Definition 
\ref{dfn:4470}.

\begin{lem} \label{lem:4406}
Let $A, B, C$ be rings; let 
$M_1, M_2 \in \dcat{D}(A)$; let $N_1 \in \dcat{D}(B)$; and let 
$N_2 \in \dcat{D}(B \ot C)$.
There is a morphism 
\[ \th : \opn{RHom}_{A}(M_1, M_2) \ot \opn{RHom}_{B}(N_1, N_2) \to 
\opn{RHom}_{A \ot B}(M_1 \ot N_1, M_2 \ot N_2) \]
in $\dcat{D}(C)$, which is functorial in these complexes. 

If $M_1$ is derived pseudo-finite over $A$, $N_1$ is derived pseudo-finite 
over $B$, and the complexes $M_2$ and $N_2$ have bounded below cohomologies, 
then $\th$ is an isomorphism. 
\end{lem}

\begin{proof}
Choose semi-free resolutions $P_1 \to M_1$ and $Q_1 \to N_1$, over $A$ and $B$ 
respectively. Then 
$P_1 \ot Q_1 \to M_1 \ot N_1$
is a semi-free resolution over $A\ot B$; see Proposition \ref{prop:4636}.
The morphism $\th$ is represented by the obvious homomorphism 
\[ \til{\th} : \opn{Hom}_{A}(P_1, M_2) \ot \opn{Hom}_{B}(Q_1, N_2) \to 
\opn{Hom}_{A \ot B}(P_1 \ot Q_1, M_2 \ot N_2) \]
in $\dcat{C}_{\mrm{str}}(C)$. 

Now assume that $\opn{H}(M_2)$ and $\opn{H}(N_2)$ are bounded below. By 
replacing  $M_2$ and $N_2$ with suitable smart truncations, we can assume these
are bounded below complexes. We fix them.

If $P_1$ and $Q_1$ are pseudo-finite semi-free, then $\til{\th}$ is an 
isomorphism in $\dcat{C}_{\mrm{str}}(C)$, by the usual calculation; 
see Lemma \ref{lem:4735}. 

If we fix $N_1$, then the complexes $M_1$ for which $\th$ is an isomorphism 
form a saturated full triangulated subcategory of $\dcat{D}(A)$. 
Likewise, if we fix $M_1$, then the complexes $N_1$ for which $\th$ is an 
isomorphism form a saturated full triangulated subcategory of 
$\dcat{D}(B)$. See Proposition \ref{prop:4603}. 
Therefore, as in the proof of Theorem \ref{thm:4320},
$\th$ is an isomorphism whenever $M_1$ is derived pseudo-finite over $A$ and 
$N_1$ is derived pseudo-finite over $B$.
\end{proof}

The inside and outside actions from Subsection \ref{subsec:RNCDC-uniq} make 
sense also in the graded setting. 
Here is the first main result of this subsection. 

\begin{thm} \label{thm:4510} 
Let $\til{B}$ be a graded ring of Laurent type, with $B := \til{B}_0$.
Assume that $B$ is a derived pseudo-finite complex over $B^{\mrm{en}}$. 
Given complexes 
$\til{M}, \til{N} \in \dcat{D}^{+}(\til{B}^{\mrm{en}},  \mrm{gr})$,
there is an isomorphism 
\[ \begin{aligned}
&
\opn{Deg}_0 \bigl( \opn{RHom}_{\til{B}^{\mrm{en, out}}}
(\til{B}, \til{M} \ot \til{N}) \bigr) 
\\ & \quad 
\iso \opn{RHom}_{B^{\mrm{en, out}}} \bigl( B, 
\opn{Deg}_0(\til{M}) \ot \opn{Deg}_0(\til{N}) \bigr)[-1]
\end{aligned} \]
in $\dcat{D}(B^{\mrm{en, in}}) = \dcat{D}(B^{\mrm{en}})$, 
which is functorial in these complexes.
\end{thm}

\begin{proof}
Say $\til{c}$ is a uniformizer of $\til{B}$. 
The ring $\til{B}^{\mrm{en}}  = \til{B} \ot \til{B}^{\mrm{op}}$ is actually
$\Z^2$-graded. Taking the degree $0$ component is for the total degree,
so the ring 
$\opn{Deg}_0(\til{B}^{\mrm{en}})$
retains a $\Z$-grading, that we shall call the {\em hidden grading}. For the 
hidden grading there is a graded ring isomorphism 
\begin{equation} \label{eqn:4440}
B^{\mrm{en}}[s, s^{-1}] = B^{\mrm{en}} \ot \K[s, s^{-1}] \iso 
\opn{Deg}_0(\til{B}^{\mrm{en}}) , 
\end{equation}
where $s$ is a variable of hidden degree $1$, and it goes to the element 
$\til{c} \ot \til{c}^{-1} \in \opn{Deg}_0(\til{B}^{\mrm{en}})$.

Similarly for complexes: given 
$\til{M}, \til{N} \in \dcat{D}(\til{B}^{\mrm{en}},  \mrm{gr})$,
there is a canonical isomorphism 
\begin{equation} \label{eqn:4441}
\opn{Deg}_0(\til{M}) \ot \opn{Deg}_0(\til{N}) \ot \K[s, s^{-1}] \iso 
\opn{Deg}_0(\til{M} \ot \til{N})
\end{equation}
in $\dcat{D} \bigl( B^{\mrm{four}}[s, s^{-1}] \bigr)$. 
We are neglecting the hidden grading from this stage onward in the proof. 

According to Lemma \ref{lem:4401} -- that applies because 
$\til{B}^{\mrm{en}}$ is a graded ring of Laurent type -- 
there is an isomorphism 
\begin{equation} \label{eqn:4442}
\opn{Deg}_0 \bigl( \opn{RHom}_{\til{B}^{\mrm{en, out}}}
(\til{B}, \til{M} \ot \til{N}) \bigr) \cong 
\opn{RHom}_{\opn{Deg}_0(\til{B}^{\mrm{en, out}})}
\bigl( B, \opn{Deg}_0(\til{M} \ot \til{N}) \bigr) 
\end{equation}
in $\dcat{D} \bigl( \opn{Deg}_0(\til{B}^{\mrm{en, in}}) \bigr)$.
After applying the restriction functor 
$\dcat{D} \bigl( \opn{Deg}_0(\til{B}^{\mrm{en, in}}) \bigr) \to 
\dcat{D}(B^{\mrm{en, in}}) =  \dcat{D}(B^{\mrm{en}})$,
(\ref{eqn:4442}) becomes an isomorphism in $\dcat{D}(B^{\mrm{en}})$. 
We then have the following isomorphisms in $\dcat{D}(B^{\mrm{en}})$~: 
\begin{equation} \label{eqn:4443}
\begin{aligned}
& 
\opn{RHom}_{\opn{Deg}_0(\til{B}^{\mrm{en, out}})}
\bigl( B, \opn{Deg}_0(\til{M} \ot \til{N}) \bigr)
\\ & \quad 
\cong^{\dag} \opn{RHom}_{B^{\mrm{en, out}} \ot \, \K[s, s^{-1}]}
\Bigl( B \ot \K, 
\bigl( \opn{Deg}_0(\til{M}) \ot \opn{Deg}_0(\til{N}) \bigr) \ot \K[s, s^{-1}] 
\Bigr)
\\ & \quad 
\cong^{\ddag} \opn{RHom}_{B^{\mrm{en, out}}} \bigl( B, 
\opn{Deg}_0(\til{M}) \ot \opn{Deg}_0(\til{N}) \bigr) \ot 
\opn{RHom}_{\K[s, s^{-1}]} \bigl( \K, \K[s, s^{-1}] \bigr) . 
\end{aligned}
\end{equation}
The isomorphism $\cong^{\dag}$ comes from formulas 
(\ref{eqn:4440}) and (\ref{eqn:4441}).
The isomorphism $\cong^{\ddag}$ is by Lemma \ref{lem:4406}; it 
justified because $B$ is a derived pseudo-finite complex over $B^{\mrm{en}}$, 
and of course $\K$ is a derived pseudo-finite complex over the noetherian ring 
$\K[s, s^{-1}]$.

The element $s - 1 \in \K[s, s^{-1}]$ is regular, so we 
have the Koszul resolution 
\[ 0 \to \K[s, s^{-1}] \xar{(s - 1) \cd (-)} \K[s, s^{-1}] \to \K \to 0 . \]
We use it to calculate
\[ \opn{RHom}_{\K[s, s^{-1}]} \bigl( \K, \K[s, s^{-1}] \bigr) \cong \K[-1] \]
in $\dcat{D}(\K)$. 
Plugging this into (\ref{eqn:4443}) finishes the proof. 
\end{proof}

From here until the end of this subsection we
consider a connected graded ring $\til{A}$, with a degree $1$
central element $\til{c} \in \til{A}$. Define the ring
$A := \til{A} / (\til{c} - 1)$,
and the canonical ring surjection $f : \til{A} \to A$.
Of course the ring $A$ is not graded. 

Let $\til{A}_{\til{c}}$ be the localization  of $\til{A}$ with respect to 
$\til{c}$, i.e.\ inverting the homogeneous multiplicatively closed set 
$\{ \til{c}^{\, i} \}_{i \in \N}$. 
Because the image of $\til{c}$ in $\til{A}_{\til{c}}$ is an invertible central 
element of degree $1$, the graded ring $\til{A}_{\til{c}}$ is a
graded ring of Laurent type, with uniformizer $\til{c}$. We let 
$\opn{inc} : \opn{Deg}_0(\til{A}_{\til{c}}) = (\til{A}_{\til{c}})_0 \to 
\til{A}_{\til{c}}$ 
be the inclusion of the degree $0$ component of $\til{A}_{\til{c}}$.
There is also the ring homomorphism 
$f_{\til{c}} : \til{A}_{\til{c}} \to A$, 
$f_{\til{c}}(\til{c}) = 1$. 

\begin{lem} \label{lem:4385}
\mbox{}
\begin{enumerate}
\item The ring homomorphism 
\[ g := f_{\til{c}} \circ \opn{inc} : 
\opn{Deg}_0(\til{A}_{\til{c}}) = (\til{A}_{\til{c}})_0 \to A \]
is an isomorphism. 

\item The ring isomorphism $g$ extends to a graded ring isomorphism 
$\til{g} : \til{A}_{\til{c}} \iso A[t, t^{-1}]$
that sends $\til{c} \mapsto t$. 
\end{enumerate}
\end{lem}

The situation is shown in this commutative diagram of rings: 
\[ \UseTips \xymatrix @C=8ex @R=6ex {
\til{A}
\ar@{->>}@(ur,ul)[rr]^{f}
\ar[r]
&
\til{A}_{\til{c}}
\ar@{->>}[r]^{f_{\til{c}}}
&
*+{A}
\\
&
*++{(\til{A}_{\til{c}})_0}
\ar@{>->}[u]^{\mrm{inc}}
\ar@{>->>}[ur]_{g}
} \]

\begin{proof}
In this proof we work with the graded ring of Laurent type 
$\til{B} := \til{A}_{\til{c}}$ and its uniformizer $\til{c}$.

\medskip \noindent 
(1) The homomorphism $f : \til{A} \to A$ is surjective, so every nonzero 
$a \in A$ can be written as 
$a = f(\til{b})$ for some nonzero $\til{b} \in \til{A}$. 
Consider the homogeneous component decomposition 
\begin{equation}  \label{eqn:4867}
\til{b} = \til{b}_{1} + \cdots + \til{b}_{r}
\end{equation}
of $\til{b}$, as in (\ref{eqn:4260}), with $r \geq 1$. Define 
$a' := \sum_{i = 1}^{r} \, \til{b}_i \cdot 
\til{c}^{\, - \opn{deg}(\til{b}_i)} \in (\til{A}_{\til{c}})_0$.
Then $a = g(a')$. We see that $g$ is surjective.

Next we look at 
\begin{equation}  \label{eqn:4868}
\opn{Ker}(g) = \opn{Ker}(f_{\til{c}}) \cap (\til{A}_{\til{c}})_0 
\sub (\til{A}_{\til{c}})_0 . 
\end{equation}
The ideal $\opn{Ker}(f_{\til{c}})$ is generated by the central element 
$\til{c} - 1$.
Consider a nonzero element $\til{b} \cd (\til{c} - 1) \in \til{A}_{\til{c}}$.
Let (\ref{eqn:4867}) be the homogeneous component decomposition of the nonzero 
element $\til{b}$. Then the homogeneous component decomposition
of $\til{b} \cd (\til{c} - 1)$ is 
$\til{b} \cd (\til{c} - 1) = -\til{b}_{1} + \cdots + \til{b}_{r} \cd \til{c}$.
This is not a homogeneous element, and hence it cannot lie inside 
$(\til{A}_{\til{c}})_0$. From formula (\ref{eqn:4868}) we conclude that 
$\opn{Ker}(g) = 0$. 
Thus $g$ is a ring isomorphism. 

\medskip \noindent 
(2) Let 
\begin{equation} \label{eqn:4476}
h : A \iso (\til{A}_{\til{c}})_0
\end{equation}
be the inverse of $g$. It extends to a graded ring homomorphism
$\til{h} : A[t, t^{-1}] \to  \til{A}_{\til{c}}$,
$t \mapsto \til{c}$. 
This is an isomorphism, because $\til{A}_{\til{c}}$ is a graded ring of Laurent 
type. The isomorphism $\til{g}$ is the inverse of $\til{h}$.
\end{proof}

\begin{exer} \label{exer:4475}
Show that the functor
$\opn{Ind}_f : \dcat{M}(\til{A}, \mrm{gr}) \to \dcat{M}(A)$,
$\opn{Ind}_f = A \ot_{\til{A}} (-)$, 
is exact.
\end{exer}

Here is the second main result of the subsection. 
Connected graded $\K$-rings were defined in Definition \ref{dfn:3703}. 

\begin{thm} \label{thm:4511}   
Let $\til{A}$ be a noetherian connected graded ring, let
$\til{c} \in \til{A}$ be a homogeneous central element of degree $1$, and let 
$A := \til{A} / (\til{c} - 1)$. Then\tup{:}
\begin{enumerate}
\item The ring $A$ is noetherian.

\item There is an isomorphism 
$A^{\mrm{en}} \ot^{\mrm{L}}_{\til{A}^{\mrm{en}}} \til{A} \cong A \oplus A[1]$
in $\dcat{D}(A^{\mrm{en}})$. 

\item The bimodule $A$ is a derived pseudo-finite complex over the ring 
$A^{\mrm{en}}$.
\end{enumerate}
\end{thm}

\begin{proof} \mbox{}

\smallskip \noindent
(1) By Theorem \ref{thm:4245} the ring $\opn{Ungr}(\til{A})$ is noetherian. 
Since there is a surjection of rings $\opn{Ungr}(\til{A}) \to A$, it follows 
that $A$ is noetherian. 

\medskip \noindent 
(2) We have these isomorphisms 
\[ \begin{aligned}
&
A^{\mrm{en}} \ot^{\mrm{L}}_{\til{A}^{\mrm{en}}} \til{A}
\cong^1 A \ot^{\mrm{L}}_{\til{A}} \til{A} \ot^{\mrm{L}}_{\til{A}} A
\\[0.3em] & \quad 
\cong^2 A \ot^{\mrm{L}}_{\til{A}} A
\cong^3 A \ot^{\mrm{L}}_{\til{A}} \til{A}_{\til{c}} 
\ot^{\mrm{L}}_{\til{A}_{\til{c}}} \! A 
\cong^4 A \ot^{\mrm{L}}_{\til{A}_{\til{c}}} \! A 
\end{aligned} \]
in $\dcat{D}(A^{\mrm{en}})$. 
The isomorphism $\cong^1$ is a rearrangement of the derived tensor factors. 
The isomorphism $\cong^2$ is from the left unitor isomorphism
$\til{A} \ot^{\mrm{L}}_{\til{A}} A \cong A$ in 
$\dcat{D}(\til{A} \ot A^{\mrm{op}})$. The isomorphism $\cong^3$
is because the ring homomorphism $\til{A} \to A$ factors through  
$\til{A}_{\til{c}}$, and 
$A \cong \til{A}_{\til{c}} \ot^{\mrm{L}}_{\til{A}_{\til{c}}} \! A$
in $\dcat{D}(\til{A} \ot A^{\mrm{op}})$. And, finally, the  
isomorphism $\cong^4$ is because 
$A \ot^{\mrm{L}}_{\til{A}} \til{A}_{\til{c}} \cong A$ in 
$\dcat{D}(\til{A} \ot \til{A}_{\til{c}}^{\mrm{op}})$

Now by Lemma \ref{lem:4505} 
we know that ${\til{c}} - 1$ is a regular central element in 
$\til{A}_{\til{c}}$. So there is a short exact sequence 
\begin{equation} \label{eqn:4510}
0 \to \til{A}_{\til{c}} \xar{\, (\til{c} - 1) \cd (-) \, } 
\til{A}_{\til{c}} \xar{ \, f_{\til{c}} \, } A \to  0
\end{equation}
in $\dcat{M}(\til{A}_{\til{c}} \ot A^{\mrm{op}})$. 
Here $\til{A}_{\til{c}}$ is an $A^{\mrm{op}}$-module via the ring homomorphism 
$h : A \to \til{A}_{\til{c}}$ from formula (\ref{eqn:4476}). We view the short 
exact sequence (\ref{eqn:4510}) as a quasi-isomorphism $\til{P} \to A$ in 
$\dcat{C}_{\mrm{str}}(\til{A}_{\til{c}} \ot A^{\mrm{op}})$,
where 
\[ \til{P} := \bigl( \cdots \to 0 \to 
\til{A}_{\til{c}} \xar{\, (\til{c} - 1) \cd (-) \, } 
\til{A}_{\til{c}} \to  0 \to \cdots \bigr) \]
is a complex concentrated in cohomological degrees $-1, 0$. 
Because $\til{P}$ is semi-free over $\til{A}_{\til{c}}$, this allows us to 
calculate: 
\[ A \ot^{\mrm{L}}_{\til{A}_{\til{c}}} \! A \cong A 
\ot_{\til{A}_{\til{c}}} \! \til{P}  \cong 
\bigl( \cdots \to 0 \to A \xar{ \,0 \, } A \to 0 \to \cdots \bigr) \cong A[1] 
\oplus A \]
in $\dcat{D}(A^{\mrm{en}})$.

\medskip \noindent 
(3) By Proposition \ref{prop:4266} there is a quasi-isomorphism 
$\til{Q} \to \til{A}$ in 
$\dcat{C}_{\mrm{str}}(\til{A}^{\mrm{en}}, \mrm{gr})$
from a nonpositive complex $\til{Q}$ of finite graded-free 
$\til{A}^{\mrm{en}}$-modules. We can forget the grading now, and just view 
$\til{Q} \to \til{A}$ as a free resolution of the module $\til{A}$ over the 
ring $\til{A}^{\mrm{en}}$. We get 
$P := A^{\mrm{en}} \ot_{\til{A}^{\mrm{en}}} \til{Q}  \cong 
A^{\mrm{en}} \ot^{\mrm{L}}_{\til{A}^{\mrm{en}}} \til{A}$ 
in $\dcat{D}(A^{\mrm{en}})$. The complex $P$ is pseudo-finite semi-free over 
$A^{\mrm{en}}$.
By item (2) the complex $A$ is a direct summand of $P$
in $\dcat{D}(A^{\mrm{en}})$, and hence it is a derived 
pseudo-finite complex over $A^{\mrm{en}}$.
\end{proof}

\mysubsection{Graded Rigid NC DC}
\label{subsec:GrR-NC-DC}

Here we make a bridge between {\em balanced dualizing complexes} and 
{\em rigid dualizing complexes}. This is done by introducing an intermediate 
kind of object: the {\em graded rigid dualizing complex}. All is in the 
noncommutative setting of course. 

The results of this subsection were originally in \cite{VdB}, with a few 
improvements later in \cite{YeZh1}. Here we give a much more detailed 
discussion, and some corrections. 

We continue with Convention \ref{conv:3670}. 
In particular $\K$ is a base field, and all rings are by default central over 
$\K$. Throughout this subsection we assume the following setup: 

\begin{setup} \label{set:4510}
We are given a noetherian connected graded central $\K$-ring $\til{A}$ and
a central homogeneous element $t \in \til{A}$ of degree $1$. 
We define the central $\K$-ring  
$A := \til{A} / (t - 1)$.
\end{setup}

Note that the ring $A$ is not graded. Also the rings
$\til{A}^{\mrm{en}}$ and $A^{\mrm{en}}$ need not be 
noetherian. However: 

\begin{prop} \label{prop:4515}
Under Setup \tup{\ref{set:4510}}, the ring $A$ is noetherian, and $A$ is a 
derived pseudo-finite complex over $A^{\mrm{en}}$.
\end{prop}

\begin{proof} 
Use Theorem \ref{thm:4511}. 
\end{proof}

The concepts of outside and inside actions from the beginning of Subsection 
\ref{subsec:RNCDC-uniq} make sense in the graded setting. Thus we have a 
connected graded ring 
$\til{A}^{\mrm{four}} = \til{A}^{\mrm{en, out}} \ot \til{A}^{\mrm{en, in}}$. 
For a complex
$\til{M} \in \dcat{D}(\til{A}^{\mrm{en}}, \mrm{gr})$,
its tensor product with itself $\til{M} \ot \til{M}$
is an object of 
$\dcat{D}(\til{A}^{\mrm{four}}, \mrm{gr})$,
and 
\begin{equation} \label{eqn:4860}
\opn{RHom}_{\til{A}^{\mrm{en, out}}} ( \til{A}, \til{M} \ot \til{M} )
\in \dcat{D}(\til{A}^{\mrm{en, in}}, \mrm{gr}) = 
\dcat{D}(\til{A}^{\mrm{en}}, \mrm{gr}) . 
\end{equation}

Recall that 
$\dcat{D}_{(\mrm{f},..)}(\til{A}^{\mrm{en}}, \mrm{gr})$
(resp.\ $\dcat{D}_{(.., \mrm{f})}(\til{A}^{\mrm{en}}, \mrm{gr})$)
is the full subcategory of \lb 
$\dcat{D}(\til{A}^{\mrm{en}}, \mrm{gr})$
on the complexes whose cohomology bimodules are finite over $\til{A}$
(resp.\  $\til{A}^{\mrm{op}}$), and 
\[ \dcat{D}_{(\mrm{f}, \mrm{f})}(\til{A}^{\mrm{en}}, \mrm{gr}) = 
\dcat{D}_{(\mrm{f},..)}(\til{A}^{\mrm{en}}, \mrm{gr}) \cap 
\dcat{D}_{(.., \mrm{f})}(\til{A}^{\mrm{en}}, \mrm{gr}) . \]

\begin{dfn} \label{dfn:4391} 
Under Setup \tup{\ref{set:4510}}, a 
{\em graded rigid NC dualizing complex}%
\index{Dualizing complex! graded rigid NC} 
over $\til{A}$ is a pair $(\til{R}, \til{\rho})$, where 
$\til{R} \in \dcat{D}^{\mrm{b}}_{(\mrm{f}, \mrm{f})}
(\til{A}^{\mrm{en}}, \mrm{gr})$
is a graded NC dualizing complex over $\til{A}$, in the sense of Definition 
\ref{dfn:3716}; and 
\[ \til{\rho} : \til{R} \iso 
\opn{RHom}_{\til{A}^{\mrm{en, out}}} ( \til{A}, \til{R} \ot \til{R} ) \]
is an isomorphism in $\dcat{D}(\til{A}^{\mrm{en}}, \mrm{gr})$,
see (\ref{eqn:4860}) with $\til{M} := \til{R}$.  
\end{dfn}

Balanced dualizing complexes were introduced in Definition \ref{dfn:3714}.

\begin{thm} \label{thm:4400}
\index{Dualizing complex! graded rigid NC}
\index{Dualizing complex! balanced}
Under Setup \tup{\ref{set:4510}}, if $\til{R}$ is a balanced NC dualizing 
complex over $\til{A}$, then $\til{R}$ is a graded rigid NC dualizing complex 
over $\til{A}$.
\end{thm}

We are a bit sloppy in stating the theorem; the proper way to state it is this:
if $(\til{R}, \til{\be})$ is a balanced dualizing complex over $\til{A}$, then 
there exists an isomorphism $\til{\rho}$ such that $(\til{R}, \til{\rho})$ is a 
graded rigid NC dualizing complex over $\til{A}$.

We need a lemma first. 
Recall (from Definition \ref{dfn:3723}) that a graded $\til{A}$-module 
$\til{M} = \bigoplus_i \til{M}_i$ is called degreewise finite over $\K$ if each 
homogeneous component $\til{M}_i$ is a finite $\K$-module. We denoted by 
$\dcat{M}_{\mrm{dwf}}(\til{A}, \mrm{gr})$
the full subcategory of $\dcat{M}(\til{A}, \mrm{gr})$ on the degreewise finite
graded modules. This is a thick abelian subcategory of 
$\dcat{M}(\til{A}, \mrm{gr})$, closed under subobjects and quotients.
Then we denoted by 
$\dcat{D}_{\mrm{dwf}}(\til{A}, \mrm{gr})$
the full subcategory of $\dcat{D}(\til{A}, \mrm{gr})$
on the complexes $\til{M}$ such that 
$\opn{H}^q(\til{M}) \in \dcat{M}_{\mrm{dwf}}(\til{A}, \mrm{gr})$
for all $q$. This is a full triangulated subcategory.

\begin{lem} \label{lem:4402}
Let 
$\til{M}, \til{N} \in 
\dcat{D}^{\mrm{b}}_{(\mrm{f}, ..)}(\til{A}^{\mrm{en}}, \mrm{gr})$.
Then\tup{:}
\begin{enumerate}
\item The complexes $\til{M}$ and $\til{M} \ot \til{N}$ satisfy 
$\til{M} \in \dcat{D}^{\mrm{b}}_{\mrm{dwf}}(\til{A}^{\mrm{en}}, \mrm{gr})$
and 
$\til{M} \ot \til{N} \in 
\dcat{D}^{\mrm{b}}_{\mrm{dwf}}(\til{A}^{\mrm{four}}, \mrm{gr})$. 

\item The canonical morphism 
$\th : \til{M}^* \ot \til{N}^* \to (\til{M} \ot \til{N})^*$ 
in $\dcat{D}^{\mrm{b}}_{\mrm{dwf}}(\til{A}^{\mrm{four}}, \mrm{gr})$
is an isomorphism. 
\end{enumerate}
\end{lem}

\begin{proof} \mbox{}

\smallskip \noindent
(1) For every $q$ the $\til{A}$-module $\opn{H}^q(\til{M})$ is finite, so
it is the image of a finite direct sum of algebraic degree shifts of $\til{A}$. 
Since $\til{A} \in \dcat{M}_{\mrm{dwf}}(\K, \mrm{gr})$,
and since the subcategory
$\dcat{M}_{\mrm{dwf}}(\K, \mrm{gr})$
is closed under quotients inside $\dcat{M}(\K, \mrm{gr})$, we see that 
$\opn{H}^q(\til{M}) \in \dcat{M}_{\mrm{dwf}}(\K, \mrm{gr})$. 
We conclude that 
\[ \opn{H}^q(\til{M}) \in \dcat{M}_{\mrm{dwf}}(\K, \mrm{gr}) \cap 
\dcat{M}(\til{A}^{\mrm{en}}, \mrm{gr}) 
= \dcat{M}_{\mrm{dwf}}(\til{A}^{\mrm{en}}, \mrm{gr}) . \]

Let $[q_0, q_1]$ be a finite integer interval that contains 
$\opn{con}(\opn{H}(\til{M}))$. 
For every $p \in \Z$ there is an isomorphism 
\[ \opn{H}^p(\til{M} \ot \til{N}) \cong 
\bigoplus_{q \in [q_0, q_1]} \, 
\bigl( \opn{H}^q(\til{M}) \ot \opn{H}^{p - q}(\til{N}) \bigr) \]
in $\dcat{M}(\K, \mrm{gr})$. 
Each $\opn{H}^q(\til{M})$ and $\opn{H}^{p - q}(\til{N})$
is a finite graded $\til{A}$-module, so it is degreewise finite and bounded 
below (in algebraic degree). Therefore 
\[ \opn{H}^q(\til{M}) \ot \opn{H}^{p - q}(\til{N}) \in 
\dcat{M}_{\mrm{dwf}}(\K, \mrm{gr}) , \]
and hence so is $\opn{H}^p(\til{M} \ot \til{N})$. 
This means that 
$\til{M} \ot \til{N} \in 
\dcat{D}^{\mrm{b}}_{\mrm{dwf}}(\til{A}^{\mrm{four}}, \mrm{gr})$. 

\medskip \noindent 
(2) We start by giving an explicit formula for the canonical morphism $\th$. 
There is a canonical homomorphism 
\begin{equation} \label{eqn:4485}
\til{\th} : \til{M}^* \ot \til{N}^* \to (\til{M} \ot \til{N})^*
\end{equation}
in $\dcat{C}_{\mrm{str}}(A^{\mrm{four}}, \mrm{gr})$; 
it is 
\[ \til{\th}(\phi \ot \psi)(m \ot n) := 
(-1)^{j \cd k} \cd \phi(m) \ot \psi(n) \in \K \]
for $\phi \in (\til{M}^*)^i$, $\psi \in (\til{N}^*)^j$, 
$m \in \til{M}^k$ and $n \in \til{N}^l$.  
Then $\th = \opn{Q}(\til{\th})$ in 
$\dcat{D}(A^{\mrm{four}}, \mrm{gr})$.
We need to show that under the given finiteness and boundedness conditions, 
$\th$ is an isomorphism. For that we can forget the $A^{\mrm{four}}$-module 
structures, and view $\th$ as a morphism in 
$\dcat{D}(\K, \mrm{gr})$.

But in $\dcat{D}(\K, \mrm{gr})$ there are canonical isomorphisms 
$\til{M} \cong \opn{H}(\til{M})$ and 
$\til{N} \cong \opn{H}(\til{N})$.
This means that we can assume that the differentials of the complexes $\til{M}$ 
and $\til{N}$ are zero. We know that for every
$\sbmat{p \\[0.1em] i} \in \sbmat{\Z \\[0.1em] \Z} = \Z^2$,
in the notation from Subsection \ref{subsec:alg-gr-mods},
the $\K$-modules $\til{M}^p_i$ and  $\til{N}^p_i$ are finite. 
Also there are uniform bounds on vanishing: there are integers 
$p_0, p_1, i_0$ such that 
$\til{M}^p_i = \til{N}^p_i = 0$
unless $p_0 \leq p \leq p_1$ and $i_0 \leq i$. 
An easy calculation shows that for these objects $\til{M}$ and $\til{N}$, the 
homomorphism $\til{\th}$ of (\ref{eqn:4485}) is bijective. 
\end{proof}

\begin{proof}[Proof of Theorem \tup{\ref{thm:4400}}]
We need to produce a graded rigidifying isomorphism $\til{\rho}$ for the graded 
NC dualizing complex $\til{R}$. 

Consider the complex 
$\til{R} \in \dcat{D}^{\mrm{b}}_{(\mrm{f}, \mrm{f})}(\til{A}^{\mrm{en}}, 
\mrm{gr})$.
According to Lemma \ref{lem:4402}(1) we know that 
\begin{equation} \label{eqn:4410}
\til{R} \ot \til{R} \in 
\dcat{D}^{\mrm{b}}_{\mrm{dwf}}(\til{A}^{\mrm{four}}, \mrm{gr})  
\end{equation}
and 
\begin{equation} \label{eqn:4412}
\til{A} \in \dcat{D}^{\mrm{b}}_{\mrm{dwf}}(\til{A}^{\mrm{en}}, \mrm{gr}) . 
\end{equation}
There is a sequence of isomorphisms in 
$\dcat{D}(\til{A}^{\mrm{en}}, \mrm{gr})$~:
\begin{equation} \label{eqn:4413}
\begin{aligned}
& \opn{RHom}_{\til{A}^{\mrm{en, out}}}(\til{A}, \til{R} \ot \til{R})
\cong^{1}
\opn{RHom}_{\til{A}^{\mrm{en, out}}} \bigl( (\til{R} \ot \til{R})^*, \til{A}^* 
\bigr)
\\
& \qquad \cong^{2}
\opn{RHom}_{\til{A}^{\mrm{en, in}}}(\til{R}^* \ot \til{R}^*, \til{A}^*)
\\
& \qquad \cong^{3}
\opn{RHom}_{\til{A}^{\mrm{en, in}}}(\til{P} \ot \til{P}, \til{A}^*) 
\\
& \qquad \cong^{4}
\opn{RHom}_{\til{A}} \bigl( \til{P} , 
\opn{RHom}_{\til{A}^{\mrm{op}}} (\til{P}, \til{A}^*) \bigr) 
\\
& \qquad \cong^{5}
\opn{RHom}_{\til{A}} (\til{P}, \til{R} ) \cong^{6} \til{R} . 
\end{aligned}
\end{equation}
The explanations for the isomorphisms are:
\begin{itemize}
\item[$\cong^{1}$~:] By Theorem \ref{thm:3800}, that applies by formulas 
(\ref{eqn:4410}) and (\ref{eqn:4412}).  

\item[$\cong^{2}$~:] This is by Lemma \ref{lem:4402}(2). Note that the outside 
action changes to an inside action. 

\item[$\cong^{3}$~:] Here 
$\til{P} := \mrm{R} \Ga_{\til{\m}}(\til{A})$, where $\til{\m}$ is the 
augmentation ideal of $\til{A}$. According to Theorems \ref{thm:3713} and 
\ref{thm:3710} there is an isomorphism 
$\til{R} \cong \til{P}^*$
in $\dcat{D}(\til{A}^{\mrm{en}}, \mrm{gr})$; dualizing we get 
an isomorphism $\til{R}^* \cong \til{P}$ in 
$\dcat{D}(\til{A}^{\mrm{en}}, \mrm{gr})$. 

\item[$\cong^{4}$~:] This is one of the standard Hom-tensor identities. 

\item[$\cong^{5}$~:] This is another instance of the isomorphism 
$\til{P}^* \cong \til{R}$. 

\item[$\cong^{6}$~:] By Lemma \ref{lem:4060}, the restriction of 
$\til{R}$ to $\dcat{D}(\til{A}, \mrm{gr})$ satisfies  
$\til{R} \in \dcat{D}_{\mrm{f}}(\til{A}, \mrm{gr}) \lb 
\sub \dcat{D}(\til{A}, \mrm{gr})_{\mrm{com}}$,
so $\til{R}$ is isomorphic, in $\dcat{D}(\til{A}^{\mrm{en}}, \mrm{gr})$,
to its abstract derived completion 
$\opn{ADC}_{\til{\m}}(\til{R}) = \opn{RHom}_{\til{A}} (\til{P}, \til{R})$.
See Subsection \ref{subsec:NC-MGM}. 
\end{itemize}

The composition of the isomorphisms in (\ref{eqn:4413}) is the graded 
rigidifying isomorphism $\til{\rho}$. 
\end{proof}

The localization homomorphism $\til{A} \to \til{A}_t$
induces a flat homomorphism of graded rings  
$\til{A}^{\mrm{en}} \to (\til{A}_t)^{\mrm{en}} =
\til{A}_t \ot \til{A}_t^{\mrm{op}}$.
This homomorphism can also be viewed as a localization of 
$\til{A}^{\mrm{en}}$ with respect to the degree $1$ central elements 
$t \ot 1$ and $1 \ot t$. 
There is a corresponding induction functor 
\begin{equation} \label{eqn:4420}
\opn{Ind}_{(\til{A}_t)^{\mrm{en}}}
: \dcat{D}(\til{A}^{\mrm{en}}, \mrm{gr}) \to 
\dcat{D} \bigl( (\til{A}_t)^{\mrm{en}}, \mrm{gr} \bigr) .
\end{equation}
Here is what it does on objects: given a complex 
$\til{M} \in \dcat{D}(\til{A}^{\mrm{en}}, \mrm{gr})$,
we have 
\[ \opn{Ind}_{\til{A}^{\mrm{en}}}(\til{M}) = 
(\til{A}_t)^{\mrm{en}} \ot_{\til{A}^{\mrm{en}}} \til{M} \cong 
\til{A}_t \ot_{\til{A}} \til{M} \ot_{\til{A}} \til{A}_t . \]
There is a ring homomorphism  
\[ \opn{inc} : (\til{A}_t)^{\mrm{en}}_0 = 
\opn{Deg}_0 \bigl( (\til{A}_t)^{\mrm{en}} \bigr)
\to (\til{A}_t)^{\mrm{en}} , \]
the inclusion of the degree $0$ component. And there is a functor 
\begin{equation} \label{eqn:4421}
\opn{Deg}_0 : \dcat{D} \bigl( (\til{A}_t)^{\mrm{en}}, \mrm{gr} \bigr) \to 
\dcat{D} \bigl( (\til{A}_t)^{\mrm{en}}_0 \bigr) . 
\end{equation}

We have a ring isomorphism 
$h : A \iso (\til{A}_t)_0$; see formula (\ref{eqn:4476}). It gives rise to 
ring homomorphisms 
\[ A^{\mrm{en}} \xar[\cong]{h^{\mrm{en}}}
(\til{A}_t)_0 \ot (\til{A}_t^{\mrm{op}})_0 
\xar[\sub]{\opn{inc}} 
(\til{A}_t \ot \til{A}_t^{\mrm{op}})_0 = (\til{A}_t)^{\mrm{en}}_0 . \]
Relative to this composed ring homomorphism we have the restriction functor 
\begin{equation} \label{eqn:4416}
\opn{Rest}_{A^{\mrm{en}}} : 
\dcat{D} \bigl( (\til{A}_t)^{\mrm{en}}_0 \bigr) \to 
\dcat{D}(A^{\mrm{en}}) . 
\end{equation}
Composing the functors from (\ref{eqn:4420}), (\ref{eqn:4421}) and 
(\ref{eqn:4416}) we get a triangulated functor 
\begin{equation} \label{eqn:4422}
\opn{Rest}_{A^{\mrm{en}}} \circ \opn{Deg}_0 \circ 
\opn{Ind}_{(\til{A}_t)^{\mrm{en}}} :
\dcat{D}(\til{A}^{\mrm{en}}, \mrm{gr}) \to \dcat{D}(A^{\mrm{en}}) . 
\end{equation}

\begin{thm} \label{thm:4395}
Under Setup \tup{\ref{set:4510}}, let 
$\til{R} \in \dcat{D}(\til{A}^{\mrm{en}}, \mrm{gr})$ be a 
graded rigid NC dualizing complex. Define the complex 
\[ R := (\opn{Rest}_{A^{\mrm{en}}} \circ \opn{Deg}_0 \circ 
\opn{Ind}_{(\til{A}_t)^{\mrm{en}}}) (\til{R})[-1] \in \dcat{D}(A^{\mrm{en}}) . 
\]
Then $R$ is a rigid NC dualizing complex over $A$.
\end{thm}

Once more we were a bit sloppy. The proper statement is this: 
we are given a graded rigid NC dualizing complex $(\til{R}, \til{\rho})$; 
then the complex $R$ defined above admits a rigidifying isomorphism $\rho$, 
such that $(R, \rho)$ is a rigid NC dualizing complex over $A$. The proof will 
come after two lemmas.  

\begin{lem} \label{lem:4393}
Assume that $\til{R}$ is a graded rigid dualizing complex over $\til{A}$.
Then for every $p$ the graded bimodule $\opn{H}^p(\til{R})$ is 
central over $\opn{Cent}(\til{A})$. 
\end{lem}

\begin{proof}
We know that the central homotheties 
\[ \opn{chm}^{}_{\til{A}, \til{A}} ,  \, 
\opn{chm}^{}_{\til{A}, \til{A}^{\mrm{op}}} :
\opn{Cent}(\til{A}) \to 
\opn{End}_{\dcat{M}(\til{A}^{\mrm{en}}, \mrm{gr})}(\til{A}) \]
are equal graded ring isomorphism.
Therefore the graded version of Theorem \ref{thm:4375} -- that is proved the 
same way -- tells us that the ring homomorphisms
\[ \opn{chm}^{\dcat{D}}_{R, \til{A}} , \,
\opn{chm}^{\dcat{D}}_{R, \til{A}^{\mrm{op}}} :
\opn{Cent}(\til{A}) 
\to \opn{End}_{\dcat{D}(\til{A}^{\mrm{en}}, \mrm{gr})}(\til{R})  \]
are equal isomorphisms.
  
The left $\opn{Cent}(\til{A})$-module structure on $\opn{H}^p(\til{R})$
coincides with the categorical action $\opn{chm}^{\dcat{D}}_{R, \til{A}}$; 
and likewise from the right. We see that the left and right 
$\opn{Cent}(\til{A})$-module structure on $\opn{H}^p(\til{R})$
are the same, so it is a central bimodule. 
\end{proof}

\begin{lem} \label{lem:4395}
Assume that $\til{R}$ is a graded rigid dualizing complex over $\til{A}$.
Then the homomorphisms 
\[ \tag{l} \til{A}_t \ot_{\til{A}} \til{R} \to 
\til{A}_t \ot_{\til{A}} \til{R} \ot_{\til{A}} \til{A}_t \]
and
\[ \tag{r} \til{R} \ot_{\til{A}} \til{A}_t \to 
\til{A}_t \ot_{\til{A}} \til{R} \ot_{\til{A}} \til{A}_t \]
in $\dcat{C}_{\mrm{str}}(\til{A}^{\mrm{en}}, \mrm{gr})$
are quasi-isomorphisms. 
\end{lem}

\begin{proof}
We shall only treat (l); the homomorphism (r) is treated similarly. 
Let $\til{C} := \opn{Cent}(\til{A})$. Then 
$\til{A}_t = \til{A} \ot_{\til{C}} \til{C}_t = 
\til{C}_t \ot_{\til{C}} \til{A}$.
Hence (l) is true if (and only if) the homomorphism
$\til{C}_t \ot_{\til{C}} \til{R} \to 
\til{C}_t \ot_{\til{C}} \til{R} \ot_{\til{C}} \til{C}_t$
is a quasi-isomorphism. 

Because of flatness of localization, for every $p$ there is a commutative 
diagram 
\begin{equation} \label{eqn:5127}
\UseTips \xymatrix @C=6ex @R=6ex {
\opn{H}^p \bigl( \til{C}_t \ot_{\til{C}} \til{R} \bigr)
\ar[r]
\ar[d]_{\cong}
&
\opn{H}^p \bigl( \til{C}_t \ot_{\til{C}} \til{R} \ot_{\til{C}} \til{C}_t \bigr)
\ar[d]_{\cong}
\\
\til{C}_t \ot_{\til{C}} \opn{H}^p (\til{R})
\ar[r]
&
\til{C}_t \ot_{\til{C}} \opn{H}^p (\til{R}) \ot_{\til{C}} \til{C}_t 
} 
\end{equation}
in $\dcat{M}(\til{A}^{\mrm{en}}, \mrm{gr})$, with vertical isomorphisms.
By Lemma \ref{lem:4393}, $\opn{H}^p (\til{R})$ is a central $\til{C}$ 
bimodule. And localization satisfies 
$\til{C}_t = \til{C}_t \ot_{\til{C}} \til{C}_t$.
This implies that the bottom arrow in diagram (\ref{eqn:5127}) is bijective. 
Hence the top arrow is also bijective. 
\end{proof}

\begin{proof}[Proof of Theorem \tup{\ref{thm:4395}}]
The proof is divided into several steps. 

\smallskip \noindent
Step 1. This is the easiest step: we prove that for every $q$ the bimodule 
$\opn{H}^q(R)$ is a finite module over $A$ and over $A^{\mrm{op}}$. 
In fact, we are going to treat the left module structure of $\opn{H}^q(R)$ 
only; the right module structure is treated the same way, by replacing $A$ 
with $A^{\mrm{op}}$.

Let us write
\[ \til{R}_t := \opn{Ind}_{(\til{A}_t)^{\mrm{en}}}(\til{R}) \cong 
\til{A}_t \ot_{\til{A}} \til{R} \ot_{\til{A}} \til{A}_t \in 
\dcat{D} \bigl( (\til{A}_t)^{\mrm{en}}, \mrm{gr} \bigr) . \]
According to Lemma \ref{lem:4395} there is an isomorphism 
$\til{R}_t \cong \til{A}_t \ot_{\til{A}} \til{R}$ 
in $\dcat{D}(\til{A}_t, \mrm{gr})$. 
Since $\til{R} \in \dcat{D}_{\mrm{f}}(\til{A}, \mrm{gr})$, 
it follows that 
$\til{R}_t \in \dcat{D}_{\mrm{f}}(\til{A}_t, \mrm{gr})$. 
We identify $(\til{A}_t)_0$ with $A$ using the canonical isomorphism $h$
from (\ref{eqn:4476}). 
By Lemma \ref{lem:4400}(2) we get 
$\opn{Deg}_0(\til{R}_t) \in \dcat{D}_{\mrm{f}} \bigl( (\til{A}_t)_0 \bigr) 
= \dcat{D}_{\mrm{f}}(A)$.
Hence 
\begin{equation} \label{eqn:4426}
R = \opn{Deg}_0(\til{R}_t)[-1] \in \dcat{D}_{\mrm{f}}(A) . 
\end{equation}
We see that $\opn{H}^q(R)$ is a finite $A$-module for all $q$. 

\medskip \noindent
Step 2. Here we prove that $R$ has finite injective dimension over $A$ and over 
$A^{\mrm{op}}$. Again, we only examine the left side; the right side is done 
similarly, just replacing $A$ with $A^{\mrm{op}}$.

Because $A$ is noetherian, it suffices to find a uniform bound for the 
cohomology of $\opn{RHom}_{A}(M, R)$, for all $M \in \dcat{M}_{\mrm{f}}(A)$. 
Take such an $A$-module $M$. 
By Lemma \ref{lem:4390}(3) there exists a graded module  
$\til{M} \in \dcat{M}_{\mrm{f}}(\til{A}, \mrm{gr})$
such that, letting
$\til{M}_t := \til{A}_t \ot_{\til{A}} M$, we have 
\begin{equation} \label{eqn:4425}
M \cong \opn{Deg}_0(\til{M}_t) 
\end{equation} 
in $\dcat{M}_{\mrm{f}}(A)$. There are isomorphisms  
\begin{equation} \label{eqn:4427}
\begin{aligned}
&
\opn{RHom}_{\til{A}}(\til{M}, \til{R}) \ot_{\til{A}} \til{A}_t
\cong^1 \opn{RHom}_{\til{A}}(\til{M}, \til{R} \ot_{\til{A}} \til{A}_t)
\\
& \quad 
\cong^2 \opn{RHom}_{\til{A}}(\til{M}, \til{R}_t)
\cong^3 \opn{RHom}_{\til{A}_t}(\til{M}_t, \til{R}_t)
\end{aligned}
\end{equation}
in $\dcat{D}(\til{A}_t^{\mrm{op}}, \mrm{gr})$.
The isomorphism $\cong^1$ is an instance of the graded tensor-evaluation 
isomorphism (Theorem \ref{thm:4535}), that is valid because $\til{A}$ is 
noetherian, $\til{M} \in \lb \dcat{M}_{\mrm{f}}(\til{A}, \mrm{gr})$, 
$\til{R} \in \dcat{D}^{\mrm{b}}(\til{A}^{\mrm{en}}, \mrm{gr})$, 
and $\til{A}_t$ is flat over $\til{A}$. 
The isomorphism $\cong^2$ is by Lemma \ref{lem:4395}. 
And the isomorphism $\cong^3$ is adjunction for the ring homomorphism 
$\til{A} \to \til{A}_t$.
 
Next, according to Lemma \ref{lem:4401}, and using the isomorphisms 
(\ref{eqn:4425}) and (\ref{eqn:4426}), there is an isomorphism 
\[ \begin{aligned}
&
\opn{Deg}_0 \bigl( \opn{RHom}_{\til{A}_t}(\til{M}_t, \til{R}_t) \bigr) \cong
\opn{RHom}_{\opn{Deg}_0(\til{A}_t)} 
\bigl( \opn{Deg}_0(\til{M}_t), \opn{Deg}_0(\til{R}_t) \bigr)
\\ & \quad 
\cong \opn{RHom}_{A}(M, R)[1] 
\end{aligned} \]
in $\dcat{D}(A^{\mrm{op}})$. 
Combining this with formula (\ref{eqn:4427}), 
we see that the injective dimension of $R$ over $A$ is at most one more than 
the graded injective dimension of $\til{R}$ over $\til{A}$; and this is known 
to be finite.

\medskip \noindent
Step 3. In this step we prove that $R$ has the derived Morita property on both 
sides. As before, we only prove this property on the $A$ side; the 
$A^{\mrm{op}}$ side is done similarly. 

The isomorphisms (\ref{eqn:4427}), but for the complex of bimodules $\til{R}$ 
instead of the for the module $\til{M}$, become isomorphisms 
\begin{equation} \label{eqn:4428}
\begin{aligned}
&
\opn{RHom}_{\til{A}}(\til{R}, \til{R}) \ot_{\til{A}} \til{A}_t
\cong^1 \opn{RHom}_{\til{A}}(\til{R}, \til{R} \ot_{\til{A}} \til{A}_t)
\\
& \quad 
\cong^2 \opn{RHom}_{\til{A}}(\til{R}, \til{R}_t)
\cong^3 \opn{RHom}_{\til{A}_t}(\til{A}_t \ot_{\til{A}} \til{R}, \til{R}_t)
\\
& \quad 
\cong^4 \opn{RHom}_{\til{A}_t}(\til{R}_t, \til{R}_t)
\end{aligned}
\end{equation}
in $\dcat{D}(\til{A}_t^{\mrm{op}}, \mrm{gr})$.
The last isomorphism $\cong^4$ is another case of Lemma \ref{lem:4395}.
Define 
$R' := \opn{Deg}_0(\til{R}_t) \in \dcat{D}_{\mrm{f}}(A)$;
so $R = R'[-1]$. 
We now apply $\opn{Deg}_0$ to the last object in (\ref{eqn:4428}), obtaining the
isomorphisms 
\begin{equation} \label{eqn:4431}
\begin{aligned}
&
\opn{Deg}_0 \bigl( \opn{RHom}_{\til{A}_t}(\til{R}_t, \til{R}_t) \bigr)
\\
& \quad 
\cong^{\tup{(i)}} \opn{RHom}_{\opn{Deg}_0(\til{A}_t)} 
\bigl( \opn{Deg}_0(\til{R}_t), \opn{Deg}_0(\til{R}_t) \bigr)
\\
& \quad 
\cong^{\tup{(ii)}} \opn{RHom}_{A}(R', R')
\cong^{\tup{(iii)}} \opn{RHom}_{A}(R, R)
\end{aligned}
\end{equation}
in $\dcat{D}(A^{\mrm{op}})$.
The isomorphism $\cong^{\tup{(i)}}$ is due to Lemma \ref{lem:4401}, with 
$\til{B} = \til{C}:= \til{A}_t$. 
The isomorphism $\cong^{\tup{(ii)}}$ is simply the definition of $R'$, and 
the isomorphism $\cong^{\tup{(iii)}}$ is because 
$R = R'[-1]$. 

We know that $\til{R}$ has the graded derived NC Morita property on the 
$\til{A}$ side. Using Lemma \ref{lem:4361} (that is true also in the graded 
situation) we see that \lb 
$\opn{H}^q \bigl( \opn{RHom}_{\til{A}}(\til{R}, \til{R}) \bigr) = 0$ 
for all $q \neq 0$. Because the functors $\opn{H}^q$ and $\opn{Deg}_0$ commute 
with each other, the isomorphisms (\ref{eqn:4428}) and (\ref{eqn:4431}) tell us 
that 
\begin{equation} \label{eqn:4433}
\opn{H}^q \bigl( \opn{RHom}_{A}(R, R) \bigr) = 0
\end{equation}
for all $q \neq 0$.

For $q = 0$ we know that 
$\opn{H}^0 \bigl( \opn{RHom}_{\til{A}}(\til{R}, \til{R}) \bigr)$
is a graded-free $\til{A}^{\mrm{op}}$-module with basis $\opn{id}_{\til{R}}$. 
This is by Lemma \ref{lem:4361}, applied in the graded situation. 
The isomorphisms (\ref{eqn:4428}) send the element
$\opn{id}_{\til{R}} \ot \, 1$ to $\opn{id}_{\til{R}_t}$; and therefore \lb 
$\opn{H}^0 \bigl( \opn{RHom}_{\til{A}_t}(\til{R}_t, \til{R}_t) \bigr)$
is a graded-free $\til{A}_t^{\mrm{op}}$-module with basis 
$\opn{id}_{\til{R}_t}$.
The element $\opn{id}_{\til{R}_t}$ has degree $0$, and hence  
\[ \opn{H}^0 \bigl(  \opn{Deg}_0 \bigl( \opn{RHom}_{\til{A}_t}(\til{R}_t, 
\til{R}_t) \bigr) \bigr) \cong 
\opn{Deg}_0 \bigl( \opn{H}^0 \bigl(   \opn{RHom}_{\til{A}_t}(\til{R}_t, 
\til{R}_t) \bigr) \bigr) \]
is a free $A^{\mrm{op}}$-module with basis $\opn{id}_{\til{R}_t}$.
The isomorphisms (\ref{eqn:4431}) send the element 
$\opn{id}_{\til{R}_t}$ to $\opn{id}_{R}$; so 
$\opn{H}^0 \bigl( \opn{RHom}_{A}(R, R) \bigr)$ 
is a free $A^{\mrm{op}}$-module with basis $\opn{id}_{R}$.
Again invoking Lemma \ref{lem:4361}, this fact, with formula (\ref{eqn:4433}), 
say that $R$ has the derived NC Morita property on the $A$ side.  

\medskip \noindent
Step 4. Here we produce a rigidifying isomorphism $\rho$ for $R$, starting with 
a given graded rigidifying isomorphism 
\[ \til{\rho} : \til{R} \iso 
\opn{RHom}_{\til{A}^{\mrm{en, out}}} ( \til{A}, \til{R} \ot \til{R} ) \]
in $\dcat{D}(\til{A}^{\mrm{en}}, \mrm{gr})$.  
Here is a list of isomorphisms in 
$\dcat{D} \bigl( (\til{A}_t)^{\mrm{en}}, \mrm{gr} \bigr)$.
\begin{equation} \label{eqn:4435}
\begin{aligned}
& \til{R}_t \cong^1 
\opn{RHom}_{\til{A}^{\mrm{en, out}}} ( \til{A}, \til{R} \ot \til{R})
\ot_{\til{A}^{\mrm{en}}} (\til{A}_t)^{\mrm{en}}
\\
& \quad \cong^{2} 
\opn{RHom}_{\til{A}^{\mrm{en, out}}} \bigl( \til{A}, 
(\til{R} \ot \til{R}) \ot_{\til{A}^{\mrm{en}}} (\til{A}_t)^{\mrm{en}} \bigr)
\\
& \quad \cong^{3} 
\opn{RHom}_{\til{A}^{\mrm{en, out}}} 
\bigl( \til{A}, (\til{R} \ot_{\til{A}} \til{A}_t) \ot
(\til{A}_t \ot_{\til{A}} \til{R}) \bigr) 
\\
& \quad \cong^{4} 
\opn{RHom}_{\til{A}^{\mrm{en, out}}} (\til{A}, \til{R}_t \ot \til{R}_t) 
\\
& \quad \cong^{5} 
\opn{RHom}_{(\til{A}_t)^{\mrm{en, out}}} (\til{A}_t, \til{R}_t \ot \til{R}_t) .
\end{aligned}
\end{equation}
The justifications are:
\begin{itemize}
\item[$\cong^{1}$~:] Here we apply the functor 
$\opn{Ind}_{(\til{A}_t)^{\mrm{en}}}$ to the isomorphism $\til{\rho}$. The 
inside action of $\til{A}^{\mrm{en}}$ on $\til{R} \ot \til{R}$ is treated as 
a right action. 

\item[$\cong^{2}$~:] This is by the graded tensor-evaluation 
isomorphism (Theorem \ref{thm:4535}). It holds because $\til{A}$ is
a derived graded-pseudo-finite complex over $\til{A}^{\mrm{en}}$,
see Corollary \ref{cor:4531}; 
$\til{R} \ot \til{R} \in \dcat{D}^{\mrm{b}}(\til{A}^{\mrm{four}}, \mrm{gr})$;
and $(\til{A}_t)^{\mrm{en}}$ is flat over $\til{A}^{\mrm{en}}$.

\item[$\cong^{3}$~:] This is a rearrangement of derived tensor factors. 

\item[$\cong^{4}$~:] By Lemma \ref{lem:4395}. 

\item[$\cong^{5}$~:] This is adjunction for the graded ring homomorphism 
$\til{A}^{\mrm{en}} \to (\til{A}_t)^{\mrm{en}}$, with the fact that 
$\til{A}_t = \til{A}_t \ot_{\til{A}} \til{A} \ot_{\til{A}} \til{A}_t$ as graded 
rings.
\end{itemize}

Now we apply $\opn{Deg}_0$ to (\ref{eqn:4435}), 
and we obtain the first isomorphism below in $\dcat{D}(A^{\mrm{en}})$.
\begin{equation} \label{eqn:4436}
\begin{aligned}
& 
R' = \opn{Deg}_0(\til{R}_t) \cong^{\tup{(i)}} 
\opn{Deg}_0 \bigl( 
\opn{RHom}_{(\til{A}_t)^{\mrm{en, out}}} (\til{A}_t, 
\til{R}_t \ot \til{R}_t) \bigr) 
\\ & \quad 
\cong^{\tup{(ii)}} \opn{RHom}_{A^{\mrm{en, out}}} (A, R' \ot R')[-1]
\end{aligned}
\end{equation}
The isomorphism $\cong^{\tup{(ii)}}$ is by Theorem \ref{thm:4510}, with 
$\til{B} := \til{A}_t$, that we are allowed to use because $A$ is a derived 
pseudo-finite complex over $A^{\mrm{en}}$; see Theorem \ref{thm:4511}. 
Plugging in $R = R'[-1]$ we get 
$R \cong \opn{RHom}_{A^{\mrm{en, out}}}(A, R \ot R)$
in $\dcat{D}(A^{\mrm{en}})$. This is the rigidifying isomorphism $\rho$.
\end{proof}

\mysubsection{Filtered Rings and Existence of Rigid NC DC}
\label{subsec:Fil-exit-R-NCDC}

In this subsection we give M. Van den Bergh's proof of the existence of NC 
rigid dualizing complexes (Theorem \ref{thm:4396}). This proof appeared in his 
seminal paper \cite{VdB}, and then there were some improvements in 
\cite{YeZh1}. The proof here is much more detailed. 

We continue with Convention \ref{conv:3670}.
Specifically, $\K$ is a base field, and all rings are central over $\K$. 

Earlier in the book (in Section \ref{sec:exist-resol}) we encountered 
filtrations of DG 
modules. Here we are interested in filtrations of rings. By a {\em filtration 
of a ring $A$} we mean a collection 
$\bigl\{ F_j(A) \bigr\}_{j \geq -1}$ of $\K$-submodules 
\[ F_{-1}(A) \sub F_{0}(A) \sub F_{1}(A) \sub \cdots \sub A \]
such that 
$F_{-1}(A) = 0$, $\bigcup_j F_j(A) = A$, 
$1_A \in F_{0}(A)$ and 
$F_{j}(A) \cd F_{k}(A) \sub F_{j + k}(A)$ 
for all $j, k$. 

Given a filtration $F = \bigl\{ F_j(A) \bigr\}_{j \geq -1}$ of the ring $A$, 
we write
$\opn{Gr}^F_j(A) := \lb F_{j}(A) / F_{j - 1}(A)$ 
for $j \geq 0$. 
The {\em associated graded ring} 
$\opn{Gr}^F(A) := \bigoplus_{j \geq 0} \opn{Gr}^F_j(A)$ 
is a nonnegative algebraically graded ring. 

A filtration $F$ on the ring $A$ gives rise to another graded ring: it is the  
{\em Rees ring} 
$\opn{Rees}^F(A) := \bigoplus_{j \geq 0} F_j(A) \cd t^j \sub A[t]$. 
Here $t$ is a central variable of degree $1$, that we call the Rees parameter. 
There is a ring isomorphism 
$A \cong \opn{Rees}^F(A) / (t - 1)$
and a graded ring isomorphism
$\opn{Gr}^F(A) \cong \opn{Rees}^F(A) / (t)$.
We often use the abbreviations 
$\bar{A} := \opn{Gr}^F(A)$ and 
$\til{A} := \opn{Rees}^F(A)$. 
This is consistent with the notation used earlier in the section: 
$\til{A}$ is a graded ring, $t$ is a degree $1$ central regular element in 
it, and $A = \til{A} / (t - 1)$.

For an element $c \in \K$ we write 
$\opn{pr}_{c} : \til{A} \to \til{A} / (t - c)$ 
for the canonical surjective ring homomorphism. 
The homomorphisms $\opn{pr}_{0}$ and $\opn{pr}_{1}$ are 
displayed in this diagram of rings:
\[ \UseTips \xymatrix @C=2ex @R=5ex {
&
\til{A} = \opn{Rees}^F(A)
\ar[dl]_(0.6){\opn{pr}_1}
\ar[dr]^(0.6){\opn{pr}_0}
\\
*+{A}
&
&
*+{\bar{A} = \opn{Gr}^F(A)}
} \] 
The filtration $F$ can be recovered from the Rees ring as follows: 
\begin{equation} \label{eqn:4506}
F_j(A) = \opn{pr}_{1} (\til{A}_j) 
= \opn{pr}_{1} \Bigl( \bigoplus\nolimits_{i = 0}^j \, \til{A}_i \Bigr) \sub A . 
\end{equation}

\begin{dfn}[\cite{YeZh1}] \label{dfn:4445}
A filtration $F$ of the ring $A$ is called a 
{\em noetherian connected filtration} 
\index{Filtration on a ring! noetherian connected}
if the graded ring $\opn{Gr}^F(A)$ is a noetherian connected graded 
ring (Definitions \ref{dfn:3703} and \ref{dfn:3706}). 
\end{dfn}

Let $\til{B}$ be a connected graded ring. Recall from Definition \ref{dfn:4532} 
that a complex $M \in \dcat{D}(\til{B}, \mrm{gr})$ is said to be 
derived graded-pseudo-finite if it belongs to the saturated full 
triangulated subcategory of $\dcat{D}(\til{B}, \mrm{gr})$
generated by the bounded above complexes $P$ such that each $P^i$ is a 
finite graded-free $\til{B}$-module. 

\begin{prop} \label{prop:4395}
Assume the ring $A$ admits a noetherian connected filtration $F$. Then\tup{:}
\begin{enumerate}
\item The Rees ring $\til{A} := \opn{Rees}^F(A)$ is a noetherian connected 
graded ring, and $\til{A}$ is a derived graded-pseudo-finite complex  
over $\til{A}^{\mrm{en}}$.

\item The ring $A$ is noetherian, and $A$ is a derived pseudo-finite complex 
over $A^{\mrm{en}}$.
\end{enumerate}
\end{prop}

\begin{proof} \mbox{}

\smallskip \noindent 
(1) Let $\bar{A} := \opn{Gr}^F(A)$, which is a noetherian connected 
graded ring. Since 
$\til{A}_0 = \opn{Rees}^F_0(A) = F_0(A) = \opn{Gr}^F_0(A) = \bar{A}_0 = \K$
and 
$\til{A}_i = \opn{Rees}^F_i(A) \cong F_i(A) \cong 
\bigoplus_{0 \leq j \leq i} \bar{A}_j$ 
as ungraded $\K$-modules, we see that $\til{A}$ is connected graded. 
Next, because $\bar{A} = \til{A} / (t)$, by Theorem \ref{thm:4240} we deduce 
that $\til{A}$ is noetherian. 
By Proposition \ref{prop:4266}, $\til{A}$ is a derived graded-pseudo-finite
complex over $\til{A}^{\mrm{en}}$.

\medskip \noindent 
(2) By item (1)  the graded ring $\til{A}$ is connected and noetherian. 
Since $A = \til{A} / (t - 1)$, it follows that $A$ is noetherian. And according 
to Theorem \ref{thm:4511}(3), $A$ is a derived pseudo-finite complex over 
$A^{\mrm{en}}$.  
\end{proof}

\begin{thm}[Van den Bergh Existence, \cite{VdB}] \label{thm:4396}
\index{Dualizing complex! balanced}
\index{Dualizing complex! noncommutative rigid}
\index{Filtration on a ring! noetherian connected}
Let $A$ be a ring. Assume $A$ admits a noetherian connected filtration 
$F$, such that the connected graded ring 
$\opn{Gr}^F(A)$ has a balanced NC dualizing complex.
Then the ring $A$ has a rigid NC dualizing complex $(R, \rho)$. Moreover, 
$(R, \rho)$ is unique up to a unique rigid isomorphism.
\end{thm}

\begin{proof}
As before we write $\til{A} := \opn{Rees}^F(A)$ and $\bar{A} := \opn{Gr}^F(A)$.
These are noetherian connected graded rings, and 
$\bar{A} = \til{A} / (t)$. Since $\bar{A}$ admits a balanced NC dualizing 
complex, Corollary \ref{cor:4585} says that $\bar{A}$ satisfies the $\chi$ 
condition and it has finite local cohomological dimension. 
By Theorem \ref{thm:4250} the graded ring $\til{A}$ also satisfies the $\chi$ 
condition and it has finite local cohomological dimension.
Using Corollary \ref{cor:4585} for $\til{A}$, we conclude that it has a 
balanced NC dualizing complex $\til{R}$. 

According to Theorem \ref{thm:4400}, $\til{R}$ is a graded rigid NC dualizing 
complex over $\til{A}$. Theorem \ref{thm:4395} then says that 
\[ R := (\opn{Rest}_{A^{\mrm{en}}} \circ \opn{Deg}_0 \circ 
\opn{Ind}_{(\til{A}_t)^{\mrm{en}}}) (\til{R})[-1] \in \dcat{D}(A^{\mrm{en}}) \]
is a rigid NC dualizing complex over $A$; i.e.\ there exists a rigidifying 
isomorphism $\rho$ for $R$. 

By Proposition \ref{prop:4395}(2), $A$ is a derived pseudo-finite complex over
$A^{\mrm{en}}$. Then, according to Theorem \ref{thm:4380}, the rigid dualizing 
complex $(R, \rho)$ is unique up to a unique rigid isomorphism.
\end{proof}

Here is a result we won't prove here. It is needed for the example that comes 
after it. 

\begin{prop} [{\cite[Proposition 6.13]{YeZh1}}]  \label{prop:4505}
Let $A$ be a commutative ring, and let $f : A \to B$ be a finite central ring 
homomorphism. If $A$ has a noetherian connected filtration $F$, then 
there exists a noetherian connected filtration $G$ of $B$, such that 
$f(F_j(A)) \sub G_j(B)$ for all $j$, and the graded ring homomorphism 
$\opn{Gr}(f) : \opn{Gr}^F(A) \to \opn{Gr}^G(B)$ 
is finite. 
\end{prop}

\begin{exa} \label{exa:4506}
Suppose the ring $B$ is finite over its center $\opn{Cent}(B)$, and 
$\opn{Cent}(B)$ is a finitely generated $\K$-ring. 
Then we can find a finite central homomorphism $f : A \to B$, where
$A =  \K[t_1, .., t_n]$ is a commutative polynomial ring. 
The grading on $A$ with $\opn{deg}(t_i) = 1$ gives rise to a noetherian 
connected filtration $F$, by the formula 
$F_j(A) := \bigoplus_{k \leq j} A_k$. Of course 
$\opn{Gr}^F(A) \cong A$. 

Now we can use Proposition \ref{prop:4505}
to conclude that $B$ has a noetherian connected filtration $G$ such that 
$\opn{Gr}^F(A) \to \opn{Gr}^G(B)$ is finite. According to Example \ref{exa:3990}
and Corollary \ref{cor:4032}, the 
graded ring $\opn{Gr}^G(B)$ has a balanced dualizing complex. 
Therefore, by Theorem \ref{thm:4396}, $B$ has a rigid NC dualizing complex. 
\end{exa}

Sometimes rings come with filtrations of the following type. 

\begin{dfn} \label{dfn:4515}
A {\em central filtration of finite type} of a ring $A$ is a filtration 
$G = \bigl\{ G_j(A) \bigr\}_{j \geq -1}$
of $A$, such that the graded ring 
$\opn{Gr}^G(A) = \bigoplus_{j \geq 0} \opn{Gr}^G_j(A)$
is finite over its center 
$C = \bigoplus_{j \geq 0} C_j  := \opn{Cent} \bigl( \opn{Gr}^G(A) \bigr)$, 
and the commutative graded ring $C$ is finitely generated over $\K$.
\end{dfn}

In \cite{YeZh2.5} such a filtration was called a ``differential 
filtration of finite type'', but this name sounds too confusing when we also 
talk about DG rings. The reason for the name ``differential'' will be apparent 
in Example \ref{exa:4515} below. 

It is not hard to see, using Theorem \ref{thm:4240}, 
that if $A$ admits a central filtration of finite type, then 
it is noetherian. However, much more is true, as the next theorem shows.  
This is a result from \cite{YeZh2.5} that we present without a 
proof. A slightly weaker version of this theorem was in the paper
\cite{McSt} of J.C. McConnell and J.T. Stafford, which 
the authors attributed to J. Bernstein. 

\begin{thm}[Two Filtrations, {\cite[Theorem 3.1]{YeZh2.5},
\cite{McSt}}] \label{thm:4505}
Let $A$ be a ring that admits a central filtration of finite type $G$.
Then $A$ admits a noetherian connected filtration $F$, such that the 
noetherian connected graded ring $\opn{Gr}^F(A)$ 
is a commutative finitely generated $\K$-ring.
\end{thm}

As explained in Example \ref{exa:3991}, the 
ring $\opn{Gr}^F(A)$ has a balanced dualizing complex. Therefore, using 
Theorems \ref{thm:4505} and \ref{thm:4396}, we conclude that:

\begin{cor} \label{cor:4515}
If $A$ admits a central filtration of finite type, then $A$ has a 
rigid NC dualizing complex $(R, \rho)$, and it is unique up to 
a unique rigid isomorphism.
\end{cor}

\begin{exa} \label{exa:4515}
Assume $\K$ has characteristic $0$. Let $C$ be a smooth commutative 
$\K$-ring of pure dimension $n$, and let $A := \mcal{D}(C)$, the ring of 
differential operators. 
Consider the filtration
$G = \bigl\{ G_j(A) \bigr\}_{j \geq -1}$
of $A$ by order of operators. For this filtration we have 
$G_0(A) = C$ and 
$G_1(A) = C \oplus \mcal{T}(C)$,
where $\mcal{T}(C)$ is the module of derivations.
Then $\opn{Gr}^G(A)$ is a commutative ring, and there is a canonical 
graded ring isomorphism 
$\opn{Gr}^G(A) \cong \opn{Sym}_{C}(\mcal{T}(C))$,
the symmetric tensor ring on the finitely generated projective $C$-module 
$\mcal{T}(C)$. This is a commutative finite type $\K$-ring, so $G$ is a central 
filtration of finite type.
Corollary \ref{cor:4515} says that $A$ has a rigid NC dualizing complex 
$(R, \rho)$. In \cite{Ye4.5} it was shown that $R \cong A[2 n]$.
\end{exa}

We end the subsection with three remarks. 

\begin{rem} \label{rem:4870}
Here is a brief history of the early days NC dualizing complexes (during the 
1990's). The first paper on this topic was \cite{Ye1} by A. Yekutieli from 
1992, in which dualizing complexes over NC connected graded rings (always over a 
base field $\K$) were introduced (see Definition \ref{dfn:3716}).
This material was part of the Ph.D. thesis of Yekutieli, and was inspired by 
ideas of his advisor M. Artin on a noncommutative analogue of Serre Duality. 
The paper \cite{Ye1} also contained the definition of balanced dualizing 
complexes (Definition \ref{dfn:3714}), a proof of NC Local Duality (Theorem 
\ref{thm:4616}), a proof of the existence of balanced dualizing complexes 
(Theorem \ref{thm:3781}), and some particular calculations (Example 
\ref{exa:3780}). 

The paper \cite{YeZh0} by Yekutieli and J.J. Zhang from 1997 used the results 
of \cite{Ye1} to establish NC Serre Duality on noncommutative projective 
schemes, in the sense of Artin and Zhang \cite{ArZh}. In the same paper we 
also find Corollary \ref{cor:4655}, saying that the existence of a balanced 
dualizing complex over the noetherian connected ring $A$ implies that $A$ 
satisfies the $\chi$ condition and has finite local cohomological dimension. 

The paper \cite{VdB} by M. Van den Bergh from 1997 begins by 
establishing the important converse to Corollary \ref{cor:4655}, namely that 
a noetherian connected ring $A$, which satisfies the $\chi$ condition and has 
finite local cohomological dimension, admits a balanced dualizing complex 
(Theorem \ref{thm:3713}). Thus we have Corollary \ref{cor:4585} on the 
``two equivalent properties''. Van den Bergh went on to define rigid dualizing 
complexes, to prove their uniqueness (Theorem \ref{thm:4380}), and to 
prove existence of rigid dualizing complexes, in the presence of suitable 
filtrations and balanced dualizing complexes (Theorem \ref{thm:4396}).
The subsequent paper \cite{VdB2} by Van den Bergh utilized the rigidity 
formalism to deduce what is now called Van den Bergh Duality 
(Theorem \ref{thm:4860} below). 

At about the same time that Van den Bergh wrote \cite{VdB}, Yekutieli worked on 
the paper \cite{Ye4}, in which the derived Picard group $\opn{DPic}(A)$ of a NC 
ring $A$ was introduced. The fact that any two NC dualizing complexes $R, R'$ 
over $A$ are related by $R' \cong T \ot^{\mrm{L}}_{A} R$ for some tilting 
complex $T$ (Theorem \ref{thm:3681}, proved in \cite{Ye4}) 
was used by Van den Bergh in his proof of the uniqueness of rigid dualizing 
complexes in \cite{VdB}.

The paper \cite{YeZh1} by Yekutieli and Zhang from 1999 built on Van den 
Bergh's ideas. The uniqueness results of Van den Bergh were sharpened, and were 
upgraded to functoriality (see Remark \ref{rem:4380}). The Auslander condition 
on dualizing complexes (see Remark \ref{rem:3821}) was also introduced and 
studied in \cite{YeZh1}.

It is hard to track the progress of the theory of NC dualizing complexes beyond 
the year 2000. This theory has merged with, or has been absorbed 
into, the more active research on Calabi-Yau rings and categories
(see next subsection), with its ties 
to birational algebraic geometry and mathematical physics. 
\end{rem}

\begin{rem} \label{rem:4380}
This remark is about the functorial properties of NC rigid dualizing 
complexes. In Remark \ref{rem:4191} we talked about the rigid trace 
homomorphism and the rigid localization homomorphism, between rigid residue 
complexes over commutative rings. These are actual homomorphisms of complexes. 
For NC rings there are usually no residue complexes (because of the pathological 
nature of injective modules in the NC setting). 
Still we can consider rigid trace morphisms and rigid localization morphisms 
between NC rigid dualizing complexes. These are morphisms in 
derived categories. Below is an outline. 

Suppose $f : A \to B$ is a finite ring homomorphism, and the rigid dualizing 
complexes $(R_A, \rho_A)$ and $(R_B, \rho_B)$ exist. A {\em rigid trace 
morphism} 
$\opn{tr}_f : (R_B, \rho_B) \to (R_A, \rho_A)$ 
is a backward morphism $\opn{tr}_f : R_B \to R_A$ in $\dcat{D}(A^{\mrm{en}})$, 
which is rigid (the obvious 
variant of Definition \ref{dfn:4372}), and is nondegenerate on both sides 
(see Definition \ref{dfn:4031}). In \cite{YeZh1} the existence and uniqueness 
of the rigid trace morphism exists was proved, when $f : A \to 
B$ is finite centralizing, and $A$ admits a noetherian connected filtration $F$ 
such that $\opn{Gr}^F(A)$ has a balanced dualizing complex. The proof is along 
the same lines as the proof of Theorem \ref{thm:4396} above, combined 
with the existence of the balanced trace morphism, established in Theorem 
\ref{thm:4030}. Among the applications of the rigid trace: (1) It was used 
in \cite{Ye4.5} to calculate the rigid dualizing complex of 
$A = \mrm{U}(\g)$, see Example \ref{exa:4505} below. (2) In 
\cite{YeZh1.5} it was used to calculate multiplicities of indecomposable 
injectives in minimal injective resolutions. 

Now consider a localization ring homomorphism $g : A \to A'$. Assume the rigid 
dualizing complexes $(R_A, \rho_A)$ and $(R_{A'}, \rho_{A'})$ exist. A {\em 
rigid localization morphism} 
$\opn{q}_g : (R_A, \rho_A) \to (R_{A'}, \rho_{A'})$ 
is a forward morphism $\opn{q}_g : R_A \to R_{A'}$ in 
$\dcat{D}(A^{\mrm{en}})$, which is rigid (it respects 
the two rigidifying isomorphisms), and is nondegenerate on both sides. 
Existence and uniqueness of the rigid localization morphism was proved in 
\cite[Theorem 6.2]{YeZh2.5}; it is related to Ore conditions. One application 
of the rigid localization morphism is in the concept of {\em homological 
transcendence degree} \cite{YeZh1.7} -- the ring $A$ is an 
integral domain which has a rigid dualizing complex (say because of the 
existence of a good filtration on it), and $D = A'$ is the division ring of 
fractions of $A$. The homological transcendence degree of $D$ is the unique 
number $n$ such that the rigid dualizing complex $R_D$ satisfies 
$\opn{H}^{-n}(R_D) \neq 0$. 
\end{rem}

\begin{rem} \label{rem:3821}
In this remark we discuss the {\em Auslander property} of NC dualizing 
complexes. Around 1990, J.E. Bj\"ork \cite{Bjo} and T. Levasseur \cite{Lev}  
had studied NC regular noetherian rings with the Auslander property. This 
property was later generalized in the papers \cite{Ye2.5}
and \cite{YeZh1} by A. Yekutieli and J.J. Zhang. Below is a quick explanation 
of the meaning and applications of this concept. 

Consider a NC noetherian ring $A$ with a NC dualizing complex $R$. 
Given a module $M \in \dcat{M}_{\mrm{f}}(A)$, its grade with respect to $R$ is 
the number 
\[ \opn{grade}_R(M) := 
\opn{\inf} \bigl( \{ q \mid \opn{Ext}^q_A(M, R) \neq 0 \} \bigr) \in 
\Z \cup  \{ \infty \} . \]
Note that $\opn{Ext}^q_A(M, R) \in \dcat{M}_{\mrm{f}}(A^{\mrm{op}})$.
The number $\opn{grade}_R(N)$ for $N \in \dcat{M}_{\mrm{f}}(A^{\mrm{op}})$
is defined similarly. The complex $R$ is called {\em Auslander} if for every 
$M \in \dcat{M}_{\mrm{f}}(A)$, every $q$, and 
every $A^{\mrm{op}}$-submodule 
$N \sub \opn{Ext}^q_A(M, R)$, the inequality $\opn{grade}_R(N) \geq q$
holds; and the same inequality holds with $A$ and $A^{\mrm{op}}$ 
exchanged. The original definition used by Bj\"ork and Levasseur, for a regular 
ring $A$, was for the particular dualizing complex $R := A$. 

In \cite{YeZh1} it was proved that in many situations the rigid dualizing 
complex $R$ from Theorem \ref{thm:4396} has the Auslander property. For 
instance, \cite[Theorem 0.7]{YeZh1} implies that if $A$ has a noetherian 
connected filtration $F$ such that $\opn{Gr}^F(A)$ is commutative, then the 
rigid dualizing complex $R$ of $A$ is Auslander. 

Having an Auslander rigid dualizing complex gives rise to the {\em canonical 
dimension} of $A$-modules. For $M \in \dcat{M}_{\mrm{f}}(A)$ its canonical 
dimension is $\opn{Cdim}(M) := - \opn{grade}_R(M)$.
This can be considered as a NC version of the rigid dimension from 
Section \ref{sec:dual-cplx-comm-rng}.
Indeed, if $A$ is a finite type commutative $\K$-ring, with commutative rigid 
dualizing complex $R$, then for a prime ideal $\p \sub A$ and the $A$-module 
$M := A / \p$ we have 
$\opn{Cdim}(M) = \opn{rig{.}dim}_{\K}(\p) = \opn{dim}(A / \p)$,
where the second number is from Definition \ref{dfn:2537}, and the 
third number is the Krull dimension of the ring $A / \p$. 

As shown in \cite{YeZh1}, several structural results for NC rings 
can be deduced from the presence of an Auslander rigid dualizing complex; for 
instance, the catenarity of $\opn{Spec}(A)$, see \cite[Theorem 0.1]{YeZh1}.
Moreover, in \cite[Section 7]{YeZh2.5} a link was made 
between Auslander rigid dualizing 
complexes and perverse t-structures: the canonical dimension of modules defines 
a perverse t-structure on $\dcat{D}^{\mrm{b}}_{\mrm{f}}(A)$, and under the 
duality functor $D := \opn{RHom}_A(-, R)$ this t-structure goes to the standard 
t-structure on $\dcat{D}^{\mrm{b}}_{\mrm{f}}(A^{\mrm{op}})$. 
\end{rem}

\mysubsection{Twisted Calabi-Yau Rings}
\label{subsec:tw-CY}

In this final subsection we study a class of regular rings that is of great 
importance in recent ring theory. 
We continue with Convention \ref{conv:3670}.
Specifically, $\K$ is a base field, and all rings are central over $\K$. 
For a ring $A$, its enveloping ring is 
$A^{\mrm{en}} = A \ot A^{\mrm{op}}$. 

\begin{dfn} \label{dfn:4860}
A ring $A$ is called {\em homologically smooth over $\K$} if $A$ is a perfect 
complex over the ring $A^{\mrm{en}}$.
\end{dfn}

This definition was proposed by M. Kontsevich, around the year 2000, cf.\ 
\cite[Erratum]{VdB2}.
Clearly if $A$ is homologically smooth over $\K$, then $A$ is a derived 
pseudo-finite complex over $A^{\mrm{en}}$. 
Recall (from Example \ref{exa:3695}) that the ring $A$ is called regular if it 
has finite global cohomological dimension. 

\begin{prop} \label{prop:4860}
If $A$ is homologically smooth over $\K$, then $A$ is a regular ring. 
\end{prop}

\begin{proof}
By Theorem \ref{thm:3415} there is a quasi-isomorphism 
$\rho : P \to A$ in $\dcat{C}_{\mrm{str}}(A^{\mrm{en}})$, where $P$ 
is a bounded complex of finite projective $A^{\mrm{en}}$-modules. 
Let $d$ be the amplitude of $P$. 
Given an $A$-module $M$, define $P_M := P \ot_A M$, and let 
$\rho_M : P_M \to M$ be the homomorphism in
$\dcat{C}_{\mrm{str}}(A)$ induced by $\rho$. Since $\rho$ is a homotopy 
equivalence in $\dcat{C}_{\mrm{str}}(A^{\mrm{op}})$, it follows that $\rho_M$ is 
a quasi-isomorphism. Since $P^i$ is a projective $A^{\mrm{en}}$-module, 
and $M$ is a free $\K$-module, we see that $P_M^i = P^i \ot_A M$ is a 
projective $A$-module (often infinitely generated; cf.\ Example \ref{exa:5125}). 
The conclusion is that the 
projective dimension of $M$ is at most $d$. 

A similar consideration shows that every $M' \in \dcat{M}(A^{\mrm{op}})$
has projective dimension at most $d$. So the global homological 
dimension of $A$ is at most $d$. 
\end{proof}

\begin{rem} \label{rem:4860}
Suppose $A$ is a finite type commutative $\K$-ring. It can be shown that $A$ is 
homologically smooth over $\K$ if and only if it is smooth over $\K$, in the 
sense of \cite{EGA-IV}. Cf.\ \cite[Proposition 3.24]{Mil}.
\end{rem}

Given a $\K$-ring automorphism $\mu$ of $A$, we define the
{\em $\mu$-twisted bimodule} $A(\mu)$ as follows: as a left $A$-module,
$A(\mu)$ is the free $A$-module $A$; and the right $A$-module action is 
twisted by $\mu$. To be explicit, let us denote by $a \cd b$ the 
multiplication of the ring $A$. Then the left action of a ring element 
$a \in A$ on a bimodule element $b \in A(\mu)$ is 
$a \, \cd^{\mu} \, b := a \cd b \in A(\mu)$; and the right action is 
$b \, \cd^{\mu} \, a := b \cd \mu(a) \in A(\mu)$.
In fancier notation,  a ring element 
$a_1 \ot \opn{op}(a_2) \in A^{\mrm{en}}$
acts on a bimodule element $b \in A(\mu)$ from the left like this:
$(a_1 \ot \opn{op}(a_2)) \, \cd^{\mu} \, b = a_1 \cd b \cd \mu(a_2) 
\in A(\mu)$.
The graded version of this twisting already appeared in Definition 
\ref{dfn:3780}. 
If $\nu$ is another automorphism of $A$, then 
$A(\mu) \ot_A A(\nu) \cong A(\mu \circ \nu)$
in $\dcat{M}(A^{\mrm{en}})$. 

It is possible to twist other bimodules: 
given $M \in \dcat{M}(A^{\mrm{en}})$ we define 
\begin{equation} \label{eqn:4896}
M(\mu) := M \ot_A A(\mu) \in \dcat{M}(A^{\mrm{en}}) . 
\end{equation}
Let us make explicit the left and right actions of the 
ring $A^{\mrm{en}}$ on the module $M(\mu)$, using formula (\ref{eqn:4878}). 
The original left and right actions of an element $a \in A$ on an element
$m \in M$ are denoted by $a \cd m$ and $m \cd a$. Then an element 
$(a_1 \ot \opn{op}(a_2)) \in A^{\mrm{en}}$ acts on an element 
$m \in M(\mu)$ from the left and the right  by 
\begin{equation} \label{eqn:4904}
\begin{aligned}
& 
(a_1 \ot \opn{op}(a_2)) \, \cd^{\mu} \, m = 
a_1 \cd m \cd \mu(a_2) ,
\\
&
m \, \cd^{\mu} \, (a_1 \ot \opn{op}(a_2)) = 
a_2 \cd m \cd \mu(a_1) . 
\end{aligned}
\end{equation}

\begin{dfn} \label{dfn:4490}
\index{Ring! Calabi-Yau}
\index{Ring! twisted Calabi-Yau}
\index{Ring! Nakayama automorphism}
A ring $A$ is called a {\em twisted CY ring} if it has these 
two properties:
\begin{itemize}
\rmitem{i} The ring $A$ is homologically smooth over $\K$. 
 
\rmitem{ii} There is a natural number $n$, called the {\em dimension of $A$}, 
and a ring automorphism $\nu$ of $A$, called the {\em Nakayama automorphism}, 
such that 
\[ \opn{Ext}^i_{A^{\mrm{en}}}(A, A^{\mrm{en}}) \cong
\begin{cases}
A(\nu) & \tup{if} \ \ i = n 
\\
0 & \tup{if} \ \ i \neq n 
\end{cases} 
\]
as $A$-bimodules. 
\end{itemize}
The ring $A$ is called a {\em CY ring} if the third property also holds:
\begin{itemize}
\rmitem{iii} The automorphism $\nu$ is inner, so that $A(\nu) \cong A$ as 
$A$-bimodules.
\end{itemize}
\end{dfn}

In property (ii) the $A^{\mrm{en}}$-module 
$\opn{Ext}^i_{A^{\mrm{en}}}(A, A^{\mrm{en}})$
should be viewed as a special case of 
$\opn{Ext}^i_{B}(M, N) \in \dcat{M}(C^{\mrm{op}})$ 
for $M \in \dcat{M}(B)$ and $N \in \dcat{M}(B \ot C^{\mrm{op}})$,
where we take $B, C := A^{\mrm{en}}$, $M := A$ and 
$N := A^{\mrm{en}}$, and we use the isomorphism (\ref{eqn:4875}) to switch
from a right $A^{\mrm{en}}$-module to a left $A^{\mrm{en}}$-module. 

The letters ``CY'' are an abbreviation for ``Calabi-Yau''. 
Most texts use the term ``CY algebra'' of course. See Remark 
\ref{rem:4505} for a discussion of the background of this definition. 
Clearly this is a relative notion: it pertains to $A$ as a central $\K$-ring. 
It is easy to see that the Nakayama automorphism 
$\nu$ is unique up to composition with an inner automorphism. 

Here is another definition, taken from \cite{Ye4.5}. 

\begin{dfn} \label{dfn:4875}
Let $A$ be a noetherian ring with rigid dualizing complex $R$. If 
$R \cong A(\mu)[n]$ in $\dcat{D}(A^{\mrm{en}})$ for some 
automorphism $\mu$ and some integer $n$, then $\mu$ is 
called the {\em dualizing automorphism} of $A$. 
\end{dfn}

Like the Nakayama automorphism, the dualizing automorphism is unique up to 
composition with an inner automorphism.  

\begin{prop} \label{prop:4900}  
Let $A$ be a noetherian ring which is homologically smooth over $\K$, let $n$ 
be a natural number, and let $\nu$ be an automorphism of $A$. 
The following two conditions are equivalent. 
\begin{enumerate}
\rmitem{i} $A$ is a twisted CY ring of dimension $n$ with Nakayama automorphism 
$\nu$. 

\rmitem{ii} The complex $R := A(\nu^{-1})[n]$ is a rigid dualizing complex over 
$A$. 
\end{enumerate}
\end{prop}

Thus the dualizing automorphism is the inverse of the Nakayama automorphism.

\begin{proof}
The key formula is this: for $M, N \in \dcat{D}(A^{\mrm{en}})$ there are 
isomorphisms
\begin{equation} \label{eqn:4890}
\begin{aligned}
& 
\opn{RHom}_{A^{\mrm{en}}}(A, M \ot N) \cong^{1}
\opn{RHom}_{A^{\mrm{en}}} \bigl( A, A^{\mrm{en}} 
\ot^{\mrm{L}}_{A^{\mrm{en}}} (M \ot N) \bigr) 
\\ & \quad
\cong^{2} \opn{RHom}_{A^{\mrm{en}}}(A, A^{\mrm{en}}) 
\ot^{\mrm{L}}_{A^{\mrm{en}}} (M \ot N)
\\ & \quad 
\cong^3 N \ot^{\mrm{L}}_{A} \opn{RHom}_{A^{\mrm{en}}}(A, A^{\mrm{en}}) 
\ot^{\mrm{L}}_{A} M 
\end{aligned}
\end{equation}
in $\dcat{D}(A^{\mrm{en}})$.
The isomorphism $\cong^1$ is the left unitor. 
The isomorphism $\cong^2$ is by derived tensor-evaluation, 
and it exists because 
$A$ is a perfect complex over $A^{\mrm{en}}$; see Theorem \ref{thm:3400}.
The isomorphism $\cong^3$ is because the left $A^{\mrm{en}}$ action on 
$M \ot N$ is the outside action, cf.\ formula (\ref{eqn:4215}). 

\medskip \noindent 
(i) $\Rightarrow$ (ii): By Example \ref{exa:3695} the complex $R$ is dualizing. 
What is needed is a rigidifying isomorphism for it. We are given that 
$\opn{RHom}_{A^{\mrm{en}}}(A, A^{\mrm{en}}) \cong A(\nu)[-n]$ 
in $\dcat{D}(A^{\mrm{en}})$. 
Plugging this and $M = N := R$ in (\ref{eqn:4890}) gives 
\[ \begin{aligned}
&
\opn{RHom}_{A^{\mrm{en}}}(A, R \ot R) \cong
A(\nu^{-1})[n] \ot^{\mrm{L}}_{A} A(\nu)[-n] \ot^{\mrm{L}}_{A} A(\nu^{-1})[n]
\\ & \quad 
\cong A(\nu^{-1} \circ \nu \circ \nu^{-1})[n - n + n] = 
A(\nu^{-1})[n] = R , 
\end{aligned} \]
as required. 

\medskip \noindent 
(ii) $\Rightarrow$ (i):
Now we are given that 
$\opn{RHom}_{A^{\mrm{en}}}(A, R \ot R) \cong R$ 
in $\dcat{D}(A^{\mrm{en}})$. Substituting 
$M = N := R$ in (\ref{eqn:4890}) gives 
\[ A(\nu^{-1})[n] \cong A(\nu^{-1})[n] \ot^{\mrm{L}}_{A}
\opn{RHom}_{A^{\mrm{en}}}(A, A^{\mrm{en}}) \ot^{\mrm{L}}_{A} A(\nu^{-1})[n] \]
in $\dcat{D}(A^{\mrm{en}})$. Therefore 
$A(\nu^{})[-n] \cong \opn{RHom}_{A^{\mrm{en}}}(A, A^{\mrm{en}})$ 
in $\dcat{D}(A^{\mrm{en}})$.
This is property (iii) of Definition \ref{dfn:4490}.
\end{proof}

\begin{exa} \label{exa:4895}
This is a continuation of Example \ref{exa:4515}: $\K$ has characteristic $0$, 
$C$ is a smooth commutative $\K$-ring of pure dimension $n$, and 
$A = \mcal{D}(C)$, the ring of differential operators.
It is known that $A$ is homologically smooth -- in fact, in this case 
$A^{\mrm{en}}$ is a regular noetherian ring, of dimension $4 n$. 
According to \cite[Theorem 2.6]{Ye4.5} the rigid dualizing complex of $A$ is
$R = A[2 n]$. Hence $A$ is a CY ring of dimension $2 n$.
\end{exa}

AS regular graded rings were introduced in Definition \ref{dfn:3782}. 

\begin{thm} \label{thm:4490}
\index{Ring! twisted Calabi-Yau}
\index{Algebraically graded ring! AS regular}
Let $A$ be a ring that admits some noetherian connected filtration $F$, such 
that the graded ring $\opn{Gr}^F(A)$ is AS regular of dimension $n$. 
Then $A$ is a noetherian twisted CY ring of dimension $n$, with a Nakayama 
automorphism $\nu$ that respects the filtration $F$. 
\end{thm}

\begin{proof}
Let $\til{A} := \opn{Rees}^F(A)$ and $\bar{A} := \opn{Gr}^F(A)$.
We are given that $\bar{A}$ is AS regular of dimension $n$. According to 
Theorem \ref{thm:4546}(2), $\til{A}$ is AS regular of dimension 
$n + 1$. 

Consider the minimal graded-free resolution 
$\til{Q} \to \til{A}$ of $\til{A}$ over $\til{A}^{\mrm{en}}$.
We know, from Theorem \ref{thm:4548}  and Proposition \ref{prop:4266}, that 
each $\til{Q}^i$ is a 
finite graded-free 
$\til{A}^{\mrm{en}}$-module, and $\til{Q}^i = 0$ for $i < -n -1$. 
Therefore the complex 
\[ P := A^{\mrm{en}} \ot_{\til{A}^{\mrm{en}}} \til{Q} \cong 
A^{\mrm{en}} \ot^{\mrm{L}}_{\til{A}^{\mrm{en}}} \til{A} \in 
\dcat{D}(A^{\mrm{en}}) \]
is perfect. But according to Theorem \ref{thm:4511}(2), $A$ is a direct summand 
of $P$ in $\dcat{D}(A^{\mrm{en}})$. Thus $A$ is a perfect complex over 
$A^{\mrm{en}}$. This is property (ii) of Definition \ref{dfn:4490}.

By Corollary \ref{cor:3780} the graded ring $\til{A}$
has a balanced dualizing complex
$\til{R} := \til{A}(\til{\mu}, l)[n + 1]$,
where $\til{\mu}$ is a graded ring automorphism of $\til{A}$, and $l$ is an 
integer. Theorem \ref{thm:4400} says that $\til{R}$ is a graded rigid dualizing 
complex over $\til{A}$. 

Lemma \ref{lem:4393} says that for every $p$ the $\til{A}$-bimodule 
$\opn{H}^{p}(\til{R})$ is a central bimodule over $\opn{Cent}(\til{A})$.
But for $p = -n - 1$ we have 
$\opn{H}^{-n - 1}(\til{R}) \cong \til{A}(\til{\mu}, l)$,
and this implies that $\til{\mu}$ acts trivially on the element 
$t \in \opn{Cent}(\til{A})$. Hence $\til{\mu}$ extends to an 
automorphism of $\til{A}_t$, and there is an isomorphism
\begin{equation} \label{eqn:4507}
\opn{Ind}_{(\til{A}_t)^{\mrm{en}}} (\til{R}) \cong 
\til{A}_t(\til{\mu}, l)[n + 1]
\end{equation}
in $\dcat{D} \bigl( (\til{A}_t)^{\mrm{en}}, \mrm{gr} \bigr)$. 
Also there is an induced automorphism $\mu$ of the ring 
$A$, and formula (\ref{eqn:4506}) shows that $\mu$ respects the filtration $F$. 

According to Theorem \ref{thm:4395} the noetherian ring $A$ 
has a rigid NC dualizing complex 
\[ R := (\opn{Rest}_{A^{\mrm{en}}} \circ \opn{Deg}_0 \circ 
\opn{Ind}_{(\til{A}_t)^{\mrm{en}}}) (\til{R})[-1] \in \dcat{D}(A^{\mrm{en}}) . 
\]
From (\ref{eqn:4507}) we see that $R \cong A(\mu)[n]$ 
in $\dcat{D}(A^{\mrm{en}})$. Now use Proposition \ref{prop:4900} with 
$\nu := \mu^{-1}$. 
\end{proof}

\begin{exa} \label{exa:4505}
Let $\g$ be a finite Lie algebra over $\K$, with 
$\opn{rank}_{\K}(\g) = n$, and let 
$A := \opn{U}(\g)$ be its universal enveloping algebra. 
The ring $A$ is noetherian twisted CY of dimension $n$. This was worked out in 
the paper \cite{Ye4.5}, and is outlined below. 

The ring $A$ has a standard noetherian connected filtration $F$, such that 
$\bar{A} = \opn{Gr}^F(A) \cong \K[t_1, \ldots, t_n]$,
the commutative polynomial ring in $n$ variables, all of degree $1$. 
The ring $\bar{A}$ is AS regular of dimension $n$, as explained in
Example \ref{exa:3990}.
The Rees ring here is the homogeneous  universal enveloping ring from Example 
\ref{exa:4075}.

Theorem \ref{thm:4490} says that $A$ is twisted CY of dimension $n$, with a 
Nakayama automorphism $\nu$ that respects the filtration $F$. 
In case $\g$ is an abelian Lie algebra (i.e.\ the Lie bracket is zero), then  
$\nu$ is the identity automorphism. When $\g$ is semi-simple, 
Van den Bergh \cite{VdB2} proved that $\nu$ is trivial too. So in these extreme 
cases the ring $A = \opn{U}(\g)$ is CY (untwisted). 

The general case was done in \cite{Ye4.5}. Consider the character (i.e.\ 
rank $1$ representation) $\bwedge^n(\g)$ of $\g$, the $n$-th exterior power of 
the adjoint representation. This can be seen as a Lie algebra 
homomorphism $\ep : \g \to \K$. Then the dualizing automorphism $\mu$ is the 
unique ring automorphism of $A = \opn{U}(\g)$ such that 
$\mu(v) = v - \ep(v) \cd 1_A$
for $v \in \g \sub F_1(A) \sub A$.
The Nakayama automorphism is $\nu = \mu^{-1}$. 

The comultiplication of $A$ makes the tensor product of $A$-bimodules into a 
bimodule. This allows us to express $A(\mu)$ like this:
$A(\mu) \cong A \ot \bwedge^n(\g)$,
where the left action of $A$ on $\bwedge^n(\g)$ is trivial, and the right 
action is the coadjoint action.
\end{exa}

Given an $A$-bimodule $M$, its {\em $i$-th Hochschild homology} is the 
$\K$-module \lb 
$\opn{HH}_i(A, M) := \opn{Tor}_i^{A^{\mrm{en}}}(A, M)$,  
and its {\em $i$-th Hochschild cohomology} is 
$\opn{HH}^i(A, M) := \lb \opn{Ext}^i_{A^{\mrm{en}}}(A, M)$. 

\begin{thm}[Van den Bergh Duality, \cite{VdB2}] \label{thm:4860}
Let $A$ be an $n$-dimensional twisted CY ring, with Nakayama automorphism 
$\nu$. For every $M \in \dcat{M}(A^{\mrm{en}})$ and 
$0 \leq i \leq n$ there is an isomorphism 
\[ \opn{HH}^i(A, M) \cong \opn{HH}_{n - i} \bigl( A, M(\nu) \bigr) \] 
in $\dcat{M}(\K)$, and it is functorial in $M$. 
\end{thm}

The theorem does not require the rings $A$ nor $A^{\mrm{en}}$ to be noetherian.
We need a lemma first. 

\begin{lem} \label{lem:4895}
For a complex $M \in \dcat{D}(A^{\mrm{en}})$, and an automorphism $\mu$ of 
$A$, there is an isomorphism 
$A \ot^{\mrm{L}}_{A^{\mrm{en}}} M(\mu) \cong 
A(\mu) \ot^{\mrm{L}}_{A^{\mrm{en}}} M$
in $\dcat{D}(\K)$. This isomorphism is functorial in $M$. 
\end{lem}

\begin{proof}
It suffices to produce a functorial isomorphism 
$\psi : A \ot_{A^{\mrm{en}}} M(\mu) \iso A(\mu) \ot_{A^{\mrm{en}}} M$
for $M \in \dcat{M}(A^{\mrm{en}})$. 
We define $\psi$ on elements $b \in A$ and $m \in M$ by the formula
$\psi(b \ot m) := \mu(b) \ot m$.
It must be verified that the relations coming from 
the action of the ring $A^{\mrm{en}}$ are respected by $\psi$. 
Take a ring element $a_1 \ot \opn{op}(a_2) \in A^{\mrm{en}}$. Then in the 
module $A \ot_{A^{\mrm{en}}} M(\mu)$ we have equalities
\begin{equation} \label{eqn:4901}
\begin{aligned}
& 
(a_2 \cd b \cd a_1) \ot m = (b \cd (a_1 \ot \opn{op}(a_2))) \ot m = 
\\ & \quad
= b \ot ((a_1 \ot \opn{op}(a_2)) \, \cd^{\mu} \, m) = 
b \ot (a_1 \cd m \cd \mu(a_2)) . 
\end{aligned}
\end{equation}
See formula (\ref{eqn:4904}). In the module $A(\mu) \ot_{A^{\mrm{en}}} M$ 
we have equalities
\begin{equation} \label{eqn:4900}
\begin{aligned} 
& 
\mu(a_2 \cd b \cd a_1) \ot m = 
(\mu(a_2) \cd \mu(b) \cd \mu(a_1)) \ot m = 
\\ & \quad
= (\mu(b) \, \cd^{\mu} \, (a_1 \ot \opn{op}(\mu(a_2)))) \ot m 
\\ & \quad
= \mu(b) \ot ((a_1 \ot \opn{op}(\mu(a_2))) \cd m) 
\\ & \quad
= \mu(b) \ot (a_1 \cd m \cd \mu(a_2))  .
\end{aligned}
\end{equation}
The isomorphism $\psi$ sends the first (resp.\ last) term in (\ref{eqn:4901}) 
to the first  (resp.\ last) term in (\ref{eqn:4900}). 
\end{proof}

\begin{proof}[Proof of Theorem \tup{\ref{thm:4860}}]
We are given an isomorphism 
\begin{equation} \label{eqn:4862}
A(\nu)[-n] \cong 
\opn{RHom}_{A^{\mrm{en}}} (A, A^{\mrm{en}})
\end{equation}
in $\dcat{D}(A^{\mrm{en}})$.
Consider these isomorphisms
\begin{equation} \label{eqn:4863}
\begin{aligned}
& A \ot^{\mrm{L}}_{A^{\mrm{en}}} M(\nu)[-n] \cong^{1} 
A(\nu)[-n] \ot^{\mrm{L}}_{A^{\mrm{en}}} M 
\\ & \quad 
\cong^{2} \opn{RHom}_{A^{\mrm{en}}} (A, A^{\mrm{en}}) 
\ot^{\mrm{L}}_{A^{\mrm{en}}} M
\\ & \quad 
\cong^{3} \opn{RHom}_{A^{\mrm{en}}} (A, A^{\mrm{en}}
\ot^{\mrm{L}}_{A^{\mrm{en}}} M)
\cong^4 \opn{RHom}_{A^{\mrm{en}}} (A, M)
\end{aligned}
\end{equation}
in $\dcat{D}(\K)$.
The first isomorphism is by Lemma \ref{lem:4895}.
The isomorphism $\cong^{2}$ comes from (\ref{eqn:4862}). 
The isomorphism $\cong^{3}$ is by derived tensor-evaluation, and it holds 
because $A$ is perfect over $A^{\mrm{en}}$, see Theorem 
\ref{thm:3400}. The isomorphism $\cong^4$ is clear. 

Taking $\opn{H}^i(-)$ in (\ref{eqn:4863}) we obtain an isomorphism 
\[ \opn{Tor}_{n - i}^{A^{\mrm{en}}} \bigl( A, M(\nu) \bigr) \cong 
\opn{Ext}^i_{A^{\mrm{en}}} (A, M) \]
in $\dcat{M}(\K)$. Functoriality is clear. 
\end{proof}

\begin{rem} \label{rem:4505}
Here is a quick discussion of the evolution of the {\em Calabi-Yau} concept. 
In algebraic geometry, an $n$-dimensional smooth proper scheme $X$ over a field 
$\K$ has the {\em canonical sheaf} $\bsym{\om}_X := \Om^n_{X / \K}$; this is a 
rank $1$ locally free $\OO_X$-module. 
In the paper \cite{BoKa0} from 1990, A.I. Bondal and M.M. Kapranov
introduced the {\em Serre functor}
\begin{equation} \label{eqn:4866}
\opn{S} := \bsym{\om}_X[n] \ot_{\OO_X} (-) . 
\end{equation}
This is an auto-equivalence of the derived category 
$\dcat{D}^{\mrm{b}}_{\mrm{c}}(X) = 
\dcat{D}^{\mrm{b}}_{\mrm{c}}(\cat{Mod} \OO_X)$
which admits a bifunctorial isomorphism 
\begin{equation} \label{eqn:4865}
\opn{Hom}_{\dcat{D}^{\mrm{b}}_{\mrm{c}}(X)}( \MM, \NN) \cong 
\opn{Hom}_{\dcat{D}^{\mrm{b}}_{\mrm{c}}(X)}( \NN, \opn{S}(\MM) )^* .
\end{equation}
Here $(-)^* := \opn{Hom}_{\K}(-, \K)$. 
This is an upgrade of the familiar Serre Duality for $X$, which is the 
special case where $\MM = \OO_X$ and $\NN$ is the translation of a locally free 
$\OO_X$-module. 

The Serre functor makes sense on any $\K$-linear triangulated category 
$\cat{D}$ that has suitable finiteness properties (``smooth'' and ``proper''). 
It is a triangulated auto-equivalence $\opn{S}$ of $\cat{D}$ admitting a 
bifunctorial isomorphism
\begin{equation} \label{eqn:4872}
\opn{Hom}_{\cat{D}}(M, N) \cong \opn{Hom}_{\cat{D}}(N, \opn{S}(M))^* 
\end{equation}
in $\dcat{M}(\K)$ for $M, N \in \cat{D}$. In \cite{BoKa0} it was proved that if 
a Serre functor exists, then it is unique. In many situations (where this makes 
sense) the Serre functor is realized by tensoring with the rigid dualizing 
complex. This is true in the geometric context above, in which the rigid 
dualizing complex of the scheme $X$ is $\bsym{\om}_X[n]$; see formula 
(\ref{eqn:4866}). It is also true for a regular finite $\K$-ring $A$, as 
explained in Example \ref{exa:4870} below. 

Back in algebraic geometry, a smooth proper scheme $X$ is called Calabi-Yau 
if the canonical sheaf $\bsym{\om}_X$ is trivial. (Examples are abelian 
varieties and K3 surfaces.) The CY condition implies that the Serre functor 
$\opn{S}$ is isomorphic to $\opn{T}^n$, the $n$-th power of the 
translation functor of $\cat{D}  = \dcat{D}^{\mrm{b}}_{\mrm{c}}(X)$. 
Thus the dimension of $X$ is precisely the ratio between the 
auto-equivalences $\opn{S}$ and $\opn{T}$. 

Around the year 1994, M. Konstevich began discussing {\em Homological Mirror 
Symmetry}, which is a vast and deep program (still not fully understood) 
relating the algebraic geometry of a CY scheme $X$ (as above) with the {\em 
symplectic geometry} of its mirror partner manifold $Y$.  See \cite{HMS}
for more information on this program. 

Among the ideas arising in the discussion 
by Kontsevich was that of a {\em Calabi-Yau category}. This is a smooth 
proper $\K$-linear triangulated category $\cat{D}$, whose translation functor 
$\opn{T}$ and Serre functor $\opn{S}$ satisfy the numerical relation 
$\opn{S} \cong \opn{T}^n$ for some integer $n$, called the {\em CY dimension} 
of $A$. More generally, if the relation $\opn{S}^m \cong \opn{T}^n$ holds for 
some $m, n \in \Z$, then $\cat{D}$ is called  {\em fractionally CY}
of dimension $(n, m)$. By abuse of notation, such a ring $A$ is sometimes 
called fractionally CY of dimension $n / m$. 

The concept of CY triangulated categories later showed up in various 
mathematical areas, including:
\begin{itemize}
\item Algebraic geometry (see \cite{Kuz}). 

\item Conformal field theories (see \cite{Cost}).

\item Representation theory of algebras and cluster categories
(see  \cite{CCS}, \cite{BMRRT}, \cite{Kel4}, \cite{Kel2}, \cite{Lad} and 
Example \ref{exa:4870} below).
\end{itemize}

In the preprint \cite{Ginz} from 2006, V. Ginzburg defined {\em Calabi-Yau 
algebras}. At about the same time K. Brown and J.J. Zhang \cite{BrZh} 
introduced {\em rigid Gorenstein algebras}. These are both precisely 
the CY rings from Definition \ref{dfn:4490} above.
Eventually the CY terminology was adopted by ring theorists. 

In case $A$ is a noetherian connected graded ring, then it is twisted CY in the 
graded sense if and only if it is AS regular; see \cite{ReRoZh}.
More results on graded CY rings can be found in \cite{ReRo} and its references.
\end{rem}

We finish the subsection (and the book) with an example of a fractionally CY 
category. 

\begin{exa} \label{exa:4870}
Suppose $A$ is a finite $\K$-ring. (Traditional texts would say that 
$A$ is a finite dimensional $\K$-algebra.)  
We claim that the rigid dualizing complex of $A$ is $R := A^{*}$, the 
$\K$-linear dual of $A$. That $R \in \dcat{D}(A^{\mrm{en}})$ is a NC dualizing 
complex over $A$ is very easy to see (cf.\ Definition \ref{dfn:3675}). 
As for rigidity: in this case it is just linear algebra. Indeed, 
$R \ot R = A^* \ot A^* \cong (A^{\mrm{en}})^*$ 
in $\dcat{M}(A^{\mrm{four}})$, and hence there are canonical isomorphisms
\[ \opn{Sq}_A(R) = \opn{RHom}_{A^{\mrm{en}}}(A, R \ot R) \cong 
\opn{Hom}_{A^{\mrm{en}}}(A, (A^{\mrm{en}})^*) \cong A^* = R \]
in $\dcat{D}(A^{\mrm{en}})$.

Now let's assume that $A$ is regular, i.e.\ it has finite global 
cohomological dimension. (If the base field $\K$ is perfect, one can show 
that $A$ is homologically smooth over it; see \cite{Hap2} and \cite{BGMS}.)
Then, according to Theorems \ref{thm:3416} and \ref{thm:3400}, 
every object of $\dcat{D}^{\mrm{b}}_{\mrm{f}}(A)$ is a compact object of 
$\dcat{D}(A)$. This includes $R$. So $R$ satisfies condition (iii) of 
Theorem \ref{thm:3455}, and therefore it is a tilting complex over $A$. 
We see that the triangulated functor 
\[ \opn{S} := R \ot^{\mrm{L}}_{A} (-) : \dcat{D}^{\mrm{b}}_{\mrm{f}}(A)
\to \dcat{D}^{\mrm{b}}_{\mrm{f}}(A) \]
is an auto-equivalence. 

It turns out that $\opn{S}$ is a Serre functor of 
$\dcat{D}^{\mrm{b}}_{\mrm{f}}(A)$. Let us prove this. 
For $M, N \in \dcat{D}^{\mrm{b}}_{\mrm{f}}(A)$
we have canonical isomorphisms
\[ \begin{aligned}
& \opn{RHom}_{A}(N, \opn{S}(M)) = \opn{RHom}_{A}(N, A^* \ot^{\mrm{L}}_{A} M)
\\ & \quad 
\cong^{\dag}  \opn{RHom}_{A}(N, A^*) \ot^{\mrm{L}}_{A} M \cong 
N^* \ot^{\mrm{L}}_{A} M 
\end{aligned} \]
in $\dcat{D}(\K)$. 
The isomorphism $\cong^{\dag}$ is by derived tensor-evaluation, and it holds 
because $N$ is a compact object of $\dcat{D}(A)$. 
Therefore we get functorial isomorphisms 
\begin{equation} \label{eqn:4874}
\begin{aligned}
& \opn{RHom}_{A}(N, \opn{S}(M))^* \cong (N^* \ot^{\mrm{L}}_{A} M)^* 
\cong \opn{RHom}_{\K}( N^* \ot^{\mrm{L}}_{A} M, \K) 
\\ & \quad 
\cong^{\heartsuit} \opn{RHom}_{A}(M, \opn{RHom}_{\K}(N^*, \K)) 
\cong \opn{RHom}_{A}(M, N) 
\end{aligned}
\end{equation}
in $\dcat{D}(\K)$. 
The isomorphism $\cong^{\heartsuit}$ is by adjunction for the ring homomorphism 
$\K \to A$. Taking $\opn{H}^0$ in (\ref{eqn:4874}) we obtain a functorial 
isomorphism 
\[ \opn{Hom}_{\dcat{D}^{\mrm{b}}_{\mrm{f}}(A)}(N, \opn{S}(M))^*
 \cong \opn{Hom}_{\dcat{D}^{\mrm{b}}_{\mrm{f}}(A)}(M, N)  \]
in $\dcat{M}(\K)$. 

Finally we specialize to the setting of Example \ref{exa:4180}.
So $A$ is the ring of upper triangular $n \times n$ matrices over $\K$ for some 
$n \geq 2$. As calculated by A. Yekutieli and E. Kreines in \cite{Ye4}, there 
is an isomorphism 
\[ \underbrace{R \ot^{\mrm{L}}_{A} \cdots \ot^{\mrm{L}}_{A} R}_{n + 1} 
\cong A[n - 1] \]
in $\dcat{D}(A^{\mrm{en}})$. Since the translation functor of
$\dcat{D}^{\mrm{b}}_{\mrm{f}}(A)$ is 
$\opn{T} \cong A[1] \ot^{\mrm{L}}_{A} (-)$, we obtain an isomorphism of 
functors $\opn{S}^{n + 1} \cong \opn{T}^{n - 1}$. 
The conclusion is that the category 
$\dcat{D}^{\mrm{b}}_{\mrm{f}}(A)$ is fractionally CY of dimension 
$(n - 1, n + 1)$.
\end{exa}

%% file: block7_190413.tex
%

\renewcommand{\thisfile}{block7\_190328}
\renewcommand{\bibname}{References}  

\cleardoublepage


%% file: main_190414.bbl
\begin{thebibliography}{155} 
\label{references}


\setlength{\parskip}{0pt}
\setlength{\itemsep}{1pt}   
\setlength{\baselineskip}{11.5pt}

\bibitem{AlJeLi}  L. Alonso, A. Jeremias and J. Lipman,
Local homology and cohomology on schemes, 
{\em Ann.\ Sci.\ ENS} {\bf 30} (1997), 1-39.
Correction   
\url{http://www.math.purdue.edu/~lipman}.

\bibitem{AlKl} A. Altman and S. Kleiman, 
``A Term of Commutative Algebra'', 
\url{http://web.mit.edu/18.705/www}.

\bibitem{AnHaKr} L. Angeleri H\"ugel, D. Happel and H. Krause (eds.),
``Handbook of Tilting Theory'', 
{\em London Math.\ Soc.\ Lecture Note Ser.} {\bf 332} (2006), 147-173.

\bibitem{AriBez} D. Arinkin and R. Bezrukavnikov, 
Perverse coherent sheaves,
{\em Mosc.\ Math.\ J.} {\bf 10} (2010), 3-29.	

\bibitem{SGA-4} 
M. Artin, A. Grothendieck, J.-L. Verdier, eds., ``S\'eminaire de 
G\'eom\'etrie Alg\'ebrique du Bois Marie -- 1963-64 -- Th\'eorie des topos et 
cohomologie \'etale des sch\'emas'', SGA 4, vol.\ 1, Lecture Notes in 
Mathematics {\bf 269}, Springer, 1972. 

\bibitem{ArSchl} M. Artin and W. Schelter,
Graded Algebras of Dimension $3$, {\em Adv.\ Math.} {\bf 66} (1987), 172-216.

\bibitem{ArSmZh} M. Artin, L.W. Small and J.J. Zhang, 
Generic flatness for strongly Noetherian algebras, 
J.\ Algebra {\bf 221} (1999), 579-610.

\bibitem{ATV} M. Artin, J. Tate and M. Van den Bergh,
Some Algebras Associated to Automorphisms of Elliptic Curves,
in ``The Grothendieck Festschrift, Vol.\ I'', Birkh\"auser, 1990, pages 33-85.

\bibitem{ATV2} M. Artin, J. Tate and M. Van den Bergh,
Modules over regular algebras of dimension 3,
Invent.\ Math.\ {\bf 106} (1991), 335-388.

\bibitem{ArZh} M. Artin and J.J. Zhang,
Noncommutative Projective Schemes, 
{\em Advances in Mathematics} {\bf 109} (1994), 228-287. 

\bibitem{AuPlRe} M. Auslander, M. Platzeck and I. Reiten,
Coxeter functors without diagrams,
{\em Trans.\ Amer.\ Math.\ Soc.} {\bf 250} (1979), 1-46.

\bibitem{Avr} L.L. Avramov,
Infinite Free Resolutions, in:
``Six Lectures on Commutative Algebra'', J. Elias et.\ al., eds., 
Birkh\"auser, 1998. 

\bibitem{AIL}  L.L. Avramov, S.B. Iyengar and J. Lipman,
Reflexivity and rigidity for complexes, I: Commutative rings, 
Algebra Number Theory {\bf 4:1} (2010), 47-86. 

\bibitem{AILN}  L.L. Avramov, S.B. Iyengar, J. Lipman and S. Nayak,
Reduction of derived Hochschild functors over commutative algebras and schemes,
Advances in Mathematics {\bf 223} (2010) 735-772.

\bibitem{Bei} A.A. Beilinson, 
Coherent sheaves on $\mbf{P}^n$ and problems in linear algebra, 
{\em Funktsional Anal.\ i Prilozhen.} {\bf 12} (1978), 68-69 (Russian).
English translation in: {\em Functional Anal.\ Appl.} {\bf 12} (1978), 214-216.

\bibitem{BBD} A.A. Beilinson, J. Bernstein and P. Deligne, 
``Faisceaux pervers'', Ast\'erisque {\bf 100}, 1980.

\bibitem{Berger} R. Berger,
A Koszul sign map, arXiv:1708.01430.

\bibitem{Ber} G.M. Bergman,
A note on abelian categories -- translating element-chasing
proofs, and exact embedding in abelian groups, 
\url{http://math.berkeley.edu/~gbergman/papers/unpub/elem-chase.pdf}.

\bibitem{BeGePo} I.N. Bernstein, I.M. Gelfand and V.A. Ponomarev, 
Coxeter functors and Gabriel's theorem, 
{\em Uspekhi Mat.\ Nauk} {\bf 28} (1973) 19-23, 
Trans.: {\em Russian Math.\ Surveys} {\bf 28} (1973), 17-32.

\bibitem{BeLu} J. Bernstein and V. Lunts,
``Equivariant Sheaves and Functors'',
Lecture Notes in Mathematics {\bf 1578}, Springer, 1994. 

\bibitem{SGA6} P. Berthelot, A. Grothendieck and L. Illusie, eds., 
``S\'eminaire de G\'eom\'etrie Alg\'ebrique du Bois Marie - 1966-67 - 
Th\'eorie des intersections et th\'eor\`eme de Riemann-Roch'', SGA 6, 
Lecture Notes in Mathematics {\bf 225}, Springer, 1971.

\bibitem{Bjo} J.E. Bj\"ork, 
The Auslander condition on noetherian rings,
in ``S\'eminaire Dubreil-Malliavin, 1987-1988'', 
LNM {\bf 1404}, Springer, 1989, pp.\ 137-173.

\bibitem{BoNe}  M. Bokstedt and A. Neeman, 
Homotopy limits in triangulated categories, 
Compositio Math.\ {\bf 86} (1993), 209-234.

\bibitem{BoKa0} A.I. Bondal and M.M. Kapranov,
Repesentable functors, Serre functors, and mutations,
Math.\ USSR Izvestia {\bf 35} (1990), 519-541.

\bibitem{BoKa} A.I. Bondal and M.M. Kapranov,
Enhanced Triangulated Categories, {\em Math.\ USSR Sbornik}, {\bf 70} (1991), 
no.\ 1, pages 93-107. 

\bibitem{BonVdB} A. Bondal and M. Van den Bergh,
Generators and Representability of Functors in Commutative and Noncommutative 
Geometry, 
Moscow Mathematical Journal {\bf 3}, Number 1 (2003), 1-36. 

\bibitem{Bor} A. Borel, 
``Algebraic D-modules'',
Academic Press, 1987.  

\bibitem{Bou}  N. Bourbaki, ``Commutative Algebra'', Hermann, Paris, 1972.

\bibitem{BrnBu} S. Brenner and M.C.R. Butler, 
Generalizations of the Bernstein-Gelfand-Ponomarev reflection functors, 
in ``Proceedings ICRA II'', 
Lectures Notes in Math.\ {\bf 832} (1980), 103-169.

\bibitem{BrSh} M.P. Brodmann and R.Y. Sharp,
``Local Cohomology: An Algebraic Introduction with Geometric Applications'',
Cambridge Studies in Advanced Mathematics {\bf 136}, 2nd Edition, 
Cambridge University Press, 2013. 

\bibitem{Brown} E. Brown, 
Cohomology theories, 
{\em Annals of Mathematics} {\bf 75} (1962), 467-484.

\bibitem{BrZh} K.A. Brown and J.J. Zhang, 
Dualising complexes and twisted Hochschild (co)homology for noetherian Hopf 
algebras,
Journal of Algebra {\bf 320}, Issue 5 (2008), 1814-1850. 

\bibitem{BMRRT} A.B. Buan, R. Marsh, M. Reineke, I. Reiten and G. Todorov, 
Tilting theory and cluster combinatorics,
{\em Adv.\ Math.} {\bf 204} (2006), no.\ 2, 572-618.

\bibitem{Buch} R.-O. Buchweitz, 
Maximal Cohen-Macaulay modules and Tate cohomology over Gorenstein rings,
(1986), 155 pp., online 
\url{http:hdl.handle.net/1807/16682}. 

\bibitem{BGMS} R.-O. Buchweitz, E.L. Green, D. Madsen and O. Solberg.
Finite Hochschild Cohomology without Finite Global Dimension,
{\em Mathematical Research Letters} {\bf 12} (2005), 805-816.

\bibitem{CCS} P. Caldero, F. Chapoton and R. Schiffler,  
Quivers with relations and cluster tilted algebras,
{\em Algebr.\ Represent.\ Theory} {\bf 9} (2006), no.\ 4, 359-376.  

\bibitem{COS} A. Canonaco, M. Ornaghi and P. Stellari, 
Localizations of the category of $\mrm{A}_{\infty}$ categories
and internal Homs, 
arXiv:1811.07830.

\bibitem{CaSt} A. Canonaco and P. Stellari, 
A tour about existence and uniqueness of dg enhancements and lifts, 
{\em J. Geom.\ Phys.} {\bf 122} (2017), 28-52.

\bibitem{CWZ} D. Chan, Q.S. Wu and J.J. Zhang,
Pre-balanced dualizing complexes,
{\em Israel Journal of Mathematics} {\bf 132} (2002), 285-314.

\bibitem{Ch} X.W. Chen,
A Note On Standard Equivalences, 
{\em Bulletin LMS} {\bf 48}, (2016), 797-801.

\bibitem{ClPaSc} E. Cline, B. Parshall and L. Scott, 
Derived categories and Morita theory, 
{\em J. Algebra} {\bf 104} (1986), 397-409.

\bibitem{Cost} K.J. Costello,
Topological conformal field theories and Calabi-Yau categories,
Advances in Mathematics {\bf 210} (2007), 165-214.

\bibitem{CrRo} D.A. Craven and R. Rouquier,
Perverse equivalences and Brou\'e's conjecture,
Advances in Mathematics {\bf 248} (2013), 1-58.

\bibitem{DoPu} A. Dold and D. Puppe, 
Homologie nicht-additiver Funktoren, 
Annales de l'Institut Fourier {\bf 11} (1961), 201-312.

\bibitem{DuSh} D. Dugger and B. Shipley,
Topological equivalences of differential graded algebras,
{\em  Advances in Mathematics} {\bf 212} (2007), 37-61.

\bibitem{Eis} D. Eisenbud,
``Commutative Algebra'', Springer, 1994.

\bibitem{Fre} P. Freyd, 
``Abelian Categories'', Harper, 1966. 

\bibitem{Fuk} K. Fukaya, 
Morse homotopy, $\mrm{A}_{\infty}$-category and Floer homologies, MSRI 
preprint No.\ 020-94, 1993.

\bibitem{GaZi} P. Gabriel and M. Zisman, 
``Calculus of Fractions and Homotopy Theory'', Ergebnisse der Mathematik und 
ihrer Grenzgebiete, Springer, 1967.

\bibitem{Gai} D. Gaitsgory,
Recent Progress in Geometric Langlands Theory,  
\url{http://www.math.harvard.edu/~gaitsgde/GL/Bourb.pdf}. 

\bibitem{GeMa} S.I. Gelfand and Y.I. Manin, 
``Methods of Homological Algebra'', Springer, 2002.

\bibitem{Ginz} V. Ginzburg,
Calabi-Yau algebras, eprint arXiv:math/0612139.

\bibitem{Gro} A. Grothendieck, 
Sur quelques points d'alg\`ebre homologique, T\^ohoku Math.\ J.\ {\bf 9}
(1957), 119-221.

\bibitem{LC} A. Grothendieck, 
``Local cohomology (lecture notes by R. Hartshorne)'', 
Lect.\ Notes Math.\ {\bf 41}, Springer, 1967.

\bibitem{EGA} A. Grothendieck and J. Dieudonn\'e, 
``\'El\'ements de g\'eometrie alg\'ebrique'', collective reference for the 
whole series. 

\bibitem{EGA-III} A. Grothendieck and J. Dieudonn\'e, 
``\'El\'ements de g\'eometrie alg\'ebrique'', 
Chapitre III, Premi\`ere partie,
Publ.\ Math.\ IHES {\bf 11}, 1961.

\bibitem{EGA-IV} A. Grothendieck and J. Dieudonn\'e,
``\'El\'ements de g\'eometrie alg\'ebrique'', 
Chapitre IV, Quatri\`eme partie,
Publ.\ Math.\ IHES {\bf 32}, 1967. 

\bibitem{Hap0} D. Happel,
On the derived category of a finite-dimensional algebra,
{\em Commentarii Mathematici Helvetici} {\bf 62} (1987), 339-389.

\bibitem{Hap} D. Happel,
``Triangulated Categories in the Representation of Finite Dimensional 
Algebras'', Cambridge University Press, 1988.

\bibitem{Hap2} D. Happel,
Hochschild cohomology of finite-dimensional algebras, 
in ``S\'eminaire d'Alg\`ebre Paul Dubreil et Marie-Paul Malliavin'',
Lecture Notes in Math.\ {\bf 1404}, Springer, 1989, pp.\ 108-126.

\bibitem{HapRi} D. Happel and C.M. Ringel, 
Tilted algebras, 
{\em Trans.\ Amer.\ Math.\ Soc.} {\bf 274} (1982), 399-443.

\bibitem{RD} R. Hartshorne, ``Residues and duality'', 
Lecture Notes in Mathematics {\bf 20}, Springer, 1966.

\bibitem{Har} R.\ Hartshorne, ``Algebraic Geometry'',
        Springer-Verlag, New-York, 1977.

\bibitem{HiVdB} L. Hille and M. Van den Bergh,
Fourier-Mukai transforms,
in: ``Handbook of Tilting Theory'', 
{\em London Math.\ Soc.\ Lecture Note Ser.} {\bf 332} (2006), pp.\ 147-173.

\bibitem{HiSt} P.J. Hilton and U. Stammbach, 
``A Course in Homological Algebra'',
Springer, 1971.

\bibitem{Hi1} V. Hinich,
Homological algebra of homotopy algebras,
{\em Comm.\ Algebra} {\bf 25} (1997), 3291-3323.
Erratum	arXiv:math/0309453. 

\bibitem{Hi2} V. Hinich,
Lectures on infinity categories, 
arXiv:1709.06271. 

\bibitem{Ho1} M. Hovey, 
``Model categories'', 
Mathematical Surveys and Monographs {\bf 63}, AMS, 1999.

\bibitem{HuSa} B. Huisgen-Zimmermann and M. Saorin,
Geometry of chain complexes and outer automorphism groups under derived 
equivalence, 
{\em Trans.\ Amer.\ Math.\ Soc.} {\bf 353} (2001), 4757-4777.

\bibitem{Huy} D. Huybrechts, 
``Fourier-Mukai Transforms in Algebraic Geometry'',
Oxford Mathematical Monographs, 2006.

\bibitem{Jac} N. Jacobson, ``Basic Algebra I'', Freeman.

\bibitem{Jor} P. J{\o}rgensen, Local Cohomology for Non-Commutative Graded 
Algebras,
{\em Comm. Algebra} {\bf 25} (1997), 575-591. 

\bibitem{Kad} T.V. Kadeisvili, 
On the theory of homology of fiber spaces, 
{\em Uspekhi Mat.\ Nauk} {\bf 35} (1980), no.\ 3(213), 183-188.


\bibitem{Ka} M. Kashiwara, 
``D-Modules and Microlocal Calculus'', 
AMS, 2003. 

\bibitem{KaSc1} M. Kashiwara and P. Schapira, 
``Sheaves on manifolds'', Springer-Verlag, 1990.

\bibitem{KaSc2} M. Kashiwara and P. Schapira, 
``Categories and sheaves'', Springer-Verlag, 2005.

\bibitem{KaSc3} M. Kashiwara and P. Schapira, 
``DQ Modules'', {\em Ast\'erisque} {\bf 345}, 2012. 

\bibitem{Kel0} B. Keller, 
Chain Complexes and Stable Categories,
{\em Manus.\ Math.} {\bf 67} (1990), 379-417.

\bibitem{Kel1} B. Keller, 
Deriving DG categories, 
{\em Ann.\ Sci.\ Ecole Norm.\ Sup.} {\bf 27}, (1994) 63-102.

\bibitem{Kel6} B. Keller, 
Introduction to A-Infinity Algebras and Modules,
{\em Homology, Homotopy and Applications}, Vol.\ {\bf 3} No.\ 1 (2001), 1-35.  
Addendum: 
{\em Homology, Homotopy And Applications} Vol.\ {\bf 4} No.\ 1 (2002), 25-28.

\bibitem{Kel3} B. Keller,
Hochschild cohomology and derived Picard groups,
{\em J. Pure Appl.\ Algebra} {\bf 190} (2004), 177-196.

\bibitem{Kel4} B. Keller,
On triangulated orbit categories, 
{\em Doc.\ Math.} {\bf 10} (2005), 551-581.

\bibitem{Kel5} B. Keller,
A-infinity algebras, modules and functor categories,
in ``Trends in representation theory of algebras and related topics'',
{\em Contemp.\ Math.} {\bf 406}, AMS, 2006, pp.\ 67-93.

\bibitem{Kel2} B. Keller,
Cluster algebras and derived categories, in:
``Derived Categories in Algebraic Geometry'', 
EMS, 2013.

\bibitem{Kelly} G.M. Kelly, 
Chain maps inducing zero homology maps, 
{\em Proc.\ Cambridge Philos.\ Soc.} {\bf 61} (1965), 847-854,

\bibitem{Kon} M. Kontsevich, 
Homological algebra of Mirror Symmetry, 
in ``Proceedings of the International Congress of Mathematicians, Zurich, 
Switzerland 1994'', Vol.\ 1, Birkhauser,  1995, 120-139.

\bibitem{KoSo} M. Kontsevich and Y. Soibelman, 
Homological mirror symmetry and torus fibrations, 
in ``Symplectic geometry and mirror symmetry'', World Sci.\ Publishing, 
2001, pp.\ 203-263.

\bibitem{HMS} M. Kontsevich et.\ al., 
Simons Collaboration on Homological Mirror Symmetry, 
\url{https://schms.math.berkeley.edu}.

\bibitem{Kosz2} J. L. Koszul, 
Sur un type d'alg\'ebres diff\'erentielles en rapport
avec la transgression, in:
``Colloque de topologie (espaces fibr\'es), 
Masson, 1951, 73-81.

\bibitem{Kr2} H. Krause, 
Localization theory for triangulated categories, 
in ``Triangulated categories'', London Math.\ Soc.\ Lecture Note 
Series {\bf 375}, 2010, pp.\ 161-253.

\bibitem{Kuz} A. Kuznetsov,
Calabi-Yau and fractional Calabi-Yau categories,
J. reine angew.\ Math.\ (online 2017). 

\bibitem{Lad} S. Ladkani,
2-CY-tilted algebras that are not Jacobian,
arXiv:1403.6814.

\bibitem{Lam} T.-Y. Lam,
``Lectures on Modules and Rings'', Springer, 1999.

\bibitem{Lef} K. Lef\`evre-Hasegawa,
``Sur les $\mrm{A}_{\infty}$-categories'', Ph.D. Thesis (2003), 
arXiv:0310337v1.

\bibitem{Lev} T. Levasseur, 
Some properties of noncommutative regular rings,
Glasgow Math.\ J. {\bf 34} (1992), 277-300. 

\bibitem{LNS} J. Lipman, S. Nayak and P. Sastry, 
Pseudofunctorial behavior of Cousin complexes on formal schemes, 
in: {\em Contemp.\ Math.} vol.\ {\bf 375}, AMS ,2005, pp.\ 31133.

\bibitem{Li1} J. Lipman,
``Dualizing sheaves, differentials and residues on algebraic varieties'',
Ast\'erisque {\bf 117} (1984).

\bibitem{Li2} J. Lipman,
``Notes on Derived Functors and Grothendieck Duality'', 
in ``Foundations of Grothendieck Duality for Diagrams of Schemes'',
LNM {\bf 1960}, Springer, 2009.

\bibitem{Lur} J. Lurie, Books and papers on derived algebraic geometry, 
\url{http://www.math.harvard.edu/~lurie}.

\bibitem{Lur2} J. Lurie, 
``Higher Algebra'', online at 
\url{http://www.math.harvard.edu/~lurie/papers/HA.pdf}.

\bibitem{Mac1} S. Maclane, ``Homology'', Springer, 1994 (reprint).

\bibitem{Mac2} S. Maclane, 
``Categories for the Working Mathematician'', 
Springer, 1978. 

\bibitem{MacSt} E. Macri and P. Stellari,
Lectures on non-commutative K3 surfaces, Bridgeland stability, and moduli 
spaces, 
arXiv:1807.06169. 

\bibitem{Matl} E. Matlis,
Injective modules over Noetherian rings,
{\em Pacific J. Math.} Volume {\bf 8}, Number 3 (1958), 511-528. 

\bibitem{Mats} H. Matsumura, 
``Commutative Ring Theory'', Cambride University Press, 1986. 

\bibitem{McRo} J.C. McConnell and J.C. Robson, 
``Noncommutative Noetherian Rings'', Wiley, Chichester, 1987.

\bibitem{McSt} J.C. McConnell and J.T. Stafford, 
Gelfand-Kirillov dimension and associated graded modules,
{\em J. Algebra} {\bf 125} (1989), 197-214.

\bibitem{Mil} J.S. Milne,
``\'Etale cohomology'', 
Princeton Univ.\ Press, 1980. 

\bibitem{MiYe} J.-I. Miyachi and A. Yekutieli,
Derived Picard Groups of Finite Dimensional Hereditary Algebras,
{\em Compositio Mathematica} {\bf 129} (2001), 341-368.

\bibitem{NaZa} D. Nadler and E. Zaslow, 
Constructible sheaves and the Fukaya category,
{\em J. Amer.\ Math.\ Soc.} {\bf 22} (2009), 233-286.

\bibitem{NV} C. N\v{a}st\v{a}sescu and F. Van Oystaeyen, 
``Graded and Filtered Rings and Modules'', LNM {\bf 758}, Springer, 1979.

\bibitem{Ne0} A. Neeman, 
The connection between the K-theory localization theorem of Thomason, Trobaugh 
and Yao and the smashing subcategories of Bousfield and Ravenel,
{\em Ann.\ Sci.\ \'Ecole Norm.\ Sup.} {\bf 25} no.\ 5 (1992), 547-566. 

\bibitem{Ne1} A. Neeman, 
``Triangulated categories'', Princeton University Press, 2001.

\bibitem{Ne2} A. Neeman, 
The Grothendieck Duality Theorem via Bousfield's Techniques and Brown 
Representability,
{\em J. AMS} {\bf 9}, Number 1, (1996), 205-236.

\bibitem{nLab} The nLab, an online source for mathematics etc., 
\url{http://ncatlab.org}.

\bibitem{Ols} M. Olsson,
``Algebraic Spaces and Stacks'', 
AMS Colloquium Publications {\bf 62}, 2016.

\bibitem{Orl} D.O. Orlov, 
Triangulated categories of singularities and D-branes in Landau-Ginzburg 
models, 
Tr.\ Mat.\ Inst.\ Steklova {\bf 246} (2004), 240-262.

\bibitem{PSY} M. Porta, L. Shaul and A. Yekutieli,
On the Homology of Completion and Torsion,
{\em Algebr.\ Repesent.\ Theory} {\bf 17} (2014), 31-67.
Erratum: {\em Algebr.\ Represent.\ Theory}, {\bf 18} (2015), 
1401-1405.

\bibitem{PSY2} M. Porta, L. Shaul and A. Yekutieli, 
Completion by Derived Double Centralizer,
{\em Algebr.\ Represent.\ Theory} {\bf 17} (2014), 481-494.

\bibitem{PSY3} M. Porta, L. Shaul and A. Yekutieli, 
Cohomologically Cofinite Complexes,
{\em Comm.\ Algebra} {\bf 43} (2015), 597-615. 

\bibitem{Pos} L. Positselski, 
Dedualizing Complexes of Bicomodules and MGM Duality over Coalgebras,
{\em Algebr.\ Represent.\ Theor.} {\bf 21} (2018), 737-767.

\bibitem{Qu1} D. Quillen, 
Higher algebraic K-theory I,
in:  ``Higher K-theories'', Lecture Notes in Math.\ {\bf 341}, Springer, 
1973, pages 85-147.

\bibitem{ReRoZh} M. Reyes, D. Rogalski and J.J. Zhang,
Skew Calabi-Yau triangulated categories and Frobenius Ext-algebras,
Trans.\ Amer.\ Math.\ Soc.\ {\bf 369} (2017), 309-340. 

\bibitem{ReRo}  M. Reyes and D. Rogalski, 
A Twisted Calabi-Yau Toolkit, 
eprint arXiv:1807.10249. 

\bibitem{Ric0} J. Rickard, 
Derived Categories and Stable Equivalence, 
{\em J. Pure Appl.\ Algebra} {\bf 61} (1989), 303-317.

\bibitem{Ric1} J. Rickard, 
Morita theory for derived categories, 
{\em J. London Math.\ Soc.} {\bf 39} (1989), 436-456.

\bibitem{Ric2} J. Rickard, 
Derived equivalences as derived functors, 
{\em J. London Math.\ Soc.} {\bf 43} (1991), 37-48.

\bibitem{Rog} D. Rogalski,
Idealizer rings and noncommutative projective geometry,
{\em J.\ Algebra} {\bf  279} (2004), 791-809.

\bibitem{Rot} J. Rotman,
``An Introduction to Homological Algebra'',
Academic Press, 1979.

\bibitem{Row} L.R. Rowen,
``Ring Theory'' (Student Edition), Academic Press, 1991.

\bibitem{Rou} R. Rouquier, 
Automorphismes, graduations et cat\'egories triangul\'ees,
{\em Journal of the Institute of Mathematics of Jussieu}, {\bf 10}
(2011) 713-751.

\bibitem{RoZi} R. Rouquier and A. Zimmermann, 
Picard groups for derived module categories, 
{\em Proc.\ London Math.\ Soc.} {\bf 87} (2003), 197-225.

\bibitem{Sao} M. Saorin, 
Dg algebras with enough idempotents, their dg modules and their derived 
categories, 
{\em Algebra and Discrete Mathematics} {\bf 23} (2017), 62-137. 

\bibitem{SKK} M. Sato, T. Kawai and M. Kashiwara,
Microfunctions and pseudo-differential equations,
pages 265-529 in 
``Hyperfunctions and pseudo-differential equations'', 
LNM {\bf 287}, Springer, 1973.

\bibitem{Sch} S. Schwede, 
Algebraic versus topological triangulated categories, 
in ``Proceedings of Conference on Triangulated Categories (Leeds 2006)'',
London Mathematical Society Lecture Note Series {\bf 375}. 

\bibitem{ShSc} S. Schwede and B. Shipley, 
Stable model categories are categories of modules, 
{\em Topology} {\bf 42} (2003), 103-153.

\bibitem{Sha1} L. Shaul, 
Hochschild cohomology commutes with adic completion,
{\em Algebra and Number Theory} {\bf 10:5} (2016), 1001-1029.

\bibitem{Sha2} L. Shaul, 
Reduction of Hochschild cohomology over algebras finite over their center,
{\em Journal of Pure and Applied Algebra} {\bf 219} (2015), 4368-4377.

\bibitem{Sha3} L. Shaul, 
Relations between Derived Hochschild Functors via Twisting,
{\em Comm.\ Algebra} {\bf 44} (2016), 2898-2907.

\bibitem{Shi} B. Shipley, 
Morita theory in stable homotopy theory, 
in: ``Handbook of Tilting Theory'', 
{\em London Math.\ Soc.\ Lecture Note Ser.} {\bf 332} (2006), 393-409.

\bibitem{Seid} P. Seidel, 
``Fukaya categories and Picard-Lefschetz theory'', 
Z\"urich Lectures in Advanced Mathematics, EMS, Z\"urich, 2008.

\bibitem{SoZa} A. Solotar and P. Zadunaisky,
Change of grading, injective dimension and dualizing complexes,
{\em Comm.\ Algebra} {\bf 46} (2018), 4414-4425.

\bibitem{Spa} N. Spaltenstein, 
Resolutions of unbounded complexes, {\em Compositio Math.} {\bf 65} (1988), no.\
2, 121-154.

\bibitem{SP} The Stacks Project, an online reference, J.A. de Jong (Editor),\\  
\url{http://stacks.math.columbia.edu}.

\bibitem{StaVdB} J. T. Stafford and M. Van den Bergh,
Noncommutative curves and noncommutative surfaces,
{\em Bull.\ Amer.\ Math.\ Soc.} {\bf 38} (2001), 171-216.

\bibitem{Sta1} J.D. Stasheff, 
Homotopy associativity of H-spaces, I,
{\em Trans.\ AMS} {\bf 108} (1963), 275-292.

\bibitem{Sta2} J D. Stasheff, 
Homotopy associativity of H-spaces, II, 
{\em Trans.\ AMS} {\bf 108} (1963), 293-312.

\bibitem{Ste} B. Stenstr\"om, ``Rings of Quotients'', Springer, 1975.

\bibitem{Ta} D. Tamarkin,
Microlocal conditions for non-displaceability,
arXiv:0809.1584 (2008). 

\bibitem{Thom} R.W. Thomason and T. Trobaugh,
Higher Algebraic K-Theory of Schemes and of Derived Categories,
in ``The Grothendieck Festschrift'', 
Progress in Mathematics {\bf 88}, Birkh\"auser, 1990.

\bibitem{To1} B. To\"en,
Lectures on DG-categories, in
``Topics in Algebraic and Topological K-Theory'', 
LNM {\bf 2008}, Springer, 2011. 

\bibitem{To2} B. To\"en,
Derived Algebraic Geometry, 
{\em EMS Surv.\ Math.\ Sci.} {\bf 1} (2014), 153-240.

\bibitem{VdB} M. Van den Bergh, 
Existence theorems for dualizing complexes over non-commutative
graded and filtered ring, {\em J. Algebra} {\bf 195} (1997), no. 2, 662-679.

\bibitem{VdB2} M. Van den Bergh, 
A relation between Hochschild homology and cohomology for Gorenstein rings,
{\em Proc.\ Amer.\ Math.\ Soc.} {\bf 126} (1998), 1345-1348.
Erratum: {\em Proc.\ Amer.\ Math.\ Soc.} {\bf 130} (2002), 2809-2810.

\bibitem{Ver0} J.-L. Verdier, 
``Cat\'egories D\'eriv\'ees, \'etat 0'', 
Lecture Notes in Mathematics {\bf 569}, Springer, 1977. 

\bibitem{Ver} J.-L. Verdier, 
Des Cat\'egories D\'eriv\'ees des Cat\'egories Ab\'eliennes, 
Ast\'erisque {\bf 239} (1996).  

\bibitem{VyYe} R. Vyas and A. Yekutieli,
Weak Proregularity, Weak Stability, and the Noncommutative MGM Equivalence,
{\em J.\ Algebra} {\bf 513} (2018), 265-325.

\bibitem{We} C. Weibel, 
``An introduction to homological algebra'',
Cambridge Studies in Advanced Math.\ {\bf 38}, 1994.

\bibitem{WuZh} Q.S. Wu and J.J. Zhang,
Some Homological Invariants of Local PI Algebras,
Journal of Algebra {\bf 225} (2000), 904-935.

\bibitem{WuZh2} Q.S. Wu and J.J. Zhang,
Dualizing Complexes over Noncommutative Local Rings,
Journal of Algebra {\bf 239} (2001), 513-548.

\bibitem{Ye1} A.\ Yekutieli, Dualizing complexes over noncommutative
graded algebras, {\em J.\ Algebra} \textbf{153} (1992), 41-84.

\bibitem{Ye2} A.\ Yekutieli, ``An Explicit Construction of the
Grothendieck Residue Complex'',
Ast\'{e}risque \textbf{208}, 1992.

\bibitem{Ye2.5} A. Yekutieli, 
The residue complex of a noncommutative graded algebra, 
J. Algebra {\bf 186} (1996), 522-543.

\bibitem{Ye3} A.\ Yekutieli,
Residues and Differential Operators on Schemes,
{\em Duke Math.\ J.} {\bf 95} (1998), 305-341. 

\bibitem{Ye4} A. Yekutieli,
Dualizing complexes, Morita equivalence and the derived Picard group of a
ring, with Appendix by E. Kreines, 
{\em J. London Math.\ Soc.} {\bf 60} (1999) 723-746.

\bibitem{Ye4.5} A. Yekutieli,
The Rigid Dualizing Complex of a Universal Enveloping Algebra,
{\em J. Pure Appl.\ Algebra} {\bf 150} (2000), 85-93.

\bibitem{Ye4.6}  A. Yekutieli,
The Derived Picard Group is a Locally Algebraic Group,
{\em Algebras and Representation Theory} {\bf 7} (2004), 53-57.

\bibitem{Ye5} A. Yekutieli,
Rigid Dualizing Complexes via Differential Graded Algebras (Survey),
in ``Proceedings of Conference on Triangulated Categories (Leeds 2006)'',
London Mathematical Society Lecture Note Series {\bf 375}. 

\bibitem{Ye5.5}  A. Yekutieli,
Continuous and Twisted L\_infinity Morphisms,
{\em J. Pure Appl.\ Algebra} {\bf 207} (2006), 575-606.

\bibitem{Ye6} A. Yekutieli,
``A Course on Derived Categories'',   
arXiv:1206.6632v2, (2012).

\bibitem{Ye7} A.\ Yekutieli, 
Central Extensions of Gerbes,
{\em Advances Math.} {\bf 225} (2010), 445-486.

\bibitem{Ye8} A.\ Yekutieli, 
Duality and Tilting for Commutative DG Rings,
arXiv:1312.6411v4 (2016).

\bibitem{Ye9} A. Yekutieli,
The Squaring Operation for Commutative DG Rings, 
{\em J. Algebra} {\bf 449} (2016), 50-107.

\bibitem{Ye10} A.\ Yekutieli,
Residues and Duality for Schemes and Stacks,
lecture notes (2013),
\url{http://www.math.bgu.ac.il/~amyekut/lectures/resid-stacks/handout.pdf}.

\bibitem{Ye11} A.\ Yekutieli,   
Rigidity, Residues and Duality for Commutative Rings, in Preperation. 

\bibitem{Ye12} A.\ Yekutieli,
Rigidity, Residues and Duality for Schemes, in Preperation. 

\bibitem{Ye13} A.\ Yekutieli,
Rigidity, Residues and Duality for DM Stacks, in Preperation.

\bibitem{Ye14} A.\ Yekutieli,
Derived Categories of Bimodules,  in Preperation. 

\bibitem{Ye15} A.\ Yekutieli,
The Derived Category of Sheaves of Commutative DG Rings, 
lecture notes (2017),
\url{http://www.math.bgu.ac.il/~amyekut/lectures/shvs-dgrings/notes.pdf}.

\bibitem{Ye16} A.\ Yekutieli,
Derived Categories of Bimodules, lecture notes (2016),
\url{http://www.math.bgu.ac.il/~amyekut/lectures%
/der-cat-bimodules/abstract.html}.

\bibitem{Ye17} A.\ Yekutieli,
Another Proof of a Theorem of Van den Bergh about Graded-Injective 
Modules, arXiv:1407.5916. 

\bibitem{YeZh0} A. Yekutieli and J.J. Zhang,
Serre duality for noncommutative projective schemes,
{\em Proc.\ Amer.\ Math.\ Soc.} {\bf 125} (1997), 697-707.

\bibitem{YeZh1} A. Yekutieli and J.J. Zhang,
Rings with Auslander Dualizing Complexes,
{\em J. Algebra} {\bf 213} (1999), 1-51. 

\bibitem{YeZh1.5}  A. Yekutieli and J. Zhang,
Multiplicities of Indecomposable Injectives,
J. London Math.\ Soc.\ (2) {\bf 71} (2005), 100-120.

\bibitem{YeZh1.7}  A. Yekutieli and J. Zhang,
Homological Transcendence Degree,
Proc.\ London Math.\ Soc.\ (3) {\bf 93} (2006) 105-137.

\bibitem{YeZh2.5} A. Yekutieli and J.J. Zhang,
Dualizing Complexes and Perverse Modules over Differential Algebras,
{\em Compositio Mathematica} {\bf 141} (2005), 620-654. 

\bibitem{YeZh2} A. Yekutieli and J. Zhang,
Dualizing Complexes and Perverse Sheaves on Noncommutative Ringed Schemes,
{\em Selecta Math.} {\bf 12} (2006), 137-177.

\bibitem{YeZh3} A. Yekutieli and J.J. Zhang,
Rigid Complexes via DG Algebras,
{\em Trans.\ AMS} {\bf 360} no.\ 6 (2008), 3211-3248.

\bibitem{YeZh4} A. Yekutieli and J.J. Zhang,
Rigid Dualizing Complexes over Commutative Rings,
{\em Algebras and Representation Theory} {\bf 12}, Number 1 (2009), 19-52.

\bibitem{YeZh5} A. Yekutieli and J.J. Zhang,
Residue Complexes over Noncommutative Rings,
{\em J. Algebra} {\bf 259} (2003), 451-493. 

\end{thebibliography}
